\documentclass{amsbook}
\usepackage{amssymb}
\usepackage{mathrsfs}
\usepackage{stmaryrd} 
\usepackage{bbm}
\usepackage{oldgerm}
\usepackage[frenchb]{babel}
\usepackage[T1]{fontenc}
\usepackage[latin1]{inputenc}
\usepackage[all]{xy}
\usepackage{hyperref}
\usepackage{upref}
\usepackage{xcolor}
\usepackage{microtype}
\usepackage[Symbol]{upgreek}

\hypersetup{
    colorlinks,
    linkcolor={red!50!black},
    citecolor={blue!50!black},
    urlcolor={blue!80!black}
}

\newtheorem{teo}[subsection]{Théorème}
\newtheorem{prop}[subsection]{Proposition}
\newtheorem{cor}[subsection]{Corollaire}
\newtheorem{lem}[subsection]{Lemme}

\theoremstyle{definition}

\newtheorem{defi}[subsection]{Définition}
\newtheorem{rema}[subsection]{Remarque}
\newtheorem{remas}[subsection]{Remarques}

\numberwithin{equation}{subsection}

\newcommand{\gtimes}{\stackrel{\leftarrow}{\times}}

\mathchardef\mhyphen="2D

\DeclareMathSymbol{\mlq}{\mathord}{operators}{``}
\DeclareMathSymbol{\mrq}{\mathord}{operators}{`'}

\newcommand{\mA}{{\mathbb A}}

\newcommand{\mD}{{\mathbb D}}
\newcommand{\mH}{{\mathbb H}}
\newcommand{\mI}{{\mathbb I}}
\newcommand{\mJ}{{\mathbb J}}

\newcommand{\mQ}{{\mathbb Q}}

\newcommand{\mN}{{\mathbb N}}

\newcommand{\mX}{{\mathbb X}}
\newcommand{\mY}{{\mathbb Y}}
\newcommand{\mZ}{{\mathbb Z}}

\newcommand{\mG}{{\mathbb G}}
\newcommand{\mK}{{\mathbb K}}
\newcommand{\mU}{{\mathbb U}}
\newcommand{\mV}{{\mathbb V}}

\newcommand{\mg}{{\mathbbm g}}

\newcommand{\bA}{{\bf A}}
\newcommand{\bB}{{\bf B}}

\newcommand{\bD}{{\bf D}}
\newcommand{\bE}{{\bf E}}

\newcommand{\bIndMod}{{\bf Ind\mhyphen Mod}}
\newcommand{\bK}{{\bf K}}
\newcommand{\bL}{{\bf L}}

\newcommand{\bN}{{\bf N}}
\newcommand{\bT}{{\bf T}}
\newcommand{\bP}{{\bf P}}
\newcommand{\bQ}{{\bf Q}}

\newcommand{\Et}{{\bf \acute{E}t}}
\newcommand{\Sch}{{\bf Sch}}
\newcommand{\Ens}{{\bf Ens}}
\newcommand{\Pt}{{\bf Pt}}

\newcommand{\bHom}{{\bf Hom}}

\newcommand{\bRep}{{\bf Rep}}

\newcommand{\bMod}{{\bf Mod}}
\newcommand{\Ind}{{\bf Ind}}

\newcommand{\bInd}{{\bf Ind}}

\newcommand{\add}{{\rm add}}

\newcommand{\colim}{{\underset{\longrightarrow}{\lim}}}
\newcommand{\indcolim}{{\mlq\mlq\colim \mrq\mrq}}
\newcommand{\indoplus}{{\mlq\mlq\bigoplus \mrq\mrq}}

\newcommand{\intern}{{\diamond}}
\newcommand{\lgg}{{\ttg}}

\newcommand{\et}{{\rm \acute{e}t}}
\newcommand{\fet}{{\rm f\acute{e}t}}

\newcommand{\ad}{{\rm ad}}

\newcommand{\MOD}{{\bf \underline{\bMod}}}
\newcommand{\INDMOD}{{\bf \underline{\bIndMod}}}
\newcommand{\INDMC}{{\bf \underline{\bIndMC}}}

\newcommand{\zar}{{\rm zar}}

\newcommand{\coh}{{\rm coh}}
\newcommand{\qcoh}{{\rm qcoh}}
\newcommand{\Tot}{{\rm Tot}}
\newcommand{\cart}{{\rm cart}}

\newcommand{\HT}{{\rm HT}}
\newcommand{\Dolb}{{\rm Dolb}}
\newcommand{\sol}{{\rm sol}}
\newcommand{\fDolb}{{\rm fDolb}}
\newcommand{\fsol}{{\rm fsol}}

\newcommand{\qsol}{{\rm qsol}}

\newcommand{\apt}{{\rm apt}}
\newcommand{\qsolnilp}{{\rm qsolnilp}}

\newcommand{\aatf}{{\text {\rm -atf}}}

\newcommand{\atf}{{\rm atf}}

\newcommand{\rf}{{\rm f}}

\newcommand{\Spec}{{\rm Spec}}

\newcommand{\Spf}{{\rm Spf}}

\newcommand{\ob}{{\rm Ob}}

\newcommand{\pr}{{\rm pr}}

\newcommand{\tor}{{\rm tor}}

\newcommand{\triv}{{\rm triv}}
\newcommand{\disc}{{\rm disc}}
\newcommand{\coker}{{\rm coker}}

\newcommand{\adj}{{\rm adj}}

\newcommand{\tot}{{\rm tot}}
\newcommand{\qpp}{{\rm qpp}}
\newcommand{\p}{{\rm p}}
\newcommand{\locp}{{\rm locp}}
\newcommand{\trqpp}{{\rm \text -qpp}}

\newcommand{\cont}{{\rm cont}}

\newcommand{\Gr}{{\rm Gr}}
\newcommand{\gp}{{\rm gp}}
\newcommand{\id}{{\rm id}}
\newcommand{\Tr}{{\rm Tr}}

\newcommand{\rb}{{\rm b}}
\newcommand{\Sym}{{\rm Sym}}

\newcommand{\Hom}{{\rm Hom}}

\newcommand{\End}{{\rm End}}

\newcommand{\Mat}{{\rm Mat}}
\newcommand{\Aut}{{\rm Aut}}

\newcommand{\Gal}{{\rm Gal}}
\newcommand{\bMH}{{\bf MH}}
\newcommand{\bIMC}{{\bf IC}}
\newcommand{\bIndMC}{{\bf Ind\mhyphen MC}}
\newcommand{\bIndMH}{{\bf Ind\mhyphen MH}}
\newcommand{\bIH}{{\bf IH}}
\newcommand{\bMC}{{\bf MC}}

\newcommand{\rC}{{\rm C}}

\newcommand{\rE}{{\rm E}}
\newcommand{\rF}{{\rm F}}
\newcommand{\rH}{{\rm H}}
\newcommand{\rI}{{\rm I}}
\newcommand{\rT}{{\rm T}}

\newcommand{\rR}{{\rm R}}
\newcommand{\rS}{{\rm S}}
\newcommand{\rV}{{\rm V}}
\newcommand{\rW}{{\rm W}}
\newcommand{\rN}{{\rm N}}

\newcommand{\rv}{{\rm v}}

\newcommand{\oF}{{\overline{F}}}
\newcommand{\oK}{{\overline{K}}}

\newcommand{\oR}{{\overline{R}}}
\newcommand{\oS}{{\overline{S}}}

\newcommand{\oU}{{\overline{U}}}

\newcommand{\oX}{{\overline{X}}}
\newcommand{\oY}{{\overline{Y}}}
\newcommand{\oZ}{{\overline{Z}}}
\newcommand{\oa}{{\overline{a}}}

\newcommand{\ogg}{{\overline{g}}}

\newcommand{\om}{{\overline{m}}}
\newcommand{\op}{{\overline{p}}}
\newcommand{\os}{{\overline{s}}}

\newcommand{\ou}{{\overline{u}}}

\newcommand{\ox}{{\overline{x}}}
\newcommand{\oy}{{\overline{y}}}
\newcommand{\oz}{{\overline{z}}}

\newcommand{\oeta}{{\overline{\eta}}}

\newcommand{\otheta}{{\overline{\theta}}}

\newcommand{\oiota}{{\overline{\iota}}}
\newcommand{\onabla}{{\overline{\nabla}}}
\newcommand{\oupmu}{{\overline{\upmu}}}

\newcommand{\ouplambda}{{\overline{\uplambda}}}

\newcommand{\ocR}{{\overline{\cR}}}

\newcommand{\ocB}{{\overline{\cB}}}

\newcommand{\ofp}{{\overline{\fp}}}

\newcommand{\uE}{{\underline{E}}}
\newcommand{\uG}{{\underline{G}}}
\newcommand{\uM}{{\underline{M}}}

\newcommand{\uX}{{\underline{X}}}

\newcommand{\ud}{{\underline{d}}}

\newcommand{\ug}{{\underline{g}}}
\newcommand{\uh}{{\underline{h}}}
\newcommand{\um}{{\underline{m}}}
\newcommand{\un}{{\underline{n}}}
\newcommand{\upp}{{\underline{p}}}
\newcommand{\uy}{{\underline{y}}}

\newcommand{\uoX}{{\underline{\oX}}}

\newcommand{\ucH}{{\underline{\cH}}}

\newcommand{\ufS}{{\underline{\fS}}}

\newcommand{\utheta}{{\underline{\theta}}}
\newcommand{\ualpha}{{\underline{\alpha}}}
\newcommand{\ubeta}{{\underline{\beta}}}
\newcommand{\uvarepsilon}{{\underline{\varepsilon}}}

\newcommand{\udelta}{{\underline{\delta}}}

\newcommand{\unabla}{{\underline{\nabla}}}
\newcommand{\upi}{{\underline{\pi}}}
\newcommand{\uPi}{{\underline{\Pi}}}

\newcommand{\ulambda}{{\underline{\lambda}}}
\newcommand{\uupgamma}{{\underline{\upgamma}}}

\newcommand{\uTheta}{{\underline{\Theta}}}

\newcommand{\uhupsigma}{{\underline{\hupsigma}}}
\newcommand{\uDelta}{{\underline{\Delta}}}

\newcommand{\umK}{{\underline{\mK}}}
\newcommand{\utmK}{{\underline{\tmK}}}
\newcommand{\uvartheta}{{\underline{\vartheta}}}

\newcommand{\ha}{{\widehat{a}}}

\newcommand{\hf}{{\widehat{f}}}
\newcommand{\hg}{{\widehat{g}}}
\newcommand{\hh}{{\widehat{h}}}
\newcommand{\hm}{{\widehat{m}}}
\newcommand{\hn}{{\widehat{n}}}

\newcommand{\hA}{{\widehat{A}}}
\newcommand{\hB}{{\widehat{B}}}

\newcommand{\hE}{{\widehat{E}}}

\newcommand{\hcA}{{\widehat{\cA}}}

\newcommand{\hRun}{{\widehat{R_1}}}
\newcommand{\hRunp}{{\widehat{R'_1}}}

\newcommand{\hRi}{{\widehat{R_\infty}}}

\newcommand{\hRinterni}{{\widehat{R^{\intern}_\infty}}}
\newcommand{\hRpi}{{\widehat{R_{p^\infty}}}}
\newcommand{\hRinternpi}{{\widehat{R^{\intern}_{p^\infty}}}}

\newcommand{\hRip}{{\widehat{R'_\infty}}}

\newcommand{\hRpip}{{\widehat{R'_{p^\infty}}}}

\newcommand{\hmZ}{{\widehat{\mZ}}}

\newcommand{\hocR}{\widehat{\ocR}}

\newcommand{\halpha}{\widehat{\alpha}}

\newcommand{\htheta}{\widehat{\theta}}

\newcommand{\hPhi}{\widehat{\Phi}}

\newcommand{\hPsi}{\widehat{\Psi}}
\newcommand{\hnabla}{\widehat{\nabla}}

\newcommand{\hupsigma}{\widehat{\upsigma}}
\newcommand{\hupnu}{\widehat{\upnu}}
\newcommand{\hupbeta}{\widehat{\upbeta}}
\newcommand{\huppi}{\widehat{\uppi}}
\newcommand{\hmg}{\widehat{\mg}}

\newcommand{\vupsigma}{\vec{\upsigma}}
\newcommand{\vuppi}{\vec{\uppi}}

\newcommand{\cA}{{\mathscr A}}
\newcommand{\cB}{{\mathscr B}}
\newcommand{\cC}{{\mathscr C}}

\newcommand{\cE}{{\mathscr E}}
\newcommand{\cF}{{\mathscr F}}
\newcommand{\cG}{{\mathscr G}}
\newcommand{\cI}{{\mathscr I}}
\newcommand{\cJ}{{\mathscr J}}
\newcommand{\cK}{{\mathscr K}}
\newcommand{\cL}{{\mathscr L}}
\newcommand{\cP}{{\mathscr P}}
\newcommand{\co}{{\mathscr O}}
\newcommand{\cR}{{\mathscr R}}
\newcommand{\cS}{{\mathscr S}}
\newcommand{\cT}{{\mathscr T}}
\newcommand{\cH}{{\mathscr H}}
\newcommand{\cM}{{\mathscr M}}
\newcommand{\cN}{{\mathscr N}}
\newcommand{\cQ}{{\mathscr Q}}
\newcommand{\cU}{{\mathscr U}}
\newcommand{\cV}{{\mathscr V}}

\newcommand{\cZ}{{\mathscr Z}}

\newcommand{\cHom}{{\mathscr Hom}}

\newcommand{\cEnd}{{\mathscr End}}

\newcommand{\fC}{{\mathfrak C}}

\newcommand{\fF}{{\mathfrak F}}
\newcommand{\fG}{{\mathfrak G}}

\newcommand{\fN}{{\mathfrak N}}

\newcommand{\fS}{{\mathfrak S}}

\newcommand{\fV}{{\mathfrak V}}

\newcommand{\fX}{{\mathfrak X}}

\newcommand{\fd}{{\mathfrak d}}
\newcommand{\fgg}{{\mathfrak g}}
\newcommand{\fm}{{\mathfrak m}}

\newcommand{\fp}{{\mathfrak p}}

\newcommand{\tta}{{\tt a}}
\newcommand{\ttb}{{\tt b}}

\newcommand{\ttg}{{\tt g}}
\newcommand{\tth}{{\tt h}}
\newcommand{\ttt}{{\tt t}}
\newcommand{\tts}{{\tt s}}

\newcommand{\ttu}{{\tt u}}
\newcommand{\ttv}{{\tt v}}

\newcommand{\hM}{{\widehat{M}}}

\newcommand{\hR}{{\widehat{R}}}
\newcommand{\hX}{{\widehat{X}}}

\newcommand{\hoR}{{\widehat{\oR}}}
\newcommand{\hoRp}{{\widehat{\oR'}}}
\newcommand{\hfS}{{\widehat{\fS}}}

\newcommand{\hcS}{{\widehat{\cS}}}

\newcommand{\hcC}{{\widehat{\cC}}}

\newcommand{\hcG}{{\widehat{\cG}}}

\newcommand{\hcR}{{\widehat{\cR}}}
\newcommand{\hmX}{{\widehat{\mX}}}
\newcommand{\hmY}{{\widehat{\mY}}}

\newcommand{\hotimes}{{\widehat{\otimes}}}

\newcommand{\htta}{{\widehat{\tta}}}

\newcommand{\bvg}{{\breve{g}}}
\newcommand{\bvu}{{\breve{u}}}

\newcommand{\bvoS}{{\breve{\oS}}}

\newcommand{\bvoX}{{\breve{\oX}}}
\newcommand{\bvocB}{{\breve{\ocB}}}

\newcommand{\bvcC}{{{\breve{\cC}}}}

\newcommand{\bvcF}{{\breve{\cF}}}
\newcommand{\bvcH}{{\breve{\cH}}}

\newcommand{\bvcP}{{\breve{\cP}}}
\newcommand{\bvtau}{{\breve{\tau}}}

\newcommand{\bvpi}{{\breve{\pi}}}

\newcommand{\bvalpha}{{\breve{\alpha}}}
\newcommand{\bvbeta}{{\breve{\beta}}}
\newcommand{\bvsigma}{{\breve{\sigma}}}
\newcommand{\bvSigma}{{\breve{\Sigma}}}

\newcommand{\bvupiota}{{\breve{\upiota}}}
\newcommand{\bvPhi}{{\breve{\Phi}}}
\newcommand{\bvpsi}{{\breve{\psi}}}
\newcommand{\bvTheta}{{\breve{\Theta}}}

\newcommand{\bvupomega}{{\breve{\upomega}}}
\newcommand{\bvupnu}{{\breve{\upnu}}}
\newcommand{\bvupgamma}{{\breve{\upgamma}}}
\newcommand{\bvuptheta}{{\breve{\uptheta}}}
\newcommand{\bvupmu}{{\breve{\upmu}}}
\newcommand{\bvmZ}{{\breve{\mZ}}}
\newcommand{\bvtta}{{\breve{\tta}}}
\newcommand{\bvogg}{{\breve{\ogg}}}
\newcommand{\bvlgg}{{\breve{\lgg}}}

\newcommand{\coS}{{\check{\oS}}}
\newcommand{\coU}{{\check{\oU}}}
\newcommand{\coX}{{\check{\oX}}}
\newcommand{\coY}{{\check{\oY}}}
\newcommand{\cog}{{\check{\ogg}}}
\newcommand{\cou}{{\check{\ou}}}
\newcommand{\calpha}{{\check{\alpha}}}
\newcommand{\cbeta}{{\check{\beta}}}

\newcommand{\coupmu}{{\check{\oupmu}}}

\newcommand{\tE}{{\widetilde{E}}}

\newcommand{\tG}{{\widetilde{G}}}
\newcommand{\tuE}{{\widetilde{\uE}}}
\newcommand{\tuG}{{\widetilde{\uG}}}

\newcommand{\tS}{{\widetilde{S}}}

\newcommand{\tU}{{\widetilde{U}}}

\newcommand{\tX}{{\widetilde{X}}}
\newcommand{\tY}{{\widetilde{Y}}}

\newcommand{\tg}{{\widetilde{g}}}
\newcommand{\tlt}{{\widetilde{t}}}

\newcommand{\tu}{{\widetilde{u}}}

\newcommand{\tx}{{\widetilde{x}}}
\newcommand{\ty}{{\widetilde{y}}}

\newcommand{\trT}{{\widetilde{\rT}}}

\newcommand{\tOmega}{{\widetilde{\Omega}}}

\newcommand{\tbeta}{{\widetilde{\beta}}}

\newcommand{\tnu}{{\widetilde{\nu}}}
\newcommand{\tupmu}{{\widetilde{\upmu}}}

\newcommand{\tpi}{{\widetilde{\pi}}}
\newcommand{\tvarphi}{{\widetilde{\varphi}}}
\newcommand{\tpsi}{{\widetilde{\psi}}}

\newcommand{\tuptau}{{\widetilde{\uptau}}}

\newcommand{\txi}{{\widetilde{\xi}}}

\newcommand{\talpha}{{\widetilde{\alpha}}}

\newcommand{\tcC}{{\widetilde{\cC}}}

\newcommand{\tmK}{{\widetilde{\mK}}}
\newcommand{\tmX}{{\widetilde{\mX}}}
\newcommand{\tmY}{{\widetilde{\mY}}}

\newcommand{\tmg}{{\widetilde{\mg}}}

\def \trace {\text{\rm tr}}
\def \oCA {\overline{\mathcal{A}}}
\def \Gal {\text{\rm {Gal}}}
\def \Aut {\text{\rm {Aut}}}
\def \os {\overline{s}}
\def \CC {\mathcal{C}}
\def \CU {\mathcal{U}}
\def \CO {\mathcal{O}}

\def \triv {\text{\rm triv}}
\def \gpt {\text{\rm gpt}}
\def \ur {\text{\rm ur}}
\def \tr {\text{\rm tr}}

\begin{document}

\title{Correspondance de Simpson $p$-adique II : fonctorialité par image directe propre et systèmes locaux de Hodge-Tate}

\author{Ahmed Abbes et Michel Gros}
\address{A.A. Laboratoire Alexander Grothendieck, UMR 9009 du CNRS, 
Institut des Hautes \'Etudes Scientifiques, 35 route de Chartres, 91440 Bures-sur-Yvette, France}
\address{M.G. Université de Rennes, CNRS, IRMAR - UMR 6625, Campus de Beaulieu, F-35042 Rennes cedex, France}
\email{abbes@ihes.fr}
\email{michel.gros@univ-rennes1.fr}

\begin{abstract}
Faltings a dégagé en 2005 un analogue $p$-adique de la correspondance de Simpson (complexe) dont la construction a été reprise par différents auteurs, selon plusieurs approches. 
Poursuivant celle que nous avons initiée dans \cite{agt}, nous développons dans la présente monographie 
de nouveaux aspects de la correspondance de Simpson $p$-adique, inspirés par notre construction
de la suite spectrale de Hodge-Tate relative \cite{ag}. Nous traitons tout d'abord du lien avec les systèmes locaux de Hodge-Tate. 
Nous établissons ensuite la fonctorialité de la correspondance de Simpson $p$-adique par image directe propre. Chemin faisant,
nous élargissons la portée de notre construction initiale.
\end{abstract}

\maketitle

\setcounter{tocdepth}{1}

\newpage{\setlength{\parindent}{0pc}}
\thispagestyle{empty}

\vspace*{13.5pc}
\begin{center}
{\Large\em À la mémoire de Pierre Berthelot}
\end{center}

\newpage\null\thispagestyle{empty}\newpage

\tableofcontents

\chapter{Un survol}

\section{Introduction}\label{intro}

\subsection{}\label{intro1}
Initiée par Faltings \cite{faltings3} et développée suivant diverses approches dont celles de Tsuji et des auteurs \cite{agt}, la
correspondance de Simpson $p$-adique fournit une équivalence de catégories entre certains {\em systèmes locaux étales $p$-adiques} 
sur une variété algébrique définie sur un corps $p$-adique et certains {\em fibrés de Higgs}. L'idée clef de sa construction vient de la
stratégie utilisée par Faltings dans son calcul de la cohomologie d'un système local $p$-adique. 

\subsection{}\label{intro3}
Soient $K$ un corps de valuation discrète complet de caractéristique $0$, à corps résiduel $k$ {\em algébriquement clos} de caractéristique $p>0$, 
$\co_K$ l'anneau de valuation de $K$, $\oK$ une clôture algébrique de $K$, $\co_\oK$ la clôture intégrale de $\co_K$ dans $\oK$, 
$G_K$ le groupe de Galois de $\oK$ sur $K$, $\co_C$ le séparé complété $p$-adique de $\co_\oK$,  $\fm_C$ son idéal maximal, 
$C$ son corps de fractions. Nous poserons $S=\Spec(\co_K)$,  $\oS=\Spec(\co_\oK)$ et nous noterons $s$ (resp.  $\eta$, resp. $\oeta$) 
le point fermé de $S$ (resp. point générique de $S$, resp. point générique de $\oS$).

Soient $X$ un $S$-schéma propre et lisse et  $L$ un faisceau localement constant constructible de
$\mZ/p^n\mZ$-modules de $X_{\oeta,\et}$ pour un entier $n\geq 0$. 
Pour calculer la cohomologie de $L$, Faltings introduit un topos annelé $(\tE,\ocB)$ équipé de deux morphismes de topos 
\begin{equation}
\xymatrix{
X_{\oeta,\et}\ar[r]^-(0.5){\psi}&{\tE}\ar[r]^-(0.5){\sigma}&X_\et}.
\end{equation}
Pour tout entier $j\geq 1$, on a $\rR^j\psi_*(L)=0$. En particulier, pour tout $i\geq 0$, on a un isomorphisme canonique
\begin{equation}\label{intro3a}
\rH^i(X_{\oeta,\et},L)\stackrel{\sim}{\rightarrow}\rH^i(\tE,\psi_*(L)).
\end{equation} 
Utilisant la théorie d'Artin-Schreier, Faltings prouve un raffinement de ce résultat, à savoir que le morphisme canonique 
\begin{equation}\label{intro3b} 
\rH^i(X_{\oeta,\et},L)\otimes_{\mZ_p}\co_C\rightarrow \rH^i(\tE,\psi_*(L)\otimes_{\mZ_p}\ocB)
\end{equation}
est un {\em presque-isomorphisme}, {\em i.e.}, son noyau et son conoyau sont annulés par $\fm_C$.  
Posant $\cM=\psi_*(L)\otimes_{\mZ_p}\ocB$, on a un isomorphisme canonique
\begin{equation}\label{intro3c}
\rR\Gamma(\tE,\cM) \stackrel{\sim}{\rightarrow}\rR\Gamma(X_\et,\rR\sigma_*(\cM)). 
\end{equation}
Ce calcul s'étend aux $\mQ_p$-systèmes locaux par passage à la limite sur $n$ puis inversion de $p$. 
Pour certains $\mQ_p$-systèmes locaux $L$, le complexe correspondant $\rR\sigma_*(\cM)$ 
est alors le complexe de Dolbeault d'un fibré de Higgs canoniquement associé à $L$ par la {\em correspondance de Simpson $p$-adique} 
(voir \ref{indsh30} pour la terminologie sur les modules de Higgs). La suite spectrale de Cartan-Leray pour $\sigma$ \eqref{intro3c} 
et le presque-isomorphisme \eqref{intro3b} conduisent, par suite, à une généralisation de la suite spectrale de Hodge-Tate \cite{ag}. 

\subsection{}
Pour décrire notre construction de la correspondance de Simpson $p$-adique, 
nous avons besoin de l'anneau $\bvocB=(\ocB/p^n\ocB)_{n\geq 0}$ du topos $\tE^{\mN^\circ}$ des sytèmes projectifs d'objets de $\tE$, indexés par 
l'ensemble ordonné $\mN$. La catégorie de coefficients avec laquelle nous travaillerons est la catégorie des $\bvocB$-modules à isogénie près, 
que nous appelons catégorie des $\bvocB_\mQ$-modules. C'est la contrepartie du côté du topos de Faltings de la catégorie des $\mQ_p$-systèmes 
locaux de $X_{\oeta,\et}$.

 \subsection{}\label{intro4}
 Reprenant le schéma classique des correspondances construites par Fontaine, on a établi  dans \cite{agt} une correspondance entre certains 
$\bvocB_\mQ$-modules et certains fibrés de Higgs dans laquelle intervient un anneau de périodes de $\tE$ que nous avons appelé l'{\em algèbre de  Higgs-Tate}. 
C'est un modèle entier de l'anneau de Hyodo $B_\HT$,  que nous avons construit en utilisant la théorie des déformations et en s'inspirant 
de l'approche originelle de Faltings et du travail  d'Ogus-Vologodsky sur la transformée de Cartier en caractéristique $p$ \cite{ov}. 
Une complétion $p$-adique ``faible'' de l'algèbre de Higgs-Tate possède de bonnes propriétés cohomologiques conduisant à une équivalence 
entre les catégories d'objets admissibles correspondants, à savoir les {\em $\bvocB_\mQ$-modules de Dolbeault} et les {\em fibrés de Higgs solubles}. 

\subsection{}\label{intro5}
Nous rappelons dans cette nouvelle monographie la construction de cette correspondance, élargissant chemin faisant sa portée, afin d'établir les nouvelles
propriétés suivantes: 
\begin{itemize}
\item[(i)] Nous caractérisons les $\bvocB_\mQ$-modules de  {\em Hodge-Tate} parmi les $\bvocB_\mQ$-modules de Dolbeault, à savoir comme ceux 
dont le  fibré de Higgs associé est nilpotent \eqref{mht1}. 
\item[(ii)] Nous prouvons la fonctorialité de la correspondance de Simpson $p$-adique par image directe propre \eqref{crindmd21}, 
ce qui conduit à une généralisation de la suite spectrale de Hodge-Tate  relative \eqref{crindmd24}.
\end{itemize}  

Une autre avancée, nouvelle elle-aussi, méritant d'être mentionnée est la descente cohomologique des modules de Dolbeault dans le cas affine petit \eqref{mdpsa18} 
conduisant à une description de la catégorie des modules de Dolbeault en terme de leurs sections globales, permettant ainsi de comparer les théories 
locales et globales (\ref{mdpsa32} et \ref{mdpsa33}).

\subsection{}\label{intro6}
La correspondance de Simpson $p$-adique requiert  l'existence d'une déformation lisse du schéma $X\otimes_{\co_K}\co_C$ 
au-dessus  du $\rW(k)$-épaississement infinitésimal $p$-adique universel  de Fontaine $\cA_2(\co_\oK)$ (\ref{eipuf9} et \ref{eipuf10}), et elle dépend
du choix d'une telle déformation, alors que la théorie des $\bvocB_\mQ$-modules de Hodge-Tate n'en dépend pas. Cet apparent paradoxe est résolu en introduisant 
une nouvelle version de notre construction de la correspondance de Simpson $p$-adique qui dépend cette fois du choix d'une déformation 
du schéma $X\otimes_{\co_K}\co_C$ au-dessus d'une version logarithmique $\cA^*_2(\co_\oK/\co_K)$ du $\co_K$-épaississement infinitésimal $p$-adique universel  
de Fontaine  $\cA_2(\co_\oK/\co_K)$ (\eqref{eipuf2f} et \eqref{epinflog7a}). 
Nous référerons à ce contexte comme étant celui du {\em cas relatif} et au précédent contexte comme étant celui du {\em cas absolu}. 
Nous traitons simultanément les deux cas, chacun ayant ses avantages et ses inconvénients.
Comme  $\cA^*_2(\co_\oK/\co_K)$ est de manière naturelle une $\co_K$-algèbre, $X\otimes_{\co_K}\co_C$  admet une $\cA^*_2(\co_\oK/\co_K)$-déformation lisse canonique, 
à savoir le changement de base de $X$. 
La théorie dans le cas relatif a néanmoins ses limites. En effet, on peut montrer que dans le cas affine, les modules de Dolbeault sont petits \eqref{pmh33}; 
grossièrement dit, ils sont ``triviaux'' modulo une puissance de $p$ prescrite. Cette puissance de $p$ est indépendante de $K$ seulement dans le cas absolu. 
Sa dépendance en $K$ limite drastiquement le champ d'applications de la théorie dans le cas relatif, 
spécialement si l'on veut étendre la correspondance à tous les $\bvocB_\mQ$-modules par descente.  
Cependant, les $\bvocB_\mQ$-modules de Hodge-Tate  sont de Dolbeault à la fois dans le cas absolu et dans le cas relatif, rendant ainsi notre
définition des $\bvocB_\mQ$-modules de Hodge-Tate  indépendante de tout choix d'une déformation. 

\subsection{}\label{intro7}
Un autre progrès apporté  par la présente monographie à \cite{agt} 
est l'extension de la théorie développée pour la catégorie des $\bvocB_\mQ$-modules à celle plus grosse des {\em ind-$\bvocB$-modules}. 
Cette dernière admet en effet des limites inductives filtrantes et a de meilleures propriétés (\ref{indsh}, \ref{indmod} et \cite{ks2}).  Cette extension qui peut donc
sembler à première vue technique s'avère en fait être très utile et nécessaire pour l'étude de la fonctorialité de la correspondance $p$-adique par image directe propre.

\subsection{}\label{intro8}
Nous donnons dans ce premier chapitre un survol détaillé du contenu de cette monographie. 
Nous traitons le cas des schémas à singularités toriques en utilisant la géométrie logarithmique mais, 
par simplicité, nous nous restreignons dans ce survol au cas lisse.  

\subsection*{Remerciements} Nous tenons à remercier chaleureusement T. Tsuji pour les nombreux conseils qu'il 
nous a très généreusement accordés tout au long de ce travail.

\section{Topos de Faltings}\label{ft}

\subsection{}\label{ft1}
Soient $X$ un $S$-schéma lisse,  $E$ la catégorie des morphismes $(V\rightarrow U)$
au-dessus du morphisme  canonique $X_\oeta\rightarrow X$, c'est-à-dire  des diagrammes commutatifs
\begin{equation}\label{ft1a}
\xymatrix{V\ar[r]\ar[d]&U\ar[d]\\
X_\oeta\ar[r]&X}
\end{equation}
tels que  $U$ soit étale sur $X$ et que le  morphisme canonique $V\rightarrow U_\oeta$ soit {\em fini  étale}. 
Il est pratique de  considérer la catégorie $E$ comme fibrée par le foncteur
\begin{equation}\label{ft1b}
\pi\colon E\rightarrow \Et_{/X}, \ \ \ (V\rightarrow U)\mapsto U,
\end{equation}
au-dessus du site étale de $X$. 

La fibre de $\pi$ au-dessus d'un objet $U$ de $\Et_{/X}$ est canoniquement équivalente à la catégorie $\Et_{\rf/U_\oeta}$ des morphismes 
finis étales au-dessus de $U_\oeta$. Nous la munissons de la topologie étale  et notons $U_{\oeta,\fet}$ le topos associé. 
Si $U_\oeta$ est connexe et si $\oy$ est un point géométrique de $U_\oeta$, le topos $U_{\oeta,\fet}$ est alors équivalent au
topos classifiant du groupe profini  $\pi_1(U_\oeta,\oy)$, {\em i.e.}, à la catégorie des ensembles discrets munis d'une action à gauche continue
de $\pi_1(U_\oeta,\oy)$.

\subsection{}\label{ft2}
Nous munissons $E$ de la topologie {\em co-évanescente} (\cite{agt} VI.10.1), c'est-à-dire de la topologie engendrée par les recouvrements 
$\{(V_i\rightarrow U_i)\rightarrow (V\rightarrow U)\}_{i\in I}$
des deux types suivants~:
\begin{itemize}
\item[(v)] $U_i=U$ pour tout  $i\in I$ et $(V_i\rightarrow V)_{i\in I}$ est un recouvrement ;
\item[(c)] $(U_i\rightarrow U)_{i\in I}$ est un recouvrement et $V_i=V\times_UU_i$ pour tout  $i\in I$. 
\end{itemize}
 
Nous notons $\tE$ le topos des faisceaux d'ensembles sur $E$. 

Se donner un faisceau $F$ sur $E$ est équivalent à se donner:
\begin{itemize}
\item[(i)] pour tout objet $U$ de $\Et_{/X}$, un faisceau $F_U$ de $U_{\oeta,\fet}$, à savoir la restriction de $F$ à la fibre
de $\pi$ au-dessus $U$;
\item[(ii)] pour tout morphisme $f\colon U'\rightarrow U$ de $\Et_{/X}$, un morphisme $\gamma_f\colon F_U\rightarrow f_{\oeta*}(F_{U'})$. 
\end{itemize}

Ces données doivent satisfaire une condition de cocycle (pour la composition des morphismes) et une condition de 
recollement (pour les recouvrements de $\Et_{/X}$). Nous écrirons $F=\{U \mapsto F_U\}$ et le verrons donc comme un faisceau sur $\Et_{/X}$ 
à valeurs dans les topos finis étales.

\subsection{}\label{ft20}
Toute application de spécialisation $\oy\rightsquigarrow \ox$ d'un point géométrique $\oy$ de $X_\oeta$ vers un point géométrique $\ox$ de $X$, 
détermine un point de $\tE$ que l'on notera par $\rho(\oy\rightsquigarrow \ox)$. {\em La collection de ces points de $\tE$ est conservative} (\cite{agt} VI.10.21).  

\subsection{}\label{ft3}
Il existe trois morphismes de topos (\cite{agt} VI.10.6 et VI.10.7)
\begin{equation}\label{ft3a}
\xymatrix{
X_{\oeta,\et}\ar[r]^-(0.5){\psi}&{\tE}\ar[r]^-(0.5){\sigma}\ar[d]^{\beta}&X_\et\\
&X_{\oeta,\fet}&}
\end{equation}
définis par 
\begin{eqnarray}
(V\mapsto U)\in \ob(E) &\mapsto& \psi^*(V\rightarrow U)=V,\label{ft3d}\\
U\in \ob(\Et_{/X}) &\mapsto& \sigma^*(U)=(U_\oeta\rightarrow U)^a,\label{ft3b}\\
V\in \ob(\Et_{\rf/X_\oeta}) &\mapsto& \beta^*(V)=(V\rightarrow X)^a,\label{ft3c}
\end{eqnarray}
où l'exposant $a$ désigne le faisceau associé.

\subsection{}\label{ft30}
Les images directes supérieures de $\sigma$ faisceautisent la cohomologie galoisienne (\cite{agt} VI.10.40): si  $F=\{U\mapsto F_U\}$ 
est un groupe abélien de $\tE$, pour tout entier $i\geq 0$, $\rR^i\sigma_*(F)$ est canoniquement isomorphe au faisceau associé au préfaisceau 
\begin{equation}\label{ft30a}
U\mapsto \rH^i(U_{\oeta,\fet}, F_U).
\end{equation}

\begin{prop}[\cite{ag} 4.4.2] \label{ft4}
Pour tout faisceau abélien localement constant constructible de torsion $F$ de $X_{\oeta,\et}$, on a $\rR^i\psi_*(F)=0$ pour tout $i\geq 1$.
\end{prop}

Ce résultat est une conséquence du fait que pour tout point géométrique $\ox$ de $X$ au-dessus de $s$, notant $\uX$ le localisé strict
de $X$ en $\ox$, $\uX_\oeta$ est un schéma $K(\pi,1)$. Cette propriété a été prouvée par Faltings (\cite{faltings1} Lemma 2.3 page 281), 
étendant ainsi des  résultats antérieurs d'Artin (\cite{sga4} XI). Elle a été généralisée plus avant encore par Achinger au cas log-lisse (\cite{achinger}, \cite{ag} 4.3.6).

\subsection{}\label{ft5}
Pour tout objet $(V\rightarrow U)$ de $E$, on notera $\oU^V$ la clôture intégrale de $\oU=U\times_S\oS$ dans $V$ et l'on posera
\begin{equation}\label{ft5a}
\ocB(V\rightarrow U)=\Gamma(\oU^V,\co_{\oU^V}).
\end{equation} 
Le  préfaisceau sur $E$ ainsi défini est en fait un faisceau. Comme pour tout faisceau, on peut écrire $\ocB=\{U\mapsto \ocB_U\}$.

Soient $U=\Spec(R)$ un $X$-schéma étale, $\oy$ un point géométrique de $U_\oeta$. La fibre $\ocB_{U,\oy}$ peut se décrire comme suit.   
Notons $(V_i)_{i\in I}$ le revêtement universel de $U_\oeta$ en $\oy$.
Pour tout  $i\in I$, soit $U_i=\Spec(R_i)$ la normalisation de $\oU$ dans $V_i$
\begin{equation}\label{ht5b}
\xymatrix{
V_i\ar[r]\ar[d]&U_i\ar[d]\\
U_\oeta\ar[r]&\oU}
\end{equation}
La fibre $\ocB_{U,\oy}$ est alors isomorphe à la $\co_\oK$-représentation suivante de $\pi_1(U_\oeta,\oy)$:
\begin{equation}\label{ht5c}
\oR=\underset{\underset{i\in I}{\longrightarrow}}{\lim} \ R_i. 
\end{equation}

\vspace{2mm}

En utilisant la théorie  d'Artin-Schreier, Faltings a démontré le raffinement suivant de \ref{ft4}:

\begin{teo}[\cite{faltings2}, \cite{ag} 4.8.13]\label{ft6}
Pour tout faisceau localement constant constructible de $(\mZ/p^n\mZ)$-modules $F$ de $X_{\oeta,\et}$, le morphisme  canonique 
\begin{equation}\label{grfmc2a} 
\rH^i(X_{\oeta,\et},F)\otimes_{\mZ_p}\co_C\rightarrow \rH^i(\tE,\psi_*(F)\otimes_{\mZ_p}\ocB)
\end{equation}
est un {\em presque-isomorphisme}, i.e., son noyau et son conoyau sont annulés par $\fm_C$. 
\end{teo}

Faltings déduit ensuite tous les théorèmes de comparaison entre la cohomologie étale $p$-adique
et les autres cohomologies $p$-adiques de ce théorème fondamental de comparaison.

\section{Théorie locale. Le torseur des déformations} \label{tordef}

\subsection{}
Nous rappelons tout d'abord la variante locale de la correspondance de Simpson $p$-adique pour les $S$-schémas affines {\em petits}.
Notre approche utilise un anneau de périodes, l'{\em algèbre de Higgs-Tate}, que nous avons introduit dans  \cite{agt}.

\subsection{}\label{tordef1}
Rappelons tout d'abord (\ref{eipo3}, \cite{fontaine3} 1.2.1, \cite{tsuji1} 1.1) que Fontaine a associé fonctoriellement à toute $\mZ_{(p)}$-algèbre $A$ telle que $A/pA\not=0$, l'anneau  
\begin{equation}\label{tordef1a}
A^\flat=\underset{\underset{\mN}{\longleftarrow}}{\lim}\ A/pA,
\end{equation} 
où les morphismes de transition sont les morphismes de Frobenius absolu de $A/pA$ ainsi que l'homomorphisme d'anneaux
\begin{equation}\label{tordef1b}
\theta\colon \rW(A^\flat)\rightarrow \hA,
\end{equation}
des vecteurs de Witt de $A^\flat$ vers la complétion $p$-adique de $A$  
défini, pour $x=(x_0,x_1,\dots)\in \rW(A^\flat)$ par 
\begin{equation}
\theta(x)=\underset{m\rightarrow +\infty}{\lim}\ (\tx_{0m}^{p^m}+p\tx_{1m}^{p^{m-1}}+\dots+p^m\tx_{mm}),
\end{equation}
où, pour tout  $n\geq 0$, on écrit $x_n=(x_{nm})_{m\geq 0}\in A^\flat$ et où, pour $x\in A/pA$, $\tx$ désigne un relèvement dans $A$. 

L'anneau $A^\flat$ est parfait de caractéristique $p$, et l'homomorphisme $\theta$ est surjectif si  le morphisme de Frobenius absolu de $A/pA$ est surjectif.

\subsection{}\label{tordef2}
L'anneau $\co_{\oK^\flat}=(\co_\oK)^\flat$  est un anneau de valuation non-discrète  complet,  de hauteur $1$. On note $\oK^\flat$ son corps de fractions.  
On choisit une suite $(p_n)_{n\geq 0}$ d'éléments de $\co_\oK$ telle que $p_0=p$ et $p_{n+1}^p=p_n$ pour tout  $n\geq 0$. 
On note $\varpi$ l'élément de $\co_{\oK^\flat}$ défini par la famille des $(p_n)$ et l'on pose
\begin{equation}\label{tordef2a}
\xi=[\varpi]-p\in \rW(\co_{\oK^\flat}).
\end{equation}
C'est un générateur du noyau de $\theta$ \eqref{tordef1b}. On pose
\begin{equation}\label{tordef2c}
\cA_2(\co_\oK)=\rW(\co_{\oK^\flat})/\ker(\theta)^2.
\end{equation}
On a alors une suite exacte \eqref{definf3e}
\begin{equation}\label{tordef2d}
0\longrightarrow \co_C\stackrel{\cdot \xi}{\longrightarrow} \cA_2(\co_\oK)
\stackrel{\theta}{\longrightarrow} \co_C \longrightarrow 0.
\end{equation}

On dispose également d'un homomorphisme canonique $\mZ_p(1)\rightarrow \co_{\oK^\flat}^\times$. 
Pour tout $\zeta\in \mZ_p(1)$, on a $\theta([\zeta]-1)=0$. On en déduit un homomorphisme de groupes
\begin{equation}\label{tordef2e}
\mZ_p(1)\rightarrow \cA_2(\co_{\oK}),\ \ \ 
\zeta\mapsto\log([\zeta])=[\zeta]-1,
\end{equation}
dont  l'image engendre l'idéal $p^{\frac{1}{p-1}}\xi\co_C$ de $\cA_2(\co_{\oK})$. Cet homomorphisme induit un isomorphisme $\co_C$-linéaire \eqref{definf17c} 
\begin{equation}\label{tordef2f}
\co_C(1)\stackrel{\sim}{\rightarrow} p^{\frac{1}{p-1}}\xi \co_C.
\end{equation}

\subsection{}\label{tordef20}
Nous introduisons dans \ref{definf4} une version relative logarithmique de l'extension $\cA_2(\co_\oK)$ \eqref{tordef2d} sur $\co_K$. 
On fixe pour cela une uniformisante $\pi$ de $\co_K$, une suite $(\pi_n)_{n\geq 0}$ 
d'éléments de $\co_\oK$ telle que $\pi_0=\pi$ et $\pi_{n+1}^p=\pi_n$ (pour tout  $n\geq 0$) et l'on note $\upi$ l'élément correspondant de $\co_{\oK^\flat}$. 
On pose
\begin{equation}\label{tordef20a}
\rW_{\co_K}(\co_{\oK^\flat})=\rW(\co_{\oK^\flat})\otimes_{\rW(k)}\co_K. 
\end{equation}
On note alors $\rW^{\ast}_{\co_K}(\co_{\oK^\flat})$ la sous-$\rW_{\co_K}(\co_{\oK^\flat})$-algèbre de $\rW_K(\co_{\oK^\flat})=\rW(\co_{\oK^\flat})\otimes_{\rW(k)}K$ 
engendrée par $[\upi]/\pi$ puis l'on pose
\begin{equation}\label{tordef20b}
\xi^{\ast}_\pi=\frac{[\upi]}{\pi}-1\in \rW^{\ast}_{\co_K}(\co_{\oK^\flat}).
\end{equation}
C'est un générateur du noyau de l'homomorphisme $\theta^{\ast}_{\co_K}\colon \rW^{\ast}_{\co_K}(\co_{\oK^\flat})\rightarrow \co_C$ induit par $\theta$ \eqref{tordef1b}.
On observera que cette algèbre dépend de $(\pi_n)_{n\geq 0}$  \eqref{definf4}. On pose
\begin{equation}\label{tordef20c}
\cA^{\ast}_2(\co_\oK/\co_K)=\rW^{\ast}_{\co_K}(\co_{\oK^\flat})/(\xi^{\ast}_\pi)^2 \rW^{\ast}_{\co_K}(\co_{\oK^\flat}).
\end{equation}
On a alors une suite exacte \eqref{definf4i}
\begin{equation}\label{tordef20d}
\xymatrix{
0\ar[r]&{\co_C}\ar[r]^-(0.5){\cdot \xi^{\ast}_\pi}&{\cA^{\ast}_2(\co_\oK/\co_K)}\ar[r]^-(0.5){\theta^{\ast}_{\co_K}}&{\co_C}\ar[r]& 0}.
\end{equation}

Notant $K_0$ le corps des fractions de $\rW(k)$ et   $\fd$ la différente de $K/K_0$, 
l'homomorphisme canonique $\cA_2(\co_\oK)\rightarrow \cA^{\ast}_2(\co_\oK/\co_K)$ induit un isomorphisme $\co_C$-linéaire \eqref{definf16}
\begin{equation}\label{tordef20e}
\xi\co_C\stackrel{\sim}{\rightarrow}\pi\fd\xi^{\ast}_\pi \co_C.
\end{equation}

\subsection{}\label{tordef3}
Soit $X=\Spec(R)$ un $S$-schéma lisse qui est  affine {\em petit}  au sens de Faltings ({\em i.e.} qui admet un $S$-morphisme étale
$X\rightarrow \mG_{m,S}^d=\Spec(\co_K[T_1^{\pm1},\dots, T_d^{\pm1}])$  pour un entier $d\geq 0$) et aussi  
tel que $X_s\not=\emptyset$. Nous fixons un point géométrique $\oy$ de $X_\oeta$ et nous notons $\oX^\star$ (resp. $X_\oeta^\star$) 
la composante connexe de $\oX=X\times_S\oS$ (resp. $X_\oeta$) contenant l'image de $\oy$; on a $X_\oeta^\star=\oX^\star\times_\oS\oeta$. 
Nous posons $\Delta=\pi_1(X_\oeta^\star,\oy)$, $\oR=\ocB_{X,\oy}$ \eqref{ht5c} et 
\begin{equation}
R_1=\Gamma(\oX^\star,\co_{\oX}).
\end{equation}
Soit $\hRun$ (resp. $\hoR$) le séparé complété $p$-adique de $R_1$ (resp. $\oR$).  Rappelons que $\Delta$ agit naturellement sur $\oR$ par des homomorphismes d'anneaux. 
Nous noterons $\bRep_{\hoR}(\Delta)$ la catégorie des $\hoR$-représentations de $\Delta$ \eqref{notconv16}.

Appliquant les constructions de \ref{tordef1}--\ref{tordef20} à l'algèbre $\oR$ (\ref{taht2} et \ref{taht3}), on pose 
\begin{eqnarray}
\cA_2(\oR)=\rW(\oR^\flat)/\xi^2\rW(\oR^\flat),\label{tordef3b}\\
\cA^*_2(\oR/\co_K)=\rW^*_{\co_K}(\oR^\flat)/(\xi^{\ast}_\pi)^2 \rW_{\co_K}(\oR^\flat).\label{tordef3bb}
\end{eqnarray}
On a donc des suites exactes
\begin{equation}\label{tordef3c}
0\longrightarrow \hoR\stackrel{\cdot \xi}{\longrightarrow} \cA_2(\oR)\longrightarrow \hoR \longrightarrow 0,
\end{equation}
\begin{equation}\label{tordef3cc}
0\longrightarrow \hoR\stackrel{\cdot \xi^{\ast}_\pi}{\longrightarrow} \cA^*_2(\oR/\co_K)\longrightarrow \hoR \longrightarrow 0.
\end{equation}

\subsection{}\label{tordef21}
Soit $\tS$ l'un des schémas $\Spec(\cA_2(\co_\oK))$ ou  $\Spec(\cA^{\ast}_2(\co_\oK/\co_K))$; dans le premier cas 
on dira qu'on est dans le  cas {\em absolu} et dans le second dans le cas {\em relatif}. 
On notera  
\begin{equation}\label{tordef21a}
i_S\colon \Spec(\co_C)\rightarrow \tS
\end{equation}
l'immersion fermée définie par l'idéal de carré nul engendré par $\txi=\xi$ dans le cas absolu et par $\txi=\xi^{\ast}_\pi$
dans le cas relatif. 

On notera que dans le cas relatif, $\tS$ est naturellement un $S$-schéma.

\subsection{}\label{tordef210}
En adéquation avec le cadre, absolu ou relatif, dans lequel on va travailler, nous notons $\tmX$ 
un des schémas $\Spec(\cA_2(\oR))$ ou $\Spec(\cA^*_2(\oR/\co_K))$. Il existe alors une immersion fermée canonique
\begin{equation}\label{tordef21b}
i_X\colon \Spec(\hoR)\rightarrow \tmX
\end{equation}
au-dessus de l'immersion fermée $i_S$  définie par l'idéal de $\co_\tmX$ engendré par $\txi$. 

Pour toute $\hRun$-algèbre $A$, on aura à considérer les $A$-modules de Higgs à coefficients dans $\txi^{-1}\Omega^1_{R/\co_K}\otimes_RA$ \eqref{MH1} 
et l'on dira abusivement qu'ils sont à coefficients dans $\txi^{-1}\Omega^1_{R/\co_K}$. 
La catégorie de ces modules sera notée $\bMH(A,\txi^{-1}\Omega^1_{R/\co_K})$. 

\subsection{}\label{tordef4}
La correspondance de Simpson $p$-adique dépend du choix d'une $\tS$-déformation lisse  $\tX$ de $X\otimes_{\co_K}\co_C$, 
\begin{equation}\label{tordef4a}
\xymatrix{
{X\otimes_{\co_K}\co_C}\ar[r]\ar[d]\ar@{}[rd]|\Box&{\tX}\ar[d]\\
{\Spec(\co_C)}\ar[r]&{\tS}}
\end{equation}
Comme $X$ est affine, une telle déformation existe toujours et est unique à isomorphisme non-unique près. 
Nous en fixons une dans la suite. 

Soient $U$ un  sous-schéma ouvert de $\Spec(\hoR)$ et $\tU$ le sous-schéma ouvert de
$\tmX$ défini par $U$ \eqref{tordef21b}.  On note $\cL(U)$ l'ensemble des morphismes représentés 
par les flèches en pointillés complétant le diagramme 
\begin{equation}\label{tordef4b}
\xymatrix{
{U}\ar[d]\ar[r]&{\tU}\ar@{.>}[d]\ar@/^2pc/[dd]\\
{X\otimes_{\co_K}\co_C}\ar[r]\ar[d]\ar@{}[rd]|\Box&{\tX}\ar[d]\\
{\Spec(\co_C)}\ar[r]&{\tS}}
\end{equation}
de telle sorte qu'il reste commutatif. Le foncteur $U\mapsto \cL(U)$ est un torseur pour la topologie de Zariski sur $\Spec(\hoR)$ 
sous le $\hoR$-module $\Hom_{R}(\Omega^1_{R/\co_K},\txi\hoR)$ \eqref{taht5}. 
Un tel torseur se décrit aisément. En effet, soit $\cF$ le $\hoR$-module des fonctions affines sur $\cL$ (\cite{agt} II.4.9). 
Ce dernier s'insère dans une suite exacte canonique 
\begin{equation}\label{tordef4c}
0\rightarrow \hoR\rightarrow \cF\rightarrow \txi^{-1}\Omega^1_{R/\co_K} \otimes_R \hoR\rightarrow 0.
\end{equation} 
Considérons la $\hoR$-algèbre
\begin{equation}\label{tordef4d}
\cC=\underset{\underset{n\geq 0}{\longrightarrow}}\lim\ \Sym^n_{\hoR}(\cF),
\end{equation}
dans laquelle les morphismes de transition sont définis en envoyant $x_1\otimes\dots \otimes x_n$ sur $1\otimes x_1\otimes\dots \otimes x_n$. Le foncteur $\cL$ est 
alors représenté par $\Spec(\cC)$ (\cite{agt} II.4.10). 

L'action naturelle de $\Delta$ sur $\oR$ induit une action sur le schéma $\tmX$, et donc une action $\hoR$-semi-linéaire sur $\cF$, 
telle que les morphismes de \eqref{tordef4c} soient $\Delta$-équivariants \eqref{taht6}. 
On en déduit une action de $\Delta$ sur $\cC$ par des homomorphismes d'anneaux. Ces actions sont continues pour la topologie $p$-adique (\cite{agt} II.12.4). 
La $\hoR$-algèbre $\cC$ munie de cette action de $\Delta$ s'appelle l'{\em algèbre de Higgs-Tate} associée à la déformation $\tX$. 

\subsection{}\label{tordef40}
L'algèbre de Higgs-Tate $\cC$ est un modèle entier de l'anneau de Hyodo (\cite{agt} II.15.6). 
Nous en introduisons une complétion $p$-adique ``faible'' qui nous sert d'anneau de périodes pour la correspondance de Simpson $p$-adique. 
Pour tout nombre rationnel $r\geq 0$, nous notons $\cF^{(r)}$ la $\hoR$-représentation de $\Delta$ déduite de $\cF$ par image inverse  
par le morphisme de multiplication par $p^r$ sur $\txi^{-1}\Omega^1_{R/\co_K}\otimes_R\hoR$, si bien qu'on a une suite exacte
\begin{equation}\label{tordef40a}
0\rightarrow \hoR\rightarrow \cF^{(r)}\rightarrow \txi^{-1}\Omega^1_{R/\co_K}\otimes_R\hoR\rightarrow 0.
\end{equation}
On considère la $\hoR$-algèbre
\begin{equation}\label{tordef40b}
\cC^{(r)}=\underset{\underset{n\geq 0}{\longrightarrow}}\lim\ \Sym^n_{\hoR}(\cF^{(r)}),
\end{equation}
où les morphismes de transition sont définis en envoyant $x_1\otimes\dots \otimes x_n$ sur $1\otimes x_1\otimes\dots \otimes x_n$. 
L'action de $\Delta$ sur $\cF^{(r)}$ induit une action sur $\cC^{(r)}$ par des
automorphismes d'anneaux, compatible avec son action sur $\hoR$.  
On notera $\hcC^{(r)}$ la complétion $p$-adique  de $\cC^{(r)}$.

Pour tous nombres rationnels $r'\geq r\geq 0$, on a un $\hoR$-homomorphisme injectif et
$\Delta$-équivariant canonique
$\alpha^{r,r'}\colon \cC^{(r')}\rightarrow \cC^{(r)}$. On vérifie aisément que l'homomorphisme induit
$\halpha^{r,r'}\colon\hcC^{(r')}\rightarrow \hcC^{(r)}$ est injectif. On pose
\begin{equation}\label{tordef40c}
\hcC^{(r+)}=\underset{\underset{t\in \mQ_{>r}}{\longrightarrow}}{\lim} \hcC^{(t)},
\end{equation}
que l'on identifie à une sous-$\hoR$-algèbre de $\hcC=\hcC^{(0)}$.  Le
groupe $\Delta$ agit naturellement sur $\hcC^{(r+)}$ par automorphismes d'anneaux de façon
compatible avec ses actions sur $\hoR$ et sur $\hcC$.

On notera  
\begin{equation}\label{tordef40d}
d_{\cC^{(r)}}\colon \cC^{(r)}\rightarrow \txi^{-1}\Omega^1_{R/\co_K}\otimes_R\cC^{(r)}
\end{equation}
la $\hoR$-dérivation universelle de $\cC^{(r)}$ \eqref{taht11} et  
\begin{equation}\label{tordef40e}
d_{\hcC^{(r)}}\colon \hcC^{(r)}\rightarrow \txi^{-1}\Omega^1_{R/\co_K}\otimes_R\hcC^{(r)}
\end{equation}
son prolongement aux complétés $p$-adiques. Ces dérivations sont clairement $\Delta$-équivariantes, 
et ce sont des  $\hoR$-champs de Higgs à coefficients dans $\txi^{-1}\Omega^1_{R/\co_K}$ d'après \ref{MH8}(i) (cf. \ref{taht11}). 

Pour tous  nombres rationnels $r'\geq r\geq 0$, on a
\begin{equation}\label{tordef40f}
p^{r'}(\id \times \alpha^{r,r'}) \circ d_{\cC^{(r')}}=p^rd_{\cC^{(r)}}\circ \alpha^{r,r'}.
\end{equation}
Les dérivations $(p^td_{\hcC^{(t)}})_{t\in \mQ_{>r}}$ induisent donc une $\hoR$-dérivation 
\begin{equation}\label{tordef40g}
d^{(r)}_{\hcC^{(r+)}}\colon \hcC^{(r+)}\rightarrow \txi^{-1}\Omega^1_{R/\co_K} \otimes_R\hcC^{(r+)},
\end{equation}
qui est aussi la restriction de $p^rd_{\hcC^{(r)}}$ à $\hcC^{(r+)}$.

\subsection{}\label{tordef41}
Pour toute $\hoR$-représentation $M$ de $\Delta$ \eqref{notconv16}, on notera $\mH(M)$ le $\hRun$-module défini  par
\begin{equation}\label{tordef41a}
\mH(M)=(M\otimes_{\hoR}\hcC^{(0+)})^\Delta.
\end{equation}
On le munit du $\hRun$-champ de Higgs à coefficients dans $\txi^{-1}\Omega^1_{R/\co_K}$ induit par $d^{(0)}_{\hcC^{(0+)}}$ \eqref{tordef40g}.   
On a ainsi défini un foncteur 
\begin{equation}\label{tordef41b}
\mH\colon \bRep_{\hoR}(\Delta) \rightarrow \bMH(\hRun,\txi^{-1}\Omega^1_{R/\co_K}).
\end{equation}

\subsection{}\label{tordef42}
Pour tout $\hRun$-module de Higgs $(N,\theta)$ à coefficients dans $\txi^{-1}\Omega^1_{R/\co_K}$, on notera  $\mV(N)$ le $\hoR$-module défini  par 
\begin{equation}\label{tordef42a}
\mV(N)=(N\otimes_{\hRun}\hcC^{(0+)})^{\theta_\tot=0},
\end{equation}
où $\theta_\tot=\theta\otimes \id+\id\otimes d^{(0)}_{\hcC^{(0+)}}$ est le $\hRun$-champ de Higgs total sur $N\otimes_{\hRun}\hcC^{(0+)}$.
On munit ce module de l'action $\hoR$-semi-linéaire de $\Delta$ induite par son action naturelle sur $\hcC^{(0+)}$. 
On a ainsi défini un foncteur 
\begin{equation}\label{tordef42b}
\mV\colon \bMH(\hRun,\txi^{-1}\Omega^1_{R/\co_K})\rightarrow \bRep_{\hoR}(\Delta).
\end{equation}

\begin{defi}[cf. \ref{repdolb8}] \label{tordef43}
On dira qu'une $\hoR[\frac 1 p]$-représentation $M$ de $\Delta$ est de {\em Dolbeault} 
si les conditions suivantes sont satisfaites:
\begin{itemize}
\item[(i)] $\mH(M)$ est un $\hRun[\frac 1 p]$-module projectif de type fini;
\item[(ii)] le morphisme  canonique 
\begin{equation}\label{tordef43a}
\mH(M) \otimes_{\hRun}\hcC^{(0+)}\rightarrow  M\otimes_{\hoR}\hcC^{(0+)}
\end{equation}
est un isomorphisme.
\end{itemize}
\end{defi}

Nous prouvons que toute $\hoR[\frac 1 p]$-représentation de $\Delta$ qui est de Dolbeault est continue pour la topologie $p$-adique \eqref{repdolb9}. 

\begin{defi}[cf. \ref{repdolb10}] \label{tordef44}
On dira qu'un $\hRun[\frac 1 p]$-module de Higgs $(N,\theta)$ à coefficients dans $\txi^{-1}\Omega^1_{R/\co_K}$
est {\em soluble} si  les conditions suivantes sont satisfaites:
\begin{itemize}
\item[(i)] $N$ est un $\hRun[\frac 1 p]$-module projectif de type fini;
\item[(ii)] le morphisme  canonique
\begin{equation}\label{tordef44a}
\mV(N) \otimes_{\hoR}\hcC^{(0+)}\rightarrow  N\otimes_{\hRun}\hcC^{(0+)}
\end{equation}
est un isomorphisme.
\end{itemize}
\end{defi}

Les notions \ref{tordef43} et \ref{tordef44} ne dépendent pas du choix de la $\tS$-déformation $\tX$ \eqref{repdolb7}. 
Néanmoins elles dépendent a priori du cadre, absolu ou relatif, dans lequel on s'est placé dans  \ref{tordef21}.
Elles sont équivalentes aux conditions de {\em petitesse}  introduites par Faltings \cite{faltings3}, {\em i.e.}, à des conditions de trivialité modulo des puissances prescrites de $p$ 
des objets considérés (\ref{pmh32}--\ref{pmh34}). 

\begin{prop}[cf. \ref{repdolb14}] \label{tordef45}
Les foncteurs $\mH$ \eqref{tordef41a} et $\mV$ \eqref{tordef42a} induisent des équivalences de catégories, quasi-inverses l'une de l'autre,
entre la catégorie des $\hoR[\frac 1 p]$-représentations de Dolbeault de $\Delta$  et celle des  
$\hRun[\frac 1 p]$-modules de Higgs solubles à coefficients dans $\txi^{-1}\Omega^1_{R/\co_K}$.
\end{prop}

\begin{prop}[cf. \ref{repdolb15}]\label{htls30}
Pour toute $\hoR[\frac 1 p]$-représentation $M$ de Dolbeault de $\Delta$,  
il existe un isomorphisme canonique fonctoriel  dans $\bD^+(\bMod(\hRun[\frac 1 p]))$ 
\begin{equation}
\rC_\cont^\bullet(\Delta, M)\stackrel{\sim}{\rightarrow} \mK^\bullet(\mH(M)),
\end{equation}
où $\rC_\cont^\bullet(\Delta, M)$ est le complexe des cochaînes continues de $\Delta$ à valeurs dans $M$
et $\mK^\bullet(\mH(M))$ est le complexe de Dolbeault du $\hRun[\frac 1 p]$-module de Higgs $\mH(M)$ associé à $M$ \eqref{tordef41b}. 
\end{prop}

\begin{prop}[cf. \ref{repht5}]\label{htls3}
Soit $M$ un $\hoR[\frac 1 p]$-module projectif de type fini muni d'une action $\hoR[\frac 1 p]$-semi-linéaire de $\Delta$. 
Les propriétés suivantes sont  alors équivalentes:
\begin{itemize}
\item[{\rm (i)}] La $\hoR[\frac 1 p]$-représentation $M$ de $\Delta$ est de Dolbeault 
et le $\hRun[\frac 1 p]$-module de Higgs associé $(\mH(M),\theta)$ \eqref{tordef41a} est nilpotent, i.e., il existe une filtration finie décroissante
$(\mH_i)_{0\leq i\leq n}$ de $\mH(M)$ par des sous-$\hRun[\frac 1 p]$-modules telle que $\mH_0=\mH(M)$, $\mH_n=0$ et pour tout  $0\leq i\leq n-1$, on a
\begin{equation}
\theta(\mH_i)\subset \txi^{-1}\Omega^1_{R/\co_K}\otimes_R\mH_{i+1}.
\end{equation}
\item[{\rm (ii)}] Il existe un $\hRun[\frac 1 p]$-module projectif de type fini $N$, un $\hRun[\frac 1 p]$-champ de Higgs $\theta$ sur $N$ à coefficients 
dans $\txi^{-1}\Omega^1_{R/\co_K}$ et un isomorphisme $\cC$-linéaire et $\Delta$-équivariant de $\hoR[\frac 1 p]$-modules  de Higgs
\begin{equation}
N \otimes_{\hRun}\cC\stackrel{\sim}{\rightarrow}  M\otimes_{\hoR}\cC.
\end{equation}
\end{itemize}
De plus, si ces conditions sont vérifiées, on a un isomorphisme de $\hRun[\frac 1 p]$-modules de Higgs
\begin{equation}
\mH(M)\stackrel{\sim}{\rightarrow}  (N,\theta).
\end{equation}
\end{prop}

\begin{defi}[cf. \ref{repdolb23}]\label{htls4}
On dit qu'une $\hoR[\frac 1 p]$-représentation $M$ de $\Delta$ est de {\em Hodge-Tate} si elle satisfait les conditions équivalentes \ref{htls3}.
\end{defi}

Cette notion ne dépend pas du choix de la $\tS$-déformation $\tX$, ni même que l'on soit dans le cas absolu ou relatif \eqref{tordef21}.

\begin{rema}
Tsuji \cite{tsuji5} a développé une version arithmétique de la correspondance de Simpson $p$-adique locale.
Il associe à une $\hoR[\frac 1 p]$-représentation $p$-adique de $\Gamma=\pi_1(X^\star_\eta,\oy)$ \eqref{tordef3} 
un champ de Higgs et un opérateur de Sen arithmétique satisfaisant une relation de compatibilité qui force
le champ de Higgs à être nilpotent (voir aussi \cite{he2}). Cela explique la relation existant entre le travail de  Liu et Zhu \cite{lz} 
et le point de vue développé ici. Le lecteur prendra garde cependant que notre notion de systèmes locaux de Hodge-Tate ne correspond pas à celle 
de Liu et Zhu puisque nous considérons des systèmes locaux géométriques alors qu'ils considèrent des systèmes locaux arithmétiques.
\end{rema}

\section{Théorie globale. Modules de Dolbeault}\label{dolb}

\subsection{}
La correspondance de Simpson $p$-adique décrite dans §~\ref{tordef} ne peut aisément se recoller pour donner une correspondance globale pour les
schémas qui ne sont pas affines petits. Pour pallier ce défaut, nous faisceautisons l'algèbre de Higgs-Tate dans le topos de Faltings et l'utilisons alors 
pour construire une correspondance de Simpson $p$-adique de façon parallèle au procédé local. 
Ensuite, nous prouvons par descente cohomologique que la correspondance globale est bien équivalente à la locale pour les schémas affines petits. 

\subsection{}\label{dolb1}
Soit $X$ un $S$-schéma lisse. 
Reprenant les notations de \ref{tordef21}, nous supposons qu'il existe une $\tS$-déformation lisse $\tX$ de $X\otimes_{\co_K}\co_C$ et la fixons dans cette section: 
\begin{equation}\label{dolb1a}
\xymatrix{
{X\otimes_{\co_K}\co_C}\ar[r]\ar[d]\ar@{}[rd]|\Box&{\tX}\ar[d]\\
{\Spec(\co_C)}\ar[r]&{\tS}}
\end{equation}
On notera que la condition est superflue dans le {\em cas relatif}; on pourra en effet prendre $\tX=X\times_S\tS$. 

Nous reprenons également les notations de la section \ref{ft}. Rappelons en particulier que l'on dispose alors de l'anneau $\ocB=\{U\mapsto \ocB_U\}$ 
du topos de Faltings $\tE$ \eqref{ft5}. Pour tout $n\geq 0$, on pose $\ocB_n=\ocB/p^n\ocB$ et pour tout $U\in \ob(\Et_{/X})$, 
$\ocB_{U,n}=\ocB_U/p^n\ocB_U$, qui est un anneau de $U_{\oeta,\fet}$. 

Pour tout  $X$-schéma étale affine petit $U$ \eqref{tordef3}, il existe une suite exacte canonique de $\ocB_{U,n}$-modules de $U_{\oeta,\fet}$ 
\begin{equation}\label{dolb1b}
0\rightarrow \ocB_{U,n}\rightarrow \cF_{U,n}\rightarrow \txi^{-1}\Omega^1_{X/S}(U)\otimes_{\co_X(U)}\ocB_{U,n} \rightarrow 0, 
\end{equation}
telle que pour tout  point géométrique $\oy$ de $U_\oeta$, on ait un isomorphisme canonique de $\ocB_{U,\oy}$-représentations de $\pi_1(U_\oeta,\oy)$
\begin{equation}
(\cF_{U,n})_\oy\stackrel{\sim}{\rightarrow}\cF^\oy_U/p^n\cF^\oy_U,
\end{equation}
où $\cF^\oy_U$ est la $\ocB_{U,\oy}$-extension de Higgs-Tate  \eqref{tordef4c} définie relativement à la restriction de $\tX$ au-dessus de $U$ \eqref{ahttf7}. 
Pour tout nombre rationnel $r\geq 0$, soit
\begin{equation}\label{dolb1c}
0\rightarrow \ocB_{U,n}\rightarrow \cF_{U,n}^{(r)}\rightarrow \txi^{-1}\Omega^1_{X/S}(U)\otimes_{\co_X(U)}\ocB_{U,n} \rightarrow 0
\end{equation}
l'extension de $\ocB_{U,n}$-modules de $U_{\oeta,\fet}$ obtenue à partir de $\cF_{U,n}$ par tirage en arrière
par la multiplication par $p^r$ sur $\txi^{-1}\Omega^1_{X/S}(U)\otimes_{\co_X(U)}\ocB_{U,n}$, et soit 
\begin{equation}\label{dolb1d}
 \cC^{(r)}_{U,n}=\underset{\underset{m\geq 0}{\longrightarrow}}\lim\ \Sym^m_{\ocB_{U,n}}(\cF^{(r)}_{U,n})
\end{equation}
la $\ocB_{U,n}$-algèbre de $U_{\oeta,\fet}$ associée, dans laquelle les morphismes de transition sont localement définis en envoyant  
 $x_1\otimes\dots \otimes x_m$ sur $1\otimes x_1\otimes\dots \otimes x_m$.  

La formation de $\cF^{(r)}_{U,n}$ étant fonctorielle en $U$, les correspondances 
\begin{equation}\label{dolb1e}
\{U\mapsto \cF^{(r)}_{U,n}\} \ \ \ {\rm et}\ \ \ \{U\mapsto \cC^{(r)}_{U,n}\}
\end{equation}
définissent des préfaisceaux sur la sous-catégorie $E^\apt$ de $E$ formée des objets $(V\rightarrow U)$ tels que $U$ soit affine petit \eqref{ft1}. 
Cette dernière est topologiquement génératrice de $E$. 
Par suite, passant aux faisceaux associés, on obtient un $\ocB_n$-module $\cF^{(r)}_n$ et une $\ocB_n$-algèbre $\cC^{(r)}_n$ de $\tE$ \eqref{ahttf37}.

\subsection{}\label{dolb3}
On pose $\cS=\Spf(\co_C)$ et on désigne par $\fX$ le schéma formel complété $p$-adique de $\oX=X\times_S\oS$. 
De même, afin de prendre en compte la topologie $p$-adique, on considère le complété formel $p$-adique du topos annelé $(\tE,\ocB)$. 
On définit d'abord la {\em fibre spéciale} $\tE_s$ de $\tE$, topos s'insérant dans le diagramme commutatif
\begin{equation}\label{dolb3a}
\xymatrix{
{\tE_s}\ar[r]^{\sigma_s}\ar[d]_{\delta}&{X_{s,\et}}\ar[d]^{\iota}\\
{\tE}\ar[r]^\sigma&{X_\et}}
\end{equation}
dans lequel $\iota$ désigne l'injection canonique  \eqref{ahttf12}. 
Concrètement, $\tE_s$ est la  sous-catégorie pleine de $\tE$ des objets $F$ tels que $F|\sigma^*(X_\eta)$ 
soit l'objet final de $\tE_{/\sigma^*(X_\eta)}$, et $\delta_*\colon \tE_s\rightarrow \tE$ le foncteur d'injection canonique.

Pour tout entier $n\geq 0$, $\ocB_n$ est un objet de $\tE_s$. Nous noterons par $\oX_n$ et $\oS_n$ les réductions de $\oX$ et $\oS$ modulo $p^n$. 
Le morphisme $\sigma_s$ est alors sous-jacent à un morphisme  canonique de topos annelés
\begin{equation}\label{dolb3b}
\sigma_n\colon (\tE_s,\ocB_n)\rightarrow (X_{s,\et},\co_{\oX_n}),
\end{equation}
dans lequel nous avons identifié les topos étales de $X_s$ et $\oX_n$, puisque $k$ est algébriquement clos.   

La complétion formelle  $p$-adique de $(\tE,\ocB)$ est le topos annelé $(\tE_s^{\mN^\circ},\bvocB)$, où $\tE_s^{\mN^\circ}$ désigne 
le topos des systèmes projectifs de $\tE_s$ indexés par l'ensemble ordonné $\mN$ (\cite{agt} III.7) et $\bvocB=(\ocB_{n})_{n\geq 0}$. 
Les morphismes $\sigma_n$ induisent un morphisme de topos 
\begin{equation}\label{dolb3e}
\hupsigma \colon (\tE_s^{\mN^\circ},\bvocB)\rightarrow (X_{s,\zar},\co_\fX).
\end{equation}

On travaille dans la catégorie $\bMod_{\mQ}(\bvocB)$ des {\em $\bvocB$-modules à isogénie près}, {\em i.e.}, la catégorie ayant 
pour objets les $\bvocB$-modules  et où, pour tous  $\bvocB$-modules $\cF$ et $\cG$, 
\begin{equation}\label{dolb3f}
\Hom_{\bMod_{\mQ}(\bvocB)}(\cF,\cG)=\Hom_{\bMod(\bvocB)}(\cF,\cG)\otimes_\mZ\mQ.
\end{equation} 
On note le foncteur de localisation $\bMod(\bvocB)\rightarrow \bMod_{\mQ}(\bvocB)$ par $\cF\mapsto \cF_\mQ$. 
On appelle {\em $\bvocB_\mQ$-modules} les  objets de $\bMod_{\mQ}(\bvocB)$. 

\subsection{}
Pour tout nombre rationnel $r\geq 0$ et tout entier $n\geq 0$, on a une suite exacte canonique localement scindée 
\begin{equation}\label{dolb3c}
0\rightarrow \ocB_n\rightarrow \cF^{(r)}_n\rightarrow 
\sigma^*_n(\txi^{-1}\Omega^1_{\oX_n/\oS_n})\rightarrow 0
\end{equation}
et un isomorphisme canonique de $\ocB_n$-algèbres  
\begin{equation}\label{dolb3d}
\cC^{(r)}_n \stackrel{\sim}{\rightarrow}\underset{\underset{m\geq 0}{\longrightarrow}}\lim\ \Sym^m_{\ocB_n}(\cF^{(r)}_n),
\end{equation}
dans lequel les morphismes de transition sont localement définis en envoyant $x_1\otimes\dots \otimes x_m$ sur 
$1\otimes x_1\otimes\dots \otimes x_m$ \eqref{ahttf38}.
On pose $\bvcF^{(r)}=(\cF^{(r)}_n)_{n\geq 0}$ qui est un $\bvocB$-module et $\bvcC^{(r)}=(\cC^{(r)}_n)_{n\geq 0}$ qui est une $\bvocB$-algèbre.
On a une suite exacte canonique de $\bvocB$-modules 
\begin{equation}\label{dolb3g}
0\rightarrow \bvocB\rightarrow \bvcF^{(r)}\rightarrow 
\hupsigma^*(\txi^{-1}\Omega^1_{\fX/\cS})\rightarrow 0.
\end{equation}
La $\bvocB$-dérivation universelle de $\bvcC^{(r)}$ peut être identifiée avec la dérivation
\begin{equation}\label{dolb3h}
d_{\bvcC^{(r)}}\colon \bvcC^{(r)}\rightarrow \hupsigma^*(\txi^{-1}\Omega^1_{\fX/\cS})\otimes_{\bvocB}\bvcC^{(r)}.
\end{equation}
C'est un $\bvocB$-champ de Higgs. On note $\mK^\bullet(\bvcC^{(r)})$ le complexe de Dolbeault du $\bvocB$-module de Higgs $(\bvcC^{(r)},p^rd_{\bvcC^{(r)}})$.

Pour tous nombres rationnels $r\geq r'\geq 0$, on a un homomorphisme canonique de $\bvocB$-algèbres $\bvcC^{(r)}\rightarrow \bvcC^{(r')}$. On observera que la  
restriction de la dérivation $p^{r'}d_{\bvcC^{(r')}}$ est $p^{r}d_{\bvcC^{(r)}}$. On a donc un morphisme de complexes 
\begin{equation}\label{dolb3i}
\mK^\bullet(\bvcC^{(r)})\rightarrow \mK^\bullet(\bvcC^{(r')}).
\end{equation}

\begin{prop}[cf. \ref{ahttf28}, \cite{agt} III.11.18]\label{dolb4}
L'homomorphisme canonique
\begin{equation}\label{dolb4a}
\co_{\fX}[\frac 1 p]\rightarrow \underset{\underset{r\in \mQ_{>0}}{\longrightarrow}}{\lim}\  \hupsigma_*(\bvcC^{(r)})[\frac 1 p]
\end{equation}
est un isomorphisme, et pour tout $q\geq 1$,
\begin{equation}\label{dolb4b}
\underset{\underset{r\in \mQ_{>0}}{\longrightarrow}}{\lim}\ \rR^q\hupsigma_*(\bvcC^{(r)})[\frac 1 p] =0.
\end{equation}
\end{prop}

Ce résultat est une faisceautisation du calcul de la cohomologie galoisienne de l'algèbre de Higgs-Tate sur un schéma affine petit 
(\cite{agt} II.12.5). Le calcul de cette cohomologie galoisienne repose sur le théorème de presque-pureté de Faltings et la faisceautisation  
requiert de prouver une version modulo $p^n$, à défaut borné près (\ref{taht14}, \cite{agt} II.12.7).

\begin{prop}[cf. \ref{ahttf33}, \cite{agt} III.11.24]\label{dolb5}
Le morphisme  canonique de $\bMod_\mQ(\bvocB)$
\begin{equation}\label{dolb5a}
\bvocB_\mQ\rightarrow \underset{\underset{r\in \mQ_{>0}}{\longrightarrow}}{\lim}\ 
\rH^0(\mK^\bullet_\mQ(\bvcC^{(r)}))
\end{equation}
est un isomorphisme, et pour tout $q\geq 1$, 
\begin{equation}\label{dolb5b}
\underset{\underset{r\in \mQ_{>0}}{\longrightarrow}}{\lim}\ \rH^q(\mK^\bullet_\mQ(\bvcC^{(r)}))=0.
\end{equation}
\end{prop}

Ce résultat est une faisceautisation du calcul de la cohomologie de de Rham de l'algèbre de Higgs-Tate sur un schéma affine petit (\ref{taht12}, \cite{agt} II.12.3). 

\subsection{}\label{dolb50}
Les limites inductives filtrantes ne sont a priori pas représentables dans la catégorie $\bMod_\mQ(\bvocB)$. Cependant, on peut 
naturellement plonger cette dernière dans la catégorie abélienne $\bIndMod(\bvocB)$ des ind-$\bvocB$-modules où les limites inductives filtrantes 
sont représentables et qui a de meilleures propriétés (\ref{indsh}, \ref{indmod}, \cite{ks2}). D'une façon analogue, on peut naturellement plonger 
la catégorie des $\co_\fX[\frac 1 p]$-modules cohérents dans la catégorie 
$\bIndMod(\co_\fX)$ des ind-$\co_\fX$-modules \eqref{ahttf50}.
Le morphisme $\hupsigma$ \eqref{dolb3e}  induit deux foncteurs adjoints
\begin{equation}\label{dolb5c}
\xymatrix{
{\bIndMod(\bvocB)}\ar@<1ex>[r]^-(0.5){\rI \hupsigma_*}&{\bIndMod(\co_\fX)}\ar@<1ex>[l]^-(0.5){\rI \hupsigma^*}}
\end{equation}
qui prolongent les foncteurs adjoints $\hupsigma^*$ et $\hupsigma_*$. 

\begin{defi}
On appelle {\em $\co_\fX[\frac 1 p]$-fibré de Higgs à coefficients dans $\txi^{-1}\Omega^1_{\fX/\cS}$} 
tout $\co_\fX[\frac 1 p]$-module localement projectif de type fini $\cN$ \eqref{notconv14}
muni d'un morphisme $\co_\fX$-linéaire $\theta\colon \cN\rightarrow \txi^{-1}\Omega^1_{\fX/\cS}\otimes_{\co_\fX}\cN$ 
tel que $\theta\wedge \theta=0$.
\end{defi}

On observera que les fibres de l'anneau $\co_\fX[\frac 1 p]$ n'étant pas nécessairement des anneaux locaux, 
les $\co_\fX[\frac 1 p]$-modules localement projectifs de type fini ne sont pas nécessairement localement libres. 

\begin{defi}[cf. \ref{indmdlb4}]\label{dolb6}
Soient $\cM$ un ind-$\bvocB$-module, $\cN$ un $\co_\fX[\frac 1 p]$-fibré de Higgs à coefficients dans $\txi^{-1}\Omega^1_{\fX/\cS}$. 
\begin{itemize}
\item[(i)] On dit que $\cM$ et  $\cN$ sont {\em $r$-associés} (pour $r\in \mQ_{>0}$) s'il existe un isomorphisme de ind-$\bvcC^{(r)}$-modules 
\begin{equation}\label{dolb6a}
\cM\otimes_{\bvocB}\bvcC^{(r)}\stackrel{\sim}{\rightarrow}\rI\hupsigma^*(\cN)\otimes_{\bvocB}\bvcC^{(r)},
\end{equation}
compatible avec les $\bvocB$-champs de Higgs totaux à coefficients dans  $\hupsigma^*(\txi^{-1}\Omega^1_{\fX/\cS})$, 
où $\cM$ est muni du champ de Higgs nul et $\bvcC^{(r)}$ du champ de Higgs $p^rd_{\bvcC^{(r)}}$. 
\item[(ii)] On dit que $\cM$ et $\cN$ sont {\em associés} s'ils sont $r$-associés pour un nombre rationnel $r>0$.
\end{itemize}
\end{defi}

En fait, \eqref{dolb6a} est un isomorphisme de {\em ind-$\bvcC^{(r)}$-modules à $p^r$-connexion relativement à l'extension $\bvcC^{(r)}/\bvocB$} 
\eqref{indmdlb1}. En particulier, pour tous nombres rationnels $r\geq r'>0$, si  $\cM$ et $\cN$ sont $r$-associés, ils sont $r'$-associés.

\begin{defi}[cf. \ref{indmdlb5}]\label{dolb7}
\
\begin{itemize}
\item[(i)] On dit qu'un ind-$\bvocB$-module est de {\em Dolbeault} s'il est associé à un  
$\co_\fX[\frac 1 p]$-fibré de Higgs à coefficients dans $\txi^{-1}\Omega^1_{\fX/\cS}$.
\item[(ii)] On dit qu'un $\co_\fX[\frac 1 p]$-fibré de Higgs à coefficients dans $\txi^{-1}\Omega^1_{\fX/\cS}$ est {\em soluble}
s'il est associé à un ind-$\bvocB$-module. 
\end{itemize}
\end{defi}

La propriété pour un ind-$\bvocB$-module d'être de Dolbeault ne dépend pas du choix de la déformation $\tX$ \eqref{dolb1} 
pourvu que l'on reste dans l'un des cadres, absolu ou bien relatif \eqref{tordef21} (voir \ref{pchp4}).
La propriété pour un $\co_\fX[\frac 1 p]$-fibré de Higgs d'être soluble dépend a priori de la déformation $\tX$ (voir toutefois \ref{pchp11}). 

La propriété d'être de Dolbeault s'applique aux $\bvocB_\mQ$-modules \eqref{dolb50}. On dit alors que le fibré de Higgs associé est  {\em rationnellement soluble}.

Dans (\cite{agt} III.12.11), nous considérions seulement des $\bvocB_\mQ$-modules et lorsqu'on les qualifiait de Dolbeault, on demandait en plus qu'ils
soient adiques de type fini. Nous requalifions ces derniers, dans \ref{aspglob1}, de $\bvocB_\mQ$-modules {\em fortement de Dolbeault} et les fibrés de Higgs  
qui leur correspondent de {\em fortement solubles}. La condition de finitude est importante pour rendre compatible les théories globales et locales pour 
les schémas affines petits \eqref{htls8}.

\begin{teo}[cf. \ref{indmdlb20}]\label{dolb8}
Il existe des équivalences explicites de catégories quasi-inverses l'une de l'autre 
\begin{equation}\label{dolb8a}
\xymatrix{
{\bIndMod^\Dolb(\bvocB)}\ar@<1ex>[r]^-(0.5){\cH}&{\bMH^\sol(\co_\fX[\frac 1 p], \txi^{-1}\Omega^1_{\fX/\cS})}
\ar@<1ex>[l]^-(0.5){\cV}}
\end{equation}
entre la catégorie des  ind-$\bvocB$-modules de Dolbeault  et la catégorie des $\co_\fX[\frac 1 p]$-fibrés de Higgs solubles à coefficients 
dans $\txi^{-1}\Omega^1_{\fX/\cS}$. 
\end{teo}

Ces foncteurs sont en effet explicitement définis dans \ref{indmdlb7} de la manière suivante. Soient $\vupsigma_*$ le foncteur composé
\begin{equation}\label{dolb8d}
\xymatrix{\bIndMod(\bvocB)\ar[r]^-(0.5){\rI \hupsigma_*}&{\bIndMod(\co_\fX)}\ar[r]^-(0.5){\kappa_{\co_\fX}}&{\bMod(\co_\fX)}},
\end{equation}
où $\rI \hupsigma_*$ est défini dans \eqref{dolb5c} et 
\begin{equation}\label{dolb8c}
\kappa_{\co_\fX}(\indcolim\alpha)= \underset{\longrightarrow}{\lim}\ \alpha.
\end{equation}
Alors, le foncteur $\cH$ est défini, pour tout ind-$\bvocB$-module $\cM$, par
\begin{equation}\label{dolb8e}
\cH(\cM)=\underset{\underset{r\in \mQ_{>0}}{\longrightarrow}}{\lim}\ \vupsigma_*(\cM\otimes_{\bvocB}\bvcC^{(r)},p^r\id\otimes d_{\bvcC^{(r)}}).
\end{equation} 
On a un procédé similaire de définition pour $\cV$. Les foncteurs $\cH$ et $\cV$ dépendent a priori de la déformation $\tX$. 

\begin{rema}
Faltings \cite{faltings3} a esquissé les grandes lignes d'une stratégie qui vise à étendre la correspondence \eqref{dolb8a}, 
au dessus d'une courbe propre et semi-stable sur $S$ qui est lisse (et stricte) sur $\eta$,  
à tous les $\bvocB_\mQ$-modules et tous les $\co_\fX[\frac 1 p]$-fibrés de Higgs. Pour ce faire, il procède par descente: 
toute représentation généralisée (resp. fibré de Higgs) devient de Dolbeault (resp. soluble), ou ce qui revient au même petite (resp. petit) 
sur un revêtement étale fini de la fibre géométrique générique de la courbe. 
La descente de la correspondance \eqref{dolb8a} définie sur le revêtement en une correspondance définie sur la courbe initiale
n'est pas immédiate car elle dépend du choix d'une déformation \eqref{dolb1a} et le revêtement étale fini ne se prolonge pas nécessairement aux déformations données des modèles entiers. 
Faltings montre qu'après torsion de l'image inverse du fibré de Higgs par l'obstruction à l'existence d'un tel prolongement,
la correspondance se descend bien. Cette image inverse tordue est étroitement lié à la {\em  fibration de Hitchin} \cite{hitchin}. 
Il y a eu récemment des développements intéressants concernant cet aspect de la théorie, dus d'une part à Heuer \cite{heuer} et d'autre part à Xu \cite{xu}. 
Nous y reviendrons dans le cadre de notre approche dans un prochain travail en commun avec T. Tsuji. 
\end{rema}

\begin{teo}[cf. \ref{indmdlb23}]\label{dolb9}
Pour tout ind-$\bvocB$-module $\cM$ de Dolbeault et tout entier $q\geq 0$, on a un isomorphisme fonctoriel canonique de $\bD^+(\bMod(\co_\fX))$
\begin{equation}\label{dolb9a}
\rR\vupsigma_*(\cM)\stackrel{\sim}{\rightarrow}\mK^\bullet(\cH(\cM)),
\end{equation}
où $\mK^\bullet(\cH(\cM))$ est le complexe de Dolbeault de $\cH(\cM)$.
\end{teo}

C'est un analogue global de \ref{htls30}.

\subsection{}\label{dolb90}
Le morphisme $\psi$ \eqref{ft3a} induit un morphisme de topos
\begin{equation}
\bvpsi\colon X_{\oeta,\et}^{\mN^\circ}\rightarrow \tE^{\mN^\circ},
\end{equation}
où $X_{\oeta,\et}^{\mN^\circ}$ est le topos des systèmes projectifs d'objets de $X_{\oeta,\et}$, indexés par l'ensemble ordonné $\mN$. 
On note $\bvmZ_p$ l'anneau $(\mZ/p^n\mZ)_{n\geq 0}$ de $X_{\oeta,\et}^{\mN^\circ}$. 

On dit qu'un $\bvmZ_p$-module $M=(M_n)_{n\in \mN}$ de $X_{\oeta,\et}^{\mN^\circ}$ est un {\em système local} si les deux conditions 
suivantes sont satisfaites:
\begin{itemize}
\item[(a)] $M$ est $p$-adique, {\em i.e.}, pour tous entiers $n\geq m\geq 0$, le morphisme $M_n/p^mM_n\rightarrow M_m$ déduit du 
morphisme de transition $M_n\rightarrow M_m$ est un isomorphisme; 
\item[(b)]  pour tout entier $n\geq 0$, le $\mZ/p^n\mZ$-module $M_n$ de $X_{\oeta,\et}$ est localement constant constructible. 
\end{itemize}

\begin{cor}[cf. \ref{sld4}]\label{dolb10}
Soient $M=(M_n)_{n\geq 0}$ un $\bvmZ_p$-système local de $X_{\oeta,\et}^{\mN^\circ}$ et $\cM=\bvpsi_*(M)\otimes_{\bvmZ_p}\bvocB$. 
Supposons que $X$ soit propre sur $S$ et que le $\bvocB_\mQ$-module $\cM_\mQ$ soit de Dolbeault.
Il existe alors une suite spectrale canonique 
\begin{equation}\label{dolb10a}
\rE_2^{i,j}=\rH^i(X_s,\rH^j(\mK^\bullet))\Rightarrow \rH^{i+j}(X_{\oeta,\et}^{\mN^\circ},M)\otimes_{\mZ_p}C,
\end{equation}
dans laquelle $\mK^\bullet$ désigne le complexe de Dolbeault de $\cH(\cM_\mQ)$. 
\end{cor}

Cela résulte de \ref{ft6} et de \ref{dolb9}.

\begin{rema}
Dans \ref{dolb10}, si l'on prend $M=\bvmZ_p$, alors $\cM=\bvocB$, le $\bvocB_\mQ$-module $\bvocB_\mQ$ est de Dolbeault et 
$\cH(\bvocB_\mQ)$ est égal à $\co_\fX[\frac 1 p]$ muni du champ de Higgs nul \eqref{aspglob12}. 
La suite spectrale \eqref{dolb10a} n'est autre que la suite spectrale de Hodge-Tate  (\cite{ag} 6.4.6). 
On notera que la construction \ref{dolb10} de cette suite spectrale montre directement qu'elle dégénère en $\rE_2$ et que 
la filtration aboutissement est scindée sans utiliser le théorème de Tate sur la cohomologie galoisienne de $C(j)$. 
Cette construction s'applique en particulier en prenant pour $\tX$ dans le cas relatif \eqref{tordef21} la déformation triviale \eqref{dolb1}. 
\end{rema}

\begin{defi}[cf. \ref{mht1}]\label{htls5}
Nous appelerons {\em $\bvocB_\mQ$-module de Hodge-Tate} tout $\bvocB_\mQ$-module de Dolbeault $\cM$ \eqref{dolb7} 
dont le $\co_\fX[\frac 1 p]$-fibré de Higgs associé $(\cH(\cM),\theta)$ \eqref{dolb8a} est nilpotent, {\em i.e.}, il existe une 
filtration finie décroissante $(\cH_i(\cM))_{0\leq i\leq n}$ de $\cH(\cM)$ par des sous-$\co_\fX[\frac 1 p]$-modules 
cohérents telle que $\cH_0(\cM)=\cH(\cM)$, $\cH_n(\cM)=0$ et telle que pour tout $0\leq i\leq n-1$, on a
\begin{equation}\label{htls5a}
\theta(\cH_i(\cM))\subset \txi^{-1}\Omega^1_{\fX/\cS}\otimes_{\co_\fX}\cH_{i+1}(\cM).
\end{equation}
\end{defi} 

Cette notion ne dépend pas du choix d'une $\tS$-déformation $\tX$, ni même du cadre absolu ou relatif dans lequel on se place \eqref{tordef21} (voir \ref{mht7}). 

\begin{prop}[cf. \ref{mht2}]\label{htls6}
Les foncteurs $\cH$ et $\cV$ \eqref{dolb8a} induisent des équivalences de catégories quasi-inverses
l'une de l'autre  
\begin{equation}\label{htls6a}
\xymatrix{
{\bMod^\HT_\mQ(\bvocB)}\ar@<1ex>[r]^-(0.5){\cH}&{\bMH^\qsolnilp(\co_\fX[\frac 1 p], \txi^{-1}\Omega^1_{\fX/\cS})}
\ar@<1ex>[l]^-(0.5){\cV}}
\end{equation} 
entre la catégorie des $\bvocB_\mQ$-modules de Hodge-Tate  et la catégorie des $\co_\fX[\frac 1 p]$-fibrés 
de Higgs à coefficients dans $\txi^{-1}\Omega^1_{\fX/\cS}$ nilpotents et rationnellement solubles \eqref{dolb7}. 
\end{prop}

\subsection{}\label{htls7}
On suppose dans la fin de cette section que $X$ est affine petit  \eqref{tordef3}, que $X_s\not=\emptyset$, et par
simplicité, que $X_\oeta$ est connexe. On fixe un point géométrique $\oy$ de $X_\oeta$. On pose 
$\Delta=\pi_1(X_\oeta,\oy)$,  on note $\bB_\Delta$ le topos classifiant de $\Delta$ et  
\begin{equation}\label{htls7a}
\nu\colon X_{\oeta,\fet} \stackrel{\sim}{\rightarrow} \bB_\Delta
\end{equation}
le foncteur fibre de $X_{\oeta,\fet}$ en $\oy$ (\cite{agt}  (VI.9.8.4)). On note $\upbeta$ le foncteur composé
\begin{equation}\label{htls7b}
\upbeta \colon \tE\rightarrow \bB_\Delta, \ \ \ F\mapsto \nu\circ (\beta_*(F)), 
\end{equation}
où $\beta$ est le morphisme de topos \eqref{ft3c}. 
On a ainsi défini un foncteur de la catégorie des faisceaux abéliens de $\tE$ dans la catégorie des $\mZ[\Delta]$-modules. 
Celui-ci étant exact à gauche, on note $\rR^q \upbeta$ $(q\geq 0)$ ses foncteurs dérivés droits. 
Pour tout  faisceau abélien $F$ de $\tE$ et tout entier $q\geq 0$, on a un isomorphisme canonique fonctoriel
\begin{equation}
\rR^q\upbeta(F)\stackrel{\sim}{\rightarrow} \nu\circ (\rR^q\beta_*(F)).
\end{equation}

Pour tout  faisceau abélien $F=(F_n)_{n\geq 0}$ de $\tE^{\mN^\circ}$, on pose
\begin{equation}\label{htls7c}
\hupbeta(F)=\underset{\underset{n\geq 0}{\longleftarrow}}{\lim}\ \upbeta(F_n).
\end{equation}
On définit ainsi un foncteur de la catégorie des faisceaux abéliens de $\tE^{\mN^\circ}$ dans la catégorie des $\mZ[\Delta]$-modules.
Celui-ci étant exact à gauche, on note abusivement $\rR^q\hupbeta(F)$ ($q\geq 0$) ses foncteurs dérivés droits. 
Grâce à (\cite{jannsen} 1.6), on a une suite exacte canonique
\begin{equation}\label{htls7d}
0\rightarrow \rR^1 \underset{\underset{n\geq 0}{\longleftarrow}}{\lim}\ \rR^{q-1}\upbeta(F_n)\rightarrow
\rR^q\hupbeta(F)\rightarrow 
\underset{\underset{n\geq 0}{\longleftarrow}}{\lim}\ \rR^q\upbeta(F_n)\rightarrow 0,
\end{equation}
où l'on pose $\rR^{-1}\upbeta(F_n)=0$ pour tout  $n\geq 0$. 

On pose $\oR=\nu(\ocB_X)$, qui n'est rien d'autre que  l'algèbre définie dans \eqref{ht5c} munie de l'action canonique de $\Delta$. 
Pour tout entier $q\geq 0$, $\rR^q\hupbeta$ induit un foncteur qu'on note également
\begin{equation}\label{htls7e}
\rR^q\hupbeta\colon \bMod(\bvocB)\rightarrow \bRep_{\hoR}(\Delta). 
\end{equation}
Ce dernier induit un foncteur que l'on note aussi 
\begin{equation}\label{htls7f}
\rR^q\hupbeta\colon \bMod_\mQ(\bvocB)\rightarrow \bRep_{\hoR[\frac 1 p]}(\Delta). 
\end{equation}

\begin{teo}[cf. \ref{mdpsa32}]\label{htls8}
On conserve les hypothèses et notations de \ref{htls7}. Alors,
\begin{itemize}
\item[{\rm (i)}] Le foncteur $\hupbeta$ \eqref{htls7f} induit une équivalence de catégories
\begin{equation}\label{htls8a}
\bMod^\fDolb_\mQ(\bvocB)\stackrel{\sim}{\rightarrow} \bRep_{\hoR[\frac 1 p]}^{\Dolb}(\Delta),
\end{equation} 
entre la catégorie des  $\bvocB_\mQ$-modules fortement de Dolbeault \eqref{dolb7} et la catégorie des 
$\hoR[\frac 1 p]$-représentations de $\Delta$ de Dolbeault \eqref{tordef43}.
\item[{\rm (ii)}] Pour tout $\bvocB_\mQ$-module fortement de Dolbeault $\cM$ et tout entier $q\geq 1$, on a
\begin{equation}\label{htls8b}
\rR^q\hupbeta(\cM)=0.
\end{equation}
\end{itemize}
\end{teo}

Cela résulte d'un résultat de descente cohomologique pour les $\bvocB_\mQ$-modules fortement de Dolbeault \ref{mdpsa18} qui,
{\it{in fine}} se réduit à un résultat de descente cohomologique pour l'anneau $\bvocB_\mQ$  (\cite{ag} 4.6.30). 

Nous prouvons que, via l'équivalence de catégories \eqref{htls8a}, les foncteurs $\cH$ \eqref{dolb8a} et $\mH$ \eqref{tordef41b}
(resp. $\cV$ \eqref{dolb8a} et $\mV$ \eqref{tordef42b}) se correspondent \eqref{mdpsa33}.

\section[Fonctorialité par image directe propre]{\texorpdfstring{Fonctorialité de la correspondance de Simpson $p$-adique par image directe propre}
{Fonctorialité de la correspondance de Simpson p-adique par image directe propre}}\label{fpscpdi}

\subsection{}\label{fpscpdi1} 
Soit $g\colon X'\rightarrow X$ un  morphisme lisse de $S$-schémas lisses. On affecte un prime $^\prime$ aux objets associés à $X'/S$. 
Par fonctorialité du topos de Faltings, $g$ induit un morphisme  canonique $\Theta$ entre les topos de Faltings correspondants qui s'insère dans 
un diagramme commutatif
\begin{equation}\label{fpscpdi1a}
\xymatrix{
{X'_{\oeta,\et}}\ar[d]_{g_\oeta}\ar[r]^{\psi'}&{\tE'}\ar[d]^{\Theta}\ar[r]^{\sigma'}&{X'_\et}\ar[d]^g\\
{X_{\oeta,\et}}\ar[r]^\psi&{\tE}\ar[r]^{\sigma}&{X_\et}}
\end{equation}
On dispose également d'un homomorphisme canonique d'anneaux
\begin{equation}\label{fpscpdi1b}
\ocB\rightarrow \Theta_*(\ocB').
\end{equation}

\begin{teo}[\cite{faltings2}, \cite{ag} 5.7.4] \label{fpscpdi2} 
Supposons que $g\colon X'\rightarrow X$ soit propre, et soit  $F'$ un faisceau localement constant constructible de 
$(\mZ/p^n\mZ)$-modules de $X'_{\oeta,\et}$ $(n\geq 1)$. Alors, pour tout entier $i\geq 0$, le morphisme  canonique 
\begin{equation}\label{fpscpdi2a} 
\psi_*(\rR^ig_{\oeta*}(F'))\otimes_{\mZ_p}\ocB\rightarrow \rR^i\Theta_*(\psi'_*(F')\otimes_{\mZ_p}\ocB')
\end{equation}
est un presque-isomorphisme. 
\end{teo}

On notera que les faisceaux $\rR^ig_{\oeta*}(F)$ sont localement constants constructibles sur $X_\oeta$ 
par les théorèmes de changement de base, propre et lisse. 

Faltings  a énoncé une {\em version relative} de son théorème principal de comparaison $p$-adique  
dans \cite{faltings2} et en a très grossièrement esquissé une preuve dans l'appendice. 
Certains arguments doivent en fait être modifiés et la preuve dans \cite{ag} requiert beaucoup plus de travail.

Dans l'énoncé (\cite{ag} 5.7.4), on requiert en fait que $g$ soit projectif. Toutefois, comme signalé dans  (\cite{ag} 5.7.6), le résultat vaut sous l'hypothèse 
plus générale que $g$ soit propre par la même preuve, en remplaçant le résultat de presque-finitude (\cite{ag} 2.8.18) par la généralisation récente de He (\cite{he3} 1.5).

\subsection{}\label{fpscpdi3}
Conservant les notations de \ref{tordef21}, nous supposons qu'il existe un diagramme commutatif à carrés cartésiens 
\begin{equation}\label{fpscpdi3b}
\xymatrix{
{X'\otimes_{\co_K}\co_C}\ar[d]_{g\otimes\id}\ar@{}[rd]|\Box\ar[r]&{\tX'}\ar[d]^{\tg}\\
{X\otimes_{\co_K}\co_C}\ar[r]\ar[d]\ar@{}[rd]|\Box&{\tX}\ar[d]\\
{\Spec(\co_C)}\ar[r]&{\tS}}
\end{equation}
où $\tX$ et $\tX'$ sont des $\tS$-schémas lisses, {\em que nous fixons dans cette section}. 
On notera que la condition est superflue dans le {\em cas relatif}; on pourra en effet prendre $\tg=g\times_S\tS$. 
On observera que $\tg$ est lisse en vertu de (\cite{ega4} 17.11.1).

Le diagramme \eqref{fpscpdi1a} induit un diagramme commutatif de morphismes de topos annelés
\begin{equation}\label{fpscpdi3a}
\xymatrix{
{(\tE'^{\mN^\circ}_s,\bvocB')}\ar[d]_{\hupsigma'}\ar[r]^{\bvuptheta}&{(\tE^{\mN^\circ}_s,\bvocB)}\ar[d]^-(0.5){\hupsigma}\\
{(X'_{s,\zar},\co_{\fX'})}\ar[r]^{\fgg}&{(X_{s,\zar},\co_{\fX})}}
\end{equation}
dans lequel les flèches horizontales sont induites par $\Theta$ et $g$, et les flèches verticales sont induites par $\sigma'$ et $\sigma$ \eqref{dolb3e}. 
Nous prouvons dans \ref{ahttfg20} que pour tout nombre rationnel $r\geq 0$, le relèvement $\tg$ \eqref{fpscpdi3b} induit un homomorphisme 
des algèbres de Higgs-Tate 
\begin{equation}\label{fpscpdi3c}
\bvuptheta^*(\bvcC^{(r)})\rightarrow \bvcC'^{(r)},
\end{equation}
dont la construction est assez délicate.

\begin{teo}[cf. \ref{crindmd21}]\label{fpscpdi4}
Supposons $g\colon X'\rightarrow X$ propre. 
Soient $\cM$ un ind-$\bvocB'$-module de Dolbeault \eqref{dolb7},  
\begin{equation}\label{fpscpdi4a}
\cH'(\cM)\rightarrow \txi^{-1}\Omega^1_{\fX'/\cS}\otimes_{\co_{\fX'}}\cH'(\cM)
\end{equation}
le fibré de Higgs associé \eqref{dolb8a}, 
\begin{equation}\label{fpscpdi4b}
\ucH'(\cM)\rightarrow \txi^{-1}\Omega^1_{\fX'/\fX}\otimes_{\co_{\fX'}}\ucH'(\cM)
\end{equation}
le fibré de Higgs relatif déduit de \eqref{fpscpdi4a}, $\umK^\bullet$ le complexe de Dolbeault de $\ucH'(\cM)$. 
Alors, pour tout entier $q\geq 0$, il existe un nombre rationnel $r>0$ et un $\bvcC^{(r)}$-isomorphisme 
\begin{equation}\label{fpscpdi4c}
\rR^q\rI\bvuptheta_*(\cM)\otimes_{\bvocB}\bvcC^{(r)}\stackrel{\sim}{\rightarrow}
\rI\hupsigma^*(\rR^q\fgg_*(\umK^\bullet))\otimes_{\bvocB}\bvcC^{(r)},
\end{equation}
compatible avec les champs de Higgs totaux, où $\bvcC^{(r)}$ est muni du champ de Higgs $p^rd_{\bvcC^{(r)}}$, $\rR^q\rI\bvuptheta_*(\cM)$ 
du champ de Higgs nul et $\rR^q\fgg_*(\umK^\bullet)$ du champ de Katz-Oda \eqref{MH96}.
\end{teo}

On observera que le $\co_\fX[\frac 1 p]$-module $\rR^q\fgg_*(\umK^\bullet)$ est cohérent, et que le foncteur 
\begin{equation}\label{fpscpdi4d}
\rI\hupsigma^*\colon \bMod^\coh(\co_{\fX}[\frac 1 p])\rightarrow \bIndMod(\bvocB)
\end{equation}
est exact \eqref{chb50}.

\begin{cor}[cf. \ref{crindmd26}]\label{fpscpdi5}
Sous les hypothèses de \ref{fpscpdi4}, si  le $\co_{\fX}[\frac 1 p]$-module $\rR^q\fgg_*(\umK^\bullet)$ est localement projectif de type fini, alors le 
ind-$\bvocB$-module $\rR^q\rI\bvuptheta_*(\cM)$ est de {\em Dolbeault}, et l'on a un isomorphisme
\begin{equation}\label{fpscpdi5a}
\cH(\rR^q\rI\bvuptheta_*(\cM))\stackrel{\sim}{\rightarrow} \rR^q\fgg_*(\umK^\bullet),
\end{equation}
où $\rR^q\fgg_*(\umK^\bullet)$ est muni du champ de Katz-Oda.
\end{cor}

\begin{cor}[cf. \ref{crindmd29}]\label{fpscpdi8}
Supposons $g\colon X'\rightarrow X$ propre. Soit $\cM^n=\bvpsi_*(\rR^n\bvg_{\oeta*}(\bvmZ_p))\otimes_{\bvmZ_p}\bvocB$ 
pour un entier $n\geq 0$.  
Alors, le $\bvocB_\mQ$-module $\cM^n_\mQ$ est de Hodge-Tate \eqref{htls5} et l'on a un isomorphisme
\begin{equation}\label{fpscpdi8a}
\cH(\cM^n_\mQ)\stackrel{\sim}{\rightarrow} \oplus_{0\leq i\leq n}\rR^ig_*(\txi^{i-n}\Omega^{n-i}_{X'/X}) \otimes_{\co_X}\co_\fX[\frac 1 p],
\end{equation}
où le champ de Higgs  à droite est induit par les applications de Kodaira-Spencer associées à $g$ 
\begin{equation}\label{fpscpdi8b}
\rR g_*(\txi^{-j}\Omega^j_{X'/X})\rightarrow \txi^{-1}\Omega^1_{X/S}\otimes_{\co_X}\rR g_*(\txi^{1-j}\Omega^{j-1}_{X'/X})[+1].
\end{equation}
\end{cor}

Cela résulte de \ref{fpscpdi2} et \ref{fpscpdi5}. En effet, $\bvpsi'_*(\bvmZ_p)=\bvmZ_p$, le $\bvocB'_\mQ$-module $\bvocB'_\mQ$ est de Dolbeault et 
$\ucH'(\bvocB'_\mQ)$ s'identifie à $\co_{\fX'}[\frac 1 p]$ muni du champ de Higgs nul \eqref{aspglob12}. 
Par conséquent, avec les notations de \ref{fpscpdi4}, pour tout $q\geq 0$, le $\co_{\fX}[\frac 1 p]$-module 
$\rR^q\fgg_*(\umK^\bullet)$ est localement libre de type fini d'après (\cite{deligne1} 5.5), ce qui complète la preuve de la première assertion. 
La seconde assertion suit aisément compte tenu de la définition du champ de Katz-Oda (\ref{MH110}, \cite{katz2} 1.2).

\begin{cor}[cf. \ref{crindmd24}]\label{fpscpdi6}
Soit $M=(M_n)_{n\geq 0}$ un $\bvmZ_p$-système local de $X'^{\mN^\circ}_{\oeta,\et}$ \eqref{dolb90}.  
On pose $\cM=\bvpsi'_*(M)\otimes_{\bvmZ_p}\bvocB'$, 
on suppose que $g\colon X'\rightarrow X$ est propre et que le $\bvocB'_\mQ$-module $\cM_\mQ$ est de Dolbeault.
Il existe alors un nombre rationnel $r>0$ et une suite spectrale
\begin{equation}\label{fpscpdi6a}
\rE_2^{i,j}=\hupsigma^*_\mQ(\rR^i\fgg_*(\rH^j(\umK^\bullet)))\otimes_{\bvocB_\mQ}\bvcC^{(r)}_\mQ
\Rightarrow \bvpsi_*(\rR^{i+j}\bvg_{\oeta*}(M))\otimes_{\bvmZ_p}\bvcC^{(r)}_\mQ,
\end{equation}
où $\umK^\bullet$ est le complexe de Dolbeault du $\co_{\fX'}[\frac 1 p]$-fibré de Higgs relatif $\ucH'(\cM_\mQ)$ \eqref{fpscpdi4b}. 
\end{cor}
 
Cela résulte de \ref{fpscpdi2} et \ref{fpscpdi4}.

\begin{cor}[cf. \ref{crindmd240}]\label{fpscpdi60}
Supposons $g\colon X'\rightarrow X$ propre.
Il existe alors un nombre rationnel $r>0$ et, pour tout entier $n\geq 0$, un isomorphisme canonique de $\bvcC^{(r)}_\mQ$-modules 
\begin{equation}\label{fpscpdi60a}
\bvpsi_*(\rR^n\bvg_{\oeta*}(\bvmZ_p))\otimes_{\bvmZ_p}\bvcC^{(r)}_\mQ
\stackrel{\sim}{\rightarrow}\oplus_{0\leq i\leq n}\sigma^*(\rR^ig_*(\Omega^{n-i}_{X'/X}))\otimes_{\sigma^*(\co_X)}\bvcC^{(r)}_\mQ(i-n).
\end{equation}
\end{cor}

Cela résulte de \ref{fpscpdi4} appliqué à $\cM=\bvocB'_\mQ$ puisque le $\bvocB'_\mQ$-module $\bvocB'_\mQ$ est de Dolbeault et que
$\ucH'(\bvocB'_\mQ)$ est le fibré trivial $\co_{\fX'}[\frac 1 p]$ muni du champ de Higgs nul.  

Ce corollaire s'applique en particulier en prenant pour $\tg$ dans le cas relatif \eqref{tordef21} la déformation triviale \eqref{fpscpdi3}. 

\begin{teo}[\cite{ag} 6.7.5]\label{fpscpdi7}
Supposons $g\colon X'\rightarrow X$ propre. On a alors une suite spectrale canonique de $\bvocB_\mQ$-modules, {\em la suite spectrale de 
Hodge-Tate relative},
\begin{equation}\label{fpscpdi7a}
\rE_2^{i,j}=\sigma^*(\rR^ig_*(\Omega^j_{X'/X}))\otimes_{\sigma^*(\co_X)}\bvocB_\mQ(-j)\Rightarrow \bvpsi_*(\rR^{i+j}\bvg_{\oeta*}(\bvmZ_p))\otimes_{\bvmZ_p}\bvocB_\mQ.
\end{equation}
\end{teo}

Cette suite spectrale ne requiert  la considération d'aucune déformation \eqref{fpscpdi3b}. 
Elle est $G_K$-équivariante pour les $G_K$-structures équivariantes sur les différents topos et objets considérés. 
Elle dégénère donc en $\rE_2$. Néanmoins, sa filtration aboutissement ne se scinde pas en général. Cependant, on peut vérifier qu'elle 
se scinde après changement de base de $\bvocB$ à $\bvcC^{(r)}$ pour un nombre rationnel $r>0$ 
et que cela correspond à la décomposition \eqref{fpscpdi60a}. 

\begin{rema}
Le théorème \ref{fpscpdi4} et son corollaire \ref{fpscpdi5} fournissent des analogues $p$-adiques de résultats de \cite{biquard} pour la correspondence de Simpson complexe, 
mais nous nous autorisons des schémas logarithmiques plus généraux. Les analogues pour la correspondence de Simpson en caractéristique $p$ 
sont établis dans \cite{ov} pour les morphismes propres et lisses. 
Dans le cadre complexe, on trouve dans \cite{dps} des énoncés qui traitent de
singularités plus générales que les singularités logarithmiques. Nous espérons établir des analogues $p$-adiques généralisant notre travail.   
\end{rema}

\subsection{}\label{fpscpdi9} 
La preuve de \ref{fpscpdi4} peut se diviser en trois étapes. Premièrement, on calcule la cohomologie galoisienne relative et la cohomologie de 
Higgs de l'algèbre de Higgs-Tate en adaptant le calcul de Faltings dans le cas absolu. Deuxièmement, on faisceautise ces calculs en considérant 
le produit fibré de topos  
\begin{equation}\label{fpscpdi9a} 
\xymatrix{
\tE'\ar[rd]^{\sigma'}\ar[d]_{\tau}&\\
{\tE\times_{X_\et}X'_\et}\ar[r]\ar[d]\ar@{}|\Box[rd]&{X'_\et}\ar[d]^g\\
{\tE}\ar[r]^-(0.5){\sigma}&{X_\et}}
\end{equation}
Le calcul, par nature local,  de la cohomologie galoisienne relative peut alors être globalisé en le calcul des faisceaux $\rR^i\tau_*(\cC'^{(r)}_n)$. 
La dernière étape est un théorème de changement de base relativement au carré cartésien ci-dessus. 

Il s'avère qu'il existe un site très naturel sous-jacent au topos $\tE\times_{X_\et}X'_\et$ qui est une variante relative du topos de Faltings.  
Sa définition nous a été inspirée par celle des produits orientés de topos (au-delà même du topos co-évanescent  qui a inspiré la définition 
du topos de Faltings). 

\section{Topos de Faltings relatif}\label{rft} 

\subsection{}\label{rft1} 
Soit $g\colon X'\rightarrow X$ un morphisme de $S$-schémas. 
On note $G$ la catégorie des morphismes $(W\rightarrow U\leftarrow V)$
au-dessus du morphisme  canonique  $X'\rightarrow X\leftarrow X_\oeta$, {\em i.e.}, les diagrammes commutatifs
\begin{equation}\label{rft1a}
\xymatrix{W\ar[r]\ar[d]&U\ar[d]&V\ar[l]\ar[d]\\
X'\ar[r]&X&X_\oeta\ar[l]}
\end{equation}
tels que $W$ soit étale sur $X'$, $U$  étale sur $X$ et  que le morphisme  canonique $V\rightarrow U_\oeta$ soit {\em fini  étale}. 
On la munit de la topologie engendrée par les recouvrements 
\[
\{(W_i\rightarrow U_i\leftarrow V_i)\rightarrow (W\rightarrow U \leftarrow V)\}_{i\in I}
\]
des trois types suivants~:
\begin{itemize}
\item[(a)] $U_i=U$, $V_i=V$ pour tout  $i\in I$ et $(W_i\rightarrow W)_{i\in I}$ est recouvrement;
\item[(b)] $W_i=W$, $U_i=U$ pour tout  $i\in I$ et $(V_i\rightarrow V)_{i\in I}$ est un recouvrement;
\item[(c)]  diagrammes
\begin{equation}\label{rft1b}
\xymatrix{
W'\ar[r]\ar@{=}[d]&U'\ar[d]\ar@{}|\Box[rd]&V'\ar[l]\ar[d]\\
W\ar[r]&U&V\ar[l]}
\end{equation}
où $U'\rightarrow U$ est un morphisme quelconque et où le carré droit est cartésien. 
\end{itemize}

On note $\tG$ le topos des faisceaux d'ensembles sur $G$ et on l'appelle {\em topos de Faltings relatif} associé au couple de morphismes $(X_\oeta\rightarrow X,X'\rightarrow X)$. 

On a un morphisme canonique de topos 
\begin{equation}\label{rft1c}
\pi\colon \tG\rightarrow X'_\et, \ \ \ W\in \ob(\Et_{/X'}) \mapsto \pi^*(W)=(W\rightarrow X\leftarrow X_\oeta)^a.
\end{equation}

\subsection{}\label{rft2} 
Si $X'=X$, le topos $\tG$ est canoniquement équivalent au topos de Faltings $\tE$. 
Par fonctorialité du topos de Faltings relatif, on obtient donc une factorisation naturelle du morphisme  canonique $\Theta\colon \tE'\rightarrow \tE$ 
s'insérant dans un diagramme commutatif
\begin{equation}\label{rft2d}
\xymatrix{
\tE'\ar[d]^{\tau}\ar[rd]^{\sigma'} \ar@/_1pc/[dd]_{\Theta}&\\
\tG\ar[d]^-(0.5){\upgamma}\ar[r]^-(0.4){\pi}\ar@{}|\Box[rd]&X'_\et\ar[d]^{g}\\
\tE\ar[r]^{\sigma}&X_\et}
\end{equation}
{\em Nous prouvons que le carré inférieur est cartésien.} 

Nous prouvons un théorème de changement de base relativement à ce carré pour tous les faisceaux
abéliens de torsion de $X'_\et$, inspiré par un théorème de changement de base pour les produits orientés dû à Gabber 
(\cite{ag} 6.5.5). Il se ramène au théorème de changement de base propre pour le topos étale. Nous prouvons 
ensuite le résultat suivant qui joue un rôle crucial à la fois dans les preuves de \ref{fpscpdi4} et de \ref{fpscpdi7}: 

\begin{teo}[\cite{ag} 6.5.31]\label{rft3} 
Soit $g\colon X'\rightarrow X$ un morphisme {\em propre} et lisse de $S$-schémas lisses. 
Il existe alors un entier $N\geq 0$ tel que pour tous entiers  $n\geq 1$ et $q\geq 0$, 
et tout $\co_{X'_n}$-module quasi-coherent $\cF$, le noyau et le conoyau du morphisme de changement de base 
\begin{equation}\label{rft3a} 
\sigma^*(\rR^q g_*(\cF))\rightarrow \rR^q\upgamma_*(\pi^*(\cF)),
\end{equation}
soient annulés par $p^N$.
\end{teo}
Dans cet énoncé, $\pi^*$ et $\sigma^*$ désignent les images inverses pour les morphismes de topos annelés 
\begin{eqnarray}
\pi\colon (\tG,\tau_*(\ocB'))&\rightarrow &(X'_\et,\co_{X'}),\label{rft3b}\\
\sigma\colon (\tE,\ocB)&\rightarrow &(X_\et,\co_{X}).\label{rft3c} 
\end{eqnarray}

\chapter{Préliminaires}

\section{Notations et conventions}

{\em Tous les anneaux considérés dans cet article possèdent un élément unité~;
les homomorphismes d'anneaux sont toujours supposés transformer l'élément unité en l'élément unité.
Nous considérons surtout des anneaux commutatifs, et lorsque nous parlons d'anneau 
sans préciser, il est sous-entendu qu'il s'agit d'un anneau commutatif~; en particulier, 
il est sous-entendu, lorsque nous parlons d'un topos annelé $(X,A)$ sans préciser, que $A$ est commutatif.

On sous-entend par {\em monoïde}  un monoïde commutatif et unitaire. 
Les homomorphismes de monoïdes sont toujours supposés transformer l'élément unité en l'élément unité.}

\subsection{}\label{notconv15}
Soit $p$ un nombre premier. On munit  $\mZ_p$ de la topologie $p$-adique, ainsi que toutes les $\mZ_p$-algèbres adiques 
({\em i.e.}, les $\mZ_p$-algèbres complètes et séparées pour la topologie $p$-adique).
Soient $A$ une $\mZ_p$-algèbre adique, $i\colon A\rightarrow A[\frac 1 p]$ l'homomorphisme canonique. 
On appelle {\em topologie $p$-adique} sur $A[\frac 1 p]$ l'unique topologie 
compatible avec sa structure de groupe additif  pour laquelle les sous-groupes $i(p^nA)$,
pour $n\in \mN$, forment un système fondamental de voisinages de $0$ (\cite{tg} chap.~III §1.2, prop.~1). 
Elle fait de $A[\frac 1 p]$ un anneau topologique. 
Soient $M$ un $A[\frac 1 p]$-module de type fini, $M^\circ$ un sous-$A$-module de type fini de $M$ 
qui l'engendre sur $A[\frac 1 p]$. On appelle {\em topologie $p$-adique} sur $M$ l'unique topologie 
compatible avec sa structure de groupe additif  pour laquelle les sous-groupes $p^nM^\circ$,
pour $n\in \mN$, forment un système fondamental de voisinages de $0$. Cette topologie ne dépend pas
du choix de $M^\circ$. En effet, si $M'$ est un autre sous-$A$-module de type fini de $M$ qui l'engendre sur $A[\frac 1 p]$,
alors il existe $m\geq 0$ tel que $p^mM^\circ\subset M'$ et $p^mM'\subset M^\circ$. 
Il est clair que $M$ est un $A[\frac 1 p]$-module topologique.

\subsection{}\label{notconv16}
Soient $G$ un groupe profini, $A$ un anneau topologique muni d'une action continue de $G$ 
par des homomorphismes d'anneaux. Une {\em $A$-représentation} de $G$ est la donnée d'un $A$-module 
$M$ et d'une action $A$-semi-linéaire de $G$ sur $M$, {\em i.e.},  telle que 
pour tous $g\in G$, $a\in A$ et $m\in M$, on ait $g(am)=g(a)g(m)$.
On dit que la $A$-représentation est {\em continue} si $M$ est un $A$-module topologique 
et si l'action de $G$ sur $M$ est continue. Soient $M$, $N$ deux $A$-représentations 
(resp. deux $A$-représentations continues) de $G$. 
Un morphisme de $M$ dans $N$ est la donnée d'un morphisme $A$-linéaire et $G$-équivariant 
(resp. $A$-linéaire, continu et $G$-équivariant) de $M$ dans $N$. 
On note $\bRep_A(G)$ (resp. $\bRep_A^\cont(G)$)
la catégorie des $A$-représentations (resp. $A$-représentations continues) de $G$.
Si $M$ et $N$ sont deux $A$-représentations de $G$, 
les $A$-modules $M\otimes_AN$ et $\Hom_A(M,N)$ sont naturellement des $A$-représentations de $G$. 

\subsection{}\label{notconv1}
Soient $A$ un anneau, $p$ un nombre premier, $n$ un entier $\geq1$.  
On désigne par $\rW(A)$ (resp. $\rW_n(A)$) l'anneau des vecteurs de Witt 
(resp. vecteurs de Witt de longueur $n$) à coefficients dans $A$ relatif à $p$. 
On a un homomorphisme d'anneaux
\begin{equation}\label{notconv1a}
\Phi_n\colon 
\begin{array}[t]{clcr}
\rW_n(A)&\rightarrow& A,\\ 
(x_0,\dots,x_{n-1})&\mapsto&x_0^{p^{n-1}}+p x_1^{p^{n-2}}+\dots+p^{n-1}x_{n-1}
\end{array}
\end{equation}
appelé $n$-ième composante fantôme. 
On dispose aussi des morphismes de restriction, de décalage et de Frobenius
\begin{eqnarray}
\rR\colon \rW_{n+1}(A)&\rightarrow& \rW_n(A),\label{notconv1b}\\
\rV\colon \rW_n(A)&\rightarrow& \rW_{n+1}(A),\label{notconv1c}\\
\rF\colon \rW_{n+1}(A)&\rightarrow& \rW_n(A).\label{notconv1d}
\end{eqnarray}
Lorsque $A$ est de caractéristique $p$, $\rF$ induit un endomorphisme de $\rW_n(A)$, encore noté $\rF$.

\subsection{} \label{notconv2}
Pour tout anneau $R$ et tout monoïde $M$, on désigne par $R[M]$ la $R$-algèbre de $M$ et par 
$e\colon M\rightarrow R[M]$ l'homomorphisme canonique, où $R[M]$ est considéré comme un
monoïde multiplicatif. Pour tout $x\in M$, on notera $e^x$ au lieu de $e(x)$. 

On désigne par $\bA_M$ le schéma $\Spec(\mZ[M])$ muni de la structure logarithmique associée à 
la structure pré-logarithmique définie par $e\colon M\rightarrow \mZ[M]$ (\cite{agt} II.5.9). 
Pour tout homomorphisme de monoïdes $\vartheta\colon M\rightarrow N$, 
on note $\bA_\vartheta\colon \bA_N\rightarrow \bA_M$ le morphisme de schémas logarithmiques associé.

\subsection{}\label{notconv4}
Pour une catégorie $\cC$, nous notons $\ob(\cC)$ l'ensemble de ses objets,
$\cC^\circ$ la catégorie opposée, et pour $X,Y\in \ob(\cC)$, 
$\Hom_\cC(X,Y)$ (ou $\Hom(X,Y)$ lorsqu'il n'y a aucune ambiguïté) 
l'ensemble des morphismes de $X$ dans $Y$. 

Si $\cC$ et $\cC'$ sont deux catégories, nous désignons par $\Hom(\cC,\cC')$ 
l'ensemble des foncteurs de $\cC$ dans $\cC'$, et  
par $\bHom(\cC,\cC')$ la catégorie des foncteurs de $\cC$ dans $\cC'$. 

Soient $I$ une catégorie, $\cC$ et $\cC'$ deux catégories sur $I$ (\cite{sga1} VI 2). 
Nous notons $\Hom_{I}(\cC,\cC')$ l'ensemble des $I$-foncteurs de $\cC$ dans $\cC'$ 
et $\Hom_{\cart/I}(\cC,\cC')$ l'ensemble des foncteurs cartésiens (\cite{sga1} VI 5.2).
Nous désignons par $\bHom_{I}(\cC,\cC')$ la catégorie des $I$-foncteurs de $\cC$ dans $\cC'$ et 
par $\bHom_{\cart/I}(\cC,\cC')$ la sous-catégorie pleine formée des foncteurs cartésiens.

\subsection{}\label{notconv18}
Pour toute catégorie abélienne $\cC$, on désigne par 
$\bD(\cC)$ sa catégorie dérivée et par $\bD^-(\cC)$, $\bD^+(\cC)$ et $\bD^\rb(\cC)$ les sous-catégories
pleines de $\bD(\cC)$ formées des complexes à cohomologie bornée
supérieurement, inférieurement et des deux côtés, respectivement  (\cite{ks2} 13.1.2). 
Sauf mention expresse du contraire, les complexes de $\cC$ sont à différentielles de degré $+1$,
le degré étant écrit en exposant. 

\subsection{}\label{notconv3}
Dans tout cette monographie, on fixe un univers $\mU$ possédant un élément de cardinal infini. 
On dit qu'un ensemble est {\em $\mU$-petit} (ou, quand aucune confusion n'en résulte, {\em petit}) s'il est isomorphe à un élément de $\mU$. 
On utilise aussi la terminologie : {\em petit groupe}, {\em petit anneau}, {\em petite catégorie}...
On dit qu'une catégorie $\cC$ est une {\em $\mU$-catégorie} si pour tous objets $X,Y$ de $\cC$, l'ensemble $\Hom_\cC(X,Y)$ est $\mU$-petit (\cite{sga4} I 1.1).
On appelle {\em catégorie des $\mU$-ensembles} et l'on note $\Ens$, 
la catégorie des ensembles qui se trouvent dans $\mU$. 
C'est un $\mU$-topos ponctuel que l'on note encore $\Pt$ (\cite{sga4} IV 2.2). 
On désigne par $\Sch$ la catégorie des schémas éléments de $\mU$. 
Sauf mention explicite du contraire, il sera sous-entendu que les anneaux et 
les schémas logarithmiques (et en particulier les schémas) envisagés dans cette monographie sont éléments de l'univers~$\mU$.

\subsection{}\label{notconv5}
Soit $\cC$ une $\mU$-catégorie \eqref{notconv3} (\cite{sga4} I 1.1). On désigne par $\hcC$ la catégorie des préfaisceaux 
de $\mU$-ensembles sur $\cC$, c'est-à-dire la catégorie des foncteurs
contravariants sur $\cC$ à valeurs dans $\Ens$ (\cite{sga4} I 1.2). 
Si $\cC$ est munie d'une topologie (\cite{sga4} II 1.1), on désigne par $\tcC$ le topos 
des faisceaux de $\mU$-ensembles sur $\cC$ (\cite{sga4} II 2.1). 

Pour $F$ un objet de $\hcC$, on note $\cC_{/F}$ la catégorie 
suivante (\cite{sga4} I 3.4.0). Les objets de $\cC_{/F}$
sont les couples formés d'un objet $X$ de $\cC$ 
et d'un morphisme $u$ de $X$ dans $F$. Si $(X,u)$ et $(Y,v)$ sont deux objets, 
un morphisme de $(X,u)$ vers $(Y,v)$ est un morphisme $g\colon X\rightarrow Y$
tel que $u=v\circ g$.

\subsection{}\label{notconv13}
Soit $X$ un $\mU$-topos (\cite{sga4} IV 1.1.2). Les systèmes projectifs d'objets de $X$ indexés par l'ensemble ordonné 
des entiers naturels $\mN$, forment un topos que l'on note $X^{\mN^\circ}$. On renvoie à  
(\cite{agt} III.7) pour des sorites utiles sur ce type de topos. On rappelle, en particulier, qu'on a un morphisme de topos 
\begin{equation}\label{notconv13a}
\uplambda\colon X^{\mN^\circ}\rightarrow X,
\end{equation}
dont le foncteur image inverse $\uplambda^*$ associe à tout objet $F$ de $X$ le système projectif constant de valeur $F$,
et dont le foncteur image directe $\uplambda_*$ associe à tout système projectif sa limite projective (\cite{agt} III.7.4).

\subsection{}\label{notconv9}
Soit $(X,A)$ un $\mU$-topos annelé (\cite{sga4} IV 11.1.1). On note $\bMod(A)$ ou $\bMod(A,X)$ 
la catégorie des $A$-modules de $X$. 
Si $M$ un $A$-module, on désigne par $\rS_A(M)$ 
(resp. $\wedge_A(M)$, resp. $\Gamma_A(M)$) l'algèbre symétrique (resp. extérieure, resp. à puissances divisées) 
de $M$ et pour tout entier $n\geq 0$, par $\rS_A^n(M)$ (resp. $\wedge_A^n(M)$,
resp. $\Gamma_A^n(M)$) sa partie homogène de degré $n$. 
Les formations de ces algèbres commutent à la localisation au-dessus d'un objet de $X$.
On omettra l'anneau $A$ des notations $\wedge_A(M)$ et $\wedge_A^n(M)$ lorsqu'il n'y a aucun risque d'ambiguïté. 

\begin{defi}[\cite{sga6} I 1.3.1]\label{notconv14}
Soit $(X,A)$ un topos annelé. 
On dit qu'un $A$-module $M$ de $X$ est {\em localement projectif de type fini} si les conditions 
équivalentes suivantes sont satisfaites :
\begin{itemize}
\item[{\rm (i)}] $M$ est de type fini et le foncteur $\cHom_A(M,\cdot)$ est exact~;
\item[{\rm (ii)}] $M$ est de type fini et tout épimorphisme de $A$-modules $N\rightarrow M$ admet localement une section~;
\item[{\rm (iii)}] $M$ est localement facteur direct d'un $A$-module libre de type fini. 
\end{itemize}
\end{defi}

Lorsque $X$ a suffisamment de points et que pour tout point $x$ de $X$, la fibre de $A$ en $x$
est un anneau local, les $A$-modules localement projectifs de type fini sont les $A$-modules 
localement libres de type fini (\cite{sga6} I 2.15.1).

\subsection{}\label{notconv10}
Pour tout schéma $X$, on désigne par $\Et_{/X}$ le {\em site étale} de $X$, 
c'est-à-dire, la sous-catégorie pleine de $\Sch_{/X}$ \eqref{notconv3} formée des schémas étales sur $X$,
munie de la topologie étale; c'est un $\mU$-site. 
On note $X_\et$ le {\em topos étale} de $X$, c'est-à-dire le topos des faisceaux de $\mU$-ensembles sur $\Et_{/X}$.

On désigne par $\Et_{\coh/X}$ la sous-catégorie pleine de $\Et_{/X}$ formée des schémas étales de présentation finie sur $X$, 
munie de la topologie induite par celle de $\Et_{/X}$; c'est un site $\mU$-petit. 
Si $X$ est quasi-séparé, le foncteur de restriction de $X_\et$ dans le topos des faisceaux de $\mU$-ensembles sur 
$\Et_{\coh/X}$ est une équivalence de catégories (\cite{sga4} VII 3.1 et 3.2).

On désigne par $\Et_{\rf/X}$ la sous-catégorie pleine de $\Et_{/X}$ formée des schémas étales finis sur $X$, 
munie de la topologie induite par celle de $\Et_{/X}$; c'est un site $\mU$-petit. 
On appelle {\em topos fini étale} de $X$ et on note  $X_\fet$, le topos des faisceaux de $\mU$-ensembles sur $\Et_{\rf/X}$ 
(cf. \cite{agt} VI.9.2). L'injection canonique $\Et_{\rf/X}\rightarrow \Et_{/X}$ induit un morphisme de topos 
\begin{equation}\label{notconv10a}
\rho_X\colon X_\et\rightarrow X_\fet.
\end{equation}

\subsection{}\label{notconv12}
Soit $X$ un schéma. On désigne par $X_\zar$ le topos de Zariski de $X$ et par 
\begin{equation}\label{notconv12b}
u_X\colon X_{\et}\rightarrow X_{\zar}
\end{equation}
le morphisme canonique (\cite{sga4} VII 4.2.2). Si $F$ est un $\co_X$-module quasi-cohérent de $X_\zar$, 
on désigne par $\iota(F)$ le faisceau de $X_\et$ défini pour tout $X$-schéma étale $U$ par (\cite{sga4} VII 2 c))
\begin{equation}\label{notconv12c}
\iota(F)(U)=\Gamma(U,F\otimes_{\co_X}\co_U).
\end{equation}
Il est commode, lorsqu'il n'y aucune risque de confusion, de désigner $\iota(F)$ abusivement par $F$. 
On notera que $\iota(\co_X)$ est un anneau de $X_\et$ et que $\iota(F)$ est un $\iota(\co_X)$-module (cf. \cite{ag} 2.1.18). 

Notant $\bMod^\qcoh(\co_X,X_\zar)$ la sous-catégorie pleine de $\bMod(\co_X,X_\zar)$ formée des $\co_X$-modules quasi-cohérents \eqref{notconv9},
la correspondance $F\mapsto \iota(F)$ définit un foncteur
\begin{equation}\label{notconv12a}
\iota\colon \bMod^\qcoh(\co_X,X_\zar)\rightarrow \bMod(\co_X,X_\et).
\end{equation}

Pour tout $\co_X$-module quasi-cohérent $F$ de $X_\zar$, on a un isomorphisme canonique 
\begin{equation}\label{notconv12d}
F\stackrel{\sim}{\rightarrow}u_{X*}(\iota(F)).
\end{equation}
Nous considérons donc $u_X$ comme un morphisme de topos annelés 
\begin{equation}\label{notconv12j}
u_X\colon (X_{\et},\co_X)\rightarrow (X_{\zar},\co_X). 
\end{equation}
Nous utilisons pour les modules la notation $u^{-1}_X$ pour désigner l'image
inverse au sens des faisceaux abéliens et nous réservons la notation 
$u^*_X$ pour l'image inverse au sens des modules. 
D'après (\cite{ag} 2.1.18), l'isomorphisme \eqref{notconv12d} induit par adjonction un isomorphisme
\begin{equation}\label{notconv12e}
u^*_X(F)\stackrel{\sim}{\rightarrow} \iota(F).
\end{equation}

\subsection{}\label{notconv17}
Soient $R$ un anneau muni d'une valuation non discrète de hauteur $1$, $\fm$ son idéal maximal, $\gamma\in R$,
$X$ un $\mU$-topos, $A$ une $R$-algèbre de $X$.  
Nous considérons les notions de $\alpha$-algèbre (ou presque-algèbre) sur $R$ introduites dans (\cite{ag} 2.6--2.9) (cf. \cite{ag} 2.10.1). 
On dit qu'un morphisme de $A$-modules est un {\em $\gamma$-isomorphisme} (resp. {\em $\alpha$-isomorphisme}) 
si son noyau et son conoyau sont annulés par $\gamma$ (resp. tout élément de $\fm$)  (\cite{ag} 2.6.2 et 2.7.2). 
On dit qu'un complexe de $A$-modules à différentielle de degré $1$, le degré étant écrit en exposant,
est {\em $\gamma$-acyclique} si ses groupes  de cohomologie sont annulés par $\gamma$, 
et qu'il est {\em $\alpha$-acyclique} s'il est $\delta$-acyclique pour tout $\delta\in \fm$.
On dit qu'une suite courte de $A$-modules est {\em $\gamma$-exacte} (resp. {\em $\alpha$-exacte}) 
si elle est $\gamma$-acyclique (resp. $\alpha$-acyclique) en tant que complexe de $A$-modules.

On désigne par $\bMod(A)$ la catégorie des $A$-modules de $X$ et par $\bD(\bMod(A))$ sa catégorie dérivée.
On dit qu'un morphisme de complexes de $\bMod(A)$ est un {\em $\gamma$-quasi-isomorphisme}  (resp. {\em $\alpha$-quasi-isomorphisme}) 
s'il induit des $\gamma$-isomorphismes (resp. $\alpha$-isomorphismes) sur les modules de cohomologie. 
De même, on dit qu'un morphisme de $\bD(\bMod(A))$ est un {\em $\gamma$-isomorphisme} 
(resp. {\em $\alpha$-isomorphisme}) s'il induit des $\gamma$-isomorphismes (resp. $\alpha$-isomorphismes) sur les modules de cohomologie. 
On désigne par $\cN$ la sous-catégorie épaisse de $\bMod(A)$ formée des $A$-modules $\alpha$-nuls
et par $\bD_\cN(\bMod(A))$ la sous-catégorie pleine de $\bD(\bMod(A))$ formée des complexes
dont les modules de cohomologie sont $\alpha$-nuls. 
Cette dernière est une sous-catégorie triangulée strictement pleine et saturée de $\bD(\bMod(A))$
(\cite{sp} \href{https://stacks.math.columbia.edu/tag/06UQ}{06UQ}). 
La famille des $\alpha$-isomorphismes de $\bD(\bMod(A))$ permet un calcul de fractions bilatéral bilatère (\cite{illusie1} I 1.4.2). 
La catégorie triangulée localisée de $\bD(\bMod(A))$ par la famille des $\alpha$-isomorphismes ``représente'' la catégorie triangulée quotient de $\bD(\bMod(A))$ 
par $\bD_\cN(\bMod(A))$ (\cite{ks2} 10.2.3).

\section{Rappel sur une construction de Fontaine-Grothendieck}\label{EIPF}

\subsection{}\label{eip1}
La construction rappelée dans cette section a été introduite indépendamment par Grothendieck (\cite{grot1} IV 3.3) et Fontaine (\cite{fontaine1} 2.2).
On fixe un nombre premier $p$. Tous les anneaux des vecteurs de Witt considérés dans ce numéro sont relatifs à $p$ (cf. \ref{notconv1}). 
Soient $A$ une $\mZ_{(p)}$-algèbre, $n$ un entier $\geq 1$. L'homomorphisme d'anneaux 
\begin{equation}\label{eip1a}
\Phi_{n+1}\colon 
\begin{array}[t]{clcr}
\rW_{n+1}(A/p^nA)&\rightarrow& A/p^nA\\
(x_0,\dots,x_{n})&\mapsto& x_0^{p^n}+p x_1^{p^{n-1}}+\dots+p^{n}x_n
\end{array}
\end{equation}
s'annule sur $\rV^n(A/p^nA)$ et induit donc par passage au quotient un homomorphisme d'anneaux
\begin{equation}\label{eip1b}
\Phi'_{n+1}\colon 
\begin{array}[t]{clcr}
\rW_{n}(A/p^nA)&\rightarrow& A/p^nA\\
(x_0,\dots,x_{n-1})&\mapsto&x_0^{p^n}+p x_1^{p^{n-1}}+\dots+p^{n-1}x_{n-1}^p.
\end{array}
\end{equation}
Ce dernier s'annule sur 
\begin{equation}\label{eip1c}
\rW_n(pA/p^nA)=\ker(\rW_n(A/p^nA)\rightarrow \rW_n(A/pA))
\end{equation}
et se factorise à son tour en un homomorphisme d'anneaux
\begin{equation}\label{eip1d}
\theta_n\colon \rW_{n}(A/pA)\rightarrow A/p^nA.
\end{equation}
Il résulte aussitôt de la définition que le diagramme 
\begin{equation}\label{eip1e}
\xymatrix{
{\rW_{n+1}(A/pA)}\ar[r]^-(0.4){\theta_{n+1}}\ar[d]_{\rR\rF}&{A/p^{n+1}A}\ar[d]\\
{\rW_{n}(A/pA)}\ar[r]^-(0.4){\theta_n}&{A/p^nA}}
\end{equation}
où $\rR$ est le morphisme de restriction \eqref{notconv1b}, 
$\rF$ est le Frobenius \eqref{notconv1d} et la flèche non libellée est l'homomorphisme canonique, est commutatif. 

Pour tout homomorphisme de $\mZ_{(p)}$-algèbres commutatives $\varphi\colon A\rightarrow B$, 
le diagramme  
\begin{equation}\label{eip1f}
\xymatrix{
{\rW_n(A/pA)}\ar[r]\ar[d]_{\theta_n}&{\rW_n(B/pB)}\ar[d]^{\theta_n}\\
{A/p^nA}\ar[r]&{B/p^nB}}
\end{equation}
où les flèches horizontales sont les morphismes induits par $\varphi$, est commutatif.

\begin{prop}[\cite{agt} II.9.2]\label{eip4}
Soit $A$ une $\mZ_{(p)}$-algèbre vérifiant les conditions suivantes~:
\begin{itemize}
\item[{\rm (i)}] $A$ est $\mZ_{(p)}$-plat.
\item[{\rm (ii)}] $A$ est intégralement clos dans $A[\frac 1 p]$.
\item[{\rm (iii)}] Le Frobenius absolu  de $A/pA$ est surjectif. 
\item[{\rm (iv)}] Il existe un entier $N\geq 1$ et une suite $(p_n)_{0\leq n\leq N}$ d'éléments de $A$
tels que $p_0=p$ et $p_{n+1}^p=p_n$ pour tout $0\leq n\leq N-1$.  
\end{itemize}
Pour tout entier $1\leq n\leq N$, on pose  
\begin{equation}\label{eip4a}
\xi_n=[\op_n]-p \in \rW_n(A/pA),
\end{equation}
où $\op_n$ est la classe de $p_n$ dans $A/pA$ et $[\ ]$ désigne le représentant multiplicatif.
Alors pour tous entiers $n\geq 1$ et $i\geq 0$ tels que $n+i\leq N$, la suite
\begin{equation}\label{eip4b}
\xymatrix{
{\rW_n(A/pA)}\ar[rr]^-(0.5){\cdot \rR^i(\xi_{n+i})}&&{\rW_n(A/pA)}
\ar[rr]^-(0.5){\theta_n\circ \rF^i}&&{A/p^nA}\ar[r]& 0}
\end{equation}
est exacte.
\end{prop}

\subsection{}\label{eipo3}
Soient $A$ une $\mZ_{(p)}$-algèbre, $\hA$ le séparé complété $p$-adique de $A$. 
On désigne par $A^\flat$ la limite projective du système projectif $(A/pA)_{\mN}$ 
dont les morphismes de transition sont les itérés de l'endomorphisme de Frobenius de $A/pA$.
\begin{equation}\label{eipo3a}
A^\flat=\underset{\underset{\mN}{\longleftarrow}}{\lim}\ A/pA.
\end{equation} 
C'est un anneau parfait de caractéristique $p$.
Pour tout entier $n\geq 1$, la projection canonique $A^\flat\rightarrow A/pA$ sur la $(n+1)$-ième composante
du système projectif $(A/pA)_{\mN}$ ({\em i.e.}, la composante d'indice $n$) induit un homomorphisme \eqref{notconv1}
\begin{equation}\label{eipo3b}
\nu_n\colon \rW(A^\flat)\rightarrow \rW_n(A/pA).
\end{equation}
Comme $\nu_n=\rF\circ\rR\circ \nu_{n+1}$, on obtient par passage à la limite projective un homomorphisme 
\begin{equation}\label{eipo3c}
\nu\colon \rW(A^\flat)\rightarrow \underset{\underset{n\geq 1}{\longleftarrow}}{\lim}\ \rW_n(A/pA),
\end{equation}
où les morphismes de transition de la limite projective sont les morphismes $\rF \rR$. On vérifie aussitôt  
qu'il est bijectif. Compte tenu de \eqref{eip1e},
les homomorphismes $\theta_n$ induisent par passage à la limite projective un 
homomorphisme 
\begin{equation}\label{eipo3d}
\theta\colon \rW(A^\flat)\rightarrow \hA.
\end{equation}
On retrouve l'homomorphisme défini par Fontaine (\cite{fontaine1} 2.2). 
Pour tout entier $r\geq 1$, on pose
\begin{equation}\label{eipo3e}
\cA_r(A)=\rW(A^\flat)/\ker(\theta)^r,
\end{equation}
et on note $\theta_r\colon \cA_r(A)\rightarrow \hA$ l'homomorphisme induit par $\theta$ (cf. \cite{fontaine3} 1.2.2). 

Pour tout homomorphisme de $\mZ_{(p)}$-algèbres commutatives $\varphi\colon A\rightarrow B$, 
le diagramme  
\begin{equation}\label{eip3f}
\xymatrix{
{\rW(A^\flat)}\ar[r]\ar[d]_{\theta}&{\rW(B^\flat)}\ar[d]^{\theta}\\
{\hA}\ar[r]&{\hB}}
\end{equation}
où les flèches horizontales sont les morphismes induits par $\varphi$, est commutatif \eqref{eip1f}.  
La correspondance $A\mapsto \cA_r(A)$ est donc fonctorielle. 

\begin{rema}\label{eip2}
Soit $k$ un corps parfait.  
La projection canonique $\rW(k)^\flat\rightarrow k$ sur la première composante ({\em i.e.}, d'indice $0$)
est un isomorphisme. Elle induit donc un isomorphisme $\rW(\rW(k)^\flat)\stackrel{\sim}{\rightarrow}\rW(k)$,
que nous utilisons pour identifier ces deux anneaux. L'homomorphisme $\theta$ s'identifie alors à
l'endomorphisme identique de $\rW(k)$.
\end{rema}

\begin{lem}\label{eip5}
Soient $A$ une $\mZ_{(p)}$-algèbre, $\hA$ le séparé complété $p$-adique de $A$, 
$A^\flat$ l'anneau défini dans \eqref{eipo3a}, $\mA$ l'ensemble des suites $(x_n)_{\mN}$ de $\hA$ 
telles que $x_{n+1}^p= x_n$ pour tout $n\geq 0$. Alors, 
\begin{itemize}
\item[{\rm (i)}] L'application
\begin{equation}\label{eip5a}
\mA\rightarrow A^\flat, \ \ \ (x_n)_{\mN}\mapsto (\ox_n)_{\mN},
\end{equation}
où $\ox_n$ est la réduction de $x_n$ modulo $p$, est un isomorphisme de monoïdes multiplicatifs. 
\item[{\rm (ii)}] 
Pour tout $(x_0,x_1,\dots)\in \rW(A^\flat)$, on a 
\begin{equation}\label{eip5b}
\theta(x_0,x_1,\dots) = \sum_{n\geq 0}p^nx_n^{(n)},
\end{equation}
où $\theta$ est l'homomorphisme \eqref{eipo3d} et pour tout $n\geq 0$, $(x_n^{(m)})_{m\geq 0}$ est l'élément de $\mA$ associé à $x_n\in A^\flat$ \eqref{eip5a}.
\end{itemize}
\end{lem}

(i) Notant $v$ la valuation de $\mZ_{(p)}$ normalisée par $v(p)=1$, pour tous entiers $m,i\geq 1$ tels que $i\leq p^m$, on a $v(\binom{p^m}{i})=m-v(i)$
et donc $v(\binom{p^m}{i})+i\geq m$. 
Soient $(x_n)_{\mN}$ et $(y_n)_{\mN}$ deux suites d'éléments de $\hA$ qui induisent par réduction modulo $p$ la même suite 
$(\ox_n)_{\mN}=(\oy_n)_{\mN}$ de $A/pA$. Pour tous $n,m\geq 0$, on a alors
\begin{equation}\label{eip5c}
x_{n+m}^{p^m}-y_{n+m}^{p^m} \in p^m \hA.
\end{equation}
En particulier, si $(x_n)_{\mN}$ et $(y_n)_{\mN}$ sont des éléments de $\mA$, 
alors $x_n=y_n$ pour tout $n\geq 0$ puisque $\hA$ est séparé pour la topologie $p$-adique.
Par suite, l'application \eqref{eip5a} est injective. Montrons qu'elle est surjective. Soient $(z_n)_{\mN}$ un élément de $A^\flat$, 
$(y_n)_{\mN}$ une suite d'éléments de $A$ qui relève $(z_n)_{\mN}$. Appliquant \eqref{eip5c} aux suites $(y_n)_{n\geq 0}$ et $(y_{n+1}^p)_{n\geq 0}$,
on voit que pour tous $n,m\geq 0$, on a 
\begin{equation}\label{eip5d}
y_{n+m+1}^{p^{m+1}}-y_{n+m}^{p^m} \in p^m \hA.
\end{equation}
Par suite, pour tout entier $n\geq 0$,
$(y_{n+m}^{p^m})_{m\geq 0}$ converge vers un élément $x_n$ de $\hA$. La suite $(x_n)_{\mN}$ appartient clairement à $\mA$ et s'envoie 
sur $(z_n)_{\mN}$ par \eqref{eip5a}; d'où la surjectivité.  

(ii) Cela résulte aussitôt des définitions.

\section{\texorpdfstring{\'Epaississements infinitésimaux $p$-adiques universels de Fontaine}
{Epaississements infinitésimaux p-adiques universels de Fontaine}}\label{eipuf}

\subsection{}\label{eipuf1}
Soient $\Lambda$ un anneau, $A$ une $\Lambda$-algèbre, $m$ un entier $\geq 0$. Un $\Lambda$-épaississement infinitésimal 
d'ordre $\leq m$ de $A$ est la donnée d'un couple $(D,\theta)$ formé d'une $\Lambda$-algèbre $D$ et d'un homomorphisme 
surjectif de $\Lambda$-algèbres $\theta\colon D\rightarrow A$ tels que $(\ker(\theta))^{m+1}=0$. 
Les $\Lambda$-épaississements infinitésimaux d'ordre $\leq m$ de $A$ forment une catégorie notée $\bE_m(A/\Lambda)$: si $(D_1,\theta_1)$ et $(D_2,\theta_2)$
sont deux $\Lambda$-épaississements infinitésimaux d'ordre $\leq m$ de $A$, un morphisme de $(D_1,\theta_1)$ vers $(D_2,\theta_2)$
est la donnée d'un homomorphisme de $\Lambda$-algèbres $f\colon D_1\rightarrow D_2$ tel que $\theta_1=\theta_2\circ f$. 
Suivant Fontaine (\cite{fontaine3} 1.1.1), si cette catégorie admet un objet initial, on l'appelle 
{\em $\Lambda$-épaississement infinitésimal universel d'ordre $\leq m$ de $A$}. 

Soit $p$ un nombre premier et supposons, de plus, que $A$ soit complet et séparé pour la topologie $p$-adique. Un objet $(D,\theta)$ de $\bE_m(A/\Lambda)$ 
est dit {\em $p$-adique} si la $\Lambda$-algèbre $D$ est complète et séparée pour la topologie $p$-adique. On note $\bE^{p}_m(A/\Lambda)$
la sous-catégorie pleine de $\bE_m(A/\Lambda)$ formée des $\Lambda$-épaississements infinitésimaux $p$-adiques d'ordre $\leq m$ de $A$. 
Suivant Fontaine (\cite{fontaine3} 1.1.3), si la catégorie $\bE^p_m(A/\Lambda)$ admet un objet initial, on l'appelle {\em $\Lambda$-épaississement infinitésimal $p$-adique universel d'ordre $\leq m$ de $A$}.

\subsection{}\label{eipuf2}
Soient $k$ un corps parfait de caractéristique $p>0$, $W$ l'anneau des vecteurs de Witt à coefficients dans $k$,
$\Lambda$ une $W$-algèbre, $A$ une $\Lambda$-algèbre.
Tous les anneaux des vecteurs de Witt considérés dans ce numéro sont relatifs à $p$ (cf. \ref{notconv1}). 
On note $\hA$ le séparé complété $p$-adique de $A$, $A^\flat$ l'anneau associé à $A$ par le foncteur défini dans \eqref{eipo3a} et 
\begin{equation}\label{eipuf2a}
\theta\colon \rW(A^\flat)\rightarrow \hA
\end{equation}
l'homomorphisme canonique \eqref{eipo3d}.
D'après \ref{eip2}, on peut canoniquement munir $A^\flat$ d'une structure de $k$-algèbre,  
et donc $\rW(A^\flat)$ d'une structure de $W$-algèbre. De plus, $\theta$ est un homomorphisme de $W$-algèbres. 
On pose 
\begin{equation}\label{eipuf2b}
\rW_\Lambda(A^\flat)=\rW(A^\flat)\otimes_W\Lambda.
\end{equation} 
On désigne par 
\begin{equation}\label{eipuf2c}
\theta_\Lambda\colon \rW_\Lambda(A^\flat)\rightarrow \hA
\end{equation}
l'homomorphisme de $\Lambda$-algèbres induit par $\theta$, et par $J_\Lambda$ son noyau. 
Pour tout entier $r\geq 1$, on pose 
\begin{equation}\label{eipuf2f}
\cA_r(A/\Lambda)=\rW_\Lambda(A^\flat)/J_\Lambda^r
\end{equation}
et on note $\hcA_r(A/\Lambda)$ son séparé complété $p$-adique,
\begin{equation}\label{eipuf2d}
\hcA_r(A/\Lambda)=\underset{\underset{n\geq 0}{\longleftarrow}}{\lim}\ \frac{\rW_\Lambda(A^\flat)}{J^r_\Lambda+p^n\rW_\Lambda(A^\flat)}.
\end{equation}
On désigne par $\theta_{\Lambda,r}\colon \cA_r(A/\Lambda)\rightarrow \hA$ l'homomorphisme induit par $\theta_\Lambda$ et par 
$\htheta_{\Lambda,r}\colon \hcA_r(A/\Lambda)\rightarrow \hA$ son prolongement aux complétés.

On désigne par $\rW'_\Lambda(A^\flat)$ le séparé complété de $\rW_\Lambda(A^\flat)$ pour la topologie définie par l'idéal 
$\theta_\Lambda^{-1}(p\hA)$,  par
\begin{equation}\label{eipuf2e}
\theta'_\Lambda\colon \rW'_\Lambda(A^\flat)\rightarrow \hA
\end{equation}
le prolongement de $\theta_\Lambda$ aux complétés, par $J'_\Lambda$ le noyau de $\theta'_\Lambda$,  
et pour tout entier $r\geq 1$, par $\hcA'_r(A/\Lambda)$ le séparé complété $p$-adique de $\rW'_\Lambda(A^\flat)/J'^r_\Lambda$.

\begin{lem}\label{eipuf4}
Conservons les hypothèses et notations de \ref{eipuf2}, supposons de plus que l'endomorphisme de Frobenius de $A/pA$ soit surjectif. Alors,
\begin{itemize}
\item[{\rm (i)}] Les homomorphismes $\theta$ et $\theta_\Lambda$ sont surjectifs.
\item[{\rm (ii)}] Pour tout $r\geq 1$, on a $\ker(\htheta_{\Lambda,r})^r=0$. 
\item[{\rm (iii)}] Si, de plus, $J_\Lambda$ est un idéal de type fini de $W_\Lambda(A^\flat)$,   
alors, pour tout entier $r\geq 1$, l'homomorphisme canonique $\rW_\Lambda(A^\flat)\rightarrow \rW'_\Lambda(A^\flat)$ induit 
un isomorphisme de $\Lambda$-algèbres
\begin{equation}\label{eipuf4a}
\hcA_r(A/\Lambda)\rightarrow \hcA'_r(A/\Lambda).
\end{equation}
\end{itemize}
\end{lem}

(i) L'homomorphisme $\theta$ se réduit modulo $p$ en le morphisme $A^\flat \rightarrow A/pA$ 
de projection sur la deuxième composante ({\em i.e.}, la composante d'indice $1$). Celui-ci est surjectif puisque l'endomorphisme de Frobenius de $A/pA$ est surjectif.
On en déduit que $\theta$ est surjectif car $\rW(A^\flat)$ et $\hA$ sont complets et séparés pour les topologies $p$-adiques (\cite{egr1} 1.8.5).  
Il en est alors de même de $\theta_\Lambda$.

(ii) En effet, compte tenu de \eqref{eipuf2d}, le noyau de $\htheta_{\Lambda,r}$ s'identifie à 
\begin{equation}\label{eipuf4k}
\underset{\underset{n\geq 0}{\longleftarrow}}{\lim}\ \frac{J_\Lambda+p^n\rW_\Lambda(A^\flat)}{J^r_\Lambda+p^n\rW_\Lambda(A^\flat)}.
\end{equation}
On a donc $\ker(\htheta_{\Lambda,r})^r=0$.

(iii) On notera d'abord que comme  $\theta_\Lambda$ et $\theta'_\Lambda$ sont surjectifs, pour tout entier $n\geq 1$, on a 
\begin{eqnarray}
\theta^{-1}_\Lambda(p^n\hA)&=&J_\Lambda+p^n \rW_\Lambda(A^\flat),\\
\theta'^{-1}_\Lambda(p^n\hA)&=&J'_\Lambda+p^n \rW'_\Lambda(A^\flat).
\end{eqnarray}
Notons
\begin{equation}\label{eipuf4b}
h_n\colon \rW'_\Lambda(A^\flat)\rightarrow \rW_\Lambda(A^\flat)/(J_\Lambda+p \rW_\Lambda(A^\flat))^n
\end{equation}
l'homomorphisme canonique. L'homomorphisme $\theta_\Lambda$ induit un homomorphisme injectif
\begin{equation}\label{eipuf4c}
\rW_\Lambda(A^\flat)/(J_\Lambda+p \rW_\Lambda(A^\flat)) \rightarrow A/pA,
\end{equation}
dont le composé avec $h_1$ est induit par $\theta'_\Lambda$. On en déduit que $\ker(h_1)=J'_\Lambda+p\rW'_\Lambda(A^\flat)$. 
D'après (\cite{ac} chap. III §~2.11 prop.~14 et cor.~1), comme $J_\Lambda$ est de type fini, 
pour tout $n\geq1$, $h_n$ est surjectif et $\ker(h_n)=(J'_\Lambda+p\rW'_\Lambda(A^\flat))^n$. 
Par suite, $h_n$ induit un isomorphisme 
\begin{equation}\label{eipuf4d}
\frac{\rW'_\Lambda(A^\flat)}{(J'_\Lambda+p\rW'_\Lambda(A^\flat))^n}\stackrel{\sim}{\rightarrow} \frac{\rW_\Lambda(A^\flat)}{(J_\Lambda+p \rW_\Lambda(A^\flat))^n}.
\end{equation}
Celui-ci s'insère dans un diagramme commutatif
\begin{equation}\label{eipuf4e}
\xymatrix{
{\rW'_\Lambda(A^\flat)/(J'_\Lambda+p\rW'_\Lambda(A^\flat))^n}\ar[r]^\sim\ar[rd]_{\theta'^{(n)}_{\Lambda}}&
{\rW_\Lambda(A^\flat)/(J_\Lambda+p \rW_\Lambda(A^\flat))^n}\ar[d]^{\theta^{(n)}_{\Lambda}}\\
&{A/p^nA}}
\end{equation}
où $\theta^{(n)}_{\Lambda}$ (resp. $\theta'^{(n)}_{\Lambda}$) est induit par $\theta_\Lambda$ (resp. $\theta'_\Lambda$). Par suite, \eqref{eipuf4d} 
induit un isomorphisme entre les noyaux de $\theta_{\Lambda,n}$ et $\theta'_{\Lambda,n}$
\begin{equation}\label{eipuf4f}
\frac{J'_\Lambda+p^n\rW'_\Lambda(A^\flat)}{(J'_\Lambda+p\rW'_\Lambda(A^\flat))^n}\stackrel{\sim}{\rightarrow} 
\frac{J_\Lambda+p^n\rW_\Lambda(A^\flat)}{(J_\Lambda+p \rW_\Lambda(A^\flat))^n}.
\end{equation}
Par ailleurs, pour tout entier $m\geq 0$, \eqref{eipuf4d}  induit un isomorphisme 
\begin{equation}\label{eipuf4g}
\frac{p^m\rW'_\Lambda(A^\flat)+(J'_\Lambda+p\rW'_\Lambda(A^\flat))^n}{(J'_\Lambda+p\rW'_\Lambda(A^\flat))^n}\stackrel{\sim}{\rightarrow} 
\frac{p^m\rW_\Lambda(A^\flat)+(J_\Lambda+p\rW_\Lambda(A^\flat))^n}{(J_\Lambda+p \rW_\Lambda(A^\flat))^n}.
\end{equation}
Remplaçons $n$ par $n+r$ et prenons $m=n$. Comme on a 
\[
(J'_\Lambda+p^{n+r}\rW'_\Lambda(A^\flat))^r+(J'_\Lambda+p\rW'_\Lambda(A^\flat))^{n+r}+ p^n\rW'_\Lambda(A^\flat)=J'^r_\Lambda+p^n\rW'_\Lambda(A^\flat),
\]
et de même pour $\rW_\Lambda(A^\flat)$ et $J_\Lambda$, on déduit de \eqref{eipuf4f} et \eqref{eipuf4g} un isomorphisme 
\begin{equation}\label{eipuf4h}
\frac{J'^r_\Lambda+p^n\rW'_\Lambda(A^\flat)}{(J'_\Lambda+p\rW'_\Lambda(A^\flat))^{n+r}}\stackrel{\sim}{\rightarrow} 
\frac{J^r_\Lambda+p^n\rW_\Lambda(A^\flat)}{(J_\Lambda+p \rW_\Lambda(A^\flat))^{n+r}}.
\end{equation}
L'isomorphisme \eqref{eipuf4d} (pour $n+r$) induit donc un isomorphisme 
\begin{equation}\label{eipuf4i}
\frac{\rW'_\Lambda(A^\flat)}{J'^r_\Lambda+p^n\rW'_\Lambda(A^\flat)}\stackrel{\sim}{\rightarrow} \frac{\rW_\Lambda(A^\flat)}{J^r_\Lambda+p^n \rW_\Lambda(A^\flat)}.
\end{equation}
Par passage à la limite projective, on obtient un isomorphisme 
\begin{equation}\label{eipuf4j}
\hcA'_r(A/\Lambda)\stackrel{\sim}{\rightarrow}\hcA_r(A/\Lambda)
\end{equation}
dont le composé avec le morphisme canonique \eqref{eipuf4a} est clairement l'identité de $\hcA_r(A/\Lambda)$. 

\begin{rema}
Sous les hypothèses de \ref{eipuf2}, Fontaine (\cite{fontaine3} 1.2.1) affirme que $\rW'_\Lambda(A^\flat)$ est complet et séparé pour la topologie $p$-adique,
mais sa preuve présente un problème. 
\end{rema}

\begin{prop}[\cite{fontaine3} 1.2.1]\label{eipuf9}
Conservons les hypothèses et notations de \ref{eipuf2}, supposons, de plus, que l'endomorphisme de Frobenius de $A/pA$ soit surjectif.
Alors, pour entier $r\geq 1$, $(\hcA_r(A/\Lambda),\htheta_{\Lambda, r})$ est un $\Lambda$-épaississement infinitésimal $p$-adique universel d'ordre $\leq r-1$ de $\hA$.
\end{prop}

On notera d'abord que $(\hcA_r(A/\Lambda),\htheta_{\Lambda, r})$ est un objet de la catégorie $\bE^p_{r-1}(\hA/\Lambda)$ \eqref{eipuf1} en vertu de \ref{eipuf4}. 
Soient $(D,\theta_D)$ un objet de $\bE^p_{r-1}(\hA/\Lambda)$, $I_D=\ker(\theta_D)$. Notons $\mD$ (resp. $\mA$) l'ensemble des suites $(x_n)_{\mN}$ d'éléments de $D$
(resp. $\hA$) telles que $x_{n+1}^p=x_n$.  Montrons que l'application 
\begin{equation}\label{eipuf9a}
\mD\mapsto \mA, \ \ \ (x_n)_{\mN}\mapsto (\theta_D(x_n))_{\mN}
\end{equation}
est bijective. On notera d'abord que $I_D+pD$ est un idéal de définition pour la topologie $p$-adique de $D$. En effet, pour tout $N\geq r$, on a 
\begin{equation}\label{eipuf9b}
p^ND\subset (I_D+pD)^N\subset p^{N-r}D.
\end{equation}    
Soient $(x_n)_{\mN}$ et $(y_n)_{\mN}$ deux suites d'éléments de $D$ telles que $\theta_D(x_n)=\theta_D(y_n)$ pour tout $n\geq 0$. 
Comme $x_n\equiv y_n$ modulo $I_D$ et donc à fortiori modulo $I_D+pD$, 
procédant comme dans la preuve de \ref{eip5}(i), on voit que pour tous $n,m\geq 0$, on a 
\begin{equation}\label{eipuf9c}
x_{n+m}^{p^m}-y_{n+m}^{p^m} \in (I_D+pD)^m.
\end{equation}
Par suite, l'application \eqref{eipuf9a} est injective. Supposons que $(\theta_D(x_n))_\mN$ appartienne à $\mA$. 
Appliquant \eqref{eipuf9c} aux suites $(x_n)_{n\geq 0}$ et $(x_{n+1}^p)_{n\geq 0}$, on en déduit que pour tout $n\geq 0$, la suite $(x_{n+m}^{p^m})_{m\geq 0}$
converge vers un élément $\tx_n$ de $D$ qui ne dépend que de $(\theta_D(x_n))_{\mN}$ mais pas de $(x_n)_{\mN}$. 
Il est clair que la suite $(\tx_n)_{\mN}$ appartient à $\mD$ et que l'on a $\theta_D(\tx_n)=\theta_D(x_n)$ pour tout $n\geq 0$; d'où la surjectivité de  \eqref{eipuf9a}.

Montrons qu'il existe un unique homomorphisme 
\begin{equation}\label{eipuf9d}
\alpha\colon \rW(A^\flat)\rightarrow D
\end{equation}
tel que $\theta=\theta_D\circ \alpha$. Supposons d'abord donné un tel morphisme et montrons qu'il est unique. 
Soient $x\in A^\flat$, $(x^{(n)})_\mN$ la suite associée de $\mA$ (cf. \ref{eip5}), $(\tx^{(n)})_\mN$ son image inverse par l'isomorphisme \eqref{eipuf9a}.
Pour tout $n\geq 0$, considérant $(x^{(n+m)})_{m\geq 0}$ comme un élément de $A^\flat$, on a 
\begin{equation}\label{eipuf9e}
\theta([(x^{(n+m)})_{m\geq 0}])=x^{(n)}.
\end{equation} 
La suite $(\alpha([(x^{(n+m)})_{m\geq 0}]))_{n\geq 0}$ étant clairement un élément de $\mD$, on en déduit que l'on a 
\begin{equation}\label{eipuf9f}
\alpha([(x^{(n+m)})_{m\geq 0}])=\tx^{(n)}.
\end{equation}
Par suite, on a 
\begin{equation}\label{eipuf9g}
\alpha(\rV^n[x])=p^n\alpha([(x^{(n+m)})_{m\geq 0}])=p^n\tx^{(n)};
\end{equation}
d'où l'unicité de $\alpha$.  

Pour construire $\alpha$, nous procédons comme dans \ref{eip1}. Pour tout entier $n\geq 1$, l'homomorphisme d'anneaux \eqref{notconv1a}
\begin{equation}\label{eipuf9h}
\Phi_{n+r+1}\colon 
\begin{array}[t]{clcr}
\rW_{n+r+1}(D/p^{n+r}D)&\rightarrow& D/p^{n+r}D\\
(x_0,\dots,x_{n+r})&\mapsto& x_0^{p^{n+r}}+p x_1^{p^{n+r-1}}+\dots+p^{n+r}x_{n+r}
\end{array}
\end{equation}
induit un homomorphisme d'anneaux
\begin{equation}\label{eipuf9i}
\Phi'_{n+r+1}\colon 
\begin{array}[t]{clcr}
\rW_n(D/p^{n+r}D)&\rightarrow& D/p^nD\\
(x_0,\dots,x_{n-1})&\mapsto&x_0^{p^{n+r}}+p x_1^{p^{n+r-1}}+\dots+p^{n-1}x_{n-1}^{p^{r+1}}.
\end{array}
\end{equation}
Pour tout $0\leq i\leq n-1$, on a  \eqref{eipuf9b}
\begin{equation}\label{eipuf9j}
p^i(I_D+pD)^{p^{n+r-i}}\subset p^{p^{n+r-i}+i-r}D \subset p^nD.
\end{equation}
L'homomorphisme $\Phi'_{n+r+1}$ s'annule donc sur 
\begin{equation}\label{eipuf9k}
\rW_n((I_D+pD)/p^{n+r}D)=\ker(\rW_n(D/p^{n+r}D)\rightarrow \rW_n(A/pA))
\end{equation}
et induit à son tour un homomorphisme d'anneaux
\begin{equation}\label{eipuf9l}
\alpha_n\colon \rW_n(A/pA)\rightarrow D/p^nD.
\end{equation}
Il résulte aussitôt de la définition que le diagramme 
\begin{equation}\label{eipuf9m}
\xymatrix{
{\rW_{n+1}(A/pA)}\ar[r]^-(0.4){\alpha_{n+1}}\ar[d]_{\rR\rF}&{D/p^{n+1}D}\ar[d]\\
{\rW_{n}(A/pA)}\ar[r]^-(0.4){\alpha_n}&{D/p^nD}}
\end{equation}
où $\rR$ est le morphisme de restriction \eqref{notconv1b}, 
$\rF$ est le Frobenius \eqref{notconv1d} et la flèche non libellée est l'homomorphisme canonique, est commutatif. 
Les homomorphismes $\alpha_n$ définissent par passage à la limite projective l'homomorphisme recherché \eqref{eipuf9d}.

Notons $\alpha_\Lambda\colon \rW(A^\flat)\otimes_W\Lambda \rightarrow D$ l'homomorphisme déduit de $\alpha$.
Comme $D$ est séparé et complet pour la topologie $p$-adique, $\alpha_\Lambda$
induit un homomorphisme $\hcA_r(A/\Lambda) \rightarrow D$ vérifiant $\htheta_{\Lambda, r} =\theta_D\circ \alpha_\Lambda$.
L'homomorphisme $\hcA_r(A/\Lambda)\rightarrow D$ de $\bE^p_{r-1}(\hA/\Lambda)$ ainsi défini est clairement unique; d'où la proposition.

\begin{lem}\label{eipuf5}
Soient $A$ un anneau, $t\in A$, $M, M', M''$ trois $A$-modules complets et séparés pour la topologie $t$-adique, 
$u\colon M'\rightarrow M$, $v\colon M\rightarrow M''$ deux morphismes $A$-linéaires tels que $v\circ u=0$. 
On note $F\mapsto F_n$ le foncteur de réduction modulo $t^n$ sur la catégorie des $R$-modules. 
On suppose que $t$ n'est diviseur de zéro d'aucun des $A$-modules $M$, $M'$ et $M''$ et que la suite de $A_1$-modules
\begin{equation}\label{eipuf5a}
0\longrightarrow M'_1\stackrel{u_1}{\longrightarrow} M_1\stackrel{v_1}{\longrightarrow} M''_1\longrightarrow 0
\end{equation}
est exacte. Alors, la suite 
\begin{equation}\label{eipuf5b}
0\longrightarrow M'\stackrel{u}{\longrightarrow} M\stackrel{v}{\longrightarrow} M''\longrightarrow 0
\end{equation}
est exacte.
\end{lem}

En effet, les colonnes du diagramme 
\begin{equation}\label{eipuf5c}
\xymatrix{
& 0\ar[d]& 0\ar[d]& 0\ar[d]& \\
0\ar[r]& M'_n\ar[r]^{u_n}\ar[d]^{\cdot t}&M_n\ar[r]^{v_n}\ar[d]^{\cdot t}&M'_n\ar[r]\ar[d]^{\cdot t}&0\\
0\ar[r]& M'_{n+1}\ar[r]^{u_{n+1}}\ar[d]&M_{n+1}\ar[r]^{v_{n+1}}\ar[d]&M'_{n+1}\ar[r]\ar[d]&0\\
0\ar[r]& M'_1\ar[r]^{u_1}\ar[d]&M_1\ar[r]^{v_1}\ar[d]&M'_1\ar[r]\ar[d]&0\\
&0 &0 &0 &}
\end{equation}
sont exactes. On en déduit par récurrence que pour tout $n\geq 1$, la suite de $A_n$-modules
\begin{equation}\label{eipuf5d}
0\longrightarrow M'_n\stackrel{u_n}{\longrightarrow} M_n\stackrel{v_n}{\longrightarrow} M''_n\longrightarrow 0
\end{equation}
est exacte.  La proposition s'ensuit en vertu de (\cite{ega3} 0.13.2.2).

\begin{prop}[\cite{agt} II.9.5; \cite{tsuji1} A.1.1 et A.2.2]\label{eipuf6}
Soient $K$ un corps de valuation discrète complet de 
caractéristique $0$, à corps résiduel parfait $k$ de caractéristique $p>0$, 
$\co_K$ l'anneau de valuation de $K$, $\pi$ une uniformisante de $\co_K$, $W$ l'anneau des vecteurs de Witt à coefficients dans $k$, 
$A$ une $\co_K$-algèbre.  Supposons les conditions suivantes vérifiées:
\begin{itemize}
\item[{\rm (i)}] $A$ est $\co_K$-plat.
\item[{\rm (ii)}] $A$ est intégralement clos dans $A[\frac 1 p]$.
\item[{\rm (iii)}] L'endomorphisme de Frobenius de $A/pA$ est surjectif. 
\item[{\rm (iv)}] Il existe une suite $(\pi_n)_{n\geq 0}$ d'éléments de $A$
telle que $\pi_0=\pi$ et $\pi_{n+1}^p=\pi_n$ pour tout $n\geq 0$.  
\end{itemize}
Reprenons les notations de \ref{eipuf2} avec $\Lambda=\co_K$. 
On désigne par $\upi$ l'élément de $A^\flat$ induit par la suite $(\pi_n)_{n\geq 0}$ et on pose  
\begin{equation}\label{eipuf6a}
\xi_\pi=[\upi]-\pi \in \rW_{\co_K}(A^\flat),
\end{equation}
où $[\ ]$ est le représentant multiplicatif. Alors, la suite
\begin{equation}\label{eipuf6b}
\xymatrix{
0\ar[r]&{\rW_{\co_K}(A^\flat)}\ar[r]^-(0.5){\cdot \xi_\pi}&{\rW_{\co_K}(A^\flat)}
\ar[r]^-(0.5){\theta_{\co_K}}&\hA\ar[r]& 0}
\end{equation}
est exacte.
\end{prop}

En effet, on a clairement $\theta_{\co_K}(\xi_\pi)=0$. 
L'algèbre $\rW_{\co_K}(A^\flat)$ étant libre de type fini sur $\rW(A^\flat)$, elle est complète et séparée pour la topologie $p$-adique, 
ou ce qui revient au-même pour la topologie $\pi$-adique. Par ailleurs, $\rW_{\co_K}(A^\flat)$ étant $\co_K$-plat, 
$\pi$ n'est pas diviseur de zéro dans  $\rW_{\co_K}(A^\flat)$. De même, pour tout $n\geq 1$, la suite 
\begin{equation}\label{eipuf6c}
0\longrightarrow A/\pi^n A\stackrel{\cdot \pi}{\longrightarrow} A/\pi^{n+1}A \longrightarrow A/\pi A\rightarrow 0
\end{equation}
est exacte, et par suite $\pi$ n'est pas diviseur de zéro dans $\hA$ (\cite{ega3} 0.13.2.2). 
Compte tenu de \ref{eipuf5}, il suffit donc de montrer que la suite 
\begin{equation}\label{eipuf6d}
\xymatrix{
0\ar[r]&{A^\flat}\ar[r]^{\cdot \upi}&{A^\flat} \ar[r]&A/\pi A\ar[r]& 0}
\end{equation}
où la troisième flèche est induite par la projection canonique $A^\flat \rightarrow A/pA$ sur le premier facteur ({\em i.e.}, d'indice $0$), est exacte.

Pour tout $n\geq 0$, il existe $q_n\in A$ tel que $p=\pi_n q_n$. Comme $\pi_n$ n'est pas diviseur de zéro dans $A$,
on a $q_{n+1}^p=p^{p-1}q_n$.  
Soit $y\in A^\flat$ tel que $\upi \cdot y=0$. Pour tout $n\geq 0$, notons $y_n$ l'image de $y$
par la projection canonique $A^\flat\rightarrow A/pA$ sur le $(n+1)$-ième facteur ({\em i.e.}, d'indice $n$), 
et soit $\ty_n$ un relèvement de $y_n$ dans $A$. On a alors $\pi_n \ty_n\in pA$ et par suite $\ty_n\in q_nA$ et 
\begin{equation}\label{eipuf6e}
y_n=y_{n+1}^p (=\ty_{n+1}^p \mod p A)= 0.
\end{equation}
On en déduit que $\upi$ n'est pas diviseur de zéro dans $A^\flat$. 

Par ailleurs, soient $n\geq 0$, $x\in A$ tels que $x^{p^n}\in \pi A$.  Comme $\pi_n^{-1}\in p^{-1} A$, on voit d'après (i) et (ii) que $x\in \pi_nA$. 
L'endomorphisme de Frobenius $\varphi$ de $A/pA$ étant surjectif d'après (iii), on en déduit que la suite 
\begin{equation}\label{eipuf6f}
\xymatrix{
0\ar[r]&{A/q_n A}\ar[r]^{\cdot \pi_n}&{A/pA} \ar[r]^{\lambda\circ \varphi^n}&A/\pi A\ar[r]& 0},
\end{equation}
où $\lambda\colon A/pA\rightarrow A/\pi A$ est le morphisme canonique, est exacte. Par suite, la suite \eqref{eipuf6d} est exacte d'après (\cite{ega3} 0.13.2.2). 

\begin{cor}\label{eipuf10}
Sous les hypothèses de \ref{eipuf6}, pour tout entier $r\geq 1$, 
les anneaux $\cA_r(A)$ \eqref{eipo3e} et $\cA_r(A/\co_K)$ \eqref{eipuf2f} sont complets et séparés pour les topologies $p$-adiques. 
\end{cor}

Montrons d'abord que $\hA$ est $\co_K$-plat. D'après (\cite{ac} chap.~III §2.11 prop.~14 et cor.~1), pour tout $n\geq 0$,  on a
\begin{equation}\label{eipuf10a}
\hA/p^n\hA\simeq A/p^nA.
\end{equation} 
Soient $x\in \hA$ tel que $px=0$, $\ox$ la classe de $x$ dans $\hA/p^n\hA$ $(n\geq 1)$. 
Comme $A$ est $\co_\oK$-plat, il résulte de \eqref{eipuf10a} que $\ox\in p^{n-1}\hA/p^n\hA$.
On en déduit que $x\in \cap_{n\geq 0}p^n\hA=\{0\}$. 
Par suite, $p$ n'est pas diviseur de zéro dans $\hA$, et donc $\hA$ est $\co_K$-plat. 

D'autre part, pour tout entier $r\geq 1$, la suite \eqref{eipuf6b} induit une suite exacte
\begin{equation}\label{eipuf10b}
\xymatrix{
0\ar[r]&{\cA_r(A/\co_K)}\ar[r]^-(0.5){\cdot \xi_\pi}&{\cA_{r+1}(A/\co_K)}
\ar[rr]^-(0.5){\theta_{\co_K,r+1}}&&\hA\ar[r]& 0}
\end{equation}
On en déduit par récurrence sur $r$ que $\cA_r(A/\co_K)$ est complet et séparé pour la topologie $p$-adique. 
Remplaçant $\co_K$ par $W$, on obtient que $\cA_r(A)$ est complet et séparé pour la topologie $p$-adique.

\section{\texorpdfstring{\'Epaississements infinitésimaux logarithmiques}
{Epaississements infinitésimaux logarithmiques}}\label{epinflog}

On renvoie à (\cite{agt} II.5) pour un lexique de géométrie logarithmique.

\subsection{}\label{epinflog1}
Soient $(X,\cM_X)$ un schéma logarithmique, $M$ un monoïde, $u\colon M\rightarrow \Gamma(X,\cM_X)$ un homomorphisme,  
$p$ un nombre premier. Supposons que le schéma $X$ soit affine d'anneau une $\mZ_{(p)}$-algèbre $A$. On note $\hA$ le séparé complété $p$-adique de $A$.
Reprenons les notations de \ref{eipo3} pour $A$. Considérons le système projectif de monoïdes multiplicatifs $(A)_{n\in \mN}$
où les morphismes de transition sont tous égaux à l'élévation à la puissance $p$-ième.
On désigne par $Q$ le produit fibré du diagramme d'homomorphismes de monoïdes 
\begin{equation}\label{epinflog1a}
\xymatrix{
&{M}\ar[d]\\
{\underset{\underset{\mN}{\longleftarrow}}{\lim}\ A}\ar[r]&A}
\end{equation}
où  la flèche horizontale est la projection sur la première composante ({\em i.e.}, d'indice $0$)
et la flèche verticale est induite par $u$. On note $\tau$ l'homomorphisme composé
\begin{equation}\label{epinflog1b}
\tau\colon Q\longrightarrow
\underset{\underset{x\mapsto x^p}{\longleftarrow}}{\lim}\ A \longrightarrow A^\flat \stackrel{[\ ]}{\longrightarrow} \rW(A^\flat),
\end{equation} 
où la première et la deuxième flèches sont les homomorphismes canoniques \eqref{eipo3a} et $[\ ]$ est le représentant multiplicatif. 
Il résulte aussitôt des définitions que le diagramme 
\begin{equation}\label{epinflog1c}
\xymatrix{
Q\ar[r]\ar[d]_{\tau}&{M}\ar[d]\\
{\rW(A^\flat)}\ar[r]^-(0.4){\theta}&{\hA}}
\end{equation}
où $\theta$ est l'homomorphisme canonique \eqref{eipo3d} et les flèches non libellées sont les morphismes canoniques, est commutatif. 

On pose $\hX=\Spec(\hA)$ que l'on munit de la structure logarithmique $\cM_{\hX}$ image inverse de $\cM_X$.  
On munit $\Spec(\rW(A^\flat))$ de la structure logarithmique $\cQ$ associée à la structure pré-logarithmique définie par 
l'homomorphisme $\tau$ \eqref{epinflog1b}. D'après \eqref{epinflog1c}, $\theta$ induit un morphisme de schémas logarithmiques
\begin{equation}\label{epinflog1d}
(\hX,\cM_\hX)\rightarrow (\Spec(\rW(A^\flat)),\cQ).
\end{equation}

\begin{prop}[\cite{agt} II.9.7]\label{epinflog2}
Conservons les hypothèses de \ref{epinflog1}, notons, de plus, $X^\circ$ l'ouvert maximal de $X$ où la structure logarithmique 
$\cM_X$ est triviale et supposons les conditions suivantes remplies~:
\begin{itemize}
\item[{\rm (a)}] $A$ est intègre et normal. 
\item[{\rm (b)}] $X^\circ$ est un $\mQ$-schéma non-vide et simplement connexe. 
\item[{\rm (c)}] $M$ est intègre et il existe un monoïde fin et saturé $M'$ et un homomorphisme $v\colon M'\rightarrow M$
tels que l'homomorphisme induit  $M'\rightarrow M/M^\times$ soit un isomorphisme.  
\end{itemize}
Alors,
\begin{itemize}
\item[{\rm (i)}] Le monoïde $Q$ est intègre et le groupe $M'^\gp$ est libre. 
\item[{\rm (ii)}] On peut compléter le diagramme \eqref{epinflog1a} en un diagramme commutatif 
\begin{equation}\label{epinflog2a}
\xymatrix{
{M'}\ar[r]^v\ar[d]_-(0.4)w&{M}\ar[d]\\
{\underset{\underset{x\mapsto x^p}{\longleftarrow}}{\lim}\ A}\ar[r]&A}
\end{equation}
Notons $\beta\colon M'\rightarrow Q$ l'homomorphisme induit. 
\item[{\rm (iii)}] La structure logarithmique $\cQ$ sur $\Spec(\rW(A^\flat))$ est associée à la structure pré-logari\-thmique 
définie par l'homomorphisme composé
\begin{equation}\label{epinflog2b}
M'\stackrel{\beta}{\longrightarrow} Q \stackrel{\tau}{\longrightarrow} \rW(A^\flat).
\end{equation}
En particulier, le schéma logarithmique $(\Spec(\rW(A^\flat)),\cQ)$ est fin et saturé.
\item[{\rm (iv)}] Si de plus, l'homomorphisme composé $u\circ v\colon M'\rightarrow \Gamma(X,\cM_X)$ est 
une carte pour $(X,\cM_X)$, alors  le morphisme \eqref{epinflog1d} est strict.
\end{itemize}
\end{prop}

\subsection{}\label{epinflog5}
Dans la suite de cette section, $K$ désigne un corps de valuation discrète complet de 
caractéristique $0$, à corps résiduel parfait $k$ de caractéristique $p>0$, $\co_K$ l'anneau de valuation de $K$, $\pi$ une uniformisante de $\co_K$,
$W$ l'anneau des vecteurs de Witt à coefficients dans $k$ (relatif à $p$) et $A$ une $\co_K$-algèbre vérifiant les conditions suivantes:
\begin{itemize}
\item[(L$_1$)] $A$ est $\co_K$-plat.
\item[(L$_2$)] $A$ est intégralement clos dans $A[\frac 1 p]$.
\item[(L$_3$)] L'endomorphisme de Frobenius de $A/pA$ est surjectif. 
\item[(L$_4$)] Il existe une suite $(\pi_n)_{n\geq 0}$ d'éléments de $A$
telle que $\pi_0=\pi$ et $\pi_{n+1}^p=\pi_n$ pour tout $n\geq 0$.  
\end{itemize}

Reprenons les notations de \ref{eipuf2} avec $\Lambda=\co_K$. On désigne par $\upi$ l'élément de $A^\flat$ induit par la suite $(\pi_n)_{n\geq 0}$ et on pose  
\begin{equation}\label{epinflog5b}
\xi_\pi=[\upi]-\pi \in \rW_{\co_K}(A^\flat),
\end{equation}
où $[\ ]$ est le représentant multiplicatif. En vertu de \ref{eipuf6}, la suite
\begin{equation}\label{epinflog5c}
\xymatrix{
0\ar[r]&{\rW_{\co_K}(A^\flat)}\ar[r]^-(0.5){\cdot \xi_\pi}&{\rW_{\co_K}(A^\flat)}
\ar[r]^-(0.5){\theta_{\co_K}}&\hA\ar[r]& 0}
\end{equation}
est exacte. En particulier, on a $\ker(\theta)\subset \xi_\pi \rW_{\co_K}(A^\flat)$. 

Il résulte de (L$_1$) que $\hA$ est $\co_K$-plat (cf. la preuve de \ref{eipuf10}). On pose 
\begin{equation}\label{epinflog5d}
\rW_{K}(A^\flat)=\rW(A^\flat)\otimes_WK
\end{equation}
et on note $\theta_K\colon \rW_{K}(A^\flat)\rightarrow \hA[\frac 1 p]$ l'homomorphisme induit par $\theta$ \eqref{eipo3d}.  
On désigne par $\rW^{\ast}_{\co_K}(A^\flat)$ la sous-$\rW_{\co_K}(A^\flat)$-algèbre de $\rW_K(A^\flat)$ engendrée par $[\upi]/\pi$ et on pose
\begin{equation}\label{epinflog5e}
\xi^{\ast}_\pi=\frac{\xi_\pi}{\pi}=\frac{[\upi]}{\pi}-1\in \rW^{\ast}_{\co_K}(A^\flat).
\end{equation}
On prendra garde que $\rW^{\ast}_{\co_K}(A^\flat)$ dépend de la suite $(\pi_n)_{n\geq 0}$. 
L'homomorphisme $\theta_K$ induit un homomorphisme 
\begin{equation}\label{epinflog5g}
\theta^{\ast}_{\co_K}\colon \rW^{\ast}_{\co_K}(A^\flat)\rightarrow \hA
\end{equation}
tel que $\theta^{\ast}_{\co_K}(\xi^{\ast}_\pi)=0$.

\begin{lem}\label{epinflog6}
{\rm (i)}\ Pour tout entier $r\geq 0$, on a 
\begin{equation}
(\xi^{\ast}_\pi)^r \rW^{\ast}_{\co_K}(A^\flat) = \rW^{\ast}_{\co_K}(A^\flat) \cap ((\xi^{\ast}_\pi)^r \rW_K(A^\flat)).
\end{equation}

{\rm (ii)}\  La suite 
\begin{equation}\label{epinflog6a}
\xymatrix{
0\ar[r]&{\rW^{\ast}_{\co_K}(A^\flat)}\ar[r]^-(0.5){\cdot \xi^{\ast}_\pi}&{\rW^{\ast}_{\co_K}(A^\flat)}
\ar[r]^-(0.5){\theta^{\ast}_{\co_K}}&\hA\ar[r]& 0}
\end{equation}
est exacte. 
\end{lem}

(i) On notera d'abord que $\rW(A^\flat)$ est $W$-plat puisque $A^\flat$ est parfait. 
Par suite, $\rW_{\co_K}(A^\flat)$ s'identifie à une sous-algèbre de $\rW_K(A^\flat)$. 
Procédons par récurrence sur $r$. L'assertion est évidente pour $r=0$. 
Supposons $r\geq 1$ et l'assertion établie pour $r-1$. Soit $x\in \rW^{\ast}_{\co_K}(A^\flat) \cap ((\xi^{\ast}_\pi)^r \rW_K(A^\flat))$. 
Par hypothèse de récurrence, il existe $a\in \rW_{\co_K}(A^\flat)$ tel que
\begin{equation}\label{epinflog6b}
x-a(\xi^{\ast}_\pi)^{r-1} \in (\xi^{\ast}_\pi)^r \rW^{\ast}_{\co_K}(A^\flat). 
\end{equation}
Comme $\xi^{\ast}_\pi$ n'est pas diviseur de zéro dans $\rW_K(A^\flat)$ \eqref{epinflog5c}, on en déduit que 
\begin{equation}\label{epinflog6c}
a\in \rW_{\co_K}(A^\flat)\cap \xi^{\ast}_\pi \rW_K(A^\flat).
\end{equation}
Par suite, $\theta_{\co_K}(a)=0$ et il existe donc $b\in  \rW_{\co_K}(A^\flat)$ tel que $a=\xi_\pi b$ en vertu de \eqref{epinflog5c}. 
Il s'ensuit que $x\in (\xi^{\ast}_\pi)^r \rW^{\ast}_{\co_K}(A^\flat)$, d'où la proposition. 

(ii) Cela résulte de (i) et de \eqref{epinflog5c}.

\subsection{}\label{epinflog7}
Pour tout entier $r\geq 1$, on pose 
\begin{equation}\label{epinflog7a}
\cA^{\ast}_r(A/\co_K)=\rW^{\ast}_{\co_K}(A^\flat)/(\xi^{\ast}_\pi)^r \rW^{\ast}_{\co_K}(A^\flat)
\end{equation}
et on note $\theta^{\ast}_{\co_K,r}\colon \cA^{\ast}_r(A/\co_K)\rightarrow \hA$ l'homomorphisme induit par $\theta^{\ast}_{\co_K}$. 
Il résulte de \ref{epinflog6}(ii) que $\cA^{\ast}_r(A/\co_K)$ est complet et séparé pour la topologie $p$-adique. 
L'homomorphisme canonique $\rW(A^\flat)\rightarrow \rW^{\ast}_{\co_K}(A^\flat)$ induit un homomorphisme \eqref{eipo3e}
\begin{equation}\label{epinflog7b}
\cA_r(A)\rightarrow \cA^{\ast}_r(A/\co_K).
\end{equation}

\subsection{}\label{epinflog3}
On pose $S=\Spec(\co_K)$ et on note $s$ (resp.  $\eta$) le point fermé (resp.  générique) de $S$ \eqref{epinflog5}.
On munit $S$ de la structure logarithmique $\cM_S$ définie par son point fermé, 
autrement dit, $\cM_S=u_*(\co_\eta^\times)\cap \co_S$, où $u\colon \eta\rightarrow S$ est l'injection canonique. 
On désigne par $\iota\colon \mN\rightarrow \Gamma(S,\cM_S)$
l'homomorphisme défini par $\iota(1)=\pi$, qui est une carte pour $(S,\cM_S)$. 
On pose $X=\Spec(A)$ que l'on suppose muni d'une structure logarithmique $\cM_X$. On se donne aussi
un morphisme $f\colon (X,\cM_X)\rightarrow (S,\cM_S)$, un monoïde $M$ et deux homomorphismes
$\mN\stackrel{\gamma}{\rightarrow}M\stackrel{u}{\rightarrow} \Gamma(X,\cM_X)$ qui s'insèrent dans un diagramme commutatif
\begin{equation}\label{epinflog3a}
\xymatrix{
\mN\ar[r]^-(0.5){\iota}\ar[d]_{\gamma}&{\Gamma(S,\cM_S)}\ar[d]^{f^\flat}\\
M\ar[r]^-(0.5)u&{\Gamma(X,\cM_X)}}
\end{equation}
où $f^\flat$ est l'homomorphisme induit par $f$. 
Reprenons les notations de \ref{epinflog1} pour le schéma logarithmique affine $(X,\cM_X)$ et l'homomorphisme $u$. 

Pour tout entier $r\geq 1$, on pose \eqref{epinflog7a}
\begin{equation}\label{epinflog3b}
\cA^{\ast}_r(X/S)=\Spec(\cA^{\ast}_r(A/\co_K)).
\end{equation}
On le munit de la structure logarithmique $\cM_{\cA^{\ast}_r(X/S)}$ associée 
à la structure pré-logarithmique définie par l'homomorphisme 
\begin{equation}\label{epinflog3c}
Q\rightarrow \cA^{\ast}_r(A/\co_K)
\end{equation} 
induit par $\tau$ \eqref{epinflog1b}. D'après \eqref{epinflog1c}, l'homomorphisme $\theta^{\ast}_{\co_K,r}$ induit un morphisme 
\begin{equation}\label{epinflog3d}
(\hX,\cM_\hX)\rightarrow (\cA^{\ast}_r(X/S),\cM_{\cA^{\ast}_r(X/S)}).
\end{equation}
Avec les notations de \ref{epinflog5}(L$_4$), on désigne par $\tpi$ l'élément de $Q$ défini par ses projections \eqref{epinflog1a}
\begin{equation}\label{epinflog3e}
(\pi_n)_{n\geq 0}\in \underset{\underset{\mN}{\longleftarrow}}{\lim}\ A \ \ \ {\rm et}\ \ \ 
\gamma(1) \in M.
\end{equation}
L'élément $(\xi^{\ast}_\pi+1)$ étant inversible dans $\cA^{\ast}_r(A/\co_K)$, l'homomorphisme 
\begin{equation}\label{epinflog3f}
\mN\rightarrow \Gamma(\cA^{\ast}_r(X/S), \cM_{\cA^{\ast}_r(X/S)}), \ \ \ 1\mapsto (\xi^{\ast}_\pi+1)^{-1}\tpi
\end{equation}
induit un morphisme de schémas logarithmiques 
\begin{equation}\label{epinflog3g}
\pr_1\colon (\cA^{\ast}_r(X/S),\cM_{\cA^{\ast}_r(X/S)})\rightarrow (S,\cM_S).
\end{equation}

On pose \eqref{eipo3e}
\begin{equation}\label{epinflog3h}
\cA_r(X)=\Spec(\cA_r(A)),
\end{equation}
que l'on munit de la structure logarithmique $\cM_{\cA_r(X)}$ image inverse de la structure logarithmique $\cQ$ sur $\Spec(\rW(A^\flat))$ (cf. \ref{epinflog1}).
On a clairement un morphisme de schémas logarithmiques 
\begin{equation}\label{epinflog3i}
\pr_2\colon (\cA^{\ast}_r(X/S),\cM_{\cA^{\ast}_r(X/S)})\rightarrow (\cA_r(X),\cM_{\cA_r(X)}).
\end{equation}

\begin{prop}\label{epinflog4}
Conservons les hypothèses et notations de \ref{epinflog5} et \ref{epinflog3}, notons, de plus, $X^\circ$ l'ouvert maximal de $X$ où la structure logarithmique 
$\cM_X$ est triviale et supposons les conditions suivantes remplies~:
\begin{itemize}
\item[{\rm (L$_5$)}] $A$ est intègre et normal. 
\item[{\rm (L$_6$)}] $X^\circ$ est un $\eta$-schéma non-vide et simplement connexe. 
\item[{\rm (L$_7$)}] $M$ est intègre et il existe un monoïde fin et saturé $M'$ et un homomorphisme $v\colon M'\rightarrow M$
tels que l'homomorphisme induit  $M'\rightarrow M/M^\times$ soit un isomorphisme et que 
l'homomorphisme composé $u\circ v\colon M'\rightarrow \Gamma(X,\cM_X)$ soit une carte pour $(X,\cM_X)$.
\end{itemize}
Alors, 
\begin{itemize}
\item[{\rm (i)}] Le monoïde $Q$ est intègre et le groupe $M'^\gp$ est libre. 
\item[{\rm (ii)}] On peut compléter le diagramme \eqref{epinflog1a} en un diagramme commutatif 
\begin{equation}\label{epinflog4a}
\xymatrix{
{M'}\ar[r]^v\ar[d]_-(0.4)w&{M}\ar[d]\\
{\underset{\underset{x\mapsto x^p}{\longleftarrow}}{\lim}\ A}\ar[r]&A}
\end{equation}
Notons $\beta\colon M'\rightarrow Q$ l'homomorphisme induit. 
\item[{\rm (iii)}] Pour tout entier $r\geq 1$, la structure logarithmique $\cM_{\cA^{\ast}_r(X/S)}$ sur $\cA^{\ast}_r(X/S)$ est associée à la structure pré-logarithmique 
définie par l'homomorphisme composé
\begin{equation}\label{epinflog4b}
M'\stackrel{\beta}{\longrightarrow} Q \longrightarrow \cA^{\ast}_r(A/\co_K),
\end{equation}
où la seconde flèche est l'homomorphisme \eqref{epinflog3c}
En particulier, le schéma logarithmique $(\cA^{\ast}_r(X/S),\cM_{\cA^{\ast}_r(X/S)})$ est fin et saturé 
et le morphisme \eqref{epinflog3d} est une immersion fermée stricte.
\end{itemize}
\end{prop}

Les deux premières propositions sont mentionnées à titre de rappel \eqref{epinflog2}. 
Soit $G$ l'image inverse de $M^\times$ par  l'homomorphisme canonique $Q\rightarrow M$. 
Il résulte aussitôt de la définition \eqref{epinflog1a} que $G$ est un sous-groupe de $Q$. 
D'autre part, l'homomorphisme composé $M'\rightarrow Q/G\rightarrow M/M^\times$,
où la première flèche est déduite de $\beta$, est un isomorphisme. 
Par suite, $M'\rightarrow Q/G$ est un isomorphisme. La proposition (iii) résulte alors de (\cite{tsuji1} 1.3.1).

\section{\texorpdfstring{Modules de Higgs et $\lambda$-connexions}{Modules de Higgs et lambda-connexions}}\label{MH}

\begin{defi}\label{MH1}
Soient $(X,A)$ un topos annelé, $E$ un $A$-module. 
\begin{itemize}
\item[(i)]
On appelle {\em $A$-module de Higgs à coefficients dans $E$}
un couple $(M,\theta)$ formé d'un $A$-module $M$ et d'un morphisme $A$-linéaire 
\begin{equation}\label{MH1a}
\theta\colon M\rightarrow M\otimes_AE
\end{equation}
tel que $\theta\wedge \theta=0$. On dit alors que $\theta$ est un {\em $A$-champ de Higgs} sur $M$ 
à coefficients dans $E$. 
\item[(ii)] Si $(M,\theta)$ et $(M',\theta')$ sont deux $A$-modules de Higgs,
un morphisme de $(M,\theta)$ dans $(M',\theta')$ est un morphisme $A$-linéaire 
$u\colon M\rightarrow M'$ tel que $(u\otimes\id_E)\circ \theta=\theta'\circ u$. 
\end{itemize}
\end{defi}

Les $A$-modules de Higgs à coefficients dans $E$ forment une catégorie que l'on note 
$\bMH(A,E)$. On peut compléter la terminologie et faire les remarques suivantes.

\addtocounter{subsubsection}{1}
\addtocounter{equation}{1}

\subsubsection{}\label{MH2a}
Soit $(M,\theta)$ un $A$-module de Higgs à coefficients dans $E$. Pour tout $i\geq 1$, on désigne par
\begin{equation}\label{MH2b}
\theta^i\colon M\otimes_A \wedge^iE \rightarrow M\otimes_A \wedge^{i+1}E
\end{equation}
le morphisme $A$-linéaire défini pour toutes sections locales 
$m$ de $M$ et $\omega$ de $\wedge^iE$ par $\theta^i(m\otimes \omega)=\theta(m)\wedge \omega$ \eqref{notconv9}.
On a $\theta^{i+1}\circ \theta^i=0$. On appelle complexe de {\em Dolbeault} de $(M,\theta)$
et l'on note $\mK^\bullet(M,\theta)$ le complexe de cochaînes de $A$-modules 
\begin{equation}\label{MH2c}
M\stackrel{\theta}{\longrightarrow}M\otimes_AE\stackrel{\theta^1}{\longrightarrow} M\otimes_A\wedge^2E \dots,
\end{equation}
où $M$ est placé en degré $0$ et les différentielles sont de degré $1$.

\addtocounter{subsubsection}{2}
\addtocounter{equation}{1}

\subsubsection{}\label{MH2j}
Soit $(M,\theta)$ un $A$-module de Higgs à coefficients dans $E$ tel que $M$ soit un $A$-module
localement libre de type fini. Considérons, pour un entier $i\geq 1$, le morphisme composé 
\begin{equation}\label{MH2k}
\xymatrix{
{\wedge^iM}\ar[r]^-(0.5){\wedge^i\theta}&{\wedge^i(M\otimes_AE)}\ar[r]&
{\wedge^iM\otimes_A\rS^iE}},
\end{equation}
où la seconde flèche est le morphisme canonique (\cite{illusie1} V 4.5).  
On appelle {\em $i$-ième invariant caractéristique} de $\theta$ et l'on note $\lambda_i(\theta)$
la trace du morphisme \eqref{MH2k} vue comme section de 
$\Gamma(X,\rS^iE)$.

\addtocounter{subsubsection}{1}
\addtocounter{equation}{1}

\subsubsection{}\label{MH2d}
Soient $(M,\theta),(M',\theta')$ deux $A$-modules de Higgs à coefficients dans $E$. 
On appelle champ de Higgs {\em total} sur $M\otimes_AM'$  
le morphisme  $A$-linéaire 
\begin{equation}\label{MH2e}
\theta_\tot\colon M\otimes_AM'\rightarrow M\otimes_AM'\otimes_AE
\end{equation} 
défini par 
\begin{equation}\label{MH2f}
\theta_\tot=\theta\otimes \id_{M'}+\id_{M}\otimes \theta'.
\end{equation}
On dit que $(M\otimes_AM',\theta_\tot)$ est le {\em produit tensoriel} de $(M,\theta)$ et $(M',\theta')$.

\addtocounter{subsubsection}{2}
\addtocounter{equation}{1}

\subsubsection{}\label{MH2g}
Supposons $E$ localement libre de type fini sur $A$ et posons $F=\cHom_A(E,A)$. 
Pour tout $A$-module $M$, le morphisme canonique 
\begin{equation}\label{MH2h}
\cEnd_A(M)\otimes_AE\rightarrow \cHom_A(M,M\otimes_AE)
\end{equation}
étant un isomorphisme, la donnée d'un $A$-champ de Higgs $\theta$ sur $M$ 
est équivalente à la donnée d'une structure de $\rS_A(F)$-module sur $M$
compatible avec sa structure de $A$-module \eqref{notconv9}. D'après (\cite{agt} II.2.8.10),  
on a alors un isomorphisme canonique de complexes de $A$-modules 
\begin{equation}\label{MH2i}
\mK^\bullet(M,\theta)\stackrel{\sim}{\rightarrow} \mK_{\rS(F)}^\bullet(M),
\end{equation}
où $\mK_{\rS(F)}^\bullet(M)$ est le complexe de Koszul du $\rS(F)$-module $M$ (\cite{agt} II.2.7.5).

\addtocounter{subsubsection}{2}
\addtocounter{equation}{1}

\subsection{}\label{MH150}
Soient $(X,A)$ un topos annelé, $E$ un $A$-module plat. 
La catégorie $\bMH(A,E)$ des $A$-modules de Higgs à coefficients dans $E$ est naturellement une catégorie abélienne. 
Soit $(M,\theta)$ un objet de cette catégorie.
La donnée d'un sous-objet de $(M,\theta)$ dans $\bMH(A,E)$ est équivalente à la donnée d'un sous-$A$-module $N$ de $M$ tel que $\theta(N)\subset N\otimes_AE$. 
Le $A$-module quotient $M/N$ muni du champ de Higgs induit par $\theta$ est un quotient de $(M,\theta)$ dans $\bMH(A,E)$. 
On dira abusivement que $N$ (resp. $M/N$) est un sous-objet (resp. quotient) de $(M,\theta)$. 

La donnée d'une filtration décroissante de $(M,\theta)$ dans $\bMH(A,E)$, indexée par $\mZ$, est équivalente à la donnée d'une filtration décroissante $(\rF^iM)_{i\in \mZ}$ 
de $M$ par des sous-$A$-modules telle que $\theta(\rF^iM)\subset \rF^iM\otimes_AE$ pour tout $i\in \mZ$.  
On dira abusivement que $\rF^\bullet M$ est une filtration de $(M,\theta)$.  
Pour tout $i\in \mZ$, le champ de Higgs $\theta$ induit sur $\Gr^i_\rF M=\rF^iM/\rF^{i+1}M$ un champ de Higgs, 
et $\Gr_\rF^\bullet M=\oplus_{i\in \mZ}\Gr^i_\rF M$ est le gradué de la filtration $\rF^\bullet M$ de $(M,\theta)$ dans $\bMH(A,E)$.

\'Etant donnés une filtration $\rF^\bullet M$ de $(M,\theta)$ dans $\bMH(A,E)$ et un sous-objet $N$ de $(M,\theta)$, 
$(N\cap \rF^iM)_{i\in \mZ}$ (resp. $((\rF^iM+N)/N)_{i\in \mZ}$) est une filtration du sous-objet $N$ (resp. quotient $M/N$)
de $(M,\theta)$, dite filtration {\em induite} (resp. {\em quotient}). On notera que pour tout $i\in \mZ$, 
on a $(N\cap \rF^iM)\otimes_AE=(N\otimes_AE)\cap (\rF^iM\otimes_AE)$.

\begin{defi}\label{MH151}
Soient $(X,A)$ un topos annelé, $E$ un $A$-module plat. 
On dit qu'un $A$-module de Higgs  $(M,\theta)$ à coefficients dans $E$ est {\em quasi-nilpotent} 
s'il existe une filtration décroissante $(\rF^iM)_{i\in \mZ}$ de $(M,\theta)$ dans $\bMH(A,E)$ \eqref{MH150} et deux entiers $m\leq n$ tels que $\rF^n M=0$, $\rF^m M=M$ et que 
le champ de Higgs de $\Gr_\rF^\bullet M$ soit nul (ou ce qui revient au même que l'on ait $\theta(\rF^iM)\subset \rF^{i+1} M \otimes_AE$ pour tout $i\in \mZ$). 
On dira alors aussi que le champ de Higgs $\theta$ est {\em quasi-nilpotent}, et que la filtration $\rF^\bullet M$ de $(M,\theta)$ est {\em quasi-nilpotente}. 
\end{defi} 

On réserve la terminologie {\em nilpotent} pour une condition un peu plus forte \eqref{MH15}.

\begin{lem}\label{MH141}
Soient $(X,A)$ un topos annelé, $E$ un $A$-module plat, $(M,\theta)$ un $A$-module de Higgs à coefficients dans $E$, 
$N$ (resp. $P$) un sous-objet (resp. quotient) de $(M,\theta)$ dans $\bMH(A,E)$, $(\rF^iM)_{i\in \mZ}$ une filtration décroissante de $(M,\theta)$ dans $\bMH(A,E)$ 
telle que $\rF^n M=0$, $\rF^m M=M$ pour deux entiers $m\leq n$ \eqref{MH150}.
\begin{itemize}
\item[{\rm (i)}] Si $(M,\theta)$ est quasi-nilpotent, il en est de même de $N$ et $P$. De plus, si la filtration $\rF^\bullet M$ de $(M,\theta)$ est quasi-nilpotente, il en est de même 
des filtrations induites sur $N$ et $P$ \eqref{MH150}. 
\item[{\rm (ii)}] Pour que $(M,\theta)$ soit quasi-nilpotent, il faut et il suffit que le gradué $\Gr_\rF^\bullet M$ de la filtration $\rF^\bullet M$ le soit. 
\end{itemize}
\end{lem} 

(i) En effet, notant $\rF^\bullet N$ (resp. $\rF^\bullet P$) la filtration induite par $\rF^\bullet M$ sur $N$ (resp. $P$) \eqref{MH150}, 
$\Gr_\rF^\bullet N$ (resp.  $\Gr_\rF^\bullet P$) s'identifie à un sous-objet (resp. quotient) de $\Gr_\rF^\bullet M$ dans $\bMH(A,E)$. Plus précisément, 
si $P=M/N$, alors on a une suite exacte canonique de $\bMH(A,E)$,
\begin{equation}
0\rightarrow \Gr_\rF^\bullet N \rightarrow \Gr_\rF^\bullet M \rightarrow \Gr_\rF^\bullet P \rightarrow 0.
\end{equation}

(ii) La condition est nécessaire d'après (i). Inversement, si $\Gr_\rF^\bullet M$ est quasi-nilpotent, il en est de même de $\Gr_\rF^i M$ pour tout $i\in \mZ$. 
Il existe donc une filtration $(\rW^i M)_{i\in \mZ}$ de $(M,\theta)$ dans $\bMH(A,E)$
qui raffine la filtration $\rF^\bullet M$ tel que le champ de Higgs du gradué $\Gr_\rW^\bullet M$ soit nul. Par suite, $(M,\theta)$ est quasi-nilpotent; la condition est donc suffisante.

\begin{prop}\label{MH14}
Soient $X$ un schéma, $\cE$ un $\co_X$-module localement libre de type fini, $\cF=\cHom_{\co_X}(\cE,\co_X)$ le dual de $\cE$, 
$\cM$ un $\co_X$-module quasi-cohérent,
$\theta\colon \cM\rightarrow \cM\otimes_{\co_X}\cE$ un $\co_X$-champ de Higgs sur $\cM$ à coefficients dans $\cE$. 
On désigne par $T^*=\Spec(\rS_{\co_X}(\cF))$ 
le fibré vectoriel sur $X$ associé à $\cF$ \eqref{notconv9}, par $\sigma\colon X\rightarrow T^*$ la section nulle, 
par $\cJ$ l'idéal de $\co_{T^*}$ correspondant à l'immersion fermée $\sigma$ 
et par $\cM^+$ le $\co_{T^*}$-module associé à $(\cM,\theta)$ \eqref{MH2g}.
Pour tout ouvert $U$ de $X$ et toute section $d\in \Gamma(U,\cF)$, on note $\theta_d$ l'endomorphisme de $\cM|U$ déduit de $\theta$ et $d$. 
Considérons les conditions suivantes:
\begin{itemize}
\item[{\rm (i)}] Le support de $\cM^+$ est contenu dans $\sigma(X)$. 
\item[{\rm (ii)}] Pour tout ouvert quasi-compact $U$ de $X$ et toute section $s\in \Gamma(U,\cM)$, il existe un entier $n\geq 1$ 
tel que pour toutes sections $d_1,\dots,d_n\in \Gamma(U,\cF)$, $\theta_{d_n}\circ \dots \circ \theta_{d_1}(s)=0$.
\item[{\rm (iii)}] Il existe une filtration décroissante finie $(\cM_i)_{0\leq i\leq n}$ de $\cM$ par des sous-$\co_X$-modules quasi-cohérents 
telle que $\cM_0=\cM$, $\cM_n=0$ et que pour tout $0\leq i\leq n-1$, on ait
\begin{equation}\label{MH14a}
\theta(\cM_i)\subset \cE\otimes_{\co_X}\cM_{i+1}.
\end{equation}
\item[{\rm (iv)}] Il existe un entier $n\geq 1$ tel que $\cJ^n\cM^+=0$. 
\end{itemize} 
Alors, on a {\rm (iv)}$\Rightarrow${\rm (iii)}$\Rightarrow${\rm (ii)}$\Leftrightarrow${\rm (i)}. Si, de plus, $X$ est quasi-compact et 
le $\co_{T^*}$-module $\cM^+$ est de type fini, les quatre conditions sont équivalentes. 
\end{prop}

Par support de $\cM^+$ on sous-entend l'ensemble (pas nécessairement fermé) des points de $T^*$ où la fibre de $\cM^+$ n'est pas nulle (\cite{ega1n} 0.3.1.5).
Il est fermé si le $\co_{T^*}$-module $\cM^+$ est de type fini (\cite{ega1n} 0.5.2.2); c'est le cas si le $\co_X$-module $\cM$ est de type fini. 
On notera que le $\co_{T^*}$-module $\cM^+$ est quasi-cohérent.

Supposons d'abord la condition (iv) remplie. Pour tout $0\leq i\leq n$, notons $\cM_i^+$ 
le plus grand sous-$\co_{T^*}$-module de $\cM^+$ annulé par $\cJ^{n-i}$. 
On a donc $\cM^+_0=\cM^+$, $\cM^+_n=0$ et pour tout $0\leq i\leq n-1$, $\cJ\cM^+_i\subset \cM^+_{i+1}$.
Notant $\pi\colon T^*\rightarrow X$ la projection canonique, la filtration de $\cM$ définie, pour tout $0\leq i\leq n$, par
$\cM_i=\pi_*(\cM^+_i)$ vérifie alors la condition (iii).  

Il est clair qu'on a (iii)$\Rightarrow$(ii). 

L'implication (ii)$\Rightarrow$(i) résulte du fait que 
pour tout ouvert $U$ de $X$, on a $\Gamma(T^*_U, \cM^+)=\Gamma(U,\cM)$, et que l'action de $\Gamma(T^*_U,\co_{T^*})$ sur $\Gamma(T^*_U, \cM^+)$  
est induite par l'application 
\begin{equation}\label{MH14b}
\Gamma(U,\cF)\times \Gamma(U,\cM)\rightarrow \Gamma(U,\cM), \ \ \ (d,m)\mapsto \theta_d(m).
\end{equation}

Soient $U$ un ouvert quasi-compact de $X$, $s\in \Gamma(T^*_U,\cM^+)$ une section de support contenu dans $\sigma(X)$. 
En vertu de (\cite{ega1n} 6.8.4) appliqué au sous-module $\co_{T^*_U} s$ de $\cM^+|T^*_U$, 
il existe un entier $n\geq 1$ tel que $\cJ^n s=0$, d'où l'implication (i)$\Rightarrow$(ii).   

Si $X$ est quasi-compact et si le $\co_{T^*}$-module $\cM^+$ est de type fini, on a (i)$\Rightarrow$(iv) en vertu de (\cite{ega1n} 6.8.4), d'où la proposition.

\begin{defi}\label{MH15}
Soient $X$ un schéma, $\cE$ un $\co_X$-module localement libre de type fini. On dit qu'un $\co_X$-module de Higgs $(\cM,\theta)$ à coefficients dans $\cE$ est 
{\em nilpotent} si le $\co_X$-module $\cM$ est quasi-cohérent et si les conditions équivalentes \ref{MH14}(i)-(ii) sont remplies. 
On dira alors aussi que le champ de Higgs $\theta$ est {\em nilpotent}. 
\end{defi}

Pour qu'un $\co_X$-module de Higgs $(\cM,\theta)$ à coefficients dans $\cE$ soit nilpotent, il faut et il suffit qu'il soit quasi-nilpotent \eqref{MH151} 
et que $(\cM,\theta)$ admette une filtration quasi-nilpotente formée de sous-$\co_X$-modules quasi-cohérents.

\begin{rema}\label{MH140}
Reprenons les hypothèses de \ref{MH14}, supposons, de plus, le schéma $X$ affine d'anneau $R$, et le $\co_X$-module $\cM$ de type fini. 
Posons $E=\Gamma(X,\cE)$, $F=\Gamma(X,\cF)$ et $M=\Gamma(X,\cM)$ qui sont des $R$-modules, 
et $M^+=\Gamma(T^*,\cM^+)$ qui est un $\rS_R(F)$-module. Notons $J$ l'idéal de l'augmentation canonique de  $\rS_R(F)$.
Alors, les propriétés suivantes sont équivalentes:
\begin{itemize}
\item[(i)] Le $\co_X$-module de Higgs $(\cM,\theta)$ est nilpotent. 
\item[(ii)] Pour tout $s\in M$, il existe un entier $n\geq 1$ tel que pour toutes sections $d_1,\dots,d_n\in F$, $\theta_{d_n}\circ \dots \circ \theta_{d_1}(s)=0$.
\item[(iii)] Il existe une filtration décroissante finie $(M_i)_{0\leq i\leq n}$ de $M$ 
telle que $M_0=M$, $M_n=0$ et que pour tout $0\leq i\leq n-1$, on ait
\begin{equation}\label{MH140a}
\theta(M_i)\subset E\otimes_RM_{i+1}.
\end{equation}
\item[{\rm (iv)}] Il existe un entier $n\geq 1$ tel que $J^n M^+=0$. 
\end{itemize} 
\end{rema}

\begin{lem}\label{MH142}
Soient $X$ un schéma, $\cE$ un $\co_X$-module localement libre de type fini, $\cM$ un $\co_X$-module quasi-cohérent,
$\theta$ un $\co_X$-champ de Higgs sur $\cM$ à coefficients dans $\cE$, 
$\cN$ (resp. $\cP$) un sous-objet (resp. quotient) de $(\cM,\theta)$ dans $\bMH(\co_X,\cE)$ tel que le $\co_X$-module sous-jacent soit quasi-cohérent, 
$(\rF^i\cM)_{i\in \mZ}$ une filtration décroissante de $(\cM,\theta)$ dans $\bMH(\co_X,\cE)$ \eqref{MH150}
telle que $\rF^n \cM=0$, $\rF^m \cM=\cM$ pour deux entiers $m\leq n$ et  que pour tout $i\in \mZ$, le $\co_X$-module $\cM_i$ soit quasi-cohérent.
\begin{itemize}
\item[{\rm (i)}] Si $(\cM,\theta)$ est nilpotent, il en est de même de $\cN$ et $\cP$. 
\item[{\rm (ii)}] Pour que $(\cM,\theta)$ soit nilpotent, il faut et il suffit que le gradué $\Gr_\rF^\bullet \cM$ de la filtration $\rF^\bullet \cM$ le soit. 
\end{itemize}
\end{lem} 

Il suffit de calquer la preuve de \ref{MH141}.

\begin{lem}\label{MH17}
Soient $X$ un schéma réduit, $\cE, \cM$ deux $\co_X$-modules localement libres de type fini, 
$\theta$ un $\co_X$-champ de Higgs nilpotent sur $\cM$ à coefficients dans $\cE$. Alors, pour tout $n\geq 1$, le $n$-ième invariant caractéristique $\lambda_n(\theta)$
de $\theta$ est nul \eqref{MH2j}.
\end{lem} 

En effet, notons $\cF=\cHom_{\co_X}(\cE,\co_X)$ le dual de $\cE$. 
Pour tout ouvert $U$ de $X$ et toute section $d\in \Gamma(U,\cF)$, 
on désigne par $\theta_d$ l'endomorphisme de $\cM|U$ déduit de $\theta$ et $d$, et par 
$e_d\colon \rS_{\co_X}(\cE)|U\rightarrow \co_X|U$ l'homomorphisme de $\co_X$-algèbres induit par $d$.  Pour tout entier $n\geq 1$, on a alors  
\begin{equation}
e_d(\lambda_n(\theta))=\Tr(\wedge^n(\theta_d)).
\end{equation}
Comme $\theta_d$ est nilpotent, on en déduit que $e_d(\lambda_n(\theta))$ est nilpotent et par suite nul car $X$ est réduit; d'où la proposition.

\begin{lem}\label{MH16}
Soient $X$ un schéma, $\cE$ un $\co_X$-module localement libre de type fini, $\cM$ un $\co_X$-module quasi-cohérent,
$\theta\colon \cM\rightarrow \cM\otimes_{\co_X}\cE$ un $\co_X$-champ de Higgs sur $\cM$ à coefficients dans $\cE$. 
Alors, il existe un sous-$\co_X$-module quasi-cohérent $\cN$ de $\cM$ vérifiant les propriétés suivantes~:
\begin{itemize}
\item[{\rm (i)}] $\theta(\cN)\subset \cN\otimes _{\co_X}\cE$ et le $\co_X$-module de Higgs $(\cN,\theta)$ à coefficients dans $\cE$ est nilpotent. 
\item[{\rm (ii)}] Pour tout morphisme de $\co_X$-modules de Higgs à coefficients dans $\cE$, $u\colon (\cM',\theta')\rightarrow (\cM,\theta)$
tel que $(\cM',\theta')$ soit nilpotent, on a $u(\cM')\subset \cN$. 
\end{itemize}
On dira que $(\cN,\theta)$ est le {\em sous-$\co_X$-module de Higgs nilpotent maximal} de $(\cM,\theta)$. 
\end{lem}

En effet, notons $\cF=\cHom_{\co_X}(\cE,\co_X)$ le dual de $\cE$, $T^*=\Spec(\rS_{\co_X}(\cF))$ le fibré vectoriel sur $X$ associé à $\cF$, 
$\pi\colon T^*\rightarrow X$ la projection canonique, $\sigma\colon X\rightarrow T^*$ la section nulle, $V$ l'ouvert complémentaire de $\sigma(X)$ dans $T^*$,
$j\colon V\rightarrow T^*$ l'injection canonique et $\cM^+$ le $\co_{T^*}$-module associé à $(\cM,\theta)$ \eqref{MH2g}. 
Le noyau $\cN^+$ du morphisme d'adjonction $\cM^+\rightarrow j_*(j^*\cM^+)$ est le sous-$\co_{T^*}$-module maximal de 
$\cM^+$ de support contenu dans $\sigma(X)$. 
Comme l'immersion $j$ est cohérente ({\em i.e.}, quasi-compacte et quasi-séparée), le $\co_{T^*}$-module $\cN^+$ est quasi-cohérent (\cite{ega1n} 6.7.1). 
On voit aussitôt que le sous-$\co_X$-module $\cN=\pi_*(\cN^+)$ de $\cM$ répond à la question.

\subsection{}\label{MH90}
Soient $(X,A)$ un topos annelé, 
\begin{equation}\label{MH90a}
0\rightarrow L\rightarrow E \rightarrow \uE \rightarrow 0
\end{equation}
une suite exacte de $A$-modules localement libres de type fini. On munit $E$ de la filtration décroissante exhaustive définie par \eqref{MH90a}, 
dont le module gradué associé est $\uE \oplus L$, concentré en degrés $[0,1]$.  
Pour tous entiers $i,j\geq 0$, on désigne par $\rW^i\wedge^jE$ l'image du morphisme canonique 
\begin{equation}\label{MH90b}
\wedge^i L\otimes_A \wedge^{j-i} E\rightarrow \wedge^jE.
\end{equation}
On définit ainsi une filtration décroissante exhaustive $(\rW^i\wedge^\bullet E)_{i\geq 0}$ de l'algèbre extérieure 
$\wedge^\bullet E$ de $E$ par des idéaux, dite {\em filtration de Koszul de $\wedge^\bullet E$ associée à la suite exacte \eqref{MH90a}}. 
On désigne par $\Gr^\bullet_\rW\wedge^\bullet E$ le gradué associé, qui est une algèbre bigraduée. 
L'injection canonique $\uE \oplus L\rightarrow \Gr^\bullet_\rW\wedge^\bullet E$
se prolonge canoniquement en un homomorphisme d'algèbres bigraduées
\begin{equation}\label{MH90c}
\wedge^\bullet(\uE \oplus L)\rightarrow \Gr^\bullet_\rW\wedge^\bullet E.
\end{equation}
C'est un isomorphisme d'après (\cite{illusie1} V 4.1.6). 
Par ailleurs, en vertu de (\cite{alg1-3} III § 7.7 prop.~10), on a un isomorphisme canonique d'algèbres bigraduées
\begin{equation}\label{MH90f}
\wedge^\bullet \uE\  {^g\otimes} \wedge^\bullet L \stackrel{\sim}{\rightarrow}\wedge^\bullet(\uE \oplus L),
\end{equation}
où le symbole ${^g\otimes}$ désigne le produit tensoriel gauche (cf. \cite{alg1-3} III § 4.7 remarques page 49). 
En particulier, pour tout entier $i\geq 0$, on a un isomorphisme canonique
\begin{equation}\label{MH90e}
\Gr^i_\rW\wedge^\bullet E\stackrel{\sim}{\rightarrow}(\wedge^i L)\otimes_A(\wedge^{\bullet -i} \uE). 
\end{equation}

Soit $(M,\theta)$ un $A$-module de Higgs à coefficients dans $E$.  
On désigne par $\utheta\colon M\rightarrow M\otimes_A\uE$ le $A$-champ de Higgs induit par $\theta$,  
et par $\mK^\bullet$ (resp. $\umK^\bullet$) le complexe de Dolbeault de $(M,\theta)$ (resp.  $(M,\utheta)$) \eqref{MH2a}. 
On observera que le complexe $\mK^\bullet$ est un module (à droite) gradué sur l'algèbre graduée $\wedge^\bullet E$ et que sa différentielle de degré 
un est $(\wedge^\bullet E)$-linéaire. Pour tous entiers $i,j\geq 0$, on pose 
\begin{equation}\label{MH90g}
\rW^i\mK^j=M\otimes_A(\rW^i\wedge^jE).
\end{equation}
On définit ainsi une filtration décroissante exhaustive du $\wedge^\bullet E$-module gradué $\mK^\bullet$ 
par des sous-$\wedge^\bullet E$-modules gradués $(\rW^i\mK^\bullet)_{i\geq 0}$, stables par la différentielle de $\mK^\bullet$,
dite {\em filtration de Koszul de $\mK^\bullet$ associée à la suite exacte \eqref{MH90a}};
en particulier, c'est une filtration du complexe $\mK^\bullet$ par des sous-complexes. On désigne par $\Gr_\rW^\bullet\mK^\bullet$ le gradué associé. 

Compte tenu de \eqref{MH90c} et \eqref{MH90f}, pour tout entier $i\geq 0$, on a un isomorphisme canonique de complexes
\begin{equation}\label{MH90h}
\Gr_\rW^i\mK^\bullet\stackrel{\sim}{\rightarrow}\wedge^i L\otimes_A \umK^{\bullet}[-i],
\end{equation}
où les différentielles de $\umK^{\bullet}[-i]$ sont celles de $\umK^\bullet$ multipliées par $(-1)^i$. 
On a donc une suite exacte de complexes de $A$-modules
\begin{equation}\label{MH90i}
0\rightarrow \wedge^{i+1} L\otimes_A \umK^{\bullet}[-i-1]\rightarrow \rW^i/\rW^{i+2}(\mK^\bullet)
\rightarrow \wedge^i L\otimes_A \umK^{\bullet}[-i]\rightarrow 0.
\end{equation}
Celle-ci induit un morphisme de la catégorie dérivée $\bD^+(\bMod(A))$, 
\begin{equation}\label{MH90j}
\partial^i\colon \wedge^i L\otimes_A \umK^{\bullet}\rightarrow \wedge^{i+1}L\otimes_A \umK^{\bullet},
\end{equation}
que l'on appellera {\em bord associé à la filtration de Koszul de $\mK^\bullet$}.

\subsection{}\label{MH98}
Conservons les hypothèses et notations de \ref{MH90}. Soit, de plus, $(N,\kappa)$ un $A$-module de Higgs à coefficients dans $L$. 
On désigne par $\kappa'$ le $A$-champ de Higgs sur $N$ à coefficients dans $E$ induit par $\kappa$, 
par $\theta'=\theta\otimes\id+\id\otimes\kappa'$ le $A$-champ de Higgs total  sur $M\otimes_A N$ à coefficients dans $E$,
et par $\utheta'$ le $A$-champ de Higgs sur $M\otimes_A N$ à coefficients dans $\uE$ induit par $\theta'$, qui n'est autre que $\utheta\otimes \id$. 
On désigne par $\mK'^\bullet$ (resp. $\umK'^\bullet$, resp. $\cK'^\bullet$) 
le complexe de Dolbeault de $(M\otimes_AN,\theta')$ (resp. $(M\otimes_A N,\utheta')$, resp. $(N,\kappa')$)
et par $\theta'^\bullet$ (resp.  $\utheta'^\bullet$, resp.  $\kappa'^\bullet$) ses différentielles. 

On munit $\mK^\bullet$ et $\mK'^\bullet$ des filtrations de Koszul associées à la suite exacte \eqref{MH90a} et on pose 
$\mG^\bullet=\rW^0/\rW^2(\mK^\bullet)$ et $\mG'^\bullet=\rW^0/\rW^2(\mK'^\bullet)$ \eqref{MH90g}.
On note encore (abusivement) les différentielles de $\mG^\bullet$ (resp. $\mG'^\bullet$) par $\theta^i$ (resp. $\theta'^i$), ce qui n'induit aucune ambiguïté.  
On a alors des suites exactes canoniques de complexes \eqref{MH90i}
\begin{eqnarray}
0\longrightarrow L\otimes_A\umK^{\bullet}[-1]\stackrel{u^\bullet}{\longrightarrow} \mG^\bullet \stackrel{v^\bullet}{\longrightarrow} \umK^\bullet \longrightarrow 0,\label{MH98b}\\
0\longrightarrow L\otimes_A\umK'^{\bullet}[-1]\stackrel{u'^\bullet}{\longrightarrow} \mG'^\bullet \stackrel{v'^\bullet}{\longrightarrow} \umK'^\bullet \longrightarrow 0,\label{MH98c}
\end{eqnarray}
dont on note les bords associés dans $\bD^+(\bMod(A))$
\begin{eqnarray}
\partial\colon \umK^\bullet&\rightarrow& L\otimes_A\umK^\bullet,\label{MH98d}\\
\partial'\colon \umK'^\bullet&\rightarrow& L\otimes_A\umK'^\bullet.\label{MH98e}
\end{eqnarray}
On désigne par $\rC^\bullet$ (resp. $\rC'^\bullet$) le cône du morphisme $u^\bullet$ (resp. $u'^\bullet$), et par $c^\bullet$ (resp. $c'^\bullet$) ses différentielles. 
Pour tout entier $i$, on a donc $\rC^i=(L\otimes_A\umK^i)\oplus \mG^i$ et $c^i$ est définie par la matrice 
\begin{equation}\label{MH98f}
\begin{pmatrix}
\id\otimes \utheta^i &0 \\
 u^{i+1} & \theta^i
\end{pmatrix}.
\end{equation}
On a une description similaire du complexe $\rC'^\bullet$.

On notera que l'on a $\mK'^i=\mK^i\otimes_A N$, $\theta'^i=\theta^i\otimes \id+\id\otimes \kappa'^i$, 
$\umK'^i=\umK^i\otimes_AN$, $\utheta'^i=\utheta^i\otimes \id$, $\mG'^i=\mG^i\otimes N$ et $\rC'^i=\rC^i\otimes N$. 
Nous utiliserons implicitement ces identifications dans la suite de cette section.

Pour tous entiers $i,j\geq 0$, le morphisme $\kappa'^i\otimes \id\colon \mK'^i\rightarrow \mK'^{i+1}$ envoie $\rF^j \mK'^i$ dans $\rF^{j+1} \mK'^{i+1}$. 
Il induit donc un morphisme
\begin{equation}\label{MH98g}
\delta'^i\colon \mG'^i \rightarrow L\otimes_A\umK'^i
\end{equation}
tel que $\delta'^i\circ u'^i=0$ \eqref{MH98c}. 

\begin{lem}\label{MH99}
Conservons les hypothèses et notations de \ref{MH98}. Alors, 
\begin{itemize}
\item[{\rm (i)}] Pour tout entier $i\geq 0$, on a $\delta'^i= (\kappa\otimes\id)\circ v'^i$ \eqref{MH98c}.
\item[{\rm (ii)}] Les morphismes $\delta'^\bullet$ définissent un morphisme de complexes $\mG'^\bullet\rightarrow L\otimes_A\umK'^\bullet$.
\item[{\rm (iii)}] On a un isomorphisme de complexes de $A$-modules
\begin{equation}\label{MH99a}
\omega^\bullet \colon \rC'^\bullet \stackrel{\sim}{\rightarrow} \rC^\bullet \otimes_A N
\end{equation}
défini en degré $i$ par l'automorphisme de $\rC'^i=(L\otimes_A\umK'^i)\oplus \mG'^i$ qui envoie $(x,y)$ sur $(x+\delta'^i(y),y)$. 
\end{itemize}
\end{lem}

La proposition (i) est immédiate et elle implique aussitôt (ii).  Pour établir (iii), il faut montrer que pour tout
élément $(x,y)$ de $\rC'^i=(L\otimes_A\umK'^i)\oplus \mG'^i$, on a 
\begin{eqnarray}
(\utheta^i\otimes\id)(x)+\delta'^{i+1}(u'^{i+1}(x)+\theta'^i(y))&=&(\utheta^i\otimes\id)(x+\delta'^i(y)),\\
u'^{i+1}(x)+\theta'^i(y)&=& u'^{i+1}(x+\delta'^i(y))+(\theta^i\otimes\id)(y).
\end{eqnarray}
La première équation résulte de (ii) et du fait que $\delta'^{i+1}\circ u'^{i+1}=0$. La seconde équation est une conséquence des 
relations $\theta'^i=\theta^i\otimes\id+ \kappa'^i\otimes\id$ et $\kappa'^i\otimes\id=u'^{i+1}\circ \delta'^i$.

\begin{prop}\label{MH100}
Sous les hypothèses et avec les notations de \ref{MH98}, identifiant $\umK'^\bullet$ avec $\umK^\bullet\otimes_A N$, on a 
\begin{equation}\label{MH100a}
\partial'=\partial\otimes \id+\id \otimes \kappa.
\end{equation}
\end{prop}

En effet, notons $\pi'^\bullet_1$ (resp. $\pi'^\bullet_2$) la projection de $\rC'^\bullet=(L\otimes_A\umK'^\bullet)\oplus \mG'^\bullet$ dans 
$L\otimes_A\umK'^\bullet$ (resp. $\mG'^\bullet$). Le composé $v'^\bullet \circ \pi'^\bullet_2\colon \rC'^\bullet\rightarrow \umK'^\bullet$ est alors un quasi-isomorphisme
\eqref{MH98c}, et $-\partial'$ \eqref{MH98e} est le composé  dans $\bD^+(\bMod(A))$ de l'inverse de $v'^\bullet \circ \pi'^\bullet_2$ et de $\pi'^\bullet_1$
(\cite{sp} \href{https://stacks.math.columbia.edu/tag/09KF}{09KF}). 
On décrit de même $\partial$ \eqref{MH98d} en terme des projections canoniques $\pi^\bullet_1$ et $\pi^\bullet_2$ de $\rC^\bullet$. 
Compte tenu de \ref{MH99}, on a 
\begin{eqnarray}
(\pi^\bullet_1\otimes \id) \circ \omega^\bullet &=& \pi'^\bullet_1 +\delta'^\bullet\circ \pi'^\bullet_2,\\
(\pi^\bullet_2\otimes \id) \circ \omega^\bullet&=&\pi'^\bullet_2.
\end{eqnarray}
On en déduit que 
\begin{eqnarray*}
\partial\otimes \id &=&\partial'-\delta'^\bullet\circ \pi'^\bullet_2\circ (v'^\bullet\circ \pi'^\bullet_2)^{-1}\\
&=& \partial'-(f^*(\kappa)\otimes\id)\circ v'^\bullet\circ \pi'^\bullet_2\circ (v'^\bullet\circ \pi'^\bullet_2)^{-1}\\
&=& \partial'-\kappa\otimes\id,
\end{eqnarray*}
où la seconde relation résulte de \ref{MH99}(ii), d'où la proposition.

\subsection{}\label{MH97}
Soient $f\colon (X,A)\rightarrow (Y,B)$ un morphisme de topos annelés, $L$ un $B$-module localement libre de type fini, 
$E$ et $\uE$ deux $A$-modules localement libres de type fini, 
\begin{equation}\label{MH97a}
0\longrightarrow f^*(L)\longrightarrow E\stackrel{u}{\longrightarrow} \uE\longrightarrow 0
\end{equation}
une suite exacte, $(M,\theta)$ un $A$-module de Higgs à coefficients dans $E$.
On désigne par $\utheta$ le $A$-champ de Higgs sur $M$ à coefficients dans $\uE$ induit par $\theta$,
et par $\mK^\bullet$ (resp. $\umK^\bullet$) le complexe de Dolbeault de $(M,\theta)$ (resp. $(M,\utheta)$). 
On munit $\mK^\bullet$ de la filtration de Koszul associée à la suite exacte \eqref{MH97a}, cf. \eqref{MH90g}. 
On considère naturellement
\begin{equation}\label{MH97c}
\oplus_{i\geq 0} \rW^i\mK^\bullet
\end{equation}
comme un module (à droite) bigradué sur la $A$-algèbre bigraduée $\oplus_{i\geq 0} \rW^i\wedge^\bullet E$.
On notera que sa différentielle de bidegré $(0,1)$ est $(\oplus_{i\geq 0} \rW^i\wedge^\bullet E)$-linéaire. Par ailleurs, on a un homomorphisme 
canonique de $A$-algèbres graduées (par $i$)
\begin{equation}\label{MH97d}
f^*(\oplus_{i\geq 0} \wedge^i L)\rightarrow \oplus_{i\geq 0} \rW^i\wedge^\bullet E.
\end{equation}

Pour tout entier $i\geq 0$, on pose 
\begin{equation}\label{MH97e}
\rE_0^{i,\bullet}=\Gr_\rW^i(\mK^\bullet). 
\end{equation}
On a une suite exacte de $(\oplus_{i\geq 0} \rF^i\wedge^\bullet E)$-modules bigradués
\begin{equation}\label{MH97f}
0\rightarrow \oplus_{i\geq 0}\rE_0^{i+1,\bullet}\rightarrow \oplus_{i\geq 0}\rW^i/\rW^{i+2}(\mK^\bullet)
\rightarrow \oplus_{i\geq 0}\rE_0^{i,\bullet}\rightarrow 0.
\end{equation}
Compte tenu de \eqref{MH97d}, on peut considérer cette suite comme une suite exacte de complexes de 
$f^*(\oplus_{i\geq 0}\wedge^iL)$-modules gradués (par $i$). 

En vertu de (\cite{sp} \href{https://stacks.math.columbia.edu/tag/015W}{015W}), il existe une suite spectrale convergente canonique
\begin{equation}\label{MH97g}
\rE_1^{i,j}=\rR^{i+j}f_*(\rE_0^{i,\bullet})\Rightarrow \rR^{i+j}f_*(\mK^\bullet),
\end{equation}
dite suite spectrale d'hypercohomologie du complexe filtré $\mK^\bullet$ (\cite{hodge2} 1.4.5).
D'après \eqref{MH90h} et la formule de projection (\cite{sp} \href{https://stacks.math.columbia.edu/tag/0B54}{0B54}), on a un isomorphisme canonique
\begin{equation}\label{MH97h}
\rE_1^{i,j}\stackrel{\sim}{\rightarrow} \rR^jf_*(\umK^\bullet) \otimes_A \wedge^i L.
\end{equation}
Pour tout entier $j\geq 0$, on obtient un complexe de $B$-modules
\begin{equation}\label{MH97i}
\rR^jf_*(\umK^\bullet)\stackrel{d_1^{0,j}}{\longrightarrow} \rR^jf_*(\umK^\bullet)\otimes_AL\stackrel{d_1^{1,j}}{\longrightarrow}
\rR^jf_*(\umK^\bullet)\otimes_A\wedge^2L\stackrel{d_1^{2,j}}{\longrightarrow} \dots
\end{equation}
dont les différentielles ne sont autres que les composantes homogènes des morphismes bord induits par la suite exacte \eqref{MH97f}. 
Comme cette dernière est une suite exacte de complexes de $f^*(\oplus_{i\geq 0}\wedge^iL)$-modules (à droite) gradués (par $i$), 
on obtient l'énoncé suivant:

\begin{lem}\label{MH96}
Sous les hypothèses de \ref{MH97} et avec les mêmes notations, pour tout entier $j\geq 0$, 
le morphisme $d_1^{0,j}$ \eqref{MH97i} est un $B$-champ de Higgs sur $\rR^jf_*(\umK^\bullet)$ à coefficients dans $L$, 
dont le complexe de Dolbeault n'est autre que le complexe $(\rE_1^{i,j})_{i\geq 0}$.
\end{lem}

On appelle {\em $B$-champ de Katz-Oda sur $\rR^jf_*(\umK^\bullet)$ à coefficients dans $L$}, le $B$-champ de Higgs $d_1^{0,j}$,
dont la construction rappelle celle de la connexion de Gauss-Manin par ces auteurs \cite{ko}. Il est aussi appelé {\em champ de Gauss-Manin} dans (\cite{ov} §~3.1).

\begin{rema}\label{MH110}
Conservons les hypothèses et notations de \ref{MH97}. Supposons de plus que le champ de Higgs $\theta$ est nul. 
La filtration de Koszul de la $A$-algèbre extérieure $\wedge^\bullet E$ associée à la suite exacte \eqref{MH97a}, 
induit pour tout entier $j\geq 0$ une suite exacte \eqref{MH90e}
\begin{equation}\label{MH110a}
0\rightarrow f^*(L)\otimes_{\co_{\fX'}} \wedge^{j-1} \uE\rightarrow \rW^0/\rW^2(\wedge^\bullet E)
\rightarrow \wedge^j \uE\rightarrow 0.
\end{equation}
Compte tenu de la formule de projection (\cite{sp} \href{https://stacks.math.columbia.edu/tag/0B54}{0B54}), on en déduit un morphisme de $\bD^+(\bMod(B))$ 
\begin{equation}\label{MH110b}
\rR f_*(M\otimes_A\wedge^j\uE)\rightarrow L\otimes_B \rR f_*(M\otimes_A\wedge^{j-1}\uE)[+1]. 
\end{equation}
Pour tout entier $q\geq 0$, on a alors un isomorphisme canonique 
\begin{equation}\label{MH110c}
\rR^q f_*(\umK^\bullet)\stackrel{\sim}{\rightarrow} \oplus_{0\leq i\leq q} \rR^i f_*(M\otimes_A\wedge^{q-i}\uE), 
\end{equation}
le champ de Katz-Oda sur $\rR^q f_*(\umK^\bullet)$ à coefficients dans $L$ \eqref{MH96} étant induit par les morphismes \eqref{MH110b}. Il est donc quasi-nilpotent \eqref{MH151}. 
\end{rema}

\subsection{}\label{MH101}
Conservons les hypothèses et notations de \ref{MH97}. Soit, de plus, $(\rF^iM)_{i\geq 0}$ une filtration décroissante de $(M,\theta)$ dans $\bMH(A,E)$ \eqref{MH150}
telle que $\rF^0 M=M$ et $\rF^nM=0$ pour un entier $n\geq 0$. C'est aussi une filtration de $(M,\utheta)$ dans $\bMH(A,\uE)$.  
On note $\rF^i\mK^\bullet = \rF^iM\otimes_A \wedge^\bullet E$ et $\rF^i\umK^\bullet= \rF^iM\otimes_A \wedge^\bullet\uE$ les complexes de Dolbeault de $\rF^i M$ vu comme objet de 
$\bMH(A,E)$ et $\bMH(A,\uE)$ respectivement. Les champs de Higgs $\theta$ et $\utheta$ induisent des champs de Higgs sur $\Gr^i_\rF M$. On note  
$\Gr^i_\rF \mK^\bullet= \Gr^i_\rF M\otimes_A \wedge^\bullet E$ et $\Gr^i_\rF \umK^\bullet = \Gr^i_\rF M\otimes_A \wedge^\bullet\uE$ les complexes de Dolbeault associés.
On obtient ainsi une filtration $(\rF^i\mK^\bullet)_{0\leq i\leq n}$ (resp. $(\rF^i\umK^\bullet)_{0\leq i\leq n}$) du complexe $\mK^\bullet$ (resp. $\umK^\bullet$)
de gradués $(\Gr^i_\rF \mK^\bullet)_{0\leq i\leq n}$ (resp. $(\Gr^i_\rF \umK^\bullet)_{0\leq i\leq n}$). 

D'après (\cite{sp} \href{https://stacks.math.columbia.edu/tag/015W}{015W}), il existe une suite spectrale convergente canonique
\begin{equation}\label{MH101a}
\rE_1^{i,j}=\rR^{i+j}f_*(\Gr^i_\rF \umK^\bullet)\Rightarrow \rR^{i+j}f_*(\umK^\bullet).
\end{equation} 

Pour tout entier $q\geq 0$, on désigne par 
\begin{equation}\label{MH101c}
\kappa^q\colon \rR^qf_*(\umK^\bullet)\rightarrow \rR^qf_*(\umK^\bullet)\otimes_BL
\end{equation}
le champ de Katz-Oda associé au $A$-module de Higgs $(M,\theta)$ \eqref{MH96}. 

Pour tous entiers $i,j\geq 0$, on désigne par 
\begin{equation}\label{MH101d}
\kappa^{i,j}_1\colon \rR^{i+j}f_*(\Gr^i_\rF \umK^\bullet)\rightarrow \rR^{i+j}f_*(\Gr^i_\rF \umK^\bullet)\otimes_BL
\end{equation}
le champ de Katz-Oda associé au $A$-module de Higgs à coefficients dans $E$ induit par $\theta$ sur $\Gr^i_\rF M$.

\begin{prop}\label{MH102}
Sous les hypothèses de \ref{MH101}, les morphismes \eqref{MH101c} et \eqref{MH101d} définissent un morphisme de la suite spectrale \eqref{MH101a} 
dans la suite spectrale qui s'en déduit par application du foncteur $-\otimes_BL$. 
\end{prop}

Soit $j$ un entier $\geq 0$. Considérons la suite exacte  
\begin{equation}\label{MH102a}
0\longrightarrow \Gr_\rW^{j+1}(\mK^\bullet) \stackrel{u}{\longrightarrow} 
\rW^j/\rW^{j+2}(\mK^\bullet) \stackrel{v}{\longrightarrow} \Gr_\rW^j(\mK^\bullet) \longrightarrow 0
\end{equation}
induite par la filtration de Koszul de $\mK^\bullet$ \eqref{MH90g}, et notons $\rC^\bullet$ le cône du morphisme $v$ et 
\begin{equation}\label{MH102b}
\Gr_\rW^j(\mK^\bullet)\stackrel{\nu}{\longrightarrow}  \rC^\bullet\stackrel{\mu}{\longleftarrow} \Gr_\rW^{j+1}(\mK^\bullet)[1]
\end{equation}
les morphismes canoniques. Alors, $\mu$ est un quasi-isomorphisme, et le morphisme $\mu^{-1}\circ \nu$ de $\bD^+(\bMod(A))$ s'identifie à l'opposé du 
bord $\partial^j$ associé à la filtration de Koszul de $\mK^\bullet$ \eqref{MH90j}. 

La filtration $\rF$ de  $\mK^\bullet$ \eqref{MH101} induit une filtration du complexe $\rW^j\mK^\bullet$ (\cite{hodge2} 1.1.8) 
qui induit à son tour des filtrations des complexes $\Gr^j_\rW \mK^\bullet$ et $\rW^j/\rW^{j+2}(\mK^\bullet)$. 
Les morphismes $u$ et $v$ sont stricts (\cite{hodge2} 1.1.9). En particulier, la filtration $\rF$ de $\mK^\bullet$ induit une filtration de $\rC^\bullet$. 

On rappelle que pour tout entier $i\geq 0$, $\rF^i\mK^\bullet$ (resp. $\Gr^i_\rF\mK^\bullet$)
est le complexe de Dolbeault de $M_i$ (resp. $\Gr^i_\rF M$). D'après (\cite{ac} I §2.6  prop.~7),
comme pour tout $\ell\geq 0$, $\wedge^\ell E/\rW^j \wedge^\ell E$ est $A$-plat \eqref{MH90e}, on a un isomorphisme canonique
\begin{equation}
\rF^i(\rW^j\mK^\bullet)\stackrel{\sim}{\rightarrow} \rW^j(\rF^i\mK^\bullet).
\end{equation}
Celui-ci induit un isomorphisme 
\begin{equation}
\Gr_\rF^i(\rW^j\mK^\bullet)\stackrel{\sim}{\rightarrow} \rW^j(\Gr_\rF^i\mK^\bullet).
\end{equation}
D'après (\cite{hodge2} 1.1.11), on en déduit pour tout entier $j'\geq j+1$, un isomorphisme canonique 
\begin{equation}
\Gr_\rF^i(\rW^j/\rW^{j'}(\mK^\bullet))\stackrel{\sim}{\rightarrow} \rW^j/\rW^{j'}(\Gr_\rF^i\mK^\bullet).
\end{equation}

Le morphisme $\mu$ \eqref{MH102b} est un quasi-isomorphisme filtré
pour les filtrations induites par la filtration $\rF$ de $\mK^\bullet$ (\cite{hodge2} 1.3.6).
 En effet, d'après ce qui précède et (\cite{hodge2} 1.1.11), l'image de la suite \eqref{MH102a} par le foncteur $\Gr^i_\rF$ s'identifie à la suite exacte canonique
\begin{equation}\label{MH102c}
0\rightarrow \Gr^{j+1}_\rW(\Gr^i_\rF(\mK^{\bullet}))\rightarrow
\rW^j/\rW^{j+2}(\Gr^i_\rF(\mK^\bullet)) \rightarrow \Gr_\rW^j(\Gr^i_\rF(\mK^{\bullet})) \rightarrow 0. 
\end{equation}
Le cône de la troisième flèche s'identifie alors à $\Gr^i_\rF\rC^\bullet$, et $\Gr^i_\rF\mu$ et $\Gr^i_\rF\nu$ s'identifient aux morphismes analogues à $\mu$ et $\nu$ pour 
la suite \eqref{MH102c}. Par suite, $\Gr^i_\rF\mu$ est un quasi-isomorphisme, d'où l'assertion recherchée. 

Prenant $j=0$, les morphismes $\mu$ et $\nu$ définissent des morphismes de suites spectrales d'hypercohomologie de complexes filtrés (par les filtrations induites par $\rF$), 
le morphisme induit par $\mu$ étant un isomorphisme. La proposition s'ensuit compte tenu de \eqref{MH90h} et de la définition du champ de Katz-Oda \eqref{MH96}.

\begin{cor}\label{MH105}
Sous les hypothèses de \ref{MH101}, la suite spectrale \eqref{MH101a} est sous-jacente à une suite spectrale dans la catégorie abélienne $\bMH(B,L)$, 
\begin{equation}\label{MH105a}
\rE_1^{i,j}=(\rR^{i+j}f_*(\Gr^i_\rF \umK^\bullet),\kappa_1^{i,j})\Rightarrow (\rR^{i+j}f_*(\umK^\bullet),\kappa^{i+j}).
\end{equation} 
\end{cor}

\begin{cor}\label{MH103}
Conservons les hypothèses de \ref{MH101} et notons pour tout entier $q\geq 0$, $(\rN^i\rR^qf_*(\umK^\bullet))_{0\leq i\leq q}$ la filtration de $\rR^qf_*(\umK^\bullet)$ 
aboutissement de la suite spectrale \eqref{MH101a}. Alors, pour tout $0\leq i\leq q$, on a 
\begin{equation}\label{MH103a}
\kappa^q(\rN^i\rR^qf_*(\umK^\bullet))\subset \rN^i\rR^qf_*(\umK^\bullet)\otimes_BL.
\end{equation}
\end{cor}
Cela résulte aussitôt de \ref{MH105}.

\begin{prop}\label{MH104}
Reprenons les hypothèses de \ref{MH97}, supposons de plus le $A$-module de Higgs $(M,\theta)$ quasi-nilpotent \eqref{MH151}. 
Alors, pour tout entier $q\geq 0$, le $B$-module de Higgs $(\rR^qf_*(\umK^\bullet),\kappa^q)$, où $\kappa^q$ est le champ de Katz-Oda \eqref{MH96}, est quasi-nilpotent. 
\end{prop}

En effet, soit $(\rF^iM)_{i\geq 0}$ une filtration décroissante quasi-nilpotente de $(M,\theta)$ dans $\bMH(A,E)$ \eqref{MH151}
telle que $\rF^0 M=M$ et $\rF^nM=0$ pour un entier $n\geq 0$. Reprenons les notations de \ref{MH101}.
Pour tous entiers $i,j\geq 0$, le champ de Higgs induit par $\utheta$ sur $\Gr_\rF^i M$ étant nul, 
le $B$-module de Higgs $(\rR^{i+j}f_*(\Gr^i_\rF \umK^\bullet),\kappa_1^{i,j})$ est quasi-nilpotent d'après \ref{MH110}. 
En vertu de \ref{MH141}(i), le gradué de la filtration de $(\rR^qf_*(\umK^\bullet),\kappa^q)$ dans $\bMH(B,L)$ aboutissement de la suite spectrale \eqref{MH105a}
est donc quasi-nilpotent. Par suite, le $B$-module de Higgs $(\rR^qf_*(\umK^\bullet),\kappa^q)$ est quasi-nilpotent d'après \ref{MH141}(ii). 

\begin{prop}\label{MH106}
Reprenons les hypothèses de \ref{MH97} et supposons de plus que $f$ soit un morphisme cohérent de schémas, que les $\co_X$-modules $E$, $\uE$ et $L$ soient quasi-cohérents, 
et que le $\co_X$-module de Higgs $(M,\theta)$ soit nilpotent \eqref{MH15}.
Alors, pour tout entier $q\geq 0$, le $\co_Y$-module de Higgs $(\rR^qf_*(\umK^\bullet),\kappa^q)$, où $\kappa^q$ est le champ de Katz-Oda \eqref{MH96}, est nilpotent.   
\end{prop}

La preuve est identique à celle de \ref{MH104} en utilisant \ref{MH142} au lieu de \ref{MH141}.

\subsection{}\label{MH9}
Soient $(X,A)$ un topos annelé, $B$ une $A$-algèbre, 
$\Omega^1_{B/A}$ le $B$-module des différentielles de Kähler de $B$ sur $A$ (\cite{illusie1} II 1.1.2), 
\begin{equation}\label{MH9a}
\Omega_{B/A}=\oplus_{n\in \mN}\Omega^n_{B/A}=\wedge_B(\Omega^1_{B/A})
\end{equation} 
l'algèbre extérieure de $\Omega^1_{B/A}$. Il existe alors une et une unique $A$-anti-dérivation 
$d\colon \Omega_{B/A}\rightarrow \Omega_{B/A}$ de degré $1$, de carré nul, qui prolonge
la $A$-dérivation universelle $d\colon B\rightarrow \Omega^1_{B/A}$. Cela résulte par exemple de 
(\cite{alg10} §2.10 prop.~13) en remarquant que $\Omega^1_{B/A}$ est le faisceau associé au préfaisceau 
$U\mapsto \Omega^1_{B(U)/A(U)}$ $(U\in \ob(X))$.

Soient $M$ un $B$-module, $\lambda\in A(X)$. 
Une {\em $\lambda$-connexion} sur $M$ relativement à l'extension $B/A$
est la donnée d'un morphisme $A$-linéaire 
\begin{equation}\label{MH9b}
\nabla\colon M\rightarrow \Omega^1_{B/A}\otimes_BM
\end{equation}
tel que pour toutes sections locales $x$ de $B$ et $s$ de $M$, on ait 
\begin{equation}\label{MH9c}
\nabla(xs)=\lambda d(x)\otimes s+x\nabla(s).
\end{equation} 
On dit aussi que $(M,\nabla)$ est un $B$-module à $\lambda$-connexion relativement à l'extension $B/A$.
On omettra l'extension $B/A$ de la terminologie lorsqu'il n'y a aucun risque de confusion. 
Le morphisme $\nabla$ se prolonge en un unique morphisme $A$-linéaire gradué de degré $1$ que l'on note aussi 
\begin{equation}\label{MH9d}
\nabla\colon \Omega_{B/A}\otimes_BM\rightarrow   \Omega_{B/A}\otimes_BM,
\end{equation} 
tel que pour toutes sections locales $\omega$ de  $\Omega^i_{B/A}$ et $s$ de $\Omega^j_{B/A}\otimes_BM$  
($i,j\in \mN$), on ait 
\begin{equation}\label{MH9e}
\nabla(\omega\wedge s)=\lambda  d(\omega)\wedge s+(-1)^i\omega \wedge \nabla(s).
\end{equation} 
Itérons cette formule~: 
\begin{equation}\label{MH9f}
\nabla \circ \nabla(\omega\wedge s)= \omega \wedge \nabla\circ \nabla(s).
\end{equation} 
On dit que $\nabla$ est {\em intégrable} si $\nabla \circ \nabla=0$.

Soient $(M,\nabla)$, $(M',\nabla')$ deux modules à $\lambda$-connexion. 
Un morphisme de $(M,\nabla)$ dans $(M',\nabla')$
est la donnée d'un morphisme $B$-linéaire $u\colon M\rightarrow M'$ 
tel que  $(\id \otimes u)\circ \nabla=\nabla'\circ u$. 

Les $1$-connexions sont classiquement appelées {\em connexions}. 
Les $0$-connexions intégrables sont les $B$-champs de Higgs à coefficients dans $\Omega^1_{B/A}$ \eqref{MH1}.

\subsection{}\label{MH10}
Soient $f\colon (X',A')\rightarrow (X,A)$ un morphisme de topos annelés, $B$ une $A$-algèbre, 
$B'$ une $A'$-algèbre, $\alpha\colon f^*(B)\rightarrow B'$ un homomorphisme de $A'$-algèbres,
$\lambda\in \Gamma(X,A)$, $(M,\nabla)$ un module à $\lambda$-connexion relativement à l'extension $B/A$ \eqref{MH9}. 
Notons $\lambda'$ l'image canonique de $\lambda$ dans $\Gamma(X',A')$, 
$d'\colon B'\rightarrow \Omega^1_{B'/A'}$ la $A'$-dérivation universelle de $B'$ et 
\begin{equation}\label{MH10a}
\gamma\colon f^*(\Omega^1_{B/A}) \rightarrow \Omega^1_{B'/A'}
\end{equation}
le morphisme $\alpha$-linéaire canonique. On voit aussitôt que 
$f^*(\nabla)$ est une $\lambda'$-connexion sur $f^*(M)$ relativement à l'extension $f^*(B)/A'$,
qui est intégrable si $\nabla$ l'est.
Par ailleurs, il existe un unique morphisme $A'$-linéaire 
\begin{equation}\label{MH10b}
\nabla'\colon B'\otimes_{f^*(B)}f^*(M)\rightarrow  \Omega^1_{B'/A'}\otimes_{f^*(B)}f^*(M)
\end{equation}
tel que pour toutes sections locales $x'$ de $B'$ et $t$ de $f^*(M)$, on ait 
\begin{equation}\label{MH10c}
\nabla'(x'\otimes t)=\lambda' d'(x')\otimes t+x'(\gamma\otimes\id)(f^*(\nabla)(t)).
\end{equation}
C'est une $\lambda'$-connexion sur $B'\otimes_{f^*(B)}f^*(M)$ relativement à l'extension 
$B'/A'$, qui est intégrable si $\nabla$ l'est.

\subsection{}\label{MH8}
Soient $(X,A)$ un topos annelé, $B$ une $A$-algèbre, $\lambda\in \Gamma(X,A)$, $(M,\nabla)$ un
module à $\lambda$-connexion relativement à l'extension $B/A$ \eqref{MH9}. 
Supposons qu'il existe un $A$-module $E$ et un $B$-isomorphisme $\gamma\colon E\otimes_AB\stackrel{\sim}{\rightarrow}\Omega^1_{B/A}$ 
tels que pour toute section locale $\omega$ de $E$, on ait $d(\gamma(\omega\otimes 1))=0$.
On note $\vartheta\colon M\rightarrow E\otimes_AM$  le morphisme induit par $\nabla$ et $\gamma$. 
\begin{itemize}
\item[(i)] Pour que la $\lambda$-connexion $\nabla$ soit intégrable, il faut et il suffit que 
$\vartheta$ soit un $A$-champ de Higgs sur $M$ à coefficients dans $E$ (\cite{agt} II.2.12).
\item[(ii)] Soit $(N,\theta)$ un $A$-module de Higgs à coefficients dans $E$. Le morphisme $A$-linéaire 
\begin{equation}\label{MH8a}
\nabla'\colon M\otimes_AN\rightarrow  \Omega^1_{B/A}\otimes_BM\otimes_AN
\end{equation}
défini par 
\begin{equation}\label{MH8b}
\nabla'=\nabla\otimes_A \id_N+ (\gamma\otimes_B\id_{M\otimes_AN})(\id_M\otimes_A \theta),
\end{equation}
est une $\lambda$-connexion sur $M\otimes_AN$ relativement à $B/A$.
Si $\nabla$ est intégrable, il en est de même de $\nabla'$. 
\item[(iii)] Soient $B'$ une $B$-algèbre, $\delta'\colon B'\rightarrow \Omega^1_{B/A}\otimes_BB'$ une $A$-dérivation telle que pour toute section locale $b$ de $B$, on ait 
$\delta'(b)=\lambda d(b)\otimes 1$. Il existe alors un unique morphisme $A$-linéaire
\begin{equation}\label{MH8c}
\nabla'\colon M\otimes_BB'\rightarrow  \Omega^1_{B/A}\otimes_BM\otimes_BB'
\end{equation}
tel que pour toutes sections locales $m$ de $M$ et $b'$ de $B'$, on ait 
\begin{equation}\label{MH8d}
\nabla'(m\otimes b')=\nabla(m)\otimes_B b'+ m\otimes_B \delta'(b').
\end{equation}
On désigne par $\vartheta'\colon M\otimes_BB'\rightarrow  E\otimes_AM\otimes_BB'$ (resp. $\theta'\colon B'\rightarrow  E\otimes_AB'$)
le morphisme induit par $\nabla'$ (resp. $\delta'$) et $\gamma$.  
Si $\nabla$ est integrable et si $\theta'$ est un $A$-champ de Higgs à coefficients dans $E$, alors $\vartheta'$ est un $A$-champ de Higgs à coefficients dans $E$. 
\end{itemize}

\section{Ind-objets d'une catégorie}\label{indsh}

\subsection{}\label{indsh1}
Soient $\cC$ une $\mU$-catégorie (\cite{sga4} I 1.1) \eqref{notconv3}, $\hcC$ la catégorie des préfaisceaux de $\mU$-ensembles sur $\cC$ \eqref{notconv5},
\begin{equation}\label{indsh1a}
\tth_\cC\colon \cC\rightarrow \hcC, \ \ \ X\mapsto (\tth_\cC(X)\colon Y\mapsto \Hom_\cC(Y,X)),
\end{equation} 
le foncteur canonique (\cite{sga4} I 1.3). Celui-ci est pleinement fidèle (\cite{sga4} I 1.4). 
On rappelle que les limites inductives dans $\hcC$ sont représentables
et qu'elles se calculent argument par argument (\cite{sga4} I 3.1). 
On notera cependant que même quand les limites inductives dans $\cC$ sont représentables,
le foncteur $\tth_\cC$ ne commute pas en général aux limites inductives. 
Nous utiliserons la notation $\indcolim$ pour désigner les limites inductives dans $\hcC$ et nous réserverons la notation 
$\colim$ pour désigner les limites inductives dans $\cC$.  

Soient $I$ et $J$ deux petites catégories, $\alpha\colon I\rightarrow \cC$ et $\beta\colon J\rightarrow \cC$ deux foncteurs. 
Pour tout objet $X$ de $\cC$, on a 
\begin{equation}\label{indsh1b}
\Hom_{\hcC}(X,\indcolim \alpha)=\underset{\underset{i\in I}{\longrightarrow}}{\lim}\ \Hom_{\cC}(X,\alpha(i)).
\end{equation}
On a des isomorphismes canoniques 
\begin{eqnarray}\label{indsh1c}
\Hom_{\hcC}(\indcolim \alpha,\indcolim \beta)&\stackrel{\sim}{\rightarrow}& 
\underset{\underset{i\in I}{\longleftarrow}}{\lim}\ \Hom_{\hcC}(\alpha(i),\indcolim \beta)\\
&\stackrel{\sim}{\rightarrow}& 
\underset{\underset{i\in I}{\longleftarrow}}{\lim}\ \underset{\underset{j\in J}{\longrightarrow}}{\lim}\ \Hom_{\cC}(\alpha(i),\beta(j)). \nonumber
\end{eqnarray}

\begin{defi}[\cite{sga4} I 8.1.8]\label{indsh2}
On dit qu'une catégorie filtrante $I$ est {\em essentiellement petite} si $I$ est une $\mU$-catégorie et 
si $\ob(I)$ admet une petite partie cofinale (\cite{sga4} I 8.1.4). 
\end{defi}

Pour qu'une petite partie $E$ de $\ob(I)$ soit cofinale, il faut et il suffit que pour tout objet $i$ de $I$, 
il existe un morphisme $i\rightarrow j$ avec $j\in E$ (\cite{sga4} 8.1.3(c)).

\begin{defi}[\cite{ks2} 6.1.1]\label{indsh3}
Soit $\cC$ une $\mU$-catégorie. 
\begin{itemize}
\item[(i)] On appelle {\em ind-objet de $\cC$} tout objet de $\hcC$
isomorphe à $\indcolim \alpha$ pour un foncteur $\alpha\colon I\rightarrow \cC$ tel que $I$ soit une petite catégorie filtrante. 
\item[(ii)] On appelle catégorie des {\em ind-objets de $\cC$} et l'on note $\Ind(\cC)$ la sous-catégorie pleine de $\hcC$ formée des ind-objets de $\cC$. 
On désigne par 
\begin{equation}\label{indsh3a}
\iota_\cC \colon \cC\rightarrow \Ind(\cC)
\end{equation}
le foncteur induit par $\tth_\cC$. 
\end{itemize}
\end{defi}

\addtocounter{subsubsection}{1}
\addtocounter{equation}{1}

On peut faire les remarques suivantes:

\subsubsection{}\label{indsh3b}
Pour tout objet $F$ de $\hcC$, notant $\alpha\colon \cC_{/F}\rightarrow \cC, (X\rightarrow F)\mapsto X$, le foncteur ``source'',  
le morphisme canonique
\begin{equation}\label{indsh3c}
\indcolim\alpha \rightarrow F
\end{equation}
est un isomorphisme (\cite{sga4} I 3.4)

\addtocounter{subsubsection}{1}
\addtocounter{equation}{1}

\subsubsection{}\label{indsh3d}
Pour qu'un objet $F$ de $\hcC$ soit un ind-objet de $\cC$, 
il faut et il suffit que la catégorie $\cC_{/F}$ \eqref{notconv5} soit filtrante et essentiellement petite (\cite{sga4} I 8.3.3, \cite{ks2} 6.1.5). 

\subsubsection{}\label{indsh3e}
La catégorie $\Ind(\cC)$ est une $\mU$-catégorie (\cite{ks2} 6.1.2).
Le foncteur canonique $\iota_\cC$ \eqref{indsh3a}
est exact à droite et petit à droite (\cite{ks2} 6.1.6) (cf. \cite{ks2} 3.3.1 et 3.3.8). 
En particulier, si les limites inductives finies dans $\cC$ sont représentables, $\iota_\cC$ commute à ces limites (\cite{ks2} 3.3.2). 

\subsubsection{}\label{indsh3f}
Les petites limites inductives filtrantes dans $\Ind(\cC)$ sont représentables et le foncteur d'injection canonique
$\Ind(\cC)\rightarrow \hcC$ commute à ces limites (\cite{ks2} 6.1.8). 

\subsubsection{}\label{indsh3g}
Si les limites projectives finies (resp. les petites limites projectives) dans $\cC$ sont représentables, alors elles le sont aussi dans $\Ind(\cC)$
et les foncteurs canoniques $\iota_\cC\colon \cC\rightarrow \Ind(\cC)$ et $\Ind(\cC)\rightarrow \hcC$ commutent à ces limites (\cite{ks2} 6.1.17).

\subsubsection{}\label{indsh3h}
Si les conoyaux des doubles flèches (resp. les sommes finis, resp. les limites inductives finies) dans $\cC$ sont représentables, 
alors les conoyaux des doubles flèches (resp. les petites sommes, resp. les petites limites inductives) dans $\Ind(\cC)$ sont représentables
(\cite{ks2} 6.1.18).

\subsubsection{}\label{indsh3i}
Si les limites projectives finies et les limites inductives finies dans $\cC$ sont représentables, 
alors les petites limites inductives filtrantes dans $\Ind(\cC)$ sont exactes (\cite{ks2} 6.1.19).

\subsection{}\label{indsh4}
Soit $\phi \colon \cC\rightarrow \cC'$ un foncteur entre $\mU$-catégories. 
D'après (\cite{sga4} I 8.6.3, \cite{ks2} 6.1.9), il existe essentiellement un unique foncteur 
\begin{equation}\label{indsh4a}
\rI\phi\colon \Ind(\cC)\rightarrow \Ind(\cC')
\end{equation}
vérifiant les deux conditions suivantes:
\begin{itemize}
\item[{\rm (i)}] le diagramme 
\begin{equation}\label{indsh4b}
\xymatrix{
\cC\ar[r]^\phi\ar[d]_{\iota_\cC}&\cC'\ar[d]^{\iota_{\cC'}}\\
{\Ind(\cC)}\ar[r]^{\rI\phi}&{\Ind(\cC')}}
\end{equation}
est commutatif à isomorphisme canonique près;
\item[{\rm (ii)}] pour toute petite catégorie filtrante $J$ et tout foncteur $\alpha\colon J\rightarrow \Ind(\cC)$, le morphisme canonique
\begin{equation}\label{indsh4c}
\rI\phi(\indcolim \alpha)\rightarrow \indcolim(\phi\circ \alpha)
\end{equation}
est un isomorphisme.
\end{itemize}

\addtocounter{subsubsection}{3}
\addtocounter{equation}{1}

On peut faire les remarques suivantes:

\subsubsection{}\label{indsh4d}
Si $\phi$ est fidèle (resp. pleinement fidèle), il en est de même de $\rI\phi$ (\cite{ks2} 6.1.10).

\addtocounter{equation}{1}

\subsubsection{}\label{indsh4e}
Si $\psi\colon \cC'\rightarrow \cC''$ est un foncteur tel que $\cC''$ soit une $\mU$-catégorie, alors on a un isomorphisme canonique entre foncteurs
de $\Ind(\cC)$ dans $\Ind(\cC'')$, 
\begin{equation}\label{indsh4f}
\rI(\psi\circ \phi)\stackrel{\sim}{\rightarrow} \rI\psi \circ \rI\phi.
\end{equation}

\addtocounter{subsubsection}{1}
\addtocounter{equation}{1}

\subsection{}\label{indsh7}
Soit $\cC$ une $\mU$-catégorie dans laquelle les petites limites inductives filtrantes sont représentables. Le foncteur canonique
$\iota_\cC\colon \cC\rightarrow \Ind(\cC)$ admet un adjoint à gauche 
\begin{equation}\label{indsh7a}
\kappa_\cC\colon \Ind(\cC)\rightarrow \cC.
\end{equation}
Pour toute petite catégorie filtrante $J$ et tout foncteur $\alpha\colon J\rightarrow \cC$, on a un isomorphisme 
\begin{equation}\label{indsh7b}
\kappa_\cC(\underset{\underset{J}{\longrightarrow}}{\mlq\mlq\lim \mrq\mrq} \alpha)
\stackrel{\sim}{\rightarrow} \underset{\underset{J}{\longrightarrow}}{\lim}\ \alpha.
\end{equation}
Le morphisme canonique $\kappa_\cC\circ \iota_\cC \rightarrow \id_\cC$ est un isomorphisme.

\subsection{}\label{indsh5}
Soient $\cC$ une $\mU$-catégorie, $S$ un {\em système multiplicatif à droite de morphismes de $\cC$} (\cite{ks2} 7.1.5). 
Pour tout objet $X$ de $\cC$, on désigne par $S^X$ la catégorie définie de la façon suivante. Les objets de $S^X$ sont les morphismes
$s\colon X\rightarrow X'$ de $\cC$ qui appartiennent à $S$. Soient $s\colon X\rightarrow X'$, $t\colon X\rightarrow X''$ deux objets de $S^X$. 
Un morphisme de $s\colon X\rightarrow X'$ dans $t\colon X\rightarrow X''$ est un morphisme $h\colon X'\rightarrow X''$ de $\cC$ 
tels que $t=h\circ s$. On notera qu'on ne demande pas que $h$ soit dans $S$. 
La catégorie $S^X$ est filtrante (\cite{ks2} 7.1.10). On désigne par 
\begin{equation}\label{indsh5a}
\alpha^X\colon S^X\rightarrow \cC, \ \ \ (s\colon X\rightarrow X')\mapsto X',
\end{equation}
le foncteur but. 

Pour tous objets $X,Y$ de $\cC$, l'application canonique
\begin{equation}\label{indsh5b}
\underset{\underset{(X\rightarrow X')\in S^X}{\longleftarrow}}{\lim}  \underset{\underset{(Y\rightarrow Y')\in S^X}{\longrightarrow}}{\lim} \Hom_\cC(X',Y')
\rightarrow \underset{\underset{(Y\rightarrow Y')\in S^X}{\longrightarrow}}{\lim} \Hom_\cC(X,Y')
\end{equation}
est bijective (cf. \cite{ks2} 7.1.5). 

On désigne par $\cC_S$ la catégorie ayant mêmes objets que $\cC$, définie par localisation par rapport à $S$  (\cite{ks2} 7.1.11) et par 
\begin{equation}\label{indsh5c}
Q\colon \cC\rightarrow \cC_S
\end{equation}
le foncteur canonique (cf. \cite{ks2} 7.1.16). 

Supposons que pour tout $X\in \ob(\cC)$, la catégorie $S^X$ soit essentiellement petite. Alors, $\cC_S$ est une $\mU$-catégorie (\cite{ks2} 7.1.14). 
Pour tous objets $X,Y$ de $\cC$, l'inverse de l'isomorphisme \eqref{indsh5b} induit un isomorphisme 
\begin{equation}\label{indsh5d}
\Hom_{\cC_S}(X,Y)\stackrel{\sim}{\rightarrow} \Hom_{\Ind(\cC)}(\indcolim \alpha^X,\indcolim \alpha^Y).
\end{equation}
Par suite, la correspondance 
\begin{equation}\label{indsh5e}
\alpha_S\colon \cC_S\rightarrow \Ind(\cC), \ \ \ X\mapsto \indcolim \alpha^X,
\end{equation}
définit un foncteur pleinement fidèle (\cite{ks2} 7.4.1). On notera que le triangle 
\begin{equation}\label{indsh5f}
\xymatrix{
\cC\ar[r]^Q\ar[rd]_{\iota_\cC}&{\cC_S}\ar[d]^{\alpha_S}\\
&{\Ind(\cC)}}
\end{equation}
n'est pas commutatif en général. Mais il existe un morphisme canonique
\begin{equation}\label{indsh5g}
\iota_\cC\rightarrow \alpha_S\circ Q,
\end{equation}
défini pour tout objet $X$ de $\cC$ par le morphisme canonique
\begin{equation}\label{indsh5h}
\iota_\cC(X)\rightarrow \indcolim \alpha^X=(\alpha_S\circ Q)(X).
\end{equation}

\subsection{}\label{indsh6}
Soit $\cC$ une $\mU$-catégorie abélienne. On désigne par $\hcC^\add$ la catégorie des foncteurs additifs de $\cC^\circ$ dans $\bMod(\mZ)$. 
C'est une catégorie abélienne et c'est une sous-catégorie pleine de $\hcC$ d'après (\cite{ks2} 8.1.12). Le foncteur canonique 
\begin{equation}\label{indsh6a}
\cC\rightarrow \hcC^\add, \ \ \ X\mapsto \Hom_\cC(-,X)
\end{equation}
est pleinement fidèle et exact à gauche, mais il n'est pas exact en général. 

Pour toute petite famille d'objets $(X_i)_{i\in I}$ de $\hcC^\add$, on désigne par 
$\underset{i\in I}{\indoplus} X_i$ l'objet $\indcolim(\underset{J}{\bigoplus} X_j)$ de $\hcC^\add$, où $J$ décrit les sous-ensembles finis de $I$. 
Pour tout objet $Y$ de $\cC$, on a donc un isomorphisme canonique
\begin{equation}\label{indsh6b}
\Hom_{\hcC^\add}(Y,\underset{i\in I}{\indoplus} X_i)\stackrel{\sim}{\rightarrow}\underset{i\in I}{\bigoplus}\ \Hom_{\hcC^\add}(Y,X_i).
\end{equation}

\addtocounter{subsubsection}{2}
\addtocounter{equation}{1}

On peut faire les remarques suivantes (\cite{ks2} 8.6.5):

\subsubsection{}\label{indsh6c}
La catégorie $\Ind(\cC)$ est abélienne. 

\subsubsection{}\label{indsh6d}
Le foncteur canonique $\iota_\cC\colon \cC\rightarrow \Ind(\cC)$ est pleinement fidèle et exact,
et le foncteur canonique $\Ind(\cC)\rightarrow \hcC^\add$ est pleinement fidèle et exact à gauche. 

\subsubsection{}\label{indsh6e}
La catégorie $\Ind(\cC)$ admet des petites limites inductives. De plus, les petites limites inductives filtrantes sont exactes. 

\subsubsection{}\label{indsh6f}
Les petites sommes dans $\Ind(\cC)$ sont représentables par $\indoplus$. 

\subsubsection{}\label{indsh6g}
Si les petites limites projectives dans $\cC$ sont représentables, elles le sont aussi dans $\Ind(\cC)$. 

\subsubsection{}\label{indsh6h} 
Si la catégorie $\cC$ est essentiellement petite, $\Ind(\cC)$ admet un générateur et est donc une catégorie de Grothendieck (\cite{ks2} 8.3.24).

\subsection{}\label{indsh8}
Soit $\cC$ une $\mU$-catégorie abélienne. On notera que la catégorie $\Ind(\cC)$ n'a pas assez d'injectifs (\cite{ks2}, 15.1.3). 
Selon (\cite{ks2}  15.2.1), on dit qu'un objet $F$ de $\Ind(\cC)$ est {\em quasi-injectif} si le foncteur 
\begin{equation}\label{indsh8a}
\cC\rightarrow \bMod(\mZ), \ \ \ X\mapsto F(X)=\Hom_{\Ind(\cC)}(X,F)
\end{equation}
est exact. 

\addtocounter{subsubsection}{1}
\addtocounter{equation}{1}

On peut faire les remarques suivantes:

\subsubsection{}\label{indsh8b}
Supposons que $\cC$ a assez d'injectifs et soit $F$ un objet de $\Ind(\cC)$. Les propriétés suivantes sont alors équivalentes (\cite{ks2}, 15.2.3):
\begin{itemize}
\item[(i)] $F$ est quasi-injectif;
\item[(ii)] il existe une petite catégorie filtrante $J$ et un foncteur $\alpha\colon J\rightarrow \cC$ tel que $F=\indcolim \alpha$
et que $\alpha(j)$ soit injectif pour tout $j\in \ob(J)$;
\item[(iii)] tout morphisme $X\rightarrow F$, où $X$ est un objet de $\cC$, se factorise à travers un objet injectif $Y$ de $\cC$.
\end{itemize}

En particulier, si $\cI$ désigne la sous-catégorie pleine des objets injectifs de $\cC$, $\Ind(\cI)$ s'identifie canoniquement 
à la sous-catégorie pleine des objets quasi-injectifs de $\Ind(\cC)$ \eqref{indsh4d}.

\subsubsection{}\label{indsh8c}
Si $\cC$ a assez d'injectifs, $\Ind(\cC)$ a assez de quasi-injectifs, c'est-à-dire que la sous-catégorie pleine 
des objets quasi-injectifs est cogénératrice (\cite{ks2}, 15.2.7).

\subsection{}\label{indsh9}
Soit $\phi\colon \cC\rightarrow \cC'$ un foncteur exact à gauche entre $\mU$-catégories abéliennes. 
Le foncteur $\rI\phi\colon \Ind(\cC)\rightarrow \Ind(\cC')$  associé à $\phi$ est exact à gauche (\cite{ks2} 8.6.8). 
Supposons qu'il existe une sous-catégorie $\phi$-injective $\cJ$ de $\cC$ dans le sens de (\cite{ks2} 13.3.4). 
En vertu de (\cite{ks2} 13.3.5) et avec les notations de \ref{notconv18}, le foncteur $\phi$ admet un foncteur dérivé à droite  
\begin{equation}\label{indsh9a}
\rR\phi\colon \bD^+(\cC)\rightarrow \bD^+(\cC'). 
\end{equation}
Pour tout entier $i$, on désigne par $\rR^i\phi\colon \cC\rightarrow \cC'$ le $i$-ème foncteur dérivé à droite de $\phi$ et par 
\begin{equation}\label{indsh9b}
\rI(\rR^i\phi)\colon \Ind(\cC)\rightarrow \Ind(\cC')
\end{equation}
le foncteur associé. D'après (\cite{ks2} 15.3.2), la catégorie $\Ind(\cJ)$ est $\rI\phi$-injective et le foncteur $\rI\phi$ admet un foncteur 
dérivé à droite 
\begin{equation}\label{indsh9c}
\rR(\rI\phi)\colon \bD^+(\Ind(\cC))\rightarrow \bD^+(\Ind(\cC')).
\end{equation}
De plus, pour tout entier $i$, on a un isomorphisme canonique
\begin{equation}\label{indsh9e}
\rR^i(\rI\phi)\stackrel{\sim}{\rightarrow}\rI(\rR^i\phi).
\end{equation}
En particulier, $\rR^i(\rI\phi)$ commute aux petites limites inductives filtrantes. 
Par ailleurs, comme le foncteur canonique $\iota_\cC\colon \cC\rightarrow \bInd(\cC)$ 
est exact (\ref{indsh3e} et \ref{indsh3g}) et que $\iota_\cC(\cJ)\subset \bInd(\cJ)$, le diagramme 
\begin{equation}\label{indsh9d}
\xymatrix{
{\bD^+(\cC)}\ar[r]^-(0.5){\rR\phi}\ar[d]_{\iota_\cC}&{\bD^+(\cC')}\ar[d]^{\iota_{\cC'}}\\
{\bD^+(\Ind(\cC))}\ar[r]^-(0.5){\rR(\rI\phi)}&{\bD^+(\Ind(\cC'))}}
\end{equation}
est commutatif à isomorphisme canonique près.

\subsection{}\label{indsh90}
Soit $\phi\colon \cC\rightarrow \cC'$ un foncteur exact à gauche entre $\mU$-catégories abéliennes tel que $\cC$ ait assez d'injectifs. 
On désigne par $\cI$ la sous-catégorie pleine de $\cC$ formée des objets injectifs.
D'après (\cite{ks2} 13.3.6(iii)), $\cI$ est $\phi$-injective. En vertu de (\cite{ks2} 13.3.5), le foncteur $\phi$ admet donc un foncteur dérivé à droite
\begin{equation}\label{indsh90a}
\rR\phi\colon \bD^+(\cC)\rightarrow \bD^+(\cC'). 
\end{equation}
La catégorie $\Ind(\cI)$ s'identifie à la sous-catégorie pleine de $\Ind(\cC)$ formée des objets quasi-injectifs \eqref{indsh8}. 
En vertu de (\cite{ks2} 15.3.2), $\Ind(\cI)$ est $\rI\phi$-injective et le foncteur $\rI\phi$ admet un foncteur dérivé à droite
\begin{equation}\label{indsh90b}
\rR(\rI\phi)\colon \bD^+(\Ind(\cC))\rightarrow \bD^+(\Ind(\cC')).
\end{equation}

\subsection{}\label{indsh37}
Soient $I$, $J$, $K$ trois catégories, $\varphi\colon I\rightarrow K$, $\psi\colon J\rightarrow K$ deux foncteurs. 
On désigne par $L$ la catégorie des triplets $(i,j,u)$ où $i\in \ob(I)$, $j\in \ob(J)$ et $u\in \Hom_K(\varphi(i),\psi(j))$, 
avec les morphismes évidents (cf. \cite{ks2} 3.4.1). 
On a alors deux foncteurs canoniques  $\tts \colon  L\rightarrow I$ et $\ttt\colon L\rightarrow J$ 
et un morphisme canonique de foncteurs $\varphi\circ \tts\rightarrow \psi\circ \ttt$.
Si les catégories $I$ et $J$ sont petites, il en est de même de $L$. 
Si les catégories $I$ et $J$ sont filtrantes et si le foncteur $\psi$ est cofinal, alors la catégorie $L$ est filtrante   
et les foncteurs $\tts$ et $\ttt$ sont cofinaux en vertu de (\cite{ks2} 3.4.5).

\section{Ind-modules}\label{indmod}

\subsection{}\label{indsh15}
Soit $(X,A)$ un $\mU$-topos annelé. 
On note $\bMod(A)$ la catégorie des $A$-modules de $X$. 
C'est une $\mU$-catégorie abélienne (\cite{sga4} II 4.11). 
Les petites limites projectives (resp. injectives) sont représentables dans $\bMod(A)$ (\cite{sga4} II 4.1),
et les petites limites inductives filtrantes sont exactes (\cite{sga4} III 4.3(4)). 
On désigne par $\bIndMod(A)$ la catégorie des {\em ind-$A$-modules}, c'est-à-dire des ind-objets de $\bMod(A)$ \eqref{indsh3}. 
C'est une $\mU$-catégorie abélienne \eqref{indsh6c}, $A(X)$-additive \eqref{indsh1c}. 
Les petites limites inductives (resp. projectives) sont représentables dans $\bIndMod(A)$, 
et les petites limites inductives filtrantes sont exactes \eqref{indsh6e}. On dispose des foncteurs
\begin{eqnarray}
\iota_A\colon \bMod(A)\rightarrow \bIndMod(A),\label{indsh15a}\\
\kappa_A \colon \bIndMod(A)\rightarrow \bMod(A),\label{indsh15b}
\end{eqnarray}
définis dans \eqref{indsh3a} et \eqref{indsh7a}, respectivement. 
Le foncteur $\iota_A$ est exact et pleinement fidèle et il commute aux petites limites projectives (cf. \ref{indsh3e} et \ref{indsh3g}).
Le foncteur $\kappa_A$ est un adjoint à gauche de $\iota_A$, et
le morphisme canonique $\kappa_A\circ \iota_A \rightarrow \id_{\bMod(A)}$ est un isomorphisme \eqref{indsh7}.

Lorsqu'il n'y a aucun risque d'ambiguïté, on identifiera $\bMod(A)$ à une sous-catégorie pleine de $\bIndMod(A)$ 
par le foncteurs $\iota_A$ qu'on omettra des notations.

Soient $B$, $C$ deux $A$-algèbres. On définit le bifoncteur
\begin{equation}\label{indsh15c}
\otimes_A\colon \bIndMod(B)\times \bIndMod(C)\rightarrow \bIndMod(B\otimes_AC)
\end{equation}
en posant pour tout ind-$B$-module $F$ et tout ind-$C$-module $G$, 
\begin{equation}\label{indsh15d}
F\otimes_A G= 
\underset{\underset{\underset{(N\rightarrow G)\in\bMod(C)_{/G}}{(M\rightarrow F)\in\bMod(B)_{/F}}}{\longrightarrow}}{\mlq\mlq\lim \mrq\mrq} M\otimes_AN, 
\end{equation}
où $M\otimes_AN$ est le produit tensoriel dans $\bMod(B\otimes_AC)$, qui est bien défini compte tenu de \ref{indsh3d}. 
Ce bifoncteur prolonge clairement le produit tensoriel des modules. 

On a un isomorphisme canonique de $\bMod(B\otimes_AC)$
\begin{equation}\label{indsh15g}
\kappa_{B\otimes_AC}(F\otimes_AG)\stackrel{\sim}{\rightarrow}\kappa_B(F)\otimes_A\kappa_C(G).
\end{equation}

Le foncteur canonique d'oubli $\bMod(B)\rightarrow \bMod(A)$ induit un foncteur additif
\begin{equation}\label{indsh15f}
\bIndMod(B)\rightarrow \bIndMod(A).
\end{equation}
On considérera tout ind-$B$-module comme un ind-$A$-module via ce foncteur, qu'on omettra des notations. 
Pour tout ind-$A$-module $F$ et tout ind-$B$-module $G$, on a un isomorphisme canonique
\begin{equation}\label{indsh15h}
\Hom_{\bIndMod(B)}(B\otimes_AF,G)\stackrel{\sim}{\rightarrow}\Hom_{\bIndMod(A)}(F,G). 
\end{equation}

D'après \ref{indsh4}(ii) appliqué au foncteur \eqref{indsh15f}, pour tous ind-$B$-modules $F$, $G$, on a un morphisme canonique de $\bIndMod(B\otimes_AB)$
\begin{equation}\label{indsh15i}
F\otimes_AG\rightarrow F\otimes_BG. 
\end{equation}

\begin{lem}\label{indsh17}
Sous les hypothèses de \ref{indsh15}, le foncteur $\kappa_A$ \eqref{indsh15b} est exact et il commute aux petites limites inductives. 
\end{lem}
En effet, $\kappa_A$ commute aux petites limites inductives puisqu'il admet un adjoint à droite. 
Soit 
\begin{equation}\label{indsh17a}
0\rightarrow F'\rightarrow F\rightarrow F''\rightarrow 0
\end{equation}
une suite exacte de $\bIndMod(A)$. D'après (\cite{ks2} 8.6.6), il existe une petite catégorie filtrante $J$ et une suite exacte de foncteurs de $J$ dans $\bMod(A)$ 
\begin{equation}\label{indsh17b}
0\rightarrow \varphi'\rightarrow \varphi\rightarrow \varphi''\rightarrow 0
\end{equation}
qui induit la suite \eqref{indsh17a} par passage à la limite inductive dans $\bIndMod(A)$ \eqref{indsh6e}. La suite 
\begin{equation}\label{indsh17c}
0\rightarrow \underset{\underset{J}{\longrightarrow}}{\lim}\ \varphi'\rightarrow 
\underset{\underset{J}{\longrightarrow}}{\lim}\ \varphi\rightarrow 
\underset{\underset{J}{\longrightarrow}}{\lim}\ \varphi''\rightarrow 0
\end{equation}
obtenue par passage à la limite inductive dans $\bMod(A)$ est exacte (\cite{sga4} III 4.3(4)). Par suite, $\kappa_A$ est exact.

\begin{lem}\label{indsh18} 
Sous les hypothèses de \ref{indsh15}, le bifoncteur \eqref{indsh15c} commute aux petites limites inductives filtrantes en chacune de ses variables. 
En particulier, pour toutes  petites catégories filtrantes $I$ et $J$ et tous foncteurs $\alpha\colon I\rightarrow \bMod(B)$ et $\beta\colon J\rightarrow \bMod(C)$, 
on a un isomorphisme canonique
\begin{equation}\label{indsh18a}
\underset{\underset{i\in I}{\longrightarrow}}{\mlq\mlq\lim \mrq\mrq} \alpha \otimes_A
\underset{\underset{j\in J}{\longrightarrow}}{\mlq\mlq\lim \mrq\mrq} \beta \stackrel{\sim}{\rightarrow} 
\underset{\underset{i\in I ; j\in J}{\longrightarrow}}{\mlq\mlq\lim \mrq\mrq} \alpha(i)\otimes_A\beta(j), 
\end{equation}
où $\alpha(i)\otimes_A\beta(j)$ est le produit tensoriel dans $\bMod(B\otimes_AC)$. 
\end{lem}

Soit $\mV$ un univers contenant $\mU$ et $\bMod(B)$. 
On désigne par $\bMod(B)^\wedge_\mU$ (resp. $\bMod(B)^\wedge_\mV$)
la catégorie des préfaisceaux de $\mU$-ensembles (resp. $\mV$-ensembles) sur $\bMod(B)$ et par $\tth^\mV_B\colon \bMod(B)\rightarrow \bMod(B)^\wedge_\mV$ 
le foncteur canonique \eqref{indsh1a}. 
Le foncteur d'injection canonique $\bIndMod(B)\rightarrow \bMod(B)^\wedge_\mU$ est pleinement fidèle et il commute aux $\mU$-petites 
limites inductives filtrantes \eqref{indsh3f}. 
Le foncteur d'injection canonique $\bMod(B)^\wedge_\mU\rightarrow \bMod(B)^\wedge_\mV$ 
est pleinement fidèle et il commute aux $\mV$-petites limites inductives (\cite{sga4} I 3.6). 

Soit $G\in \bIndMod(C)$. Considérons le foncteur
\begin{equation}
\varphi\colon \bIndMod(B)\rightarrow \bIndMod(B\otimes_AC),\ \ \ F\mapsto F\otimes_AG,
\end{equation}
et notons $\phi$ le foncteur composé
\begin{equation}
\bMod(B)\rightarrow \bIndMod(B)\stackrel{\varphi}{\rightarrow} \bIndMod(B\otimes_AC)\rightarrow \bMod(B\otimes_AC)^\wedge_\mV,
\end{equation}
où la première et la dernière flèche sont les foncteurs d'injection canoniques. 
D'après (\cite{ks2} 2.7.1), il existe un foncteur  $\Phi\colon \bMod(B)^\wedge_\mV\rightarrow \bMod(B\otimes_AC)^\wedge_\mV$, 
commutant aux $\mV$-petites limites inductives, 
et un isomorphisme $\Phi\circ \tth^\mV_B\stackrel{\sim}{\rightarrow} \phi$. Par ailleurs, il résulte de la preuve de 
{\em loc. cit.} et de \ref{indsh3d} que le diagramme 
\begin{equation}
\xymatrix{
{\bIndMod(B)}\ar[r]^-(0.5){\varphi}\ar[d]&{\bIndMod(B\otimes_AC)}\ar[d]\\
{\bMod(B)^\wedge_\mV}\ar[r]^-(0.5){\Phi}&{\bMod(B\otimes_AC)^\wedge_\mV}}
\end{equation}
où les flèches verticales sont les foncteurs d'injection canoniques, est commutatif à isomorphisme près. Par suite, 
$\varphi$ commute aux $\mU$-petites limites inductives filtrantes. La proposition s'ensuit puisque le produit tensoriel \eqref{indsh15c} est symétrique. 

\begin{lem}\label{indsh19}
Sous les hypothèses de \ref{indsh15}, le bifoncteur \eqref{indsh15c} est exact à droite.
\end{lem}

Cela résulte de \ref{indsh6e}, \ref{indsh18} et (\cite{ks2} 8.6.6(a)).

\subsection{}\label{indsh27}
Reprenons les hypothèses de \ref{indsh15}. 
En vertu de \ref{indsh18}, pour tous ind-$B$-modules $F$ et $F'$ et tous ind-$C$-modules $G$ et $G'$, 
on a des isomorphismes canoniques de $\bIndMod(B\otimes_AC)$
\begin{eqnarray}
(F\otimes_AG)\otimes_CG'&\stackrel{\sim}{\rightarrow}&F\otimes_A(G\otimes_CG'),\label{indsh27a}\\
(F'\otimes_BF)\otimes_AG&\stackrel{\sim}{\rightarrow}&F'\otimes_B(F\otimes_AG).\label{indsh27b}
\end{eqnarray}

D'après \ref{indsh4}(ii) et \ref{indsh18}, pour tout ind-$B$-module $F$ et tout ind-$C$-module $G$, 
l'image de $F\otimes_AG$ par le foncteur canonique $\bIndMod(B\otimes_AC)\rightarrow \bIndMod(A)$
n'est autre que le produit tensoriel de $F$ et $G$ dans $\bIndMod(A)$. 

Lorsque $A=B=C$, le produit tensoriel \eqref{indsh15c} fait de $\bIndMod(A)$ une catégorie monoïdale symétrique, ayant $A$ pour objet unité.

\begin{lem}\label{indsh33}
Soient $(X,A)$ un $\mU$-topos annelé, $B$ une $A$-algèbre, $M$, $N$ deux ind-$B$-modules, 
$u\colon M\rightarrow N$ un morphisme de $\bIndMod(A)$ \eqref{indsh15f} tels que le diagramme
\begin{equation}\label{indsh33a}
\xymatrix{
{B\otimes_AM}\ar[r]^-(0.5){\id\otimes u}\ar[d]&{B\otimes_AN}\ar[d]\\
M\ar[r]^-(0.5)u&N}
\end{equation}
où les flèches verticales sont les morphismes canoniques \eqref{indsh15i}, soit commutatif. 
Alors, il existe une petite catégorie filtrante $J$, deux foncteurs $\alpha,\beta\colon K\rightrightarrows \bMod(B)$, 
un morphisme de foncteurs $\sigma\colon \alpha\rightarrow \beta$ et deux isomorphismes 
\begin{equation}\label{indsh33b}
M\stackrel{\sim}{\rightarrow} \indcolim \alpha \ \ \ {\rm et}\ \ \ N\stackrel{\sim}{\rightarrow} \indcolim \beta, 
\end{equation}
tels que l'image de la limite inductive de $\sigma$ par le foncteur canonique $\bIndMod(B)\rightarrow \bIndMod(A)$ soit $u$. 
\end{lem}

\'Ecrivons $M=\indcolim \varphi$ et $N=\indcolim \psi$ où $\varphi \colon I\rightarrow \bMod(B)$ et $\psi \colon J\rightarrow \bMod(B)$
sont deux foncteurs tels que les catégories $I$ et $J$ soient petites et filtrantes. 
Notons $\tvarphi \colon I\rightarrow \bMod(A)_{/M}$, $\tpsi\colon J\rightarrow \bMod(A)_{/N}$ et $\tu\colon \bMod(A)_{/M} \rightarrow \bMod(A)_{/N}$
les foncteurs induits par $\varphi$, $\psi$ et $u$. 
On désigne par $K'$ la catégorie des triplets $(i,j,v)$ où $i\in \ob(I)$, $j\in \ob(J)$ et $v\in \Hom_{\bMod(A)_{/N}}(\tu\circ \tvarphi(i),\tpsi(j))$, 
avec les morphismes évidents \eqref{indsh37}. Concrètement, les objets de $K'$ sont les triplets $(i,j,v)$ où $i\in \ob(I)$, $j\in \ob(J)$ 
et $v\colon \varphi(i)\rightarrow \psi(j)$ est un morphisme $A$-linéaire tel que le diagramme de $\bIndMod(A)$
\begin{equation}\label{indsh33c}
\xymatrix{
{\varphi(i)}\ar[r]^-(0.5)v\ar[d]&{\psi(j)}\ar[d]\\
M\ar[r]^-(0.5)u&N}
\end{equation}
où les flèches verticales sont les morphismes canoniques, soit commutatif.

Pour tout $A$-module $F$, on a  
\begin{equation}\label{indsh33d}
\Hom_{\bIndMod(A)}(F,N)=\underset{\underset{j\in J}{\longrightarrow}}{\lim}\ \Hom_{\bMod(A)}(F,\psi(j)).
\end{equation}
Comme la catégorie $J$ est filtrante, on en déduit que le foncteur $\tpsi$ est cofinal d'après (\cite{sga4} I 8.1.3(b)). 
En vertu de (\cite{ks2} 3.4.5), la catégorie $K'$ est filtrante et les projections canoniques $\alpha'\colon K'\rightarrow I$ 
et $\beta'\colon K'\rightarrow J$ sont cofinales. De plus, $K'$ est clairement petite.

On désigne par $K$ la sous-catégorie pleine de $K'$ formée des triplets $(i,j,v)$ où $i\in \ob(I)$, $j\in \ob(J)$ 
et $v\colon \varphi(i)\rightarrow \psi(j)$ est un morphisme $B$-linéaire. 
Il résulte de \eqref{indsh33a}, \eqref{indsh33d} et (\cite{sga4} I 8.1.3(c)) que le foncteur d'injection $\iota\colon K\rightarrow K'$ est cofinal et que $K$ est filtrante. 
On prend $\alpha=\varphi\circ \alpha'\circ \iota$ et $\beta=\psi\circ \beta'\circ \iota$. On a un morphisme canonique $\sigma\colon \beta\rightarrow \alpha$
dont la limite inductive s'identifie à $u$ en vertu de \eqref{indsh33c}, d'où la proposition. 

\begin{rema}\label{indsh51}
Sous les hypothèses de \ref{indsh36}, le foncteur canonique $\bIndMod(B)\rightarrow \bIndMod(A)$ étant fidèle d'après \ref{indsh4d}, 
on peut considérer $u$ canoniquement comme un morphisme de $\bIndMod(B)$.
\end{rema}

\begin{lem}\label{indsh45}
Soient $(X,A)$ un $\mU$-topos annelé, $F$ un ind-$A$-module. Alors, pour que le foncteur
\begin{equation}\label{indsh45a}
\bMod(A)\rightarrow \bIndMod(A), \ \ \ M\mapsto M\otimes_AF
\end{equation}
soit exact, il faut et il suffit que le foncteur 
\begin{equation}\label{indsh45b}
\bIndMod(A)\rightarrow \bIndMod(A), \ \ \ G\mapsto G\otimes_AF
\end{equation}
soit exact.
\end{lem}
En effet, si le foncteur \eqref{indsh45b} est exact, il en est de même du foncteur \eqref{indsh45a} puisque le foncteur canonique 
$\bMod(A)\rightarrow \bIndMod(A)$ est exact \eqref{indsh6d}. L'implication inverse résulte de \ref{indsh6e}, \ref{indsh18} et (\cite{ks2} 8.6.6(1)).

\begin{defi}\label{indsh46}
Soit $(X,A)$ un $\mU$-topos annelé. On dit qu'un ind-$A$-module est {\em plat} (ou {\em $A$-plat}) s'il vérifie les conditions équivalentes de \ref{indsh45}.
\end{defi}

Cette notion prolonge celle de platitude pour les $A$-modules (\cite{sga4} V 1.1). 
On peut faire les remarques suivantes:

\subsubsection{}\label{indsh46a}
Une petite limite inductive filtrante de ind-$A$-modules plats est un ind-$A$-module plat d'après \ref{indsh6e} et \ref{indsh18}. 

\subsubsection{}\label{indsh46b}
Soient $B$ une $A$-algèbre, $F$ un ind-$A$-module plat. Alors, le ind-$B$-module $B\otimes_AF$ est plat. 
En effet, le foncteur canonique $\bIndMod(B)\rightarrow \bIndMod(A)$ est exact et fidèle d'après \ref{indsh4d} et (\cite{ks2} 8.6.8), et pour tout $B$-module $N$, 
on a un isomorphisme canonique $N\otimes_B(B\otimes_AF)\stackrel{\sim}{\rightarrow}N\otimes_AF$ de $\bIndMod(B)$ \eqref{indsh27b}. 

\addtocounter{equation}{3}

\subsubsection{}\label{indsh46c}
Soit $B$ une $A$-algèbre telle que le foncteur 
\begin{equation}\label{indsh46d}
\bMod(A)\rightarrow \bMod(B), \ \ \ M\mapsto M\otimes_AB,
\end{equation}
est exact et fidèle; en particulier $B$ est une $A$-algèbre plate.  Le foncteur 
\begin{equation}\label{indsh46e}
\bIndMod(A)\rightarrow \bIndMod(B), \ \ \ F\mapsto F\otimes_AB,
\end{equation}
est alors exact et fidèle d'après \ref{indsh4d} et (\cite{ks2} 8.6.8).  

Pour qu'un ind-$A$-module $F$ soit plat, il faut et il suffit que le ind-$B$-module $B\otimes_AF$ soit plat. 
En effet, la condition est suffisante d'après \ref{indsh46b} et elle est nécessaire compte tenu de l'hypothèse \eqref{indsh46e} 
et du fait que pour tout $A$-module $M$, on a un isomorphisme canonique de $\bIndMod(B)$
\begin{equation}
B\otimes_A(M\otimes_AF) \stackrel{\sim}{\rightarrow} (B\otimes_AM)\otimes_B(B\otimes_AF).
\end{equation}

\subsection{}\label{indsh21}
Soit $f\colon (X,A)\rightarrow (Y,B)$ un morphisme de $\mU$-topos annelés. 
En vertu de \ref{indsh4}, les foncteurs adjoints $f^*$ et $f_*$ entre les catégories $\bMod(B)$ et $\bMod(A)$ 
induisent des foncteurs additifs 
\begin{eqnarray}
\rI f^*\colon \bIndMod(B) \rightarrow \bIndMod(A),\label{indsh21a}\\
\rI f_*\colon \bIndMod(A) \rightarrow \bIndMod(B).\label{indsh21b}
\end{eqnarray}
D'après (\cite{ks2} 8.6.8), le foncteur $\rI f^*$ (resp. $\rI f_*$) 
est exact à droite (resp. gauche). Il résulte aussitôt de \eqref{indsh1c} que le foncteur $\rI f^*$ est un adjoint à gauche du foncteur $\rI f_*$. 

Compte tenu de \ref{indsh4}(ii) et \ref{indsh18}, pour tous ind-$B$-modules $F$, $F'$, on a un isomorphisme canonique
\begin{equation}\label{indsh21h}
\rI f^*(F\otimes_B F')\stackrel{\sim}{\rightarrow} \rI f^*(F)\otimes_A\rI f^*(F').
\end{equation}

D'après \ref{indsh90}, le foncteur $\rI f_*$ admet un foncteur dérivé à droite
\begin{equation}\label{indsh21c}
\rR(\rI f_*)\colon \bD^+(\bIndMod(A))\rightarrow \bD^+(\bIndMod(B)).
\end{equation}
Pour tout entier $q$, on a un isomorphisme canonique \eqref{indsh9e}
\begin{equation}\label{indsh21g}
\rR^q(\rI f_*)\stackrel{\sim}{\rightarrow}\rI(\rR^qf_*).
\end{equation}
Le diagramme \eqref{indsh9d}
\begin{equation}\label{indsh21d}
\xymatrix{
{\bD^+(\bMod(A))}\ar[r]^-(0.5){\rR f_*}\ar[d]_{\iota_A}&{\bD^+(\bMod(B))}\ar[d]^{\iota_{B}}\\
{\bD^+(\bIndMod(A))}\ar[r]^-(0.5){\rR(\rI f_*)}&{\bD^+(\bIndMod(B))}}
\end{equation}
est commutatif à isomorphisme canonique près.

\subsection{}\label{indsh212}
Soient $f\colon (X,A)\rightarrow (Y,B)$ un morphisme de $\mU$-topos annelés,
$F$, $G$ deux ind-$A$-modules, $p$, $q$ deux entiers $\geq 0$. Il existe alors un morphisme canonique de cup-produit
\begin{equation}\label{indsh212b}
\rR^p\rI f_*(F)\otimes_B \rR^q\rI f_*(G)\rightarrow \rR^{p+q}\rI f_*(F\otimes_AG).
\end{equation}
En effet, écrivons $F=\indcolim \alpha$ et $G=\indcolim \beta$, où 
$\alpha \colon I\rightarrow \bMod(A)$ et $\beta \colon J\rightarrow \bMod(A)$
sont deux foncteurs tels que les catégories $I$ et $J$ soient petites et filtrantes. 
Compte tenu de \ref{indsh18}, \eqref{indsh21g} et \eqref{indsh4c}, on prend pour morphisme \eqref{indsh212b} la limite inductive sur $\alpha$ et $\beta$ 
des morphismes de cup-produit 
\begin{equation}\label{indsh212a}
\rR^pf_*(\alpha(i))\otimes_B\rR^qf_*(\beta(j))\rightarrow \rR^{p+q}f_*(\alpha(i)\otimes_A \beta(j)), \ \ \ (i\in \ob(I), \ j\in \ob(J)). 
\end{equation}
Les foncteurs $\talpha\colon  I\rightarrow \bMod(A)_{/F}$ et $\tbeta \colon J\rightarrow \bMod(A)_{/G}$ induits par $\alpha$
et $\beta$ étant cofinaux d'après \eqref{indsh1b} et (\cite{sga4} I 8.1.3(b)), remplaçant $I$ (resp. $J$) par $\bMod(A)_{/F}$ (resp. $\bMod(A)_{/G}$),  
on en déduit que \eqref{indsh212b} ne dépend pas du choix de $\alpha$ et $\beta$. 

Il résulte aussitôt de l'associativité du cup-produit des modules (\cite{sp} \href{https://stacks.math.columbia.edu/tag/0FP4}{0FP4}) 
que pour tous ind-$A$-modules $F$, $G$, $H$ et tous entiers $p,q,r\geq 0$, le diagramme 
\begin{equation}\label{indsh212c}
\xymatrix{
{\rR^p\rI f_*(F) \otimes_B \rR^q\rI f_*(G) \otimes_B \rR^r\rI f_*(H)}\ar[r]\ar[d]&{\rR^{p+q}\rI f_*(F\otimes_A G) \otimes_B \rR^r\rI f_*(H)}\ar[d]\\
{\rR^p\rI f_*(F) \otimes_B \rR^{q+r}\rI f_*(G\otimes_A H)}\ar[r]&{\rR^{p+q+r}\rI f_*(F\otimes_A G\otimes_AH)}}
\end{equation}
défini par \eqref{indsh27a} et les morphismes de cup-produit \eqref{indsh212b}, est commutatif. 

\subsection{}\label{indsh211}
Soit $f\colon (X,A)\rightarrow (Y,B)$ un morphisme de $\mU$-topos annelés. 
Notons $\cF_A$ (resp. $\cF_B$) la sous-catégorie pleine de $\bMod(A)$ (resp. $\bMod(B)$) formée 
des modules flasques (\cite{sga4} V 4.1). 
Les $A$-modules injectifs étant flasques, la catégorie $\cF_A$ est cogénératrice de $\bIndMod(A)$. 
D'après (\cite{sga4} V 5.2), les objets de $\cF_A$ sont acycliques pour $f_*$.
Par suite, en vertu de (\cite{ks2} 13.3.8), la catégorie $\cF_A$ est $f_*$-injective dans le sens de (\cite{ks2} 13.3.4).
On en déduit que la catégorie $\Ind(\cF_A)$ est $\rI f_*$-injective d'après (\cite{ks2} 15.3.2). 
Par ailleurs, comme $f_*(\cF_A)\subset \cF_B$ (\cite{sga4} V 4.9), on a $\rI f_*(\Ind(\cF_A))\subset \Ind(\cF_B)$. 

Soit $g\colon (Y,B)\rightarrow (Z,C)$ un morphisme de $\mU$-topos annelés. Posons $h=gf\colon (X,A)\rightarrow (Z,C)$.
En vertu de ce qui précède et de (\cite{ks2} 13.3.13), on a un isomorphisme canonique
\begin{equation}\label{indsh211a}
\rR (\rI h_*)\stackrel{\sim}{\rightarrow} \rR (\rI g_*)\circ \rR(\rI f_*).
\end{equation}

\subsection{}\label{indsh42}
Soit 
\begin{equation}\label{indsh42a}
\xymatrix{
{(X',A')}\ar[r]^-(0.5){g'}\ar[d]_{f'}&{(X,A)}\ar[d]^f\\
{(Y',B')}\ar[r]^-(0.5)g&{(Y,B)}}
\end{equation}
un diagramme de morphismes de topos annelés, 
commutatif à isomorphisme canonique près; autrement dit, on a un isomorphisme
\begin{equation}\label{indsh42b}
f_*g'_*\stackrel{\sim}{\rightarrow}g_*f'_*
\end{equation}
et le diagramme 
\begin{equation}\label{indsh42c}
\xymatrix{
{f_*(A)}\ar[d]_{f_*(g'^\#)}&B\ar[r]^-(0.5){g^\#}\ar[l]_-(0.5){f^\#}&{g_*(B')}\ar[d]^{g_*(f'^\#)}\\
{f_*(g'_*(A'))}\ar[rr]&&{g_*(f'_*(A'))}}
\end{equation}
où la flèche non labellisée est l'isomorphisme \eqref{indsh42b}, est commutatif. 
Si $F$ est un ind-$A$-module, on a, pour tout $q\geq 0$, un morphisme canonique de $\bIndMod(B')$
\begin{equation}\label{indsh42d}
\rI g^*(\rR^q\rI f_*(F))\rightarrow \rR^q\rI f'_*(\rI g'^*(F)),
\end{equation}
dit {\em morphisme de changement de base}. En effet, cela revient à donner un morphisme
\begin{equation}\label{indsh42e}
\rR^q\rI f_*(F)\rightarrow \rI g_*(\rR^q\rI f'_*(\rI g'^*(F))),
\end{equation}
et on prend le morphisme composé 
\begin{eqnarray}\label{indsh42f}
\lefteqn{\rR^q\rI f_*(F)\rightarrow \rR^q\rI f_*(\rI g'_*(\rI g'^*(F)))\rightarrow} \\
&&\rR^q\rI (fg')_*(\rI g'^*(F))\stackrel{\sim}{\rightarrow} \rR^q\rI (g f')_*(\rI g'^*(F))
\rightarrow \rI g_*(\rR^q \rI f'_*(g'^*(F))),\nonumber
\end{eqnarray}
où la première flèche provient du morphisme d'adjonction
$\id\rightarrow \rI g'_*\rI g'^*$, la deuxième et la dernière flèche de la suite spectrale de 
Cartan-Leray \eqref{indsh211a} et la troisième flèche de l'isomorphisme \eqref{indsh42b}.

Compte tenu de \eqref{indsh21g}, on vérifie aussitôt que le morphisme de changement de base \eqref{indsh42d}
est la limite inductive des morphismes de changement de base pour les modules. 
Par suite, \eqref{indsh42d} est un isomorphisme lorsque le morphisme de changement de base pour les modules est un isomorphisme. 
En particulier, \eqref{indsh42d} est un isomorphisme si $g$ est le morphisme de localisation de $Y$
en un objet $V$, $B'=g^{-1}(B)$, $g'$ est le morphisme de localisation de $X$ en $U=f^*(V)$, $A'=g'^{-1}(A)$ et $f'=f_{|V}$ (\cite{sga4} 5.1(3)).

\begin{lem}\label{indsh47}
Soient $f\colon (X,A)\rightarrow (Y,B)$ un morphisme de $\mU$-topos annelés, $E$ un $B$-module localement libre de type fini, $\cF$ un objet de $\bD^+(\bIndMod(A))$.
Il existe alors un isomorphisme canonique fonctoriel en $\cF$ de  $\bD^+(\bIndMod(B))$
\begin{equation}\label{indsh47a}
E\otimes_A\rR \rI f_*(\cF)\rightarrow \rR \rI f_*(f^*(E)\otimes_B\cF).
\end{equation}
\end{lem}

On notera d'abord que le foncteur $E\otimes_A-$ étant exact sur la catégorie des ind-$A$-modules \eqref{indsh45}, 
il se prolonge à la catégorie $\bD^+(\bIndMod(A))$. On définit de même  le foncteur $f^*(E)\otimes_B-$ sur la catégorie $\bD^+(\bIndMod(B))$. 

On désigne par $\cI$ la catégorie des $A$-modules injectifs, 
et par $\bK^+(\bIndMod(A))$ (resp. $\bK^+(\Ind(\cI))$) la catégorie des complexes bornés inférieurement de ind-$A$-modules
(resp. injectifs) à homotopie près (\cite{ks2} 11.3.7). 
On notera que la catégorie $\Ind(\cI)$ s'identifie à la sous-catégorie pleine de $\bIndMod(A)$ formée des objets quasi-injectifs \eqref{indsh8}. 
Comme $\cI$ est $f_*$-injective d'après (\cite{ks2} 13.3.6(iii)), $\Ind(\cI)$ est $(\rI f_*)$-injective en vertu de (\cite{ks2} 15.3.2), autrement dit, on a 
les propriétés suivantes (cf. \cite{ks2} 10.3.2 et 13.3.4):
\begin{itemize}
\item[{\rm (i)}] pour tout complexe borné inférieurement de ind-$A$-modules $\cF$, il existe un objet $\cG$ de $\bK^+(\Ind(\cI))$
et un quasi-isomorphisme $\cF\rightarrow \cG$;
\item[{\rm (ii)}] pour tout complexe exact $\cG$ de $\bK^+(\Ind(\cI))$, le complexe $\rI f_*(\cG)$ est exact. 
\end{itemize}
Notant $\bN^+$ le système ``nul'' des complexes exacts de $\bK^+(\bIndMod(A))$ (cf. \cite{ks2} 13.1.2), le foncteur $\rR \rI f_*$ est alors défini par le diagramme suivant:
\[
\xymatrix{
{\bK^+(\bIndMod(A))}\ar[r]&{\bD^+(\bIndMod(A))}\ar[rd]^{\rR \rI f_*}&\\
{\bK^+(\Ind(\cI))}\ar[r]\ar[u]\ar@/_2pc/[rr]_{\rI f_*}&{\bK^+(\Ind(\cI))/(\bK^+(\Ind(\cI))\cap \bN^+)}\ar[r]\ar[u]^-(0.4)*[@]{\sim}&{\bD^+(\bIndMod(B))}}
\]

Pour tout objet $\cG$ de $\bK^+(\Ind(\cI))$, $f^*(E)\otimes_A\cG$ est un objet de $\bK^+(\Ind(\cI))$ d'après 
(\cite{sp} \href{https://stacks.math.columbia.edu/tag/01E7}{01E7}). 
Par ailleurs, il résulte de \eqref{indsh4c} qu'on a un isomorphisme canonique, fonctoriel en $\cG$,
\begin{equation}\label{indsh47b}
E\otimes_A \rI f_*(\cG)\rightarrow \rI f_*(f^*(E)\otimes_B\cG).
\end{equation}
La proposition s'ensuit.

\subsection{}\label{indsh36}
Soient $(X,A)$ un topos annelé, $U$ un objet de $X$. 
On désigne par $j_U\colon X_{/U}\rightarrow X$ le morphisme de localisation de $X$ en $U$.
Pour tout $F\in \ob(X)$, le faisceau $j_U^*(F)$ sera aussi noté $F|U$. Le topos $X_{/U}$ sera annelé par $A|U$. 
D'après (\cite{ks2} 8.6.8),  $j_U^*$ induit un foncteur exact 
\begin{equation}\label{indsh36a}
\rI j_U^*\colon \bIndMod(A)\rightarrow  \bIndMod(A|U).
\end{equation} 
Pour tout ind-$A$-module $F$, le ind-$(A|U)$-module $\rI j_U^*(F)$ sera aussi noté $F|U$. 

Le foncteur prolongement par zéro $j_{U!}\colon \bMod(A|U)\rightarrow \bMod(A)$ est exact et fidèle (\cite{sga4} IV 11.3.1). 
D'après \ref{indsh4d} et (\cite{ks2} 8.6.8), il induit un foncteur exact et fidèle 
\begin{equation}\label{indsh36b}
\rI j_{U!}\colon \bIndMod(A|U)\rightarrow \bIndMod(A).
\end{equation}
C'est un adjoint à gauche du foncteur $\rI j_U^*$.

\begin{lem}\label{indsh48}
Soient $(X,A)$ un topos annelé, $U$ un objet de $X$.  Alors,
\begin{itemize}
\item[{\rm (i)}] Pour tout ind-$(A|U)$-module plat $P$ \eqref{indsh46}, le ind-$A$-module $\rI j_{U!}(P)$ est plat. 
\item[{\rm (ii)}] Pour tout ind-$A$-module plat $M$, le ind-$(A|U)$-module  $\rI j_U^*(M)$ est plat. 
\end{itemize}
\end{lem}

En effet, il résulte aussitôt de (\cite{sga4} IV 12.11) que pour tout ind-$A$-module $M$
et tout ind-$(A|U)$-module $P$, on a un isomorphisme canonique fonctoriel
\begin{equation}\label{indsh48a}
\rI j_{U!}(P\otimes_{(A|U)}\rI j_U^*(M))\stackrel{\sim}{\rightarrow} \rI j_{U!}(P)\otimes_AM.
\end{equation}

(i) Cela résulte de \eqref{indsh48a} et du fait que les foncteurs $\rI j_U^*$ et $\rI j_{U!}$ sont exacts. 

(ii) Il résulte de \eqref{indsh48a} et du fait que le foncteur $\rI j_{U!}$ est exact que le foncteur 
\begin{equation}
\bIndMod(A|U)\rightarrow \bIndMod(A), \ \ \ P\mapsto \rI j_{U!}(P\otimes_{(A_\mQ|U)}\rI j_U^*(M))
\end{equation}
est exact. Comme le foncteur $\rI j_{U!}$ est de plus fidèle \eqref{indsh36}, on en déduit que
le foncteur $P\mapsto P\otimes_{(A_\mQ|U)}\rI j_U^*(M)$ sur la catégorie $\bIndMod(A|U)$ est exact; d'où la proposition.

\begin{lem}\label{indsh39}
Soient $(X,A)$ un topos annelé, $(U_p)_{1\leq p\leq n}$ un recouvrement 
fini de l'objet final de $X$, $F$, $G$ deux ind-$A$-modules. Pour tous $1\leq p,q\leq n$, on pose $U_{pq}=U_p\times U_q$. 
Reprenons les notations de \ref{indsh36}. Alors,
\begin{itemize}
\item[{\rm (i)}] Le diagramme d'applications d'ensembles
\begin{eqnarray}\label{indsh39a}
&&\Hom_{\bIndMod(A)}(F,G)\rightarrow\prod_{1\leq p\leq n}\Hom_{\bIndMod(A|U_p)}(F|U_p,G|U_p)\\
&&\rightrightarrows \prod_{1\leq p,q\leq n}\Hom_{\bIndMod(A|U_{p,q})}(F|U_{pq},G|U_{pq})\nonumber
\end{eqnarray}
est exact.
\item[{\rm (ii)}] Pour que $F$ soit nul, il faut et il suffit que pour tout $1\leq p\leq n$, $F|U_p$ soit nul.
\item[{\rm (iii)}] Pour qu'un morphisme $u\colon F\rightarrow G$ de $\bIndMod(A)$ soit un isomorphisme, 
il faut et il suffit que pour tout $1\leq p\leq n$, $u|U_p$ soit un isomorphisme. 
\item[{\rm (iv)}] Pour que $F$ soit $A$-plat \eqref{indsh46}, il faut et il suffit que pour tout $1\leq p\leq n$, $F|U_p$ soit $(A|U_p)$-plat. 
\end{itemize}
\end{lem}

(i) \'Ecrivons $F=\indcolim \alpha$ et $G=\indcolim \beta$ où 
$\alpha \colon I\rightarrow \bMod(A)$ et $\beta \colon J\rightarrow \bMod(A)$
sont deux foncteurs tels que les catégories $I$ et $J$ soient petites et filtrantes. 
Soit $u\colon F\rightarrow G$ un morphisme de $\bIndMod(A)$. 
Compte tenu de \eqref{indsh1c}, pour tout $i\in \ob(I)$, il existe $j\in \ob(J)$ et un morphisme $u_{ij}\colon \alpha(i)\rightarrow \beta(j)$ 
qui représente l'image de $u$ dans 
\begin{equation}
\underset{\underset{j\in J}{\longrightarrow}}{\lim}\ \Hom_{\bMod(A)}(\alpha(i),\beta(j)).
\end{equation}
On dira que $u_{ij}$ {\em représente} $u$. 

Soient $u,v\colon F\rightarrow G$ deux morphismes de $\bIndMod(A)$ tels que pour tout $1\leq p\leq n$, on ait $u|U_p=v|U_p$. 
Compte tenu de \eqref{indsh1b}, 
pour tout $i\in \ob(I)$, il existe $j\in \ob(J)$ et des représentants  $u_{ij},v_{ij}\colon \alpha(i)\rightarrow \beta(j)$ de $u$ et $v$, respectivement, 
tels que $u_{ij}|U_p=v_{ij}|U_p$ pour tout $1\leq p\leq n$. On en déduit que $u_{ij}=v_{ij}$ et par suite que $u=v$ \eqref{indsh1c}.  
La première application de \eqref{indsh39a} est donc injective. 

Soient, pour tout $1\leq p\leq n$, $u_p\colon F|U_p\rightarrow G|U_p$ un morphisme $\bIndMod(A|U_p)$
tels que $(u_p)_{1\leq p\leq n}$ soit dans le noyau de la double flèche de \eqref{indsh39a}. 
Compte tenu de \eqref{indsh1b}, pour tout $i\in \ob(I)$, 
il existe $j\in \ob(J)$ et pour tout $1\leq p\leq n$ un représentant $u_{ijp}\colon \alpha(i)|U_p\rightarrow \beta(j)|U_p$ de $u_p$,
tels que $u_{ijp}|U_{pq}=u_{ijq}|U_{pq}$ pour tous $1\leq p,q\leq n$. Il existe donc $u_{ij}\colon \alpha(i)\rightarrow \beta(j)$
tel que $u_{ij}|U_p=u_{ijp}$ pour tout $1\leq p\leq n$. Pour tout morphisme $h\colon i'\rightarrow i$ de $I$, l'image de $u_{ij}$ par l'application
\begin{equation}
\underset{\underset{j\in J}{\longrightarrow}}{\lim}\ \Hom_{\bMod(A)}(\alpha(i),\beta(j))
\rightarrow \underset{\underset{j\in J}{\longrightarrow}}{\lim}\ \Hom_{\bMod(A)}(\alpha(i'),\beta(j))
\end{equation}
induite par $h$ coïncide avec celle de $u_{i'j}$; on le vérifie par restriction aux $U_p$ pour $1\leq p\leq n$. 
Les morphismes $u_{ij}$ définissent alors un morphisme $u\colon F\rightarrow G$ de $\bIndMod(A)$ \eqref{indsh1b} d'image $(u_p)_{1\leq p\leq n}$, 
d'où l'exactitude au centre de \eqref{indsh39a}. 

(ii) En effet, $F$ est nul si et seulement si $\id_F=0$. L'assertion résulte donc de (i). 

(iii) Supposons que pour tout $1\leq p\leq n$, $u|U_p$ soit un isomorphisme et soit $v_p\colon G|U_p\rightarrow F|U_p$ son inverse. 
Comme $(v_p)_{1\leq p\leq n}$ est clairement dans le noyau de la double flèche de \eqref{indsh39a}, il est l'image d'un unique morphisme
$u\colon G\rightarrow F$ de $\bIndMod(A)$. Compte tenu de l'injectivité de la première flèche de \eqref{indsh39a}, $v$ est l'inverse de $u$.  

(iv) En effet, la condition est nécessaire en vertu de \ref{indsh48}(ii) et elle est suffisante compte tenu de (iii) et \eqref{indsh21h}. 

\subsection{}\label{indsh40}
Soit $f\colon (X,A)\rightarrow (Y,B)$ un morphisme de $\mU$-topos annelés, $F$ un ind-$A$-module, $G$ un ind-$B$-module, $q$ un entier $\geq 0$. 
Le morphisme d'adjonction $G\rightarrow \rI f_*(\rI f^*(G))$ et le cup-produit \eqref{indsh212b} induisent un morphisme bifonctoriel
\begin{equation}\label{indsh40a}
G\otimes_A\rR^q(\rI f_*)(F)\rightarrow \rR^q(\rI f_*)(\rI f^*(G)\otimes_AF). 
\end{equation}
On peut faire les remarques suivantes~:
\begin{itemize}
\item[(i)] Pour tout ind-$B$-module $G'$, le composé 
\begin{equation}\label{indsh40b}
\xymatrix{
{G\otimes_BG'\otimes_B\rR^q(\rI f_*)(F)}\ar[r]\ar[rd]&
{G\otimes_B\rR^q(\rI f_*)(\rI f^*(G')\otimes_AF)}\ar[d]\\ 
&{\rR^q(\rI f_*)(\rI f^*(G\otimes_BG')\otimes_AF)}}
\end{equation}
des morphismes induits par les morphismes \eqref{indsh40a} relatifs à $G$ et $G'$, 
n'est autre que le morphisme \eqref{indsh40a} relatif à $B\otimes_BG'$. Cela résulte aussitôt de \eqref{indsh212c}. 
\item[(ii)] Lorsque $q=0$, le morphisme \eqref{indsh40a} est l'adjoint du morphisme composé 
\begin{equation}\label{indsh40c}
\rI f^*(G\otimes_B\rI f_*(F))\stackrel{\sim}{\rightarrow} \rI f^*(G)\otimes_A\rI f^*(\rI f_*(F))\rightarrow \rI f^*(G)\otimes_AF,
\end{equation}
où la première flèche est l'isomorphisme \eqref{indsh21h} et 
la seconde flèche est induite par le morphisme canonique $\rI f^*(\rI f_*(F))\rightarrow F$.
Cela résulte de l'énoncé analogue pour les modules (\cite{sp} \href{https://stacks.math.columbia.edu/tag/0B68}{0B68}).
\end{itemize}

\begin{lem}\label{indsh41}
Soient $f\colon (X,A)\rightarrow (Y,B)$ un morphisme de $\mU$-topos annelés, $N$ un $B$-module
localement projectif de type fini \eqref{notconv14}, $F$ un ind-$A$-module, $q$ un entier $\geq 0$.
Supposons l'objet final de $Y$ quasi-compact. Alors, le morphisme canonique \eqref{indsh40a} 
\begin{equation}\label{indsh41a}
N\otimes_B\rR^q\rI f_*(F)\rightarrow\rR^q(\rI f_*)(\rI f^*(N)\otimes_AF)
\end{equation}
est un isomorphisme.
\end{lem}
Il existe un recouvrement fini $(U_i)_{0\leq i\leq n}$ de l'objet final de $Y$ tel que pour tout $0\leq i\leq n$, 
le  $(B|U_i)$-module $N|U_i$ soit un facteur direct d'un $(B|U_i)$-module libre de type fini. 
Compte tenu de \ref{indsh39}(iii) et du fait que $\rR^q\rI f_*$ commute aux localisations \eqref{indsh42}, 
on peut alors se borner au cas où $N$ est un facteur direct d'un $B$-module libre de type fini, et même au cas où $N$ est un 
$B$-module libre de type fini, auquel cas l'assertion est évidente.

\section{\texorpdfstring{Ind-modules de Higgs et $\lambda$-connexions}{Ind-modules de Higgs et lambda-connexions}}

\begin{defi}\label{indsh30}
Soient $(X,A)$ un topos annelé, $E$ un $A$-module. 
\begin{itemize}
\item[(i)]
On appelle {\em ind-$A$-module de Higgs à coefficients dans $E$}
un couple $(M,\theta)$ formé d'un ind-$A$-module $M$ et d'un morphisme de $\bIndMod(A)$ 
\begin{equation}\label{indsh30a}
\theta\colon M\rightarrow M\otimes_AE
\end{equation}
tel que le morphisme composé 
\begin{equation}\label{indsh30b}
\xymatrix{M\ar[r]^-(0.5)\theta&{M\otimes_AE}\ar[rr]^-(0.5){\theta\otimes_A\id_E}&&
{M\otimes_AE\otimes_AE}\ar[rr]^-(0.5){\id_M\otimes_Aw}&&{M\otimes_A\wedge^2E}},
\end{equation}
où $w\colon E\otimes_AE\rightarrow \wedge^2E$ est le produit extérieur \eqref{notconv9}, soit nul. 
On dit alors que $\theta$ est un {\em $A$-champ de Higgs} sur $M$ à coefficients dans $E$. 
\item[(ii)] Si $(M,\theta)$ et $(M',\theta')$ sont deux ind-$A$-modules de Higgs,
un morphisme de $(M,\theta)$ dans $(M',\theta')$ est un morphisme  
$u\colon M\rightarrow M'$ de $\bIndMod(A)$ tel que $(u\otimes\id_E)\circ \theta=\theta'\circ u$. 
\end{itemize}
\end{defi}

Les ind-$A$-modules de Higgs à coefficients dans $E$ forment une catégorie notée $\bIndMH(A,E)$. 
Le lecteur prendra garde que malgré la notation, cette catégorie n'est pas la catégorie des ind-objets d'une autre catégorie \eqref{indsh3}. 
On peut compléter la terminologie et faire les remarques suivantes.

\addtocounter{subsubsection}{2}
\addtocounter{equation}{1}

\subsubsection{}\label{indsh30c}
Soit $(M,\theta)$ un ind-$A$-module de Higgs à coefficients dans $E$. Pour tout $i\geq 1$, on désigne par
\begin{equation}\label{indsh30d}
\theta_i\colon M\otimes_A \wedge^iE \rightarrow M\otimes_A \wedge^{i+1}E
\end{equation}
le morphisme composé de $\bIndMod(A)$
\begin{equation}\label{indsh30e}
\xymatrix{{M\otimes_A \wedge^iE}\ar[rr]^-(0.5){\theta\otimes_A\id_{\wedge^iE}}&&
{M\otimes_AE\otimes_A\wedge^i E}\ar[rr]^-(0.5){\id_M\otimes_Aw^i}&&{M\otimes_A\wedge^{i+1}E}}
\end{equation}
où $w^i\colon E\otimes_A\wedge^i E\rightarrow \wedge^{i+1}E$ est le produit extérieur \eqref{notconv9}.
On a $\theta_{i+1}\circ \theta_i=0$. On appelle complexe de {\em Dolbeault} de $(M,\theta)$
et l'on note $\mK^\bullet(M,\theta)$ le complexe de cochaînes de $\bIndMod(A)$
\begin{equation}\label{indsh30f}
M\stackrel{\theta}{\longrightarrow}M\otimes_AE\stackrel{\theta_1}{\longrightarrow} M\otimes_A\wedge^2E \dots,
\end{equation}
où $M$ est placé en degré $0$ et les différentielles sont de degré $1$.

\addtocounter{subsubsection}{3}
\addtocounter{equation}{1}

\subsubsection{}\label{indsh30g}
Soient $(M,\theta),(M',\theta')$ deux ind-$A$-modules de Higgs à coefficients dans $E$. 
On appelle champ de Higgs {\em total} sur $M\otimes_AM'$  
le morphisme  $A$-linéaire 
\begin{equation}\label{indsh30h}
\theta_\tot\colon M\otimes_AM'\rightarrow M\otimes_AM'\otimes_AE
\end{equation} 
défini par 
\begin{equation}\label{indsh30i}
\theta_\tot=\theta\otimes \id_{M'}+\id_{M}\otimes \theta'.
\end{equation}
On dit que $(M\otimes_AM',\theta_\tot)$ est le {\em produit tensoriel} de $(M,\theta)$ et $(M',\theta')$.

\addtocounter{subsubsection}{2}
\addtocounter{equation}{1}

\subsubsection{}\label{indsh30j}
Pour tout ind-$A$-module de Higgs $M$ à coefficients dans $E$, $\kappa_A(\theta)$ s'identifie à un $A$-champ de Higgs 
$\kappa_A(M)\rightarrow \kappa_A(M)\otimes_AE$ \eqref{indsh15g}. Notant $\mK^\bullet(M,\theta)$ le complexe de Dolbeault de $(M,\theta)$, 
$\kappa_A(\mK^\bullet(M,\theta))$ s'identifie au complexe de Dolbeault de $(\kappa_A(M),\kappa_A(\theta))$.

\subsection{}\label{indsh60}
Soient $(X,A)$ un topos annelé, 
\begin{equation}\label{indsh60a}
0\rightarrow L\rightarrow E\rightarrow \uE\rightarrow 0
\end{equation}
une suite exacte de $A$-modules localement libres de type fini. 
On munit l'algèbre extérieure $\wedge^\bullet E$ de $E$ de la filtration de Koszul $\rW^\bullet\wedge^\bullet E$ associée à cette suite (cf. \ref{MH90}). 

Soit $(M,\theta)$ un ind-$A$-module de Higgs à coefficients dans $E$.  
On désigna par $\utheta\colon M\rightarrow M\otimes_A\uE$ le $A$-champ de Higgs induit par $\theta$,
et par $\mK^\bullet$ (resp. $\umK^\bullet$) le complexe de Dolbeault de $(M,\theta)$ (resp. $(M,\utheta)$) \eqref{indsh30c}. 
Pour tous entiers $i,j\geq 0$, on pose 
\begin{equation}\label{indsh60c}
\rF^i\mK^j=M\otimes_A(\rW^i\wedge^jE).
\end{equation}
On définit ainsi une filtration décroissante exhaustive du $\wedge^\bullet E$-module gradué $\mK^\bullet$ 
par des sous-$\wedge^\bullet E$-modules gradués $(\rF^i\mK^\bullet)_{i\geq 0}$, stables par la différentielle de $\mK^\bullet$,
dite {\em filtration de Koszul de $\mK^\bullet$ associée à la suite exacte \eqref{indsh60a}};
en particulier, c'est une filtration du complexe $\mK^\bullet$ par des sous-complexes.

Compte tenu de \eqref{MH90c} et \eqref{MH90f}, pour tout entier $i\geq 0$, on a un isomorphisme canonique de complexes
\begin{equation}\label{indsh60e}
\rW^i\mK^\bullet/\rW^{i+1}\mK^\bullet\stackrel{\sim}{\rightarrow}\wedge^i L\otimes_A \umK^{\bullet}[-i],
\end{equation}
où les différentielles de $\umK^{\bullet}[-i]$ sont celles de $\umK^\bullet$ multipliées par $(-1)^i$.
On a donc une suite exacte canonique de ind-$A$-modules
\begin{equation}\label{indsh60f}
0\rightarrow \wedge^{i+1} L\otimes_A \umK^{\bullet}[-i-1]\rightarrow \rW^i\mK^\bullet/\rW^{i+2}\mK^\bullet
\rightarrow \wedge^i L\otimes_A \umK^{\bullet}[-i]\rightarrow 0.
\end{equation}
Cette dernière induit un morphisme de la catégorie dérivée $\bD^+(\bIndMod(A))$, 
\begin{equation}\label{indsh60g}
\partial^i\colon \wedge^i L\otimes_A \umK^{\bullet}\rightarrow \wedge^{i+1}L\otimes_A \umK^{\bullet},
\end{equation}
que l'on appellera {\em bord associé à la filtration de Koszul de $\mK^\bullet$}.

\begin{prop}\label{indsh62}
Soient $(X,A)$ un topos annelé,  
\begin{equation}\label{indsh62a}
0\rightarrow L\rightarrow E\rightarrow \uE\rightarrow 0
\end{equation}
une suite exacte de $A$-modules localement libres de type fini, 
$(M,\theta)$ un ind-$A$-module de Higgs à coefficients dans $E$, $(N,\kappa)$ un ind-$A$-module de Higgs à coefficients dans $L$. 
On désigne par $\utheta$ le $A$-champ de Higgs sur $M$ à coefficients dans $\uE$ induit par $\theta$,
par $\kappa'$ le $A$-champ de Higgs sur $N$ à coefficients dans $E$ induit par $\kappa$, 
par $\theta'=\theta\otimes\id+\id\otimes\kappa'$ le $A$-champ de Higgs total sur $M\otimes_A N$ à coefficients dans $E$,
et par $\utheta'$ le $A$-champ de Higgs sur $M\otimes_A N$ à coefficients dans $\uE$ induit par $\theta'$, qui n'est autre que $\utheta\otimes \id$. 
On note $\mK^\bullet$ (resp. $\umK^\bullet$, resp. $\mK'^\bullet$, resp. $\umK'^\bullet$) 
le complexe de Dolbeault de $(M,\theta)$ (resp. $(M,\utheta)$, resp. $(M\otimes_A N,\theta')$, resp. $(M\otimes_A N,\utheta')$) \eqref{indsh30c}, et 
\begin{eqnarray}
\partial\colon \umK^\bullet&\rightarrow& L\otimes_A\umK^\bullet,\label{indsh62b}\\
\partial'\colon \umK'^\bullet&\rightarrow& L\otimes_A\umK'^\bullet,\label{indsh62c}
\end{eqnarray}
les bords dans $\bD^+(\bIndMod(A))$ associés aux filtrations de Koszul de $\mK^\bullet$ et $\mK'^\bullet$ relatives à la suite exacte \eqref{indsh62a}; cf. \eqref{indsh60g}. 
Identifiant $\umK'^\bullet$ avec $\umK^\bullet\otimes_A N$, on a alors 
\begin{equation}\label{indsh62d}
\partial'=\partial\otimes \id-\id \otimes \kappa.
\end{equation}
\end{prop} 

Il suffit de calquer la preuve de \ref{MH100}.

\begin{defi}\label{indsh28}
Soient $(X,A)$ un topos annelé, $B$ une $A$-algèbre, $\lambda\in A(X)$,
$\Omega^1_{B/A}$ le $B$-module des differentielles de Kähler de $B/A$, $d\colon B\rightarrow \Omega^1_{B/A}$ la $A$-dérivation universelle \eqref{MH9},
$M$ un ind-$B$-module,
\begin{equation}\label{indsh28a}
\nabla\colon M\rightarrow \Omega^1_{B/A}\otimes_BM
\end{equation}
un morphisme de $\bIndMod(A)$ \eqref{indsh15f}.
On dit que $\nabla$ est une {\em $\lambda$-connexion sur $M$  
relativement à l'extension $B/A$} (et que $(M,\nabla)$ est un {\em ind-$B$-module à $\lambda$-connexion}) 
si le morphisme composé de $\bIndMod(A)$ 
\begin{equation}\label{indsh28b}
\xymatrix{{B\otimes_AM}\ar[r]&M\ar[r]^-(0.5)\nabla&{\Omega^1_{B/A}\otimes_BM}},
\end{equation}
où la première flèche est le morphisme canonique \eqref{indsh15i}, est la somme des deux morphismes composés de $\bIndMod(A)$
\begin{eqnarray}
\xymatrix{
{B\otimes_AM}\ar[rr]^-(0.5){\id_B\otimes_A\nabla}&&{B\otimes_A\Omega^1_{B/A}\otimes_BM}\ar[r]&{\Omega^1_{B/A}\otimes_BM}},\label{indsh28c}\\
\xymatrix{
{B\otimes_AM}\ar[rr]^-(0.5){\lambda d\otimes_A\id_M}&&{\Omega^1_{B/A}\otimes_AM}\ar[r]&{\Omega^1_{B/A}\otimes_BM}},\label{indsh28d}
\end{eqnarray}
où les secondes flèches sont les morphismes canoniques \eqref{indsh27a}.
\end{defi}

Soient $(M,\nabla)$, $(M',\nabla')$ deux ind-$B$-modules à $\lambda$-connexion. 
Un morphisme de $(M,\nabla)$ dans $(M',\nabla')$
est la donnée d'un morphisme $u\colon M\rightarrow M'$ de $\bIndMod(B)$ tel que  $(\id \otimes_B u)\circ \nabla=\nabla'\circ u$.

\begin{lem}\label{indsh29}
Soient $(X,A)$ un topos annelé, $B$ une $A$-algèbre, $\lambda\in A(X)$, $J$ une petite catégorie filtrante, 
$\alpha\colon J\rightarrow \bMod(B)$ un foncteur, $M=\indcolim \alpha$ sa limite inductive dans $\bIndMod(B)$, 
$\nabla$ une $\lambda$-connexion sur $M$ relativement à l'extension $B/A$. 
Notons $\iota\colon \bMod(B)\rightarrow \bMod(A)$
le foncteur d'oubli, et pour tout $j\in \ob(J)$, posons $M_j=\alpha(j)$. 
Alors, il existe une petite catégorie filtrante $K$, deux foncteurs cofinaux $\beta,\gamma \colon K\rightrightarrows J$
et deux morphismes de foncteurs $\sigma\colon \beta\rightarrow \gamma$ et 
$V\colon \iota\circ \alpha \circ \beta \rightarrow  \iota\circ (\Omega^1_{B/A}\otimes_B \alpha\circ \gamma)$, 
de sorte que pour tout $k\in \ob(K)$, on a un morphisme $B$-linéaire, fonctoriel en $k$, 
\begin{equation}\label{indsh29a}
\alpha(\sigma_k)\colon M_{\beta(k)}\rightarrow M_{\gamma(k)}
\end{equation}
et un morphisme $A$-linéaire, fonctoriel en $k$,
\begin{equation}\label{indsh29b}
V_k\colon M_{\beta(k)}\rightarrow \Omega^1_{B/A}\otimes_BM_{\gamma(k)},
\end{equation}
tels que les deux propriétés suivantes soient remplies:
\begin{itemize}
\item[{\rm (i)}] pour tout $k\in \ob(K)$ et pour toutes sections locales $b$ de $B$ et $m$ de $M_{\beta(k)}$, on a 
\begin{equation}\label{indsh29c}
V_k(bm)=\lambda d(b)\otimes_B \alpha(\sigma_k)(m)+bV_k(m).
\end{equation} 
\item[{\rm (ii)}] la limite inductive de $\alpha(\sigma)$ est l'identité de $M$ et la limite inductive de $V$ est $\nabla$. 
\end{itemize}
\end{lem}

Notons $\cC$ la catégorie filtrante $\bMod(A)_{/\Omega^1_{B/A}\otimes_BM}$ \eqref{indsh3d} et considérons les foncteurs 
\begin{eqnarray}
\psi\colon J\rightarrow \cC,\ \ \ j\mapsto (\Omega^1_{B/A}\otimes_B M_j\rightarrow \Omega^1_{B/A}\otimes_B M),\label{indsh29d}\\
\varphi\colon J\rightarrow \cC,\ \ \ i\mapsto (M_i\rightarrow M\stackrel{\nabla}{\rightarrow}\Omega^1_{B/A}\otimes_B M),\label{indsh29e}
\end{eqnarray}
induits par la relation $M=\indcolim \alpha$. Pour tout $A$-module $N$, on a 
\begin{equation}\label{indsh29f}
\Hom_{\bIndMod(A)}(N,\Omega^1_{B/A}\otimes_B M)=\underset{\underset{j\in J}{\longrightarrow}}{\lim}\ \Hom_{\bMod(A)}(N,\Omega^1_{B/A}\otimes_BM_j).
\end{equation}
Comme $J$ est filtrante, on en déduit que le foncteur $\psi$ est cofinal d'après (\cite{sga4} I 8.1.3(b)). 

Considérons les foncteurs 
\begin{eqnarray}
\psi'\colon J\rightarrow J\times \cC,\ \ \ j\mapsto (j,\psi(j)),\label{indsh29g}\\
\varphi'\colon J\rightarrow J\times \cC,\ \ \ i\mapsto (i,\varphi(i)).\label{indsh29h} 
\end{eqnarray}
Il résulte de ce qui précède et de (\cite{sga4} I 8.1.3(b)) que $\psi'$ est cofinal. 
On désigne par $K'$ la catégorie des triplets $(i,j,u)$ où $i,j\in \ob(J)$ et $u\in \Hom_{J\times \cC}(\varphi'(i),\psi'(j))$, avec les morphismes évidents 
\eqref{indsh37}. Plus explicitement, les objets de $K'$ sont les quadruplets 
$(i,j,s,v)$ où $i,j\in \ob(J)$ , $s\colon i\rightarrow j$ est un morphisme de $J$ et 
\begin{equation}\label{indsh29i}
v\colon M_i\rightarrow \Omega^1_{B/A}\otimes_AM_j
\end{equation}
est un morphisme $A$-linéaire qui s'insère dans le diagramme commutatif
\begin{equation}\label{indsh29j}
\xymatrix{
{M_i}\ar[r]^-(0.5)v\ar[d]&{\Omega^1_{B/A}\otimes_AM_j}\ar[d]\\
M\ar[r]^-(0.5)\nabla&{\Omega^1_{B/A}\otimes_AM}}
\end{equation}
où les flèches verticales sont induites par la relation $M=\indcolim \alpha$. 
D'après (\cite{ks2} 3.4.5), la catégorie $K'$ est filtrante et la première et la deuxième projections canoniques $\beta',\gamma'\colon K'\rightrightarrows J$, 
respectivement, sont cofinales. De plus, $K'$ est clairement petite. 

Soit $K$ la sous-catégorie pleine de $K'$ formée des quadruplets $(i,j,s,v)$ tels que 
pour toutes sections locales $b$ de $B$ et $m$ de $M_{i}$, on ait 
\begin{equation}\label{indsh29k}
v(bm)=\lambda d(b)\otimes_B \alpha(s)(m)+bv(m).
\end{equation} 
Il résulte de \eqref{indsh29f} et (\cite{sga4} I 8.1.3(c)) que le foncteur d'injection $K\rightarrow K'$ est cofinal et que $K$ est filtrante. 
On note $\beta,\gamma \colon K\rightrightarrows J$ les foncteurs induits par $\beta'$ et $\gamma'$. 
Pour tout $k\in \ob(K)$, écrivons 
\begin{equation}\label{indsh29l}
k=(\beta(k),\gamma(k),\sigma_k\colon \beta(k)\rightarrow \gamma(k), V_k \colon M_{\beta(k)}\rightarrow \Omega^1_{B/A}\otimes_BM_{\gamma(k)}). 
\end{equation}
On a donc deux morphismes de foncteurs $\sigma\colon \beta\rightarrow \gamma$ et 
$V\colon \iota\circ \alpha \circ \beta \rightarrow  \iota\circ (\Omega^1_{B/A}\otimes_B \alpha\circ \gamma)$. 
Il résulte aussitôt de \eqref{indsh29j} que la limite inductive de $V$ est le morphisme $\nabla$. De même, pour tout $k\in \ob(K)$, on a 
un diagramme commutatif 
\begin{equation}\label{indsh29m}
\xymatrix{
M_{\beta(k)}\ar[r]\ar[rd]&{M_{\gamma(k)}}\ar[d]\\
&M}
\end{equation}
On en déduit que la limite inductive de $\alpha(\sigma)$ est l'identité de $M$, d'où la proposition.

\begin{prop}\label{indsh31}
Soient $(X,A)$ un topos annelé, $B$ une $A$-algèbre, $\lambda\in A(X)$, $(M,\nabla)$ un ind-$B$-module
à $\lambda$-connexion relativement à l'extension $B/A$. 
Alors, $\nabla$ se prolonge uniquement en un morphisme gradué de degré $1$ de $\bIndMod(A)$ que l'on note aussi
\begin{equation}\label{indsh31a}
\nabla\colon \Omega_{B/A}\otimes_BM\rightarrow \Omega_{B/A}\otimes_BM,
\end{equation}
tel que pour tous entiers $p,q\geq 0$, le morphisme composé de $\bIndMod(A)$ 
\begin{equation}\label{indsh31b}
\xymatrix{{\Omega^p_{B/A}\otimes_A\Omega^q_{B/A}\otimes_BM}\ar[r]&{\Omega^{p+q}_{B/A}\otimes_BM}\ar[r]^-(0.5)\nabla&{\Omega^{p+q+1}_{B/A}\otimes_BM}}
\end{equation}
soit la somme des deux morphismes composés de $\bIndMod(A)$
\begin{equation}\label{indsh31c}
\xymatrix{
{\Omega^p_{B/A}\otimes_A\Omega^q_{B/A}\otimes_BM}\ar[rr]^-(0.5){(-1)^p\id\otimes_A\nabla}&&
{\Omega^p_{B/A}\otimes_A\Omega^{q+1}_{B/A}\otimes_BM}\ar[r]&{\Omega^{p+q+1}_{B/A}\otimes_BM}},
\end{equation}
\begin{equation}\label{indsh31d}
\xymatrix{
{\Omega^p_{B/A}\otimes_A\Omega^q_{B/A}\otimes_BM}\ar[rr]^-(0.5){\lambda d\otimes_A\id}&&
{\Omega^{p+1}_{B/A}\otimes_A\Omega^q_{B/A}\otimes_BM}\ar[r]&{\Omega^{p+q+1}_{B/A}\otimes_BM}},
\end{equation}
où les flèches non labellisées sont induites par la multiplication de $\Omega_{B/A}$.
De plus, pour tous entiers $p,q\geq 0$, le diagramme 
\begin{equation}\label{indsh31e}
\xymatrix{
{\Omega^p_{B/A}\otimes_A\Omega^q_{B/A}\otimes_BM}\ar[r]\ar[d]_{\id\otimes_A\nabla\circ \nabla}&{\Omega^{p+q}_{B/A}\otimes_BM}\ar[d]^{\nabla\circ \nabla}\\
{\Omega^p_{B/A}\otimes_A\Omega^{q+2}_{B/A}\otimes_BM}\ar[r]&{\Omega^{p+q+2}_{B/A}\otimes_BM}}
\end{equation}
où les flèches horizontales sont induites par la multiplication de $\Omega_{B/A}$, est commutatif. 
\end{prop}

En effet, l'unicité de $\nabla$ \eqref{indsh31a} 
résulte de la surjectivité des morphismes de multiplication $\Omega^p_{B/A}\otimes_A\Omega^q_{B/A}\rightarrow \Omega^{p+q}_{B/A}$ ($p,q\geq 0$).
Pour établir son existence, considérons une petite catégorie filtrante  $J$ et un foncteur  $\alpha\colon J\rightarrow \bMod(B)$
tels que $M=\indcolim \alpha$ dans $\bIndMod(B)$. 
Utilisons alors \ref{indsh29} et reprenons les mêmes notations. Pour tout $k\in \ob(K)$, on note aussi 
\begin{equation}\label{indsh32a}
\alpha(\sigma_k)\colon \Omega_{B/A}\otimes_BM_{\beta(k)}\rightarrow \Omega_{B/A}\otimes_BM_{\gamma(k)}
\end{equation}
le morphisme $B$-linéaire induit par $\alpha(\sigma_k)$ \eqref{indsh29a}. 
Le morphisme $V_k$ se prolonge en un unique morphisme $A$-linéaire gradué de degré $1$ que l'on note aussi 
\begin{equation}\label{indsh32b} 
V_k\colon \Omega_{B/A}\otimes_BM_{\beta(k)}\rightarrow   \Omega_{B/A}\otimes_BM_{\gamma(k)},
\end{equation} 
tel que pour toutes sections locales $\omega$ de  $\Omega^p_{B/A}$ et $m$ de $\Omega^q_{B/A}\otimes_BM_{\beta(k)}$  
($p,q\in \mN$), on ait 
\begin{equation}\label{indsh32c} 
V_k(\omega\wedge m)=\lambda  d(\omega)\wedge \alpha(\sigma_k)(m)+(-1)^p\omega \wedge V_k(m).
\end{equation} 
On obtient ainsi un morphisme de foncteurs que l'on note aussi
\begin{equation}\label{indsh32d} 
V\colon \iota\circ (\Omega_{B/A}\otimes_B \alpha \circ \beta) \rightarrow  \iota\circ (\Omega_{B/A}\otimes_B \alpha\circ \gamma).
\end{equation} 
On prend pour $\nabla$ \eqref{indsh31a} la limite inductive de $V$ dans $\bIndMod(A)$, qui vérifie la propriété requise puisque la limite inductive 
du morphisme de foncteurs $\alpha\circ \sigma$ est l'identité de $M$.  

D'après \ref{indsh37}, il existe une petite catégorie filtrante $L$, deux foncteurs cofinaux
$\tts,\ttt\colon L\rightrightarrows K$ et un morphisme de foncteurs $\tau\colon \gamma\circ \tts \rightarrow \beta\circ \ttt$. 
Pour tout $\ell\in \ob(L)$, on a donc les morphismes de $J$ 
\begin{equation}\label{indsh32h} 
\xymatrix{
{\beta\circ \tts(\ell)}\ar[r]^-(0.5){\sigma_{\tts(\ell)}}&{\gamma\circ \tts(\ell)}\ar[r]^-(0.5){\tau_\ell}&
{\beta \circ \tts (\ell)}\ar[r]^-(0.5){\sigma_{\ttt(\ell)}}&{\gamma\circ \ttt(\ell)}}.
\end{equation}
On désigne par $W_\ell$ le morphisme $A$-linéaire de degré $2$, défini par le composé
\begin{equation}\label{indsh32e} 
\xymatrix{
{\Omega_{B/A}\otimes_BM_{\beta(\tts(\ell))}}\ar[d]_-(0.5){V_{\tts(\ell)}}\ar[rr]^{W_\ell}&&{\Omega_{B/A}\otimes_BM_{\gamma(\ttt(\ell))}}\\
{\Omega_{B/A}\otimes_BM_{\gamma(\tts(\ell))}}\ar[rr]^-(0.5){\id\otimes\alpha(\tau_\ell)}&&{\Omega_{B/A}\otimes_BM_{\beta(\ttt(\ell))}}\ar[u]_-(0.5){V_{\ttt(\ell)}}}
\end{equation}
On obtient ainsi un morphisme de foncteurs que l'on note 
\begin{equation}\label{indsh32g} 
W\colon \iota\circ (\Omega_{B/A}\otimes_B \alpha \circ \beta \circ \tts) \rightarrow  \iota\circ (\Omega_{B/A}\otimes_B \alpha\circ \gamma\circ \ttt).
\end{equation} 
Comme la limite inductive du morphisme de foncteurs $\alpha\circ \tau$ est l'identité de $M$, 
la limite inductive de $W$  est le foncteur composé $\nabla\circ \nabla$. 

Pour tout $\ell\in \ob(L)$ et toutes sections locales $\omega$ de  $\Omega^p_{B/A}$ et $m$ de $\Omega^q_{B/A}\otimes_BM_{\beta(\tts(\ell))}$, on a 
\[ 
W_\ell(\omega\wedge m)=\omega \wedge W_\ell(m)+(-1)^p\lambda  d(\omega)\wedge [(\alpha(\sigma_{\ttt(\ell)}\circ \tau_\ell)\otimes \id)\circ V_{\tts(\ell)}(m)-V_{\ttt(\ell)}\circ 
\alpha(\tau_\ell\circ\sigma_{\tts(\ell)})(m)].
\]
La commutativité du diagramme \eqref{indsh31e} s'ensuit.

\begin{cor}\label{indsh34}
Sous les hypothèses de \ref{indsh31}, le morphisme $\nabla\circ \nabla$ \eqref{indsh31a} est l'image d'un morphisme gradué de degré $2$ de $\bIndMod(B)$ 
\begin{equation}\label{indsh31f}
\Omega_{B/A}\otimes_BM\rightarrow \Omega_{B/A}\otimes_BM.
\end{equation}
\end{cor}

Cela résulte de \ref{indsh33} et \eqref{indsh31e}.

\begin{defi}\label{indsh35}
Soient $(X,A)$ un topos annelé, $B$ une $A$-algèbre, $\lambda\in A(X)$, $(M,\nabla)$ un ind-$B$-module
à $\lambda$-connexion relativement à l'extension $B/A$. 
On dit que $\nabla$ est {\em intégrable} si les conditions équivalentes suivantes sont satisfaites:
\begin{itemize}
\item[(i)] le morphisme composé $\nabla\circ \nabla\colon \Omega_{B/A}\otimes_BM\rightarrow \Omega_{B/A}\otimes_BM$ est nul (cf. \ref{indsh31});
\item[(ii)] le morphisme composé
\begin{equation}
M\stackrel{\nabla}{\longrightarrow} \Omega^1_{B/A}\otimes_BM\stackrel{\nabla}{\longrightarrow} \Omega^2_{B/A}\otimes_BM
\end{equation}
est nul. 
\end{itemize}
\end{defi}

\begin{lem}\label{indsh43} 
Soient $(X,A)$ un topos annelé, $B$ une $A$-algèbre, $\lambda\in A(X)$, $(M,\nabla)$ un ind-$B$-module à 
$\lambda$-connexion relativement à l'extension $B/A$, 
$B'$ une $B$-algèbre,  $d'\colon B'\rightarrow \Omega^1_{B'/A}$ sa $A$-dérivation universelle. On désigne par 
\begin{equation}\label{indsh43a}
u\colon B'\otimes_B\Omega_{B/A}\rightarrow \Omega_{B'/A}
\end{equation}
le morphisme canonique de $B'$-algèbres graduées \eqref{MH9a}. Alors, il existe un unique morphisme de $\bIndMod(A)$
\begin{equation}\label{indsh43b}
\nabla'\colon B'\otimes_BM\rightarrow \Omega^1_{B'/A}\otimes_BM,
\end{equation}
dont le composé avec le morphisme canonique $B'\otimes_AM\rightarrow B'\otimes_BM$ est la somme des morphismes composés $v_1$ et $v_2$ suivants:
\begin{equation}\label{indsh43c}
\xymatrix{
{M\otimes_AB'}\ar[rr]^-(0.5){\nabla\otimes_A\id_{B'}}\ar[rrrd]_-(0.5){v_1}&&{(\Omega^1_{B/A}\otimes_BM)\otimes_AB'}\ar[r]&{\Omega^1_{B/A}\otimes_BM\otimes_BB'}
\ar[d]^-(0.5){u^1\otimes_B\id_M}\\
&&&{\Omega^1_{B'/A}\otimes_BM}\\
{B'\otimes_AM}\ar[rr]_-(0.5){\lambda d'\otimes_A\id_M}\ar[rrru]^-(0.5){v_2}&&{(\Omega^1_{B/A}\otimes_BB')\otimes_AM}\ar[r]&
{\Omega^1_{B/A}\otimes_BB'\otimes_BM}\ar[u]_-(0.5){u^1\otimes_B\id_M}}
\end{equation}
où les flèches non labellisées sont les morphismes canoniques. C'est une $\lambda$-connexion sur $B'\otimes_BM$ relativement à l'extension $B'/A$. 
Si, de plus, la $\lambda$-connexion $\nabla$ est intégrable, il en est de même de~$\nabla'$.
\end{lem}

En effet, l'unicité de $\nabla'$ \eqref{indsh43a} 
résulte de la surjectivité du morphisme canonique $B'\otimes_AM\rightarrow B'\otimes_BM$.
Pour établir son existence, considérons une petite catégorie filtrante  $J$ et un foncteur  $\alpha\colon J\rightarrow \bMod(B)$
tels que $M=\indcolim \alpha$ dans $\bIndMod(B)$. 
Utilisons alors \ref{indsh29} et reprenons les mêmes notations.
Pour tout $k\in \ob(K)$, il existe un unique morphisme $A$-linéaire 
\begin{equation}\label{indsh43d}
V'_k\colon B'\otimes_{B}M_{\beta(k)}\rightarrow  \Omega^1_{B'/A}\otimes_{B}M_{\gamma(k)}
\end{equation}
tel que pour toutes sections locales $b'$ de $B'$ et $m$ de $M_{\beta(k)}$, on ait 
\begin{equation}\label{indsh43e}
V'_k(b'\otimes m)=\lambda d'(b')\otimes \alpha(\sigma_k)(m)+(u\otimes\id_{M_{\gamma(k)}})(b'\otimes V_k(m)).
\end{equation}
On obtient un morphisme de foncteurs 
\begin{equation}\label{indsh43f}
V'\colon \iota\circ (B'\otimes_B \alpha \circ \beta) \rightarrow  \iota\circ (\Omega^1_{B'/A}\otimes_B \alpha\circ \gamma).
\end{equation} 
Comme la limite inductive du morphisme de foncteurs $\alpha\circ \sigma$ est l'identité de $M$, 
en passant à limite inductive, on obtient un morphisme de $\bIndMod(A)$
\begin{equation}\label{indsh43g}
\nabla'\colon B'\otimes_BM\rightarrow \Omega^1_{B'/A}\otimes_BM.
\end{equation}
Il résulte aussitôt de \eqref{indsh43e} que $\nabla'$ est une $\lambda$-connexion sur $B'\otimes_BM$ relativement à l'extension $B'/A$. 

Pour tout $k \in \ob(K)$, $V_k$ se prolonge en un unique morphisme $A$-linéaire gradué de degré $1$ que l'on note aussi 
\begin{equation}\label{indsh43h} 
V_k\colon \Omega_{B/A}\otimes_BM_{\beta(k)}\rightarrow   \Omega_{B/A}\otimes_BM_{\gamma(k)},
\end{equation} 
tel que pour toutes sections locales $\omega$ de  $\Omega^p_{B/A}$ et $m$ de $\Omega^q_{B/A}\otimes_BM_{\beta(k)}$  
($p,q\in \mN$), on ait 
\begin{equation}\label{indsh43i} 
V_k(\omega\wedge m)=\lambda  d(\omega)\wedge \alpha(\sigma_k)(m)+(-1)^p\omega \wedge V_k(m).
\end{equation} 
De même, pour tout $k \in \ob(K)$, $V'_k$ se prolonge en un unique morphisme $A$-linéaire gradué de degré $1$ que l'on note aussi 
\begin{equation}\label{indsh43j} 
V'_k\colon \Omega_{B'/A}\otimes_BM_{\beta(k)}\rightarrow   \Omega_{B'/A}\otimes_BM_{\gamma(k)},
\end{equation} 
tel que pour toutes sections locales $\omega'$ de  $\Omega^p_{B'/A}$ et $m'$ de $\Omega^q_{B'/A}\otimes_BM_{\beta(k)}$  
($p,q\in \mN$), on ait 
\begin{equation}\label{indsh43k} 
V'_k(\omega'\wedge m')=\lambda  d'(\omega')\wedge \alpha(\sigma_k)(m')+(-1)^p\omega' \wedge V'_k(m').
\end{equation}

D'après \ref{indsh37}, il existe une petite catégorie filtrante $L$, deux foncteurs cofinaux
$\tts,\ttt\colon L\rightrightarrows K$ et un morphisme de foncteurs $\tau\colon \gamma\circ \tts \rightarrow \beta\circ \ttt$. 
Pour tout $\ell\in \ob(L)$, on a donc les morphismes de $J$ 
\begin{equation}\label{indsh43l} 
\xymatrix{
{\beta\circ \tts(\ell)}\ar[r]^-(0.5){\sigma_{\tts(\ell)}}&{\gamma\circ \tts(\ell)}\ar[r]^-(0.5){\tau_\ell}&
{\beta \circ \ttt (\ell)}\ar[r]^-(0.5){\sigma_{\ttt(\ell)}}&{\gamma\circ \ttt(\ell)}}.
\end{equation}
On désigne par $W_\ell$ et $W'_\ell$ les morphismes $A$-linéaires de degré $2$, définis par les composés
\begin{equation}\label{indsh43m} 
\xymatrix{
{\Omega_{B/A}\otimes_BM_{\beta(\tts(\ell))}}\ar[d]_-(0.5){V_{\tts(\ell)}}\ar[rr]^{W_\ell}&&{\Omega_{B/A}\otimes_BM_{\gamma(\ttt(\ell))}}\\
{\Omega_{B/A}\otimes_BM_{\gamma(\tts(\ell))}}\ar[rr]^-(0.5){\id\otimes\alpha(\tau_\ell)}&&{\Omega_{B/A}\otimes_BM_{\beta(\ttt(\ell))}}\ar[u]_-(0.5){V_{\ttt(\ell)}}}
\end{equation}
\begin{equation}\label{indsh43n} 
\xymatrix{
{\Omega_{B'/A}\otimes_BM_{\beta(\tts(\ell))}}\ar[d]_-(0.5){V'_{\tts(\ell)}}\ar[rr]^{W'_\ell}&&{\Omega_{B'/A}\otimes_BM_{\gamma(\ttt(\ell))}}\\
{\Omega_{B'/A}\otimes_BM_{\gamma(\tts(\ell))}}\ar[rr]^-(0.5){\id\otimes\alpha(\tau_\ell)}&&{\Omega_{B'/A}\otimes_BM_{\beta(\ttt(\ell))}}\ar[u]_-(0.5){V'_{\ttt(\ell)}}}
\end{equation}
Pour toute section locale $m$ de $M_{\beta(\tts(\ell))}$, on a
\begin{eqnarray}\label{indsh43p}  
W'_\ell(m)&=&V'_{\ttt(\ell)}\circ (\id\otimes\alpha(\tau_\ell))\circ V'_{\tts(\ell)}(m)\\
&=&V'_{\ttt(\ell)}\circ (u\otimes\alpha(\tau_\ell))\circ V_{\tts(\ell)}(m)\nonumber\\
&=&(u\otimes\id_{M_{\gamma(\ttt(\ell))}})\circ V_{\ttt(\ell)}\circ(\id\otimes\alpha(\tau_\ell))\circ V_{\tts(\ell)}(m)\nonumber\\
&=&(u\otimes\id_{M_{\gamma(\ttt(\ell))}})\circ W_\ell(m). \nonumber
\end{eqnarray} 
On en déduit que si la $\lambda$-connexion $\nabla$ est intégrable, il en est de même de $\nabla'$ (cf. \ref{indsh34}).

\subsection{}\label{indsh44}  
Soient $f\colon (X',A')\rightarrow (X,A)$ un morphisme de topos annelés, $B$ une $A$-algèbre, 
$B'$ une $A'$-algèbre, $\alpha\colon f^*(B)\rightarrow B'$ un homomorphisme de $A'$-algèbres, $\lambda\in \Gamma(X,A)$,
$M$ un ind-$B$-module, $\nabla$ une $\lambda$-connexion sur $M$ relativement à l'extension $B/A$ \eqref{indsh28}. 
Notons  $\lambda'$ l'image canonique de $\lambda$ dans $\Gamma(X',A')$. Le foncteur $f^*$ induit un foncteur \eqref{indsh21a}
\begin{equation}\label{indsh44a}  
\rI f^* \colon \bIndMod(B)\rightarrow \bIndMod(f^*(B)),
\end{equation}
compatible avec le foncteur $\rI f^* \colon \bIndMod(A)\rightarrow \bIndMod(A')$ via les foncteurs d'oubli \eqref{indsh15f}. 
On vérifie aussitôt \eqref{indsh21h} que le morphisme de $\bIndMod(A')$
\begin{equation}\label{indsh44b}  
\rI f^*(\nabla)\colon \rI f^*(M)\rightarrow \rI f^*(M)\otimes_{f^*(B)} \Omega^1_{f^*(B)/A'}
\end{equation}
est une $\lambda'$-connexion sur le ind-$f^*(B)$-module $\rI f^*(M)$ relativement à l'extension $f^*(B)/A'$. 
Il résulte de \ref{indsh31} que $(\rI f^*(M), \rI f^*(\nabla))$ est intégrable si $(M,\nabla)$ l'est.

D'après \ref{indsh43}, $(\rI f^*(M), \rI f^*(\nabla))$ définit grâce à $\alpha$ 
une $\lambda'$-connexion sur le ind-$B'$-module $\rI f^*(M)\otimes_{f^*(B)} B'$ relativement à l'extension $B'/A'$,
\begin{equation}\label{indsh44c}  
\nabla'\colon \rI f^*(M)\otimes_{f^*(B)} B'\rightarrow \rI f^*(M)\otimes_{f^*(B)} \Omega^1_{B'/A'},
\end{equation}
qui est intégrable si $(M,\nabla)$ l'est.

\subsection{}\label{indsh38}
Soient $(X,A)$ un topos annelé, $B$ une $A$-algèbre, $\lambda\in \Gamma(X,A)$, $(M,\nabla)$ un
ind-$B$-module à $\lambda$-connexion relativement à l'extension $B/A$ \eqref{indsh28}. 
Supposons qu'il existe un $A$-module $E$ et un $B$-isomorphisme $\tau\colon E\otimes_AB\stackrel{\sim}{\rightarrow}\Omega^1_{B/A}$ 
tels que pour toute section locale $\omega$ de $E$, on ait $d(\tau(\omega\otimes 1))=0$.
On note $\vartheta\colon M\rightarrow E\otimes_AM$  le morphisme induit par $\nabla$ et $\tau$. 
\begin{itemize}
\item[(i)] Pour que la $\lambda$-connexion $\nabla$ soit intégrable \eqref{indsh35}, il faut et il suffit que 
$\vartheta$ soit un $A$-champ de Higgs sur $M$ à coefficients dans $E$ \eqref{indsh30}.  En effet, le diagramme 
\begin{equation}
\xymatrix{
{E\otimes_AM}\ar[r]^-(0.5){-\id\otimes\vartheta}\ar[d]&{E\otimes_AE\otimes_AM}\ar[r]&{(\wedge^2E)\otimes_AM}\ar[d]\\
{\Omega^1_{B/A}\otimes_BM}\ar[rr]^-(0.5)\nabla&&{\Omega^2_{B/A}\otimes_BM}}
\end{equation}
où les flèches verticales sont les isomorphismes induits par $\tau$, est clairement commutatif. 

\item[(ii)] Soit $(N,\theta)$ un ind-$A$-module de Higgs à coefficients dans $E$. Le morphisme de $\bIndMod(A)$
\begin{equation}\label{indsh38a}
\nabla'\colon M\otimes_AN\rightarrow  \Omega^1_{B/A}\otimes_BM\otimes_AN
\end{equation}
défini par
\begin{equation}\label{indsh38b}
\nabla'=\nabla\otimes_A \id_N+ (\tau\otimes_B\id_{M\otimes_AN})(\id_M\otimes_A \theta),
\end{equation}
est une $\lambda$-connexion sur $M\otimes_AN$ relativement à $B/A$.
Si $\nabla$ est intégrable, il en est de même de $\nabla'$. 
\end{itemize}

\begin{lem}\label{indsh380}
Soient $(X,A)$ un topos annelé, $B$ une $A$-algèbre, $\lambda\in \Gamma(X,A)$, $(M,\nabla)$ un
ind-$B$-module à $\lambda$-connexion relativement à l'extension $B/A$, 
$B'$ une $B$-algèbre, $\delta'\colon B'\rightarrow \Omega^1_{B/A}\otimes_BB'$ une $A$-dérivation 
telle que pour toute section locale $b$ de $B$, on ait $\delta'(b)=\lambda d(b)\otimes 1$. Alors, 
\begin{itemize}
\item[{\rm (i)}] Il existe un unique morphisme de $\bIndMod(A)$
\begin{equation}\label{indsh380a}
\nabla'\colon B'\otimes_BM\rightarrow  \Omega^1_{B/A}\otimes_BB'\otimes_BM,
\end{equation}
tel que le composé avec le morphisme canonique $B'\otimes_AM\rightarrow B'\otimes_BM$ soit la somme des morphismes
\begin{equation}\label{indsh380b}
\xymatrix{
{M\otimes_AB'}\ar[rr]^-(0.5){\nabla\otimes_A\id_{B'}}&&{(\Omega^1_{B/A}\otimes_BM)\otimes_AB'}\ar[r]&{\Omega^1_{B/A}\otimes_BM\otimes_BB'}},
\end{equation}
\begin{equation}\label{indsh380c}
\xymatrix{
{B'\otimes_AM}\ar[rr]^-(0.5){\delta'\otimes_A\id_M}&&{(\Omega^1_{B/A}\otimes_BB')\otimes_AM}\ar[r]&{\Omega^1_{B/A}\otimes_BB'\otimes_BM}},
\end{equation}
où les secondes flèches sont les morphismes canoniques.
\item[{\rm (ii)}] Supposons, de plus, qu'il existe un $A$-module $E$ et un $B$-isomorphisme $\tau\colon E\otimes_AB\stackrel{\sim}{\rightarrow}\Omega^1_{B/A}$ 
tels que pour toute section locale $\omega$ de $E$, on ait $d(\tau(\omega\otimes 1))=0$.
On désigne par $\vartheta'\colon B'\otimes_BM\rightarrow  E\otimes_AB'\otimes_BM$ (resp. $\theta'\colon B'\rightarrow  E\otimes_AB'$)
le morphisme induit par $\nabla'$ (resp. $\delta'$) et $\tau$.  
Si $\nabla$ est integrable et si $\theta'$ est un $A$-champ de Higgs à coefficients dans $E$, alors $\vartheta'$ est un $A$-champ de Higgs à coefficients dans $E$. 
\end{itemize}
\end{lem}

(i) En effet, l'unicité de $\nabla'$ \eqref{indsh380a} 
résulte de la surjectivité du morphisme canonique $M\otimes_AB'\rightarrow M\otimes_BB'$.
Pour établir son existence, considérons une petite catégorie filtrante  $J$ et un foncteur  $\alpha\colon J\rightarrow \bMod(B)$
tels que $M=\indcolim \alpha$ dans $\bIndMod(B)$. 
Utilisons alors \ref{indsh29} et reprenons les mêmes notations.
Pour tout $k\in \ob(K)$, il existe un unique morphisme $A$-linéaire
\begin{equation}\label{indsh380d}
V'_k\colon B'\otimes_BM_{\beta(k)}\rightarrow \Omega^1_{B/A}\otimes_BB'\otimes_BM_{\gamma(k)},
\end{equation}
tel que pour toutes sections locales $b'$ de $B'$ et $m$ de $M_{\beta(k)}$, on ait
\begin{equation}\label{indsh380e}
V'_k(b'\otimes_Bm)=\delta'(b')\otimes_B \alpha(\sigma_k)(m)+b'\otimes_BV_k(m).
\end{equation}
On définit ainsi un morphisme de foncteurs 
\begin{equation}\label{indsh380f}
V'\colon \iota\circ (B'\otimes_B\alpha \circ \beta)\rightarrow  \iota\circ (\Omega^1_{B/A}\otimes_BB'\otimes_B \alpha\circ \gamma).
\end{equation}
On prend pour $\nabla'$ \eqref{indsh380a} la limite inductive de $V'$, qui vérifie clairement la propriété requise. 

(ii) Conservons les notations précédentes. Pour tout $k \in \ob(K)$, on désigne par 
\begin{eqnarray}
\vartheta_k\colon M_{\beta(k)}&\rightarrow& M_{\gamma(k)} \otimes_A E,\\
\vartheta'_k\colon M_{\beta(k)}\otimes_BB'&\rightarrow& M_{\gamma(k)} \otimes_BB'\otimes_A E,
\end{eqnarray} 
les morphismes induits par $V_k$, $V'_k$ et $\tau$. 
Le morphisme $\vartheta_k$ induit un unique morphisme $A$-linéaire que l'on note aussi 
\begin{equation}
\vartheta_k\colon M_{\beta(k)}\otimes_A E\rightarrow M_{\gamma(k)} \otimes_A \wedge^2 E,
\end{equation} 
tel que pour toutes sections locales 
$m$ de $M_{\beta(k)}\otimes_A \wedge^p E$ et $\omega$ de $\wedge^qE$, on ait $\vartheta_k(m\otimes \omega)=\vartheta_k(m)\wedge \omega$. 
De même, $\theta'$ et $\vartheta'_k$ induisent canoniquement des morphismes $A$-linéaires que l'on note aussi 
\begin{eqnarray}
\theta'\colon B'\otimes_A E&\rightarrow& B' \otimes_A \wedge^2 E,\\
\vartheta'_k\colon M_{\beta(k)}\otimes_BB'\otimes_A E&\rightarrow& M_{\gamma(k)} \otimes_BB'\otimes_A \wedge^2 E.
\end{eqnarray} 

D'après \ref{indsh37}, il existe une petite catégorie filtrante $L$, deux foncteurs cofinaux
$\tts,\ttt\colon L\rightrightarrows K$ et un morphisme de foncteurs $\tau\colon \gamma\circ \tts \rightarrow \beta\circ \ttt$. 
Pour tout $\ell\in \ob(L)$, on a donc les morphismes de $J$ 
\begin{equation}
\xymatrix{
{\beta\circ \tts(\ell)}\ar[r]^-(0.5){\sigma_{\tts(\ell)}}&{\gamma\circ \tts(\ell)}\ar[r]^-(0.5){\tau_\ell}&
{\beta \circ \ttt (\ell)}\ar[r]^-(0.5){\sigma_{\ttt(\ell)}}&{\gamma\circ \ttt(\ell)}}.
\end{equation}
On désigne par $W_\ell$ et $W'_\ell$ les morphismes $A$-linéaires définis par les composés
\begin{equation}\label{indsh380g}
\xymatrix{
{M_{\beta(\tts(\ell))}}\ar[d]_-(0.5){\vartheta_{\tts(\ell)}}\ar[rr]^{W_\ell}&&{M_{\gamma(\ttt(\ell))}\otimes_A\wedge^2E}\\
{M_{\gamma(\tts(\ell))}\otimes_AE}\ar[rr]^-(0.5){\alpha(\tau_\ell)\otimes\id}&&{M_{\beta(\ttt(\ell))}\otimes_AE}\ar[u]_-(0.5){\vartheta_{\ttt(\ell)}}}
\end{equation}
\begin{equation}\label{indsh380h}
\xymatrix{
{M_{\beta(\tts(\ell))}\otimes_BB'}\ar[d]_-(0.5){\vartheta'_{\tts(\ell)}}\ar[rr]^{W'_\ell}&&{M_{\gamma(\ttt(\ell))}\otimes_BB'\otimes_A\wedge^2E}\\
{M_{\gamma(\tts(\ell))}\otimes_BB'\otimes_A E}\ar[rr]^-(0.5){\alpha(\tau_\ell)\otimes\id}&&{M_{\beta(\ttt(\ell))}\otimes_BB'\otimes_AE}\ar[u]_-(0.5){\vartheta'_{\ttt(\ell)}}}
\end{equation}
Ce sont des morphismes de foncteurs. 
Comme la limite inductive du morphisme de foncteurs $\alpha\circ \tau$ est l'identité de $M$, on voit que la limite inductive de $\ell\mapsto W_\ell$ est nulle.

Par ailleurs, pour toutes sections locales $m$ de $M_{\beta(\tts(\ell))}$ et $b'$ de $B'$, comme $\theta'^2 (b')=0$, on a 
\[
W'_{\ell}(m\otimes_Bb')=W_\ell(m)\otimes b'-(\alpha(\sigma_{\ttt(\ell)}\circ \tau_\ell)\otimes\id)\circ \vartheta_{\tts(\ell)}(m)\wedge \theta'(b')+
\vartheta_{\ttt(\ell)}\circ \alpha(\tau_\ell\circ \sigma_{\tts(\ell)})(m)\wedge \theta'(b').
\]
Comme la limite inductive du morphisme de foncteurs $\ell\mapsto \alpha(\tau_\ell\circ \sigma_{\tts(\ell)})$ 
(resp. $\ell\mapsto \alpha(\sigma_{\ttt(\ell)}\circ \tau_\ell)$) est l'identité de $M$, on en déduit que la limite inductive du morphisme de foncteurs 
$\ell\mapsto W'_\ell$ est nulle, ce qui démontre que $\vartheta'$ est un $A$-champ de Higgs \eqref{indsh380h}.

\section{Modules à isogénie près}

\begin{defi}\label{caip1}
Soit $\cC$ une catégorie additive. 
\begin{itemize}
\item[(i)] Un morphisme $u\colon M\rightarrow N$ de $\cC$ est appelé {\em isogénie} s'il existe
un entier $n\not=0$ et un morphisme $v\colon N\rightarrow M$ de $\cC$ tels que $v\circ u=n\cdot \id_M$ 
et $u\circ v= n\cdot \id_N$. 
\item[(ii)] Un objet $M$ de $\cC$ est dit {\em d'exposant fini} s'il existe un entier $n\not=0$ tel que $n\cdot \id_M=0$. 
\end{itemize}
\end{defi}
On peut compléter la terminologie et faire les remarques suivantes~:

\subsubsection{}\label{caip1a}
\addtocounter{equation}{1}
La famille des isogénies de $\cC$ permet un calcul de fractions bilatéral (\cite{illusie1} I 1.4.2).
On appelle {\em catégorie des objets de $\cC$ à isogénie près}, et l'on note $\cC_\mQ$,  
la catégorie localisée de $\cC$ par rapport aux isogénies. On désigne par
\begin{equation}\label{caip1aa}
Q\colon \cC\rightarrow \cC_\mQ, \ \ \ X\mapsto X_\mQ,
\end{equation}
le foncteur de localisation. On vérifie aisément que pour tous $X,Y\in \ob(\cC)$, on a  
\begin{equation}\label{caip1ab}
\Hom_{\cC_\mQ}(X_\mQ,Y_\mQ)=\Hom_{\cC}(X,Y)\otimes_{\mZ}\mQ.
\end{equation}
En particulier, la catégorie $\cC_\mQ$ est additive et le foncteur de localisation est additif. 
Pour qu'un objet $X$ de $\cC$ soit d'exposant fini, il faut et il suffit que $X_\mQ$ soit nul. 
Pour qu'un morphisme $f$ de $\cC$ soit une isogénie, il faut et il suffit que $f_\mQ$ soit un isomorphisme de $\cC_\mQ$. 

\addtocounter{subsubsection}{2}

\subsubsection{}\label{caip1b}
\addtocounter{equation}{1}
Si $\cC$ est une catégorie abélienne, la catégorie $\cC_\mQ$ est abélienne et le foncteur de localisation 
$Q\colon \cC\rightarrow \cC_\mQ$ est exact. En fait, $\cC_\mQ$ s'identifie canoniquement
à la catégorie quotient de $\cC$ par la sous-catégorie épaisse $\cE$ des objets d'exposant fini.
En effet, notons $\cC/\cE$ la catégorie quotient de $\cC$ par $\cE$ et $T\colon \cC\rightarrow \cC/\cE$ 
le foncteur canonique  (\cite{gabriel} III §1). Pour tout $X\in \ob(\cE)$, on a $Q(X)=0$. 
Par suite, il existe un et un unique foncteur $Q'\colon \cC/\cE\rightarrow \cC_\mQ$ tel que $Q=Q'\circ T$. 
Par ailleurs, pour tout $X\in \ob(\cC)$ et tout entier $n\not=0$, $T(n\cdot \id_X)$ est un isomorphisme. 
Il existe donc un et un unique foncteur $T'\colon \cC_\mQ\rightarrow \cC/\cE$ tel que $T=T'\circ Q$. 
On voit aussitôt que $T'$ et $Q'$ sont des équivalences de catégories quasi-inverses l'une de l'autre. 
Le foncteur $Q$, qui s'identifie au foncteur $T$, est donc exact (\cite{gabriel} III prop.~1).

\subsubsection{}\label{caip1c}
\addtocounter{equation}{1}
Tout foncteur additif (resp. exact)  entre catégories additives (resp. abéliennes) $\cC\rightarrow \cC'$ 
s'étend de manière unique en un foncteur additif (resp. exact) $\cC_\mQ\rightarrow \cC'_\mQ$, compatible aux foncteurs 
de localisation.

\subsection{}\label{indsh11}
Soient $\cC$ une $\mU$-catégorie abélienne \eqref{notconv3}, $I$ le système multiplicatif bilatéral des isogénies de $\cC$, 
$\cC_\mQ$ la catégorie des objets de $\cC$ à isogénie près, qui est une $\mU$-catégorie abélienne, 
\begin{equation}\label{indsh11a}
Q\colon \cC\rightarrow \cC_\mQ, \ \ \ X\mapsto X_\mQ,
\end{equation}
le foncteur canonique. 
Pour tout objet $X$ de $\cC$, on désigne par $I^X$ la catégorie des isogénies de source $X$ et par 
\begin{equation}\label{indsh11b}
\alpha^X\colon I^X\rightarrow \cC, \ \ \ (X\rightarrow X')\mapsto X',
\end{equation}
le foncteur but (cf. \ref{indsh5}). La catégorie $I^X$ est filtrante (\cite{ks2} 7.1.10) et essentiellement petite  \eqref{indsh2}. 
En effet, les endomorphismes de $X$ définis par la multiplication par un entier non-nul 
forment une petite partie cofinale de $\ob(I^X)$. On peut donc considérer le foncteur pleinement fidèle \eqref{indsh5e}
\begin{equation}\label{indsh11c}
\upalpha\colon \cC_\mQ\rightarrow \Ind(\cC), \ \ \ X\mapsto \indcolim \alpha^X. 
\end{equation}
On rappelle que le morphisme canonique \eqref{indsh5g}
\begin{equation}\label{indsh11d}
\iota_\cC\rightarrow \upalpha\circ Q
\end{equation}
n'est pas un isomorphisme en général. 

\begin{lem}\label{indsh16}
Conservons les hypothèses et notations de \ref{indsh11}. Notons, de plus, $M$ le monoïde multiplicatif $\mZ-\{0\}$ et $\uM$ 
la catégorie filtrante dont les objets sont les 
éléments de $M$ et les morphismes sont déterminés par la relation de divisibilité dans $M$. Alors,
\begin{itemize}
\item[{\rm (i)}] Pour tout objet $X$ de $\cC$, le foncteur $\upmu^X\colon \uM\rightarrow I^X$ qui envoie un entier non nul $n$ 
sur l'isogénie de multiplication par $n$ dans $X$, est cofinal. 
\item[{\rm (ii)}] Le foncteur $\upalpha$ \eqref{indsh11c} est exact. 
\item[{\rm (iii)}] Pour tout objet injectif $X$ de $\cC$, $Q(X)$ est injectif et $\upalpha(Q(X))$ est quasi-injectif. 
\end{itemize}
\end{lem}

(i) En effet, d'après (\cite{sga4} I 8.1.3(b)), 
comme la catégorie $\uM$ est filtrante, il suffit de montrer que le foncteur $\upmu^X$ satisfait les propriétés F1) et F2) de (\cite{sga4} I 8.1.3). 
Soit $u\colon X\rightarrow X'$ une isogénie de $\cC$. Il existe un morphisme $v\colon X'\rightarrow X$ de $\cC$ 
et un élément $n$ de $M$ tels que $v\circ u=n \cdot \id_X$ 
et $u\circ v=n\cdot \id_{X'}$. Par suite, $v$ définit un morphisme de $(u\colon X\rightarrow X')$ vers $\upmu^X([n])$ dans $I^{X}$, d'où la  propriété F1). 

Soient $m\in M$ et $v_1,v_2 \colon X'\rightrightarrows X$ deux morphismes de $\cC$ tels que  $m \cdot \id_X=v_1\circ u=v_2\circ u$.
On a $m \cdot v =n \cdot v_1=n \cdot v_2$ dans $\cC$. Par suite, on a $\upmu^X([m]\rightarrow [mn])\circ v_1=\upmu^X([m]\rightarrow [mn])\circ v_2$ dans $I^X$, 
d'où la  propriété F2). 

(ii) En effet, il résulte de (i) et \ref{indsh6e} que le foncteur $\upalpha\circ Q$ est exact. 
La proposition s'ensuit compte tenu de (\cite{gabriel} III cor.~1 à prop.~1). 

(iii) La première assertion résulte de (\cite{gabriel} III cor.~1 à prop.~1) et la seconde assertion est une conséquence de (i) et \ref{indsh8b}.

\subsection{}\label{indsh12}
Soit $\phi\colon \cC\rightarrow \cC'$ un foncteur additif entre $\mU$-catégories abéliennes. 
On reprend les notations de \ref{indsh11} pour $\cC$ et on considère les notations analogues pour $\cC'$ que l'on équipe d'un exposant $^\prime$. 
Pour tout objet $X$ de $\cC$, $\phi$ induit un foncteur 
\begin{equation}\label{indsh12a}
\phi^X\colon I^X\rightarrow I'^{\phi(X)}
\end{equation}
qui s'insère dans un diagramme strictement commutatif
\begin{equation}\label{indsh12b}
\xymatrix{
I^X\ar[r]^-(0.5){\phi^X}\ar[d]_-(0.5){\alpha^X}&{I'^{\phi(X)}}\ar[d]^-(0.5){\alpha'^{\phi(X)}}\\
{\cC}\ar[r]^-(0.5)\phi&{\cC'}}
\end{equation}

\begin{lem}\label{indsh13}
Conservons les hypothèses de \ref{indsh12}. Alors,
\begin{itemize}
\item[{\rm (i)}] Pour tout objet $X$ de $\cC$, le foncteur $\phi^X$ est cofinal. 
\item[{\rm (ii)}] Le diagramme 
\begin{equation}\label{indsh13a}
\xymatrix{
{\cC_\mQ}\ar[r]^{\phi_\mQ}\ar[d]_{\upalpha}&{\cC'_\mQ}\ar[d]^{\upalpha'}\\
{\Ind(\cC)}\ar[r]^{\rI\phi}&{\Ind(\cC')}}
\end{equation}
où $\phi_\mQ$ et $\rI\phi$ sont les foncteurs additifs induits par $\phi$, et 
$\upalpha$ et $\upalpha'$ sont les foncteurs canoniques \eqref{indsh11b}, est commutatif à
isomorphisme canonique près. 
\end{itemize}
\end{lem}

(i) En effet, d'après (\cite{sga4} I 8.1.3(b)), comme la catégorie $I^X$ est filtrante, il suffit de montrer que le foncteur $\phi^X$
satisfait les propriétés F1) et F2) de (\cite{sga4} I 8.1.3). Soit $u\colon \phi(X)\rightarrow Y$ une isogénie de $\cC'$. 
Il existe un morphisme $v\colon Y\rightarrow \phi(X)$ et un entier $n$ non nul tels que $v\circ u=n \cdot \id_{\phi(X)}$ et $u\circ v=n\cdot \id_Y$. 
Par suite, $v$ définit un morphisme de $u$ vers $\phi^X(n\cdot \id_X)$ dans $I'^{\phi(X)}$, d'où la  propriété F1). 

Soient $w\colon X\rightarrow X'$ une isogénie de $\cC$ et $v_1,v_2\colon Y\rightrightarrows \phi(X')$ deux morphismes de $\cC'$ 
tels que $\phi(w)=v_1\circ u=v_2\circ u$. On a $\phi(w)\circ v=n \cdot v_1=n \cdot v_2$. Par suite, $\phi^X(n\cdot \id_{X'})\circ v_1=\phi^X(n\cdot \id_{X'})\circ v_2$, 
d'où la  propriété F2). 

(ii) En effet, on a des isomorphismes canoniques
\begin{equation}\label{indsh13b}
\rI\phi(\indcolim \alpha^X)\stackrel{\sim}{\rightarrow} \indcolim (\phi\circ \alpha^X) \stackrel{\sim}{\rightarrow} \indcolim (\alpha'^{\phi(X)}\circ \phi^X)
\stackrel{\sim}{\rightarrow} \indcolim \alpha'^{\phi(X)},
\end{equation}
le premier est \eqref{indsh4c}, le deuxième est sous-jacent au diagramme commutatif \eqref{indsh12b} et le troisième est le morphisme
canonique, qui est un isomorphisme par (i).

\subsection{}\label{indsh14}
Soit $\phi\colon \cC\rightarrow \cC'$ un foncteur exact à gauche entre $\mU$-catégories abéliennes tel que $\cC$ ait assez d'injectifs. 
D'après (\cite{gabriel} III cor.~1 à prop.~1), le foncteur additif $\phi_\mQ\colon \cC_\mQ\rightarrow \cC'_\mQ$, induit par $\phi$  \eqref{indsh11}, 
est exact à gauche. D'après \ref{indsh16}(iii), le foncteur canonique $Q\colon \cC\rightarrow \cC_\mQ$ \eqref{indsh11a} transforme 
les objets injectifs en des objets injectifs. 
On désigne par $\cI$ la sous-catégorie pleine de $\cC$ formée des objets injectifs et par $\cI_\mQ$ la sous-catégorie pleine de $\cC_\mQ$
formée des images par $Q$ des objets injectifs de $\cC$. 
Alors, $\cI$ (resp. $\cI_\mQ$) est cogénératrice de $\cC$ (resp. $\cC_\mQ$), en particulier $\cC_\mQ$ a assez d'injectifs. 
Par ailleurs,  $\cI$ est $\phi$-injective d'après (\cite{ks2} 13.3.6(iii)), et  $\cI_\mQ$ est $\phi_\mQ$-injective
en vertu de (\cite{ks2} 13.3.7) et (\cite{gabriel} III cor.~1 à prop.~1). 
D'après (\cite{ks2} 13.3.5), les foncteurs $\phi$ et $\phi_\mQ$ admettent donc des foncteurs dérivés à droite
\begin{eqnarray}
\rR\phi\colon \bD^+(\cC)\rightarrow \bD^+(\cC'),\label{indsh14a}\\
\rR(\phi_\mQ)\colon \bD^+(\cC_\mQ)\rightarrow\bD^+(\cC'_\mQ).\label{indsh14b}
\end{eqnarray}
Pour justifier l'existence de $\rR(\phi_\mQ)$, on peut plus simplement évoquer le fait que $\cC_\mQ$ a assez injectifs. 
Mais l'introduction de $\cI_\mQ$ est nécessaire pour \eqref{indsh14g} ci-dessous. 

Comme les foncteurs $Q$ et $Q'$ sont exacts et qu'ils transforment les objets injectifs en des objets injectifs,  on a des isomorphismes canoniques 
\begin{equation}\label{indsh14c}
\rR(\phi_\mQ) \circ Q' \stackrel{\sim}{\rightarrow}\rR(\phi_\mQ\circ Q') \stackrel{\sim}{\rightarrow} 
\rR(Q\circ \phi) \stackrel{\sim}{\rightarrow}  Q\circ\rR(\phi),
\end{equation}
où le premier et le dernier isomorphismes résultent de (\cite{ks2} 13.3.13) et l'isomorphisme central est induit par la définition de $\phi_\mQ$. 
Le diagramme 
\begin{equation}\label{indsh14d}
\xymatrix{
{\bD^+(\cC)}\ar[r]^{\rR(\phi)}\ar[d]_{Q}&{\bD^+(\cC')}\ar[d]^{Q'}\\
{\bD^+(\cC_\mQ)}\ar[r]^{\rR(\phi_\mQ)}&{\bD^+(\cC'_\mQ)}}
\end{equation}
est donc commutatif à isomorphisme canonique près. 
En particulier, pour tout entier $i$, le foncteur $\rR^i(\phi_\mQ)$ est déduit par localisation du foncteur $\rR^i\phi$.  

La catégorie $\Ind(\cI)$ s'identifie à la sous-catégorie pleine de $\Ind(\cC)$ formée des objets quasi-injectifs \eqref{indsh8}. 
En vertu de (\cite{ks2} 15.3.2), $\Ind(\cI)$ est $\rI\phi$-injective et le foncteur $\rI\phi$ admet un foncteur dérivé à droite
\begin{equation}\label{indsh14e}
\rR(\rI\phi)\colon \bD^+(\Ind(\cC))\rightarrow \bD^+(\Ind(\cC')).
\end{equation}
Il résulte de \ref{indsh9} que pour tout entier $i$, on a un isomorphisme canonique $\rR^i(\rI\phi)\stackrel{\sim}{\rightarrow}\rI(\rR^i\phi)$.

D'après \ref{indsh16}, le foncteur $\upalpha\colon \cC_\mQ\rightarrow \Ind(\cC)$ \eqref{indsh11c} est exact et on a $\upalpha(\cJ_\mQ)\subset \Ind(\cI)$. 
Par suite, on a des isomorphismes canoniques 
\begin{equation}\label{indsh14f}
\rR(\rI\phi) \circ \upalpha \stackrel{\sim}{\rightarrow}\rR((\rI\phi) \circ \upalpha) \stackrel{\sim}{\rightarrow} 
\rR(\upalpha'\circ \phi_\mQ) \stackrel{\sim}{\rightarrow} \upalpha'\circ\rR(\phi_\mQ),
\end{equation}
où le premier et le dernier isomorphisme résultent de (\cite{ks2} 13.3.13) et l'isomorphisme central est induit par l'isomorphisme sous-jacent à \eqref{indsh13a}.  
Le diagramme
\begin{equation}\label{indsh14g}
\xymatrix{
{\bD^+(\cC_\mQ)}\ar[r]^{\rR(\phi_\mQ)}\ar[d]_{\upalpha}&{\bD^+(\cC'_\mQ)}\ar[d]^{\upalpha'}\\
{\bD^+(\Ind(\cC))}\ar[r]^{\rR(\rI\phi)}&{\bD^+(\Ind(\cC'))}}
\end{equation}
est donc commutatif à un isomorphisme canonique près. 
En particulier, pour tout entier $i$, le diagramme
\begin{equation}\label{indsh14h}
\xymatrix{
{\cC_\mQ}\ar[r]^{\rR^i(\phi_\mQ)}\ar[d]_{\upalpha}&{\cC'_\mQ}\ar[d]^{\upalpha'}\\
{\Ind(\cC)}\ar[r]^{\rR^i(\rI\phi)}&{\Ind(\cC')}}
\end{equation}
est commutatif à un isomorphisme canonique près.

\subsection{}\label{indsh20}
Soit $(X,A)$ un $\mU$-topos annelé.  
On désigne par $\bMod_{\mQ}(A)$ la catégorie des $A$-modules de $X$ à isogénie près \eqref{indsh11} et par  
\begin{equation}\label{indsh20a}
Q_A\colon \bMod(A)\rightarrow \bMod_\mQ(A),\ \ \ M\mapsto M_\mQ,
\end{equation}
le foncteur canonique \eqref{indsh11a}. Le produit tensoriel de $\bMod(A)$ induit un bifoncteur 
\begin{equation}\label{indsh20b}
\bMod_{\mQ}(A)\times \bMod_{\mQ}(A)\rightarrow \bMod_{\mQ}(A),\ \ \ (M,N)\mapsto M\otimes_{A_\mQ}N,
\end{equation} 
faisant de $\bMod_{\mQ}(A)$ une catégorie monoïdale symétrique, ayant $A_\mQ$ pour objet unité. 
Les objets de $\bMod_{\mQ}(A)$ seront aussi appelés des {\em $A_\mQ$-modules}. Cette terminologie
se justifie en considérant $A_\mQ$ comme un monoïde de $\bMod_{\mQ}(A)$. 
Le bifoncteur \eqref{indsh20b} est exact à droite. 

On a un foncteur canonique \eqref{indsh11c}
\begin{equation}\label{indsh20c}
\upalpha_{A}\colon \bMod_\mQ(A)\rightarrow \bIndMod(A), 
\end{equation}
qui est exacte d'après \ref{indsh16}(ii).
Compte tenu de \ref{indsh16}(i), pour tous $A$-modules $M$ et $N$, on a des isomorphismes canoniques
\begin{equation}\label{indsh20d}
\upalpha_{A}(M_\mQ)\otimes_A N\stackrel{\sim}{\rightarrow} \upalpha_{A}(M_\mQ)\otimes_A\upalpha_{A}(N_\mQ)
\stackrel{\sim}{\rightarrow} \upalpha_{A}(M_\mQ\otimes_{A_\mQ}N_\mQ).
\end{equation}

\subsection{}\label{indsh22}
Soit $f\colon (Y,B)\rightarrow (X,A)$ un morphisme de $\mU$-topos annelés. 
Les foncteurs adjoints $f^*$ et $f_*$ entre les catégories $\bMod(A)$ et $\bMod(B)$ induisent des foncteurs additifs adjoints  
\begin{eqnarray}
f^*_\mQ\colon \bMod_\mQ(A) \rightarrow \bMod_{\mQ}(B),\label{indsh22a}\\
f_{\mQ*}\colon \bMod_{\mQ}(B) \rightarrow \bMod_\mQ(A),\label{indsh22b}
\end{eqnarray}
dont le premier (resp. second) est exact à droite (resp. gauche). 
D'après \ref{indsh14}, le foncteur $f_{\mQ *}$ admet un foncteur dérivé à droite
\begin{equation}\label{indsh22c}
\rR f_{\mQ*}\colon \bD^+(\bMod_\mQ(B))\rightarrow \bD^+(\bMod_\mQ(A)).
\end{equation}
Le diagramme 
\begin{equation}\label{indsh22d}
\xymatrix{
{\bD^+(\bMod(A))}\ar[r]^{\rR f_*}\ar[d]_{Q_A}&{\bD^+(\bMod(B))}\ar[d]^{Q_B}\\
{\bD^+(\bMod_\mQ(A))}\ar[r]^{\rR f_{\mQ*}}&{\bD^+(\bMod_\mQ(B))}}
\end{equation}
où $Q_A$ et $Q_B$ sont les foncteurs canoniques \eqref{indsh20a}, est commutatif à isomorphisme canonique près. 

D'après \ref{indsh13}(ii), le diagramme 
\begin{equation}\label{indsh22f}
\xymatrix{
{\bMod_\mQ(A)}\ar[r]^{f^*_{\mQ}}\ar[d]_{\upalpha_{A}}&{\bMod_\mQ(B)}\ar[d]^{\upalpha_{B}}\\
{\bIndMod(A)}\ar[r]^{\rI f^*}&{\bIndMod(B)}}
\end{equation}
où $\upalpha_A$ et $\upalpha_B$ sont les foncteurs canoniques \eqref{indsh20c}, est commutatif à isomorphisme canonique près. 
En vertu de \eqref{indsh14g}, le diagramme
\begin{equation}\label{indsh22e}
\xymatrix{
{\bD^+(\bMod_\mQ(B))}\ar[r]^-(0.5){\rR(f_{\mQ*})}\ar[d]_{\upalpha_{B}}&{\bD^+(\bMod_\mQ(A))}\ar[d]^{\upalpha_{A}}\\
{\bD^+(\bIndMod(B))}\ar[r]^-(0.5){\rR(\rI f_*)}&{\bD^+(\bIndMod(A))}}
\end{equation}
est commutatif à un isomorphisme canonique près.

\subsection{}\label{indsh23}
Soient $(X,A)$ un topos annelé, $E$ un $A$-module. 
On appelle {\em $A$-isogénie de Higgs à coefficients dans $E$} la donnée d'un quadruplet
\begin{equation}\label{indsh23b}
(M,N,u\colon M\rightarrow N,\theta\colon M\rightarrow N\otimes_AE)
\end{equation}
formé de deux $A$-modules $M$ et $N$ et de deux morphismes $A$-linéaires $u$ et $\theta$  
vérifiant la propriété suivante~: il existe un entier $n\not=0$ et un morphisme $A$-linéaire $v\colon N\rightarrow M$ tels que 
$v\circ u=n\cdot \id_M$, $u\circ v=n\cdot \id_N$ et que $(M,(v\otimes \id_E)\circ \theta)$ et $(N,\theta\circ v)$ 
soient des $A$-modules de Higgs à coefficients dans $E$. 
On notera que $u$ induit une isogénie de modules de Higgs de $(M,(v\otimes \id_E)\circ \theta)$
dans $(N,\theta\circ v)$ \eqref{caip1}, d'où la terminologie.  
Soient $(M,N,u,\theta)$, $(M',N',u',\theta')$ deux $A$-isogénies de Higgs à coefficients dans $E$. 
Un morphisme de $(M,N,u,\theta)$ dans $(M',N',u',\theta')$ est la donnée de deux morphismes $A$-linéaires
$\alpha\colon M\rightarrow M'$ et $\beta\colon N\rightarrow N'$ tels que $\beta\circ u=u'\circ \alpha$
et $(\beta\otimes \id_E)\circ \theta=\theta'\circ \alpha$. 
On désigne par $\bIH(A,E)$ la catégorie des $A$-isogénies de Higgs à coefficients dans $E$. C'est une catégorie 
additive. On note $\bIH_\mQ(A,E)$ la catégorie des objets de $\bIH(A,E)$ à isogénie près \eqref{caip1}. 

Soit $(M,N,u,\theta)$ une $A$-isogénie de Higgs à coefficients dans $E$. Pour tout $i\geq 1$, on désigne par
\begin{equation}\label{indsh23c}
\theta_i\colon M\otimes_A \wedge^iE \rightarrow N\otimes_A \wedge^{i+1}E
\end{equation}
le morphisme $A$-linéaire défini pour toutes sections locales 
$m$ de $M$ et $\omega$ de $\wedge^iE$ par $\theta_i(m\otimes \omega)=\theta(m)\wedge \omega$.
On note 
\begin{equation}\label{indsh23d}
\otheta_i= (u^{-1}_\mQ\otimes \id_{ \wedge^{i+1}E})\circ \theta_i
\colon M_\mQ\otimes_{A_\mQ} (\wedge^iE)_\mQ \rightarrow M_\mQ\otimes_{A_\mQ} (\wedge^{i+1}E)_\mQ.
\end{equation}
On montre aussitôt que $\otheta_{i+1}\circ \otheta_i=0$ (\cite{agt} III.6.9). 
On appelle complexe de {\em Dolbeault} de $(M,N,u,\theta)$
et l'on note $\mK^\bullet(M,N,u,\theta)$ le complexe de cochaînes de $\bMod_\mQ(A)$ 
\begin{equation}\label{indsh23e}
M_\mQ\stackrel{\otheta_0}{\longrightarrow}M_\mQ\otimes_{A_\mQ}E_\mQ\stackrel{\otheta_1}{\longrightarrow} 
M_\mQ\otimes_{A_\mQ}(\wedge^2E)_\mQ\rightarrow \dots,
\end{equation}
où $M_\mQ$ est placé en degré $0$ et les différentielles sont de degré $1$. On obtient ainsi un foncteur
de la catégorie $\bIH(A,E)$ dans la catégorie des complexes de $\bMod_\mQ(A)$.
Toute isogénie de $\bIH(A,E)$ induit un isomorphisme des complexes de Dolbeault associés. Le foncteur ``complexe de Dolbeault''
induit donc un foncteur de $\bIH_\mQ(A,E)$ dans la catégorie des complexes de $\bMod_\mQ(A)$.

Soient $(M,N,u,\theta)$, $(M',N',u',\theta')$ deux $A$-isogénies de Higgs à coefficients dans $E$. Posant 
\begin{equation}\label{indsh23f}
\theta_\tot=\theta\otimes_Au'+u\otimes_A\theta'\colon M\otimes_AM'\rightarrow N\otimes_AN'\otimes_AE,
\end{equation}
on voit aussitôt que $(M\otimes_AM',N\otimes_AN',u\otimes_Au',\theta_\tot)$ est une $A$-isogénie de Higgs à coefficients dans $E$. 
On définit ainsi un foncteur bi-additif que l'on note 
\begin{equation}\label{indsh23g}
\bIH(A,E)\times \bIH(A,E)\rightarrow \bIH(A,E), \ \ \ (I,I')\mapsto I\otimes_AI'.
\end{equation}
Celui-ci s'étend naturellement en un foncteur bi-additif que l'on note encore
\begin{equation}\label{indsh23h}
\bIH_\mQ(A,E)\times \bIH_\mQ(A,E)\rightarrow \bIH_\mQ(A,E), \ \ \ (I,I')\mapsto I\otimes_{A_\mQ}I'.
\end{equation}

\subsection{}\label{indsh49}
Soient $(X,A)$ un topos annelé, $E$ un $A$-module, $(M,N,u,\theta)$, $(M',N',u',\theta')$ deux $A$-isogénies de Higgs à coefficients dans $E$. Comme 
$u_\mQ$ est un isomorphisme de $\bMod_\mQ(A)$ \eqref{indsh20}, on peut considérer le composé 
\begin{equation}\label{indsh49a}
\otheta=(u^{-1}_\mQ\otimes \id_{E})\circ \theta_\mQ\colon M_\mQ\rightarrow M_\mQ\otimes_{A_\mQ}E_\mQ.
\end{equation} 
On définit de même le morphisme $\otheta'\colon M'_\mQ\rightarrow M'_\mQ\otimes_{A_\mQ}E_\mQ$ de $\bMod_\mQ(A)$. 
Le morphisme 
\begin{equation}\label{indsh49b}
\begin{array}[t]{clcr}
\Hom_{\bIH(A,E)}((M,N,u,\theta),(M',N',u',\theta'))&\rightarrow &\Hom_A(M,M')\\
(\alpha,\beta)&\mapsto &\alpha
\end{array}
\end{equation}
induit alors un isomorphisme
\begin{eqnarray}\label{indsh49c}
\lefteqn{\Hom_{\bIH_\mQ(A,E)}((M,N,u,\theta),(M',N',u',\theta'))\stackrel{\sim}{\rightarrow}}\\
&& \{\alpha\in \Hom_A(M,M')\otimes_\mZ\mQ, \ (\alpha\otimes \id_E)\circ \otheta=\otheta'\circ \alpha\}.\nonumber
\end{eqnarray}
En effet, comme $u_\mQ$ et $u'_\mQ$ sont des isomorphismes de $\bMod_\mQ(A)$, le morphisme \eqref{indsh49b} induit un morphisme injectif
\begin{equation}\label{indsh49d}
\Hom_{\bIH_\mQ(A,E)}((M,N,u,\theta),(M',N',u',\theta'))\rightarrow \Hom_A(M,M')\otimes_\mZ\mQ.
\end{equation}
Soit $\alpha\in \Hom_A(M,M')\otimes_\mZ\mQ$ tel que $(\alpha\otimes \id_E)\circ \otheta=\otheta'\circ \alpha$. 
Montrons que $\alpha$ est dans l'image du morphisme \eqref{indsh49c}. Comme ce dernier est $\mQ$-linéaire,
on peut supposer que $\alpha \in \Hom_A(M,M')$. Ils existent deux entiers non-nuls $n$ et $n'$ et deux morphismes $A$-linéaires 
$v\colon N\rightarrow M$ et $v'\colon N'\rightarrow M'$ tels que $v\circ u=n\cdot \id_M$, $u\circ v=n\cdot \id_N$,  $v'\circ u'=n'\cdot \id_{M'}$, 
$u'\circ v'=n'\cdot \id_{N'}$ et que $(M,(v\otimes \id_E)\circ \theta)$ et $(M',(v'\otimes \id_E)\circ \theta')$ soient des $A$-modules de Higgs à coefficients dans $E$. 
Quitte à multiplier $n$ et $n'$, on peut supposer que $n=n'$ et que le diagramme 
\begin{equation}
\xymatrix{
M\ar[r]^-(0.5)\theta\ar[d]_\alpha&{N\otimes_AE}\ar[r]^-(0.5){v\otimes \id_E}&{M\otimes_AE}\ar[d]^{\alpha\otimes\id_E}\\
M'\ar[r]^-(0.5){\theta'}&{N'\otimes_AE}\ar[r]^-(0.5){v'\otimes \id_E}&{M'\otimes_AE}}
\end{equation} 
est commutatif. Par suite, $(\alpha,\alpha)$ définit un morphisme de $\bIH(A,E)$ de $(M,M,\id,(v\otimes\id_E)\circ \theta)$ dans 
$(M',M',\id,(v'\otimes\id_E)\circ \theta')$. Comme les morphismes  
\begin{eqnarray}
(\id_M,v)\colon (M,N,u,\theta)&\rightarrow& (M,M,\id,(v\otimes\id_E)\circ \theta),\\
(\id_{M'},v')\colon (M',N',u',\theta')&\rightarrow& (M',M',\id,(v'\otimes\id_E)\circ \theta'),
\end{eqnarray}
sont des isogénies de $\bIH(A,E)$, on en déduit que le morphisme \eqref{indsh49c} est surjectif et par suite bijectif.

\subsection{}\label{indsh24}
Soient $(X,A)$ un topos annelé, $B$ une $A$-algèbre, $d\colon B\rightarrow \Omega^1_{B/A}$ la $A$-dérivation universelle, $\lambda\in \Gamma(X,A)$. 
On appelle {\em $\lambda$-isoconnexion relativement à l'extension $B/A$}
(ou simplement {\em $\lambda$-isoconnexion} lorsqu'il n'y a aucun risque de confusion) la donnée d'un quadruplet   
\begin{equation}\label{indsh24a}
(M,N,u\colon M\rightarrow N,\nabla\colon M\rightarrow \Omega^1_{B/A}\otimes_BN)
\end{equation}
où $M$ et $N$ sont des $B$-modules, $u$ est une isogénie de $B$-modules \eqref{caip1}
et $\nabla$ est un morphisme $A$-linéaire tel que pour toutes sections locales 
$x$ de $B$ et $t$ de $M$, on ait 
\begin{equation}\label{indsh24b}
\nabla(xt)=\lambda d(x) \otimes u(t)+x\nabla(t).
\end{equation} 
Pour tout morphisme $B$-linéaire $v\colon N\rightarrow M$ pour lequel il existe un entier $n$ 
tel que  $u\circ v=n\cdot \id_N$ et $v\circ u=n\cdot \id_M$, 
les couples $(M,(\id\otimes v)\circ \nabla)$ et $(N,\nabla\circ v)$ sont des modules à $(n\lambda)$-connexions \eqref{MH9},
et $u$ est un morphisme de $(M,(\id\otimes v)\circ \nabla)$ dans $(N,\nabla\circ v)$. On dit que la $\lambda$-isoconnexion 
$(M,N,u,\nabla)$ est {\em intégrable} s'il existe un morphisme $B$-linéaire $v\colon N\rightarrow M$ et un entier $n\not= 0$ 
tels que  $u\circ v=n\cdot \id_N$, $v\circ u=n\cdot \id_M$ et que les $(n\lambda)$-connexions $(\id\otimes v)\circ \nabla$ sur $M$ et $\nabla\circ v$ 
sur $N$ soient intégrables. 

Soient $(M,N,u,\nabla)$, $(M',N',u',\nabla')$ deux $\lambda$-isoconnexions. 
Un morphisme de $(M,N,u,\nabla)$ dans $(M',N',u',\nabla')$ est la donnée de 
deux morphismes $B$-linéaires $\alpha\colon M\rightarrow M'$ et $\beta\colon N\rightarrow N'$
tels que  $\beta\circ u=u'\circ \alpha$ et $(\id \otimes \beta)\circ \nabla=\nabla'\circ \alpha$. 

On désigne par $\bIMC^\lambda(B/A)$ la catégorie des $\lambda$-isoconnexions intégrables relativement à l'extension $B/A$. C'est une catégorie 
additive. On note $\bIMC^\lambda_\mQ(B/A)$ la catégorie des objets de $\bIMC^\lambda(B/A)$ à isogénie près \eqref{caip1}.

\subsection{}\label{indsh50}
Soient $(X,A)$ un topos annelé, $B$ une $A$-algèbre, $\lambda\in \Gamma(X,A)$
$(M,N,u,\nabla)$ et $(M',N',u',\nabla')$ deux $\lambda$-isoconnexions relativement à l'extension $B/A$. 
Comme $u_\mQ$ est un isomorphisme de $\bMod_\mQ(B)$ \eqref{indsh20}, on peut considérer le morphisme composé de $\bMod_\mQ(A)$
\begin{equation}\label{indsh50a}
\onabla=(\id\otimes u^{-1}_\mQ)\circ \nabla\colon M_\mQ\rightarrow (\Omega^1_{B/A})_\mQ\otimes_{A_\mQ}M_\mQ.
\end{equation} 
On définit de même le morphisme $\onabla'\colon M'_\mQ\rightarrow (\Omega^1_{B/A})_\mQ\otimes_{A_\mQ}M'_\mQ$ de $\bMod_\mQ(A)$. 
Calquant la preuve de \ref{indsh49}, on montre que le morphisme 
\begin{equation}\label{indsh50b}
\begin{array}[t]{clcr}
\Hom_{\bIMC^\lambda(B/A)}((M,N,u,\nabla),(M',N',u',\nabla'))&\rightarrow &\Hom_B(M,M')\\
(\alpha,\beta)&\mapsto &\alpha
\end{array}
\end{equation}
induit un isomorphisme
\begin{eqnarray}\label{indsh50c}
\lefteqn{\Hom_{\bIMC^\lambda_\mQ(B/A)}((M,N,u,\nabla),(M',N',u',\nabla'))\stackrel{\sim}{\rightarrow}}\\
&& \{\alpha\in \Hom_B(M,M')\otimes_\mZ\mQ, \ (\id \otimes \alpha)\circ \onabla=\onabla'\circ \alpha \ {\rm dans}\ \bMod_\mQ(A)\}.\nonumber
\end{eqnarray}

\subsection{}\label{indsh25}
Soient $f\colon (X',A')\rightarrow (X,A)$ un morphisme de topos annelés, $B$ une $A$-algèbre, 
$B'$ une $A'$-algèbre, $\alpha\colon f^*(B)\rightarrow B'$ un homomorphisme de $A'$-algèbres, $\lambda\in \Gamma(X,A)$,
$(M,N,u,\nabla)$ une $\lambda$-isoconnexion relativement à l'extension $B/A$. 
Notons  $\lambda'$ l'image canonique de $\lambda$ dans $\Gamma(X',A')$, 
$d'\colon B'\rightarrow \Omega^1_{B'/A'}$ la $A'$-dérivation universelle de $B'$ et 
\begin{equation}\label{indsh25a}
\gamma\colon f^*(\Omega^1_{B/A}) \rightarrow \Omega^1_{B'/A'}
\end{equation}
le morphisme $\alpha$-linéaire canonique. On voit aussitôt que 
$(f^*(M),f^*(N),f^*(u),f^*(\nabla))$ est une $\lambda'$-isoconnexion relativement à l'extension $f^*(B)/A'$,
qui est intégrable si $(M,N,u,\nabla)$ l'est.

Il existe un unique morphisme $A'$-linéaire 
\begin{equation}\label{indsh25b}
\nabla'\colon B'\otimes_{f^*(B)}f^*(M)\rightarrow  \Omega^1_{B'/A'}\otimes_{f^*(B)}f^*(N)
\end{equation}
tel que pour toutes sections locales $x'$ de $B'$ et $t$ de $f^*(M)$, on ait 
\begin{equation}\label{indsh25c}
\nabla'(x'\otimes t)=\lambda'd'(x')\otimes f^*(u)(t)+x'(\gamma\otimes\id_{f^*(N)})(f^*(\nabla)(t)).
\end{equation}
Le quadruplet $(B'\otimes_{f^*(B)}f^*(M),B'\otimes_{f^*(B)}f^*(N),\id_{B'}\otimes_{f^*(B)} f^*(u),\nabla')$ 
est une $\lambda'$-isoconnexion relativement à l'extension $B'/A'$, qui est intégrable si $(M,N,u,\nabla)$ l'est.

\subsection{}\label{indsh26}
Soient $(X,A)$ un topos annelé, $B$ une $A$-algèbre, $\lambda\in \Gamma(X,A)$, $(M,N,u,\nabla)$ une $\lambda$-isoconnexion intégrable relativement à $B/A$. 
Supposons qu'il existe un $A$-module $E$ et un $B$-isomorphisme $\gamma\colon E\otimes_AB\stackrel{\sim}{\rightarrow}\Omega^1_{B/A}$ 
tels que pour toute section locale $\omega$ de $E$, on ait $d(\gamma(\omega\otimes 1))=0$.
Notons $\vartheta\colon M\rightarrow E\otimes_AN$ le morphisme induit par $\nabla$ et $\gamma$. Alors,
\begin{itemize}
\item[(i)] Le quadruplet $(M,N,u,\vartheta)$ est une $A$-isogénie de Higgs à coefficients dans $E$.  
\item[(ii)] Pour toute $A$-isogénie de Higgs $(M',N',u',\theta')$ à coefficients dans $E$, il existe un unique morphisme $A$-linéaire 
\begin{equation}\label{indsh26a}
\nabla'\colon M\otimes_AM'\rightarrow  \Omega^1_{B/A}\otimes_BN\otimes_AN'
\end{equation}
tel que pour toutes sections locales $t$ de $M$ et $t'$ de $M'$, on ait 
\begin{equation}\label{indsh26b}
\nabla'(t\otimes t')=\nabla(t)\otimes_Au'(t')+ (\gamma\otimes_B\id_{N\otimes_A N'})(u(t)\otimes_A \theta'(t')).
\end{equation}
Le quadruplet $(M\otimes_AM',N\otimes_AN',u\otimes u',\nabla')$ est une $\lambda$-isoconnexion intégrable. 
\end{itemize}

\section{Complément sur la fonctorialité des topos co-évanescents généralisés}\label{cftf}

\subsection{}\label{cftf1}
Commençons par rappeler et préciser la terminologie introduite dans (\cite{agt} § VI.5) sur les topos co-évanescents généralisés. 
On appelle {\em site fibré co-évanescent} 
la donnée d'un $\mU$-site $I$ et d'une catégorie fibrée, clivée et normalisée au-dessus de la catégorie sous-jacente à $I$ (\cite{sga1} VI 7.1)
\begin{equation}\label{cftf1a}
\pi\colon E\rightarrow I,
\end{equation}
vérifiant les conditions suivantes:
\begin{itemize}
\item[(i)] les produits fibrés sont représentables dans $I$; 
\item[(ii)]  pour tout $i\in \ob(I)$, la catégorie fibre $E_i$ de $E$ au-dessus de $i$ 
est munie d'une topologie faisant de celle-ci un $\mU$-site, 
et les limites projectives finies sont représentables dans $E_i$. 
On note $\tE_i$ le topos des faisceaux de $\mU$-ensembles sur $E_i$. 
\item[(iii)] pour tout morphisme $f\colon i\rightarrow j$ de $I$, le foncteur image inverse $f^+\colon E_j\rightarrow E_i$
est continu et exact à gauche. Il définit donc un morphisme de topos que l'on note aussi (abusivement) 
$f\colon \tE_i\rightarrow \tE_j$ (\cite{sga4} IV 4.9.2). 
\end{itemize}

Pour tout $i\in \ob(I)$, on note
\begin{equation}\label{cftf1b}
\alpha_i\colon E_i\rightarrow E
\end{equation}
le foncteur d'inclusion canonique. On prendra garde que ce foncteur était noté $\alpha_{i!}$ dans (\cite{agt} (VI.5.1.2)).

Le foncteur $\pi$ est en fait un $\mU$-site fibré  (\cite{sga4} VI 7.2.1 et 7.2.4). On désigne par 
\begin{equation}\label{cftf1c}
\cF\rightarrow I
\end{equation}
le $\mU$-topos fibré associé à $\pi$ (\cite{sga4} VI 7.2.6). La catégorie fibre de $\cF$ au-dessus de tout $i\in \ob(I)$
est canoniquement équivalente au topos $\tE_i$, et le foncteur image inverse par tout  
morphisme $f\colon i\rightarrow j$ de $I$ s'identifie au foncteur image inverse $f^*\colon \tE_j\rightarrow \tE_i$ 
du morphisme de topos $f\colon \tE_i\rightarrow \tE_j$. On désigne par
\begin{equation}\label{cftf1d}
\cF^\vee\rightarrow I^\circ
\end{equation}
la catégorie fibrée obtenue en associant à tout $i\in \ob(I)$ la catégorie $\tE_i$, et à tout morphisme 
$f\colon i\rightarrow j$ de $I$ le foncteur $f_*\colon \tE_i\rightarrow \tE_j$ image directe par le morphisme 
de topos $f\colon \tE_i\rightarrow \tE_j$. On désigne par
\begin{equation}\label{cftf1e}
\cP^\vee\rightarrow I^\circ
\end{equation}
la catégorie fibrée obtenue en associant à tout $i\in \ob(I)$ la catégorie $\hE_i$ des préfaisceaux de $\mU$-ensembles
sur $E_i$, et à tout morphisme $f\colon i\rightarrow j$ de $I$ le foncteur $\hf^*\colon \hE_i\rightarrow \hE_j$ obtenu en composant 
avec le foncteur image inverse $f^+\colon E_j\rightarrow E_i$.  On notera que le $I^\circ$-foncteur canonique 
$\cF^\vee\rightarrow \cP^\vee$ est compatible aux foncteurs images inverses.

On fixe dans la suite de cette section le site fibré co-évanescent $\pi\colon E\rightarrow I$ 
et les catégories fibrés associées $\cF/I$, $\cF^\vee/I^\circ$ et $\cP^\vee/I^\circ$.

\subsection{}\label{cftf2}
On note que $E$ est une $\mU$-catégorie. 
On désigne par $\hE$ la catégorie des préfaisceaux de $\mU$-ensembles sur $E$. 
D'après (\cite{agt} VI.5.2) et avec les notations de \ref{notconv4}, 
on a une équivalence de catégories
\begin{eqnarray}\label{cftf2a}
\hE&\stackrel{\sim}{\rightarrow} & \bHom_{I^\circ}(I^\circ,\cP^\vee)\\
F&\mapsto& \{i\mapsto F\circ \alpha_i\},\nonumber
\end{eqnarray}
où $\alpha_i$ est le foncteur \eqref{cftf1b}.  
On identifiera dans la suite $F$ à la section $\{i\mapsto F\circ \alpha_i\}$ qui lui est associée par cette équivalence.

\begin{defi}\label{cftf3}
On appelle {\em v-préfaisceau} de $\mU$-ensembles sur $E$ tout préfaisceau $F$ de $\mU$-ensembles sur $E$ 
tel que pour tout $i\in \ob(I)$, $F\circ \alpha_i$ soit un faisceau sur $E_i$ \eqref{cftf1b}. 
\end{defi}

On désigne par $\hE_{\rv}$ la catégorie des v-préfaisceaux, 
c'est-à-dire la sous-catégorie pleine de $\hE$ formée des v-préfaisceaux. 
On a alors une équivalence de catégories
\begin{eqnarray}\label{cftf3a}
\hE_\rv&\stackrel{\sim}{\rightarrow} & \bHom_{I^\circ}(I^\circ,\cF^\vee)\\
F&\mapsto& \{i\mapsto F\circ \alpha_i\}.\nonumber
\end{eqnarray}

On notera que cette notion n'utilise pas la condition \ref{cftf1}(i).

\subsection{}\label{cftf4}
Suivant (\cite{agt} VI.5.3), on appelle topologie {\em co-évanescente} de $E$ 
la topologie engendrée par les familles de recouvrements $(V_n\rightarrow V)_{n\in \Sigma}$
des deux types suivants:
\begin{itemize}
\item[(v)] Il existe $i\in \ob(I)$ tel que $(V_n\rightarrow V)_{n\in \Sigma}$ soit 
une famille couvrante de $E_i$.
\item[(c)] Il existe une famille couvrante de morphismes $(f_n\colon i_n\rightarrow i)_{n\in \Sigma}$ de $I$ 
telle que $\pi(V)=i$ et pour tout $n\in \Sigma$, $V_n$ soit isomorphe à $f_n^+(V)$.
\end{itemize}
Les recouvrements du type (v) sont dits {\em verticaux}, et ceux du type (c) sont dits {\em cartésiens}. 
Le site ainsi défini est appelé {\em site co-évanescent} associé au site fibré co-évanescent $\pi$ \eqref{cftf1a}; c'est un $\mU$-site. 
On appelle topos {\em co-évanescent} associé au site fibré co-évanescent $\pi$, et l'on note $\tE$, 
le topos des faisceaux de $\mU$-ensembles sur $E$. Pour tout préfaisceau $F$ de $\mU$-ensembles sur $E$, on note $F^a$ le faisceau associé.

\subsection{}\label{cftf10}
Soient $J$ une sous-catégorie pleine et topologiquement génératrice de $I$, $\psi \colon J\rightarrow I$ le foncteur canonique, 
qu'on omettra des notations quand il n'y a aucun risque de confusion.  On désigne par 
\begin{equation}\label{cftf10a}
\pi_J\colon E_J\rightarrow J
\end{equation}
la catégorie fibrée déduite de $\pi$ par changement de base par $\psi$, et par 
\begin{equation}\label{cftf10b}
\Psi\colon E_J\rightarrow E
\end{equation}
la projection canonique (\cite{sga1} VI § 3).  
Pour tout $j\in \ob(J)$, la catégorie fibre $(E_J)_j$ est canoniquement équivalente à la catégorie fibre $E_j$; on les identifiera dans la suite. 
Le topos fibré $\cF_J\rightarrow J$ associé à $\pi_J$ est canoniquement $J$-équivalent 
au topos fibré déduit de $\cF/I$ \eqref{cftf1c} par changement de base par $\psi$. On désigne par  
\begin{eqnarray}
\cF_J^\vee&\rightarrow &J^\circ,\label{cftf10c}\\
\cP_J^\vee&\rightarrow &J^\circ,\label{cftf10d}
\end{eqnarray}
les catégories fibrées déduites de $\cF^\vee/I^\circ$ \eqref{cftf1d} et $\cP^\vee/I^\circ$ \eqref{cftf1e} par changement de base par $\psi^\circ$.
On note $\hE_J$ la catégorie des préfaisceaux de $\mU$-ensembles sur $E_J$.
On vérifie aussitôt que le diagramme de foncteurs 
\begin{equation}\label{cftf10e}
\xymatrix{
{\hE}\ar[r]^-(0.5)\sim\ar[d]_{\hPsi^*}&{\bHom_{I^\circ}(I^\circ,\cP^\vee)}\ar[d]\\
{\hE_J}\ar[r]^-(0.5)\sim&{\bHom_{J^\circ}(J^\circ,\cP_J^\vee)}}
\end{equation}
où les flèches horizontales sont les équivalences de catégories \eqref{cftf2a}, 
$\hPsi^*$ est le foncteur défini par la composition avec $\Psi$ 
et la flèche verticale de droite est le foncteur canonique (\cite{sga1} VI § 3), est commutatif
à isomorphisme canonique près. Par suite, pour tout préfaisceau $F=\{i\mapsto F_i\}$ sur $E$, on a 
\begin{equation}\label{cftf10f}
\hPsi^*(F)=\{j\in J^\circ\mapsto F_{j}\}.
\end{equation}

On munit $E_J$ de la topologie induite par la topologie co-évanescente sur $E$ au moyen du foncteur $\Psi$ (\cite{sga4} III §3),
et on note $\tE_J$ le topos des faisceaux de $\mU$-ensembles sur $E_J$. 
On observera que le foncteur $\Psi$ est pleinement fidèle et que la famille $E_J$ est topologiquement génératrice de $E$.  
Supposons que les catégories $J$ et $E_J$ soient $\mU$-petites. 
D'après (\cite{sga4} III 4.1) et sa preuve, le foncteur $\Psi$ est continu et cocontinu, et le foncteur $\hPsi^*$ induit une équivalence de catégories 
\begin{equation}\label{cftf10g}
\Psi_s\colon \tE\stackrel{\sim}{\rightarrow} \tE_J.
\end{equation}

\begin{lem}\label{cftf11}
Conservons les hypothèses et notations de \ref{cftf10}. Alors, 
\begin{itemize}
\item[{\rm (i)}] Pour tout préfaisceau $G$ sur $E$, on a un isomorphisme canonique 
\begin{equation}\label{cftf11a}
\Psi_s(G^a)\stackrel{\sim}{\rightarrow} \hPsi^*(G)^a,
\end{equation}
où l'exposant $^a$ désigne les faisceaux associés. 
\item[{\rm (ii)}] Pour tout préfaisceau $F=\{j\in J^\circ\mapsto F_j\}$ sur $E_J$,  
il existe un préfaisceau canonique $G$ sur $E$ et un isomorphisme $F\stackrel{\sim}{\rightarrow} \hPsi^*(G)$.
Ce dernier induit un isomorphisme 
\begin{equation}\label{cftf11b}
\Psi_s(G^a)\stackrel{\sim}{\rightarrow} F^a,
\end{equation}
où l'exposant $^a$ désigne les faisceaux associés.
\end{itemize}
\end{lem}

(i) Cela résulte de (\cite{sga4} III 2.3) puisque le foncteur $\hPsi^*$ est continu et cocontinu. 

(ii)  Pour tout objet $i$ de $I$, on définit la catégorie $J_{/i}$ de la façon suivante.  
Les objets de $J_{/i}$ sont les couples $(j,f)$ formés d'un objet $j$ de $J$ et d'un morphisme $f\colon j\rightarrow i$ de $I$. 
Soient $(j,f)$ et $(j',f')$ deux objets de $J_{/i}$. Un morphisme de $(j,f)$ dans $(j',f')$ 
est la donnée d'un morphisme $g\colon j\rightarrow j'$ de $J$ tel que $f=f'\circ g$ dans $I$. 
Pour tout morphisme $h\colon i\rightarrow i'$ de $I$, on a le foncteur 
\begin{equation}\label{cftf11c}
J_{/h}\colon J_{/i}\rightarrow J_{/i'},\ \ \ (j,f)\mapsto (j,h\circ f). 
\end{equation}

Pour tout $i\in \ob(I)$, les préfaisceaux $\hf^*(F_j)$, pour $(j,f)\in \ob(J^\circ_{/i})$, forment naturellement un système projectif. On pose 
\begin{equation}\label{cftf11d}
G_i=\underset{\underset{(j,f)\in J^\circ_{/i}}{\longleftarrow}}{\lim} \hf^*(F_j).
\end{equation}
Pour tout morphisme $h\colon i\rightarrow i'$ de $I$, comme le foncteur $\hh^*\colon \hE_i\rightarrow \hE_{i'}$ admettant un adjoint à gauche, 
il commute aux limites projectives. On en déduit, compte tenu de \eqref{cftf11c}, un morphisme $G_{i'}\rightarrow \hh^*(G_i)$. 
On vérifie aussitôt que la collection $G=\{i\in I^\circ\mapsto G_i\}$ forment un préfaisceau sur $E$. 

Pour tout $j\in \ob(J)$, $(j,\id_j)$ est l'objet final de $J_{/j}$, ce qui implique que $\hPsi^*(G)=F$.

\subsection{}\label{cftf5}
Soient $\pi'\colon E'\rightarrow I'$ un site fibré co-évanescent, 
\begin{equation}\label{cftf5a}
\xymatrix{
E\ar[r]^{\pi}\ar[d]_{\Phi}&I\ar[d]^{\varphi}\\
E'\ar[r]^{\pi'}&I'}
\end{equation}
un diagramme de foncteurs strictement commutatif, {\em i.e.} tel que $\varphi \circ \pi=\pi'\circ \Phi$. On désigne par 
\begin{equation}\label{cftf5b}
\pi'_{I}\colon E'_{I}\rightarrow I
\end{equation} 
la catégorie fibrée déduite de $\pi'$ par changement de base par $\varphi$ (\cite{sga1} VI § 3), et par 
\begin{equation}\label{cftf5c}
\upphi\colon E\rightarrow E'_{I}
\end{equation}
le $I$-foncteur induit par $\Phi$ \eqref{cftf5a}. On suppose le foncteur $\upphi$ {\em cartésien}. 

On désigne par $\cF'\rightarrow I'$ le topos fibré associé au site fibré $\pi'$ et par 
\begin{eqnarray}
\cF'^\vee&\rightarrow &I'^\circ,\label{cftf5f}\\
\cP'^\vee&\rightarrow &I'^\circ,\label{cftf5g}
\end{eqnarray}
les catégories fibrées associées au site fibré $\pi'$, définies dans \eqref{cftf1d} et \eqref{cftf1e}, respectivement.

Pour tout $i\in \ob(I)$, la catégorie fibre $(E'_I)_i$ est canoniquement équivalente à la catégorie fibre $E'_{\varphi(i)}$; on les identifiera dans la suite. 
On voit donc que $\pi'_I$ est un site fibré co-évanescent. Le topos fibré $\cF'_I\rightarrow I$ associé à $\pi'_I$ est canoniquement $I$-équivalent 
au topos fibré déduit de $\cF'/I'$ par changement de base par $\varphi$. On note
\begin{eqnarray}
\cF'^\vee_I&\rightarrow &I^\circ,\label{cftf5d}\\
\cP'^\vee_I&\rightarrow &I^\circ,\label{cftf5e}
\end{eqnarray}
les catégories fibrées associées à $\pi'_I$, définies dans \eqref{cftf1d} et \eqref{cftf1e}, respectivement. 
Elles sont canoniquement $I^\circ$-équivalentes  aux catégories fibrées déduites de 
$\cF'^\vee/I'^\circ$ et $\cP'^\vee/I'^\circ$ par changement de base par $\varphi^\circ$.

Pour tout $i\in \ob(I)$, on désigne par $\Phi_i\colon E_i\rightarrow E'_{\varphi(i)}$ le foncteur induit par $\Phi$,
par $ \hE_i$ la catégorie des préfaisceaux de $\mU$-ensembles sur $E_i$, 
et par $\hPhi_i^*\colon \hE'_{\varphi(i)}\rightarrow \hE_i$ le foncteur défini par composition avec $\Phi_i$. 
Soient $f\colon j\rightarrow i$ un morphisme de $I$, $f'=\varphi(f)\colon \varphi(j)\rightarrow \varphi(i)$.
On a un isomorphisme de foncteurs 
\begin{equation}\label{cftf5h}
\Phi_j \circ f^+ \stackrel{\sim}{\rightarrow} f'^+\circ \Phi_i,
\end{equation}
où $f^+\colon E_i\rightarrow E_j$ (resp. $f'^+\colon E'_{\varphi(i)}\rightarrow E'_{\varphi(j)}$) est le foncteur image inverse par $f$ de $E$ 
(resp. $f'$ de $E'$). Il induit un isomorphisme de foncteurs 
\begin{equation}\label{cftf5i}
\hPhi_i^* \circ \hf'^* \stackrel{\sim}{\rightarrow}  \hf^* \circ \hPhi_j^*, 
\end{equation}
où $\hf^*\colon \hE_j\rightarrow \hE_i$ (resp. $\hf'^*\colon \hE'_{\varphi(j)}\rightarrow \hE'_{\varphi(i)}$) 
est le foncteur image inverse par $f^\circ$ de $\cP^\vee$ (resp. $f'^\circ$ de $\cP'^\vee$) \eqref{cftf1e}. 
Les isomorphismes \eqref{cftf5h} vérifient une relation de cocycle du type (\cite{egr1} (1.1.2.2)), 
qui induit une relation analogue pour les isomorphismes \eqref{cftf5i}. 
D'après (\cite{sga1} VI 12; cf. aussi \cite{egr1} 1.1.2), 
les foncteurs $\hPhi_i^*$ définissent donc un $I^\circ$-foncteur cartésien
\begin{equation}\label{cftf5j}
\cP'^\vee_I\rightarrow \cP^\vee.
\end{equation}

On désigne par 
\begin{equation}\label{cftf5k}
\hPhi^*\colon \hE'\rightarrow \hE
\end{equation}
le foncteur défini par la composition avec $\Phi$. On vérifie aussitôt que le diagramme de foncteurs 
\begin{equation}\label{cftf5l}
\xymatrix{
{\hE'}\ar[r]^-(0.5)\sim\ar[dd]_{\hPhi^*}&{\bHom_{I'^\circ}(I'^\circ,\cP'^\vee)}\ar[d]^u\\
&{\bHom_{I^\circ}(I^\circ,\cP'^\vee_I)}\ar[d]^v\\
{\hE}\ar[r]^-(0.5)\sim&{\bHom_{I^\circ}(I^\circ,\cP^\vee)}}
\end{equation}
où $u$ est le foncteur canonique (\cite{sga1} VI § 3),
$v$ est le foncteur défini par composition avec le foncteur \eqref{cftf5j},  
et les flèches horizontales sont les équivalences de catégories \eqref{cftf2a}, est commutatif
à isomorphisme canonique près. Par suite, pour tout préfaisceau $F'=\{i'\mapsto F'_{i'}\}$ sur $E'$, on a 
\begin{equation}\label{cftf5m}
\hPhi^*(F')=\{i\mapsto \hPhi_i^*(F'_{\varphi(i)})\}.
\end{equation}

\subsection{}\label{cftf6}
Conservons les hypothèses et notations de \ref{cftf5}; supposons, de plus, que les catégories $E$, $I$ et $(E'_{i'})_{i'\in \ob(I')}$ soient $\mU$-petites.   
D'après (\cite{sga4} I 5.1), le foncteur $\hPhi^*$ \eqref{cftf5k} admet un adjoint à gauche 
\begin{equation}\label{cftf6a}
\hPhi_!\colon \hE\rightarrow \hE'.
\end{equation}
De même, pour tout $i\in \ob(I)$, la catégorie $E_i$ étant $\mU$-petite, 
le foncteur $\hPhi_i^* \colon \hE'_{\varphi(i)}\rightarrow \hE_i$ admet un adjoint à gauche 
\begin{equation}\label{cftf6b}
\hPhi_{i!}\colon \hE_i\rightarrow \hE'_{\varphi(i)}.
\end{equation}
Pour tout morphisme $f\colon j'\rightarrow i'$ de $I'$,  
le foncteur $\hf^*\colon \hE'_{j'}\rightarrow \hE'_{i'}$  image inverse par $f^\circ$ de $\cP'^\vee$ 
admet un adjoint à gauche 
\begin{equation}\label{cftf6c}
\hf_!\colon \hE'_{i'}\rightarrow \hE'_{j'}.
\end{equation}

Pour tous morphismes composables $h\colon k'\rightarrow j'$ et $f\colon j'\rightarrow i'$ de $I'$, posant $g=f\circ h\colon k'\rightarrow i'$, on a  un
isomorphisme canonique $\hg^*\stackrel{\sim}{\rightarrow} \hf^*\circ \hh^*$ qui induit par adjonction un isomorphisme
\begin{equation}\label{cftf6d}
\hg_!\stackrel{\sim}{\rightarrow} \hh_!\circ \hf_!.
\end{equation}

Pour tout objet $i'$ de $I'$, on définit la catégorie $I^{i'}_\varphi$ de la façon suivante. Les objets de $I^{i'}_\varphi$ sont les couples $(i,f)$ 
formés d'un objet $i$ de $I$ et d'un morphisme $f\colon i'\rightarrow \varphi(i)$ de $I'$. Soient $(i_1,f_1)$ et $(i_2,f_2)$ deux objets de $I^{i'}_\varphi$. 
Un morphisme de $(i_1,f_1)$ dans $(i_2,f_2)$ est la donnée d'un morphisme $g\colon i_1\rightarrow i_2$ de $I$ tel que $f_2=\varphi(g)\circ f_1$. 

Pour tout morphisme $h\colon i'\rightarrow j'$ de $I'$, on a le foncteur 
\begin{equation}\label{cftf6e}
I^h_\varphi\colon I^{j'}_\varphi\rightarrow I^{i'}_\varphi,\ \ \ (j,g)\mapsto (j,g\circ h). 
\end{equation}

\begin{lem}\label{cftf7}
Conservons les hypothèses et notations de \ref{cftf5} et \ref{cftf6}. Soit, de plus, $F=\{i\mapsto F_i\}$ un préfaisceau sur $E$. 
Alors, pour tout $i'\in \ob(I')$, on a un isomorphisme 
\begin{equation}\label{cftf7a}
\hPhi_!(F)\circ \alpha'_{i'} \stackrel{\sim}{\rightarrow} \underset{\underset{(i,f)\in (I_\varphi^{i'})^\circ}{\longrightarrow}}{\lim}\ \hf_!(\hPhi_{i!}(F_i)),
\end{equation}
où $\alpha'_{i'}\colon E'_{i'}\rightarrow E'$ est l'injection canonique. 
De plus, pour tout morphisme $h\colon i'\rightarrow j'$ de $I'$ et tout objet $(j,g)$ de $I^{j'}_\varphi$, posant $f=g\circ h\colon i'\rightarrow \varphi(j)$
de sorte que $(j,f)=I^h_\varphi(j,g)$ \eqref{cftf6e}, le diagramme 
\begin{equation}\label{cftf7b}
\xymatrix{
{\hh_!(\hg_!(\hPhi_{j!}(F_j)))}\ar[r]\ar[d]&{\hf_!(\hPhi_{j!}(F_j))}\ar[d]\\
{\hh_!(\hPhi_!(F)\circ \alpha'_{j'})}\ar[r]&{\hPhi_!(F)\circ \alpha'_{i'}}}
\end{equation}
où les flèches verticales sont induites par les morphismes \eqref{cftf7a}, la flèche horizontale supérieure est l'isomorphisme canonique \eqref{cftf6d}
et la flèche horizontale inférieure est le morphisme adjoint du morphisme canonique
\begin{equation}\label{cftf7c}
\hPhi_!(F)\circ \alpha'_{j'}\rightarrow \hh^*(\hPhi_!(F)\circ \alpha'_{i'}),
\end{equation}
est commutatif. 
\end{lem}

Soient $U'\in \ob(E')$, $i'=\pi'(U')$. On définit la catégorie $I^{U'}_\Phi$ de la façon suivante. 
Les objets de $I^{U'}_\Phi$ sont les couples $(U,m)$ 
formés d'un objet $U$ de $E$ et d'un morphisme $m\colon U'\rightarrow \Phi(U)$ de $E'$. 
Soient $(U_1,m_1)$ et $(U_2,m_2)$ deux objets de $I^{U'}_\Phi$. 
Un morphisme de $(U_1,m_1)$ dans $(U_2,m_2)$ est la donnée d'un morphisme $n\colon U_1\rightarrow U_2$ de $E$ 
tel que $m_2=\Phi(n)\circ m_1$. D'après la preuve de (\cite{sga4} I 5.1), on a un isomorphisme 
\begin{equation}\label{cftf7d}
\hPhi_!(F)(U')\stackrel{\sim}{\rightarrow} \underset{\underset{(U,m)\in (I^{U'}_\Phi)^\circ}{\longrightarrow}}{\lim}\ F(U).
\end{equation}

Pour tout $(i,f)\in \ob(I^{i'}_\varphi)$, on a  un
isomorphisme canonique $(f^+\circ \Phi_i)^*\stackrel{\sim}{\rightarrow} \hPhi_i^*\circ \hf^*$ qui induit par adjonction un isomorphisme
\begin{equation}\label{cftf7e}
(f^+\circ \Phi_i)_!\stackrel{\sim}{\rightarrow} \hf_!\circ \hPhi_{i!}.
\end{equation}
On définit la catégorie $I^{U'}_{\Phi_i}$ de la façon suivante. 
Les objets de $I^{U'}_{\Phi_i}$ sont les couples $(U,m)$ 
formés d'un objet $U$ de $E_i$ et d'un morphisme $m\colon U'\rightarrow f^+(\Phi_i(U))$ de $E'_{i'}$. 
Soient $(U_1,m_1)$ et $(U_2,m_2)$ deux objets de $I^{U'}_{\Phi_i}$. 
Un morphisme de $(U_1,m_1)$ dans $(U_2,m_2)$ est la donnée d'un morphisme $n\colon U_1\rightarrow U_2$ de $E_i$ 
tel que $m_2=f^+(\Phi_i(n))\circ m_1$. D'après \eqref{cftf7e}, on a un isomorphisme 
\begin{equation}\label{cftf7f}
\hf_!(\hPhi_{i!}(F_i))(U')\stackrel{\sim}{\rightarrow} \underset{\underset{(U,m)\in (I^{U'}_{\Phi_i})^\circ}{\longrightarrow}}{\lim}\ F(U).
\end{equation}

Par ailleurs, on a le foncteur 
\begin{equation}\label{cftf7g}
\varpi\colon I^{U'}_\Phi \rightarrow I^{i'}_\varphi,\ \ \ (U,m)\mapsto (\pi(U),\pi'(m)).
\end{equation}
En effet, on a $\pi'(m)\colon i'\rightarrow \pi'(\Phi(U))=\varphi(\pi(U))$. Comme le foncteur $\upphi$ \eqref{cftf5c} est cartésien, 
pour tout $(i,f)\in \ob(I^{i'}_\varphi)$, la catégorie fibre de $\varpi$ 
au-dessus de $(i,f)$ s'identifie canoniquement à la catégorie $I^{U'}_{\Phi_i}$. 
La proposition s'ensuit compte tenu de \eqref{cftf7e} et \eqref{cftf7f}. 

\begin{rema}\label{cftf15}
Conservons les hypothèses et notations de \ref{cftf5} et \ref{cftf6}. 
\begin{itemize}
\item[(i)] Pour tout $i\in \ob(I)$, le diagramme 
\begin{equation}\label{cftf15a}
\xymatrix{
E_i\ar[r]^-(0.5){\alpha_i}\ar[d]_{\Phi_i}&E\ar[d]^\Phi\\
E'_{\varphi(i)}\ar[r]^{\alpha'_{\varphi(i)}}&{E'}}
\end{equation}
est strictement commutatif. On en déduit que pour tout préfaisceau $F'$ sur $E'$, on a un isomorphisme canonique 
\begin{equation}\label{cftf15b}
\hPhi_i^*(F'\circ \alpha'_{\varphi(i)})\stackrel{\sim}{\rightarrow} \alpha_i\circ \hPhi^*(F').
\end{equation}
\item[(ii)] Avec les hypothèses et notations de \ref{cftf7}, le diagramme  
\begin{equation}\label{cftf15c}
\xymatrix{
{F_i}\ar[rr]^-(0.5){a}\ar[d]&&{F\circ \alpha_i}\ar[d]\\
{\hPhi_i^*(\hPhi_{i!}(F_i))}\ar[r]^-(0.5){\hPhi_i^*(b)}&{\hPhi^*_i(\hPhi_!(F)\circ \alpha'_{\varphi(i)})}\ar[r]^-(0.5){c}&{\hPhi^*(\hPhi_!(F))\circ \alpha_i}}
\end{equation}
où les flèches verticales sont induites par les morphismes d'adjonction, $a$ est l'isomorphisme canonique,
$b\colon \hPhi_{i!}(F_i)\rightarrow \hPhi_!(F)\circ \alpha'_{\varphi(i)}$ est le morphisme induit par l'isomorphisme \eqref{cftf7a} 
et l'objet $(\varphi(i),\id)$ de $I^{\varphi(i)}_\varphi$ et $c$ est l'isomorphisme \eqref{cftf15b}, est commutatif. Cela résulte aussitôt de la preuve de 
\ref{cftf7}, en particulier de la description des catégories fibre de $\varpi$ \eqref{cftf7g}. 
\end{itemize}
\end{rema}

\subsection{}\label{cftf8}
Conservons les hypothèses et notations de \ref{cftf5} et \ref{cftf6}. 
Supposons, de plus, satisfaites les hypothèses suivantes:
\begin{itemize}
\item[(i)] le foncteur $\upphi$ \eqref{cftf5c} est cartésien; 
\item[(ii)] les catégories $E$, $I$ et $(E'_{i'})_{i'\in \ob(I')}$ sont $\mU$-petites;
\item[(iii)] le foncteur $\varphi$ est continu; 
\item[(iv)] pour tout $i\in \ob(I)$, le foncteur $\Phi_i\colon E_i\rightarrow E'_{\varphi(i)}$ est continu. 
\end{itemize}
Les hypothèses (i)-(ii) sont mentionnées pour rappel. 

Le foncteur $\Phi$ \eqref{cftf5a} 
est continu pour les topologies co-évanescentes sur $E$ et $E'$, en vertu de \eqref{cftf5l} et (\cite{agt} VI 5.10). 
Notant $\tE'$ le topos des faisceaux de $\mU$-ensembles sur $E'$, le foncteur $\hPhi^*$ induit un foncteur 
\begin{equation}\label{cftf8a}
\Phi_s\colon \tE'\rightarrow \tE.
\end{equation}
D'après (\cite{sga4} III 1.3), le foncteur $\Phi_s$ admet un adjoint à gauche 
\begin{equation}\label{cftf8b}
\Phi^s\colon \tE\rightarrow \tE'
\end{equation}
tel que pour tout préfaisceau $F$ sur $E$, on ait un isomorphisme canonique fonctoriel
\begin{equation}\label{cftf8c}
\Phi^s(F^a)=(\hPhi_!(F))^a,
\end{equation} 
où l'exposant $^a$ désigne les faisceaux associés. 

Pour tout $i'\in \ob(I')$, notons $\tE'_{i'}$ le topos de faisceaux de $\mU$-ensembles sur $E'_{i'}$. 
Le foncteur $\hPhi^*_i$ induit un foncteur 
\begin{equation}\label{cftf8d}
\Phi_{i,s}\colon \tE'_{\varphi(i)}\rightarrow \tE_i.
\end{equation}
Celui-ci admet un adjoint à gauche 
\begin{equation}\label{cftf8e}
\Phi^s_i\colon \tE_i\rightarrow \tE'_{\varphi(i)}
\end{equation}
tel que pour tout préfaisceau $G$ sur $E_i$, on ait un isomorphisme canonique fonctoriel
\begin{equation}\label{cftf8f}
\Phi^s_i(G^a)=(\hPhi_{i!}(G))^a,
\end{equation} 
où l'exposant $^a$ désigne les faisceaux associés.

\begin{lem}\label{cftf9}
Conservons les hypothèses et notations de \ref{cftf8}. Soit, de plus, $F=\{i\mapsto F_i\}$ un v-préfaisceau sur $E$ \eqref{cftf3}. 
Alors,
\begin{itemize}
\item[{\rm (i)}] Pour tout $i'\in \ob(I')$, les faisceaux $f^*(\Phi_i^s(F_i))$, pour $(i,f)\in \ob((I^{i'}_\varphi)^\circ)$, 
forment naturellement un système inductif de $\tE'_{i'}$. Posons 
\begin{equation}\label{cftf9a}
F'_{i'}=\underset{\underset{(i,f)\in (I^{i'}_\varphi)^\circ}{\longrightarrow}}{\lim}\ f^*(\Phi_i^s(F_i)).
\end{equation}
\item[{\rm (ii)}] La collection $F'=\{i'\in I'^\circ\mapsto F'_{i'}\}$ forme naturellement un v-préfaisceau sur $E'$. 
Pour tout morphisme $h\colon i'\rightarrow j'$ de $I'$ et tout objet $(j,g)$ de $I^{j'}_\varphi$, posant $f=g\circ h\colon i'\rightarrow \varphi(j)$
de sorte que $(j,f)=I^h_\varphi(j,g)$ \eqref{cftf6e}, le diagramme 
\begin{equation}\label{cftf9b}
\xymatrix{
{h^*(g^*(\Phi_j^s(F_i)))}\ar[r]\ar[d]&{f^*(\Phi_j^s(F_i))}\ar[d]\\
{h^*(F'_{j'})}\ar[r]&{F'_{i'}}}
\end{equation}
où les flèches verticales sont les morphismes canoniques \eqref{cftf9a}, la flèche horizontale supérieure est l'isomorphisme canonique 
et la flèche horizontale inférieure est le morphisme adjoint du morphisme $F'_{j'}\rightarrow h_*(F'_{i'})$ définissant la structure de préfaisceau sur
$\{i'\mapsto F'_{i'}\}$, est commutatif. 
\item[{\rm (iii)}] On a un isomorphisme canonique fonctoriel 
\begin{equation}\label{cftf9c}
\Phi^s(F^a)\stackrel{\sim}{\rightarrow} F'^a,
\end{equation} 
où l'exposant $^a$ désigne les faisceaux associés. 
\end{itemize}
\end{lem}

Les propositions (i)-(ii) sont immédiates et la proposition (iii) résulte de \ref{cftf7}, \eqref{cftf8c}, \eqref{cftf8f},
(\cite{agt} VI 5.17) et le fait que le foncteur faisceau associé commute aux limites inductives (\cite{sga4} II 4.1).

\subsection{}\label{cftf12}
Reprenons les hypothèses et notations de \ref{cftf5}. On désigne par $\mI$ la catégorie suivante. 
Les objets de $\mI$ sont les triplets $(i,i',f)$ formés d'un objet $i$ de $I$, d'un objet $i'$ de $I'$ et 
d'un morphisme $f\colon i'\rightarrow \varphi(i)$ de $I'$. Soient $(i_1,i'_1,f_1)$ et $(i_2,i'_2,f_2)$ deux objets de $\mI$. 
Un morphisme de $(i_1,i'_1,f_1)$ dans $(i_2,i'_2,f_2)$ est la donnée d'un morphisme $m\colon i_1\rightarrow i_2$ de $I$
et d'un morphisme $m'\colon i'_1\rightarrow i'_2$ de $I'$ tels que $\varphi(m)\circ f_1=f_2\circ m'$. 
Considérons les foncteurs 
\begin{eqnarray}
\ttb\colon \mI&\rightarrow& I, \ \ \ (i,i',f) \mapsto i, \label{cftf12a}\\
\tts\colon \mI&\rightarrow& I', \ \  (i,i',f) \mapsto i'. \label{cftf12b}
\end{eqnarray}
On notera que pour tout $i'\in \ob(I')$, la catégorie $I^{i'}_\varphi$ \eqref{cftf6} s'identifie à la catégorie fibre du foncteur $\tts$ au-dessus de $i'$. 

Soient $\mJ$ une sous-catégorie pleine de $\mI$, $J$ et $J'$ ses images essentielles par les foncteurs $\ttb$ et $\tts$ respectivement. 

Pour tout $i\in \ob(I)$, on désigne par $J_{/i}$ la catégorie des couples $(j,m)$ formés d'un objet $j$ de $J$
et d'un morphisme $m\colon j\rightarrow i$ de $I$, les morphismes de $J_{/i}$ étant naturellement induits par ceux de $J$.
Pour tout $i'\in \ob(I')$, on définit de même la catégorie $J'_{/i'}$. 
Pour tout $(i,i',f)\in \ob(\mI)$, on désigne par $\mJ_{/(i,i',f)}$ la catégorie des quintuples $(j,j',g,m,m')$
formés d'un objet $(j,j',g)$ de $\mJ$ et d'un morphisme $(m,m')\colon (j,j',g)\rightarrow (i,i',f)$ de $\mI$, {\em i.e.} un diagramme commutatif
\begin{equation}\label{cftf12c}
\xymatrix{
j'\ar[d]_{m'}\ar[r]^g&{\varphi(j)}\ar[d]^{\varphi(m)}\\
i'\ar[r]^f&{\varphi(i)}}
\end{equation}
les morphismes de $\mJ_{/(i,i',f)}$ étant naturellement induits par ceux de $\mJ$. 

On désigne par $\pi_J\colon E_J\rightarrow J$ la catégorie fibrée déduite de $\pi$ par changement de base par le foncteur canonique $J\rightarrow I$, 
par $\Psi\colon E_J\rightarrow E$ la projection canonique (\cite{sga1} VI § 3) et par $\hPsi^*\colon \hE_J\rightarrow \hE$ le foncteur
défini par composition avec $\Psi$. De même, on désigne par $\pi'_{J'}\colon E'_{J'}\rightarrow J'$ la catégorie fibrée déduite de 
$\pi'$ par changement de base par le foncteur canonique $J'\rightarrow I'$, par 
$\Psi'\colon E'_{J'}\rightarrow E'$ la projection canonique et par $\hPsi'^*\colon \hE'_{J'}\rightarrow \hE'$ le foncteur
défini par composition avec $\Psi'$.

\begin{defi}\label{cftf120}
Conservons les hypothèses et notations de \ref{cftf12}. 
Soient $F=\{j\in J^\circ\mapsto F_j\}$ un préfaisceau sur $E_J$, $F'=\{j'\in J'^\circ\mapsto F'_{j'}\}$ un préfaisceau sur $E'_{J'}$.
On appelle {\em $\mJ$-système de $\Phi$-morphismes compatibles de $F$ dans $F'$} la donnée, pour tout objet $(j,j',g)$ de $\mJ$, d'un morphisme 
\begin{equation}\label{cftf120a}
u_{(j,j',g)}\colon F_j\rightarrow \hPhi_j^*(\hg^*(F'_{j'}))
\end{equation}
de $\hE_j$, tel que pour tout morphisme $(m,m')\colon (j_1,j'_1,g_1)\rightarrow (j_2,j'_2,g_2)$ de $\mJ$, 
\begin{equation}\label{cftf120b}
\xymatrix{
{j'_1}\ar[r]^{g_1}\ar[d]_{m'}&{\varphi(j_1)}\ar[d]^{\varphi(m)}\\
{j'_2}\ar[r]_{g_2}&{\varphi(j_2)}}
\end{equation}
posant $n=\varphi(m)$, le diagramme 
\begin{equation}\label{cftf120c}
\xymatrix{
{F_{j_2}}\ar[rr]^-(0.5){u_{(j_2,j'_2,g_2)}}\ar[d]_{a}&&{\hPhi_{j_2}^*(\hg^*_2(F'_{j'_2}))}\ar[rr]^-(0.5){\hPhi_{j_2}^*(\hg^*_2(a'))}&&{\hPhi_{j_2}^*(\hg^*_2(\hm'^*(F'_{j'_1})))}\ar[d]^{v_2}\\
{\hm^*(F_{j_1})}\ar[rr]^-(0.5){\hm^*(u_{(j_1,j'_1,g_1)})}&&{\hm^*(\hPhi_{j_1}^*(\hg^*_1(F'_{j'_1})))}\ar[rr]^{v_1}&&{\hPhi_{j_2}^*(\hn^*(\hg^*_1(F'_{j'_1})))}}
\end{equation}
où $a$ (resp. $a'\colon F'_{j'_2}\rightarrow \hm'^*(F'_{j'_1})$) est le morphisme de transition du faisceau $F$ (resp. $F'$), 
$v_1$ est l'isomorphisme \eqref{cftf5i} et $v_2$ est l'isomorphisme induit par le diagramme commutatif \eqref{cftf120a}, soit commutatif. 
\end{defi}

\subsection{}\label{cftf122}
Conservons les hypothèses et notations de \ref{cftf12} et reprenons, de plus, celles de \ref{cftf8}. 
Soient $F=\{j\in J^\circ\mapsto F_j\}$ un v-préfaisceau sur $E_J$ \eqref{cftf3}, 
$F'=\{j'\in J'^\circ\mapsto F'_{j'}\}$ un v-préfaisceau sur $E'_{J'}$. 
La donnée d'un $\mJ$-système de $\Phi$-morphismes compatibles de $F$ dans $F'$ \eqref{cftf120} est équivalente par adjonction 
à la donnée, pour tout objet $(j,j',g)$ de $\mJ$, d'un morphisme de $\tE'_j$, 
\begin{equation}\label{cftf122a}
u'_{(j,j',g)}\colon g^*(\Phi_j^s(F_j))\rightarrow F'_{j'},
\end{equation}
tel que pour tout morphisme $(m,m')\colon (j_1,j'_1,g_1)\rightarrow (j_2,j'_2,g_2)$ de $\mJ$, 
\begin{equation}\label{cftf122b}
\xymatrix{
{j'_1}\ar[r]^{g_1}\ar[d]_{m'}&{\varphi(j_1)}\ar[d]^{\varphi(m)}\\
{j'_2}\ar[r]_{g_2}&{\varphi(j_2)}}
\end{equation}
posant $n=\varphi(m)$, le diagramme 
\begin{equation}\label{cftf122c}
\xymatrix{
{g_1^*(n^*(\Phi_{j_2}^s(F_{j_2})))}\ar[rr]^-(0.5){v_2}\ar[d]_{v_1}&&{m'^*(g_2^*(\Phi_{j_2}^s(F_{j_2})))}\ar[rr]^-(0.5){m'^*(u'_{(j_2,j'_2,f_2)})}&&{m'^*(F'_{j'_2})}\ar[d]^{a'}\\
{g_1^*(\Phi_{j_1}^s(m^*(F_{j_2})))}\ar[rr]^-(0.5){g_1^*(\Phi_{j_1}^*(a))}&&{g_1^*(\Phi_{j_1}^s(F_{j_1}))}\ar[rr]^-(0.5){u'_{(j_1,j'_1,f_1)}}&&{F'_{j'_1}}}
\end{equation}
où $a\colon m^*(F_{j_2})\rightarrow F_{j_1}$ (resp. $a'$) est le morphisme de transition du v-préfaisceau $F$ (resp. $F'$), 
$v_1$ est l'isomorphisme adjoint de \eqref{cftf5i} et $v_2$ est l'isomorphisme induit par le diagramme commutatif \eqref{cftf122b}, soit commutatif.

\begin{lem}\label{cftf123}
Conservons les hypothèses et notations de \ref{cftf12}; supposons de plus que les conditions suivantes soient remplies:
\begin{itemize}
\item[{\rm (i)}] les produits fibrés sont représentables dans $I$ et $I'$ et $\varphi$ commute à ces produits; 
\item[{\rm (ii)}] pour tout objet $(i,i',f)$ de $\mI$, tout objet $\ell'$ de $J'$ et tout morphisme $u\colon \ell'\rightarrow i'$ de $I'$, 
il existe un objet $(j,j',g)$ de $\mJ$, un morphisme $v\colon \ell'\rightarrow j'$ de $I'$ et un morphisme $(m,m')\colon (j,j',g)\rightarrow (i,i',f)$ de $\mI$
tels que $u=m'\circ v$;
\begin{equation}\label{cftf123a}
\xymatrix{
\ell'\ar[r]_v\ar@/^1pc/[rr]^u&j'\ar[r]_{m'}\ar[d]_g&i'\ar[d]^f\\
&{\varphi(j)}\ar[r]^{\varphi(m)}&{\varphi(i)}}
\end{equation}
\end{itemize}
Soient $C$ une catégorie où les limites projectives sont représentables, $(i,i',f)$ un objet de $\mI$, $F\colon J'^\circ_{/i'}\rightarrow C$  
un foncteur. Alors, le morphisme canonique
\begin{equation}\label{cftf123b}
\underset{\underset{J'^\circ_{/i'}}{\longleftarrow}}{\lim} \ F\rightarrow \underset{\underset{\mJ^\circ_{/(i,i',f)}}{\longleftarrow}}{\lim} F\circ \tts^\circ,
\end{equation}
où on a encore noté $\tts\colon \mJ_{/(i,i',f)} \rightarrow J'_{/i'}$ le foncteur induit par $\tts$ \eqref{cftf12b}, est un isomorphisme. 
\end{lem} 

On notera d'abord qu'il résulte de (i) que la limite projective d'un diagramme 
\begin{equation}\label{cftf123c}
\xymatrix{
&{(i_1,i'_1,f_1)}\ar[d]\\
{(i_2,i'_2,f_2)}\ar[r]&{(i,i',f)}}
\end{equation}
de morphismes de $\mI$ est représentable par l'objet $(i_1\times_ii_2,i'_1\times_{i'} i'_2, f_1\times f_2)$ de $\mI$.  

D'après (ii), pour tout objet $\ell'$ de $J'_{/i'}$, il existe un objet $(j,j',g)$ de $\mJ_{/(i,i',f)}$ et un $i'$-morphisme $v\colon \ell'\rightarrow j'$. 
Le morphisme composé 
\begin{equation}\label{cftf123d}
\underset{\underset{\mJ^\circ_{/(i,i',f)}}{\longleftarrow}}{\lim} F\circ \tts^\circ \rightarrow F(j')\rightarrow F(\ell'),
\end{equation}
où la première flèche est le morphisme canonique et la seconde flèche est induite par $v$, est indépendant des choix de $(j,j',g)$ et $v$. 
En effet, si $(j_1,j'_1,g_1)$ et  $(j_2,j'_2,g_2)$ sont deux objets de $\mJ_{/(i,i',f)}$ 
et si $v_1\colon \ell'\rightarrow j'_1$ et $v_2\colon \ell'\rightarrow j'_2$ sont deux $i'$-morphismes, 
il existe, d'après (ii), un objet $(j_3,j'_3,g_3)$ de $\mJ$, un morphisme $v_3\colon \ell'\rightarrow j'_3$ de $I'$ 
et un morphisme $(m,m')\colon (j_3,j'_3,g_3)\rightarrow (j_1\times_ij_2,j'_1\times_{i'} j'_2, g_1\times g_2)$ de $\mI$ 
tels que $v_1\times v_2=m'\circ v_3$;
\begin{equation}\label{cftf123e}
\xymatrix{
\ell'\ar[r]_{v_3}\ar@/^1pc/[rr]^-(0.5){v_1\times v_2}&j'_3\ar[r]_-(0.5){m'}\ar[d]_{g_3}&{j'_1\times_{i'} j'_2}\ar[d]^{g_1\times g_2}\\
&{\varphi(j_3)}\ar[r]^-(0.5){\varphi(m)}&{\varphi(j_1\times_ij_2)}}
\end{equation}
Il s'ensuit aussitôt que le morphisme \eqref{cftf123b} est un isomorphisme.

\begin{lem}\label{cftf121}
Conservons les hypothèses et notations de \ref{cftf12}; supposons, de plus, l'une des deux conditions suivantes remplie:
\begin{itemize}
\item[{\rm (a)}] pour tout $(i,i',f)\in \ob(\mI)$, le foncteur $\mJ_{/(i,i',f)}\rightarrow J'_{/i'}$ induit par le foncteur $\tts$ \eqref{cftf12b}
est cofinal {\rm (\cite{sga4} I 8.1.1)}; ou 
\item[{\rm (b)}] les conditions \ref{cftf123}{\rm (i)-(ii)} sont satisfaites. 
\end{itemize}
Soient $F=\{j\in J^\circ\mapsto F_j\}$ un préfaisceau sur $E_J$, $F'=\{j'\in J'^\circ\mapsto F'_{j'}\}$ un préfaisceau sur $E'_{J'}$, 
\begin{equation}\label{cftf121a}
u_{(j,j',g)}\colon F_j\rightarrow \hPhi_j^*(\hg^*(F'_{j'})), \ \ \ (j,j',g)\in \ob(\mJ),
\end{equation}
un $\mJ$-système de $\Phi$-morphismes compatibles de $F$ dans $F'$. Alors, il existe un préfaisceau canonique $G=\{i\in I^\circ\mapsto G_i\}$ 
sur $E$, un préfaisceau canonique $G'=\{i'\in I'^\circ\mapsto G'_{i'}\}$ sur $E'$, 
un $\mI$-système canonique de $\Phi$-morphismes compatibles de $G$ dans $G'$, 
\begin{equation}\label{cftf121b}
v_{(i,i',f)}\colon G_i\rightarrow \hPhi_i^*(\hf^*(G'_{i'})), \ \ \ (i,i',f)\in \ob(\mI),
\end{equation}
et des isomorphismes $h\colon \hPsi^*(G)\stackrel{\sim}{\rightarrow}F$
et $h'\colon \hPsi'^*(G')\stackrel{\sim}{\rightarrow}F'$ tels que pour tout $(j,j',g)\in \ob(\mJ)$, on ait 
\begin{equation}\label{cftf121c}
u_{(j,j',f)}\circ h_j = \hPhi_i^*(\hg^*(h'_{j'}))\circ v_{(j,j',g)}.
\end{equation}
\end{lem}

Pour tout objet $i$ de $I$, posons 
\begin{equation}\label{cftf121d}
G_i=\underset{\underset{(j,m)\in J^\circ_{/i}}{\longleftarrow}}{\lim} \hm^*(F_j).
\end{equation}
La collection $G=\{i\in I^\circ\mapsto G_i\}_{i\in I^\circ}$ forme alors un préfaisceau sur $E$ et on a un isomorphisme canonique 
$h\colon \hPsi^*(G) \stackrel{\sim}{\rightarrow} F$. Pour tout objet $i'$ de $I'$, posons 
\begin{equation}\label{cftf121e}
G'_{i'}=\underset{\underset{(j',m')\in J'^\circ_{/i'}}{\longleftarrow}}{\lim} \hm'^*(F'_{j'}).
\end{equation}
La collection $G'=\{i'\in I'^\circ\mapsto G'_{i'}\}_{i'\in I'^\circ}$ forme alors un préfaisceau sur $E'$ et on a un isomorphisme canonique 
$h'\colon \hPsi'^*(G')\stackrel{\sim}{\rightarrow} F'$.

Soient $(i,i',f)$ un objet de $\mI$, $(j',j',g,m,m')$ un objet de $\mJ_{/(i,i',f)}$;
\begin{equation}\label{cftf121f}
\xymatrix{
{j'}\ar[r]^-(0.5)g\ar[d]_{m'}&{\varphi(j)}\ar[d]^{\varphi(m)}\\
{i'}\ar[r]^-(0.5)f&{\varphi(i)}}
\end{equation}
Posant $n=\varphi(m)$, on définit le morphisme $v_{(j,j',g,m,m')}$ par le diagramme commutatif
\begin{equation}\label{cftf121g}
\xymatrix{
{G_i}\ar[rrr]^-(0.5){v_{(j,j',g,m,m')}}\ar[d]_{a}&&&{\hPhi_i^*(\hf^*(\hm'^*(F'_{j'})))}\\
{\hm^*(F_j)}\ar[rr]^-(0.5){\hm^*(u_{(j,j',g)})}&&{\hm^*(\hPhi_j^*(\hg^*(F'_{j'})))}\ar[r]^{w_1}&{\hPhi_i^*(\hn^*(\hg^*(F'_{j'})))}\ar[u]_{w_2}}
\end{equation}
où $a$ est le morphisme canonique \eqref{cftf121d}, 
$w_1$ est l'isomorphisme \eqref{cftf5i} et $w_2$ est l'isomorphisme induit par le diagramme commutatif \eqref{cftf121f}. 
D'après les relations \eqref{cftf120c}, les morphismes $v_{(j,j',g,m,m')}$ sont compatibles pour $(j,j',g,m,m')\in \ob(\mJ_{/(i,i',f)})$.  
Compte tenu des hypothèses et de \ref{cftf123}, passant à la limite projective sur $\mJ_{/(i,i',f)}^\circ$, on obtient un morphisme 
\begin{equation}\label{cftf121h}
v_{(i,i',f)}\colon G_i\rightarrow \hPhi_i^*(\hf^*(G'_{i'})).
\end{equation}
En effet, les foncteurs $\hPhi_i^*$ et $\hf^*$ commutent aux limites projectives puisqu'ils admettent des adjoints à gauche. 
On vérifie péniblement que les morphismes $v_{(i,i',f)}$ pour $(i,i',f)\in \ob(\mI)$ forment 
un $\mI$-système de $\Phi$-morphismes compatibles de $G$ dans $G'$ qui vérifie les conditions requises. 

\subsection{}\label{cftf13}
Conservons les hypothèses et notations de \ref{cftf121} et reprenons, de plus, celles de  \ref{cftf6}.
Pour tout objet $(i,i',f)$ de $\mI$, on désigne par 
\begin{equation}\label{cftf13a}
v'_{(i,i',f)}\colon \hf_!(\Phi_{i!}(G_i))\rightarrow G'_{i'}
\end{equation}
le morphisme adjoint de $v_{(i,i',f)}$ \eqref{cftf121b}. 
D'après \eqref{cftf120c}, pour tout $i'\in \ob(I')$, les morphismes $v'_{(i,i',f)}$, pour $(i,f)\in I^{i'}_{\varphi}$, 
sont compatibles. Ils induisent donc un morphisme 
\begin{equation}\label{cftf13b}
v'_{i'}\colon \underset{\underset{(i,f)\in (I^{i'}_{\varphi})^\circ}{\longrightarrow}}{\lim}\hf_!(\Phi_{i!}(G_i))\rightarrow G'_{i'}.
\end{equation}
Compte tenu de \ref{cftf7} et de \eqref{cftf120c}, les morphismes $v'_{i'}$, pour $i'\in \ob(I')$, définissent un morphisme de $\hE'$,
\begin{equation}\label{cftf13c}
v'\colon \Phi_!(G)\rightarrow G'.
\end{equation}

\subsection{}\label{cftf14}
Conservons les hypothèses et notations de \ref{cftf121} et reprenons, de plus, celles de  \ref{cftf8}. 
Supposons aussi que la sous-catégorie $J$ (resp. $J'$) de $I$ (resp. $I'$) soit topologiquement génératrice. 
D'après \ref{cftf10}, le foncteur $\Psi\colon E_J\rightarrow E$ est continu et cocontinu et le foncteur $\hPsi^*$ induit une équivalence de catégories 
\begin{equation}\label{cftf14a}
\Psi_s\colon \tE\stackrel{\sim}{\rightarrow} \tE_J.
\end{equation}
De même, le foncteur $\Psi'\colon E'_{J'}\rightarrow E'$ est continu et cocontinu et le foncteur $\hPsi'^*$ induit une équivalence de catégories 
\begin{equation}\label{cftf14b}
\Psi'_s\colon \tE'\stackrel{\sim}{\rightarrow} \tE'_{J'}.
\end{equation}

Soient $F=\{j\in J^\circ\mapsto F_j\}$ un v-préfaisceau sur $E_J$ \eqref{cftf3}, 
$F'=\{j'\in J'^\circ\mapsto F'_{j'}\}$ un v-préfaisceau sur $E'_{J'}$, 
\begin{equation}\label{cftf14c}
u'_{(j,j',g)}\colon g^*(\Phi_j^s(F_j))\rightarrow F'_{j'}, \ \ \ (j,j',g)\in \ob(\mJ), 
\end{equation}
un $\mJ$-système de $\Phi$-morphismes compatibles de $F$ dans $F'$ \eqref{cftf122}.
D'après \ref{cftf121} et sa preuve, il existe un v-préfaisceau canonique $G=\{i\in I^\circ\mapsto G_i\}$ 
sur $E$, un v-préfaisceau canonique $G'=\{i'\in I'^\circ\mapsto G'_{i'}\}$ sur $E'$, 
un $\mI$-système canonique de $\Phi$-morphismes compatibles de $G$ dans $G'$, 
\begin{equation}\label{cftf14f}
v'_{(i,i',f)}\colon f^*(\Phi^s_i(G_i))\rightarrow G'_{i'}, \ \ \ (i,i',f)\in \ob(\mI),
\end{equation}
et des isomorphismes $h\colon \hPsi^*(G)\stackrel{\sim}{\rightarrow}F$
et $h'\colon \hPsi'^*(G')\stackrel{\sim}{\rightarrow}F'$ tels que pour tout $(j,j',g)\in \ob(\mJ)$, on ait 
\begin{equation}\label{cftf14g}
u'_{(j,j',f)}\circ g^*(\Phi_j^s(h_j)) = h'_{j'}\circ v'_{(i,i',f)}.
\end{equation}

D'après \ref{cftf11}, on a des isomorphismes canoniques $\Psi_s(G^a)\stackrel{\sim}{\rightarrow}F^a$ et 
$\Psi'_s(G'^a)\stackrel{\sim}{\rightarrow}F'^a$, où l'exposant~$^a$ désigne les faisceaux associés. 
Pour tout $i'\in \ob(I')$, les morphismes $v'_{(i,i',f)}$, pour $(i,f)\in I^{i'}_{\varphi}$, 
sont compatibles \eqref{cftf122c}. Ils induisent donc un morphisme 
\begin{equation}\label{cftf14h}
v'_{i'}\colon \underset{\underset{(i,f)\in (I^{i'}_{\varphi})^\circ}{\longrightarrow}}{\lim}f^*(\Phi_i^s(G_i))\rightarrow G'_{i'}.
\end{equation}
Compte tenu de \ref{cftf9} et \eqref{cftf122c}, les morphismes $v'_{i'}$, pour $i'\in \ob(I')$, induisent un morphisme de $\tE'$,
\begin{equation}\label{cftf14i}
v'\colon \Phi^s(G^a)\rightarrow G'^a.
\end{equation}

\chapter[\'Etude locale]{\texorpdfstring{Correspondance de Simpson $p$-adique et modules de Hodge-Tate. \'Etude locale}
{Correspondance de Simpson p-adique et modules de Hodge-Tate. \'Etude locale}}\label{cspel}

\section{\texorpdfstring{Hypothèses et notations. Déformations infinitésimales $p$-adiques}{Hypothèses et notations. Déformations infinitésimales p-adiques}}
\label{definf}

\subsection{}\label{definf1}
Dans ce chapitre, $K$ désigne un corps de valuation discrète complet de 
caractéristique $0$, à corps résiduel parfait $k$ de caractéristique $p>0$,  
$\co_K$ l'anneau de valuation de $K$, $W$ l'anneau des vecteurs de Witt à coefficients dans $k$ relatif à $p$, 
$\oK$ une clôture algébrique de $K$, $\co_\oK$ la clôture intégrale de $\co_K$ dans $\oK$,
$\fm_\oK$ l'idéal maximal de $\co_\oK$ et $G_K$ le groupe de Galois de $\oK$ sur $K$.
On note $\co_C$ le séparé complété $p$-adique de $\co_\oK$, $\fm_C$ son idéal maximal,
$C$ son corps des fractions et $v$ sa valuation, normalisée par $v(p)=1$. 
On désigne par $\hmZ(1)$ et $\mZ_p(1)$ 
les $\mZ[G_K]$-modules 
\begin{eqnarray}
\hmZ(1)&=&\underset{\underset{n\geq 1}{\longleftarrow}}{\lim}\ \mu_{n}(\co_{\oK}),\label{definf1aa}\\
\mZ_p(1)&=&\underset{\underset{n\geq 0}{\longleftarrow}}{\lim}\ \mu_{p^n}(\co_{\oK}),\label{definf1a}
\end{eqnarray}  
où $\mu_{n}(\co_\oK)$ désigne le sous-groupe des racines $n$-ièmes de l'unité dans $\co_\oK$. 
Pour tout $\mZ_p[G_K]$-module $M$ et tout entier $n$, on pose $M(n)=M\otimes_{\mZ_p}\mZ_p(1)^{\otimes n}$.

Pour tout $\mZ_p$-module $A$, on note $\hA$ son séparé complété $p$-adique.

\subsection{} 
On pose $S=\Spec(\co_K)$, $\oS=\Spec(\co_\oK)$ et $\coS=\Spec(\co_C)$. 
On note $s$ (resp.  $\eta$, resp. $\oeta$) le point fermé de $S$ (resp.  générique de $S$, resp. générique de $\oS$).
Pour tout entier $n\geq 1$, on pose $S_n=\Spec(\co_K/p^n\co_K)$. Pour tout $S$-schéma $X$, on pose 
\begin{equation}\label{definf1c}
\oX=X\times_S\oS,  \ \ \ \coX=X\times_S\coS \ \ \ {\rm et}\ \ \  X_n=X\times_SS_n.
\end{equation} 

On munit $S$ de la structure logarithmique $\cM_S$ définie par son point fermé, 
et $\oS$ et $\coS$ des structures logarithmiques $\cM_\oS$ et $\cM_\coS$ images inverses de $\cM_S$.

\subsection{}\label{definf2}
Comme $\co_\oK$ est un anneau de valuation non discrète de hauteur $1$, 
il est loisible de développer la $\alpha$-algèbre (ou presque-algèbre) sur cet anneau (\cite{ag} 2.10.1) (cf. \cite{ag} 2.6-2.10).   
On choisit un système compatible $(\beta_n)_{n>0}$ 
de racines $n$-ièmes de $p$ dans $\co_\oK$. Pour tout nombre rationnel $\varepsilon>0$, 
on pose $p^\varepsilon=(\beta_n)^{\varepsilon n}$, où $n$ est un entier $>0$ tel que $\varepsilon n$ soit entier.

\subsection{}\label{definf3}
Tous les anneaux des vecteurs de Witt considérés dans ce chapitre sont relatifs à $p$ (cf. \ref{notconv1}). 
On désigne par $\co_{\oK^\flat}$ la limite projective du système projectif $(\co_\oK/p\co_\oK)_{\mN}$ 
dont les morphismes de transition sont les itérés de l'endomorphisme de Frobenius absolu de $\co_\oK/p\co_\oK$; 
\begin{equation}\label{definf3a}
\co_{\oK^\flat}= \underset{\underset{\mN}{\longleftarrow}}{\lim}\ \co_\oK/p\co_\oK.
\end{equation}
C'est un anneau de valuation non-discrète, de hauteur $1$, complet et 
parfait de caractéristique $p$ (\cite{ag} 4.8.1 et 4.8.2). On note $\oK^\flat$ son corps des fractions
et $\fm_{\oK^\flat}$ son idéal maximal. 

On fixe une suite $(p_n)_{n\geq 0}$ d'éléments de $\co_\oK$ telle que $p_0=p$ et $p_{n+1}^p=p_n$ (pour tout $n\geq 0$) et 
on note $\varpi$ l'élément associé de $\co_{\oK^\flat}$. On pose  
\begin{equation}\label{definf3b}
\xi=[\varpi]-p \in \rW(\co_{\oK^\flat}),
\end{equation}
où $[\ ]$ est le représentant multiplicatif. Reprenons les notations de \ref{eipo3} pour $A=\co_\oK$, en particulier, posons 
\begin{equation}\label{definf3d}
\cA_2(\co_\oK)=\rW(\co_{\oK^\flat})/\ker(\theta)^2,
\end{equation}
et notons $\theta_2\colon \cA_2(\co_{\oK})\rightarrow \co_C$ l'homomorphisme induit par $\theta$ \eqref{eipo3d}.
D'après \ref{eipuf6}, la suite 
\begin{equation}\label{definf3c}
0\longrightarrow \rW(\co_{\oK^\flat})\stackrel{\cdot \xi}{\longrightarrow} \rW(\co_{\oK^\flat})
\stackrel{\theta}{\longrightarrow} \co_C \longrightarrow 0
\end{equation}
est exacte. Elle induit donc une suite exacte 
\begin{equation}\label{definf3e}
0\longrightarrow \co_C\stackrel{\cdot \xi}{\longrightarrow} \cA_2(\co_\oK)
\stackrel{\theta_2}{\longrightarrow} \co_C \longrightarrow 0,
\end{equation}
où on a encore noté $\cdot \xi$ le morphisme induit par la multiplication par $\xi$ dans $\cA_2(\co_\oK)$. 
L'idéal de carré nul $\ker(\theta_2)$ de $\cA_2(\co_\oK)$ est un $\co_C$-module libre de base $\xi$. Il sera noté $\xi\co_C$. 
Contrairement à $\xi$, il ne dépend pas du choix de la suite $(p_n)_{n\geq 0}$. 
On note $\xi^{-1}\co_C$ le $\co_C$-module dual de $\xi\co_C$. 
Pour tout $\co_C$-module $M$, on désigne les $\co_C$-modules $M\otimes_{\co_C}(\xi \co_C)$ 
et $M\otimes_{\co_C}(\xi^{-1} \co_C)$ simplement par $\xi M$ et $\xi^{-1} M$, respectivement.

Le groupe de Galois $G_K$ agit naturellement sur $\rW(\co_{\oK^\flat})$ par des automorphismes d'anneaux,
et l'homomorphisme $\theta$ est $G_K$-équivariant. On en déduit une action de $G_K$ 
sur $\cA_2(\co_{\oK})$ par des automorphismes d'anneaux tel que l'homomorphisme $\theta_2$ soit 
$G_K$-équivariant. 

\subsection{}\label{definf17}
On a un homomorphisme canonique 
\begin{equation}\label{definf17a}
\mZ_p(1)\rightarrow \co_{\oK^\flat}^\times.
\end{equation} 
Pour tout $\zeta\in \mZ_p(1)$, on note encore $\zeta$ son image dans $\co_{\oK^\flat}^\times$.  
Comme $\theta([\zeta]-1)=0$, on obtient un homomorphisme de groupes
\begin{equation}\label{definf17b}
\mZ_p(1)\rightarrow \cA_2(\co_\oK),\ \ \ 
\zeta\mapsto\log([\zeta])=[\zeta]-1,
\end{equation}
dont l'image est contenue dans $\ker(\theta_2)=\xi\co_C$. Celui-ci est clairement $\mZ_p$-linéaire. D'après (\cite{agt} II.9.18),
son image engendre l'idéal $p^{\frac{1}{p-1}}\xi\co_C$ de $\cA_2(\co_\oK)$, et le morphisme $\co_C$-linéaire induit
\begin{equation}\label{definf17c}
\co_C(1)\rightarrow p^{\frac{1}{p-1}}\xi \co_C
\end{equation}
est un isomorphisme.

\subsection{}\label{definf4}
L'anneau $\rW(\co_{\oK^\flat})$ étant naturellement muni d'une structure de $W$-algèbre \eqref{eip2}, on pose 
\begin{equation}\label{definf4a}
\rW_{\co_K}(\co_{\oK^\flat})=\rW(\co_{\oK^\flat})\otimes_W\co_K 
\end{equation}
et on note $\theta_{\co_K}\colon \rW_{\co_K}(\co_{\oK^\flat})\rightarrow \co_C$ l'homomorphisme induit par $\theta$ \eqref{definf3c}.
On fixe une uniformisante $\pi$ de $\co_K$ et une suite $(\pi_n)_{n\geq 0}$ 
d'éléments de $\co_\oK$ telle que $\pi_0=\pi$ et $\pi_{n+1}^p=\pi_n$ (pour tout $n\geq 0$) et 
on note $\upi$ l'élément associé de $\co_{\oK^\flat}$. On pose  
\begin{equation}\label{definf4b}
\xi_\pi=[\upi]-\pi \in \rW_{\co_K}(\co_{\oK^\flat}). 
\end{equation}
D'après \ref{eipuf6}, la suite 
\begin{equation}\label{definf4c}
\xymatrix{
0\ar[r] &{\rW_{\co_K}(\co_{\oK^\flat})}\ar[r]^-(0.5){\cdot \xi_\pi}& {\rW_{\co_K}(\co_{\oK^\flat})}\ar[r]^-(0.5){\theta_{\co_K}}&
{\co_C}\ar[r]& 0}
\end{equation}
est exacte. En particulier, on a $\xi\in \xi_\pi \rW_{\co_K}(\co_{\oK^\flat})$ \eqref{definf3b}.  On pose 
\begin{equation}\label{definf4d}
\rW_{K}(\co_{\oK^\flat})=\rW(\co_{\oK^\flat})\otimes_WK
\end{equation}
et on note $\theta_K\colon \rW_{K}(\co_{\oK^\flat})\rightarrow C$ l'homomorphisme induit par $\theta$ \eqref{definf3c}.  
On désigne par $\rW^{\ast}_{\co_K}(\co_{\oK^\flat})$ la sous-$\rW_{\co_K}(\co_{\oK^\flat})$-algèbre de $\rW_K(\co_{\oK^\flat})$ engendrée par $[\upi]/\pi$ et on pose 
\begin{equation}\label{definf4e}
\xi^{\ast}_\pi=\frac{\xi_\pi}{\pi}=\frac{[\upi]}{\pi}-1\in \rW^{\ast}_{\co_K}(\co_{\oK^\flat}).
\end{equation}
On prendra garde que $\rW^{\ast}_{\co_K}(\co_{\oK^\flat})$ dépend de la suite $(\pi_n)_{n\geq 0}$. 
L'homomorphisme $\theta_K$ induit un homomorphisme 
\begin{equation}\label{definf4f}
\theta^{\ast}_{\co_K}\colon \rW^{\ast}_{\co_K}(\co_{\oK^\flat})\rightarrow \co_C
\end{equation}
tel que $\theta^{\ast}_{\co_K}(\xi^{\ast}_\pi)=0$. D'après \ref{epinflog6}(ii), la suite 
\begin{equation}\label{definf4g}
\xymatrix{
0\ar[r]&{\rW^{\ast}_{\co_K}(\co_{\oK^\flat})}\ar[r]^-(0.5){\cdot \xi^{\ast}_\pi}&{\rW^{\ast}_{\co_K}(\co_{\oK^\flat})}
\ar[r]^-(0.5){\theta^{\ast}_{\co_K}}&\co_C\ar[r]& 0}
\end{equation}
est exacte. On pose 
\begin{equation}\label{definf4h}
\cA^{\ast}_2(\co_\oK/\co_K)=\rW^{\ast}_{\co_K}(\co_{\oK^\flat})/(\xi^{\ast}_\pi)^2 \rW^{\ast}_{\co_K}(\co_{\oK^\flat})
\end{equation}
et on note $\theta^{\ast}_{\co_K,2}\colon \cA^{\ast}_2(\co_\oK/\co_K)\rightarrow \co_C$ l'homomorphisme induit par $\theta^{\ast}_{\co_K}$. 
La suite \eqref{definf4g} induit une suite exacte 
\begin{equation}\label{definf4i}
\xymatrix{
0\ar[r]&{\co_C}\ar[r]^-(0.5){\cdot \xi^{\ast}_\pi}&{\cA^{\ast}_2(\co_\oK/\co_K)}\ar[r]^-(0.5){\theta^{\ast}_{\co_K,2}}&{\co_C}\ar[r]& 0},
\end{equation}
où on a encore noté $\cdot \xi^{\ast}_\pi$ le morphisme induit par la multiplication par $\xi^{\ast}_\pi$ dans $\cA^{\ast}_2(\co_\oK/\co_K)$. 
L'idéal de carré nul $\ker(\theta^{\ast}_{\co_K,2})$ de $\cA^{\ast}_2(\co_\oK/\co_K)$ 
est un $\co_C$-module libre de base $\xi^{\ast}_\pi$. On le note $\xi^{\ast}_\pi\co_C$ et 
on note $(\xi^{\ast}_\pi)^{-1}\co_C$ son $\co_C$-module dual. 
Pour tout $\co_C$-module $M$, on désigne les $\co_C$-modules $M\otimes_{\co_C}(\xi^{\ast}_\pi \co_C)$ 
et $M\otimes_{\co_C}((\xi^{\ast}_\pi)^{-1} \co_C)$ simplement par $\xi^{\ast}_\pi M$ et $(\xi^{\ast}_\pi)^{-1} M$, respectivement.

\subsection{}\label{definf9}
Pour tout entier $n\geq 1$, on désigne par $\chi_n\colon G_K\rightarrow \mu_{p^n}(\co_\oK)$ l'homomorphisme défini pour tout $g\in G_K$ par
\begin{equation}\label{definf9a}
g(\pi_n)=\chi_n(g)\pi_n,
\end{equation}
et par 
\begin{equation}\label{definf9b}
\chi\colon G_K\rightarrow \mZ_p(1)
\end{equation}
la limite projective des $\chi_n$ \eqref{definf1a}.
On note encore $\chi\colon G_K\rightarrow \co_{\oK^\flat}^\times$ le composé de \eqref{definf9b} et de l'homomorphisme canonique 
$\mZ_p(1)\rightarrow \co_{\oK^\flat}^\times$ \eqref{definf17a}. 
Pour tout $g\in G_K$, on a alors $g([\upi])=[\chi(g)].[\upi]$. Par suite, l'action naturelle de $G_K$ sur $\rW_K(\co_{\oK^\flat})$ préserve $\rW^{\ast}_{\co_K}(\co_{\oK^\flat})$, 
et l'homomorphisme $\theta^{\ast}_{\co_K}$ est $G_K$-équivariant. On en déduit une action de $G_K$ 
sur $\cA^{\ast}_2(\co_{\oK}/\co_K)$ par des automorphismes d'anneaux tel que l'homomorphisme $\theta^{\ast}_{\co_K,2}$ soit $G_K$-équivariant.

\subsection{}\label{definf14}
Soient $K'$ une extension finie de $K$ contenue dans $\oK$, $\co_{K'}$ la clôture intégrale de $\co_K$ dans $K'$,
$K'_0$ l'extension non ramifiée maximale de $K$ contenue dans $K'$, $\co_{K'_0}$ la clôture intégrale de $\co_K$ dans $K'_0$, 
$e$ le degré de l'extension $K'/K'_0$. On fixe une uniformisante $\pi'$ de $\co_{K'}$ et une suite $(\pi'_n)_{n\geq 0}$ 
d'éléments de $\co_\oK$ telle que $\pi'_0=\pi'$ et $\pi'^p_{n+1}=\pi'_n$ (pour tout $n\geq 0$), et 
on note $\upi'$ l'élément associé de $\co_{\oK^\flat}$. On pose 
\begin{equation}\label{definf14a}
\xi_{\pi'}=[\upi']-\pi' \in \rW_{\co_{K'}}(\co_{\oK^\flat}). 
\end{equation}
Soient $f\in \co_{K'_0}[T]$ un polynôme d'Eisenstein de degré $e$ tel que $f(\pi')=0$, $f'\in \co_{K'_0}[T]$ sa dérivée, 
$g\in \co_{K'}[T]$ le polynôme défini par la relation $f(T)=(T-\pi')g(T)$. On a alors
\begin{equation}\label{definf14b}
\theta_{\co_{K'}}(g([\upi']))=f'(\pi')\in \co_C.
\end{equation}

\begin{prop}\label{definf15}
Sous les hypothèses de \ref{definf14}, il existe une unité $u$ de $\rW_{\co_{K'}}(\co_{\oK^\flat})$ telle que 
\begin{equation}\label{definf15a}
\xi_\pi= u \cdot g([\upi'])\cdot \xi_{\pi'} \in \rW_{\co_{K'}}(\co_{\oK^\flat}).
\end{equation}
\end{prop}
En effet, si $K'=K'_0$, on a $g=1$ et la proposition est immédiate puisque $\xi_\pi$ et $\xi_{\pi'}$ sont deux générateurs de $\ker(\theta_{\co_{K'}})$ d'après \eqref{definf4c}.  
On peut donc se borner au cas où $K=K'_0$.  On a $\theta_{\co_K}(f([\upi']))=0$ et $f([\upi']) \mod \pi = \upi'^e$. 
Il existe donc une unité $v$ de $\co_{\oK^\flat}$ telle que $f([\upi']) \mod \pi = v \upi$. Par ailleurs, d'après \eqref{eipuf6d}, la suite 
\begin{equation}\label{definf15b}
\xymatrix{
0\ar[r]&{\co_{\oK^\flat}}\ar[r]^{\cdot \upi}&{\co_{\oK^\flat}} \ar[r]&\co_\oK/\pi\co_\oK\ar[r]& 0}
\end{equation}
où la troisième flèche est induite par la projection canonique $\co_{\oK^\flat} \rightarrow \co_{\oK}/p\co_{\oK}$ sur le premier facteur du système projectif \eqref{definf3a}, 
est exacte. Il résulte alors de \ref{eipuf5} que $f([\upi'])$ engendre $\ker(\theta_{\co_K})$. La proposition s'ensuit compte tenu de \eqref{definf4c}. 

\begin{cor}\label{definf16}
Notant $K_0$ le corps des fractions de $W$ \eqref{definf1} et $\fd$ la différente de l'extension $K/K_0$, 
l'homomorphisme canonique $\cA_2(\co_\oK)\rightarrow \cA^{\ast}_2(\co_\oK/\co_K)$ induit un isomorphisme $\co_C$-linéaire
\begin{equation}\label{definf16a}
\xi\co_C\stackrel{\sim}{\rightarrow}\pi\fd\xi^{\ast}_\pi \co_C.
\end{equation}
\end{cor}
Cela résulte de \ref{definf15} et \eqref{definf14b}.

\subsection{}\label{definf5}
Considérons le système projectif de monoïdes multiplicatifs $(\co_\oK)_{n\in \mN}$, 
où les morphismes de transition sont tous égaux à l'élévation à la puissance $p$-ième.
On note $Q_S$ le produit fibré du diagramme d'homomorphismes de monoïdes 
\begin{equation}\label{definf5a}
\xymatrix{
&{\co_K-\{0\}}\ar[d]\\
{\underset{\underset{x\mapsto x^p}{\longleftarrow}}{\lim}\ \co_\oK}\ar[r]&\co_\oK}
\end{equation}
où  la flèche verticale est l'homomorphisme canonique
et la flèche horizontale est la projection sur la première composante ({\em i.e.}, d'indice $0$).
On désigne par $\uptau_S$ l'homomorphisme composé
\begin{equation}\label{definf5b}
\uptau_S \colon Q_S\longrightarrow
\underset{\underset{x\mapsto x^p}{\longleftarrow}}{\lim}\ \co_\oK \longrightarrow \co_{\oK^\flat} \stackrel{[\ ]}{\longrightarrow} 
\rW(\co_{\oK^\flat}),
\end{equation} 
où $[\ ]$ est le représentant multiplicatif et les autres flèches sont les morphismes canoniques. 
Il résulte aussitôt des définitions que le diagramme 
\begin{equation}\label{definf5c}
\xymatrix{
Q_S\ar[r]\ar[d]_{\uptau_S}&{\co_K-\{0\}}\ar[d]\\
{\rW(\co_{\oK^\flat})}\ar[r]^-(0.4)\theta&{\co_C}}
\end{equation}
où les flèches non libellées sont les morphismes canoniques, est commutatif. 
Par ailleurs, le groupe de Galois $G_K$ agit naturellement sur le monoïde $Q_S$,
et l'homomorphisme $\uptau_S$ est $G_K$-équivariant.

\subsection{}\label{definf6}
On pose 
\begin{equation}\label{definf6a}
\cA_2(\oS)=\Spec(\cA_2(\co_\oK))
\end{equation} 
que l'on munit de la structure logarithmique $\cM_{\cA_2(\oS)}$
associée à la structure pré-logarithmique définie par l'homomorphisme $Q_S\rightarrow \cA_2(\co_\oK)$ induit par $\uptau_S$ \eqref{definf5b}. 
En vertu de \ref{epinflog2}, $\cM_{\cA_2(\oS)}$ est la structure logarithmique sur 
$\cA_2(\oS)$ associée à la structure pré-logarithmique définie par l'homomorphisme 
\begin{equation}\label{definf6b}
\mN\rightarrow \cA_2(\co_\oK),\ \ \ 1\mapsto [\upi].
\end{equation}
En effet, notant $\tpi$ l'élément de $Q_S$ défini par ses projections \eqref{definf4}
\begin{equation}\label{definf6c}
(\pi_n)_{n\in \mN}\in \underset{\underset{\mN}{\longleftarrow}}{\lim}\ \co_\oK \ \ \ {\rm et}\ \ \ 
\pi \in \co_K-\{0\},
\end{equation}
on a $\uptau_S(\tpi)=[\upi]\in \rW(\co_{\oK^\flat})$. 
Le schéma logarithmique $(\cA_2(\oS),\cM_{\cA_2(\oS)})$ est donc fin et saturé, 
et $\theta_2$ \eqref{definf3d} induit une immersion fermée exacte \eqref{definf1}
\begin{equation}\label{definf6d}
(\coS,\cM_\coS)\rightarrow (\cA_2(\oS),\cM_{\cA_2(\oS)}).
\end{equation}

\subsection{}\label{definf7}
On pose 
\begin{equation}\label{definf7a}
\cA^{\ast}_2(\oS/S)=\Spec(\cA^{\ast}_2(\co_\oK/\co_K)),
\end{equation} 
que l'on munit de la structure logarithmique $\cM_{\cA^{\ast}_2(\oS/S)}$
associée à la structure pré-logarithmique définie par l'homomorphisme $Q_S\rightarrow \cA^{\ast}_2(\co_\oK/\co_K)$ induit par $\uptau_S$ \eqref{definf5b}.
En vertu de \ref{epinflog4}, $\cM_{\cA^{\ast}_2(\oS/S)}$ est la structure logarithmique sur 
$\cA^{\ast}_2(\oS/S)$ associée à la structure pré-logarithmique définie par l'homomorphisme 
\begin{equation}\label{definf7b}
\mN\rightarrow \cA^{\ast}_2(\co_\oK/\co_K),\ \ \ 1\mapsto [\upi].
\end{equation}
Le schéma logarithmique $(\cA^{\ast}_2(\oS/S),\cM_{\cA^{\ast}_2(\oS/S)})$ est donc fin et saturé, 
et $\theta^{\ast}_{\co_K,2}$ \eqref{definf4h} induit une immersion fermée exacte \eqref{definf1}
\begin{equation}\label{definf7d}
(\coS,\cM_\coS)\rightarrow (\cA^{\ast}_2(\oS/S),\cM_{\cA^{\ast}_2(\oS/S)}).
\end{equation}

L'élément $(\xi^{\ast}_\pi+1)$ étant inversible dans $\cA^{\ast}_2(\co_\oK/\co_K)$, l'homomorphisme 
\begin{equation}\label{definf7e}
\mN\rightarrow \Gamma(\cA^{\ast}_2(\oS/S), \cM_{\cA^{\ast}_2(\oS/S)}), \ \ \ 1\mapsto (\xi^{\ast}_\pi+1)^{-1}\tpi
\end{equation}
induit un morphisme strict de schémas logarithmiques 
\begin{equation}\label{definf7f}
\pr_1\colon (\cA^{\ast}_2(\oS/S),\cM_{\cA^{\ast}_2(\oS/S)})\rightarrow (S,\cM_S).
\end{equation}
Par ailleurs, on a clairement un morphisme strict de schémas logarithmiques 
\begin{equation}\label{definf7g}
\pr_2\colon (\cA^{\ast}_2(\oS/S),\cM_{\cA^{\ast}_2(\oS/S)})\rightarrow (\cA_2(\oS),\cM_{\cA_2(\oS)}).
\end{equation}

\subsection{}\label{definf10}
Dans la suite de ce chapitre, $(\tS,\cM_\tS)$ désigne l'un des deux schémas logarithmiques 
\begin{equation}\label{definf10a}
(\cA_2(\oS),\cM_{\cA_2(\oS)})\ \ \ {\rm ou} \ \ \ (\cA^{\ast}_2(\oS/S),\cM_{\cA^{\ast}_2(\oS/S)});
\end{equation}
le premier cas sera dit {\em absolu} et le second sera dit {\em relatif}. On note 
\begin{equation}\label{definf10b}
i_S\colon (\coS,\cM_\coS)\rightarrow (\tS,\cM_{\tS})
\end{equation}
l'immersion fermée exacte canonique \eqref{definf6d} ou \eqref{definf7d}. On pose
\begin{equation}\label{definf10c}
\txi=\xi \ \ \ {\rm ou} \ \ \  \txi=\xi^{\ast}_\pi
\end{equation}
selon que l'on est dans le cas absolu ou relatif. Nous traitons simultanément les deux cas, chacun ayant ses avantages et ses inconvénients.
On notera que dans le cas relatif, $\tS$ est naturellement un $S$-schéma \eqref{definf7f}. 

L'immersion fermée $\coS\rightarrow \tS$ est définie par l'idéal de carré nul $\txi\co_\tS$ de $\co_\tS$, 
associé au $\co_C$-module $\txi \co_C$ (cf. \ref{definf3} et \ref{definf4}). On note $\txi^{-1}\co_C$ le $\co_C$-module dual de $\txi \co_C$. 
Pour tout $\co_C$-module $M$ et tout entier $i\geq 1$, on désigne les $\co_C$-modules $M\otimes_{\co_C}(\txi \co_C)^{\otimes i}$ 
et $M\otimes_{\co_C}(\txi^{-1} \co_C)^{\otimes i}$ simplement par $\txi^i M$ et $\txi^{-i} M$, respectivement. 

Le groupe de Galois $G_K$ agit naturellement à gauche sur le schéma logarithmique 
$(\tS,\cM_{\tS})$ et l'immersion fermée \eqref{definf10b} est $G_K$-équivariante.

\subsection{}\label{cad1}
Soit $f\colon (X,\cM_X)\rightarrow (S,\cM_S)$ un morphisme {\em adéquat} de schémas logarithmiques (\cite{agt} III.4.7), 
ayant une carte adéquate (\cite{agt} III.4.4), tel que $X=\Spec(R)$ soit affine et que $X_s$ soit non-vide. 
On désigne par $X^\circ$ le sous-schéma ouvert maximal de $X$
où la structure logarithmique $\cM_X$ est triviale~; c'est un sous-schéma ouvert de $X_\eta$.
On note $j\colon X^\circ\rightarrow X$ l'injection canonique. Pour tout $X$-schéma $U$, on pose  
\begin{equation}\label{definf11a}
U^\circ=U\times_XX^\circ.
\end{equation} 
Pour alléger les notations, on pose
\begin{equation}\label{definf11b}
\tOmega^1_{X/S}=\Omega^1_{(X,\cM_X)/(S,\cM_S)},
\end{equation}
que l'on considère comme un faisceau de $X_\zar$ ou $X_\et$, selon le contexte (cf. \ref{notconv12}), et 
\begin{equation}\label{cad1i}
\tOmega^1_{R/\co_K}=\tOmega^1_{X/S}(X).
\end{equation}

On se donne, de plus, une carte adéquate $((P,\gamma),(\mN,\iota),\vartheta)$ pour $f$, autrement dit une carte 
$(P,\gamma)$ pour  $(X,\cM_X)$ (\cite{agt} II.5.13), une carte $(\mN,\iota)$ pour $(S,\cM_S)$ et un homomorphisme 
de monoïdes $\vartheta\colon \mN\rightarrow P$ tels que les conditions suivantes soient remplies:  
\begin{itemize}
\item[(i)] Le diagramme d'homomorphismes de monoïdes
\begin{equation}\label{cad1a}
\xymatrix{
P\ar[r]^-(0.5)\gamma&{\Gamma(X,\cM_X)}\\
\mN\ar[r]^-(0.5){\iota}\ar[u]^\vartheta&{\Gamma(S,\cM_S)}\ar[u]_{f^\flat}}
\end{equation}
est commutatif, ou ce qui revient au même (avec les notations de \ref{notconv2}),
le diagramme associé de morphismes de schémas logarithmiques 
\begin{equation}\label{cad1b}
\xymatrix{
{(X,\cM_X)}\ar[r]^-(0.5){\gamma^a}\ar[d]_f&{\bA_P}\ar[d]^{\bA_\vartheta}\\
{(S,\cM_S)}\ar[r]^-(0.5){\iota^a}&{\bA_\mN}}
\end{equation}
est commutatif.
\item[(ii)] Le monoïde $P$ est torique, {\em i.e.}, $P$ est fin et saturé et $P^\gp$ est un $\mZ$-module libre (\cite{agt} II.5.1).
\item[(iii)] L'homomorphisme $\vartheta$ est saturé (\cite{agt} II.5.2).
\item[(iv)] L'homomorphisme $\vartheta^\gp\colon \mZ\rightarrow P^\gp$ est injectif, 
le sous-groupe de torsion de $\coker(\vartheta^\gp)$ est d'ordre premier à $p$ et le morphisme de schémas usuels
\begin{equation}\label{cad1c}
X\rightarrow S\times_{\bA_\mN}\bA_P
\end{equation}
déduit de \eqref{cad1b} est étale.  
\item[(v)] Posons $\lambda=\vartheta(1)\in P$, 
\begin{eqnarray}
L&=&\Hom_{\mZ}(P^\gp,\mZ),\label{cad1d}\\
\rH(P)&=&\Hom(P,\mN).\label{cad1e}
\end{eqnarray} 
On notera que $\rH(P)$ est un monoïde fin, saturé et affûté et que l'homomorphisme canonique 
$\rH(P)^\gp\rightarrow \Hom((P^\sharp)^\gp,\mZ)$ est un isomorphisme (\cite{ogus} I 2.2.3). 
On suppose qu'il existe $h_1,\dots,h_r\in \rH(P)$, qui sont $\mZ$-linéairement indépendants dans $L$, tels que  
\begin{equation}\label{cad1f}
\ker(\lambda)\cap \rH(P)=\{\sum_{i=1}^ra_ih_i | \ (a_1,\dots,a_r)\in \mN^r\},
\end{equation}
où l'on considère $\lambda$ comme un homomorphisme $L\rightarrow \mZ$. 
\end{itemize}

\vspace{2mm}

On note $\alpha\colon P\rightarrow R$ l'homomorphisme induit par la carte $(P,\gamma)$. 
On pose $\pi=\iota(1)$ qui est une uniformisante de $\co_K$. On a alors \eqref{notconv2}
\begin{equation}\label{cad1h}
S\times_{\bA_\mN}\bA_P=\Spec(\co_K[P]/(\pi-e^\lambda)).
\end{equation}
On désigne par $P^\gp$ le groupe associé à $P$. 
En vertu de (\cite{ag} 4.2.2), le $\mZ$-module $P^\gp/\mZ\lambda$ est libre de type fini. 
D'après \ref{cad1}(iv) et (\cite{kato1} 1.8) ou (\cite{ogus} IV 1.1.4), on a un isomorphisme $R$-linéaire canonique
\begin{equation}\label{cad1g}
(P^\gp/\mZ \lambda)\otimes_\mZ R \stackrel{\sim}{\rightarrow} \tOmega^1_{R/\co_K}. 
\end{equation}
Ce $R$-module est donc libre de rang $d=\dim(X/S)$. 

\begin{rema}\label{cad100}
Les hypothèses de \ref{cad1} correspondent aux conditions fixées dans (\cite{agt} II.6.2), à l'exception de la connexité de $X$. 
On se ramène au cas où cette dernière condition est satisfaite en remplaçant $X$ par ses composantes connexes.
\end{rema}

\subsection{}\label{definf12}
On munit  $\coX=X\times_S\coS$ \eqref{definf1c} de la structure logarithmique $\cM_\coX$ image inverse de $\cM_X$. 
On a alors un isomorphisme canonique
\begin{equation}\label{definf12a}
(\coX,\cM_{\coX})\stackrel{\sim}{\rightarrow}(X,\cM_X)\times_{(S,\cM_S)}(\coS,\cM_\coS),
\end{equation} 
le produit fibré étant indifféremment pris dans la catégorie des schémas logarithmiques ou 
dans celle des schémas logarithmiques fins. 

{\em On se donne dans la suite de ce chapitre une $(\tS,\cM_{\tS})$-déformation lisse $(\tX,\cM_\tX)$ de $(\coX,\cM_{\coX})$}, 
autrement dit, un morphisme lisse de schémas logarithmiques fins $(\tX,\cM_\tX)\rightarrow (\tS, \cM_{\tS})$
et un $(\coS,\cM_{\coS})$-isomorphisme 
\begin{equation}\label{definf12b}
(\coX,\cM_{\coX})\stackrel{\sim}{\rightarrow}
(\tX,\cM_\tX)\times_{(\tS, \cM_{\tS})}(\coS,\cM_{\coS}).
\end{equation}
Une telle déformation existe et est unique à isomorphisme près en vertu de (\cite{kato1}, 3.14).

\section{Torseurs et algèbres de Higgs-Tate}\label{taht}

\subsection{}\label{cad3}
Pour tout entier $n\geq 1$, on pose
\begin{eqnarray}\label{cad3a}
\co_{K_n}=\co_K[\zeta]/(\zeta^{n}-\pi),
\end{eqnarray}
qui est un anneau de valuation discrète. On note $K_n$ le corps des fractions de $\co_{K_n}$
et $\pi_n$ la classe de $\zeta$ dans $\co_{K_n}$, qui est une uniformisante de $\co_{K_n}$.  
On pose $S^{(n)}=\Spec(\co_{K_n})$
que l'on munit de la structure logarithmique $\cM_{S^{(n)}}$ définie par son point fermé. 
On désigne par $\iota_n\colon \mN\rightarrow \Gamma(S^{(n)},\cM_{S^{(n)}})$
l'homomorphisme défini par $\iota_n(1)=\pi_n$; c'est une carte pour $(S^{(n)},\cM_{S^{(n)}})$.  

Considérons le système inductif de monoïdes $(\mN^{(n)})_{n\geq 1}$, 
indexé par l'ensemble $\mZ_{\geq 1}$ ordonné par la relation de divisibilité, 
défini par $\mN^{(n)}=\mN$ pour tout $n\geq 1$ et dont l'homomorphisme de transition
$\mN^{(n)}\rightarrow \mN^{(mn)}$ (pour $m,n\geq 1$) est l'endomorphisme de Frobenius d'ordre $m$ 
de $\mN$ ({\em i.e.}, l'élévation à la puissance $m$-ième). On notera $\mN^{(1)}$ simplement $\mN$. Les  schémas logarithmiques
$(S^{(n)},\cM_{S^{(n)}})_{n\geq 1}$ forment naturellement un système projectif.
Pour tous entiers $m, n\geq 1$, avec les notations de  \ref{notconv2}, on a un diagramme cartésien de morphismes de schémas logarithmiques
\begin{equation}\label{cad3b}
\xymatrix{
{(S^{(mn)},\cM_{S^{(mn)}})}\ar[r]^-(0.5){\iota^a_{mn}}\ar[d]&{\bA_{\mN^{(mn)}}}\ar[d]\\
{(S^{(n)},\cM_{S^{(n)}})}\ar[r]^-(0.5){\iota^a_n}&{\bA_{\mN^{(n)}}}}
\end{equation}
où $\iota^a_n$ (resp. $\iota^a_{mn}$) est le morphisme associé à $\iota_n$ (resp. $\iota_{mn}$) (\cite{agt} II.5.13). 

\subsection{}\label{cad4}
Considérons le système inductif de monoïdes $(P^{(n)})_{n\geq 1}$, 
indexé par l'ensemble $\mZ_{\geq 1}$ ordonné par la relation de divisibilité, 
défini par $P^{(n)}=P$ pour tout $n\geq 1$ et dont l'homomorphisme de transition
$i_{n,mn}\colon P^{(n)}\rightarrow P^{(mn)}$ (pour $m, n\geq 1$) est l'endomorphisme 
de Frobenius d'ordre $m$ de $P$ ({\em i.e.}, l'élévation à la puissance $m$-ième) (cf. \ref{cad1}). Pour tout $n\geq 1$, on note
\begin{equation}\label{cad4a}
P\stackrel{\sim}{\rightarrow}P^{(n)}, \ \ \ t\mapsto t^{(n)},
\end{equation} 
l'isomorphisme canonique. Pour tout $t\in P$ et tous $m, n\geq 1$, on a donc
\begin{equation}\label{cad4b}
i_{n,mn}(t^{(n)})=(t^{(mn)})^m. 
\end{equation}
On notera $P^{(1)}$ simplement $P$.

Pour tout entier $n\geq 1$, on pose (avec les notations de \ref{notconv2})
\begin{equation}\label{cad4c}
(X^{(n)},\cM_{X^{(n)}})=(X,\cM_{X})\times_{\bA_P}\bA_{P^{(n)}}.
\end{equation}
On notera que la projection canonique $(X^{(n)},\cM_{X^{(n)}})\rightarrow \bA_{P^{(n)}}$ est stricte. 
Comme le diagramme \eqref{cad3b} est cartésien, il existe un unique morphisme
\begin{equation}\label{cad4d}
f^{(n)}\colon (X^{(n)},\cM_{X^{(n)}})\rightarrow (S^{(n)},\cM_{S^{(n)}}),
\end{equation} 
qui s'insère dans le diagramme commutatif
\begin{equation}\label{cad4e}
\xymatrix{
{(X^{(n)},\cM_{X^{(n)}})}\ar[rrr]\ar[ddd]\ar[rd]_{f^{(n)}}\ar@{}[drrr]|*+[o][F-]{1}&&&
{\bA_{P^{(n)}}}\ar[ld]^{\bA_{\vartheta}}\ar[ddd]\\
&{(S^{(n)},\cM_{S^{(n)}})}\ar[r]^-(0.5){\iota^a_n}\ar[d]&{\bA_{\mN^{(n)}}}\ar[d]&\\
&{(S,\cM_{S})}\ar[r]^-(0.5){\iota^a}&{\bA_\mN}&\\
{(X,\cM_{X})}\ar[rrr]\ar[ru]^f&&&{\bA_P}\ar[lu]_{\bA_\vartheta}}
\end{equation}
D'après (\cite{ag} 4.2.7(i)), la face $\xymatrix{\ar@{}|*+[o][F-]{1}}$ du diagramme \eqref{cad4e} est une carte adéquate pour $f^{(n)}$; 
en particulier, $f^{(n)}$ est lisse et saturé.

\subsection{}\label{cad5}
Soit $\oy$ un point géométrique de $\oX^\circ$ (cf. \eqref{definf1c} et \eqref{definf11a}). 
Le schéma $\oX$ étant localement irréductible d'après (\cite{ag} 4.2.7(iii)),  
il est la somme des schémas induits sur ses composantes irréductibles. On note $\oX^\star$
la composante irréductible de $\oX$ contenant $\oy$. 
De même, $\oX^\circ$ est la somme des schémas induits sur ses composantes irréductibles
et $\oX^{\star \circ}=\oX^\star\times_{X}X^\circ$ est la composante irréductible de $\oX^\circ$ contenant $\oy$. 
On désigne par $\Delta$ le groupe profini $\pi_1(\oX^{\star \circ},\oy)$ et par $(V_i)_{i\in I}$ le revêtement universel normalisé de 
$\oX^{\star \circ}$ en $\oy$ (\cite{ag} 2.1.20). Pour chaque $i\in I$, on note $\oX_i$ la fermeture intégrale de $\oX$ dans $V_i$.
Les schémas $(\oX_i)_{i\in I}$ forment alors un système projectif filtrant. On pose   
\begin{equation}\label{taht1d}
\oR=\underset{\underset{i\in I}{\longrightarrow}}{\lim}\ \Gamma(\oX_i,\co_{\oX_i}).
\end{equation} 
C'est un anneau intègre et normal (\cite{ag} 4.1.10), sur lequel agit naturellement $\Delta$ par des homomorphismes d'anneaux et l'action est discrète. 

On pose 
\begin{equation}\label{taht1e}
\mX=\Spec(\oR)\ \ \ {\rm et} \ \ \ \hmX=\Spec(\hoR)
\end{equation}
que l'on munit des structures logarithmiques images inverses de $\cM_X$, 
notées respectivement $\cM_\mX$ et $\cM_\hmX$. 
Les actions de $\Delta$ sur $\oR$ et $\hoR$ induisent des actions à gauche sur 
les schémas logarithmiques $(\mX,\cM_\mX)$ et $(\hmX,\cM_\hmX)$. 
Munissant $(\coX,\cM_\coX)$ de l'action triviale de $\Delta$ (cf. \ref{definf12}), on a un morphisme canonique $\Delta$-équivariant 
\begin{equation}\label{taht1f}
(\hmX,\cM_\hmX)\rightarrow (\coX,\cM_\coX).
\end{equation}

\subsection{}\label{tpcg4}
Pour tous entiers $m,n\geq 1$, le morphisme canonique $X^{(mn)}\rightarrow X^{(n)}$ 
est fini et surjectif. D'après (\cite{ega4} 8.3.8(i)), il existe alors un $X$-morphisme
\begin{equation}\label{tpcg4a}
\oy\rightarrow \underset{\underset{n\geq 1}{\longleftarrow}}{\lim}\ X^{(n)},
\end{equation} 
la limite projective étant indexée par l'ensemble $\mZ_{\geq 1}$ ordonné par la relation de divisibilité.
On fixe un tel morphisme dans toute la suite de cette section. Celui-ci induit un $S$-morphisme
\begin{equation}\label{tpcg4b}
\oS\rightarrow \underset{\underset{n\geq 1}{\longleftarrow}}{\lim}\ S^{(n)}. 
\end{equation}
Pour tout entier $n\geq 1$, on pose
\begin{equation}\label{tpcg4c}
\oX^{(n)}= X^{(n)}\times_{S^{(n)}}\oS \ \ \ {\rm et}\ \ \ \oX^{(n)\circ}=\oX^{(n)}\times_XX^\circ.
\end{equation} 
On en déduit un $\oX$-morphisme 
\begin{equation}\label{tpcg4d}
\oy\rightarrow \underset{\underset{n\geq 1}{\longleftarrow}}{\lim}\ \oX^{(n)}.
\end{equation} 

Pour tout entier $n\geq 1$, le schéma $\oX^{(n)}$ étant normal et localement irréductible d'après (\cite{ag} 4.2.7(iii)), 
il est la somme des schémas induits sur ses composantes irréductibles. 
On note $\oX^{(n)\star}$ la composante irréductible de $\oX^{(n)}$ contenant l'image de $\oy$ \eqref{tpcg4d}.
De même, $\oX^{(n)\circ}$ est la somme des schémas induits sur ses composantes irréductibles
et $\oX^{(n)\star\circ}=\oX^{(n)\star}\times_XX^\circ$  est la composante irréductible de $\oX^{(n)\circ}$ 
contenant l'image de $\oy$. On notera que $\oX^{(n)}$ étant fini sur $\oX$ \eqref{cad4e},  
$\oX^{(n)\star}$ est la fermeture intégrale de $\oX^{\star}$ dans $\oX^{(n)\star\circ}$. 
On pose 
\begin{equation}\label{tpcg4f}
R_n=\Gamma(\oX^{(n)\star},\co_{\oX^{(n)}}).
\end{equation}

D'après (\cite{ag} 4.2.7(v)), le morphisme $\oX^{(n)\star\circ}\rightarrow \oX^{\star\circ}$ 
est étale fini. Il résulte de la preuve de (\cite{agt} II.6.8(iv)) que $\oX^{(n)\star \circ}$ 
est en fait un revêtement étale fini et galoisien de $\oX^{\star \circ}$ de groupe $\Delta_n$ 
canoniquement isomorphe à un sous-groupe de $\Hom_\mZ(P^\gp/\mZ\lambda,\mu_{n}(\oK))$.  
Le groupe $\Delta_n$ agit naturellement sur $R_n$. 

Si $n$ est une puissance de $p$, le morphisme canonique
\begin{equation}
\oX^{(n)\star}\rightarrow\oX^{(n)}\times_\oX\oX^\star
\end{equation}
est un isomorphisme en vertu de (\cite{agt} II.6.6(v)), et on a donc $\Delta_n\simeq \Hom_\mZ(P^\gp/\mZ\lambda,\mu_{n}(\oK))$.

Les anneaux $(R_n)_{n\geq 1}$ forment naturellement un système inductif. On pose 
\begin{eqnarray}
R_\infty&=&\underset{\underset{n\geq 1}{\longrightarrow}}{\lim}\ R_n,\label{tpcg4g}\\
R_{p^\infty}&=&\underset{\underset{n\geq 0}{\longrightarrow}}{\lim}\ R_{p^n}.\label{tpcg4gg}
\end{eqnarray}
Ce sont des anneaux normaux et intègres d'après (\cite{ega1n} 0.6.1.6(i) et 0.6.5.12(ii)).

Le morphisme \eqref{tpcg4d} induit un $\oX^\star$-morphisme 
\begin{equation}\label{tpcg4e}
\Spec(\oR)\rightarrow \underset{\underset{n\geq 1}{\longleftarrow}}{\lim}\ \oX^{(n)\star},
\end{equation} 
et par suite des homomorphismes injectifs 
\begin{equation}\label{tpcg4h}
R_{p^\infty}\rightarrow R_\infty\rightarrow \oR.
\end{equation}
On notera que les $\co_C$-modules $\hRun$ et $\hoR$ sont plats et que l'homomorphisme canonique 
$\hRun\rightarrow \hoR$ est injectif (\cite{agt} II.6.14). 

Les groupes $(\Delta_n)_{n\geq 1}$ forment naturellement un système projectif. On pose 
\begin{eqnarray}
\Delta_\infty&=&\underset{\underset{n\geq 1}{\longleftarrow}}{\lim}\ \Delta_n,\label{tpcg4k}\\
\Delta_{p^\infty}&=&\underset{\underset{n\geq 0}{\longleftarrow}}{\lim}\ \Delta_{p^n}.\label{tpcg4kk}
\end{eqnarray}
On a des homomorphismes canoniques 
\begin{equation}\label{tpcg4l}
\xymatrix{
{\Delta_{\infty}}\ar@{->>}[d]\ar@{^(->}[r]&{\Hom_{\mZ}(P^\gp/\mZ\lambda,\hmZ(1))}\ar[d]\\
{\Delta_{p^\infty}}\ar[r]^-(0.5){\sim}&{\Hom_{\mZ}(P^\gp/\mZ \lambda,\mZ_p(1))}}
\end{equation}
Le noyau $\Sigma_0$ de l'homomorphisme canonique $\Delta_\infty\rightarrow \Delta_{p^\infty}$
est un groupe profini d'ordre premier à $p$. Par ailleurs, le morphisme \eqref{tpcg4d} détermine un homomorphisme 
surjectif $\Delta\rightarrow \Delta_\infty$. On note $\Sigma$ son noyau. Les homomorphismes \eqref{tpcg4h} sont alors $\Delta$-équivariants. 
On résume les constructions introduites dans ce numéro dans le diagramme suivant:
\begin{equation}\label{tpcg4m}
\xymatrix{
R_1\ar[rr]^{\Delta_{p^\infty}}\ar@/^2pc/[rrrr]|{\Delta_\infty}\ar@/_2pc/[rrrrrr]|{\Delta}&&
{R_{p^\infty}}\ar[rr]^{\Sigma_0}&&{R_\infty}\ar[rr]^\Sigma&&\oR}
\end{equation}

\subsection{}\label{cad6}
Considérons le système projectif de monoïdes multiplicatifs $(\oR)_{n\in \mN}$, 
où les morphismes de transition sont les itérés de l'application d'élévation à la puissance $p$-ième de $\oR$.
On note $Q_X$ le produit fibré du diagramme d'homomorphismes de monoïdes 
\begin{equation}\label{cad6a}
\xymatrix{
&{\Gamma(X,\cM_X)}\ar[d]\\
{\underset{\underset{\mN}{\longleftarrow}}{\lim}\ \oR}\ar[r]&\oR}
\end{equation}
où  la flèche verticale est l'homomorphisme canonique (qui se factorise par $R$)
et la flèche horizontale est la projection sur la première composante ({\em i.e.}, d'indice $0$).
On désigne par $\uptau_X$ l'homomorphisme composé
\begin{equation}\label{cad6b}
\uptau_X \colon Q_X\longrightarrow
\underset{\underset{\mN}{\longleftarrow}}{\lim}\ \oR \longrightarrow \oR^\flat \stackrel{[\ ]}{\longrightarrow} 
\rW(\oR^\flat),
\end{equation} 
où $\oR^\flat$ est l'anneau défini dans \eqref{eipo3a}, $[\ ]$ est le représentant multiplicatif et les autres flèches sont les morphismes canoniques. 
Il résulte aussitôt des définitions que le diagramme 
\begin{equation}\label{cad6c}
\xymatrix{
Q_X\ar[r]\ar[d]_{\uptau_X}&{\Gamma(X,\cM_X)}\ar[d]\\
{\rW(\oR^\flat)}\ar[r]^-(0.4)\theta&{\hoR}}
\end{equation}
où les flèches non libellées sont les morphismes canoniques, est commutatif. 
Le groupe $\Delta$ agit naturellement sur le monoïde $Q_X$, et l'homomorphisme $\uptau_X$ est $\Delta$-équivariant.  

On a un homomorphisme canonique 
$Q_S\rightarrow Q_X$ \eqref{definf5} qui s'insère dans un diagramme commutatif 
\begin{equation}\label{cad6d}
\xymatrix{
Q_S\ar[r]\ar[d]_{\uptau_S}&Q_X\ar[d]^{\uptau_X}\\
{\rW(\co_{\oK^\flat})}\ar[r]&{\rW(\oR^\flat)}}
\end{equation}

\subsection{}\label{pmh6}
Posons \eqref{cad4}
\begin{equation}\label{mph6a}
P_\infty=\underset{\underset{n\geq 1}{\longrightarrow}}{\lim}\ P^{(n)}.
\end{equation} 
Pour tout $n\geq 1$, on désigne par $\alpha_{n}\colon P^{(n)}\rightarrow R_{n}$ l'homomorphisme induit par 
le morphisme strict canonique $(X_{n},\cM_{X_{n}})\rightarrow \bA_{P^{(n)}}$ \eqref{cad4c}. 
Pour tous entiers $m,n\geq 1$, le diagramme 
\begin{equation}\label{pmh6b}
\xymatrix{
{P^{(n)}}\ar[r]^{\alpha_{n}}\ar[d]_{i_{n,mn}}&R_{n}\ar[d]\\
{P^{(mn)}}\ar[r]^-(0.5){\alpha_{mn}}&R_{mn}}
\end{equation}
est commutatif. Les $\alpha_{n}$ définissent donc par passage à la limite inductive un homomorphisme 
\begin{equation}\label{pmh6c}
\alpha_\infty\colon P_\infty\rightarrow R_\infty.
\end{equation}
On note encore $\alpha_\infty\colon P_\infty\rightarrow \oR$ le composé de $\alpha_\infty$ et de l'injection canonique 
$R_\infty\rightarrow \oR$.  

Pour tout $t\in P$, on désigne par $\tlt$ l'élément  de  $Q_X$ \eqref{cad6a} défini par ses projections
\begin{equation}\label{pmh6d}
(\alpha_\infty(t^{(p^n)}))_{n\in \mN}\in \underset{\underset{\mN}{\longleftarrow}}{\lim}\ \oR \ \ \ {\rm et}\ \ \ 
\gamma(t) \in \Gamma(X,\cM_X),
\end{equation}
où $t^{(p^n)}$ est l'image de $t$ dans $P^{(p^n)}$ par l'isomorphisme \eqref{cad4a}. L'application 
\begin{equation}\label{pmh6e}
\upnu\colon P\rightarrow Q_X, \ \ \ t\mapsto \tlt
\end{equation}
ainsi définie est un morphisme de monoïdes.

\subsection{}\label{taht201}
Il existe une unique application 
\begin{equation}\label{taht201a}
\langle \ ,\ \rangle\colon  \Delta_\infty\times P_\infty \rightarrow 
\mu_\infty(\co_{\oK})=\underset{\underset{n\geq 1}{\longrightarrow}}{\lim}\ \mu_n(\co_{\oK}), 
\end{equation}
où la limite inductive est indexée par l'ensemble $\mZ_{\geq 1}$ ordonné par la relation de divisibilité,
telle que pour tout $g\in \Delta_\infty$ et tout $x\in P^{(n)}$ $(n\geq 1)$, 
on ait $\langle g,x\rangle\in \mu_n(\co_{\oK})$ et 
\begin{equation}\label{taht201b}
g(\alpha_\infty(x))=\langle g,x\rangle\cdot  \alpha_\infty(x),
\end{equation}
où $\alpha_\infty$ est l'homomorphisme \eqref{pmh6c}. 
En effet, $R_\infty$ est intègre, et on a $\alpha_\infty(x)^n\in \alpha(P)\subset R$ qui est inversible sur $X^\circ$ \eqref{pmh6}; 
donc $\alpha_\infty(x)\not=0$ et $\alpha_\infty(x)^n$ est invariant par $\Delta_\infty$. 

L'accouplement \eqref{taht201a} est multiplicatif en chacun de ses facteurs. 

Soit $n$ un entier $\geq 1$. Rappelons qu'on a un isomorphisme canonique $P^{(n)}\stackrel{\sim}{\rightarrow} P$ 
\eqref{cad4a} et que $\Delta_\infty$ est canoniquement isomorphe 
à un sous-groupe de $\Hom_{\mZ}(P^\gp/\mZ\lambda,\hmZ(1))$ \eqref{tpcg4l}. 
On a donc un homomorphisme canonique $\Delta_\infty\rightarrow \Hom_{\mZ}(P^\gp/\mZ\lambda,\mu_n(\co_{\oK}))$.
D'après (\cite{agt} II.6.6(vi)), le diagramme 
\begin{equation}\label{taht201d}
\xymatrix{
{\Delta_\infty\times P^{(n)}}\ar[rr]^-(0.5){\langle \ ,\ \rangle}\ar[d]&&{\mu_n(\co_{\oK})}\\
{\Hom_{\mZ}(P^\gp/\mZ\lambda,\mu_n(\co_{\oK}))\times P}\ar[rru]&&}
\end{equation}
où la flèche oblique est l'accouplement canonique, est commutatif.

Soit $t\in P$. On désigne par  
\begin{equation}\label{taht202a}
\chi_t\colon \Delta_\infty\rightarrow \mZ_p(1) 
\end{equation}
l'application qui à tout $g\in \Delta_\infty$ associe l'élément 
\begin{equation}\label{taht202b}
\chi_t(g)=\underset{\underset{n\geq 0}{\longleftarrow}}{\lim}\ 
\langle g,t^{(p^n)}\rangle,
\end{equation}
où $t^{(p^n)}$ est l'image de $t$ dans $P^{(p^n)}$ par l'isomorphisme \eqref{cad4a}
et $\langle g,t^{(p^n)}\rangle\in \mu_{p^n}(\co_\oK)$ est défini dans \eqref{taht201a}. 
Il est clair que $\chi_t$ est un morphisme de groupes.

On a clairement $\chi_0=1$, et pour tous $t,t'\in P$,  
\begin{equation}\label{taht202d}
\chi_{tt'}=\chi_t\cdot\chi_{t'}.
\end{equation}
Par suite, l'application $P\rightarrow \Hom(\Delta_\infty,\mZ_p(1))$ définie par  $t\mapsto \chi_t$
est un morphisme de monoïdes. Elle induit donc un homomorphisme que l'on note encore 
\begin{equation}\label{taht202e}
P^\gp\rightarrow \Hom(\Delta_\infty,\mZ_p(1)), \ \ \ t\mapsto \chi_t.
\end{equation}
Comme $\chi_\lambda=1$, on en déduit un homomorphisme 
\begin{equation}\label{taht202f}
P^\gp/\mZ\lambda\rightarrow \Hom(\Delta_\infty,\mZ_p(1)).
\end{equation}
En vertu de (\cite{agt} (II.7.19.9)), ce dernier induit un isomorphisme 
\begin{equation}\label{taht202g}
(P^\gp/\mZ \lambda)\otimes_\mZ\mZ_p\stackrel{\sim}{\rightarrow} \Hom(\Delta_\infty,\mZ_p(1)).
\end{equation}
On en déduit, compte tenu de  \eqref{cad1g} et (\cite{agt} (II.6.12.2)), un isomorphisme $\hRun$-linéaire 
\begin{equation}\label{taht202h}
\tOmega^1_{R/\co_K}\otimes_R\hRun\stackrel{\sim}{\rightarrow} \Hom(\Delta_\infty,\hRun(1)).
\end{equation}

\subsection{}\label{taht2}
Reprenons les notations de \ref{eipo3} pour $A=\oR$: posons 
\begin{equation}\label{taht2a}
\cA_2(\oR)=\rW(\oR^\flat)/\ker (\theta)^2,
\end{equation}
et notons $\theta_2\colon \cA_2(\oR)\rightarrow \hoR$ l'homomorphisme induit par $\theta$ \eqref{eipo3d}.
D'après \ref{eipuf6} et (\cite{agt} II.9.10), la suite 
\begin{equation}\label{taht2b}
0\longrightarrow \rW(\oR^\flat)\stackrel{\cdot \xi}{\longrightarrow} \rW(\oR^\flat)
\stackrel{\theta}{\longrightarrow} \hoR \longrightarrow 0,
\end{equation}
où $\xi\in \rW(\co_{\oK^\flat})$ est l'élément défini dans \eqref{definf3b}, est exacte. Elle induit donc une suite exacte 
\begin{equation}\label{taht2c}
0\longrightarrow \hoR\stackrel{\cdot \xi}{\longrightarrow} \cA_2(\oR)
\stackrel{\theta_2}{\longrightarrow} \hoR \longrightarrow 0,
\end{equation}
où on a encore noté $\cdot \xi$ le morphisme induit par la multiplication par $\xi$ dans $\cA_2(\oR)$. 

Le groupe $\Delta$ agit naturellement sur $\rW(\oR^\flat)$ par des automorphismes d'anneaux,
et l'homomorphisme $\theta$ est $\Delta$-équivariant. On en déduit une action de $\Delta$ 
sur $\cA_2(\oR)$ par des automorphismes d'anneaux telle que l'homomorphisme $\theta_2$ soit 
$\Delta$-équivariant. 

On pose 
\begin{equation}\label{taht2g}
\cA_2(\mX)=\Spec(\cA_2(\oR))
\end{equation} 
que l'on munit de la structure logarithmique $\cM'_{\cA_2(\mX)}$ associée à la structure pré-logarithmique définie 
par l'homomorphisme $Q_X\rightarrow \cA_2(\oR)$ induit par $\uptau_X$ \eqref{cad6b}. 
L'homomorphisme $\theta_2$ induit alors un morphisme 
\begin{equation}\label{taht2h}
(\hmX,\cM_\hmX)\rightarrow (\cA_2(\mX),\cM'_{\cA_2(\mX)}),
\end{equation}
où $(\hmX,\cM_\hmX)$ est le schéma logarithmique défini dans \ref{cad5}.  

Les actions de $\Delta$ sur $\cA_2(\oR)$ et $Q_X$ induisent une action à gauche sur
le schéma logarithmique $(\cA_2(\mX),\cM'_{\cA_2(\mX)})$. Le morphisme \eqref{taht2h} est $\Delta$-équivariant. 

Sans les hypothèses de \ref{epinflog2}, on ne sait pas si 
le schéma logarithmique $(\cA_2(\mX),\cM'_{\cA_2(\mX)})$ est fin et saturé et 
si \eqref{taht2h} est une immersion fermée exacte. C'est pourquoi nous équipons $\cA_2(\mX)$ 
d'une autre structure logarithmique, à savoir la structure logarithmique $\cM_{\cA_2(\mX)}$ associée à la structure pré-logarithmique définie 
par l'homomorphisme composé
\begin{equation}\label{taht2e}
\tuptau_X\colon P\stackrel{\upnu}{\longrightarrow} Q_X\stackrel{\uptau_X}{\longrightarrow} \rW(\oR^\flat) \longrightarrow  \cA_2(\oR),
\end{equation}
où la première flèche est l'homomorphisme \eqref{pmh6e} et seconde flèche est l'homomorphisme \eqref{cad6b}. 

L'homomorphisme $\upnu\colon P\rightarrow Q_X$ \eqref{pmh6e} induit un morphisme de structures logarithmiques sur $\cA_2(\mX)$ 
\begin{equation}\label{taht2f}
\cM_{\cA_2(\mX)}\rightarrow \cM'_{\cA_2(\mX)}.
\end{equation}
Il est clair que l'homomorphisme composé $\theta_2\circ \tuptau_X\colon P\rightarrow\hoR$ 
est induit par $\alpha$ (cf. \ref{cad1}). Par suite, $\theta_2$ induit une immersion fermée exacte
\begin{equation}\label{taht2i}
(\hmX,\cM_{\hmX})\rightarrow (\cA_2(\mX),\cM_{\cA_2(\mX)}),
\end{equation}
qui se factorise à travers le morphisme \eqref{taht2h}.

\subsection{}\label{taht200}
L'homomorphisme canonique $\mZ_p(1)\rightarrow (\oR^\flat,\times)$ et l'homomorphisme
trivial $\mZ_p(1)\rightarrow \Gamma(X,\cM_X)$  (de valeur $1$) induisent un homomorphisme
\begin{equation}\label{taht200aa}
\mZ_p(1)\rightarrow Q_X.
\end{equation} 
Pour tous $g\in \Delta$ et $t\in P$, on a dans $Q_X$
\begin{equation}\label{taht200ab}
g(\upnu(t))=\chi_t(g)\cdot \upnu(t),
\end{equation}
où $\upnu$ est l'application \eqref{pmh6e} et l'on a (abusivement) noté $\chi_t\colon \Delta\rightarrow \mZ_p(1)$ l'application déduite de \eqref{taht202a}.
On en déduit la relation suivante dans $\cA_2(\oR)$
\begin{equation}\label{taht200a}
g(\tuptau_X(t))=[\chi_t(g)]\cdot \tuptau_X(t),
\end{equation}
où $[\chi_t(g)]$ désigne l'image de $\chi_t(g)$ par l'application composée 
\[
\mZ_p(1)\longrightarrow \oR^\flat\stackrel{[\ ]}{\longrightarrow}\rW(\oR^\flat)\longrightarrow \cA_2(\oR).
\]

Pour tout $g\in \Delta$, notons $\gamma_g$ l'automorphisme de $\cA_2(\mX)$ induit par l'action de $g$ sur $\cA_2(\oR)$. 
La structure logarithmique $\gamma_g^*(\cM_{\cA_2(\mX)})$ sur $\cA_2(\mX)$ 
est associée à la structure pré-logarithmique définie par l'homomorphisme
composé $g\circ \tuptau_X\colon P\rightarrow \cA_2(\oR)$ \eqref{taht2e}. L'application
\begin{equation}\label{taht200b}
P\rightarrow \Gamma(\cA_2(\mX),\cM_{\cA_2(\mX)}), \ \ \ t\mapsto [\chi_t(g)]\cdot t,
\end{equation}
est un morphisme de monoïdes \eqref{taht202d}. 
Elle induit donc un morphisme de structures logarithmiques sur $\cA_2(\mX)$
\begin{equation}\label{taht200c}
a_g\colon \gamma_g^*(\cM_{\cA_2(\mX)})\rightarrow \cM_{\cA_2(\mX)}.
\end{equation}
De même, le morphisme de monoïdes  
\begin{equation}\label{taht200d}
P\rightarrow \Gamma(\cA_2(\mX),\gamma_g^*(\cM_{\cA_2(\mX)})), \ \ \ t\mapsto [\chi_t(g^{-1})] \cdot t,
\end{equation}
induit un morphisme de structures logarithmiques sur $\cA_2(\mX)$
\begin{equation}\label{taht200e}
b_g\colon \cM_{\cA_2(\mX)}\rightarrow \gamma_g^*(\cM_{\cA_2(\mX)}).
\end{equation}
On voit aussitôt que $a_g$ et $b_g$ sont des isomorphismes inverses l'un de l'autre, et que 
l'application $g\mapsto  (\gamma_{g^{-1}},a_{g^{-1}})$ est une action à gauche de $\Delta$ sur le schéma logarithmique
$(\cA_2(\mX),\cM_{\cA_2(\mX)})$. 

On vérifie aussitôt que l'immersion \eqref{taht2i} et le morphisme canonique 
\begin{equation}\label{taht200ac}
(\cA_2(\mX),\cM_{\cA_2(\mX)})\rightarrow (\cA_2(\oS),\cM_{\cA_2(\oS)})
\end{equation}
sont $\Delta$-équivariants. Par ailleurs, pour tout $g\in \Delta$, le diagramme 
\begin{equation}\label{taht200ad}
\xymatrix{
{\gamma_g^*(\cM_{\cA_2(\mX)})}\ar[r]\ar[d]_{a_g}&{\gamma_g^*(\cM'_{\cA_2(\mX)})}\ar[d]^{a'_g}\\
{\cM_{\cA_2(\mX)}}\ar[r]&{\cM'_{\cA_2(\mX)}}}
\end{equation}
où les flèches horizontales sont induites par l'homomorphisme \eqref{taht2f} et 
$a'_g$ est l'automorphisme de structures logarithmiques sur $\cA_2(\mX)$ induit par l'action de $g$ sur $Q_X$, 
est commutatif.

\subsection{}\label{taht3}
L'anneau $\rW(\oR^\flat)$ étant naturellement muni d'une structure de $W$-algèbre \eqref{eip2}, on pose 
\begin{equation}\label{taht3a}
\rW_{\co_K}(\oR^\flat)=\rW(\oR^\flat)\otimes_W\co_K 
\end{equation}
et on note $\theta_{\co_K}\colon \rW_{\co_K}(\oR^\flat)\rightarrow \hoR$ l'homomorphisme induit par $\theta$ \eqref{taht2b}.
D'après \ref{eipuf6} et (\cite{agt} II.9.10), la suite 
\begin{equation}\label{taht3b}
\xymatrix{
0\ar[r] &{\rW_{\co_K}(\oR^\flat)}\ar[r]^-(0.5){\cdot \xi_\pi}& {\rW_{\co_K}(\oR^\flat)}\ar[r]^-(0.5){\theta_{\co_K}}&
{\hoR}\ar[r]& 0},
\end{equation}
où $\xi_\pi\in \rW_{\co_K}(\co_{\oK^\flat})$ est l'élément défini dans \eqref{definf4b}, est exacte. On pose 
\begin{equation}\label{taht3c}
\rW_{K}(\oR^\flat)=\rW(\oR^\flat)\otimes_WK
\end{equation}
et on note $\theta_K\colon \rW_{K}(\oR^\flat)\rightarrow \hoR[\frac 1 p]$ l'homomorphisme induit par $\theta$ \eqref{taht2b}.  
On désigne par $\rW^{\ast}_{\co_K}(\oR^\flat)$ la sous-$\rW_{\co_K}(\oR^\flat)$-algèbre de $\rW_K(\oR^\flat)$ engendrée par $\xi^{\ast}_\pi=\xi_\pi/\pi$ \eqref{definf4e}.
L'homomorphisme $\theta_K$ induit un homomorphisme 
\begin{equation}\label{taht3d}
\theta^{\ast}_{\co_K}\colon \rW^{\ast}_{\co_K}(\oR^\flat)\rightarrow \hoR
\end{equation}
tel que $\theta^{\ast}_{\co_K}(\xi^{\ast}_\pi)=0$. D'après \ref{epinflog6}(ii), la suite 
\begin{equation}\label{taht3e}
\xymatrix{
0\ar[r]&{\rW^{\ast}_{\co_K}(\oR^\flat)}\ar[r]^-(0.5){\cdot \xi^{\ast}_\pi}&{\rW^{\ast}_{\co_K}(\oR^\flat)}
\ar[r]^-(0.5){\theta^{\ast}_{\co_K}}&\hoR\ar[r]& 0}
\end{equation}
est exacte. On pose 
\begin{equation}\label{taht3f}
\cA^{\ast}_2(\oR/\co_K)=\rW^{\ast}_{\co_K}(\oR^\flat)/(\xi^{\ast}_\pi)^2 \rW^{\ast}_{\co_K}(\oR^\flat)
\end{equation}
et on note $\theta^{\ast}_{\co_K,2}\colon \cA^{\ast}_2(\oR/\co_K)\rightarrow \hoR$ l'homomorphisme induit par $\theta^{\ast}_{\co_K}$. 
La suite \eqref{taht3e} induit une suite exacte 
\begin{equation}\label{taht3g}
\xymatrix{
0\ar[r]&{\hoR}\ar[r]^-(0.5){\cdot \xi^{\ast}_\pi}&{\cA^{\ast}_2(\oR/\co_K)}\ar[r]^-(0.5){\theta^{\ast}_{\co_K,2}}&{\hoR}\ar[r]& 0},
\end{equation}
où on a encore noté $\cdot \xi^{\ast}_\pi$ le morphisme induit par la multiplication par $\xi^{\ast}_\pi$ dans $\cA^{\ast}_2(\oR/\co_K)$. 

Le groupe $\Delta$ agit naturellement sur $\rW^{\ast}_{\co_K}(\oR^\flat)$ par des automorphismes d'anneaux \eqref{definf9},
et l'homomorphisme $\theta^{\ast}_{\co_K}$ est $\Delta$-équivariant. On en déduit une action de $\Delta$ 
sur $\cA^{\ast}_2(\oR/\co_K)$ par des automorphismes d'anneaux telle que l'homomorphisme $\theta^{\ast}_{\co_K,2}$ soit 
$\Delta$-équivariant.  

On pose 
\begin{equation}\label{taht3h}
\cA^{\ast}_2(\mX/S)=\Spec(\cA^{\ast}_2(\oR/\co_K))
\end{equation} 
que l'on munit de la structure logarithmique $\cM'_{\cA^{\ast}_2(\mX/S)}$ associée à la structure pré-logarithmique définie 
par l'homomorphisme $Q_X\rightarrow \cA^{\ast}_2(\oR/\co_K)$ induit par $\uptau_X$ \eqref{cad6b}. 
L'homomorphisme $\theta^{\ast}_{\co_K,2}$ induit un morphisme
\begin{equation}\label{taht3i}
(\hmX,\cM_\hmX)\rightarrow (\cA^{\ast}_2(\mX/S),\cM'_{\cA^{\ast}_2(\mX/S)}),
\end{equation}
où $(\hmX,\cM_\hmX)$ est le schéma logarithmique défini dans \ref{cad5}.  

Les actions de $\Delta$ sur $\cA^{\ast}_2(\oR/\co_K)$ et $Q_X$ induisent une action à gauche sur
le schéma logarithmique $(\cA^{\ast}_2(\mX/S),\cM'_{\cA^{\ast}_2(\mX/S)})$. Le morphisme \eqref{taht3i} est $\Delta$-équivariant.

Sans les hypothèses de \ref{epinflog4}, on ne sait pas si 
le schéma logarithmique $(\cA^{\ast}_2(\mX/S),\cM'_{\cA^{\ast}_2(\mX/S)})$ est fin et saturé et 
si \eqref{taht3i} est une immersion fermée exacte. C'est pourquoi nous équipons $\cA^{\ast}_2(\mX/S)$ 
d'une autre structure logarithmique, à savoir la structure logarithmique $\cM_{\cA^{\ast}_2(\mX/S)}$ associée à la structure pré-logarithmique définie 
par l'homomorphisme composé
\begin{equation}\label{taht3j}
\tuptau^*_X\colon P\stackrel{\upnu}{\longrightarrow} Q_X\stackrel{\uptau_X}{\longrightarrow} \rW(\oR^\flat) \longrightarrow  \cA^{\ast}_2(\oR/\co_K),
\end{equation}
où la première flèche est l'homomorphisme \eqref{pmh6e} et seconde flèche est l'homomorphisme \eqref{cad6b}.

L'homomorphisme $\upnu\colon P\rightarrow Q_X$ \eqref{pmh6e} induit un morphisme de structures logarithmiques sur $\cA^{\ast}_2(\mX/S)$ 
\begin{equation}\label{taht3k}
\cM_{\cA^{\ast}_2(\mX/S)}\rightarrow \cM'_{\cA^{\ast}_2(\mX/S)}.
\end{equation}
Il est clair que l'homomorphisme composé $\theta^{\ast}_{\co_K,2}\circ \tuptau_X\colon P\rightarrow\hoR$ 
est induit par $\alpha$ (cf. \ref{cad1}). Par suite, $\theta^{\ast}_{\co_K,2}$ induit une immersion fermée exacte
\begin{equation}\label{taht3l}
(\hmX,\cM_{\hmX})\rightarrow (\cA^{\ast}_2(\mX/S),\cM_{\cA^{\ast}_2(\mX/S)}),
\end{equation}
qui se factorise à travers le morphisme \eqref{taht3i}. 

On voit aussitôt que les morphismes \eqref{taht200c} induisent une action à gauche de $\Delta$ sur le schéma logarithmique
$(\cA^{\ast}_2(\mX/S),\cM_{\cA^{\ast}_2(\mX/S)})$, compatible avec le morphisme de structures logarithmiques \eqref{taht3k}. 
L'immersion fermée \eqref{taht3l} est $\Delta$-équivariante.

\subsection{}\label{taht4}
Dans la suite de cette section, selon que l'on est dans le cas absolu ou relatif \eqref{definf10}, 
on désigne par $\tmX$ l'un des deux schémas 
\begin{equation}\label{taht4a}
\cA_2(\mX)\ \ \ {\rm ou} \ \ \ \cA^{\ast}_2(\mX/S),
\end{equation}
et par $\cM_{\tmX}$ (resp. $\cM'_{\tmX}$) l'une des deux structures logarithmiques sur $\tmX$, 
\begin{equation}\label{taht4e}
\cM_{\cA_2(\mX)} \ \ \ {\rm ou} \ \ \ \cM_{\cA^{\ast}_2(\mX/S)} \ \ \ ({\rm resp.} \ \ \ \cM'_{\cA_2(\mX)} \ \ \  {\rm ou} \ \ \ \cM'_{\cA^{\ast}_2(\mX/S)}).
\end{equation}
Le schéma logarithmique $(\tmX,\cM_{\tmX})$ ainsi défini est fin et saturé. 
L'homomorphisme $\upnu\colon P\rightarrow Q_X$ \eqref{pmh6e} induit un morphisme de structures logarithmiques sur $\tmX$ 
\begin{equation}\label{taht4g}
\cM_{\tmX}\rightarrow \cM'_{\tmX}.
\end{equation}
On a un morphisme canonique \eqref{taht2h} ou \eqref{taht3i}
\begin{equation}\label{taht4f}
i'_X\colon (\hmX,\cM_\hmX)\rightarrow (\tmX,\cM'_{\tmX}),
\end{equation}
et une immersion fermée exacte canonique \eqref{taht2i} ou \eqref{taht3l}
\begin{equation}\label{taht4b}
i_X\colon (\hmX,\cM_\hmX)\rightarrow (\tmX,\cM_{\tmX}),
\end{equation}
qui se factorise à travers $i'_X$. 
L'immersion fermée $\hmX\rightarrow \tmX$ est définie par l'idéal de carré nul $\txi\co_\tmX$ de $\co_\tmX$, associé au $\hoR$-module $\txi\hoR$ 
(voir \ref{definf10} pour les notations).

D'après \eqref{cad6d}, on a un morphisme canonique
\begin{equation}\label{taht4c}
(\tmX,\cM_{\tmX})\rightarrow (\tS,\cM_\tS)
\end{equation}
qui s'insère dans un diagramme commutatif \eqref{definf10b}
\begin{equation}\label{taht4d}
\xymatrix{
{(\hmX,\cM_\hmX)}\ar[r]^{i_X}\ar[d]&{(\tmX,\cM_{\tmX})}\ar[d]\\
{(\coS,\cM_{\coS})}\ar[r]^{i_S}&{(\tS,\cM_\tS)}}
\end{equation}

Le groupe $\Delta$ agit à gauche sur les schémas logarithmique $(\tmX,\cM_{\tmX})$ et $(\tmX,\cM'_{\tmX})$, et les morphismes 
$i_X$, $i'_X$ et \eqref{taht4c} sont $\Delta$-équivariants. 

\begin{prop}[\cite{agt} II.9.13]\label{pmh7}
Supposons qu'il existe une carte fine et saturée $M\rightarrow \Gamma(X,\cM_X)$ pour $(X,\cM_X)$ induisant un isomorphisme 
\begin{equation}\label{pmha}
M\stackrel{\sim}{\rightarrow} \Gamma(X,\cM_X)/\Gamma(X,\co^\times_X).
\end{equation}
Alors, le morphisme $\cM_{\tmX}\rightarrow \cM'_{\tmX}$ \eqref{taht4g} est un isomorphisme. 
En particulier, l'homomorphisme $P\rightarrow \Gamma(\tmX,\cM'_\tmX)$ induit par \eqref{pmh6e} est une carte pour $(\tmX,\cM'_\tmX)$, 
et le morphisme $i'_X$ \eqref{taht4f} est une immersion fermée exacte.
\end{prop}

En vertu de \ref{epinflog2} et \ref{epinflog4}, le schéma logarithmique $(\tmX,\cM'_{\tmX})$ est fin et saturé et 
le morphisme $i'_X$ \eqref{taht4f} est une immersion fermée exacte.
Pour tout point géométrique $\oz$ de $\hmX$, notant encore  $\oz$ le point géométrique $i'_X(\oz)=i_X(\oz)$
de $\tmX$, l'homomorphisme
\begin{equation}
\cM_{\tmX,\oz}/\co^\times_{\tmX,\oz} \rightarrow \cM'_{\tmX,\ox}/\co^\times_{\tmX,\oz}
\end{equation}
est un isomorphisme. Comme $\cM'_{\tmX,\oz}$ est intègre, 
on en déduit que le morphisme $\cM_{\tmX,\oz}\rightarrow \cM'_{\tmX,\oz}$ est un isomorphisme; d'où la proposition.

\begin{rema}\label{pmh70}
Contrairement à $\cM_\tmX$, la structure logarithmique $\cM'_\tmX$ sur $\tmX$ ne dépend pas de la carte adéquate choisie sur $(X,\cM_X)$.
Mais sans l'hypothèse de \ref{pmh7}, on ne sait pas si le schéma logarithmique $(\tmX,\cM'_{\tmX})$ 
est fin et saturé et si le morphisme $i'_X$ \eqref{taht4f} est une immersion fermée exacte. 
On notera toutefois que la condition de \ref{pmh7} est remplie sur un recouvrement ouvert affine de $X$ d'après (\cite{agt} II.5.17).
\end{rema}

\subsection{}\label{taht5}
On pose \eqref{cad1i}
\begin{equation}\label{taht5a}
\rT=\Hom_{\hoR}(\tOmega^1_{R/\co_K}\otimes_R\hoR,\txi\hoR).
\end{equation} 
On identifie le $\hoR$-module dual à $\txi^{-1}\tOmega^1_{R/\co_K}\otimes_R\hoR$ (cf. \ref{definf10})
et on note $\cG$ la $\hoR$-algèbre symétrique associée \eqref{notconv9}
\begin{equation}\label{taht5b}
\cG=\rS_{\hoR}(\txi^{-1}\tOmega^1_{R/\co_K}\otimes_R\hoR).
\end{equation}
On désigne par $\hmX_\zar$ le topos de Zariski de $\hmX=\Spec(\hoR)$ \eqref{taht1e}, par $\trT$ le $\co_\hmX$-module associé à $\rT$
et par $\bT$ le $\hmX$-fibré vectoriel associé à son dual, autrement dit,  
\begin{equation}\label{taht5c}
\bT=\Spec(\cG).
\end{equation}
 
Soient $U$ un ouvert de Zariski de $\hmX$, $\tU$ l'ouvert correspondant de $\tmX$ (cf. \ref{taht4}). 
On désigne par $\cL(U)$ 
l'ensemble des morphismes représentés par des flèches pointillées qui complètent  le diagramme canonique
\begin{equation}\label{taht5d}
\xymatrix{
{(U,\cM_\hmX|U)}\ar[rr]^-(0.5){i_X\times_{\tmX}\tU}\ar[d]&&{(\tU,\cM_{\tmX}|\tU)}\ar@{.>}[d]\ar@/^2pc/[dd]\\
{(\coX,\cM_{\coX})}\ar[rr]\ar[d]&&{(\tX,\cM_\tX)}\ar[d]\\
{(\coS,\cM_\coS)}\ar[rr]^-(0.5){i_S}&&{(\tS,\cM_{\tS})}}
\end{equation}
de façon à le laisser commutatif \eqref{definf12b}. D'après (\cite{agt} II.5.23),
le foncteur $U\mapsto \cL(U)$ est un $\trT$-torseur de $\hmX_\zar$.  
On l'appelle {\em torseur de Higgs-Tate} associé à $(\tX,\cM_\tX)$.  
On désigne par $\cF$ le $\hoR$-module des fonctions affines sur $\cL$ (\cite{agt} II.4.9). 
Celui-ci s'insère dans une suite exacte canonique 
\begin{equation}\label{taht5e}
0\rightarrow \hoR\rightarrow \cF\rightarrow \txi^{-1}\tOmega^1_{R/\co_K} \otimes_R \hoR\rightarrow 0.
\end{equation} 
D'après (\cite{illusie1} I 4.3.1.7), cette suite induit pour tout entier $n\geq 1$, une suite exacte \eqref{notconv9}
\begin{equation}\label{taht5f}
0\rightarrow \rS^{n-1}_{\hoR}(\cF)\rightarrow \rS^{n}_{\hoR}(\cF)\rightarrow \rS^n_{\hoR}(\txi^{-1}\tOmega^1_{R/\co_K}
\otimes_R\hoR)\rightarrow 0.
\end{equation}
Les $\hoR$-modules $(\rS^{n}_{\hoR}(\cF))_{n\in \mN}$ forment donc un système inductif filtrant, 
dont la limite inductive 
\begin{equation}\label{taht5g}
\cC=\underset{\underset{n\geq 0}{\longrightarrow}}\lim\ \rS^n_{\hoR}(\cF)
\end{equation}
est naturellement munie d'une structure de $\hoR$-algèbre. D'après (\cite{agt} II.4.10), le $\hmX$-schéma 
\begin{equation}\label{taht5h}
\bL=\Spec(\cC)
\end{equation}
est naturellement un $\bT$-fibré principal homogène sur $\hmX$ qui représente canoniquement $\cL$. 
On prendra garde que $\cL$, $\cF$, $\cC$ et $\bL$ dépendent de $(\tX,\cM_\tX)$.

\subsection{}\label{taht6}
On munit $\hmX$ de l'action naturelle à gauche de $\Delta$~; pour tout $g\in \Delta$, 
l'automorphisme de $\hmX$ défini par $g$, que l'on note aussi $g$, est induit par l'automorphisme $g^{-1}$ de $\hoR$. 
On considère $\trT$ comme un $\co_\hmX$-module $\Delta$-équivariant au moyen 
de la donnée de descente correspondant au $\hRun$-module 
$\Hom_{\hRun}(\tOmega^1_{R/\co_K}\otimes_R\hRun,\txi\hRun)$ (cf. \cite{agt} II.4.18). Pour tout $g\in \Delta$, on a donc 
un isomorphisme canonique de $\co_\hmX$-modules
\begin{equation}\label{taht6a}
\tau_g^\trT\colon \trT\stackrel{\sim}{\rightarrow} g^*(\trT).
\end{equation}
Celui-ci induit un isomorphisme de $\hmX$-schémas en groupes
\begin{equation}\label{taht6b}
\tau_g^\bT\colon \bT\stackrel{\sim}{\rightarrow} g^\bullet(\bT),
\end{equation}
où $g^\bullet$ désigne le foncteur de changement de base par l'automorphisme $g$ de $\hmX$.
On obtient ainsi une structure $\Delta$-équivariante sur le $\hmX$-schéma en groupes $\bT$ (cf. \cite{agt} II.4.17)
et par suite une action à gauche de $\Delta$ sur $\bT$ compatible avec son action sur $\hmX$; 
l'automorphisme de $\bT$ défini par un élément $g$ de $\Delta$ est le composé de $\tau^\bT_g$ 
et de la projection canonique $g^\bullet(\bT)\rightarrow \bT$. 
On en déduit une action de $\Delta$ sur $\cG$ \eqref{taht5b} par des automorphismes 
d'anneaux, compatible avec son action sur $\hoR$, que l'on appelle {\em action canonique}.
Cette dernière est concrètement  induite par l'action triviale sur 
$\rS_{\hRun}(\txi^{-1}\tOmega^1_{R/\co_K}\otimes_R\hRun)$.

L'action naturelle à gauche de $\Delta$ sur le schéma logarithmique $(\tmX,\cM_{\tmX})$ (\ref{taht2} ou \ref{taht3})
induit sur le $\trT$-torseur $\cL$ une structure $\Delta$-équivariante (cf. \cite{agt} II.4.18), 
autrement dit, elle induit pour tout $g\in \Delta$, un isomorphisme $\tau_g^\trT$-équivariant
\begin{equation}\label{taht6c}
\tau^{\cL}_g\colon \cL\stackrel{\sim}{\rightarrow} g^*(\cL);
\end{equation}
ces isomorphismes étant soumis à des relations de compatibilité (cf. \cite{agt} II.4.16). 
En effet, pour tout ouvert de Zariski $U$ de $\hmX$, on prend pour 
\begin{equation}\label{taht6d}
\tau^{\cL}_g(U)\colon \cL(U)\stackrel{\sim}{\rightarrow} \cL(g(U))
\end{equation}
l'isomorphisme défini de la façon suivante. 
Soient $\tU$ l'ouvert de $\tmX$ correspondant à $U$, $\mu\in \cL(U)$ que l'on considère comme un morphisme
\begin{equation}\label{taht6e}
\mu\colon (\tU,\cM_{\tmX}|\tU)\rightarrow (\tX,\cM_\tX).
\end{equation}
Comme $i_X$ \eqref{taht4b} et le morphisme \eqref{taht1f} sont $\Delta$-équivariants, le morphisme composé
\begin{equation}\label{taht6f}
(g(\tU),\cM_{\tmX}|g(\tU))\stackrel{g^{-1}}{\longrightarrow} (\tU,\cM_{\tmX}|\tU)
\stackrel{\mu}{\longrightarrow} (\tX,\cM_\tX)
\end{equation}
prolonge le morphisme canonique $(g(U),\cM_\hmX|g(U))\rightarrow (\tX,\cM_\tX)$. 
Il correspond à l'image de $\mu$ par $\tau^{\cL}_g(U)$. On vérifie aussitôt que 
le morphisme $\tau^{\cL}_g$ ainsi défini est un isomorphisme $\tau_g^\trT$-équivariant et que 
ces isomorphismes vérifient les relations de compatibilité requises dans (\cite{agt} II.4.16).

D'après (\cite{agt} II.4.21), les structures $\Delta$-équivariantes sur $\trT$ et $\cL$ induisent une structure
$\Delta$-équivariante sur le $\co_\hmX$-module associé à $\cF$, ou, ce qui revient au même,
une action $\hoR$-semi-linéaire de $\Delta$ sur $\cF$, telle que les morphismes de la suite \eqref{taht5e} soient 
$\Delta$-équivariants. 
On en déduit sur $\bL$ \eqref{taht5h} une structure de $\bT$-fibré principal homogène $\Delta$-équivariant 
sur $\hmX$ (cf. \cite{agt} II.4.20). Pour tout $g\in \Delta$, on a donc un isomorphisme $\tau_g^\bT$-équivariant
\begin{equation}\label{taht6g}
\tau_g^{\bL}\colon \bL\stackrel{\sim}{\rightarrow} g^\bullet(\bL).
\end{equation}
Cette structure détermine une action à gauche de $\Delta$ sur $\bL$ compatible avec son action sur $\hmX$; 
l'automorphisme de $\bL$ défini par un élément $g$ de $\Delta$ est le composé de $\tau^{\bL}_g$ 
et de la projection canonique $g^\bullet(\bL)\rightarrow \bL$. On obtient ainsi une action de $\Delta$ sur 
$\cC$ \eqref{taht5g} par des automorphismes d'anneaux, compatible avec son action sur $\hoR$, 
que l'on appelle {\em action canonique}. Cette dernière est concrètement induite par l'action de $\Delta$ sur $\cF$. 

Pour tout $g\in \Delta$, on désigne par 
\begin{equation}\label{taht6h}
\bL(\hmX)\stackrel{\sim}{\rightarrow}\bL(\hmX), \ \ \ \mu\mapsto {^g\mu} 
\end{equation}
le composé des isomorphismes 
\begin{eqnarray}
\tau_g^{\bL}\colon \bL(\hmX)&\stackrel{\sim}{\rightarrow}& g^\bullet(\bL)(\hmX),\\
g^\bullet(\bL)(\hmX)&\stackrel{\sim}{\rightarrow}&\bL(\hmX), \ \ \ \mu\mapsto \pr\circ \mu\circ g^{-1}, 
\end{eqnarray}
où $g^{-1}$ agit sur $\hmX$ et $\pr\colon g^\bullet(\bL)\rightarrow \bL$ est la projection canonique,
de sorte que le diagramme 
\begin{equation}\label{taht6i}
\xymatrix{
{\bL}\ar[r]^g&{\bL}\\
\hmX\ar[r]^g\ar[u]^{\mu}&\hmX\ar[u]_{{^g\mu}}}
\end{equation}
est commutatif. En particulier, pour tous $\mu\in \cL(\hmX)$ et $\beta\in \cF$, on a 
\begin{equation}\label{taht6j}
(g^{-1}(\beta))(\mu)=g^{-1}(\beta({^g\mu})).
\end{equation}
Par ailleurs, ${^g\mu}$ est défini par le morphisme composé 
\begin{equation}\label{taht6k}
(\tmX,\cM_{\tmX})\stackrel{g^{-1}}{\longrightarrow} (\tmX,\cM_{\tmX})
\stackrel{\mu}{\longrightarrow} (\tX,\cM_\tX).
\end{equation}

\begin{defi}\label{taht7}
La $\hoR$-algèbre $\cC$ \eqref{taht5g}, munie de l'action canonique de $\Delta$ 
\eqref{taht6}, est appelée {\em l'algèbre de Higgs-Tate} associée à $(\tX,\cM_\tX)$. 
La $\hoR$-représentation $\cF$ \eqref{taht5e} de $\Delta$ est appelée 
l'{\em extension de Higgs-Tate} associée à $(\tX,\cM_\tX)$.
\end{defi}

\subsection{}\label{taht10}
Pour tout nombre rationnel $r\geq 0$, on note $\cF^{(r)}$ la $\hoR$-représentation de $\Delta$ 
déduite de $\cF$ \eqref{taht5e} par image inverse par la multiplication par $p^r$ sur 
$\txi^{-1}\tOmega^1_{R/\co_K}\otimes_R\hoR$, de sorte qu'on a une suite exacte scindée de $\hoR$-modules
\begin{equation}\label{taht10a}
0\rightarrow \hoR\longrightarrow \cF^{(r)}\rightarrow 
\txi^{-1}\tOmega^1_{R/\co_K}\otimes_R\hoR\rightarrow 0.
\end{equation}
D'après (\cite{illusie1} I 4.3.1.7), cette suite induit pour tout entier $n\geq 1$, une suite exacte 
\begin{equation}\label{taht10b}
0\rightarrow \rS^{n-1}_{\hoR}(\cF^{(r)})\rightarrow \rS^{n}_{\hoR}(\cF^{(r)})\rightarrow 
\rS^n_{\hoR}(\txi^{-1}\tOmega^1_{R/\co_K}
\otimes_R\hoR)\rightarrow 0.
\end{equation}
Les $\hoR$-modules $(\rS^{n}_{\hoR}(\cF^{(r)}))_{n\in \mN}$ forment donc un système inductif filtrant, 
dont la limite inductive 
\begin{equation}\label{taht10c}
\cC^{(r)}=\underset{\underset{n\geq 0}{\longrightarrow}}\lim\ \rS^n_{\hoR}(\cF^{(r)})
\end{equation}
est naturellement munie d'une structure de $\hoR$-algèbre. 
L'action de $\Delta$ sur $\cF^{(r)}$ induit une action sur $\cC^{(r)}$ 
par des automorphismes d'anneaux, compatible avec son action sur $\hoR$, que l'on appelle {\em action canonique}.
La $\hoR$-algèbre $\cC^{(r)}$ munie de cette action
est appelée {\em l'algèbre de Higgs-Tate d'épaisseur $r$ associée à $(\tX,\cM_\tX)$}. 

On note $\hcC^{(r)}$ (resp. $\hcC$) le séparé complété $p$-adique de $\cC^{(r)}$ (resp. $\cC$) que l'on suppose toujours muni de 
la topologie $p$-adique. On munit $\hcC^{(r)}\otimes_{\mZ_p}\mQ_p$ de la topologie $p$-adique \eqref{notconv15}.
On vérifie aussitôt que $\cC^{(r)}$ et $\hcC^{(r)}$ sont $\co_C$-plats.
Pour tous nombres rationnels $r'\geq r\geq 0$, on a un $\hoR$-homomorphisme canonique injectif et $\Delta$-équivariant
$\alpha^{r,r'}\colon \cC^{(r')}\rightarrow \cC^{(r)}$. L'homomorphisme induit
$\halpha^{r,r'}\colon\hcC^{(r')}\rightarrow \hcC^{(r)}$ est aussi injectif. On pose 
\begin{equation}\label{taht10d}
\hcC^{(r+)}=\underset{\underset{t\in \mQ_{>r}}{\longrightarrow}}{\lim} \hcC^{(t)},
\end{equation}
que l'on identifie à une sous-$\hoR$-algèbre de $\hcC^{(r)}$ par la limite inductive des homomorphismes 
$(\halpha^{r,t})_{t\in \mQ_{>r}}$. Les actions de $\Delta$ sur les anneaux $(\hcC^{(t)})_{t\in \mQ_{>r}}$ 
induisent une action sur $\hcC^{(r+)}$ par des automorphismes d'anneaux, 
compatible avec ses actions sur $\hoR$ et sur $\hcC$. L'algèbre $\hcC^{(0+)}$ a été notée $\cC^\dagger$ dans (\cite{agt} (II.12.1.6)).

\begin{rema}\label{taht18}
Pour une $\hRun$-algèbre $A$, on considère dans la suite de ce chapitre des $A$-modules 
de Higgs à coefficients dans $\txi^{-1}\tOmega^1_{R/\co_K}\otimes_RA$ (cf. \ref{MH1} et \ref{definf10}). 
On dira abusivement qu'ils sont à coefficients dans $\txi^{-1}\tOmega^1_{R/\co_K}$. 
La catégorie de ces modules sera notée $\bMH(A,\txi^{-1}\tOmega^1_{R/\co_K})$. 
\end{rema}

\subsection{}\label{taht11}
On a un isomorphisme $\cG$-linéaire canonique \eqref{taht5b}
\begin{equation}\label{taht11a}
\Omega^1_{\cG/\hoR}\stackrel{\sim}{\rightarrow} \txi^{-1}\tOmega^1_{R/\co_K} \otimes_R\cG.
\end{equation}
Celui-ci induit un isomorphisme 
\begin{equation}\label{taht11b}
\Omega^1_{\cC/\hoR}\stackrel{\sim}{\rightarrow} \txi^{-1}\tOmega^1_{R/\co_K} \otimes_R\cC.
\end{equation}
On désigne par 
\begin{equation}\label{taht11c}
d_{\cC}\colon \cC\rightarrow \txi^{-1}\tOmega^1_{R/\co_K} \otimes_R\cC
\end{equation}
la $\hoR$-dérivation universelle de $\cC$ et par 
\begin{equation}\label{taht11d}
d_{\hcC}\colon \hcC\rightarrow \txi^{-1}\tOmega^1_{R/\co_K}\otimes_R\hcC
\end{equation}
son prolongement aux complétés (on notera que le $R$-module $\tOmega^1_{R/\co_K}$ est libre de type fini). 
Pour tout $x\in \cF$, $d_{\cC}(x)$ est l'image canonique 
de $x$ dans $\txi^{-1}\tOmega^1_{R/\co_K} \otimes_R\hoR$ \eqref{taht5e}.

De même, pour tout nombre rationnel $r\geq 0$, on désigne par  
\begin{equation}\label{taht11e}
d_{\cC^{(r)}}\colon \cC^{(r)}\rightarrow \txi^{-1}\tOmega^1_{R/\co_K}\otimes_R\cC^{(r)}
\end{equation}
la $\hoR$-dérivation universelle de $\cC^{(r)}$ et par 
\begin{equation}\label{taht11f}
d_{\hcC^{(r)}}\colon \hcC^{(r)}\rightarrow \txi^{-1}\tOmega^1_{R/\co_K}\otimes_R\hcC^{(r)}
\end{equation}
son prolongement aux complétés. On voit aussitôt que les dérivations $d_{\cC^{(r)}}$ et $d_{\hcC^{(r)}}$ sont $\Delta$-équivariantes. 
Comme $\txi^{-1}\tOmega^1_{R/\co_K}\otimes_R \hoR=d_{\cC^{(r)}}(\cF^{(r)})\subset d_{\cC^{(r)}}(\cC^{(r)})$,
$d_{\cC^{(r)}}$ et $d_{\hcC^{(r)}}$ sont également des $\hoR$-champs 
de Higgs à coefficients dans $\txi^{-1}\tOmega^1_{R/\co_K}$ d'après \ref{MH8}(i). 
On désigne par $\mK^\bullet(\hcC^{(r)})$
le complexe de Dolbeault de $(\hcC^{(r)},p^rd_{\hcC^{(r)}})$  \eqref{MH2a} et par $\tmK^\bullet(\hcC^{(r)})$ 
le complexe de Dolbeault augmenté 
\begin{equation}\label{taht11g}
\hoR\rightarrow \mK^0(\hcC^{(r)})\rightarrow \mK^1(\hcC^{(r)})\rightarrow \dots
\rightarrow \mK^n(\hcC^{(r)})\rightarrow \dots,
\end{equation}
où $\hoR$ est placé en degré $-1$ et la différentielle $\hoR\rightarrow\hcC^{(r)}$ est l'homomorphisme canonique. 

Pour tous nombres rationnels $r'\geq r\geq 0$, on a 
\begin{equation}\label{taht11h}
p^{r'}(\id \times \alpha^{r,r'}) \circ d_{\cC^{(r')}}=p^rd_{\cC^{(r)}}\circ \alpha^{r,r'}.
\end{equation}
Par suite, $\halpha^{r,r'}$ induit un morphisme de complexes 
\begin{equation}\label{taht11i}
\upiota^{r,r'}\colon \tmK^\bullet(\hcC^{(r')})\rightarrow \tmK^\bullet(\hcC^{(r)}).
\end{equation}

D'après \eqref{taht11h}, les dérivations $(p^td_{\hcC^{(t)}})_{t\in \mQ_{>r}}$ induisent une $\hoR$-dérivation 
\begin{equation}\label{taht11j}
d^{(r)}_{\hcC^{(r+)}}\colon \hcC^{(r+)}\rightarrow \txi^{-1}\tOmega^1_{R/\co_K} \otimes_R\hcC^{(r+)},
\end{equation}
qui n'est autre que la restriction de $p^rd_{\hcC^{(r)}}$ à $\hcC^{(r+)}$. 
C'est également un $\hoR$-champ de Higgs à coefficients dans $\txi^{-1}\tOmega^1_{R/\co_K}$. 
On note $\mK^\bullet(\hcC^{(r+)})$ 
le complexe de Dolbeault de $(\hcC^{(r+)},d_{\hcC^{(r+)}}^{(r)})$. 
Comme $\hoR$ est $\co_C$-plat, pour tout nombre rationnel $r\geq 0$, on a 
\begin{equation}\label{taht11k}
\ker(d_{\hcC^{(r+)}}^{(r)})=\ker(d_{\hcC^{(r)}})=\hoR.
\end{equation}

\subsection{}\label{taht8}
Considérons une deuxième $(\tS,\cM_{\tS})$-déformation lisse $(\tX',\cM_{\tX'})$ de $(\coX,\cM_{\coX})$
et affectons d'un exposant $'$ les objets associés \eqref{taht5}. 
D'après (\cite{kato1} 3.14), il existe un isomorphisme de $(\tS,\cM_{\tS})$-déformations
\begin{equation}\label{taht8a}
h\colon (\tX,\cM_{\tX})\stackrel{\sim}{\rightarrow} (\tX',\cM_{\tX'}).
\end{equation}
L'isomorphisme de $\trT$-torseurs $\cL\stackrel{\sim}{\rightarrow} \cL'$, 
$\psi\mapsto h\circ \psi$ \eqref{taht5d} induit un isomorphisme $\hoR$-linéaire et $\Delta$-équivariant 
\begin{equation}\label{taht8b}
\cF'\stackrel{\sim}{\rightarrow}\cF,
\end{equation}
qui s'insère dans un diagramme commutatif \eqref{taht5e}
\begin{equation}\label{taht8c}
\xymatrix{
0\ar[r]&{\hoR}\ar[r]\ar@{=}[d]&{\cF'}\ar[r]\ar[d]&{\txi^{-1}\tOmega^1_{R/\co_K} \otimes_R \hoR}\ar[r]\ar@{=}[d]&0\\
0\ar[r]&{\hoR}\ar[r]&{\cF}\ar[r]&{\txi^{-1}\tOmega^1_{R/\co_K} \otimes_R \hoR}\ar[r] & 0}
\end{equation} 
On en déduit un $\hoR$-isomorphisme $\Delta$-équivariant
\begin{equation}\label{taht8d}
\cC'\stackrel{\sim}{\rightarrow} \cC.
\end{equation}

\subsection{}\label{taht9}
Considérons le cas absolu, {\em i.e.}, $(\tS,\cM_\tS)=(\cA_2(\oS),\cM_{\cA_2(\oS)})$ \eqref{definf10} et posons \eqref{definf12}
\begin{equation}\label{taht9a}
(\tX',\cM_{\tX'})=(\tX,\cM_{\tX})\times_{(\cA_2(\oS),\cM_{\cA_2(\oS)})}(\cA^{\ast}_2(\oS/S),\cM_{\cA^{\ast}_2(\oS/S)})
\end{equation}
où le changement de base est défini par le morphisme $\pr_2$ \eqref{definf7g}, 
le produit étant indifféremment pris dans la catégorie des schémas logarithmiques ou dans celle des schémas logarithmiques fins. 
On affecte d'un exposant $'$ les objets associés à la $(\cA^{\ast}_2(\oS/S),\cM_{\cA^{\ast}_2(\oS/S)})$-déformation $(\tX',\cM_{\tX'})$ \eqref{taht5}. 
On note $K_0$ le corps des fractions de $W$ \eqref{definf1} et $\fd$ la différente de l'extension $K/K_0$ et on pose $\rho=v(\pi\fd)$.  
On désigne par 
\begin{equation}\label{taht9b}
\iota\colon \xi\co_C \stackrel{\sim}{\rightarrow} \xi^{\ast}_\pi \co_C
\end{equation}
l'isomorphisme $\co_C$-linéaire tel que le composé 
\begin{equation}\label{taht9c}
\xi\co_C \stackrel{\iota}{\longrightarrow} \xi^{\ast}_\pi \co_C \stackrel{\cdot p^\rho}{\longrightarrow} p^\rho \xi^{\ast}_\pi \co_C
\end{equation}
coïncide avec l'isomorphisme \eqref{definf16a}. L'isomorphisme $\iota\otimes_{\co_C}\hoR$ induit un isomorphisme $\hoR$-linéaire 
$u\colon \rT\stackrel{\sim}{\rightarrow} \rT'$ \eqref{taht5a}.

Une chasse au diagramme commutatif 
\begin{equation}\label{taht9d}
\xymatrix{
{(\hmX,\cM_{\hmX})}\ar[r]\ar[d]&{(\cA^{\ast}_2(\mX/S),\cM_{\cA^{\ast}_2(\mX/S)})}\ar[r]\ar@/^3pc/[dd]|\hole&{(\cA_2(\mX),\cM_{\cA_2(\mX)})}\ar@/^2pc/[dd]\\
{(\coX,\cM_{\coX})}\ar[r]\ar[d]&{(\tX',\cM_{\tX'})}\ar[r]\ar[d]\ar@{}[rd]|\Box&{(\tX,\cM_{\tX})}\ar[d]\\
{(\coS,\cM_{\coS})}\ar[r]&{(\cA^{\ast}_2(\oS/S),\cM_{\cA^{\ast}_2(\oS/S)})}\ar[r]&{(\cA_2(\oS),\cM_{\cA_2(\oS)})}}
\end{equation}
induit un morphisme $(p^\rho u)$-équivariant et $\Delta$-équivariant $\cL\rightarrow \cL'$ de $\hmX_\zar$. D'après (\cite{agt} II.4.12), on en déduit un morphisme 
$\hoR$-linéaire et $\Delta$-équivariant
\begin{equation}\label{taht9e}
v\colon \cF'\rightarrow \cF
\end{equation}
qui s'insère dans un diagramme commutatif \eqref{taht5e}
\begin{equation}\label{taht9f}
\xymatrix{
0\ar[r]&{\hoR}\ar[r]\ar@{=}[d]&{\cF'}\ar[r]\ar[d]^v&{(\xi^{\ast}_\pi)^{-1}\tOmega^1_{R/\co_K} \otimes_R \hoR}\ar[r]\ar[d]^{p^\rho u^\vee}&0\\
0\ar[r]&{\hoR}\ar[r]&{\cF}\ar[r]&{\xi^{-1}\tOmega^1_{R/\co_K} \otimes_R \hoR}\ar[r] & 0}
\end{equation}  
où $u^\vee$ est le morphisme dual de $u$. On en déduit un $\hoR$-morphisme $\Delta$-équivariant
\begin{equation}\label{taht9g}
w\colon \cC'\rightarrow \cC.
\end{equation}
Il est clair que $v\otimes_{\mZ_p}\mQ_p$ et $w\otimes_{\mZ_p}\mQ_p$ sont des isomorphismes \eqref{taht9b}.

D'après \eqref{taht9f}, pour tout nombre rationnel $r\geq 0$, le morphisme $v$ \eqref{taht9e} induit un morphisme  
$\hoR$-linéaire et $\Delta$-équivariant 
\begin{equation}\label{taht9i}
v^{(r)}\colon\cF'^{(r)}\rightarrow \cF^{(r)}
\end{equation} 
qui s'insère dans un diagramme commutatif \eqref{taht10a}
\begin{equation}\label{taht9j}
\xymatrix{
0\ar[r]&{\hoR}\ar[r]\ar@{=}[d]&{\cF'^{(r)}}\ar[r]\ar[d]^{v^{(r)}}&{(\xi^{\ast}_\pi)^{-1}\tOmega^1_{R/\co_K} \otimes_R \hoR}\ar[r]\ar[d]^{p^\rho u^\vee}&0\\
0\ar[r]&{\hoR}\ar[r]&{\cF^{(r)}}\ar[r]&{\xi^{-1}\tOmega^1_{R/\co_K} \otimes_R \hoR}\ar[r] & 0}
\end{equation} 
Comme $u^\vee$ est un isomorphisme, on en déduit un isomorphisme canonique 
$\hoR$-linéaire et $\Delta$-équi\-variant
\begin{equation}\label{taht9k}
\cF'^{(r)}\stackrel{\sim}{\rightarrow} \cF^{(r+\rho)}
\end{equation}
qui s'insère dans un diagramme commutatif 
\begin{equation}
\xymatrix{
0\ar[r]&{\hoR}\ar[r]\ar@{=}[d]&{\cF'^{(r)}}\ar[r]\ar[d]&{(\xi^{\ast}_\pi)^{-1}\tOmega^1_{R/\co_K} \otimes_R \hoR}\ar[r]\ar[d]^{u^\vee}&0\\
0\ar[r]&{\hoR}\ar[r]&{\cF^{(r+\rho)}}\ar[r]&{\xi^{-1}\tOmega^1_{R/\co_K} \otimes_R \hoR}\ar[r] & 0}
\end{equation} 
On en déduit un isomorphisme $\Delta$-équivariant de $\hoR$-algèbres 
\begin{equation}\label{taht9m}
\cC'^{(r)}\stackrel{\sim}{\rightarrow} \cC^{(r+\rho)}.
\end{equation} 
On a donc un isomorphisme canonique $\Delta$-équivariant de $\hoR$-algèbres
\begin{equation}\label{taht9l}
\hcC'^{(r+)}\stackrel{\sim}{\rightarrow}\hcC^{(r+\rho+)}.
\end{equation}
Il identifie $\hcC'^{(r+)}$ à une sous-$\hoR$-algèbre de $\hcC^{(r+)}$.

Pour tout nombre rationnel $r\geq 0$, le diagramme 
\begin{equation}\label{taht9n}
\xymatrix{
{\cC'^{(r)}}\ar[rr]^-(0.5){d_{\cC'^{(r)}}}\ar[d]&&{(\xi^{\ast}_\pi)^{-1}\tOmega^1_{R/\co_K} \otimes_R \cC'^{(r)}}\ar[d]^{u^\vee}\\
{\cC^{(r+\rho)}}\ar[rr]^-(0.5){d_{\cC^{(r+\rho)}}}&&{\xi^{-1}\tOmega^1_{R/\co_K} \otimes_R \cC^{(r+\rho)}}}
\end{equation} 
est commutatif \eqref{taht11e}. On en déduit que le diagramme  
\begin{equation}\label{taht9p}
\xymatrix{
{\hcC'^{(r+)}}\ar[rr]^-(0.5){d^{(r)}_{\hcC'^{(r+)}}}\ar[d]&&{(\xi^{\ast}_\pi)^{-1}\tOmega^1_{R/\co_K} \otimes_R \hcC'^{(r+)}}\ar[d]^{p^\rho u^\vee}\\
{\hcC^{(r+)}}\ar[rr]^-(0.5){d^{(r)}_{\hcC^{(r+)}}}&&{\xi^{-1}\tOmega^1_{R/\co_K} \otimes_R \hcC^{(r+)}}}
\end{equation} 
est commutatif \eqref{taht11j}.

\begin{prop}\label{taht12}
Soient $r,r'$ deux nombres rationnels tels que $r'>r>0$. Alors,  
\begin{itemize}
\item[{\rm (i)}] Il existe un nombre rationnel $\alpha\geq 0$ dépendant de $r$ et $r'$ mais pas du morphisme $f$  
vérifiant les conditions de \ref{cad1} ni de la carte adéquate, tel que
\begin{equation}\label{taht12a}
p^\alpha\upiota^{r,r'}\colon \tmK^\bullet(\hcC^{(r')})\rightarrow 
\tmK^\bullet(\hcC^{(r)}),
\end{equation}
où $\upiota^{r,r'}$ est le morphisme \eqref{taht11i}, soit homotope à $0$ par une homotopie $\hoR$-linéaire. 
\item[{\rm (ii)}] Le morphisme canonique
\begin{equation}\label{taht12b}
\upiota^{r,r'}\otimes_{\mZ_p}\mQ_p\colon \tmK^\bullet(\hcC^{(r')})\otimes_{\mZ_p}\mQ_p\rightarrow 
\tmK^\bullet(\hcC^{(r)})\otimes_{\mZ_p}\mQ_p
\end{equation}
est homotope à $0$ par une homotopie continue. 
\item[{\rm (iii)}] Le complexe $\mK^\bullet(\hcC^{(0+)})\otimes_{\mZ_p}\mQ_p$ est une résolution de $\hoR[\frac 1 p]$ \eqref{taht11j}.
\end{itemize}
\end{prop}

En effet, une section $\psi\in \cL(\hmX)$ induit un isomorphisme de $\hoR$-algèbres $\cG\stackrel{\sim}{\rightarrow}\cC$ \eqref{taht5b}.
Il suffit alors de calquer les preuves de (\cite{agt} II.11.2, II.11.3 et II.11.4), qui correspondent au cas absolu \eqref{definf10} (cf. \cite{agt} II.12.3).

\begin{prop}\label{taht13}
Pour tout nombre rationnel $r\geq 0$, les actions de $\Delta$ sur $\cF^{(r)}$, $\cC^{(r)}$ et $\hcC^{(r)}$ sont 
continues pour les topologies $p$-adiques. 
\end{prop}

Le cas absolu \eqref{taht4} a été démontré dans (\cite{agt} II.12.4) et le cas relatif résulte du cas absolu,  
compte tenu de \ref{taht8} et \ref{taht9}, en particulier de \eqref{taht9m}.

\begin{prop}\label{taht14}
Soient $r,r'$ deux nombres rationnels tels que $r'>r>0$. Alors,
\begin{itemize}
\item[{\rm (i)}] Pour tout entier $n\geq 1$, l'homomorphisme canonique 
\begin{equation}\label{taht14a}
R_1/p^nR_1\rightarrow (\cC^{(r)}/p^n\cC^{(r)})^{\Delta}
\end{equation}
est $\alpha$-injectif \eqref{definf2}. Notons $\cH^{(r)}_n$ son conoyau. 
\item[{\rm (ii)}] Il existe un entier $\alpha\geq 0$, dépendant de $r$, $r'$ et $d=\dim(X/S)$, 
mais pas du morphisme $f$ vérifiant les conditions de \ref{cad1} ni de la carte adéquate, tel que pour tout entier $n\geq 1$, 
le morphisme canonique $\cH^{(r')}_n\rightarrow \cH^{(r)}_n$ soit annulé par $p^\alpha$. 
\item[{\rm (iii)}] Il existe un entier $\gamma\geq 0$, dépendant de $r$, $r'$ et $d$, 
mais pas du morphisme $f$ vérifiant les conditions de \ref{cad1} ni de la carte adéquate, tel que pour tous entiers $n,q\geq 1$, 
le morphisme canonique
\begin{equation}\label{taht14b}
\rH^q(\Delta,\cC^{(r')}/p^n\cC^{(r')})\rightarrow \rH^q(\Delta,\cC^{(r)}/p^n\cC^{(r)})
\end{equation}
soit annulé par $p^\gamma$. 
\end{itemize}
\end{prop}

Le cas absolu \eqref{taht4} a été démontré dans (\cite{agt} II.12.7) et le cas relatif résulte du cas absolu,  
compte tenu de \ref{taht8} et \ref{taht9}, en particulier de \eqref{taht9m}.

\begin{prop}\label{taht16}
Soit $r$ un nombre rationnel $>0$. Alors~:
\begin{itemize}
\item[{\rm (i)}] Le morphisme canonique 
\begin{equation}\label{taht16a}
\hRun\otimes_{\mZ_p}\mQ_p\rightarrow (\hcC^{(r)}\otimes_{\mZ_p}\mQ_p)^\Delta
\end{equation}
est un isomorphisme. 
\item[{\rm (ii)}] Pour tout entier $i\geq 1$, on a 
\begin{equation}\label{taht16b}
\underset{\underset{r\in \mQ_{>0}}{\longrightarrow}}{\lim}\ 
\rH^i_\cont(\Delta,\hcC^{(r)}\otimes_{\mZ_p}\mQ_p)=0.
\end{equation}
\end{itemize}
\end{prop}

(i) Le cas absolu \eqref{taht4} a été démontré dans (\cite{agt} II.12.5(i)) et le cas relatif résulte du cas absolu,  
compte tenu de \ref{taht8} et \ref{taht9}, en particulier de \eqref{taht9m}. 

(ii) Pour tout $\co_C$-$\Delta$-module topologique $M$, notons $\rC^\bullet_\cont(\Delta,M)$ le complexe des cochaînes non homogènes continues
de $\Delta$ à valeurs dans $M$ (\cite{agt} II.3.8). Comme $\Delta$ est compact, pour tout nombre rationnel $r>0$, le morphisme canonique 
\begin{equation}\label{taht16c}
\rC^\bullet_\cont(\Delta,\hcC^{(r)})\otimes_{\mZ_p}\mQ_p\rightarrow 
\rC^\bullet_\cont(\Delta,\hcC^{(r)}\otimes_{\mZ_p}\mQ_p)
\end{equation}
est un isomorphisme. Par ailleurs, en vertu de (\cite{agt} (II.3.10.4) et (II.3.10.5)), on a une suite exacte 
\begin{equation}\label{taht16d}
0\rightarrow \rR^1\underset{\underset{n}{\longleftarrow}}{\lim} \ \rH^{i-1}(\Delta,\cC^{(r)}/p^n\cC^{(r)})\rightarrow 
\rH^i_\cont(\Delta,\hcC^{(r)})\rightarrow 
\underset{\underset{n}{\longleftarrow}}{\lim} \ \rH^{i}(\Delta,\cC^{(r)}/p^n\cC^{(r)})\rightarrow 0.
\end{equation}
Il résulte de \ref{taht14}(iii) que 
\begin{equation}\label{taht16e}
\underset{\underset{r\in \mQ_{>0}}{\longrightarrow}}{\lim}\ 
(\underset{\underset{n\in \mN}{\longleftarrow}}{\lim}\ \rH^{i}(\Delta,\cC^{(r)}/p^n\cC^{(r)}))\otimes_{\mZ_p}\mQ_p =0.
\end{equation}
Si $i\geq 2$, on a, de même,
\begin{equation}\label{taht16f}
\underset{\underset{r\in \mQ_{>0}}{\longrightarrow}}{\lim}\ 
(\rR^1 \underset{\underset{n\in \mN}{\longleftarrow}}{\lim}\ \rH^{i-1}(\Delta,\cC^{(r)}/p^n\cC^{(r)}))\otimes_{\mZ_p}\mQ_p =0.
\end{equation}
Pour tout nombre rationnel $r>0$, notons $\cH^{(r)}_n$ le conoyau du morphisme canonique $R_1/p^nR_1\rightarrow (\cC^{(r)}/p^n\cC^{(r)})^{\Delta}$. 
En vertu de (\cite{jannsen} 1.15) et du fait que $\rR^2 \underset{\underset{n\in \mN}{\longleftarrow}}{\lim} =0$, le morphisme
\begin{equation}\label{taht16g}
\rR^1 \underset{\underset{n\in \mN}{\longleftarrow}}{\lim}\ (\cC^{(r)}/p^n\cC^{(r)})^{\Delta} \rightarrow \rR^1 \underset{\underset{n\in \mN}{\longleftarrow}}{\lim}\ \cH^{(r)}_n
\end{equation}
est un isomorphisme. Il résulte alors de \ref{taht14}(ii) que  
\begin{equation}\label{taht16h}
\underset{\underset{r\in \mQ_{>0}}{\longrightarrow}}{\lim}\ 
(\rR^1 \underset{\underset{n\in \mN}{\longleftarrow}}{\lim}\ (\cC^{(r)}/p^n\cC^{(r)})^{\Delta} )\otimes_{\mZ_p}\mQ_p=0.
\end{equation}
La proposition résulte alors de \eqref{taht16d}.

\begin{cor}\label{taht17}
Pour tout nombre rationnel $r> 0$, on a $(\hcC^{(0+)})^\Delta=(\hcC^{(r)})^\Delta=\hRun$.
\end{cor}

Cela résulte de \ref{taht16} en calquant la preuve de (\cite{agt} II.11.8).

\section{Représentations de Dolbeault} \label{repdolb}

\subsection{}\label{repdolb1}
Pour tout nombre rationnel $r\geq 0$, on désigne par $\bMC^r$ la catégorie des $\cC^{(r)}$-modules à $p^r$-connexion 
intégrable relativement à l'extension $\cC^{(r)}/\hoR$ (cf. \ref{MH9}).
Comme $\txi^{-1}\tOmega^1_{R/\co_K}\otimes_R \hoR=d_{\cC^{(r)}}(\cF^{(r)})\subset d_{\cC^{(r)}}(\cC^{(r)})$ \eqref{taht11},
tout objet de $\bMC^r$ définit un $\hoR$-module de Higgs à coefficients dans $\txi^{-1}\tOmega^1_{R/\co_K}$ en vertu de \ref{MH8}(i). 
On a le foncteur
\begin{equation}\label{repdolb1a}
\bMod(\hoR)\rightarrow \bMC^r, \ \ \ M\mapsto (\cC^{(r)}\otimes_{\hoR}M, p^r d_{\cC^{(r)}}\otimes \id_M).
\end{equation}
Compte tenu de \ref{MH8}(ii) et avec les notations de \ref{taht18}, on a le foncteur 
\begin{equation}\label{repdolb1b}
\bMH(\hRun,\txi^{-1}\tOmega^1_{R/\co_K})\rightarrow \bMC^r, \ \ \ 
(N,\theta)\mapsto (\cC^{(r)}\otimes_{\hRun} N,p^r d_{\cC^{(r)}}\otimes \id_N+\id_{\cC^{(r)}}\otimes \theta).
\end{equation} 

Soient $r'$ un nombre rationnel tel que  $r\geq r'\geq 0$, $(M,\nabla)$ un objet de $\bMC^r$. 
Compte tenu de \eqref{taht11h}, il existe un et un unique morphisme $\hoR$-linéaire 
\begin{equation}\label{repdolb1c}
\nabla'\colon \cC^{(r')}\otimes_{\cC^{(r)}}M\rightarrow \txi^{-1}\tOmega^1_{R/\co_K} \otimes_R \cC^{(r')}\otimes_{\cC^{(r)}}M
\end{equation}
tel que pour tous $t\in \cC^{(r')}$ et $x\in M$, on ait
\begin{equation}\label{repdolb1d}
\nabla'(t\otimes_{\cC^{(r)}}x)=p^{r'} d_{\cC^{(r')}}(t)\otimes_{\cC^{(r)}} x+t\otimes_{\cC^{(r)}}\nabla(x).
\end{equation} 
Par suite, $\nabla'$ est une $p^{r'}$-connexion intégrable sur $\cC^{(r')}\otimes_{\cC^{(r)}}M$ relativement à l'extension $\cC^{(r')}/\hoR$. On notera que le morphisme 
canonique $M\rightarrow \cC^{(r')}\otimes_{\cC^{(r)}}M$ est un morphisme de $\hoR$-modules de Higgs à coefficients dans $\txi^{-1}\tOmega^1_{R/\co_K}$. 
On obtient ainsi un foncteur 
\begin{equation}\label{repdolb1e}
\bMC^r\rightarrow \bMC^{r'},\ \ \ (M,\nabla)\mapsto (\cC^{(r')}\otimes_{\cC^{(r)}}M,\nabla').
\end{equation}

\begin{rema}\label{repdolb20}
Les objets de $\bMC^0$ sont les $\cC$-modules à connexion intégrable relativement à l'extension $\cC/\hoR$. 
On définit de même la catégorie des $\cG$-modules à connexion intégrable relativement à l'extension $\cG/\hoR$ \eqref{taht5b} et les foncteurs analogues 
à \eqref{repdolb1a} et \eqref{repdolb1b}. Tout objet de cette catégorie définit un $\hoR$-module de Higgs à coefficients dans 
$\txi^{-1}\tOmega^1_{R/\co_K}$ d'après \eqref{taht11a} et \ref{MH8}(i).
\end{rema}

\subsection{}\label{repdolb2}
Soient $r$ un nombre rationnel $\geq 0$, $\lambda\in \hoR$, $M$ un $\hcC^{(r)}$-module. 
On appelle {\em $\lambda$-connexion $p$-adique sur $M$ relativement à l'extension $\hcC^{(r)}/\hoR$} la donnée d'un morphisme $\hoR$-linéaire 
\begin{equation}\label{repdolb2a}
\nabla\colon M\rightarrow \txi^{-1}\tOmega^1_{R/\co_K} \otimes_R M
\end{equation}
tel que pour tous $t\in \hcC^{(r)}$ et $x\in M$, on ait 
\begin{equation}\label{repdolb2b}
\nabla(tx)=\lambda d_{\hcC^{(r)}}(t)\otimes x+t\nabla(x).
\end{equation} 
On dit que $\nabla$ est {\em intégrable} si elle est un $\hoR$-champ de Higgs à coefficients dans $\txi^{-1}\tOmega^1_{R/\co_K}$ \eqref{MH1}. 
Si $M$ est complet et séparé pour la topologie $p$-adique, on retrouve la notion introduite dans  (\cite{agt} II.2.14) (cf. \cite{agt} II.2.16). 

\subsection{}\label{repdolb3}
Soit $r$ un nombre rationnel $\geq 0$. 
On désigne par $\bMC^r_p$ la catégorie des $\hcC^{(r)}$-modules à $p^r$-connexion 
$p$-adique intégrable relativement à l'extension $\hcC^{(r)}/\hoR$ \eqref{repdolb2}.
On a le foncteur
\begin{equation}\label{repdolb3a}
\bMod(\hoR)\rightarrow \bMC^r_p, \ \ \ M\mapsto (\hcC^{(r)}\otimes_{\hoR}M, p^r d_{\hcC^{(r)}}\otimes \id_M).
\end{equation}
Soient $(M,\nabla)$ un objet de $\bMC^r_p$, $(N,\theta)$ un $\hRun$-module de Higgs à coefficients dans  $\txi^{-1}\tOmega^1_{R/\co_K}$ (cf. \ref{taht18}).
Il existe un et un unique morphisme $\hoR$-linéaire 
\begin{equation}\label{repdolb3b}
\nabla'\colon M\otimes_{\hRun}N\rightarrow  \txi^{-1}\tOmega^1_{R/\co_K}\otimes_RM\otimes_{\hRun}N
\end{equation}
tel que pour tous $x\in M$ et $y\in N$, on ait 
\begin{equation}\label{repdolb3c}
\nabla'(x\otimes_{\hRun} y)=\nabla(x)\otimes_{\hRun} y+ x\otimes_{\hRun} \theta(y).
\end{equation}
C'est une $p^r$-connexion $p$-adique intégrable sur $M\otimes_{\hRun}N$ relativement à l'extension $\hcC^{(r)}/\hoR$  \eqref{MH2d}.
En particulier, on a un foncteur 
\begin{equation}\label{repdolb3d}
\bMH(\hRun,\txi^{-1}\tOmega^1_{R/\co_K})\rightarrow \bMC^r_p, \ \ \ 
(N,\theta)\mapsto (\hcC^{(r)}\otimes_{\hRun} N,p^r d_{\hcC^{(r)}}\otimes \id_N+\id_{\hcC^{(r)}}\otimes \theta).
\end{equation}

Pour tout objet $(M,\nabla)$ de $\bMC^r$ \eqref{repdolb1}, il existe un et un unique morphisme $\hoR$-linéaire 
\begin{equation}\label{repdolb3i}
\hnabla\colon \hcC^{(r)}\otimes_{\cC^{(r)}}M\rightarrow \txi^{-1}\tOmega^1_{R/\co_K} \otimes_R \hcC^{(r)}\otimes_{\cC^{(r)}}M
\end{equation}
tel que pour tous $t\in \hcC^{(r)}$ et $x\in M$, on ait
\begin{equation}\label{repdolb3j}
\hnabla(t\otimes_{\cC^{(r)}}x)=p^r d_{\hcC^{(r)}}(t)\otimes_{\cC^{(r)}} x+t\otimes_{\cC^{(r)}}\nabla(x).
\end{equation} 
C'est une $p^{r}$-connexion $p$-adique intégrable sur $\hcC^{(r)}\otimes_{\cC^{(r)}}M$ relativement à l'extension $\hcC^{(r)}/\hoR$.
On obtient ainsi un foncteur 
\begin{equation}\label{repdolb3k}
\bMC^r\rightarrow \bMC^r_p,\ \ \ (M,\nabla)\mapsto (\hcC^{(r)}\otimes_{\cC^{(r)}}M,\hnabla).
\end{equation}

Soient $r'$ un nombre rationnel tel que  $r\geq r'\geq 0$, $(M,\nabla)$ un objet de $\bMC^r_p$. 
Compte tenu de \eqref{taht11h}, il existe un et un unique morphisme $\hoR$-linéaire 
\begin{equation}\label{repdolb3e}
\nabla'\colon \hcC^{(r')}\otimes_{\hcC^{(r)}}M\rightarrow \txi^{-1}\tOmega^1_{R/\co_K} \otimes_R \hcC^{(r')}\otimes_{\hcC^{(r)}}M
\end{equation}
tel que pour tous $t\in \hcC^{(r')}$ et $x\in M$, on ait
\begin{equation}\label{repdolb3f}
\nabla'(t\otimes_{\hcC^{(r)}}x)=p^{r'} d_{\hcC^{(r')}}(t)\otimes_{\hcC^{(r)}} x+t\otimes_{\hcC^{(r)}}\nabla(x).
\end{equation} 
C'est une $p^{r'}$-connexion $p$-adique intégrable sur $\hcC^{(r')}\otimes_{\hcC^{(r)}}M$ relativement à l'extension $\hcC^{(r')}/\hoR$. On notera que le morphisme 
canonique $M\rightarrow \hcC^{(r')}\otimes_{\hcC^{(r)}}M$ est un morphisme de $\hoR$-modules de Higgs à coefficients dans $\txi^{-1}\tOmega^1_{R/\co_K}$. 
On obtient ainsi un foncteur 
\begin{equation}\label{repdolb3g}
\bMC^r_p\rightarrow \bMC^{r'}_p,\ \ \ (M,\nabla)\mapsto (\hcC^{(r')}\otimes_{\hcC^{(r)}}M,\nabla').
\end{equation}

Supposons $r>0$ et soit $(M,\nabla)$ un objet de $\bMC^{r}_p$. 
On obtient par passage à la limite inductive des morphismes \eqref{repdolb3e} 
sur les nombres rationnels $0< r'\leq r$ tendant vers $0$, un $\hoR$-champ de Higgs 
\begin{equation}\label{repdolb3h}
\nabla^{(0+)}\colon \hcC^{(0+)}\otimes_{\hcC^{(r)}}M \rightarrow \txi^{-1}\tOmega^1_{R/\co_K}\otimes_R \hcC^{(0+)}\otimes_{\hcC^{(r)}}M.
\end{equation}

\begin{lem}\label{repdolb4}
Soient $M$ un $\hoR[\frac 1 p]$-module projectif de type fini, $r$ un nombre rationnel $\geq 0$. 
On note $\cC^{(r)}\otimes_{\hoR}M$ le $\cC^{(r)}$-module à $p^r$-connexion intégrable relativement à l'extension $\cC^{(r)}/\hoR$ associé à $M$ \eqref{repdolb1} et 
$\hcC^{(r)}\otimes_{\hoR}M$ le $\hcC^{(r)}$-module à $p^r$-connexion $p$-adique intégrable relativement à l'extension $\hcC^{(r)}/\hoR$ associé à $M$ \eqref{repdolb3}. 
Alors, $\cC^{(r)}\otimes_{\hoR}M$ est le sous-$\hoR[\frac 1 p]$-module de Higgs nilpotent maximal de $\hcC^{(r)}\otimes_{\hoR}M$ \eqref{MH16}.
\end{lem}

En effet, $\cC^{(r)}\otimes_{\hoR}M$ est nilpotent puisque la restriction de $d_{\cC^{(r)}}$ \eqref{taht11e} à $\cF^{(r)}$ est le morphisme canonique $\cF^{(r)}\rightarrow
\txi^{-1}\tOmega^1_{R/\co_K}\otimes_R\hoR$ \eqref{taht10a}.
Comme $M$ est un facteur direct d'un $\hoR[\frac 1 p]$-module libre de type fini, on peut se borner au cas où $M=\hoR[\frac 1 p]$ (cf. la preuve de \ref{MH16}),
auquel cas l'assertion est immédiate.

\subsection{}\label{repdolb5}
Pour toute $\hoR$-représentation $M$ de $\Delta$ \eqref{notconv16}, on note $\mH(M)$ le $\hRun$-module défini par 
\begin{equation}\label{repdolb5a}
\mH(M)=(M\otimes_{\hoR}\hcC^{(0+)})^\Delta.
\end{equation}
On le munit du $\hRun$-champ de Higgs 
à coefficients dans $\txi^{-1}\tOmega^1_{R/\co_K}$ induit par $d^{(0)}_{\hcC^{(0+)}}$ \eqref{taht11j} (cf. \ref{taht18}).   
On définit ainsi un foncteur 
\begin{equation}\label{repdolb5b}
\mH\colon \bRep_{\hoR}(\Delta) \rightarrow \bMH(\hRun,\txi^{-1}\tOmega^1_{R/\co_K}).
\end{equation}

\subsection{}\label{repdolb6}
Pour tout $\hRun$-module de Higgs $(N,\theta)$ à coefficients dans $\txi^{-1}\tOmega^1_{R/\co_K}$ \eqref{taht18}, 
on note $\mV(N)$ le $\hoR$-module défini par 
\begin{equation}\label{repdolb6a}
\mV(N)=(N\otimes_{\hRun}\hcC^{(0+)})^{\theta_\tot=0},
\end{equation}
où $\theta_\tot=\theta\otimes \id+\id\otimes d^{(0)}_{\hcC^{(0+)}}$ est le $\hRun$-champ de Higgs total sur 
$N\otimes_{\hRun}\hcC^{(0+)}$ \eqref{MH2d}.
On le munit de l'action $\hoR$-semi-linéaire de $\Delta$ induite
par son action naturelle sur $\hcC^{(0+)}$. 
On définit ainsi un foncteur 
\begin{equation}\label{repdolb6b}
\mV\colon \bMH(\hRun,\txi^{-1}\tOmega^1_{R/\co_K})\rightarrow \bRep_{\hoR}(\Delta).
\end{equation}

\begin{remas}\label{repdolb7}\
\begin{itemize}
\item[(i)] Il résulte de \ref{taht8} que les foncteurs $\mH$ et $\mV$ ne dépendent pas du choix 
de la $(\tS,\cM_\tS)$-déformation $(\tX,\cM_\tX)$ \eqref{definf12}, à isomorphisme {\em non-canonique} près.
\item[(ii)] Pour toute $\hoR$-représentation $M$ de $\Delta$, le morphisme canonique 
\begin{equation}\label{repdolb7a}
\mH(M)\otimes_{\hRun}\hRun[\frac 1 p]\rightarrow \mH(M\otimes_{\hoR}\hoR[\frac 1 p])
\end{equation}
est un isomorphisme.   
\item[(iii)] Pour tout $\hRun$-module de Higgs $(N,\theta)$ à coefficients dans $\txi^{-1}\tOmega^1_{R/\co_K}$, 
le morphisme canonique 
\begin{equation}\label{repdolb7b}
\mV(N)\otimes_{\hoR}\hoR[\frac 1 p]\rightarrow \mV(N\otimes_{\hRun}\hRun[\frac 1 p])
\end{equation}
est un isomorphisme.   
\end{itemize}
\end{remas}

\begin{defi}[\cite{agt} II.12.11] \label{repdolb8}
On dit qu'une $\hoR[\frac 1 p]$-représentation $M$ de $\Delta$ est {\em de Dolbeault} 
si les conditions suivantes sont remplies~:
\begin{itemize}
\item[(i)] $\mH(M)$ est un $\hRun[\frac 1 p]$-module projectif de type fini;
\item[(ii)] le morphisme canonique 
\begin{equation}\label{repdolb8a}
\mH(M) \otimes_{\hRun}\hcC^{(0+)}\rightarrow  M\otimes_{\hoR}\hcC^{(0+)}
\end{equation}
est un isomorphisme.
\end{itemize}
\end{defi}

Cette notion ne dépend pas du choix de la $(\tS,\cM_\tS)$-déformation $(\tX,\cM_{\tX})$ d'après \ref{repdolb7}(i).
Mais elle dépend a priori du cas relatif ou absolu considéré dans \ref{definf10}.
Lorsqu'il y a lieu de préciser, on dira que la $\hoR[\frac 1 p]$-représentation $M$ est {\em de Dolbeault absolue} ou 
{\em de Dolbeault relativement à $\co_K$}.

Cette définition est en fait équivalente à (\cite{agt} II.12.11), même si elle semble plus générale, cf. \ref{repdolb9}.

\begin{defi}[\cite{agt} II.12.12] \label{repdolb10}
On dit qu'un $\hRun[\frac 1 p]$-module de Higgs $(N,\theta)$ à coefficients dans $\txi^{-1}\tOmega^1_{R/\co_K}$
est {\em soluble} si les conditions suivantes sont remplies~:
\begin{itemize}
\item[(i)] $N$ est un $\hRun[\frac 1 p]$-module projectif de type fini~;
\item[(ii)] le morphisme canonique 
\begin{equation}\label{repdolb10a}
\mV(N) \otimes_{\hoR}\hcC^{(0+)}\rightarrow  N\otimes_{\hRun}\hcC^{(0+)}
\end{equation}
est un isomorphisme.
\end{itemize}
\end{defi}

Cette notion ne dépend pas du choix de la $(\tS,\cM_\tS)$-déformation $(\tX,\cM_{\tX})$ d'après \ref{repdolb7}(i).
Mais elle dépend a priori du cas relatif ou absolu considéré dans \ref{definf10}.
Les coefficients du module de Higgs le précise. Il n'y a donc pas besoin de le préciser dans la terminologie. 

Cette définition est en fait équivalente à (\cite{agt} II.12.12), même si elle semble plus générale, cf. \ref{repdolb11}. 

\begin{prop}\label{repdolb17}
Soient $M$ une $\hoR[\frac 1 p]$-représentation de $\Delta$, 
$N$ un $\hRun[\frac 1 p]$-module projectif de type fini, muni  d'un $\hRun$-champ de Higgs $\theta$ 
à coefficients dans $\txi^{-1}\tOmega^1_{R/\co_K}$, $r$ un nombre rationnel $>0$,
\begin{equation}\label{repdolb17a}
N \otimes_{\hRun}\hcC^{(r)}\stackrel{\sim}{\rightarrow}  M\otimes_{\hoR}\hcC^{(r)}
\end{equation}
un isomorphisme $\hcC^{(r)}$-linéaire
et $\Delta$-équivariant de $\hoR[\frac 1 p]$-modules de Higgs à coefficients dans $\txi^{-1}\tOmega^1_{R/\co_K}$,
où $N$ est muni de l'action triviale de $\Delta$, $M$ est muni du champ de Higgs nul et 
$\hcC^{(r)}$ est muni de l'action canonique de $\Delta$ et du champ de Higgs $p^rd_{\hcC^{(r)}}$.
Alors, 
\begin{itemize}
\item[{\rm (i)}] Le $\hoR[\frac 1 p]$-module $M$ est projectif de type fini, et l'action de $\Delta$ sur $M$ est continue pour la topologie $p$-adique \eqref{notconv15}. 
\item[{\rm (ii)}] La $\hoR[\frac 1 p]$-représentation $M$ de $\Delta$ est de Dolbeault, et on a un isomorphisme de $\hRun[\frac 1 p]$-modules de Higgs
à coefficients dans $\txi^{-1}\tOmega^1_{R/\co_K}$
\begin{equation}\label{repdolb17b}
\mH(M)\stackrel{\sim}{\rightarrow}  (N,\theta). 
\end{equation}
\item[{\rm (iii)}] Le $\hRun[\frac 1 p]$-module de Higgs $(N,\theta)$ est soluble, et on a un isomorphisme de $\hoR[\frac 1 p]$-représenta\-tions de $\Delta$
\begin{equation}\label{repdolb17c}
\mV(N,\theta)\stackrel{\sim}{\rightarrow}  M.
\end{equation}
\end{itemize}
\end{prop}

Toute section $\cF^{(r)}\rightarrow \hoR$ de l'extension \eqref{taht10a}, en tant que suite exacte de $\hoR$-modules sans actions de $\Delta$, 
définit une rétraction de la $\hoR$-algèbre $\hcC^{(r)}$. 
L'isomorphisme \eqref{repdolb17a} implique alors que le $\hoR[\frac 1 p]$-module $M$ est projectif de type fini. 
L'isomorphisme \eqref{repdolb17a} est en fait un isomorphisme de $\hcC^{(r)}$-modules à $p^r$-connexion $p$-adique intégrable 
relativement à l'extension $\hcC^{(r)}/\hoR$ \eqref{repdolb3}. On en déduit 
un isomorphisme $\hcC^{(0+)}$-linéaire et $\Delta$-équivariant de $\hoR[\frac 1 p]$-modules de Higgs à coefficients dans 
$\txi^{-1}\tOmega^1_{R/\co_K}$ \eqref{repdolb3h}
\begin{equation}\label{repdolb17d}
N \otimes_{\hRun}\hcC^{(0+)}\stackrel{\sim}{\rightarrow}  M\otimes_{\hoR}\hcC^{(0+)}.
\end{equation}
Comme $N$ est un facteur direct d'un $\hRun[\frac 1 p]$-module libre de type fini, on a 
$(N\otimes_{\hRun}\hcC^{(0+)})^\Delta=N$ en vertu de \ref{taht17}. On obtient alors de \eqref{repdolb17d} un isomorphisme
de $\hRun[\frac 1 p]$-modules de Higgs $\mH(M)\stackrel{\sim}{\rightarrow} (N,\theta)$. De même, compte tenu de \eqref{taht11k}, 
on obtient de \eqref{repdolb17d} un isomorphisme de $\hoR[\frac 1 p]$-représentations de $\Delta$, $\mV(N,\theta) \stackrel{\sim}{\rightarrow} M$. 
On déduit que le $\hRun[\frac 1 p]$-module de Higgs $(N,\theta)$ est soluble et que la $\hoR[\frac 1 p]$-représentation $M$ de $\Delta$ est de Dolbeault. 
Par ailleurs, la preuve de (\cite{agt} II.12.22) montre que la $\hoR[\frac 1 p]$-représentation $M$ de $\Delta$ est continue, d'où la proposition.  
Le lecteur attentif notera que (\cite{agt} II.12.22) se place dans le cas absolu \eqref{taht4}, mais la même preuve vaut pour le cas relatif.

\begin{cor}\label{repdolb9}
Pour qu'une $\hoR[\frac 1 p]$-représentation $M$ de $\Delta$ 
soit de Dolbeault, il faut et il suffit qu'il existe un $\hRun[\frac 1 p]$-module projectif de type fini $N$, un $\hRun$-champ de Higgs 
$\theta$ sur $N$ à coefficients dans $\txi^{-1}\tOmega^1_{R/\co_K}$, un nombre rationnel $r>0$ et un isomorphisme $\hcC^{(r)}$-linéaire
et $\Delta$-équivariant de $\hoR[\frac 1 p]$-modules de Higgs à coefficients dans $\txi^{-1}\tOmega^1_{R/\co_K}$
\begin{equation}\label{repdolb9a}
N \otimes_{\hRun}\hcC^{(r)}\stackrel{\sim}{\rightarrow}  M\otimes_{\hoR}\hcC^{(r)},
\end{equation}
où $N$ est muni de l'action triviale de $\Delta$, $M$ est muni du champ de Higgs nul et 
$\hcC^{(r)}$ est muni de l'action canonique de $\Delta$ et du champ de Higgs $p^rd_{\hcC^{(r)}}$.
De plus, dans ce cas, on a les propriétés suivantes:
\begin{itemize}
\item[{\rm (i)}] Le $\hoR[\frac 1 p]$-module $M$ est projectif de type fini, et l'action de $\Delta$ sur $M$ est continue pour la topologie $p$-adique \eqref{notconv15}. 
\item[{\rm (ii)}] On a un isomorphisme de $\hRun[\frac 1 p]$-modules de Higgs
à coefficients dans $\txi^{-1}\tOmega^1_{R/\co_K}$
\begin{equation}\label{repdolb9e}
\mH(M)\stackrel{\sim}{\rightarrow}  (N,\theta).
\end{equation}
\end{itemize}
\end{cor}

En effet, la condition est suffisante en vertu de \ref{repdolb17}. 
Montrons qu'elle est nécessaire. Supposons $M$ de Dolbeault.  
Toute section $\cF\rightarrow \hoR$ de l'extension \eqref{taht5e}, en tant que suite exacte de $\hoR$-modules sans actions de $\Delta$, 
définit une rétraction de la $\hoR$-algèbre $\hcC=\hcC^{(0)}$, et par suite une rétraction de la sous-$\hoR$-algèbre $\hcC^{(0+)}$ \eqref{taht10d}. 
L'isomorphisme \eqref{repdolb8a} implique alors que le $\hoR[\frac 1 p]$-module $M$ est projectif de type fini. 
Par suite, pour tout nombre rationnel $r> 0$, le morphisme canonique 
$M\otimes_\hoR\hcC^{(r)}\rightarrow M\otimes_\hoR\hcC^{(0+)}$ est injectif.
Le $\hRun[\frac 1 p]$-module $\mH(M)$ étant de type fini,
il existe un nombre rationnel $r>0$ tel que $\mH(M)$ soit contenu dans 
$M\otimes_{\hoR}\hcC^{(r)}$. Comme $M$ est de type fini sur $\hoR[\frac 1 p]$, 
quitte à diminuer $r$, le morphisme canonique
\begin{equation}\label{repdolb9b}
\mH(M) \otimes_{\hRun}\hcC^{(r)}\rightarrow  M\otimes_{\hoR}\hcC^{(r)}
\end{equation}
est surjectif. Par ailleurs, $\mH(M)$ étant $\hRun$-plat, le morphisme canonique 
\begin{equation}\label{repdolb9c}
\mH(M)\otimes_\hRun\hcC^{(r)}\rightarrow \mH(M)\otimes_\hRun\hcC^{(0+)}
\end{equation} 
est injectif. On en déduit que \eqref{repdolb9b} est un isomorphisme; d'où la condition recherchée. 
Les propriétés (i) et (ii) résultent de \ref{repdolb17}.

\begin{cor}\label{repdolb11}
Pour qu'un $\hRun[\frac 1 p]$-module de Higgs $(N,\theta)$ à coefficients dans $\txi^{-1}\tOmega^1_{R/\co_K}$
soit soluble, il faut et il suffit qu'il vérifie les conditions suivantes~:
\begin{itemize}
\item[{\rm (i)}] le $\hRun[\frac 1 p]$-module $N$ est projectif de type fini;
\item[{\rm (ii)}] il existe une $\hoR[\frac 1 p]$-représentation $M$ de $\Delta$, 
un nombre rationnel $r>0$ et un isomorphisme $\hcC^{(r)}$-linéaire et $\Delta$-équivariant de $\hoR[\frac 1 p]$-modules de Higgs
à coefficients dans $\txi^{-1}\tOmega^1_{R/\co_K}$
\begin{equation}\label{repdolb11a}
M \otimes_{\hoR}\hcC^{(r)}\stackrel{\sim}{\rightarrow}  N\otimes_{\hRun}\hcC^{(r)},
\end{equation}
où $N$ est muni de l'action triviale de $\Delta$, $M$ est muni du champ de Higgs nul et 
$\hcC^{(r)}$ est muni de l'action canonique de $\Delta$ et du champ de Higgs $p^rd_{\hcC^{(r)}}$. 
\end{itemize}
De plus, dans ce cas, le $\hoR[\frac 1 p]$-module $M$ est projectif de type fini, et on a un isomorphisme de $\hoR[\frac 1 p]$-représentations de $\Delta$
\begin{equation}\label{repdolb11b}
\mV(N,\theta)\stackrel{\sim}{\rightarrow}  M.
\end{equation}
\end{cor}

La preuve est similaire à celle de \ref{repdolb9} et est laissée au lecteur.

\subsection{}\label{repdolb12}
Reprenons les hypothèses et notations de \ref{taht9}, et notons de plus 
\begin{equation}\label{repdolb12c}
\nu\colon (\xi_\pi^{\ast})^{-1}\hRun\otimes_R\tOmega^1_{R/\co_K}\stackrel{\sim}{\rightarrow} \xi^{-1}\hRun\otimes_R\tOmega^1_{R/\co_K}
\end{equation}
l'isomorphisme induit par \eqref{taht9b}. On notera que $\nu\otimes_{\hRun}\hoR=u^\vee$ est le morphisme dual de $u$ \eqref{taht9}. On note 
\begin{equation}
\varepsilon \colon \bMH(\hRun,(\xi_\pi^{\ast})^{-1}\tOmega^1_{R/\co_K})\rightarrow \bMH(\hRun,\xi^{-1}\tOmega^1_{R/\co_K})
\end{equation}
le foncteur induit par $p^\rho \nu$. On désigne par 
\begin{equation}\label{repdolb12a}
\xymatrix{
{\bRep_{\hoR}(\Delta)} \ar@<1ex>[r]^-(0.5){\mH}& {\bMH(\hRun,\xi^{-1}\tOmega^1_{R/\co_K})} \ar@<1ex>[l]^-(0.5){\mV}}
\end{equation}
\begin{equation}\label{repdolb12b}
\xymatrix{
{\bRep_{\hoR}(\Delta)} \ar@<1ex>[r]^-(0.5){\mH'}& {\bMH(\hRun,(\xi_\pi^{\ast})^{-1}\tOmega^1_{R/\co_K})} \ar@<1ex>[l]^-(0.5){\mV'}}
\end{equation}
les foncteurs \eqref{repdolb5b} et \eqref{repdolb6b} associés aux déformations $(\tX,\cM_{\tX})$ et $(\tX',\cM_{\tX'})$, respectivement.

Pour toute $\hoR$-représentation $M$ de $\Delta$, l'homomorphisme canonique $\hcC'^{(0+)}\rightarrow \hcC^{(0+)}$ \eqref{taht9l} 
induit un morphisme $\hRun$-linéaire 
\begin{equation}\label{repdolb12d}
a_M\colon \mH'(M)\rightarrow \mH(M).
\end{equation}
Compte tenu de \eqref{taht9p}, celui-ci s'insère dans un diagramme commutatif
\begin{equation}\label{repdolb12e}
\xymatrix{
{\mH'(M)}\ar[r]^-(0.5){\theta'}\ar[d]_{a_M}&{(\xi_\pi^{\ast})^{-1}\tOmega^1_{R/\co_K}\otimes_R \mH'(M)}\ar[d]^{p^\rho \nu \otimes_{\hRun} a_M}\\
{\mH(M)}\ar[r]^-(0.5){\theta}&{\xi^{-1}\tOmega^1_{R/\co_K}\otimes_R \mH(M)}}
\end{equation}
où $\theta$ et $\theta'$ sont les champs de Higgs canoniques. Par suite, $a_M$ définit un morphisme de foncteurs 
\begin{equation}\label{repdolb12f}
\varepsilon\circ \mH'\rightarrow \mH.
\end{equation}

De même, le diagramme \eqref{taht9p} induit un morphisme canonique de foncteurs 
\begin{equation}\label{repdolb12g}
\mV'\rightarrow \mV \circ \varepsilon.
\end{equation}

\begin{prop}\label{repdolb13}
Les hypothèses étant celles de \ref{repdolb12}, soient, de plus, $M$ une $\hoR[\frac 1 p]$-représen\-tation de Dolbeault de $\Delta$  relativement à $\co_K$ \eqref{repdolb8},
$N$ un $\hRun[\frac 1 p]$-module de Higgs soluble à coefficients dans $(\xi_\pi^{\ast})^{-1}\tOmega^1_{R/\co_K}$ \eqref{repdolb10}. 
Alors,
\begin{itemize}
\item[{\rm (i)}] La $\hoR[\frac 1 p]$-représentation $M$ de $\Delta$ est de Dolbeault absolue et le morphisme canonique 
\begin{equation}\label{repdolb13a}
\varepsilon(\mH'(M))\rightarrow \mH(M)
\end{equation}
est un isomorphisme.
\item[{\rm (ii)}] Le $\hRun[\frac 1 p]$-module de Higgs $\varepsilon(N)$ à coefficients dans $\xi^{-1}\tOmega^1_{R/\co_K}$ est soluble et le 
morphisme canonique  
\begin{equation}\label{repdolb13b}
\mV'(N)\rightarrow \mV(\varepsilon(N))
\end{equation}
est un isomorphisme.
\end{itemize}
\end{prop}

Supposons qu'il existe un nombre rationnel $r>0$ et un isomorphisme $\hcC'^{(r)}$-linéaire
et $\Delta$-équivariant de $\hoR[\frac 1 p]$-modules de Higgs à coefficients dans $(\xi^{\ast}_\pi)^{-1}\tOmega^1_{R/\co_K}$
\begin{equation}\label{repdolb13c}
N \otimes_{\hRun}\hcC'^{(r)}\stackrel{\sim}{\rightarrow}  M\otimes_{\hoR}\hcC'^{(r)},
\end{equation}
où $N$ est muni de l'action triviale de $\Delta$, $M$ est muni du champ de Higgs nul et 
$\hcC'^{(r)}$ est muni du champ de Higgs $p^rd_{\hcC'^{(r)}}$ et de l'action canonique de $\Delta$.
Celui-ci induit un isomorphisme de $\hRun$-modules de Higgs à coefficients dans $(\xi^{\ast}_\pi)^{-1}\tOmega^1_{R/\co_K}$ \eqref{repdolb9e}
\begin{equation}\label{repdolb13d}
N\stackrel{\sim}{\rightarrow}\mH'(M),
\end{equation}
et un isomorphisme $\hoR$-linéaire et $\Delta$-équivariant \eqref{repdolb11b}
\begin{equation}\label{repdolb13e}
\mV'(N)\stackrel{\sim}{\rightarrow} M.
\end{equation}
Compte tenu de \eqref{taht9m} et \eqref{taht9n}, on déduit de \eqref{repdolb13c} un isomorphisme $\hcC^{(r+\rho)}$-linéaire
et $\Delta$-équivariant de $\hoR[\frac 1 p]$-modules de Higgs à coefficients dans $\xi^{-1}\tOmega^1_{R/\co_K}$
\begin{equation}\label{repdolb13f}
\varepsilon(N) \otimes_{\hRun}\hcC^{(r+\rho)}\stackrel{\sim}{\rightarrow}  M\otimes_{\hoR}\hcC^{(r+\rho)},
\end{equation}
où $\varepsilon(N)$ est muni de l'action triviale de $\Delta$, $M$ est muni du champ de Higgs nul et 
$\hcC^{(r+\rho)}$ est muni du champ de Higgs $p^{r+\rho}d_{\hcC^{(r+\rho)}}$ et de l'action canonique de $\Delta$.
En vertu de \ref{repdolb17}, $M$ est de Dolbeault absolue, $\varepsilon(N)$ est soluble,
et on a un isomorphisme de $\hRun$-modules de Higgs à coefficients dans $\xi^{-1}\tOmega^1_{R/\co_K}$ 
\begin{equation}\label{repdolb13g}
\varepsilon(N)\stackrel{\sim}{\rightarrow}\mH(M),
\end{equation}
et un isomorphisme $\hoR$-linéaire et $\Delta$-équivariant 
\begin{equation}\label{repdolb13h}
\mV(\varepsilon(N))\stackrel{\sim}{\rightarrow} M.
\end{equation}
On vérifie aussitôt que les isomorphismes \eqref{repdolb13d} et \eqref{repdolb13g} sont compatibles via le morphisme \eqref{repdolb12f}, 
et les isomorphismes \eqref{repdolb13e} et \eqref{repdolb13h} sont compatibles via le morphisme \eqref{repdolb12g}, d'où la proposition.

\begin{prop}\label{repdolb14}
Les foncteurs $\mH$ \eqref{repdolb5b} et $\mV$ \eqref{repdolb6b} induisent des équivalences de catégories quasi-inverses l'une de l'autre,
entre la catégorie des $\hoR[\frac 1 p]$-représentations de Dolbeault de $\Delta$ et celle des 
$\hRun[\frac 1 p]$-modules de Higgs solubles à coefficients dans $\txi^{-1}\tOmega^1_{R/\co_K}$.
\end{prop}

Le cas absolu \eqref{definf10} a été démontré dans (\cite{agt} II.12.24). Le cas relatif s'ensuit compte tenu de \ref{taht8} et \ref{repdolb13}.

\begin{prop}\label{repdolb15}
Soient $M$ une $\hoR[\frac 1 p]$-représentation de Dolbeault de $\Delta$,  
$(\mH(M),\theta)$ le $\hRun[\frac 1 p]$-module de Higgs à coefficients 
dans $\txi^{-1}\tOmega^1_{R/\co_K}$ associé. 
On a alors un isomorphisme canonique fonctoriel dans $\bD^+(\bMod(\hRun[\frac 1 p]))$ 
\begin{equation}\label{repdolb15a}
\rC_\cont^\bullet(\Delta, M)\stackrel{\sim}{\rightarrow} \mK^\bullet(\mH(M),\theta),
\end{equation}
où $\rC_\cont^\bullet(\Delta, M)$ est le complexe de cochaînes continues de $\Delta$ à valeurs dans $M$
et $\mK^\bullet(\mH(M),\theta)$ est le complexe de Dolbeault de $(\mH(M),\theta)$ \eqref{MH2a}. 
\end{prop}

Le cas absolu \eqref{definf10} a été démontré dans (\cite{agt} II.12.26) et le cas relatif s'ensuit compte tenu de \ref{taht8} et \ref{repdolb13}.

\section{Petits modules de Higgs}\label{pmh}

\subsection{}\label{pmh100}
On note $K_0$ le corps des fractions de $W$ \eqref{definf1} et $\fd$ la différente de l'extension $K/K_0$.
On pose $\rho=0$ dans le cas absolu \eqref{definf10} et $\rho=v(\pi\fd)$ dans le cas relatif. 
D'après \eqref{definf16a} et \eqref{definf17c}, on a un isomorphisme $\co_C$-linéaire canonique
\begin{equation}\label{mph1g}
\co_C(1)\stackrel{\sim}{\rightarrow} p^{\rho+\frac{1}{p-1}}\txi\co_C.
\end{equation}

\begin{defi}[\cite{agt} II.13.1]\label{pmh3}
Soient $G$ un groupe topologique, $A$ une $\co_C$-algèbre complète et séparée pour la topologie $p$-adique,
munie d'une action continue de $G$ (par des homomorphismes de $\co_C$-algèbres), 
$\alpha$ un nombre rationnel $>0$, 
$M$ une $A$-représentation continue de $G$, munie de la topologie $p$-adique.  
\begin{itemize}
\item[(i)] On dit que $M$ est {\em $\alpha$-quasi-petite} si le $A$-module $M$ est complet et séparé pour 
la topologie $p$-adique, et est engendré par un nombre fini d'éléments $G$-invariants modulo $p^{\alpha}M$. 
\item[(ii)] On dit que $M$ est {\em quasi-petite} si elle est $\alpha'$-quasi-petite pour un nombre rationnel $\alpha'>\rho+\frac{2}{p-1}$ \eqref{pmh100}.
\end{itemize}
\end{defi} 

On désigne par $\bRep^{\alpha\trqpp}_{A}(G)$ (resp. $\bRep^{\qpp}_{A}(G)$) la sous-catégorie pleine de 
$\bRep^\cont_A(G)$ formée des $A$-représentations $\alpha$-quasi-petites (resp. quasi-petites) de $G$ 
dont le $A$-module sous-jacent est $\co_C$-plat.

\begin{defi}[\cite{agt} II.13.4]\label{pmh4}
Soient $\varepsilon$ un nombre rationnel $>0$, $(N,\theta)$ un $\hRun$-module de Higgs  
à coefficients dans $\txi^{-1}\tOmega^1_{R/\co_K}$ \eqref{taht18}.
\begin{itemize}
\item[(i)] On dit que $(N,\theta)$ est {\em $\varepsilon$-quasi-petit} si $N$ est de type fini sur $\hRun$
et si $\theta$ est un multiple de  $p^{\varepsilon}$ dans $\txi^{-1}\End_{\hRun}(N)\otimes_R\tOmega^1_{R/\co_K}$
\eqref{MH2g}. On dit alors aussi que le $\hRun$-champ de Higgs $\theta$
est {\em $\varepsilon$-quasi-petit}. 
\item[{\rm (ii)}] On dit que $(N,\theta)$ est {\em quasi-petit} s'il est $\varepsilon'$-quasi-petit pour un nombre rationnel $\varepsilon'>\frac{1}{p-1}$. 
On dit alors aussi que le $\hRun$-champ de Higgs $\theta$ est {\em quasi-petit}. 
\end{itemize}
\end{defi}

On note $\bMH^{\varepsilon\trqpp}(\hRun,\txi^{-1}\tOmega^1_{R/\co_K})$ (resp. $\bMH^{\qpp}(\hRun,\txi^{-1}\tOmega^1_{R/\co_K})$)
la sous-catégorie pleine de $\bMH(\hRun,\txi^{-1}\tOmega^1_{R/\co_K})$
formée des $\hRun$-modules de Higgs $\varepsilon$-quasi-petits (resp. quasi-petits) dont le $\hRun$-module sous-jacent est $\co_C$-plat.

\subsection{}\label{pmh1}
Reprenons les notations de \ref{cad1}. On rappelle que le $\mZ$-module $P^\gp/\mZ\lambda$ est libre de type fini en vertu de (\cite{ag} 4.2.2). 
Soient $t_1,\dots,t_d$ des éléments de $P^\gp$ tels que leurs images dans $P^\gp/\mZ\lambda$ forment une $\mZ$-base. 
Pour tout $1\leq i\leq d$, notons $d\log(t_i)$ l'image de $t_i$ par le morphisme \eqref{cad1g}
et posons $y_i=\txi^{-1}d\log(t_i)\in \txi^{-1}\tOmega^1_{R/\co_K}\otimes_R\hRun$, de sorte que $(y_i)_{1\leq i\leq d}$ forme une $\hRun$-base de  
$\txi^{-1}\tOmega^1_{R/\co_K}\otimes_R\hRun$.
Notons $\chi_{t_i}$ l'image de $t_i$ par l'homomorphisme \eqref{taht202e}, et $\chi_i$ l'homomorphisme composé 
\begin{equation}\label{pmh1h}
\xymatrix{
{\Delta_\infty}\ar[r]^-(0.5){\chi_{t_i}}& 
{\mZ_p(1)}\ar[r]^-(0.5){\log([\ ])}&{p^{\rho+\frac{1}{p-1}}\txi\co_C}},
\end{equation}
où  la seconde flèche est induite par l'isomorphisme \eqref{mph1g}. On note encore $\chi_i\colon \Delta\rightarrow p^{\rho+\frac{1}{p-1}}\txi\co_C$ 
l'homomorphisme induit.

\subsection{}\label{pmh5}
Soient $M$ un $\hRun$-module de type fini et $\co_C$-plat, $\varepsilon$ un nombre rationnel $>0$, 
$\alpha=\varepsilon+\rho+\frac{1}{p-1}$ \eqref{pmh100}. 
On désigne par $\Psi_M$ l'isomorphisme composé
\begin{equation}\label{pmh5a}
\xymatrix{
{\Hom_{\mZ}(\Delta_\infty, p^{\alpha}\End_{\hRun}(M))}\ar[r]^-(0.5)\sim\ar[dr]_{\Psi_M}&
{p^{\varepsilon}\txi^{-1}\End_{\hRun}(M)\otimes_{\hRun}\Hom_{\mZ}(\Delta_\infty,\hRun(1))}
\ar[d]\\
&{p^{\varepsilon}\txi^{-1}\End_{\hRun}(M)\otimes_{R}\tOmega^1_{R/\co_K}}}
\end{equation}
où l'isomorphisme vertical est induit par l'isomorphisme \eqref{taht202h} et l'isomorphisme horizontal provient de l'isomorphisme \eqref{mph1g} et de (\cite{agt} (II.6.12.2)). 
On notera que le $\hRun$-module $\End_{\hRun}(M)$ est complet et séparé pour la topologie $p$-adique 
et $\co_C$-plat (\cite{agt} II.13.9).

Pour tout homomorphisme
$\phi\colon \Delta_\infty \rightarrow p^{\alpha}\End_{\hRun}(M)$, avec les notations de \ref{pmh1}, 
il existe $\phi_i\in p^{\varepsilon}\End_{\hRun}(M)$ $(1\leq i\leq d)$ tels que 
\begin{equation}\label{pmh5b}
\phi=\sum_{i=1}^d\txi^{-1}\phi_i\otimes \chi_i.
\end{equation}
On a alors (\cite{agt} II.13.10)
\begin{equation}\label{pmh5c}
\Psi_M(\phi)=\sum_{i=1}^d\txi^{-1}\phi_i \otimes d\log(t_i).
\end{equation}

Soit $\varphi$ une $\hRun$-représentation $\alpha$-quasi-petite de $\Delta_\infty$ sur $M$ \eqref{pmh3}.
Comme $\Delta_\infty$ agit trivialement sur $\hRun$, $\varphi$ est un homomorphisme 
\begin{equation}\label{pmh5d}
\varphi\colon \Delta_\infty\rightarrow \Aut_{\hRun}(M)
\end{equation}
d'image contenue dans le sous-groupe $\id+p^{\alpha}\End_{\hRun}(M)$ de $\Aut_{\hRun}(M)$.
Comme $\Delta_\infty$ est abélien, on peut définir l'homomorphisme (\cite{agt} II.13.9)
\begin{equation}\label{pmh5e}
\log(\varphi)\colon \Delta_\infty\rightarrow p^{\alpha}\End_{\hRun}(M).
\end{equation}
On voit aussitôt que $\Psi_M(\log(\varphi))\wedge \Psi_M(\log(\varphi))=0$ (\cite{agt} II.13.10), 
autrement dit, $\Psi_M(\log(\varphi))$ est un $\hRun$-champ de Higgs 
$\varepsilon$-quasi-petit sur $M$ à coefficients dans $\txi^{-1}\tOmega^1_{R/\co_K}$ \eqref{pmh4}. 
On obtient ainsi un foncteur 
\begin{equation}\label{pmh5f}
\begin{array}[t]{clcr}
\bRep_{\hRun}^{\alpha\trqpp}(\Delta_\infty)&\rightarrow &\bMH^{\varepsilon\trqpp}(\hRun,\txi^{-1}\tOmega^1_{R/\co_K})\\
(M,\varphi)&\mapsto& (M,\Psi_M(\log(\varphi))).
\end{array}
\end{equation}

Soit $\theta$ un $\hRun$-champ de Higgs $\varepsilon$-quasi-petit sur $M$ à coefficients dans 
$\txi^{-1}\tOmega^1_{R/\co_K}$. Comme $\theta\wedge \theta=0$,
l'image de l'homomorphisme $\Psi^{-1}_M(\theta)\colon \Delta_\infty \rightarrow p^\alpha\End_{\hRun}(M)$ 
est formée d'endomorphismes de $M$ qui commutent deux à deux  (cf. \cite{agt} 13.10). 
On peut donc définir l'homomorphisme (\cite{agt} II.13.9)
\begin{equation}\label{pmh5g}
\exp(\Psi^{-1}_M(\theta))\colon \Delta_\infty \rightarrow \Aut_{\hRun}(M),
\end{equation}
qui est clairement une $\hRun$-représentation $\alpha$-quasi-petite de $\Delta_\infty$ sur $M$. 
On définit ainsi un foncteur 
\begin{equation}\label{pmh5h}
\begin{array}[t]{clcr}
\bMH^{\varepsilon\trqpp}(\hRun,\txi^{-1}\tOmega^1_{R/\co_K})&\rightarrow& \bRep_{\hRun}^{\alpha\trqpp}(\Delta_\infty)\\ 
(M,\theta)&\mapsto& (M,\exp(\Psi^{-1}_M(\theta))). 
\end{array}
\end{equation}

\subsection{}\label{pmh8}
D'après la condition \ref{cad1}(iv), il existe essentiellement un unique morphisme étale 
\begin{equation}\label{pmh8a}
(\tX_0,\cM_{\tX_0})\rightarrow (\tS,\cM_{\tS}) \times_{\bA_\mN}\bA_P
\end{equation}
qui s'insère dans un diagramme commutatif à carrés cartésiens
\begin{equation}\label{pmh8b}
\xymatrix{
{(\coX,\cM_\coX)}
\ar[r]\ar[d]&{(\coS,\cM_{\coS})\times_{\bA_\mN}\bA_P}\ar[d]\ar[r]&{(\coS,\cM_{\coS})}\ar[d]^{i_S}\\
{(\tX_0,\cM_{\tX_0})}\ar[r]&{(\tS,\cM_{\tS})
\times_{\bA_\mN}\bA_P}\ar[r]\ar[d]&{(\tS,\cM_{\tS})}\ar[d]^a\\
&{\bA_P}\ar[r]&{\bA_\mN}}
\end{equation}
où le morphisme $a$ est défini par la carte 
$\mN\rightarrow \Gamma(\tS,\cM_{\tS}), 1\mapsto [\upi]$ (\eqref{definf6b} ou \eqref{definf7b}). 
On dit que {\em $(\tX_0,\cM_{\tX_0})$ est la $(\tS,\cM_{\tS})$-déformation lisse de 
$(\coX,\cM_\coX)$ définie par la carte adéquate $((P,\gamma),(\mN,\iota),\vartheta)$}. 
On désigne par $\cL_0$ le torseur de Higgs-Tate associé à $(\tX_0,\cM_{\tX_0})$ \eqref{taht5}, par $\cC_0$ la $\hoR$-algèbre de Higgs-Tate 
et par $\cF_0$ la $\hoR$-extension de Higgs-Tate associées à $(\tX_0,\cM_{\tX_0})$ \eqref{taht7}.

On vérifie aussitôt que le diagramme 
\begin{equation}\label{pmh8c}
\xymatrix{
{(\hmX,\cM_\hmX)}\ar[r]^-(0.4){i_X}\ar[d]&{(\tmX,\cM_{\tmX})}\ar[r]^-(0.4)b\ar[d]&{\bA_P}\ar[d]\\
{(\coS,\cM_\coS)}\ar[r]^-(0.4){i_S}&{(\tS,\cM_{\tS})}\ar[r]^-(0.4)a&{\bA_\mN}}
\end{equation}
où le morphisme $b$ est l'homomorphisme canonique \eqref{taht4}, est commutatif. 
On en déduit un morphisme $\phi_0$ qui s'insère dans le diagramme commutatif (sans la flèche pointillée) 
\begin{equation}\label{pmh8d}
\xymatrix{
{(\hmX,\cM_{\hmX})}\ar[r]\ar[d]_{i_X}&{(\coX,\cM_\coX)}\ar[r]&{(\tX_0,\cM_{\tX_0})}\ar[d]\\
{(\tmX,\cM_{\tmX})}\ar[rr]^{\phi_0}\ar@{.>}[rru]^{\psi_0}&&{(\tS,\cM_{\tS})
\times_{\bA_\mN}\bA_P}}
\end{equation}
On peut compléter ce dernier par une unique flèche pointillée  
$\psi_0\in \cL_0(\hmX)$ de façon à le laisser commutatif. 
On dit que {\em $\psi_0$ est la section de $\cL_0(\hmX)$ définie par la carte $((P,\gamma),(\mN,\iota),\vartheta)$}.

\subsection{}\label{pmh9}
Distinguons provisoirement dans ce numéro le cas absolu du cas relatif \eqref{definf10}:
notons $(\tX_0,\cM_{\tX_0})$ (resp. $(\tX'_0,\cM_{\tX'_0})$) la déformation 
lisse de $(\coX,\cM_\coX)$ au-dessus de 
\[
(\cA_2(\oS),\cM_{\cA_2(\oS)})\ \ \ ({\rm resp.}\  (\cA^*_2(\oS/S),\cM_{\cA^*_2(\oS/S)}))
\] 
définie par la carte adéquate $((P,\gamma),(\mN,\iota),\vartheta)$ \eqref{pmh8}, $\cL_0$ (resp. $\cL'_0$) le torseur de Higgs-Tate associé \eqref{taht5}
et $\psi_0\in \cL_0(\hmX)$ (resp. $\psi'_0\in \cL'_0(\hmX)$) la section définie par la même carte adéquate. 
On vérifie aussitôt que le diagramme \eqref{taht4}
\begin{equation}\label{pmh9a}
\xymatrix{
{(\cA_2^*(\mX/S),\cM_{\cA_2^*(\mX/S)})}\ar[r]\ar[d]_{\psi_0'}&{(\cA_2(\mX),\cM_{\cA_2(\mX)})}\ar[d]^{\psi_0}\\
{(\tX'_0,\cM_{\tX'_0})}\ar[r]\ar[d]\ar@{}[rd]|\Box&{(\tX_0,\cM_{\tX_0})}\ar[d]\\
{(\cA_2^*(\oS/S),\cM_{\cA_2^*(\oS/S)})}\ar[r]&{(\cA_2(\oS),\cM_{\cA_2(\oS)})}}
\end{equation}
est commutatif et que le carré inférieur est cartésien. En particulier, $\psi'_0$ est induite par $\psi_0$.

\begin{prop}\label{pmh90}
Conservons les notations de \ref{pmh8}. Pour tous $t\in P^\gp$ et $g\in \Delta$, on a 
\begin{equation}\label{pmh90a}
(\psi_0-{^g\psi_0})(d\log(t))=-\log([\chi_t(g)]),
\end{equation}  
où ${^g\psi_0}$ est l'image de $\psi_0$ par l'isomorphisme \eqref{taht6h}, $\psi_0-{^g\psi_0}$ est vu comme un élément de $\bT(\hmX)=\rT$ \eqref{taht5a},
$d\log(t)$ est l'image canonique de $t$ dans $\tOmega^1_{R/\co_K}$ et 
$\log([\chi_t])$ désigne l'homomorphisme composé 
\begin{equation}\label{pmh90b}
\xymatrix{
\Delta\ar[r]&{\Delta_\infty}\ar[r]^-(0.5){\chi_t}& 
{\mZ_p(1)}\ar[r]^-(0.5){\log([\ ])}&{p^{\rho+\frac{1}{p-1}}\txi\hoR}\ar[r]&{\txi\hoR}},
\end{equation}
où la première et la dernière flèches sont les morphismes canoniques, $\chi_t$ est l'image de $t$ par l'homomorphisme \eqref{taht202e} 
et la troisième flèche est induite par l'isomorphisme \eqref{mph1g}. 
\end{prop}
 
Comme les deux membres de l'équation \eqref{pmh90a} sont des homomorphismes de $P^\gp$ dans $\txi\hoR$, on peut se borner au cas où $t\in P$.   
Les morphismes $\phi_0$ et $\phi_0\circ g^{-1}$, où $\phi_0$ est le morphisme défini dans \eqref{pmh8d} 
et $g^{-1}$ agit sur $(\tmX,\cM_{\tmX})$, prolongent le même morphisme 
\[
(\hmX,\cM_\hmX)\rightarrow (\tS,\cM_{\tS})\times_{\bA_\mN}\bA_P.
\]
D'après les définitions et la condition \ref{cad1}(iv), la différence $\phi_0-\phi_0\circ g^{-1}$ 
correspond au morphisme $\psi_0-{^g\psi_0}\in \rT$. 
D'autre part, on a $g(\upnu(t))=[\chi_t(g)]\cdot \upnu(t)$ dans $\Gamma(\tmX,\cM_{\tmX})$ \eqref{taht200ab}.
La proposition s'ensuit compte tenu de \eqref{definf17b} et (\cite{agt} II.5.23).

\subsection{}\label{pmh10}
Posons $\bL_0=\Spec(\cC_0)$ \eqref{taht5h} qui est naturellement un $\bT$-fibré principal homogène sur $\hmX$ qui représente $\cL_0$. 
Considérons l'isomorphisme de $\bT$-fibrés principaux homogènes sur $\hmX$ 
\begin{equation}\label{pmh10a}
\ttt_0\colon \bT\stackrel{\sim}{\rightarrow} \bL_0,\ \ \ v\mapsto v+\psi_0.
\end{equation}
La structure de $\bT$-fibré principal homogène $\Delta$-équivariant sur $\bL_0$ \eqref{taht6g} se transporte 
par $\ttt_0$ en une structure de $\bT$-fibré principal homogène $\Delta$-équivariant sur $\bT$ (cf. \cite{agt} II.4.20).
Pour tout $g\in \Delta$, on a donc un isomorphisme $\tau_g^\bT$-équivariant
\begin{equation}\label{pmh10b}
\tau_{g,\psi_0}^{\bT}\colon \bT\stackrel{\sim}{\rightarrow} g^\bullet(\bT).
\end{equation}
Cette structure détermine une action à gauche de $\Delta$ sur $\bT$ compatible avec son action sur $\hmX$. 
On en déduit une action 
\begin{equation}\label{pmh10c}
\varphi_0\colon \Delta\rightarrow \Aut_\hRun(\cG)
\end{equation}
de $\Delta$ sur $\cG$ \eqref{taht5b} par des automorphismes d'anneaux, 
compatible avec son action sur $\hoR$~;  pour tout $g\in \Delta$, 
$\varphi_0(g)$ est induit par l'automorphisme de $\bT$ défini par $g^{-1}$. 

\begin{prop}\label{pmh11}
Avec les notations de \ref{pmh1}, pour tout $g\in \Delta$, on a \eqref{pmh10c}
\begin{equation}\label{pmh11a}
\varphi_0(g)=\exp(-\sum_{i=1}^d\txi^{-1}\frac{\partial}{\partial y_i}\otimes \chi_i(g)) \circ g.
\end{equation}
\end{prop}

Cela résulte de \ref{pmh90} et (\cite{agt} (II.10.11.6)), cf. (\cite{agt} II.10.17) où le cas absolu \eqref{definf10} a été démontré. 
On peut aussi déduire le cas relatif du cas absolu comme suit. 
Reprenons les notations de \ref{pmh9} et posons 
\begin{eqnarray}
\cG&=&\rS_{\hoR}(\xi^{-1}\tOmega^1_{R/\co_K}\otimes_R\hoR),\label{pmh11b}\\
\cG'&=&\rS_{\hoR}((\xi^*_\pi)^{-1}\tOmega^1_{R/\co_K}\otimes_R\hoR).\label{pmh11c}
\end{eqnarray}
On note $\varphi_0$ (resp. $\varphi'_0$) l'action de $\Delta$ sur $\cG$ (resp. $\cG'$) induite par la section $\psi_0$ (resp. $\psi'_0$). 
D'après \ref{definf16}, on a un isomorphisme canonique 
\begin{equation}\label{pmh11d}
(\xi^*_\pi)^{-1}\co_C\stackrel{\sim}{\rightarrow}p^\rho\xi^{-1}\co_C.
\end{equation}
On en déduit un homomorphisme injectif de $\hoR$-algèbres
\begin{equation}\label{pmh11e}
\cG'\rightarrow \cG.
\end{equation}
Posons $\bT=\Spec(\cG)$ et $\bT'=\Spec(\cG')$ et notons $h\colon \bT\rightarrow \bT'$ le morphisme 
de $\hmX$-fibrés vectoriels déduit de \eqref{pmh11e}. 
D'après \eqref{taht9d}, on a un morphisme canonique $h$-équivariant et $\Delta$-équivariant $\cL_0\rightarrow \cL'_0$. 
Celui-ci transforme $\psi_0$ en $\psi'_0$ compte tenu de \eqref{pmh9a}. Il s'ensuit que l'homomorphisme \eqref{pmh11e} est $\Delta$-équivariant
lorsque l'on munit $\cG$ et $\cG'$ des actions $\varphi_0$ et $\varphi'_0$ respectivement. La proposition dans le cas relatif résulte donc du cas absolu.

\subsection{}\label{pmh12}
Pour tout nombre rationnel $r\geq 0$, on désigne par $\cG^{(r)}$ la sous-$\hoR$-algèbre de $\cG$ \eqref{taht5b} définie par \eqref{notconv9}
\begin{equation}\label{pmh12a}
\cG^{(r)}=\rS_{\hoR}(p^r\txi^{-1}\tOmega^1_{R/\co_K}\otimes_R\hoR)
\end{equation}
et par $\hcG^{(r)}$ son séparé complété $p$-adique.
Compte tenu de (\cite{agt} II.6.14) et de sa preuve, $\cG^{(r)}$ et $\hcG^{(r)}$ sont $\co_C$-plats. 
Pour tous nombres rationnels $r'\geq r\geq 0$, on a un homomorphisme injectif canonique 
$a^{r,r'}\colon \cG^{(r')}\rightarrow \cG^{(r)}$. On vérifie aussitôt que l'homomorphisme 
induit $\ha^{r,r'}\colon\hcG^{(r')}\rightarrow \hcG^{(r)}$ est injectif. 
On a un $\cG^{(r)}$-isomorphisme canonique 
\begin{equation}\label{pmh12c}
\Omega^1_{\cG^{(r)}/\hoR}\stackrel{\sim}{\rightarrow}\txi^{-1}\tOmega^1_{R/\co_K}\otimes_R\cG^{(r)}.
\end{equation} 
On désigne par
\begin{equation}\label{pmh12d}
d_{\cG^{(r)}}\colon \cG^{(r)}\rightarrow \txi^{-1}\tOmega^1_{R/\co_K}\otimes_R\cG^{(r)}
\end{equation}
la $\hoR$-dérivation universelle de $\cG^{(r)}$ et par 
\begin{equation}\label{pmh12e}
d_{\hcG^{(r)}}\colon \hcG^{(r)}\rightarrow \txi^{-1}\tOmega^1_{R/\co_K}\otimes_R\hcG^{(r)}
\end{equation}
son prolongement aux complétés (on notera que le $R$-module $\tOmega^1_{R/\co_K}$ est libre de type fini). 
Comme $\txi^{-1}\tOmega^1_{R/\co_K}\otimes_R \hoR\subset d_{\cG^{(r)}}(\cG^{(r)})$, 
$d_{\cG^{(r)}}$ et $d_{\hcG^{(r)}}$ sont également des $\hoR$-champs de Higgs à coefficients dans 
$\txi^{-1}\tOmega^1_{R/\co_K}$ d'après \ref{MH8}(i). Pour tous nombres rationnels $r'\geq r\geq 0$, on a 
\begin{equation}\label{pmh12f}
p^{r'-r}(\id \times a^{r,r'}) \circ d_{\cG^{(r')}}=d_{\cG^{(r)}}\circ a^{r,r'}.
\end{equation}

La section $\psi_0\in \cL_0(\hmX)$ définie par la carte adéquate $((P,\gamma),(\mN,\iota),\vartheta)$ \eqref{pmh8}
induit un isomorphisme de $\hoR$-algèbres 
\begin{equation}\label{pmh12g}
\cG\stackrel{\sim}{\rightarrow}\cC_0
\end{equation}
qui est $\Delta$-équivariant lorsque l'on munit $\cG$ de l'action $\varphi_0$  \eqref{pmh10c} et $\cC_0$ de l'action canonique. 

D'après \ref{pmh120} ci-dessous, la sous-$\hoR$-algèbre $\cG^{(r)}$ de $\cG$ est stable par l'action $\varphi_0$ de $\Delta$ sur $\cG$. 
Notant $\cC_0^{(r)}$ la $\hoR$-algèbre de Higgs-Tate d'épaisseur $r$ associée à $(\tX_0,\cM_{\tX_0})$ \eqref{taht10c}, 
on démontre que $\Spec(\cC_0^{(r)})$ est naturellement un $\Spec(\cG^{(r)})$-fibré principal homogène sur $\hmX$ (cf. \cite{agt} II.12.1). 
La section $\psi_0$ induit un $\hoR$-homomorphisme $\cC_0^{(r)}\rightarrow \hoR$ et par suite un isomorphisme de $\hoR$-algèbres 
\begin{equation}\label{pmh12h}
\cG^{(r)}\stackrel{\sim}{\rightarrow}\cC_0^{(r)}.
\end{equation}
Celui-ci est $\Delta$-équivariant lorsque l'on munit $\cG^{(r)}$ de l'action induite par $\varphi_0$  et $\cC_0^{(r)}$ de l'action canonique, 
et est compatible aux dérivations $d_{\cG^{(r)}}$ et $d_{\cC_0^{(r)}}$.

\begin{prop}[\cite{agt} II.11.6]\label{pmh120}
Pour tout nombre rationnel $r\geq 0$, la sous-$\hoR$-algèbre $\cG^{(r)}$ \eqref{pmh12a} de $\cG$ 
est stable par l'action $\varphi_0$ de $\Delta$ sur $\cG$ \eqref{pmh10c}, et les actions induites  
de $\Delta$ sur $\cG^{(r)}$ et $\hcG^{(r)}$ sont continues pour les topologies $p$-adiques. 
\end{prop}

Reprenons les hypothèses et notations de \ref{pmh1}. D'après \ref{pmh11}, pour tout $g\in \Delta$ et tout $1\leq i\leq d$, on a 
\begin{equation}\label{pmh120a}
\varphi_{0}(g)(y_i)=y_i-\txi^{-1}\chi_i(g).
\end{equation} 
Comme $\txi^{-1}\chi_i(g)\in p^{\rho+\frac{1}{p-1}}\co_C$ \eqref{pmh1h}, $\cG^{(r)}$ est stable par $\varphi_0(g)$.
Soit $\zeta$ un générateur de $\mZ_p(1)$. Il existe $a_g\in \mZ_p$ tel que $\chi_i(g)=[\zeta^{a_g}]-1\in \cA_2(\co_\oK)$. 
Par linéarité, on a $\log([\zeta^{a_g}])\in p^{v_p(a_g)}\xi \co_C$, et par suite $\varphi_{0}(g)(y_i)- y_i \in p^{\rho+v_p(a_g)}\cG$ \eqref{definf16a}. 
Pour tout entier $n\geq 0$, l'ensemble des $g\in \Delta$ tels que $\rho+v_p(a_g)\geq n$ 
étant un sous-groupe ouvert de $\Delta$, on en déduit que le stabilisateur de 
la classe de $p^ry_i$ dans $\cG^{(r)}/p^n\cG^{(r)}$ est ouvert dans $\Delta$.
La seconde assertion s'ensuit car l'action de $\Delta$ sur $\oR/p^n\oR$ est continue pour la topologie discrète.

\subsection{}\label{pmh16}
Reprenons les notations de \ref{pmh9} et pour tout nombre rationnel $r\geq 0$, posons 
\begin{eqnarray}
\cG^{(r)}&=&\rS_{\hoR}(p^r\xi^{-1}\tOmega^1_{R/\co_K}\otimes_R\hoR),\label{pmh16a}\\
\cG'^{(r)}&=&\rS_{\hoR}(p^r(\xi^*_\pi)^{-1}\tOmega^1_{R/\co_K}\otimes_R\hoR).\label{pmh16b}
\end{eqnarray}
On a alors un isomorphisme canonique 
\begin{equation}\label{pmh16c}
\cG'^{(r)}\stackrel{\sim}{\rightarrow}\cG^{(r+\rho)}.
\end{equation}
Celui-ci est compatible aux actions de $\Delta$ définies par les sections $\psi_0$ et $\psi'_0$ d'après \ref{pmh11}. De plus, le diagramme 
\begin{equation}\label{pmh16d}
\xymatrix{
{\cG'^{(r)}}\ar[r]^-(0.5){d_{\cG'^{(r)}}}\ar[d]&{(\xi^{\ast}_\pi)^{-1}\tOmega^1_{R/\co_K} \otimes_R \cG'^{(r)}}\ar[d]^{\nu}\\
{\cG^{(r+\rho)}}\ar[r]^-(0.5){d_{\cG^{(r+\rho)}}}&{\xi^{-1}\tOmega^1_{R/\co_K} \otimes_R \cG^{(r+\rho)}}}
\end{equation} 
où $d_{\cG'^{(r)}}$ et $d_{\cG^{(r+\rho)}}$ sont les $\hRun$-dérivations universelles \eqref{pmh12d} et $\nu$ est l'isomorphisme induit par \eqref{taht9b},
est commutatif

\subsection{}\label{pmh18}
Pour tout nombre rationnel $r\geq 0$, on désigne par $\fS^{(r)}$ la sous-$\hRun$-algèbre de $\cG^{(r)}$ \eqref{pmh12a}
définie par \eqref{notconv9}
\begin{equation}\label{pmh18a}
\fS^{(r)}=\rS_{\hRun}(p^r\txi^{-1}\tOmega^1_{R/\co_K}\otimes_R\hRun),
\end{equation}
et par $\hfS^{(r)}$ son séparé complété $p$-adique. On pose $\fS=\fS^{(0)}$ et $\hfS=\hfS^{(0)}$. 
On notera que $\hfS^{(r)}$ est $\hRun$-plat en vertu de (\cite{egr1} 1.12.4) et est donc $\co_C$-plat (\cite{agt} II.6.14). 
Pour tous nombres rationnels $r'\geq r\geq 0$, on a un homomorphisme injectif canonique 
$\tta^{r,r'}\colon \fS^{(r')}\rightarrow \fS^{(r)}$. On vérifie aussitôt que l'homomorphisme 
induit $\htta^{r,r'}\colon\hfS^{(r')}\rightarrow \hfS^{(r)}$ est injectif. On a un $\fS^{(r)}$-isomorphisme canonique 
\begin{equation}\label{pmh18b}
\Omega^1_{\fS^{(r)}/\hRun}\stackrel{\sim}{\rightarrow}\txi^{-1}\tOmega^1_{R/\co_K}\otimes_R\fS^{(r)}.
\end{equation} 
On note
\begin{equation}\label{pmh18c}
d_{\fS^{(r)}}\colon \fS^{(r)}\rightarrow \txi^{-1}\tOmega^1_{R/\co_K}\otimes_R\fS^{(r)}
\end{equation}
la $\hRun$-dérivation universelle de $\fS^{(r)}$ et par 
\begin{equation}\label{pmh18d}
d_{\hfS^{(r)}}\colon \hfS^{(r)}\rightarrow \txi^{-1}\tOmega^1_{R/\co_K}\otimes_R\hfS^{(r)}
\end{equation}
son prolongement aux complétés. 
Comme $\txi^{-1}\tOmega^1_{R/\co_K}\otimes_R \hRun\subset d_{\fS^{(r)}}(\fS^{(r)})$, 
$d_{\fS^{(r)}}$ et $d_{\hfS^{(r)}}$ sont également des $\hRun$-champs de Higgs à coefficients dans 
$\txi^{-1}\tOmega^1_{R/\co_K}$ d'après \ref{MH8}(i). Pour tous nombres rationnels $r'\geq r\geq 0$, on a 
\begin{equation}\label{pmh18e}
p^{r'-r}(\id \times \tta^{r,r'}) \circ d_{\fS^{(r')}}=d_{\fS^{(r)}}\circ \tta^{r,r'}.
\end{equation}

Calquant la preuve de \ref{pmh120}, on montre que pour tout nombre rationnel $r\geq 0$, 
l'action $\varphi_0$ de $\Delta$ sur $\cG$ \eqref{pmh10c} préserve la sous-$\hRun$-algèbre $\fS^{(r)}$, 
que l'action induite de $\Delta$ sur $\fS^{(r)}$ se factorise à travers $\Delta_{\infty}$ et que  l'action de $\Delta_\infty$ sur $\fS^{(r)}$ 
ainsi définie est continue pour la topologie $p$-adique.

\subsection{}\label{pmh15}
Soient $r$ un nombre rationnel $\geq 0$, $\lambda\in \hfS^{(r)}$, $M$ un $\hfS^{(r)}$-module. 
On appelle {\em $\lambda$-connexion $p$-adique sur $M$ relativement à l'extension $\hfS^{(r)}/\hRun$} la donnée d'un morphisme $\hRun$-linéaire 
\begin{equation}\label{pmh15a}
\nabla\colon M\rightarrow \txi^{-1}\tOmega^1_{R/\co_K} \otimes_R M
\end{equation}
tel que pour tous $t\in \hfS^{(r)}$ et $x\in M$, on ait 
\begin{equation}\label{pmh15b}
\nabla(tx)=\lambda d_{\hfS^{(r)}}(t)\otimes x+t\nabla(x).
\end{equation} 
On dit que $\nabla$ est {\em intégrable} si elle est un $\hRun$-champ de Higgs à coefficients dans $\txi^{-1}\tOmega^1_{R/\co_K}$. 
Si $M$ est complet et séparé pour la topologie $p$-adique, on retrouve la notion introduite dans  (\cite{agt} II.2.14) (cf. \cite{agt} II.2.16). 

Soient $\nabla$ une $\lambda$-connexion $p$-adique intégrable sur $M$ relativement à l'extension $\hfS^{(r)}/\hRun$, 
$(N,\theta)$ un $\hRun$-module de Higgs à coefficients dans  $\txi^{-1}\tOmega^1_{R/\co_K}$.
Il existe un et un unique morphisme $\hRun$-linéaire 
\begin{equation}\label{pmh15c}
\nabla'\colon M\otimes_{\hRun}N\rightarrow  \txi^{-1}\tOmega^1_{R/\co_K}\otimes_RM\otimes_{\hRun}N
\end{equation}
tel que pour tous $x\in M$ et $y\in N$, on ait 
\begin{equation}\label{pmh15d}
\nabla'(x\otimes_{\hRun} y)=\nabla(x)\otimes_{\hRun} y+ x\otimes_{\hRun} \theta(y).
\end{equation}
C'est une $\lambda$-connexion $p$-adique intégrable sur $M\otimes_{\hRun}N$ relativement à l'extension $\hfS^{(r)}/\hRun$. 

On définit de même la notion de $\lambda$-connexion $p$-adique intégrable sur un $\hcG^{(r)}$-module relativement à l'extension $\hcG^{(r)}/\hoR$. 
Si $\nabla$ est une $\lambda$-connexion $p$-adique intégrable sur $M$ relativement à l'extension $\hfS^{(r)}/\hRun$, 
il existe un et un unique morphisme $\hoR$-linéaire 
\begin{equation}\label{pmh15e}
\nabla'\colon \hcG^{(r)}\otimes_{\hfS^{(r)}}M\rightarrow \txi^{-1}\tOmega^1_{R/\co_K} \otimes_R \hcG^{(r)}\otimes_{\hfS^{(r)}}M
\end{equation}
tel que pour tous $t\in \hcG^{(r)}$ et $x\in M$, on ait
\begin{equation}\label{pmh15f}
\nabla'(t\otimes_{\hfS^{(r)}}x)=\lambda d_{\hcG^{(r)}}(t)\otimes_{\hfS^{(r)}} x+t\otimes_{\hfS^{(r)}}\nabla(x).
\end{equation} 
C'est une $\lambda$-connexion $p$-adique intégrable sur $\hcG^{(r)}\otimes_{\hfS^{(r)}}M$ relativement à l'extension $\hcG^{(r)}/\hoR$.

\subsection{}\label{pmh13}
Soient $r$, $\varepsilon$ deux nombres rationnels tels que $r\geq 0$ et $\varepsilon>r+\frac{1}{p-1}$, 
$(N,\theta)$ un $\hRun$-module de Higgs $\varepsilon$-quasi-petit à coefficients dans $\txi^{-1}\tOmega^1_{R/\co_K}$ \eqref{pmh4}
tel que $N$ soit $\co_C$-plat. On peut écrire de manière unique \eqref{pmh1}
\begin{equation}\label{pmh13a}
\theta=\sum_{i=1}^d\theta_i \otimes y_i,
\end{equation}
où les $\theta_i$ sont des endomorphismes de $N$ appartenant à $p^{\varepsilon}\End_{\hRun}(N)$  
et commutant deux à deux. Pour tout $\un=(n_1,\dots,n_d)\in \mN^d$, 
posons $|\un|=\sum_{i=1}^d n_i$, $\un!=\prod_{i=1}^d n_i!$, $\utheta^{\un}=\prod_{i=1}^d\theta_i^{n_i} \in \End_{\hRun}(N)$
et $\uy^{\un}=\prod_{i=1}^dy_i^{n_i}\in \fS$. 
On notera que $N\otimes_{\hRun}\hfS^{(r)}$ est complet
et séparé pour la topologie $p$-adique (\cite{egr1} 1.10.2), et qu'il est $\co_C$-plat puisque $\hfS^{(r)}$ est $\hRun$-plat. 
Par suite, pour tout $z\in N\otimes_{\hRun}\hfS^{(r)}$, la série 
\begin{equation}\label{pmh13b}
\sum_{\un\in \mN^d} \frac{1}{\un!}(\utheta^{\un}\otimes\uy^{\un})(z)
\end{equation}
converge dans $N\otimes_{\hRun}\hfS^{(r)}$, et définit un endomorphisme $\hfS^{(r)}$-linéaire 
de $N\otimes_{\hRun} \hfS^{(r)}$, que l'on note 
\begin{equation}\label{pmh13c}
\exp_r(\theta)\colon N\otimes_{\hRun} \hfS^{(r)}\rightarrow N\otimes_{\hRun} \hfS^{(r)}.
\end{equation}
Pour tout nombre rationnel $r'$ tel que $0\leq r'\leq r$, le diagramme 
\begin{equation}
\xymatrix{
{N\otimes_{\hRun} \hfS^{(r)}}\ar[rr]^-(0.5){\exp_r(\theta)}\ar[d]_-(0.4){\id \otimes \htta^{r',r}}&&
{N\otimes_{\hRun} \hfS^{(r)}}\ar[d]^-(0.4){\id \otimes \htta^{r',r}}\\
{N\otimes_{\hRun} \hfS^{(r')}}\ar[rr]^-(0.5){\exp_{r'}(\theta)}&&{N\otimes_{\hRun} \hfS^{(r')}}}
\end{equation}
est commutatif. On peut donc se permettre d'omettre l'indice $r$ de la notation $\exp_r(\theta)$ sans risque d'ambiguïté.

\begin{prop}\label{pmh14}
Soient $r$, $\varepsilon$ deux nombres rationnels tels que $r\geq 0$ et $\varepsilon>r+\frac{1}{p-1}$, 
$N$ un $\hRun$-module de type fini et $\co_C$-plat, 
$\theta$ un $\hRun$-champ de Higgs $\varepsilon$-quasi-petit sur $N$ à coefficients dans $\txi^{-1}\tOmega^1_{R/\co_K}$,
$\varphi$ la $\hRun$-représentation quasi-petite de $\Delta_\infty$ sur $N$ associée à $\theta$ par le foncteur \eqref{pmh5h}. 
Alors, l'endomorphisme \eqref{pmh13c}
\begin{equation}\label{pmh14a}
\exp_r(\theta)\colon N\otimes_{\hRun} \hfS^{(r)}\rightarrow N\otimes_{\hRun} \hfS^{(r)}
\end{equation}
est un isomorphisme $\Delta_\infty$-équivariant de $\hfS^{(r)}$-modules à $p^r$-connexion $p$-adique intégrable relativement à l'extension 
$\hfS^{(r)}/\hRun$ \eqref{pmh15}, où $\hfS^{(r)}$ est muni de l'action de $\Delta_{\infty}$ induite par $\varphi_0$, 
le module $N$ de la source est muni de l'action triviale de $\Delta_{\infty}$ 
et du $\hRun$-champ de Higgs $\theta$, et le module $N$ du but est muni de l'action $\varphi$
de $\Delta_{\infty}$ et du $\hRun$-champ de Higgs nul. 
En particulier, $\exp_r(\theta)$ est un isomorphisme de $\hRun$-modules de Higgs à coefficients dans
$\txi^{-1}\tOmega^1_{R/\co_K}$. 
\end{prop}

Le cas absolu \eqref{definf10} a été démontré dans (\cite{agt} II.13.15) et le cas relatif s'en déduit compte tenu de \ref{pmh16} et du fait que $N\otimes_{\hRun}\hfS^{(r)}$ est 
$\co_C$-plat. 

\begin{cor}\label{pmh17}
Sous les hypothèses de \ref{pmh14}, on a
un isomorphisme fonctoriel et $\Delta$-équiva\-riant de $\bMC^r_p$ \eqref{repdolb3}
\begin{equation}\label{pmh17a}
N\otimes_{\hRun} \hcC^{(r)}\stackrel{\sim}{\rightarrow} N\otimes_{\hRun} \hcC^{(r)},
\end{equation}
où $\hcC^{(r)}$ est muni de l'action canonique de $\Delta$,
le module $N$ de la source est muni de l'action triviale de $\Delta_{\infty}$ 
et du $\hRun$-champ de Higgs $\theta$, et le module $N$ du but est muni de l'action $\varphi$
de $\Delta_{\infty}$ et du $\hRun$-champ de Higgs nul \eqref{repdolb3d}. 
Si, de plus, la déformation $(\tX,\cM_\tX)$  est définie par la carte adéquate $((P,\gamma),(\mN,\iota),\vartheta)$ \eqref{pmh8}, l'isomorphisme est canonique. 
\end{cor}

En effet, d'après \ref{pmh14} et compte tenu de \ref{pmh15}, $\exp_r(\theta)$ induit un isomorphisme fonctoriel et $\Delta$-équivariant de $\hcG^{(r)}$-modules à 
$p^r$-connexion $p$-adique intégrable relativement à l'extension $\hcG^{(r)}/\hoR$ 
\begin{equation}\label{pmh17b}
N\otimes_{\hRun} \hcG^{(r)}\stackrel{\sim}{\rightarrow} N\otimes_{\hRun} \hcG^{(r)},
\end{equation}
où $\hcG^{(r)}$ est muni de l'action de $\Delta$ induite par $\varphi_0$ \eqref{pmh12},
le module $N$ de la source est muni de l'action triviale de $\Delta_{\infty}$ 
et du $\hRun$-champ de Higgs $\theta$, et le module $N$ du but est muni de l'action $\varphi$
de $\Delta_{\infty}$ et du $\hRun$-champ de Higgs nul. 
La proposition s'ensuit compte tenu de \ref{taht8} et \eqref{pmh12h}.

\begin{defi}\label{pmh21}
Soient $G$ un groupe topologique, $A$ une $\co_C$-algèbre complète et séparée pour la topologie $p$-adique,
munie d'une action continue de $G$ (par des homomorphismes de $\co_C$-algèbres).
On munit $A[\frac 1 p]$ de la topologie $p$-adique \eqref{notconv15}.
On dit qu'une $A[\frac 1 p]$-représentation continue $M$ de $G$ est {\em petite} si
les conditions suivantes sont remplies~: 
\begin{itemize}
\item[(i)] $M$ est un $A[\frac 1 p]$-module projectif de type fini, muni de la topologie $p$-adique \eqref{notconv15}~;
\item[(ii)] il existe un nombre rationnel $\alpha>\rho+\frac{2}{p-1}$ \eqref{pmh100}
et un sous-$A$-module de type fini $M^\circ$ de $M$, stable par $G$, 
engendré par un nombre fini d'éléments $G$-invariants modulo $p^{\alpha}M^\circ$,
et qui engendre $M$ sur $A[\frac 1 p]$. 
\end{itemize}
\end{defi}

On désigne par $\bRep^{\p}_{A[\frac 1 p]}(G)$
la sous-catégorie pleine de $\bRep^\cont_{A[\frac 1 p]}(G)$ formée des $A$-représenta\-tions petites de $G$ \eqref{notconv16}. 
On notera que cette définition correspond à celle donnée dans (\cite{agt} II.13.2) dans le cas absolu \eqref{definf10}. 

\begin{remas}\label{pmh22}
Soient $G$ un groupe topologique, $A$ une $\co_C$-algèbre complète et séparée pour la topologie $p$-adique,
munie d'une action continue de $G$ (par des homomorphismes de $\co_C$-algèbres).
\begin{itemize}
\item[{\rm (i)}] Soient $M$ un $A[\frac 1 p]$-module projectif de type fini, $M^\circ$ un sous-$A$-module de type fini
de $M$. Alors $M^\circ$ est complet et séparé pour la topologie $p$-adique. En effet, $M^\circ$ est complet en vertu de 
(\cite{ac} chap. III §2.12 cor.~1 de prop.~16). D'autre part, quitte à ajouter à $M$ un facteur direct, 
on peut le supposer libre de type fini sur $A[\frac 1 p]$. Par suite, il existe un entier $m\geq 0$ tel que $p^m M^\circ$
soit contenu dans un $A$-module libre de type fini $N$. Donc $\cap_{n\geq 0}p^nM^\circ
\subset \cap_{n\geq 0}p^nN=0$. 
\item[{\rm (ii)}] Soient $M$ une $A[\frac 1 p]$-représentation petite de $G$, 
$M^\circ$ un sous-$A$-module de type fini de $M$ vérifiant la condition \ref{pmh21}(ii). 
Il résulte alors de (i) que $M^\circ$ est une $A$-représentation quasi-petite de $G$ \eqref{pmh3}. 
\end{itemize}
\end{remas}

\subsection{}\label{pmh23}
Soient $G$ un groupe topologique, $A$ une $\co_C$-algèbre complète et séparée pour la topologie $p$-adique,
munie d'une action continue de $G$ (par des homomorphismes de $\co_C$-algèbres). On désigne par $\bRep'^{\qpp}_{A}(G)$ la sous-catégorie 
pleine de la catégorie $\bRep^{\qpp}_{A}(G)$ \eqref{pmh3} formée des $A$-représentations $M$ de $G$ telle que le $A[\frac 1 p]$-module sous-jacent à 
$M[\frac 1 p]$ soit projectif de type fini. C'est une catégorie additive. On note $\bRep'^{\qpp}_{A,\mQ}(G)$ la catégorie des objets 
de $\bRep'^{\qpp}_{A}(G)$ à isogénie près. Le foncteur 
\begin{equation}\label{pmh23a}
\bRep'^{\qpp}_{A}(G)\rightarrow \bRep^{\p}_{A[\frac 1 p]}(G),\ \ \ M\mapsto M[\frac 1 p]
\end{equation}
induit alors un foncteur 
\begin{equation}\label{pmh23b}
\bRep'^{\qpp}_{A,\mQ}(G)\rightarrow \bRep^{\p}_{A[\frac 1 p]}(G).
\end{equation}

\begin{lem}\label{pmh24}
Le foncteur \eqref{pmh23b} est une équivalence de catégories.
\end{lem}
En effet, ce foncteur est essentiellement surjectif d'après \ref{pmh22}. Soient $M,N$ deux $A$-modules de type fini et $\co_C$-plats. 
Alors, le morphisme canonique 
\begin{equation}\label{pmh24a}
\Hom_A(M,N)\otimes_{\mZ_p}\mQ_p\rightarrow \Hom_{A[\frac 1 p]}(M[\frac 1 p],N[\frac 1 p])
\end{equation}
est un isomorphisme. Supposons $M$ et $N$ munis d'actions $A$-semi-linéaires de $G$. 
Le morphisme canonique 
\begin{equation}\label{pmh24b}
\Hom_{A\langle G\rangle}(M,N)\otimes_{\mZ_p}\mQ_p\rightarrow \Hom_{A[\frac 1 p]\langle G\rangle}(M[\frac 1 p],N[\frac 1 p])
\end{equation}
où la source (resp. le but) désigne l'ensemble des morphismes de $A$-représentations (resp. $A[\frac 1 p]$-représentations) de $G$ est alors un isomorphisme:
l'injectivité résulte aussitôt de celle de \eqref{pmh24a} et la surjectivité de celle de \eqref{pmh24a} et du fait que $M$ et $N$ sont $\co_C$-plats; d'où la proposition.

\begin{defi}\label{pmh25}
On dit qu'un $\hRun[\frac 1 p]$-module de Higgs $(N,\theta)$ à coefficients dans $\txi^{-1}\tOmega^1_{R/\co_K}$ 
est {\em petit} si les conditions suivantes sont remplies~: 
\begin{itemize}
\item[(i)] $N$ est un $\hRun[\frac 1 p]$-module projectif de type fini~; 
\item[(ii)] il existe un nombre rationnel $\varepsilon>\frac{1}{p-1}$ et 
un sous-$\hRun$-module de type fini $N^\circ$ de $N$, qui l'engendre sur $\hRun[\frac 1 p]$, tels que l'on ait 
\begin{equation}\label{pmh25a}
\theta(N^\circ)\subset p^{\varepsilon}\txi^{-1}N^\circ \otimes_R\tOmega^1_{R/\co_K}.
\end{equation}
\end{itemize}
\end{defi}
On désigne par $\bMH^{\p}(\hRun[\frac 1 p],\txi^{-1}\tOmega^1_{R/\co_K})$ la sous-catégorie pleine de 
$\bMH(\hRun[\frac 1 p],\txi^{-1}\tOmega^1_{R/\co_K})$ formée des $\hRun[\frac 1 p]$-modules de Higgs petits. 
On notera que cette définition correspond à celle donnée dans (\cite{agt} II.13.5) dans le cas absolu \eqref{definf10}.

\begin{lem}[\cite{agt} II.13.7]\label{pmh26}
Soit  $(N,\theta)$ un $\hRun[\frac 1 p]$-module de Higgs à coefficients dans $\txi^{-1}\tOmega^1_{R/\co_K}$ 
tel que les conditions suivantes soient remplies~: 
\begin{itemize}
\item[{\rm (i)}] $N$ est un $\hRun[\frac 1 p]$-module projectif de type fini~; 
\item[{\rm (ii)}] il existe un nombre rationnel $\varepsilon>\frac{1}{p-1}$ tel que pour tout $i\geq 1$, 
le $i$-ième invariant caractéristique de $\theta$ appartienne à 
$p^{i\varepsilon}\txi^{-i}\rS_R^i(\tOmega^1_{R/\co_K})\otimes_R\hRun$ \eqref{MH2j}. 
\end{itemize}
Alors $(N,\theta)$ est petit. 
\end{lem}

\subsection{}\label{pmh27}
On note $\bMH'^\qpp(\hRun,\txi^{-1}\tOmega^1_{R/\co_K})$ la sous-catégorie pleine de
$\bMH^\qpp(\hRun,\txi^{-1}\tOmega^1_{R/\co_K})$ \eqref{pmh4} 
formée des $\hRun$-modules de Higgs $(N,\theta)$ tels que le $\hRun[\frac 1 p]$-module $N[\frac 1 p]$ soit projectif de type fini. 
C'est une catégorie additive. On désigne par $\bMH'^\qpp_\mQ(\hRun,\txi^{-1}\tOmega^1_{R/\co_K})$
la catégorie des objets de $\bMH'^\qpp(\hRun,\txi^{-1}\tOmega^1_{R/\co_K})$ à isogénie près. Le foncteur 
\begin{equation}\label{pmh27a}
\begin{array}[t]{clcr}
\bMH'^\qpp(\hRun,\txi^{-1}\tOmega^1_{R/\co_K})&\rightarrow& \bMH^{\p}(\hRun[\frac 1 p],\txi^{-1}\tOmega^1_{R/\co_K}),\\ 
(N,\theta)&\mapsto& (N\otimes_{\mZ_p}\mQ_p, \theta\otimes_{\mZ_p}\mQ_p)
\end{array}
\end{equation}
induit alors un foncteur 
\begin{equation}\label{pmh27b}
\bMH'^\qpp_\mQ(\hRun,\txi^{-1}\tOmega^1_{R/\co_K})\rightarrow \bMH^{\p}(\hRun[\frac 1 p],\txi^{-1}\tOmega^1_{R/\co_K}).
\end{equation}

\begin{lem}\label{pmh28}
Le foncteur \eqref{pmh27b} est une équivalence de catégories.
\end{lem}
En effet, il résulte aussitôt des définitions que ce foncteur est essentiellement surjectif. 
Soient $M,N$ deux $\hRun$-modules de type fini et $\co_C$-plats. 
Alors, le morphisme canonique 
\begin{equation}\label{pmh28a}
\Hom_\hRun(M,N)\otimes_{\mZ_p}\mQ_p\rightarrow \Hom_{\hRun[\frac 1 p]}(M\otimes_{\mZ_p}\mQ_p,N\otimes_{\mZ_p}\mQ_p)
\end{equation}
est un isomorphisme. Supposons que $M$ et $N$ soient munis de $\hRun$-champs de Higgs à coefficients dans $\txi^{-1}\tOmega^1_{R/\co_K}$. 
Le morphisme canonique 
\begin{equation}\label{pmh28b}
\Hom_{\bMH}(M,N)\otimes_{\mZ_p}\mQ_p\rightarrow \Hom_{\bMH}(M\otimes_{\mZ_p}\mQ_p,N\otimes_{\mZ_p}\mQ_p)
\end{equation}
où la source (resp. le but) désigne l'ensemble des morphismes de $\hRun$-modules ($\hRun[\frac 1 p]$-modules) de Higgs est alors un isomorphisme:
l'injectivité résulte aussitôt de celle de \eqref{pmh28a} et la surjectivité de celle de \eqref{pmh28a} et du fait que $M$ et $N$ sont $\co_C$-plats et que
$\tOmega^1_{R/\co_K}$ est $R$-plat; d'où la proposition.

\subsection{}\label{pmh29}
Compte tenu de \ref{pmh24} et \ref{pmh28}, les foncteurs \eqref{pmh5f} induisent un foncteur 
\begin{equation}\label{pmh29a}
\bRep_{\hRun[\frac 1 p]}^{\p}(\Delta_\infty)\rightarrow \bMH^{\p}(\hRun[\frac 1 p],\txi^{-1}\tOmega^1_{R/\co_K}),
\end{equation}
et les foncteurs \eqref{pmh5h} induisent un foncteur 
\begin{equation}\label{pmh29b}
\bMH^{\p}(\hRun[\frac 1 p],\txi^{-1}\tOmega^1_{R/\co_K})\rightarrow \bRep_{\hRun[\frac 1 p]}^{\p}(\Delta_\infty). 
\end{equation}

\begin{prop}\label{pmh30}
Soit $(N,\theta)$ un petit $\hRun[\frac 1 p]$-module de Higgs à coefficients dans $\txi^{-1}\tOmega^1_{R/\co_K}$ \eqref{pmh25} et soit
$\varphi$ la petite $\hRun[\frac 1 p]$-représentation de $\Delta_\infty$ sur $N$ associée à $\theta$ par le foncteur \eqref{pmh29b}. 
Alors,
\begin{itemize}
\item[{\rm (i)}] On a un $\hcC^{(0+)}$-isomorphisme fonctoriel 
$\Delta$-équivariant de $\hoR$-modules de Higgs à coefficients dans $\txi^{-1}\tOmega^1_{R/\co_K}$
\begin{equation}\label{pmh30a}
N \otimes_{\hRun} \hcC^{(0+)}\stackrel{\sim}{\rightarrow} N \otimes_{\hRun} \hcC^{(0+)},
\end{equation}
où $\hcC^{(0+)}$ est muni de l'action canonique de $\Delta$ et du $\hoR$-champ de Higgs $d^{(0)}_{\hcC^{(0+)}}$ 
\eqref{taht11j}, le module $N$ de la source est muni de l'action triviale de $\Delta_{\infty}$ 
et du $\hRun$-champ de Higgs $\theta$, et le module $N$ du but est muni de l'action $\varphi$
de $\Delta_{\infty}$ et du $\hRun$-champ de Higgs nul. 
Si, de plus, la déformation $(\tX,\cM_\tX)$  est définie par la carte adéquate $((P,\gamma),(\mN,\iota),\vartheta)$ \eqref{pmh8}, 
l'isomorphisme est canonique.
\item[{\rm (ii)}] Le $\hRun[\frac 1 p]$-module de Higgs $(N,\theta)$ est soluble \eqref{repdolb10},
et on a un $\hoR[\frac 1 p]$-isomorphisme $\Delta$-équivariant fonctoriel 
\begin{equation}\label{pmh30b}
\mV(N)\stackrel{\sim}{\rightarrow} (N,\varphi)\otimes_{\hRun}\hoR,
\end{equation}
où $\mV$ est le foncteur \eqref{repdolb6b}. 
Si, de plus, la déformation $(\tX,\cM_\tX)$  est définie par la carte adéquate $((P,\gamma),(\mN,\iota),\vartheta)$ \eqref{pmh8}, 
l'isomorphisme est canonique.

\item[{\rm (iii)}]  La $\hoR[\frac 1 p]$-représentation $\mV(N)$ de $\Delta$ est petite et de Dolbeault \eqref{repdolb8},
et on a un isomorphisme fonctoriel de $\hRun[\frac 1 p]$-modules de Higgs
\begin{equation}\label{pmh30c}
\mH(\mV(N))\stackrel{\sim}{\rightarrow}(N,\theta),
\end{equation}
où $\mH$ est le foncteur \eqref{repdolb5b}. 
Si, de plus, la déformation $(\tX,\cM_\tX)$  est définie par la carte adéquate $((P,\gamma),(\mN,\iota),\vartheta)$ \eqref{pmh8}, 
l'isomorphisme est canonique.
\end{itemize}
\end{prop}

L'isomorphisme \eqref{pmh30a} résulte de \ref{pmh17}. Les autres assertions s'en déduisent compte tenu du fait que
$\ker(d^{(0)}_{\hcC^{(0+)}})=\hoR$ et $(\hcC^{(0+)})^{\Delta}=\hRun$ \eqref{taht17}.

\begin{prop}[\cite{agt} II.13.24, IV.5.3.10]\label{pmh31}
Soient $N$ un $\hRun[\frac 1 p]$-module projectif de type fini, 
$\theta$ un $\hRun[\frac 1 p]$-champ de Higgs sur $N$ à coefficients dans $\txi^{-1}\tOmega^1_{R/\co_K}$,
$M$ un $\hoR[\frac 1 p]$-module muni du champ de Higgs nul, $r$ un nombre rationnel $>0$, 
\begin{equation}\label{pmh31a}
N\otimes_{\hRun}\hcC^{(r)}\stackrel{\sim}{\rightarrow}M\otimes_{\hoR}\hcC^{(r)}
\end{equation}
un isomorphisme de $\bMC^r_p$ \eqref{repdolb3}.
Alors, $(N,\theta)$ est un petit $\hRun[\frac 1 p]$-module de Higgs à coefficients 
dans $\txi^{-1}\tOmega^1_{R/\co_K}$ \eqref{pmh25}.
\end{prop}
En fait, les propositions (\cite{agt} II.13.24 et IV.5.3.10) sont formulées dans le cas absolu \eqref{definf10}, mais la preuve s'applique
{\em mutatis mutandis} au cas relatif.

\begin{cor}\label{pmh32}
Pour qu'un $\hRun[\frac 1 p]$-module de Higgs à coefficients dans $\txi^{-1}\tOmega^1_{R/\co_K}$ 
soit soluble \eqref{repdolb10}, il faut et il suffit qu'il soit petit \eqref{pmh25}.  
\end{cor}

Cela résulte de \ref{pmh30}(ii) et \ref{pmh31}

\begin{cor}\label{pmh33}
Toute $\hoR[\frac 1 p]$-représentation de Dolbeault de $\Delta$ \eqref{repdolb8} est petite \eqref{pmh21}. 
\end{cor}

Cela résulte de \ref{repdolb14}, \ref{pmh32} et \ref{pmh30}(iii)

\begin{prop}\label{pmh34}
Dans le cas absolu \eqref{definf10}, pour qu'une $\hoR[\frac 1 p]$-représentation de $\Delta$ soit de Dolbeault \eqref{repdolb8},  
il faut et il suffit qu'elle soit petite \eqref{pmh21}. 
\end{prop}

Cela résulte de (\cite{agt} II.14.8) et (\cite{tsuji5} 13.7). 
On notera que les hypothèses de (\cite{tsuji5} §~2 et §~13) sont satisfaites compte tenu de \ref{cad1}, (\cite{ag} 4.2.2(ii)) et (\cite{tsuji4} I.5.1).

\section{Représentations de Hodge-Tate}\label{repht}

Les hypothèses et notations de \ref{taht}, \ref{repdolb} et \ref{pmh} sont en vigueur dans cette section. 

\subsection{}\label{repht1}
Reprenons les notations de \ref{pmh1}.
Soient $N$ un $\hRun[\frac 1 p]$-module projectif de type fini, 
$\theta$ un $\hRun[\frac 1 p]$-champ de Higgs {\em nilpotent} à coefficients dans $\txi^{-1}\tOmega^1_{R/\co_K}$ \eqref{MH15}. 
Il existe des endomorphismes $\hRun[\frac 1 p]$-linéaires $\theta_1,\dots,\theta_d$ de $N$
qui commutent deux à deux, tels que 
\begin{equation}\label{repht1a}
\theta=\sum_{i=1}^d\txi^{-1}\theta_i\otimes d\log(t_i).
\end{equation}
Les $\theta_i$ étant nilpotents, on peut définir une représentation $\hRun[\frac 1 p]$-linéaire de $\Delta_\infty$ sur $N$ par 
la formule
\begin{equation}\label{repht1b} 
\varphi=\exp\left(\sum_{i=1}^d\txi^{-1}\theta_i\otimes \chi_i\right),
\end{equation}
où $\chi_i$ est défini dans \eqref{pmh1h}. 
Par ailleurs, le champ de Higgs $\theta$ induit un endomorphisme $\cG[\frac 1 p]$-linéaire de $N\otimes_{\hRun}\cG$ \eqref{taht5b} 
que l'on note encore $\theta$, défini pour tous $h\in \cG$ et $x\in N$ par 
\begin{equation}\label{repht1c} 
\theta(x\otimes h)=\sum_{1\leq i\leq d}\theta_i(x)\otimes y_ih.
\end{equation}
Celui-ci étant clairement nilpotent, on peut définir son exponentielle
\begin{equation}\label{repht1d} 
\exp(\theta)\colon N\otimes_{\hRun} \cG\rightarrow N\otimes_{\hRun} \cG.
\end{equation}

\begin{prop}\label{repht2}
Conservons les hypothèses de \ref{repht1}. 
\begin{itemize}
\item[{\em (i)}] Pour tout nombre rationnel $\varepsilon>\frac{1}{p-1}$, il existe 
un sous-$\hRun$-module de type fini $N^\circ$ de $N$, qui l'engendre sur $\hRun[\frac 1 p]$, tels que l'on ait 
\begin{equation}\label{repht2a}
\theta(N^\circ)\subset p^{\varepsilon}\txi^{-1}N^\circ \otimes_R\tOmega^1_{R/\co_K}.
\end{equation}
En particulier, le $\hRun[\frac 1 p]$-module de Higgs $(N,\theta)$ est petit \eqref{pmh25}. 
\item[{\em (ii)}] La $\hRun[\frac 1 p]$-représentation $(N,\varphi)$ de $\Delta_\infty$  \eqref{repht1b} est l'image du $\hRun[\frac 1 p]$-module de Higgs $(N,\theta)$
par le foncteur \eqref{pmh29b}.  
\item[{\em (iii)}] Soient $r$ un nombre rationnel $\geq 0$, $\varepsilon$ nombre rationnel  $>r+\frac{1}{p-1}$, 
$N^\circ$ un sous-$\hRun$-module de type fini de $N$ vérifiant les conclusions de {\rm (i)}, $\hfS^{(r)}$ l'anneau défini dans \ref{pmh18}, 
$\exp_r(\theta)$ l'endomorphisme $\hfS^{(r)}$-linéaire de 
$N^\circ\otimes_{\hRun}\hfS^{(r)}$ défini dans \eqref{pmh13c}. 
Alors, les endomorphismes $\hcG$-linéaires sur $N\otimes_{\hRun}\hcG$ induits par $\exp_r(\theta)$ et $\exp(\theta)$ \eqref{repht1d} coïncident. 
\item[{\em (iv)}] L'endomorphisme $\exp(\theta)$ \eqref{repht1d} est un isomorphisme $\Delta$-équivariant 
de $\cG$-modules à connexion intégrable relativement à l'extension $\cG/\hoR$ \eqref{repdolb20}, 
où le module $N$ de la source est muni de l'action triviale de $\Delta$ 
et du $\hRun$-champ de Higgs $\theta$, le module $N$ du but est muni de l'action $\varphi$
de $\Delta$ et du $\hRun$-champ de Higgs nul et $\cG$ est muni de l'action $\varphi_0$ de $\Delta$ \eqref{pmh10c}.
\end{itemize}
\end{prop}

(i) Comme l'anneau $\hRun$ est normal d'après (\cite{agt} II.6.15) et donc réduit, pour tout $n\geq 1$, le $n$-ième invariant caractéristique $\lambda_n(\theta)$
de $\theta$ est nul en vertu de \ref{MH17}. La proposition résulte alors de (\cite{agt} II.13.7).  

(ii) Cela résulte aussitôt de (i) et \ref{pmh5}.

(iii) Cela résulte aussitôt des définitions (cf. \ref{pmh13}). 

(iv) Cela résulte de (ii), (iii) et \ref{pmh14} puisque le morphisme canonique $N\otimes_{\hRun}\cG\rightarrow N\otimes_{\hRun}\hcG$
est injectif. 

\begin{cor}\label{repht3}
Sous les hypothèses \ref{repht1}, il existe un isomorphisme $\Delta$-équivariant 
de $\cC$-modules à connexion intégrable relativement à l'extension $\cC/\hoR$ \eqref{repdolb20}
\begin{equation}\label{repht3a} 
N\otimes_{\hRun} \cC\rightarrow N\otimes_{\hRun} \cC,
\end{equation}
où $\cC$ est muni de l'action canonique de $\Delta$ \eqref{taht7}, le module $N$ de la source est muni de l'action triviale de $\Delta$ 
et du $\hRun$-champ de Higgs $\theta$ et le module $N$ du but est muni de l'action $\varphi$
de $\Delta$ \eqref{repht1b} et du $\hRun$-champ de Higgs nul.
\end{cor}

On peut supposer que la déformation $(\tX,\cM_{\tX})$ est définie par la carte adéquate $((P,\gamma),(\mN,\iota),\vartheta)$ \eqref{pmh8}, 
auquel cas la proposition résulte de \ref{repht2}(iv) et \eqref{pmh12g}.

\begin{cor}\label{repht4}
Soient $N$ un $\hRun[\frac 1 p]$-module projectif de type fini, $\theta$ un $\hRun[\frac 1 p]$-champ de Higgs à coefficients dans 
$\txi^{-1}\tOmega^1_{R/\co_K}$. Alors, les conditions suivantes sont équivalentes:
\begin{itemize}
\item[{\rm (i)}] Le champ de Higgs $\theta$ est nilpotent \eqref{MH15}.
\item[{\rm (ii)}] Il existe une $\hoR[\frac 1 p]$-représentation $M$ de $\Delta$ et 
un isomorphisme $\cC$-linéaire et $\Delta$-équivariant de $\hoR[\frac 1 p]$-modules de Higgs à coefficients dans $\txi^{-1}\tOmega^1_{R/\co_K}$
\begin{equation}\label{repht4a}
N \otimes_{\hRun}\cC\stackrel{\sim}{\rightarrow}  M\otimes_{\hoR}\cC,
\end{equation}
où $\cC$ est muni de l'action canonique de $\Delta$ et du champ de Higgs $d_{\cC}$, $N$ est muni de l'action triviale de $\Delta$ et $M$ est muni du champ de Higgs nul.
\end{itemize}
De plus, dans ce cas, le $\hoR[\frac 1 p]$-module $M$ est projectif de type fini, le $\hRun[\frac 1 p]$-module de Higgs $(N,\theta)$ est soluble et on a un isomorphisme de 
$\hoR[\frac 1 p]$-représentations de $\Delta$,
\begin{equation}\label{repht4b}
\mV(N,\theta)\stackrel{\sim}{\rightarrow}  M.
\end{equation}
\end{cor}

L'implication (i)$\Rightarrow$(ii) résulte de \ref{repht3}. Montrons l'implication inverse. Supposons la condition (ii) satisfaite. 
Toute section $\cF\rightarrow \hoR$ de l'extension \eqref{taht5e}, en tant que suite exacte de $\hoR$-modules sans actions de $\Delta$, 
définit une rétraction de la $\hoR$-algèbre $\cC$. 
L'isomorphisme \eqref{repht4a} implique alors que le $\hoR[\frac 1 p]$-module $M$ est projectif de type fini. 
Le $\hRun[\frac 1 p]$-module $N$ étant projectif de type fini, il existe un entier $n\geq 1$ tel que l'isomorphisme \eqref{repht4a} 
envoie $N$ dans $M\otimes_{\hoR} \rS^n_{\hoR}(\cF)$. On en déduit aussitôt que $\theta$ est nilpotent.
L'isomorphisme \eqref{repht4a} est en fait un isomorphisme de $\cC$-modules à connexion intégrable 
relativement à l'extension $\cC/\hoR$ \eqref{repdolb20}. 
Pour tout nombre rationnel $r\geq 0$, le morphisme canonique 
$\alpha^{0,r}\colon \cC^{(r)}\rightarrow \cC$ \eqref{taht10} induit un isomorphisme de $\hoR$-algèbres
\begin{equation}\label{repht4c}
\cC^{(r)}[\frac 1 p]\stackrel{\sim}{\rightarrow} \cC[\frac 1 p].
\end{equation}
Compte tenu de \eqref{taht11h} et \eqref{repdolb3k}, l'isomorphisme \eqref{repht4a} induit alors un isomorphisme 
$\hcC^{(r)}$-linéaire et $\Delta$-équivariant de $\hoR[\frac 1 p]$-modules de Higgs à coefficients dans $\txi^{-1}\tOmega^1_{R/\co_K}$
\begin{equation}\label{repht4d}
N \otimes_{\hRun}\hcC^{(r)}\stackrel{\sim}{\rightarrow}  M\otimes_{\hoR}\hcC^{(r)},
\end{equation}
où $\hcC^{(r)}$ est muni de l'action canonique de $\Delta$ et du champ de Higgs $p^rd_{\hcC^{(r)}}$, $N$ est muni de l'action triviale de $\Delta$
et $M$ est muni du champ de Higgs nul.
Il résulte alors de \ref{repdolb17}  que le $\hRun[\frac 1 p]$-module de Higgs $(N,\theta)$ est soluble 
et qu'on a un isomorphisme de $\hoR[\frac 1 p]$-représentations de $\Delta$
\begin{equation}\label{repht4e}
\mV(N,\theta)\stackrel{\sim}{\rightarrow}  M.
\end{equation}

\begin{cor}\label{repht5}
Pour toute $\hoR[\frac 1 p]$-représentation $M$ de $\Delta$, les conditions suivantes sont équivalentes:
\begin{itemize}
\item[{\rm (i)}] La $\hoR[\frac 1 p]$-représentation $M$ de $\Delta$ est de Dolbeault et le $\hRun[\frac 1 p]$-module de Higgs 
à coefficients dans $\txi^{-1}\tOmega^1_{R/\co_K}$ associé, $\mH(M)$, est nilpotent \eqref{MH15}.
\item[{\rm (ii)}] Il existe un $\hRun[\frac 1 p]$-module projectif de type fini $N$, un $\hRun[\frac 1 p]$-champ de
Higgs $\theta$ sur $N$ à coefficients dans $\txi^{-1}\tOmega^1_{R/\co_K}$ et un isomorphisme $\cC$-linéaire
et $\Delta$-équivariant de $\hoR[\frac 1 p]$-modules de Higgs à coefficients dans $\txi^{-1}\tOmega^1_{R/\co_K}$
\begin{equation}\label{repht5a}
N \otimes_{\hRun}\cC\stackrel{\sim}{\rightarrow}  M\otimes_{\hoR}\cC,
\end{equation}
où $\cC$ est muni de l'action canonique de $\Delta$ et du champ de Higgs $d_{\cC}$,
$N$ est muni de l'action triviale de $\Delta$ et $M$ est muni du champ de Higgs nul.
\end{itemize}
De plus, dans ce cas, on a un isomorphisme de $\hRun[\frac 1 p]$-modules de Higgs
à coefficients dans $\txi^{-1}\tOmega^1_{R/\co_K}$
\begin{equation}\label{repht5b}
\mH(M)\stackrel{\sim}{\rightarrow}  (N,\theta).
\end{equation}
\end{cor}

L'implication (i)$\Rightarrow$(ii) résulte de \ref{repdolb14} et \ref{repht4} appliqué à $\mH(M)$. Montrons l'implication inverse. 
Supposons la condition (ii) satisfaite. Pour tout nombre rationnel $r\geq 0$, l'isomorphisme \eqref{repht5a} induit un isomorphisme 
$\hcC^{(r)}$-linéaire et $\Delta$-équivariant de $\hoR[\frac 1 p]$-modules de Higgs à coefficients dans $\txi^{-1}\tOmega^1_{R/\co_K}$
\begin{equation}\label{repht5c}
N \otimes_{\hRun}\hcC^{(r)}\stackrel{\sim}{\rightarrow}  M\otimes_{\hoR}\hcC^{(r)},
\end{equation}
où $\hcC^{(r)}$ est muni de l'action canonique de $\Delta$ et du champ de Higgs $p^rd_{\hcC^{(r)}}$,
$N$ est muni de l'action triviale de $\Delta$ et $M$ est muni du champ de Higgs nul (cf. la preuve de \ref{repht4}).
Il résulte alors de \ref{repdolb9} que $M$ est de Dolbeault et qu'on a un isomorphisme de $\hRun[\frac 1 p]$-modules de Higgs
à coefficients dans $\txi^{-1}\tOmega^1_{R/\co_K}$
\begin{equation}\label{repht5d}
\mH(M)\stackrel{\sim}{\rightarrow}  (N,\theta).
\end{equation}
Le $\hRun[\frac 1 p]$-module $N$ (resp. $\hoR[\frac 1 p]$-module $M$) étant projectif de type fini, il existe un entier $n\geq 1$ tel que l'isomorphisme \eqref{repht5a} envoie $N$ dans 
$M\otimes_{\hoR} \rS^n_{\hoR}(\cF)$. On en déduit aussitôt que $\theta$ est nilpotent, d'où la condition~(i).

\begin{defi}\label{repdolb23}
On dit qu'une $\hoR[\frac 1 p]$-représentation $M$ de $\Delta$ est {\em de Hodge-Tate} si elle vérifie les conditions équivalentes de \ref{repht5}.
\end{defi}

Cette notion ne dépend pas du choix de la $(\tS,\cM_\tS)$-déformation $(\tX,\cM_{\tX})$ d'après \ref{repdolb7}(i),
ni du cas absolu ou relatif \eqref{definf10} en vertu de \ref{taht9}. 
Lorsqu'il y a lieu de préciser, on dira que la $\hoR[\frac 1 p]$-représentation $M$ de $\Delta$ est {\em géométriquement de Hodge-Tate} pour éviter toute confusion
avec la notion de représentation de Hodge-Tate considérée dans (\cite{hyodo1} 2.1) (cf. \cite{agt} II.15).

\begin{prop}\label{repdolb24}
Les foncteurs $\mH$  \eqref{repdolb5b} et $\mV$ \eqref{repdolb6b} induisent des équivalences de catégories quasi-inverses l'une de l'autre,
entre la catégorie des $\hoR[\frac 1 p]$-représentations de Hodge-Tate de $\Delta$ et celle des 
$\hRun[\frac 1 p]$-modules de Higgs nilpotents à coefficients dans $\txi^{-1}\tOmega^1_{R/\co_K}$ dont le $\hRun[\frac 1 p]$-module sous-jacent est projectif de type fini.
\end{prop}

Cela résulte de \ref{repdolb14}, \ref{repht4} et \ref{repht5}.

\chapter[\'Etude globale]{\texorpdfstring{Correspondance de Simpson $p$-adique et modules de Hodge-Tate. \'Etude globale}
{Correspondance de Simpson p-adique et modules de Hodge-Tate. \'Etude globale}}

\section{Hypothèses et notations}

\subsection{}\label{defing1}
Dans ce chapitre, $K$ désigne un corps de valuation discrète complet de 
caractéristique $0$, à corps résiduel {\em algébriquement clos} $k$ de caractéristique $p>0$,  
$\co_K$ l'anneau de valuation de $K$, $W$ l'anneau des vecteurs de Witt à coefficients dans $k$ relatif à $p$, 
$\oK$ une clôture algébrique de $K$, $\co_\oK$ la clôture intégrale de $\co_K$ dans $\oK$,
$\fm_\oK$ l'idéal maximal de $\co_\oK$ et $G_K$ le groupe de Galois de $\oK$ sur $K$.
On note $\co_C$ le séparé complété $p$-adique de $\co_\oK$, $\fm_C$ son idéal maximal,
$C$ son corps des fractions et $v$ sa valuation, normalisée par $v(p)=1$. 
On désigne par $\mZ_p(1)$ le $\mZ[G_K]$-module 
\begin{equation}\label{defing1a}
\mZ_p(1)=\underset{\underset{n\geq 1}{\longleftarrow}}{\lim}\ \mu_{p^n}(\co_\oK),
\end{equation}  
où $\mu_{p^n}(\co_\oK)$ désigne le sous-groupe des racines $p^n$-ièmes de l'unité dans $\co_\oK$. 
Pour tout $\mZ_p[G_K]$-module $M$ et tout entier $n$, on pose $M(n)=M\otimes_{\mZ_p}\mZ_p(1)^{\otimes n}$.

On pose $S=\Spec(\co_K)$, $\oS=\Spec(\co_\oK)$ et $\coS=\Spec(\co_C)$. 
On note $s$ (resp.  $\eta$, resp. $\oeta$) le point fermé de $S$ (resp.  générique de $S$, resp. générique de $\oS$).
Pour tout entier $n\geq 1$, on pose $S_n=\Spec(\co_K/p^n\co_K)$. Pour tout $S$-schéma $X$, on pose 
\begin{equation}\label{defing1c}
\oX=X\times_S\oS,  \ \ \ \coX=X\times_S\coS \ \ \ {\rm et}\ \ \  X_n=X\times_SS_n.
\end{equation} 

On munit $S$ de la structure logarithmique $\cM_S$ définie par son point fermé, 
et $\oS$ et $\coS$ des structures logarithmiques $\cM_\oS$ et $\cM_\coS$ images inverses de $\cM_S$. 

\subsection{}\label{defing2}
Comme $\co_\oK$ est un anneau de valuation non discrète de hauteur $1$, 
il est loisible de développer la $\alpha$-algèbre (ou presque-algèbre) sur cet anneau (\cite{ag} 2.10.1) (cf. \cite{ag} 2.6-2.10).   
On choisit un système compatible $(\beta_n)_{n>0}$ 
de racines $n$-ièmes de $p$ dans $\co_\oK$. Pour tout nombre rationnel $\varepsilon>0$, 
on pose $p^\varepsilon=(\beta_n)^{\varepsilon n}$, où $n$ est un entier $>0$ tel que $\varepsilon n$ soit entier.

\subsection{}\label{defing11}
Soit $f\colon (X,\cM_X)\rightarrow (S,\cM_S)$ un morphisme {\em adéquat} de schémas logarithmiques (\cite{agt} III.4.7). 
On désigne par $X^\circ$ le sous-schéma ouvert maximal de $X$
où la structure logarithmique $\cM_X$ est triviale~; c'est un sous-schéma ouvert de $X_\eta$.
On note $j\colon X^\circ\rightarrow X$ l'injection canonique. Pour tout $X$-schéma $U$, on pose  
\begin{equation}\label{defing11a}
U^\circ=U\times_XX^\circ.
\end{equation} 
On note $\hbar\colon \oX\rightarrow X$ et $h\colon \oX^\circ\rightarrow X$ les morphismes canoniques \eqref{defing1c}, de sorte que 
l'on a $h=\hbar\circ j_\oX$. Pour alléger les notations, on pose
\begin{equation}\label{defing11b}
\tOmega^1_{X/S}=\Omega^1_{(X,\cM_X)/(S,\cM_S)},
\end{equation}
que l'on considère comme un faisceau de $X_\zar$ ou $X_\et$, selon le contexte (cf. \ref{notconv12}). 

\subsection{}\label{ahttf1}
Pour tout entier $n\geq 1$, on note $a\colon X_s\rightarrow X$, $a_n\colon X_s\rightarrow X_n$, 
$\iota_n\colon X_n\rightarrow X$ et $\oiota_n\colon \oX_n\rightarrow \oX$ les injections canoniques \eqref{defing1c}. 
Comme $k$ est algébriquement clos, il existe un unique $S$-morphisme $s\rightarrow \oS$. 
Celui-ci induit des immersions fermées $\oa\colon X_s\rightarrow \oX$ et $\oa_n\colon X_s\rightarrow \oX_n$
qui relèvent $a$ et $a_n$, respectivement. 
\begin{equation}\label{ahttf1a}
\xymatrix{
{X_s}\ar[r]_{\oa_n}\ar@{=}[d]\ar@/^1pc/[rr]^{\oa}&{\oX_n}\ar[r]_{\oiota_n}\ar[d]^{\hbar_n}&\oX\ar[d]^\hbar\\
{X_s}\ar[r]^{a_n}\ar@/_1pc/[rr]_{a}&{X_n}\ar[r]^{\iota_n}&X}
\end{equation}
Comme $\hbar$ est entier et que $\hbar_n$ est un homéomorphisme
universel, pour tout faisceau $\cF$ de $\oX_\et$, le morphisme de changement de base 
\begin{equation}
a^*(\hbar_*(\cF))\rightarrow \oa^*(\cF)
\end{equation}
est un isomorphisme (\cite{sga4} VIII 5.6). Par ailleurs, $\oa_n$ étant un homéomorphisme universel, on peut considérer $\co_{\oX_n}$
comme un faisceau de $X_{s,\zar}$ ou $X_{s,\et}$, selon le contexte (cf. \ref{notconv12}). 

On pose \eqref{defing11b}
\begin{equation}\label{ahttf1b}
\tOmega^1_{\oX_n/\oS_n}=\tOmega^1_{X/S}\otimes_{\co_X}\co_{\oX_n},
\end{equation}
que l'on considère aussi comme un faisceau de $X_{s,\zar}$ ou $X_{s,\et}$, selon le contexte.

\subsection{}\label{defing12}
Reprenons les notations de \ref{definf10}. On pose 
\begin{equation}\label{ahttf1c}
\txi^{-1}\tOmega^1_{\oX_n/\oS_n}=\txi^{-1}\co_C\otimes_{\co_C}\tOmega^1_{\oX_n/\oS_n}.
\end{equation}
On munit  $\coX=X\times_S\coS$ \eqref{defing1c} de la structure logarithmique $\cM_\coX$ image inverse de $\cM_X$. 
On a alors un isomorphisme canonique
\begin{equation}\label{defing12a}
(\coX,\cM_{\coX})\stackrel{\sim}{\rightarrow}(X,\cM_X)\times_{(S,\cM_S)}(\coS,\cM_\coS),
\end{equation} 
le produit fibré étant indifféremment pris dans la catégorie des schémas logarithmiques ou 
dans celle des schémas logarithmiques fins. 

{\em On suppose dans ce chapitre qu'il existe une $(\tS,\cM_{\tS})$-déformation lisse $(\tX,\cM_\tX)$ de $(\coX,\cM_{\coX})$ que l'on fixe}; 
autrement dit, un morphisme lisse de schémas logarithmiques fins $(\tX,\cM_\tX)\rightarrow (\tS, \cM_{\tS})$ et un $(\coS,\cM_{\coS})$-isomorphisme 
\begin{equation}\label{defing12b}
(\coX,\cM_{\coX})\stackrel{\sim}{\rightarrow}
(\tX,\cM_\tX)\times_{(\tS, \cM_{\tS})}(\coS,\cM_{\coS}).
\end{equation}

\begin{rema}\label{defing13}
Dans le cas relatif \eqref{definf10}, il existe une $(\tS,\cM_{\tS})$-déformation lisse canonique de $(\coX,\cM_{\coX})$, à savoir
\begin{equation}\label{defing13a}
(X,\cM_{X})\times_{(S,\cM_S)}(\tS,\cM_\tS)
\end{equation}
où l'on considère $(\tS,\cM_\tS)$ comme un schéma logarithmique au-dessus de $(S,\cM_S)$ via $\pr_1$ \eqref{definf7f}, 
le produit fibré étant indifféremment pris dans la catégorie des schémas logarithmiques ou dans celle des schémas logarithmiques fins. 
\end{rema}

\section{Modules de Higgs}\label{nrmh}

\subsection{}\label{definf18}
On pose $\cS=\Spf(\co_C)$ et on désigne par $\fX$ le schéma formel complété $p$-adique de $\oX$ \eqref{defing11}. 
C'est un $\cS$-schéma formel de présentation finie (\cite{egr1} 2.3.15). Il est donc idyllique (\cite{egr1} 2.6.13). 
On désigne par $\bMod(\co_{\fX})$ (resp. $\bMod(\co_{\fX}[\frac 1 p])$) la catégorie des $\co_{\fX}$-modules 
(resp. $\co_{\fX}[\frac 1 p]$-modules) de $X_{s,\zar}$, par 
$\bMod^\coh(\co_{\fX})$ (resp. $\bMod^\coh(\co_{\fX}[\frac 1 p])$) 
la sous-catégorie pleine formée des $\co_{\fX}$-modules (resp. $\co_{\fX}[\frac 1 p]$-modules) cohérents 
et par $\bMod^{\coh}_\mQ(\co_{\fX})$ la catégorie des $\co_\fX$-modules cohérents à isogénie près (\cite{agt} III.6.1.1). 
D'après (\cite{agt} III.6.16), le foncteur canonique
\begin{equation}\label{aspglob1c}
\bMod^{\coh}(\co_{\fX})\rightarrow \bMod^{\coh}(\co_{\fX}[\frac 1 p]),\ \ \ \cN\mapsto \cN_{\mQ_p}
\end{equation}
induit une équivalence de catégories abéliennes
\begin{equation}\label{aspglob1d}
\bMod^{\coh}_\mQ(\co_{\fX})\stackrel{\sim}{\rightarrow} \bMod^{\coh}(\co_{\fX}[\frac 1 p]). 
\end{equation}

\subsection{}\label{definf19}
Avec les conventions et notations de \ref{definf10} et \eqref{defing11b}, 
on désigne par $\txi^{-1}\tOmega^1_{\fX/\cS}$ le complété $p$-adique du $\co_\coX$-module \eqref{defing1c} (\cite{egr1} 2.5.1)
\begin{equation}\label{definf19a}
\txi^{-1}\tOmega^1_{\coX/\coS}=\txi^{-1}\tOmega^1_{X/S}\otimes_{\co_X}\co_{\coX}.
\end{equation} 
On note $\bMH(\co_\fX,\txi^{-1}\tOmega^1_{\fX/\cS})$ la catégorie des $\co_\fX$-modules de Higgs à coefficients dans 
$\txi^{-1}\tOmega^1_{\fX/\cS}$  \eqref{MH1}. On dit qu'un tel module de Higgs est {\em cohérent} si le $\co_\fX$-module sous-jacent est cohérent.
On désigne par $\bMH^\coh(\co_\fX, \txi^{-1}\tOmega^1_{\fX/\cS})$ la sous-catégorie pleine de $\bMH(\co_\fX,\txi^{-1}\tOmega^1_{\fX/\cS})$
formée des modules de Higgs cohérents.

On sous-entend par {\em $\co_\fX[\frac 1 p]$-module de Higgs à coefficients dans $\txi^{-1}\tOmega^1_{\fX/\cS}$}, 
un $\co_\fX[\frac 1 p]$-module de Higgs à coefficients dans $\txi^{-1}\tOmega^1_{\fX/\cS}[\frac 1 p]$.
On dit qu'un tel module de Higgs est {\em cohérent} si le $\co_\fX[\frac 1 p]$-module sous-jacent est cohérent. 
On désigne par $\bMH(\co_\fX[\frac 1 p], \txi^{-1}\tOmega^1_{\fX/\cS})$ la catégorie abélienne des 
$\co_\fX[\frac 1 p]$-modules de Higgs à coefficients dans $\txi^{-1}\tOmega^1_{\fX/\cS}$
et par $\bMH^\coh(\co_\fX[\frac 1 p], \txi^{-1}\tOmega^1_{\fX/\cS})$ la sous-catégorie pleine formée des modules de Higgs cohérents; c'est une sous-catégorie abélienne.  

Dans la suite de ce chapitre, on omettra le champ de Higgs de la notation d'un module de Higgs lorsqu'on n'en a pas explicitement besoin.

On désigne par $\bIH(\co_\fX,\txi^{-1}\tOmega^1_{\fX/\cS})$ la catégorie des $\co_\fX$-isogénies  
de Higgs à coefficients dans $\txi^{-1}\tOmega^1_{\fX/\cS}$ \eqref{indsh23} et par
$\bIH^\coh(\co_\fX,\txi^{-1}\tOmega^1_{\fX/\cS})$ la sous-catégorie pleine 
formée des quadruplets $(\cM,\cN,u,\theta)$ tels que $\cM$ et $\cN$ soient des $\co_\fX$-modules cohérents. Ce sont des catégories additives. 
On note $\bIH_\mQ(\co_\fX,\txi^{-1}\tOmega^1_{\fX/\cS})$ (resp. $\bIH^\coh_\mQ(\co_\fX,\txi^{-1}\tOmega^1_{\fX/\cS})$)
la catégorie des objets de $\bIH(\co_\fX,\txi^{-1}\tOmega^1_{\fX/\cS})$ 
(resp. $\bIH^\coh(\co_\fX,\txi^{-1}\tOmega^1_{\fX/\cS})$) à isogénie près (\cite{agt} III.6.1.1). Le foncteur 
\begin{equation}\label{definf19b}
\begin{array}[t]{clcr}
\bIH(\co_\fX,\txi^{-1}\tOmega^1_{\fX/\cS})&\rightarrow& \bMH(\co_\fX[\frac 1 p], \txi^{-1}\tOmega^1_{\fX/\cS})\\
(\cM,\cN,u,\theta)&\mapsto& (\cM_{\mQ_p}, (\id \otimes u_{\mQ_p}^{-1})\circ\theta_{\mQ_p})
\end{array}
\end{equation}
induit un foncteur 
\begin{equation}\label{definf19c}
\bIH_\mQ(\co_\fX,\txi^{-1}\tOmega^1_{\fX/\cS})\rightarrow \bMH(\co_\fX[\frac 1 p], \txi^{-1}\tOmega^1_{\fX/\cS}).
\end{equation}
D'après (\cite{agt} III.6.21), celui-ci induit une équivalence de catégories 
\begin{equation}\label{definf19d}
\bIH^\coh_\mQ(\co_\fX,\txi^{-1}\tOmega^1_{\fX/\cS})\stackrel{\sim}{\rightarrow} 
\bMH^\coh(\co_\fX[\frac 1 p], \txi^{-1}\tOmega^1_{\fX/\cS}).
\end{equation}

\begin{defi}\label{definf20}
On appelle {\em $\co_\fX[\frac 1 p]$-fibré de Higgs à coefficients dans $\txi^{-1}\tOmega^1_{\fX/\cS}$} 
tout $\co_\fX[\frac 1 p]$-module de Higgs à coefficients dans $\txi^{-1}\tOmega^1_{\fX/\cS}$ 
dont le $\co_\fX[\frac 1 p]$-module sous-jacent est localement projectif de type fini \eqref{notconv14}. 
\end{defi}

\begin{lem}\label{definf21}
Soit $(\cM,\theta)$ un $\co_\fX[\frac 1 p]$-module de Higgs cohérent  à coefficients dans $\txi^{-1}\tOmega^1_{\fX/\cS}$. 
On pose $\cT=\cHom_{\co_\fX}(\tOmega^1_{\fX/\cS},\txi\co_\fX)$
et on munit $\cM$ de la structure de $\rS_{\co_\fX}(\cT)$-module induite par $\theta$ \eqref{MH2g}.
On désigne par $\cJ$ le noyau de l'augmentation canonique $\rS_{\co_\fX}(\cT)\rightarrow \co_\fX$. 
Pour tout ouvert $U$ de $\fX$ et toute section $d\in \Gamma(U,\cT)$, on note $\theta_d$ l'endomorphisme de $\cM|U$ déduit de $\theta$ et $d$. 
Alors, les conditions suivantes sont équivalentes:
\begin{itemize}
\item[{\rm (i)}] Il existe un entier $n\geq 1$ tel que $\cJ^n\cM=0$. 
\item[{\rm (ii)}] Il existe un entier $n\geq 1$ tel que pour tout ouvert $U$ de $\fX$ et pour toutes sections $s\in \Gamma(U,\cM)$ et 
$d_1,\dots,d_n\in \Gamma(U,\cT)$, on ait $\theta_{d_n}\circ \dots \circ \theta_{d_1}(s)=0$.
\item[{\rm (iii)}] Il existe une filtration décroissante finie $(\cM_i)_{0\leq i\leq n}$ de $\cM$ par des sous-$\co_\fX[\frac 1 p]$-modules cohérents 
telle que $\cM_0=\cM$, $\cM_n=0$ et que pour tout $0\leq i\leq n-1$, on ait
\begin{equation}\label{definf21a}
\theta(\cM_i)\subset \txi^{-1}\tOmega^1_{\fX/\cS}\otimes_{\co_\fX}\cM_{i+1}.
\end{equation}
\end{itemize} 
\end{lem}

En effet, on a clairement (iii)$\Rightarrow$(ii)$\Leftrightarrow$(i). Supposons la condition (i) remplie. Pour tout $0\leq i\leq n$, on désigne par $\cM_i$ 
le plus grand sous-$\rS_{\co_\fX}(\cT)$-module de $\cM$ annulé par $\cJ^{n-i}$. On a donc $\cM_0=\cM$, $\cM_n=0$ et pour tout
$0\leq i\leq n-1$, $\cJ\cM_i\subset \cM_{i+1}$. Comme le $\co_\fX$-module $\cT$ est de type fini, les $\co_\fX[\frac 1 p]$-modules $\cM_i$ sont cohérents 
(\cite{ega1n} 0.5.3.4). La filtration de $\cM$ ainsi définie vérifie donc la condition (iii), d'où la proposition.    

\begin{defi}\label{definf22}
On dit qu'un $\co_\fX[\frac 1 p]$-module de Higgs cohérent $(\cM,\theta)$ à coefficients dans $\txi^{-1}\tOmega^1_{\fX/\cS}$ 
est {\em nilpotent} s'il vérifie les conditions équivalentes de \ref{definf21}. On dit alors aussi que le champ de Higgs $\theta$ est {\em nilpotent}.
\end{defi}

Pour qu'un $\co_\fX[\frac 1 p]$-module de Higgs $(\cM,\theta)$ à coefficients dans $\txi^{-1}\tOmega^1_{\fX/\cS}$ soit nilpotent, il faut et il suffit qu'il soit quasi-nilpotent 
et que $(\cM,\theta)$ admette une filtration quasi-nilpotente formée de sous-$\co_\fX[\frac 1 p]$-modules cohérents \eqref{MH151}.

\begin{rema}\label{definf23}
Conservons les hypothèses de \ref{definf21}, supposons, de plus, $\fX$ affine d'anneau $A$. Posons $T=\cT(\fX)$,
$M=\cM(\fX)$ qui est un $A[\frac 1 p]$-module de type fini et notons $\uptheta$ le $A[\frac 1 p]$-champ de Higgs sur $M$ à coefficients dans 
$\txi^{-1}\tOmega^1_{\fX/\cS}(\fX)$ induit par $\theta$ (\cite{egr1} 2.7.2.3 et 2.10.24). 
Alors, les propriétés suivantes sont équivalentes:
\begin{itemize}
\item[(a)] Le $\co_\fX[\frac 1 p]$-module de Higgs $(\cM,\theta)$ est nilpotent dans le sens de \ref{definf22}.
\item[(b)] Il existe une filtration décroissante finie $(M_i)_{0\leq i\leq n}$ de $M$ par des sous-$A[\frac 1 p]$-modules 
telle que $M_0=M$, $M_n=0$ et que pour tout $0\leq i\leq n-1$, on ait
\begin{equation}
\uptheta(M_i)\subset \txi^{-1}\tOmega^1_{\fX/\cS}(\fX)\otimes_AM_{i+1}.
\end{equation}
\item[(c)] Pour tout $s\in M$, il existe un entier $n\geq 1$ tel que pour toutes sections $d_1,\dots,d_n\in T$, $\theta_{d_n}\circ \dots \circ \theta_{d_1}(s)=0$.
\item[(d)] Le $A[\frac 1 p]$-module de Higgs $(M,\uptheta)$ est nilpotent dans le sens de \ref{MH15}.
\end{itemize}
En effet, l'anneau $A[\frac 1 p]$ étant noethérien (\cite{egr1} 1.10.12(i)), les conditions (b) et \ref{definf21}(iii) sont équivalentes en vertu de (\cite{egr1} 2.7.2.3 et 2.10.24).  
D'autre part, les conditions (b), (c) et (d) sont équivalentes d'après \ref{MH140}. 
\end{rema}

\begin{rema}\label{definf230}
Sous les hypothèses de \ref{definf23}, pour que le $\co_\fX[\frac 1 p]$-module de Higgs $(\cM,\theta)$ soit nilpotent, 
il faut et il suffit que la condition suivante soit remplie:
\begin{itemize}
\item[{\rm (ii')}] Pour tout ouvert $U$ de $\fX$, il existe un entier $n\geq 1$ tel que pour toutes sections $s\in \Gamma(U,\cM)$ et 
$d_1,\dots,d_n\in \Gamma(U,\cT)$, on ait $\theta_{d_n}\circ \dots \circ \theta_{d_1}(s)=0$.
\end{itemize}
De plus, il suffit qu'elle le soit pour des ouverts affines $U$ couvrant $\fX$. 
\end{rema}

\begin{lem}\label{definf231}
Soient $(\cM,\theta)$ un $\co_\fX[\frac 1 p]$-module de Higgs à coefficients dans $\txi^{-1}\tOmega^1_{\fX/\cS}$,
$\cN$ (resp. $\cP$) un sous-objet (resp. quotient) de $(\cM,\theta)$ dans $\bMH^\coh(\co_\fX[\frac 1 p], \txi^{-1}\tOmega^1_{\fX/\cS})$, 
$(\rF^i\cM)_{i\in \mZ}$ une filtration décroissante de $(\cM,\theta)$ dans $\bMH^\coh(\co_\fX[\frac 1 p], \txi^{-1}\tOmega^1_{\fX/\cS})$ \eqref{MH150}
telle que $\rF^n \cM=0$, $\rF^m \cM=\cM$ pour deux entiers $m\leq n$. Alors,
\begin{itemize}
\item[{\rm (i)}] Si $(\cM,\theta)$ est nilpotent, il en est de même de $\cN$ et $\cP$. 
\item[{\rm (ii)}] Pour que $(\cM,\theta)$ soit nilpotent, il faut et il suffit que le gradué $\Gr_\rF^\bullet \cM$ de la filtration $\rF^\bullet \cM$ le soit. 
\end{itemize}
\end{lem} 

Il suffit de calquer la preuve de \ref{MH141}.

\begin{defi}\label{definf24}
Soient $\varepsilon$ un nombre rationnel $>0$, $(\cN,\theta)$ un $\co_\fX$-module de Higgs à coefficients dans 
$\txi^{-1}\tOmega^1_{\fX/\cS}$.
\begin{itemize}
\item[(i)] On dit que $(\cN,\theta)$ est {\em $\varepsilon$-quasi-petit} si le $\co_\fX$-module $\cN$ est cohérent et si 
$\theta$ est un multiple de  $p^{\varepsilon}$ en tant que section de $\txi^{-1}\cEnd_{\co_\fX}(\cN)\otimes_R\tOmega^1_{\fX/\cS}$
\eqref{MH2g}. On dit alors aussi que le $\co_\fX$-champ de Higgs $\theta$ est {\em $\varepsilon$-quasi-petit}. 
\item[{\rm (ii)}] On dit que $(\cN,\theta)$ est {\em quasi-petit} s'il est $\varepsilon'$-quasi-petit pour un nombre rationnel $\varepsilon'>\frac{1}{p-1}$. 
On dit alors aussi que le $\co_\fX$-champ de Higgs $\theta$ est {\em quasi-petit}. 
\end{itemize}
\end{defi}

On désigne par $\bMH^{\varepsilon\trqpp}(\co_\fX,\txi^{-1}\tOmega^1_{\fX/\cS})$ (resp. $\bMH^{\qpp}(\co_\fX,\txi^{-1}\tOmega^1_{\fX/\cS})$)
la sous-catégorie pleine de $\bMH(\co_\fX,\txi^{-1}\tOmega^1_{\fX/\cS})$ 
formée des $\co_\fX$-modules de Higgs $\varepsilon$-quasi-petits (resp. quasi-petits) dont le $\co_\fX$-module sous-jacent est $\cS$-plat.

\begin{rema}\label{definf25}
Soient $\varepsilon$ un nombre rationnel $>0$, $\cN$ un $\co_\fX$-module cohérent et $\cS$-plat, $\theta$ un $\co_\fX$-champ de Higgs à coefficients dans 
$\txi^{-1}\tOmega^1_{\fX/\cS}$. Pour que $(\cN,\theta)$ soit $\varepsilon$-quasi-petit, il faut et il suffit que l'on ait
\begin{equation}\label{definf25a}
\theta(\cN)\subset  p^\varepsilon \txi^{-1}\tOmega^1_{\fX/\cS}\otimes_{\co_\fX} \cN.
\end{equation}
\end{rema}

\begin{defi}\label{definf26}
Soit $(\cN,\theta)$ un $\co_\fX[\frac 1 p]$-fibré de Higgs  à coefficients dans $\txi^{-1}\tOmega^1_{\fX/\cS}$ \eqref{definf20}.
\begin{itemize}
\item[(i)] On dit que $(\cN,\theta)$ est {\em petit} s'il existe un sous-$\co_\fX$-module cohérent $\fN$ de $\cN$
qui l'engendre sur $\co_\fX[\frac 1 p]$ et un nombre rationnel $\varepsilon>\frac{1}{p-1}$ tels que 
\begin{equation}\label{definf26a}
\theta(\fN)\subset  p^\varepsilon \txi^{-1}\tOmega^1_{\fX/\cS}\otimes_{\co_\fX} \fN.
\end{equation}
\item[(ii)] On dit que $(\cN,\theta)$ est {\em localement petit} s'il existe un recouvrement ouvert $(U_i)_{i\in I}$ de $X_s$
tel que pour tout $i\in I$, $(\cN|U_i,\theta|U_i)$ soit petit. 
\end{itemize}
\end{defi}

On désigne par $\bMH^{\p}(\co_\fX[\frac 1 p],\txi^{-1}\tOmega^1_{\fX/\cS})$ (resp. $\bMH^{\locp}(\co_\fX[\frac 1 p],\txi^{-1}\tOmega^1_{\fX/\cS})$) 
la sous-catégorie pleine de $\bMH(\co_\fX[\frac 1 p],\txi^{-1}\tOmega^1_{\fX/\cS})$ formée des $\co_\fX[\frac 1 p]$-fibrés de Higgs petits (resp. localement petits).

\begin{rema}\label{definf260}
Supposons $X$ affine et le $\co_X$-module $\tOmega^1_{X/S}$ libre de type fini. 
Posons $R=\Gamma(X,\co_X)$ et $R_1=R\otimes_{\co_K}\co_\oK$
et notons $\hR$ et $\hRun$ leurs séparés complétés $p$-adiques.
Soit  $(\cN,\theta)$ un $\co_\fX[\frac 1 p]$-fibré de Higgs à coefficients dans $\txi^{-1}\tOmega^1_{\fX/\cS}$.
Posons $N=\Gamma(\fX,\cN)$ qui est un $\hRun[\frac 1 p]$-module projectif de type fini d'après (\cite{agt} III.6.17), 
et notons $\uptheta$ le $\hRun[\frac 1 p]$-champ de Higgs sur $N$ à coefficients dans 
$\txi^{-1}\tOmega^1_{X/S}(X)$ induit par $\theta$ (\cite{egr1} 2.7.2.3 et 2.10.24). 
Pour que le $\co_\fX[\frac 1 p]$-fibré de Higgs $(\cN,\theta)$ soit petit dans le sens de \ref{definf26}, 
il faut et il suffit que le $\hRun[\frac 1 p]$-module de Higgs $(N,\uptheta)$ soit petit dans le sens de \ref{pmh25}.
En effet, la condition est nécessaire en vertu de (\cite{egr1} (2.10.5.1)) et elle est suffisante compte tenu de  (\cite{egr1} 1.10.2).
\end{rema}

\begin{prop}\label{definf27}
Supposons que $X$ soit affine, que $X_s$ soit non-vide et que $f$ \eqref{defing11} admette une carte adéquate \eqref{cad1}. 
Alors, tout $\co_\fX[\frac 1 p]$-fibré de Higgs nilpotent \eqref{definf22} à coefficients dans $\txi^{-1}\tOmega^1_{\fX/\cS}$ est petit \eqref{definf26}. 
\end{prop}

Posons $R=\Gamma(X,\co_X)$ et $R_1=R\otimes_{\co_K}\co_\oK$ et notons $\hRun$ le séparé complété $p$-adique de $R_1$. 
Soit $(\cN,\theta)$ un $\co_\fX[\frac 1 p]$-fibré de Higgs nilpotent à coefficients dans $\txi^{-1}\tOmega^1_{\fX/\cS}$.
Posons $N=\Gamma(\fX,\cN)$ qui est un $\hRun[\frac 1 p]$-module projectif de type fini d'après (\cite{agt} III.6.17), 
et notons $\uptheta$ le $\hRun[\frac 1 p]$-champ de Higgs sur $N$ à coefficients dans 
$\txi^{-1}\tOmega^1_{X/S}(X)$ induit par $\theta$. 
D'après \ref{definf23}, le $\hRun[\frac 1 p]$-module de Higgs $(N,\uptheta)$ est nilpotent dans le sens de \ref{MH15}. 
Il est donc petit dans le sens de \ref{pmh25}, en vertu de \ref{repht2}(i). Par suite, le $\co_\fX[\frac 1 p]$-fibré de Higgs $(\cN,\theta)$ est petit d'après \ref{definf260}.

\begin{cor}\label{definf270}
Tout $\co_\fX[\frac 1 p]$-fibré de Higgs nilpotent \eqref{definf22} à coefficients dans $\txi^{-1}\tOmega^1_{\fX/\cS}$ est localement petit \eqref{definf26}. 
\end{cor}

\subsection{}\label{definf28}
On note $\bMH'^{\qpp}(\co_\fX,\txi^{-1}\tOmega^1_{\fX/\cS})$ la sous-catégorie pleine de $\bMH^{\qpp}(\co_\fX,\txi^{-1}\tOmega^1_{\fX/\cS})$ formée des 
$\co_\fX$-modules de Higgs $(\cN,\theta)$ tels que le $\co_\fX[\frac 1 p]$-module $\cN_{\mQ_p}$ soit localement projectif de type fini. C'est une catégorie additive.
On note $\bMH'^{\qpp}_\mQ(\co_\fX,\txi^{-1}\tOmega^1_{\fX/\cS})$ la catégorie des objets de $\bMH'^{\qpp}(\co_\fX,\txi^{-1}\tOmega^1_{\fX/\cS})$ à isogénie près
(\cite{agt} III.6.1.1). Le foncteur 
\begin{equation}\label{definf28a}
\begin{array}{clcr}
\bMH'^{\qpp}(\co_\fX,\txi^{-1}\tOmega^1_{\fX/\cS})&\rightarrow &\bMH^{\p}(\co_\fX[\frac 1 p],\txi^{-1}\tOmega^1_{\fX/\cS})\\
(\cN,\theta)&\mapsto&(\cN_{\mQ_p},\theta)
\end{array}
\end{equation}
induit un foncteur 
\begin{equation}\label{definf28b}
\bMH'^{\qpp}_{\mQ}(\co_\fX,\txi^{-1}\tOmega^1_{\fX/\cS})\rightarrow \bMH^{\p}(\co_\fX[\frac 1 p],\txi^{-1}\tOmega^1_{\fX/\cS}).
\end{equation}

\begin{lem}\label{definf29}
Le foncteur \eqref{definf28b} est une équivalence de catégories.
\end{lem}

En effet, il résulte aussitôt des définitions que ce foncteur est essentiellement surjectif. 
Soient $\cM,\cN$ deux $\co_\fX$-modules cohérents et $\cS$-plats. 
D'après \eqref{aspglob1d}, le morphisme canonique 
\begin{equation}\label{definf29a}
\Hom_{\co_\fX}(\cM,\cN)\otimes_{\mZ_p}\mQ_p\rightarrow \Hom_{\co_\fX[\frac 1 p]}(\cM\otimes_{\mZ_p}\mQ_p,\cN\otimes_{\mZ_p}\mQ_p)
\end{equation}
est un isomorphisme. Supposons que $\cM$ et $\cN$ soient munis de $\co_\fX$-champs de Higgs à coefficients dans $\txi^{-1}\tOmega^1_{\fX/\cS}$. 
Le morphisme canonique 
\begin{equation}\label{definf29b}
\Hom_{\bMH}(\cM,\cN)\otimes_{\mZ_p}\mQ_p\rightarrow \Hom_{\bMH}(\cM\otimes_{\mZ_p}\mQ_p,\cN\otimes_{\mZ_p}\mQ_p)
\end{equation}
où la source (resp. le but) désigne le module des morphismes de $\co_\fX$-modules ($\co_\fX[\frac 1 p]$-modules) de Higgs est alors un isomorphisme:
l'injectivité résulte aussitôt de celle de \eqref{definf29a} et la surjectivité de celle de \eqref{definf29a} et du fait que $\cM$ et $\cN$ sont $\cS$-plats et que
$\tOmega^1_{\fX/\cS}$ est $\fX$-plat; d'où la proposition.

\section{Topos de Faltings}\label{ahttf}

\subsection{}\label{tf0}
On désigne par $\cR$ la catégorie des revêtements étales ({\em i.e.}, la sous-catégorie pleine de la catégorie 
des morphismes de schémas, formée des revêtements étales) et par 
\begin{equation}\label{tf0a}
\cR\rightarrow \Sch
\end{equation}
le ``foncteur but'', qui fait de $\cR$ une catégorie fibrée clivée et normalisée au-dessus de $\Sch$ \eqref{notconv3} (\cite{sga1} VI):  
la catégorie fibre au-dessus de tout schéma $X$ est canoniquement équivalente à
la catégorie $\Et_{\rf/X}$ des schémas étales finis sur $X$, et pour tout morphisme de schémas 
$f\colon Y\rightarrow X$, le foncteur image inverse $f^+\colon \Et_{\rf/X}\rightarrow \Et_{\rf/Y}$
n'est autre que le foncteur changement de base par $f$ \eqref{notconv10}. 
Munissant chaque fibre de la topologie étale, $\cR/\Sch$ devient un $\mU$-site fibré (\cite{sga4} VI 7.2.1). 

\subsection{}\label{ahttf2}
Avec les notations de \ref{notconv10}, on désigne par
\begin{equation}\label{ahttf2a}
\pi\colon E\rightarrow \Et_{/X}
\end{equation}
le $\mU$-site fibré de Faltings associé au morphisme $h\colon \oX^\circ\rightarrow X$ \eqref{defing11a} (\cite{agt} VI.10.1),
c'est-à-dire le $\mU$-site fibré déduit du site fibré des revêtements étales $\cR/\Sch$ \eqref{tf0a} 
par changement de base par le foncteur 
\begin{equation}
\Et_{/X}\rightarrow \Sch, \ \ \ U\mapsto \oU^\circ=U\times_X\oX^\circ.
\end{equation} 
Pour tout $U\in \ob(\Et_{/X})$, on note 
\begin{equation}\label{ahttf2b}
\iota_U\colon \Et_{\rf/\oU^\circ}\rightarrow E
\end{equation} 
le foncteur canonique. On prendra garde que ce foncteur était noté $\alpha_{U!}$ dans (\cite{agt} (VI.5.1.2)). 

On munit $E$ de la topologie co-évanescente associée à $\pi$ (\cite{agt} VI.5.3), autrement dit,  
la topologie engendrée par les recouvrements $\{(V_i\rightarrow U_i)\rightarrow (V\rightarrow U)\}_{i\in I}$
des deux types suivants: 
\begin{itemize}
\item[(v)] $U_i=U$ pour tout $i\in I$, et $(V_i\rightarrow V)_{i\in I}$ est un recouvrement étale. 
\item[(c)] $(U_i\rightarrow U)_{i\in I}$ est un recouvrement étale et $V_i=U_i\times_UV$ pour tout $i\in I$. 
\end{itemize}
Le site co-évanescent $E$ ainsi défini est encore appelé {\em site de Faltings} associé à $h$;  c'est un $\mU$-site.  
On désigne par $\hE$ (resp. $\tE$) la catégorie des préfaisceaux (resp. le topos des faisceaux) 
de $\mU$-ensembles sur $E$. On appelle $\tE$ le {\em topos de Faltings} associé à $h$ (\cite{agt} VI.10.1).

On désigne par 
\begin{equation}\label{ahttf2g}
\fF\rightarrow  \Et_{/X}
\end{equation}
le $\mU$-topos fibré associé à $\pi$. La catégorie fibre de $\fF$ au-dessus de tout $U\in \ob(\Et_{/X})$ est 
canoniquement équivalente au topos fini étale $\oU^\circ_\fet$ de $\oU^\circ$
et le foncteur image inverse pour tout morphisme $\upmu\colon U'\rightarrow U$ de $\Et_{/X}$ s'identifie au foncteur 
$\oupmu^{\circ*}_{\fet}\colon \oU^\circ_\fet\rightarrow \oU'^\circ_\fet$ image inverse par le morphisme 
de topos $\oupmu^\circ_\fet\colon \oU'^\circ_\fet\rightarrow \oU^\circ_\fet$ (\cite{agt} VI.9.3). On désigne par
\begin{equation}\label{ahttf2h}
\fF^\vee\rightarrow (\Et_{/X})^\circ
\end{equation}
la catégorie fibrée obtenue en associant à tout $U\in \ob(\Et_{/X})$ la catégorie $\oU^\circ_\fet$, et à tout morphisme 
$\upmu\colon U'\rightarrow U$ de $\Et_{/X}$ le foncteur 
$\oupmu^\circ_{\fet*}\colon \oU'^\circ_\fet\rightarrow \oU^\circ_\fet$ 
image directe par le morphisme de topos $\oupmu^\circ_\fet$. On désigne par
\begin{equation}\label{ahttf2i}
\cP^\vee\rightarrow (\Et_{/X})^\circ
\end{equation}
la catégorie fibrée obtenue en associant à tout $U\in \ob(\Et_{/X})$ la catégorie $(\Et_{\rf/\oU^\circ})^\wedge$ 
des préfaisceaux de $\mU$-ensembles sur $\Et_{\rf/\oU^\circ}$, et à tout morphisme $\upmu\colon U'\rightarrow U$ de $\Et_{/X}$ 
le foncteur 
\begin{equation}\label{ahttf2j}
\oupmu^\circ_{\fet*}\colon (\Et_{\rf/\oU'^\circ})^\wedge\rightarrow (\Et_{\rf/\oU^\circ})^\wedge
\end{equation} 
obtenu en composant avec le foncteur image inverse $\oupmu^{\circ+}\colon \Et_{\rf/\oU^\circ}\rightarrow \Et_{\rf/\oU'^\circ}$.  

On a une équivalence de catégories (\cite{agt} VI.5.2)
\begin{eqnarray}\label{ahttf2k}
\hE&\rightarrow& \bHom_{(\Et_{/X})^\circ}((\Et_{/X})^\circ,\cP^\vee)\\
F&\mapsto &\{U\mapsto F\circ \iota_U\}.\nonumber
\end{eqnarray}
On identifiera dans la suite $F$ à la section $\{U\mapsto F\circ \iota_U\}$ qui lui est associée par cette équivalence.

D'après (\cite{agt} VI.5.11), le foncteur \eqref{ahttf2k} induit un foncteur pleinement fidèle 
\begin{equation}\label{ahttf2l}
\tE\rightarrow \bHom_{(\Et_{/X})^\circ}((\Et_{/X})^\circ,\fF^\vee),
\end{equation}
d'image essentielle les sections $\{U\mapsto F_U\}$ vérifiant une condition de recollement.

Les foncteurs 
\begin{eqnarray}
\sigma^+\colon \Et_{/X}\rightarrow E,&& U\mapsto (\oU^\circ\rightarrow U),\label{ahttf2e}\\
\iota_{X}\colon \Et_{\rf/\oX^\circ}\rightarrow E,&& V\mapsto (V\rightarrow X),\label{ahttf2f}
\end{eqnarray}
sont continus et exacts à gauche (\cite{agt} VI.10.6). Ils définissent donc deux morphismes de topos 
\begin{eqnarray}
\sigma\colon \tE \rightarrow X_\et,\label{ahttf2c}\\
\beta\colon \tE \rightarrow \oX^\circ_\fet.\label{ahttf2d}
\end{eqnarray}
Pour tout faisceau $F=\{U\mapsto F_U\}$ sur $E$, on a $\beta_*(F)=F_X$. 

Le foncteur 
\begin{equation}\label{ahttf2m}
\psi^+\colon E\rightarrow \Et_{/\oX^\circ},\ \ \ (V\rightarrow U)\mapsto V
\end{equation}
est continu et exact à gauche (\cite{agt} VI.10.7). Il définit donc un morphisme de topos 
\begin{equation}\label{ahttf2n}
\psi\colon \oX^\circ_\et\rightarrow \tE.
\end{equation}
Nous changeons ici de notations par rapport à {\em loc. cit.}

\subsection{}\label{ahttf12}
Comme $X_\eta$ est un ouvert de $X_\et$, {\em i.e.}, un sous-objet de l'objet final $X$ (\cite{sga4} IV 8.3),
$\sigma^*(X_\eta)$ est un ouvert de $\tE$. On note 
\begin{equation}
\gamma\colon \tE_{/\sigma^*(X_\eta)}\rightarrow \tE
\end{equation}
le morphisme de localisation de $\tE$ en $\sigma^*(X_\eta)$ (\cite{sga4} IV 5.2). 
On désigne par $\tE_s$ le sous-topos fermé de $\tE$ complémentaire de l'ouvert $\sigma^*(X_\eta)$, 
c'est-à-dire la sous-catégorie pleine de $\tE$ formée des faisceaux $F$ tels que $\gamma^*(F)$
soit un objet final de $\tE_{/\sigma^*(X_\eta)}$  (\cite{sga4} IV 9.3.5), par 
\begin{equation}\label{ahttf12a}
\delta\colon \tE_s\rightarrow \tE
\end{equation} 
le plongement canonique et par
\begin{equation}\label{ahttf12b}
\sigma_s\colon \tE_s\rightarrow X_{s,\et}
\end{equation} 
le morphisme de topos induit par $\sigma$ (\cite{agt} (III.9.8.3)). Le diagramme de morphismes de topos 
\begin{equation}\label{ahttf12c}
\xymatrix{
{\tE_s}\ar[r]^-(0.5){\sigma_s}\ar[d]_{\delta}&{X_{s,\et}}\ar[d]^{a}\\
{\tE}\ar[r]^-(0.5)\sigma&{X_\et}}
\end{equation}
est commutatif à isomorphisme près.

\subsection{}\label{TFA9}
On désigne par $X_\et\gtimes_{X_\et}\oX^\circ_\et$ le topos co-évanescent du morphisme $f_\et\colon \oX^\circ_\et\rightarrow X_\et$ (\cite{agt} VI.3.12) et par
\begin{equation}\label{TFA9a}
\rho\colon X_\et\gtimes_{X_\et}\oX^\circ_\et\rightarrow \tE
\end{equation}
le morphisme canonique (\cite{agt} VI.10.15). 
D'après (\cite{agt} VI.4.20) et (\cite{sga4} VIII 7.9), la donnée d'un point de $X_\et\gtimes_{X_\et}\oX^\circ_\et$ 
est équivalente à la donnée d'une paire de points géométriques $\ox$ de $X$ et $\oy$ de $\oX^\circ$
et d'une flèche de spécialisation de $\hbar(\oy)$ vers $\ox$, c'est-à-dire, d'un $X$-morphisme 
$\oy\rightarrow X_{(\ox)}$, où $X_{(\ox)}$ désigne le localisé strict de $X$ en $\ox$. 
Un tel point sera noté $(\oy\rightsquigarrow \ox)$, et son image par $\rho$ sera notée 
$\rho(\oy\rightsquigarrow \ox)$, qui est donc un point de $\tE$.  
La famille des foncteurs fibres de $\tE$ associés à ces points est conservative d'après (\cite{agt} VI.10.21).

\subsection{}\label{TFA14}
Soient $\ox$ un point géométrique de $X$, $\uX$ le localisé strict de $X$ en $\ox$. On pose $\uoX=\uX\times_S\oS$ \eqref{defing1c} et 
$\uoX^\circ=\uoX\times_XX^\circ$ \eqref{defing11a}. 
On désigne par $\tuE$ le topos de Faltings associé au morphisme canonique $\uoX^\circ\rightarrow \uX$, par 
\begin{equation}\label{TFA14b}
\Phi\colon \tuE\rightarrow \tE
\end{equation}
le morphisme de fonctorialité induit par le morphisme canonique $\uX\rightarrow X$, par 
\begin{equation}\label{TFA14c}
\ubeta\colon \tuE\rightarrow \uoX^\circ_\fet
\end{equation}
le morphisme canonique \eqref{ahttf2d} et par 
\begin{equation}\label{TFA14d}
\theta\colon \uoX^\circ_\fet\rightarrow \tuE
\end{equation}
la section canonique de $\ubeta$ définie dans (\cite{agt} VI.10.23). On note 
\begin{equation}\label{TFA14e}
\varphi_\ox\colon \tE\rightarrow \uoX^\circ_\fet
\end{equation}
le foncteur composé $\theta^*\circ \Phi^*$ (\cite{agt} VI.10.29).  

Le foncteur composé $\varphi_\ox\circ\beta^*$ est canoniquement isomorphe au foncteur image inverse par 
le morphisme canonique $\uoX^\circ_\fet\rightarrow \oX^\circ_\fet$, d'après (\cite{agt} (VI.10.24.3) et (VI.10.12.6)). 
Pour tout objet $F$ de $X_\et$, on a un isomorphisme canonique et fonctoriel 
\begin{equation}\label{TFA14f}
\varphi_\ox(\sigma^*(F))\stackrel{\sim}{\rightarrow} F_\ox,
\end{equation}
de but le faisceau constant de $\uoX^\circ_\fet$ de valeur $F_\ox$, en vertu de (\cite{agt} VI.10.24 et (VI.10.12.6)). 

On désigne par $\fV_\ox$ la catégorie des $X$-schémas étales $\ox$-pointés, ou ce qui revient au même,
la catégorie des voisinages du point de $X_\et$ associé à $\ox$ dans le site $\Et_{/X}$ (\cite{sga4} IV 6.8.2 et VIII 3.9).
Pour tout objet $(U,\fp\colon \ox\rightarrow U)$ de $\fV_\ox$, 
on note encore $\fp\colon \uX\rightarrow U$ le morphisme déduit de 
$\fp$ (\cite{sga4} VIII 7.3) et on pose
\begin{equation}\label{TFA14a}
\ofp^\circ=\fp\times_X\oX^\circ\colon \uoX^\circ \rightarrow \oU^\circ.
\end{equation}

Pour tout objet $F=\{U\mapsto F_U\}$ de $\hE$ \eqref{ahttf2k},
on note $F^\tta$ le faisceau de $\tE$ associé à $F$, et pour tout $U\in \ob(\Et_{/X})$, 
$F_U^\tta$ le faisceau de $\oU^\circ_\fet$ associé à $F_U$. 
D'après (\cite{agt} VI.10.37),  on a un isomorphisme canonique et fonctoriel
\begin{equation}\label{TFA14g}
\varphi_\ox(F^\tta) \stackrel{\sim}{\rightarrow}\underset{\underset{(U,\fp)\in \fV_\ox^\circ}{\longrightarrow}}{\lim}\ (\ofp^\circ)_\fet^*(F^\tta_U).
\end{equation}

En vertu de (\cite{agt} VI.10.30), pour tout groupe abélien $F$ de $\tE$
et tout entier $q\geq 0$, on a un isomorphisme canonique et fonctoriel
\begin{equation}\label{TFA14h}
\rR^q\sigma_*(F)_\ox\stackrel{\sim}{\rightarrow}\rH^q(\uoX^\circ_\fet,\varphi_\ox(F)). 
\end{equation}

\subsection{}\label{TFA2}
Le schéma $\oX$ est normal et localement irréductible d'après (\cite{agt} III.4.2(iii)). 
Par ailleurs, l'immersion $j\colon X^\circ\rightarrow X$ est quasi-compacte puisque $X$ est noethérien.  
Pour tout objet $(V\rightarrow U)$ de $E$, on note $\oU^V$ la fermeture intégrale de $\oU$ dans $V$. 
Pour tout morphisme $(V'\rightarrow U')\rightarrow (V\rightarrow U)$ de $E$, on a un morphisme canonique 
$\oU'^{V'}\rightarrow \oU^V$ qui s'insère dans un diagramme commutatif 
\begin{equation}\label{TFA2a}
\xymatrix{
V'\ar[r]\ar[d]&{\oU'^{V'}}\ar[d]\ar[r]&\oU'\ar[r]\ar[d]&U'\ar[d]\\
V\ar[r]&{\oU^V}\ar[r]&\oU\ar[r]&U}
\end{equation} 
On désigne par $\ocB$ le préfaisceau sur $E$ défini pour tout $(V\rightarrow U)\in \ob(E)$, par 
\begin{equation}\label{TFA2b}
\ocB((V\rightarrow U))=\Gamma(\oU^V,\co_{\oU^V}).
\end{equation}
C'est un faisceau pour la topologie co-évanescente sur $E$ en vertu de (\cite{agt} III.8.16). 
Pour tout $U\in \ob(\Et_{/X})$, on pose 
\begin{equation}\label{TFA2d}
\ocB_{U}=\ocB\circ \iota_{U}.
\end{equation} 

D'après (\cite{agt} III.8.17), on a un homomorphisme canonique 
\begin{equation}\label{TFA2c}
\sigma^*(\hbar_*(\co_\oX))\rightarrow \ocB.
\end{equation}
Sauf mention explicite du contraire, on considère $\sigma$ \eqref{ahttf2c}
comme un morphisme de topos annelés
\begin{equation}\label{TFA2e}
\sigma\colon (\tE,\ocB)\rightarrow (X_\et,\hbar_*(\co_\oX)).
\end{equation}

Notant encore $\co_\oK$ le faisceau constant de $\oX^\circ_\fet$ de valeur $\co_\oK$. 
On considèrera aussi $\beta$ \eqref{ahttf2d} comme un morphisme de topos annelés
\begin{equation}\label{TFA2g}
\beta\colon (\tE,\ocB)\rightarrow (\oX^\circ_\fet,\co_\oK).
\end{equation}

\subsection{}\label{ahttf44}
Soient $U\in \ob(\Et_{/X})$, $\oy$ un point géométrique de $\oU^\circ$ \eqref{defing11a}. 
Le schéma $\oU$ étant localement irréductible d'après (\cite{agt} III.3.3 et III.4.2(iii)),  
il est la somme des schémas induits sur ses composantes irréductibles. On note $\oU^\star$
la composante irréductible de $\oU$ contenant $\oy$. 
De même, $\oU^\circ$ est la somme des schémas induits sur ses composantes irréductibles
et $\oU^{\star \circ}=\oU^\star\times_{X}X^\circ$ est la composante irréductible de $\oU^\circ$ contenant $\oy$. 
On note $\bB_{\pi_1(\oU^{\star \circ},\oy)}$ le topos classifiant du groupe profini $\pi_1(\oU^{\star \circ},\oy)$ et
\begin{equation}\label{ahttf44a}
\nu_\oy\colon \oU^{\star \circ}_\fet \stackrel{\sim}{\rightarrow}\bB_{\pi_1(\oU^{\star \circ},\oy)}
\end{equation}
le foncteur fibre  de $\oU^{\star \circ}_\fet$ en $\oy$ (\cite{agt} VI.9.8). On pose
\begin{equation}\label{ahttf44b}
\oR^\oy_U=\nu_\oy(\ocB_U|\oU^{\star \circ}).
\end{equation}

\subsection{}\label{TFA11}
Soient $(\oy\rightsquigarrow \ox)$ un point de $X_\et\gtimes_{X_\et}\oX^\circ_\et$ \eqref{TFA9}, 
$\uX$ le localisé strict de $X$ en $\ox$, $\fV_\ox$ la catégorie des $X$-schémas étales $\ox$-pointés (cf. \ref{TFA14}). 
On pose $\uoX=\uX\times_S\oS$ \eqref{defing1c} et $\uoX^\circ=\uoX\times_XX^\circ$ \eqref{defing11a}.
Le $X$-morphisme $u\colon \oy\rightarrow \uX$
définissant $(\oy\rightsquigarrow \ox)$ se relève en un $\oX^\circ$-morphisme $v\colon \oy\rightarrow \uoX^\circ$ et 
induit donc un point géométrique de $\uoX^\circ$ que l'on note aussi (abusivement) $\oy$.
Pour tout objet $(U,\fp\colon \ox\rightarrow U)$ de $\fV_\ox$, 
on note encore $\fp\colon \uX\rightarrow U$ le morphisme déduit de $\fp$.
On note aussi (abusivement) $\oy$ le point géométrique $\ofp^\circ(v(\oy))$ de $\oU^\circ$ \eqref{TFA14a}.

Pour tout objet $F=\{U\mapsto F_U\}$ de $\hE$ \eqref{ahttf2k}, on note
$F^\tta$ le faisceau de $\tE$ associé à $F$, et pour tout $U\in \ob(\Et_{/X})$, 
$F_U^\tta$ le faisceau de $\oU^\circ_\fet$ associé à $F_U$. 
D'après (\cite{agt} VI.10.36 et (VI.9.3.4)), on a un isomorphisme canonique et fonctoriel 
\begin{equation}\label{TFA11a}
(F^\tta)_{\rho(\oy\rightsquigarrow \ox)} \stackrel{\sim}{\rightarrow} 
\underset{\underset{(U,\fp)\in \fV_\ox^\circ}{\longrightarrow}}{\lim}\ (F^\tta_U)_{\rho_{\oU^\circ}(\oy)},
\end{equation}
où $\rho$ est le morphisme \eqref{TFA9a} et $\rho_{\oU^\circ}\colon \oU^\circ_\et\rightarrow \oU^\circ_\fet$ 
est le morphisme canonique \eqref{notconv10a}. 
Compte tenu de \eqref{ahttf44b} et (\cite{agt} VI.9.9), on en déduit un isomorphisme canonique de $\Gamma(\uoX,\co_{\uoX})$-algèbres
\begin{equation}\label{TFA11c}
\ocB_{\rho(\oy\rightsquigarrow \ox)}\stackrel{\sim}{\rightarrow} 
\underset{\underset{(U,\fp)\in \fV_\ox^\circ}{\longrightarrow}}{\lim}\ \oR^{\oy}_U.
\end{equation}

\subsection{}\label{TFA12}
Conservons les hypothèses et notations de \ref{TFA11}; supposons, de plus, que {\em $\ox$ soit au-dessus de $s$}. 
D'après (\cite{agt} III.3.7), $\uoX$ est normal et strictement local (et en particulier intègre). 
Pour tout objet $(U,\fp\colon \ox\rightarrow U)$ de $\fV_\ox$, on désigne par $\oU^\star$  
la composante irréductible de $\oU$ contenant $\oy$ 
et on pose $\oU^{\star\circ}=\oU^\star\times_XX^\circ$ qui est la composante irréductible de $\oU^\circ$ contenant $\oy$ (cf. \ref{ahttf44}).
Le morphisme $\ofp^\circ \colon \uoX^\circ\rightarrow \oU^\circ$
\eqref{TFA14a} se factorise donc à travers $\oU^{\star \circ}$. 
On désigne par $\bB_{\pi_1(\uoX^{\circ},\oy)}$ le topos classifiant du groupe profini 
$\pi_1(\uoX^{\circ},\oy)$ et par
\begin{equation}\label{TFA12d}
\nu_\oy\colon \uoX^\circ_\fet \stackrel{\sim}{\rightarrow}\bB_{\pi_1(\uoX^\circ,\oy)}
\end{equation}
le foncteur fibre de $\uoX^\circ_\fet$ en $\oy$  (\cite{agt}  (VI.9.8.4)). 
D'après (\cite{agt} VI.10.31 et VI.9.9), le foncteur composé
\begin{equation}\label{TFA12e}
\xymatrix{
{\tE}\ar[r]^-(0.5){\varphi_\ox}&{\uoX^\circ_\fet}\ar[r]^-(0.5){\nu_\oy}&{\bB_{\pi_1(\uoX^\circ,\oy)}}\ar[r]&\Ens},
\end{equation}
où $\varphi_\ox$ est le foncteur canonique \eqref{TFA14e} et la dernière flèche est le foncteur d'oubli de l'action de $\pi_1(\uoX^\circ,\oy)$, 
est canoniquement isomorphe au foncteur fibre associé au point $\rho(\oy\rightsquigarrow \ox)$ de $\tE$ \eqref{TFA9a}.  

On définit la $\Gamma(\uoX,\co_{\uoX})$-algèbre $\oR_{\uX}^\oy$ de $\bB_{\pi_1(\uoX^\circ,\oy)}$ par la formule
\begin{equation}\label{TFA12f}
\oR_{\uX}^\oy=\underset{\underset{(U,\fp)\in \fV_\ox^\circ}{\longrightarrow}}{\lim}\ \oR^{\oy}_U, 
\end{equation}
où l'on considère $\oR^{\oy}_U$ comme une $\Gamma(\oU,\co_\oU)$-algèbre de $\bB_{\pi_1(\oU^{\star\circ},\oy)}$ \eqref{ahttf44b}.  
L'isomorphisme \eqref{TFA14g} induit un isomorphisme de $\Gamma(\uoX,\co_{\uoX})$-algèbres de 
$\bB_{\pi_1(\uoX^\circ,\oy)}$
\begin{equation}\label{TFA12g}
\nu_\oy(\varphi_\ox(\ocB))\stackrel{\sim}{\rightarrow} \oR_{\uX}^\oy,
\end{equation} 
dont l'isomorphisme de $\Gamma(\uoX,\co_{\uoX})$-algèbres sous-jacent est \eqref{TFA11c}.

\subsection{}\label{ahttf42}
Pour tout entier $n\geq 0$, on pose
\begin{equation}\label{ahttf3b}
\ocB_n=\ocB/p^n\ocB.
\end{equation}
D'après (\cite{agt} III.9.7), $\ocB_n$ est un anneau de $\tE_s$. Pour tout $U\in \ob(\Et_{/X})$, on pose \eqref{TFA2d}
\begin{equation}\label{ahttf3e}
\ocB_{U,n}=\ocB_U/p^n\ocB_U.
\end{equation}

On a un homomorphisme canonique $\sigma_s^*(\co_{\oX_n})\rightarrow \ocB_n$ (\cite{agt} (III.9.9.3)).
Le morphisme $\sigma_s$ \eqref{ahttf12b} est donc sous-jacent à un morphisme de topos annelés que l'on note
\begin{equation}\label{ahttf42a}
\sigma_n\colon (\tE_s,\ocB_n)\rightarrow (X_{s,\et},\co_{\oX_n}).
\end{equation}
Nous utilisons pour les modules la notation $\sigma_n^{-1}$ (ou $\sigma_s^*$)
pour désigner l'image inverse au sens des faisceaux abéliens et nous réservons la notation 
$\sigma_n^*$ pour l'image inverse au sens des modules.

On a un homomorphisme canonique $\ocB_{X,n}\rightarrow \beta_*(\ocB_n)$, qui n'est pas un isomorphisme en général (\cite{ag} 4.1.8).
Le morphisme composé $\beta\circ \delta$ est donc sous-jacent à un morphisme de topos annelés que l'on note
\begin{equation}\label{ahttf42b}
\beta_n\colon (\tE_s,\ocB_n) \rightarrow (\oX^\circ_\fet,\ocB_{X,n}).
\end{equation}

On désigne par 
\begin{equation}\label{ahttf42c}
\Sigma_n\colon (\tE_s,\ocB_n)\rightarrow (X_{s,\zar},\co_{\oX_n})
\end{equation}
le composé du morphisme $\sigma_n$ et du morphisme canonique \eqref{notconv12j}
\begin{equation}\label{ahttf42d}
u_n\colon (X_{s,\et},\co_{\oX_n})\rightarrow (X_{s,\zar},\co_{\oX_n}).
\end{equation}
Nous utilisons pour les modules la notation $\Sigma_n^{-1}$ 
pour désigner l'image inverse au sens des faisceaux abéliens et nous réservons la notation 
$\Sigma_n^*$ pour l'image inverse au sens des modules.

\subsection{}\label{ahttf13}
Pour tout $\mU$-topos $T$, on note $T^{\mN^\circ}$ le topos des systèmes projectifs de $T$, 
indexés par l'ensemble ordonné $\mN$ des entiers naturels \eqref{notconv13}. 
On désigne par $\bvocB$ l'anneau $(\ocB_{n+1})_{n\in \mN}$ de $\tE_s^{\mN^\circ}$ \eqref{ahttf3b},
par $\co_{\bvoX}$ l'anneau $(\co_{\oX_{n+1}})_{n\in \mN}$ de $X_{s,\et}^{\mN^\circ}$ ou $X_{s,\zar}^{\mN^\circ}$, selon le contexte,
et par $\txi^{-1}\tOmega^1_{\bvoX/\bvoS}$ le $\co_{\bvoX}$-module $(\txi^{-1}\tOmega^1_{\oX_{n+1}/\oS_{n+1}})_{n\in \mN}$ \eqref{ahttf1c}. 
On prendra garde de ne pas confondre $\co_{\bvoX}$ et $\co_{\coX}$ \eqref{defing1c}.

On note
\begin{equation}\label{ahttf13a}
\bvsigma\colon (\tE_s^{\mN^\circ},\bvocB)\rightarrow(X_{s,\et}^{\mN^\circ},\co_{\bvoX})
\end{equation}
le morphisme de topos annelés induit par les $(\sigma_{n+1})_{n\in \mN}$ \eqref{ahttf42a}. 
Nous utilisons pour les modules la notation $\bvsigma^{-1}$ pour désigner l'image
inverse au sens des faisceaux abéliens et nous réservons la notation 
$\bvsigma^*$ pour l'image inverse au sens des modules.

On note
\begin{equation}\label{ahttf13c}
\bvSigma\colon (\tE_s^{\mN^\circ},\bvocB)\rightarrow (X_{s,\zar}^{\mN^\circ},\co_{\bvoX})
\end{equation}
le morphisme de topos annelés défini par les $(\Sigma_{n+1})_{n\in \mN}$ \eqref{ahttf42c}. Nous utilisons pour les modules la notation 
$\bvSigma^{-1}$ pour désigner l'image inverse au sens des faisceaux abéliens et nous réservons la notation 
$\bvSigma^*$ pour l'image inverse au sens des modules.

On rappelle que $\fX$ désigne le schéma formel complété $p$-adique de $\oX$ \eqref{definf18}. On note
\begin{equation}\label{ahttf13d}
\uplambda\colon (X_{s,\zar}^{\mN^\circ},\co_{\bvoX})\rightarrow (X_{s,\zar}, \co_\fX)
\end{equation}
le morphisme de topos annelés pour lequel  le foncteur $\uplambda_*$ est le foncteur limite projective \eqref{notconv13a}. On désigne par
\begin{equation}\label{ahttf13e}
\hupsigma\colon (\tE_s^{\mN^\circ},\bvocB)\rightarrow (X_{s,\zar},\co_{\fX})
\end{equation}
le morphisme composé $\uplambda\circ \bvSigma$. Nous utilisons pour les modules la notation $\hupsigma^{-1}$ pour désigner l'image
inverse au sens des faisceaux abéliens et nous réservons la notation $\hupsigma^*$ pour l'image inverse au sens des modules. 
On observera que le morphisme $\hupsigma$ a été noté $\top$ dans (\cite{agt} (III.11.1.11)).

Pour tout $\co_\fX$-module $\cF$ de $X_{s,\zar}$, on a un isomorphisme canonique 
\begin{equation}\label{ahttf17c}
\hupsigma^*(\cF)\stackrel{\sim}{\rightarrow}\bvSigma^*((\cF/p^{n+1}\cF)_{n\in \mN}). 
\end{equation}
En particulier, $\hupsigma^*(\cF)$ est adique (\cite{agt} III.7.18).

\subsection{}\label{ahttf43}
On désigne par $\bMod(\bvocB)$ la catégorie des $\bvocB$-modules de $\tE_s^{\mN^\circ}$,
par $\bIndMod(\bvocB)$ la catégorie des {\em ind-$\bvocB$-modules} \eqref{indsh15} et par 
\begin{equation}\label{ahttf43a}
\iota_\bvocB\colon \bMod(\bvocB)\rightarrow \bIndMod(\bvocB)
\end{equation}
le foncteur canonique, qui est exact et pleinement fidèle \eqref{indsh15}. 
On identifiera $\bMod(\bvocB)$ à une sous-catégorie pleine de $\bIndMod(\bvocB)$ 
par le foncteur $\iota_\bvocB$ qu'on omettra des notations. 

On désigne par $\bMod(\co_{\fX})$ la catégorie des $\co_{\fX}$-modules de $X_{s,\zar}$ \eqref{definf18}, 
par $\bIndMod(\co_\fX)$ la catégorie des ind-$\co_\fX$-modules et par 
\begin{equation}\label{ahttf43b}
\iota_{\co_\fX}\colon \bMod(\co_\fX)\rightarrow \bIndMod(\co_\fX)
\end{equation}
le foncteur canonique, qui est exact et pleinement fidèle. 
On identifiera $\bMod(\co_\fX)$ à une sous-catégorie pleine de $\bIndMod(\co_\fX)$ 
par le foncteur $\iota_{\co_\fX}$ qu'on omettra des notations.

D'après \ref{indsh21}, le morphisme de topos annelés $\hupsigma$ \eqref{ahttf13e} induit deux foncteurs additifs adjoints
\begin{eqnarray}
\rI \hupsigma^*\colon \bIndMod(\co_\fX) \rightarrow \bIndMod(\bvocB),\label{ahttf43c}\\
\rI \hupsigma_*\colon \bIndMod(\bvocB) \rightarrow \bIndMod(\co_\fX).\label{ahttf43d}
\end{eqnarray}
Le foncteur $\rI \hupsigma^*$ (resp. $\rI \hupsigma_*$) est exact à droite (resp. gauche). 
Le foncteur $\rI \hupsigma_*$ admet un foncteur dérivé à droite
\begin{equation}\label{ahttf43e}
\rR\rI \hupsigma_*\colon \bD^+(\bIndMod(\bvocB))\rightarrow \bD^+(\bIndMod(\co_\fX)).
\end{equation}
Le diagramme 
\begin{equation}\label{ahttf43f}
\xymatrix{
{\bD^+(\bMod(\bvocB))}\ar[r]^-(0.5){\rR \hupsigma_*}\ar[d]_{\iota_{\bvocB}}&{\bD^+(\bMod(\co_\fX))}\ar[d]^{\iota_{\co_\fX}}\\
{\bD^+(\bIndMod(\bvocB))}\ar[r]^-(0.5){\rR\rI \hupsigma_*}&{\bD^+(\bIndMod(\co_\fX))}}
\end{equation}
est commutatif à isomorphisme canonique près \eqref{indsh21d}. 

D'après \ref{indsh7}, Le foncteur canonique $\iota_{\co_\fX}$ admet un adjoint à gauche 
\begin{equation}\label{ahttf43g}
\kappa_{\co_\fX}\colon \bIndMod(\co_\fX) \rightarrow  \bMod(\co_\fX), 
\end{equation}
tel que pour toute petite catégorie filtrante $J$ et tout foncteur $\alpha\colon J\rightarrow \bMod(\co_\fX)$, on a un isomorphisme 
\begin{equation}\label{ahttf43h}
\kappa_{\co_\fX}(\underset{\underset{J}{\longrightarrow}}{\mlq\mlq\lim \mrq\mrq} \alpha)\stackrel{\sim}{\rightarrow} \underset{\underset{J}{\longrightarrow}}{\lim}\ \alpha.
\end{equation}
Le morphisme canonique $\kappa_{\co_\fX}\circ \iota_{\co_\fX} \rightarrow \id_{\bMod(\co_\fX)}$ est un isomorphisme. 
Le foncteur $\kappa_{\co_\fX}$ est exact d'après \ref{indsh17}.

On désigne par $\vupsigma_*$ le foncteur composé
\begin{equation}\label{ahttf43i}
\vupsigma_*=\kappa_{\co_\fX}\circ \rI \hupsigma_*\colon \bIndMod(\bvocB) \rightarrow \bMod(\co_\fX).
\end{equation}
Celui-ci est exact à gauche. Il admet un foncteur dérivé à droite
\begin{equation}\label{ahttf43j}
\rR\vupsigma_*\colon \bD^+(\bIndMod(\bvocB))\rightarrow \bD^+(\bMod(\co_\fX)),
\end{equation}
canoniquement isomorphe à $\kappa_{\co_\fX}\circ \rR\rI \hupsigma_*$ (\cite{ks2} 13.3.13).

\subsection{}\label{ahttf40}
On désigne par $\bMod_{\mQ}(\co_\fX)$ (resp. $\bMod_{\mQ}(\bvocB)$) la catégorie des $\co_\fX$-modules (resp. $\bvocB$-modules) 
à isogénie près \eqref{indsh20} et par
\begin{eqnarray}
\bMod(\co_\fX)\rightarrow \bMod_\mQ(\co_\fX),\ \ \ \cF\mapsto \cF_{\mQ},\label{ahttf40a}\\
\bMod(\bvocB)\rightarrow \bMod_\mQ(\bvocB),\ \ \ \cM\mapsto \cM_{\mQ},\label{ahttf40b}
\end{eqnarray}
les foncteurs canoniques. Ce sont deux catégories abéliennes et monoïdales symétriques, ayant $\co_{\fX,\mQ}$ (resp. $\bvocB_\mQ$) pour objet unité.
Les objets de $\bMod_{\mQ}(\bvocB)$ seront aussi appelés des {\em $\bvocB_\mQ$-modules}. 

D'après \ref{indsh11}, on a des foncteurs canoniques \eqref{indsh11c}
\begin{eqnarray}
\upalpha_{\co_\fX}\colon \bMod_\mQ(\co_\fX)&\rightarrow& \bIndMod(\co_\fX),\label{ahttf40c}\\
\upalpha_{\bvocB}\colon \bMod_\mQ(\bvocB)&\rightarrow& \bIndMod(\bvocB).\label{ahttf40d} 
\end{eqnarray}
Ces foncteurs sont pleinement fidèles \eqref{indsh5e} et exacts \eqref{indsh16}.
On identifiera $\bMod_\mQ(\co_\fX)$ (resp. $\bMod_\mQ(\bvocB)$) à une sous-catégorie pleine de $\bIndMod(\co_\fX)$ (resp. $\bIndMod(\bvocB)$)
par le foncteur $\upalpha_{\co_\fX}$ (resp. $\upalpha_{\bvocB}$) qu'on omettra des notations. On considérera donc tout
$\bvocB_\mQ$-module comme un ind-$\bvocB$-module.

D'après \ref{indsh20}, le morphisme de topos annelés $\hupsigma$ \eqref{ahttf13e} induit deux foncteurs additifs adjoints 
\begin{eqnarray}
\hupsigma^*_\mQ\colon \bMod_\mQ(\co_\fX) \rightarrow \bMod_\mQ(\bvocB),\label{ahttf40e}\\
\hupsigma_{\mQ*}\colon \bMod_\mQ(\bvocB) \rightarrow \bMod_\mQ(\co_\fX).\label{ahttf40f}
\end{eqnarray}
Le foncteur $\hupsigma^*_\mQ$ (resp. $\hupsigma_{\mQ*}$) est exact à droite (resp. gauche). 
D'après \ref{indsh13}(ii), les diagrammes
\begin{equation}\label{ahttf40fg}
\xymatrix{
{\bMod_\mQ(\bvocB)}\ar[r]^-(0.5){\hupsigma_{\mQ*}}\ar[d]_{\upalpha_{\bvocB}}&{\bMod_\mQ(\co_\fX)}\ar[d]^{\upalpha_{\co_\fX}}\\
{\bIndMod(\bvocB)}\ar[r]^{\rI \hupsigma_*}&{\bIndMod(\co_\fX)}}
\end{equation}
\begin{equation}\label{ahttf40g}
\xymatrix{
{\bMod_\mQ(\co_\fX)}\ar[r]^-(0.5){\hupsigma^*_\mQ}\ar[d]_{\upalpha_{\co_\fX}}&{\bMod_\mQ(\bvocB)}\ar[d]^{\upalpha_{\bvocB}}\\
{\bIndMod(\co_\fX)}\ar[r]^{\rI \hupsigma^*}&{\bIndMod(\bvocB)}}
\end{equation}
sont commutatifs à isomorphismes canoniques près. 

Le foncteur canonique
\begin{equation}\label{ahttf40j}
\bMod(\co_{\fX})\rightarrow \bMod(\co_{\fX}[\frac 1 p]),\ \ \ \cF\mapsto \cF[\frac 1 p]
\end{equation}
induit un foncteur exact
\begin{equation}\label{ahttf40k}
\bMod_\mQ(\co_{\fX})\rightarrow \bMod(\co_{\fX}[\frac 1 p]). 
\end{equation} 
On note encore 
\begin{equation}\label{ahttf40l}
\hupsigma_{\mQ*}\colon \bMod_\mQ(\bvocB)\rightarrow \bMod(\co_{\fX}[\frac 1 p])
\end{equation}
le composé du foncteur $\hupsigma_{\mQ*}$ \eqref{ahttf40f} et du foncteur \eqref{ahttf40k}.

D'après \ref{indsh14}, le foncteur $\hupsigma_{\mQ*}$ \eqref{ahttf40f} admet un foncteur dérivé à droite
\begin{equation}\label{ahttf40h}
\rR \hupsigma_{\mQ*}\colon \bD^+(\bMod_\mQ(\bvocB))\rightarrow \bD^+(\bMod_\mQ(\co_\fX)).
\end{equation}
En vertu de \eqref{indsh14g}, le diagramme
\begin{equation}\label{ahttf40i}
\xymatrix{
{\bD^+(\bMod_\mQ(\bvocB))}\ar[r]^-(0.5){\rR\hupsigma_{\mQ*}}\ar[d]_{\upalpha_{\bvocB}}&{\bD^+(\bMod_\mQ(\co_\fX))}\ar[d]^{\upalpha_{\co_\fX}}\\
{\bD^+(\bIndMod(\bvocB))}\ar[r]^-(0.5){\rR\rI \hupsigma_*}&{\bD^+(\bIndMod(\co_\fX))}}
\end{equation}
est commutatif à un isomorphisme canonique près. 

De même, le foncteur $\hupsigma_{\mQ*}$ \eqref{ahttf40l} admet un foncteur dérivé à droite
\begin{equation}\label{ahttf40m}
\rR \hupsigma_{\mQ*}\colon \bD^+(\bMod_\mQ(\bvocB))\rightarrow \bD^+(\bMod(\co_\fX[\frac 1 p])),
\end{equation}
qui n'est autre que le composé du foncteur $\rR \hupsigma_{\mQ*}$ \eqref{ahttf40h} et du foncteur exact \eqref{ahttf40k}. 
Cet abus de notation n'induit aucune confusion. Il résulte de \eqref{ahttf40i} que le diagramme
\begin{equation}\label{ahttf40n}
\xymatrix{
{\bD^+(\bMod_\mQ(\bvocB))}\ar[r]^-(0.5){\rR\hupsigma_{\mQ*}}\ar[d]_{\upalpha_{\bvocB}}&{\bD^+(\bMod(\co_\fX[\frac 1 p]))}\ar[d]\\
{\bD^+(\bIndMod(\bvocB))}\ar[r]^-(0.5){\rR\vupsigma_*}&{\bD^+(\bMod(\co_\fX))}}
\end{equation}
où $\rR\vupsigma_*$ est le foncteur \eqref{ahttf43j} et la flèche non libellée est le foncteur canonique, est commutatif à un isomorphisme canonique près. 

\begin{defi}\label{ahttf60}
On dit qu'un $\bvocB_\mQ$-module est {\em adique de type fini} s'il est isomorphe à $\cM_\mQ$ \eqref{ahttf40b}, 
où $\cM$ est un $\bvocB$-module adique de type fini (\cite{agt} III.7.16). 
\end{defi}

On désigne par $\bMod^\atf(\bvocB)$ la sous-catégorie pleine de $\bMod(\bvocB)$ formée des $\bvocB$-modules adiques de type fini 
et par $\bMod^\atf_{\mQ}(\bvocB)$ la catégorie des objets de $\bMod^\atf(\bvocB)$ à isogénie près \eqref{caip1}. Le foncteur canonique 
\begin{equation}\label{ahttf60a}
\bMod^\atf_{\mQ}(\bvocB)\rightarrow\bMod_{\mQ}(\bvocB)
\end{equation} 
étant pleinement fidèle, pour qu'un $\bvocB_\mQ$-module soit adique de type fini, il faut et il suffit qu'il soit dans l'image essentielle de ce foncteur. 

\begin{defi}\label{ahttf49}
On dit qu'un ind-$\bvocB$-module $\cF$ est {\em rationnel} si la multiplication par $p$ sur $\cF$ est un isomorphisme. 
\end{defi}

On notera que l'image essentielle de $\upalpha_{\bvocB}$ est formée d'ind-$\bvocB$-modules rationnels. 

\subsection{}\label{ahttf50}
On désigne par $\bMod^\coh(\co_{\fX})$ (resp. $\bMod^\coh(\co_{\fX}[\frac 1 p])$) 
la catégorie des $\co_{\fX}$-modules (resp. $\co_{\fX}[\frac 1 p]$-modules) cohérents de $X_{s,\zar}$
et par $\bMod^{\coh}_\mQ(\co_{\fX})$ la catégorie des $\co_\fX$-modules cohérents à isogénie près (\cite{agt} III.6.1.1). 
D'après (\cite{agt} III.6.16), le foncteur canonique
\begin{equation}\label{ahttf50a}
\bMod^{\coh}(\co_{\fX})\rightarrow \bMod^{\coh}(\co_{\fX}[\frac 1 p]),\ \ \ \cF\mapsto \cF[\frac 1 p]
\end{equation}
induit une équivalence de catégories abéliennes
\begin{equation}\label{ahttf50b}
\bMod^{\coh}_\mQ(\co_{\fX})\stackrel{\sim}{\rightarrow} \bMod^{\coh}(\co_{\fX}[\frac 1 p]). 
\end{equation} 

Considérons le foncteur composé 
\begin{equation}\label{ahttf50c}
\bMod^{\coh}(\co_{\fX}[\frac 1 p])\stackrel{\sim}{\longrightarrow} \bMod^{\coh}_\mQ(\co_{\fX}) \longrightarrow \bMod_\mQ(\co_\fX)
\stackrel{\upalpha_{\co_\fX}}{\longrightarrow} \bIndMod(\co_\fX), 
\end{equation}
où la première flèche est un quasi-inverse de l'équivalence de catégories \eqref{ahttf50b} 
et la seconde flèche est le foncteur pleinement fidèle canonique. 
On identifiera $\bMod^{\coh}(\co_{\fX}[\frac 1 p])$ à une sous-catégorie pleine de $\bIndMod(\co_\fX)$ par ce foncteur composé. 
On considérera donc tout $\co_\fX[\frac 1 p]$-module cohérent aussi comme un ind-$\co_\fX$-module. 

Pour tout $\co_\fX[\frac 1 p]$-module cohérent $\cN$, considéré comme un ind-$\co_\fX$-module, on a un isomorphisme canonique fonctoriel 
\begin{equation}\label{ahttf50d}
\kappa_{\co_\fX}(\cN)\stackrel{\sim}{\rightarrow} \cN.
\end{equation}
Cela résulte aussitôt des définitions et de \ref{indsh16}(i). 

Soient $\cN$ un $\co_{\fX}[\frac 1 p]$-module cohérent, $\cM$ un ind-$\bvocB$-module. Les foncteurs adjoints $\rI\hupsigma^*$ et $\rI\hupsigma_*$
et l'isomorphisme \eqref{ahttf50d} induisent une application 
\begin{equation}\label{ahttf50g}
\Hom_{\bIndMod(\bvocB)}(\rI\hupsigma^*(\cN),\cM)\rightarrow \Hom_{\bMod(\co_\fX)}(\cN,\vupsigma_*(\cM)).
\end{equation}
Celle-ci est injective d'après \ref{ahttf53} ci-dessous. 
Compte tenu de \eqref{indsh15g} et \eqref{ahttf50d}, pour tout entier $q\geq 0$, le morphisme canonique \eqref{indsh40a} 
\begin{equation}\label{ahttf50e}
\cN\otimes_{\co_\fX}\rR^q\rI \hupsigma_*(\cM)\rightarrow\rR^q\rI \hupsigma_*(\rI\hupsigma^*(\cN)\otimes_{\bvocB}\cM)
\end{equation}
induit par application du foncteur $\kappa_{\co_\fX}$ un morphisme 
\begin{equation}\label{ahttf50f}
\cN\otimes_{\co_\fX}\rR^q \vupsigma_*(\cM)\rightarrow\rR^q\vupsigma_*(\rI\hupsigma^*(\cN)\otimes_{\bvocB}\cM).
\end{equation}

\begin{lem}\label{ahttf53}
Soient $\cN$ un $\co_{\fX}[\frac 1 p]$-module cohérent, $\cF$ un ind-$\co_\fX$-module. Alors, le foncteur $\kappa_{\co_\fX}$ induit une application
injective 
\begin{equation}
\Hom_{\bIndMod(\co_\fX)}(\cN,\cF)\rightarrow \Hom_{\bMod(\co_\fX)}(\cN,\kappa_{\co_\fX}(\cF)). 
\end{equation}
\end{lem}

Soient $\cN^\circ$ un $\co_\fX$-module cohérent tel que $\cN^\circ[\frac 1p]=\cN$ \eqref{ahttf50a},
$\alpha\colon J\rightarrow \bMod(\co_\fX)$ un foncteur avec $J$ une petite catégorie filtrante tel que $\cF=\indcolim\alpha$. 
Compte tenu de \eqref{indsh1c}, il s'agit de montrer que l'application canonique 
\begin{equation}
\underset{\underset{\mN}{\longleftarrow}}{\lim}\ \underset{\underset{j\in J}{\longrightarrow}}{\lim}\ \Hom_{\co_\fX}(\cN^\circ,\alpha(j))
\rightarrow \Hom_{\co_\fX}(\cN, \underset{\underset{j\in J}{\longrightarrow}}{\lim}\ \alpha(j)),
\end{equation}
où les morphismes de transition de la limite projective sont induits par la multiplication par les puissances de $p$, 
est injective. La question étant locale, on peut supposer $\cN^\circ$ engendré par ses sections globales, et par suite se borner au cas où
$\cN^\circ$ est libre de type fini et même au cas où $\cN^\circ=\co_\fX$. D'après (\cite{sga4} VI 5.3), le morphisme canonique 
\begin{equation}
\underset{\underset{j\in J}{\longrightarrow}}{\lim}\ \Hom_{\co_\fX}(\co_\fX,\alpha(j))
\rightarrow \Hom_{\co_\fX}(\co_\fX, \underset{\underset{j\in J}{\longrightarrow}}{\lim}\ \alpha(j))
\end{equation}
est un isomorphisme. Par ailleurs, le morphisme canonique 
\begin{equation}
\underset{\underset{\mN}{\longleftarrow}}{\lim}\ \Hom_{\co_\fX}(\co_\fX, \underset{\underset{j\in J}{\longrightarrow}}{\lim}\ \alpha(j))
\rightarrow \Hom_{\co_\fX}(\co_\fX[\frac 1 p], \underset{\underset{j\in J}{\longrightarrow}}{\lim}\ \alpha(j)),
\end{equation}
où les morphismes de transition de la limite projective sont induits par la multiplication par les puissances de $p$, 
est un isomorphisme, d'où l'assertion recherchée.

\begin{lem}\label{ahttf51}
Soient $\cN$ un $\co_{\fX}[\frac 1 p]$-module localement projectif de type fini \eqref{notconv14}, $\cM$ un ind-$\bvocB$-module, $q$ un entier $\geq 0$.
Alors, le morphisme canonique \eqref{ahttf50e} 
\begin{equation}\label{ahttf51a}
\cN\otimes_{\co_\fX}\rR^q\rI \hupsigma_*(\cM)\rightarrow\rR^q\rI \hupsigma_*(\rI\hupsigma^*(\cN)\otimes_{\bvocB}\cM)
\end{equation}
est un isomorphisme. Il en est de même du morphisme induit \eqref{ahttf50f}
\begin{equation}\label{ahttf51b}
\cN\otimes_{\co_\fX}\rR^q\vupsigma_*(\cM)\rightarrow\rR^q\vupsigma_*(\rI\hupsigma^*(\cN)\otimes_{\bvocB}\cM).
\end{equation}
\end{lem}
Il existe un recouvrement ouvert de Zariski $(U_i)_{0\leq i\leq n}$ de $X$ tel que pour tout $0\leq i\leq n$, la restriction 
de $\cN$ à $(U_i)_s$ soit un facteur direct d'un $(\co_\fX|U_i)[\frac 1 p]$-module libre de type fini. 
Compte tenu de \ref{indsh13}(ii), \ref{indsh39}(iii) et du fait que $\rR^q\rI \hupsigma_*$ commute aux localisations \eqref{indsh42}, 
on peut alors se borner au cas où $\cN$ est un facteur direct
d'un $\co_\fX[\frac 1 p]$-module libre de type fini, et même au cas où $\cN$ est un 
$\co_\fX[\frac 1 p]$-module libre de type fini. Comme $\rR^q\rI \hupsigma_*$ commute aux limites inductives \eqref{indsh21g}, 
le morphisme \eqref{ahttf51a} est un isomorphisme. Il s'ensuit aussitôt que le morphisme \eqref{ahttf51b} est un isomorphisme.

\section{Algèbres de Higgs-Tate dans le topos de Faltings}\label{ahttff}

\subsection{}\label{ahttf4}
On désigne par $\bP$ la sous-catégorie pleine de $\Et_{/X}$ formée des schémas affines $U$ tels que l'une des deux conditions suivantes 
soit remplie:
\begin{itemize}
\item[(i)] le schéma $U_s$ est vide; ou 
\item[(ii)] le morphisme $(U,\cM_X|U)\rightarrow (S,\cM_S)$ induit par $f$ \eqref{defing11} admet une carte adéquate \eqref{cad1}. 
\end{itemize}

On munit $\bP$ de la topologie induite par celle de $\Et_{/X}$. Comme $X$ est noethérien et donc quasi-séparé, 
tout objet de $\bP$ est cohérent sur $X$. Par suite, $\bP$ est 
une famille $\mU$-petite, topologiquement génératrice du site $\Et_{/X}$ et est stable par produits fibrés. 

On désigne par 
\begin{equation}\label{ahttf4a}
\pi_\bP\colon E_\bP\rightarrow \bP
\end{equation} 
le site fibré déduit de $\pi$ \eqref{ahttf2a} par changement de base par le foncteur d'injection canonique $\bP\rightarrow \Et_{/X}$, par 
\begin{equation}\label{ahttf4c}
\cP^\vee_\bP\rightarrow \bP^\circ
\end{equation}
la catégorie fibrée sur $\bP^\circ$ déduite de la catégorie fibrée $\cP^\vee$ \eqref{ahttf2i} 
par changement de base par le foncteur d'injection canonique $\bP\rightarrow \Et_{/X}$, 
et par $\hE_\bP$ la catégorie des préfaisceaux de $\mU$-ensembles sur $E_\bP$.
Pour tout $U\in \ob(\bP)$, on note $\iota_U\colon \Et_{\rf/\oU^\circ}\rightarrow E_\bP$ 
le foncteur canonique \eqref{ahttf2b}. On a une équivalence de catégories 
\begin{eqnarray}\label{ahttf4d}
\hE_\bP&\stackrel{\sim}{\rightarrow}& \bHom_{\bP^\circ}(\bP^\circ,\cP^\vee_\bP)\\
F&\mapsto &\{U\mapsto F\circ \iota_U\}.\nonumber
\end{eqnarray}
On identifiera dans la suite $F$ à la section $\{U\mapsto F\circ \iota_U\}$ qui lui est associée par cette équivalence.

On munit $E_\bP$ de la topologie co-évanescente définie par $\pi_\bP$ (\cite{agt} VI.5.3)  et on note $\tE_\bP$ le topos des faisceaux de $\mU$-ensembles 
sur $E_\bP$. D'après (\cite{agt} VI.5.21 et VI.5.22), la topologie de $E_\bP$ est induite par celle de $E$ au moyen du foncteur de projection canonique 
$E_\bP\rightarrow E$, et celui-ci induit par restriction une équivalence de catégories 
\begin{equation}\label{ahttf4b}
\tE\stackrel{\sim}{\rightarrow}\tE_\bP.
\end{equation}

\subsection{}\label{ahttf5}
On désigne par $\bQ$ la sous-catégorie pleine de $\bP$ formée des schémas affines $U$ tels que l'une des conditions suivantes soit remplie:
\begin{itemize}
\item[(i)] le schéma $U_s$ est vide; ou 
\item[(ii)] il existe une carte fine et saturée $M\rightarrow \Gamma(U,\cM_X)$ pour $(U,\cM_X|U)$
induisant un isomorphisme 
\begin{equation}\label{ahttf5a}
M\stackrel{\sim}{\rightarrow} \Gamma(U,\cM_X)/\Gamma(U,\co^\times_X).
\end{equation}
Cette carte est a priori indépendante de la carte adéquate requise dans \ref{ahttf4}(ii). 
\end{itemize}

On munit $\bQ$ de la topologie induite par celle de $\Et_{/X}$. 
Il résulte de (\cite{agt} II.5.17) que $\bQ$ est une sous-catégorie 
topologiquement génératrice de $\Et_{/X}$. 

On désigne par 
\begin{equation}\label{ahttf5b}
\pi_\bQ\colon E_\bQ\rightarrow \bQ
\end{equation}
le site fibré déduit de $\pi$ \eqref{ahttf2a}
par changement de base par le foncteur d'injection canonique $\bQ\rightarrow \Et_{/X}$, par 
\begin{equation}\label{ahttf5c}
\cP^\vee_\bQ\rightarrow \bQ^\circ
\end{equation}
la catégorie fibrée sur $\bQ^\circ$ déduite de la catégorie fibrée $\cP^\vee$ \eqref{ahttf2i} 
par changement de base par le foncteur d'injection canonique $\bQ\rightarrow \Et_{/X}$, 
et par $\hE_\bQ$ la catégorie des préfaisceaux de $\mU$-ensembles sur $E_\bQ$.
Pour tout $U\in \ob(\bQ)$, on note $\iota_U\colon \Et_{\rf/\oU^\circ}\rightarrow E_\bQ$ 
le foncteur canonique \eqref{ahttf2b}. On a une équivalence de catégories 
\begin{eqnarray}\label{ahttf5d}
\hE_\bQ&\stackrel{\sim}{\rightarrow}& \bHom_{\bQ^\circ}(\bQ^\circ,\cP^\vee_\bQ)\\
F&\mapsto &\{U\mapsto F\circ \iota_U\}.\nonumber
\end{eqnarray}
On identifiera dans la suite $F$ à la section $\{U\mapsto F\circ \iota_U\}$ qui lui est associée par cette équivalence.

Le foncteur de projection canonique $E_\bQ\rightarrow E$ est pleinement fidèle et 
la catégorie $E_\bQ$ est $\mU$-petite et topologiquement génératrice du site $E$. On munit $E_\bQ$ de la 
topologie induite par celle de $E$. Par restriction, le topos $\tE$ est alors équivalent à la catégorie 
des faisceaux de $\mU$-ensembles sur $E_\bQ$ (\cite{sga4} III 4.1). 
On prendra garde qu'en général, $\bQ$ n'étant pas stable par produits fibrés, 
on ne peut pas parler de topologie co-évanescente sur $E_\bQ$ associée à $\pi_\bQ$, 
et encore moins appliquer  (\cite{agt} VI.5.21 et VI.5.22). 

\begin{rema}\label{ahttf45}
\
\begin{itemize}
\item[(i)] Les sous-catégories $\bP$ et $\bQ$ de $\Et_{/X}$ ont été introduites dans (\cite{agt} III.10.3 et III.10.5), mais en omettant d'y inclure 
les objets $U$ tels que $U_s$ soit vide, ce qui est toutefois nécessaire pour que ces sous-catégories soient topologiquement génératrices. 
Nous réparons cette omission dans \ref{ahttf4} et \ref{ahttf5}. Cet erratum n'a aucune conséquence sur la suite de (\cite{agt} III). 
\item[(ii)] Contrairement à (\cite{agt} III.10.5), nous ne supposons pas les schémas dans $\bQ$ connexes. 
\end{itemize}
\end{rema}

\subsection{}\label{ahttf6}
Soient $Y$ un objet de $\bP$ tel que $Y_s$ soit non vide, $((P,\gamma),(\mN,\iota),\vartheta)$
une carte adéquate pour le morphisme $f|Y\colon (Y,\cM_X|Y)\rightarrow (S,\cM_S)$ induit par $f$, 
$\oy$ un point géométrique de $\oY^\circ$. On notera que les hypothèses de \ref{cad1} sont remplies.
Le schéma $\oY$ étant localement irréductible d'après (\cite{ag} 4.2.7(iii) et \cite{agt} III.3.3),  
il est la somme des schémas induits sur ses composantes irréductibles. On note $\oY^\star$
la composante irréductible de $\oY$ contenant $\oy$. 
De même, $\oY^\circ$ est la somme des schémas induits sur ses composantes irréductibles
et $\oY^{\star \circ}=\oY^\star\times_{X}X^\circ$ est la composante irréductible de $\oY^\circ$ contenant $\oy$. 
On désigne par $\oR^\oy_Y$ la représentation discrète de $\pi_1(\oY^{\star\circ},\oy)$ définie dans \eqref{ahttf44b}
et par $\hoR^\oy_Y$ son séparé complété $p$-adique. On pose 
\begin{equation}\label{ahttf6a}
\mY^\oy=\Spec(\oR^\oy_Y)\ \ \ {\rm et} \ \ \ \hmY^\oy=\Spec(\hoR^\oy_Y)
\end{equation}
que l'on munit des structures logarithmiques images inverses de $\cM_X$, 
notées respectivement $\cM_{\mY^\oy}$ et $\cM_{\hmY^\oy}$. Selon le cas considéré dans \ref{definf10}, 
on désigne par $(\tmY^\oy,\cM_{\tmY^\oy})$ l'un ou l'autre des schémas logarithmiques 
\begin{equation}\label{ahttf6b}
(\cA_2(\mY^\oy),\cM_{\cA_2(\mY^\oy)})\ \ \ {\rm ou} \ \ \ (\cA^{\ast}_2(\mY^\oy/S),\cM_{\cA^{\ast}_2(\mY^\oy/S)})
\end{equation}
associés au morphisme $f|Y$ et à la carte adéquate $((P,\gamma),(\mN,\iota),\vartheta)$ dans \ref{taht2} et \ref{taht3}. 
On a une immersion fermée exacte canonique \eqref{taht4b}
\begin{equation}\label{ahttf6c}
i^\oy_Y\colon (\hmY^\oy,\cM_{\hmY^\oy})\rightarrow (\tmY^\oy,\cM_{\tmY^\oy}).
\end{equation}

On pose \eqref{defing11b}
\begin{equation}\label{ahttf6d}
\rT^\oy_Y=\Hom_{\hoR^\oy_Y}(\tOmega^1_{X/S}(Y)\otimes_{\co_X(Y)}\hoR^\oy_Y,\txi\hoR^\oy_Y).
\end{equation} 
On identifie le $\hoR^\oy_Y$-module dual à $\txi^{-1}\tOmega^1_{X/S}(Y)\otimes_{\co_X(Y)}\hoR^\oy_Y$ (cf. \ref{definf10}).
On désigne par $\hmY^\oy_\zar$ le topos de Zariski de $\hmY^\oy$ et par $\trT^\oy_Y$ le $\co_{\hmY^\oy}$-module associé à $\rT^\oy_Y$. 
Soient $U$ un ouvert de Zariski de $\hmY^\oy$, $\tU$ l'ouvert correspondant de $\tmY^\oy$. On désigne par $\cL^\oy_Y(U)$ 
l'ensemble des morphismes représentés par des flèches pointillées qui complètent  le diagramme canonique \eqref{defing12b}
\begin{equation}\label{ahttf6e}
\xymatrix{
{(U,\cM_{\hmY^\oy}|U)}\ar[rr]^-(0.5){i^\oy_Y\times_{\tmY^\oy}\tU}\ar[d]&&{(\tU,\cM_{\tmY^\oy}|\tU)}\ar@{.>}[d]\ar@/^2pc/[dd]\\
{(\coX,\cM_{\coX})}\ar[rr]\ar[d]&&{(\tX,\cM_\tX)}\ar[d]\\
{(\coS,\cM_\coS)}\ar[rr]^-(0.5){i_S}&&{(\tS,\cM_{\tS})}}
\end{equation}
de façon à le laisser commutatif. D'après (\cite{agt} II.5.23), le foncteur $U\mapsto \cL^\oy_Y(U)$ est un $\trT^\oy_Y$-torseur de $\hmY_\zar$. 
Notant $\tY\rightarrow \tX$ l'unique morphisme étale qui relève $\coY\rightarrow \coX$ (\cite{ega4} 18.1.2) (cf. \ref{definf12}),  
$\cL^\oy_Y$ s'identifie au torseur de Higgs-Tate associé à $(\tY,\cM_\tX|\tY)$ \eqref{taht5}. 
On note $\cF^\oy_Y$ le $\hoR^\oy_Y$-module des fonctions affines sur $\cL^\oy_Y$ (cf. \cite{agt} II.4.9).
Celui-ci s'insère dans une suite exacte canonique \eqref{taht5e}
\begin{equation}\label{ahttf6f}
0\rightarrow \hoR^\oy_Y\rightarrow \cF^\oy_Y\rightarrow \txi^{-1}\tOmega^1_{X/S}(Y) \otimes_{\co_X(Y)}\hoR^\oy_Y\rightarrow 0.
\end{equation} 
On désigne par $\cC^\oy_Y$ la $\hoR^\oy_Y$-algèbre \eqref{taht5g} 
\begin{equation}\label{ahttf6g}
\cC^\oy_Y=\underset{\underset{n\geq 0}{\longrightarrow}}\lim\ \rS^n_{\hoR^\oy_Y}(\cF^\oy_Y).
\end{equation}

D'après \ref{taht6}, le $\hoR^\oy_Y$-module  $\cF^\oy_Y$ est canoniquement muni d'une action $\hoR^\oy_Y$-semi-linéaire de $\pi_1(\oY^{\star\circ},\oy)$, 
qui est continue pour la topologie $p$-adique. Les morphismes de la suite \eqref{ahttf6f} 
sont $\pi_1(\oY^{\star\circ},\oy)$-équivariants; on retrouve ainsi l'extension de Higgs-Tate associée à $(\tY,\cM_\tX|\tY)$ \eqref{taht7}.
On en déduit une action de $\pi_1(\oY^{\star\circ},\oy)$ sur  $\cC^\oy_Y$ par des homomorphismes d'anneaux, 
qui est continue pour la topologie $p$-adique et qui prolonge l'action canonique sur $\hoR^\oy_Y$; on retrouve ainsi l'algèbre 
de Higgs-Tate associée à $(\tY,\cM_\tX|\tY)$ \eqref{taht7}. On notera que ces représentations dépendent de 
la carte adéquate $((P,\gamma),(\mN,\iota),\vartheta)$. Toutefois, ils n'en dépendent pas si $Y$ est un objet de $\bQ$ d'après \ref{pmh7}.

\subsection{}\label{ahttf7}
Soient $Y$ un objet de $\bP$ tel que $Y_s$ soit non vide, $((P,\gamma),(\mN,\iota),\vartheta)$
une carte adéquate pour le morphisme $f|Y\colon (Y,\cM_X|Y)\rightarrow (S,\cM_S)$ induit par $f$. 
Si $A$ est un anneau et $M$ un $A$-module, on note encore $A$ (resp. $M$) le faisceau constant de valeur $A$ (resp. $M$) de $\oY^\circ_\fet$. 
On rappelle que $\oY^\circ$ est la somme des schémas induits sur ses composantes irréductibles \eqref{ahttf6}.
Soient $W$ une composante irréductible de $\oY^\circ$, $\Pi(W)$ son groupoïde fondamental (\cite{agt} VI.9.10). 
Compte tenu de \eqref{ahttf44b} et (\cite{agt} VI.9.11), le faisceau $\ocB_Y|W$ \eqref{TFA2d} de $W_\fet$ définit un foncteur 
\begin{equation}\label{ahttf7a}
\Pi(W)\rightarrow  \Ens, \ \ \ \oy\mapsto \oR^\oy_Y.
\end{equation}
On en déduit un foncteur 
\begin{equation}\label{ahttf7b}
\Pi(W)\rightarrow \Ens, \ \ \ \oy\mapsto \cF^\oy_Y/p^n\cF^\oy_Y.
\end{equation}
Pour tout point géométrique $\oy$ de $W$, $\cF^\oy_Y/p^n\cF^\oy_Y$ est une représentation discrète et continue de $\pi_1(W,\oy)$.
Par suite, en vertu de (\cite{agt} VI.9.11), le foncteur \eqref{ahttf7b} définit un $(\ocB_{Y,n}|W)$-module 
$\cF_{W,n}$ de $W_\fet$, unique à isomorphisme canonique près, où $\ocB_{Y,n}=\ocB_Y/p^n\ocB_Y$ \eqref{ahttf3e}. 
Par descente (\cite{giraud2} II 3.4.4), il existe un $\ocB_{Y,n}$-module $\cF_{Y,n}$ de $\oY^\circ_\fet$, unique 
à isomorphisme canonique près, tel que pour toute composante irréductible $W$ de $\oY^\circ$, on ait $\cF_{Y,n}|W=\cF_{W,n}$.

La suite exacte \eqref{ahttf6f} induit une suite exacte de $\ocB_{Y,n}$-modules 
\begin{equation}\label{ahttf7c}
0\rightarrow \ocB_{Y,n}\rightarrow \cF_{Y,n}\rightarrow 
\txi^{-1}\tOmega^1_{X/S}(Y)\otimes_{\co_X(Y)}\ocB_{Y,n} \rightarrow 0. 
\end{equation}
D'après (\cite{illusie1} I 4.3.1.7), celle-ci induit pour tout entier $m\geq 1$, une suite exacte
\[
0\rightarrow \rS^{m-1}_{\ocB_{Y,n}}(\cF_{Y,n})\rightarrow \rS^{m}_{\ocB_{Y,n}}(\cF_{Y,n})\rightarrow 
\rS^m_{\ocB_{Y,n}}(\txi^{-1}\tOmega^1_{X/S}(Y)\otimes_{\co_X(Y)} \ocB_{Y,n})\rightarrow 0.
\]
Les $\ocB_{Y,n}$-modules $(\rS^{m}_{\ocB_{Y,n}}(\cF_{Y,n}))_{m\in \mN}$ forment donc un système inductif 
dont la limite inductive 
\begin{equation}\label{ahttf7d}
 \cC_{Y,n}=\underset{\underset{m\geq 0}{\longrightarrow}}\lim\ \rS^m_{\ocB_{Y,n}}(\cF_{Y,n})
\end{equation}
est naturellement munie d'une structure de $\ocB_{Y,n}$-algèbre de $\oY^\circ_\fet$. 

On notera que $\cF_{Y,n}$ et $\cC_{Y,n}$ dépendent du choix la carte adéquate $((P,\gamma),(\mN,\iota),\vartheta)$. 
Toutefois, ils n'en dépendent pas si $Y$ est un objet de $\bQ$ d'après \ref{pmh7}.
Par ailleurs, $\cF_{Y,n}$ et $\cC_{Y,n}$ dépendent du choix de la déformation $(\tX,\cM_{\tX})$ fixée dans \ref{defing12}.

\subsection{}\label{ahttf8}
Soient $g\colon Y\rightarrow Z$ un morphisme de $\bP$ tel que $Y_s$ soit non vide, 
$((P,\gamma),(\mN,\iota),\vartheta)$ une carte adéquate pour le morphisme $f|Z\colon (Z,\cM_X|Z)\rightarrow (S,\cM_S)$ induit par $f$, 
$\oy$ un point géométrique de $\oY^\circ$, $\oz=\ogg(\oy)$. 
On munit le morphisme $f|Y\colon (Y,\cM_X|Y)\rightarrow (S,\cM_S)$ induit par $f$ de la même carte adéquate $((P,\gamma),(\mN,\iota),\vartheta)$. 
On rappelle que $\oY$ et $\oZ$ sont sommes des schémas induits sur leurs composantes irréductibles. 
On note $\oY^\star$ la composante irréductible de $\oY$ contenant $\oy$ et 
$\oZ^\star$ la composante irréductible de $\oZ$ contenant $\oz$, de sorte que $\ogg(\oY^\star)\subset \oZ^\star$.
Le morphisme $\ogg^\circ\colon \oY^\circ\rightarrow \oZ^\circ$ 
induit un homomorphisme de groupes $\pi_1(\oY^{\star \circ},\oy) \rightarrow\pi_1(\oZ^{\star \circ},\oz)$. 
Le morphisme canonique $(\ogg^\circ)^*_\fet(\ocB_Z)\rightarrow \ocB_Y$ 
induit un homomorphisme d'anneaux $\pi_1(\oY^{\star \circ},\oy)$-équivariant 
\begin{equation}\label{ahttf8a}
\oR^\oz_Z\rightarrow \oR^{\oy}_Y,
\end{equation} 
et par suite un morphisme $\pi_1(\oY^{\star \circ},\oy)$-équivariant de schémas 
$h\colon \hmY^{\oy}\rightarrow \hmZ^\oz$. Comme $g$ est étale, on a un morphisme canonique 
$\co_{\hmZ^\oz}$-linéaire et $\pi_1(\oY^{\star \circ},\oy)$-équivariant $u\colon \trT^\oz_Z\rightarrow h_*(\trT^\oy_Y)$
tel que le morphisme adjoint $u^\sharp\colon h^*(\trT^\oz_Z)\rightarrow \trT^\oy_Y$ soit un isomorphisme.
Il résulte aussitôt des définitions \eqref{ahttf6e}
qu'on a un morphisme canonique $u$-équivariant et $\pi_1(\oY^{\star \circ},\oy)$-équivariant
\begin{equation}\label{ahttf8k}
v\colon \cL_Z^\oz\rightarrow h_*(\cL_Y^\oy).
\end{equation}
D'après (\cite{agt} II.4.22), le couple $(u,v)$ induit un isomorphisme $\hoR^\oz_Y$-linéaire et $\pi_1(\oY^{\star \circ},\oy)$-équivariant  
\begin{equation}\label{ahttf8b}
\cF^\oy_Y\stackrel{\sim}{\rightarrow} \cF^\oz_Z\otimes_{\hoR^\oz_Z}\hoR^\oy_Y,
\end{equation}
et par suite, un morphisme $\hoR^\oz_Z$-linéaire et $\pi_1(\oY^{\star \circ},\oy)$-équivariant 
\begin{equation}\label{ahttf8c}
\cF^\oz_Z\rightarrow \cF^\oy_Y
\end{equation}
qui s'insère dans un diagramme commutatif 
\begin{equation}\label{ahttf8d}
\xymatrix{
0\ar[r]&{\hoR^\oz_Z}\ar[r]\ar[d]&{\cF^\oz_Z}\ar[r]\ar[d]&
{\txi^{-1}\tOmega^1_{X/S}(Z)\otimes_{\co_X(Z)}\hoR^\oz_Z}\ar[r]\ar[d]&0\\
0\ar[r]&{\hoR^\oy_Y}\ar[r]&{\cF^\oy_Y}\ar[r]&
{\txi^{-1}\tOmega^1_{X/S}(Y)\otimes_{\co_X(Y)}\hoR^\oy_Y}\ar[r]&0}
\end{equation}
On en déduit un homomorphisme $\pi_1(\oY^{\star \circ},\oy)$-équivariant de $\hoR^\oz_Z$-algèbres
\begin{equation}\label{ahttf8e}
\cC^{\oz}_Z\rightarrow\cC^\oy_Y.
\end{equation}

On désigne par $\Pi(\oY^{\star \circ})$ et $\Pi(\oZ^{\star \circ})$ les groupoïdes fondamentaux de $\oY^{\star \circ}$ et 
$\oZ^{\star \circ}$ et par
\begin{equation}\label{ahttf8l}
\gamma\colon \Pi(\oY^{\star \circ})\rightarrow \Pi(\oZ^{\star \circ})
\end{equation}
le foncteur induit par le foncteur image inverse $\Et_{\rf/\oZ^{\star \circ}}\rightarrow \Et_{\rf/\oY^{\star \circ}}$.
Pour tout entier $n\geq 0$, on note $F_{Y,n}\colon \Pi(\oY^{\star \circ})\rightarrow \Ens$ et 
$F_{Z,n}\colon \Pi(\oZ^{\star \circ})\rightarrow \Ens$ les foncteurs associés par 
(\cite{agt} VI.9.11) aux objets $\cF_{Y,n}|\oY^{\star \circ}$ de $\oY^{\star \circ}_\fet$ et 
$\cF_{Z,n}|\oZ^{\star \circ}$ de $\oZ^{\star \circ}_\fet$, respectivement. 
Le morphisme \eqref{ahttf8c} induit clairement un morphisme de foncteurs
\begin{equation}\label{ahttf8m}
F_{Z,n}\circ \gamma \rightarrow F_{Y,n}. 
\end{equation} 
On en déduit par (\cite{agt} VI.9.11) un morphisme $(\ogg^{\circ})_\fet^*(\ocB_{Z,n})$-linéaire 
\begin{equation}\label{ahttf8f}
(\ogg^{\circ})^*_\fet(\cF_{Z,n}) \rightarrow\cF_{Y,n},
\end{equation} 
et donc par adjonction, un morphisme $\ocB_{Z,n}$-linéaire 
\begin{equation}\label{ahttf8g}
\cF_{Z,n} \rightarrow \ogg^{\circ}_{\fet *}(\cF_{Y,n}). 
\end{equation} 
Il résulte de \eqref{ahttf8d} que le diagramme 
\begin{equation}\label{ahttf8h}
\xymatrix{
0\ar[r]&{(\ogg^{\circ})^*_\fet(\ocB_{Z,n})}\ar[r]\ar[d]&{(\ogg^{\circ})^*_\fet(\cF_{Z,n})}\ar[r]\ar[d]&
{\txi^{-1}\tOmega^1_{X/S}(Z)\otimes_{\co_X(Z)}(\ogg^{\circ})^*_\fet(\ocB_{Z,n})}\ar[r]\ar[d]&0\\
0\ar[r]&{\ocB_{Y,n}}\ar[r]&{\cF_{Y,n}}\ar[r]&
{\txi^{-1}\tOmega^1_{X/S}(Y)\otimes_{\co_X(Y)}\ocB_{Y,n}}\ar[r]&0}
\end{equation}
est commutatif. 
On en déduit un homomorphisme de $(\ogg^{\circ})_\fet^*(\ocB_{Z,n})$-algèbres 
\begin{equation}\label{ahttf8i}
(\ogg^{\circ})^*_\fet(\cC_{Z,n}) \rightarrow\cC_{Y,n},
\end{equation} 
et donc par adjonction un homomorphisme de $\ocB_{Z,n}$-algèbres
\begin{equation}\label{ahttf8j}
\cC_{Z,n} \rightarrow\ogg^{\circ}_{\fet *}(\cC_{Y,n}). 
\end{equation} 

On notera que les morphismes \eqref{ahttf8g} et \eqref{ahttf8j} dépendent du choix la carte adéquate $((P,\gamma),(\mN,\iota),\vartheta)$. 
Toutefois, ils n'en dépendent pas si $Y$ et $Z$ sont des objets de $\bQ$ d'après \ref{pmh7}.

\subsection{}\label{ahttf35}
Reprenons les hypothèses et notations de \ref{ahttf7}. Soient, de plus, $r$ un nombre rationnel $\geq 0$, $n$ un entier $\geq 0$.   
On désigne par $\cF^{(r)}_{Y,n}$ l'extension de $\ocB_{Y,n}$-modules de $\oY^\circ_\fet$ déduite de 
$\cF_{Y,n}$ \eqref{ahttf7c} par image inverse 
par le morphisme de multiplication par $p^r$ sur $\txi^{-1}\tOmega^1_{X/S}(Y)\otimes_{\co_X(Y)}\ocB_{Y,n}$,
de sorte qu'on a une suite exacte canonique de $\ocB_{Y,n}$-modules
\begin{equation}\label{ahttf35a}
0\rightarrow \ocB_{Y,n}\rightarrow \cF^{(r)}_{Y,n}\rightarrow 
\txi^{-1}\tOmega^1_{X/S}(Y)\otimes_{\co_X(Y)}\ocB_{Y,n} \rightarrow 0. 
\end{equation}
Celle-ci induit pour tout entier $m\geq 1$, une suite exacte de $\ocB_{Y,n}$-modules
\[
0\rightarrow \rS^{m-1}_{\ocB_{Y,n}}(\cF^{(r)}_{Y,n})\rightarrow \rS^{m}_{\ocB_{Y,n}}(\cF^{(r)}_{Y,n})\rightarrow 
\rS^m_{\ocB_{Y,n}}(\txi^{-1}\tOmega^1_{X/S}(Y)\otimes_{\co_X(Y)} \ocB_{Y,n})\rightarrow 0.
\]
Les $\ocB_{Y,n}$-modules $(\rS^{m}_{\ocB_{Y,n}}(\cF^{(r)}_{Y,n}))_{m\in \mN}$ forment donc un système inductif 
dont la limite inductive 
\begin{equation}\label{ahttf35b}
 \cC^{(r)}_{Y,n}=\underset{\underset{m\geq 0}{\longrightarrow}}\lim\ \rS^m_{\ocB_{Y,n}}(\cF^{(r)}_{Y,n})
\end{equation}
est naturellement munie d'une structure de $\ocB_{Y,n}$-algèbre de $\oY^\circ_\fet$. 

Pour tous nombres rationnels $r\geq r'\geq 0$, on a un morphisme $\ocB_{Y,n}$-linéaire canonique 
\begin{equation}\label{ahttf35c}
\tta_{Y,n}^{r,r'}\colon\cF_{Y,n}^{(r)}\rightarrow \cF_{Y,n}^{(r')}
\end{equation} 
qui relève la multiplication par $p^{r-r'}$ sur  
$\txi^{-1}\tOmega^1_{X/S}(Y)\otimes_{\co_X(Y)}\ocB_{Y,n}$ et qui étend l'identité sur $\ocB_{Y,n}$ \eqref{ahttf35a}. 
Il induit un homomorphisme de $\ocB_{Y,n}$-algèbres 
\begin{equation}\label{ahttf35d}
\alpha_{Y,n}^{r,r'}\colon \cC_{Y,n}^{(r)}\rightarrow \cC_{Y,n}^{(r')}.
\end{equation}

\subsection{}\label{ahttf36}
Reprenons les hypothèses et notations de \ref{ahttf8}. 
Soient, de plus, $r$ un nombre rationnel $\geq 0$, $n$ un entier $\geq 0$.
Le diagramme \eqref{ahttf8h} induit un morphisme $(\ogg^{\circ})^*_\fet(\ocB_{Z,n})$-linéaire
\begin{equation}\label{ahttf36a}
(\ogg^{\circ})^*_\fet(\cF^{(r)}_{Z,n}) \rightarrow\cF^{(r)}_{Y,n}
\end{equation} 
qui s'insère dans un diagramme commutatif
\begin{equation}\label{ahttf36b}
\xymatrix{
0\ar[r]&{(\ogg^{\circ})^*_\fet(\ocB_{Z,n})}\ar[r]\ar[d]&{(\ogg^{\circ})^*_\fet(\cF^{(r)}_{Z,n})}\ar[r]\ar[d]&
{\txi^{-1}\tOmega^1_{X/S}(Z)\otimes_{\co_X(Z)}(\ogg^{\circ})^*_\fet(\ocB_{Z,n})}\ar[r]\ar[d]&0\\
0\ar[r]&{\ocB_{Y,n}}\ar[r]&{\cF^{(r)}_{Y,n}}\ar[r]&
{\txi^{-1}\tOmega^1_{X/S}(Y)\otimes_{\co_X(Y)}\ocB_{Y,n}}\ar[r]&0}
\end{equation}
On en déduit par adjonction un morphisme $\ocB_{Z,n}$-linéaire
\begin{equation}\label{ahttf36c}
\cF^{(r)}_{Z,n} \rightarrow(\ogg^{\circ})_{\fet*}(\cF^{(r)}_{Y,n}).
\end{equation} 
On en déduit aussi un morphisme de $(\ogg^{\circ})^*_\fet(\ocB_{Z,n})$-algèbres
\begin{equation}\label{ahttf36d}
(\ogg^{\circ})^*_\fet(\cC^{(r)}_{Z,n}) \rightarrow\cC^{(r)}_{Y,n},
\end{equation} 
et donc par adjonction un morphisme de $\ocB_{Z,n}$-algèbres
\begin{equation}\label{ahttf36f}
\cC^{(r)}_{Z,n} \rightarrow(\ogg^{\circ})_{\fet*}(\cC^{(r)}_{Y,n}). 
\end{equation} 
Les morphismes \eqref{ahttf36c} et \eqref{ahttf36f} vérifient des relations de cocycles 
du type (\cite{egr1} (1.1.2.2)).

\subsection{}\label{ahttf9}
Pour tout nombre rationnel $r\geq 0$,  tout entier $n\geq 0$ et tout objet $Y$ de $\bP$ tel que $Y_s$ soit vide, 
on pose $\cC^{(r)}_{Y,n}=\cF^{(r)}_{Y,n}=0$.  La suite exacte \eqref{ahttf35a} vaut encore dans ce cas, puisque $\ocB_Y$ est une $\oK$-algèbre.
Les morphismes \eqref{ahttf36c} et \eqref{ahttf36f} sont alors définis pour tout morphisme 
de $\bP$, et ils vérifient des relations de cocycles du type (\cite{egr1} (1.1.2.2)).

\subsection{}\label{ahttf37}
Soient $r$ un nombre rationnel $\geq 0$, $n$ un entier $\geq 0$. 
Les correspondances $\{U\mapsto p^n\ocB_U\}$ et $\{U\mapsto \ocB_{U,n}\}$
forment naturellement des préfaisceaux sur $E$ \eqref{ahttf2k}, et les morphismes canoniques 
\begin{eqnarray}
\{U\mapsto p^n\ocB_U\}^\tta&\rightarrow&p^n\ocB,\label{TFA8c}\\
\{U\mapsto \ocB_{U,n}\}^\tta&\rightarrow&\ocB_n,\label{TFA8d}
\end{eqnarray}
où  les termes de gauche désignent les faisceaux associés dans $\tE$, sont des isomorphismes
en vertu de (\cite{agt} VI.8.2 et VI.8.9). D'après \ref{ahttf36}, les correspondances 
\begin{equation}\label{ahttf37a}
\{Y\in \bQ^\circ\mapsto \cF^{(r)}_{Y,n} \} \ \ \ {\rm et}\ \ \ \{Y\in \bQ^\circ\mapsto \cC^{(r)}_{Y,n}\}
\end{equation} 
définissent des préfaisceaux  sur $E_\bQ$ \eqref{ahttf5d} de modules et d'algèbres, respectivement, 
relativement à l'anneau $\{Y\in \bQ^\circ\mapsto \ocB_{Y,n}\}$. On pose 
\begin{eqnarray}
\cF^{(r)}_n&=&\{Y\in \bQ^\circ\mapsto \cF^{(r)}_{Y,n}\}^a,\label{ahttf37b}\\
\cC^{(r)}_n&=&\{Y\in \bQ^\circ\mapsto \cC^{(r)}_{Y,n}\}^a, \label{ahttf37c}
\end{eqnarray}
les faisceaux associés dans $\tE$ (cf. \ref{cftf10} et \ref{ahttf5}). D'après \eqref{TFA8d} et (\cite{agt} (III.10.6.5)),  
$\cF^{(r)}_n$ est un $\ocB_n$-module~; on l'appelle la {\em $\ocB_n$-extension de Higgs-Tate d'épaisseur $r$} 
associée à $(f,\tX,\cM_\tX)$.
De même, $\cC^{(r)}_n$ est une $\ocB_n$-algèbre~; on l'appelle la {\em $\ocB_n$-algèbre de Higgs-Tate d'épaisseur $r$} 
associée à $(f,\tX,\cM_\tX)$. On pose $\cF_n=\cF_n^{(0)}$ et $\cC_n=\cC_n^{(0)}$, et  
on les appelle la {\em $\ocB_n$-extension de Higgs-Tate} et  la {\em $\ocB_n$-algèbre de Higgs-Tate}, respectivement,
associées à $(f,\tX,\cM_\tX)$.

Pour tous nombres rationnels $r\geq r'\geq 0$, les morphismes \eqref{ahttf35c} induisent un morphisme $\ocB_n$-linéaire 
\begin{equation}\label{ahttf37d}
\tta_n^{r,r'}\colon \cF^{(r)}_n\rightarrow \cF_n^{(r')}.
\end{equation}
Les homomorphismes \eqref{ahttf35d} induisent un homomorphisme de $\ocB_n$-algèbres 
\begin{equation}\label{ahttf37e}
\alpha_n^{r,r'}\colon \cC_n^{(r)}\rightarrow \cC_n^{(r')}.
\end{equation}
Pour tous nombres rationnels $r\geq r'\geq r''\geq 0$, on a
\begin{equation}\label{ahttf37f}
\tta_n^{r,r''}=\tta_n^{r',r''} \circ \tta_n^{r,r'} \ \ \ {\rm et}\ \ \ \alpha_n^{r,r''}=\alpha_n^{r',r''} \circ \alpha_n^{r,r'}.
\end{equation}

\begin{prop}[\cite{agt} III.10.22]\label{ahttf38}
Soient $r$ un nombre rationnel $\geq 0$, $n$ un entier $\geq 1$. Alors~:
\begin{itemize}
\item[{\rm (i)}] Les faisceaux $\cF^{(r)}_n$ et $\cC^{(r)}_n$ sont des objets de $\tE_s$.  
\item[{\rm (ii)}] Avec les notations de \eqref{ahttf1c} et \eqref{ahttf42a}, on a une suite exacte localement scindée canonique de $\ocB_n$-modules 
\begin{equation}\label{ahttf38a}
0\rightarrow \ocB_n\rightarrow \cF^{(r)}_n\rightarrow 
\sigma_n^*(\txi^{-1}\tOmega^1_{\oX_n/\oS_n})\rightarrow 0.
\end{equation}
Elle induit pour tout entier $m\geq 1$, une suite exacte de $\ocB_n$-modules 
\begin{equation}\label{ahttf38b}
0\rightarrow \rS^{m-1}_{\ocB_n}(\cF^{(r)}_n)\rightarrow \rS^m_{\ocB_n}(\cF^{(r)}_n)\rightarrow 
\sigma_n^*(\rS^m_{\co_{\oX_n}}(\txi^{-1}\tOmega^1_{\oX_n/\oS_n}))\rightarrow 0.
\end{equation}
En particulier, les $\ocB_n$-modules $(\rS^m_{\ocB_n}(\cF^{(r)}_n))_{m\in \mN}$ forment un système inductif filtrant. 
\item[{\rm (iii)}] On a un isomorphisme canonique de $\ocB_n$-algèbres  
\begin{equation}\label{ahttf38c}
\cC^{(r)}_n \stackrel{\sim}{\rightarrow}\underset{\underset{m\geq 0}{\longrightarrow}}\lim\ \rS^m_{\ocB_n}(\cF^{(r)}_n).
\end{equation}
\item[{\rm (iv)}] Pour tous nombres rationnels $r\geq r'\geq 0$, le diagramme 
\begin{equation}\label{ahttf38d}
\xymatrix{
0\ar[r]&{\ocB_n}\ar[r]\ar@{=}[d]&
{\cF^{(r)}_n}\ar[r]\ar[d]^{\tta_n^{r,r'}}&{\sigma_n^*(\txi^{-1}\tOmega^1_{\oX_n/\oS_n})}\ar[r]\ar[d]^{\cdot p^{r-r'}}& 0\\
0\ar[r]&{\ocB_n}\ar[r]&{\cF^{(r')}_n}\ar[r]&{\sigma_n^*(\txi^{-1}\tOmega^1_{\oX_n/\oS_n})}\ar[r]& 0}
\end{equation}
où les lignes horizontales sont les suites exactes \eqref{ahttf38a}
et la flèche verticale de droite désigne la multiplication par $p^{r-r'}$, est commutatif. 
De plus, les morphismes $\tta_n^{r,r'}$ et $\alpha_n^{r,r'}$ sont compatibles avec les isomorphismes 
\eqref{ahttf38c} pour $r$ et $r'$. 
\end{itemize}
\end{prop}

\subsection{}\label{ahttf39}
Supposons que le morphisme $f\colon (X,\cM_X)\rightarrow (S,\cM_S)$ 
admette une carte adéquate $((P,\gamma),(\mN,\iota),\vartheta)$ que l'on fixe. 
Soient $r$ un nombre rationnel $\geq 0$, $n$ un entier $\geq 0$. 
D'après \ref{ahttf36}, les correspondances 
\begin{equation}\label{ahttf39a}
\{Y\in \bP^\circ\mapsto \cF^{(r)}_{Y,n}\}\ \ \ {\rm et}\ \ \ \{Y\in \bP^\circ\mapsto \cC^{(r)}_{Y,n}\}
\end{equation} 
définissent alors des préfaisceaux  sur $E_\bP$ \eqref{ahttf4d} de modules et d'algèbres, respectivement, 
relativement à l'anneau $\{Y\mapsto \ocB_{Y,n}\}$. Ces préfaisceaux dépendent de la carte
adéquate $((P,\gamma),(\mN,\iota),\vartheta)$, mais les faisceaux associés n'en dépendent pas. 
En effet, en vertu de \ref{cftf11}(i), on a des isomorphismes canoniques 
\begin{eqnarray}
\cF^{(r)}_n&\stackrel{\sim}{\rightarrow}&\{Y\in \bP^\circ\mapsto \cF^{(r)}_{Y,n}\}^a,\label{ahttf39b}\\
\cC^{(r)}_n&\stackrel{\sim}{\rightarrow}&\{Y\in \bP^\circ\mapsto \cC^{(r)}_{Y,n}\}^a, \label{ahttf39c}
\end{eqnarray}
où le $\ocB_n$-module $\cF^{(r)}_n$ et la $\ocB_n$-algèbre $\cC^{(r)}_n$ sont définis dans \ref{ahttf37}.

\subsection{}\label{ahttf55}
Soient $Y$ un objet de $\bP$ tel que $Y_s$ soit non vide, $((P,\gamma),(\mN,\iota),\vartheta)$
une carte adéquate pour le morphisme $f|Y\colon (Y,\cM_X|Y)\rightarrow (S,\cM_S)$ induit par $f$, 
$\oy$ un point géométrique de $\oY^\circ$. Reprenons les notations de \ref{ahttf6}.
Pour tout nombre rationnel $r\geq 0$, on désigne par $\cF_Y^{\oy,(r)}$ 
l'extension de $\hoR_Y^\oy$-modules déduite de $\cF_Y^\oy$ 
\eqref{ahttf6f} par image inverse par le morphisme de multiplication par $p^r$ sur 
$\txi^{-1}\tOmega^1_{X/S}(Y)\otimes_{\co_X(Y)} \hoR^\oy_Y$,
de sorte qu'on a une suite exacte de $\hoR^\oy_Y$-modules 
\begin{equation}\label{ahttf55a}
0\rightarrow \hoR^\oy_Y\rightarrow \cF_Y^{\oy,(r)}\rightarrow \txi^{-1}\tOmega^1_{X/S}(Y)\otimes_{\co_X(Y)} \hoR^\oy_Y
\rightarrow 0.
\end{equation}
On désigne par $\cC_Y^{\oy,(r)}$ la $\hoR^\oy_Y$-algèbre \eqref{taht10c}
\begin{equation}\label{ahttf55b}
\cC_Y^{\oy,(r)}= \underset{\underset{m\geq 0}{\longrightarrow}}\lim\ \rS^m_{\hoR^\oy_Y}(\cF_Y^{\oy,(r)}).
\end{equation}

Il résulte de \ref{ahttf8} que les formations de $\cF_Y^{\oy,(r)}$ et $\cC_Y^{\oy,(r)}$ sont fonctorielles 
en la paire $(Y,\oy)$. Plus précisément, soient $g\colon Z\rightarrow Y$ un morphisme de $\bP$, $\oz$ un point géométrie de $\oZ^\circ$
d'image $\oy$ par le morphisme $\ogg^\circ\colon \oZ^\circ\rightarrow \oY^\circ$. 
On munit le morphisme $f|Z\colon (Z,\cM_X|Z)\rightarrow (S,\cM_S)$ induit par $f$ de la même carte adéquate $((P,\gamma),(\mN,\iota),\vartheta)$. 
Il résulte aussitôt de \ref{ahttf8} que le diagramme canonique
\begin{equation}\label{ahttf55i}
\xymatrix{
0\ar[r]&{\hoR^\oy_Y}\ar[r]\ar[d]&{\cF^{\oy,(r)}_Y}\ar[r]\ar[d]&
{\txi^{-1}\tOmega^1_{X/S}(Y)\otimes_{\co_X(Y)}\hoR^\oy_Y}\ar[r]\ar[d]&0\\
0\ar[r]&{\hoR^\oz_Z}\ar[r]&{\cF^{\oz,(r)}_Z}\ar[r]&
{\txi^{-1}\tOmega^1_{X/S}(Z)\otimes_{\co_X(Z)}\hoR^\oz_Z}\ar[r]&0}
\end{equation}
est commutatif. Le morphisme canonique 
\begin{equation}\label{ahttf55j}
\tOmega^1_{X/S}(Y)\otimes_{\co_X(Y)}\co_X(Z)\rightarrow \tOmega^1_{X/S}(Z)
\end{equation} 
étant un isomorphisme, on en déduit que les morphismes canoniques
\begin{eqnarray}
\cF^{\oy,(r)}_Y\otimes_{\hoR_Y^\oy}\hoR_Z^\oz&\rightarrow&\cF^{\oz,(r)}_Z,\label{ahttf55k}\\
\cC^{\oy,(r)}_Y\otimes_{\hoR_Y^\oy}\hoR_Z^\oz&\rightarrow&\cC^{\oz,(r)}_Z,\label{ahttf55l}
\end{eqnarray} 
sont des isomorphismes. 

Comme $\oY$ est localement irréductible \eqref{ahttf6}, 
il est la somme des schémas induits sur ses composantes irréductibles. On note  
$\oY^\star$ la composante irréductible de $\oY$ contenant $\oy$. De même, $\oY^\circ$
est la somme des schémas induits sur ses composantes irréductibles, et $\oY^{\star \circ}=\oY^\star\times_XX^\circ$ 
est la composante irréductible de $\oY^\circ$ contenant $\oy$. On note $\bB_{\pi_1(\oY^{\star \circ},\oy)}$ 
le topos classifiant du groupe profini $\pi_1(\oY^{\star \circ},\oy)$ et 
\begin{equation}\label{ahttf55e}
\nu_\oy\colon \oY^{\star \circ}_\fet \stackrel{\sim}{\rightarrow} \bB_{\pi_1(\oY^{\star \circ},\oy)}
\end{equation}
le foncteur fibre en $\oy$ (\cite{agt} VI.9.8). On a alors \eqref{ahttf44b} 
\begin{equation}\label{ahttf55f}
\oR^\oy_Y=\nu_\oy(\ocB_Y|\oY^{\star \circ}).
\end{equation}
Comme $\nu_\oy$ est exact et qu'il commute aux limites inductives, 
pour tout entier $n\geq 0$, on a des isomorphismes canoniques 
de $\oR^\oy_Y$-modules et de $\oR^\oy_Y$-algèbres, respectivement (cf. \ref{ahttf35}),
\begin{eqnarray}
\nu_\oy(\cF^{(r)}_{Y,n}|\oY^{\star \circ})&\stackrel{\sim}{\rightarrow}& \cF^{\oy,(r)}_Y/p^n \cF^{\oy,(r)}_Y,\label{ahttf55g}\\
\nu_\oy(\cC^{(r)}_{Y,n}|\oY^{\star \circ})&\stackrel{\sim}{\rightarrow}& \cC^{\oy,(r)}_Y/p^n \cC^{\oy,(r)}_Y.\label{ahttf55h}
\end{eqnarray}

\subsection{}\label{ahttf56}
Soient $(\oy\rightsquigarrow \ox)$ un point de $X_\et\gtimes_{X_\et}\oX^\circ_\et$ \eqref{TFA9} tel que $\ox$ soit
au-dessus de $s$, $\uX$ le localisé strict de $X$ en $\ox$. 
On rappelle que la donnée d'un voisinage du point de $X_\et$ associé à $\ox$ 
dans le site $\Et_{/X}$ (resp. $\bP$ \eqref{ahttf4}, resp. $\bQ$ \eqref{ahttf5})
est équivalente à la donnée d'un $X$-schéma étale $\ox$-pointé (resp. de $\bP$, resp. de $\bQ$) (\cite{sga4} IV 6.8.2). 
Ces objets forment naturellement une catégorie cofiltrante, que l'on note $\fV_\ox$ (resp. $\bP_\ox$,
resp. $\bQ_\ox$). Les catégories $\bP_\ox$ et $\bQ_\ox$ sont $\mU$-petites, et 
les foncteurs d'injection canoniques $\bQ\rightarrow \bP\rightarrow \Et_{/X}$ induisent des foncteurs pleinement 
fidèles et cofinaux $\bQ_\ox\rightarrow \bP_\ox\rightarrow \fV_\ox$.

Reprenons les notations de \ref{TFA12}, en particulier l'algèbre
\begin{equation}\label{ahttf56d}
\oR^\oy_{\uX}=\underset{\underset{(U,\fp)\in \fV_\ox^\circ}{\longrightarrow}}{\lim}\ \oR^{\oy}_U,
\end{equation}
où $\oR^{\oy}_U$ est l'anneau défini dans \eqref{ahttf55f}.  
On note $\hoR^\oy_{\uX}$ le complété séparé $p$-adique de $\oR^\oy_{\uX}$. 
On a un isomorphisme canonique \eqref{TFA12g}
\begin{equation}\label{ahttf56m}
\nu_{\oy}(\varphi_\ox(\ocB)) \stackrel{\sim}{\rightarrow} \oR^{\oy}_{\uX}. 
\end{equation} 

Pour tout nombre rationnel $r\geq 0$, on pose 
\begin{eqnarray}
\cF^{\oy,(r)}_{\uX}&=&\underset{\underset{(U,\fp)\in \bQ_\ox^\circ}{\longrightarrow}}{\lim}\ 
\cF^{\oy,(r)}_U\otimes_{\hoR_U^\oy}\hoR_{\uX}^\oy,\label{ahttf56e}\\
\cC^{\oy,(r)}_{\uX}&=&\underset{\underset{(U,\fp)\in \bQ_\ox^\circ}{\longrightarrow}}{\lim}\ 
\cC^{\oy,(r)}_U\otimes_{\hoR_U^\oy}\hoR_{\uX}^\oy,\label{ahttf56f}
\end{eqnarray}
où $\cF^{\oy,(r)}_U$ est le $\hoR^{\oy}_U$-module 
défini dans \eqref{ahttf55a} et $\cC^{\oy,(r)}_U$ est la $\hoR^{\oy}_U$-algèbre définie dans \eqref{ahttf55b}. 
D'après \eqref{ahttf55k} et \eqref{ahttf55l}, pour tout objet $(U,\fp)$ de $\bQ_\ox$, les morphismes canoniques
\begin{eqnarray}
\cF^{\oy,(r)}_U\otimes_{\hoR_U^\oy}\hoR_{\uX}^\oy&\rightarrow& \cF^{\oy,(r)}_{\uX},\label{ahttf56k}\\
\cC^{\oy,(r)}_U\otimes_{\hoR_U^\oy}\hoR_{\uX}^\oy&\rightarrow& \cC^{\oy,(r)}_{\uX},\label{ahttf56l}
\end{eqnarray}
sont des isomorphismes.  

\begin{rema}\label{ahttf57}
Sous les hypothèses de \ref{ahttf56}, pour tout entier $n\geq 0$, les morphismes canoniques 
\begin{eqnarray}
\underset{\underset{(U,\fp)\in \bQ_\ox^\circ}{\longrightarrow}}{\lim}\ 
\cF^{\oy,(r)}_U/p^n\cF^{\oy,(r)}_U&\rightarrow&\cF^{\oy,(r)}_{\uX}/p^n\cF^{\oy,(r)}_{\uX},\label{ahttf57b}\\
\underset{\underset{(U,\fp)\in \bQ_\ox^\circ}{\longrightarrow}}{\lim}\ 
\cC^{\oy,(r)}_U/p^n\cC^{\oy,(r)}_U&\rightarrow&\cC^{\oy,(r)}_{\uX}/p^n\cC^{\oy,(r)}_{\uX},\label{ahttf57a}
\end{eqnarray}
sont des isomorphismes.
\end{rema}

\subsection{}\label{ahttf41}
Soient $g\colon X'\rightarrow X$ un morphisme étale de type fini, $r$ un nombre rationnel $\geq 0$, $n$ un entier $\geq 0$. 
On munit $X'$ de la structure logarithmique $\cM_{X'}$ image inverse de $\cM_X$ et on note $f'\colon (X',\cM_{X'})\rightarrow (S,\cM_S)$ 
le morphisme induit par $f$ et $g$. On observera que $f'$ est adéquat (\cite{agt} III.4.7) et que $X'^\circ=X^\circ\times_XX'$ 
est le sous-schéma ouvert maximal de $X'$ où la structure logarithmique $\cM_{X'}$ est triviale.
On munit $\oX'$ et $\coX'$ \eqref{defing1c} des structures logarithmiques $\cM_{\oX'}$ et $\cM_{\coX'}$ 
images inverses de $\cM_{X'}$. 
Il existe essentiellement un unique morphisme étale $\tg\colon \tX'\rightarrow \tX$ 
qui s'insère dans un diagramme cartésien \eqref{defing12}
\begin{equation}\label{ahttf41a}
\xymatrix{
{\coX'}\ar[r]\ar[d]_{\cog}&{\tX'}\ar[d]^{\tg}\\
{\coX}\ar[r]&{\tX}}
\end{equation}
On munit $\tX'$ de la structure logarithmique $\cM_{\tX'}$ image inverse de $\cM_{\tX}$, de sorte que 
$(\tX',\cM_{\tX'})$ est une $(\tS,\cM_{\tS})$-déformation lisse de $(\coX',\cM_{\coX'})$.

On associe à $(f',\tX',\cM_{\tX'})$ des objets analogues à ceux associés à $(f,\tX,\cM_{\tX})$, 
qu'on note par les mêmes symboles affectés d'un exposant $^\prime$.

Tout $X'$-schéma étale est naturellement un $X$-schéma étale. 
On définit ainsi un foncteur 
\begin{equation}\label{ahttf41b}
\Et_{/X'}\rightarrow \Et_{/X},
\end{equation}
qui se factorise à travers une équivalence de catégories $\Et_{/X'} \stackrel{\sim}{\rightarrow} (\Et_{/X})_{/X'}$. 
Pour qu'un objet $U'$ de $\Et_{/X'}$ soit un objet de $\bP'$ (resp. $\bQ'$), 
il faut et il suffit que $U'$ soit un objet de $\bP$ (resp. $\bQ$). 

Tout objet de $E'$ est naturellement un objet de $E$. On définit ainsi un foncteur 
\begin{equation}\label{ahttf41c}
\Phi\colon E'\rightarrow E.
\end{equation}
Celui-ci se factorise à travers une équivalence de catégories $E'\stackrel{\sim}{\rightarrow} E_{/(\oX'^\circ\rightarrow X')}$. 
On désigne par 
\begin{equation}\label{ahttf41d}
\hPhi^*\colon \hE\rightarrow \hE'
\end{equation}
le foncteur défini par la composition avec $\Phi$. 

D'après (\cite{agt} VI.5.38), la topologie co-évanescente de $E'$ est induite par celle de $E$ au moyen du foncteur $\Phi$. 
Par suite, $\Phi$ est continu et cocontinu (\cite{sga4} III 5.2). Il définit donc une suite de trois foncteurs adjoints:
\begin{equation}\label{ahttf41e}
\Phi_!\colon \tE'\rightarrow \tE, \ \ \ \Phi^*\colon \tE\rightarrow \tE', \ \ \ \Phi_*\colon \tE'\rightarrow \tE,
\end{equation}
dans le sens que pour deux foncteurs consécutifs de la suite, celui de droite est
adjoint à droite de l'autre. D'après (\cite{sga4} III 5.4), le foncteur $\Phi_!$ se factorise à travers une équivalence de catégories 
\begin{equation}\label{ahttf41f}
\tE'\stackrel{\sim}{\rightarrow} \tE_{/\sigma^*(X')},
\end{equation}
où $\sigma^*(X')=(\oX'^\circ\rightarrow X')^a$ \eqref{ahttf2e}.
Comme $\Phi\colon E'\rightarrow E$ est un adjoint à gauche du foncteur 
\begin{equation}\label{ahttf41g}
E\rightarrow E',\ \ \ (V\rightarrow U)\mapsto (V\times_XX'\rightarrow U\times_XX'), 
\end{equation}
le morphisme de localisation $(\Phi^*,\Phi_*)\colon \tE'\rightarrow \tE$ de $\tE$ en $\sigma^*(X')$
s'identifie au morphisme de fonctorialité induit par $g$ (\cite{agt} (VI.10.12)), en vertu de (\cite{sga4} III 2.5). 

On a un homomorphisme canonique $\Phi^*(\ocB)\rightarrow \ocB'$ (\cite{agt} (III.8.20.6)), 
qui est un isomorphisme en vertu de (\cite{agt} III.8.21(i)). 

Soient $n$ un entier $\geq 0$, $Y'$ un objet de $\bP'$, $((P',\gamma'),(\mN,\iota),\vartheta')$ une carte adéquate pour le morphisme 
\begin{equation}\label{ahttf41h}
f'|Y'\colon (Y',\cM_{X'}|Y')\rightarrow (S,\cM_S)
\end{equation} 
induit par $f'$. On peut alors considérer le faisceau $\cF'^{(r)}_{Y',n}$ de $\oY'^\circ_\fet$ \eqref{ahttf35}. 
Comme $\cM_X|Y'=\cM_{X'}|Y'$, considérant $Y'$ comme objet de $\bP$,
on peut aussi considérer le faisceau $\cF^{(r)}_{Y',n}$ de $\oY'^\circ_\fet$. 
Le morphisme $\tg$ étant étale, on a clairement $\cF^{(r)}_{Y',n}=\cF'^{(r)}_{Y',n}$ (cf. \cite{agt} III.14.3).  

D'après \ref{cftf11}(ii), il existe un préfaisceau $F=\{U\in \Et^\circ_{/X}\mapsto F_U\}$ sur $E$ tel que 
pour tout $U\in \ob(\bQ)$, on ait $F_U=\cF^{(r)}_{U,n}$, de sorte qu'on a un isomorphisme canonique $\cF^{(r)}_n\stackrel{\sim}{\rightarrow} F^a$. 
Par ailleurs, on a 
\begin{equation}\label{ahttf41i}
\hPhi^*(F)=\{U'\in \Et^\circ_{/X'}\mapsto F_{U'}\}.
\end{equation}
On en déduit, compte tenu de \ref{cftf11}(i) et (\cite{sga4} III 2.3(2)), un isomorphisme $\ocB'_n$-linéaire canonique 
\begin{equation}\label{ahttf41j}
\Phi^*(\cF^{(r)}_n)\stackrel{\sim}{\rightarrow}\cF'^{(r)}_n.
\end{equation}

\subsection{}\label{ahttf46}
Soit $U$ un objet de $\Et_{/X}$. D'après (\cite{agt} VI.10.14)  (cf. aussi \ref{ahttf41}), le topos  $\tE_{/\sigma^*(U)}$, localisé de $\tE$ en $\sigma^*(U)$, 
est canoniquement équivalent au topos de Faltings associé au morphisme $\oU^\circ\rightarrow U$ \eqref{ahttf2e}. 
On désigne par
\begin{equation}\label{ahttf46a}
\jmath_U\colon \tE_{/\sigma^*(U)}\rightarrow \tE
\end{equation}
le morphisme de localisation de $\tE$ en $\sigma^*(U)$, qui s'identifie au morphisme de fonctorialité induit par le morphisme canonique 
$U\rightarrow X$  (\cite{agt} (VI.10.12)), et par 
\begin{equation}\label{ahttf46b}
\beta_U\colon \tE_{/\sigma^*(U)} \rightarrow \oU^\circ_\fet
\end{equation}
le morphisme canonique \eqref{ahttf2d}. 

Soit $\upmu \colon U'\rightarrow U$ un morphisme de $\Et_{/X}$. On désigne par 
\begin{equation}\label{ahttf46c}
\Phi_\upmu\colon \tE_{/\sigma^*(U')}\rightarrow \tE_{/\sigma^*(U)}
\end{equation}
le morphisme de localisation associé au morphisme $\sigma^*(U)\rightarrow \sigma^*(U)$  (\cite{sga4} IV 5.5),
qui s'identifie au morphisme de fonctorialité induit par $\upmu$. Les diagrammes 
\begin{equation}\label{ahttf46d}
\xymatrix{
{\tE_{/\sigma^*(U')}}\ar[r]^-(0.5){\Phi_\upmu}\ar[rd]_{\jmath_{U'}}&{\tE_{/\sigma^*(U)}}\ar[d]^{\jmath_U}\\
&\tE}
\end{equation}
\begin{equation}\label{ahttf46e}
\xymatrix{
{\tE_{/\sigma^*(U')}}\ar[r]^-(0.5){\Phi_\upmu}\ar[d]_{\beta_{U'}}&{\tE_{/\sigma^*(U)}}\ar[d]^{\beta_U}\\
{\oU'^\circ_\fet}\ar[r]^-(0.5){\oupmu^\circ}&{\oU^\circ_\fet}}
\end{equation}
sont commutatifs à isomorphismes canoniques près (\cite{agt} (VI.10.12.6)). 

\begin{lem}\label{ahttfg9}
Soient $U$ un objet de $\bP$, $((P,\gamma),(\mN,\iota),\vartheta)$ une carte adéquate pour le morphisme 
$f|U\colon (U,\cM_X|U)\rightarrow (S,\cM_S)$ induit par $f$, $r$ un nombre rationnel $\geq 0$, $n$ un entier $\geq 0$. 
Alors, avec les notations de \ref{ahttf46}, on a un morphisme canonique 
\begin{equation}\label{ahttfg9a}
\cF^{(r)}_{U,n}\rightarrow \beta_{U*}(\jmath_U^*(\cF_n^{(r)})).
\end{equation}
Celui-ci est indépendant de la carte adéquate $((P,\gamma),(\mN,\iota),\vartheta)$ si $U$ est un objet de $\bQ$.  

Soit $\upmu\colon U'\rightarrow U$ un morphisme de $\bP$. Munissons le morphisme 
$f|U'\colon (U',\cM_X|U')\rightarrow (S,\cM_S)$ induit par $f$ de la carte adéquate induite par $((P,\gamma),(\mN,\iota),\vartheta)$. 
Alors, le diagramme 
\begin{equation}\label{ahttfg9b}
\xymatrix{
{\beta_{U'}^*(\oupmu^{\circ*}(\cF^{(r)}_{U,n}))}\ar[d]_{\beta_{U'}^*(a)}\ar[r]^b&{\Phi_\upmu^*(\beta_U^*(\cF^{(r)}_{U,n}))}\ar[r]&
{\Phi_\upmu^*(\jmath_{U}^*(\cF_n^{(r)}))}\ar[d]^c\\
{\beta_{U'}^*(\cF^{(r)}_{U',n})}\ar[rr]&&{\jmath_{U'}^*(\cF_n^{(r)})}}
\end{equation}
où $a\colon \oupmu^{\circ*}(\cF^{(r)}_{U,n})\rightarrow \cF^{(r)}_{U',n}$ (resp. $b$, resp. $c$) est le morphisme canonique \eqref{ahttf36a} 
(resp. l'isomorphisme sous-jacent à \eqref{ahttf46e}, resp. l'isomorphisme sous-jacent à \eqref{ahttf46d}), 
et les deux autres flèches sont induites par l'adjoint du morphisme \eqref{ahttfg9a}, est commutatif. 
\end{lem}

En effet, compte tenu de \eqref{ahttf41j}, on peut se réduire au cas où $X=U$. Pour tout objet $Y$ de $\bP$, 
on peut donc définir le faisceau $\cF^{(r)}_{Y,n}$ de $\oY^\circ_\fet$ relativement à la carte adéquate $((P,\gamma),(\mN,\iota),\vartheta)$. 
En vertu de \ref{cftf11}(i), on a un isomorphisme canonique 
\begin{equation}\label{ahttfg9c}
\cF^{(r)}_n\stackrel{\sim}{\rightarrow}\{Y\in \bP^\circ\mapsto \cF^{(r)}_{Y,n}\}^a.
\end{equation}
D'après \ref{cftf11}(ii), il existe un préfaisceau canonique $F=\{Y\in \Et^\circ_{/X}\mapsto F_Y\}$ sur $E$ tel que 
pour tout $Y\in \ob(\bP)$, on ait $F_Y=\cF^{(r)}_{Y,n}$, de sorte qu'on a un isomorphisme canonique $\cF^{(r)}_n\stackrel{\sim}{\rightarrow} F^a$. 
Le morphisme canonique $F\rightarrow \cF^{(r)}_n$ induit, pour tout $Y\in \ob(\Et^\circ_{/X})$, un morphisme 
\begin{equation}\label{ahttfg9d}
F_Y\rightarrow \beta_{Y*}(\jmath_Y^*(\cF_n^{(r)})),
\end{equation}
d'où le morphisme \eqref{ahttfg9a}. 

Si $U$ est un objet de $\bQ$, appliquant le même argument au préfaisceau $\{Y\in \bQ^\circ\mapsto \cF^{(r)}_{Y,n}\}$ sur $E_\bQ$,
on voit que le morphisme \eqref{ahttfg9a} ne dépend pas de la carte adéquate $((P,\gamma),(\mN,\iota),\vartheta)$. 

Pour tout morphisme $\upmu\colon U'\rightarrow U$ de $\Et^\circ_{/X}$, le diagramme 
\begin{equation}\label{ahttfg9e}
\xymatrix{
{F_U}\ar[rr]^-(0.5){t}\ar[d]&&{\oupmu^{\circ}_*(F_{U'})}\ar[d]\\
{\beta_{U*}(\jmath_U^*(\cF_n^{(r)}))}\ar[r]^-(0.5){\ad}&
{\beta_{U*}(\Phi_{\upmu*}(\Phi_{\upmu}^*(\jmath_U^*(\cF_n^{(r)}))))}\ar[r]^-(0.5){\lambda}&{\oupmu^{\circ}_*(\beta_{U'*}(\jmath_{U'}^*(\cF_n^{(r)})))}}
\end{equation}
où $t$ est le morphisme de transition de $F$, $\ad$ est induit par le morphisme d'adjonction $\id\rightarrow \Phi_{\upmu*}\Phi_{\upmu}^*$,
$\lambda$ est induit par les diagrammes commutatifs \eqref{ahttf46d} et \eqref{ahttf46e}, et les flèches verticales sont induites par \eqref{ahttfg9d},
est commutatif. On en déduit aussitôt que le diagramme \eqref{ahttfg9b} est commutatif.

\subsection{}\label{ahttf47}
Soient $\ox$ un point géométrique de $X$ au-dessus de $s$, $\uX$ le localisé strict de $X$ en $\ox$, $n$ un entier $\geq 0$.
On note $\fV_\ox$, $\bP_\ox$ et $\bQ_\ox$ les catégories des voisinages de $\ox$ dans les sites $\Et_{/X}$, $\bP$ et $\bQ$, respectivement (cf. \ref{ahttf56}). 
Pour tout objet $(U,\iota\colon \ox\rightarrow U)$ de $\fV_\ox$, on note encore $\iota\colon \uX\rightarrow U$ le $X$-morphisme déduit de $\iota$ (\cite{sga4} VIII 7.3),
et on note abusivement $\oiota\colon \uoX^\circ\rightarrow \oU^\circ$ le morphisme induit.

On désigne par $\tuE$ le topos de Faltings associé au morphisme canonique $\uoX^\circ\rightarrow \uX$, par 
\begin{equation}\label{ahttf47a}
\Phi\colon \tuE\rightarrow \tE
\end{equation}
le morphisme de fonctorialité induit par le morphisme canonique $\uX\rightarrow X$, par 
\begin{equation}\label{ahttf47b}
\ubeta\colon \tuE\rightarrow \uoX^\circ_\fet
\end{equation}
le morphisme canonique \eqref{ahttf2d} et par 
\begin{equation}\label{ahttf47c}
\theta\colon \uoX^\circ_\fet\rightarrow \tuE
\end{equation}
la section canonique de $\ubeta$ définie dans (\cite{agt} VI.10.23). On note 
\begin{equation}\label{ahttf47d}
\varphi_\ox\colon \tE\rightarrow \uoX^\circ_\fet
\end{equation}
le foncteur composé $\theta^*\circ \Phi^*$ \eqref{TFA14e}. 

Soit $(U,\iota\colon \ox\rightarrow U)$ un objet de $\fV_\ox$. Reprenons les notations de \ref{ahttf46}. 
Le morphisme $\iota\colon \uX\rightarrow U$ induit par fonctorialité un morphisme 
\begin{equation}\label{ahttf47e}
\Phi_\iota\colon \tuE\rightarrow \tE_{/\sigma^*(U)}
\end{equation}  
qui s'insère dans un diagramme commutatif à isomorphisme canonique près 
\begin{equation}\label{ahttf47f}
\xymatrix{
{\tuE}\ar[r]^-(0.5){\Phi_\iota}\ar[d]_{\ubeta}&{\tE_{/\sigma^*(U)}}\ar[d]^{\beta_U}\\
{\uoX^\circ_\fet}\ar[r]^-(0.5){\oiota}&{\oU^\circ_\fet}}
\end{equation}

Supposons que $(U,\iota\colon \ox\rightarrow U)$ soit un objet de $\bP_\ox$
et soit $((P,\gamma),(\mN,\iota),\vartheta)$ une carte adéquate pour le morphisme $f|U\colon (U,\cM_X|U)\rightarrow (S,\cM_S)$ induit par $f$. 
Appliquant le foncteur composé $\theta^*\circ \Phi^*_\iota$ au morphisme $\beta^*_{U}(\cF^{(r)}_{U,n})\rightarrow \jmath_U^*(\cF_n^{(r)})$
adjoint de \eqref{ahttfg9a} et tenant compte du diagramme commutatif \eqref{ahttf47f}, on obtient un morphisme 
\begin{equation}\label{ahttf47g}
\oiota^*(\cF^{(r)}_{U,n})\rightarrow \varphi_\ox(\cF^{(r)}_n).
\end{equation}
Celui-ci est indépendant de la carte adéquate $((P,\gamma),(\mN,\iota),\vartheta)$ si $U$ est un objet de~$\bQ$.  

Le morphisme $\beta^*_{U}(\cF^{(r)}_{U,n})\rightarrow \jmath_U^*(\cF_n^{(r)})$ adjoint de \eqref{ahttfg9a} est le composé 
\begin{equation}\label{ahttf47i}
\beta_U^*(\cF^{(r)}_{U,n})\rightarrow \beta_U^*(\beta_{U*}(\jmath_U^*(\cF^{(r)}_n)))\rightarrow \jmath_U^*(\cF^{(r)}_n),
\end{equation}
où la première flèche est l'image par le foncteur $\beta_U^*$ de \eqref{ahttfg9a} et la seconde flèche est le morphisme d'adjonction. 
On en déduit une factorisation du morphisme \eqref{ahttf47g} en
\begin{equation}\label{ahttf47j}
\oiota^*(\cF^{(r)}_{U,n})\rightarrow \oiota^*(\beta_{U*}(\jmath_U^*(\cF^{(r)}_n)))\rightarrow \varphi_\ox(\cF^{(r)}_n).
\end{equation}
D'après (\cite{egr1} 1.2.4(i)), le second morphisme est le composé 
\begin{equation}\label{ahttf47k}
\oiota^*(\beta_{U*}(\jmath_U^*(\cF^{(r)}_n)))\rightarrow \ubeta_*(\Phi_\iota^*(\jmath_U^*(\cF^{(r)}_n)))\rightarrow \theta^*(\Phi^*(\cF^{(r)}_n)),
\end{equation}
où la première flèche est le morphisme de changement de base relativement au diagramme \eqref{ahttf47f} et la seconde flèche est induite par l'isomorphisme 
de changement de base $\ubeta_*\stackrel{\sim}{\rightarrow}\theta^*$ (\cite{agt} VI.10.27). 

Soient $\upmu \colon U'\rightarrow U$ un morphisme de $\bP$, $\iota'\colon \ox\rightarrow U'$ un $X$-morphisme, 
$\iota=\upmu\circ \iota'\colon \ox\rightarrow U$. 
Considérons d'abord une carte adéquate  $((P,\gamma),(\mN,\iota),\vartheta)$ pour le morphisme 
$f|U\colon (U,\cM_X|U)\rightarrow (S,\cM_S)$ induit par $f$ et munissons le morphisme 
$f|U'\colon (U',\cM_X|U')\rightarrow (S,\cM_S)$ induit par $f$ de la carte induite par $((P,\gamma),(\mN,\iota),\vartheta)$. 
Il résulte alors de \eqref{ahttfg9b} que le diagramme 
\begin{equation}\label{ahttf47h}
\xymatrix{
{\oiota^*(\cF^{(r)}_{U,n})}\ar[d]\ar[rd]&\\
{\oiota'^*(\cF^{(r)}_{U',n})}\ar[r]&{\varphi_\ox(\cF^{(r)}_n)}}
\end{equation}
où la flèche verticale est induite par le morphisme canonique $\oupmu^{\circ*}(\cF^{(r)}_{U,n})\rightarrow \cF^{(r)}_{U',n}$ \eqref{ahttf36a} 
et les deux autres flèches sont les morphismes \eqref{ahttf47g}, est commutatif. 

On déduit de ce qui précède que si $\upmu$ est un morphisme de $\bQ$, alors le diagramme \eqref{ahttf47h} est aussi commutatif.

\begin{prop}\label{ahttf48}
Soient $\ox$ un point géométrique de $X$, $n$ un entier $\geq 0$. Reprenons les notations de \ref{ahttf47}. 
Alors, 
\begin{itemize}
\item[{\rm (i)}] Les morphismes \eqref{ahttf47g} pour $(U,\iota)\in \ob(\bQ_\ox)$, induisent des isomorphismes 
\begin{eqnarray}
\underset{\underset{(U,\iota)\in \bQ^\circ_\ox}{\longrightarrow}}{\lim}\  \oiota^*(\cF^{(r)}_{U,n})
&\stackrel{\sim}{\rightarrow}& \varphi_\ox(\cF^{(r)}_n), \label{ahttf48a}\\
\underset{\underset{(U,\iota)\in \bQ^\circ_\ox}{\longrightarrow}}{\lim}\  \oiota^*(\cC^{(r)}_{U,n})
&\stackrel{\sim}{\rightarrow}& \varphi_\ox(\cC^{(r)}_n). \label{ahttf48aa}
\end{eqnarray}
\item[{\rm (ii)}] Supposons que le morphisme $f\colon (X,\cM_X)\rightarrow (S,\cM_S)$ admette une carte adéquate. 
Alors, les morphismes \eqref{ahttf47g} définis relativement à une carte adéquate de $f$, induisent des isomorphismes
\begin{eqnarray}
\underset{\underset{(U,\iota)\in \bP^\circ_\ox}{\longrightarrow}}{\lim}\  \oiota^*(\cF^{(r)}_{U,n})
&\stackrel{\sim}{\rightarrow}& \varphi_\ox(\cF^{(r)}_n), \label{ahttf48b}\\
\underset{\underset{(U,\iota)\in \bP^\circ_\ox}{\longrightarrow}}{\lim}\  \oiota^*(\cC^{(r)}_{U,n})
&\stackrel{\sim}{\rightarrow}& \varphi_\ox(\cC^{(r)}_n). \label{ahttf48bb}
\end{eqnarray}
\end{itemize}
\end{prop}

(i) Cela résulte aussitôt de (ii) puisque le foncteur canonique $\bQ_\ox\rightarrow \bP_\ox$ est cofinal.

(ii) D'après \ref{cftf11}(ii), il existe un préfaisceau canonique $F=\{U\in \Et^\circ_{/X}\mapsto F_U\}$ sur $E$ tel que pour tout $U\in \ob(\bP)$, 
$F_U$ soit le faisceau $\cF^{(r)}_{U,n}$ défini relativement à la carte adéquate donnée de $f$, 
de sorte qu'on a un isomorphisme canonique $\cF^{(r)}_n\stackrel{\sim}{\rightarrow} F^a$.
En vertu de (\cite{agt} VI.10.37), on a un isomorphisme canonique 
\begin{equation}\label{ahttf48c}
\underset{\underset{(U,\iota)\in \fV^\circ_\ox}{\longrightarrow}}{\lim}\  \oiota^*(F^a_U)
\stackrel{\sim}{\rightarrow} \varphi_\ox(F^a).
\end{equation}
Pour tout $(U,\iota)\in \ob(\bP_\ox)$, le morphisme 
$\oiota^*(F^a_U)\rightarrow \varphi_\ox(F^a)$ s'identifie au morphisme \eqref{ahttf47g} 
d'après la description \eqref{ahttf47j} (cf. \cite{agt} VI.10.34). 
L'isomorphisme \eqref{ahttf48b} s'ensuit puisque le foncteur canonique $\bP_\ox\rightarrow \fV_\ox$ est cofinal.
Il induit l'isomorphisme \eqref{ahttf48bb}.

\subsection{}\label{ahttf23}
Soit $r$ un nombre rationnel $\geq 0$.   D'après \ref{ahttf38}, on a un isomorphisme canonique $\cC_n^{(r)}$-linéaire
\begin{equation}\label{ahttf23a}
\Omega^1_{\cC_n^{(r)}/\ocB_n}\stackrel{\sim}{\rightarrow} 
\sigma_n^*(\txi^{-1}\tOmega^1_{\oX_n/\oS_n})\otimes_{\ocB_n}\cC_n^{(r)}.
\end{equation}
La $\ocB_n$-dérivation universelle de $\cC_n^{(r)}$ correspond via cet isomorphisme à l'unique $\ocB_n$-dérivation 
\begin{equation}\label{ahttf23b}
d_{\cC_n^{(r)}}\colon \cC_n^{(r)}\rightarrow \sigma_n^*(\txi^{-1}\tOmega^1_{\oX_n/\oS_n})\otimes_{\ocB_n}\cC_n^{(r)}
\end{equation}
qui prolonge le morphisme canonique $\cF_n^{(r)}\rightarrow \sigma_n^*(\txi^{-1}\tOmega^1_{\oX_n/\oS_n})$ \eqref{ahttf38a}. 
Il résulte de \ref{ahttf38}(iv) que pour tous nombres rationnels $r\geq r'\geq 0$, on a 
\begin{equation}\label{ahttf23c}
p^{r-r'}(\id \otimes \alpha^{r,r'}_n) \circ d_{\cC_n^{(r)}}=d_{\cC_n^{(r')}}\circ \alpha^{r,r'}_n.
\end{equation}

\subsection{}\label{ahttf14}
Soit $r$ un nombre rationnel $\geq 0$.
Pour tous entiers $m\geq n\geq 1$, on a un morphisme $\ocB_m$-linéaire canonique
$\cF^{(r)}_m\rightarrow \cF^{(r)}_n$ compatible avec la suite exacte \eqref{ahttf38a}  
et un homomorphisme canonique de $\ocB_m$-algèbres $\cC^{(r)}_m\rightarrow \cC^{(r)}_n$ tels que 
les morphismes induits 
\begin{equation}\label{ahttf14a}
\cF^{(r)}_m\otimes_{\ocB_m}\ocB_n\rightarrow \cF^{(r)}_n\ \ \ {\rm et}\ \ \ 
\cC^{(r)}_m\otimes_{\ocB_m}\ocB_n\rightarrow \cC^{(r)}_n
\end{equation}
soient des isomorphismes. Ces morphismes forment des systèmes compatibles lorsque $m$ et $n$ varient, 
de sorte que $(\cF^{(r)}_{n+1})_{n\in \mN}$ et $(\cC^{(r)}_{n+1})_{n\in \mN}$ sont des systèmes projectifs. 
Avec les notations de \ref{ahttf13}, on appelle {\em $\bvocB$-extension de Higgs-Tate d'épaisseur $r$} associée à $(f,\tX,\cM_\tX)$, 
et l'on note $\bvcF^{(r)}$, le $\bvocB$-module $(\cF^{(r)}_{n+1})_{n\in \mN}$ de $\tE_s^{\mN^\circ}$. 
On appelle {\em $\bvocB$-algèbre de Higgs-Tate d'épaisseur $r$} associée à $(f,\tX,\cM_\tX)$,
et l'on note $\bvcC^{(r)}$, la $\bvocB$-algèbre $(\cC^{(r)}_{n+1})_{n\in \mN}$ de $\tE_s^{\mN^\circ}$.
Ce sont des $\bvocB$-modules adiques (\cite{agt} III.7.16). 
D'après (\cite{agt} III.7.3(i), (III.7.5.4) et (III.7.12.1)), la suite exacte \eqref{ahttf38a}
induit une suite exacte de $\bvocB$-modules 
\begin{equation}\label{ahttf14b}
0\rightarrow \bvocB\rightarrow \bvcF^{(r)}\rightarrow 
\bvsigma^*(\txi^{-1}\tOmega^1_{\bvoX/\bvoS})\rightarrow 0.
\end{equation}
Comme le $\co_X$-module $\tOmega^1_{X/S}$ est localement libre de type fini, 
le $\bvocB$-module $\bvsigma^*(\txi^{-1}\tOmega^1_{\bvoX/\bvoS})$ est localement libre de type fini
et la suite \eqref{ahttf14b} est localement scindée. 
D'après (\cite{illusie1} I 4.3.1.7), elle induit pour tout entier $m\geq 1$, une suite exacte de $\bvocB$-modules 
\begin{equation}\label{ahttf14c}
0\rightarrow \rS^{m-1}_{\bvocB}(\bvcF^{(r)})\rightarrow \rS^m_{\bvocB}(\bvcF^{(r)})\rightarrow 
\bvsigma^*(\rS^m_{\co_{\bvoX}}(\txi^{-1}\tOmega^1_{\bvoX/\bvoS}))\rightarrow 0.
\end{equation}
En particulier, les $\bvocB$-modules $(\rS^m_{\bvocB}(\bvcF^{(r)}))_{m\in \mN}$ forment un système inductif filtrant. 
D'après (\cite{agt} III.7.3(i) et (III.7.12.3)), on a un isomorphisme canonique de $\bvocB$-algèbres
\begin{equation}\label{ahttf14d}
\bvcC^{(r)}\stackrel{\sim}{\rightarrow}\underset{\underset{m\geq 0}{\longrightarrow}}\lim\ \rS^m_{\bvocB}(\bvcF^{(r)}). 
\end{equation}

On pose $\bvcF=\bvcF^{(0)}$ et $\bvcC=\bvcC^{(0)}$, et  
on les appelle la {\em $\bvocB$-extension de Higgs-Tate} et  la {\em $\bvocB$-algèbre de Higgs-Tate}, respectivement,
associées à $(f,\tX,\cM_\tX)$.
Pour tous nombres rationnels $r\geq r'\geq 0$, les morphismes $(\tta_n^{r,r'})_{n\in \mN}$ \eqref{ahttf37d}
induisent un morphisme $\bvocB$-linéaire 
\begin{equation}\label{ahttf14e}
\bvtta^{r,r'}\colon \bvcF^{(r)}\rightarrow \bvcF^{(r')}.
\end{equation}
Les homomorphismes $(\alpha_n^{r,r'})_{n\in \mN}$ \eqref{ahttf37e}
induisent un homomorphisme de $\bvocB$-algèbres 
\begin{equation}\label{ahttf14f}
\bvalpha^{r,r'}\colon \bvcC^{(r)}\rightarrow \bvcC^{(r')}.
\end{equation}
Pour tous nombres rationnels $r\geq r'\geq r''\geq 0$, on a
\begin{equation}\label{ahttf14g}
\bvtta^{r,r''}=\bvtta^{r',r''} \circ \bvtta^{r,r'} \ \ \ {\rm et}\ \ \ \bvalpha^{r,r''}=\bvalpha^{r',r''} \circ \bvalpha^{r,r'}.
\end{equation}

Les dérivations $(d_{n+1}^{(r)})_{n\in \mN}$ \eqref{ahttf23b} définissent un morphisme
\begin{equation}\label{ahttf14h}
d_{\bvcC^{(r)}}\colon \bvcC^{(r)}\rightarrow \bvsigma^*(\txi^{-1}\tOmega^1_{\bvoX/\bvoS})\otimes_{\bvocB}\bvcC^{(r)},
\end{equation}
qui n'est autre que la $\bvocB$-dérivation universelle de $\bvcC^{(r)}$. Elle prolonge le morphisme canonique 
$\bvcF^{(r)}\rightarrow \bvsigma^*(\xi^{-1}\tOmega^1_{\bvoX/\bvoS})$.
Pour tous nombres rationnels $r\geq r'\geq 0$, on a 
\begin{equation}\label{ahttf14i}
p^{r-r'}(\id \otimes \bvalpha^{r,r'}) \circ d_{\bvcC^{(r)}}=d_{\bvcC^{(r')}}\circ \bvalpha^{r,r'}.
\end{equation}

\begin{remas}\label{ahttf15}
Soient $r$ un nombre rationnel $\geq 0$, $n$ un  entier $\geq 1$. 
\begin{itemize}
\item[{\rm (i)}] Pour tout entier $m\geq 0$, les morphismes canoniques $\rS^m_{\ocB_n}(\cF^{(r)}_n)\rightarrow \cC_n^{(r)}$
et $\rS^m_{\bvocB}(\bvcF^{(r)})\rightarrow \bvcC^{(r)}$ sont injectifs. En effet, pour tout entier $m'\geq m$, le morphisme canonique 
$\rS^m_{\ocB_n}(\cF^{(r)}_n)\rightarrow  \rS^{m'}_{\ocB_n}(\cF^{(r)}_n)$
est injectif \eqref{ahttf14c}. Comme les limites injectives filtrantes commutent 
aux limites projectives finies dans $\tE_s^{\mN^\circ}$,
$\rS^m_{\ocB_n}(\cF^{(r)}_n)\rightarrow \cC_n^{(r)}$ est injectif. 
La seconde assertion se déduit de la première par (\cite{agt} III.7.3(i)). 
\item[{\rm (ii)}]  On a $\sigma_n^*(\txi^{-1}\tOmega^1_{\oX_n/\oS_n})=d_{\cC_n^{(r)}}(\cF_n^{(r)}) \subset d_{\cC_n^{(r)}}(\cC_n^{(r)})$
\eqref{ahttf23b}. Par suite, la dérivation $d_{\cC_n^{(r)}}$ est un $\ocB_n$-champ de Higgs à coefficients dans 
$\sigma_n^*(\txi^{-1}\tOmega^1_{\oX_n/\oS_n})$ d'après \ref{MH8}(i).
\item[{\rm (iii)}] On a $\bvsigma^*(\txi^{-1}\tOmega^1_{\bvoX/\bvoS})=d_{\bvcC^{(r)}}(\bvcF^{(r)}) \subset d_{\bvcC^{(r)}}(\bvcC^{(r)})$
\eqref{ahttf14h}. Par suite, la dérivation $d_{\bvcC^{(r)}}$ est un $\bvocB$-champ de Higgs à coefficients dans 
$\bvsigma^*(\txi^{-1}\tOmega^1_{\bvoX/\bvoS})$.  
\end{itemize}
\end{remas}

\begin{prop}\label{ahttf52}
Pour tout nombre rationnel $r\geq 0$, le foncteur 
\begin{equation}\label{ahttf52a}
\bMod(\bvocB)\rightarrow \bMod(\bvcC^{(r)}), \ \ \ M\mapsto M\otimes_{\bvocB} \bvcC^{(r)}
\end{equation}
est exact et fidèle; en particulier, $\bvcC^{(r)}$ est $\bvocB$-plat. 
\end{prop}
Comme le $\co_X$-module $\tOmega^1_{X/S}$ est localement libre de type fini, 
la suite exacte \eqref{ahttf14b} est localement scindée. Un scindage local de cette suite 
induit, pour tout entier $m\geq 0$, un scindage local de la suite exacte \eqref{ahttf14c}.
On en déduit que $\rS^m_{\bvocB}(\bvcF^{(r)})$ est $\bvocB$-plat et que l'homomorphisme 
canonique $\bvocB\rightarrow \bvcC^{(r)}$ admet localement des sections. 
La proposition s'ensuit compte tenu de \eqref{ahttf14d}.

\subsection{}\label{ahttf29}
Considérons le cas absolu, {\em i.e.}, $(\tS,\cM_\tS)=(\cA_2(\oS),\cM_{\cA_2(\oS)})$ \eqref{definf10} et posons \eqref{defing12}
\begin{equation}\label{ahttf29a}
(\tX',\cM_{\tX'})=(\tX,\cM_{\tX})\times_{(\cA_2(\oS),\cM_{\cA_2(\oS)})}(\cA^{\ast}_2(\oS/S),\cM_{\cA^{\ast}_2(\oS/S)})
\end{equation}
où le changement de base est défini par le morphisme $\pr_2$ \eqref{definf7g}. 
On affecte d'un exposant $'$ les objets associés à la $(\cA^{\ast}_2(\oS/S),\cM_{\cA^{\ast}_2(\oS/S)})$-déformation $(\tX',\cM_{\tX'})$ \eqref{ahttf14}. 
On note $K_0$ le corps des fractions de $W$ \eqref{defing1} et $\fd$ la différente de l'extension $K/K_0$ et on pose $\rho=v(\pi\fd)$ \eqref{defing1}.
On désigne par 
\begin{equation}\label{ahttf29b}
\iota\colon \xi\co_C \stackrel{\sim}{\rightarrow} \xi^{\ast}_\pi \co_C
\end{equation}
l'isomorphisme $\co_C$-linéaire tel que le composé 
\begin{equation}\label{ahttf29c}
\xi\co_C \stackrel{\iota}{\longrightarrow} \xi^{\ast}_\pi \co_C \stackrel{\cdot p^\rho}{\longrightarrow} p^\rho \xi^{\ast}_\pi \co_C
\end{equation}
coïncide avec l'isomorphisme \eqref{definf16a}. 
L'isomorphisme $\iota$ induit un isomorphisme $\co_{\bvoX}$-linéaire
\begin{equation}\label{ahttf29d}
\nu\colon (\xi_\pi^{\ast})^{-1}\tOmega^1_{\bvoX/\bvoS}\stackrel{\sim}{\rightarrow}  \xi^{-1}\tOmega^1_{\bvoX/\bvoS}.
\end{equation}
Il résulte aussitôt de \ref{taht9} que pour tout nombre rationnel $r\geq 0$, on a un isomorphisme $\bvocB$-linéaire 
\begin{equation}\label{ahttf29e}
\bvcF'^{(r)}\stackrel{\sim}{\rightarrow} \bvcF^{(r+\rho)}
\end{equation}
qui s'insère dans un diagramme commutatif
\begin{equation}\label{ahttf29f}
\xymatrix{
0\ar[r]&{\bvocB}\ar[r]\ar@{=}[d]&{\bvcF'^{(r)}}\ar[r]\ar[d]&{\bvsigma^*((\xi^{\ast}_\pi)^{-1}\tOmega^1_{\bvoX/\bvoS})}\ar[r]\ar[d]^{\nu}&0\\
0\ar[r]&{\bvocB}\ar[r]&{\bvcF^{(r+\rho)}}\ar[r]&{\bvsigma^*(\xi^{-1}\tOmega^1_{\bvoX/\bvoS})}\ar[r] & 0}
\end{equation}  
On en déduit un isomorphisme de $\bvocB$-algèbres 
\begin{equation}\label{ahttf29g}
\bvcC'^{(r)}\stackrel{\sim}{\rightarrow} \bvcC^{(r+\rho)}.
\end{equation} 
D'après \eqref{taht9n}, le diagramme
\begin{equation}\label{ahttf29h}
\xymatrix{
{\bvcC'^{(r)}}\ar[rr]^-(0.5){d_{\bvcC'^{(r)}}}\ar[d]&&{\bvsigma^*((\xi^{\ast}_\pi)^{-1}\tOmega^1_{\bvoX/\bvoS}) \otimes_{\bvocB} \bvcC'^{(r)}}\ar[d]^{\nu}\\
{\bvcC^{(r+\rho)}}\ar[rr]^-(0.5){d_{\bvcC^{(r+\rho)}}}&&{\bvsigma^*(\xi^{-1}\tOmega^1_{\bvoX/\bvoS}) \otimes_{\bvocB} \bvcC^{(r+\rho)}}}
\end{equation} 
est commutatif.

\subsection{}\label{ahttf17}
On rappelle que $\fX$ désigne le schéma formel complété $p$-adique de $\oX$ \eqref{definf18} et $\txi^{-1}\tOmega^1_{\fX/\cS}$ 
le complété $p$-adique du $\co_\oX$-module $\txi^{-1}\tOmega^1_{\oX/\oS}$ \eqref{definf19a}.
D'après (\cite{ag} 2.1.18.6) et \eqref{ahttf17c}, on a un isomorphisme $\bvocB$-linéaire canonique 
\begin{equation}\label{ahttf17d}
\hupsigma^*(\txi^{-1}\tOmega^1_{\fX/\cS})\stackrel{\sim}{\rightarrow} \bvsigma^*(\txi^{-1}\tOmega^1_{\bvoX/\bvoS}),
\end{equation}
où $\hupsigma$ est le morphisme \eqref{ahttf13e}. 
On désigne par
\begin{equation}\label{ahttf17e}
\txi^{-1}\tOmega^1_{\fX/\cS}\rightarrow \hupsigma_*(\bvsigma^*(\txi^{-1}\tOmega^1_{\bvoX/\bvoS}))
\end{equation}
le morphisme adjoint et par 
\begin{equation}\label{ahttf17f}
\updelta\colon \txi^{-1}\tOmega^1_{\fX/\cS}\rightarrow \rR^1\hupsigma_*(\bvocB)
\end{equation}
le composé de \eqref{ahttf17e} et du morphisme bord de la suite exacte longue de cohomologie déduite
de la suite exacte \eqref{ahttf14b}.

\begin{prop}\label{ahttf20}
Il existe un et un unique isomorphisme de $\co_\fX[\frac 1 p]$-algèbres graduées 
\begin{equation}\label{ahttf20a}
\wedge(\txi^{-1}\tOmega^1_{\fX/\cS}[\frac 1 p])\stackrel{\sim}{\rightarrow} \oplus_{i\geq 0}\rR^i\hupsigma_*(\bvocB)[\frac 1 p]
\end{equation}
dont la composante en degré un est le morphisme $\updelta\otimes_{\mZ_p}\mQ_p$ \eqref{ahttf17f}.
\end{prop} 

Le cas absolu \eqref{definf10} a été démontré dans (\cite{agt} III.11.8). Pour le cas relatif, 
la question étant locale pour la topologie de Zariski de $X$, on peut supposer $X$ affine. La proposition résulte alors du cas absolu,  
compte tenu de (\cite{kato1} 3.14) et \ref{ahttf29}, en particulier de \eqref{ahttf29f}.

\begin{prop}\label{ahttf24}
Soient $r,r'$ deux nombres rationnels tels que $r>r'>0$. Alors,
\begin{itemize}
\item[{\rm (i)}] Pour tout entier $n\geq 1$, l'homomorphisme canonique \eqref{ahttf37c}
\begin{equation}\label{ahttf24a}
\co_{\oX_n}\rightarrow \sigma_{n*}(\cC_n^{(r)})
\end{equation}
est presque-injectif \eqref{defing2}. Notons $\cH^{(r)}_n$ son conoyau. 
\item[{\rm (ii)}] Il existe un nombre rationnel $a>0$ tel que pour tout entier $n\geq 1$, 
le morphisme 
\begin{equation}\label{ahttf24b}
\cH^{(r)}_n\rightarrow \cH^{(r')}_n
\end{equation} 
induit par l'homomorphisme 
$\alpha_n^{r,r'}\colon \cC_n^{(r)}\rightarrow \cC_n^{(r')}$ \eqref{ahttf37e} soit annulé par $p^a$. 
\item[{\rm (iii)}] Il existe un nombre rationnel $b>0$ tel que pour tous entiers $n,q\geq 1$, 
le morphisme canonique
\begin{equation}\label{ahttf24c}
\rR^q\sigma_{n*}(\cC_n^{(r)})\rightarrow \rR^q\sigma_{n*}(\cC_n^{(r')})
\end{equation}
soit annulé par $p^b$.
\end{itemize}
\end{prop}
Le cas absolu \eqref{definf10} a été démontré dans (\cite{agt} III.11.13). Pour le cas relatif, 
la question étant locale pour la topologie de Zariski de $X$, on peut supposer $X$ affine. La proposition résulte alors du cas absolu,  
compte tenu de (\cite{kato1} 3.14) et \ref{ahttf29}.

\begin{cor}\label{ahttf25}
Soient $r,r'$ deux nombres rationnels tels que $r>r'>0$. Alors,
\begin{itemize}
\item[{\rm (i)}] L'homomorphisme canonique de $X_{s,\et}^{\mN^\circ}$
\begin{equation}\label{ahttf25a}
\co_{\bvoX}\rightarrow \bvsigma_*(\bvcC^{(r)})
\end{equation}
est presque-injectif. Notons $\bvcH^{(r)}$ son conoyau. 
\item[{\rm (ii)}] Il existe un nombre rationnel $a>0$ tel que le morphisme 
\begin{equation}\label{ahttf25b}
\bvcH^{(r)}\rightarrow \bvcH^{(r')}
\end{equation} 
induit par l'homomorphisme canonique $\bvalpha^{r,r'}\colon \bvcC^{(r)}\rightarrow \bvcC^{(r')}$ \eqref{ahttf14f}
soit annulé par $p^a$. 
\item[{\rm (iii)}] Il existe un nombre rationnel $b>0$ tel que pour tout entier $q\geq 1$, 
le morphisme canonique de $X_{s,\et}^{\mN^\circ}$
\begin{equation}\label{ahttf25c}
\rR^q\bvsigma_*(\bvcC^{(r)})\rightarrow \rR^q\bvsigma_*(\bvcC^{(r')})
\end{equation}
soit annulé par $p^b$.
\end{itemize}
\end{cor}

Cela résulte de \ref{ahttf24} et (\cite{agt} III.7.3(i) et (III.7.5.5)).

\begin{prop}\label{ahttf26}
Soient $r,r'$ deux nombres rationnels tels que $r>r'>0$. Alors,
\begin{itemize}
\item[{\rm (i)}] L'homomorphisme canonique 
\begin{equation}\label{ahttf26a}
\co_{\fX}\rightarrow \hupsigma_*(\bvcC^{(r)})
\end{equation}
est injectif. Notons $\cL^{(r)}$ son conoyau. 
\item[{\rm (ii)}] Il existe un nombre rationnel $a>0$ tel que le morphisme 
\begin{equation}\label{ahttf26b}
\cL^{(r)}\rightarrow \cL^{(r')}
\end{equation} 
induit par l'homomorphisme canonique $\bvalpha^{r,r'}\colon \bvcC^{(r)}\rightarrow \bvcC^{(r')}$ \eqref{ahttf14f}
soit annulé par $p^a$. 
\item[{\rm (iii)}] Pour tout entier $q\geq 1$, il existe un nombre rationnel $b>0$ tel que  
le morphisme canonique 
\begin{equation}\label{ahttf26c}
\rR^q\hupsigma_*(\bvcC^{(r)})\rightarrow \rR^q\hupsigma_*(\bvcC^{(r')})
\end{equation}
soit annulé par $p^b$.
\end{itemize}
\end{prop}
Le cas absolu \eqref{definf10} a été démontré dans (\cite{agt} III.11.16). Pour le cas relatif, 
la question étant locale pour la topologie de Zariski de $X$, on peut supposer $X$ affine. La proposition résulte alors du cas absolu,  
compte tenu de (\cite{kato1} 3.14) et \ref{ahttf29}.

\begin{cor}\label{ahttf27}
Soient $r$, $r'$ deux nombres rationnels tels que $r>r'>0$. Alors,
\begin{itemize}
\item[{\rm (i)}] L'homomorphisme canonique 
\begin{equation}\label{ahttf27a}
u^r\colon \co_{\fX}[\frac 1 p]\rightarrow \hupsigma_*(\bvcC^{(r)})[\frac 1 p]
\end{equation}
admet (en tant que morphisme $\co_\fX[\frac 1 p]$-linéaire) un inverse à gauche canonique
\begin{equation}
v^r\colon \hupsigma_*(\bvcC^{(r)})[\frac 1 p]\rightarrow \co_{\fX}[\frac 1 p].
\end{equation}
\item[{\rm (ii)}] Le composé 
\begin{equation}\label{ahttf27b}
\hupsigma_*(\bvcC^{(r)})[\frac 1 p]\stackrel{v^r}{\longrightarrow} \co_{\fX}[\frac 1 p]
\stackrel{u^{r'}}{\longrightarrow} \hupsigma_*(\bvcC^{(r')})[\frac 1 p]
\end{equation}
est l'homomorphisme canonique.
\item[{\rm (iii)}] Pour tout entier $q\geq 1$, le morphisme canonique 
\begin{equation}\label{ahttf27c}
\rR^q\hupsigma_*(\bvcC^{(r)})[\frac 1 p]\rightarrow \rR^q\hupsigma_*(\bvcC^{(r')})[\frac 1 p]
\end{equation}
est nul. 
\end{itemize}
\end{cor}

Cela résulte de \ref{ahttf26} (cf. \cite{agt} III.11.17).

\begin{cor}[\cite{agt} III.11.18]\label{ahttf28}
L'homomorphisme canonique 
\begin{equation}\label{ahttf28a}
\co_{\fX}[\frac 1 p]\rightarrow \underset{\underset{r\in \mQ_{>0}}{\longrightarrow}}{\lim}\  \hupsigma_*(\bvcC^{(r)})[\frac 1 p]
\end{equation}
est un isomorphisme, et pour tout entier $q\geq 1$,
\begin{equation}\label{ahttf28b}
\underset{\underset{r\in \mQ_{>0}}{\longrightarrow}}{\lim}\ \rR^q\hupsigma_*(\bvcC^{(r)})[\frac 1 p] =0.
\end{equation}
\end{cor}

\subsection{}\label{ahttf30}
Pour tout nombre rationnel $r\geq 0$, on note encore 
\begin{equation}\label{ahttf30a}
d_{\bvcC^{(r)}}\colon \bvcC^{(r)}\rightarrow \hupsigma^*(\txi^{-1}\tOmega^1_{\fX/\cS})\otimes_{\bvocB}\bvcC^{(r)}
\end{equation}
la $\bvocB$-dérivation induite par $d_{\bvcC^{(r)}}$ \eqref{ahttf14h} et l'isomorphisme \eqref{ahttf17d},
que l'on identifie à la $\bvocB$-dérivation universelle de $\bvcC^{(r)}$.
C'est un $\bvocB$-champ de Higgs à coefficients dans $\hupsigma^*(\txi^{-1}\tOmega^1_{\fX/\cS})$ \eqref{ahttf15}.
On note $\mK^\bullet(\bvcC^{(r)})$ le complexe de Dolbeault
du $\bvocB$-module de Higgs $(\bvcC^{(r)},p^rd_{\bvcC^{(r)}})$
et $\tmK^\bullet(\bvcC^{(r)})$ le complexe de Dolbeault augmenté
\begin{equation}\label{ahttf30b}
\bvocB\rightarrow \mK^0(\bvcC^{(r)})\rightarrow \mK^1(\bvcC^{(r)})\rightarrow \dots 
\rightarrow \mK^n(\bvcC^{(r)})\rightarrow \dots,
\end{equation}
où $\bvocB$ est placé en degré $-1$ et la différentielle $\bvocB\rightarrow \bvcC^{(r)}$ est l'homomorphisme canonique. 

Pour tous nombres rationnels $r\geq r'\geq 0$, on a \eqref{ahttf14i} 
\begin{equation}\label{ahttf30c}
p^r(\id \otimes \bvalpha^{r,r'}) \circ d_{\bvcC^{(r)}}=p^{r'}d_{\bvcC^{(r')}}\circ \bvalpha^{r,r'},
\end{equation}
où $\bvalpha^{r,r'}\colon \bvcC^{(r)}\rightarrow \bvcC^{(r')}$ est l'homomorphisme \eqref{ahttf14f}. Par suite, 
$\bvalpha^{r,r'}$ induit un morphisme de complexes 
\begin{equation}\label{ahttf30d}
\bvupiota^{r,r'}\colon \tmK^\bullet(\bvcC^{(r)})\rightarrow \tmK^\bullet(\bvcC^{(r')}).
\end{equation}

On note $\mK^\bullet_\mQ(\bvcC^{(r)})$ et $\tmK^\bullet_\mQ(\bvcC^{(r)})$ 
les images de $\mK^\bullet(\bvcC^{(r)})$ et $\tmK^\bullet(\bvcC^{(r)})$ dans $\bMod_{\mQ}(\bvocB)$ \eqref{ahttf40a}.
On considèrera ces complexes aussi comme des complexes de $\bIndMod(\bvocB)$ via le foncteur $\alpha_{\bvocB}$ \eqref{ahttf40d}.

\begin{prop}\label{ahttf31}
Pour tous nombres rationnels  $r> r'> 0$ et tout entier $q$, le morphisme canonique \eqref{ahttf30d} 
\begin{equation}\label{ahttf31a}
\rH^q(\bvupiota^{r,r'}_\mQ)\colon \rH^q(\tmK^\bullet_\mQ(\bvcC^{(r)}))\rightarrow 
\rH^q(\tmK^\bullet_\mQ(\bvcC^{(r')}))
\end{equation}
est nul.
\end{prop}

Le cas absolu \eqref{definf10} a été démontré dans (\cite{agt} III.11.22). Pour le cas relatif, 
la question étant locale pour la topologie de Zariski de $X$ (\cite{agt} III.6.7), on peut supposer $X$ affine. La proposition résulte alors du cas absolu,  
compte tenu de (\cite{kato1} 3.14) et \ref{ahttf29}, en particulier de \eqref{ahttf29h}.

\begin{cor}\label{ahttf32}
Soient $r$, $r'$ deux nombres rationnels tels que $r>r'>0$. Alors,
\begin{itemize}
\item[{\rm (i)}] Le morphisme canonique 
\begin{equation}\label{ahttf32a}
u^r\colon \bvocB_\mQ\rightarrow \rH^0(\mK^\bullet_\mQ(\bvcC^{(r)}))
\end{equation}
admet un inverse à gauche canonique 
\begin{equation}\label{ahttf32b}
v^r\colon  \rH^0(\mK^\bullet_\mQ(\bvcC^{(r)})) \rightarrow \bvocB_\mQ.
\end{equation}
\item[{\rm (ii)}] Le composé 
\begin{equation}\label{ahttf32c}
\rH^0(\mK^\bullet_\mQ(\bvcC^{(r)})) \stackrel{v^r}{\longrightarrow} \bvocB_\mQ
\stackrel{u^{r'}}{\longrightarrow} \rH^0(\mK^\bullet_\mQ(\bvcC^{(r')}))
\end{equation}
est le morphisme canonique.
\item[{\rm (iii)}] Pour tout entier $q\geq 1$, le morphisme canonique 
\begin{equation}\label{ahttf32d}
\rH^q(\mK^\bullet_\mQ(\bvcC^{(r)}))\rightarrow 
\rH^q(\mK^\bullet_\mQ(\bvcC^{(r')}))
\end{equation}
est nul.
\end{itemize}
\end{cor}

Cela résulte de \ref{ahttf31} (cf. \cite{agt} III.11.23).

\begin{cor}[\cite{agt} III.11.24]\label{ahttf33}
Le morphisme canonique de complexes de ind-$\bvocB$-modules 
\begin{equation}\label{ahttf33a}
\bvocB_\mQ[0]\rightarrow \underset{\underset{r\in \mQ_{>0}}{\longrightarrow}}{\mlq\mlq\lim \mrq\mrq}\ \mK^\bullet_\mQ(\bvcC^{(r)})
\end{equation}
est un quasi-isomorphisme.
\end{cor}

On rappelle d'abord que $\bIndMod(\bvocB)$ admet des petites limites inductives et que les petites limites inductives filtrantes sont exactes \eqref{indsh6e}. 
La proposition est donc équivalente au fait que le morphisme canonique 
\begin{equation}\label{ahttf33b}
\bvocB_\mQ\rightarrow \underset{\underset{r\in \mQ_{>0}}{\longrightarrow}}{\mlq\mlq\lim \mrq\mrq}\ 
\rH^0(\mK^\bullet_\mQ(\bvcC^{(r)}))
\end{equation}
est un isomorphisme, et que pour tout entier $q\geq 1$, 
\begin{equation}\label{ahttf33c}
\underset{\underset{r\in \mQ_{>0}}{\longrightarrow}}{\mlq\mlq\lim \mrq\mrq}\ \rH^q(\mK^\bullet_\mQ(\bvcC^{(r)}))=0.
\end{equation}
Ces énoncés résultent aussitôt de \ref{ahttf32}.

\begin{remas}\label{ahttf34}
Bien que les limites inductives filtrantes ne soient pas a priori représentables dans la catégorie $\bMod_\mQ(\bvocB)$, 
il résulte de \ref{ahttf32} que dans cette catégorie, le morphisme canonique 
\begin{equation}\label{ahttf34a}
\bvocB_\mQ\rightarrow \underset{\underset{r\in \mQ_{>0}}{\longrightarrow}}{\lim}\ 
\rH^0(\mK^\bullet_\mQ(\bvcC^{(r)}))
\end{equation}
est un isomorphisme, et pour tout entier $q\geq 1$, 
\begin{equation}\label{ahttf34b}
\underset{\underset{r\in \mQ_{>0}}{\longrightarrow}}{\lim}\ \rH^q(\mK^\bullet_\mQ(\bvcC^{(r)}))=0.
\end{equation}
\end{remas}

\section{Ind-modules de Dolbeault}\label{indmdlb}

\subsection{}\label{indmdlb1}
Soit $r$ un nombre rationnel $\geq 0$. On rappelle que $\hupsigma$ désigne le morphisme \eqref{ahttf13e} et 
\begin{equation}\label{indmdlb1a}
d_{\bvcC^{(r)}}\colon \bvcC^{(r)}\rightarrow \hupsigma^*(\txi^{-1}\tOmega^1_{\fX/\cS})\otimes_{\bvocB}\bvcC^{(r)}
\end{equation}
la $\bvocB$-dérivation universelle de $\bvcC^{(r)}$ \eqref{ahttf30a}, 
qui est un $\bvocB$-champ de Higgs à coefficients dans $\hupsigma^*(\txi^{-1}\tOmega^1_{\fX/\cS})$.
On désigne par $\bIndMC(\bvcC^{(r)}/\bvocB)$ la catégorie des ind-$\bvcC^{(r)}$-modules à $p^r$-connexion intégrable relativement à l'extension
$\bvcC^{(r)}/\bvocB$ (cf. \ref{indsh28} et \ref{indsh35}). 
Le lecteur prendra garde que malgré la notation, cette catégorie n'est pas la catégorie des ind-objets d'une autre catégorie \eqref{indsh3}. 
Chaque objet de $\bIndMC(\bvcC^{(r)}/\bvocB)$ est un ind-$\bvocB$-module de Higgs
à coefficients dans $\hupsigma^*(\txi^{-1}\tOmega^1_{\fX/\cS})$ d'après \ref{indsh38}(i). 
On peut donc lui associer un complexe de Dolbeault dans $\bIndMod(\bvocB)$ \eqref{indsh30c}.  

Considérons les foncteurs 
\begin{equation}\label{indmdlb1b}
\rI\fS^{(r)}\colon
\begin{array}[t]{clcr}
\bIndMod(\bvocB)&\rightarrow &\bIndMC(\bvcC^{(r)}/\bvocB)\\
\cM&\mapsto& (\bvcC^{(r)}\otimes_{\bvocB}\cM,p^rd_{\bvcC^{(r)}}\otimes_{\bvocB}\id_{\cM}),
\end{array}
\end{equation}
\begin{equation}\label{indmdlb1c}
\rI\cK^{(r)}\colon 
\begin{array}[t]{clcr}
\bIndMC(\bvcC^{(r)}/\bvocB)&\rightarrow &\bIndMod(\bvocB)\\
(\cF,\nabla)&\mapsto& \ker(\nabla). 
\end{array}
\end{equation}
Il est clair que $\rI\fS^{(r)}$ est un adjoint à gauche de $\rI \cK^{(r)}$ \eqref{indsh15h}. 

On désigne par $\bIndMH(\co_\fX,\txi^{-1}\tOmega^1_{\fX/\cS})$ 
(resp. $\bIndMH(\bvocB,\hupsigma^*(\txi^{-1}\tOmega^1_{\fX/\cS}))$) la catégorie des ind-$\co_\fX$-modules de Higgs (resp. ind-$\bvocB$-modules de Higgs)
à coefficients dans $\txi^{-1}\tOmega^1_{\fX/\cS}$ (resp. $\hupsigma^*(\txi^{-1}\tOmega^1_{\fX/\cS})$) \eqref{indsh30}. 
Le foncteur $\rI\hupsigma^*$ \eqref{ahttf43c} induit un foncteur que l'on note encore
\begin{equation}\label{indmdlb1d}
\begin{array}[t]{clcr}
\rI\hupsigma^*\colon \bIndMH(\co_\fX,\txi^{-1}\tOmega^1_{\fX/\cS})&\rightarrow& \bIndMH(\bvocB, \hupsigma^*(\txi^{-1}\tOmega^1_{\fX/\cS}))\\
(\cN,\theta)&\mapsto& (\rI\hupsigma^*(\cN),\rI\hupsigma^*(\theta)).
\end{array}
\end{equation}
Compte tenu de \ref{indsh41} et \ref{indsh40}(i), le foncteur $\rI\hupsigma_*$ \eqref{ahttf43d} induit un foncteur que l'on note encore
\begin{equation}\label{indmdlb1e}
\begin{array}[t]{clcr}
\rI\hupsigma_*\colon \bIndMH(\bvocB, \hupsigma^*(\txi^{-1}\tOmega^1_{\fX/\cS})) &\rightarrow& \bIndMH(\co_\fX,\txi^{-1}\tOmega^1_{\fX/\cS})\\ 
(\cM,\theta)&\mapsto& (\rI\hupsigma_*(\cM),\rI\hupsigma_*(\theta)).
\end{array}
\end{equation}
Il résulte de \ref{indsh40}(ii) que le foncteur $\rI\hupsigma^*$ \eqref{indmdlb1d} est un adjoint à gauche de $\rI\hupsigma_*$ \eqref{indmdlb1e}.

Compte tenu de \ref{indsh38}(ii), le foncteur \eqref{indmdlb1d} induit un foncteur 
\begin{equation}\label{indmdlb1f}
\begin{array}[t]{clcr}
\rI\hupsigma^{(r)*}\colon\bIndMH(\co_\fX,\txi^{-1}\tOmega^1_{\fX/\cS})&\rightarrow& \bIndMC(\bvcC^{(r)}/\bvocB)\\
(\cN,\theta)&\mapsto&(\bvcC^{(r)}\otimes_{\bvocB}\rI\hupsigma^*(\cN),p^rd_{\bvcC^{(r)}} \otimes_{\bvocB}\id+\id \otimes_{\bvocB} \rI\hupsigma^*(\theta)).
\end{array}
\end{equation}
Le foncteur $\rI\hupsigma_*$ \eqref{indmdlb1e} induit un foncteur 
\begin{equation}\label{indmdlb1g}
\rI\hupsigma^{(r)}_*\colon \bIndMC(\bvcC^{(r)}/\bvocB) \rightarrow \bIndMH(\co_\fX,\txi^{-1}\tOmega^1_{\fX/\cS}), \ \ \ (\cF,\nabla)\mapsto (\rI\hupsigma_*(\cF),
\rI\hupsigma_*(\nabla)).
\end{equation}
Le foncteur $\rI\hupsigma^{(r)*}$ \eqref{indmdlb1f} est un adjoint à gauche de $\rI\hupsigma^{(r)}_*$ \eqref{indmdlb1g}. 

Compte tenu de \eqref{indsh15g}, le foncteur $\kappa_{\co_\fX}$ \eqref{ahttf43g} induit un foncteur que l'on note encore
\begin{equation}\label{indmdlb1h}
\kappa_{\co_\fX}\colon \bIndMH(\co_\fX,\txi^{-1}\tOmega^1_{\fX/\cS})\rightarrow \bMH(\co_\fX,\txi^{-1}\tOmega^1_{\fX/\cS}).
\end{equation}
On désigne par $\rI\vupsigma^{(r)}_*$ le foncteur composé 
\begin{equation}\label{indmdlb1i}
\vupsigma^{(r)}_*=\kappa_{\co_\fX} \circ \rI\hupsigma^{(r)}_*\colon \bIndMC(\bvcC^{(r)}/\bvocB) \rightarrow \bMH(\co_\fX,\txi^{-1}\tOmega^1_{\fX/\cS}).
\end{equation}

\subsection{}\label{indmdlb2}
Reprenons les notations de \ref{definf19}. 
Compte tenu de \eqref{indsh20d}, le foncteur $\upalpha_{\co_\fX}$ \eqref{ahttf40c} induit un foncteur 
\begin{equation}\label{indmdlb2d}
\begin{array}[t]{clcr}
\bIH_\mQ(\co_\fX,\txi^{-1}\tOmega^1_{\fX/\cS})&\rightarrow& \bIndMH(\co_\fX,\txi^{-1}\tOmega^1_{\fX/\cS})\\
(\cM,\cN,u,\theta)&\mapsto& (\upalpha_{\co_\fX}(\cM_\mQ),(\id \otimes \upalpha_{\co_\fX}(u_\mQ)^{-1})\circ\upalpha_{\co_\fX}(\theta_\mQ)).
\end{array}
\end{equation}
Celui-ci est pleinement fidèle d'après \eqref{indsh5d} et \eqref{indsh49c}. Composant avec les foncteurs 
\begin{equation} \label{indmdlb2e}
\bMH^\coh(\co_\fX[\frac 1 p], \txi^{-1}\tOmega^1_{\fX/\cS})\stackrel{\sim}{\rightarrow}
\bIH^\coh_\mQ(\co_\fX,\txi^{-1}\tOmega^1_{\fX/\cS})\rightarrow \bIH_\mQ(\co_\fX,\txi^{-1}\tOmega^1_{\fX/\cS}),
\end{equation}
où la première flèche est un quasi-inverse de l'équivalence de catégories \eqref{definf19d} et la seconde flèche est l'injection canonique, 
on obtient un foncteur pleinement fidèle
\begin{equation}\label{indmdlb2f} 
\bMH^\coh(\co_\fX[\frac 1 p], \txi^{-1}\tOmega^1_{\fX/\cS})\rightarrow \bIndMH(\co_\fX,\txi^{-1}\tOmega^1_{\fX/\cS}).
\end{equation}
On identifiera $\bMH^\coh(\co_\fX[\frac 1 p], \txi^{-1}\tOmega^1_{\fX/\cS})$ à une sous-catégorie pleine de $\bIndMH(\co_\fX,\txi^{-1}\tOmega^1_{\fX/\cS})$
par ce foncteur qu'on omettra des notations. On considérera donc tout
$\co_\fX[\frac 1 p]$-module de Higgs cohérent à coefficients dans $\txi^{-1}\tOmega^1_{\fX/\cS}$ comme 
un ind-$\co_\fX$-module de Higgs à coefficients dans $\txi^{-1}\tOmega^1_{\fX/\cS}$.

Pour tous $(\cN,\theta)\in \ob(\bMH^\coh(\co_\fX[\frac 1 p], \txi^{-1}\tOmega^1_{\fX/\cS}))$ et $(\cF,\nabla)\in \ob(\bIndMC(\bvcC^{(r)}/\bvocB))$, on a un
application canonique bifonctorielle
\begin{eqnarray}\label{indmdlb2g} 
\lefteqn{\Hom_{\bIndMC(\bvcC^{(r)}/\bvocB)}(\rI \hupsigma^{(r)*}(\cN,\theta),(\cF,\nabla))\rightarrow}\\
&&\Hom_{\bMH^\coh(\co_\fX, \txi^{-1}\tOmega^1_{\fX/\cS})}((\cN,\theta),\vupsigma^{(r)}_*(\cF,\nabla)),\nonumber
\end{eqnarray}
induite par l'isomorphisme d'adjonction entre $\rI \hupsigma^{(r)*}$ et $\rI \hupsigma^{(r)}_*$ et l'isomorphisme \eqref{ahttf50d}. 
On appelle abusivement {\em l'adjoint} d'un morphisme $\rI \hupsigma^{(r)*}(\cN,\theta)\rightarrow (\cF,\nabla)$ de $\bIndMC(\bvcC^{(r)}/\bvocB)$,
son image par l'application \eqref{indmdlb2g}. 
Il résulte aussitôt de l'injectivité de l'application \eqref{ahttf50g} que l'application \eqref{indmdlb2g} est injective.

\subsection{}\label{indmdlb3}
Soient $r$, $r'$ deux nombres rationnels tels que $r\geq r'\geq 0$, $(\cF,\nabla)$ un ind-$\bvcC^{(r)}$-module 
à $p^r$-connexion intégrable relativement à l'extension $\bvcC^{(r)}/\bvocB$. 
D'après \ref{indsh43} et \eqref{ahttf14i}, l'ind-$\bvcC^{(r')}$-module $\bvcC^{(r')}\otimes_{\bvcC^{(r)}}\cF$ est alors canoniquement muni d'une $p^{r'}$-connexion intégrable 
$\nabla'$ relativement à l'extension $\bvcC^{(r')}/\bvocB$. On définit ainsi un foncteur  
\begin{equation}\label{indmdlb3a}
\rI\varepsilon^{r,r'}\colon 
\begin{array}[t]{clcr}
\bIndMC(\bvcC^{(r)}/\bvocB)&\rightarrow& \bIndMC(\bvcC^{(r')}/\bvocB)\\
(\cF,\nabla)&\mapsto& (\bvcC^{(r')}\otimes_{\bvcC^{(r)}}\cF,\nabla').
\end{array}
\end{equation}
On a un isomorphisme canonique de foncteurs de $\bIndMod(\bvocB)$ dans $\bIndMC(\bvcC^{(r')}/\bvocB)$ 
\begin{equation}\label{indmdlb3b}
\rI\varepsilon^{r,r'}\circ \rI\fS^{(r)}\stackrel{\sim}{\longrightarrow} \rI\fS^{(r')}.
\end{equation}
On a un isomorphisme canonique de foncteurs de $\bIndMH(\co_\fX,\txi^{-1}\tOmega^1_{\fX/\cS})$ dans $\bIndMC(\bvcC^{(r')}/\bvocB)$ 
\begin{equation}\label{indmdlb3c}
\rI\varepsilon^{r,r'}\circ \rI\hupsigma^{(r)*}\stackrel{\sim}{\longrightarrow} \rI\hupsigma^{(r')*}.
\end{equation}

Le diagramme 
\begin{equation}\label{indmdlb3d}
\xymatrix{
{\cF}\ar[r]^-(0.5){\nabla}\ar[d]_{\bvalpha^{r,r'}\otimes_{\bvcC^{(r)}} \id}&
{\hupsigma^*(\txi^{-1}\tOmega^1_{\fX/\cS})\otimes_{\bvocB}\cF}\ar[d]^{\id\otimes_{\bvocB}\bvalpha^{r,r'}\otimes_{\bvcC^{(r)}} \id}\\
{\bvcC^{(r')}\otimes_{\bvcC^{(r)}}\cF}\ar[r]^-(0.5){\nabla'}&
{\hupsigma^*(\txi^{-1}\tOmega^1_{\fX/\cS})\otimes_{\bvocB}\bvcC^{(r')}\otimes_{\bvcC^{(r)}}\cF}}
\end{equation}
est commutatif \eqref{indsh43}. On en déduit un morphisme canonique de foncteurs de $\bIndMC(\bvcC^{(r)}/\bvocB)$ dans $\bIndMod(\bvocB)$
\begin{equation}\label{indmdlb3e}
\rI\cK^{(r)}\rightarrow \rI\cK^{(r')}\circ \rI\varepsilon^{r,r'},
\end{equation}
et un morphisme canonique de foncteurs de $\bIndMC(\bvcC^{(r)}/\bvocB)$ dans $\bIndMH(\co_\fX,\txi^{-1}\tOmega^1_{\fX/\cS})$ 
\begin{equation}\label{indmdlb3f}
\rI\hupsigma^{(r)}_*\longrightarrow \rI\hupsigma^{(r')}_*\circ \rI\varepsilon^{r,r'}.
\end{equation}

Pour tout nombre rationnel $r''$ tel $r'\geq r''\geq 0$, on a un isomorphisme canonique de foncteurs 
de $\bIndMC(\bvcC^{(r)}/\bvocB)$ dans $\bIndMC(\bvcC^{(r'')}/\bvocB)$ 
\begin{equation}\label{indmdlb3g}
\rI\varepsilon^{r',r''}\circ \rI\varepsilon^{r,r'}\stackrel{\sim}{\rightarrow}\rI\varepsilon^{r,r''}.
\end{equation}

\begin{defi}\label{indmdlb4}
Soient $\cM$ un ind-$\bvocB$-module, 
$\cN$ un $\co_\fX[\frac 1 p]$-fibré de Higgs à coefficients dans $\txi^{-1}\tOmega^1_{\fX/\cS}$ \eqref{definf20},
que l'on considère aussi comme un ind-$\co_\fX$-module de Higgs à coefficients dans $\txi^{-1}\tOmega^1_{\fX/\cS}$ \eqref{indmdlb2f}.
\begin{itemize}
\item[(i)] Soit $r$ un nombre rationnel $>0$.  On dit que $\cM$ et  $\cN$ sont {\em $r$-associés} 
s'il existe un isomorphisme de $\bIndMC(\bvcC^{(r)}/\bvocB)$ 
\begin{equation}\label{indmdlb4a}
\alpha\colon \rI\hupsigma^{(r)*}(\cN) \stackrel{\sim}{\rightarrow}\rI\fS^{(r)}(\cM).
\end{equation}
On dit alors aussi que le triplet $(\cM,\cN,\alpha)$ est {\em $r$-admissible}. 
\item[(ii)] On dit que $\cM$ et  $\cN$ sont {\em associés} s'il existe un nombre rationnel $r>0$ tel que 
$\cM$ et  $\cN$ soient $r$-associés.
\end{itemize}
\end{defi}

On notera que pour tous nombres rationnels $r\geq r'>0$, 
si $\cM$ et  $\cN$ sont $r$-associés, ils sont $r'$-associés, compte tenu de \eqref{indmdlb3b} et \eqref{indmdlb3c}. 

\begin{defi}\label{indmdlb5}
\
\begin{itemize}
\item[(i)] On dit qu'un ind-$\bvocB$-module est {\em de Dolbeault} s'il est associé à  
un $\co_\fX[\frac 1 p]$-fibré de Higgs à coefficients dans $\txi^{-1}\tOmega^1_{\fX/\cS}$.
\item[(ii)] On dit qu'un $\co_\fX[\frac 1 p]$-fibré de Higgs à coefficients dans $\txi^{-1}\tOmega^1_{\fX/\cS}$ est {\em soluble}
s'il est associé à un ind-$\bvocB$-module. 
\end{itemize}
\end{defi}

Ces notions dépendent a priori de la déformation $(\tX,\cM_\tX)$ fixée dans \ref{defing12}.
On désigne par $\bIndMod^\Dolb(\bvocB)$ la sous-catégorie pleine de $\bIndMod(\bvocB)$
formée des ind-$\bvocB$-modules de Dolbeault, et par $\bMH^\sol(\co_\fX[\frac 1 p], \txi^{-1}\tOmega^1_{\fX/\cS})$ 
la sous-catégorie pleine de  $\bMH(\co_\fX[\frac 1 p], \txi^{-1}\tOmega^1_{\fX/\cS})$
formée des $\co_\fX[\frac 1 p]$-fibrés de Higgs solubles à coefficients dans $\txi^{-1}\tOmega^1_{\fX/\cS}$.  

\begin{prop}\label{indmdlb6}
Tout ind-$\bvocB$-module de Dolbeault est rationnel \eqref{ahttf49} et plat \eqref{indsh46}.
\end{prop}

En effet, soient $\cM$ un ind-$\bvocB$-module de Dolbeault, 
$\cN$ un $\co_\fX[\frac 1 p]$-fibré de Higgs à coefficients dans 
$\txi^{-1}\tOmega^1_{\fX/\cS}$, $r$ est un nombre rationnel $>0$,
\begin{equation}\label{indmdlb6a}
\rI \hupsigma^{(r)*}(\cN) \stackrel{\sim}{\rightarrow}\rI\fS^{(r)}(\cM)
\end{equation} 
un isomorphisme de $\bIndMC(\bvcC^{(r)}/\bvocB)$. 

Montrons d'abord que la multiplication par $p$ sur $\cM$ est un isomorphisme. 
La question étant locale pour la topologie étale de $\fX$ d'après \ref{indsh39}(iii), 
on peut supposer que l'homomorphisme canonique $\bvocB\rightarrow \bvcC^{(r)}$ admet des sections (cf. la preuve de \ref{ahttf52}). 
Comme la multiplication par $p$ sur $\cM\otimes_{\bvocB}\bvcC^{(r)}$ est un isomorphisme \eqref{indmdlb6a}, 
la multiplication par $p$ sur $\cM$ est aussi un isomorphisme. 

Montrons ensuite que le ind-$\bvocB$-module $\rI \hupsigma^*(\cN)$ est plat. Compte tenu de \ref{indsh39}(iv), 
on peut supposer que $\cN$ est libre de type fini sur $\co_\fX[\frac 1 p]$. 
L'assertion recherchée résulte alors du fait que $\rI \hupsigma^*(\co_\fX[\frac 1 p])=\bvocB_\mQ$ \eqref{indsh13a}
est un ind-$\bvocB$-module plat \eqref{indsh46a}.
Par suite, le ind-$\bvcC^{(r)}$-module $\rI \hupsigma^*(\cN)\otimes_{\bvocB}\bvcC^{(r)}$ est plat \eqref{indsh46b}
et il en est alors de même de $\cM\otimes_{\bvocB}\bvcC^{(r)}$ \eqref{indmdlb6a}. 
On en déduit que le ind-$\bvocB$-module $\cM$ est plat en vertu de \ref{indsh46c} et \ref{ahttf52}.

\subsection{}\label{indmdlb7}
Pour tout ind-$\bvocB$-module $\cM$ et tous nombres rationnels $r\geq r'\geq 0$, 
le morphisme \eqref{indmdlb3f} et l'isomorphisme \eqref{indmdlb3b} induisent un morphisme de 
$\bIndMH(\co_\fX, \txi^{-1}\tOmega^1_{\fX/\cS})$
\begin{equation}\label{indmdlb7a}
\rI\hupsigma^{(r)}_*(\rI\fS^{(r)}(\cM))\rightarrow \rI\hupsigma^{(r')}_*(\rI\fS^{(r')}(\cM)).
\end{equation}
On obtient ainsi un petit système inductif filtrant $(\rI\hupsigma^{(r)}_*(\rI\fS^{(r)}(\cM)))_{r\in \mQ_{\geq 0}}$. 
On désigne par $\rI\cH$ le foncteur 
\begin{equation}\label{indmdlb7b}
\rI\cH\colon \bIndMod(\bvocB)\rightarrow \bIndMH(\co_\fX, \txi^{-1}\tOmega^1_{\fX/\cS}), \ \ \ \cM\mapsto 
\underset{\underset{r\in \mQ_{>0}}{\longrightarrow}}{\mlq\mlq\lim \mrq\mrq}\ \rI\hupsigma^{(r)}_*(\rI\fS^{(r)}(\cM)). 
\end{equation}
Composant avec le foncteur $\kappa_{\co_\fX}$ \eqref{indmdlb1h}, on obtient le foncteur
\begin{equation}\label{indmdlb7c}
\cH=\kappa_{\co_\fX}\circ \rI\cH\colon \bIndMod(\bvocB)\rightarrow \bMH(\co_\fX, \txi^{-1}\tOmega^1_{\fX/\cS}). 
\end{equation}
Comme le foncteur $\kappa_{\co_\fX}$ commute aux petites limites inductives filtrantes (\cite{ks2} 6.3.1), pour tout ind-$\bvocB$-module $\cM$, on a
\begin{equation}\label{indmdlb7d}
\cH(\cM)= \underset{\underset{r\in \mQ_{>0}}{\longrightarrow}}{\lim}\ \vupsigma^{(r)}_*(\rI\fS^{(r)}(\cM)),
\end{equation}
où le système inductif $(\vupsigma^{(r)}_*(\rI\fS^{(r)}(\cM)))_{r\in \mQ_{\geq 0}}$ est induit par \eqref{indmdlb7a}. Les foncteurs $\rI\cH$ et $\cH$
s'appliquent en particulier aux $\bvocB_\mQ$-modules \eqref{ahttf40d}. 

Pour tout objet $\cN$ de $\bIndMH(\co_\fX, \txi^{-1}\tOmega^1_{\fX/\cS})$ 
et tous nombres rationnels $r\geq r'\geq 0$, 
le morphisme \eqref{indmdlb3e} et l'isomorphisme \eqref{indmdlb3c} induisent un morphisme de $\bIndMod(\bvocB)$
\begin{equation}\label{indmdlb7e}
\rI\cK^{(r)}(\rI\hupsigma^{(r)*}(\cN))\rightarrow \rI\cK^{(r')}(\rI\hupsigma^{(r')*}(\cN)).
\end{equation}
On obtient ainsi un petit système inductif filtrant $(\rI\cK^{(r)}(\rI\hupsigma^{(r)*}(\cN)))_{r\geq  0}$. 
On désigne par $\rI\cV$ le foncteur 
\begin{equation}\label{indmdlb7f}
\rI\cV\colon \bIndMH(\co_\fX, \txi^{-1}\tOmega^1_{\fX/\cS})\rightarrow \bIndMod(\bvocB), \ \ \ \cN\mapsto 
\underset{\underset{r\in \mQ_{>0}}{\longrightarrow}}{\mlq\mlq\lim \mrq\mrq}\ \rI\cK^{(r)}(\rI\hupsigma^{(r)*}(\cN)). 
\end{equation}
On désigne par $\cV$ le composé des foncteurs $\rI\cV$ et \eqref{indmdlb2f},
\begin{equation}\label{indmdlb7g}
\cV\colon \bMH^\coh(\co_\fX[\frac 1 p], \txi^{-1}\tOmega^1_{\fX/\cS})\rightarrow \bIndMH(\co_\fX, \txi^{-1}\tOmega^1_{\fX/\cS})
\rightarrow \bIndMod(\bvocB).
\end{equation}

\begin{lem}\label{indmdlb8}
On a un isomorphisme canonique de $\co_\fX$-modules de Higgs à coefficients dans $\txi^{-1}\tOmega^1_{\fX/\cS}$,
\begin{equation}\label{indmdlb8a}
(\co_\fX[\frac 1 p],0)\stackrel{\sim}{\rightarrow}\cH(\bvocB_\mQ).
\end{equation}
\end{lem}

En effet, pour tout nombre rationnel $r\geq 0$, on a un isomorphisme canonique \eqref{ahttf43i}
\begin{equation}
\vupsigma_*(\bvcC^{(r)}_\mQ)\stackrel{\sim}{\rightarrow} \hupsigma_*(\bvcC^{(r)})[\frac 1 p],
\end{equation}
où $\bvcC^{(r)}_\mQ$ est considéré comme un ind-$\bvocB$-module \eqref{ahttf40d}. 
La proposition résulte alors de \ref{ahttf28}.

\begin{lem}\label{indmdlb11}
Soient $\cN$ un $\co_\fX[\frac 1 p]$-fibré de Higgs à coefficients dans $\txi^{-1}\tOmega^1_{\fX/\cS}$ \eqref{definf20},
$r$ un nombre rationnel $\geq 0$.   
On a alors un isomorphisme canonique de $\bIndMH(\co_\fX, \txi^{-1}\tOmega^1_{\fX/\cS})$
\begin{equation}\label{indmdlb11a}
\cN\otimes_{\co_\fX}\rI\hupsigma^{(r)}_*(\rI\fS^{(r)}(\bvocB))\stackrel{\sim}{\rightarrow}
\rI\hupsigma^{(r)}_*(\rI\hupsigma^{(r)*}(\cN)),
\end{equation}
où le membre de gauche est le produit tensoriel des ind-modules de Higgs \eqref{indsh30g}. 
Appliquant le foncteur $\kappa_{\co_\fX}$ \eqref{indmdlb1h}, on obtient un isomorphisme de $\bMH(\co_\fX, \txi^{-1}\tOmega^1_{\fX/\cS})$
\begin{equation}\label{indmdlb11b}
\gamma^{(r)}\colon \cN\otimes_{\co_\fX}\vupsigma^{(r)}_*(\rI\fS^{(r)}(\bvocB))\stackrel{\sim}{\rightarrow}\vupsigma^{(r)}_*(\rI\hupsigma^{(r)*}(\cN)),
\end{equation}
où le membre de gauche est le produit tensoriel des ind-modules de Higgs \eqref{MH2d}. 
De plus, on a les propriétés suivantes:
\begin{itemize}
\item[{\rm (i)}] Le morphisme 
\begin{equation}\label{indmdlb11c}
\cN\rightarrow \vupsigma^{(r)}_*(\rI\hupsigma^{(r)*}(\cN))
\end{equation}
induit par \eqref{indmdlb11b} et le morphisme canonique $\co_\fX\rightarrow \vupsigma^{(r)}_*(\rI\fS^{(r)}(\bvocB))$,
est l'adjoint du morphisme identité de $\rI\hupsigma^{(r)*}(\cN)$ \eqref{indmdlb2g}.
\item[{\rm (ii)}] Pour tout nombre rationnel $r'$ tel que $r\geq r'\geq 0$,  le diagramme 
\begin{equation}\label{indmdlb11cd}
\xymatrix{
{\cN\otimes_{\co_\fX}\vupsigma^{(r)}_*(\rI\fS^{(r)}(\bvocB))}\ar[r]^-(0.5){\gamma^{(r)}}\ar[d]&
{\vupsigma^{(r)}_*(\rI\hupsigma^{(r)*}(\cN))}\ar[d]\\
{\cN\otimes_{\co_\fX}\vupsigma^{(r')}_*(\rI\fS^{(r')}(\bvocB))}\ar[r]^-(0.5){\gamma^{(r')}}&
{\vupsigma^{(r')}_*(\rI\hupsigma^{(r')*}(\cN))}}
\end{equation}
où les flèches verticales sont induites par le morphisme \eqref{indmdlb3f} et par les isomorphismes \eqref{indmdlb3b} 
et \eqref{indmdlb3c}, est commutatif. 
\end{itemize}
\end{lem}

En effet, d'après \ref{ahttf51}, on a des isomorphismes canoniques de ind-$\co_\fX$-modules 
\begin{eqnarray}
\cN\otimes_{\co_\fX}\rI\hupsigma_*(\bvcC^{(r)})&\stackrel{\sim}{\rightarrow}&\rI\hupsigma_*(\rI\hupsigma^*(\cN)\otimes_{\bvocB}
\bvcC^{(r)}),\label{indmdlb11d}\\
\ \ \ \ \ \ \txi^{-1}\tOmega^1_{\fX/\cS}\otimes_{\co_\fX}\rI\hupsigma_*(\rI\hupsigma^*(\cN)\otimes_{\bvocB} \bvcC^{(r)})
&\stackrel{\sim}{\rightarrow}&\rI\hupsigma_*(\rI\hupsigma^*(\txi^{-1}\tOmega^1_{\fX/\cS}\otimes_{\co_\fX}\cN)\otimes_{\bvocB}
\bvcC^{(r)}),\label{indmdlb11e}\\
\txi^{-1}\tOmega^1_{\fX/\cS}\otimes_{\co_\fX}\cN\otimes_{\co_\fX}\rI\hupsigma_*(\bvcC^{(r)})
&\stackrel{\sim}{\rightarrow}&\rI\hupsigma_*(\rI\hupsigma^*(\txi^{-1}\tOmega^1_{\fX/\cS}\otimes_{\co_\fX}\cN)\otimes_{\bvocB}
\bvcC^{(r)}) \label{indmdlb11f}
\end{eqnarray}
où le dernier l'isomorphisme est induit par les deux premiers d'après \ref{indsh40}(i). 
De plus, compte tenu du caractère bifonctoriel de l'isomorphisme \eqref{ahttf51a}, le diagramme 
\[
\xymatrix{
{\cN\otimes_{\co_\fX}\rI\hupsigma_*(\bvcC^{(r)})}\ar[r]\ar[d]_{\theta\otimes \id+p^r\id\otimes \rI\hupsigma_*(d_{\bvcC^{(r)}})}&
{\rI\hupsigma_*(\rI\hupsigma^*(\cN)\otimes_{\bvocB} \bvcC^{(r)})}\ar[d]^{\rI\hupsigma_*(\rI\hupsigma^*(\theta)\otimes \id+p^r\id\otimes d_{\bvcC^{(r)}})}\\
{\txi^{-1}\tOmega^1_{\fX/\cS}\otimes_{\co_\fX}\cN\otimes_{\co_\fX}\rI\hupsigma_*(\bvcC^{(r)})}\ar[r]&
{\rI\hupsigma_*(\rI\hupsigma^*(\txi^{-1}\tOmega^1_{\fX/\cS}\otimes_{\co_\fX}\cN)\otimes_{\bvocB}\bvcC^{(r)})}}
\]
où $\theta$ est le champ de Higgs de $\cN$, est commutatif. 
On prend alors pour morphisme \eqref{indmdlb11a} l'isomorphisme \eqref{indmdlb11d}. 
Compte tenu de \eqref{indsh15g} et \eqref{ahttf50d}, l'image de \eqref{indmdlb11a} par $\kappa_{\co_\fX}$
s'identifie à l'isomorphisme \eqref{indmdlb11b}.

La proposition (i) résulte de \ref{indsh40}(ii).
La proposition (ii) est une conséquence du caractère bifonctoriel de l'isomorphisme \eqref{indmdlb11b}.

\subsection{}\label{indmdlb9}
Soient $r$ un nombre rationnel $>0$, $(\cM,\cN,\alpha)$ un triplet $r$-admissible \eqref{indmdlb4}. 
Pour tout nombre rationnel $r'$ tel que $0< r'\leq r$, on désigne par 
\begin{equation}\label{indmdlb9a}
\alpha^{(r')}\colon \rI\hupsigma^{(r')*}(\cN) \stackrel{\sim}{\rightarrow}\rI\fS^{(r')}(\cM)
\end{equation}
l'isomorphisme de $\bIndMC(\bvcC^{(r')}/\bvocB)$ induit par $\rI\varepsilon^{r,r'}(\alpha)$ et les isomorphismes \eqref{indmdlb3b} et \eqref{indmdlb3c}, 
et par 
\begin{equation}\label{indmdlb9b}
\beta^{(r')}\colon \cN \rightarrow \vupsigma^{(r')}_*(\rI\fS^{(r')}(\cM))
\end{equation}
son adjoint \eqref{indmdlb2g}.

\begin{prop}\label{indmdlb10}
Les hypothèses étant celles de \ref{indmdlb9}, soient, de plus, $r'$, $r''$ deux nombres rationnels 
tels que $0<r''< r'\leq r$. Alors,
\begin{itemize}
\item[{\rm (i)}] Le morphisme composé 
\begin{equation}\label{indmdlb10a}
\cN\stackrel{\beta^{(r')}}{\longrightarrow} \vupsigma^{(r')}_*(\rI\fS^{(r')}(\cM))\longrightarrow \cH(\cM),
\end{equation}
où la seconde flèche est le morphisme canonique \eqref{indmdlb7d}, est un isomorphisme, indépendant de $r'$. 
\item[{\rm (ii)}] Le morphisme composé 
\begin{equation}\label{indmdlb10b}
\vupsigma^{(r')}_*(\rI\fS^{(r')}(\cM))\longrightarrow \cH(\cM)\stackrel{\sim}{\longrightarrow}\cN\stackrel{\beta^{(r'')}}{\longrightarrow} 
\vupsigma^{(r'')}_*(\rI\fS^{(r'')}(\cM)),
\end{equation}
où la première flèche est le morphisme canonique \eqref{indmdlb7d} et la deuxième flèche est l'isomorphisme inverse de 
\eqref{indmdlb10a}, est le morphisme canonique déduit de \eqref{indmdlb7a}.
\end{itemize}
\end{prop}

(i) Pour tout nombre rationnel $0<t\leq r$, l'isomorphisme \eqref{indmdlb11b} induit un isomorphisme de $\bMH(\co_\fX, \txi^{-1}\tOmega^1_{\fX/\cS})$
que l'on note encore
\begin{equation}\label{indmdlb10c}
\gamma^{(t)}\colon \cN\otimes_{\co_\fX}\vupsigma^{(t)}_*(\rI\fS^{(t)}(\bvocB_\mQ)) 
\stackrel{\sim}{\rightarrow} \vupsigma^{(t)}_*(\rI\hupsigma^{(t)*}(\cN)).
\end{equation} 
En effet, pour tout $\co_\fX$-module $\cM$, le morphisme canonique 
$\cN\otimes_{\co_\fX}\cM\rightarrow \cN\otimes_{\co_\fX}\cM_\mQ$ est un isomorphisme de $\bIndMod(\co_\fX)$ \eqref{indsh5g}. 
On désigne par 
\begin{equation}\label{indmdlb10d}
\delta^{(t)}\colon \cN\otimes_{\co_\fX}\vupsigma^{(t)}_*(\rI\fS^{(t)}(\bvocB_\mQ)) 
\stackrel{\sim}{\rightarrow} \vupsigma^{(t)}_*(\rI\fS^{(t)}(\cM))
\end{equation} 
le composé $\vupsigma^{(t)}_*(\alpha^{(t)})\circ \gamma^{(t)}$. Le diagramme 
\begin{equation}\label{indmdlb10e}
\xymatrix{
{\cN\otimes_{\co_\fX}\vupsigma^{(r')}_*(\rI\fS^{(r')}(\bvocB_\mQ))}\ar[r]^-(0.5){\delta^{(r')}}\ar[d]&{\vupsigma^{(r')}_*(\rI\fS^{(r')}(\cM))}\ar[d]\\
{\cN\otimes_{\co_\fX}\vupsigma^{(r'')}_*(\rI\fS^{(r'')}(\bvocB_\mQ))}\ar[r]^-(0.5){\delta^{(r'')}}&{\vupsigma^{(r'')}_*(\rI\fS^{(r'')}(\cM))}}
\end{equation}
où les flèches verticales sont les morphismes canoniques \eqref{indmdlb7a}, 
est commutatif en vertu de \ref{indmdlb11}(ii). Les isomorphismes $(\delta^{(t)})_{0<t\leq r}$ 
induisent par passage à la limite inductive un isomorphisme de $\bMH(\co_\fX, \txi^{-1}\tOmega^1_{\fX/\cS})$
\begin{equation}\label{indmdlb10f}
\delta\colon \cN\otimes_{\co_\fX}\cH(\bvocB_\mQ) \stackrel{\sim}{\rightarrow} \cH(\cM).
\end{equation}
Considérons le diagramme commutatif
\begin{equation}\label{indmdlb10g}
\xymatrix{
{\cN}\ar[r]^-(0.5){\iota^{(r')}}\ar[rd]&
{\cN\otimes_{\co_\fX}\vupsigma^{(r')}_*(\rI\fS^{(r')}(\bvocB_\mQ))}\ar[r]^-(0.5){\delta^{(r')}}
\ar[d]&{\vupsigma^{(r')}_*(\rI\fS^{(r')}(\cM))}\ar[d]\\
&{\cN\otimes_{\co_\fX}\cH(\bvocB_\mQ)}\ar[r]^-(0.5){\delta}&{\cH(\cM)}}
\end{equation}
où $\iota^{(r')}$ est induit par le morphisme canonique $\co_\fX\rightarrow \vupsigma^{(r')}_*(\rI\fS^{(r')}(\bvocB_\mQ))$ 
et les flèches verticales sont les morphismes canoniques. D'après \ref{indmdlb11}(i), on a 
\begin{equation}\label{indmdlb10h}
\delta^{(r')}\circ \iota^{(r')}=\vupsigma^{(r')}_*(\alpha^{(r')})\circ\gamma^{(r')}\circ \iota^{(r')}=\beta^{(r')}.
\end{equation} 
La proposition s'ensuit en vertu de \ref{indmdlb8}. 

(ii) Cela résulte de \eqref{indmdlb10e}, \eqref{indmdlb10g} et \ref{ahttf27}(ii).

\begin{cor}\label{indmdlb12}
Pour tout ind-$\bvocB$-module de Dolbeault $\cM$, 
$\cH(\cM)$ \eqref{indmdlb7d} est un $\co_\fX[\frac 1 p]$-fibré de Higgs soluble, associé à $\cM$.  
En particulier, $\cH$ induit un foncteur que l'on note encore
\begin{equation}\label{indmdlb12a}
\cH\colon 
\begin{array}[t]{clcr}
\bIndMod^\Dolb(\bvocB)&\rightarrow& \bMH^\sol(\co_\fX[\frac 1 p], \txi^{-1}\tOmega^1_{\fX/\cS})\\
\cM&\mapsto& \cH(\cM).
\end{array}
\end{equation}
\end{cor}

\begin{cor}\label{indmdlb13}
Pour tout ind-$\bvocB$-module de Dolbeault $\cM$, il existe un nombre rationnel $r>0$
et un isomorphisme de $\bIndMC(\bvcC^{(r)}/\bvocB)$,
\begin{equation}\label{indmdlb13a}
\alpha\colon \rI\hupsigma^{(r)*}(\cH(\cM))\stackrel{\sim}{\rightarrow} \rI\fS^{(r)}(\cM)
\end{equation}
vérifiant les propriétés suivantes. Pour tout nombre rationnel $r'$ tel que $0< r'\leq r$, on désigne par
\begin{equation}\label{indmdlb13b}
\alpha^{(r')}\colon \rI\hupsigma^{(r')*}(\cH(\cM)) \stackrel{\sim}{\rightarrow}\rI\fS^{(r')}(\cM)
\end{equation}
l'isomorphisme de $\bIndMC(\bvcC^{(r')}/\bvocB)$ induit par $\rI\varepsilon^{r,r'}(\alpha)$ et 
les isomorphismes \eqref{indmdlb3b} et \eqref{indmdlb3c}, et par
\begin{equation}\label{indmdlb13c}
\beta^{(r')}\colon \cH(\cM)\rightarrow\vupsigma^{(r')}_*(\rI\fS^{(r')}(\cM))
\end{equation}
son adjoint \eqref{indmdlb2g}. Alors,
\begin{itemize}
\item[{\rm (i)}] Pour tout nombre rationnel $r'$ tel que $0< r'\leq r$, 
le morphisme $\beta^{(r')}$ est un inverse à droite du morphisme canonique 
$\varpi^{(r')}\colon\vupsigma^{(r')}_*(\rI\fS^{(r')}(\cM))\rightarrow \cH(\cM)$. 
\item[{\rm (ii)}] Pour tous nombres rationnels $r'$ et $r''$ tels que $0<r''<r'\leq r$, le composé
\begin{equation}\label{indmdlb13d}
\xymatrix{
{\vupsigma^{(r')}_*(\rI\fS^{(r')}(\cM))}\ar[r]^-(0.5){\varpi^{(r')}}&
{\cH(\cM)}\ar[r]^-(0.5){\beta^{(r'')}}&{\vupsigma^{(r'')}_*(\rI\fS^{(r'')}(\cM))}}
\end{equation}
est le morphisme canonique. 
\end{itemize}
\end{cor}

\begin{rema}\label{indmdlb16}
Sous les hypothèses de \ref{indmdlb13}, l'isomorphisme $\alpha$ n'est a priori pas uniquement déterminé par $(\cM,r)$, 
mais pour tout nombre rationnel $0<r'<r$, le morphisme $\alpha^{(r')}$ \eqref{indmdlb13b} ne dépend que de $\cM$,
et il en dépend fonctoriellement (cf. la preuve de \ref{indmdlb20}).
\end{rema}

\subsection{}\label{indmdlb14}
Soient $r$ un nombre rationnel $>0$, $(\cM,\cN,\alpha)$ un triplet $r$-admissible \eqref{indmdlb4}. 
Pour éviter toute ambiguïté avec \eqref{indmdlb9a}, notons 
\begin{equation}\label{indmdlb14a}
\calpha\colon \rI\fS^{(r)}(\cM)\rightarrow \rI\hupsigma^{(r)*}(\cN)
\end{equation}
l'inverse de $\alpha$ dans $\bIndMC(\bvcC^{(r)}/\bvocB)$. 
Pour tout nombre rationnel $r'$ tel que $0< r'\leq r$, on désigne par 
\begin{equation}\label{indmdlb14b}
\calpha^{(r')}\colon \rI\fS^{(r')}(\cM) \stackrel{\sim}{\rightarrow}\rI\hupsigma^{(r')*}(\cN)
\end{equation}
l'isomorphisme de $\bIndMC(\bvcC^{(r')}/\bvocB)$ induit par $\rI\varepsilon^{r,r'}(\calpha)$ et les isomorphismes \eqref{indmdlb3b} et \eqref{indmdlb3b}, et par
\begin{equation}\label{indmdlb14c}
\cbeta^{(r')}\colon \cM \rightarrow\rI\cK^{(r')}(\rI\hupsigma^{(r')*}(\cN))
\end{equation}
le morphisme adjoint \eqref{indmdlb1c}. 

\begin{prop}\label{indmdlb15}
Les hypothèses étant celles de \ref{indmdlb14}, soient, de plus, $r'$, $r''$ deux nombres rationnels tels que $0<r''< r'\leq r$. Alors,
\begin{itemize} 
\item[{\rm (i)}] Le morphisme composé 
\begin{equation}\label{indmdlb15a}
\xymatrix{
{\cM}\ar[r]^-(0.5){\cbeta^{(r')}}&{\rI\cK^{(r')}(\rI\hupsigma^{(r')*}(\cN))}\ar[r]&{\cV(\cN)}},
\end{equation}
où la seconde flèche est le morphisme canonique \eqref{indmdlb7f}, est un isomorphisme, indépendant de $r'$. 
\item[{\rm (ii)}] Le morphisme composé 
\begin{equation}\label{indmdlb15b}
\xymatrix{
{\rI\cK^{(r')}(\rI\hupsigma^{(r')*}(\cN))}\ar[r]&{\cV(\cN)}\ar[r]^{\sim}&{\cM}\ar[r]^-(0.5){\cbeta^{(r'')}}&{\rI\cK^{(r'')}(\rI\hupsigma^{(r'')*}(\cN))}},
\end{equation}
où la première flèche est le morphisme canonique et la seconde flèche est l'isomorphisme inverse de \eqref{indmdlb15a},
est le morphisme canonique \eqref{indmdlb7e}.
\end{itemize}
\end{prop}

(i) Comme l'ind-$\bvocB$-module $\cM$ est rationnel et plat d'après \ref{indmdlb6}, pour tout nombre rationnel $t\geq 0$, 
on a des isomorphismes canoniques de $\bIndMod(\bvocB)$
\begin{equation}\label{indmdlb15c}
\cM\otimes_{\bvocB}\rI\cK^{(t)}(\rI\fS^{(t)}(\bvocB_\mQ))
\stackrel{\sim}{\rightarrow}\cM\otimes_{\bvocB}\rI\cK^{(t)}(\rI\fS^{(t)}(\bvocB))
\stackrel{\sim}{\rightarrow}\rI\cK^{(t)}(\rI\fS^{(t)}(\cM)).
\end{equation}
On note $\gamma^{(t)}$ l'isomorphisme composé. 
Pour tout nombre rationnel $0<t\leq r$, on désigne par 
\begin{equation}\label{indmdlb15d}
\delta^{(t)}\colon \cM\otimes_{\bvocB}\rI\cK^{(t)}(\rI\fS^{(t)}(\bvocB_\mQ))
\stackrel{\sim}{\rightarrow}\rI\cK^{(t)}(\rI\hupsigma^{(t)*}(\cN))
\end{equation}
l'isomorphisme $\cK^{(t)}(\calpha^{(t)})\circ \gamma^{(t)}$. Le diagramme 
\begin{equation}\label{indmdlb15e}
\xymatrix{
{\cM\otimes_{\bvocB}\rI\cK^{(r')}(\rI\fS^{(r')}(\bvocB_\mQ))}\ar[d]\ar[r]^-(0.5){\delta^{(r')}}&
{\rI\cK^{(r')}(\rI\hupsigma^{(r')*}(\cN))}\ar[d]\\
{\cM\otimes_{\bvocB}\rI\cK^{(r'')}(\rI\fS^{(r'')}(\bvocB_\mQ))}\ar[r]^-(0.5){\delta^{(r'')}}&
{\rI\cK^{(r'')}(\rI\hupsigma^{(r'')*}(\cN))}}
\end{equation}
où les flèches verticales sont induites par le morphisme \eqref{indmdlb3e} et les isomorphismes 
\eqref{indmdlb3b} et \eqref{indmdlb3c}, est clairement commutatif. 

D'après \ref{ahttf32}, le morphisme canonique 
\begin{equation}\label{indmdlb15f}
\cM\rightarrow \underset{\underset{t\in \mQ_{>0}}{\longrightarrow}}{\lim}\ 
\cM\otimes_{\bvocB} \rI\cK^{(t)}(\rI\fS^{(t)}(\bvocB_\mQ))
\end{equation}
est un isomorphisme. Les isomorphismes $(\delta^{(t)})_{0<t\leq r}$ 
induisent alors par passage à la limite inductive un isomorphisme 
\begin{equation}\label{indmdlb15g}
\delta\colon \cM\stackrel{\sim}{\rightarrow} \cV(\cN).
\end{equation}
Il résulte aussitôt des définitions que le diagramme 
\begin{equation}\label{indmdlb15h}
\xymatrix{
{\cM\otimes_{\bvocB} \rI\cK^{(r')}(\rI\fS^{(r')}(\bvocB_\mQ))}
\ar[r]^-(0.5){\delta^{(r')}}&{\rI\cK^{(r')}(\rI\hupsigma^{(r')*}(\cN))}\ar[d]\\
{\cM}\ar[r]^-(0.5){\delta}\ar[u]^-(0.5){\iota^{(r')}}&{\cV(\cN)}}
\end{equation}
où $\iota^{(r')}$ est induit par le morphisme canonique $\bvocB\rightarrow \rI\cK^{(r')}(\rI\fS^{(r')}(\bvocB_\mQ))$
et la flèche non-libellée est le morphisme canonique, est commutatif. On vérifie aussitôt qu'on a 
\begin{equation}\label{indmdlb15i}
\delta^{(r')}\circ \iota^{(r')}=\cK^{(r')}(\calpha^{(r')})\circ \gamma^{(r')}\circ \iota^{(r')}=\cbeta^{(r')}.
\end{equation}
La proposition s'ensuit. 

(ii) Cela résulte de \eqref{indmdlb15e} et \ref{ahttf32}(ii).

\begin{cor}\label{indmdlb17}
Pour tout $\co_\fX[\frac 1 p]$-fibré de Higgs $\cN$ à coefficients dans $\txi^{-1}\tOmega^1_{\fX/\cS}$, 
le  ind-$\bvocB$-module $\cV(\cN)$ \eqref{indmdlb7f} est de Dolbeault, associé à $\cN$.  
En particulier, $\cV$ induit un foncteur que l'on note encore
\begin{equation}\label{indmdlb17a}
\cV\colon 
\begin{array}[t]{clcr}
\bMH^\sol(\co_\fX[\frac 1 p], \txi^{-1}\tOmega^1_{\fX/\cS})&\rightarrow &\bIndMod^\Dolb(\bvocB)\\
\cN&\mapsto & \cV(\cN).
\end{array}
\end{equation}
\end{cor}

\begin{cor}\label{indmdlb18}
Pour tout $\co_\fX[\frac 1 p]$-fibré de Higgs soluble $\cN$ à coefficients dans $\txi^{-1}\tOmega^1_{\fX/\cS}$,
il existe un nombre rationnel $r>0$ et un isomorphisme de $\bIndMC(\bvcC^{(r)}/\bvocB)$
\begin{equation}\label{indmdlb18a}
\calpha\colon \rI\fS^{(r)}(\cV(\cN))\stackrel{\sim}{\rightarrow}\rI\hupsigma^{(r)*}(\cN)
\end{equation}
vérifiant les propriétés suivantes. Pour tout nombre rationnel $r'$ tel que $0< r'\leq r$, notons
\begin{equation}\label{indmdlb18b}
\calpha^{(r')}\colon\rI\fS^{(r')}(\cV(\cN))\stackrel{\sim}{\rightarrow}\rI\hupsigma^{(r')*}(\cN)
\end{equation}
l'isomorphisme de $\bIndMC(\bvcC^{(r')}/\bvocB)$ induit par $\rI\varepsilon^{r,r'}(\calpha)$ et 
les isomorphismes \eqref{indmdlb3b} et \eqref{indmdlb3c}, et 
\begin{equation}\label{indmdlb18c}
\cbeta^{(r')}\colon \cV(\cN)\rightarrow\rI\cK^{(r')}(\rI\hupsigma^{(r')*}(\cN))
\end{equation}
son adjoint. Alors,
\begin{itemize}
\item[{\rm (i)}] Pour tout nombre rationnel $r'$ tel que $0<r'\leq r$,
le morphisme $\cbeta^{(r')}$ est un inverse à droite du morphisme canonique
$\varpi^{(r')}\colon \rI\cK^{(r')}(\rI\hupsigma^{(r')*}(\cN))\rightarrow \cV(\cN)$.  
\item[{\rm (ii)}] Pour tous nombres rationnels $r'$ et $r''$ tels que $0<r''<r'\leq r$, le composé
\begin{equation}\label{indmdlb18d}
\xymatrix{
{\rI\cK^{(r')}(\hupsigma^{(r')*}(\cN))}\ar[r]^-(0.5){\varpi^{(r')}}&
{\cV(\cN)}\ar[r]^-(0.5){\cbeta^{(r'')}}&{\rI\cK^{(r'')}(\hupsigma^{(r'')*}(\cN))}}
\end{equation}
est le morphisme canonique. 
\end{itemize}
\end{cor}

\begin{rema}\label{indmdlb19}
Sous les hypothèses de \ref{indmdlb18}, l'isomorphisme $\calpha$ n'est a priori pas uniquement déterminé par $(\cN,r)$, 
mais pour tout nombre rationnel $0<r'<r$, 
le morphisme $\calpha^{(r')}$ \eqref{indmdlb18b} ne dépend que de $\cN$ et il en dépend fonctoriellement (cf. la preuve
de \ref{indmdlb20}). 
\end{rema}

\begin{teo}\label{indmdlb20}
Les foncteurs \eqref{indmdlb12a} et \eqref{indmdlb17a} 
\begin{equation}
\xymatrix{
{\bIndMod^\Dolb(\bvocB)}\ar@<1ex>[r]^-(0.5){\cH}&{\bMH^\sol(\co_\fX[\frac 1 p], \txi^{-1}\tOmega^1_{\fX/\cS})}
\ar@<1ex>[l]^-(0.5){\cV}}
\end{equation}
sont des équivalences de catégories quasi-inverses l'une de l'autre. 
\end{teo}

Pour tout objet $\cM$ de $\bIndMod^\Dolb(\bvocB)$,
$\cH(\cM)$ est un $\co_\fX[\frac 1 p]$-fibré de Higgs soluble associé à $\cM$, en vertu de \ref{indmdlb12}. 
Choisissons un nombre rationnel $r_\cM>0$ et un isomorphisme de $\bIndMC(\bvcC^{(r_\cM)}/\bvocB)$ 
\begin{equation}\label{indmdlb20a}
\alpha_\cM\colon \rI\hupsigma^{(r_\cM)*}(\cH(\cM))\stackrel{\sim}{\rightarrow} \rI\fS^{(r_\cM)}(\cM)
\end{equation}
vérifiant les propriétés du \ref{indmdlb13}. Pour tout nombre rationnel $r$ tel que $0< r\leq r_\cM$, on désigne par
\begin{equation}\label{indmdlb20b}
\alpha^{(r)}_\cM\colon \rI\hupsigma^{(r)*}(\cH(\cM)) \stackrel{\sim}{\rightarrow}\rI\fS^{(r)}(\cM)
\end{equation}
l'isomorphisme de $\bIndMC(\bvcC^{(r)}/\bvocB)$ induit par $\rI\varepsilon^{r_\cM,r}(\alpha_\cM)$ et 
les isomorphismes \eqref{indmdlb3b} et \eqref{indmdlb3c}, par 
\begin{eqnarray}
\calpha_\cM\colon \rI\fS^{(r_\cM)}(\cM)&\stackrel{\sim}{\rightarrow}& \rI\hupsigma^{(r_\cM)*}(\cH(\cM)),\label{indmdlb20c}\\
\calpha^{(r)}_\cM\colon \rI\fS^{(r)}(\cM) &\stackrel{\sim}{\rightarrow}&\rI\hupsigma^{(r)*}(\cH(\cM)),\label{indmdlb20d}
\end{eqnarray}
les inverses de $\alpha_\cM$ et $\alpha_\cM^{(r)}$, respectivement, et par
\begin{eqnarray}
\beta^{(r)}_\cM\colon \cH(\cM)&\rightarrow&\vupsigma^{(r)}_*(\rI\fS^{(r)}(\cM)),\label{indmdlb20e}\\
\cbeta^{(r)}_\cM\colon \cM&\rightarrow&\rI\cK^{(r)}(\rI\hupsigma^{(r)*}(\cH(\cM))),\label{indmdlb20f}
\end{eqnarray}
les morphismes adjoints de $\alpha^{(r)}_\cM$ et $\calpha^{(r)}_\cM$, respectivement. On notera que $\calpha^{(r)}_\cM$ est
induit par $\rI\varepsilon^{r_\cM,r}(\calpha_\cM)$ et les isomorphismes \eqref{indmdlb3b} et \eqref{indmdlb3c}.
D'après \ref{indmdlb15}(i), le morphisme composé
\begin{equation}\label{indmdlb20g}
\xymatrix{
{\cM}\ar[r]^-(0.5){\cbeta^{(r)}_\cM}&{\rI\cK^{(r)}(\rI\hupsigma^{(r)*}(\cH(\cM)))}\ar[r]&{\cV(\cH(\cM))}},
\end{equation}
où la seconde flèche est le morphisme canonique, 
est un isomorphisme, qui dépend a priori de $\alpha_\cM$ mais pas de $r$.  Montrons que cet isomorphisme
ne dépend que de $\cM$ (mais pas du choix de $\alpha_\cM$) et qu'il en dépend fonctoriellement. 
Il suffit de montrer que pour tout morphisme $u\colon \cM\rightarrow \cM'$ de $\bIndMod^\Dolb(\bvocB)$
et tout nombre rationnel $0<r<\inf(r_\cM,r_{\cM'})$, le diagramme de $\bIndMC(\bvcC^{(r)}/\bvocB)$
\begin{equation}\label{indmdlb20h}
\xymatrix{
{\rI\hupsigma^{(r)*}(\cH(\cM))}\ar[r]^-(0.5){\alpha^{(r)}_\cM}\ar[d]_{\rI\hupsigma^{(r)*}(\cH(u))}&{\rI\fS^r(\cM)}\ar[d]^{\rI\fS^r(u)}\\
{\rI\hupsigma^{(r)*}(\cH(\cM'))}\ar[r]^-(0.5){\alpha^{(r)}_{\cM'}}&{\rI\fS^r(\cM')}}
\end{equation}
est commutatif. Soient $r$, $r'$ deux nombres rationnels tels que $0<r<r'<\inf(r_\cM,r_{\cM'})$. 
Considérons le diagramme 
\begin{equation}
\xymatrix{
{\vupsigma^{(r')}_*(\rI\fS^{(r')}(\cM))}\ar[r]^-(0.4){\varpi^{(r')}_\cM}\ar[d]_{\vupsigma^{(r')}_*(\rI\fS^{r'}(u))}&
{\cH(\cM)}\ar[r]^-(0.5){\beta^{(r)}_\cM}\ar[d]^{\cH(u)}&{\vupsigma^{(r)}_*(\rI\fS^{(r)}(\cM))}\ar[d]^{\vupsigma^{(r)}_*(\rI\fS^{r}(u))}\\
{\vupsigma^{(r')}_*(\rI\fS^{(r')}(\cM'))}\ar[r]^-(0.4){\varpi^{(r')}_{\cM'}}&{\cH(\cM')}\ar[r]^-(0.5){\beta^{(r)}_{\cM'}}&{\vupsigma^{(r)}_*(\rI\fS^{r}(\cM'))}}
\end{equation}
où $\varpi^{(r')}_\cM$ et $\varpi^{(r')}_{\cM'}$ sont les morphismes canoniques. 
Il résulte de \ref{indmdlb13}(ii) que le grand rectangle est commutatif. 
Comme le carré de gauche est commutatif et que $\varpi^{(r')}_{\cM'}$ est surjectif d'après \ref{indmdlb13}(i),
le carré de droite est aussi commutatif. L'assertion recherchée s'ensuit compte tenu de l'injectivité de \eqref{indmdlb2g}. 

De même, pour tout objet $\cN$ de $\bMH^\sol(\co_\fX[\frac 1 p], \xi^{-1}\tOmega^1_{\fX/\cS})$,
$\cV(\cN)$ est un ind-$\bvocB$-module de Dolbeault associé à $\cN$, en vertu de \ref{indmdlb17}.
Choisissons un nombre rationnel $r_\cN>0$ et un isomorphisme de $\bIndMC(\bvcC^{(r_\cN)}/\bvocB)$ 
\begin{equation}\label{indmdlb20aa}
\calpha_\cN\colon \rI\fS^{(r_\cN)}(\cV(\cN))\stackrel{\sim}{\rightarrow} \rI\hupsigma^{(r_\cN)*}(\cN)
\end{equation}
vérifiant les propriétés de \ref{indmdlb18}. Pour tout nombre rationnel $r$ tel que $0< r\leq r_\cN$, on désigne par
\begin{equation}\label{indmdlb20bb}
\calpha^{(r)}_\cN\colon \rI\fS^{(r)}(\cV(\cN))\stackrel{\sim}{\rightarrow} \rI\hupsigma^{(r)*}(\cN)
\end{equation}
l'isomorphisme de $\bIndMC(\bvcC^{(r)}/\bvocB)$ induit par $\rI\varepsilon^{r_\cN,r}(\calpha_\cN)$ et 
les isomorphismes \eqref{indmdlb3b} et \eqref{indmdlb3c}, par 
\begin{eqnarray}
\alpha_\cN\colon \rI\hupsigma^{(r_\cN)*}(\cN) &\stackrel{\sim}{\rightarrow}&\rI\fS^{(r_\cN)}(\cV(\cN)),\label{indmdlb20cc}\\
\alpha^{(r)}_\cN\colon \rI\hupsigma^{(r)*}(\cN)&\stackrel{\sim}{\rightarrow} &\rI\fS^{(r)}(\cV(\cN)),\label{indmdlb20dd}
\end{eqnarray}
les inverses de $\calpha_\cM$ et $\calpha_\cN^{(r)}$, respectivement, et par
\begin{eqnarray}
\cbeta^{(r)}_\cN\colon \cV(\cN)&\rightarrow&\rI\cK^{(r)}(\rI\hupsigma^{(r)*}(\cN)),\label{indmdlb20ee}\\
\beta^{(r)}_\cN\colon \cN&\rightarrow&\vupsigma^{(r)}_*(\rI\fS^r(\cV(\cN))),\label{indmdlb20ff}
\end{eqnarray}
les morphismes adjoints de $\calpha^{(r)}_\cN$ et $\alpha^{(r)}_\cN$, respectivement. D'après \ref{indmdlb10}(i), 
le morphisme composé 
\begin{equation}
\xymatrix{
{\cN}\ar[r]^-(0.5){\beta^{(r)}_\cN}&{\vupsigma^{(r)}_*(\rI\fS^{(r)}(\cV(\cN)))}\ar[r]&{\cH(\cV(\cN))}},
\end{equation}
où la seconde flèche est le morphisme canonique, est un isomorphisme, qui dépend a priori de $\calpha_\cN$
mais pas de $r$. Montrons que cet isomorphisme
ne dépend que de $\cN$ (mais pas du choix de $\calpha_\cN$) et qu'il en dépend fonctoriellement. 
Il suffit de montrer que pour tout morphisme $v\colon \cN\rightarrow \cN'$ de 
$\bMH^\sol(\co_\fX[\frac 1 p], \xi^{-1}\tOmega^1_{\fX/\cS})$ et tout nombre rationnel 
$0<r<\inf(r_\cN,r_{\cN'})$, le diagramme de $\bIndMC(\bvcC^{(r)}/\bvocB)$
\begin{equation}
\xymatrix{
{\rI\fS^{(r)}(\cV(\cN))}\ar[r]^-(0.5){\calpha^{(r)}_\cN}\ar[d]_{\rI\fS^{(r)}(\cV(v))}&{\rI\hupsigma^{(r)*}(\cN)}\ar[d]^{\rI\hupsigma^{(r)*}(v)}\\
{\rI\fS^{(r)}(\cV(\cN'))}\ar[r]^-(0.5){\calpha^{(r)}_{\cN'}}&{\rI\hupsigma^{(r)*}(\cN')}}
\end{equation}
est commutatif. Soient $r$, $r'$ deux nombres rationnels tels que $0<r<r'<\inf(r_\cN,r_{\cN'})$. 
Considérons le diagramme de $\bIndMod(\bvocB)$
\begin{equation}
\xymatrix{
{\rI\cK^{(r')}(\rI\hupsigma^{(r')*}(\cN))}\ar[r]^-(0.5){\varpi^{(r')}_\cN}\ar[d]_{\rI\cK^{(r')}(\rI\hupsigma^{(r')*}(v))}&
{\cV(\cN)}\ar[r]^-(0.5){\cbeta^{(r)}_\cM}\ar[d]^{\cV(v)}&{\rI\cK^{(r)}(\rI\hupsigma^{(r)*}(\cN))}\ar[d]^{\rI\cK^{(r)}(\rI\hupsigma^{(r)*}(v))}\\
{\rI\cK^{(r')}(\rI\hupsigma^{(r')*}(\cN'))}\ar[r]^-(0.5){\varpi^{(r')}_{\cN'}}&{\cV(\cN')}\ar[r]^-(0.5){\cbeta^{(r)}_{\cN'}}&{\rI\cK^{(r)}(\rI\hupsigma^{(r)*}(\cN'))}}
\end{equation}
où $\varpi^{(r')}_\cN$ et $\varpi^{(r')}_{\cN'}$ sont les morphismes canoniques. 
Il résulte de \ref{indmdlb18}(ii) que le grand rectangle est commutatif. 
Comme le carré de gauche est commutatif et que $\varpi_\cN^{(r')}$ est inversible à droite d'après \ref{indmdlb18}(i),
le carré de droite est aussi commutatif; d'où l'assertion recherchée.

\subsection{}\label{indmdlb25}
Reprenons les hypothèses et notations de \ref{ahttf29}. On rappelle qu'on affecte d'un exposant $'$ les objets associés à la 
$(\cA^{\ast}_2(\oS/S),\cM_{\cA^{\ast}_2(\oS/S)})$-déformation $(\tX',\cM_{\tX'})$ \eqref{ahttf14}. 
On dispose des foncteurs \eqref{indmdlb12a} et \eqref{indmdlb17a}
\begin{equation}\label{indmdlb25a}
\xymatrix{
{\bIndMod^\Dolb(\bvocB)} \ar@<1ex>[r]^-(0.5){\cH}& {\bMH^\sol(\co_\fX[\frac 1 p],\xi^{-1}\tOmega^1_{\fX/\cS})} \ar@<1ex>[l]^-(0.5){\cV}}
\end{equation}
\begin{equation}\label{indmdlb25b}
\xymatrix{
{\bIndMod^{\Dolb'}(\bvocB)} \ar@<1ex>[r]^-(0.5){\cH'}& {\bMH^{\sol'}(\co_\fX[\frac 1 p],(\xi_\pi^{\ast})^{-1}\tOmega^1_{\fX/\cS})} \ar@<1ex>[l]^-(0.5){\cV'}}
\end{equation}
associés aux déformations $(\tX,\cM_{\tX})$ et $(\tX',\cM_{\tX'})$, respectivement. 
L'isomorphisme $\nu$ \eqref{ahttf29d} induit un  isomorphisme $\co_{\fX}$-linéaire
\begin{equation}\label{indmdlb25c}
\uplambda_*(\nu)\colon (\xi_\pi^{\ast})^{-1}\tOmega^1_{\fX/\cS}\stackrel{\sim}{\rightarrow}  \xi^{-1}\tOmega^1_{\fX/\cS}.
\end{equation}
On note 
\begin{equation}\label{indmdlb25d}
\tau \colon \bMH(\co_\fX[\frac 1 p],(\xi_\pi^{\ast})^{-1}\tOmega^1_{\fX/\cS})\rightarrow \bMH(\co_\fX[\frac 1 p],\xi^{-1}\tOmega^1_{\fX/\cS})
\end{equation}
le foncteur induit par $p^\rho \uplambda_*(\nu)$. 

Compte tenu de \eqref{ahttf29h}, pour tout nombre rationnel $r\geq 0$, on a des isomorphismes canoniques 
\begin{eqnarray}
\tau \circ \vupsigma'^{(r)}_*\circ \rI \fS'^{(r)} &\stackrel{\sim}{\rightarrow}& \vupsigma^{(r+\rho)}_*\circ \rI \fS^{(r+\rho)},\label{indmdlb25e}\\
\rI \cK'^{(r)} \circ \rI \hupsigma'^{(r)*}&\stackrel{\sim}{\rightarrow}&  \rI \cK^{(r+\rho)} \circ \rI \hupsigma^{(r+\rho)*} \circ \tau.\label{indmdlb25f}
\end{eqnarray}

\begin{prop}\label{indmdlb26}
Les hypothèses étant celles de \ref{indmdlb25}, soient, de plus, $\cM$ un ind-$\ocB$-module de Dolbeault  relativement à la déformation $(\tX',\cM_{\tX'})$,
$\cN$ un $\co_\fX[\frac 1 p]$-fibré de Higgs à coefficients dans $(\xi_\pi^{\ast})^{-1}\tOmega^1_{\fX/\cS}$ soluble relativement à la déformation $(\tX',\cM_{\tX'})$. 
Alors,
\begin{itemize}
\item[{\rm (i)}] Le ind-$\ocB$-module $\cM$ est de Dolbeault relativement à la déformation $(\tX,\cM_{\tX})$, et les morphismes canoniques \eqref{indmdlb25e} 
induisent un isomorphisme fonctoriel
\begin{equation}\label{indmdlb26a}
\tau(\cH'(\cM))\stackrel{\sim}{\rightarrow} \cH(\cM).
\end{equation}
\item[{\rm (ii)}] Le $\co_\fX[\frac 1 p]$-module de Higgs $\tau(\cN)$ à coefficients dans $\xi^{-1}\tOmega^1_{\fX/\cS}$ est soluble 
relativement à la déformation $(\tX,\cM_{\tX})$ et les morphismes canoniques \eqref{indmdlb25f} induisent un isomorphisme fonctoriel
\begin{equation}\label{indmdlb26b}
\cV'(\cN)\stackrel{\sim}{\rightarrow} \cV(\tau(\cN)).
\end{equation}
\end{itemize}
\end{prop}

(i) D'après \ref{indmdlb13}, il existe un nombre rationnel $r>0$ et un isomorphisme de 
$\bIndMC(\bvcC'^{(r)}/\bvocB)$
\begin{equation}\label{indmdlb26c}
\alpha'\colon \rI\hupsigma'^{(r)*}(\cH'(\cM)) \stackrel{\sim}{\rightarrow}\rI\fS'^{(r)}(\cM).
\end{equation}
D'après \ref{indmdlb16}, quitte à diminuer $r$, cet isomorphisme est canonique et fonctoriel en $\cM$.  
Compte tenu de \eqref{ahttf29h}, $\alpha'$ induit un isomorphisme fonctoriel de $\bIndMC(\bvcC^{(r+\rho)}/\bvocB)$
\begin{equation}\label{indmdlb26d}
\alpha\colon \rI\hupsigma^{(r+\rho)*}(\tau(\cH'(\cM))) \stackrel{\sim}{\rightarrow}\rI\fS^{(r+\rho)}(\cM).
\end{equation}
On en déduit, en vertu de \ref{indmdlb10}(i), un isomorphisme fonctoriel
\begin{equation}\label{indmdlb26e}
\tau(\cH'(\cM))\stackrel{\sim}{\rightarrow} \cH(\cM).
\end{equation}

(ii) Il suffit de calquer la preuve de (i) en remplaçant \ref{indmdlb13} et \ref{indmdlb10}(i) par \ref{indmdlb18} et \ref{indmdlb15}(i).

\section{\texorpdfstring{$\bvocB_\mQ$-modules de Dolbeault}{BQ-modules de Dolbeault}}\label{aspglob}

\subsection{}\label{aspglob6}
Soit $r$ un nombre rationnel $\geq 0$. 
On notera simplement $\bIMC(\bvcC^{(r)}/\bvocB)$ la catégorie des $p^r$-isoconnexions intégrables relativement à l'extension $\bvcC^{(r)}/\bvocB$ (cf. \ref{indsh24}); 
on omet donc l'exposant $p^r$ de la notation  introduite dans \ref{indsh24} considérant qu'il est redondant avec l'exposant $r$ de $\bvcC^{(r)}$. 
C'est une catégorie additive.
On désigne par $\bIMC_\mQ(\bvcC^{(r)}/\bvocB)$ la catégorie des objets de $\bIMC(\bvcC^{(r)}/\bvocB)$ à isogénie près (\cite{agt} III.6.1.1). 
D'après \ref{indsh26}(i) et \ref{ahttf15}(iii), tout objet de $\bIMC(\bvcC^{(r)}/\bvocB)$ est une $\bvocB$-isogénie
de Higgs à coefficients dans $\hupsigma^*(\txi^{-1}\tOmega^1_{\fX/\cS})$ \eqref{indmdlb1}. 
En particulier, on peut associer fonctoriellement à tout objet de $\bIMC_\mQ(\bvcC^{(r)}/\bvocB)$ un complexe de Dolbeault dans 
$\bMod_{\mQ}(\bvocB)$ (cf. \ref{indsh23}).

Considérons le foncteur 
\begin{equation}\label{aspglob6a}
\begin{array}[t]{clcr}
\fS^{(r)}\colon\bMod(\bvocB)&\rightarrow &\bIMC(\bvcC^{(r)}/\bvocB)\\ 
\cM&\mapsto& (\bvcC^{(r)}\otimes_{\bvocB}\cM,\bvcC^{(r)}\otimes_{\bvocB}\cM,
\id,p^rd_{\bvcC^{(r)}}\otimes \id)
\end{array}
\end{equation}
et notons encore
\begin{equation}\label{aspglob6b}
\fS^{(r)}\colon \bMod_\mQ(\bvocB)\rightarrow \bIMC_\mQ(\bvcC^{(r)}/\bvocB)
\end{equation}
le foncteur induit. Considérons par ailleurs le foncteur 
\begin{equation}\label{aspglob6c}
\cK^{(r)}\colon 
\begin{array}[t]{clcr}
\bIMC(\bvcC^{(r)}/\bvocB)&\rightarrow &\bMod(\bvocB)\\ 
(\cF,\cG,u,\nabla)&\mapsto& \ker(\nabla)
\end{array}
\end{equation}
et notons encore
\begin{equation}\label{aspglob6d}
\cK^{(r)}\colon \bIMC_\mQ(\bvcC^{(r)}/\bvocB)\rightarrow \bMod_\mQ(\bvocB)
\end{equation}
le foncteur induit. 
Il est clair que le foncteur \eqref{aspglob6a} est un adjoint à gauche du foncteur \eqref{aspglob6c}. 
Par suite, le foncteur \eqref{aspglob6b} est un adjoint à gauche du foncteur \eqref{aspglob6d}.

D'après \ref{indsh26}(ii), si $(\cN,\cN',v,\theta)$ est une $\co_\fX$-isogénie de Higgs
à coefficients dans $\txi^{-1}\tOmega^1_{\fX/\cS}$ \eqref{indsh23},
\begin{equation}\label{aspglob6e}
(\bvcC^{(r)}\otimes_{\bvocB}\hupsigma^*(\cN),\bvcC^{(r)}\otimes_{\bvocB}\hupsigma^*(\cN'),\id \otimes_{\bvocB}\hupsigma^*(v),
p^rd_{\bvcC^{(r)}} \otimes\hupsigma^*(v)+\id \otimes \hupsigma^*(\theta))
\end{equation}
est un objet de $\bIMC(\bvcC^{(r)}/\bvocB)$. On obtient ainsi un foncteur \eqref{definf19}
\begin{equation}\label{aspglob6f}
\hupsigma^{(r)*}\colon \bIH(\co_\fX,\txi^{-1}\tOmega^1_{\fX/\cS})\rightarrow \bIMC(\bvcC^{(r)}/\bvocB).
\end{equation}
Compte tenu de \eqref{indmdlb2e}, celui-ci induit un foncteur que l'on note encore
\begin{equation}\label{aspglob6g}
\hupsigma^{(r)*}\colon \bMH^\coh(\co_\fX[\frac 1 p], \txi^{-1}\tOmega^1_{\fX/\cS})\rightarrow \bIMC_\mQ(\bvcC^{(r)}/\bvocB).
\end{equation}

Soit $(\cM,\cM',u,\nabla)$ un objet de $\bIMC(\bvcC^{(r)}/\bvocB)$.
Compte tenu de (\cite{agt} III.12.4(i)), $\nabla$ induit un morphisme $\co_\fX$-linéaire
\begin{equation}\label{aspglob6h}
\hupsigma_*(\nabla)\colon \hupsigma_*(\cM)\rightarrow \txi^{-1}\tOmega^1_{\fX/\cS}\otimes_{\co_\fX}\hupsigma_*(\cM').
\end{equation}
On vérifie facilement que $(\hupsigma_*(\cM),\hupsigma_*(\cM'),\hupsigma_*(u),\hupsigma_*(\nabla))$ est une $\co_\fX$-isogénie 
de Higgs à coefficients dans $\txi^{-1}\tOmega^1_{\fX/\cS}$ (cf. \cite{agt} III.12.3(i)). On obtient ainsi un foncteur
\begin{equation}\label{aspglob6i}
\hupsigma^{(r)}_*\colon \bIMC(\bvcC^{(r)}/\bvocB)\rightarrow \bIH(\co_\fX,\txi^{-1}\tOmega^1_{\fX/\cS}).
\end{equation}
Le composé des foncteurs \eqref{aspglob6i} et \eqref{definf19b} induit un foncteur  que l'on note encore
\begin{equation}\label{aspglob6j}
\hupsigma^{(r)}_*\colon \bIMC_\mQ(\bvcC^{(r)}/\bvocB)\rightarrow \bMH(\co_\fX[\frac 1 p], \txi^{-1}\tOmega^1_{\fX/\cS}).
\end{equation}

Il est clair que  le foncteur \eqref{aspglob6f} est un adjoint à gauche du foncteur \eqref{aspglob6i}. 
On en déduit que pour tous $\cN\in \ob(\bMH^\coh(\co_\fX[\frac 1 p], \txi^{-1}\tOmega^1_{\fX/\cS}))$ 
et $\cA\in \ob(\bIMC_\mQ(\bvcC^{(r)}/\bvocB))$, on a un homomorphisme canonique bifonctoriel 
\begin{equation}\label{aspglob6k}
\Hom_{\bIMC_\mQ(\bvcC^{(r)}/\bvocB)}(\hupsigma^{(r)*}(\cN),\cA)\rightarrow
\Hom_{\bMH(\co_\fX[\frac 1 p], \txi^{-1}\tOmega^1_{\fX/\cS})}(\cN,\hupsigma^{(r)}_*(\cA)),
\end{equation}
qui est injectif d'après (\cite{agt} III.6.20 et III.6.21). 
On appelle abusivement {\em adjoint} d'un morphisme $\hupsigma^{(r)*}(\cN)\rightarrow \cA$ de $\bIMC_\mQ(\bvcC^{(r)}/\bvocB)$,
son image par l'homomorphisme \eqref{aspglob6k}.

\subsection{}\label{aspglob7}
Soient $r$, $r'$ deux nombres rationnels tels que $r\geq r'\geq 0$, $(\cF,\cG,u,\nabla)$ une $p^r$-isoconnexion 
intégrable relativement à l'extension $\bvcC^{(r)}/\bvocB$. 
D'après \eqref{ahttf14i}, il existe un et un unique morphisme $\bvocB$-linéaire
\begin{equation}
\nabla'\colon \bvcC^{(r')}\otimes_{\bvcC^{(r)}}\cF\rightarrow \hupsigma^*(\xi^{-1}\tOmega^1_{\fX/\cS}) \otimes_{\bvocB}
\bvcC^{(r')}\otimes_{\bvcC^{(r)}} \cG
\end{equation}
tel que pour toutes sections locales $x'$ de $\bvcC^{(r')}$ et $s$ de $\cF$, on ait 
\begin{equation}
\nabla'(x'\otimes_{\bvcC^{(r)}} s)=p^{r'} d_{\bvcC^{(r')}}(x')\otimes_{\bvcC^{(r)}}u(s) + x'\otimes_{\bvcC^{(r)}}\nabla(s).
\end{equation}
Le quadruplet $(\bvcC^{(r')}\otimes_{\bvcC^{(r)}}\cF,\bvcC^{(r')}\otimes_{\bvcC^{(r)}}\cG,
\id\otimes_{\bvcC^{(r)}}u,\nabla')$ est une $p^{r'}$-isoconnexion intégrable relativement à l'extension 
$\bvcC^{(r')}/\bvocB$. On obtient ainsi un foncteur 
\begin{equation}\label{aspglob7a}
\varepsilon^{r,r'}\colon \bIMC(\bvcC^{(r)}/\bvocB)\rightarrow \bIMC(\bvcC^{(r')}/\bvocB).
\end{equation}
Celui-ci induit un foncteur que l'on note encore 
\begin{equation}\label{aspglob7b}
\varepsilon^{r,r'}\colon \bIMC_\mQ(\bvcC^{(r)}/\bvocB)\rightarrow \bIMC_\mQ(\bvcC^{(r')}/\bvocB).
\end{equation}
On a un isomorphisme canonique de foncteurs de $\bMod(\bvocB)$ dans $\bIMC(\bvcC^{(r')}/\bvocB)$ 
(resp.  de $\bMod_\mQ(\bvocB)$ dans $\bIMC_\mQ(\bvcC^{(r')}/\bvocB)$) 
\begin{equation}\label{aspglob7c}
\varepsilon^{r,r'}\circ \fS^{(r)}\stackrel{\sim}{\longrightarrow} \fS^{(r')},
\end{equation}
et un isomorphisme canonique de foncteurs de $\bIH(\co_\fX,\xi^{-1}\tOmega^1_{\fX/\cS})$ dans $\bIMC(\bvcC^{(r')}/\bvocB)$ 
(resp.  de $\bMH^\coh(\co_\fX[\frac 1 p], \txi^{-1}\tOmega^1_{\fX/\cS})$ dans $\bIMC_\mQ(\bvcC^{(r')}/\bvocB)$) 
\begin{equation}\label{aspglob7d}
\varepsilon^{r,r'}\circ \hupsigma^{(r)*}\stackrel{\sim}{\longrightarrow} \hupsigma^{(r')*}.
\end{equation}
Par ailleurs, on a un morphisme canonique de foncteurs de $\bIMC(\bvcC^{(r)}/\bvocB)$ dans $\bMod(\bvocB)$ (resp. 
de $\bIMC_\mQ(\bvcC^{(r)}/\bvocB)$ dans $\bMod_\mQ(\bvocB)$)
\begin{equation}\label{aspglob7e}
\cK^{(r)}\rightarrow \cK^{(r')}\circ \varepsilon^{r,r'},
\end{equation}
et un morphisme canonique de foncteurs de $\bIMC(\bvcC^{(r)}/\bvocB)$ dans $\bIH(\co_\fX,\txi^{-1}\tOmega^1_{\fX/\cS})$ 
(resp.  de $\bIMC_\mQ(\bvcC^{(r)}/\bvocB)$ dans $\bMH(\co_\fX[\frac 1 p], \txi^{-1}\tOmega^1_{\fX/\cS})$) 
\begin{equation}\label{aspglob7f}
\hupsigma^{(r)}_*\longrightarrow \hupsigma^{(r')}_*\circ \varepsilon^{r,r'}.
\end{equation}

Pour tout nombre rationnel $r''$ tel $r'\geq r''\geq 0$, on a un isomorphisme canonique de foncteurs 
de $\bIMC(\bvcC^{(r)}/\bvocB)$ dans $\bIMC(\bvcC^{(r'')}/\bvocB)$ (resp. de $\bIMC_\mQ(\bvcC^{(r)}/\bvocB)$ dans $\bIMC_\mQ(\bvcC^{(r'')}/\bvocB)$)
\begin{equation}\label{aspglob7g}
\varepsilon^{r',r''}\circ \varepsilon^{r,r'}\stackrel{\sim}{\rightarrow}\varepsilon^{r,r''}.
\end{equation}

\subsection{}\label{dolbff1}
Soit $r$ un nombre rationnel $\geq 0$. 
On désigne par $\bMod_{\mQ}(\bvcC^{(r)})$ la catégorie des $\bvcC^{(r)}$-modules à isogénie près \eqref{indsh20} et 
par  $\bIndMod(\bvcC^{(r)})$ la catégorie des ind-$\bvcC^{(r)}$-modules \eqref{indsh15}.
D'après \ref{indsh11}, on a un foncteur pleinement fidèle canonique \eqref{indsh11c}
\begin{equation}\label{dolbff1a}
\upalpha_{\bvcC^{(r)}}\colon \bMod_\mQ(\bvcC^{(r)})\rightarrow \bIndMod(\bvcC^{(r)}),
\end{equation}
compatible avec le foncteur $\upalpha_{\bvocB}$ \eqref{ahttf40d}, via les foncteurs d'oubli \eqref{indsh15f}. 
On vérifie aussitôt, en tenant compte de \eqref{indsh20d}, que pour tout  objet $(\cM,\cM',u,\nabla)$ de $\bIMC(\bvcC^{(r)}/\bvocB)$, le morphisme de $\bIndMod(\bvocB)$
\begin{equation}\label{dolbff1b}
(\id \otimes \upalpha_{\bvocB}(u_\mQ)^{-1})\circ\upalpha_{\bvocB}(\nabla_\mQ)\colon \upalpha_{\bvcC^{(r)}}(\cM_\mQ)\rightarrow 
\hupsigma^*(\txi^{-1}\tOmega^1_{\fX/\cS})\otimes_{\bvocB}\upalpha_{\bvcC^{(r)}}(\cM_\mQ)
\end{equation}
est une $p^r$-connexion intégrable sur $\upalpha_{\bvcC^{(r)}}(\cM_\mQ)$  
relativement à l'extension $\bvcC^{(r)}/\bvocB$ \eqref{indmdlb1}.
On définit ainsi un foncteur 
\begin{equation}\label{dolbff1c}
\upalpha_{\bvcC^{(r)}/\bvocB}\colon \bIMC_\mQ(\bvcC^{(r)}/\bvocB)\rightarrow \bIndMC(\bvcC^{(r)}/\bvocB).
\end{equation}
Celui-ci est pleinement fidèle d'après \eqref{indsh5d} et \eqref{indsh50c}. 

Le diagramme 
\begin{equation}\label{dolbff1d}
\xymatrix{
{\bMod_\mQ(\bvocB)}\ar[r]^-(0.5){\fS^{(r)}}\ar[d]_{\upalpha_{\bvocB}}&{\bIMC_\mQ(\bvcC^{(r)}/\bvocB)}\ar[d]^{\upalpha_{\bvcC^{(r)}/\bvocB}}\\
{\bIndMod(\bvocB)}\ar[r]^-(0.5){\rI\fS^{(r)}}&{\bIndMC(\bvcC^{(r)}/\bvocB)}}
\end{equation}
où $\fS^{(r)}$ (resp. $\rI\fS^{(r)}$) est le foncteur \eqref{aspglob6b} (resp. \eqref{indmdlb1b}), est commutatif à isomorphisme canonique près. 
De même, le diagramme 
\begin{equation}\label{dolbff1e}
\xymatrix{
{\bMH^\coh(\co_\fX[\frac 1 p], \txi^{-1}\tOmega^1_{\fX/\cS})}\ar[r]^-(0.5){\hupsigma^{(r)*}}\ar[d]_{u}&
{\bIMC_\mQ(\bvcC^{(r)}/\bvocB)}\ar[d]^{\upalpha_{\bvcC^{(r)}/\bvocB}}\\
{\bIndMH(\co_\fX,\txi^{-1}\tOmega^1_{\fX/\cS})}\ar[r]^-(0.5){\rI\hupsigma^{(r)*}}&{\bIndMC(\bvcC^{(r)}/\bvocB)}}
\end{equation}
où $\hupsigma^{(r)*}$ (resp. $\rI\hupsigma^{(r)*}$) est le foncteur \eqref{aspglob6g} (resp. \eqref{indmdlb1f}) et $u$ est le foncteur \eqref{indmdlb2f}, 
est commutatif à isomorphisme canonique près.

Comme le foncteur $\upalpha_{\bvocB}$ est exact d'après \ref{indsh16}(ii), le diagramme 
\begin{equation}\label{dolbff1f}
\xymatrix{
{\bIMC_\mQ(\bvcC^{(r)}/\bvocB)}\ar[d]_{\upalpha_{\bvcC^{(r)}/\bvocB}}\ar[r]^-(0.5){\cK^{(r)}}&{\bMod_\mQ(\bvocB)}\ar[d]^{\upalpha_{\bvocB}}\\
{\bIndMC(\bvcC^{(r)}/\bvocB)}\ar[r]^-(0.5){\rI\cK^{(r)}}&{\bIndMod(\bvocB)}}
\end{equation}
où $\cK^{(r)}$ (resp. $\rI\cK^{(r)}$) est le foncteur \eqref{aspglob6d} (resp. \eqref{indmdlb1c}), est commutatif à isomorphisme canonique près. 

D'après \ref{indsh13}(ii), le diagramme 
\begin{equation}\label{dolbff1g}
\xymatrix{
{\bIMC_\mQ(\bvcC^{(r)}/\bvocB)}\ar[d]_{\upalpha_{\bvcC^{(r)}/\bvocB}}\ar[r]^-(0.5){(\hupsigma^{(r)}_*)_\mQ}&
{\bIH_\mQ(\co_\fX, \txi^{-1}\tOmega^1_{\fX/\cS})}\ar[d]^{v}\\
{\bIndMC(\bvcC^{(r)}/\bvocB)}\ar[r]^-(0.5){\rI\hupsigma^{(r)}_*}&{\bIndMH(\co_\fX, \txi^{-1}\tOmega^1_{\fX/\cS})}}
\end{equation}
où $\hupsigma^{(r)}_*$ (resp. $\rI\hupsigma^{(r)}_*$) est le foncteur \eqref{aspglob6i} (resp. \eqref{indmdlb1g}) et $v$ est le foncteur \eqref{indmdlb2d}, 
est commutatif à isomorphisme canonique près. On en déduit que le diagramme 
\begin{equation}\label{dolbff1h}
\xymatrix{
{\bIMC_\mQ(\bvcC^{(r)}/\bvocB)}\ar[d]_{\upalpha_{\bvcC^{(r)}/\bvocB}}\ar[r]^-(0.5){\hupsigma^{(r)}_*}&
{\bMH(\co_\fX[\frac 1 p], \txi^{-1}\tOmega^1_{\fX/\cS})}\ar[d]^{w}\\
{\bIndMC(\bvcC^{(r)}/\bvocB)}\ar[r]^-(0.5){\vupsigma^{(r)}_*}&{\bMH(\co_\fX, \txi^{-1}\tOmega^1_{\fX/\cS})}}
\end{equation}
où $\hupsigma^{(r)}_*$ (resp. $\vupsigma^{(r)}_*$) est le foncteur \eqref{aspglob6j} (resp. \eqref{indmdlb1i}) et $w$ est le foncteur canonique, 
est commutatif à isomorphisme canonique près. 

Pour tous nombres rationnels $r\geq r'\geq 0$, le diagramme 
\begin{equation}\label{dolbff1i}
\xymatrix{
{\bIMC_\mQ(\bvcC^{(r)}/\bvocB)}\ar[d]_{\upalpha_{\bvcC^{(r)}/\bvocB}}\ar[r]^-(0.5){\varepsilon^{r,r'}}&
{\bIMC_\mQ(\bvcC^{(r')}/\bvocB)}\ar[d]^{\upalpha_{\bvcC^{(r')}/\bvocB}}\\
{\bIndMC(\bvcC^{(r)}/\bvocB)}\ar[r]^-(0.5){\rI\varepsilon^{r,r'}}&{\bIndMC(\bvcC^{(r')}/\bvocB)}}
\end{equation}
où $\varepsilon^{r,r'}$ (resp. $\rI\varepsilon^{r,r'}$) est le foncteur \eqref{aspglob7b} (resp. \eqref{indmdlb3a}), 
est commutatif à isomorphisme canonique près. Les isomorphismes \eqref{aspglob7c} et \eqref{indmdlb3b} (resp. \eqref{aspglob7d} et \eqref{indmdlb3c}) 
sont compatibles. Les morphismes \eqref{aspglob7e} et \eqref{indmdlb3e} (resp. \eqref{aspglob7f} et \eqref{indmdlb3f}) sont compatibles.

\begin{defi}\label{aspglob9}
Soient $\cM$ un $\bvocB_\mQ$-module, 
$\cN$ un $\co_\fX[\frac 1 p]$-fibré de Higgs à coefficients dans $\txi^{-1}\tOmega^1_{\fX/\cS}$ \eqref{definf20}.
\begin{itemize}
\item[(i)] Soit $r$ un nombre rationnel $>0$.  On dit que $\cM$ et  $\cN$ sont {\em $r$-associés} 
s'il existe un isomorphisme de $\bIMC_\mQ(\bvcC^{(r)}/\bvocB)$ 
\begin{equation}\label{aspglob9a}
\alpha\colon \hupsigma^{(r)*}(\cN) \stackrel{\sim}{\rightarrow}\fS^{(r)}(\cM).
\end{equation}
On dit alors aussi que le triplet $(\cM,\cN,\alpha)$ est {\em $r$-admissible}. 
\item[(ii)] On dit que $\cM$ et  $\cN$ sont {\em associés} s'il existe un nombre rationnel $r>0$ tel que 
$\cM$ et  $\cN$ soient $r$-associés.
\end{itemize}
\end{defi}

On notera que pour tous nombres rationnels $r\geq r'>0$, 
si $\cM$ et  $\cN$ sont $r$-associés, ils sont $r'$-associés, compte tenu de \eqref{aspglob7c} et \eqref{aspglob7d}.

\begin{lem}\label{dolbff2}
Soient $\cM$ un $\bvocB_\mQ$-module,  
$\cN$ un $\co_\fX[\frac 1 p]$-fibré de Higgs à coefficients dans $\txi^{-1}\tOmega^1_{\fX/\cS}$,
$r$ un nombre rationnel $>0$. 
Alors, pour que $\cM$ et $\cN$ soient $r$-associés dans le sens de \ref{aspglob9}, 
il faut et il suffit que $\upalpha_{\bvocB}(\cM)$ \eqref{ahttf40d} et $\cN$ soient $r$-associés dans le sens de \ref{indmdlb4}.
\end{lem}

Cela résulte aussitôt de \eqref{dolbff1d}, \eqref{dolbff1e} et de la pleine fidélité du foncteur $\upalpha_{\bvcC^{(r)}/\bvocB}$ \eqref{dolbff1c}.

\begin{defi}\label{aspglob1}\
\begin{itemize}
\item[(i)] On dit qu'un $\bvocB_\mQ$-module est {\em de Dolbeault} s'il est associé à  
un $\co_\fX[\frac 1 p]$-fibré de Higgs à coefficients dans $\txi^{-1}\tOmega^1_{\fX/\cS}$ \eqref{aspglob9}. 
\item[(ii)] On dit qu'un $\bvocB_\mQ$-module est {\em fortement de Dolbeault} s'il est  de Dolbeault et adique de type fini \eqref{ahttf60}. 
\item[(iii)] On dit qu'un $\co_\fX[\frac 1 p]$-fibré de Higgs à coefficients dans $\txi^{-1}\tOmega^1_{\fX/\cS}$ \eqref{definf20} 
est {\em rationnellement soluble} (resp. {\em fortement soluble})
s'il est associé à un $\bvocB_\mQ$-module (resp. adique de type fini).
\end{itemize}
\end{defi}

Ces notions dépendent a priori de la déformation $(\tX,\cM_\tX)$ fixée dans \ref{defing12}.
On désigne par $\bMod_\mQ^\Dolb(\bvocB)$ (resp. $\bMod_\mQ^\fDolb(\bvocB)$) la sous-catégorie pleine de $\bMod_\mQ(\bvocB)$  
formée des $\bvocB_\mQ$-modules de Dolbeault (resp. fortement de Dolbeault), et par $\bMH^\qsol(\co_\fX[\frac 1 p], \txi^{-1}\tOmega^1_{\fX/\cS})$ 
(resp. $\bMH^\fsol(\co_\fX[\frac 1 p], \txi^{-1}\tOmega^1_{\fX/\cS})$)
la sous-catégorie pleine de $\bMH(\co_\fX[\frac 1 p], \txi^{-1}\tOmega^1_{\fX/\cS})$
formée des $\co_\fX[\frac 1 p]$-fibrés de Higgs rationnellement (resp. fortement) solubles à coefficients dans $\txi^{-1}\tOmega^1_{\fX/\cS}$. 

\begin{rema}\label{aspglob10}
D'après \ref{dolbff2}, pour qu'un $\bvocB_\mQ$-module $\cM$ soit de Dolbeault, il faut et il suffit que $\upalpha_{\bvocB}(\cM)$
soit de Dolbeault \eqref{indmdlb5}. 
\end{rema}

\begin{rema}
On prendra garde que les $\bvocB_\mQ$-modules fortement de Dolbeault
(resp. les $\co_\fX[\frac 1 p]$-fibrés de Higgs fortement solubles) avaient été baptisés {\em $\bvocB_\mQ$-module de Dolbeault}
(resp. {\em $\co_\fX[\frac 1 p]$-fibrés de Higgs solubles}) dans (\cite{agt} III.12.11). 
\end{rema}

\subsection{}\label{aspglob11}
Pour tout $\bvocB_\mQ$-module $\cM$ et tous nombres rationnels $r\geq r'\geq 0$, 
le morphisme \eqref{aspglob7f} et l'isomorphisme \eqref{aspglob7c} induisent un morphisme de 
$\bMH(\co_\fX[\frac 1 p], \txi^{-1}\tOmega^1_{\fX/\cS})$
\begin{equation}\label{aspglob11a}
\hupsigma^{(r)}_*(\fS^{(r)}(\cM))\rightarrow \hupsigma^{(r')}_*(\fS^{(r')}(\cM)).
\end{equation}
On obtient ainsi un système inductif filtrant $(\hupsigma^r_+(\fS^{(r)}(\cM)))_{r\in \mQ_{\geq 0}}$. 
On désigne par $\cH_\mQ$ le foncteur 
\begin{equation}\label{aspglob11b}
\cH_\mQ\colon \bMod_\mQ(\bvocB)\rightarrow \bMH(\co_\fX[\frac 1 p], \txi^{-1}\tOmega^1_{\fX/\cS}), \ \ \ \cM\mapsto 
\underset{\underset{r\in \mQ_{>0}}{\longrightarrow}}{\lim}\ \hupsigma^{(r)}_*(\fS^{(r)}(\cM)). 
\end{equation}

Compte tenu de \eqref{dolbff1d} et \eqref{dolbff1h}, le diagramme 
\begin{equation}\label{aspglob11d}
\xymatrix{
{\bMod_\mQ(\bvocB)}\ar[r]^-(0.5){\cH_\mQ}\ar[d]_{\upalpha_{\bvocB}}&{\bMH(\co_\fX[\frac 1 p], \txi^{-1}\tOmega^1_{\fX/\cS})}\ar[d]^w\\
{\bIndMod(\bvocB)}\ar[r]^-(0.5){\cH}&{\bMH(\co_\fX, \txi^{-1}\tOmega^1_{\fX/\cS})}}
\end{equation}
où $\cH$ est le foncteur \eqref{indmdlb7c} et $w$ est le foncteur canonique, est commutatif à isomorphisme canonique près. 

Pour tout objet $\cN$ de $\bMH(\co_\fX[\frac 1 p], \txi^{-1}\tOmega^1_{\fX/\cS})$ 
et tous nombres rationnels $r\geq r'\geq 0$, 
le morphisme \eqref{aspglob7e} et l'isomorphisme \eqref{aspglob7d} induisent un morphisme de $\bMod_{\mQ}(\bvocB)$
\begin{equation}\label{aspglob11c}
\cK^{(r)}(\hupsigma^{(r)*}(\cN))\rightarrow \cK^{(r')}(\hupsigma^{(r')*}(\cN)).
\end{equation}
On obtient ainsi un système inductif filtrant $(\cK^{(r)}(\hupsigma^{(r)*}(\cN)))_{r\geq  0}$. 
On rappelle \eqref{ahttf34} que les limites inductives filtrantes ne sont pas nécessairement représentables dans la catégorie 
$\bMod_\mQ(\bvocB)$.

Dans le reste de cette section, on établit pour les $\bvocB_\mQ$-modules des énoncés analogues à ceux établis dans \ref{indmdlb} pour les ind-$\bvocB$-modules,
suivant le point de vue adopté dans (\cite{agt} § III.12).

\begin{lem}\label{aspglob12}
Le $\bvocB_\mQ$-module $\bvocB_\mQ$ est de Dolbeault et 
on a un isomorphisme canonique de $\co_\fX[\frac 1 p]$-modules de Higgs à coefficients dans $\txi^{-1}\tOmega^1_{\fX/\cS}$
\begin{equation}\label{aspglob12a}
(\co_\fX[\frac 1 p],0)\stackrel{\sim}{\rightarrow}\cH_\mQ(\bvocB_\mQ).
\end{equation}
\end{lem}
En effet, la première assertion résulte aussitôt des définitions et la seconde de \ref{ahttf28}, ou ce qui revient au même de \ref{indmdlb8} et \eqref{aspglob11d}.

\subsection{}\label{aspglob13}
Soient $r$ un nombre rationnel $>0$, $\cM$ un $\bvocB_\mQ$-module,  
$\cN$ un $\co_\fX[\frac 1 p]$-fibré de Higgs à coefficients dans $\txi^{-1}\tOmega^1_{\fX/\cS}$ tels que 
le triplet $(\cM,\cN,\alpha)$ soit $r$-admissible \eqref{aspglob9}. 
Pour tout nombre rationnel $r'$ tel que $0< r'\leq r$, on désigne par 
\begin{equation}\label{aspglob13a}
\alpha^{(r')}\colon \hupsigma^{(r')*}(\cN) \stackrel{\sim}{\rightarrow}\fS^{(r')}(\cM)
\end{equation}
l'isomorphisme de $\bIMC_\mQ(\bvcC^{(r')}/\bvocB)$ induit par $\varepsilon^{r,r'}(\alpha)$ et les isomorphismes 
\eqref{aspglob7c} et \eqref{aspglob7d}, et par 
\begin{equation}\label{aspglob13b}
\beta^{(r')}\colon \cN \rightarrow\hupsigma^{(r')}_*(\fS^{(r')}(\cM))
\end{equation}
son adjoint \eqref{aspglob6k}.

\begin{prop}\label{aspglob14}
Les hypothèses étant celles de \ref{aspglob13}, soient, de plus, $r'$, $r''$ deux nombres rationnels 
tels que $0<r''< r'\leq r$. Alors,
\begin{itemize}
\item[{\rm (i)}] Le morphisme composé 
\begin{equation}\label{aspglob14a}
\cN\stackrel{\beta^{(r')}}{\longrightarrow} \hupsigma^{(r')}_*(\fS^{(r')}(\cM))\longrightarrow \cH_\mQ(\cM),
\end{equation}
où la seconde flèche est le morphisme canonique \eqref{aspglob11b}, est un isomorphisme, indépendant de $r'$. 
\item[{\rm (ii)}] Le morphisme composé 
\begin{equation}\label{aspglob14b}
\hupsigma^{(r')}_*(\fS^{(r')}(\cM))\longrightarrow \cH_\mQ(\cM)\stackrel{\sim}{\longrightarrow}\cN\stackrel{\beta^{(r'')}}{\longrightarrow} 
\hupsigma^{(r'')}_*(\fS^{(r'')}(\cM))
\end{equation}
où la première flèche est le morphisme canonique \eqref{aspglob11b} et la deuxième flèche est l'isomorphisme inverse de 
\eqref{aspglob14a}, est le morphisme canonique \eqref{aspglob11a}.
\end{itemize}
\end{prop}

Il suffit de calquer la preuve de (\cite{agt} III.12.17) en tenant compte de \ref{ahttf27} et \ref{aspglob12}.

\begin{cor}\label{aspglob15}
Pour tout $\bvocB_\mQ$-module de Dolbeault $\cM$, 
$\cH_\mQ(\cM)$ \eqref{aspglob11b} est un $\co_\fX[\frac 1 p]$-fibré de Higgs rationnellement soluble, associé à $\cM$.  
En particulier, $\cH_\mQ$ induit un foncteur que l'on note encore
\begin{equation}\label{aspglob15a}
\cH_\mQ\colon \bMod^\Dolb_\mQ(\bvocB)\rightarrow \bMH^\qsol(\co_\fX[\frac 1 p], \txi^{-1}\tOmega^1_{\fX/\cS}), 
\ \ \ \cM\mapsto \cH_\mQ(\cM).
\end{equation}
\end{cor}

Il résulte aussitôt des définitions que le foncteur \eqref{aspglob15a} induit un foncteur que l'on note encore
\begin{equation}\label{aspglob15b}
\cH_\mQ\colon \bMod^\fDolb_\mQ(\bvocB)\rightarrow \bMH^\fsol(\co_\fX[\frac 1 p], \txi^{-1}\tOmega^1_{\fX/\cS}).
\end{equation}

\begin{cor}\label{aspglob16}
Pour tout $\bvocB_\mQ$-module de Dolbeault $\cM$, il existe un nombre rationnel $r>0$
et un isomorphisme de $\bIMC_\mQ(\bvcC^{(r)}/\bvocB)$
\begin{equation}\label{aspglob16a}
\alpha\colon \hupsigma^{(r)*}(\cH_\mQ(\cM))\stackrel{\sim}{\rightarrow} \fS^{(r)}(\cM)
\end{equation}
vérifiant les propriétés suivantes. Pour tout nombre rationnel $r'$ tel que $0< r'\leq r$, notons
\begin{equation}\label{aspglob16b}
\alpha^{(r')}\colon \hupsigma^{(r')*}(\cH_\mQ(\cM)) \stackrel{\sim}{\rightarrow}\fS^{(r')}(\cM)
\end{equation}
l'isomorphisme de $\bIMC_\mQ(\bvcC^{(r')}/\bvocB)$ induit par $\varepsilon^{r,r'}(\alpha)$ et 
les isomorphismes \eqref{aspglob7c} et \eqref{aspglob7d}, et 
\begin{equation}\label{aspglob16c}
\beta^{(r')}\colon \cH_\mQ(\cM)\rightarrow\hupsigma^{(r')}_*(\fS^{(r')}(\cM))
\end{equation}
son adjoint \eqref{aspglob6k}. Alors,
\begin{itemize}
\item[{\rm (i)}] Pour tout nombre rationnel $r'$ tel que $0< r'\leq r$, 
le morphisme $\beta^{(r')}$ est un inverse à droite du morphisme canonique 
$\varpi^{(r')}\colon\hupsigma^{(r')}_*(\fS^{(r')}(\cM))\rightarrow \cH_\mQ(\cM)$. 
\item[{\rm (ii)}] Pour tous nombres rationnels $r'$ et $r''$ tels que $0<r''<r'\leq r$, le composé
\begin{equation}\label{aspglob16d}
\hupsigma^{(r')}_*(\fS^{(r')}(\cM))\stackrel{\varpi^{(r')}}{\longrightarrow} 
\cH_\mQ(\cM)\stackrel{\beta^{(r'')}}{\longrightarrow} \hupsigma^{(r'')}_*(\fS^{(r'')}(\cM))
\end{equation}
est le morphisme canonique. 
\end{itemize}
\end{cor}

\begin{rema}\label{aspglob160}
Sous les hypothèses de \ref{aspglob16}, l'isomorphisme $\alpha$ n'est a priori pas uniquement déterminé par $(\cM,r)$, 
mais pour tout nombre rationnel $0<r'<r$, le morphisme $\alpha^{r'}$ \eqref{aspglob16b} ne dépend que de $\cM$,
et il en dépend fonctoriellement (cf. la preuve de \cite{agt} III.12.26).
\end{rema}

\subsection{}\label{aspglob17}
Soient $r$ un nombre rationnel $>0$, $\cM$ un $\bvocB_\mQ$-module,  
$\cN$ un $\co_\fX[\frac 1 p]$-fibré de Higgs à coefficients dans $\txi^{-1}\tOmega^1_{\fX/\cS}$ tels que 
le triplet $(\cM,\cN,\alpha)$ soit $r$-admissible \eqref{aspglob9}. 
Pour éviter toute ambiguïté avec \eqref{aspglob13a}, notons 
\begin{equation}\label{aspglob17a}
\calpha\colon \fS^{(r)}(\cM)\rightarrow \hupsigma^{(r)*}(\cN)
\end{equation}
l'inverse de $\alpha$ dans $\bIMC_\mQ(\bvcC^{(r)}/\bvocB)$. 
Pour tout nombre rationnel $r'$ tel que $0< r'\leq r$, on désigne par 
\begin{equation}\label{aspglob17b}
\calpha^{(r')}\colon \fS^{(r')}(\cM) \stackrel{\sim}{\rightarrow}\hupsigma^{(r')*}(\cN)
\end{equation}
l'isomorphisme de $\bIMC_\mQ(\bvcC^{(r')}/\bvocB)$ induit par $\varepsilon^{r,r'}(\calpha)$ et les isomorphismes \eqref{aspglob7c} et \eqref{aspglob7d}, et par
\begin{equation}\label{aspglob17c}
\cbeta^{(r')}\colon \cM \rightarrow\cK^{(r')}(\hupsigma^{(r')*}(\cN))
\end{equation}
le morphisme adjoint. 

\begin{prop}\label{aspglob18}
Les hypothèses étant celles de \ref{aspglob17}, soient, de plus, $r'$, $r''$ deux nombres rationnels tels que $0<r''< r'\leq r$. Alors,
\begin{itemize} 
\item[{\rm (i)}] La limite inductive $\cV_\mQ(\cN)$ du système inductif $(\cK^{(t)}(\hupsigma^{(t)*}(\cN)))_{t\in \mQ_{>0}}$ \eqref{aspglob11c}
est représentable dans $\bMod_\mQ(\bvocB)$.
\item[{\rm (ii)}] Le morphisme composé 
\begin{equation}\label{aspglob18a}
\cM\stackrel{\cbeta^{(r')}}{\longrightarrow} \cK^{(r')}(\hupsigma^{(r')*}(\cN))\longrightarrow \cV_\mQ(\cN),
\end{equation}
où la seconde flèche est le morphisme canonique, est un isomorphisme, indépendant de $r'$. 
\item[{\rm (iii)}] Le morphisme composé 
\begin{equation}\label{aspglob18b}
\cK^{(r')}(\hupsigma^{(r')*}(\cN))\longrightarrow \cV_\mQ(\cN)\stackrel{\sim}{\longrightarrow} \cM 
\stackrel{\cbeta^{(r'')}}{\longrightarrow} \cK^{(r'')}(\hupsigma^{(r'')*}(\cN))
\end{equation}
où la première flèche est le morphisme canonique et la seconde flèche est l'isomorphisme inverse de \eqref{aspglob18a},
est le morphisme canonique \eqref{aspglob11c}.
\end{itemize}
\end{prop}

Il suffit de calquer la preuve de (\cite{agt} III.12.22) en tenant compte de \ref{ahttf32}.

\begin{cor}\label{aspglob19}
On a un foncteur
\begin{equation}\label{aspglob19a}
\cV_\mQ\colon \bMH^\qsol(\co_\fX[\frac 1 p], \txi^{-1}\tOmega^1_{\fX/\cS})\rightarrow \bMod^\Dolb_\mQ(\bvocB), 
\ \ \ \cN\mapsto \underset{\underset{r\in \mQ_{>0}}{\longrightarrow}}{\lim}\ \cK^{(r)}(\hupsigma^{(r)*}(\cN))).
\end{equation}
De plus, pour tout objet $\cN$ de $\bMH^\qsol(\co_\fX[\frac 1 p], \txi^{-1}\tOmega^1_{\fX/\cS})$, 
$\cV_\mQ(\cN)$ est associé à $\cN$.
\end{cor} 

Il résulte aussitôt des définitions que le foncteur \eqref{aspglob19a} induit un foncteur que l'on note encore 
\begin{equation}\label{aspglob19b}
\cV_\mQ\colon \bMH^\fsol(\co_\fX[\frac 1 p], \txi^{-1}\tOmega^1_{\fX/\cS})\rightarrow \bMod^\fDolb_\mQ(\bvocB).
\end{equation}

\begin{cor}\label{aspglob20}
Pour tout $\co_\fX[\frac 1 p]$-fibré de Higgs rationnellement soluble $\cN$ à coefficients dans $\txi^{-1}\tOmega^1_{\fX/\cS}$,
il existe un nombre rationnel $r>0$ et un isomorphisme de $\bIMC_\mQ(\bvcC^{(r)}/\bvocB)$
\begin{equation}\label{aspglob20a}
\calpha\colon \fS^{(r)}(\cV_\mQ(\cN))\stackrel{\sim}{\rightarrow}\hupsigma^{(r)*}(\cN)
\end{equation}
vérifiant les propriétés suivantes. Pour tout nombre rationnel $r'$ tel que $0< r'\leq r$, notons
\begin{equation}\label{aspglob20b}
\calpha^{(r')}\colon\fS^{(r')}(\cV_\mQ(\cN))\stackrel{\sim}{\rightarrow}\hupsigma^{(r')*}(\cN)
\end{equation}
l'isomorphisme de $\bIMC_\mQ(\bvcC^{(r')}/\bvocB)$ induit par $\varepsilon^{r,r'}(\calpha)$ et 
les isomorphismes \eqref{aspglob7c} et \eqref{aspglob7d}, et 
\begin{equation}\label{aspglob20c}
\cbeta^{(r')}\colon \cV_\mQ(\cN)\rightarrow\cK^{(r')}(\hupsigma^{(r')*}(\cN))
\end{equation}
son adjoint. Alors,
\begin{itemize}
\item[{\rm (i)}] Pour tout nombre rationnel $r'$ tel que $0<r'\leq r$,
le morphisme $\cbeta^{(r')}$ est un inverse à droite du morphisme canonique
$\varpi^{(r')}\colon \cK^{(r')}(\hupsigma^{(r')*}(\cN))\rightarrow \cV_\mQ(\cN)$.  
\item[{\rm (ii)}] Pour tous nombres rationnels $r'$ et $r''$ tels que $0<r''<r'\leq r$, le composé
\begin{equation}\label{aspglob20d}
\cK^{(r')}(\hupsigma^{(r')*}(\cN))\stackrel{\varpi^{(r')}}{\longrightarrow} 
\cV_\mQ(\cN)\stackrel{\cbeta^{(r'')}}{\longrightarrow} \cK^{(r'')}(\hupsigma^{(r'')*}(\cN))
\end{equation}
est le morphisme canonique. 
\end{itemize}
\end{cor}

\begin{rema}\label{aspglob200}
Sous les hypothèses de \ref{aspglob20}, l'isomorphisme $\calpha$ n'est a priori pas uniquement déterminé par $(\cN,r)$, 
mais pour tout nombre rationnel $0<r'<r$, 
le morphisme $\calpha^{r'}$ \eqref{aspglob20b} ne dépend que de $\cN$ et il en dépend fonctoriellement (cf. la preuve
de \cite{agt} III.12.26). 
\end{rema}

\begin{lem}\label{aspglob201}
Le diagramme 
\begin{equation}\label{aspglob201a}
\xymatrix{
{\bMH^\qsol(\co_\fX[\frac 1 p], \txi^{-1}\tOmega^1_{\fX/\cS})}\ar[r]^-(0.5){\cV_\mQ}\ar[d]&{\bMod^\Dolb_\mQ(\bvocB)}\ar[d]^{\upalpha_{\bvocB}}\\
{\bMH^\sol(\co_\fX[\frac 1 p], \txi^{-1}\tOmega^1_{\fX/\cS})}\ar[r]^-(0.5){\cV}&{\bIndMod^\Dolb(\bvocB)}}
\end{equation}
où $\cV$ est le foncteur \eqref{indmdlb17a} et les flèches verticales sont les foncteurs canoniques \eqref{dolbff2}, 
est commutatif à isomorphisme canonique près. 
\end{lem}

En effet, soit $\cN$ un objet de $\bMH^\qsol(\co_\fX[\frac 1 p], \txi^{-1}\tOmega^1_{\fX/\cS})$. 
En vertu de \ref{aspglob20}, il existe un nombre rationnel $r>0$ et un isomorphisme de $\bIMC_\mQ(\bvcC^{(r)}/\bvocB)$
\begin{equation}
\lambda\colon \fS^{(r)}(\cV_\mQ(\cN))\stackrel{\sim}{\rightarrow}\hupsigma^{(r)*}(\cN).
\end{equation}
D'après \ref{aspglob200}, quitte à diminuer $r$, cet isomorphisme est canonique et fonctoriel en $\cN$.  
Compte tenu de \eqref{dolbff1d} et \eqref{dolbff1e}, $\upalpha_{\cC^{(r)}/\bvocB}(\lambda)$ induit un isomorphisme de $\bIndMC(\bvcC^{(r)}/\bvocB)$
\begin{equation}
\upalpha_{\cC^{(r)}/\bvocB}(\lambda)\colon \rI\fS^{(r)}(\upalpha_{\bvocB}(\cV_\mQ(\cN)))\stackrel{\sim}{\rightarrow}\rI\hupsigma^{(r)*}(\cN).
\end{equation}
En vertu de \ref{indmdlb15}(i), on en déduit un isomorphisme fonctoriel
\begin{equation}
\upalpha_{\bvocB}(\cV_\mQ(\cN))\stackrel{\sim}{\rightarrow}\cV(\cN).
\end{equation}

\begin{teo}\label{aspglob21}
Les foncteurs \eqref{aspglob15a} et \eqref{aspglob19a} 
\begin{equation}\label{aspglob21a}
\xymatrix{
{\bMod^\Dolb_\mQ(\bvocB)}\ar@<1ex>[r]^-(0.5){\cH_\mQ}&{\bMH^\qsol(\co_\fX[\frac 1 p], \txi^{-1}\tOmega^1_{\fX/\cS})}
\ar@<1ex>[l]^-(0.5){\cV_\mQ}}
\end{equation}
sont des équivalences de catégories quasi-inverses l'une de l'autre. 
\end{teo}

Il suffit de calquer la preuve de (\cite{agt} III.12.26) en tenant compte de \ref{aspglob14}, \ref{aspglob15}, 
\ref{aspglob16}, \ref{aspglob18}, \ref{aspglob19} et \ref{aspglob20}.

\begin{cor}\label{aspglob210}
Les foncteurs \eqref{aspglob15b} et \eqref{aspglob19b} 
\begin{equation}\label{aspglob210a}
\xymatrix{
{\bMod^\fDolb_\mQ(\bvocB)}\ar@<1ex>[r]^-(0.5){\cH_\mQ}&{\bMH^\fsol(\co_\fX[\frac 1 p], \txi^{-1}\tOmega^1_{\fX/\cS})}
\ar@<1ex>[l]^-(0.5){\cV_\mQ}}
\end{equation}
sont des équivalences de catégories quasi-inverses l'une de l'autre. 
\end{cor}

\subsection{}\label{aspglob24}
Reprenons les hypothèses et notations de \ref{ahttf29}. On rappelle qu'on affecte d'un exposant $'$ les objets associés à la 
$(\cA^{\ast}_2(\oS/S),\cM_{\cA^{\ast}_2(\oS/S)})$-déformation $(\tX',\cM_{\tX'})$ \eqref{ahttf14}. 
On dispose des foncteurs  \eqref{aspglob15a} et \eqref{aspglob19a}
\begin{equation}\label{aspglob24a}
\xymatrix{
{\bMod^\Dolb_\mQ(\bvocB)} \ar@<1ex>[r]^-(0.5){\cH_\mQ}& {\bMH^\qsol(\co_\fX[\frac 1 p],\xi^{-1}\tOmega^1_{\fX/\cS})} \ar@<1ex>[l]^-(0.5){\cV_\mQ}}
\end{equation}
\begin{equation}\label{aspglob24b}
\xymatrix{
{\bMod^{\Dolb'}_\mQ(\bvocB)} \ar@<1ex>[r]^-(0.5){\cH'_\mQ}& {\bMH^{\qsol'}(\co_\fX[\frac 1 p],(\xi_\pi^{\ast})^{-1}\tOmega^1_{\fX/\cS})} \ar@<1ex>[l]^-(0.5){\cV'_\mQ}}
\end{equation}
associés aux déformations $(\tX,\cM_{\tX})$ et $(\tX',\cM_{\tX'})$, respectivement. 
L'isomorphisme $\nu$ \eqref{ahttf29d} induit un  isomorphisme $\co_{\fX}$-linéaire
\begin{equation}\label{aspglob24c}
\uplambda_*(\nu)\colon (\xi_\pi^{\ast})^{-1}\tOmega^1_{\fX/\cS}\stackrel{\sim}{\rightarrow}  \xi^{-1}\tOmega^1_{\fX/\cS}.
\end{equation}
On note 
\begin{equation}\label{aspglob24d}
\tau \colon \bMH(\co_\fX[\frac 1 p],(\xi_\pi^{\ast})^{-1}\tOmega^1_{\fX/\cS})\rightarrow \bMH(\co_\fX[\frac 1 p],\xi^{-1}\tOmega^1_{\fX/\cS})
\end{equation}
le foncteur induit par $p^\rho \uplambda_*(\nu)$. 

Compte tenu de \eqref{ahttf29h}, pour tout nombre rationnel $r\geq 0$, on a des isomorphismes canoniques 
\begin{eqnarray}
\tau \circ \hupsigma'^{(r)}_*\circ \fS'^{(r)} &\stackrel{\sim}{\rightarrow}& \hupsigma^{(r+\rho)}_*\circ \fS^{(r+\rho)},\label{aspglob24e}\\
\cK'^{(r)} \circ \hupsigma'^{(r)*}&\stackrel{\sim}{\rightarrow}&  \cK^{(r+\rho)} \circ \hupsigma^{(r+\rho)*} \circ \tau.\label{aspglob24f}
\end{eqnarray}

\begin{prop}\label{aspglob25}
Les hypothèses étant celles de \ref{aspglob24}, soient, de plus, $\cM$ un $\ocB_\mQ$-module
de Dolbeault relativement à la déformation $(\tX',\cM_{\tX'})$,
$\cN$ un $\co_\fX[\frac 1 p]$-fibré de Higgs à coefficients dans 
$(\xi_\pi^{\ast})^{-1}\tOmega^1_{\fX/\cS}$ rationnellement soluble relativement à la déformation $(\tX',\cM_{\tX'})$. 
Alors,
\begin{itemize}
\item[{\rm (i)}] Le $\ocB_\mQ$-module $\cM$ est de Dolbeault relativement à la déformation $(\tX,\cM_{\tX})$, 
et les morphismes canoniques \eqref{aspglob24e} 
induisent un isomorphisme fonctoriel
\begin{equation}\label{aspglob25a}
\tau(\cH'_\mQ(\cM))\stackrel{\sim}{\rightarrow} \cH_\mQ(\cM).
\end{equation}
\item[{\rm (ii)}] Le $\co_\fX[\frac 1 p]$-module de Higgs $\tau(\cN)$ à coefficients dans $\xi^{-1}\tOmega^1_{\fX/\cS}$ est rationnellement soluble 
relativement à la déformation $(\tX,\cM_{\tX})$ et les morphismes canoniques \eqref{aspglob24f} induisent un isomorphisme fonctoriel
\begin{equation}\label{aspglob25b}
\cV'_\mQ(\cN)\stackrel{\sim}{\rightarrow} \cV_\mQ(\tau(\cN)).
\end{equation}
\end{itemize}
\end{prop}

(i) En vertu \ref{aspglob16}, il existe un nombre rationnel $r>0$ et un isomorphisme de $\bIMC_\mQ(\bvcC'^{(r)}/\bvocB)$
\begin{equation}\label{aspglob25c}
\alpha'\colon \hupsigma'^{(r)*}(\cH'_\mQ(\cM)) \stackrel{\sim}{\rightarrow}\fS'^{(r)}(\cM).
\end{equation}
D'après \ref{aspglob160}, quitte à diminuer $r$, cet isomorphisme est canonique et fonctoriel en $\cM$.  
Compte tenu de \eqref{ahttf29h}, $\alpha'$ induit un isomorphisme de $\bIMC(\bvcC^{(r+\rho)}/\bvocB)$
\begin{equation}\label{aspglob25d}
\alpha\colon \hupsigma^{(r+\rho)*}(\tau(\cH'_\mQ(\cM))) \stackrel{\sim}{\rightarrow}\fS^{(r+\rho)}(\cM).
\end{equation}
On en déduit, en vertu de \ref{aspglob14}(i), un isomorphisme fonctoriel
\begin{equation}\label{aspglob25e}
\tau(\cH'_\mQ(\cM))\stackrel{\sim}{\rightarrow} \cH_\mQ(\cM).
\end{equation}

(ii) Il suffit de calquer la preuve de (i) en remplaçant \ref{aspglob14}(i) par \ref{aspglob18}(ii).

\begin{prop}[\cite{agt} III.15.8]\label{aspglob31}
Supposons le schéma $X$ affine et le $\co_X$-module $\tOmega^1_{X/S}$ libre. 
Alors, tout $\co_\fX[\frac 1 p]$-fibré de Higgs fortement soluble $(\cN,\theta)$ à coefficients dans $\txi^{-1}\tOmega^1_{\fX/\cS}$ est petit \eqref{definf26}.
\end{prop}

En fait, la proposition (\cite{agt} III.15.8) est formulée dans le cas absolu \eqref{definf10}, mais la preuve s'applique mutatis mutandis au cas relatif. 

\begin{rema}\label{aspglob32}
La preuve de la proposition \ref{aspglob31} dans (\cite{agt} III.15.8) utilise d'une façon cruciale la condition que $(\cN,\theta)$ est {\em fortement soluble}. 
En effet, cette preuve repose sur (\cite{agt} III.12.31) qui requiert que le $\bvocB_\mQ$-module $\cM$ soit adique de type fini, 
condition nécessaire pour avoir l'isomorphisme (\cite{agt} (III.12.29.9)). 
\end{rema}

\section{Cohomologie des ind-modules de Dolbeault}

\begin{lem}\label{indmdlb21}
Pour tout ind-$\bvocB$-module rationnel \eqref{ahttf49} et plat \eqref{indsh46} 
$\cM$ et tout entier $q\geq 0$, les morphismes canoniques, le premier de ind-$\co_\fX$-modules et le second de $\co_\fX$-modules,
\begin{eqnarray}
\rR^q\rI\hupsigma_*(\cM) &\rightarrow &
\underset{\underset{r\in \mQ_{>0}}{\longrightarrow}}{\mlq\mlq\lim \mrq\mrq}\
\rR^q\rI\hupsigma_*(\cM\otimes_{\bvocB}\mK^\bullet(\bvcC^{(r)})),\label{indmdlb21a}\\
\rR^q\vupsigma_*(\cM)  &\rightarrow &
\underset{\underset{r\in \mQ_{>0}}{\longrightarrow}}{\lim}\
\rR^q\vupsigma_*(\cM\otimes_{\bvocB}\mK^\bullet(\bvcC^{(r)})),\label{indmdlb21b}
\end{eqnarray}
où $\mK^\bullet(\bvcC^{(r)})$ est le complexe de Dolbeault du $\bvocB$-module de Higgs $(\bvcC^{(r)},p^rd_{\bvcC^{(r)}})$ \eqref{ahttf30}
et $\vupsigma_*$ est le foncteur \eqref{ahttf43i}, sont des isomorphismes.
\end{lem}

En effet, en vertu de \ref{ahttf33}, le morphisme canonique de complexes de ind-$\bvocB$-modules 
\begin{equation}\label{indmdlb21c}
\bvocB_\mQ[0]\rightarrow \underset{\underset{r\in \mQ_{>0}}{\longrightarrow}}{\mlq\mlq\lim \mrq\mrq}\ \mK^\bullet_\mQ(\bvcC^{(r)})
\end{equation}
est un quasi-isomorphisme. Par ailleurs, $\cM$ étant rationnel, pour tout $\bvocB$-module $\cF$, le morphisme canonique 
$\cM\otimes_{\bvocB}\cF\rightarrow \cM\otimes_{\bvocB}\cF_\mQ$ \eqref{indsh5g} est un isomorphisme. Comme $\cM$ est plat, 
on en déduit que le morphisme canonique de complexes de ind-$\bvocB$-modules 
\begin{equation}\label{indmdlb21d}
\cM[0]\rightarrow \underset{\underset{r\in \mQ_{>0}}{\longrightarrow}}{\mlq\mlq\lim \mrq\mrq}\ \cM\otimes_{\bvocB}\mK^\bullet(\bvcC^{(r)})
\end{equation}
est un quasi-isomorphisme. 
Comme $\rR^q\rI\hupsigma_*$ commute aux petites limites inductives filtrantes \eqref{indsh9e}, on en déduit que \eqref{indmdlb21a} est un isomorphisme
puis que \eqref{indmdlb21b} est un isomorphisme \eqref{indsh7b}.

\begin{lem}\label{indmdlb220}
Soit $\cN$ un $\co_\fX[\frac 1 p]$-fibré de Higgs à coefficients dans $\txi^{-1}\tOmega^1_{\fX/\cS}$. 
Notons $\mK^\bullet(\cN)$ le complexe de Dolbeault de $\cN$ \eqref{MH2a} et 
pour tout nombre rationnel $r\geq 0$, $\hupsigma^{(r)*}(\cN)$ l'objet de $\bIMC_\mQ(\bvcC^{(r)}/\bvocB)$ associé à $\cN$ \eqref{aspglob6g}
et $\mK^\bullet(\hupsigma^{(r)*}(\cN))$ son complexe de Dolbeault dans $\bMod_\mQ(\bvocB)$ \eqref{aspglob6}. 
On a alors un morphisme canonique fonctoriel de $\bD^+(\bMod(\co_\fX[\frac 1 p]))$ 
\begin{equation}\label{indmdlb220ab}
\mK^\bullet(\cN)\rightarrow 
\rR \hupsigma_{\mQ*}(\mK^\bullet(\hupsigma^{(r)*}(\cN))),
\end{equation}
où $\hupsigma_{\mQ*}$ est le foncteur \eqref{ahttf40l}. C'est en fait un morphisme de systèmes inductifs indexés par $r\in \mQ_{>0}$, 
où les morphismes de transition du but sont induits par les homomorphismes 
$\bvalpha^{r,r'}$ \eqref{ahttf14f} pour $r\geq r'>0$. 
De plus, pour tout entier $q\geq 0$, le morphisme induit 
\begin{equation}\label{indmdlb220a}
\rH^q(\mK^\bullet(\cN))
\rightarrow \underset{\underset{t\in \mQ_{>r}}{\longrightarrow}}{\lim}\
\rR^q\hupsigma_{\mQ*}(\mK^\bullet(\hupsigma^{(r)*}(\cN)))
\end{equation}
est un isomorphisme. 
\end{lem}

Soit $r$ un nombre rationnel $>0$.
On a un morphisme canonique de complexes de  $\co_\fX[\frac 1 p]$-modules 
\begin{equation}\label{indmdlb220i}
\mK^\bullet(\cN)\rightarrow 
\hupsigma_{\mQ*}(\mK^\bullet(\hupsigma^{(r)*}(\cN))).
\end{equation}
Comme la catégorie abélienne $\bMod_{\mQ}(\bvocB)$ a assez d'injectifs \eqref{indsh14},  
il existe un complexe borné inférieurement de $\bvocB_\mQ$-modules injectifs $\cL^\bullet$
et un quasi-isomorphisme $u\colon \mK^\bullet(\uhupsigma'^{(r)*}(\cN))\rightarrow \cL^\bullet$ (\cite{sp} \href{https://stacks.math.columbia.edu/tag/013K}{013K}).
Ce dernier induit un morphisme de $\bD^+(\bMod(\co_\fX[\frac 1 p]))$ 
\begin{equation}\label{indmdlb220j}
\hupsigma_{\mQ*}(\mK^\bullet(\hupsigma^{(r)*}(\cN)))\rightarrow 
\rR \hupsigma_{\mQ*}(\mK^\bullet(\hupsigma^{(r)*}(\cN))),
\end{equation}
qui ne depend que de $\mK^\bullet(\hupsigma^{(r)*}(\cN))$, mais pas de $u$
(\cite{sp} \href{https://stacks.math.columbia.edu/tag/05TG}{05TG}). 
On prend pour morphisme \eqref{indmdlb220ab} le composé de \eqref{indmdlb220i} et \eqref{indmdlb220j}. 

Compte tenu de (\cite{agt} III.12.4(i) et III.12.2(i)), pour tout entier $j\geq 0$, 
$d_{\bvcC^{(r)}}$ \eqref{indmdlb1a} induit un morphisme $\co_\fX$-linéaire
\begin{equation}\label{indmdlb220b}
\delta^{j,(r)}\colon \rR^j\hupsigma_*(\bvcC^{(r)})\rightarrow \txi^{-1}\tOmega^1_{\fX/\cS}\otimes_{\co_\fX}\rR^j\hupsigma_*(\bvcC^{(r)}),
\end{equation}
qui est un $\co_\fX$-champ de Higgs sur $\rR^j\hupsigma_*(\bvcC^{(r)})$ à coefficients dans 
$\txi^{-1}\tOmega^1_{\fX/\cS}$. On note $\theta$ le $\co_\fX[\frac 1 p]$-champ de Higgs sur $\cN$
et $\vartheta^{j,(r)}_\tot=\theta\otimes \id+p^r\id \otimes \delta^{j,(r)}$ le 
$\co_\fX[\frac 1 p]$-champ de Higgs total sur $\cN\otimes_{\co_\fX}\rR^j\hupsigma_*(\hcC^{(r)})$ \eqref{MH2d}.
D'après (\cite{agt} III.12.4(ii)), pour tout entier $i\geq 0$, on a un isomorphisme canonique $\co_\fX[\frac 1 p]$-linéaire
\begin{equation}\label{indmdlb220c}
\rR^j\hupsigma_{\mQ*}(\mK^i(\hupsigma^{(r)*}(\cN)))\stackrel{\sim}{\rightarrow} 
\mK^i(\cN\otimes_{\co_\fX}\rR^j\hupsigma_*(\bvcC^{(r)}),\vartheta_\tot^{j,(r)}),
\end{equation} 
compatible avec les morphismes induits par les différentielles des deux complexes de Dolbeault.

Par ailleurs, la catégorie abélienne $\bMod_{\mQ}(\bvocB)$ ayant assez d'injectifs,
on a une suite spectrale canonique fonctorielle
\begin{equation}\label{indmdlb220d}
{^r\rE}_1^{i,j}=\rR^j\hupsigma_{\mQ*}(\mK^i(\hupsigma^{(r)*}(\cN)))\Rightarrow \rR^{i+j}\hupsigma_{\mQ*}(\mK^\bullet(\hupsigma^{(r)*}(\cN))).
\end{equation}
Celle-ci induit, pour tout entier $i\geq 0$, un morphisme canonique
\begin{equation}
\rH^i({^t\rE}_1^{\bullet,0})\rightarrow 
\rR^i\hupsigma_{\mQ*}(\umK^\bullet(\uhupsigma^{(r)*}(\cN))),
\end{equation}
qui n'est autre que le morphisme induit par \eqref{indmdlb220j} d'après (\cite{ega3} (0.11.3.4.2)).

En vertu de \ref{ahttf28} et \eqref{indmdlb220c}, pour tout $i\geq 0$, on a un isomorphisme canonique 
\begin{equation}\label{indmdlb220e}
\underset{\underset{r\in \mQ_{>0}}{\longrightarrow}}{\lim}\ {^r\rE}_1^{i,0}\stackrel{\sim}{\rightarrow}\mK^i(\cN,\theta),
\end{equation}
et pour tout $j\geq 1$, on a 
\begin{equation}\label{indmdlb220f}
\underset{\underset{r\in \mQ_{>0}}{\longrightarrow}}{\lim}\ {^r\rE}_1^{i,j}=0.
\end{equation}
De plus, les isomorphismes \eqref{indmdlb220e} (pour $i\in \mN$) forment un isomorphisme de complexes. 
Comme les petites limites inductives filtrantes existent et sont exactes dans $\bMod(\co_\fX[\frac 1 p])$ (\cite{sga4} II 4.3), 
on en déduit que les morphismes \eqref{indmdlb220a} sont des isomorphismes.

\begin{lem}\label{indmdlb22}
Soit $\cN$ un $\co_\fX[\frac 1 p]$-fibré de Higgs à coefficients dans $\txi^{-1}\tOmega^1_{\fX/\cS}$. 
Notons $\mK^\bullet(\cN)$ le complexe de Dolbeault de $\cN$ et 
pour tout nombre rationnel $r\geq 0$, $\rI\hupsigma^{(r)*}(\cN)$ l'objet de $\bIndMC(\bvcC^{(r)}/\bvocB)$ associé à $\cN$ \eqref{indmdlb1f}
et $\mK^\bullet(\rI\hupsigma^{(r)*}(\cN))$ son complexe de Dolbeault \eqref{indsh30c}. 
On a alors un morphisme canonique fonctoriel de $\bD^+(\bMod(\co_\fX))$ 
\begin{equation}\label{indmdlb22b}
\mK^\bullet(\cN)\rightarrow 
\rR \vupsigma_*(\mK^\bullet(\rI\hupsigma^{(r)*}(\cN))),
\end{equation}
où $\vupsigma_*$ est le foncteur \eqref{ahttf43i}. C'est en fait un morphisme de systèmes inductifs indexés par $r\in \mQ_{>0}$, 
où les morphismes de transition du but sont induits par les homomorphismes 
$\bvalpha^{r,r'}$ \eqref{ahttf14f} pour $r\geq r'>0$. 
De plus, le morphisme induit de $\bD^+(\bMod(\co_\fX))$ 
\begin{equation}\label{indmdlb22c}
\mK^\bullet(\cN)\rightarrow 
\rR\vupsigma_*(\underset{\underset{t\in \mQ_{>r}}{\longrightarrow}}{\mlq\mlq\lim \mrq\mrq}\ \mK^\bullet(\rI\hupsigma^{(r)*}(\cN)))
\end{equation}
est un isomorphisme. 
\end{lem}

Compte tenu de \eqref{dolbff1e} et \eqref{ahttf40n}, on prend pour morphisme \eqref{indmdlb22b} l'image canonique du morphisme \eqref{indmdlb220ab}.
D'après \ref{indsh17}, pour tout entier $q\geq 0$, on a  $\rR^q\vupsigma_*=\kappa_{\co_\fX}\circ \rR^q\rI\hupsigma_*$ \eqref{ahttf43i} 
et ce foncteur commute aux petites limites inductives filtrantes \eqref{indsh9e}.
Il résulte alors de \ref{indmdlb220} que \eqref{indmdlb22c} est un isomorphisme.

\begin{teo}\label{indmdlb23}
Soient $\cM$ un ind-$\bvocB$-module de Dolbeault \eqref{indmdlb5}, $q$ un entier $\geq 0$. 
Notons $\mK^\bullet(\cH(\cM))$ le complexe de Dolbeault du $\co_\fX[\frac 1 p]$-fibré de Higgs 
$\cH(\cM)$ \eqref{indmdlb12}. On a alors un isomorphisme canonique fonctoriel de $\bD^+(\bMod(\co_\fX))$
\begin{equation}\label{indmdlb23a}
\rR\vupsigma_*(\cM)\stackrel{\sim}{\rightarrow}\mK^\bullet(\cH(\cM)),
\end{equation}
où $\vupsigma_*$ est le foncteur \eqref{ahttf43i}.
\end{teo}

En effet, en vertu de \ref{indmdlb13}, il existe un nombre rationnel $r_\cM>0$ et un isomorphisme de $\bIndMC(\bvcC^{(r_\cM)}/\bvocB)$, 
\begin{equation}\label{indmdlb23b}
\alpha_\cM\colon \rI\hupsigma^{(r_\cM)*}(\cH(\cM))\stackrel{\sim}{\rightarrow} \rI\fS^{(r_\cM)}(\cM),
\end{equation}
vérifiant les propriétés \ref{indmdlb13}(i)-(ii). Pour tout nombre rationnel $r$ tel que $0< r< r_\cM$, on désigne par
\begin{equation}\label{indmdlb23c}
\alpha^{(r)}_\cM\colon \rI\hupsigma^{(r)*}(\cH(\cM))\stackrel{\sim}{\rightarrow} \rI\fS^{(r)}(\cM)
\end{equation}
l'isomorphisme de $\bIndMC(\bvcC^{(r)}/\bvocB)$ induit par $\rI\varepsilon^{r_\cM,r}(\alpha_\cM)$ et 
les isomorphismes \eqref{indmdlb3b} et \eqref{indmdlb3c}. D'après la preuve de \ref{indmdlb20}, 
$\alpha^{(r)}_\cM$ ne dépend que de $\cM$ (mais pas de $\alpha_\cM$)
et il en dépend fonctoriellement. On note $\mK^\bullet(\rI\hupsigma^{(r)*}(\cH(\cM)))$ 
le complexe de Dolbeault de $\rI\hupsigma^{(r)*}(\cH(\cM))$ \eqref{indsh30c}.
L'isomorphisme $\alpha^{(r)}_\cM$ induit un isomorphisme \eqref{ahttf30}
\begin{equation}
\mK^\bullet(\rI\hupsigma^{(r)*}(\cH(\cM))) \stackrel{\sim}{\rightarrow}
\cM\otimes_{\bvocB}\mK^\bullet(\bvcC^{(r)}),
\end{equation}
où $\mK^\bullet(\bvcC^{(r)})$ est le complexe de Dolbeault du $\bvocB$-module de Higgs $(\bvcC^{(r)},p^rd_{\bvcC^{(r)}})$ \eqref{ahttf30}.
Ces isomorphismes forment un isomorphisme de systèmes inductifs (pour $0<r<r_\cM$).
On en déduit un isomorphisme canonique fonctoriel de complexes de ind-$\bvocB$-modules
\begin{equation}
\underset{\underset{r\in \mQ_{>0}}{\longrightarrow}}{\mlq\mlq\lim \mrq\mrq}\ 
\mK^\bullet(\rI\hupsigma^{(r)*}(\cH(\cM))) \stackrel{\sim}{\rightarrow}
\underset{\underset{r\in \mQ_{>0}}{\longrightarrow}}{\mlq\mlq\lim \mrq\mrq}\ 
\cM\otimes_{\bvocB}\mK^\bullet(\bvcC^{(r)}).
\end{equation}
Comme le ind-$\bvocB$-module $\cM$ est rationnel et plat d'après \ref{indmdlb6}, le théorème s'ensuit compte tenu de \eqref{indmdlb21d} et \ref{indmdlb22}.

\begin{cor}\label{aspglob22}
Soient $\cM$ un $\bvocB_\mQ$-module de Dolbeault \eqref{aspglob1},
$\mK^\bullet(\cH_\mQ(\cM))$ le complexe de Dolbeault du $\co_\fX[\frac 1 p]$-fibré de Higgs $\cH_\mQ(\cM)$ \eqref{aspglob11b}. 
Alors, on a un isomorphisme canonique fonctoriel de $\bD^+(\bMod(\co_\fX))$
\begin{equation}\label{aspglob22a}
\rR\hupsigma_{\mQ*}(\cM)\stackrel{\sim}{\rightarrow}\mK^\bullet(\cH_\mQ(\cM)),
\end{equation}
où l'on a encore noté $\rR\hupsigma_\mQ$ le composé du foncteur \eqref{ahttf40m} et du foncteur 
canonique de $\bD^+(\bMod(\co_\fX[\frac 1 p]))$ dans $\bD^+(\bMod(\co_\fX))$.
En particulier, pour tout entier $q\geq 0$, on a un isomorphisme canonique de $\co_\fX[\frac 1 p]$-modules
\begin{equation}\label{aspglob22b}
\rR^q\hupsigma_{\mQ*}(\cM)\stackrel{\sim}{\rightarrow}\rH^q(\mK^\bullet(\cH_\mQ(\cM))). 
\end{equation}
\end{cor}

Cela résulte de \ref{indmdlb23}, \eqref{aspglob11d} et \eqref{ahttf40n}

\begin{rema}
Les isomorphismes \eqref{aspglob22b} ont été démontrés dans (\cite{agt} III.12.34) dans le cas absolu \eqref{definf10}. 
La preuve est essentiellement la même que celle de \ref{indmdlb23}. 
Le passage aux ind-modules permet d'obtenir l'isomorphisme \eqref{aspglob22a}. 
\end{rema}

\section{Modules de Dolbeault sur un schéma affine petit}\label{mdpsa}

\subsection{}\label{mdpsa1}
On suppose dans cette section que $X$ est un objet de $\bP$ \eqref{ahttf4} et que $X_s$ est non-vide. 
On se donne une carte adéquate $((P,\gamma),(\mN,\iota),\vartheta)$ pour $f$ \eqref{cad1}. 
On pose $R=\Gamma(X,\co_X)$, $R_1=R\otimes_{\co_K}\co_\oK$ et 
\begin{equation}\label{mdpsa1b}
\tOmega^1_{R/\co_K}=\Gamma(X,\tOmega^1_{X/S}).
\end{equation}

Si $A$ est un anneau et $M$ un $A$-module, on note encore $A$ (resp. $M$) le faisceau constant 
de valeur $A$ (resp. $M$) de $\oX^\circ_\fet$ ou $(\oX^\circ_\fet)^{\mN^\circ}$, selon le contexte. 

On désigne par $\bvocB_X$ l'anneau $(\ocB_{X,n+1})_{n\in \mN}$ de $(\oX^\circ_\fet)^{\mN^\circ}$ \eqref{ahttf3e},
par $\bMod(\bvocB_X)$ la catégorie des $\bvocB_X$-modules de $(\oX^\circ_\fet)^{\mN^\circ}$ 
et par $\bMod_{\mQ}(\bvocB_X)$ la catégorie des objets de $\bMod(\bvocB_X)$ à isogénie près \eqref{caip1a}. 
On note
\begin{equation}\label{mdpsa13a}
\bMod(\bvocB_X)\rightarrow \bMod_\mQ(\bvocB_X),\ \ \ \cM\mapsto \cM_{\mQ}
\end{equation}
le foncteur canonique. La catégorie $\bMod_\mQ(\bvocB_X)$ est  abélienne et monoïdale symétrique, ayant $\bvocB_{X,\mQ}$ pour objet unité. 
Les objets de $\bMod_{\mQ}(\bvocB_X)$ seront aussi appelés des {\em $\bvocB_{X,\mQ}$-modules}.

On note 
\begin{equation}\label{mdpsa13c}
\bvbeta\colon (\tE_s^{\mN^\circ},\bvocB) \rightarrow ((\oX^\circ_\fet)^{\mN^\circ},\bvocB_X)
\end{equation}
le morphisme de topos annelés  induit par les $(\beta_{n+1})_{n\in \mN}$ \eqref{ahttf42b}. 
Nous utilisons pour les modules la notation $\bvbeta^{-1}$ pour désigner l'image
inverse au sens des faisceaux abéliens et nous réservons la notation 
$\bvbeta^*$ pour l'image inverse au sens des modules.

Le foncteur $\bvbeta_*$ induit un foncteur additif et exact à gauche 
\begin{equation}\label{mdpsa13d}
\bvbeta_{\mQ *}\colon \bMod_{\mQ}(\bvocB) \rightarrow \bMod_{\mQ}(\bvocB_X).
\end{equation}
Le foncteur $\bvbeta^*$ induit un foncteur additif 
\begin{equation}\label{mdpsa13e}
\bvbeta^*_\mQ\colon \bMod_{\mQ}(\bvocB_X) \rightarrow \bMod_{\mQ}(\bvocB). 
\end{equation}

\begin{defi}\label{mdpsa13}
On dit qu'un $\bvocB_{X,\mQ}$-module est {\em adique de type fini} s'il est isomorphe à $\cM_\mQ$ \eqref{mdpsa13a}, 
où $\cM$ est un $\bvocB_X$-module adique de type fini (\cite{agt} III.7.16). 
\end{defi}

On désigne par $\bMod^\atf(\bvocB_X)$ la sous-catégorie pleine de $\bMod(\bvocB_X)$ formée des $\bvocB_X$-modules adiques de type fini 
et par $\bMod^\atf_{\mQ}(\bvocB_X)$ la catégorie des objets de $\bMod^\atf(\bvocB_X)$ à isogénie près \eqref{caip1}. Le foncteur canonique 
\begin{equation}\label{mdpsa13b}
\bMod^\atf_{\mQ}(\bvocB_X)\rightarrow\bMod_{\mQ}(\bvocB_X)
\end{equation} 
étant pleinement fidèle, pour qu'un $\bvocB_{X,\mQ}$-module soit adique de type fini, il faut et il suffit qu'il soit dans l'image essentielle de ce foncteur.

\begin{prop}[\cite{agt} III.13.2]\label{mdpsa2}
Pour tout $\co_\fX$-module cohérent $\cN$, on a un isomorphisme $\bvocB$-linéaire canonique et fonctoriel
\begin{equation}\label{mdpsa2a}
\bvbeta^*(\cN(\fX)\otimes_{\hRun}\bvocB_X)\stackrel{\sim}{\rightarrow}\hupsigma^*(\cN),
\end{equation}
où $\hupsigma$ est le morphisme de topos annelés \eqref{ahttf13e}.
\end{prop}

\subsection{}\label{mdpsa3}
Soit $r$ un nombre rationnel $\geq 0$. Avec les notations de \ref{ahttf35}, 
on désigne par $\bvcF_X^{(r)}$ le $\bvocB_X$-module $(\cF_{X,n+1}^{(r)})_{n\in \mN}$ 
et par $\bvcC_X^{(r)}$ la $\bvocB_X$-algèbre $(\cC_{X,n+1}^{(r)})_{n\in \mN}$.
D'après (\cite{agt} III.7.3(i) et (III.7.12.1)), on a une suite exacte de $\bvocB_X$-modules 
\begin{equation}\label{mdpsa3a}
0\rightarrow \bvocB_X\rightarrow \bvcF^{(r)}_X\rightarrow 
\xi^{-1}\tOmega^1_{R/\co_K}\otimes_{R}\bvocB_X\rightarrow 0.
\end{equation} 
Compte tenu de (\cite{agt} III.7.3(i) et (III.7.12.3)), on a un isomorphisme canonique de $\bvocB_X$-algèbres
\begin{equation}\label{mdpsa3b}
\bvcC^{(r)}_X\stackrel{\sim}{\rightarrow} 
\underset{\underset{m\geq 0}{\longrightarrow}}\lim\ \rS^m_{\bvocB_X}(\bvcF^{(r)}_X).
\end{equation}

Pour tous nombres rationnels $r\geq r'\geq 0$, les morphismes 
$(\tta_{X,n+1}^{r,r'})_{n\in \mN}$ \eqref{ahttf35c} induisent un morphisme $\bvocB_X$-linéaire 
\begin{equation}\label{mdpsa3c}
\bvtta^{r,r'}_X\colon \bvcF^{(r)}_X\rightarrow \bvcF^{(r')}_X.
\end{equation}
Les homomorphismes 
$(\alpha_{X,n+1}^{r,r'})_{n\in \mN}$ \eqref{ahttf35d} induisent un homomorphisme de $\bvocB_X$-algèbres 
\begin{equation}\label{mdpsa3d}
\bvalpha^{r,r'}_X\colon \bvcC^{(r)}_X\rightarrow \bvcC^{(r')}_X.
\end{equation}
Pour tous nombres rationnels $r\geq r'\geq r''\geq 0$, on a
\begin{equation}\label{mdpsa3e}
\bvtta^{r,r''}_X=\bvtta^{r',r''}_X \circ \bvtta^{r,r'}_X\ \ \ {\rm et}\ \ \  
\bvalpha^{r,r''}_X=\bvalpha^{r',r''}_X \circ \bvalpha^{r,r'}_X.
\end{equation}

On a un isomorphisme canonique $\bvcC_X^{(r)}$-linéaire
\begin{equation}\label{mdpsa3f}
\Omega^1_{\bvcC_X^{(r)}/\bvocB_X}\stackrel{\sim}{\rightarrow} 
\txi^{-1}\tOmega^1_{R/\co_K}\otimes_R\bvcC_X^{(r)}.
\end{equation}
La $\bvocB_X$-dérivation universelle de $\bvcC_X^{(r)}$ correspond via cet isomorphisme à l'unique 
$\bvocB_X$-dérivation 
\begin{equation}\label{mdpsa3g}
d_{\bvcC^{(r)}_X}\colon \bvcC_X^{(r)}\rightarrow \txi^{-1}\tOmega^1_{R/\co_K}\otimes_R\bvcC_X^{(r)}
\end{equation}
qui prolonge le morphisme canonique $\bvcF^{(r)}_X\rightarrow \txi^{-1}\tOmega^1_{R/\co_K}\otimes_{R}\bvocB_X$ 
\eqref{mdpsa3a}. Comme 
\begin{equation}\label{mdpsa3h}
\txi^{-1}\tOmega^1_{R/\co_K}\otimes_R\bvocB_X =d_{\bvcC^{(r)}_X} (\bvcF_X^{(r)}) 
\subset d_{\bvcC^{(r)}_X}(\bvcC_X^{(r)}),
\end{equation} 
la dérivation $d_{\bvcC^{(r)}_X}$ est un $\bvocB_X$-champ de Higgs à coefficients dans 
$\txi^{-1}\tOmega^1_{R/\co_K}$ d'après \ref{MH8}(i). 
Pour tous nombres rationnels $r\geq r'\geq 0$, on a 
\begin{equation}\label{mdpsa3i}
p^{r-r'}(\id \otimes \bvalpha^{r,r'}_X) \circ d_{\bvcC^{(r)}_X}=d_{\bvcC^{(r')}_X}\circ \bvalpha^{r,r'}_X.
\end{equation}

\begin{prop}\label{mdpsa6}
Pour tout nombre rationnel $r\geq 0$, les morphismes canoniques
\begin{eqnarray}
\bvbeta^*(\bvcF^{(r)}_X)&\rightarrow&\bvcF^{(r)},\label{mdpsa6a}\\
\bvbeta^*(\bvcC^{(r)}_X)&\rightarrow&\bvcC^{(r)},\label{mdpsa6b}
\end{eqnarray}
sont des isomorphismes. De plus, pour tous nombres rationnels 
$r\geq r'\geq 0$, les morphismes $\bvbeta^*(\bvtta_X^{r,r'})$ et $\bvbeta^*(\bvalpha_X^{r,r'})$ 
s'identifient aux morphismes $\bvtta^{r,r'}$ \eqref{ahttf14e} et $\bvalpha^{r,r'}$ \eqref{ahttf14f}, respectivement.
\end{prop}

Le cas absolu \eqref{definf10}  où $X$ est un objet de $\bQ$ \eqref{ahttf5} a été démontré dans (\cite{agt} III.13.4).
La condition que $X$ est un objet de $\bQ$ est en fait superflue (cf. \ref{ahttf39}). Le cas relatif se traite de même.

\subsection{}\label{mdpsa4}
Soit $r$ un nombre rationnel $\geq 0$. 
On notera simplement $\bIMC(\bvcC_X^{(r)}/\bvocB_X)$ la catégorie des $p^r$-isoconnexions intégrables relativement à l'extension $\bvcC_X^{(r)}/\bvocB_X$ (cf. \ref{indsh24}) ; on omet donc l'exposant $p^r$ de la notation  introduite dans \ref{indsh24} considérant qu'il est redondant avec l'exposant $r$ de $\bvcC_X^{(r)}$. 
C'est une catégorie additive.
On désigne par $\bIMC_\mQ(\bvcC_X^{(r)}/\bvocB_X)$ la catégorie des objets de $\bIMC(\bvcC_X^{(r)}/\bvocB_X)$ à isogénie près.

Considérons le foncteur 
\begin{equation}\label{mdpsa4a}
\fS^{(r)}_X\colon
\begin{array}[t]{clcr}
\bMod(\bvocB_X)&\rightarrow& \bIMC(\bvcC_X^{(r)}/\bvocB_X)\\
\cM&\mapsto& (\bvcC^{(r)}_X\otimes_{\bvocB_X}\cM,\bvcC^{(r)}_X\otimes_{\bvocB_X}\cM,
\id,p^rd_{\bvcC^{(r)}_X}\otimes \id)
\end{array}
\end{equation}
et notons encore 
\begin{equation}\label{mdpsa4d}
\fS^{(r)}_X\colon\bMod_\mQ(\bvocB_X)\rightarrow \bIMC_\mQ(\bvcC_X^{(r)}/\bvocB_X)
\end{equation}
le foncteur induit.

Reprenons les notations de \ref{definf19}.
D'après \ref{indsh26}, si $(\cN,\cN',v,\theta)$ est une $\co_\fX$-isogénie de Higgs
à coefficients dans $\txi^{-1}\tOmega^1_{\fX/\cS}$,
\begin{equation}\label{mdpsa4b}
(\bvcC^{(r)}_X\otimes_{\hRun}\cN(\fX),\bvcC^{(r)}_X\otimes_{\hRun}\cN'(\fX),\id \otimes_{\hRun}v,
p^r d_{\bvcC^{(r)}_X} \otimes v+\id \otimes \theta)
\end{equation}
est un objet de $\bIMC(\bvcC_X^{(r)}/\bvocB_X)$. On obtient ainsi un foncteur 
\begin{equation}\label{mdpsa4c}
\hupsigma^{(r)*}_X\colon\bIH(\co_\fX,\txi^{-1}\tOmega^1_{\fX/\cS})\rightarrow \bIMC(\bvcC_X^{(r)}/\bvocB_X).
\end{equation}
D'après \eqref{definf19d}, celui-ci induit un foncteur que l'on note encore
\begin{equation}\label{mdpsa4e}
\hupsigma^{(r)*}_X\colon\bMH^\coh(\co_\fX[\frac 1 p], \txi^{-1}\tOmega^1_{\fX/\cS})\rightarrow \bIMC_\mQ(\bvcC_X^{(r)}/\bvocB_X).
\end{equation}

Pour tous nombres rationnels $r\geq r'\geq 0$, on dispose d'un foncteur canonique
\begin{equation}\label{mdpsa4f}
\varepsilon_X^{r,r'}\colon \bIMC(\bvcC_X^{(r)}/\bvocB_X)\rightarrow \bIMC(\bvcC_X^{(r')}/\bvocB_X).
\end{equation}
analogue du foncteur \eqref{aspglob7a}. Celui-ci induit un foncteur que l'on note encore 
\begin{equation}\label{mdpsa4g}
\varepsilon_X^{r,r'}\colon \bIMC_\mQ(\bvcC_X^{(r)}/\bvocB_X)\rightarrow \bIMC_\mQ(\bvcC_X^{(r')}/\bvocB_X).
\end{equation}
On a un isomorphisme canonique de foncteurs de $\bMod(\bvocB_X)$ dans $\bIMC(\bvcC_X^{(r')}/\bvocB_X)$ 
(resp.  de $\bMod_\mQ(\bvocB_X)$ dans $\bIMC_\mQ(\bvcC_X^{(r')}/\bvocB_X)$) 
\begin{equation}\label{mdpsa4h}
\varepsilon_X^{r,r'}\circ \fS_X^{(r)}\stackrel{\sim}{\longrightarrow} \fS_X^{(r')},
\end{equation}
et un isomorphisme canonique de foncteurs de $\bIH(\co_\fX,\xi^{-1}\tOmega^1_{\fX/\cS})$ dans $\bIMC(\bvcC_X^{(r')}/\bvocB_X)$ 
(resp.  de $\bMH^\coh(\co_\fX[\frac 1 p], \txi^{-1}\tOmega^1_{\fX/\cS})$ dans $\bIMC_\mQ(\bvcC_X^{(r')}/\bvocB_X)$) 
\begin{equation}\label{mdpsa4i}
\varepsilon_X^{r,r'}\circ \hupsigma_X^{(r)*}\stackrel{\sim}{\longrightarrow} \hupsigma_X^{(r')*}.
\end{equation}

\begin{prop}[\cite{agt} III.13.7]\label{mdpsa5}
Pour tout nombre rationnel $r\geq 0$, les diagrammes de foncteurs
\begin{equation}\label{mdpsa5a}
\xymatrix{
{\bMod(\bvocB_X)}\ar[r]^-(0.5){\fS^{(r)}_X}\ar[d]_{\bvbeta^*}&{\bIMC(\bvcC_X^{(r)}/\bvocB_X)}\ar[d]^{\bvbeta^*}\\
{\bMod(\bvocB)}\ar[r]^-(0.5){\fS^{(r)}}&{\bIMC(\bvcC^{(r)}/\bvocB)}}
\end{equation}
\begin{equation}\label{mdpsa5b}
\xymatrix{
{\bIH^\coh(\co_\fX,\txi^{-1}\tOmega^1_{\fX/\cS})}\ar[r]^-(0.5){\hupsigma^{(r)*}_X}\ar[rd]_-(0.5){\hupsigma^{(r)*}}&{\bIMC(\bvcC_X^{(r)}/\bvocB_X)}\ar[d]^{\bvbeta^*}\\
&{\bIMC(\bvcC^{(r)}/\bvocB)}}
\end{equation}
où  l'image inverse $\bvbeta^*$ pour les $p^r$-isoconnexions est définie dans \ref{indsh25}, sont commutatifs à isomorphismes
canoniques près.  
\end{prop}

\begin{defi}\label{mdpsa30}
Soient $\cM$ un $\bvocB_{X,\mQ}$-module \eqref{mdpsa1}, 
$\cN$ un $\co_\fX[\frac 1 p]$-fibré de Higgs à coefficients dans $\txi^{-1}\tOmega^1_{\fX/\cS}$ \eqref{definf20}.
\begin{itemize}
\item[(i)] Soit $r$ un nombre rationnel $>0$.  On dit que $\cM$ et  $\cN$ sont {\em $r$-associés} 
s'il existe un isomorphisme de $\bIMC_\mQ(\bvcC^{(r)}_X/\bvocB_X)$ 
\begin{equation}\label{mdpsa30a}
\alpha\colon \hupsigma_X^{(r)*}(\cN) \stackrel{\sim}{\rightarrow}\fS_X^{(r)}(\cM).
\end{equation}
On dit alors aussi que le triplet $(\cM,\cN,\alpha)$ est {\em $r$-admissible}. 
\item[(ii)] On dit que $\cM$ et  $\cN$ sont {\em associés} s'il existe un nombre rationnel $r>0$ tel que 
$\cM$ et  $\cN$ soient $r$-associés.
\end{itemize}
\end{defi}

On notera que pour tous nombres rationnels $r\geq r'>0$, 
si $\cM$ et  $\cN$ sont $r$-associés, ils sont $r'$-associés, compte tenu de \eqref{mdpsa4h} et \eqref{mdpsa4i}.

\begin{rema}\label{mdpsa31}
On notera que tout $\co_\fX[\frac 1 p]$-module localement projectif de type fini \eqref{notconv14} est en fait projectif de type fini, 
{\em i.e.} est facteur direct d'un $\co_\fX[\frac 1 p]$-module libre de type fini, d'après (\cite{agt} III.6.17).
\end{rema}

\begin{lem}\label{mdpsa21}
Soient $\cM$ un $\bvocB_{X,\mQ}$-module, 
$\cN$ un $\co_\fX[\frac 1 p]$-fibré de Higgs à coefficients dans $\txi^{-1}\tOmega^1_{\fX/\cS}$, $r$ un nombre rationnel $>0$.
Si $\cM$ et $\cN$ sont $r$-associés, alors le $\bvocB_\mQ$-module $\bvbeta^*_\mQ(\cM)$ et 
le $\co_\fX[\frac 1 p]$-fibré de Higgs $\cN$ sont $r$-associés dans le sens de \ref{aspglob9}. 
\end{lem}

Cela résulte aussitôt de \ref{mdpsa5}.

\subsection{}\label{mdpsa12}
Le schéma $\oX$ étant localement irréductible d'après (\cite{ag} 4.2.7(iii)),  
il est la somme des schémas induits sur ses composantes irréductibles que l'on note $\oX_1,\dots,\oX_c$. 
Pour tout $1\leq i\leq c$, on pose $\cR_i=\Gamma(\oX_i,\co_{\oX})$ et on note $\hcR_i$ le séparé complété $p$-adique de $\cR_i$.
On se donne un point géométrique $\oy_i$ de $\oX^\circ_i=\oX_i\times_XX^\circ$ et on pose $\Delta_i=\pi_1(\oX^{\circ}_i,\oy_i)$.
On désigne par $\bB_{\Delta_i}$ le topos classifiant de $\Delta_i$, par
\begin{equation}\label{mdpsa12a}
\nu_i\colon \oX^{\circ}_{i,\fet} \stackrel{\sim}{\rightarrow} \bB_{\Delta_i}
\end{equation}
le foncteur fibre  de $\oX^{\circ}_{i,\fet}$ en $\oy_i$ (\cite{agt}  (VI.9.8.4)) et par
\begin{equation}\label{mdpsa12b}
\mu_i\colon \bB_{\Delta_i}\rightarrow \oX^\circ_{i,\fet}
\end{equation}
le foncteur quasi-inverse défini dans (\cite{agt}  (VI.9.8.3)). 

On notera que l'anneau $R_1=R\otimes_{\co_K}\co_\oK$ s'identifie au produit des anneaux $\cR_i$ pour $1\leq i\leq c$. 

\subsection{}\label{mdpsa8}
Conservons les notations de \ref{mdpsa12} et reprenons celles de \ref{taht}. 
Pour tout entier $n\geq 1$, on considère le schéma logarithmique $(X^{(n)},\cM_{X^{(n)}})$ défini par la formule \eqref{cad4c}.
Pour tout $1\leq i\leq c$, on se donne un $X$-morphisme 
\begin{equation}
\oy_i\rightarrow \underset{\underset{n\geq 1}{\longleftarrow}}{\lim}\ X^{(n)}.
\end{equation} 
On peut alors appliquer les constructions de \ref{pmh}, en particulier celles de \ref{pmh5}. 

Soient $\cN$ un $\co_\fX$-module cohérent et $\cS$-plat, $\theta$ un $\co_\fX$-champ de Higgs {\em quasi-petit} sur $\cN$ à coefficients dans 
$\txi^{-1}\tOmega^1_{\fX/\cS}$ \eqref{definf24}. 
Pour tout $1\leq i\leq c$, on pose $N_i=\Gamma(\oX_{i,s},\cN)$ et on note
\begin{equation}\label{mdpsa8c}
\theta_i \colon N_i\rightarrow \txi^{-1}\tOmega^1_{R/\co_K}\otimes_RN_i
\end{equation}
le $\hcR_i$-champ de Higgs induit par $\theta$. Le $\hcR_i$-module de Higgs $(N_i,\theta_i)$ est quasi-petit \eqref{pmh4},
et $N_i$ est $\co_C$-plat.  On note $\varphi_i$  la $\hcR_i$-représentation quasi-petite de $\Delta_i$ sur $N_i$ associée à $(N_i,\theta_i)$
par le foncteur \eqref{pmh5h}. Pour tout entier $n\geq 1$, on note $\cP_{n,i}$ le $\cR_i$-module $\mu_i(N_i/p^nN_i,\varphi_i)$ de $\oX^{\circ}_{i,\fet}$ \eqref{mdpsa12b}. 
Il existe alors un $R_1$-module $\cP_n$ de $\oX^{\circ}_\fet$, unique à isomorphisme unique près, tel que pour tout $1\leq i\leq c$, on ait $\cP_n|\oX^{\circ}_i=\cP_{n,i}$. 
D'après (\cite{agt} III.2.11), $\cP_n$ est de type fini sur $R_1$. Les $(\cP_{n+1})_{n\in \mN}$ forment naturellement un système projectif.
Pour tous entiers $n\geq m\geq 1$, le morphisme $\cP_n/p^m\cP_n\rightarrow\cP_m$ induit par le morphisme
de transition $\cP_n\rightarrow \cP_m$ est un isomorphisme. Le $R_1$-module $\bvcP=(\cP_{n+1})_{n\in \mN}$ de $(\oX^{\circ}_\fet)^{\mN^\circ}$ est adique de type fini 
d'après (\cite{agt} III.7.14). 
On obtient ainsi un foncteur \eqref{definf24}
\begin{equation}\label{mdpsa8d}
\begin{array}[t]{clcr}
\bMH^{\qpp}(\co_\fX,\txi^{-1}\tOmega^1_{\fX/\cS})&\rightarrow& \bMod^\atf((\oX^{\circ}_\fet)^{\mN^\circ},R_1)\\
(\cN,\theta)&\mapsto &\bvcP=(\cP_{n+1})_{n\in \mN}
\end{array}
\end{equation}

\begin{prop}\label{mdpsa9}
Soient $r, \varepsilon$ deux nombres rationnels tels que $r\geq 0$ et $\varepsilon>r+\frac{1}{p-1}$, 
$\cN$ un $\co_\fX$-module cohérent et $\cS$-plat, 
$\theta$ un $\co_\fX$-champ de Higgs $\varepsilon$-quasi-petit sur $\cN$ à coefficients dans $\txi^{-1}\tOmega^1_{\fX/\cS}$ \eqref{definf24}. 
Notons $\bvcP$ le $R_1$-module de $(\oX^{\circ}_\fet)^{\mN^\circ}$
associé à $(\cN,\theta)$ par le foncteur \eqref{mdpsa8d} et notons encore $\cN$ la $\co_\fX$-isogénie de Higgs $(\cN,\cN,\id,\theta)$
à coefficients dans $\txi^{-1}\tOmega^1_{\fX/\cS}$ \eqref{definf19}. 
Alors, il existe un isomorphisme fonctoriel de $\bIMC(\bvcC_X^{(r)}/\bvocB_X)$ \eqref{mdpsa4}
\begin{equation}\label{mdpsa9a}
\hupsigma^{(r)*}_X(\cN)\stackrel{\sim}{\rightarrow} \fS^{(r)}_X(\bvcP\otimes_{R_1}\bvocB_X).
\end{equation}
\end{prop}
Il suffit de calquer la preuve de (\cite{agt} III.13.8) qui correspond au cas absolu \eqref{definf10}, en utilisant dans le cas général \ref{pmh17} au lieu de (\cite{agt} II.13.17).

\begin{cor}\label{mdpsa10}
Les hypothèses étant celles de \ref{mdpsa9}, supposons, de plus, que le $\co_\fX[\frac 1 p]$-module $\cN[\frac 1 p]$ soit projectif de type fini \eqref{mdpsa31}.
Alors, 
\begin{itemize}
\item[{\rm (i)}] Le $\bvocB_{X,\mQ}$-module $(\bvcP\otimes_{R_1}\bvocB_X)_\mQ$ et le $\co_\fX[\frac 1 p]$-fibré de Higgs $\cN[\frac 1 p]$ 
sont associés dans le sens de \ref{mdpsa30}. 
\item[{\rm (ii)}] Le $\co_\fX[\frac 1 p]$-fibré de Higgs $\cN[\frac 1 p]$ est fortement soluble, 
le $\bvocB_\mQ$-module $(\bvbeta^*(\bvcP\otimes_{R_1}\bvocB_X))_\mQ$ est fortement de Dolbeault,
et on a un isomorphisme fonctoriel de $\bvocB_\mQ$-modules 
\begin{equation}\label{mdpsa10a}
\cV_\mQ(\cN[\frac 1 p])\stackrel{\sim}{\rightarrow} \bvbeta^*(\bvcP\otimes_{R_1}\bvocB_X)_\mQ,
\end{equation}
où $\cV_\mQ$ est le foncteur \eqref{aspglob19a}.
\end{itemize}
\end{cor}

(i) Cela résulte de \ref{mdpsa9}. 

(ii) Cela résulte de (i), \ref{mdpsa21} et \ref{aspglob18}(ii).  On notera que le $\bvocB$-module $\bvbeta^*(\bvcP\otimes_{R_1}\bvocB_X)$ est adique de type fini 
d'après \eqref{mdpsa8d} et (\cite{agt} III.7.5).

\begin{cor}\label{mdpsa15}
Pour qu'un $\co_\fX[\frac 1 p]$-fibré de Higgs à coefficients dans $\txi^{-1}\tOmega^1_{\fX/\cS}$ soit petit \eqref{definf26}, il faut et il suffit qu'il soit fortement soluble \eqref{aspglob1}. 
\end{cor}
Cela résulte de \ref{aspglob31} et \ref{mdpsa10}.

\begin{prop}\label{mdpsa18}
Reprenons les hypothèses de \ref{mdpsa8}, supposons, de plus, que le $\co_\fX[\frac 1 p]$-module $\cN[\frac 1 p]$ soit projectif de type fini \eqref{mdpsa31}. 
Alors, 
\begin{itemize}
\item[{\rm (i)}] Le morphisme d'adjonction
\begin{equation}\label{mdpsa18aa}
\bvcP\otimes_{R_1}\bvocB_X \rightarrow \bvbeta_*(\bvbeta^*(\bvcP\otimes_{R_1}\bvocB_X)), 
\end{equation}
où $\bvbeta$ désigne le morphisme de topos annelés \eqref{mdpsa13c}, est une isogénie. 
\item[{\rm (ii)}]  Pour tout entier $q\geq 1$, $\rR^q\bvbeta_*(\bvbeta^*(\bvcP\otimes_{R_1}\bvocB_X))$ est d'exposant fini \eqref{caip1}. 
\end{itemize}
\end{prop} 

En effet, en vertu de (\cite{agt} III.6.17 et \cite{egr1} 1.10.2(iii)), il existe un $\co_\fX$-module cohérent $\cN'$, deux entiers $t,m\geq 0$ et un $p^t$-isomorphisme $\co_\fX$-linéaire 
$u\colon \cN \oplus \cN'\rightarrow \co_\fX^m$ \eqref{notconv17}. 

(i) Compte tenu de (\cite{ag} 2.6.3 et \cite{agt} III.7.5), il suffit de montrer que pour tout entier $n\geq 0$, le morphisme canonique
\begin{equation}\label{mdpsa18b}
a\colon \cP_n\otimes_{R_1}\ocB_{X,n} \rightarrow \beta_*(\beta^{-1}(\cP_n)\otimes_{R_1}\ocB_n) 
\end{equation}
où $\beta$ est le morphisme de topos \eqref{ahttf2d}, est un  $p^{3t}$-isomorphisme, 
ou encore que pour tout $1\leq i\leq c$, la fibre $\nu_i(a|\oX_i^\circ)$ \eqref{mdpsa12a} est un  $p^{3t}$-isomorphisme

Pour tout objet $V$ de $\Et_{\rf/\oX_i^\circ}$ tel que le faisceau $\cP_n|V$ soit constant, le morphisme canonique
\begin{equation}
a_V\colon (\cP_n\otimes_{R_1}\ocB_{X,n})(V) \rightarrow (\beta^{-1}(\cP_n)\otimes_{R_1}\ocB_n)(V\rightarrow X)
\end{equation}
est un $p^{3t}$-isomorphisme. En effet, on se ramène par $u$ au cas où 
$\cP_n$ est le faisceau constant de valeur $\cR_i/p^n\cR_i$ \eqref{mdpsa8}, {\em i.e.}, correspondant à la représentation triviale de $\Delta_i$ sur ce $R_1$-module, 
auquel cas $a_V$ est un $\alpha$-isomorphisme en vertu de (\cite{ag} 4.6.29).  

Le morphisme de $R_1$-modules sous-jacent à $\nu_i(a|\oX_i^\circ)$ est la limite inductive des morphismes $a_{V_j}$, où $(V_j)_{j\in J_i}$ 
est le revêtement universel normalisé de $\oX^\circ_i$ en $\oy_i$ (\cite{agt} VI.9.8). 
La représentation $\nu_i(\cP_n)$ étant discrète, le sous-ensemble $J'_i$ de $J_i$ formé des éléments $j$ tels que 
$\cP_n|V_j$ soit constant, est cofinal dans $J_i$, d'où l'assertion recherchée.

(ii) En vertu de (\cite{ag} 4.6.30(ii)), pour tous entiers $1\leq i\leq c$ et $n\geq 0$, la limite inductive 
\begin{equation}
\underset{\underset{j\in J_i}{\longrightarrow}}{\lim}\ \rH^q((V_j\rightarrow X),\ocB_n),
\end{equation}
où $(V_j)_{j\in J_i}$ est défini plus haut, est $\alpha$-nul. On en déduit par $u$ que   
\begin{equation}
p^{3t}\underset{\underset{j\in J'_i}{\longrightarrow}}{\lim}\ \rH^q((V_j\rightarrow X),\beta^{-1}(\cP_n)\otimes_{R_1}\ocB_n))=0,
\end{equation}
où $J'_i$ est l'ensemble ordonné défini plus haut. Par suite, 
$\nu_i(\rR^q\beta_*(\beta^{-1}(\cP_n)\otimes_{R_1}\ocB_n)|\oX_i^\circ)$ est nul compte tenu de (\cite{sga4} V 5.1 et IV (6.3.3)), d'où la proposition (\cite{agt} III.7.5).

\begin{cor}\label{mdpsa16}
Sous les hypothèses de \ref{mdpsa10}, on a un isomorphisme fonctoriel de $\bvocB_{X,\mQ}$-modules
\begin{equation}
\bvbeta_{\mQ *}(\cV_\mQ(\cN[\frac 1 p]))\stackrel{\sim}{\rightarrow}(\bvcP\otimes_{R_1}\bvocB_X)_\mQ, 
\end{equation}
où $\bvbeta_\mQ$ désigne le morphisme de topos annelés \eqref{mdpsa13d}.
\end{cor}
Cela résulte de \ref{mdpsa10}(ii) et \ref{mdpsa18}(i). 

\begin{prop}\label{mdpsa14}
Soient $\cM$ un $\bvocB_\mQ$-module fortement de Dolbeault, $(\cN,\theta)$ un sous-$\co_\fX$-module de Higgs quasi-petit  
à coefficients dans $\txi^{-1}\tOmega^1_{\fX/\cS}$ de $\cH_\mQ(\cM)$ \eqref{definf24}
tel que le $\co_\fX$-module $\cN$ soit cohérent et $\cS$-plat et qu'il engendre $\cH_\mQ(\cM)$ sur $\co_\fX[\frac 1 p]$,
$\bvcP$ le $R_1$-module adique de type fini de $(\oX^{\circ}_\fet)^{\mN^\circ}$ associé à $(\cN,\theta)$ par le foncteur \eqref{mdpsa8d}. 
Alors, il existe des isomorphismes adjoints 
\begin{eqnarray}
\cM\stackrel{\sim}{\rightarrow}(\bvbeta^*(\bvcP\otimes_{R_1}\bvocB_X))_\mQ,\label{mdpsa14a}\\
\bvbeta_{\mQ *}(\cM)\stackrel{\sim}{\rightarrow}(\bvcP\otimes_{R_1}\bvocB_X)_\mQ.\label{mdpsa14b}
\end{eqnarray}
\end{prop}

En effet, d'après \ref{aspglob21} et \ref{mdpsa10}(ii), on a un isomorphisme de $\bvocB_\mQ$-modules
\begin{equation}
\cM\stackrel{\sim}{\rightarrow}(\bvbeta^*(\bvcP\otimes_{R_1}\bvocB_X))_\mQ.
\end{equation}
En vertu de \ref{mdpsa18}(i), on en déduit par adjonction un isomorphisme 
\begin{equation}
\bvbeta_{\mQ *}(\cM)\stackrel{\sim}{\rightarrow}(\bvcP\otimes_{R_1}\bvocB_X)_\mQ.
\end{equation}

\begin{rema}\label{mdpsa141}
D'après \ref{aspglob15} et \ref{mdpsa15}, pour tout $\bvocB_\mQ$-module fortement de Dolbeault $\cM$, 
il existe un sous-$\co_\fX$-module de Higgs quasi-petit $(\cN,\theta)$ de $\cH_\mQ(\cM)$ vérifiant les propriétés requises dans \ref{mdpsa14}.
\end{rema}

\begin{cor}\label{mdpsa140}
Soit $\cM$ un $\bvocB_\mQ$-module fortement de Dolbeault. Alors,
\begin{itemize}
\item[{\rm (i)}] Le $\bvocB_{X,\mQ}$-module $\bvbeta_{\mQ *}(\cM)$ est adique de type \eqref{mdpsa13} et il est associé au $\co_\fX[\frac 1 p]$-fibré de Higgs $\cH_\mQ(\cM)$
dans le sens de \ref{mdpsa30}.  
\item[{\rm (ii)}] Le morphisme d'adjonction $\bvbeta^*_\mQ(\bvbeta_{\mQ *}(\cM))\rightarrow \cM$ est un isomorphisme.
\item[{\rm (iii)}] Le morphisme d'adjonction $\bvbeta_{\mQ *}(\cM) \rightarrow \bvbeta_{\mQ *}(\bvbeta^*_\mQ(\bvbeta_{\mQ *}(\cM)))$ est un isomorphisme.
\end{itemize}
\end{cor}

Les assertions (i) et (ii) résultent de \ref{mdpsa14}, \ref{mdpsa141} et \ref{mdpsa10}(i).
L'assertion (iii) résulte de (ii) compte tenu du diagramme commutatif
\begin{equation}
\xymatrix{
&{\bvbeta_*}\\
{\bvbeta_*}\ar[r]_-(0.5){\adj}\ar[ru]^{\id}&{\bvbeta_*\bvbeta^*\bvbeta_*}\ar[u]_-(0.5){\bvbeta_*(\adj)}}
\end{equation}
où les flèches notées ``$\adj$'' désignent les morphismes d'adjonction.

\begin{defi}\label{mdpsa19}
On dit qu'un $\bvocB_{X,\mQ}$-module est {\em fortement de Dolbeault} s'il est isomorphisme à $\bvbeta_{\mQ *}(\cM)$ \eqref{mdpsa13d}, 
pour un $\bvocB_\mQ$-module fortement de Dolbeault $\cM$.
\end{defi}

Il résulte de \ref{mdpsa14} que tout  $\bvocB_{X,\mQ}$-module fortement de Dolbeault est adique de type fini et qu'il est associé à un $\co_\fX[\frac 1 p]$-fibré de Higgs fortement soluble,
dans le sens de \ref{mdpsa30}. 

On désigne par $\bMod_\mQ^\fDolb(\bvocB_X)$ la sous-catégorie pleine de $\bMod_\mQ(\bvocB_X)$  
formée des $\bvocB_{X,\mQ}$-modules fortement de Dolbeault. 

\begin{prop}\label{mdpsa17}
Les foncteurs \eqref{mdpsa13d} et \eqref{mdpsa13e}
\begin{equation}
\xymatrix{
{\bMod^\fDolb_{\mQ}(\bvocB_X)}\ar@<1ex>[r]^-(0.5){\bvbeta^*_{\mQ}}&{\bMod^\fDolb_{\mQ}(\bvocB)}
\ar@<1ex>[l]^-(0.5){\bvbeta_{\mQ*}}}
\end{equation}
induisent des équivalences de catégories quasi-inverses l'une de l'autre entre la catégorie des $\bvocB_{X,\mQ}$-modules fortement de Dolbeault
et celle des $\bvocB_\mQ$-modules fortement de Dolbeault.
\end{prop}

Cela résulte de \ref{mdpsa14}.

\subsection{}\label{mdpsa22}
Reprenons les notations de \ref{mdpsa12}. Pour tout entier $1\leq i\leq c$, posons
\begin{equation}\label{mdpsa22a}
\ocR_i=\nu_i(\ocB_X|\oX^{\circ}_i),
\end{equation}
et notons $\hocR_i$ son séparé complété $p$-adique. 
On munit les anneaux $\hocR_i$ et $\hocR_i[\frac 1 p]$ des topologies $p$-adiques \eqref{notconv15}.
On dit qu'une $\hocR_i$-représentation continue de $\Delta_i$ \eqref{notconv16} est {\em $p$-adique de type fini} si le $\hocR_i$-module sous-jacent 
est muni de la topologie $p$-adique et s'il est séparé de type fini. 
On notera que tout $\hocR_i$-module de type fini est complet pour la topologie $p$-adique (\cite{ac} chap. III § 2.11 cor.~1 de prop.~16). 
On désigne par $\bRep_{\ocR_i}^\disc(\Delta_i)$ la catégorie des $\ocR_i$-représentations discrètes de $\Delta_i$, 
par $\bRep_{\hocR_i}^{p\aatf}(\Delta_i)$ la sous-catégorie pleine de 
$\bRep_{\hocR_i}(\Delta_i)$ formée des $\hocR_i$-représentations 
$p$-adiques de type fini de $\Delta_i$ et par $\bRep_{\hocR_i[\frac 1 p]}^{\Dolb}(\Delta_i)$ la sous-catégorie pleine de 
$\bRep_{\hocR_i[\frac 1 p]}(\Delta_i)$ formée des $\hocR_i[\frac 1 p]$-représentations de Dolbeault de $\Delta_i$ \eqref{repdolb8}.

Le foncteur $\nu_i$ \eqref{mdpsa12a} induit un foncteur  
\begin{equation}\label{mdpsa22b}
\upnu_i\colon \bMod(\ocB_X)\rightarrow \bRep_{\ocR_i}^\disc(\Delta_i).
\end{equation}
On en déduit par passage à la limite projective un foncteur
\begin{equation}\label{mdpsa22c}
\hupnu_i\colon \bMod(\bvocB_X)\rightarrow \bRep_{\hocR_i}(\Delta_i). 
\end{equation}
Celui-ci induit un foncteur que l'on note encore
\begin{equation}\label{mdpsa22e}
\hupnu_i\colon \bMod_{\mQ}(\bvocB_X)\rightarrow \bRep_{\hocR_i[\frac 1 p]}(\Delta_i). 
\end{equation}

\begin{lem}\label{mdpsa23}
Avec les notations de \ref{mdpsa22}, les foncteurs $(\hupnu_i)_{1\leq i\leq c}$ \eqref{mdpsa22c} induisent une équivalence de catégories 
\begin{equation}\label{mdpsa23a}
\bMod^\atf(\bvocB_X)\stackrel{\sim}{\rightarrow} \prod_{1\leq i\leq c} \bRep_{\hocR_i}^{p\aatf}(\Delta_i).
\end{equation}
\end{lem}
En effet, notant 
\begin{equation}\label{mdpsa23b}
\uplambda\colon (\oX^\circ_\fet)^{\mN^\circ}\rightarrow \oX^\circ_\fet
\end{equation}
le morphisme canonique \eqref{notconv13a}, les $\uplambda^*(\oX^\circ_i)$ $(1\leq i\leq c)$ forment un recouvrement disjoint de l'objet final de $(\oX^\circ_\fet)^{\mN^\circ}$.
La proposition résulte donc de (\cite{agt} III.7.6 et III.7.21).

\begin{lem}\label{mdpsa24}
Soient $N$ un $\hRun$-module cohérent tel que le $\hRun[\frac 1 p]$-module $N[\frac 1 p]$ soit projectif, $A$ une $\hRun$-algèbre, $\hA$ son séparé complété $p$-adique. 
Alors, le morphisme $\hA$-linéaire canonique
\begin{equation}\label{mdpsa24a}
N\otimes_{\hRun}\hA\rightarrow  N\hotimes_{\hRun}A,
\end{equation}
où le produit tensoriel $\hotimes$ est complété pour la topologie $p$-adique, est une isogénie. 
\end{lem}

En effet, d'après (\cite{egr1} 1.10.2(iii)), il existe un $\hRun$-module cohérent $N'$, deux entiers $t,n\geq 0$ et un $p^t$-isomorphisme $u\colon N\oplus N'\rightarrow \hRun^n$. 
Considérons le diagramme commutatif canonique
\begin{equation}
\xymatrix{
{(N\oplus N')\otimes_{\hRun}\hA}\ar[r]^-(0.5)v\ar[d]_{u\otimes\id_\hA}&{(N\oplus N')\hotimes_{\hRun}A}\ar[d]^{u\hotimes \id_A}\\
{\hA^n}\ar@{=}[r]&{\hA^n}}
\end{equation}
Comme $u\otimes\id_\hA$ et $u\hotimes \id_A$ sont des $p^t$-isomorphismes, 
$v$ est un $p^{2t}$-isomorphisme, d'où la proposition (\cite{ag} 2.6.3).

\begin{lem}\label{mdpsa27}
Reprenons les notations de \ref{mdpsa22}. Soient $i$ un entier tel que $1\leq i\leq c$, $(N,\theta)$, $(N',\theta')$ deux $\hcR_i$-modules de Higgs 
à coefficients dans $\txi^{-1}\Omega^1_{R/\co_K}$ quasi-petits \eqref{pmh4}  tels que les $\hcR_i$-modules $N$ et $N'$ soient cohérents et $\co_C$-plats et 
que les $\hcR_i[\frac 1 p]$-modules $N[\frac 1 p]$ et $N'[\frac 1 p]$ soient projectifs. 
On note $\varphi$ (resp. $\varphi'$) la $\hcR_i$-représentation de $\Delta_i$ 
sur $N$ (resp. $N'$) associée à $(N,\theta)$ (resp. $(N',\theta')$) par le foncteur \eqref{pmh5h}. Alors, le morphisme canonique
\begin{eqnarray}\label{mdpsa27a}
\lefteqn{
\Hom_{\hocR_i[\Delta_i]}((N,\varphi)\hotimes_{\hcR_i}\hocR_i,(N',\varphi')\hotimes_{\hcR_i}\hocR_i)[\frac 1 p]}\\
&\rightarrow & \Hom_{\hocR_i[\frac 1 p][\Delta_i]}((N,\varphi)\hotimes_{\hcR_i}\hocR_i[\frac 1 p],(N',\varphi')\hotimes_{\hcR_i}\hocR_i[\frac 1 p])\nonumber
\end{eqnarray}
est un isomorphisme. 
\end{lem}

On notera d'abord que le $\hcR_i$-module $N$ étant de présentation finie, le morphisme canonique
\begin{equation}\label{mdpsa27b}
\Hom_{\hocR_i}(N\otimes_{\hcR_i}\hocR_i,N'\otimes_{\hcR_i}\hocR_i)[\frac 1 p]\rightarrow \Hom_{\hocR_i[\frac 1 p]}(N\otimes_{\hcR_i}\hocR_i[\frac 1 p],N'\otimes_{\hcR_i}\hocR_i[\frac 1 p])
\end{equation}
est un isomorphisme. Par suite, le morphisme canonique
\begin{eqnarray}\label{mdpsa27c}
\lefteqn{
j\colon \Hom_{\hocR_i[\Delta_i]}((N,\varphi)\otimes_{\hcR_i}\hocR_i,(N',\varphi')\otimes_{\hcR_i}\hocR_i)[\frac 1 p]}\\
&\rightarrow & \Hom_{\hocR_i[\frac 1 p][\Delta_i]}((N,\varphi)\otimes_{\hcR_i}\hocR_i[\frac 1 p],(N',\varphi')\otimes_{\hcR_i}\hocR_i[\frac 1 p])\nonumber
\end{eqnarray}
est injectif. 

De même, le morphisme canonique 
\begin{equation}
\Hom_{\hcR_i}(N,N')[\frac 1 p] \rightarrow \Hom_{\hcR_i[\frac 1 p]}(N[\frac 1 p],N'[\frac 1 p]),
\end{equation}
est un isomorphisme. Comme $N'$ est $\co_C$-plat, on en déduit que le morphisme canonique 
\begin{equation}
i\colon \Hom_{\bMH(\hcR_i)}((N,\theta),(N',\theta'))[\frac 1 p] \rightarrow \Hom_{\bMH(\hcR_i[\frac 1 p])}((N[\frac 1 p],\theta),(N'[\frac 1 p],\theta')),
\end{equation}
où on a abusivement noté $\theta$ (resp. $\theta'$) le champ de Higgs induit par $\theta$ (resp. $\theta'$) sur $N[\frac 1 p]$ (resp. $N'[\frac 1 p]$), est un isomorphisme. 

Considérons le diagramme commutatif
\[
\xymatrix{
{\Hom_{\bMH(\hcR_i)}((N,\theta),(N',\theta'))[\frac 1 p]}\ar[r]^-(0.5)u\ar[d]_i&{\Hom_{\hocR_i[\Delta_i]}((N,\varphi)\otimes_{\hcR_i}\hocR_i,(N',\varphi')\otimes_{\hcR_i}\hocR_i)[\frac 1 p]}\ar[d]^j\\
{\Hom_{\bMH(\hcR_i[\frac 1 p])}((N[\frac 1 p],\theta),(N'[\frac 1 p],\theta'))}\ar[r]^-(0.5)v&
{\Hom_{\hocR_i[\frac 1 p][\Delta_i]}((N,\varphi)\otimes_{\hcR_i}\hocR_i[\frac 1 p],(N',\varphi')\otimes_{\hcR_i}\hocR_i[\frac 1 p])}}
\]
où $u$ est induit par le foncteur \eqref{pmh5h} et $v$ par le foncteur \eqref{pmh29b}. 
En vertu de \ref{pmh30}, les $\hcR_i[\frac 1 p]$-modules de Higgs $(N[\frac 1 p],\theta)$ et $(N'[\frac 1 p],\theta')$ sont solubles et le morphisme $v$ s'identifie au morphisme 
\begin{equation}
\Hom_{\bMH(\hcR_i[\frac 1 p])}((N[\frac 1 p],\theta),(N'[\frac 1 p],\theta'))\rightarrow 
\Hom_{\hocR_i[\frac 1 p][\Delta_i]}(\mV(N[\frac 1 p],\theta)),\mV(N[\frac 1 p],\theta))
\end{equation}
induit par le foncteur $\mV$ \eqref{repdolb6b}. D'après \ref{repdolb14}, $v$ est donc un isomorphisme. Comme $i$ est un isomorphisme et $j$ est injectif, on en déduit
que $j$ est bijectif. Pour conclure la preuve de la proposition, il suffit d'observer que le morphisme canonique de $\hocR_i$-représentations de $\Delta_i$
\begin{equation}
(N,\varphi)\otimes_{\hcR_i}\hocR_i \rightarrow (N,\varphi)\hotimes_{\hcR_i}\hocR_i
\end{equation}
est une isogénie d'après \ref{mdpsa24} et (\cite{ag} 2.6.3), et de même pour $(N',\varphi')$.

\begin{lem}\label{mdpsa25}
Avec les notations de \ref{mdpsa22}, pour tout $\bvocB_{X,\mQ}$-module fortement de Dolbeault $\cM$ \eqref{mdpsa19} 
et pour tout entier $1\leq i\leq c$, $\hupnu_i(\cM)$ est une $\hocR_i[\frac 1 p]$-représentation de Dolbeault de $\Delta_i$ \eqref{repdolb8}.
\end{lem}

En effet, d'après \ref{mdpsa14} et \ref{mdpsa141}, 
il existe un $\co_\fX$-module cohérent et $\cS$-plat $\cN$ tel que le $\co_\fX[\frac 1 p]$-module $\cN[\frac 1 p]$ soit projectif de type fini \eqref{mdpsa31}, 
un $\co_\fX$-champ de Higgs quasi-petit $\theta$ sur $\cN$ à coefficients dans $\txi^{-1}\tOmega^1_{\fX/\cS}$ \eqref{definf24} et un isomorphisme 
\begin{equation}\label{mdpsa25a}
\cM\stackrel{\sim}{\rightarrow}(\bvcP\otimes_{R_1}\bvocB_X)_\mQ,
\end{equation}
où $\bvcP$ est le $R_1$-module adique de type fini de $(\oX^{\circ}_\fet)^{\mN^\circ}$ associé à $(\cN,\theta)$ par le foncteur \eqref{mdpsa8d}.
Posons $N_i=\Gamma(\oX_{i,s},\cN)$ qui est un $\hcR_i$-module cohérent, et notons $\varphi_i$ la $\hcR_i$-représentation quasi-petite de $\Delta_i$ sur $N_i$
associée dans \ref{mdpsa8} à $\theta$.
D'après \ref{mdpsa23}, comme le $\hcR_i[\frac 1 p]$-module $N_i[\frac 1 p]$ est projectif, 
l'isomorphisme \eqref{mdpsa25a} induit alors un isomorphisme de $\hocR_i[\frac 1 p]$-représentations de $\Delta_i$
\begin{equation}\label{mdpsa25b}
\hupnu_i(\cM) \stackrel{\sim}{\rightarrow} (N_i,\varphi_i)\otimes_{\hcR_i}\hocR_i[\frac 1 p]. 
\end{equation}
On en déduit que la $\hocR_i[\frac 1 p]$-représentation $\hupnu_i(\cM)$ de $\Delta_i$ est de Dolbeault en vertu de \ref{pmh30}.

\begin{lem}\label{mdpsa26}
Avec les notations de \ref{mdpsa22}, tout entier $1\leq i\leq c$ et toute $\hocR_i[\frac 1 p]$-représentation de Dolbeault $M$ de $\Delta_i$, 
il existe un $\bvocB_{X,\mQ}$-module fortement de Dolbeault $\cM$ tel que 
$\hupnu_i(\cM)=M$ et $\hupnu_j(\cM)=0$ pour tout $1\leq j\leq c$ avec $j\not=i$. 
\end{lem}

En effet, d'après \ref{repdolb14} et \ref{pmh32}, le $\hcR_i[\frac 1 p]$-module de Higgs $\mH(M)$ est petit. 
Il existe donc un sous-$\hcR_i$-module cohérent $N$ de $\mH(M)$ qui l'engendre sur $\hcR_i[\frac 1 p]$, 
tel que le $\hcR_i[\frac 1 p]$-champ de Higgs de $\mH(M)$
induit sur $N$ un $\hcR_i$-champ de Higgs quasi-petit $\theta$ à coefficients dans $\txi^{-1}\tOmega^1_{R/\co_K}$. 
On note $\varphi$  la $\hcR_i$-représentation quasi-petite de $\Delta_i$ sur $N$ associée à $(N,\theta)$ par le foncteur \eqref{pmh5h}. 
En vertu de \ref{repdolb14} et \ref{pmh30}(ii), on a un isomorphisme $\hocR_i$-linéaire et $\Delta_i$-équivariant 
\begin{equation}\label{mdpsa26a}
M\stackrel{\sim}{\rightarrow}(N,\varphi)\otimes_{\hcR_i}\hocR_i[\frac 1 p]. 
\end{equation}
Notons $\cN$ le $\co_\fX$-module cohérent tel que $\Gamma(\oX_{i,s},\cN)=N$ et $\Gamma(\oX_{j,s},\cN)=0$ pour tout $1\leq j\leq c$ avec $j\not=i$. 
On note encore $\theta$ le $\co_\fX$-champ de Higgs quasi-petit induit sur $\cN$ par le $\hcR_i$-champ de Higgs $\theta$ sur $N$. 
On désigne par $\bvcP$ le $R_1$-module adique de type fini de $(\oX^{\circ}_\fet)^{\mN^\circ}$ associé à $(\cN,\theta)$ par le foncteur \eqref{mdpsa8d}. 
Il résulte de \ref{mdpsa10} et \ref{mdpsa16} que le $\bvocB_{X,\mQ}$-module $\cM=(\bvcP\otimes_{R_1}\bvocB_X)_\mQ$ est fortement de Dolbeault. 
D'après la preuve de \ref{mdpsa25}, on a un isomorphisme de $\hocR_i[\frac 1 p]$-représentations de $\Delta_i$
\begin{equation}\label{mdpsa26b}
\hupnu_i(\cM) \stackrel{\sim}{\rightarrow} (N,\varphi)\otimes_{\hcR_i}\hocR_i[\frac 1 p],
\end{equation}
et $\hupnu_j(\cM)=0$ pour tout $1\leq j\leq c$ avec $j\not=i$, d'où la proposition.

\begin{prop}\label{mdpsa28}
Avec les notations de \ref{mdpsa22}, les foncteurs $(\hupnu_i)_{1\leq i\leq c}$ \eqref{mdpsa22e} induisent une équivalence de catégories 
\begin{equation}\label{mdpsa28a}
\bMod^\fDolb_\mQ(\bvocB_X)\stackrel{\sim}{\rightarrow} \prod_{1\leq i\leq c} \bRep_{\hocR_i[\frac 1 p]}^{\Dolb}(\Delta_i).
\end{equation} 
\end{prop}

En effet, le foncteur \eqref{mdpsa28a} est bien défini en vertu de \ref{mdpsa25}. Il est pleinement fidèle d'après \ref{mdpsa14}, \ref{mdpsa23} et \ref{mdpsa27}, 
et il est essentiellement surjectif d'après \ref{mdpsa26}.

\subsection{}\label{mdpsa29}
Reprenons les notations de \ref{mdpsa12}. Pour tout entier $1\leq i\leq c$, on note $\upbeta_i$ le foncteur composé 
\begin{equation}\label{mdpsa29a}
\upbeta_i \colon \tE\rightarrow \bB_{\Delta_i}, \ \ \ F\mapsto \nu_i\circ (\beta_*(F)|\oX^\circ_i), 
\end{equation}
où $\beta$ est le morphisme de topos \eqref{ahttf2d}. 
On définit ainsi un foncteur de la catégorie des faisceaux abéliens de $\tE$ dans celle des $\mZ[\Delta_i]$-modules. 
Celui-ci étant exact à gauche, on désigne par $\rR^q \upbeta_i$ $(q\geq 0)$ ses foncteurs dérivés à droite. 
Pour tout faisceau abélien $F$ de $\tE$ et tout entier $q\geq 0$, on a un isomorphisme canonique fonctoriel
\begin{equation}\label{mdpsa29b}
\rR^q\upbeta_i(F)\stackrel{\sim}{\rightarrow} \nu_i\circ (\rR^q\beta_*(F)|\oX^\circ_i).
\end{equation}
Oubliant l'action de $\Delta_i$, on a un isomorphisme canonique fonctoriel
\begin{equation}\label{mdpsa29c}
\rR^q\upbeta_i(F)\stackrel{\sim}{\rightarrow} \underset{\underset{j\in J_i}{\longrightarrow}}{\lim}\ \rH^q((V_j\rightarrow X),F),
\end{equation}
où $(V_j)_{j\in J_i}$ est le revêtement universel normalisé de $\oX^\circ_i$ en $\oy_i$  (\cite{ag} 2.1.20).  

Pour tout faisceau abélien $F=(F_n)_{n\geq 0}$ de $\tE^{\mN^\circ}$, on désigne par $\hupbeta_i(F)$ 
la limite projective des $\upbeta_i(F_n)$,
\begin{equation}\label{mdpsa29d}
\hupbeta_i(F)=\underset{\underset{n\geq 0}{\longleftarrow}}{\lim}\ \upbeta_i(F_n).
\end{equation}
On définit ainsi un foncteur de la catégorie des faisceaux abéliens de $\tE^{\mN^\circ}$ dans celle des $\mZ[\Delta_i]$-modules.
Celui-ci étant exact à gauche, on désigne abusivement par $\rR^q\hupbeta_i(F)$ ($q\geq 0$) ses foncteurs dérivés à droite. 
D'après (\cite{jannsen} 1.6), on a une suite exacte canonique 
\begin{equation}\label{mdpsa29e}
0\rightarrow \rR^1 \underset{\underset{n\geq 0}{\longleftarrow}}{\lim}\ \rR^{q-1}\upbeta_i(F_n)\rightarrow
\rR^q\hupbeta_i(F)\rightarrow 
\underset{\underset{n\geq 0}{\longleftarrow}}{\lim}\ \rR^q\upbeta_i(F_n)\rightarrow 0,
\end{equation}
où on a posé $\rR^{-1}\upbeta_i(F_n)=0$ pour tout $n\geq 0$. 

Pour tout entier $q\geq 0$, le foncteur $\rR^q\hupbeta_i$ induit un foncteur qu'on note aussi 
\begin{equation}\label{mdpsa29f}
\rR^q\hupbeta_i\colon \bMod(\bvocB)\rightarrow \bRep_{\hocR_i}(\Delta_i),
\end{equation}
dans la catégorie des $\hocR_i$-représentations de $\Delta_i$, où $\hocR_i$ est le séparé complété $p$-adique de l'anneau $\ocR_i$ \eqref{mdpsa22a}. 
Celui-ci induit un foncteur qu'on note aussi 
\begin{equation}\label{mdpsa29g}
\rR^q\hupbeta_i\colon \bMod_\mQ(\bvocB)\rightarrow \bRep_{\hocR_i[\frac 1 p]}(\Delta_i). 
\end{equation}

Il résulte facilement de \eqref{mdpsa29d}, \eqref{mdpsa29e} et (\cite{sga4} V 3.5) que les foncteurs $\rR^q\hupbeta_i$ $(q\geq 0)$ \eqref{mdpsa29f}
sont les foncteurs dérivés à droite du foncteur $\hupbeta_i=\rR^0\hupbeta_i$ sur la catégorie des $\bvocB$-modules. 
Par suite, les foncteurs $\rR^q\hupbeta_i$ $(q\geq 0)$ \eqref{mdpsa29g}
sont les foncteurs dérivés à droite du foncteur $\hupbeta_i=\rR^0\hupbeta_i$ sur la catégorie des $\bvocB_\mQ$-modules \eqref{indsh14}.

\begin{teo}\label{mdpsa32}
Conservons les notations de \ref{mdpsa29}. Alors,
\begin{itemize}
\item[{\rm (i)}] Les foncteurs $(\hupbeta_i)_{1\leq i\leq c}$ \eqref{mdpsa29g} induisent une équivalence de catégories 
\begin{equation}\label{mdpsa32a}
\bMod^\fDolb_\mQ(\bvocB)\stackrel{\sim}{\rightarrow} \prod_{1\leq i\leq c} \bRep_{\hocR_i[\frac 1 p]}^{\Dolb}(\Delta_i).
\end{equation} 
\item[{\rm (ii)}] Pour tout $\bvocB_\mQ$-module fortement de Dolbeault $\cM$ et tout entier $q\geq 1$, on a 
\begin{equation}
\rR^q\hupbeta_i(\cM)=0.
\end{equation}
\end{itemize}
\end{teo}

(i) En effet, pour tout $1\leq i\leq c$, on a $\hupbeta_i=\hupnu_i\circ \bvbeta_*$, où $\hupbeta_i=\rR^0\hupbeta_i$ est le foncteur \eqref{mdpsa29f}, 
$\hupnu_i$ est le foncteur \eqref{mdpsa22c} et $\bvbeta$ est le morphisme de topos annelés \eqref{mdpsa13c}. 
La proposition résulte alors de \ref{mdpsa17} et \ref{mdpsa28}. 

(ii) En effet, d'après \ref{mdpsa14}, il existe un $\bvocB_X$-module adique de type fini $\bvcF$ et deux isomorphismes adjoints
\begin{eqnarray}
\cM\stackrel{\sim}{\rightarrow}(\bvbeta^*(\bvcF))_\mQ,\\
\bvbeta_{\mQ *}(\cM)\stackrel{\sim}{\rightarrow}\bvcF_\mQ.
\end{eqnarray}
De plus, en vertu de \ref{mdpsa18}(ii), le $\bvocB_X$-module $\rR^q\bvbeta_*(\bvbeta^*(\bvcF))$ est d'exposant fini. 
La proposition résulte alors de \eqref{mdpsa29e} appliqué au faisceau $\bvbeta^*(\bvcF)$, \eqref{mdpsa29b}, (\cite{agt} III.7.5) et (\cite{jannsen} 1.15).

\begin{prop}\label{mdpsa33}
Conservons les notations de \ref{mdpsa29}. Alors,
\begin{itemize}
\item[{\rm (i)}] Pour tout $\co_\fX[\frac 1 p]$-fibré fortement soluble $\cN$ à coefficients dans $\txi^{-1}\tOmega^1_{\fX/\cS}$ \eqref{aspglob1}, 
et tout entier $1\leq i\leq c$, le $\hcR_i[\frac 1 p]$-module de Higgs $N_i=\Gamma(\oX_{i,s},\cN)$ à coefficients dans $\txi^{-1}\tOmega^1_{\fX/\cS}(\oX_{i,s})$ 
est soluble dans le sens de \ref{repdolb10}. 
\item[{\rm (ii)}] Pour tout entier $1\leq i\leq c$, le diagramme de foncteurs 
\begin{equation}\label{mdpsa33a}
\xymatrix{
{\bMH^{\fsol}(\co_\fX[\frac 1 p],\txi^{-1}\tOmega^1_{\fX/\cS})}\ar[d]_{\Gamma(\oX_{i,s},-)}\ar[r]^-(0.5){\cV_\mQ}&{\bMod^\fDolb_\mQ(\bvocB)}\ar[d]^{\hupbeta_i}\\
{\bMH^\sol(\hcR_i[\frac 1 p],\txi^{-1}\tOmega^1_{\fX/\cS}(\oX_{i,s}))}\ar[r]\ar[r]^-(0.5){\mV}&{\bRep_{\hocR_i[\frac 1 p]}^{\Dolb}(\Delta_i)}}
\end{equation}
où $\cV_\mQ$ est le foncteur \eqref{aspglob19b} et $\mV$ est le foncteur \eqref{repdolb6b}, est commutatif à isomorphisme canonique près.
\item[{\rm (iii)}] Pour tout entier $1\leq i\leq c$, le diagramme de foncteurs 
\begin{equation}\label{mdpsa33b}
\xymatrix{
{\bMod^\fDolb_\mQ(\bvocB)}\ar[r]^-(0.5){\cH_\mQ}\ar[d]_{\hupbeta_i}&{\bMH^{\fsol}(\co_\fX[\frac 1 p],\txi^{-1}\tOmega^1_{\fX/\cS})}\ar[d]^{\Gamma(\oX_{i,s},-)}\\
{\bRep_{\hocR_i[\frac 1 p]}^{\Dolb}(\Delta_i)}\ar[r]\ar[r]^-(0.5){\mH}&{\bMH^\sol(\hcR_i[\frac 1 p],\txi^{-1}\tOmega^1_{\fX/\cS}(\oX_{i,s}))}}
\end{equation}
où $\cH_\mQ$ est le foncteur \eqref{aspglob15b} et $\mH$ est le foncteur \eqref{repdolb5b}, est commutatif à isomorphisme canonique près.
\end{itemize}
\end{prop}

(i) En effet, comme $\cN$ est petit en vertu de \ref{aspglob31}, $N_i$ est petit d'après \ref{definf260} et donc soluble en vertu de \ref{pmh32}. 

(ii)  Cela résulte de \ref{pmh30}(ii), \eqref{mdpsa10a} et \eqref{mdpsa14b}, en tenant compte de \ref{mdpsa15} et \ref{definf29}. 
 
(iii) Cela résulte de (ii), \ref{repdolb14} et \ref{aspglob210}.

\section{Image inverse d'un ind-module de Dolbeault par un morphisme étale}

\subsection{}\label{pchamp1}
Soit $g\colon X'\rightarrow X$ un morphisme étale de type fini. On munit $X'$ de la structure logarithmique $\cM_{X'}$ 
image inverse de $\cM_X$ et on note $f'\colon (X',\cM_{X'})\rightarrow (S,\cM_S)$ le morphisme induit par $f$ et $g$. 
On observera que $f'$ est adéquat (\cite{agt} III.4.7) et que $X'^\circ=X^\circ\times_XX'$ 
est le sous-schéma ouvert maximal de $X'$ où la structure logarithmique $\cM_{X'}$ est triviale.
On munit $\oX'$ et $\coX'$ \eqref{defing1c} des structures logarithmiques $\cM_{\oX'}$ et $\cM_{\coX'}$ 
images inverses de $\cM_{X'}$. 
Il existe essentiellement un unique morphisme étale $\tg\colon \tX'\rightarrow \tX$ 
qui s'insère dans un diagramme cartésien \eqref{defing12}
\begin{equation}\label{pchamp1a}
\xymatrix{
{\coX'}\ar[r]\ar[d]_{\cog}&{\tX'}\ar[d]^{\tg}\\
{\coX}\ar[r]&{\tX}}
\end{equation}
On munit $\tX'$ de la structure logarithmique $\cM_{\tX'}$ image inverse de $\cM_{\tX}$, de sorte que 
$(\tX',\cM_{\tX'})$ est une $(\tS,\cM_{\tS})$-déformation lisse de $(\coX',\cM_{\coX'})$.

On associe à $(f',\tX',\cM_{\tX'})$ des objets analogues à ceux définis dans \ref{ahttf}--\ref{aspglob} pour $(f,\tX,\cM_{\tX})$, 
qu'on note par les mêmes symboles affectés d'un exposant $^\prime$.  
On désigne par
\begin{eqnarray}
\Phi\colon \tE'&\rightarrow& \tE,\label{pchamp1b}\\
\Phi_s\colon \tE'_s&\rightarrow& \tE_s,\label{pchamp1c}
\end{eqnarray}
les morphismes de topos induits par fonctorialité par $g$ (\cite{agt} (III.8.5.3) et (III.9.118)).
D'après (\cite{agt} VI.10.14), $\Phi$ s'identifie au morphisme de localisation de $\tE$ en $\sigma^*(X')$. 
De plus, on a un homomorphisme canonique $\Phi^*(\ocB)\rightarrow \ocB'$ (\cite{agt} (III.8.20.6)), 
qui est un isomorphisme en vertu de (\cite{agt} III.8.21(i)). 
Pour tout entier $n\geq 1$, $\Phi_s$ est sous-jacent à un morphisme canonique de topos annelés (\cite{agt} (III.9.11.11)) 
\begin{equation}\label{pchamp1d}
\Phi_n\colon (\tE'_s,\ocB'_n)\rightarrow (\tE_s,\ocB_n).
\end{equation}
L'homomorphisme $\Phi_s^*(\ocB_n)\rightarrow \ocB'_n$ étant un isomorphisme d'après (\cite{agt} III.9.13),
il n'y a pas de différence pour les $\ocB_n$-modules entre l'image
inverse par $\Phi_s$ au sens des faisceaux abéliens et l'image inverse par $\Phi_n$ au sens des modules.
Le diagramme de morphismes de topos annelés (\cite{agt} (III.9.11.12))
\begin{equation}\label{pchamp1e}
\xymatrix{
{(\tE'_s,\ocB'_n)}\ar[r]^{\Phi_n}\ar[d]_{\sigma'_n}&{(\tE_s,\ocB_n)}\ar[d]^{\sigma_n}\\
{(X'_{s,\et},\co_{\oX'_n})}\ar[r]^{\ogg_n}&{(X_{s,\et},\co_{\oX_n})}}
\end{equation}
où $\ogg_n$ est le morphisme induit par $g$, est commutatif à isomorphisme canonique près. 

Pour tout nombre rationnel $r\geq 0$, on a un isomorphisme $\ocB'_n$-linéaire canonique \eqref{ahttf41j}
\begin{equation}\label{pchamp1f}
\upnu^{(r)}_n\colon \Phi^*_n(\cF^{(r)}_n)\stackrel{\sim}{\rightarrow}\cF'^{(r)}_n,
\end{equation}
où $\cF^{(r)}_n$ et $\cF'^{(r)}_n$ sont les modules définis dans \eqref{ahttf37b}. 
On a un isomorphisme canonique \eqref{ahttf12c}
\begin{equation}\label{pchamp1g}
\sigma^*_n(\txi^{-1}\tOmega^1_{\oX_n/\oS_n})\stackrel{\sim}{\rightarrow}
\sigma^{-1}(\txi^{-1}\tOmega^1_{X/S})\otimes_{\sigma^{-1}(\co_X)}\ocB_n.
\end{equation}
Donc en vertu de (\cite{agt} VI.5.34(ii), VI.8.9 et VI.5.17), $\sigma^*_n(\xi^{-1}\tOmega^1_{\oX_n/\oS_n})$ est 
le faisceau de $\tE$ associé au préfaisceau sur $E$ défini par la correspondance
\begin{equation}\label{pchamp1h}
\{U\in \Et^\circ_{/X} \mapsto \xi^{-1}\tOmega^1_{X/S}(U)\otimes_{\co_X(U)}\ocB_{U,n}\}.
\end{equation}
On en déduit que le diagramme 
\begin{equation}\label{pchamp1i}
\xymatrix{
0\ar[r]&{\ocB'_n}\ar@{=}[d]\ar[r]&{\Phi^*_n(\cF^{(r)}_n)}\ar[d]_{\upnu^{(r)}_n}\ar[r]&
{\Phi^*_n(\sigma^*_n(\txi^{-1}\tOmega^1_{\oX_n/\oS_n}))}\ar[r]\ar[d]&0\\
0\ar[r]&{\ocB'_n}\ar[r]&{\cF'^{(r)}_n}\ar[r]&{\sigma'^*_n(\txi^{-1}\tOmega^1_{\oX'_n/\oS_n})}\ar[r]&0}
\end{equation}
où les lignes horizontales sont induites par les suites exactes \eqref{ahttf38a} et la flèche verticale de droite est l'isomorphisme induit par l'isomorphisme canonique
\begin{equation}\label{pchamp1j}
\ogg_n^*(\txi^{-1}\tOmega^1_{\oX_n/\oS_n})\stackrel{\sim}{\rightarrow} \txi^{-1}\tOmega^1_{\oX'_n/\oS_n}
\end{equation}
et le diagramme commutatif \eqref{pchamp1e}, est commutatif.

\subsection{}\label{pchamp15}
On désigne par $\fX'$ le schéma formel complété $p$-adique de $\oX'$, par
\begin{equation}\label{pchamp15a}
\fgg\colon \fX'\rightarrow \fX
\end{equation} 
le prolongement de $\ogg\colon \oX'\rightarrow \oX$ aux complétés et par
\begin{equation}\label{pchamp15b}
\bvPhi\colon (\tE'^{\mN^\circ}_s,\bvocB')\rightarrow (\tE^{\mN^\circ}_s,\bvocB)
\end{equation} 
le morphisme de topos annelés induit par les morphismes $(\Phi_n)_{n\geq 1}$ \eqref{pchamp1d} (cf. \cite{agt} III.7.5).  
D'après (\cite{agt} III.9.14), $\bvPhi$ est canoniquement isomorphe 
au morphisme de localisation du topos annelé $(\tE_s^{\mN^\circ},\bvocB)$ en $\uplambda^*(\sigma^*_s(X'_s))$,
où $\uplambda\colon \tE_s^{\mN^\circ}\rightarrow \tE_s$ est le morphisme canonique de topos défini dans \eqref{notconv13a}.
Par suite, il n'y a pas de différence pour les $\bvocB$-modules entre l'image
inverse par $\bvPhi$ au sens des faisceaux abéliens et l'image inverse au sens des modules.

Il résulte aussitôt de \eqref{pchamp1e} que diagramme de morphismes de topos annelés 
\begin{equation}\label{pchamp15g}
\xymatrix{
{(\tE'^{\mN^\circ}_s,\bvocB')}\ar[r]^{\bvPhi}\ar[d]_{\hupsigma'}&{(\tE^{\mN^\circ}_s,\bvocB)}\ar[d]^\hupsigma\\
{(X'_{s,\zar},\co_{\fX'})}\ar[r]^{\fgg}&{(X_{s,\zar},\co_{\fX})}}
\end{equation}
où $\hupsigma$ et $\hupsigma'$ sont les morphismes définis dans \eqref{ahttf13e},
est commutatif à isomorphisme canonique près.

Pour tout nombre rationnel $r\geq 0$, les isomorphismes \eqref{pchamp1f} induisent un isomorphisme 
\begin{equation}\label{pchamp15i}
\bvupnu^{(r)}\colon \bvPhi^*(\bvcF^{(r)})\stackrel{\sim}{\rightarrow}\bvcF'^{(r)},
\end{equation}
où $\bvcF^{(r)}$ et $\bvcF'^{(r)}$ sont les modules définis dans \eqref{ahttf14b}.  
Les diagrammes commutatifs \eqref{pchamp1i} induisent un diagramme commutatif
\begin{equation}\label{pchamp15j}
\xymatrix{
0\ar[r]&{\bvPhi^*(\bvocB)}\ar[r]\ar@{=}[d]&{\Phi^*(\bvcF^{(r)})}\ar[r]
\ar[d]^{\bvupnu^{(r)}}&{\bvPhi^*(\hupsigma^*(\txi^{-1}\tOmega^1_{\oX/\oS}))}\ar[r]\ar[d]^{\updelta}&0\\
0\ar[r]&{\bvocB'}\ar[r]&{\bvcF'^{(r)}}\ar[r]&{\hupsigma'^*(\txi^{-1}\tOmega^1_{\oX'/\oS})}\ar[r]&0}
\end{equation}
où $\delta$ est l'isomorphisme induit par l'isomorphisme canonique 
\begin{equation}\label{pchamp15e}
\fgg^*(\txi^{-1}\tOmega^1_{\fX/\cS})\stackrel{\sim}{\rightarrow} \txi^{-1}\tOmega^1_{\fX'/\cS}
\end{equation} 
et le diagramme commutatif \eqref{pchamp15g}. 

Compte tenu de \eqref{ahttf14d}, l'isomorphisme \eqref{pchamp15i} induit un isomorphisme de $\bvocB'$-algèbres
\begin{equation}\label{pchamp15k}
\bvupmu^{(r)}\colon \bvPhi^*(\bvcC^{(r)})\stackrel{\sim}{\rightarrow}\bvcC'^{(r)}.
\end{equation}
Il résulte aussitôt de \eqref{pchamp15j} que le diagramme 
\begin{equation}\label{pchamp15m}
\xymatrix{
{\bvPhi^*(\bvcC^{(r)})}\ar[r]^-(0.5){\bvupmu^{(r)}}\ar[d]_-(0.5){\bvPhi^*(d_{\bvcC^(r)})}&{\bvcC'^{(r)}}\ar[d]^{d_{\bvcC'(r)}}\\
{\bvPhi^*(\hupsigma^*(\txi^{-1}\tOmega^1_{\fX/\cS})\otimes_{\bvocB}\bvcC^{(r)})}\ar[r]^-(0.5){\updelta\otimes \bvupnu^{(r)}}&
{\hupsigma'^*(\txi^{-1}\tOmega^1_{\fX'/\cS})\otimes_{\bvocB'}\bvcC'^{(r)}}}
\end{equation}
où $d_{\bvcC^{(r)}}$ et $d_{\bvcC'^{(r)}}$ sont les dérivations \eqref{indmdlb1a}, est commutatif. 
Pour tous nombres rationnels $r\geq r'\geq 0$, le diagramme 
\begin{equation}\label{pchamp15l}
\xymatrix{
{\bvPhi^*(\bvcC^{(r)})}\ar[r]^-(0.5){\bvupnu^{(r)}}\ar[d]_{\bvPhi^*(\bvalpha^{r,r'})}&{\bvcC'^{(r)}}\ar[d]^-(0.5){\bvalpha'^{r,r'}}\\
{\bvPhi^*(\bvcC^{(r')})}\ar[r]^-(0.5){\bvupnu^{(r')}}&{\bvcC'^{(r')}}}
\end{equation}
où $\bvalpha^{r,r'}$ et  $\bvalpha'^{r,r'}$ sont les homomorphismes canoniques \eqref{ahttf14f}, est commutatif.

\subsection{}\label{pchamp16}
Compte tenu de l'isomorphisme \eqref{pchamp15e}, le foncteur $\fgg^*$ induit des foncteurs 
\begin{eqnarray}
\fgg^*\colon \bMH(\co_\fX,\txi^{-1}\tOmega^1_{\fX/\cS})&\rightarrow& \bMH(\co_{\fX'},\txi^{-1}\tOmega^1_{\fX'/\cS}),\label{pchamp16a}\\
\fgg^*\colon \bIH(\co_\fX,\txi^{-1}\tOmega^1_{\fX/\cS})&\rightarrow& \bIH(\co_{\fX'},\txi^{-1}\tOmega^1_{\fX'/\cS}),\label{pchamp16b}\\
\rI\fgg^*\colon \bIndMH(\co_\fX,\txi^{-1}\tOmega^1_{\fX/\cS})&\rightarrow&\bIndMH(\co_{\fX'},\txi^{-1}\tOmega^1_{\fX'/\cS}).\label{pchamp16c}
\end{eqnarray}

Le foncteur $\bvPhi^*$ \eqref{pchamp15b} induit un foncteur \eqref{indsh21a}
\begin{equation} \label{pchamp16d}
\rI\bvPhi^*\colon \bIndMod(\bvocB)\rightarrow \bIndMod(\bvocB').
\end{equation}

Soit $r$ un nombre rationnel $\geq 0$. D'après \ref{indsh44}, le foncteur $\bvPhi^*$ \eqref{pchamp15b} et l'isomorphisme \eqref{pchamp15k} 
induisent un foncteur que l'on note encore \eqref{indmdlb1}
\begin{equation}\label{pchamp16e}
\rI\bvPhi^*\colon \bIndMC(\bvcC^{(r)}/\bvocB)\rightarrow \bIndMC(\bvcC'^{(r)}/\bvocB').
\end{equation}
Les diagrammes de foncteurs
\begin{equation}\label{pchamp16f}
\xymatrix{
{\bIndMod(\bvocB)}\ar[r]^-(0.5){\rI\fS^{(r)}}\ar[d]_{\rI\bvPhi^*}&{\bIndMC(\bvcC^{(r)}/\bvocB)}\ar[d]^{\rI\bvPhi^*}\\
{\bIndMod(\bvocB')}\ar[r]^-(0.5){\rI\fS'^{(r)}}&{\bIndMC(\bvcC'^{(r)}/\bvocB')}}
\end{equation}
où les flèches horizontales sont les foncteurs \eqref{indmdlb1b}, et 
\begin{equation}\label{pchamp16g}
\xymatrix{
{\bIndMH(\co_\fX,\txi^{-1}\tOmega^1_{\fX/\cS})}\ar[r]^-(0.5){\rI\hupsigma^{(r)*}}\ar[d]_{\rI\fgg^*}&{\bIndMC(\bvcC^{(r)}/\bvocB)}\ar[d]^{\rI\bvPhi^*}\\
{\bIndMH(\co_{\fX'},\txi^{-1}\tOmega^1_{\fX'/\cS})}\ar[r]^-(0.5){\rI\hupsigma'^{(r)*}}&{\bIndMC(\bvcC'^{(r)}/\bvocB')}}
\end{equation}
où les flèches horizontales sont les foncteurs \eqref{indmdlb1f}, 
sont clairement commutatifs à isomorphismes canoniques près.  

D'après (\cite{agt} III.9.14) et \ref{indsh36}, le diagramme de foncteurs
\begin{equation}\label{pchamp16h}
\xymatrix{
{\bIndMC(\bvcC^{(r)}/\bvocB)}\ar[r]^-(0.5){\rI\cK^{(r)}}\ar[d]_{\rI\bvPhi^*}&{\bMod(\bvocB)}\ar[d]^{\rI\bvPhi^*}\\
{\bIndMC(\bvcC'^{(r)}/\bvocB')}\ar[r]^-(0.5){\rI\cK'^{(r)}}&{\bMod(\bvocB')}}
\end{equation}
où les flèches horizontales sont les foncteurs \eqref{indmdlb1c}, est commutatif à isomorphisme canonique près. 

Le morphisme de changement de base \eqref{indsh42d} relatif au diagramme \eqref{pchamp15g} induit un morphisme de foncteurs
de $\bIndMC(\bvcC^{(r)}/\bvocB)$ dans $\bIndMH(\co_{\fX'},\txi^{-1}\tOmega^1_{\fX'/\cS})$ 
\begin{equation}\label{pchamp16i}
\rI\fgg^*\circ \rI\hupsigma^{(r)}_*\rightarrow \rI\hupsigma'^{(r)}_*\circ \rI\bvPhi^*,
\end{equation}
où $\rI\hupsigma^{(r)}_*$ et $\rI\hupsigma'^{(r)}_*$ sont les foncteurs \eqref{indmdlb1g}. 
Comme le foncteur $\fgg^*$ commute aux limites inductives, le foncteur $\rI\fgg^*$ 
commute au foncteur $\kappa_{\co_\fX}$ \eqref{indmdlb1h}. On en déduit  un morphisme de foncteurs de $\bIndMC(\bvcC^{(r)}/\bvocB)$ dans 
$\bMH(\co_{\fX'},\txi^{-1}\tOmega^1_{\fX'/\cS})$
\begin{equation}\label{pchamp16ij}
\fgg^*\circ \vupsigma^{(r)}_*\rightarrow \vupsigma'^{(r)}_*\circ \rI\bvPhi^*,
\end{equation}
où $\vupsigma^{(r)}_*$ et $\vupsigma'^{(r)}_*$ sont les foncteurs \eqref{indmdlb1i}.

D'après (\cite{sga4} XVII 2.1.3), le morphisme \eqref{pchamp16i} est l'adjoint du morphisme composé
\begin{equation}\label{pchamp16j}
\rI\hupsigma'^{(r)*}\circ \rI\fgg^*\circ \rI\hupsigma^{(r)}_*\stackrel{\sim}{\rightarrow} \rI\bvPhi^*\circ \rI\hupsigma^{(r)*} \circ \rI\hupsigma^{(r)}_* \rightarrow \rI\bvPhi^*,
\end{equation} 
où la première flèche est l'isomorphisme sous-jacent au diagramme \eqref{pchamp16g} 
et la seconde flèche est le morphisme d'adjonction. 
Par suite, pour tout objet $\cN$ de $\bIndMH(\co_\fX,\txi^{-1}\tOmega^1_{\fX/\cS})$ et tout objet $\cF$ de $\bIndMC(\bvcC^{(r)}/\bvocB)$, 
le diagramme d'applications d'ensembles
\begin{equation}\label{pchamp16k}
\text{\tiny {\xymatrix{
{\Hom_{\bIndMC(\bvcC^{(r)}/\bvocB)}(\rI\hupsigma^{(r)*}(\cN),\cF)}\ar[d]\ar[r]^-(0.5)a&
{\Hom_{\bIndMC(\bvcC'^{(r)}/\bvocB')}(\rI\hupsigma'^{(r)*}(\rI\fgg^*(\cN)),\rI\bvPhi^*(\cF))}\ar[d]\\
{\Hom_{\bIndMH(\co_\fX,\txi^{-1}\tOmega^1_{\fX/\cS})}(\cN,\rI\hupsigma^{(r)}_*(\cF))}\ar[r]^-(0.5)b&
{\Hom_{\bIndMH(\co_{\fX'},\txi^{-1}\tOmega^1_{\fX'/\cS})}(\rI\fgg^*(\cN),\rI\hupsigma'^{(r)}_*(\rI\bvPhi^*(\cF)))}}}}
\end{equation}
où les flèches verticales sont les isomorphismes d'adjonction, $a$ est induit par 
le foncteur $\rI\bvPhi^*$ et l'isomorphisme sous-jacent au diagramme \eqref{pchamp16g}, 
et $b$ est induit par le foncteur $\rI\fgg^*$ et le morphisme \eqref{pchamp16i}, est commutatif. 

Pour tous nombres rationnels $r\geq r'\geq 0$, le diagramme de foncteurs 
\begin{equation}\label{pchamp16l}
\xymatrix{
{\bIndMC(\bvcC^{(r)}/\bvocB)}\ar[r]^-(0.5){\rI\varepsilon^{r,r'}}\ar[d]_{\rI\bvPhi^*}&{\bIndMC(\bvcC^{(r')}/\bvocB)}\ar[d]^{\rI\bvPhi^*}\\
{\bIndMC(\bvcC'^{(r)}/\bvocB')}\ar[r]^-(0.5){\rI\varepsilon'^{r,r'}}&{\bIndMC(\bvcC'^{(r')}/\bvocB')}}
\end{equation}
où les flèches horizontales sont les foncteurs \eqref{indmdlb3a}, est  
commutatif à isomorphisme canonique près.
Il résulte aussitôt de \eqref{pchamp15l} que le diagramme de morphismes de foncteurs 
\begin{equation}\label{pchamp16m}
\xymatrix{
{\rI\fgg^*\circ \rI\hupsigma^{(r)}_*}\ar[rr]\ar[d]&&{\rI\fgg^*\circ \rI\hupsigma^{(r')}_*\circ \rI\varepsilon^{r,r'}}\ar[d]\\
{\rI\hupsigma'^{(r)}_*\circ \rI\bvPhi^*}\ar[r]&{\rI\hupsigma'^{(r')}_*\circ \rI\varepsilon^{r,r'}\circ \rI\bvPhi^*}\ar@{=}[r]&
{\rI\hupsigma'^{(r')}_*\circ \rI\bvPhi^*\circ \rI\varepsilon^{r,r'}}}
\end{equation}
où  les flèches horizontales sont induites par le morphisme \eqref{indmdlb3f}, 
les flèches verticales sont induites par le morphisme \eqref{pchamp16i}
et l'identification notée avec un symbole $=$ provient du diagramme \eqref{pchamp16l}, est commutatif.
Par suite, le morphisme composé 
\begin{equation}\label{pchamp16n}
\rI\fgg^*\circ \rI\hupsigma^{(r)}_*\circ \rI\fS^{(r)}\rightarrow \rI\hupsigma'^{(r)}_*\circ \rI\bvPhi^*\circ \rI\fS^{(r)} 
\stackrel{\sim}{\rightarrow} \rI\hupsigma'^{(r)}_*\circ \rI\fS'^{(r)}\circ \rI\bvPhi^*,
\end{equation}
où la première flèche est induite par \eqref{pchamp16i} et la seconde flèche est l'isomorphisme sous-jacent au 
diagramme \eqref{pchamp16f},
induit par passage à la limite inductive, pour $r\in \mQ_{>0}$, un morphisme de foncteurs de $\bIndMod(\bvocB)$
dans $\bIndMH(\co_{\fX'},\txi^{-1}\tOmega^1_{\fX'/\cS})$
\begin{equation}\label{pchamp16p}
\rI\fgg^*\circ \rI\cH\rightarrow \rI\cH'\circ \rI\bvPhi^*,
\end{equation} 
où $\cH$ et $\cH'$ sont les foncteurs \eqref{indmdlb7b}. Comme le foncteur $\fgg^*$ commute aux limites inductives, le foncteur $\rI\fgg^*$ 
commute au foncteur $\kappa_{\co_\fX}$ \eqref{indmdlb1h}. On en déduit  un morphisme de foncteurs de $\bIndMod(\bvocB)$
dans $\bMH(\co_{\fX'},\txi^{-1}\tOmega^1_{\fX'/\cS})$
\begin{equation}\label{pchamp16q}
\fgg^*\circ \cH\rightarrow \cH'\circ \rI\bvPhi^*,
\end{equation} 
où $\cH$ et $\cH'$ sont les foncteurs \eqref{indmdlb7c}.

\begin{prop}\label{pchamp17}
Supposons que $g$ soit une immersion ouverte. Alors,
\begin{itemize}
\item[{\rm (i)}] Pour tout nombre rationnel $r\geq 0$, le morphisme \eqref{pchamp16i} est un isomorphisme. 
Il rend commutatif le diagramme de foncteurs
\begin{equation}\label{pchamp17a}
\xymatrix{
{\bIndMC(\bvcC^{(r)}/\bvocB)}\ar[d]_{\rI\bvPhi^*}\ar[r]^-(0.5){\rI\hupsigma^{(r)}_*}&{\bIndMH(\co_\fX,\xi^{-1}\tOmega^1_{\fX/\cS})}\ar[d]^{\rI\fgg^*}\\
{\bIndMC(\bvcC'^{(r)}/\bvocB')}\ar[r]^-(0.5){\rI\hupsigma'^{(r)}_*}&{\bIndMH(\co_{\fX'},\xi^{-1}\tOmega^1_{\fX'/\cS})}}
\end{equation}
\item[{\rm (ii)}]  Le morphisme \eqref{pchamp16q} est un isomorphisme. Il rend commutatif le diagramme de foncteurs  
\begin{equation}\label{pchamp17b}
\xymatrix{
{\bIndMod(\bvocB)}\ar[r]^-(0.5)\cH\ar[d]_{\rI\bvPhi^*}&
{\bMH(\co_\fX,\txi^{-1}\tOmega^1_{\fX/\cS})}\ar[d]^{\fgg^*}\\
{\bIndMod(\bvocB')}\ar[r]^-(0.5){\cH'}&{\bMH(\co_{\fX'},\txi^{-1}\tOmega^1_{\fX'/\cS})}}
\end{equation}
\end{itemize}
\end{prop}

(i) Cela résulte de (\cite{agt} III.9.15) et \ref{indsh42}. 

(ii) Cela résulte de (i) et des définitions.

\begin{prop}\label{pchamp18}
Soient $\cM$ un ind-$\bvocB$-module de Dolbeault, 
$\cN$ un $\co_\fX[\frac 1 p]$-fibré de Higgs soluble à coefficients dans $\txi^{-1}\tOmega^1_{\fX/\cS}$.
Alors, $\rI\bvPhi^*(\cM)$ est un ind-$\bvocB'$-module de Dolbeault et $\fgg^*(\cN)$ est un 
$\co_{\fX'}[\frac 1 p]$-fibré de Higgs soluble à coefficients dans $\txi^{-1}\tOmega^1_{\fX'/\cS}$.
Si de plus, $\cM$ et $\cN$ sont associés, $\rI\bvPhi^*(\cM)$  et $\fgg^*(\cN)$ sont associés. 
\end{prop}

Supposons qu'il existe un nombre rationnel $r>0$ et un isomorphisme de $\bIndMC(\bvcC^{(r)}/\bvocB)$
\begin{equation}\label{pchamp18a}
\alpha\colon \rI\hupsigma^{(r)*}(\cN)\stackrel{\sim}{\rightarrow}\rI\fS^{(r)}(\cM).
\end{equation}
D'après \eqref{pchamp16f}, \eqref{pchamp16g} et \eqref{indsh13a}, $\rI\bvPhi^*(\alpha)$ induit un isomorphisme de $\bIndMC(\bvcC'^{(r)}/\bvocB')$
\begin{equation}\label{pchamp18b}
\alpha'\colon \rI\hupsigma'^{(r)*}(\fgg^*(\cN))\stackrel{\sim}{\rightarrow}\rI\fS'^{(r)}(\rI\bvPhi^*(\cM));
\end{equation}
d'où la proposition.

\subsection{}\label{pchamp19}
D'après \ref{pchamp18}, $\rI\bvPhi^*$ induit un foncteur 
\begin{equation}\label{pchamp19a}
\rI\bvPhi^*\colon \bIndMod^\Dolb(\bvocB)\rightarrow \bIndMod^\Dolb(\bvocB'),
\end{equation}
et $\fgg^*$ induit un foncteur 
\begin{equation}\label{pchamp19b}
\fgg^*\colon \bMH^\sol(\co_\fX[\frac 1 p],\txi^{-1}\tOmega^1_{\fX/\cS})\rightarrow \bMH^\sol(\co_{\fX'}[\frac 1 p],
\txi^{-1}\tOmega^1_{\fX'/\cS}).
\end{equation}

\begin{prop}\label{pchamp20}\
\begin{itemize}
\item[{\rm (i)}] Le diagramme de foncteurs  
\begin{equation}\label{pchamp20a}
\xymatrix{
{\bIndMod^\Dolb(\bvocB)}\ar[r]^-(0.5)\cH\ar[d]_{\rI\bvPhi^*}&
{\bMH^\sol(\co_\fX[\frac 1 p],\txi^{-1}\tOmega^1_{\fX/\cS})}\ar[d]^{\fgg^*}\\
{\bIndMod^\Dolb(\bvocB')}\ar[r]^-(0.5){\cH'}&{\bMH^\sol(\co_{\fX'}[\frac 1 p],\txi^{-1}\tOmega^1_{\fX'/\cS})}}
\end{equation}
où $\cH$ et $\cH'$ sont les foncteurs \eqref{indmdlb12} est commutatif à isomorphisme canonique près.
\item[{\rm (ii)}] Le diagramme de foncteurs  
\begin{equation}\label{pchamp20b}
\xymatrix{
{\bMH^\sol(\co_\fX[\frac 1 p],\txi^{-1}\tOmega^1_{\fX/\cS})}\ar[d]_{\fgg^*}
\ar[r]^-(0.5)\cV&{\bIndMod^\Dolb(\bvocB)}\ar[d]^{\rI\bvPhi^*}\\
{\bMH^\sol(\co_{\fX'}[\frac 1 p],\txi^{-1}\tOmega^1_{\fX'/\cS})}\ar[r]^-(0.5){\cV'}&{\bIndMod^\Dolb(\bvocB')}}
\end{equation}
où $\cV$ et $\cV'$ sont les foncteurs \eqref{indmdlb17a} est commutatif à isomorphisme canonique près.
\end{itemize}
\end{prop}

(i) D'après \ref{indmdlb13}, il existe un nombre rationnel $r_\cM>0$ et un isomorphisme de \\
$\bIndMC(\bvcC^{(r_\cM)}/\bvocB)$ 
\begin{equation}\label{pchamp20c}
\alpha_\cM\colon \rI\hupsigma^{(r_\cM)*}(\cH(\cM))\stackrel{\sim}{\rightarrow} \rI\fS^{(r_\cM)}(\cM)
\end{equation}
vérifiant les propriétés (i) et (ii) de {\em loc.cit.} Pour tout nombre rationnel $r$ tel que $0< r\leq r_\cM$, on désigne par
\begin{equation}\label{pchamp20d}
\alpha^{(r)}_\cM\colon \rI\hupsigma^{(r)*}(\cH(\cM)) \stackrel{\sim}{\rightarrow}\rI\fS^{(r)}(\cM)
\end{equation}
l'isomorphisme de $\bIndMC(\bvcC^{(r)}/\bvocB)$ induit par $\rI\varepsilon^{r_\cM,r}(\alpha_\cM)$ \eqref{indmdlb3a} et 
les isomorphismes \eqref{indmdlb3b} et \eqref{indmdlb3c}. Compte tenu de \eqref{pchamp16f}, \eqref{pchamp16g} et \eqref{indsh13a}, 
$\rI\bvPhi^*(\alpha_\cM)$ induit un isomorphisme de $\bIndMC(\bvcC'^{(r_\cM)}/\bvocB')$
\begin{equation}\label{pchamp20e}
\alpha'_\cM\colon \rI\hupsigma'^{(r_\cM)*}(\fgg^*(\cH(\cM)))\stackrel{\sim}{\rightarrow}\rI\fS'^{(r_\cM)}(\rI\bvPhi^*(\cM)).
\end{equation}
De même, $\bvPhi^*(\alpha^{(r)}_\cM)$ induit un isomorphisme de $\bIndMC(\bvcC'^{(r)}/\bvocB')$
\begin{equation}\label{pchamp20f}
\alpha'^{(r)}_\cM\colon \rI\hupsigma'^{(r)*}(\fgg^*(\cH(\cM)))\stackrel{\sim}{\rightarrow}\rI\fS'^{(r)}(\rI\bvPhi^*(\cM)),
\end{equation}
qu'on peut aussi déduire de $\rI\varepsilon'^{r_\cM,r}(\alpha'_\cM)$ par \eqref{pchamp16l}.
On désigne par  
\begin{equation}\label{pchamp20g}
\beta'^{(r)}_\cM\colon \fgg^*(\cH(\cM))\rightarrow\vupsigma'^{(r)}_*(\rI\fS'^{(r)}(\rI\bvPhi^*(\cM))
\end{equation}
son adjoint \eqref{indmdlb2g}. 
D'après \ref{pchamp18}, $\rI\bvPhi^*(\cM)$ est un ind-$\bvocB'$-module de Dolbeault et 
$\fgg^*(\cH(\cM))$ est un $\co_{\fX'}[\frac 1 p]$-fibré de Higgs soluble à coefficients 
dans $\txi^{-1}\tOmega^1_{\fX'/\cS}$, associé à $\rI\bvPhi^*(\cM)$. 
Par suite, en vertu de \ref{indmdlb10}(i), le morphisme composé 
\begin{equation}\label{pchamp20h}
\fgg^*(\cH(\cM))\stackrel{\beta'^{(r)}_\cM}{\longrightarrow} \vupsigma'^{(r)}_*(\rI\fS'^{(r)}(\rI\bvPhi^*(\cM))\longrightarrow \cH'(\rI\bvPhi^*(\cM)),
\end{equation}
où la seconde flèche est le morphisme canonique \eqref{indmdlb7d}, est un isomorphisme 
qui dépend a priori de $\alpha_\cM$ mais pas de $r$.   
D'après la preuve de \ref{indmdlb20}, pour tout morphisme $u\colon \cM\rightarrow \cM'$ de $\bIndMod^\Dolb(\bvocB)$ 
et tout nombre rationnel $r$ tel que $0<r<\inf(r_\cM,r_{\cM'})$, le diagramme de $\bIndMC(\bvcC^{(r)}/\bvocB)$
\begin{equation}
\xymatrix{
{\rI\hupsigma^{(r)*}(\cH(\cM))}\ar[r]^-(0.5){\alpha^{(r)}_\cM}\ar[d]_{\rI\hupsigma^{(r)*}(\cH(u))}&{\rI\fS^{(r)}(\cM)}\ar[d]^{\rI\fS^{(r)}(u)}\\
{\rI\hupsigma^{(r)*}(\cH(\cM'))}\ar[r]^-(0.5){\alpha^{(r)}_{\cM'}}&{\rI\fS^{(r)}(\cM')}}
\end{equation}
est commutatif. On en déduit que l'isomorphisme composé \eqref{pchamp20h}
\begin{equation}\label{pchamp20i}
\fgg^*(\cH(\cM))\stackrel{\sim}{\rightarrow} \cH'(\rI\bvPhi^*(\cM))
\end{equation}
ne dépend que de $\cM$ (mais pas du choix de $\alpha_\cM$) et qu'il en dépend fonctoriellement; d'où la proposition. 
 
(ii) La preuve est similaire à celle de (i) et est laissée au lecteur.

\begin{remas}\label{pchamp21}
Soit $\cM$ un ind-$\bvocB$-module de Dolbeault.
\begin{itemize}
\item[(i)] Le morphisme canonique \eqref{pchamp16q}
\begin{equation}\label{pchamp21a}
\fgg^*(\cH(\cM))\rightarrow \cH'(\bvPhi^*(\cM))
\end{equation}
est un isomorphisme~; c'est l'isomorphisme sous-jacent au diagramme commutatif \eqref{pchamp20a}.
En effet, reprenons les notations de la preuve de \ref{pchamp20}(i) et notons, de plus,
\begin{equation}\label{pchamp21b}
\beta^{(r)}_\cM\colon \cH(\cM)\rightarrow \vupsigma^{(r)}_*(\fS^r(\cM))
\end{equation}
le morphisme adjoint de $\alpha^{(r)}_\cM$. Il résulte de \eqref{pchamp16k}
que le morphisme $\beta'^{(r)}_\cM$ \eqref{pchamp20g} est égal au composé 
\begin{equation}\label{pchamp21c}
\xymatrix{
{\fgg^*(\cH(\cM))}\ar[r]^-(0.5){\fgg^*(\beta^{(r)}_\cM)}&{\fgg^*(\vupsigma^{(r)}_*(\rI\fS^{(r)}(\cM)))}\ar[r]&
{\vupsigma'^{(r)}_*(\rI\bvPhi^*(\rI\fS^{(r)}(\cM)))}\ar[r]^-(0.5)\sim&{\vupsigma'^{(r)}_*(\fS'^{(r)}(\rI\bvPhi^*(\cM)))}},
\end{equation}
où la deuxième flèche est le morphisme \eqref{pchamp16ij} et la dernière flèche est l'isomorphisme sous-jacent 
au diagramme \eqref{pchamp16f}. Par ailleurs, la limite inductive des morphismes $\beta^r_\cM$, pour $r\in \mQ_{>0}$, 
est l'identité, et la limite inductive des morphismes $\beta'^{(r)}_\cM$, pour $r\in \mQ_{>0}$,
est égale à l'isomorphisme composé \eqref{pchamp20h}, sous-jacent au diagramme commutatif \eqref{pchamp20a}. 
\item[(ii)] Soient $r$ un nombre rationnel $>0$, 
\begin{equation}\label{pchamp21d}
\alpha\colon \rI\hupsigma^{(r)*}(\cH(\cM))\stackrel{\sim}{\rightarrow}\rI\fS^{(r)}(\cM)
\end{equation}
un isomorphisme de $\bIndMC(\bvcC^{(r)}/\bvocB)$ vérifiant les propriétés de \ref{indmdlb13}. 
Compte tenu de (i), \eqref{pchamp16f} et \eqref{pchamp16g}, on peut identifier $\rI\bvPhi^*(\alpha)$ à un isomorphisme
\begin{equation}\label{pchamp21e}
\alpha'\colon \rI\hupsigma'^{(r)*}(\cH'(\bvPhi^*(\cM)))\stackrel{\sim}{\rightarrow}\rI\fS'^{(r)}(\rI\bvPhi^*(\cM)).
\end{equation}
Par ailleurs, $\rI\bvPhi^*(\cM)$ est un ind-$\bvocB'$-module de Dolbeault d'après \ref{pchamp18}.
Il résulte aussitôt de \ref{indmdlb10} que $\alpha'$ vérifie les propriétés de \ref{indmdlb13}. 
\end{itemize}
\end{remas}

\subsection{}\label{pchamp22}
Le foncteur $\bvPhi^*$ \eqref{pchamp15b} induit un foncteur \eqref{indsh22a}
\begin{equation} \label{pchamp22a}
\bvPhi^*_\mQ\colon \bMod_\mQ(\bvocB)\rightarrow \bMod_\mQ(\bvocB').
\end{equation}
On développe dans (\cite{agt} III.14.7 et III.14.8) les compatibilités pour les $\bvocB_\mQ$-modules analogues à celles développées dans \ref{pchamp16} 
et \ref{pchamp17} pour les ind-$\bvocB$-modules. On en déduit la proposition suivante. 

\begin{prop}[\cite{agt} III.14.9]\label{pchamp2}
Soient $\cM$ un $\bvocB_\mQ$-module de Dolbeault (resp. fortement de Dolbeault), 
$\cN$ un $\co_\fX[\frac 1 p]$-fibré de Higgs rationnellement soluble (resp. fortement soluble) à coefficients dans $\txi^{-1}\tOmega^1_{\fX/\cS}$.
Alors, $\bvPhi^*_\mQ(\cM)$ est un $\bvocB'_\mQ$-module de Dolbeault (resp. fortement de Dolbeault) et $\fgg^*(\cN)$ est un 
$\co_{\fX'}[\frac 1 p]$-fibré de Higgs rationnellement soluble (resp. fortement soluble) à coefficients dans $\txi^{-1}\tOmega^1_{\fX'/\cS}$.
Si, de plus, $\cM$ et $\cN$ sont associés, $\bvPhi^*_\mQ(\cM)$ et $\fgg^*(\cN)$ sont associés. 
\end{prop}

\begin{prop}\label{pchamp23}
Tout $\co_\fX[\frac 1 p]$-fibré de Higgs fortement soluble $(\cN,\theta)$ à coefficients dans $\txi^{-1}\tOmega^1_{\fX/\cS}$ \eqref{aspglob1} 
est localement petit \eqref{definf26}.
\end{prop}

Cela résulte de \ref{aspglob31} et \ref{pchamp2}.

\subsection{}\label{pchamp3}
D'après \ref{pchamp2}, $\bvPhi^*$ induit un foncteur 
\begin{equation}\label{pchamp3a}
\bvPhi^*\colon \bMod_\mQ^\Dolb(\bvocB)\rightarrow \bMod_\mQ^\Dolb(\bvocB'),
\end{equation}
et $\fgg^*$ induit un foncteur 
\begin{equation}\label{pchamp3b}
\fgg^*\colon \bMH^\qsol(\co_\fX[\frac 1 p],\txi^{-1}\tOmega^1_{\fX/\cS})\rightarrow \bMH^\qsol(\co_{\fX'}[\frac 1 p],
\txi^{-1}\tOmega^1_{\fX'/\cS}).
\end{equation}

\begin{prop}[\cite{agt} III.14.11]\label{pchamp4}\
\begin{itemize}
\item[{\rm (i)}] Le diagramme de foncteurs  
\begin{equation}\label{pchamp4a}
\xymatrix{
{\bMod_\mQ^\Dolb(\bvocB)}\ar[r]^-(0.5){\cH_\mQ}\ar[d]_{\bvPhi^*}&
{\bMH^\qsol(\co_\fX[\frac 1 p],\txi^{-1}\tOmega^1_{\fX/\cS})}\ar[d]^{\fgg^*}\\
{\bMod_\mQ^\Dolb(\bvocB')}\ar[r]^-(0.5){\cH'_\mQ}&{\bMH^\qsol(\co_{\fX'}[\frac 1 p],\txi^{-1}\tOmega^1_{\fX'/\cS})}}
\end{equation}
où $\cH_\mQ$ et $\cH'_\mQ$ sont les foncteurs \eqref{aspglob15a} est commutatif à isomorphisme canonique près.
\item[{\rm (ii)}] Le diagramme de foncteurs  
\begin{equation}\label{pchamp4b}
\xymatrix{
{\bMH^\qsol(\co_\fX[\frac 1 p],\txi^{-1}\tOmega^1_{\fX/\cS})}\ar[d]_{\fgg^*}
\ar[r]^-(0.5){\cV_\mQ}&{\bMod_\mQ^\Dolb(\bvocB)}\ar[d]^{\bvPhi^*}\\
{\bMH^\qsol(\co_{\fX'}[\frac 1 p],\txi^{-1}\tOmega^1_{\fX'/\cS})}\ar[r]^-(0.5){\cV'_\mQ}&{\bMod_\mQ^\Dolb(\bvocB')}}
\end{equation}
où $\cV_\mQ$ et $\cV'_\mQ$ sont les foncteurs \eqref{aspglob19a} est commutatif à isomorphisme canonique près.
\end{itemize}
\end{prop}
La preuve est identique à celle de \ref{pchamp20}.

\section{Propriétés champêtres des modules de Dolbeault}\label{pchp}

\subsection{}\label{pchp1}
On note $\uppsi$ le morphisme composé
\begin{equation}\label{pchp1a}
\uppsi\colon \tE_s^{\mN^\circ}\stackrel{\uplambda}{\longrightarrow} \tE_s \stackrel{\sigma_s}{\longrightarrow} X_{s,\et}
\stackrel{a_\et}{\longrightarrow} X_\et,
\end{equation}
où $\uplambda$ est le morphisme canonique de topos défini dans \eqref{notconv13a}, $\sigma_s$ est le morphisme canonique de topos \eqref{ahttf12b} 
et $a\colon X_s\rightarrow X$ est l'injection canonique \eqref{ahttf1a}.
Pour tout objet $U$ de $\Et_{/X}$, on désigne par $f_U\colon (U,\cM_X|U)\rightarrow (S,\cM_S)$ 
le morphisme induit par $f$, et par $\tU\rightarrow \tX$ l'unique morphisme étale 
qui relève $\coU\rightarrow \coX$, de sorte que $(\tU,\cM_{\tX}|\tU)$ est une 
$(\tS,\cM_{\tS})$-déformation lisse de $(\coU,\cM_{\coX}|\coU)$. 
Le localisé du topos annelé $(\tE_s^{\mN^\circ},\bvocB)$ en $\uppsi^*(U)$  
est canoniquement équivalent au topos annelé analogue associé à $f_U$ en vertu de (\cite{agt} III.9.14). 
Pour tout nombre rationnel $r\geq 0$, $\bvcC^{(r)}|\uppsi^*(U)$ s'identifie à la $(\bvocB|\uppsi^*(U))$-algèbre de Higgs-Tate d'épaisseur 
$r$ associée à la déformation $(\tU,\cM_\tX|\tU)$ d'après \eqref{pchamp15k}. 
On désigne par $\INDMOD_U(\bvocB)$ la catégorie des ind-$(\bvocB|\uppsi^*(U))$-modules,
par $\INDMOD^\Dolb_U(\bvocB)$ la sous-catégorie des ind-$(\bvocB|\uppsi^*(U))$-modules de 
Dolbeault relativement à la déformation $(\tU,\cM_\tX|\tU)$, 
et par $\INDMC_U(\bvcC^{(r)}/\bvocB)$ la catégorie des ind-$(\bvcC^{(r)}|\uppsi^*(U))$-modules à 
$p^r$-connexion intégrable relativement à l'extension $(\bvcC^{(r)}|\uppsi^*(U))/(\bvocB|\uppsi^*(U))$ \eqref{indmdlb1}.

Pour tout morphisme $g\colon U'\rightarrow U$ de $\Et_{/X}$, on note 
\begin{equation}\label{pchp1b}
\begin{array}[t]{clcr}
\INDMOD_U(\bvocB)&\rightarrow &\INDMOD_{U'}(\bvocB)\\
\cM&\mapsto& \cM|\uppsi^*(U')
\end{array}
\end{equation}
le foncteur de restriction \eqref{indsh36}. Celui-ci induit, d'après \ref{pchamp18}, un foncteur 
\begin{equation}\label{pchp1c}
\begin{array}[t]{clcr}
\INDMOD^\Dolb_U(\bvocB)&\rightarrow &\INDMOD^\Dolb_{U'}(\bvocB)\\
\cM&\mapsto& \cM|\uppsi^*(U').
\end{array}
\end{equation}

Pour tout nombre rationnel $r\geq 0$, on note  
\begin{equation}\label{pchp1d}
\begin{array}[t]{clcr}
\INDMC_U(\bvcC^{(r)}/\bvocB)&\rightarrow &\INDMC_{U'}(\bvcC^{(r)}/\bvocB)\\
\cF&\mapsto& \cF|\uppsi^*(U')
\end{array}
\end{equation}
le foncteur de restriction. 

On désigne par $\Et_{\coh/X}$ la sous-catégorie pleine de $\Et_{/X}$ formée des schémas étales de présentation finie sur $X$ et par 
\begin{equation}\label{pchp1e}
\INDMOD^\Dolb(\bvocB)\rightarrow \Et_{\coh/X}
\end{equation}
la catégorie fibrée clivée normalisée dont la fibre au-dessus d'un objet $U$ de $\Et_{\coh/X}$ est la catégorie $\INDMOD^\Dolb_U(\bvocB)$
et le foncteur image inverse par un morphisme $U'\rightarrow U$ de $\Et_{\coh/X}$ 
est le foncteur de restriction \eqref{pchp1c}.

\begin{lem}\label{pchp2}
Soient $r$ un nombre rationnel $\geq 0$, 
$\cF, \cF'$ deux objets de $\bIndMC(\bvcC^{(r)}/\bvocB)$, $(U_i)_{i\in I}$ un recouvrement étale de $X$. 
Pour tous $(i,j)\in I^2$, on pose $U_{ij}=U_i\times_XU_j$. Alors, le diagramme d'applications d'ensembles
\begin{eqnarray}\label{pchp2a}
&&\Hom_{\INDMC_X(\bvcC^{(r)}/\bvocB)}(\cF,\cF')\rightarrow \prod_{i\in I}\Hom_{\INDMC_{U_i}(\bvcC^{(r)}/\bvocB)}(\cF|\uppsi^*(U_i),\cF'|\uppsi^*(U_i))\\
&&\rightrightarrows
\prod_{(i,j)\in I^2}\Hom_{\INDMC_{U_{ij}}(\bvcC^{(r)}/\bvocB)}(\cF|\uppsi^*(U_{ij}),\cF'|\uppsi^*(U_{ij}))\nonumber
\end{eqnarray}
est exact.
\end{lem}
En effet, comme $X$ est quasi-compact, on peut supposer $I$ fini, auquel cas l'assertion résulte facilement de \ref{indsh39}.

\begin{prop}\label{pchp3}
Soient $\cM$ un ind-$\bvocB$-module, $(U_i)_{i\in I}$ un recouvrement de $\Et_{\coh/X}$. 
Pour que $\cM$ soit de Dolbeault, il faut et il suffit que pour tout $i\in I$, le ind-$(\bvocB|\uppsi^*(U_i))$-module
$\cM|\uppsi^*(U_i)$ soit de Dolbeault.
\end{prop}

En effet, la condition est nécessaire en vertu de \ref{pchamp18}. Supposons que pour tout $i\in I$, $\cM|\uppsi^*(U_i)$ 
soit de Dolbeault et montrons que $\cM$ est de Dolbeault. 
Comme $X$ est quasi-compact, on peut supposer $I$ fini. Pour tout $i\in I$, notons $\fX_i$ le schéma formel
complété $p$-adique de $\oU_i$. Pour tout $(i,j)\in I^2$, posons $U_{ij}=U_i\times_XU_j$ et notons $\fX_{ij}$
le schéma formel complété $p$-adique de $\oU_{ij}$.  
En vertu de \ref{pchamp20}(i), on a un isomorphisme canonique de $\co_{\fX_{ij}}[\frac 1 p]$-modules de Higgs
à coefficients dans $\xi^{-1}\tOmega^1_{\fX_{ij}/\cS}$
\begin{equation}\label{pchp3a}
\cH_i(\cM|\uppsi^*(U_i))\otimes_{\co_{\fX_i}}\co_{\fX_{ij}}\stackrel{\sim}{\rightarrow}\cH_{ij}(\cM|\uppsi^*(U_{ij})),
\end{equation}
où $\cH_i$ et $\cH_{ij}$ sont les foncteurs \eqref{indmdlb12a} associés à $(f_{U_i},\tU_i,\cM_\tX|\tU_i)$ et 
$(f_{U_{ij}},\tU_{ij},\cM_\tX|\tU_{ij})$, respectivement.
On en déduit une donnée de descente $\delta$ sur les modules de Higgs $(\cH_i(\cM|\uppsi^*(U_i)))_{i\in I}$ relativement 
au recouvrement étale $(\fX_i\rightarrow \fX)_{i\in I}$. Celle-ci étant effective d'après (\cite{agt} III.6.22), 
il existe un $\co_\fX[\frac 1 p]$-fibré de Higgs $\cN$ à coefficients dans $\txi^{-1}\tOmega^1_{\fX/\cS}$ 
et pour tout $i\in I$, un isomorphisme de 
$\co_{\fX_i}[\frac 1 p]$-modules de Higgs 
\begin{equation}\label{pchp3b}
\cN\otimes_{\co_\fX}\co_{\fX_i}\stackrel{\sim}{\rightarrow} \cH_i(\cM|\uppsi^*(U_i)),
\end{equation}
qui induisent la donnée de descente $\delta$.

Pour tout $(i,j)\in I^2$ et tout nombre rationnel $r>0$, on note $\rI\hupsigma^{(r)*}_i$ et $\rI\fS^{(r)}_i$ 
(resp. $\rI\hupsigma^{(r)*}_{ij}$ et $\rI\fS^{(r)}_{ij}$) les foncteurs 
\eqref{indmdlb1d} et \eqref{indmdlb1b} associés à $(f_{U_i},\tU_i,\cM_\tX|\tU_i)$ (resp. $(f_{U_{ij}},\tU_{ij},\cM_\tX|\tU_{ij})$). 
D'après \ref{indmdlb13}, pour tout $i\in I$, il existe un nombre rationnel $r_i>0$ et un isomorphisme de $\INDMC_{U_i}(\bvcC^{(r)}/\bvocB)$
\begin{equation}\label{pchp3c}
\alpha_i\colon \rI\hupsigma^{(r_i)*}_i(\cH_i(\cM|\uppsi^*(U_i)))\stackrel{\sim}{\rightarrow} \rI\fS^{(r_i)}_i(\cM|\uppsi^*(U_i))
\end{equation}
vérifiant les propriétés (i) et (ii) de {\em loc.cit.} Pour tout $(i,j)\in I^2$, $\cM|\uppsi^*(U_{ij})$ est de Dolbeault en vertu de \ref{pchamp18}.  
D'après \eqref{pchamp16f}, \eqref{pchamp16g} et \ref{pchamp20}(i), $\alpha_i|\uppsi^*(U_{ij})$  s'identifie à un  isomorphisme 
\begin{equation}\label{pchp3d}
\alpha_i|\uppsi^*(U_{ij})\colon \rI\hupsigma^{(r_i)*}_{ij}(\cH_{ij}(\cM|\uppsi^*(U_{ij})))\stackrel{\sim}{\rightarrow} \rI\fS^{(r_i)}_{ij}(\cM|\uppsi^*(U_{ij})).
\end{equation}
Celui-ci vérifie les propriétés de \ref{indmdlb13}, compte tenu de \ref{pchamp21}.
Pour tout nombre rationnel $r$ tel que $0< r\leq r_i$, on désigne par 
$\rI\varepsilon_i^{r_i,r}\colon \INDMC_{U_i}(\bvcC^{(r_i)}/\bvocB)\rightarrow \INDMC_{U_i}(\bvcC^{(r)}/\bvocB)$
le foncteur \eqref{indmdlb3a} associé à $(f_{U_i},\tU_i,\cM_\tX|\tU_i)$ et par 
\begin{equation}\label{pchp3e}
\alpha^{(r)}_i\colon \rI\hupsigma^{(r)*}_i(\cH_i(\cM|\uppsi^*(U_i)))\stackrel{\sim}{\rightarrow} \rI\fS^{(r)}_i(\cM|\uppsi^*(U_i))
\end{equation}
l'isomorphisme de $\INDMC_{U_i}(\bvcC^{(r)}/\bvocB)$ induit par $\rI\varepsilon_i^{r_i,r}(\alpha_i)$  et 
les isomorphismes \eqref{indmdlb3b} et \eqref{indmdlb3c}. 
D'après \eqref{pchamp16f} et \eqref{pchamp16g}, on peut identifier $\alpha^{(r)}_i$ à un isomorphisme 
\begin{equation}\label{pchp3f}
\alpha^{(r)}_i\colon \rI\hupsigma^{(r)*}(\cN)|\uppsi^*(U_i)\stackrel{\sim}{\rightarrow} \rI\fS^{(r)}(\cM)|\uppsi^*(U_i).
\end{equation}
Il résulte de la preuve de \ref{indmdlb20} que pour tout nombre rationnel $0<r<\inf(r_i,r_j)$, on a dans $\INDMC_{U_{ij}}(\bvcC^{(r)}/\bvocB)$
\begin{equation}\label{pchp3g}
\alpha^{(r)}_i|\uppsi^*(U_{ij})=\alpha^{(r)}_j|\uppsi^*(U_{ij}).
\end{equation}
En vertu de \ref{pchp2}, pour tout nombre rationnel $0<r<\inf(r_i,i\in I)$, 
les isomorphismes $(\alpha^{(r)}_i)_{i\in I}$ se recollent en un isomorphisme de $\bIndMC(\bvcC^{(r)}/\bvocB)$
\begin{equation}\label{pchp3h}
\alpha^{(r)}\colon \rI\hupsigma^{(r)*}(\cN)\stackrel{\sim}{\rightarrow}  \rI\fS^{(r)}(\cM).
\end{equation}
Par suite, $\cM$ est de Dolbeault.

\begin{cor}\label{pchp4}
La propriété pour un ind-$\bvocB$-module d'être de Dolbeault ne dépend pas du choix de la déformation $(\tX,\cM_\tX)$ \eqref{defing12} 
pourvu que l'on reste dans l'un des cadres, absolu ou bien relatif \eqref{definf10}.
\end{cor}

En effet, la question est locale d'après \ref{pchp3}, et 
si $X$ est affine, toutes les $(\tS,\cM_\tS)$-déformations lisses de $(\coX,\cM_\coX)$ sont isomorphes en vertu de (\cite{kato1} 3.14).

\begin{prop}\label{pchp5}
Les conditions suivantes sont équivalentes~:
\begin{itemize}
\item[{\rm (i)}] La catégorie fibrée \eqref{pchp1e}
\begin{equation}\label{pchp5a}
\INDMOD^\Dolb(\bvocB)\rightarrow \Et_{\coh/X}
\end{equation}
est un champ {\rm (\cite{giraud2} II 1.2.1)}.
\item[{\rm (ii)}] Pour tout recouvrement $(U_i\rightarrow U)_{i\in I}$ de $\Et_{\coh/X}$, 
notant $\cU$ (resp. pour tout $i\in I$, $\cU_i$) le schéma formel complété $p$-adique de $\oU$ (resp. $\oU_i$), 
pour qu'un $\co_\cU[\frac 1 p]$-fibré de Higgs $\cN$ à coefficients dans $\txi^{-1}\tOmega^1_{\cU/\cS}$
soit soluble, il faut et il suffit que pour tout $i\in I$, le $\co_{\cU_i}[\frac 1 p]$-fibré de Higgs $\cN\otimes_{\co_\cU}\co_{\cU_i}$ 
à coefficients dans $\txi^{-1}\tOmega^1_{\cU_i/\cS}$ soit soluble.
\end{itemize}
\end{prop}

Soit $(U_i\rightarrow U)_{i\in I}$ un recouvrement de $\Et_{\coh/X}$. Pour tout $(i,j)\in I^2$, posons $U_{ij}=U_i\times_XU_j$. 
Notons $\cU$ le schéma formel complété $p$-adique de $\oU$ et 
$\cH^\star$ et $\cV^\star$ les foncteurs \eqref{indmdlb12a} et \eqref{indmdlb17a} associés à $(f_{U},\tU,\cM_\tX|\tU)$. 
Pour tout $i\in I$, notons $\cU_i$ le schéma formel complété $p$-adique de $\oU_i$ et 
$\cH_i$ et $\cV_i$ les foncteurs \eqref{indmdlb12a} et \eqref{indmdlb17a} associés à $(f_{U_i},\tU_i,\cM_\tX|\tU_i)$. 

Montrons d'abord (i)$\Rightarrow$(ii). Soit $\cN$ un $\co_\cU[\frac 1 p]$-fibré de Higgs à coefficients dans 
$\txi^{-1}\tOmega^1_{\cU/\cS}$. 
Si $\cN$ est soluble, pour tout $i\in I$, $\cN\otimes_{\co_\cU}\co_{\cU_i}$ est soluble d'après \ref{pchamp18}. 
Inversement, supposons que pour tout $i\in I$, $\cN\otimes_{\co_\cU}\co_{\cU_i}$ soit soluble et montrons que $\cN$ est soluble. 
Pour tout $i\in I$, $\cM_i=\cV_i(\cN\otimes_{\co_\cU}\co_{\cU_i})$ est un ind-$\bvocB|\uppsi^*(U_i)$-module de Dolbeault.
D'après \ref{pchamp20}(ii), la donnée de descente canonique sur les fibrés de Higgs $(\cN\otimes_{\co_\cU}\co_{\cU_i})_{i\in I}$ 
relativement au recouvrement étale $(\cU_i\rightarrow \cU)_{i\in I}$
induit une donnée de descente $\delta$ sur les modules de Dolbeault $(\cM_i)_{i\in I}$
relativement au recouvrement $(U_i\rightarrow U)_{i\in I}$.
Cette dernière étant effective d'après (i), il existe un ind-$\bvocB|\uppsi^*(U)$-module de Dolbeault 
$\cM$ et pour tout $i\in I$, un isomorphisme de ind-$\bvocB|\uppsi^*(U_i)$-modules 
\begin{equation}
\cM|\uppsi^*(U_i)\stackrel{\sim}{\rightarrow}\cM_i
\end{equation}
qui induisent la donnée de descente $\delta$. 
En vertu de \ref{indmdlb20} et \ref{pchamp20}(i), on a un isomorphisme canonique de $\co_\cU[\frac 1 p]$-fibrés de Higgs 
$\cH^\star(\cM)\stackrel{\sim}{\rightarrow}\cN$. Par suite, $\cN$ est soluble.

Montrons ensuite (ii)$\Rightarrow$(i).  Pour tous ind-$\bvocB|\uppsi^*(U)$-modules $\cM$ et $\cM'$,  
le diagramme d'applications d'ensembles
\begin{eqnarray}
\lefteqn{\Hom_{\INDMOD_{U}(\bvocB)}(\cM,\cM')\rightarrow \prod_{i\in I}
\Hom_{\INDMOD_{U_i}(\bvocB)}(\cM|\uppsi^*(U_i),\cM'|\uppsi^*(U_i))}\\
&&\rightrightarrows 
\prod_{(i,j)\in I^2}\Hom_{\INDMOD_{U_{ij}}(\bvocB))}(\cM|\uppsi^*(U_{ij}),\cM'|\uppsi^*(U_{ij}))\nonumber
\end{eqnarray}
est exact. En effet, comme $U$ est quasi-compact, on peut supposer $I$ fini, auquel cas l'assertion résulte de \ref{indsh39}.

Pour tout $i\in I$, soit $\cM_i$ un ind-$\bvocB_\mQ|\uppsi^*(U_i)$-module de Dolbeault et soit $\delta$ une donnée de descente
sur $(\cM_i)_{i\in I}$ relativement au recouvrement $(U_i\rightarrow U)_{i\in I}$. Montrons que $\delta$ est effective. 
Par hypothèse, pour tout $i\in I$, 
$\cN_i=\cH_i(\cM_i)$ est un $\co_{\cU_i}[\frac 1 p]$-fibré de Higgs soluble à coefficients dans $\txi^{-1}\tOmega^1_{\cU_i/\cS}$. 
Compte tenu de \ref{pchamp20}(i), $\delta$ induit une donnée de descente $\gamma$ 
sur les fibrés de Higgs $(\cN_i)_{i\in I}$ relativement
au recouvrement étale $(\cU_i\rightarrow \cU)_{\in I}$. Celle-ci étant effective d'après (\cite{agt} III.6.22), 
il existe un $\co_\cU[\frac 1 p]$-fibré de Higgs $\cN$ et pour tout $i\in I$, un isomorphisme de 
$\co_{\cU_i}[\frac 1 p]$-modules de Higgs 
\begin{equation}
\cN\otimes_{\co_\cU}\co_{\cU_i}\stackrel{\sim}{\rightarrow} \cN_i,
\end{equation}
qui induisent la donnée de descente $\gamma$. D'après (ii), $\cN$ est soluble. Par suite, 
$\cM=\cV^\star(\cN)$ est un ind-$\bvocB_\mQ|\uppsi^*(U)$-module de Dolbeault.
D'après \ref{indmdlb20} et \ref{pchamp20}(ii), pour tout $i\in I$, on a un isomorphisme canonique de ind-$\bvocB_\mQ|\uppsi^*(U_i)$-modules
$\cM|\uppsi^*(U_i)\stackrel{\sim}{\rightarrow}\cM_i$, lesquels induisent la donnée de descente $\delta$, ce qui prouve l'assertion.

\subsection{}\label{pchp6}
Pour tout objet $U$ de $\Et_{/X}$, on désigne par $\MOD_U(\bvocB)$ la catégorie des $(\bvocB|\uppsi^*(U))$-modules,
par $\MOD_{\mQ,U}(\bvocB)$ la catégorie des $(\bvocB|\uppsi^*(U))_\mQ$-modules et par  $\MOD_{\mQ,U}^\Dolb(\bvocB)$ \\ 
(resp. $\MOD_{\mQ,U}^\fDolb(\bvocB)$) la sous-catégorie des $(\bvocB|\uppsi^*(U))_\mQ$-modules de 
Dolbeault (resp. fortement de Dolbeault) relativement à la déformation $(\tU,\cM_\tX|\tU)$ \eqref{pchp1}. 
Pour tout morphisme $g\colon U'\rightarrow U$ de $\Et_{/X}$, le foncteur de restriction
\begin{equation}\label{pchp6a}
\begin{array}[t]{clcr}
\MOD_U(\bvocB)&\rightarrow& \MOD_{U'}(\bvocB)\\
\cM&\mapsto& \cM|\uppsi^*(U')
\end{array}
\end{equation}
induit un foncteur
\begin{equation}\label{pchp6b}
\begin{array}[t]{clcr}
\MOD_{\mQ,U}(\bvocB)&\rightarrow &\MOD_{\mQ,U'}(\bvocB)\\
\cM&\mapsto& \cM|\uppsi^*(U').
\end{array}
\end{equation}
D'après \ref{pchamp2}, ce dernier induit deux foncteurs
\begin{equation}\label{pchp6c}
\begin{array}[t]{clcr}
\MOD_{\mQ,U}^\Dolb(\bvocB)&\rightarrow& \MOD_{\mQ,U'}^\Dolb(\bvocB)\\
\cM&\mapsto& \cM|\uppsi^*(U').
\end{array}
\end{equation}
\begin{equation}\label{pchp6cd}
\begin{array}[t]{clcr}
\MOD_{\mQ,U}^\fDolb(\bvocB)&\rightarrow& \MOD_{\mQ,U'}^\fDolb(\bvocB)\\
\cM&\mapsto& \cM|\uppsi^*(U').
\end{array}
\end{equation}

On note ${\bf MOD}(\bvocB)$ la $(\tE_s^{\mN^\circ})$-catégorie fibrée (et même scindée \cite{sga1} VI § 9) 
des $\bvocB$-modules sur $\tE_s^{\mN^\circ}$ (\cite{giraud2} II 3.4.1). 
C'est un champ au-dessus de $\tE_s^{\mN^\circ}$ d'après (\cite{giraud2} II 3.4.4). On désigne par 
\begin{equation}\label{pchp6d}
\MOD(\bvocB)\rightarrow \Et_{\coh/X}
\end{equation}
le changement de base de ${\bf MOD}(\bvocB)$ (\cite{sga1} VI § 3) par $\uppsi^*\circ \varepsilon$,  où
$\uppsi$ est le morphisme \eqref{pchp1a} et $\varepsilon \colon \Et_{\coh/X}\rightarrow X_{\et}$ est le foncteur canonique. 
C'est aussi un champ d'après (\cite{giraud2} II 3.1.1). 
On en déduit une catégorie fibrée 
\begin{equation}\label{pchp6e}
\MOD_\mQ(\bvocB)\rightarrow \Et_{\coh/X},
\end{equation}
dont la fibre au-dessus d'un objet $U$ de $\Et_{\coh/X}$ est la catégorie $\MOD_{\mQ,U}(\bvocB)$
et le foncteur image inverse par un morphisme 
$U'\rightarrow U$ de $\Et_{\coh/X}$ est le foncteur de restriction \eqref{pchp6b}. 
On notera que ce n'est a priori pas un champ. Elle induit deux catégoriex fibrées 
\begin{eqnarray}
\MOD_\mQ^\Dolb(\bvocB)&\rightarrow& \Et_{\coh/X},\label{pchp6f}\\
\MOD_\mQ^\fDolb(\bvocB)&\rightarrow& \Et_{\coh/X},\label{pchp6g}
\end{eqnarray}
dont les fibres au-dessus d'un objet $U$ de $\Et_{\coh/X}$ sont les catégorie $\MOD_{\mQ,U}^\Dolb(\bvocB)$ et $\MOD_{\mQ,U}^\fDolb(\bvocB)$, 
respectivement, et les foncteurs image inverse par un morphisme $U'\rightarrow U$ de $\Et_{\coh/X}$ 
sont les foncteurs de restriction \eqref{pchp6c} et \eqref{pchp6cd}.

\begin{prop}\label{pchp7}
Soient $\cM$ un $\bvocB_\mQ$-module (resp. adique de type fini), $(U_i)_{i\in I}$ un recouvrement de $\Et_{\coh/X}$. 
Pour que $\cM$ soit de Dolbeault (resp. fortement de Dolbeault), il faut et il suffit que pour tout $i\in I$, le $(\bvocB|\uppsi^*(U_i))_\mQ$-module
$\cM|\uppsi^*(U_i)$ soit de Dolbeault (resp. fortement de Dolbeault) \eqref{pchp6b}.
\end{prop}
Il suffit de calquer la preuve de (\cite{agt} III.15.4). 

\begin{cor}\label{pchp8}
La propriété pour un $\bvocB_\mQ$-module d'être de Dolbeault (resp. fortement de Dolbeault) 
ne dépend pas du choix de la déformation $(\tX,\cM_\tX)$ \eqref{defing12} 
pourvu que l'on reste dans l'un des cadres, absolu ou bien relatif \eqref{definf10}.
\end{cor}

En effet, le fait qu'un $\bvocB_\mQ$-module soit de Dolbeault (resp. fortement de Dolbeault) est une question locale d'après \ref{pchp7}, et 
si $X$ est affine, toutes les $(\tS,\cM_\tS)$-déformations lisses de $(\coX,\cM_\coX)$ sont isomorphes en vertu de (\cite{kato1} 3.14).

\begin{prop}\label{pchp9}
Les conditions suivantes sont équivalentes:
\begin{itemize}
\item[{\rm (i)}] La catégorie fibrée \eqref{pchp6f}
\begin{equation}\label{pchp9a}
\MOD_\mQ^\Dolb(\bvocB)\rightarrow \Et_{\coh/X}
\end{equation}
est un champ {\rm (\cite{giraud2} II 1.2.1)}.
\item[{\rm (ii)}] Pour tout recouvrement $(U_i\rightarrow U)_{i\in I}$ de $\Et_{\coh/X}$, 
notant $\cU$ (resp. pour tout $i\in I$, $\cU_i$) le schéma formel complété $p$-adique de $\oU$ (resp. $\oU_i$), 
pour qu'un $\co_\cU[\frac 1 p]$-fibré de Higgs $\cN$ à coefficients dans $\txi^{-1}\tOmega^1_{\cU/\cS}$
soit rationnellement soluble \eqref{aspglob1}, il faut et il suffit que pour tout $i\in I$, le $\co_{\cU_i}[\frac 1 p]$-fibré de Higgs $\cN\otimes_{\co_\cU}\co_{\cU_i}$ 
à coefficients dans $\txi^{-1}\tOmega^1_{\cU_i/\cS}$ soit rationnellement soluble.
\end{itemize}
\end{prop}
Il suffit de calquer la preuve de (\cite{agt} III.15.5).

\begin{prop}\label{pchp90}
Les conditions suivantes sont équivalentes:
\begin{itemize}
\item[{\rm (i)}] La catégorie fibrée \eqref{pchp6g}
\begin{equation}\label{pchp90a}
\MOD_\mQ^\fDolb(\bvocB)\rightarrow \Et_{\coh/X}
\end{equation}
est un champ.
\item[{\rm (ii)}] Pour tout recouvrement $(U_i\rightarrow U)_{i\in I}$ de $\Et_{\coh/X}$, 
notant $\cU$ (resp. pour tout $i\in I$, $\cU_i$) le schéma formel complété $p$-adique de $\oU$ (resp. $\oU_i$), 
pour qu'un $\co_\cU[\frac 1 p]$-fibré de Higgs $\cN$ à coefficients dans $\txi^{-1}\tOmega^1_{\cU/\cS}$
soit fortement soluble \eqref{aspglob1}, il faut et il suffit que pour tout $i\in I$, le $\co_{\cU_i}[\frac 1 p]$-fibré de Higgs $\cN\otimes_{\co_\cU}\co_{\cU_i}$ 
à coefficients dans $\txi^{-1}\tOmega^1_{\cU_i/\cS}$ soit fortement soluble.
\end{itemize}
\end{prop}

Il suffit de calquer la preuve de (\cite{agt} III.15.5), qui correspond au cas absolu \eqref{definf10}.

\begin{prop}\label{pchp10}\
Sous les conditions de \ref{pchp90},  
pour qu'un $\co_\fX[\frac 1 p]$-fibré de Higgs à coefficients dans $\txi^{-1}\tOmega^1_{\fX/\cS}$ soit localement petit \eqref{definf26}, 
il faut et il suffit qu'il soit fortement soluble \eqref{aspglob1}. 
\end{prop}

Cela résulte de \ref{mdpsa15}. 

\begin{rema}\label{pchp11}
Sous les conditions de \ref{pchp5} (resp. \ref{pchp9}, resp. \ref{pchp90}), la propriété pour un $\co_\fX[\frac 1 p]$-fibré de Higgs à coefficients dans $\txi^{-1}\tOmega^1_{\fX/\cS}$ 
d'être soluble (resp. rationnellement soluble, resp. fortement soluble) ne dépend pas du choix de la déformation $(\tX,\cM_\tX)$ \eqref{defing12} 
pourvu que l'on reste dans l'un des cadres, absolu ou bien relatif \eqref{definf10}. En effet, la question est locale par hypothèse, et 
si $X$ est affine, toutes les $(\tS,\cM_\tS)$-déformations lisses de $(\coX,\cM_\coX)$ sont isomorphes. 
\end{rema}

\section{Modules de Hodge-Tate}\label{mht}

\begin{defi}\label{mht1}
On appelle {\em $\bvocB_\mQ$-module de Hodge-Tate} tout $\bvocB_\mQ$-module de Dolbeault $\cM$ \eqref{aspglob1} 
dont le $\co_\fX[\frac 1 p]$-fibré de Higgs associé $\cH_\mQ(\cM)$ \eqref{aspglob15a} est nilpotent \eqref{definf22}. 
\end{defi} 

On désigne par $\bMod_\mQ^\HT(\bvocB)$ la sous-catégorie pleine de $\bMod_\mQ(\bvocB)$ 
formée des $\bvocB_\mQ$-modules de Hodge-Tate, et par $\bMH^\qsolnilp(\co_\fX[\frac 1 p], \txi^{-1}\tOmega^1_{\fX/\cS})$ 
la sous-catégorie pleine de  $\bMH(\co_\fX[\frac 1 p], \txi^{-1}\tOmega^1_{\fX/\cS})$
formée des $\co_\fX[\frac 1 p]$-fibrés de Higgs rationnellement solubles et nilpotents à coefficients dans $\txi^{-1}\tOmega^1_{\fX/\cS}$ \eqref{aspglob1}.

\begin{prop}\label{mht2}
Les foncteurs $\cH_\mQ$ \eqref{aspglob15a} et $\cV_\mQ$ \eqref{aspglob19a} induisent des équivalences de catégories quasi-inverses l'une de l'autre
\begin{equation}\label{mht2a}
\xymatrix{
{\bMod^\HT_\mQ(\bvocB)}\ar@<1ex>[r]^-(0.5){\cH_\mQ}&{\bMH^\qsolnilp(\co_\fX[\frac 1 p], \txi^{-1}\tOmega^1_{\fX/\cS})}
\ar@<1ex>[l]^-(0.5){\cV_\mQ}}
\end{equation} 
\end{prop}

Cela résulte aussitôt de \ref{aspglob21}.

\begin{prop}\label{mht3}
Supposons que $X$ soit un objet de $\bP$ \eqref{ahttf4}.
Alors, tout $\co_\fX[\frac 1 p]$-fibré de Higgs nilpotent à coefficients dans $\txi^{-1}\tOmega^1_{\fX/\cS}$ est fortement soluble \eqref{aspglob1}.
\end{prop}

En effet, $(\cN,\theta)$ est petit d'après \ref{definf27}. Il est donc fortement soluble en vertu de \ref{mdpsa15}. 

\begin{cor}\label{mht30}
Supposons que $X$ soit un objet de $\bP$ \eqref{ahttf4}.
Alors, tout $\bvocB_\mQ$-module de Hodge-Tate est fortement de Dolbeault \eqref{aspglob1}, et est en particulier adique de type fini \eqref{ahttf60}.
\end{cor}

Cela résulte de \ref{mht3} et \ref{aspglob21}. 

\begin{cor}\label{mht31}
Sous les conditions de \ref{pchp90}, tout $\co_\fX[\frac 1 p]$-fibré de Higgs nilpotent à coefficients dans $\txi^{-1}\tOmega^1_{\fX/\cS}$ est fortement soluble,
et tout $\bvocB_\mQ$-module de Hodge-Tate est fortement de Dolbeault \eqref{aspglob1}, et est en particulier adique de type fini.
\end{cor}

Cela résulte de \ref{mht3} et \ref{aspglob21}. 

\begin{teo}\label{mht32}
Supposons que $X$ soit un objet de $\bP$ \eqref{ahttf4} et reprenons les notations de \ref{mdpsa29}. 
Pour tout entier $1\leq i\leq c$, notons $\bRep_{\hocR_i[\frac 1 p]}^{\HT}(\Delta_i)$ la sous-catégorie de $\bRep_{\hocR_i[\frac 1 p]}(\Delta_i)$ formée des 
$\hocR_i[\frac 1 p]$-représentations de Hodge-Tate de $\Delta_i$ \eqref{repdolb23}. 
Alors, les foncteurs $(\hupbeta_i)_{1\leq i\leq c}$ \eqref{mdpsa29g} induisent une équivalence de catégories 
\begin{equation}\label{mht32a}
\bMod^\HT_\mQ(\bvocB)\stackrel{\sim}{\rightarrow} \prod_{1\leq i\leq c} \bRep_{\hocR_i[\frac 1 p]}^{\HT}(\Delta_i).
\end{equation} 
\end{teo}

Cela résulte de \ref{mdpsa32}(i), \ref{mdpsa33}(iii) et \ref{definf23}.

\subsection{}\label{mht4}
Reprenons les notations de \ref{pchp6}. Pour tout objet $U$ de $\Et_{/X}$, on désigne par 
\begin{equation}\label{mht4a}
\MOD_{U,\mQ}^\HT(\bvocB)
\end{equation}
la catégorie des $(\bvocB|\uppsi^*(U))_\mQ$-modules de Hodge-Tate relativement à la déformation $(\tU,\cM_\tX|\tU)$.
D'après \ref{pchamp2}, pour tout morphisme $g\colon U'\rightarrow U$ de $\Et_{/X}$, le foncteur de restriction \eqref{pchp6b}
induit un foncteur 
\begin{equation}\label{mht4b}
\begin{array}[t]{clcr}
\MOD_{U,\mQ}^\HT(\bvocB)&\rightarrow& \MOD_{U',\mQ}^\HT(\bvocB)\\
\cM&\mapsto&\cM|\uppsi^*(U_i).
\end{array}
\end{equation}
La catégorie fibrée \eqref{pchp6f} induit une catégorie fibrée 
\begin{equation}\label{mht4c}
\MOD_\mQ^\HT(\bvocB)\rightarrow \Et_{\coh/X}
\end{equation}
dont la fibre au-dessus d'un objet $U$ de $\Et_{\coh/X}$ est la catégorie $\MOD_{U,\mQ}^\HT(\bvocB)$
et le foncteur image inverse par un morphisme $U'\rightarrow U$ de $\Et_{\coh/X}$ 
est le foncteur de restriction \eqref{mht4b}.

\begin{prop}\label{mht6}
Soient $\cM$ un $\bvocB_\mQ$-module, $(U_i)_{i\in I}$ un recouvrement de $\Et_{\coh/X}$. 
Pour que $\cM$ soit de Hodge-Tate, il faut et il suffit que pour tout $i\in I$, le $(\bvocB|\uppsi^*(U_i))_\mQ$-module
$\cM|\uppsi^*(U_i)$ soit de Hodge-Tate.
\end{prop}
Cela résulte de \ref{pchp7} et \ref{pchamp4}(i).

\begin{cor}\label{mht7}
La propriété pour un $\bvocB_\mQ$-module d'être de Hodge-Tate ne dépend pas du choix de la déformation $(\tX,\cM_\tX)$ \eqref{defing12} 
ni même du cadre absolu ou relatif \eqref{definf10}.
\end{cor}

En effet, le fait qu'un $\bvocB_\mQ$-module soit de Hodge-Tate étant une question locale d'après \ref{mht6}, 
on peut supposer $X$ affine vérifiant les hypothèses de \ref{mdpsa1}.
Toutes les $(\tS,\cM_\tS)$-déformations lisses de $(\coX,\cM_\coX)$ sont alors isomorphes en vertu de (\cite{kato1} 3.14).
La propriété pour un $\bvocB_\mQ$-module d'être de Hodge-Tate ne dépend donc pas du choix de la 
$(\tS,\cM_\tS)$-déformation $(\tX,\cM_\tX)$ pourvu que l'on reste dans l'un des cadres, absolu ou bien relatif. 

Pour montrer que cette propriété ne dépend pas non plus du cadre absolu ou relatif, 
reprenons les hypothèses et notations de \ref{ahttf29} et \ref{aspglob24}. Compte tenu de \ref{aspglob25}, il suffit de montrer que si $\cM$ est 
un $\ocB_\mQ$-module de Hodge-Tate relativement à la déformation $(\tX,\cM_{\tX})$, il est de Hodge-Tate relativement à la déformation $(\tX',\cM_{\tX'})$.
Supposons donc le $\ocB_\mQ$-module $\cM$ de Hodge-Tate relativement à la déformation $(\tX,\cM_{\tX})$,
et notons $\theta$ le champ de Higgs canonique sur $\cH_\mQ(\cM)$ à coefficients dans $\xi^{-1}\tOmega^1_{\fX/\cS}$ \eqref{aspglob24a}. 
D'après les preuves de \ref{definf27} et \ref{repht2}(i), pour tout $\varepsilon \geq 0$, il existe un sous-$\co_\fX$-module cohérent $\cN_\varepsilon$ de $\cH_\mQ(\cM)$ 
qui l'engendre sur $\co_\fX[\frac 1 p]$ tel que  
\begin{equation}\label{mht7a}
\theta(\cN_\varepsilon)\subset  p^\varepsilon \xi^{-1}\tOmega^1_{\fX/\cS}\otimes_{\co_\fX} \cN_\varepsilon.
\end{equation}
Prenons $\varepsilon >\rho+\frac{1}{p-1}$ \eqref{ahttf29}, on a alors \eqref{definf16a}
\begin{equation}\label{mht7b}
\theta(\cN_\varepsilon)\subset  p^{\varepsilon-\rho} (\xi^*_\pi)^{-1}\tOmega^1_{\fX/\cS}\otimes_{\co_\fX} \cN_\varepsilon.
\end{equation}
Par suite, $\theta$ induit sur $\cN_\varepsilon$ un champ de Higgs à coefficients dans $(\xi^*_\pi)^{-1}\tOmega^1_{\fX/\cS}$ que l'on note 
$\theta'$ pour éviter toute confusion.  
Reprenons les hypothèses et notations de \ref{mdpsa8} et notons $\bvcP$ le $R_1$-module de $(\oX^{\circ}_\fet)^{\mN^\circ}$ associé à $(\cN_\varepsilon,\theta)$ 
par le foncteur  
\begin{equation}\label{mht7c}
\bMH^{\qpp}(\co_\fX,\xi^{-1}\tOmega^1_{\fX/\cS})\rightarrow \bMod((\oX^{\circ}_\fet)^{\mN^\circ},R_1)
\end{equation}
défini dans \eqref{mdpsa8d}. Il résulte aussitôt de \ref{pmh5} que $\bvcP$ est aussi le  $R_1$-module de $(\oX^{\circ}_\fet)^{\mN^\circ}$ associé à $(\cN_\varepsilon,\theta')$ 
par le foncteur analogue
\begin{equation}\label{mht7d}
\bMH^{\qpp}(\co_\fX,(\xi^*_\pi)^{-1}\tOmega^1_{\fX/\cS})\rightarrow \bMod((\oX^{\circ}_\fet)^{\mN^\circ},R_1).
\end{equation}
En vertu de \ref{aspglob21} et \ref{mdpsa10} appliqué à $(\cN_\varepsilon,\theta)$, on a un isomorphisme 
\begin{equation}\label{mht7e}
\cM\stackrel{\sim}{\rightarrow}\bvbeta^*(\bvcP\otimes_{R_1}\bvocB_X)_\mQ.
\end{equation}
Il résulte alors de \ref{mdpsa10} appliqué à $(\cN_\varepsilon,\theta')$ que $\cM$ est fortement de Dolbeault relativement à la déformation $(\tX',\cM_{\tX'})$ et qu'on a un isomorphisme 
\begin{equation}\label{mht7f}
\cH'_\mQ(\cM)\stackrel{\sim}{\rightarrow}(\cN_\varepsilon\otimes_{\mZ_p}\mQ_p,\theta').
\end{equation}
Par suite, $\cM$ est de Hodge-Tate relativement à la déformation $(\tX',\cM_{\tX'})$, d'où la proposition.

\begin{prop}\label{mht5}
Les conditions suivantes sont équivalentes:
\begin{itemize}
\item[{\rm (i)}] La catégorie fibrée \eqref{mht4c}
\begin{equation}\label{mht5a}
\MOD_\mQ^\HT(\bvocB)\rightarrow \Et_{\coh/X}
\end{equation}
est un champ {\rm (\cite{giraud2} II 1.2.1)}.
\item[{\rm (ii)}] 
Tout $\co_\fX[\frac 1 p]$-fibré de Higgs nilpotent à coefficients dans $\txi^{-1}\tOmega^1_{\fX/\cS}$ est rationnellement soluble.
\end{itemize}

De plus, ces conditions sont vérifiées si les conditions de \ref{pchp9} le sont. 
\end{prop}

En effet, calquant la preuve de (\cite{agt} III.15.5), on montre que la condition (i) est équivalente à la condition suivante~:
\begin{itemize}
\item[{\rm (ii')}] Pour tout recouvrement $(U_i\rightarrow U)_{i\in I}$ de $\Et_{\coh/X}$, 
notant $\cU$ (resp. pour tout $i\in I$, $\cU_i$) le schéma formel complété $p$-adique de $\oU$ (resp. $\oU_i$), 
pour qu'un $\co_\cU[\frac 1 p]$-fibré de Higgs {\em nilpotent} $\cN$ à coefficients dans $\txi^{-1}\tOmega^1_{\cU/\cS}$
soit rationnellement soluble, il faut et il suffit que pour tout $i\in I$, le $\co_{\cU_i}[\frac 1 p]$-fibré de Higgs $\cN\otimes_{\co_\cU}\co_{\cU_i}$ 
à coefficients dans $\txi^{-1}\tOmega^1_{\cU_i/\cS}$ soit rationnellement soluble.
\end{itemize}
On a clairement (ii)$\Rightarrow$(ii'), et on a (ii')$\Rightarrow$(ii) en vertu de \ref{mht3}. Enfin, la condition \ref{pchp9}(ii) implique clairement (ii').

\section{Systèmes locaux de Dolbeault et de Hodge-Tate}\label{sld}

\subsection{}\label{sld1}
Pour tout $\mU$-topos $T$, on désigne par $\bvmZ_p$ la $\mZ_p$-algèbre $(\mZ/p^n\mZ)_{n\geq 1}$ de $T^{\mN^\circ}$ \eqref{notconv13} et par
\begin{equation}\label{sld1a}
Q_{\bvmZ_p}\colon \bMod(\bvmZ_p,T^{\mN^\circ}) \rightarrow \bMod_\mQ(\bvmZ_p,T^{\mN^\circ}), \ \ \ M\mapsto M_\mQ,
\end{equation}
le foncteur canonique \eqref{indsh20a}. On rappelle que les objets de $\bMod_\mQ(\bvmZ_p,T^{\mN^\circ})$ sont aussi appelés des 
{\em $\bvmZ_{p,\mQ}$-modules}. On note 
\begin{equation}\label{sld1b}
\upalpha_{\bvmZ_p}\colon \bMod_\mQ(\bvmZ_p,T^{\mN^\circ}) \rightarrow \bIndMod(\bvmZ_p,T^{\mN^\circ}),
\end{equation}
le foncteur canonique \eqref{indsh20c}.

\subsection{}\label{sld2}
On désigne par 
\begin{equation}\label{sld2a}
\bvpsi\colon (\oX^\circ_\et)^{\mN^\circ}\rightarrow \tE^{\mN^\circ}
\end{equation}
le morphisme de topos induit par $\psi$ \eqref{ahttf2n}, que l'on considérera naturellement comme un morphisme de topos annelés par les anneaux $\bvmZ_p$.
Il résulte de (\cite{agt} III.7.5 et VI.10.9(iii)) que l'homomorphisme canonique $\bvmZ_p\rightarrow \bvpsi_*(\bvmZ_p)$ est un isomorphisme. 
On utilisera pour $\bvpsi$ les notations introduites dans \ref{indsh21} et \ref{indsh22}. 

Pour tout entier $n\geq 0$, tout $\ocB_n$-module $\cM_n$ de $\tE$ est naturellement un objet de $\tE_s$ (\cite{agt} III.9.7). 
De plus, le foncteur $\delta_*$ \eqref{ahttf12a} étant exact, on a un isomorphisme canonique 
\begin{equation}\label{sld2b}
\rR\Gamma(\tE,\cM_n)\stackrel{\sim}{\rightarrow}\rR\Gamma(\tE_s,\cM_n). 
\end{equation}
Par suite, compte tenu de (\cite{agt} III.7.5 et III.7.11), tout $\bvocB$-module $\cM$ de $\tE^{\mN^\circ}$ est naturellement un objet de $\tE_s^{\mN^\circ}$,
et on a un isomorphisme canonique 
\begin{equation}\label{sld2c}
\rR\Gamma(\tE^{\mN^\circ},\cM)\stackrel{\sim}{\rightarrow}\rR\Gamma(\tE_s^{\mN^\circ},\cM). 
\end{equation}
Il s'ensuit que tout ind-$\bvocB$-module de $\tE^{\mN^\circ}$ est naturellement un ind-$\bvocB$-module de $\tE_s^{\mN^\circ}$ \eqref{indsh4d}.

\begin{defi}\label{sld3}
\
\begin{itemize}
\item[(i)] On dit qu'un $\bvmZ_p$-module $M=(M_n)_{n\in \mN}$ de $(\oX^\circ_\et)^{\mN^\circ}$ est un {\em système local} 
(ou que $M$ est un {\em $\bvmZ_p$-système local}) si les deux conditions suivantes sont satisfaites:
\begin{itemize}
\item[(a)] $M$ est $p$-adique, autrement dit, pour tous entiers $n\geq m\geq 0$, 
le morphisme $M_n/p^mM_n\rightarrow M_m$ déduit du morphisme de transition $M_n\rightarrow M_m$, est un isomorphisme; 
\item[(b)]  pour tout entier $n\geq 0$, le $\mZ/p^n\mZ$-module $M_n$ de $\oX^\circ_\et$ est localement constant constructible. 
\end{itemize}
\item[(ii)] On dit qu'un $\bvmZ_{p,\mQ}$-module de $(\oX^\circ_\et)^{\mN^\circ}$ est un {\em système local} s'il est isomorphe 
à $M_\mQ$ pour un $\bvmZ_p$-système local $M$ \eqref{sld1a}.
\end{itemize}
\end{defi}

\begin{defi}\label{sld5}
\
\begin{itemize}
\item[(i)] On dit qu'un ind-$\bvmZ_p$-module $M$ de $(\oX^\circ_\et)^{\mN^\circ}$ est {\em de Dolbeault} si le ind-$\bvocB$-module 
$\rI\bvpsi_*(M)\otimes_{\bvmZ_p}\bvocB$ est de Dolbeault \eqref{indmdlb5}.
\item[(ii)] On dit qu'un $\bvmZ_{p,\mQ}$-système local $M$ de $(\oX^\circ_\et)^{\mN^\circ}$ est {\em de Dolbeault} (resp. {\em de Hodge-Tate}) si le $\bvocB_\mQ$-module 
$\bvpsi_{\mQ*}(M)\otimes_{\bvmZ_{p,\mQ}}\bvocB_\mQ$ est de Dolbeault \eqref{aspglob1} (resp. de Hodge-Tate \eqref{mht1}).
\end{itemize}
\end{defi}

D'après \eqref{indsh13a} \eqref{indsh20d} et \ref{aspglob10}, pour qu'un $\bvmZ_{p,\mQ}$-système local $M$ de $(\oX^\circ_\et)^{\mN^\circ}$ soit de Dolbeault, 
il faut et il suffit que le ind-$\bvmZ_p$-module $\upalpha_{\bvmZ_p}(M)$ soit de Dolbeault.

\begin{rema}
Tout $\bvmZ_{p,\mQ}$-système local de Dolbeault $M$ de $(\oX^\circ_\et)^{\mN^\circ}$ est fortement de Dolbeault \eqref{aspglob1}. 
En effet, le $\bvocB_\mQ$-module $\bvpsi_{\mQ*}(M)\otimes_{\bvmZ_{p,\mQ}}\bvocB_\mQ$ est adique de type fini \eqref{ahttf60} 
d'après (\cite{agt} III.7.5, VI.9.20 et VI.10.9(iii)). 
\end{rema}

\begin{teo}\label{sld4}
Soient $M=(M_n)_{n\geq 0}$ un $\bvmZ_p$-système local de $(\oX^\circ_\et)^{\mN^\circ}$, $\cM=\bvpsi_*(M)\otimes_{\bvmZ_p}\bvocB$, 
$\mK^\bullet$ le complexe de Dolbeault du $\co_\fX[\frac 1 p]$-module de Higgs $\cH_\mQ(\cM_\mQ)$ \eqref{aspglob11b}. 
Supposons le morphisme $f\colon X\rightarrow S$ propre et le $\bvmZ_{p,\mQ}$-système local $M_\mQ$ de Dolbeault, i.e., le $\bvocB_\mQ$-module $\cM_\mQ$ de Dolbeault.
Alors, il existe une suite spectrale canonique
\begin{equation}\label{sld4a}
\rE_2^{i,j}=\rH^i(X_s,\rH^j(\mK^\bullet))\Rightarrow \rH^{i+j}((\oX^\circ_\et)^{\mN^\circ},M)\otimes_{\mZ_p}C.
\end{equation}
\end{teo}

En effet, en vertu de \ref{aspglob22} et \eqref{sld2c}, la suite spectrale de Cartan-Leray 
\begin{equation}\label{sld4b}
\rE_2^{i,j}=\rH^i(X_s,\rR^j\hupsigma_*(\cM))\Rightarrow \rH^{i+j}(\tE_s^{\mN^\circ},\cM)
\end{equation}
induit une suite spectrale 
\begin{equation}\label{sld4c}
\rE_2^{i,j}=\rH^i(X_s,\rH^j(\mK^\bullet))\Rightarrow \rH^{i+j}(\tE^{\mN^\circ},\cM)\otimes_\mZ\mQ.
\end{equation}

Pour tout $n\geq 0$, posons $\cM_n=\psi_*(M_n)\otimes_{\mZ_p}\ocB$, de sorte que $\cM=(\cM_n)_{n\geq 0}$ (\cite{agt} VI.7.5). 
D'après (\cite{agt} VI.7.10), pour tout entier $q\geq 0$ on a une suite exacte
\begin{equation}\label{sld4d}
0\rightarrow \rR^1\underset{\underset{n\geq 1}{\longleftarrow}}\lim\ \rH^{q-1}(\tE,\cM_n)\rightarrow 
\rH^q(\tE^{\mN^\circ},\cM)\rightarrow \underset{\underset{n\geq 1}{\longleftarrow}}\lim\ \rH^q(\tE,\cM_n)\rightarrow 0.
\end{equation}
En vertu de (\cite{ag} 4.8.13), pour tout $n\geq 0$, on a un morphisme canonique 
\begin{equation}\label{sld4e}
u^q_n\colon \rH^q(\oX^\circ_{\et},M_n)\otimes_{\mZ_p}\co_\oK\rightarrow \rH^q(\tE,\cM_n),
\end{equation}
qui est un $\alpha$-isomorphisme. Les morphismes 
\begin{equation}\label{sld4f}
\underset{\underset{n\geq 1}{\longleftarrow}}\lim\ u^q_n \ \ \ {\rm et}\ \ \ \rR^1\underset{\underset{n\geq 1}{\longleftarrow}}\lim\ u^q_n
\end{equation}
sont donc des $\alpha$-isomorphismes (\cite{gr} 2.4.2(ii)). 

De même, on a une suite exacte
\begin{equation}\label{sld4g}
0\rightarrow \rR^1\underset{\underset{n\geq 0}{\longleftarrow}}\lim\ \rH^{q-1}(\oX^\circ_\et,M_n)\rightarrow 
\rH^q((\oX^\circ_\et)^{\mN^\circ},M)\rightarrow \underset{\underset{n\geq 0}{\longleftarrow}}\lim\ \rH^q(\oX^\circ_\et,M_n)\rightarrow 0.
\end{equation}
Les groupes $\rH^q(\oX^\circ_{\et},M_n)$ étant finis en vertu de (\cite{sga45} Th.finitude 1.1), 
le système projectif $(\rH^q(\oX^\circ_{\et},M_n))_{n\geq 1}$ vérifie la condition de Mittag-Leffler. 
Le morphisme canonique 
\begin{equation}\label{sld4h}
\rH^q((\oX^\circ_\et)^{\mN^\circ},M)\rightarrow \underset{\underset{n\geq 0}{\longleftarrow}}\lim\ \rH^q(\oX^\circ_\et,M_n)
\end{equation}
est donc un isomorphisme d'après (\cite{jannsen} 1.15).
On en déduit aussi que le morphisme canonique
\begin{equation}\label{sld4i}
\rH^q(\tE_s^{\mN^\circ},\cM)\rightarrow \underset{\underset{n\geq 1}{\longleftarrow}}\lim\ \rH^q(\tE_s,\cM_n)
\end{equation}
est un $\alpha$-isomorphisme compte tenu de (\cite{roos} théo.~1). 

D'après (\cite{sga5} VI 2.2.2 et la remarque après 2.2.3), le système projectif $\rH^q=(\rH^q(\oX^\circ_{\et},M_n))_{n\geq 0}$ est AR-$p$-adique constructible, 
autrement dit, il existe un système projectif $p$-adique noethérien de groupes abéliens $A^q=(A^q_n)_{n\geq 0}$ et un AR-isomorphisme 
$A^q\rightarrow \rH^q$ (\cite{sga5} V 3.2.2). On en déduit un isomorphisme
\begin{equation}
\underset{\underset{n\geq 0}{\longleftarrow}}\lim\ A^q_n\stackrel{\sim}{\rightarrow} \underset{\underset{n\geq 0}{\longleftarrow}}\lim\ 
\rH^q(\oX^\circ_{\et},M_n).
\end{equation}
En particulier, le $\mZ_p$-module $\rH^q((\oX^\circ_\et)^{\mN^\circ},M)$ est de type fini. 
Le morphisme de systèmes projectifs $(A^q_n\otimes_{\mZ_p}\co_C)_{n\geq 0}\rightarrow (\rH^q(\oX^\circ_{\et},M_n)\otimes_{\mZ_p}\co_C)_{n\geq 0}$ est aussi 
un AR-isomorphisme. Comme la source est un système $p$-adique, on en déduit que le morphisme canonique 
\begin{equation}\label{sshtrl300b}
\rH^q((\oX^\circ_\et)^{\mN^\circ},M)\otimes_{\mZ_p}\co_C \rightarrow \underset{\underset{n\geq 0}{\longleftarrow}}\lim\  
\rH^q(\oX^\circ_{\et},M_n) \otimes_{\mZ_p}\co_C
\end{equation}
est un isomorphisme. La proposition résulte alors de \eqref{sld4c}, \eqref{sld4e}, \eqref{sld4h}, \eqref{sld4i} et \eqref{sshtrl300b}.

\begin{remas}\label{sld6}\
Dans \ref{sld4}, si l'on prend $M=\bvmZ_p$, alors $\cM=\bvocB$, le $\bvocB_\mQ$-module $\bvocB_\mQ$ est de Dolbeault et 
$\cH(\bvocB_\mQ)$ est égal à $\co_\fX[\frac 1 p]$ muni du champ de Higgs nul \eqref{aspglob12}. 
La suite spectrale \eqref{dolb10a} n'est autre que la suite spectrale de Hodge-Tate  (\cite{ag} 6.4.6). 
On notera que la construction \ref{sld4} de cette suite spectrale montre directement qu'elle dégénère en $\rE_2$ et que 
la filtration aboutissement est scindée sans utiliser le théorème de Tate sur la cohomologie galoisienne de $C(j)$. 
Cette construction s'applique en particulier en prenant pour $(\tX,\cM_\tX)$ dans le cas relatif \eqref{definf10} la déformation triviale \eqref{defing13}. 
\end{remas}

\chapter{Cohomologies relatives des algèbres de Higgs-Tate. \'Etude locale}

\section{Hypothèses et notations. Rappels de cohomologie galoisienne relative}\label{hypmdef}

\subsection{}\label{hypmdef1}
Dans ce chapitre, $K$ désigne un corps de valuation discrète complet de 
caractéristique $0$, à corps résiduel {\em algébriquement clos} $k$ de caractéristique $p>0$,  
$\co_K$ l'anneau de valuation de $K$, 
$\oK$ une clôture algébrique de $K$, $\co_\oK$ la clôture intégrale de $\co_K$ dans $\oK$,
$\fm_\oK$ l'idéal maximal de $\co_\oK$ et $G_K$ le groupe de Galois de $\oK$ sur $K$.
On note $\co_C$ le séparé complété $p$-adique de $\co_\oK$, $\fm_C$ son idéal maximal,
$C$ son corps des fractions et $v$ sa valuation, normalisée par $v(p)=1$. 
On désigne par $\hmZ(1)$ et $\mZ_p(1)$ les $\mZ[G_K]$-modules 
\begin{eqnarray}
\hmZ(1)&=&\underset{\underset{n\geq 1}{\longleftarrow}}{\lim}\ \mu_{n}(\co_{\oK}),\label{hypmdef1a}\\
\mZ_p(1)&=&\underset{\underset{n\geq 0}{\longleftarrow}}{\lim}\ \mu_{p^n}(\co_{\oK}),\label{hypmdef1c}\\
\mu_{p^\infty}(\co_\oK)&=&\underset{\underset{n\geq 0}{\longrightarrow}}{\lim}\ \mu_{p^n}(\co_{\oK}),
\end{eqnarray}  
où $\mu_n(\co_{\oK})$ désigne le sous-groupe des racines $n$-ièmes de l'unité dans $\co_\oK$. 
Pour tout $\mZ_p[G_K]$-module $M$ et tout entier $n$, on pose $M(n)=M\otimes_{\mZ_p}\mZ_p(1)^{\otimes n}$.  

Pour tout $\mZ_p$-module $A$, on note $\hA$ son séparé complété $p$-adique.

On pose $S=\Spec(\co_K)$, $\oS=\Spec(\co_\oK)$ et $\coS=\Spec(\co_C)$. 
On note $s$ (resp.  $\eta$, resp. $\oeta$) le point fermé de $S$ (resp.  générique de $S$, resp. générique de $\oS$).
Pour tout entier $n\geq 1$, on pose $S_n=\Spec(\co_K/p^n\co_K)$. Pour tout $S$-schéma $X$, on pose 
\begin{equation}\label{hypmdef1b}
\oX=X\times_S\oS,  \ \ \ \coX=X\times_S\coS \ \ \ {\rm et}\ \ \  X_n=X\times_SS_n.
\end{equation} 

On munit $S$ de la structure logarithmique $\cM_S$ définie par son point fermé, 
et $\oS$ et $\coS$ des structures logarithmiques $\cM_\oS$ et $\cM_\coS$ images inverses de $\cM_S$.

\subsection{}\label{hypmdef100}
Comme $\co_\oK$ est un anneau de valuation non discrète de hauteur $1$, 
il est loisible de développer la $\alpha$-algèbre (ou presque-algèbre) sur cet anneau (\cite{ag} 2.10.1) (cf. \cite{ag} 2.6-2.10).   
On choisit un système compatible $(\beta_n)_{n>0}$ 
de racines $n$-ièmes de $p$ dans $\co_\oK$. Pour tout nombre rationnel $\varepsilon>0$, 
on pose $p^\varepsilon=(\beta_n)^{\varepsilon n}$, où $n$ est un entier $>0$ tel que $\varepsilon n$ soit entier.

\subsection{}\label{hypmdef2}
Dans ce chapitre, $f\colon (X,\cM_X)\rightarrow (S,\cM_S)$ et $f'\colon (X',\cM_{X'})\rightarrow (S,\cM_S)$
désignent des morphismes {\em adéquats} de schémas logarithmiques (\cite{agt} III.4.7) et 
\begin{equation}\label{hypmdef2a}
g\colon (X',\cM_{X'})\rightarrow (X,\cM_X)
\end{equation} 
un $(S,\cM_S)$-morphisme lisse et saturé. On suppose que les schémas $X=\Spec(R)$ et $X'=\Spec(R')$ sont affines et que $X'_s$ est non-vide. 
De plus, sauf dans \ref{mtht101}, on suppose que le morphisme $g$ admet une carte relativement adéquate (\cite{ag} 5.1.11)  
\begin{equation}\label{hypmdef2d}
((P',\gamma'),(P,\gamma),(\mN,\iota),\vartheta\colon \mN\rightarrow P, h\colon P\rightarrow P'),
\end{equation} 
que l'on fixe. On pose $\pi=\iota(1)$ qui est une uniformisante de $\co_K$.  

On désigne par $X^\circ$ le sous-schéma ouvert maximal de $X$
où la structure logarithmique $\cM_X$ est triviale; c'est un sous-schéma ouvert de $X_\eta$. 
L'immersion $j\colon X^\circ\rightarrow X$ est schématiquement dominante (\cite{agt} III.4.2(iv)). 
Pour tout $X$-schéma $U$, on pose  
\begin{equation}\label{hypmdef2b}
U^\circ=U\times_XX^\circ.
\end{equation} 

On désigne par $X'^\rhd$ le sous-schéma ouvert maximal de $X'$
où la structure logarithmique $\cM_{X'}$ est triviale; c'est un sous-schéma ouvert de $X'^\circ=X'\times_XX^\circ$.
L'immersion $j'\colon X'^\rhd\rightarrow X'$ est schématiquement dominante. 
Pour tout $X'$-schéma $U'$, on pose  
\begin{equation}\label{hypmdef2c}
U'^\rhd=U'\times_{X'}X'^\rhd.
\end{equation}

\subsection{}\label{hypmdef4}
On pose
\begin{equation}\label{hypmdef4a}
\tOmega^1_{R/\co_K}=\Gamma(X,\Omega^1_{(X,\cM_X)/(S,\cM_S)}),
\end{equation}
qui est un $R$-module libre de type fini (\cite{ag} (4.5.9.2)). De même, on pose
\begin{eqnarray}
\tOmega^1_{R'/\co_K}&=&\Gamma(X',\Omega^1_{(X',\cM_{X'})/(S,\cM_S)}),\label{hypmdef4b}\\
\tOmega^1_{R'/R}&=&\Gamma(X',\Omega^1_{(X',\cM_{X'})/(X,\cM_X)}), \label{hypmdef4c}
\end{eqnarray}
qui sont des $R'$-modules libres de type fini. On a alors une suite exacte scindée de $R'$-modules 
\begin{equation}\label{hypmdef4d}
0\rightarrow \tOmega^1_{R/\co_K}\otimes_RR'\rightarrow \tOmega^1_{R'/\co_K}\rightarrow \tOmega^1_{R'/R}\rightarrow 0.
\end{equation}

\subsection{}\label{eccr20}
Pour tout entier $n\geq 1$, on pose
\begin{eqnarray}\label{eccr20a}
\co_{K_n}=\co_K[\zeta]/(\zeta^{n}-\pi),
\end{eqnarray}
qui est un anneau de valuation discrète \eqref{hypmdef2}. On note $K_n$ le corps des fractions de $\co_{K_n}$
et $\pi_n$ la classe de $\zeta$ dans $\co_{K_n}$, qui est une uniformisante de $\co_{K_n}$.  
On pose $S^{(n)}=\Spec(\co_{K_n})$
que l'on munit de la structure logarithmique $\cM_{S^{(n)}}$ définie par son point fermé. 
On désigne par $\iota_n\colon \mN\rightarrow \Gamma(S^{(n)},\cM_{S^{(n)}})$
l'homomorphisme défini par $\iota_n(1)=\pi_n$; c'est une carte pour $(S^{(n)},\cM_{S^{(n)}})$.  

Considérons le système inductif de monoïdes $(\mN^{(n)})_{n\geq 1}$, 
indexé par l'ensemble $\mZ_{\geq 1}$ ordonné par la relation de divisibilité, 
défini par $\mN^{(n)}=\mN$ pour tout $n\geq 1$ et dont l'homomorphisme de transition
$\mN^{(n)}\rightarrow \mN^{(mn)}$ (pour $m,n\geq 1$) est l'endomorphisme de Frobenius d'ordre $m$ 
de $\mN$ ({\em i.e.}, l'élévation à la puissance $m$-ième). On notera $\mN^{(1)}$ simplement $\mN$. Les  schémas logarithmiques
$(S^{(n)},\cM_{S^{(n)}})_{n\geq 1}$ forment naturellement un système projectif.
Pour tous entiers $m, n\geq 1$, avec les notations de  \ref{notconv2}, on a un diagramme cartésien de morphismes de schémas logarithmiques
\begin{equation}\label{eccr20b}
\xymatrix{
{(S^{(mn)},\cM_{S^{(mn)}})}\ar[r]^-(0.5){\iota^a_{mn}}\ar[d]&{\bA_{\mN^{(mn)}}}\ar[d]\\
{(S^{(n)},\cM_{S^{(n)}})}\ar[r]^-(0.5){\iota^a_n}&{\bA_{\mN^{(n)}}}}
\end{equation}
où $\iota^a_n$ (resp. $\iota^a_{mn}$) est le morphisme associé à $\iota_n$ (resp. $\iota_{mn}$) et $\bA_{\mN^{(n)}}$ est défini dans \ref{notconv2} (cf. \cite{agt} II.5.13). 

\subsection{}\label{eccr2}
On désigne par $(P^{(n)})_{n\geq 1}$ le système inductif de monoïdes,
indexé par l'ensemble $\mZ_{\geq 1}$ ordonné par la relation de divisibilité, 
défini par $P^{(n)}=P$ pour tout $n\geq 1$ et dont l'homomorphisme de transition
$i_{n,mn}\colon P^{(n)}\rightarrow P^{(mn)}$ (pour $m, n\geq 1$) est l'endomorphisme 
de Frobenius d'ordre $m$ de $P$ ({\em i.e.}, l'élévation à la puissance $m$-ième). On note $P^{(1)}$ simplement $P$.

Pour tout $n\geq 1$, on pose (avec les notations de \ref{notconv2})
\begin{equation}\label{eccr2a}
(X^{(n)},\cM_{X^{(n)}})=(X,\cM_{X})\times_{\bA_P}\bA_{P^{(n)}}.
\end{equation}
On rappelle \eqref{cad4e} qu'il existe un unique morphisme 
\begin{equation}\label{eccr2b}
f^{(n)}\colon (X^{(n)},\cM_{X^{(n)}})\rightarrow (S^{(n)},\cM_{S^{(n)}})
\end{equation} 
qui s'insère dans le diagramme commutatif
\begin{equation}\label{eccr2c}
\xymatrix{
{(X^{(n)},\cM_{X^{(n)}})}\ar[rrr]\ar[ddd]\ar[rd]_{f^{(n)}}&&&
{\bA_{P^{(n)}}}\ar[ld]^{\bA_{\vartheta}}\ar[ddd]\\
&{(S^{(n)},\cM_{S^{(n)}})}\ar[r]^-(0.5){\iota^a_n}\ar[d]&{\bA_{\mN^{(n)}}}\ar[d]&\\
&{(S,\cM_{S})}\ar[r]^-(0.5){\iota^a}&{\bA_\mN}&\\
{(X,\cM_{X})}\ar[rrr]\ar[ru]^f&&&{\bA_P}\ar[lu]_{\bA_\vartheta}}
\end{equation}

\subsection{}\label{eccr3}
On désigne par $(P'^{[n]})_{n\geq 1}$ le système inductif de monoïdes 
indexé par l'ensemble $\mZ_{\geq 1}$ ordonné par la relation de divisibilité, 
défini par $P'^{[n]}=P'$ pour tout $n\geq 1$ et dont l'homomorphisme de transition
$P'^{[n]}\rightarrow P'^{[mn]}$ (pour $m, n\geq 1$) est l'endomorphisme 
de Frobenius d'ordre $m$ de $P'$. On note $P'^{[1]}$ simplement $P'$. 

Pour tout $n\geq 1$, on pose (avec les notations de \ref{notconv2})
\begin{equation}\label{eccr3a}
(X'^{[n]},\cM_{X'^{[n]}})=(X',\cM_{X'})\times_{\bA_{P'}}\bA_{P'^{[n]}}.
\end{equation}
Il existe alors un unique morphisme
\begin{equation}\label{eccr3b}
f'^{[n]}\colon (X'^{[n]},\cM_{X'^{[n]}})\rightarrow (S^{(n)},\cM_{S^{(n)}}) 
\end{equation} 
qui s'insère dans le diagramme commutatif 
\begin{equation}\label{eccr3c}
\xymatrix{
{(X'^{[n]},\cM_{X'^{[n]}})}\ar[rrr]\ar[ddd]\ar[rd]_{f'^{[n]}}&&&
{\bA_{P'^{[n]}}}\ar[ld]^{\bA_{h\circ \vartheta}}\ar[ddd]\\
&{(S^{(n)},\cM_{S^{(n)}})}\ar[r]\ar[d]&{\bA_{\mN^{(n)}}}\ar[d]&\\
&{(S,\cM_{S})}\ar[r]&{\bA_\mN}&\\
{(X',\cM_{X'})}\ar[rrr]\ar[ru]^{f'}&&&{\bA_{P'}}\ar[lu]_{\bA_{h\circ \vartheta}}}
\end{equation}

\subsection{}\label{eccr4}
Pour tout $n\geq 1$, on pose 
\begin{equation}\label{eccr4a}
(X'^{(n)},\cM_{X'^{(n)}})=(X^{(n)},\cM_{X^{(n)}})\times_{(X,\cM_X)}(X',\cM_{X'}),
\end{equation}
le produit étant indifféremment pris dans la catégorie des schémas logarithmiques ou dans celle des schémas logarithmiques saturés (\cite{agt} II.5.20). 
On note  
\begin{equation}\label{eccr4b}
f'^{(n)}\colon (X'^{(n)},\cM_{X'^{(n)}})\rightarrow (S^{(n)},\cM_{S^{(n)}})
\end{equation}
le morphisme déduit de $f^{(n)}$ \eqref{eccr2b}. 

Il existe un unique $(S^{(n)},\cM_{S^{(n)}})$-morphisme
\begin{equation}\label{eccr4g}
g^{(n)}\colon (X'^{[n]},\cM_{X'^{[n]}})\rightarrow (X^{(n)},\cM_{X^{(n)}})
\end{equation}
au-dessus de $g$ et du morphisme $\bA_{h}\colon \bA_{P'^{[n]}}\rightarrow \bA_{P^{(n)}}$. Celui-ci induit un morphisme canonique 
\begin{equation}\label{eccr4h}
(X'^{[n]},\cM_{X'^{[n]}})\rightarrow (X'^{(n)},\cM_{X'^{(n)}})
\end{equation}
au-dessus de $(X',\cM_{X'})$ et de $(X^{(n)},\cM_{X^{(n)}})$.

\subsection{}\label{eccr40}
Soient $\oy'$ un point géométrique de $\oX'^\rhd$ et $\oy$ son image dans $\oX^\circ$ (cf. \ref{hypmdef1} et \ref{hypmdef2}). 
Les schémas $\oX$ et $\oX'$ étant localement irréductibles d'après (\cite{ag} 4.2.7(iii)),  
ils sont les sommes des schémas induits sur leurs composantes irréductibles. 
On note $\oX^\star$ (resp. $\oX'^\star$)
la composante irréductible de $\oX$ (resp. $\oX'$) contenant $\oy$ (resp. $\oy'$). 
De même, $\oX^\circ$ (resp. $\oX'^\rhd$) est la somme des schémas induits sur ses composantes irréductibles
et $\oX^{\star \circ}=\oX^\star\times_{X}X^\circ$ (resp. $\oX'^{\star \rhd}=\oX'^\star\times_{X'}X'^\rhd$) est la composante irréductible de 
$\oX^\circ$ (resp. $\oX'^\rhd$) contenant $\oy$ (resp. $\oy'$).

On désigne par $\Delta$ le groupe profini $\pi_1(\oX^{\star \circ},\oy)$ et par $(V_i)_{i\in I}$ le revêtement universel normalisé de 
$\oX^{\star \circ}$ en $\oy$ (\cite{agt} VI.9.8). Pour tout $i\in I$, on désigne par $\oX^{V_i}$ la fermeture intégrale de $\oX$ dans $V_i$. 
On pose 
\begin{equation}\label{eccr40b}
\oR=\underset{\underset{i\in I}{\longrightarrow}}{\lim}\  \Gamma(\oX^{V_i},\co_{\oX^{V_i}}),
\end{equation}
qui est naturellement muni d'une action discrète de $\Delta$ par des automorphismes d'anneaux. 
On notera que $\oR$ est la représentation discrète $\oR^{\oy}_{X}$ de $\Delta$ définie dans \eqref{ahttf44b} (\cite{ag} (4.1.9.3)).

On désigne par $\Delta'$ le groupe profini $\pi_1(\oX'^{\star \rhd},\oy')$ et par 
$(W_j)_{j\in J}$ le revêtement universel normalisé de $\oX'^{\star \rhd}$ en $\oy'$. 
Pour tout $j\in J$, on désigne par $\oX'^{W_j}$ la fermeture intégrale de $\oX'$ dans $W_j$. 
On pose  
\begin{equation}\label{eccr40a}
\oR'=\underset{\underset{j\in J}{\longrightarrow}}{\lim}\  \Gamma(\oX'^{W_j},\co_{\oX'^{W_j}}),
\end{equation}
qui est naturellement muni d'une action discrète de $\Delta'$ par des automorphismes d'anneaux. 
On notera que $\oR'$ est la représentation discrète $\oR'^{\oy'}_{X'}$ de $\Delta'$ définie dans \eqref{ahttf44b}. 

Pour tout $i\in I$, on a un $\oX^\circ$-morphisme canonique $\oy\rightarrow V_i$.
On en déduit un $\oX'^\rhd$-morphisme $\oy'\rightarrow V_i\times_{\oX^\circ}\oX'^\rhd$.
Le schéma $V_i\times_{\oX^\circ}\oX'^\rhd$ étant localement irréductible,  
il est la somme des schémas induits sur ses composantes irréductibles. 
On désigne par $V'_i$ la composante irréductible de $V_i\times_{\oX^\circ}\oX'^\rhd$ contenant l'image de $\oy'$
et par $\oX'^{V'_i}$ la fermeture intégrale de $\oX'$ dans $V'_i$.
Les schémas $(V'_i)_{i\in I}$ forment naturellement un système projectif de revêtements étales finis connexes $\oy'$-pointés de $\oX'^{\star \rhd}$.
Pour tout $i\in I$, on note $\Pi_i$ le sous-groupe ouvert de $\Delta'$ correspondant à $V'_i$, autrement dit le noyau de l'action canonique 
de $\Delta'$ sur la fibre de $V'_i$ au-dessus de $\oy'$ (\cite{ag} (2.1.20.2)). 
On désigne par $\Pi$ le sous-groupe fermé de $\Delta'$ défini par
\begin{equation}\label{eccr40c}
\Pi=\cap_{i\in I}\Pi_i.
\end{equation}
Posant $\Delta^\intern=\Delta'/\Pi$, on a un homomorphisme canonique $\Delta'\rightarrow \Delta$ qui se factorise à travers un homomorphisme 
injectif $\Delta^\intern\rightarrow \Delta$.

Considérons l'anneau 
\begin{equation}\label{eccr40e}
\oR^\intern=\underset{\underset{i\in I}{\longrightarrow}}{\lim}\  \Gamma(\oX'^{V'_i},\co_{\oX'^{V'_i}}),
\end{equation}
qui est naturellement muni d'une action discrète de $\Delta^\intern$ par des automorphismes d'anneaux. 
On a un homomorphisme canonique $\Delta^\intern$-équivariant $\oR\rightarrow \oR^\intern$. 
 
Pour tous $i\in I$ et $j\in J$, il existe au plus un morphisme de $\oX'^{\star \rhd}$-schémas pointés $W_j\rightarrow V'_i$. 
De plus, pour tout $i\in I$, il existe $j\in J$ et un morphisme de $\oX'^{\star \rhd}$-schémas pointés $W_j\rightarrow V'_i$. 
On a donc un homomorphisme canonique
\begin{equation}\label{eccr40d}
\underset{\underset{i\in I}{\longrightarrow}}{\lim}\  \Gamma(\oX'^{V'_i},\co_{\oX'^{V'_i}})\rightarrow 
\underset{\underset{j\in J}{\longrightarrow}}{\lim}\  \Gamma(\oX'^{W_j},\co_{\oX'^{W_j}}).
\end{equation}
On en déduit un homomorphisme canonique $\Delta'$-équivariant de $R'_1$-algèbres
\begin{equation}\label{eccr40f}
\oR^\intern\rightarrow \oR'.
\end{equation}

\subsection{}\label{eccr41}
D'après \eqref{eccr4h}, on a des morphismes canoniques de systèmes projectifs de schémas logarithmiques, indexés par l'ensemble 
$\mZ_{\geq 1}$ ordonné par la relation de divisibilité, 
\begin{equation}\label{eccr41a}
(X'^{[n]}, \cM_{X'^{[n]}})\rightarrow (X'^{(n)}, \cM_{X'^{(n)}})\rightarrow (X^{(n)}, \cM_{X^{(n)}})\rightarrow (S^{(n)},\cM_{S^{(n)}}).
\end{equation}
Pour tous entiers $m,n\geq 1$, le morphisme canonique $X'^{[mn]}\rightarrow X'^{[n]}$ est fini et surjectif.
Il existe donc un $X'$-morphisme (\cite{ega4} 8.3.8(i))
\begin{equation}\label{eccr41b}
\oy'\rightarrow \underset{\underset{n\geq 1}{\longleftarrow}}{\lim}\ X'^{[n]}.
\end{equation} 
On fixe un tel morphisme dans la suite de cette section.
Celui-ci induit un $S$-morphisme
\begin{equation}\label{eccr41c}
\oS\rightarrow \underset{\underset{n\geq 1}{\longleftarrow}}{\lim}\ S^{(n)}.
\end{equation}
Pour tout entier $n\geq 1$, on pose
\begin{eqnarray}
\oX'^{[n]}= X'^{[n]}\times_{S^{(n)}}\oS,&& \oX'^{[n]\rhd}=\oX'^{[n]}\times_{X'}X'^\rhd,\label{eccr41d1}\\
\oX'^{(n)}= X'^{(n)}\times_{S^{(n)}}\oS,&& \oX'^{(n)\rhd}=\oX'^{(n)}\times_{X'}X'^\rhd,\label{eccr41d2}\\
\oX^{(n)}= X^{(n)}\times_{S^{(n)}}\oS,&& \oX^{(n)\circ}=\oX^{(n)}\times_{X}X^\circ.\label{eccr41d3}
\end{eqnarray} 
Les morphismes \eqref{eccr41a} induisent des morphismes de systèmes projectifs de $\oS$-schémas
\begin{equation}\label{eccr41e}
\oX'^{[n]}\rightarrow \oX'^{(n)}\rightarrow \oX^{(n)}.
\end{equation}
On déduit de \eqref{eccr41a} et \eqref{eccr41b} un $\oX'$-morphisme 
\begin{equation}\label{eccr41f}
\oy'\rightarrow \underset{\underset{n\geq 1}{\longleftarrow}}{\lim}\ \oX'^{[n]}.
\end{equation}

\subsection{}\label{eccr42}
Pour tout entier $n\geq 1$, le schéma $\oX'^{[n]}$ est normal et localement irréductible d'après (\cite{ag} 4.2.7(iii)). 
Il est donc la somme des schémas induits sur ses composantes irréductibles. 
On note $\oX'^{[n]\star}$ la composante irréductible de $\oX'^{[n]}$ contenant l'image de $\oy'$ \eqref{eccr41f}.
De même, $\oX'^{[n]\rhd}$ est la somme des schémas induits sur ses composantes irréductibles
et $\oX'^{[n]\star\rhd}=\oX'^{[n]\star}\times_{X'}X'^\rhd$ est la composante irréductible de $\oX'^{[n]\rhd}$ contenant l'image de $\oy'$.
On pose 
\begin{equation}\label{eccr42a}
R'_n=\Gamma(\oX'^{[n]\star},\co_{\oX'^{[n]}}).
\end{equation}
Les anneaux $(R'_n)_{n\geq 1}$ forment naturellement un système inductif. On pose 
\begin{eqnarray}
R'_\infty&=&\underset{\underset{n\geq 1}{\longrightarrow}}{\lim}\ R'_n,\label{eccr42b1}\\
R'_{p^\infty}&=&\underset{\underset{n\geq 0}{\longrightarrow}}{\lim}\ R'_{p^n}.\label{eccr42b2}
\end{eqnarray}
Le morphisme \eqref{eccr41f} induit des homomorphismes injectifs de $R'_1$-algèbres 
\begin{equation}\label{eccr42c}
R'_{p^\infty}\rightarrow R'_\infty\rightarrow \oR'.
\end{equation}

De même, pour tout entier $n\geq 1$, le schéma $\oX^{(n)}$ est la somme des schémas induits sur ses composantes irréductibles. 
On note $\oX^{(n)\star}$ la composante irréductible de $\oX^{(n)}$ contenant l'image de $\oy'$ \eqref{eccr41f}.
Le schéma $\oX^{(n)\circ}$ est la somme des schémas induits sur ses composantes irréductibles
et $\oX^{(n)\star\circ}=\oX^{(n)\star}\times_XX^\circ$ est la composante irréductible de $\oX^{(n)\circ}$ contenant l'image de $\oy'$.
On pose 
\begin{equation}\label{eccr42d}
R_n=\Gamma(\oX^{(n)\star},\co_{\oX^{(n)}}).
\end{equation}
Les anneaux $(R_n)_{n\geq 1}$ forment naturellement un système inductif. On pose 
\begin{eqnarray}
R_\infty&=&\underset{\underset{n\geq 1}{\longrightarrow}}{\lim}\ R_n,\label{eccr42e1}\\
R_{p^\infty}&=&\underset{\underset{n\geq 0}{\longrightarrow}}{\lim}\ R_{p^n}.\label{eccr42e2}
\end{eqnarray}
Les morphismes \eqref{eccr41e} et \eqref{eccr41f} induisent des homomorphismes injectifs de $R_1$-algèbres 
\begin{equation}\label{eccr42f}
R_{p^\infty}\rightarrow R_\infty\rightarrow \oR.
\end{equation}

Pour tout entier $n\geq 1$, le schéma $\oX'^{(n)}$ est normal et localement irréductible d'après \eqref{eccr4a} et (\cite{agt} III.4.2(iii)). 
Il est donc la somme des schémas induits sur ses composantes irréductibles. 
Par ailleurs, l'immersion $\oX'^{(n)\rhd}\rightarrow \oX'^{(n)}$ est schématiquement dominante (\cite{ega4} 11.10.5). 
On note $\oX'^{(n)\star}$ la composante irréductible de $\oX'^{(n)}$ contenant l'image de $\oy'$ \eqref{eccr41f}.
Le schéma $\oX'^{(n)\rhd}$ est aussi la somme des schémas induits sur ses composantes irréductibles
et $\oX'^{(n)\star\rhd}=\oX'^{(n)\star}\times_{X'}X'^\rhd$  est la composante irréductible de $\oX'^{(n)\rhd}$ 
contenant l'image de $\oy'$. On pose 
\begin{equation}\label{eccr42g}
R^\intern_n=\Gamma(\oX'^{(n)\star},\co_{\oX'^{(n)}}).
\end{equation}
Les anneaux $(R^\intern_n)_{n\geq 1}$ forment naturellement un système inductif. On pose 
\begin{eqnarray}
R^\intern_\infty&=&\underset{\underset{n\geq 1}{\longrightarrow}}{\lim}\ R^\intern_n,\label{eccr42h1}\\
R^\intern_{p^\infty}&=&\underset{\underset{n\geq 0}{\longrightarrow}}{\lim}\ R^\intern_{p^n}.\label{eccr42h2}
\end{eqnarray}
Ce sont des anneaux normaux et intègres d'après (\cite{ega1n} 0.6.1.6(i) et 0.6.5.12(ii)). 
Compte tenu de la définition \eqref{eccr40e}, les morphismes \eqref{eccr41e} et \eqref{eccr41f} 
induisent des homomorphismes injectifs de $R^\intern_1$-algèbres 
\begin{equation}\label{eccr42i}
R^\intern_{p^\infty}\rightarrow R^\intern_\infty\rightarrow \oR^\intern.
\end{equation}

Les morphismes de \eqref{eccr41e} induisent des homomorphismes de systèmes inductifs d'algèbres 
\begin{equation}\label{eccr42j}
(R_n)_{n\geq 1}\rightarrow (R^\intern_n)_{n\geq 1}\rightarrow (R'_n)_{n\geq 1}.
\end{equation} 
L'homomorphisme $R^\intern_1\rightarrow R'_1$ est un isomorphisme.

\subsection{}\label{eccr43}
Soit $n$ un entier $\geq 1$. Il résulte de  (\cite{ag} 4.2.7(v)) et de la preuve de (\cite{agt} II.6.8(iv)) 
que le morphisme $\oX'^{[n]\star\rhd}\rightarrow \oX'^{\star\rhd}$ est un revêtement étale fini et galoisien de groupe $\Delta'_n$, 
un sous-groupe de $\Hom_\mZ(P'^\gp/h^\gp(\vartheta^\gp(\mZ)),\mu_{n}(\oK))$.  
Le groupe $\Delta'_n$ agit naturellement sur $R'_n$. Si $n$ est une puissance de $p$, le morphisme canonique
\begin{equation}\label{eccr43a}
\oX'^{[n]\star}\rightarrow\oX'^{[n]}\times_{\oX'}\oX'^\star
\end{equation}
est un isomorphisme en vertu de (\cite{agt} II.6.6(v)), et on a donc 
\begin{equation}
\Delta'_n\stackrel{\sim}{\rightarrow} \Hom_\mZ(P'^\gp/h^\gp(\vartheta^\gp(\mZ)),\mu_{n}(\oK)).
\end{equation}

Les groupes $(\Delta'_n)_{n\geq 1}$ forment naturellement un système projectif. 
On pose 
\begin{eqnarray}
\Delta'_\infty&=&\underset{\underset{n\geq 1}{\longleftarrow}}{\lim}\ \Delta'_n,\label{eccr43b1}\\
\Delta'_{p^\infty}&=&\underset{\underset{n\geq 0}{\longleftarrow}}{\lim}\ \Delta'_{p^n}.\label{eccr43b2}
\end{eqnarray}
On identifie $\Delta'_\infty$ (resp. $\Delta'_{p^\infty}$) au groupe de Galois de l'extension des corps de fractions de $R'_\infty$ sur $R'_1$
(resp. $R'_{p^\infty}$ sur $R'_1$). On a des homomorphismes canoniques \eqref{hypmdef1a}
\begin{equation}\label{eccr43b3}
\xymatrix{
{\Delta'_{\infty}}\ar@{->>}[d]\ar@{^(->}[r]&{\Hom(P'^\gp/h^\gp(\vartheta^\gp(\mZ)),\hmZ(1))}\ar[d]\\
{\Delta'_{p^\infty}}\ar[r]^-(0.5)\sim&{\Hom(P'^\gp/h^\gp(\vartheta^\gp(\mZ)),\mZ_p(1))}}
\end{equation}
Le noyau $\Sigma'_0$ de l'homomorphisme canonique $\Delta'_\infty\rightarrow \Delta'_{p^\infty}$
est un groupe profini d'ordre premier à $p$. 
Par ailleurs, le morphisme \eqref{eccr41f} détermine un homomorphisme surjectif $\Delta'\rightarrow \Delta'_\infty$ \eqref{eccr40}. 
On note $\Sigma'$ son noyau. Les homomorphismes \eqref{eccr42c} sont $\Delta'$-équivariants.

De même, le morphisme $\oX^{(n)\star\circ}\rightarrow \oX^{\star\circ}$ 
est un revêtement étale fini et galoisien de groupe $\Delta_n$, 
un sous-groupe de $\Hom_\mZ(P^\gp/\vartheta^\gp(\mZ),\mu_{n}(\oK))$.  
Le groupe $\Delta_n$ agit naturellement sur $R_n$. Si $n$ est une puissance de $p$, le morphisme canonique
\begin{equation}\label{eccr43c}
\oX^{(n)\star}\rightarrow\oX^{(n)}\times_{\oX}\oX^\star
\end{equation}
est un isomorphisme, et on a donc $\Delta_n\simeq \Hom_\mZ(P^\gp/\vartheta^\gp(\mZ),\mu_{n}(\oK))$.

Les groupes $(\Delta_n)_{n\geq 1}$ forment naturellement un système projectif. 
On pose 
\begin{eqnarray}
\Delta_\infty&=&\underset{\underset{n\geq 1}{\longleftarrow}}{\lim}\ \Delta_n,\label{eccr43d1}\\
\Delta_{p^\infty}&=&\underset{\underset{n\geq 0}{\longleftarrow}}{\lim}\ \Delta_{p^n}.\label{eccr43d2}
\end{eqnarray}
On identifie $\Delta_\infty$ (resp. $\Delta_{p^\infty}$) au groupe de Galois de l'extension des corps de fractions de $R_\infty$ sur $R_1$
(resp. $R_{p^\infty}$ sur $R_1$). On a des homomorphismes canoniques \eqref{hypmdef1a}
\begin{equation}\label{eccr43d3}
\xymatrix{
{\Delta_{\infty}}\ar@{->>}[d]\ar@{^(->}[r]&{\Hom(P^\gp/\vartheta^\gp(\mZ),\hmZ(1))}\ar[d]\\
{\Delta_{p^\infty}}\ar[r]^-(0.5)\sim&{\Hom(P^\gp/\vartheta^\gp(\mZ),\mZ_p(1))}}
\end{equation}
Le noyau $\Sigma_0$ de l'homomorphisme canonique $\Delta_\infty\rightarrow \Delta_{p^\infty}$
est un groupe profini d'ordre premier à $p$. 
Par ailleurs, les morphismes \eqref{eccr41e} et \eqref{eccr41f} déterminent un homomorphisme surjectif $\Delta\rightarrow \Delta_\infty$ \eqref{eccr40}.
On note $\Sigma$ son noyau. Les homomorphismes \eqref{eccr42f} sont $\Delta$-équivariants.

Il résulte de (\cite{ag} 5.2.8(iii)) que le morphisme $\oX'^{[n]\star\rhd}\rightarrow \oX'^{(n)\star\rhd}$ est un revêtement étale, 
fini et galoisien de groupe $\fS_n$, un sous-groupe de $\Hom_\mZ(P'^{\gp}/h^\gp(P^{\gp}),\mu_n(\oK))$. 
Il résulte de \ref{tpcg4} et \eqref{eccr4a} que le morphisme $\oX'^{(n)\star\rhd}\rightarrow \oX'^{\star\rhd}$ 
est un revêtement étale, fini et galoisien de groupe $\Delta^\intern_n$, un sous-groupe de $\Delta_n$. 
Le groupe $\Delta^\intern_n$ agit naturellement sur $R^\intern_n$. On a une suite exacte canonique
\begin{equation}\label{eccr43e}
0\rightarrow \fS_n\rightarrow \Delta'_n\rightarrow \Delta^\intern_n \rightarrow 0.
\end{equation}
Par ailleurs, le diagramme 
\begin{equation}\label{eccr43f}
\xymatrix{
\fS_n\ar[r]\ar[d]&{\Hom_\mZ(P'^{\gp}/h^\gp(P^{\gp}),\mu_n(\oK))}\ar[d]\\
{\Delta'_n}\ar[r]&{\Hom_\mZ(P'^\gp/h^\gp(\vartheta^\gp(\mZ)),\mu_{n}(\oK))}}
\end{equation}
où la flèche verticale de droite est l'homomorphisme canonique, est commutatif.
En vertu de (\cite{ag} 5.2.14(iv)), si $n$ est une puissance de $p$, le morphisme canonique
\begin{equation}\label{eccr43g}
\oX'^{(n)\star}\rightarrow \oX^{(n)}\times_{\oX}\oX'^\star
\end{equation}
est un isomorphisme, et on a 
\begin{eqnarray}
\fS_{n}&=&\Hom(P'^{\gp}/h^\gp(P^{\gp}),\mu_{n}(\oK)),\label{eccr43h1}\\
\Delta^\intern_{n}&=&\Hom(P^{\gp}/\vartheta^\gp(\mZ),\mu_{n}(\oK)).\label{eccr43h2}
\end{eqnarray}

Les groupes $(\fS_n)_{n\geq 1}$ et $(\Delta^\intern_n)_{n\geq 1}$ forment naturellement des systèmes projectifs. On pose 
\begin{eqnarray}
\fS_\infty&=&\underset{\underset{n\geq 1}{\longleftarrow}}{\lim}\ \fS_n,\label{eccr43i1}\\
\fS_{p^\infty}&=&\underset{\underset{n\geq 0}{\longleftarrow}}{\lim}\ \fS_{p^n}, \label{eccr43i2}\\
\Delta^\intern_\infty&=&\underset{\underset{n\geq 1}{\longleftarrow}}{\lim}\ \Delta^\intern_n,\label{eccr43i3}\\
\Delta^\intern_{p^\infty}&=&\underset{\underset{n\geq 0}{\longleftarrow}}{\lim}\ \Delta^\intern_{p^n}.\label{eccr43i4}
\end{eqnarray}
On identifie $\Delta^\intern_\infty$ (resp. $\Delta^\intern_{p^\infty}$) 
au groupe de Galois de l'extension des corps des fractions de $R^\intern_\infty$ (resp. $R^\intern_{p^\infty}$) sur $R'_1$. 
D'après la théorie de Galois, les suites 
\begin{equation}\label{eccr43j}
0\rightarrow \fS_\infty\rightarrow \Delta'_\infty\rightarrow \Delta^\intern_\infty\rightarrow 0,
\end{equation}
\begin{equation}\label{eccr43k}
0\rightarrow \fS_{p^\infty}\rightarrow \Delta'_{p^\infty}\rightarrow \Delta^\intern_{p^\infty}\rightarrow 0,
\end{equation}
déduites de \eqref{eccr43e} sont exactes. Le groupe $\fS_\infty$ (resp. $\fS_{p^\infty}$) s'identifie donc au groupe de Galois
de l'extension des corps de fractions de $R'_\infty$ sur $R^\intern_\infty$ (resp. $R'_{p^\infty}$ sur $R^\intern_{p^\infty}$).  
Comme le noyau de l'homomorphisme canonique $\Delta^\intern_\infty\rightarrow \Delta^\intern_{p^\infty}$ 
est un groupe profini d'ordre premier à $p$, 
on en déduit que l'homomorphisme canonique $\fS_\infty\rightarrow \fS_{p^\infty}$ est surjectif. 
On désigne par $\fN$ son noyau, qui est un groupe profini d'ordre premier à $p$, de sorte qu'on a la suite exacte
\begin{equation}\label{eccr43l}
0\rightarrow \fN\rightarrow \fS_\infty \rightarrow \fS_{p^\infty}\rightarrow 0.
\end{equation}

On désigne par $\fS$ le noyau de l'homomorphisme canonique $\Delta'\rightarrow \Delta^\intern_\infty$. 
On a alors un isomorphisme canonique 
\begin{equation}\label{eccr43m}
\fS\stackrel{\sim}{\rightarrow}\underset{\underset{n\geq 1}{\longleftarrow}}{\lim} \ \pi_1(\oX'^{(n)\star\rhd},\oy').
\end{equation}

\subsection{}\label{eccr34}
On peut résumer les principales constructions de cette section dans le diagramme commutatif suivant d'extensions d'anneaux intègres et normaux
\begin{equation}\label{eccr34a}
\xymatrix{
\oR\ar[r]&\oR^\intern\ar[rr]^-(0.5)\Pi&&\oR'\\
R_\infty\ar[r]\ar[u]^{\Sigma}&{R^\intern_\infty}\ar[r]_-(0.5){\fS_{p^\infty}}\ar[u]^{\Sigma^\intern}
\ar@/^1pc/[rr]^-(0.5){\fS_\infty}\ar@/^1pc/[rru]^-(0.3){\fS}&
{R^\intern_\infty\otimes_{R^\intern_{p^\infty}}R'_{p^\infty}}\ar[r]_-(0.4){\fN}&{R'_\infty}\ar[u]_{\Sigma'}\\
R_{p^\infty}\ar[r]\ar[u]^{\Sigma_0}&R^\intern_{p^\infty}\ar[r]^-(0.5){\fS_{p^\infty}}\ar[u]&R'_{p^\infty}\ar[ru]_{\Sigma'_0}\ar[u]&\\
R_1\ar[u]^{\Delta_{p^\infty}}\ar[r]&R^\intern_1\ar@{=}[r]\ar[u]^{\Delta^\intern_{p^\infty}}&R'_1\ar[u]_{\Delta'_{p^\infty}}&}
\end{equation}
où pour certaines extensions entières, on a noté le groupe de Galois de l'extension associée des corps de fractions. 
On observera que l'homomorphisme $R^\intern_1\rightarrow \oR^\intern$ induit une extension galoisienne des corps de fractions, 
et il en est donc de même de l'homomorphisme $R^\intern_\infty\rightarrow \oR^\intern$. 
On note $\oF^\intern$ (resp. $F^\intern_\infty$) le corps de fractions
de $\oR^\intern$ (resp. $R^\intern_\infty$) et $\Sigma^\intern$ le groupe de Galois de l'extension $\oF^\intern/F^\intern_\infty$.

\begin{prop}[\cite{ag} 5.2.15]\label{eccr9}
\
\begin{itemize}
\item[{\rm (i)}] L'anneau $R^\intern_\infty\otimes_{R^\intern_{p^\infty}}R'_{p^\infty}$ est normal et intègre, 
et on a un isomorphisme canonique \eqref{eccr43l}
\begin{equation}\label{eccr9d}
R^\intern_\infty\otimes_{R^\intern_{p^\infty}}R'_{p^\infty}\stackrel{\sim}{\rightarrow}(R'_\infty)^\fN. 
\end{equation}
\item[{\rm (ii)}] Pour tout $a\in \co_\oK$, l'homomorphisme canonique
\begin{equation}\label{eccr9e}
(R^\intern_\infty/aR^\intern_\infty)\otimes_{R^\intern_{p^\infty}}R'_{p^\infty}\rightarrow(R'_\infty/aR'_\infty)^\fN 
\end{equation}
est un isomorphisme. 
\end{itemize}
\end{prop}

\begin{prop}\label{eccr99}
La $R^\intern_\infty$-algèbre $\oR^\intern$ est $\alpha$-fidèlement plate.
\end{prop}
Cela résulte de (\cite{ag} 2.9.12 et 5.2.15) et (\cite{agt} V.12.4 et V.12.9).

\subsection{}\label{eccr12}
Posons 
\begin{eqnarray}
\Xi_{p^n}&=&\Hom(\Delta_{p^\infty},\mu_{p^n}(\co_\oK)),\label{eccr12a}\\
\Xi_{p^\infty}&=&\Hom(\Delta_{p^\infty},\mu_{p^\infty}(\co_\oK)),\label{eccr12b}\\
\Xi'_{p^n}&=&\Hom(\Delta'_{p^\infty},\mu_{p^n}(\co_\oK)),\label{eccr12c}\\
\Xi'_{p^\infty}&=&\Hom(\Delta'_{p^\infty},\mu_{p^\infty}(\co_\oK)).\label{eccr12d}
\end{eqnarray}
On identifie $\Xi_{p^n}$ (resp. $\Xi'_{p^n}$) à un sous-groupe de $\Xi_{p^\infty}$ (resp. $\Xi'_{p^\infty}$). D'après (\cite{agt} II.8.9), 
il existe une décomposition canonique de $R_{p^\infty}$ en somme directe de $R_1$-modules de présentation finie,
stables sous l'action de $\Delta_{p^\infty}$,
\begin{equation}\label{eccr12e}
R_{p^\infty}=\bigoplus_{\lambda\in \Xi_{p^\infty}}R^{(\lambda)}_{p^\infty},
\end{equation}
telle que l'action de $\Delta_{p^\infty}$ sur le facteur $R^{(\lambda)}_{p^\infty}$ soit donnée par le caractère $\lambda$. 
De plus, pour tout $n\geq 0$, on a 
\begin{equation}\label{eccr12f}
R_{p^n}=\bigoplus_{\lambda\in \Xi_{p^n}}R^{(\lambda)}_{p^\infty}. 
\end{equation}
De même, il existe une décomposition canonique de $R'_{p^\infty}$ en somme directe de $R'_1$-modules de présentation finie,
stables sous l'action de $\Delta'_{p^\infty}$,
\begin{equation}\label{eccr12g}
R'_{p^\infty}=\bigoplus_{\lambda\in \Xi'_{p^\infty}}R'^{(\lambda)}_{p^\infty},
\end{equation}
telle que l'action de $\Delta'_{p^\infty}$ sur le facteur $R'^{(\lambda)}_{p^\infty}$ soit donnée par le caractère $\lambda$. 
De plus, pour tout $n\geq 0$, on a 
\begin{equation}\label{eccr12h}
R'_{p^n}=\bigoplus_{\lambda\in \Xi'_{p^n}}R'^{(\lambda)}_{p^\infty}. 
\end{equation}

Posons 
\begin{eqnarray}
\upchi_{p^n}&=&\Hom(\fS_{p^\infty},\mu_{p^n}(\co_\oK)),\label{eccr12i}\\
\upchi_{p^\infty}&=&\Hom(\fS_{p^\infty},\mu_{p^\infty}(\co_\oK)).\label{eccr12j}
\end{eqnarray}
On identifie $\upchi_{p^n}$ à un sous-groupe de $\upchi_{p^\infty}$. 
En vertu de  \eqref{eccr43k} et (\cite{ag} 5.2.14(iv)), la suite canonique de $\mZ_p$-modules  
\begin{equation}\label{eccr12k}
0\rightarrow \fS_{p^\infty}\rightarrow \Delta'_{p^\infty}\rightarrow \Delta_{p^\infty} \rightarrow 0
\end{equation}
est exacte et scindée. On en déduit une suite exacte 
\begin{equation}\label{eccr12l}
0\rightarrow \Xi_{p^\infty}\rightarrow \Xi'_{p^\infty}\rightarrow \upchi_{p^\infty} \rightarrow 0.
\end{equation}

\begin{lem}[\cite{ag} 5.2.24]\label{eccr13}
L'homomorphisme canonique $R^\intern_{p^\infty}\rightarrow R'_{p^\infty}$ induit un isomorphisme
\begin{equation}\label{eccr13a}
R^\intern_{p^\infty}\stackrel{\sim}{\rightarrow}\bigoplus_{\lambda\in \Xi_{p^\infty}}R'^{(\lambda)}_{p^\infty}. 
\end{equation}
\end{lem}

\subsection{}\label{eccr16}
Pour tout $\nu\in \upchi_{p^\infty}$, fixons un relèvement $\tnu\in \Xi'_{p^\infty}$ \eqref{eccr12l} et posons 
\begin{equation}\label{eccr16a}
R'^{(\nu)}_{p^\infty}=\bigoplus_{\lambda\in \Xi_{p^\infty}}R'^{(\tnu+\lambda)}_{p^\infty},
\end{equation}
qui est naturellement un $R^\intern_{p^\infty}$-module \eqref{eccr13}. 
La décomposition \eqref{eccr12g} induit alors une décomposition en $R^\intern_{p^\infty}$-modules
\begin{equation}\label{eccr16b}
R'_{p^\infty}=\bigoplus_{\nu\in \upchi_{p^\infty}}R'^{(\nu)}_{p^\infty},
\end{equation}
telle que l'action de $\fS_{p^\infty}$ sur le facteur $R'^{(\nu)}_{p^\infty}$ soit donnée par le caractère $\nu$.

\begin{lem}[\cite{ag} 5.2.26]\label{eccr17}
Pour tout $n\geq 0$, l'homomorphisme canonique $R^\intern_{p^\infty}\otimes_{R^\intern_{p^n}}R'_{p^n}\rightarrow R'_{p^\infty}$ induit un isomorphisme
\begin{equation}\label{eccr17a}
R^\intern_{p^\infty}\otimes_{R^\intern_{p^n}}R'_{p^n}\stackrel{\sim}{\rightarrow}\bigoplus_{\nu\in \upchi_{p^n}}R'^{(\nu)}_{p^\infty}. 
\end{equation}
En particulier, pour tout $\nu\in \upchi_{p^n}$, le $R^\intern_{p^\infty}$-module $R'^{(\nu)}_{p^\infty}$ est de présentation finie.
\end{lem}

\subsection{}\label{eccr32}
Soit $a$ un élément non nul de $\co_\oK$. 
Comme le sous-groupe de torsion de $P'^\gp/h^\gp(P^\gp)$ est d'ordre premier à $p$,  
on a d'après (\cite{ag} (5.2.12.3)), un isomorphisme canonique 
\begin{equation}\label{eccr32a}
(P'^\gp/h^\gp(P^\gp))\otimes_\mZ\mZ_p(-1)\stackrel{\sim}{\rightarrow}
\Hom_{\mZ_p}(\fS_{p^\infty}, \mZ_p).
\end{equation}
On en déduit, compte tenu de (\cite{ogus} IV 1.1.4), un isomorphisme $R^\intern_\infty$-linéaire \eqref{hypmdef4c}
\begin{equation}\label{eccr32b}
\tOmega^1_{R'/R}\otimes_{R'}(R^\intern_\infty/aR^\intern_\infty)(-1)\stackrel{\sim}{\rightarrow} 
\Hom_{\mZ}(\fS_{p^\infty},R^\intern_\infty/aR^\intern_\infty).
\end{equation}
On interprètera dans la suite le but de ce morphisme comme un groupe de cohomologie. Posons 
\begin{equation}\label{eccr32c}
\delta=\dim_{\mQ}((P'^\gp/h^\gp(P^\gp))\otimes_{\mZ}\mQ).
\end{equation}

\begin{prop}[\cite{ag} 5.2.32]\label{eccr51}
Soit $a$ un élément non nul de $\co_\oK$. 
\begin{itemize}
\item[{\rm (i)}] Il existe un et un unique homomorphisme de $\oR^\intern$-algèbres graduées
\begin{equation}\label{eccr51a}
\wedge(\tOmega^1_{R'/R}\otimes_{R'}(\oR^\intern/a\oR^\intern)(-1))\rightarrow \rH^*(\Pi,\oR'/a\oR')
\end{equation}
dont la composante en degré un est induite par \eqref{eccr32b} (cf. \ref{eccr34}). Celui-ci 
est $\alpha$-injectif et son conoyau est annulé par $p^{\frac{1}{p-1}}\fm_\oK$.  
\item[{\rm (ii)}] Le $\oR^\intern$-module $\rH^i(\Pi,\oR'/a \oR')$ est de présentation $\alpha$-finie 
pour tout $i\geq 0$, et est $\alpha$-nul pour tout $i\geq \delta+1$ \eqref{eccr32c}. 
\item[{\rm (iii)}] Pour tous entiers $r'\geq r\geq 0$, notons
\begin{equation}\label{eccr51b}
\hbar_{r,r'}\colon \rH^*(\Pi,\oR'/p^{r'} \oR')\rightarrow \rH^*(\Pi,\oR'/p^r \oR')
\end{equation}
le morphisme canonique. Alors pour tout entier $r\geq 1$,
il existe un entier $r'\geq r$, dépendant seulement de $\delta$ mais pas des autres données dans \ref{hypmdef2}, 
tel que pour tout entier $r''\geq r'$, les images de $\hbar_{r,r'}$ et $\hbar_{r,r''}$ soient $\alpha$-isomorphes.  
\end{itemize}
\end{prop}

\section{Cohomologies galoisienne relative des algèbres de Higgs-Tate}\label{mtht}

\subsection{}\label{hypmdef3}
Les hypothèses et notations de \ref{hypmdef} sont en vigueur dans cette section, en particulier celles de \ref{hypmdef2}.
Reprenons de plus les notations introduites dans \ref{definf10}. 
On munit  $\coX=X\times_S\coS$ et $\coX'=X'\times_S\coS$ \eqref{hypmdef1b} des structures logarithmiques $\cM_\coX$ et $\cM_{\coX'}$ 
images inverses respectivement de $\cM_X$ et $\cM_{X'}$. 
On suppose qu'il existe des $(\tS,\cM_{\tS})$-déformations lisses $(\tX,\cM_\tX)$ de $(\coX,\cM_{\coX})$ 
et $(\tX',\cM_{\tX'})$ de $(\coX',\cM_{\coX'})$ \eqref{defing12b} et un $(\tS,\cM_{\tS})$-morphisme
\begin{equation}\label{hypmdef3a}
\tg\colon (\tX',\cM_{\tX'})\rightarrow (\tX,\cM_\tX)
\end{equation}
qui s'insère dans un diagramme commutatif (à carrés cartésiens)
\begin{equation}\label{hypmdef3b}
\xymatrix{
{(\coX',\cM_{\coX'})}\ar[r]\ar[d]\ar@{}[rd]|{\Box}&{(\tX',\cM_{\tX'})}\ar[d]^{\tg}\\
{(\coX,\cM_\coX)}\ar[r]\ar[d]\ar@{}[rd]|{\Box}&{(\tX,\cM_\tX)}\ar[d]\\
{(\coS,\cM_{\coS})}\ar[r]&{(\tS,\cM_{\tS})}}
\end{equation}
On notera que les carrés sont cartésiens aussi bien dans la catégorie des schémas logarithmiques que 
dans celle des schémas logarithmiques fins. 
{\em On fixe dans la suite de ce chapitre les déformations et le $(\tS,\cM_{\tS})$-morphisme $\tg$ \eqref{hypmdef3a}.}  

\begin{lem}\label{hypmdef30}
Le morphisme $\tg$ \eqref{hypmdef3a} est lisse. 
\end{lem}

En effet, en vertu de (\cite{kato1} 3.12 ou \cite{ogus} IV 3.2.3), il suffit de montrer que le morphisme canonique 
\begin{equation}
\tg^*(\Omega^1_{(\tX,\cM_\tX)/(\tS,\cM_{\tS})})\rightarrow \Omega^1_{(\tX',\cM_{\tX'})/(\tS,\cM_{\tS})}
\end{equation}
est localement inversible à gauche, autrement dit qu'il induit localement un isomorphisme de la source sur un facteur direct du but. 
Il revient au même de dire que pour tout point $x'$ de $\tX'$ (ou ce qui revient au même de $\coX'$), d'image $x$ de $\tX$, le morphisme canonique 
\begin{equation}
\tg^*(\Omega^1_{(\tX,\cM_\tX)/(\tS,\cM_{\tS})})_{x'}\rightarrow \Omega^1_{(\tX',\cM_{\tX'})/(\tS,\cM_{\tS}),x'}
\end{equation}
est inversible à gauche. 
Comme la source et le but sont des $\co_{\tX',x'}$-modules libres de type fini, 
cette dernière condition est équivalente au fait que le morphisme canonique 
\begin{equation}
\Omega^1_{(\tX,\cM_\tX)/(\tS,\cM_{\tS}),x}\otimes_{\co_{\tX,x}}\kappa(x')\rightarrow \Omega^1_{(\tX',\cM_{\tX'})/(\tS,\cM_{\tS}),x'}\otimes_{\co_{\tX',x'}}\kappa(x')
\end{equation}
soit injectif en vertu de (\cite{ega4} 0.19.1.12). Celui-ci s'identifie au morphisme 
\begin{equation}
\Omega^1_{(\coX,\cM_\coX)/(\coS,\cM_{\coS}),x}\otimes_{\co_{\coX,x}}\kappa(x')\rightarrow \Omega^1_{(\coX',\cM_{\coX'})/(\coS,\cM_{\coS}),x'}\otimes_{\co_{\coX',x'}}\kappa(x')
\end{equation}
induit par le morphisme lisse $g\times_S\coS$ \eqref{hypmdef2a}. Il est donc injectif; d'où la proposition.  

\subsection{}\label{mtht3}
On pose 
\begin{equation}\label{mtht3a}
\mX=\Spec(\oR) \ \ \ {\rm et} \ \ \  \hmX=\Spec(\hoR),
\end{equation}
que l'on munit des structures logarithmiques images inverses de $\cM_X$, 
notées respectivement $\cM_\mX$ et $\cM_\hmX$. 
Reprenant les notations de \ref{taht4}, on désigne par $(\tmX,\cM_\tmX)$ l'un des deux schémas logarithmiques 
\begin{equation}\label{mtht3b}
(\cA_2(\mX),\cM_{\cA_2(\mX)})\ \ \ {\rm ou} \ \ \ (\cA^{\ast}_2(\mX/S),\cM_{\cA^{\ast}_2(\mX/S)}),
\end{equation}
selon que l'on est dans le cas absolu ou relatif \eqref{definf10}, et par 
\begin{equation}\label{mtht3c}
i_X\colon (\hmX,\cM_\hmX)\rightarrow (\tmX,\cM_{\tmX})
\end{equation}
l'immersion fermée exacte canonique. 

De même, on pose  
\begin{equation}\label{mtht3d}
\mX'=\Spec(\oR')\ \ \ {\rm et} \ \ \ \hmX'=\Spec(\hoRp),
\end{equation}
que l'on munit des structures logarithmiques images inverses de $\cM_{X'}$, 
notées respectivement $\cM_{\mX'}$ et $\cM_{\hmX'}$. 
On désigne par $(\tmX',\cM_{\tmX'})$ l'un des deux schémas logarithmiques 
\begin{equation}\label{mtht3e}
(\cA_2(\mX'),\cM_{\cA_2(\mX')})\ \ \ {\rm ou} \ \ \ (\cA^{\ast}_2(\mX'/S),\cM_{\cA^{\ast}_2(\mX'/S)}),
\end{equation}
selon que l'on est dans le cas absolu ou relatif \eqref{definf10}, et par 
\begin{equation}\label{mtht3f}
i_{X'}\colon (\hmX',\cM_{\hmX'})\rightarrow (\tmX',\cM_{\tmX'})
\end{equation}
l'immersion fermée exacte canonique (cf. \ref{taht4}). 

On a un morphisme canonique d'algèbres $\oR\rightarrow \oR'$ \eqref{eccr40} et un morphisme canonique de monoïdes 
$Q_X\rightarrow Q_{X'}$, où $Q_X$ et $Q_{X'}$ sont les monoïdes définis dans \ref{cad6} relativement aux schémas logarithmiques 
$(X,\cM_X)$ et $(X',\cM_{X'})$. Le diagramme d'homomorphisme de monoïdes 
\begin{equation}\label{mtht3i}
\xymatrix{
P\ar[r]^-(0.5){\upnu}\ar[d]_h&Q_{X}\ar[r]^-(0.5){\uptau_X}\ar[d]&{\rW(\oR^\flat)}\ar[d]\\
P'\ar[r]^-(0.5){\upnu'}&Q_{X'}\ar[r]^-(0.5){\uptau_{X'}}&{\rW(\oR'^\flat)}} 
\end{equation}
où $\uptau_X$ et $\uptau_{X'}$ sont les homomorphismes \eqref{cad6b} et $\upnu$ et $\upnu'$ 
sont les homomorphismes \eqref{pmh6e}, est commutatif. 

On déduit de ce qui précède un morphisme de $(\coS,\cM_{\coS})$-schémas logarithmiques
\begin{equation}\label{mtht3h}
\hmg\colon (\hmX',\cM_{\hmX'})\rightarrow (\hmX,\cM_{\hmX}),
\end{equation}
et un morphisme de $(\tS,\cM_{\tS})$-schémas logarithmiques
\begin{equation}\label{mtht3hh}
\tmg\colon (\tmX',\cM_{\tmX'})\rightarrow (\tmX,\cM_{\tmX}),
\end{equation}
compatible avec $\hmg$. 

Soient $U$ un ouvert de Zariski de $\hmX$ et $U'$ un ouvert de Zariski de $\hmX'$ tels que $\hmg(U')\subset U$. 
On note $\tU$ l'ouvert de $\tmX$ correspondant à $U$ et $\tU'$ l'ouvert de $\tmX'$ correspondant à $U'$. 
\begin{equation}\label{mtht3g}
\xymatrix{
&{(U',\cM_{\hmX'}|U')}\ar[ld]\ar[rr]\ar@{->}'[d]^-(0.5){\hmg|U'}[dd]&&{(\tU',\cM_{\tmX'}|\tU')}\ar@{.>}[ld]\ar[dd]^-(0.3){\tmg|\tU'}\\
{(\coX',\cM_{\coX'})}\ar[rr]\ar[dd]&&{(\tX',\cM_{\tX'})}\ar[dd]&\\
&{(U,\cM_\hmX|U)}\ar[ld]\ar@{->}'[r][rr]\ar[ddl]|\hole&&{(\tU,\cM_{\tmX}|\tU)}\ar@{.>}[ld]\ar[ldd]\\
{(\coX,\cM_\coX)}\ar[rr]\ar[d]&&{(\tX,\cM_\tX)}\ar[d]\\
{(\coS,\cM_{\coS})}\ar[rr]&&{(\tS,\cM_{\tS})}}
\end{equation}

\subsection{}\label{mtht4}
On pose \eqref{cad1i}
\begin{equation}\label{mtht4a}
\rT=\Hom_{\hoR}(\tOmega^1_{R/\co_K}\otimes_R\hoR,\txi\hoR).
\end{equation} 
On identifie le $\hoR$-module dual à $\txi^{-1}\tOmega^1_{R/\co_K}\otimes_R\hoR$ (cf. \ref{definf10})
et on note $\cG$ la $\hoR$-algèbre symétrique associée \eqref{notconv9}
\begin{equation}\label{mtht4b}
\cG=\rS_{\hoR}(\txi^{-1}\tOmega^1_{R/\co_K}\otimes_R\hoR).
\end{equation}
On désigne par $\hmX_\zar$ le topos de Zariski de $\hmX$, par $\trT$ le $\co_\hmX$-module associé à $\rT$
et par $\bT$ le $\hmX$-fibré vectoriel associé à son dual, autrement dit,  
\begin{equation}\label{mtht4c}
\bT=\Spec(\cG).
\end{equation}

On pose, de même,
\begin{equation}\label{mtht4d}
\rT'=\Hom_{\hoRp}(\tOmega^1_{R'/\co_K}\otimes_{R'}\hoRp,\txi\hoRp).
\end{equation} 
On identifie le $\hoRp$-module dual  à $\txi^{-1}\tOmega^1_{R'/\co_K}\otimes_{R'}\hoRp$, 
et on note $\cG'$ la $\hoRp$-algèbre symétrique associée 
\begin{equation}\label{mtht4e}
\cG'=\rS_{\hoRp}(\txi^{-1}\tOmega^1_{R'/\co_K}\otimes_{R'}\hoRp).
\end{equation}
On désigne par $\hmX'_\zar$ le topos de Zariski de $\hmX'$, par $\trT'$ le $\co_{\hmX'}$-module associé à $\rT'$
et par $\bT'$ le $\hmX'$-fibré vectoriel associé à son dual, autrement dit,  
\begin{equation}\label{mtht4f}
\bT'=\Spec(\cG').
\end{equation}

On pose, enfin, 
\begin{equation}\label{mtht4i}
\rT_{R'/R}=\Hom_{\hoRp}(\tOmega^1_{R'/R}\otimes_{R'}\hoRp,\txi\hoRp).
\end{equation} 
On identifie le $\hoRp$-module dual  à $\txi^{-1}\tOmega^1_{R'/R}\otimes_{R'}\hoRp$, 
et on note $\cG_{R'/R}$ la $\hoRp$-algèbre symétrique associée \eqref{notconv9}
\begin{equation}\label{mtht4j}
\cG_{R'/R}=\rS_{\hoRp}(\txi^{-1}\tOmega^1_{R'/R}\otimes_{R'}\hoRp).
\end{equation}
On désigne par $\trT_{X'/X}$ le $\co_{\hmX'}$-module associé à $\rT_{R'/R}$
et par $\bT_{X'/X}$ le $\hmX'$-fibré vectoriel associé à son dual, autrement dit,  
\begin{equation}\label{mtht4k}
\bT_{X'/X}=\Spec(\cG_{R'/R}).
\end{equation}

La suite exacte localement scindée \eqref{hypmdef4d} induit une suite exacte de $\co_{\hmX'}$-modules 
\begin{equation}\label{mtht4g}
\xymatrix{
0\ar[r]&{\trT_{X'/X}}\ar[r]^-(0.4){\iota}&{\trT'}\ar[r]^-(0.4)u&{\hmg^*(\trT)}\ar[r]&0}.
\end{equation}
On en déduit une suite exacte de $\hmX'$-fibrés vectoriels
\begin{equation}\label{mtht4h}
\xymatrix{
0\ar[r]&{\bT_{X'/X}}\ar[r]&{\bT'}\ar[r]&{\bT\times_{\hmX}\hmX'}\ar[r]&0}.
\end{equation}

\subsection{}\label{mtht5}
Soient $U$ un ouvert de Zariski de $\hmX$, $\tU$ l'ouvert correspondant de $\tmX$ \eqref{mtht3}. 
On désigne par $\cL(U)$ l'ensemble des $\tS$-morphismes représentés par des flèches pointillées qui complètent  
la face inférieure du parallélépipède du diagramme \eqref{mtht3g} de façon à la laisser commutative. 
D'après \ref{taht5}, le foncteur $U\mapsto \cL(U)$ est un $\trT$-torseur de $\hmX_\zar$; 
c'est le torseur de Higgs-Tate associé à $(\tX,\cM_\tX)$.  
Le $\trT$-torseur $\cL$ est naturellement muni d'une structure $\Delta$-équivariante (cf. \ref{taht6} et \ref{eccr40}).
On désigne par $\cF$ le $\hoR$-module des fonctions affines sur $\cL$ (cf. \cite{agt} II.4.9). 
Celui-ci s'insère dans une suite exacte canonique 
\begin{equation}\label{mtht5a}
0\rightarrow \hoR\rightarrow \cF\rightarrow \txi^{-1}\tOmega^1_{R/\co_K} \otimes_R \hoR\rightarrow 0.
\end{equation} 
Le $\hoR$-module $\cF$ est naturellement muni d'une action $\hoR$-semi-linéaire de $\Delta$ 
telle que les morphismes de la suite \eqref{mtht5a} soient 
$\Delta$-équivariants (cf. \ref{taht6}). 
D'après (\cite{illusie1} I 4.3.1.7), cette suite induit pour tout entier $n\geq 1$, une suite exacte \eqref{notconv9}
\begin{equation}\label{mtht5b}
0\rightarrow \rS^{n-1}_{\hoR}(\cF)\rightarrow \rS^{n}_{\hoR}(\cF)\rightarrow \rS^n_{\hoR}(\txi^{-1}\tOmega^1_{R/\co_K}
\otimes_R\hoR)\rightarrow 0.
\end{equation}
Les $\hoR$-modules $(\rS^{n}_{\hoR}(\cF))_{n\in \mN}$ forment donc un système inductif filtrant, 
dont la limite inductive 
\begin{equation}\label{mtht5c}
\cC=\underset{\underset{n\geq 0}{\longrightarrow}}\lim\ \rS^n_{\hoR}(\cF)
\end{equation}
est naturellement munie d'une structure de $\hoR$-algèbre et d'une action de $\Delta$
par des automorphismes d'anneaux, compatible avec son action sur $\hoR$; 
c'est l'algèbre de Higgs-Tate associée à $(\tX,\cM_\tX)$ \eqref{taht7}. D'après (\cite{agt} II.4.10), le $\hmX$-schéma 
\begin{equation}\label{mtht5d}
\bL=\Spec(\cC)
\end{equation}
est naturellement un $\bT$-fibré principal homogène sur $\hmX$ qui représente canoniquement $\cL$. 
Celui-ci est muni d'une structure $\Delta$-équivariante (cf. \ref{taht6}). 
Cette structure détermine une action à gauche de $\Delta$ sur $\bL$ compatible avec son action sur $\hmX$.

\subsection{}\label{mtht6}
Soient $U'$ un ouvert de Zariski de $\hmX'$, $\tU'$ l'ouvert correspondant de $\tmX'$. 
On désigne par $\cL'(U')$ l'ensemble des $\tS$-morphismes représentés par des flèches pointillées qui complètent  
la face supérieure du parallélépipède du diagramme \eqref{mtht3g} de façon à la laisser commutative. 
Le foncteur $U'\mapsto \cL'(U')$ est un $\trT'$-torseur de $\hmX'_\zar$; 
c'est le torseur de Higgs-Tate associé à $(\tX',\cM_{\tX'})$.  
Le $\trT'$-torseur $\cL'$ est naturellement muni d'une structure $\Delta'$-équivariante (cf. \ref{taht6} et \ref{eccr40}).
On désigne par $\cF'$ le $\hoRp$-module des fonctions affines sur $\cL'$. 
Celui-ci s'insère dans une suite exacte canonique 
\begin{equation}\label{mtht6a}
0\rightarrow \hoRp\rightarrow \cF'\rightarrow \txi^{-1}\tOmega^1_{R'/\co_K} \otimes_{R'} \hoRp\rightarrow 0.
\end{equation} 
Le $\hoRp$-module $\cF'$ est naturellement muni d'une action $\hoRp$-semi-linéaire de $\Delta'$ 
telle que les morphismes de la suite \eqref{mtht6a} soient 
$\Delta'$-équivariants (cf. \ref{taht6}). 
Cette suite induit pour tout entier $n\geq 1$, une suite exacte \eqref{notconv9}
\begin{equation}\label{mtht6b}
0\rightarrow \rS^{n-1}_{\hoRp}(\cF')\rightarrow \rS^{n}_{\hoRp}(\cF')\rightarrow \rS^n_{\hoRp}(\txi^{-1}\tOmega^1_{R'/\co_K}
\otimes_{R'}\hoRp)\rightarrow 0.
\end{equation}
Les $\hoR$-modules $(\rS^{n}_{\hoRp}(\cF'))_{n\in \mN}$ forment donc un système inductif filtrant, 
dont la limite inductive 
\begin{equation}\label{mtht6c}
\cC'=\underset{\underset{n\geq 0}{\longrightarrow}}\lim\ \rS^n_{\hoRp}(\cF')
\end{equation}
est naturellement munie d'une structure de $\hoRp$-algèbre et d'une action de $\Delta'$
par des automorphismes d'anneaux, compatible avec son action sur $\hoRp$; 
c'est l'algèbre de Higgs-Tate associée à $(\tX',\cM_{\tX'})$ \eqref{taht7}. 
Le $\hmX'$-schéma 
\begin{equation}\label{mtht6d}
\bL'=\Spec(\cC')
\end{equation}
est naturellement un $\bT'$-fibré principal homogène sur $\hmX'$ qui représente canoniquement $\cL'$. On a donc 
un $\hmX'$-morphisme canonique
\begin{equation}\label{mtht6e}
\bT'\times_{\hmX'}\bL'\rightarrow \bL'.
\end{equation}
Le $\bT'$-fibré principal homogène $\bL'$ est muni d'une structure $\Delta'$-équivariante (cf. \ref{taht6}). 
Cette structure détermine une action à gauche de $\Delta'$ sur $\bL'$ compatible avec son action sur $\hmX'$.

\subsection{}\label{mtht7}
On désigne par $\cL^+$ l'image inverse affine $\hmg^+(\cL)$ du $\trT$-torseur $\cL$ par le morphisme $\hmg\colon \hmX'\rightarrow \hmX$ 
(\cite{agt} II.4.5), autrement dit le $\hmg^*(\trT)$-torseur de $\hmX'_\zar$ déduit du $\hmg^{-1}(\trT)$-torseur $\hmg^*(\cL)$
par extension de son groupe structural par l'homomorphisme canonique $\hmg^{-1}(\trT)\rightarrow \hmg^*(\trT)$~:
\begin{equation}\label{mtht7a}
\hmg^+(\cL)=\hmg^*(\cL)\wedge^{\hmg^{-1}(\trT)}\hmg^*(\trT).
\end{equation}

Comme $\hmX$ est affine, le $\trT$-torseur $\cL$ est trivial; soit $\sigma \in \cL(\hmX)$. 
D'autre part, le morphisme $(\tX',\cM_{\tX'})\rightarrow (\tX,\cM_\tX)$ étant lisse et $\hmX'$ étant affine, 
il existe $\sigma' \in \cL'(\hmX')$ qui complète la face 
latérale droite du parallélépipède du diagramme \eqref{mtht3g} (pour $U=\hmX$ et $U'=\hmX'$) de façon à la laisser commutative.
On dira dans la suite que $\sigma'$ relève $\sigma$, ou que $\sigma$ et $\sigma'$ sont {\em compatibles}. 
Il existe alors un unique morphisme $u$-équivariant \eqref{mtht4g}
\begin{equation}\label{mtht7b}
v\colon \cL'\rightarrow \cL^+
\end{equation}
qui envoie $\sigma'$ sur l'image canonique de $\sigma$ dans $\cL^+(\hmX')$. 
Il résulte aussitôt de \ref{mtht8} ci-dessous que le morphisme $v$ ne dépend pas du choix du couple de sections compatibles $(\sigma,\sigma')$.

D'après (\cite{agt} II.4.12 et II.4.13), le morphisme $v$ induit un morphisme $\hoRp$-linéaire 
\begin{equation}\label{mtht7d}
\cF\otimes_{\hoR}\hoRp\rightarrow \cF'
\end{equation}
qui s'insère dans un diagramme commutatif 
\begin{equation}\label{mtht7e}
\xymatrix{
0\ar[r]&{\hoRp}\ar@{=}[d]\ar[r]&{\cF\otimes_{\hoR}\hoRp}\ar[d]\ar[r]&{\txi^{-1}\tOmega^1_{R/\co_K}\otimes_R\hoRp}\ar[r]\ar[d]&0\\
0\ar[r]&{\hoRp}\ar[r]&{\cF'}\ar[r]&{\txi^{-1}\tOmega^1_{R'/\co_K}\otimes_{R'}\hoRp}\ar[r]&0}
\end{equation}
où les lignes sont induites par les suites exactes \eqref{mtht5a} et \eqref{mtht6a} 
et  la troisième flèche verticale est induite par le morphisme canonique \eqref{hypmdef4d}. 

D'après (\cite{agt} II.4.6), le morphisme $v$ \eqref{mtht7b} induit un morphisme $\bT'$-équivariant \eqref{mtht4h}
\begin{equation}\label{mtht7c}
\nu\colon \bL'\rightarrow \bL\times_{\hmX}\hmX', 
\end{equation}
et donc un homomorphisme de $\hoR$-algèbres 
\begin{equation}\label{mtht7f}
\cC\rightarrow \cC'.
\end{equation}
Ce dernier prolonge le morphisme $\hoR$-linéaire $\cF\rightarrow \cF'$ induit par \eqref{mtht7d} 
et peut en être déduit compte tenu de \eqref{mtht7e} et des définitions \eqref{mtht5c} et \eqref{mtht6c}.

Le morphisme $\hmg$ étant $\Delta'$-équivariant \eqref{eccr40}, le $(\bT\times_{\hmX}\hmX')$-fibré principal homogène 
$\bL\times_{\hmX}\hmX'$ sur $\hmX'$ est naturellement muni d'une structure $\Delta'$-équivariante.

\begin{lem}\label{mtht8}
Soient $(\sigma,\sigma'), (\sigma_1,\sigma'_1) \in \cL(\hmX)\times \cL'(\hmX')$ deux couples de sections compatibles \eqref{mtht7}. 
Alors, on a 
\begin{equation}\label{mtht8a}
u(\sigma'_1-\sigma')=\sigma_1-\sigma\in \rT\otimes_{\hoR}\hoRp,
\end{equation}
où $u$ est le morphisme défini dans \eqref{mtht4g}. En particulier, 
le morphisme $v$ \eqref{mtht7b} ne dépend pas du choix du couple de sections compatibles qui sert à le définir. 
\end{lem}

Cela résulte aussitôt par fonctorialité du diagramme \eqref{mtht3g}.

\begin{lem}\label{mtht9}
Le morphisme $\bT_{X'/X}\times_{\hmX'}\bL'\rightarrow \bL'$ induit par les morphismes \eqref{mtht6e} et 
$\bT_{X'/X}\rightarrow \bT'$ \eqref{mtht4h} fait de $\bL'$ 
un $\bL\times_{\hmX}\bT_{X'/X}$-fibré principal homogène sur $\bL\times_{\hmX}\hmX'$ \eqref{mtht7c}.
\end{lem}

Cela résulte aussitôt de la suite exacte de $\hmX'$-fibrés vectoriels \eqref{mtht4h} et du fait que le morphisme 
$\nu\colon \bL'\rightarrow \bL\times_{\hmX}\hmX'$ \eqref{mtht7c} est $\bT'$-équivariant.

\begin{lem}\label{mtht90}
Le morphisme $\nu\colon \bL'\rightarrow \bL\times_{\hmX}\hmX'$ \eqref{mtht7c}  est $\Delta'$-équivariant. 
En particulier, l'homomorphisme  $\cC\rightarrow \cC'$ \eqref{mtht7f} est $\Delta'$-équivariant.
\end{lem}

Soient $(\sigma,\sigma') \in \cL(\hmX)\times \cL'(\hmX')$ un couple de sections compatibles \eqref{mtht7}, $c\in \Delta'$. 
On désigne par ${^c\sigma}$ la section de $\cL(\hmX)$ définie par le morphisme composé \eqref{taht6k}
\begin{equation}\label{mtht90a}
(\tmX,\cM_{\tmX})\stackrel{c^{-1}}{\longrightarrow} (\tmX,\cM_{\tmX})
\stackrel{\sigma}{\longrightarrow} (\tX,\cM_\tX)
\end{equation}
et par ${^c\sigma'}$ la section de $\cL'(\hmX')$ définie par le morphisme composé 
\begin{equation}\label{mtht90b}
(\tmX',\cM_{\tmX'})\stackrel{c^{-1}}{\longrightarrow} (\tmX',\cM_{\tmX'})
\stackrel{\sigma'}{\longrightarrow} (\tX',\cM_{\tX'}).
\end{equation}
Comme le morphisme $\tmg\colon (\tmX',\cM_{\tmX'})\rightarrow (\tmX,\cM_{\tmX})$ est $\Delta'$-équivariant, 
les sections ${^c\sigma}$ et ${^c\sigma'}$ sont compatibles \eqref{mtht3g}.
Par suite, $\nu({^c\sigma'})={^c\sigma}$ et donc $\nu(c(\sigma'))=c(\sigma)$ en vertu de \eqref{taht6i}. 
Comme l'homomorphisme canonique $\mu\colon \bT'\rightarrow \bT\times_{\hmX}\hmX'$ \eqref{mtht4h} est $\Delta'$-équivariant 
et que $\nu$ est $\mu$-équivariant, on en déduit que $\nu$ est $\Delta'$-équivariant.

\subsection{}\label{mtht100}
Considérons des $(\tS,\cM_{\tS})$-déformations lisses $(\tX^\natural,\cM_{\tX^\natural})$ de $(\coX,\cM_{\coX})$
et $(\tX'^\natural,\cM_{\tX'^\natural})$ de $(\coX',\cM_{\coX'})$ \eqref{defing12} et un $(\tS,\cM_{\tS})$-morphisme lisse
\begin{equation}\label{mtht100a}
\tg^\natural\colon (\tX'^\natural,\cM_{\tX'^\natural})\rightarrow (\tX^\natural,\cM_{\tX^\natural})
\end{equation}
qui s'insère dans un diagramme commutatif 
\begin{equation}\label{mtht100b}
\xymatrix{
{(\coX',\cM_{\coX'})}\ar[r]\ar[d]&{(\tX'^\natural,\cM_{\tX'^\natural})}\ar[d]^{\tg^\natural}\\
{(\coX,\cM_\coX)}\ar[r]\ar[d]&{(\tX^\natural,\cM_{\tX^\natural})}\ar[d]\\
{(\coS,\cM_{\coS})}\ar[r]&{(\tS,\cM_{\tS})}}
\end{equation}
D'après (\cite{kato1} 3.14), il existe un isomorphisme de $(\tS,\cM_{\tS})$-déformations
\begin{equation}\label{mtht100c}
h\colon (\tX,\cM_{\tX})\stackrel{\sim}{\rightarrow} (\tX^\natural,\cM_{\tX^\natural}).
\end{equation}
Considérons $(\tX',\cM_{\tX'})$ comme une $(\tX,\cM_{\tX})$-déformation lisse de $(\coX',\cM_{\coX'})$
et $(\tX'^\natural,\cM_{\tX'^\natural})$ comme une $(\tX^\natural,\cM_{\tX^\natural})$-déformation lisse de $(\coX',\cM_{\coX'})$. 
D'après (\cite{kato1} 3.14), il existe un isomorphisme de $(\tS,\cM_{\tS})$-déformations
\begin{equation}\label{mtht100d}
h'\colon (\tX',\cM_{\tX'})\stackrel{\sim}{\rightarrow} (\tX'^\natural,\cM_{\tX'^\natural})
\end{equation}
qui s'insère dans un diagramme commutatif 
\begin{equation}\label{mtht100e}
\xymatrix{
{(\tX',\cM_{\tX'})}\ar[r]^{h'}\ar[d]_{\tg}&{(\tX'^\natural,\cM_{\tX'^\natural})}\ar[d]^{\tg^\natural}\\
{(\tX,\cM_{\tX})}\ar[r]^h&{(\tX^\natural,\cM_{\tX^\natural})}}
\end{equation}

Nous affectons d'un exposant $^\natural$ les objets associés à ces déformations (cf. \ref{mtht5} et \ref{mtht6}). 
L'isomorphisme de $\trT$-torseurs $\cL\stackrel{\sim}{\rightarrow} \cL^\natural$, 
$\psi\mapsto h\circ \psi$ \eqref{mtht3g} induit un isomorphisme $\hoR$-linéaire et $\Delta$-équivariant 
\begin{equation}\label{mtht100f}
\cF^\natural\stackrel{\sim}{\rightarrow}\cF,
\end{equation}
qui s'insère dans un diagramme commutatif \eqref{mtht5a}
\begin{equation}\label{mtht100g}
\xymatrix{
0\ar[r]&{\hoR}\ar[r]\ar@{=}[d]&{\cF^\natural}\ar[r]\ar[d]&{\txi^{-1}\tOmega^1_{R/\co_K} \otimes_R \hoR}\ar[r]\ar@{=}[d]&0\\
0\ar[r]&{\hoR}\ar[r]&{\cF}\ar[r]&{\txi^{-1}\tOmega^1_{R/\co_K} \otimes_R \hoR}\ar[r] & 0}
\end{equation} 
On en déduit un $\hoR$-isomorphisme $\Delta$-équivariant
\begin{equation}\label{mtht100h}
\cC^\natural\stackrel{\sim}{\rightarrow} \cC.
\end{equation}

L'isomorphisme de $\trT'$-torseurs $\cL'\stackrel{\sim}{\rightarrow} \cL'^\natural$, 
$\psi'\mapsto h'\circ \psi'$ \eqref{mtht3g} induit un isomorphisme $\hoRp$-linéaire et $\Delta'$-équivariant 
\begin{equation}\label{mtht100i}
\cF'^\natural\stackrel{\sim}{\rightarrow}\cF',
\end{equation}
qui s'insère dans un diagramme commutatif \eqref{mtht6a}
\begin{equation}\label{mtht100j}
\xymatrix{
0\ar[r]&{\hoRp}\ar[r]\ar@{=}[d]&{\cF'^\natural}\ar[r]\ar[d]&{\txi^{-1}\tOmega^1_{R'/\co_K} \otimes_{R'}\hoRp}\ar[r]\ar@{=}[d]&0\\
0\ar[r]&{\hoRp}\ar[r]&{\cF'}\ar[r]&{\txi^{-1}\tOmega^1_{R'/\co_K} \otimes_{R'}\hoRp}\ar[r] & 0}
\end{equation} 
On en déduit un $\hoRp$-isomorphisme $\Delta'$-équivariant
\begin{equation}\label{mtht100k}
\cC'^\natural\stackrel{\sim}{\rightarrow} \cC'.
\end{equation}

On vérifie aussitôt \eqref{mtht7} que le diagramme 
\begin{equation}\label{mtht100l}
\xymatrix{
\cC\ar[r]\ar[d]&{\cC'}\ar[d]\\
{\cC^\natural}\ar[r]&{\cC'^\natural}}
\end{equation}
où les flèches verticales sont les isomorphismes \eqref{mtht100h} et \eqref{mtht100k} et les flèches horizontales sont les homomorphismes \eqref{mtht7f},
est commutatif.

\subsection{}\label{mtht101}
Dans ce numéro, on remplace exceptionnellement l'hypothèse que $g$ admet une carte relativement adéquate \eqref{hypmdef2d} par les hypothèses 
suivantes : les morphismes $f$ et $f'$ admettent des cartes adéquates et le schéma logarithmique 
$(X',\cM_{X'})$ admet une carte fine et saturée $M'\rightarrow \Gamma(X',\cM_{X'})$ induisant un isomorphisme 
\begin{equation}\label{mtht101a}
M'\stackrel{\sim}{\rightarrow}\Gamma(X',\cM_{X'})/(X',\co_{X'}^\times).
\end{equation}
Ces trois cartes sont a priori indépendantes les unes des autres. 
De plus, on se donne une carte adéquate $((P,\gamma),(\mN,\iota),\vartheta)$ pour $f$
\eqref{cad1}. On peut alors considérer le diagramme commutatif 
\begin{equation}\label{mtht101b}
\xymatrix{
{(\hmX',\cM_{\hmX'})}\ar[r]^-(0.5){i'_{X'}}\ar[d]_{\hmg}&{(\tmX',\cM'_{\tmX'})}\ar[d]\\
{(\hmX,\cM_{\hmX})}\ar[r]^-(0.5){i_X}&{(\tmX,\cM_{\tmX})}}
\end{equation}
où $\hmg$ est le morphisme défini dans \eqref{mtht3h},
$\cM'_{\tmX'}$ (resp. $\cM_{\tmX}$) est la structure logarithmique sur $\tmX'$ (resp. $\tmX$) définie dans \ref{taht4}
et $i'_{X'}$ (resp. $i_X$) est l'immersion fermée exacte \eqref{taht4f} (resp. \eqref{taht4b}). On notera ici la différence entre $\cM'_{\tmX'}$
et $\cM_{\tmX}$ (cf. \ref{pmh70}). 
On dispose alors des représentations $\cF$ et $\cC$ de $\Delta$ définies dans  \ref{mtht5} relativement à $i_X$.  
Calquant \ref{mtht6}, on définit les représentations $\cF'$ et $\cC'$ de $\Delta'$ relativement à $i'_{X'}$. 
On voit aussitôt que les résultats \ref{mtht7}--\ref{mtht90} valent encore. 
En particulier, il existe un morphisme $\hoR$-linéaire et $\Delta'$-équivariant canonique 
\begin{equation}\label{mtht101c}
\cF\rightarrow \cF',
\end{equation}
et un homomorphisme $\Delta'$-équivariant de $\hoR$-algèbres 
\begin{equation}\label{mtht101d}
\cC\rightarrow \cC'
\end{equation}
qui le prolonge.

\subsection{}\label{mtht10}
Pour tout nombre rationnel $r\geq 0$, on note $\cF^{(r)}$ la $\hoR$-représentation de $\Delta$ 
déduite de $\cF$ \eqref{mtht5a} par image inverse par la multiplication par $p^r$ sur 
$\txi^{-1}\tOmega^1_{R/\co_K}\otimes_R\hoR$, de sorte qu'on a une suite exacte localement scindée de $\hoR$-modules
\begin{equation}\label{mtht10a}
0\rightarrow \hoR\longrightarrow \cF^{(r)}\rightarrow 
\txi^{-1}\tOmega^1_{R/\co_K}\otimes_R\hoR\rightarrow 0.
\end{equation}
Celle-ci induit pour tout entier $n\geq 1$, une suite exacte 
\begin{equation}\label{mtht10b}
0\rightarrow \rS^{n-1}_{\hoR}(\cF^{(r)})\rightarrow \rS^{n}_{\hoR}(\cF^{(r)})\rightarrow 
\rS^n_{\hoR}(\txi^{-1}\tOmega^1_{R/\co_K}
\otimes_R\hoR)\rightarrow 0.
\end{equation}
Les $\hoR$-modules $(\rS^{n}_{\hoR}(\cF^{(r)}))_{n\in \mN}$ forment donc un système inductif filtrant, 
dont la limite inductive 
\begin{equation}\label{mtht10c}
\cC^{(r)}=\underset{\underset{n\geq 0}{\longrightarrow}}\lim\ \rS^n_{\hoR}(\cF^{(r)})
\end{equation}
est naturellement munie d'une structure de $\hoR$-algèbre et d'une action de $\Delta$
par des automorphismes d'anneaux, compatible avec son action sur $\hoR$;  c'est l'algèbre de Higgs-Tate d'épaisseur $r$ associée 
à $(\tX,\cM_\tX)$ \eqref{taht10}. On note $\hcC^{(r)}$ le séparé complété $p$-adique de $\cC^{(r)}$ que l'on suppose toujours muni de 
la topologie $p$-adique. On munit $\hcC^{(r)}[\frac 1 p]$ de la topologie $p$-adique \eqref{notconv15}.

Pour tous nombres rationnels $r'\geq r\geq 0$, on a un $\hoR$-homomorphisme canonique injectif et $\Delta$-équivariant
$\alpha^{r,r'}\colon \cC^{(r')}\rightarrow \cC^{(r)}$.

\subsection{}\label{mtht11}
Pour tout nombre rationnel $r\geq 0$, on note $\cF'^{(r)}$ la $\hoRp$-représentation de $\Delta'$ 
déduite de $\cF'$ \eqref{mtht6a} par image inverse par la multiplication par $p^r$ sur 
$\txi^{-1}\tOmega^1_{R'/\co_K}\otimes_{R'}\hoRp$, de sorte qu'on a une suite exacte localement scindée de $\hoRp$-modules
\begin{equation}\label{mtht11a}
0\rightarrow \hoRp\longrightarrow \cF'^{(r)}\rightarrow 
\txi^{-1}\tOmega^1_{R'/\co_K}\otimes_{R'}\hoRp\rightarrow 0.
\end{equation}
Celle-ci induit pour tout entier $n\geq 1$, une suite exacte 
\begin{equation}\label{mtht11b}
0\rightarrow \rS^{n-1}_{\hoRp}(\cF'^{(r)})\rightarrow \rS^{n}_{\hoRp}(\cF'^{(r)})\rightarrow 
\rS^n_{\hoRp}(\txi^{-1}\tOmega^1_{R'/\co_K}
\otimes_{R'}\hoRp)\rightarrow 0.
\end{equation}
Les $\hoRp$-modules $(\rS^{n}_{\hoRp}(\cF'^{(r)}))_{n\in \mN}$ forment donc un système inductif filtrant, 
dont la limite inductive 
\begin{equation}\label{mtht11c}
\cC'^{(r)}=\underset{\underset{n\geq 0}{\longrightarrow}}\lim\ \rS^n_{\hoRp}(\cF'^{(r)})
\end{equation}
est naturellement munie d'une structure de $\hoRp$-algèbre et d'une action de $\Delta'$
par des automorphismes d'anneaux, compatible avec son action sur $\hoRp$;  c'est l'algèbre de Higgs-Tate d'épaisseur $r$ associée 
à $(\tX',\cM_{\tX'})$ \eqref{taht10}. 

Pour tous nombres rationnels $r'\geq r\geq 0$, on a un $\hoRp$-homomorphisme canonique injectif et $\Delta'$-équivariant
$\alpha'^{r,r'}\colon \cC'^{(r')}\rightarrow \cC'^{(r)}$. 

Le morphisme \eqref{mtht7d} induit un morphisme $\hoRp$-linéaire $\Delta'$-équivariant $\cF^{(r)}\otimes_{\hoR}\hoRp\rightarrow \cF'^{(r)}$ 
et par suite un homomorphisme de $\hoR$-algèbres $\Delta'$-équivariant
\begin{equation}\label{mtht11d}
\cC^{(r)}\rightarrow \cC'^{(r)}.
\end{equation}

Pour tout nombre rationnel $t\geq r$, on considère la $\cC^{(r)}$-algèbre
\begin{equation}\label{mtht11e}
\cC'^{(t,r)}=\cC'^{(t)}\otimes_{\cC^{(t)}}\cC^{(r)}
\end{equation}
déduite de $\cC'^{(t)}$ par changement de base par l'homomorphisme $\alpha^{t,r}\colon \cC^{(t)}\rightarrow \cC^{(r)}$ \eqref{mtht10}. 
On note $\hcC'^{(t,r)}$ le séparé complété $p$-adique de $\cC'^{(t,r)}$ que l'on suppose toujours muni de 
la topologie $p$-adique. On munit $\hcC'^{(t,r)}[\frac 1 p]$ de la topologie $p$-adique \eqref{notconv15}.
Les actions de $\Delta'$ sur les anneaux $\cC'^{(t)}$ et $\cC^{(r)}$
induisent une actions sur $\cC'^{(t,r)}$ par des automorphismes d'anneaux, compatible avec son action sur $\cC^{(r)}$.

Pour tous nombres rationnels  $t,t',r'$  tels que $t'\geq t\geq r$ et $t'\geq r'\geq r$, le diagramme 
\begin{equation}\label{mtht11f}
\xymatrix{
{\cC^{(r')}}\ar[d]_{\alpha^{r',r}}&&{\cC^{(t')}}\ar[ll]_{\alpha^{t',r'}}\ar[r]\ar[d]_{\alpha^{t',t}}&{\cC'^{(t')}}\ar[d]^{\alpha'^{t',t}}\\
{\cC^{(r)}}&&{\cC^{(t)}}\ar[ll]^{\alpha^{t,r}}\ar[r]&{\cC'^{(t)}}}
\end{equation}
où les flèches non-labellisées sont les homomorphismes \eqref{mtht11d}, est commutatif. On en déduit un homomorphisme canonique de $\cC^{(r')}$-algèbres
\begin{equation}\label{mtht11g}
\cC'^{(t',r')}\rightarrow \cC'^{(t,r)}.
\end{equation}

\begin{prop}\label{mtht20}
Soient $r,t,t'$ trois nombres rationnels tels que $t'>t>r\geq 0$, $\oR^\intern$ l'anneau défini dans \eqref{eccr40e}, $\Pi$ le groupe défini dans \eqref{eccr40c}. Alors,
\begin{itemize}
\item[{\rm (i)}] Pour tout entier $n\geq 1$, l'homomorphisme canonique 
\begin{equation}\label{mtht20a}
(\cC^{(r)}/p^n\cC^{(r)})\otimes_{\oR}\oR^\intern\rightarrow (\cC'^{(t,r)}/p^n\cC'^{(t,r)})^{\Pi}
\end{equation}
est $\alpha$-injectif \eqref{hypmdef100}. Notons $\cH^{(t,r)}_n$ son conoyau. 
\item[{\rm (ii)}] Il existe un entier $a\geq 0$, dépendant de $t$, $t'$ et $\ell=\dim(X'/X)$, 
mais pas des morphismes $f$, $f'$ et $g$ vérifiant les conditions de \ref{hypmdef2} et \ref{hypmdef3}, tel que pour tout entier $n\geq 1$, 
le morphisme canonique $\cH^{(t',r)}_n\rightarrow \cH^{(t,r)}_n$ soit annulé par $p^a$. 
\item[{\rm (iii)}] Il existe un entier $b\geq 0$, dépendant de $t$, $t'$ et $\ell$, 
mais pas des morphismes $f$, $f'$ et $g$ vérifiant les conditions de \ref{hypmdef2} et \ref{hypmdef3}, tel que pour tous entiers $n,q\geq 1$, 
le morphisme canonique
\begin{equation}\label{mtht20b}
\rH^q(\Pi,\cC'^{(t',r)}/p^n\cC'^{(t',r)})\rightarrow \rH^q(\Pi,\cC'^{(t,r)}/p^n\cC'^{(t,r)})
\end{equation}
soit annulé par $p^b$. 
\end{itemize}
\end{prop}

Le reste de cette section est consacré à la preuve de cet énoncé qui sera déduit dans \ref{mtht42} de \ref{mtht41}.

\begin{cor}\label{mtht21}
Soit $r$ un nombre rationnel $\geq 0$. 
\begin{itemize}
\item[{\rm (i)}] Le morphisme canonique 
\begin{equation}\label{mtht21a}
(\cC^{(r)}\hotimes_{\oR}\oR^\intern)[\frac 1 p]\rightarrow 
\underset{\underset{t\in \mQ_{>r}}{\longrightarrow}}{\lim} (\hcC'^{(t,r)}[\frac 1 p])^\Pi,
\end{equation}
où le produit tensoriel $\hotimes$ est complété pour la topologie $p$-adique, est un isomorphisme. 
\item[{\rm (ii)}] Pour tout entier $i\geq 1$, on a 
\begin{equation}\label{mtht21b}
\underset{\underset{t\in \mQ_{>r}}{\longrightarrow}}{\lim}\ 
\rH^i_\cont(\Pi,\hcC'^{(t,r)}[\frac 1 p])=0.
\end{equation}
\end{itemize}
\end{cor}

Reprenons les notations de \ref{mtht20}. Il résulte de \ref{mtht20}(i) et (\cite{gr} 2.4.2(ii)) que pour tout nombre rationnel $t>r$, la suite canonique
\begin{equation}\label{mtht21c}
0\rightarrow \cC^{(r)}\hotimes_{\oR}\oR^\intern\rightarrow (\hcC'^{(t,r)})^{\Pi}\rightarrow  \underset{\underset{n\geq 0}{\longleftarrow}}{\lim}\ \cH_n^{(t,r)}
\end{equation}
est $\alpha$-exacte. La proposition (i) s'ensuit en inversant $p$ puis en 
passant à la limite inductive sur les rationnels $t>r$,  compte tenu de \ref{mtht20}(ii).

Comme on a (\cite{agt} (II.3.10.2))
\begin{equation}\label{mtht21d}
\rR^1\underset{\underset{n}{\longleftarrow}}{\lim}\ (\cC^{(r)}/p^n\cC^{(r)})\otimes_{\oR}\oR^\intern=0,
\end{equation}
il résulte de \ref{mtht20}(i) et (\cite{gr} 2.4.2(ii)) que le morphisme canonique
\begin{equation}\label{mtht21e}
\rR^1\underset{\underset{n\geq 0}{\longleftarrow}}{\lim}\ (\cC'^{(t,r)}/p^n\cC'^{(t,r)})^\Pi\rightarrow \rR^1\underset{\underset{n\geq 0}{\longleftarrow}}{\lim}\ \cH_n^{(t,r)}
\end{equation}
est un $\alpha$-isomorphisme.

D'après (\cite{agt} (II.3.10.4) et (II.3.10.5)), pour tout $i \geq 1$, on a une suite exacte canonique
\[
0\rightarrow \rR^1\underset{\underset{n}{\longleftarrow}}{\lim}\ \rH^{i-1}(\Pi,\cC'^{(t,r)}/p^n\cC'^{(t,r)})\rightarrow 
\rH^i_\cont(\Pi,\hcC'^{(t,r)})\rightarrow \underset{\underset{n}{\longleftarrow}}{\lim}\ \rH^{i}(\Pi,\cC'^{(t,r)}/p^n\cC'^{(t,r)})\rightarrow 0.
\]
On en déduit, compte tenu de \eqref{mtht21e} et \ref{mtht20}(ii)-(iii), que pour tous entiers $t'>t>r$, il existe un nombre rationnel $\gamma \geq 0$
tel que le morphisme canonique
\begin{equation}\label{mtht21g}
\rH^i_\cont(\Pi,\hcC'^{(t',r)})\rightarrow \rH^i_\cont(\Pi,\hcC'^{(t,r)})
\end{equation}
soit annulé par $p^\gamma$. La proposition (ii) s'ensuit en passant à la limite inductive sur les rationnels $t>r$.

\subsection{}\label{mtht22}
Pour tout nombre rationnel $r\geq 0$, on désigne par $\cG^{(r)}$ la sous-$\hoR$-algèbre de $\cG$ \eqref{mtht4b} définie par \eqref{notconv9}
\begin{equation}\label{mtht22a}
\cG^{(r)}=\rS_{\hoR}(p^r\txi^{-1}\tOmega^1_{R/\co_K}\otimes_R\hoR)
\end{equation}
et par $\hcG^{(r)}$ son séparé complété $p$-adique que l'on suppose toujours muni de la topologie $p$-adique.
Pour tous nombres rationnels $r'\geq r\geq 0$, on a un homomorphisme injectif canonique 
$a^{r,r'}\colon \cG^{(r')}\rightarrow \cG^{(r)}$.

\subsection{}\label{mtht23}
Soient $(\tX_0,\cM_{\tX_0})$ la $(\cA_2(\oS),\cM_{\cA_2(\oS)})$-déformation lisse de 
$(\coX,\cM_\coX)$ définie par la carte $(P,\gamma)$ \eqref{hypmdef2d} (cf. \ref{pmh8}), 
$\cL_0$ le torseur de Higgs-Tate associé, $\psi_0\in \cL_0(\hmX)$ la section définie par la même carte \eqref{pmh8d}.
On désigne par $\bL_0$ le $\bT$-fibré principal homogène sur $\hmX$ associé à $\cL_0$ \eqref{taht5}
et par $\ttt$ l'isomorphisme de $\bT$-fibrés principaux homogènes sur $\hmX$
\begin{equation}\label{mtht23a}
\ttt\colon \bT\stackrel{\sim}{\rightarrow} \bL_0,\ \ \ v\mapsto v+\psi_0.
\end{equation}
La structure canonique de $\bT$-fibré principal homogène $\Delta$-équivariant sur $\bL_0$ \eqref{taht6}
se transporte par $\ttt$ en une structure de $\bT$-fibré principal homogène $\Delta$-équivariant sur $\bT$.
Cette dernière détermine une action à gauche de $\Delta$ sur $\bT$ compatible avec son action sur $\hmX$. 
On en déduit une action 
\begin{equation}\label{mtht23b}
\varphi\colon \Delta\rightarrow \Aut_\hRun(\cG)
\end{equation}
de $\Delta$ sur $\cG$ \eqref{mtht4b} par des automorphismes d'anneaux, 
compatible avec son action sur $\hoR$~;  pour tout $c\in \Delta$, 
$\varphi(c)$ est induit par l'automorphisme de $\bT$ défini par $c^{-1}$.

D'après \ref{pmh120}, pour tout nombre rationnel $r\geq 0$, la sous-$\hoR$-algèbre $\cG^{(r)}$ \eqref{mtht22a} de $\cG$ 
est stable par l'action $\varphi$ de $\Delta$ sur $\cG$, et les actions induites  
de $\Delta$ sur $\cG^{(r)}$ et $\hcG^{(r)}$ sont continues pour les topologies $p$-adiques. 
Sauf mention expresse du contraire, on munit $\cG^{(r)}$ et $\hcG^{(r)}$ des actions de $\Delta$ induites par $\varphi$.

\subsection{}\label{mtht27}
On reprend les notations de \ref{eccr42}. 
On désigne par $\fG$ la sous-$\hRun$-algèbre de $\cG$ \eqref{mtht4b} définie par \eqref{notconv9}
\begin{equation}\label{mtht27a}
\fG=\rS_{\hRun}(\txi^{-1}\tOmega^1_{R/\co_K}\otimes_R\hRun). 
\end{equation} 
D'après \ref{pmh11}, l'action $\varphi$ de $\Delta$ sur $\cG$ préserve $\fG$, 
et l'action induite sur $\fG$ se factorise à travers $\Delta_{p^\infty}$ \eqref{eccr43d2}. 
On note encore $\varphi$ l'action de $\Delta_{p^\infty}$ sur $\fG$ ainsi définie. 
En calquant la preuve de \ref{pmh120}, on montre que celle-ci est continue pour la topologie $p$-adique sur $\fG$. On pose
\begin{eqnarray}
\cG_\infty&=&\fG\otimes_{\hRun}\hRi,\label{mtht27b}\\
\cG_{p^\infty}&=&\fG\otimes_{\hRun}\hRpi.\label{mtht27c}
\end{eqnarray}
L'action $\varphi$ de $\Delta_{p^\infty}$ sur $\fG$ induit des actions de $\Delta_{p^\infty}$ sur $\cG_{p^\infty}$ et
de  $\Delta_{\infty}$ sur $\cG_{\infty}$.

Pour tout nombre rationnel $r\geq 0$, on désigne par $\cG^{(r)}_\infty$ 
la sous-$\hRi$-algèbre de $\cG_\infty$ définie par 
\begin{equation}\label{mtht27d}
\cG^{(r)}_\infty=\rS_{\hRi}(p^r\txi^{-1}\tOmega^1_{R/\co_K}\otimes_R\hRi),
\end{equation}
et par $\hcG^{(r)}_\infty$ son séparé complété $p$-adique.
Pour tout nombre rationnel $r'\geq r$, on a un homomorphisme injectif canonique 
$\cG^{(r')}_\infty\rightarrow \cG^{(r)}_\infty$. D'après la preuve de \ref{pmh120},
$\cG^{(r)}_\infty$ est stable par l'action de $\Delta_\infty$ sur $\cG_\infty=\cG^{(0)}_\infty$, 
et les actions induites de $\Delta_\infty$ sur $\cG^{(r)}_\infty$ et $\hcG^{(r)}_\infty$ 
sont continues pour les topologies $p$-adiques. Sauf mention expresse du contraire, 
on munit $\cG^{(r)}_\infty$ et $\hcG^{(r)}_\infty$ de ces actions et des topologies $p$-adiques.

On désigne par $\cG^{(r)}_{p^\infty}$ la sous-$\hRpi$-algèbre de $\cG_{p^\infty}$ définie par
\begin{equation}\label{mtht27e}
\cG^{(r)}_{p^\infty}=\rS_{\hRpi}(p^r\txi^{-1}\tOmega^1_{R/\co_K}\otimes_R\hRpi),
\end{equation}
et par $\hcG^{(r)}_{p^\infty}$ son séparé complété $p$-adique.
L'algèbre $\cG^{(r)}_{p^\infty}$ vérifie des propriétés analogues à celles vérifiées par $\cG^{(r)}_{\infty}$. 
En particulier, $\cG^{(r)}_{p^\infty}$ est stable par l'action de $\Delta_{p^\infty}$ sur 
$\cG_{p^\infty}=\cG^{(0)}_{p^\infty}$, et les actions induites de $\Delta_{p^\infty}$
sur $\cG^{(r)}_{p^\infty}$ et $\hcG^{(r)}_{p^\infty}$ sont continues pour les topologies $p$-adiques.  
Sauf mention expresse du contraire, on munit $\cG^{(r)}_{p^\infty}$ et $\hcG^{(r)}_{p^\infty}$ de ces actions
et des topologies $p$-adiques.

\subsection{}\label{mtht24}
Pour tout nombre rationnel $r\geq 0$, on désigne par $\cG'^{(r)}$ la sous-$\hoRp$-algèbre de $\cG'$ \eqref{mtht4e} définie par \eqref{notconv9}
\begin{equation}\label{mtht24a}
\cG'^{(r)}=\rS_{\hoRp}(p^r\txi^{-1}\tOmega^1_{R'/\co_K}\otimes_{R'}\hoRp).
\end{equation}
Pour tous nombres rationnels $r'\geq r\geq 0$, on a un homomorphisme injectif canonique 
$a'^{r,r'}\colon \cG'^{(r')}\rightarrow \cG'^{(r)}$. 

Le morphisme canonique $\tOmega^1_{R/\co_K}\otimes_RR'\rightarrow \tOmega^1_{R'/\co_K}$ induit un homomorphisme 
de $\hoR$-algèbres $\cG^{(r)}\rightarrow \cG'^{(r)}$ \eqref{mtht22a}. 
Pour tout nombre rationnel $t\geq r$, on considère la $\cG^{(r)}$-algèbre
\begin{equation}\label{mtht24b}
\cG'^{(t,r)}=\cG'^{(t)}\otimes_{\cG^{(t)}}\cG^{(r)}
\end{equation}
déduite de $\cG'^{(t)}$ par changement de base par l'homomorphisme canonique $a^{t,r}\colon \cG^{(t)}\rightarrow \cG^{(r)}$. 
On note $\hcG'^{(t,r)}$ le séparé complété $p$-adique de $\cG'^{(t,r)}$ que l'on suppose toujours muni de 
la topologie $p$-adique. 

Pour tous nombres rationnels  $t,t',r'$  tels que $t'\geq t\geq r$ et $t'\geq r'\geq r$, le diagramme 
\begin{equation}\label{mtht24c}
\xymatrix{
{\cG^{(r')}}\ar[d]_{a^{r',r}}&&{\cG^{(t')}}\ar[ll]_{a^{t',r'}}\ar[r]\ar[d]_{a^{t',t}}&{\cG'^{(t')}}\ar[d]^{a'^{t',t}}\\
{\cG^{(r)}}&&{\cG^{(t)}}\ar[ll]^{a^{t,r}}\ar[r]&{\cG'^{(t)}}}
\end{equation}
est commutatif. On en déduit un homomorphisme canonique de $\cG^{(r')}$-algèbres
\begin{equation}\label{mtht24d}
a'^{t',t,r',r}\colon \cG'^{(t',r')}\rightarrow \cG'^{(t,r)}.
\end{equation}

\subsection{}\label{mtht25}
Soient $(\tX'_0,\cM_{\tX'_0})$ la $(\cA_2(\oS),\cM_{\cA_2(\oS)})$-déformation lisse de 
$(\coX',\cM_{\coX'})$ définie par la carte $(P',\gamma')$ \eqref{hypmdef2d} (cf. \ref{pmh8}), 
$\cL'_0$ le torseur de Higgs-Tate associé, $\psi'_0\in \cL'_0(\hmX')$ la section définie par la même carte \eqref{pmh8d}.
On désigne par $\bL'_0$ le $\bT'$-fibré principal homogène sur $\hmX'$ associé à $\cL'_0$ \eqref{taht5}
et par $\ttt'$ l'isomorphisme de $\bT'$-fibrés principaux homogènes sur $\hmX'$
\begin{equation}\label{mtht25a}
\ttt'\colon \bT'\stackrel{\sim}{\rightarrow} \bL'_0,\ \ \ v\mapsto v+\psi'_0.
\end{equation}
La structure canonique de $\bT'$-fibré principal homogène $\Delta'$-équivariant sur $\bL'_0$ \eqref{taht6}
se transporte par $\ttt'$ en une structure de $\bT'$-fibré principal homogène $\Delta'$-équivariant sur $\bT'$.
Cette dernière détermine une action à gauche de $\Delta'$ sur $\bT'$ compatible avec son action sur $\hmX'$. 
On en déduit une action 
\begin{equation}\label{mtht25b}
\varphi'\colon \Delta'\rightarrow \Aut_{\hRun'}(\cG')
\end{equation}
de $\Delta'$ sur $\cG'$ \eqref{mtht4b} par des automorphismes d'anneaux, 
compatible avec son action sur $\hoRp$;  pour tout $\gamma\in \Delta'$, 
$\varphi'(\gamma)$ est induit par l'automorphisme de $\bT'$ défini par $\gamma^{-1}$.

D'après \ref{pmh120}, pour tout nombre rationnel $r\geq 0$, la sous-$\hoRp$-algèbre $\cG'^{(r)}$ \eqref{mtht24a} de $\cG'$ 
est stable par l'action $\varphi'$ de $\Delta'$ sur $\cG'$, et l'action induite  
de $\Delta'$ sur $\cG'^{(r)}$ est continue pour la topologie $p$-adique. 
Sauf mention expresse du contraire, on munit $\cG'^{(r)}$ de l'action de $\Delta'$ induite par $\varphi'$.

\subsection{}\label{mtht26}
D'après \eqref{hypmdef2d} et (\cite{ag} 5.1.11), avec les notations de \ref{notconv2}, on a un morphisme étale strict canonique 
\begin{equation}\label{mtht26a}
(\coX',\cM_{\coX'})\rightarrow (\coX,\cM_{\coX})\times_{\bA_P}\bA_{P'}, 
\end{equation}
le produit étant indifféremment pris dans la catégorie des schémas logarithmiques ou 
dans celle des schémas logarithmiques fins (\cite{agt} I.5.19). 
Par ailleurs, avec les notations de \ref{mtht23} et \ref{mtht25}, on a deux $(\tS,\cM_\tS)$-morphismes étales stricts canoniques \eqref{pmh8a}
\begin{equation}\label{mtht26b}
(\tX'_0,\cM_{\tX'_0})\rightarrow (\tS,\cM_\tS)\times_{\bA_\mN}\bA_{P'}\ \ \ {\rm et}\ \ \ (\tX_0,\cM_{\tX_0})\times_{\bA_P}\bA_{P'}\rightarrow (\tS,\cM_{\tS})\times_{\bA_\mN}\bA_{P'},
\end{equation}
qui induisent par réduction au-dessus de $(\coS,\cM_{\coS})$ les morphismes étales stricts canoniques
\begin{equation}\label{mtht26c}
(\coX',\cM_{\coX'})\rightarrow (\coS,\cM_\coS)\times_{\bA_\mN}\bA_{P'}\ \ \ {\rm et}\ \ \ (\coX,\cM_{\coX})\times_{\bA_P}\bA_{P'}\rightarrow (\coS,\cM_{\coS})\times_{\bA_\mN}\bA_{P'}.
\end{equation}
Le morphisme \eqref{mtht26a} se prolonge donc en un unique morphisme étale strict 
\begin{equation}\label{mtht26d}
(\tX'_0,\cM_{\tX'_0})\rightarrow (\tX_0,\cM_{\tX_0})\times_{\bA_P}\bA_{P'}
\end{equation}
au-dessus de $(\tS,\cM_\tS)\times_{\bA_\mN}\bA_{P'}$. Notons $\tg_0\colon (\tX'_0,\cM_{\tX'_0})\rightarrow (\tX_0,\cM_{\tX_0})$ le $(\tS,\cM_{\tS})$-morphisme induit. Celui-ci induit un morphisme $\bT'$-équivariant \eqref{mtht7c}
\begin{equation}\label{mtht26g}
\bL'_0\rightarrow \bL_0\times_{\hmX}\hmX'. 
\end{equation}

Le diagramme 
\begin{equation}\label{mtht26f}
\xymatrix{
{(\tmX',\cM_{\tmX'})}\ar[r]^{\tmg}\ar[d]_{\phi'_0}&{(\tmX,\cM_{\tmX})}\ar[d]^{\phi_0}\\
{(\tS,\cM_\tS)\times_{\bA_\mN}\bA_{P'}}\ar[r]&{(\tS,\cM_\tS)\times_{\bA_\mN}\bA_{P}}}
\end{equation}
où $\phi_0$ et $\phi'_0$ sont les morphismes définis par les flèches horizontales de \eqref{mtht3i}, est commutatif. 
Par ailleurs, le diagramme de morphismes canoniques
\begin{equation}
\xymatrix{
{(\hmX',\cM_{\hmX'})}\ar[r]^{\hmg}\ar[d]&{(\hmX,\cM_{\hmX})}\ar[d]\\
{(\coX',\cM_{\coX'})}\ar[r]\ar[d]&{(\coX,\cM_{\coX})}\ar[d]\\
{(\coS,\cM_\coS)\times_{\bA_\mN}\bA_{P'}}\ar[r]&{(\coS,\cM_\coS)\times_{\bA_\mN}\bA_{P}}}
\end{equation}
est commutatif. On en déduit que le diagramme 
\begin{equation}
\xymatrix{
{(\tmX',\cM_{\tmX'})}\ar[r]^{\tmg}\ar[d]^{\psi'_0}\ar@/_3pc/[dd]_{\phi'_0}&{(\tmX,\cM_{\tmX})}\ar[d]_{\psi_0}\ar@/^3pc/[dd]^{\phi_0}\\
{(\tX'_0,\cM_{\tX'_0})}\ar[r]^{\tg_0}\ar[d]&{(\tX_0,\cM_{\tX_0})}\ar[d]\\
{(\tS,\cM_\tS)\times_{\bA_\mN}\bA_{P'}}\ar[r]&{(\tS,\cM_\tS)\times_{\bA_\mN}\bA_{P}}}
\end{equation}
est commutatif \eqref{pmh8d}. Les sections $\psi_0$ et $\psi'_0$ sont donc compatibles \eqref{mtht7}. 
Compte tenu de \ref{mtht90}, on en déduit que le $\hoR$-homomorphisme 
canonique $\cG\rightarrow \cG'$ est $\Delta'$-équivariant pour les actions $\varphi$ de $\Delta$ sur $\cG$ \eqref{mtht23b} 
et $\varphi'$ de $\Delta'$ sur $\cG'$ \eqref{mtht25b}. 
Par suite, pour tout nombre rationnel $r\geq 0$, l'homomorphisme canonique $\cG^{(r)}\rightarrow \cG'^{(r)}$ est $\Delta'$-équivariant \eqref{mtht25}.
On en déduit, pour tous nombres rationnels $t\geq r\geq 0$, une action de $\Delta'$ sur $\cG'^{(t,r)}$ \eqref{mtht24b} 
par des automorphismes d'anneaux, compatible avec l'action de $\Delta$ sur $\cG^{(r)}$.  Cette action ainsi que 
l'action induite de $\Delta'$ sur $\hcG'^{(t,r)}$ sont continues pour les topologies $p$-adiques.

\subsection{}\label{mtht28}
On reprend les notations de \ref{eccr42}. 
On désigne par $\fG'$ la sous-$\hRunp$-algèbre de $\cG'$ \eqref{mtht4e} définie par \eqref{notconv9}
\begin{equation}\label{mtht28a}
\fG'=\rS_{\hRunp}(\txi^{-1}\tOmega^1_{R'/\co_K}\otimes_{R'}\hRunp). 
\end{equation} 
D'après \ref{pmh11}, l'action $\varphi'$ de $\Delta'$ sur $\cG'$ préserve $\fG'$, 
et l'action induite sur $\fG'$ se factorise à travers $\Delta'_{p^\infty}$ \eqref{eccr43b2}. 
On note encore $\varphi'$ l'action de $\Delta'_{p^\infty}$ sur $\fG'$ ainsi définie. 
Celle-ci est continue pour la topologie $p$-adique sur $\fG'$ (cf. \cite{agt} II.11.9). On pose
\begin{eqnarray}
\cG'_\infty&=&\fG'\otimes_{\hRunp}\hRip,\label{mtht28b}\\
\cG'_{p^\infty}&=&\fG'\otimes_{\hRunp}\hRpip.\label{mtht28c}
\end{eqnarray}
L'action $\varphi'$ de $\Delta'_{p^\infty}$ sur $\fG'$ induit des actions de $\Delta'_{p^\infty}$ sur $\cG'_{p^\infty}$ et
de  $\Delta'_{\infty}$ sur $\cG'_{\infty}$.

Soit $r$ un nombre rationnel $\geq 0$. On désigne par $\cG'^{(r)}_\infty$ 
la sous-$\hRip$-algèbre de $\cG'_\infty$ définie par 
\begin{equation}\label{mtht28d}
\cG'^{(r)}_\infty=\rS_{\hRip}(p^r\txi^{-1}\tOmega^1_{R'/\co_K}\otimes_{R'}\hRip).
\end{equation}
Pour tout nombre rationnel $r'\geq r$, on a un homomorphisme injectif canonique 
$\cG'^{(r')}_\infty\rightarrow \cG'^{(r)}_\infty$. D'après la preuve de \ref{pmh120},
$\cG'^{(r)}_\infty$ est stable par l'action de $\Delta'_\infty$ sur $\cG'_\infty=\cG'^{(0)}_\infty$, 
et l'action induite de $\Delta'_\infty$ sur $\cG'^{(r)}_\infty$ est continue pour la topologie $p$-adique. 

Le morphisme canonique $\tOmega^1_{R/\co_K}\otimes_RR'\rightarrow \tOmega^1_{R'/\co_K}$ induit un homomorphisme 
$\Delta'_\infty$-équivariant de $\hRi$-algèbres $\cG^{(r)}_\infty\rightarrow \cG'^{(r)}_\infty$ \eqref{mtht27d}. 
Pour tout nombre rationnel $t\geq r$, on considère la $\cG^{(r)}_\infty$-algèbre
\begin{equation}
\cG'^{(t,r)}_\infty=\cG'^{(t)}_\infty\otimes_{\cG^{(t)}_\infty}\cG^{(r)}_\infty
\end{equation}
déduite de $\cG'^{(t)}_\infty$ par changement de base par l'homomorphisme canonique $\cG^{(t)}_\infty\rightarrow \cG^{(r)}_\infty$. 
On note $\hcG'^{(t,r)}_\infty$ le séparé complété $p$-adique de $\cG'^{(t,r)}_\infty$ que l'on suppose toujours muni de 
la topologie $p$-adique. 
L'action de $\Delta'_\infty$ sur $\cG'^{(t)}_\infty$ induit une action de $\Delta'_\infty$ sur $\cG'^{(t,r)}_\infty$ 
par des automorphismes d'anneaux, compatible avec l'action de $\Delta_\infty$ sur $\cG^{(r)}_\infty$. 
Cette action ainsi que l'action induite de $\Delta'_\infty$ sur $\hcG'^{(t,r)}_\infty$ sont continues pour les topologies $p$-adiques.

On désigne par $\cG'^{(r)}_{p^\infty}$ la sous-$\hRpip$-algèbre de $\cG'_{p^\infty}$ définie par
\begin{equation}\label{mtht28e}
\cG'^{(r)}_{p^\infty}=\rS_{\hRpip}(p^r\txi^{-1}\tOmega^1_{R'/\co_K}\otimes_{R'}\hRpip).
\end{equation}
L'algèbre $\cG'^{(r)}_{p^\infty}$ vérifie des propriétés analogues à celles vérifiées par $\cG'^{(r)}_{\infty}$. 
En particulier, $\cG'^{(r)}_{p^\infty}$ est stable par l'action de $\Delta'_{p^\infty}$ sur 
$\cG'_{p^\infty}=\cG'^{(0)}_{p^\infty}$, et l'action induite de $\Delta'_{p^\infty}$
sur $\cG'^{(r)}_{p^\infty}$ est continue pour la topologie $p$-adique.

Le morphisme canonique $\tOmega^1_{R/\co_K}\otimes_RR'\rightarrow \tOmega^1_{R'/\co_K}$ induit un homomorphisme 
$\Delta'_{p^\infty}$-équivariant de $\hRpi$-algèbres $\cG^{(r)}_{p^\infty}\rightarrow \cG'^{(r)}_{p^\infty}$ \eqref{mtht27e}. 
Pour tout nombre rationnel $t\geq r$, on considère la $\cG^{(r)}_{p^\infty}$-algèbre
\begin{equation}
\cG'^{(t,r)}_{p^\infty}=\cG'^{(t)}_{p^\infty}\otimes_{\cG^{(t)}_{p^\infty}}\cG^{(r)}_{p^\infty}
\end{equation}
déduite de $\cG'^{(t)}_{p^\infty}$ par changement de base par l'homomorphisme canonique $\cG^{(t)}_{p^\infty}\rightarrow \cG^{(r)}_{p^\infty}$. 
On note $\hcG'^{(t,r)}_{p^\infty}$ le séparé complété $p$-adique de $\cG'^{(t,r)}_{p^\infty}$ que l'on suppose toujours muni de la topologie $p$-adique. 
L'action de $\Delta'_{p^\infty}$ sur $\cG'^{(t)}_{p^\infty}$ induit une action de $\Delta'_{p^\infty}$ sur $\cG'^{(t,r)}_{p^\infty}$ 
par des automorphismes d'anneaux, compatible avec l'action de $\Delta_{p^\infty}$ sur $\cG^{(r)}_{p^\infty}$. 
Cette action ainsi que l'action induite de $\Delta'_{p^\infty}$ sur $\hcG'^{(t,r)}_{p^\infty}$ sont continues pour les topologies $p$-adiques.

Nous considérons principalement les actions des sous-groupes 
$\fS\subset \Delta'$ \eqref{eccr43m}, $\fS_\infty \subset \Delta'_\infty$ \eqref{eccr43i1} et $\fS_{p^\infty} \subset \Delta'_{p^\infty}$ \eqref{eccr43i2} sur 
sur les algèbres introduites plus haut (cf. \ref{eccr34}).

\begin{lem}\label{mtht29}
Soient $r, t$ deux nombres rationnels tels que $t\geq r\geq 0$, $i$ un entier $\geq 0$. Alors,  
\begin{itemize}
\item[{\rm (i)}] Le morphisme canonique 
\begin{equation}\label{mtht29a}
\rH^i_\cont(\fS_{p^\infty},\cG'^{(t,r)}_{p^\infty}\hotimes_{R^\intern_{p^\infty}}R^\intern_\infty)\rightarrow \rH^i_\cont(\fS_\infty,\hcG'^{(t,r)}_\infty),
\end{equation}
où les anneaux $R^\intern_{\infty}$ et $R^\intern_{p^\infty}$ sont définis dans \eqref{eccr42h1} et \eqref{eccr42h2} 
et le produit tensoriel $\hotimes$ est complété pour la topologie $p$-adique, est un isomorphisme.
\item[{\rm (ii)}] Le morphisme canonique 
\begin{equation}\label{mtht29b}
\rH^i_\cont(\fS_\infty,\hcG'^{(t,r)}_\infty)\rightarrow \rH^i_\cont(\fS,\hcG'^{(t,r)})
\end{equation}
est un $\alpha$-isomorphisme.
\end{itemize}
\end{lem} 

(i) Soit $n$ un entier $\geq 0$. L'homomorphisme canonique
\begin{equation}\label{mtht29c}
\cG'^{(t,r)}_{p^\infty}\otimes_{R^\intern_{p^\infty}}(R^\intern_\infty/p^n R^\intern_\infty) \rightarrow (\cG'^{(t,r)}_\infty/p^n\cG'^{(t,r)}_\infty)^{\fN}
\end{equation}
est un isomorphisme en vertu de \ref{eccr9}(ii). Comme la $p$-dimension cohomologique de $\fN$ est nulle \eqref{eccr43l}
(\cite{serre1} I cor.~2 de prop.~14), le morphisme
canonique 
\begin{equation}\label{mtht29d}
\rH^i(\fS_{p^\infty},\cG'^{(t,r)}_{p^\infty}\otimes_{R^\intern_{p^\infty}}(R^\intern_\infty/p^n R^\intern_\infty))
\rightarrow \rH^i(\fS_\infty,\cG'^{(t,r)}_\infty/p^n\cG'^{(t,r)}_\infty)
\end{equation}
est un isomorphisme. On en déduit par (\cite{agt} (II.3.10.4) et (II.3.10.5)) que le morphisme \eqref{mtht29a} est un isomorphisme. 

(ii) Soit $n$ un entier $\geq 0$.  L'homomorphisme canonique
\begin{equation}\label{mtht29e}
\cG'^{(t,r)}_\infty/p^n\cG'^{(t,r)}_\infty\rightarrow (\cG'^{(t,r)}/p^n\cG'^{(t,r)})^{\Sigma'}
\end{equation}
est un $\alpha$-isomorphisme en vertu de (\cite{agt} II.6.22). 
On en déduit par (\cite{agt} II.6.20) que le morphisme canonique 
\begin{equation}\label{mtht29f}
\psi_n\colon \rH^i(\fS_\infty,\cG'^{(t,r)}_\infty/p^n\cG'^{(t,r)}_\infty)\rightarrow \rH^i(\fS,\cG'^{(t,r)}/p^n\cG'^{(t,r)})
\end{equation}
est un $\alpha$-isomorphisme. Notons $A_n$ (resp. $C_n$) le noyau (resp. conoyau) de $\psi_n$. 
Alors les $\co_\oK$-modules
\[
\underset{\underset{n\geq 0}{\longleftarrow}}{\lim}\ A_n, \ \ \ \underset{\underset{n\geq 0}{\longleftarrow}}{\lim}\ C_n, 
\ \ \ \rR^1\underset{\underset{n\geq 0}{\longleftarrow}}{\lim}\ A_n, \ \ \ 
\rR^1\underset{\underset{n\geq 0}{\longleftarrow}}{\lim}\ C_n
\]
sont $\alpha$-nuls en vertu de (\cite{gr} 2.4.2(ii)). Par suite, les morphismes 
\[
\underset{\underset{n\geq 0}{\longleftarrow}}{\lim}\ \psi_n
\ \ \ {\rm et} \ \ \ \rR^1\underset{\underset{n\geq 0}{\longleftarrow}}{\lim}\ \psi_n
\]
sont des $\alpha$-isomorphismes. On en déduit par (\cite{agt} (II.3.10.4) et (II.3.10.5)) que le morphisme 
\eqref{mtht29b} est un $\alpha$-isomorphisme.

\subsection{}\label{mtht30}
Posons $\lambda=\vartheta(1)$ \eqref{hypmdef2d}. 
L'homomorphisme $h^\gp\colon P^\gp\rightarrow P'^\gp$ est injectif, et le sous-groupe de torsion de $P'^\gp/h^\gp(P^\gp)$ est d'ordre premier à $p$
(\cite{ag} 5.1.11). Considérons la suite exacte canonique
\begin{equation}\label{mtht30a}
0\longrightarrow P^\gp/\lambda\mZ\stackrel{h^\gp}{\longrightarrow} P'^\gp/h^\gp(\lambda)\mZ \longrightarrow P'^{\gp}/h^\gp(P^{\gp}) \longrightarrow 0.
\end{equation}
Les $\mZ$-modules $P^\gp/\mZ\lambda$ et $P'^\gp/h^\gp(\lambda)\mZ$ sont libres de type fini d'après (\cite{ag} 4.2.2); notons $d$ et $d'$ leurs rangs
respectifs et posons $\ell=d'-d$. En vertu de (\cite{ag} (5.2.12.3)), on a un isomorphisme canonique 
\begin{equation}\label{mtht30b}
\fS_{p^\infty}\stackrel{\sim}{\rightarrow}\Hom(P'^{\gp}/h^\gp(P^{\gp}),\mZ_p(1)).
\end{equation}
On en déduit un isomorphisme 
\begin{equation}\label{mtht30c}
(P'^{\gp}/h^\gp(P^{\gp})) \otimes_\mZ\mZ_p\stackrel{\sim}{\rightarrow}\Hom(\fS_{p^\infty},\mZ_p(1)).
\end{equation}
D'après (\cite{ogus} IV 1.1.4), on a des isomorphismes linéaires canoniques \eqref{hypmdef4}
\begin{eqnarray}
\tOmega^1_{R'/R}&\stackrel{\sim}{\rightarrow}& (P'^{\gp}/h^\gp(P^{\gp}))\otimes_\mZ R',\label{mtht30d1}\\
\tOmega^1_{R'/\co_K}&\stackrel{\sim}{\rightarrow}& (P'^{\gp}/h^\gp(\lambda)\mZ)\otimes_\mZ R',\label{mtht30d2}\\
\tOmega^1_{R/\co_K}&\stackrel{\sim}{\rightarrow}& (P^{\gp}/\lambda\mZ) \otimes_\mZ R.\label{mtht30d3}
\end{eqnarray}

Soient $(t_i)_{1\leq i\leq \ell}$ des éléments de $P'^\gp$, $(t_i)_{\ell+1\leq i\leq d'}$ des éléments de $P^\gp$, 
tels que les images des $(t_i)_{1\leq i\leq \ell}$ dans $P'^\gp/h^\gp(P^\gp)\otimes_\mZ\mZ_p$ forment une $\mZ_p$-base,
que les images des $(t_i)_{\ell+1\leq i\leq d'}$ dans $P^\gp/\lambda\mZ$ forment une $\mZ$-base,
et que les images des $t_{1},\dots,t_{\ell},h^\gp(t_{\ell+1}),\dots,h^\gp(t_{d'})$ dans $P'^\gp/h^\gp(\lambda)\mZ$ forment une $\mZ$-base.
Pour alléger, on notera abusivement $(h^\gp(t_i))_{\ell+1\leq i\leq d'}$ simplement par $(t_i)_{\ell+1\leq i\leq d'}$, ce qui n'induit 
aucune ambiguïté puisque $h^\gp$ est injectif.

Les $(d\log(t_i))_{1\leq i\leq d'}$ forment une $R'$-base de $\tOmega^1_{R'/\co_K}$ \eqref{mtht30d2}. 
Pour tout $1\leq i\leq d'$ et tout $\un=(n_1,\dots,n_{d'})\in \mN^{d'}$, 
posons $y_i=\txi^{-1}d\log(t_i)\in \txi^{-1}\tOmega^1_{R'/\co_K}\subset \cG'_{p^\infty}$ \eqref{mtht28c},
$|\un|=\sum_{i=1}^{d'} n_i$ et $\uy^{\un}=\prod_{i=1}^{d'} y_i^{n_i}\in \cG'_{p^\infty}$.

On note $W$ l'anneau des vecteurs de Witt à coefficients dans $k$ relatif à $p$,
$K_0$ le corps des fractions de $W$ et $\fd$ la différente de l'extension $K/K_0$.
On pose $\rho=0$ dans le cas absolu \eqref{definf10} et $\rho=v(\pi\fd)$ dans le cas relatif. 
D'après \eqref{definf16a} et \eqref{definf17c}, on a un isomorphisme $\co_C$-linéaire canonique
\begin{equation}\label{mtht30e}
\co_C(1)\stackrel{\sim}{\rightarrow} p^{\rho+\frac{1}{p-1}}\txi\co_C.
\end{equation}
Pour tout $1\leq i\leq d'$, on désigne par $\chi_{t_i}$ l'image de $t_i$ par l'homomorphisme \eqref{mtht30c}, et par $\chi_i$ l'homomorphisme composé 
\begin{equation}\label{mtht30f}
\xymatrix{
{\fS_{p^\infty}}\ar[r]^-(0.5){\chi_{t_i}}& 
{\mZ_p(1)}\ar[r]^-(0.5){\log([\ ])}&{p^{\rho+\frac{1}{p-1}}\txi\co_C}},
\end{equation}
où  la seconde flèche est induite par l'isomorphisme \eqref{mtht30e}. On notera que $\chi_{t_i}=1$ pour tout $\ell+1\leq i\leq d'$.

D'après \ref{pmh11}, pour tout $\gamma\in \fS_{p^\infty}$ et tout $1\leq i\leq d'$, on a 
\begin{equation}\label{mtht30g}
\varphi'(\gamma)(y_i)=y_i-\txi^{-1}\chi_i(\gamma).
\end{equation} 

Soit $\zeta$ une $\mZ_p$-base de $\mZ_p(1)$.
Comme les $(\chi_{t_i})_{1\leq i\leq \ell}$ forment une $\mZ_p$-base de $\Hom(\fS_{p^\infty},\mZ_p(1))$,   
il existe une unique $\mZ_p$-base $(\gamma_i)_{1\leq i\leq \ell}$  de $\fS_{p^\infty}$ 
telle que  $\chi_{t_i}(\gamma_j)=\delta_{ij}\zeta$ pour tous $1\leq i,j\leq \ell$. 
Pour tout entier $1\leq i\leq \ell$, on désigne par $\upchi^{> i}_{p^\infty}$ le sous-groupe de \eqref{eccr12j}
\begin{equation}\label{mtht30h}
{\upchi}_{p^\infty}=\Hom(\fS_{p^\infty},\mu_{p^\infty}(\co_\oK))
\end{equation}
formé des homomorphismes $\nu\colon \fS_{p^\infty}\rightarrow \mu_{p^\infty}(\co_\oK)$ tels que $\nu(\gamma_j)=1$ pour tout $1\leq j\leq i$.
On a donc $\upchi^{>0}_{p^\infty}=\upchi_{p^\infty}$ et $\upchi^{>\ell}_{p^\infty}=0$.

On observera que $R'_{p^\infty}$ est séparé pour la topologie $p$-adique 
et s'identifie donc à un sous-anneau de $\hRpip$; cela résulte par exemple de la preuve de (\cite{agt} II.8.9), plus précisément,
avec les notations de {\em loc. cit.}, de (\cite{agt} (II.8.9.8) et (II.8.9.18)) et du fait que $\Spec(R_1)$ est un ouvert de $\Spec(C_1)$.

Soient $t,r$ deux nombres rationnels tels que $t\geq r\geq 0$. Pour tout élément $\un=(n_1,\dots,n_{d'})$ de $\mN^{d'}$, on pose
\begin{equation}
|\un|^{(t,r)} = t\sum_{i=1}^\ell n_i+r\sum_{i=\ell+1}^{d'} n_i.
\end{equation}
Pour tout entier $0\leq i\leq \ell$ et tout $\nu\in \upchi^{>i}_{p^\infty}$, 
tenant compte de \eqref{eccr16b} et avec les mêmes notations, 
on désigne par $\cG'^{(t,r), >i}_{p^\infty}(\nu)$ et $\cG'^{(t,r),>i}_{p^\infty}$
les sous-$R^\intern_{p^\infty}$-modules de $\cG'^{(t,r)}_{p^\infty}$ définis par 
\begin{eqnarray}
\cG'^{(t,r),>i}_{p^\infty}(\nu)&=&\bigoplus_{\un\in J_i}p^{|\un|^{(t,r)}} R'^{(\nu)}_{p^\infty}\uy^\un,\label{mtht30i}\\
\cG'^{(t,r),>i}_{p^\infty}&=&\bigoplus_{\nu'\in \upchi^{>i}_{p^\infty}}\cG'^{(t,r),>i}_{p^\infty}(\nu'),\label{mtht30j}
\end{eqnarray}
où $J_{i}$ est le sous-ensemble de $\mN^{d'}$ formé des éléments $\un=(n_1,\dots,n_{d'})$ tels que $n_{1}=\dots=n_i=0$.  
On note $\hcG'^{(t,r),>i}_{p^\infty}$ le séparé complété $p$-adique de $\cG'^{(t,r),>i}_{p^\infty}$. 
Comme $R'_{p^\infty}$ est séparé pour la topologie $p$-adique, il en est de même de $R'^{(\nu)}_{p^\infty}$, 
$\cG'^{(t,r),>i}_{p^\infty}(\nu)$ et $\cG'^{(t,r),>i}_{p^\infty}$. Par suite, $\cG'^{(t,r),>i}_{p^\infty}$ s'identifie à un 
sous-$R^\intern_{p^\infty}$-module de $\hcG'^{(t,r),>i}_{p^\infty}$.

La topologie $p$-adique de $R'_{p^\infty}$ étant induite par la topologie $p$-adique de $\hRpip$, 
on déduit facilement de \eqref{eccr16b} 
que la topologie $p$-adique de $\cG'^{(t,r),>i}_{p^\infty}$ est induite par la topologie $p$-adique de $\cG'^{(t,r)}_{p^\infty}$.
Par suite, $\hcG'^{(t,r),>i}_{p^\infty}$ est l'adhérence de $\cG'^{(t,r),>i}_{p^\infty}$ dans $\hcG'^{(t,r)}_{p^\infty}$. 

Il résulte de \eqref{eccr16b} et \eqref{mtht30g} que pour tout $1\leq j\leq \ell$ et tout $\nu\in \upchi^{>i}_{p^\infty}$, 
$\gamma_{j}$ préserve $\cG'^{(t,r),>i}_{p^\infty}(\nu)$ et donc aussi $\cG'^{(t,r),>i}_{p^\infty}$. 
Si $1\leq j\leq i$, $\gamma_{j}$ fixe $\cG'^{(t,r),>i}_{p^\infty}(\nu)$ et  $\cG'^{(t,r),>i}_{p^\infty}$.

\begin{prop}\label{mtht31}
Les hypothèses étant celles de \ref{mtht30}, soient, de plus, $i$ un entier tel que $1\leq i\leq \ell=d'-d$,
$t,t',r$ trois nombres rationnels tels que $t'>t>r\geq 0$. Alors,
\begin{itemize}
\item[{\rm (i)}] On a $\hcG'^{(t,r),>0}_{p^\infty}=\hcG'^{(t,r)}_{p^\infty}$ et 
$\hcG'^{(t,r),>\ell}_{p^\infty}=\cG^{(r)}_{p^\infty}\hotimes_{R_{p^\infty}}R^\intern_{p^\infty}$, 
où le produit tensoriel $\hotimes$ est complété pour la topologie $p$-adique.
\item[{\rm (ii)}] La suite 
\begin{equation}\label{mtht31a}
\xymatrix{
0\ar[r]&{\hcG'^{(t,r),>i}_{p^\infty}}\ar[r]^-(0.45){u}&{\hcG'^{(t,r),>i-1}_{p^\infty}}\ar[r]^{\gamma_i-\id}&{\hcG'^{(t,r),>i-1}_{p^\infty}}},
\end{equation}
où $u$ est le morphisme canonique, est exacte. 
\item[{\rm (iii)}] Il existe un entier $\alpha\geq 0$, dépendant de $t$ et $t'$ 
mais pas des morphismes $f$, $f'$ et $g$ vérifiant les conditions de \ref{hypmdef2} et \ref{hypmdef3}, tel que l'on ait
\begin{equation}\label{mtht31b}
p^\alpha\cdot \hcG'^{(t',r),>i-1}_{p^\infty}
\subset (\gamma_i-\id)(\hcG'^{(t,r),>i-1}_{p^\infty}).
\end{equation}
\end{itemize}
\end{prop}

(i) Notant $e$ le caractère trivial de $\upchi_{p^\infty}$, on a  $R'^{(e)}_{p^\infty}=R^\intern_{p^\infty}$ d'après \ref{eccr13}.  
On en déduit que $\hcG'^{(t,r),>\ell}_{p^\infty}=\cG^{(r)}_{p^\infty}\hotimes_{R_{p^\infty}}R^\intern_{p^\infty}$. 
Posons
\begin{eqnarray}
\rS'^{(r)}&=&\rS_{R'_{p^\infty}}(p^r\txi^{-1}\tOmega^1_{R'/\co_K}\otimes_{R'}R'_{p^\infty}),\\
\rS^{(r)}&=&\rS_{R'_{p^\infty}}(p^r\txi^{-1}\tOmega^1_{R/\co_K}\otimes_{R'}R'_{p^\infty}).
\end{eqnarray}
Compte tenu de \eqref{eccr16b}, on a  
\begin{equation}
\cG'^{(t,r),>0}_{p^\infty}=\rS'^{(t)}\otimes_{\rS^{(t)}} \rS^{(r)},
\end{equation}
ce qui implique que $\hcG'^{(t,r),>0}_{p^\infty}=\hcG'^{(t,r)}_{p^\infty}$.

(ii) Comme pour tout entier $n\geq 0$, on a clairement 
\begin{equation}\label{mtht31d}
p^n \cdot\cG'^{(t,r),>i}_{p^\infty}=\cG'^{(t,r),>i}_{p^\infty}\cap (p^n\cdot\cG'^{(t,r),>i-1}_{p^\infty}),
\end{equation} 
le morphisme canonique $u\colon \hcG'^{(t,r),>i}_{p^\infty}\rightarrow \hcG'^{(t,r),>i-1}_{p^\infty}$ est injectif. 

Pour tout $\nu\in \upchi^{>i-1}_{p^\infty}$, notons $(R'^{(\nu)}_{p^\infty})^\wedge$ et $(\cG'^{(t,r),>i-1}_{p^\infty}(\nu))^\wedge$
les séparés complétés $p$-adiques de $R'^{(\nu)}_{p^\infty}$ et $\cG'^{(t,r),>i-1}_{p^\infty}(\nu)$, respectivement. 
On observera que $(\cG'^{(t,r),>i-1}_{p^\infty}(\nu))^\wedge$ s'identifie à un sous-$R^\intern_{p^\infty}$-module de $\hcG'^{(t,r),>i-1}_{p^\infty}$, stable par $\gamma_i$. 
Compte tenu de \eqref{mtht30j}, tout élément $x$ de $\hcG'^{(t,r),>i-1}_{p^\infty}$ peut s'écrire comme la somme d'une série 
\begin{equation}\label{mtht31dd}
\sum_{\nu\in \upchi^{>i-1}_{p^\infty}}x_\nu,
\end{equation}
où $x_\nu\in (\cG'^{(t,r),>i-1}_{p^\infty}(\nu))^\wedge$ et, pour tout entier $n\geq 0$, sauf pour un nombre fini de $\nu\in \upchi^{>i-1}_{p^\infty}$, 
$x_\nu\in p^n (\cG'^{(t,r),>i-1}_{p^\infty}(\nu))^\wedge$.
Pour qu'un tel élément $x$ soit nul, il faut et il suffit que $x_\nu$ soit nul pour tout $\nu\in \upchi^{>i-1}_{p^\infty}$.

Soit $\nu\in \upchi^{>i-1}_{p^\infty}$. Pour tout
\begin{equation}\label{mtht31f}
z=\sum_{\un\in J_{i-1}}p^{|\un|^{(t,r)}}a_\un \uy^\un\in \cG'^{(t,r),>i-1}_{p^\infty}(\nu),
\end{equation}
on a \eqref{mtht30g}
\begin{equation}\label{mtht31g}
(\gamma_i-\id)(z)=\sum_{\un\in J_{i-1}}p^{|\un|^{(t,r)}}b_\un \uy^\un\in \cG'^{(t,r),>i-1}_{p^\infty}(\nu),
\end{equation}
où pour tout $\un=(n_1,\dots,n_{d'})\in J_{i-1}$, 
\begin{equation}\label{mtht31h}
b_\un=(\nu(\gamma_i)-1)a_\un+\sum_{\um=(m_1,\dots,m_{d'})
\in J_{i-1}(\un)} p^{t(m_i-n_i)}\binom{m_i}{n_i} \nu(\gamma_i)a_{\um}w^{m_i-n_i},
\end{equation}
$J_{i-1}(\un)$ désigne le sous-ensemble de $J_{i-1}$ formé des éléments $\um=(m_1,\dots,m_{d'})$ tels que 
$m_j=n_j$ pour $j\not=i$ et $m_i>n_i$, 
et $w=-\txi^{-1}\log([\zeta])$ est un élément de valuation $\rho+\frac{1}{p-1}$ de $\co_C$ \eqref{mtht30e}.

Soit $z$ un élément de $(\cG'^{(t,r),>i-1}_{p^\infty}(\nu))^\wedge$. Alors, $z$ peut s'écrire comme somme d'une série 
\begin{equation}
z=\sum_{\un\in J_{i-1}}p^{|\un|^{(t,r)}}a_\un \uy^\un,
\end{equation}
où $a_\un \in (R'^{(\nu)}_{p^\infty})^\wedge$ et $a_\un$ tend vers $0$ quand $|\un|$ tend vers l'infini. 
Comme $R'^{(\nu)}_{p^\infty}$ est $\co_C$-plat, il en est de même de $(R'^{(\nu)}_{p^\infty})^\wedge$ (cf. la preuve de \cite{agt} II.6.14). 
Par suite, pour que $z$ soit nul, il faut et il suffit que $a_\un$ soit nul pour tout $\un\in J_{i-1}$.
On voit aussitôt que $(\gamma_i-\id)(z)$ est encore donné par la formule \eqref{mtht31g}.

Supposons que $\gamma_i(z)=z$ et $\nu(\gamma_i)\not=1$.  
Comme $(R'^{(\nu)}_{p^\infty})^\wedge$ est $\co_C$-plat et que $v(\nu(\gamma_i)-1)\leq \frac{1}{p-1}$, 
pour tout  $\un=(n_1,\dots,n_{d'})\in J_{i-1}$, on a 
\[
a_\un=-(\nu(\gamma_i)-1)^{-1}\sum_{\um=(m_1,\dots,m_{d'})
\in J_{i-1}(\un)} p^{t(m_i-n_i)}\binom{m_i}{n_i} \nu(\gamma_i)a_{\um}w^{m_i-n_i}.
\]
On en déduit que pour tout $\alpha\in \mN$ et tout $\un\in J_{i-1}$,
on a $a_\un\in p^{t\alpha} (R'^{(\nu)}_{p^\infty})^\wedge$ (on le prouve par récurrence sur $\alpha$);
donc $z=0$ puisque $(R'^{(\nu)}_{p^\infty})^\wedge$ est séparé pour la topologie $p$-adique. 
Par suite, $\gamma_i-\id$ est injectif sur $(\cG'^{(t,r),>i-1}_{p^\infty}(\nu))^\wedge$. 

Supposons que $\gamma_i(z)=z$ et $\nu(\gamma_i)=1$, de sorte que $\nu\in \upchi^{>i}_{p^\infty}$. 
Alors, pour tout $\un=(n_1,\dots,n_{d'})\in J_{i-1}$, si on pose 
$\un'=(n'_1,\dots,n'_{d'})\in J_{i-1}(\un)$ avec $n'_i=n_i+1$, on a 
\[
(n_i+1)!a_{\un'}=-\sum_{\um=(m_1,\dots,m_{d'})
\in J_{i-1}(\un')} p^{t(m_i-n_i-1)}m_i! a_{\um}\frac{w^{m_i-n_i-1}}{(m_i-n_i)!}.
\]
On a $w^{m-1}/m!\in \co_C$ pour tout entier $m\geq 1$. On en déduit que 
pour tout $\alpha\in \mN$ et tout $\un=(n_1,\dots,n_{d'})\in J_{i-1}$ tel que $n_i\geq 1$, 
on a $n_i!a_{\un}\in p^{t\alpha} (R'^{(\nu)}_{p^\infty})^\wedge$ (on le prouve par récurrence sur $\alpha$); donc $a_\un=0$. 
Par suite $z \in (\cG'^{(t,r),>i}_{p^\infty}(\nu))^\wedge$. On en déduit que la suite 
\begin{equation}
\xymatrix{
0\ar[r]&{\hcG'^{(t,r),>i}_{p^\infty}(\nu)}\ar[r]^-(0.45){u_\nu}&{\hcG'^{(t,r),>i-1}_{p^\infty}(\nu)}\ar[r]^{\gamma_i-\id}&{\hcG'^{(t,r),>i-1}_{p^\infty}(\nu)}},
\end{equation}
où $u_\nu$ est le morphisme canonique, est exacte.
 
Il résulte de ce qui précède que la suite \eqref{mtht31a} est exacte. 

(iii) Compte tenu de \eqref{mtht31dd}, il suffit de montrer qu'il existe un entier $\alpha\geq 0$ ne dépendant que de $t$ et $t'$
tel que pour tout $\nu\in \upchi^{>i-1}_{p^\infty}$, on ait
\begin{equation}\label{mtht31aa}
p^\alpha\cdot (\cG'^{(t',r),>i-1}_{p^\infty}(\nu))^\wedge\subset  (\gamma_i-\id)((\cG'^{(t,r),>i-1}_{p^\infty}(\nu))^\wedge).
\end{equation}

Supposons $\nu(\gamma_i)\not=1$. D'après \eqref{mtht31g} appliquée aux éléments de $(\cG'^{(t,r),>i-1}_{p^\infty}(\nu))^\wedge$, on a 
\[
(\nu(\gamma_i)-1)(\cG'^{(t,r),>i-1}_{p^\infty}(\nu))^\wedge\subset  (\gamma_i-\id)((\cG'^{(t,r),>i-1}_{p^\infty}(\nu))^\wedge)
+(\nu(\gamma_i)-1)p^t (\cG'^{(t,r),>i-1}_{p^\infty}(\nu))^\wedge.
\]
On en déduit que 
\begin{equation}\label{mtht31ad}
(\nu(\gamma_i)-1)(\cG'^{(t,r),>i-1}_{p^\infty}(\nu))^\wedge\subset  (\gamma_i-\id)((\cG'^{(t,r),>i-1}_{p^\infty}(\nu))^\wedge).
\end{equation}
On obtient l'inclusion recherchée \eqref{mtht31aa} avec $\alpha=1$ car $v(\nu(\gamma_i)-1)\leq \frac{1}{p-1}$. 

Supposons que $\nu(\gamma_i)=1$, de sorte que $\nu\in \upchi^{>i}_{p^\infty}$. 
Compte tenu de \eqref{mtht31g}, pour établir l'inclusion recherchée \eqref{mtht31aa}, 
il suffit de montrer qu'il existe un entier $\alpha\geq 0$ ne dépendant que de $t$ et $t'$ tel que pour tout élément
\begin{equation}\label{mtht31ae}
\sum_{\un\in J_{i-1}}p^{|\un|^{(t',r)}}b_\un \uy^\un\in \hcG'^{(t',r),>i-1}_{p^\infty}(\nu),
\end{equation}
où $b_\un \in (R'^{(\nu)}_{p^\infty})^\wedge$ et $b_\un$ tend vers $0$ quand $|\un|$ tend vers l'infini,
le système d'équations linéaires définies, pour $\un=(n_1,\dots,n_{d'})\in J_{i-1}$, par 
\begin{equation}\label{mtht31af}
p^\alpha p^{(t'-t)\sum_{j=1}^\ell n_j}n_i!b_\un=\sum_{\um=(m_1,\dots,m_{d'})
\in J_{i-1}(\un)} p^{t(m_i-n_i)}m_i! a_{\um}\frac{w^{m_i-n_i}}{(m_i-n_i)!},
\end{equation}
admette une solution $a_\um\in (R'^{(\nu)}_{p^\infty})^\wedge$ pour $\um\in J_{i-1}$ telle que 
$a_\um$ tend vers $0$ quand $|\um|$ tend vers l'infini. Pour $\un=(n_1,\dots,n_{d'}) \in J_{i-1}$, posons
\begin{eqnarray}
a'_\un&=&n_i!p^{-\frac{n_i}{p-1}}a_\un,\label{mtht31ba}\\
b'_\un&=&p^\alpha p^{(t'-t)\sum_{j=1}^\ell n_j} n_i!p^{-\frac{n_i}{p-1}}b_\un,\label{mtht31bb}
\end{eqnarray}
de sorte que l'équation \eqref{mtht31af} devient 
\begin{equation}\label{mtht31ag}
b'_\un=p^{t+\frac{1}{p-1}}w\sum_{\um=(m_1,\dots,m_{d'})
\in J_{i-1}(\un)} p^{(t+\frac{1}{p-1})(m_i-n_i-1)}a'_{\um}\frac{w^{(m_i-n_i-1)}}{(m_i-n_i)!}. 
\end{equation}
Considérons l'endomorphisme $R^\intern_{p^\infty}$-linéaire $\Phi$ de $(\oplus_{\un\in J_{i-1}}R'^{(\nu)}_{p^\infty})^\wedge$ défini,  
pour une suite $(x_\un)_{\un\in J_{i-1}}$ d'éléments de $(R'^{(\nu)}_{p^\infty})^\wedge$ tendant vers $0$ quand
$|\un|$ tend vers l'infini, par
\begin{equation}\label{mtht31ak}
\Phi(\sum_{\un\in J_{i-1}}x_\un)= \sum_{\un\in J_{i-1}}z_\un,
\end{equation}
où pour $\un=(n_1,\dots,n_{d'})\in J_{i-1}$, 
\begin{equation}
z_\un=\sum_{\um=(m_1,\dots,m_{d'})
\in \{\un\}\cup J_{i-1}(\un)} p^{(t+\frac{1}{p-1})(m_i-n_i)}x_{\um}\frac{w^{(m_i-n_i)}}{(m_i-n_i+1)!}.
\end{equation}
Comme $\Phi$ est congru à l'identité modulo $p^{t+\frac{1}{p-1}}$, il est surjectif en vertu (\cite{egr1} 1.8.5). 
Par suite, pour toute suite $b'_\un\in p^{t+\frac{1}{p-1}}w (R'^{(\nu)}_{p^\infty})^\wedge$ pour $\un \in J_{i-1}$ 
tendant vers $0$ quand $|\un|$ tend vers l'infini, 
l'équation \eqref{mtht31ag} admet une solution $a'_\um\in (R'^{(\nu)}_{p^\infty})^\wedge$ pour $\um\in J_{i-1}$ 
tendant vers $0$ quand $|\um|$ tend vers l'infini; il suffit de prendre $a'_\un=x_\un$ un antécédent par
le morphisme $\Phi$ \eqref{mtht31ak} de la suite $z_\un=p^{-t-\frac{1}{p-1}}w^{-1} b'_\un$.  
D'autre part, $v(w)=\rho+\frac{1}{p-1}$ et il existe un entier $\alpha\geq 0$ tel que pour tout $n\in \mN$, on ait 
\begin{equation}
(t'-t)n+v(n!)-\frac{n}{p-1}+\alpha\geq t+\rho+\frac{2}{p-1}. 
\end{equation}
L'assertion recherchée s'ensuit en prenant pour $b'_\un$ pour
$\un\in J_{i-1}$ les éléments définis par \eqref{mtht31bb}.
On notera que $v(n!)\leq \frac{n}{p-1}$ pour tout entier $n\geq 0$ \eqref{mtht31ba}.

\subsection{}\label{mtht33}
Reprenons les notations de \ref{mtht30}. Soient $t,r$ deux nombres rationnels tels que $t\geq r\geq  0$, $i$ un entier tel que $0\leq i \leq \ell$. 
Pour alléger les notations, posons \eqref{mtht30j}
\begin{equation}
G'^{(t,r)}=\cG'^{(t,r),>0}_{p^\infty}, 
\end{equation}
qui est une $R'_{p^\infty}$-algèbre d'après \eqref{eccr16b}. On a $\cG'^{(t,r)}_{p^\infty}=G'^{(t,r)}\otimes_{R'_{p^\infty}}\hRpip$ \eqref{mtht28c}. 
Comme l'anneau $R'_{p^\infty}$ est séparé pour la topologie $p$-adique \eqref{mtht30}, 
$G'^{(t,r)}$ s'identifie à une sous-$R'_{p^\infty}$-algèbre de $\hcG'^{(t,r)}_{p^\infty}$, stable par l'action de $\fS_{p^\infty}$ d'après \eqref{mtht30g}. 
Le séparé complété $p$-adique de $G'^{(t,r)}$ s'identifie à $\hcG'^{(t,r)}_{p^\infty}$. Par ailleurs, avec les notations de \eqref{eccr34a}, le morphisme canonique 
\begin{equation}
G'^{(t,r)}\hotimes_{R^\intern_{p^\infty}}R^\intern_\infty\rightarrow \cG'^{(t,r)}_{p^\infty}\hotimes_{R^\intern_{p^\infty}}R^\intern_\infty,
\end{equation}
où le produit tensoriel $\hotimes$ est complété pour la topologie $p$-adique, est un isomorphisme.  
Pour tout entier $0\leq i\leq \ell$, le $R^\intern_\infty$-module
$\cG'^{(t,r),>i}_{p^\infty}\otimes_{R^\intern_{p^\infty}}R^\intern_\infty$ étant facteur direct de $G'^{(t,r)}\otimes_{R^\intern_{p^\infty}}R^\intern_\infty$,
$\cG'^{(t,r),>i}_{p^\infty}\hotimes_{R^\intern_{p^\infty}}R^\intern_\infty$ est un sous-$\hRinterni$-module de 
$\cG'^{(t,r)}_{p^\infty}\hotimes_{R^\intern_{p^\infty}}R^\intern_\infty$.

\begin{prop}\label{mtht32}
Les hypothèses étant celles de \ref{mtht33}, soient, de plus, $i$ un entier tel que $1\leq i\leq \ell=d'-d$,
$t,t',r$ trois nombres rationnels tels que $t'>t>r\geq 0$. Alors,
\begin{itemize}
\item[{\rm (i)}] On a $\cG'^{(t,r),>0}_{p^\infty}\hotimes_{R^\intern_{p^\infty}}R^\intern_\infty=\cG'^{(t,r)}_{p^\infty}\hotimes_{R^\intern_{p^\infty}}R^\intern_\infty$ et 
$\cG'^{(t,r),>\ell}_{p^\infty}\hotimes_{R^\intern_{p^\infty}}R^\intern_\infty=\cG^{(t,r)}_{p^\infty}\hotimes_{R_{p^\infty}}R^\intern_\infty$, 
où le produit tensoriel $\hotimes$ est complété pour la topologie $p$-adique.
\item[{\rm (ii)}] La suite 
\begin{equation}\label{mtht32a}
\xymatrix{
0\ar[r]&{\cG'^{(t,r),>i}_{p^\infty}\hotimes_{R^\intern_{p^\infty}}R^\intern_\infty}\ar[r]^-(0.45){u}&
{\cG'^{(t,r),>i-1}_{p^\infty}\hotimes_{R^\intern_{p^\infty}}R^\intern_\infty}\ar[r]^{\gamma_i-\id}&{\cG'^{(t,r),>i-1}_{p^\infty}\hotimes_{R^\intern_{p^\infty}}R^\intern_\infty}},
\end{equation}
où $u$ est le morphisme canonique, est exacte. 
\item[{\rm (iii)}] Il existe un entier $\alpha\geq 0$, dépendant de $t$ et $t'$ 
mais pas des morphismes $f$, $f'$ et $g$ vérifiant les conditions de \ref{hypmdef2} et \ref{hypmdef3}, tel que l'on ait
\begin{equation}\label{mtht32b}
p^\alpha(\cG'^{(t',r),>i-1}_{p^\infty}\hotimes_{R^\intern_{p^\infty}}R^\intern_\infty)
\subset (\gamma_i-\id)(\cG'^{(t,r),>i-1}_{p^\infty}\hotimes_{R^\intern_{p^\infty}}R^\intern_\infty).
\end{equation}
\end{itemize}
\end{prop}

(i) Cela résulte ausitôt de \ref{mtht31}(i). 

(ii) Comme pour tout entier $n\geq 0$, on a clairement 
\begin{equation}
p^n (\cG'^{(t,r),>i}_{p^\infty}\otimes_{R^\intern_{p^\infty}}R^\intern_\infty)=(\cG'^{(t,r),>i}_{p^\infty}\otimes_{R^\intern_{p^\infty}}R^\intern_\infty)
\cap p^n(\cG'^{(t,r),>i-1}_{p^\infty}\otimes_{R^\intern_{p^\infty}}R^\intern_\infty),
\end{equation} 
le morphisme canonique $u\colon \cG'^{(t,r),>i}_{p^\infty}\hotimes_{R^\intern_{p^\infty}}R^\intern_\infty\rightarrow 
\cG'^{(t,r),>i-1}_{p^\infty}\hotimes_{R^\intern_{p^\infty}}R^\intern_\infty$ est injectif. 

Pour tout $\nu\in \upchi^{>i-1}_{p^\infty}$, notons $R'^{(\nu)}_{p^\infty}\hotimes_{R^\intern_{p^\infty}}R^\intern_\infty$ et 
$\cG'^{(t,r),>i-1}_{p^\infty}(\nu)\hotimes_{R^\intern_{p^\infty}}R^\intern_\infty$
les séparés complétés $p$-adiques de $R'^{(\nu)}_{p^\infty}\otimes_{R^\intern_{p^\infty}}R^\intern_\infty$ et 
$\cG'^{(t,r),>i-1}_{p^\infty}(\nu)\otimes_{R^\intern_{p^\infty}}R^\intern_\infty$, respectivement. 
On observera que $\cG'^{(t,r),>i-1}_{p^\infty}(\nu)\hotimes_{R^\intern_{p^\infty}}R^\intern_\infty$ 
s'identifie à un sous-$\hRinternpi$-module de $\cG'^{(t,r),>i-1}_{p^\infty}\hotimes_{R^\intern_{p^\infty}}R^\intern_\infty$, stable par $\gamma_i$. 
Compte tenu de \eqref{mtht30j}, tout élément $x$ de $\cG'^{(t,r),>i-1}_{p^\infty}\hotimes_{R^\intern_{p^\infty}}R^\intern_\infty$ 
peut s'écrire comme la somme d'une série 
\begin{equation}
\sum_{\nu\in \upchi^{>i-1}_{p^\infty}}x_\nu,
\end{equation}
où $x_\nu\in \cG'^{(t,r),>i-1}_{p^\infty}(\nu)\hotimes_{R^\intern_{p^\infty}}R^\intern_\infty$ 
et, pour tout entier $n\geq 0$, sauf pour un nombre fini de $\nu\in \upchi^{>i-1}_{p^\infty}$, 
$x_\nu\in p^n (\cG'^{(t,r),>i-1}_{p^\infty}(\nu)\hotimes_{R^\intern_{p^\infty}}R^\intern_\infty)$.
Pour qu'un tel élément $x$ soit nul, il faut et il suffit que $x_\nu$ soit nul pour tout $\nu\in \upchi^{>i-1}_{p^\infty}$.

Soit $\nu\in \upchi^{>i-1}_{p^\infty}$. D'après \eqref{mtht31g} et avec les mêmes notations, pour tout
\begin{equation}\label{mtht32f}
z=\sum_{\un\in J_{i-1}}p^{|\un|^{(t,r)}}a_\un \uy^\un\in \cG'^{(t,r),>i-1}_{p^\infty}(\nu)\otimes_{R^\intern_{p^\infty}}R^\intern_\infty,
\end{equation}
où $a_\un\in R'^{(\nu)}_{p^\infty}\otimes_{R^\intern_{p^\infty}}R^\intern_\infty$, on a \eqref{mtht30g}
\begin{equation}\label{mtht32g}
(\gamma_i-\id)(z)=\sum_{\un\in J_{i-1}}p^{|\un|^{(t,r)}}b_\un \uy^\un\in \cG'^{(t,r),>i-1}_{p^\infty}(\nu)\otimes_{R^\intern_{p^\infty}}R^\intern_\infty,
\end{equation}
où pour tout $\un=(n_1,\dots,n_{d'})\in J_{i-1}$, 
\begin{equation}\label{mtht32h}
b_\un=(\nu(\gamma_i)-1)a_\un+\sum_{\um=(m_1,\dots,m_{d'})
\in J_{i-1}(\un)} p^{t(m_i-n_i)}\binom{m_i}{n_i} \nu(\gamma_i)a_{\um}w^{m_i-n_i}.
\end{equation}

Soit $z$ un élément de $\cG'^{(t,r),>i-1}_{p^\infty}(\nu)\hotimes_{R^\intern_{p^\infty}}R^\intern_\infty$. Alors, $z$ peut s'écrire comme somme d'une série 
\begin{equation}
z=\sum_{\un\in J_{i-1}}p^{|\un|^{(t,r)}}a_\un \uy^\un,
\end{equation}
où $a_\un \in R'^{(\nu)}_{p^\infty}\hotimes_{R^\intern_{p^\infty}}R^\intern_\infty$ et $a_\un$ tend vers $0$ quand $|\un|$ tend vers l'infini. 
Comme $R'^{(\nu)}_{p^\infty}\otimes_{R^\intern_{p^\infty}}R^\intern_\infty$ 
est facteur direct de $R'_{p^\infty}\otimes_{R^\intern_{p^\infty}}R^\intern_\infty$ \eqref{eccr16b} et que 
$R'_{p^\infty}\otimes_{R^\intern_{p^\infty}}R^\intern_\infty$ est un sous-anneau de $R'_\infty$ \eqref{eccr9d},
on voit que $R'^{(\nu)}_{p^\infty}\otimes_{R^\intern_{p^\infty}}R^\intern_\infty$ est $\co_C$-plat.
Il en est donc de même de $R'^{(\nu)}_{p^\infty}\hotimes_{R^\intern_{p^\infty}}R^\intern_\infty$ (cf. la preuve de \cite{agt} II.6.14). 
Par suite, pour que $z$ soit nul, il faut et il suffit que $a_\un$ soit nul pour tout $\un\in J_{i-1}$.
On voit aussitôt que $(\gamma_i-\id)(z)$ est encore donné par la formule \eqref{mtht32g}.

Supposons que $\gamma_i(z)=z$ et $\nu(\gamma_i)\not=1$.  
Comme $R'^{(\nu)}_{p^\infty}\hotimes_{R^\intern_{p^\infty}}R^\intern_\infty$ est $\co_C$-plat et que $v(\nu(\gamma_i)-1)\leq \frac{1}{p-1}$, 
pour tout  $\un=(n_1,\dots,n_{d'})\in J_{i-1}$, on a 
\[
a_\un=-(\nu(\gamma_i)-1)^{-1}\sum_{\um=(m_1,\dots,m_{d'})
\in J_{i-1}(\un)} p^{t(m_i-n_i)}\binom{m_i}{n_i} \nu(\gamma_i)a_{\um}w^{m_i-n_i}.
\]
On en déduit que pour tout $\alpha\in \mN$ et tout $\un\in J_{i-1}$,
on a $a_\un\in p^{t\alpha} R'^{(\nu)}_{p^\infty}\hotimes_{R^\intern_{p^\infty}}R^\intern_\infty$ (on le prouve par récurrence sur $\alpha$);
donc $z=0$. Par suite, $\gamma_i-\id$ est injectif sur $\cG'^{(t,r),>i-1}_{p^\infty}(\nu)\hotimes_{R^\intern_{p^\infty}}R^\intern_\infty$. 

Supposons que $\gamma_i(z)=z$ et $\nu(\gamma_i)=1$, de sorte que $\nu\in \upchi^{>i}_{p^\infty}$. 
Alors, pour tout $\un=(n_1,\dots,n_{d'})\in J_{i-1}$, si on pose 
$\un'=(n'_1,\dots,n'_{d'})\in J_{i-1}(\un)$ avec $n'_i=n_i+1$, on a 
\[
(n_i+1)!a_{\un'}=-\sum_{\um=(m_1,\dots,m_{d'})
\in J_{i-1}(\un')} p^{t(m_i-n_i-1)}m_i! a_{\um}\frac{w^{m_i-n_i-1}}{(m_i-n_i)!}.
\]
On a $w^{m-1}/m!\in \co_C$ pour tout entier $m\geq 1$. On en déduit que 
pour tout $\alpha\in \mN$ et tout $\un=(n_1,\dots,n_{d'})\in J_{i-1}$ tel que $n_i\geq 1$, 
on a $n_i!a_{\un}\in p^{t\alpha} R'^{(\nu)}_{p^\infty}\hotimes_{R^\intern_{p^\infty}}R^\intern_\infty$ (on le prouve par récurrence sur $\alpha$); donc $a_\un=0$. 
Par suite $z \in \cG'^{(t,r),>i}_{p^\infty}(\nu)\hotimes_{R^\intern_{p^\infty}}R^\intern_\infty$. On en déduit que la suite 
\begin{equation}
\xymatrix{
0\ar[r]&{\cG'^{(t,r),>i}_{p^\infty}(\nu)\hotimes_{R^\intern_{p^\infty}}R^\intern_\infty}\ar[r]^-(0.45){u_\nu}&
{\cG'^{(t,r),>i-1}_{p^\infty}(\nu)\hotimes_{R^\intern_{p^\infty}}R^\intern_\infty}\ar[r]^{\gamma_i-\id}&{\cG'^{(t,r),>i-1}_{p^\infty}(\nu)\hotimes_{R^\intern_{p^\infty}}R^\intern_\infty}},
\end{equation}
où $u_\nu$ est le morphisme canonique, est exacte.
 
Il résulte de ce qui précède que la suite \eqref{mtht32a} est exacte.

(iii) En vertu de \ref{mtht31}(iii), il existe un entier $\alpha\geq 0$, dépendant de $t$ et $t'$ 
mais pas des données \ref{hypmdef2}, tel que l'on ait
\begin{equation}
p^\alpha(\hcG'^{(t',r),>i-1}_{p^\infty}\otimes_{R^\intern_{p^\infty}}R^\intern_\infty)
\subset (\gamma_i-\id)(\hcG'^{(t,r),>i-1}_{p^\infty}\otimes_{R^\intern_{p^\infty}}R^\intern_\infty).
\end{equation}
Posant $M=(\gamma_i-\id)(\hcG'^{(t,r),>i-1}_{p^\infty}\otimes_{R^\intern_{p^\infty}}R^\intern_\infty)$, on en déduit 
par complétion $p$-adique un diagramme commutatif 
\begin{equation}\label{mtht32d}
\xymatrix{
{p^\alpha(\cG'^{(t',r),>i-1}_{p^\infty}\hotimes_{R^\intern_{p^\infty}}R^\intern_\infty)}\ar[r]&\hM\ar[r]&
{\cG'^{(t,r),>i-1}_{p^\infty}\hotimes_{R^\intern_{p^\infty}}R^\intern_\infty}\\
&{\cG'^{(t,r),>i-1}_{p^\infty}\hotimes_{R^\intern_{p^\infty}}R^\intern_\infty}\ar[ru]_{\gamma_i-\id}\ar@{->>}[u]&}
\end{equation}
où la flèche verticale est surjective en vertu de (\cite{egr1} 1.8.5), d'où la proposition recherchée.

\subsection{}\label{mtht35}
Conservons les hypothèses de \ref{mtht33}. Pour tous nombres rationnels $t\geq r \geq 0$, 
on définit par récurrence, pour tout entier $0\leq i\leq \ell$, un complexe $\mK_i^{(t,r),\bullet}$ de $\hRinterni$-représentations
continues de $\fS_{p^\infty}$ en posant $\mK_0^{(t,r),\bullet}=\cG'^{(t,r)}_{p^\infty}\hotimes_{R^\intern_{p^\infty}}R^\intern_\infty[0]$ 
et pour tout $1\leq i\leq \ell$, $\mK_i^{(t,r),\bullet}$ est la fibre du morphisme 
\begin{equation}
\gamma_i-\id\colon \mK_{i-1}^{(t,r),\bullet}\rightarrow \mK_{i-1}^{(t,r),\bullet}.
\end{equation} 
Il résulte de (\cite{agt} II.3.25 et (II.2.7.9)) qu'on a un isomorphisme canonique dans $\bD^+(\bMod(\hRinterni))$
\begin{equation}\label{mtht35a}
\rC_{\cont}^\bullet(\fS_{p^\infty},\cG'^{(t,r)}_{p^\infty}\hotimes_{R^\intern_{p^\infty}}R^\intern_\infty)\stackrel{\sim}{\rightarrow}\mK_\ell^{(t,r),\bullet}.
\end{equation}
Pour tous nombres rationnels $t'\geq t \geq r\geq 0$, l'homomorphisme canonique $\cG'^{(t',r)}_{p^\infty}\rightarrow \cG'^{(t,r)}_{p^\infty}$
induit pour tout entier $0\leq i\leq \ell$ un morphisme $\mK_i^{(t',r),\bullet}\rightarrow \mK_i^{(t,r),\bullet}$
de complexes de $\hRinterni$-représentations continues de $\fS_{p^\infty}$.

\begin{prop}\label{mtht36}
Les hypothèses étant celles de \ref{mtht33}, soient, de plus, $i$ un entier tel que $0\leq i\leq \ell=d'-d$,
$t,t',r$ trois nombres rationnels tels que $t'>t>r\geq 0$. Alors,
\begin{itemize}
\item[{\rm (i)}] On a un isomorphisme canonique $\hRinterni$-linéaire et $\fS_{p^\infty}$-équivariant 
\begin{equation}\label{mtht36a}
\cG'^{(t,r),>i}_{p^\infty}\hotimes_{R^\intern_{p^\infty}}R^\intern_\infty\stackrel{\sim}{\rightarrow}\rH^0(\mK_i^{(t,r),\bullet}).
\end{equation}
\item[{\rm (ii)}] Il existe un entier $\alpha_i\geq 0$, dépendant de $t$, $t'$ et $i$,
mais pas des morphismes $f$, $f'$ et $g$ vérifiant les conditions de \ref{hypmdef2} et \ref{hypmdef3}, tel que 
pour tout entier $j\geq 1$, le morphisme canonique 
\begin{equation}\label{mtht36b}
\rH^j(\mK_i^{(t',r),\bullet}) \rightarrow\rH^j(\mK_i^{(t,r),\bullet})
\end{equation}
soit annulé par $p^{\alpha_i}$. 
\end{itemize}
\end{prop}

Procédons par récurrence sur $i$. La proposition est immédiate pour $i=0$ d'après \ref{mtht32}(i). 
Supposons $i\geq 1$ et la proposition démontrée pour $i-1$. Le triangle distingué 
\begin{equation}\label{mtht36c}
\xymatrix{
{\mK_i^{(t,r),\bullet}}\ar[r]&{\mK_{i-1}^{(t,r),\bullet}}\ar[r]^-(0.4){\gamma_i-\id}& 
{\mK_{i-1}^{(t,r),\bullet}}\ar[r]^-(0.5){+1}&}
\end{equation}
et l'hypothèse de récurrence induisent une suite exacte
\begin{equation}\label{mtht36d}
\xymatrix{
0\ar[r]&{\rH^0(\mK_i^{(t,r),\bullet})}\ar[r]&{\cG'^{(t,r),>i-1}_{p^\infty}\hotimes_{R^\intern_{p^\infty}}R^\intern_\infty}\ar[r]^-(0.4){\gamma_i-\id}& 
{\cG'^{(t,r),>i-1}_{p^\infty}\hotimes_{R^\intern_{p^\infty}}R^\intern_\infty}},
\end{equation}
qui implique proposition (i) en vertu de \ref{mtht32}(ii).

Pour tout entier $j\geq 1$, le triangle distingué \eqref{mtht36c} induit une suite exacte de $\hRinterni$-modules
\begin{equation}\label{mtht36e}
0\rightarrow C_j^{(t,r)}\rightarrow \rH^j(\mK_i^{(t,r),\bullet}) \rightarrow D_j^{(t,r)} \rightarrow 0,
\end{equation}
où $C_j^{(t,r)}$ est un quotient de $\rH^{j-1}(\mK_{i-1}^{(t,r),\bullet})$ et 
$D_j^{(t,r)}$ est un sous-module de $\rH^j(\mK_{i-1}^{(t,r),\bullet})$. De plus, d'après l'hypothèse de récurrence, 
on a un isomorphisme canonique 
\begin{equation}\label{mtht36f}
C_1^{(t,r)}\stackrel{\sim}{\rightarrow}(\cG'^{(t,r),>i-1}_{p^\infty}\hotimes_{R^\intern_{p^\infty}}R^\intern_\infty)/(\gamma_i-1)(\cG'^{(t,r),>i-1}_{p^\infty}\hotimes_{R^\intern_{p^\infty}}R^\intern_\infty).
\end{equation}
Les morphismes canoniques $\mK_{i-1}^{(t',r),\bullet}\rightarrow \mK_{i-1}^{(t,r),\bullet}$
et $\mK_i^{(t',r),\bullet}\rightarrow \mK_i^{(t,r),\bullet}$ induisent des morphismes $C_j^{(t',r)}\rightarrow C_j^{(t,r)}$ et 
$D_j^{(t',r)}\rightarrow D_j^{(t,r)}$ qui s'insèrent dans un diagramme commutatif 
\begin{equation}\label{mtht36g}
\xymatrix{
0\ar[r]&{C_j^{(t',r)}}\ar[r]\ar[d]&{\rH^j(\mK_i^{(t',r),\bullet})}\ar[r]\ar[d]&{D_j^{(t',r)}}\ar[r]\ar[d]&0\\
0\ar[r]&{C_j^{(t,r)}}\ar[r]&{\rH^j(\mK_i^{(t,r),\bullet})}\ar[r]&{D_j^{(t,r)}}\ar[r]&0}
\end{equation}
où la flèche verticale au centre est le morphisme canonique. 

Posons $t''=(t+t')/2$. D'après l'hypothèse de récurrence, il existe un entier 
$\alpha_{i-1}\geq 0$, dépendant seulement de $t$, $t'$ et $i-1$, 
tel que pour tout entier $j\geq 1$, le morphisme $D_j^{(t',r)}\rightarrow D_j^{(t'',r)}$ soit annulé par $p^{\alpha_{i-1}}$. 
D'autre part, compte tenu de l'hypothèse de récurrence et en vertu de \eqref{mtht36f} et \ref{mtht32}(iii), il existe un 
entier $\alpha'_{i-1}\geq 0$, dépendant seulement de $t$, $t'$ et $i-1$, tel que pour tout entier $j\geq 1$, 
le morphisme $C_j^{(t'',r)}\rightarrow C_j^{(t,r)}$ soit annulé par $p^{\alpha'_{i-1}}$. La proposition (ii) s'ensuit en prenant 
$\alpha_i=\alpha_{i-1}+\alpha'_{i-1}$. 

\begin{cor}\label{mtht37}
Soient $t,t',r$ trois nombres rationnels tels que $t'>t>r\geq 0$. Alors,
\begin{itemize}
\item[{\rm (i)}] L'homomorphisme canonique 
\begin{equation}\label{mtht37a}
\cG^{(r)}_{p^\infty}\hotimes_{R_{p^\infty}}R^\intern_\infty \rightarrow (\cG'^{(t,r)}_{p^\infty}\hotimes_{R^\intern_{p^\infty}}R^\intern_\infty)^{\fS_{p^\infty}}
\end{equation}
est un isomorphisme.
\item[{\rm (ii)}] Il existe un entier $\alpha\geq 0$, dépendant de $t$, $t'$ et $\ell=d'-d$, 
mais pas des morphismes $f$, $f'$ et $g$ vérifiant les conditions de \ref{hypmdef2} et \ref{hypmdef3}, tel que 
pour tout entier $j\geq 1$, le morphisme canonique 
\begin{equation}\label{mtht37b}
\rH^j_\cont(\fS_{p^\infty},\cG'^{(t',r)}_{p^\infty}\hotimes_{R^\intern_{p^\infty}}R^\intern_\infty)\rightarrow 
\rH^j_\cont(\fS_{p^\infty},\cG'^{(t,r)}_{p^\infty}\hotimes_{R^\intern_{p^\infty}}R^\intern_\infty)
\end{equation}
soit annulé par $p^\alpha$. 
\end{itemize}
\end{cor}

Cela résulte de \eqref{mtht35a} et \ref{mtht36}.

\begin{cor}\label{mtht38}
Soient $t,t',r$ trois nombres rationnels tels que $t'>t>r\geq 0$, $\fS$ le groupe défini dans \eqref{eccr43m}. Alors,
\begin{itemize}
\item[{\rm (i)}] L'homomorphisme canonique 
\begin{equation}\label{mtht38a}
\cG^{(r)}_{p^\infty}\hotimes_{R_{p^\infty}}R^\intern_\infty  \rightarrow (\hcG'^{(t,r)})^\fS
\end{equation}
est un $\alpha$-isomorphisme.
\item[{\rm (ii)}] Il existe un entier $\alpha\geq 0$, dépendant de $t$, $t'$ et $\ell=d'-d$, 
mais pas des données \ref{hypmdef2}, tel que 
pour tout entier $j\geq 1$, le morphisme canonique 
\begin{equation}\label{mtht38b}
\rH^j_\cont(\fS,\hcG'^{(t',r)})\rightarrow \rH^j_\cont(\fS,\hcG'^{(t,r)})
\end{equation}
soit annulé par $p^\alpha$. 
\end{itemize}
\end{cor}

Cela résulte de \ref{mtht29} et \ref{mtht37} (cf. \ref{eccr34}).

\begin{prop}\label{mtht40}
Soient $t,t',r$ trois nombres rationnels tels que $t'>t>r\geq 0$. Alors,
\begin{itemize}
\item[{\rm (i)}] Pour tout entier $n\geq 1$, l'homomorphisme canonique 
\begin{equation}\label{mtht40a}
(\cG^{(r)}_{p^\infty}/p^n\cG^{(r)}_{p^\infty})\otimes_{R_{p^\infty}}R^\intern_{\infty}\rightarrow (\cG'^{(t,r)}/p^n\cG'^{(r)})^{\fS}
\end{equation}
est $\alpha$-injectif. Notons $\cK^{(r)}_n$ son conoyau. 
\item[{\rm (ii)}] Il existe un entier $\alpha\geq 0$, dépendant de $t$, $t'$ et $\ell=d'-d$, 
mais pas des morphismes $f$, $f'$ et $g$ vérifiant les conditions de \ref{hypmdef2} et \ref{hypmdef3}, tel que pour tout entier $n\geq 1$, 
le morphisme canonique $\cK^{(t',r)}_n\rightarrow \cK^{(t,r)}_n$ soit annulé par $p^\alpha$. 
\item[{\rm (iii)}] Il existe un entier $\gamma\geq 0$, dépendant de $t$, $t'$ et $\ell$, 
mais pas des morphismes $f$, $f'$ et $g$ vérifiant les conditions de \ref{hypmdef2} et \ref{hypmdef3}, tel que pour tous entiers $n,q\geq 1$, 
le morphisme canonique
\begin{equation}\label{mtht40b}
\rH^q(\fS,\cG'^{(t',r)}/p^n\cG'^{(t',r)})\rightarrow \rH^q(\fS,\cG'^{(t,r)}/p^n\cG'^{(t,r)})
\end{equation}
soit annulé par $p^\gamma$. 
\end{itemize}
\end{prop}

(i) Cela résulte de \ref{mtht38}(i) et de 
la suite exacte longue de cohomologie associée à la suite exacte courte de $\mZ_p$-représentations de $\fS$ 
\begin{equation}\label{mtht40c}
0\longrightarrow \hcG'^{(t,r)}\stackrel{\cdot p^n}{\longrightarrow} \hcG'^{(t,r)}\longrightarrow \cG'^{(r)}/p^n\cG'^{(t,r)}\longrightarrow 0.
\end{equation}
On en déduit aussi un morphisme $\hRinterni$-linéaire $\alpha$-injectif
\begin{equation}\label{mtht40d}
\cK^{(t,r)}_n\rightarrow \rH^1_{\cont}(\fS,\hcG'^{(t,r)}).
\end{equation}

(ii) Cela résulte de \eqref{mtht40d} et \ref{mtht38}(ii).

(iii) Pour tous entiers $n,q\geq 1$, la suite exacte longue de cohomologie déduite de \eqref{mtht40c} 
fournit une suite exacte de $\hRun$-modules
\begin{equation}
0\rightarrow \rH^q_{\cont}(\fS,\hcG'^{(t,r)})/p^n\rH^q_{\cont}(\fS,\hcG'^{(t,r)})
\rightarrow \rH^q(\fS,\cG'^{(t,r)}/p^n\cG'^{(t,r)}) \rightarrow T_n^{(t,r),q}\rightarrow 0,
\end{equation}
où $T_n^{(t,r),q}$ est un sous-module de $p^n$-torsion de $\rH^{q+1}_{\cont}(\fS,\hcG'^{(t,r)})$.
Posons $t''=(t+t')/2$. D'après \ref{mtht38}(ii), il existe un entier $\beta'>0$, dépendant 
seulement de $t$, $t'$ et $\ell$, tel que
pour tout entier $q\geq 1$, les morphismes canoniques 
\[
\rH^q_{\cont}(\fS,\hcG'^{(t',r)})\rightarrow \rH^q_{\cont}(\fS,\hcG'^{(t'',r)}) \ \ \ {\rm et} \ \ \ 
\rH^q_{\cont}(\fS,\hcG'^{(t'',r)})\rightarrow \rH^q_{\cont}(\fS,\hcG'^{(t,r)})
\] 
soient annulés par $p^{\beta'}$. La proposition s'ensuit en prenant $\beta=2\beta'$.

\begin{cor}\label{mtht41}
Soient $t,t',r$ trois nombres rationnels tels que $t'>t>r\geq 0$, $\Pi$ le groupe défini dans \eqref{eccr40c}. Alors,
\begin{itemize}
\item[{\rm (i)}] Pour tout entier $n\geq 1$, l'homomorphisme canonique 
\begin{equation}\label{mtht41a}
(\cG^{(r)}/p^n\cG^{(r)})\otimes_{\oR}\oR^\intern\rightarrow (\cG'^{(t,r)}/p^n\cG'^{(t,r)})^{\Pi}
\end{equation}
est $\alpha$-injectif. Notons $\cH^{(t,r)}_n$ son conoyau. 
\item[{\rm (ii)}] Il existe un entier $\alpha\geq 0$, dépendant de $t$, $t'$ et $\ell=d'-d$, 
mais pas des morphismes $f$, $f'$ et $g$ vérifiant les conditions de \ref{hypmdef2} et \ref{hypmdef3}, tel que pour tout entier $n\geq 1$, 
le morphisme canonique $\cH^{(t',r)}_n\rightarrow \cH^{(t,r)}_n$ soit annulé par $p^\alpha$. 
\item[{\rm (iii)}] Il existe un entier $\gamma\geq 0$, dépendant de $t$, $t'$ et $\ell$, 
mais pas des morphismes $f$, $f'$ et $g$ vérifiant les conditions de \ref{hypmdef2} et \ref{hypmdef3}, tel que pour tous entiers $n,q\geq 1$, 
le morphisme canonique
\begin{equation}\label{mtht41b}
\rH^q(\Pi,\cG'^{(t',r)}/p^n\cG'^{(t',r)})\rightarrow \rH^q(\Pi,\cG'^{(t,r)}/p^n\cG'^{(t,r)})
\end{equation}
soit annulé par $p^\gamma$. 
\end{itemize}
\end{cor}

Nous utiliserons les notations de \ref{eccr34}. Il résulte de (\cite{ag} 5.2.17) et de la suite spectrale 
\begin{equation}
\rE_2^{ij}=\rH^i(\Sigma^\intern,\rH^j(\Pi,\hcG'^{(t,r)}/p^n\hcG'^{(t,r)}))\Rightarrow \rH^{i+j}(\fS,\hcG'^{(t,r)}/p^n\hcG'^{(t,r)})
\end{equation}
que pour tout entier $i\geq 0$, le morphisme canonique
\begin{equation}
\rH^i(\fS,\cG'^{(t,r)}/p^n\cG'^{(t,r)})\rightarrow \rH^i(\Pi,\cG'^{(t,r)}/p^n\cG'^{(t,r)})^{\Sigma^\intern}
\end{equation}
est un $\alpha$-isomorphisme, et que le morphisme $\oR^\intern$-linéaire canonique 
\begin{equation}
\rH^i(\Pi,\cG'^{(t,r)}/p^n\cG'^{(t,r)})^{\Sigma^\intern}\otimes_{R^\intern_\infty}\oR^\intern\rightarrow \rH^i(\Pi,\cG'^{(t,r)}/p^n\cG'^{(t,r)})
\end{equation}
est un $\alpha$-isomorphisme. On notera que les $\oR^\intern$-représentations $\rH^j(\Pi,\hcG'^{(t,r)}/p^n\hcG'^{(t,r)})$ de $\Sigma^\intern$ 
sont discrètes (\cite{ribzal} § 7.2 page 257). Par suite, pour tout entier $i\geq 0$, le morphisme $\oR^\intern$-linéaire canonique 
\begin{equation}
\rH^i(\fS,\cG'^{(t,r)}/p^n\cG'^{(t,r)})\otimes_{R^\intern_\infty}\oR^\intern\rightarrow \rH^i(\Pi,\cG'^{(t,r)}/p^n\cG'^{(t,r)})
\end{equation}
est un $\alpha$-isomorphisme. Par ailleurs, la $R^\intern_\infty$-algèbre $\oR^\intern$ est $\alpha$-plate d'après \ref{eccr99}. 
La proposition résulte alors de \ref{mtht40}.

\subsection{}\label{mtht42}
On peut maintenant démontrer la proposition \ref{mtht20}. Compte tenu de \ref{mtht100}, \ref{mtht23} et \ref{mtht25}, on peut se réduire 
au cas où $\tg$ \eqref{hypmdef3a} est le morphisme $\tg_0$ considéré dans \ref{mtht26}. 
Les sections $(\psi_0,\psi'_0)\in \cL_0(\hmX)\times \cL'_0(\hmX')$ sont alors compatibles \eqref{mtht7}. 
Elles définissent un isomorphisme $\Delta$-équivariant de $\hoR$-algèbres 
$\cG\stackrel{\sim}{\rightarrow} \cC$ \eqref{mtht4b} et un isomorphisme $\Delta'$-équivariant de $\hoRp$-algèbres 
$\cG'\stackrel{\sim}{\rightarrow} \cC'$ \eqref{mtht4e} qui s'insèrent dans un diagramme commutatif
\begin{equation}
\xymatrix{
\cG\ar[r]\ar[d]&\cC\ar[d]\\
\cG'\ar[r]&\cC'}
\end{equation}
où les flèches verticales sont les homomorphismes canoniques \eqref{mtht7f}. 
D'après (\cite{agt} II.12.1), elles définissent aussi un isomorphisme $\Delta$-équivariant  de $\hoR$-algèbres $\cG^{(r)}\stackrel{\sim}{\rightarrow} \cC^{(r)}$ 
\eqref{mtht22a} et un isomorphisme $\Delta'$-équivariant de $\hoRp$-algèbres $\cG'^{(r)}\stackrel{\sim}{\rightarrow} \cC'^{(r)}$ \eqref{mtht24a} 
qui s'insèrent dans un diagramme commutatif
\begin{equation}
\xymatrix{
\cG^{(r)}\ar[r]\ar[d]&\cC^{(r)}\ar[d]\\
\cG'^{(r)}\ar[r]&\cC'^{(r)}}
\end{equation} 
où les flèches verticales sont les homomorphismes canoniques \eqref{mtht11d}. 
On en déduit un isomorphisme $\Delta'$-équivariant de $\hoRp$-algèbres $\cG'^{(t,r)}\stackrel{\sim}{\rightarrow}  \cC'^{(t,r)}$ \eqref{mtht24b}
compatible avec l'isomorphisme $\cG^{(r)}\stackrel{\sim}{\rightarrow} \cC^{(r)}$. La proposition \ref{mtht20} résulte alors de \ref{mtht41}.

\section{Cohomologie de Dolbeault relative des algèbres de Higgs-Tate}\label{cdlbr}

\subsection{}\label{cdlbr1}
Les hypothèses et notations de \ref{hypmdef} et \ref{mtht} sont en vigueur dans cette section. 
Pour tout nombre rationnel $r\geq 0$, on désigne par  
\begin{equation}\label{cdlbr1a}
d_{\cC'^{(r)}}\colon \cC'^{(r)}\rightarrow \txi^{-1}\tOmega^1_{R'/\co_K}\otimes_{R'}\cC'^{(r)}
\end{equation}
la $\hoRp$-dérivation universelle de l'algèbre $\cC'^{(r)}$ définie dans \eqref{mtht11c} (cf. \ref{taht11}). 
Le morphisme
\begin{equation}\label{cdlbr1b}
\ud_{\cC'^{(r)}}\colon \cC'^{(r)}\rightarrow \txi^{-1}\tOmega^1_{R'/R}\otimes_{R'}\cC'^{(r)}
\end{equation}
induit par $d_{\cC'^{(r)}}$ s'identifie à la $(\cC^{(r)}\otimes_{\hoR}\hoRp)$-dérivation universelle de $\cC'^{(r)}$. 
Cela résulte de \ref{mtht9} lorsque $r=0$ et la preuve dans le cas général est similaire. 
Comme $\txi^{-1}\tOmega^1_{R'/R}\otimes_{R'}\hoRp=\ud_{\cC'^{(r)}}(\cF'^{(r)})\subset \ud_{\cC'^{(r)}}(\cC'^{(r)})$,
$\ud_{\cC'^{(r)}}$ est un $(\cC^{(r)}\otimes_{\hoR}\hoRp)$-champ de Higgs à coefficients dans $\txi^{-1}\tOmega^1_{R'/R}$ d'après \ref{MH8}(i).

Pour tous nombres rationnels $r'\geq r\geq 0$, on a 
\begin{equation}\label{cdlbr1c}
p^{r'}(\id \times \alpha'^{r,r'}) \circ d_{\cC'^{(r')}}=p^rd_{\cC'^{(r)}}\circ \alpha'^{r,r'},
\end{equation}
où $\alpha'^{r,r'}\colon \cC'^{(r')}\rightarrow \cC'^{(r)}$ est l'homomorphisme canonique. 

Pour tous nombres rationnels $t\geq r\geq 0$, le morphisme 
\begin{equation}\label{cdlbr1d}
\ud_{\cC'^{(t,r)}}\colon \cC'^{(t,r)}\rightarrow \txi^{-1}\tOmega^1_{R'/R}\otimes_{R'}\cC'^{(t,r)}
\end{equation}
induit par $\ud_{\cC'^{(t)}}$ s'identifie à la $(\cC^{(r)}\otimes_{\hoR}\hoRp)$-dérivation universelle de l'algèbre $\cC'^{(t,r)}$ définie dans \eqref{mtht11e}.
On désigne par 
\begin{equation}\label{cdlbr1g}
\ud_{\hcC'^{(t,r)}}\colon \hcC'^{(t,r)}\rightarrow \txi^{-1}\tOmega^1_{R'/R}\otimes_{R'}\hcC'^{(t,r)}
\end{equation}
son prolongement aux complétés (on notera que le $R'$-module $\tOmega^1_{R'/R}$ est libre de type fini). 
Les dérivations $\ud_{\cC'^{(t,r)}}$ et $\ud_{\hcC'^{(t,r)}}$ sont des $(\cC^{(r)}\otimes_{\hoR}\hoRp)$-champs de Higgs à coefficients dans 
$\txi^{-1}\tOmega^1_{R'/R}$. On désigne par $\umK^\bullet(\hcC'^{(t,r)})$ le complexe de Dolbeault  de $(\hcC'^{(t,r)},p^t\ud_{\hcC'^{(t,r)}})$
et par $\utmK^\bullet(\hcC'^{(t,r)})$ le complexe de Dolbeault augmenté 
\begin{equation}\label{cdlbr1e}
\cC^{(r)}\hotimes_{\oR}\oR'\rightarrow \umK^0(\hcC'^{(t,r)})\rightarrow \umK^1(\hcC'^{(t,r)})\rightarrow \dots
\rightarrow \umK^n(\hcC'^{(t,r)})\rightarrow \dots,
\end{equation}
où $\cC^{(r)}\hotimes_{\oR}\oR'$ est placé en degré $-1$, le produit tensoriel $\hotimes$ est complété pour la topologie $p$-adique 
et la différentielle $\cC^{(r)}\hotimes_{\oR}\oR'\rightarrow\hcC'^{(t,r)}$ est l'homomorphisme canonique. 
  
Compte tenu de \eqref{cdlbr1c}, pour tous nombres rationnels $t'\geq t\geq r\geq 0$, 
le morphisme $\halpha'^{t,t'}$ induit un morphisme de complexes 
\begin{equation}\label{cdlbr1f}
\upiota^{t',t,r}\colon \utmK^\bullet(\hcC'^{(t',r)})\rightarrow \utmK^\bullet(\hcC'^{(t,r)}).
\end{equation}

\begin{prop}\label{cdlbr2}
Soient $t',t,r$ trois nombres rationnels tels que $t'>t>r\geq 0$. Alors,  
\begin{itemize}
\item[{\rm (i)}] Il existe un nombre rationnel $\alpha\geq 0$ dépendant de $t$ et $t'$ mais pas des morphismes $f$, $f'$ et $g$
vérifiant les conditions de \ref{hypmdef2} et \ref{hypmdef3}, tel que
\begin{equation}\label{cdlbr2a}
p^\alpha\upiota^{t',t,r}\colon \utmK^\bullet(\hcC'^{(t',r)})\rightarrow \utmK^\bullet(\hcC'^{(t,r)}),
\end{equation}
où $\upiota^{t',t,r}$ est le morphisme \eqref{cdlbr1f}, soit homotope à $0$ par une homotopie $\hoRp$-linéaire. 
\item[{\rm (ii)}] Le morphisme canonique
\begin{equation}\label{cdlbr2b}
\upiota^{t',t,r}\otimes_{\mZ_p}\mQ_p\colon \utmK^\bullet(\hcC'^{(t',r)})\otimes_{\mZ_p}\mQ_p\rightarrow 
\utmK^\bullet(\hcC'^{(t,r)})\otimes_{\mZ_p}\mQ_p
\end{equation}
est homotope à $0$ par une homotopie continue. 
\end{itemize}
\end{prop}

Soit $(\sigma,\sigma')\in \cL(\hmX)\times \cL'(\hmX')$ un couple de sections compatibles \eqref{mtht7}. 
D'après (\cite{agt} II.12.1), il définit un isomorphisme de $\hoR$-algèbres $\cG^{(r)}\stackrel{\sim}{\rightarrow} \cC^{(r)}$ \eqref{mtht22a} et  
un isomorphisme de $\hoRp$-algèbres $\cG'^{(r)}\stackrel{\sim}{\rightarrow} \cC'^{(r)}$ \eqref{mtht24a} qui s'insèrent dans un diagramme commutatif
\begin{equation}
\xymatrix{
\cG^{(r)}\ar[r]\ar[d]&\cC^{(r)}\ar[d]\\
\cG'^{(r)}\ar[r]&\cC'^{(r)}}
\end{equation} 
où les flèches verticales sont les homomorphismes canoniques \eqref{mtht11d}. 
On en déduit un isomorphisme de $\hoRp$-algèbres $\cG'^{(t,r)}\stackrel{\sim}{\rightarrow}  \cC'^{(t,r)}$ \eqref{mtht24b}
compatible avec l'isomorphisme $\cG^{(r)}\stackrel{\sim}{\rightarrow} \cC^{(r)}$. 
On identifie les $\hoR$-algèbres $\cG^{(r)}$ et $\cC^{(r)}$ 
et les $\hoR$-algèbres $\cG'^{(t,r)}$ et $\cC'^{(t,r)}$ par ces isomorphismes, et on adapte les notations de \ref{cdlbr1} en conséquence. 

Considérons la suite exacte canonique \eqref{hypmdef2d}
\begin{equation}
0\longrightarrow P^\gp/\lambda\mZ\stackrel{h^\gp}{\longrightarrow} P'^\gp/h^\gp(\lambda)\mZ \longrightarrow P'^{\gp}/h^\gp(P^{\gp}) \longrightarrow 0.
\end{equation}
Les $\mZ$-modules $P^\gp/\mZ\lambda$ et $P'^\gp/h^\gp(\lambda)\mZ$ sont libres de type fini d'après (\cite{ag} 4.2.2); notons $d$ et $d'$ leurs rangs
respectifs et posons $\ell=d'-d$. Soient $(t_i)_{1\leq i\leq \ell}$ des éléments de $P'^\gp$, $(t_i)_{\ell+1\leq i\leq d'}$ des éléments de $P^\gp$, 
tels que les images des $(t_i)_{1\leq i\leq \ell}$ dans $P'^\gp/h^\gp(P^\gp)\otimes_\mZ\mZ_p$ forment une $\mZ_p$-base,
que les images des $(t_i)_{\ell+1\leq i\leq d'}$ dans $P^\gp/\lambda\mZ$ forment une $\mZ$-base,
et que les images des $t_{1},\dots,t_{\ell},h^\gp(t_{\ell+1}),\dots,h^\gp(t_{d'})$ dans $P'^\gp/h^\gp(\lambda)\mZ$ forment une $\mZ$-base.
Pour alléger les notations, on notera abusivement $(h^\gp(t_i))_{\ell+1\leq i\leq d'}$ simplement par $(t_i)_{\ell+1\leq i\leq d'}$, ce qui n'induit 
aucune ambiguïté puisque $h^\gp$ est injectif. Les $(d\log(t_i))_{1\leq i\leq d'}$ forment une $R'$-base de $\tOmega^1_{R'/\co_K}$ \eqref{mtht30d2}. 
Pour tout $1\leq i\leq d'$, posons $y_i=\txi^{-1}d\log(t_i)\in \txi^{-1}\tOmega^1_{R'/\co_K}\subset \cG'$.
Pour tout $1\leq i\leq \ell$, on note $J_i$ le sous-ensemble de $\mN^{d'}$ formé des éléments $\un=(n_1,\dots,n_{d'})$
tels que $n_1=\dots=n_{i}=0$.

On désigne par
\begin{equation}
\uptau^{-1}\colon \hcG'^{(t',r)}\otimes_{\mZ_p}\mQ_p\rightarrow (\cG^{(r)}\hotimes_{\oR}\oR') \otimes_{\mZ_p}\mQ_p
\end{equation}
le morphisme $\hoRp$-linéaire défini par 
\[
\uptau^{-1}(\sum_{\un=(n_1,\dots,n_{d'})\in \mN^{d'}}a_\un \prod_{1\leq i\leq d'}y_i^{n_i})=
\sum_{\un=(n_1,\dots,n_{d'})\in J_\ell}a_{\un}\prod_{\ell+1\leq i\leq d'}y_{i}^{n_i}.
\]
Pour tout entier $m\geq 0$, il existe un et un unique morphisme $\hoRp$-linéaire
\begin{equation}
\uptau^{m}\colon \txi^{-m-1}\tOmega^{m+1}_{R'/R}\otimes_{R'}\hcG'^{(t',r)}\otimes_{\mZ_p}\mQ_p\rightarrow 
\txi^{-m}\tOmega^{m}_{R'/R}\otimes_{R'}\hcG'^{(t,r)}\otimes_{\mZ_p}\mQ_p
\end{equation}
tel que pour tout $1\leq i_1<\dots<i_{m+1}\leq \ell$, on ait 
\begin{eqnarray*}
\ \ \ \ \ \lefteqn{\uptau^{m}(\sum_{\un=(n_1,\dots,n_{d'})\in \mN^{d'}}a_{\un}\prod_{1\leq i\leq d'}y_i^{n_i}\otimes 
\txi^{-1}d\log(t_{i_1})\wedge \dots\wedge \txi^{-1}d\log(t_{i_{m+1}}))}\\
&=&\sum_{\un=(n_1,\dots,n_{d'})\in J_{i_1-1}}
\frac{a_{\un}}{n_{i_1}+1}\prod_{1\leq i\leq d'}y_i^{n_i+\delta_{ii_1}}\otimes 
\txi^{-1}d\log(t_{i_2})\wedge \dots\wedge \txi^{-1}d\log(t_{i_{m+1}}).\nonumber
\end{eqnarray*}

Soit $\alpha$ un nombre rationnel tel que 
\begin{equation}
\alpha\geq \sup_{x\in \mQ_{\geq 0}}(\log_p(x+1)+(x+1)t-xt'),
\end{equation}
où $\log_p$ est la fonction logarithme de base $p$. Pour tout entier $m\geq 0$, on a clairement
\begin{equation}
p^\alpha \uptau^m(\txi^{-m-1}\tOmega^{m+1}_{R'/R}\otimes_{R'}\hcG'^{(t',r)})\subset 
\txi^{-m}\tOmega^{m}_{R'/R}\otimes_{R'}\hcG'^{(t,r)}. 
\end{equation}
On vérifie aussitôt que les morphismes $(p^\alpha \uptau^m)_{m\geq -1}$ définissent une homotopie reliant $0$ au morphisme 
$p^\alpha \upiota^{t',t,r}$, d'où la proposition.

\begin{cor}\label{cdlbr3}
Pour tout nombre rationnel $r\geq 0$, le morphisme canonique de complexes
\begin{equation}
(\cC^{(r)}\hotimes_{\oR}\oR')\otimes_{\mZ_p}\mQ_p[0]\rightarrow \underset{\underset{t\in \mQ_{>r}}{\longrightarrow}}{\lim} 
\umK^\bullet(\hcC'^{(t,r)})\otimes_{\mZ_p}\mQ_p
\end{equation} 
est un quasi-isomorphisme.
\end{cor}

\chapter{Cohomologie relative des modules de Dolbeault}

\section{Hypothèses et notations; topos de Faltings relatif}\label{hmdf}

\subsection{}\label{hmdf1}
Dans ce chapitre, $K$ désigne un corps de valuation discrète complet de 
caractéristique $0$, à corps résiduel {\em algébriquement clos} $k$ de caractéristique $p>0$,  
$\co_K$ l'anneau de valuation de $K$, $\oK$ une clôture algébrique de $K$, $\co_\oK$ la clôture intégrale de $\co_K$ dans $\oK$,
$\fm_\oK$ l'idéal maximal de $\co_\oK$ et $G_K$ le groupe de Galois de $\oK$ sur $K$.
On note $\co_C$ le séparé complété $p$-adique de $\co_\oK$, $\fm_C$ son idéal maximal,
$C$ son corps des fractions et $v$ sa valuation, normalisée par $v(p)=1$. 
On désigne par $\mZ_p(1)$ le $\mZ[G_K]$-module
\begin{equation}\label{hmdf1a}
\mZ_p(1)=\underset{\underset{n\geq 0}{\longleftarrow}}{\lim}\ \mu_{p^n}(\co_{\oK}),
\end{equation}  
où $\mu_{p^n}(\co_{\oK})$ désigne le sous-groupe des racines $p^n$-ièmes de l'unité dans $\co_\oK$. 
Pour tout $\mZ_p[G_K]$-module $M$ et tout entier $n$, on pose $M(n)=M\otimes_{\mZ_p}\mZ_p(1)^{\otimes n}$. 

Pour tout groupe abélien $A$, on note $\hA$ son séparé complété $p$-adique.

On pose $S=\Spec(\co_K)$, $\oS=\Spec(\co_\oK)$ et $\coS=\Spec(\co_C)$. 
On note $s$ (resp.  $\eta$, resp. $\oeta$) le point fermé de $S$ (resp.  générique de $S$, resp. générique de $\oS$).
Pour tout entier $n\geq 1$, on pose $S_n=\Spec(\co_K/p^n\co_K)$. Pour tout $S$-schéma $X$, on pose 
\begin{equation}\label{hmdf1b}
\oX=X\times_S\oS,  \ \ \ \coX=X\times_S\coS \ \ \ {\rm et}\ \ \  X_n=X\times_SS_n.
\end{equation} 

On munit $S$ de la structure logarithmique $\cM_S$ définie par son point fermé, 
et $\oS$ et $\coS$ des structures logarithmiques $\cM_\oS$ et $\cM_\coS$ images inverses de $\cM_S$. 

\subsection{}\label{hmdf2}
Comme $\co_\oK$ est un anneau de valuation non discrète de hauteur $1$, 
il est loisible de développer la $\alpha$-algèbre (ou presque-algèbre) sur cet anneau (\cite{ag} 2.10.1) (cf. \cite{ag} 2.6-2.10).   
On choisit un système compatible $(\beta_n)_{n>0}$ 
de racines $n$-ièmes de $p$ dans $\co_\oK$. Pour tout nombre rationnel $\varepsilon>0$, 
on pose $p^\varepsilon=(\beta_n)^{\varepsilon n}$, où $n$ est un entier $>0$ tel que $\varepsilon n$ soit entier.

\subsection{}\label{hmdf3}
Dans ce chapitre, $f\colon (X,\cM_X)\rightarrow (S,\cM_S)$ et $f'\colon (X',\cM_{X'})\rightarrow (S,\cM_S)$
désignent des morphismes {\em adéquats} de schémas logarithmiques (\cite{agt} III.4.7) et 
\begin{equation}\label{hmdf3a}
g\colon (X',\cM_{X'})\rightarrow (X,\cM_X)
\end{equation} 
un $(S,\cM_S)$-morphisme lisse et saturé. On désigne par $X^\circ$ le sous-schéma ouvert maximal de $X$
où la structure logarithmique $\cM_X$ est triviale~; c'est un sous-schéma ouvert de $X_\eta$.
On note $j\colon X^\circ\rightarrow X$ l'injection canonique. Pour tout $X$-schéma $U$, on pose  
\begin{equation}\label{hmdf3b}
U^\circ=U\times_XX^\circ.
\end{equation} 
On note $\hbar\colon \oX\rightarrow X$ et $h\colon \oX^\circ\rightarrow X$ les morphismes canoniques \eqref{hmdf1b}, de sorte que 
l'on a $h=\hbar\circ j_\oX$. 

On désigne par $X'^\rhd$ le sous-schéma ouvert maximal de $X'$
où la structure logarithmique $\cM_{X'}$ est triviale~; c'est un sous-schéma ouvert de $X'^\circ=X'\times_XX^\circ$.
On note $j'\colon X'^\rhd\rightarrow X'$ l'injection canonique. Pour tout $X'$-schéma $U'$, on pose  
\begin{equation}\label{hmdf3c}
U'^\rhd=U'\times_{X'}X'^\rhd.
\end{equation} 
On note $\hbar'\colon \oX'\rightarrow X'$ et $h'\colon \oX'^\rhd\rightarrow X'$ les morphismes canoniques, de sorte que 
l'on a $h'=\hbar'\circ j'_{\oX'}$.

\subsection{}\label{hmdf5}
Reprenons les notations introduites dans \ref{definf10}. 
On munit  $\coX=X\times_S\coS$ et $\coX'=X'\times_S\coS$ \eqref{hmdf1b} des structures logarithmiques $\cM_\coX$ et $\cM_{\coX'}$ 
images inverses respectivement de $\cM_X$ et $\cM_{X'}$. 
On suppose qu'il existe des $(\tS,\cM_{\tS})$-déformations lisses $(\tX,\cM_\tX)$ de $(\coX,\cM_{\coX})$ 
et $(\tX',\cM_{\tX'})$ de $(\coX',\cM_{\coX'})$ \eqref{defing12} et un $(\tS,\cM_{\tS})$-morphisme
\begin{equation}\label{hmdf5a}
\tg\colon (\tX',\cM_{\tX'})\rightarrow (\tX,\cM_\tX)
\end{equation}
qui s'insère dans un diagramme commutatif (à carrés cartésiens) 
\begin{equation}\label{hmdf5b}
\xymatrix{
{(\coX',\cM_{\coX'})}\ar[r]\ar[d]_{\cog}\ar@{}[rd]|{\Box}&{(\tX',\cM_{\tX'})}\ar[d]^{\tg}\\
{(\coX,\cM_\coX)}\ar[r]\ar[d]\ar@{}[rd]|{\Box}&{(\tX,\cM_\tX)}\ar[d]\\
{(\coS,\cM_{\coS})}\ar[r]&{(\tS,\cM_{\tS})}}
\end{equation}
On notera que les carrés sont cartésiens aussi bien dans la catégorie des schémas logarithmiques que 
dans celle des schémas logarithmiques fins. 
{\em On fixe dans ce chapitre les déformations et le $(\tS,\cM_{\tS})$-morphisme $\tg$ \eqref{hmdf5a}.}  

On observera que $\tg$ est lisse en vertu de \ref{hypmdef30}. 

\begin{rema}\label{hmdf50}
Dans le cas relatif \eqref{definf10}, il existe une $(\tS,\cM_{\tS})$-déformation lisse canonique $\tg$ de $g$, à savoir 
\begin{equation}
\tg=g\times_{(S,\cM_S)}(\tS,\cM_\tS)\colon (X',\cM_{X'})\times_{(S,\cM_S)}(\tS,\cM_\tS)\rightarrow (X,\cM_X)\times_{(S,\cM_S)}(\tS,\cM_\tS),
\end{equation}
où l'on considère $(\tS,\cM_\tS)$ comme un schéma logarithmique au-dessus de $(S,\cM_S)$ via $\pr_1$ \eqref{definf7f}, 
le produit étant indifféremment pris dans la catégorie des schémas logarithmiques ou dans celle des schémas logarithmiques fins. 
\end{rema}

\subsection{}\label{hmdf4}
Pour tout entier $n\geq 1$, on note $a\colon X_s\rightarrow X$, $a_n\colon X_s\rightarrow X_n$, 
$a'\colon X'_s\rightarrow X'$ et $a'_n\colon X'_s\rightarrow X'_n$ 
les injections canoniques \eqref{hmdf1b}. 
Le corps résiduel de $\co_K$ étant algébriquement clos, il existe un unique $S$-morphisme $s\rightarrow \oS$. 
Celui-ci induit des immersions fermées $\oa\colon X_s\rightarrow \oX$, $\oa_n\colon X_s\rightarrow \oX_n$, 
$\oa'\colon X'_s\rightarrow \oX'$ et $\oa'_n\colon X'_s\rightarrow \oX'_n$ qui relèvent $a$, $a_n$, $a'$ et $a'_n$, respectivement. 
Comme $\oa_n$ (resp. $\oa'_n$) est un homéomorphisme universel, on peut considérer $\co_{\oX_n}$ (resp. $\co_{\oX'_n}$)
comme un faisceau de $X_{s,\zar}$ ou $X_{s,\et}$ (resp. $X'_{s,\zar}$ ou $X'_{s,\et}$), selon le contexte (cf. \ref{notconv12}). 

\subsection{}\label{hmdf6}
Pour alléger les notations, on pose
\begin{equation}\label{hmdf6a}
\tOmega^1_{X/S}=\Omega^1_{(X,\cM_X)/(S,\cM_S)},
\end{equation}
que l'on considère comme un faisceau de $X_\zar$ ou $X_\et$, selon le contexte \eqref{notconv12}. 
C'est un $\co_X$-module localement libre de type fini.
De même, on pose
\begin{equation}\label{hmdf6b}
\tOmega^1_{X'/S}=\Omega^1_{(X',\cM_{X'})/(S,\cM_S)} \ \ \ {\rm et} \ \ \
\tOmega^1_{X'/X}=\Omega^1_{(X',\cM_{X'})/(X,\cM_X)},
\end{equation}
que l'on considère comme un faisceau de $X'_\zar$ ou $X'_\et$, selon le contexte. 
Ce sont des $\co_{X'}$-modules  localement libres de type fini. 
De plus, on a une suite exacte localement scindée canonique de $\co_{X'}$-modules 
\begin{equation}\label{hmdf6c}
0\rightarrow g^*(\tOmega^1_{X/S})\rightarrow \tOmega^1_{X'/S}\rightarrow \tOmega^1_{X'/X}\rightarrow 0.
\end{equation}

Suivant les conventions de \ref{hmdf4}, pour tout entier $n\geq 1$, on pose 
\begin{equation}\label{hmdf6d}
\tOmega^1_{\oX_n/\oS_n}=\tOmega^1_{X/S}\otimes_{\co_X}\co_{\oX_n},
\end{equation}
que l'on considère comme un faisceau de $X_{s,\zar}$ ou $X_{s,\et}$, selon le contexte, et 
\begin{equation}\label{hmdf6e}
\tOmega^1_{\oX'_n/\oS_n}=\tOmega^1_{X'/S}\otimes_{\co_{X'}}\co_{\oX'_n}\ \ \ 
({\rm resp.} \ \ \ \tOmega^1_{\oX'_n/\oX_n}=\tOmega^1_{X'/X}\otimes_{\co_{X'}}\co_{\oX'_n}),
\end{equation}
que l'on considère comme un faisceau de $X'_{s,\zar}$ ou $X'_{s,\et}$, selon le contexte.

\subsection{}\label{hmdf7}
On désigne par 
\begin{equation}\label{hmdf7a}
\pi\colon E\rightarrow \Et_{/X}
\end{equation}
le $\mU$-site fibré de Faltings associé au morphisme $h\colon \oX^\circ\rightarrow X$ \eqref{hmdf3b} (cf. \ref{ahttf2}).
On munit $E$ de la topologie co-évanescente définie par $\pi$ \eqref{ahttf2}  
et on note $\tE$ le topos des faisceaux de $\mU$-ensembles sur $E$, dit topos de Faltings associé à $h$. 

On note $\Et_{\coh/X}$ la sous-catégorie pleine de $\Et_{/X}$ formée des schémas étales 
de présentation finie sur $X$, munie de la topologie induite par celle de $\Et_{/X}$ \eqref{notconv10}. 
Comme $X$ est noethérien et donc quasi-séparé, $\Et_{\coh/X}$ est une famille $\mU$-petite, topologiquement génératrice du site $\Et_{/X}$ 
et est stable par produits fibrés. On désigne par
\begin{equation}\label{hmdf7f}
\pi_\coh\colon E_\coh\rightarrow \Et_{\coh/X}
\end{equation}
le site fibré déduit de $\pi$ par changement de base par le foncteur d'injection canonique $\Et_{\coh/X}\rightarrow \Et_{/X}$. 
On munit $E_\coh$ de la topologie co-évanescente définie par $\pi_\coh$.
D'après (\cite{agt} VI.10.4), la projection canonique $E_\coh\rightarrow E$ 
induit par restriction une équivalence entre le topos $\tE$ et le topos des faisceaux de $\mU$-ensembles sur $E_\coh$. 
De plus, la topologie co-évanescente de $E_\coh$ est induite par celle de $E$.

On désigne par 
\begin{eqnarray}
\sigma\colon \tE \rightarrow X_\et,\label{hmdf7b}\\
\beta\colon \tE \rightarrow \oX^\circ_\fet,\label{hmdf7c}\\
\psi\colon \oX^\circ_\et\rightarrow \tE,\label{hmdf7h}
\end{eqnarray}
les morphismes canoniques \eqref{ahttf2c}, \eqref{ahttf2d} et \eqref{ahttf2n}.

On désigne par $X_\et\gtimes_{X_\et}\oX^\circ_\et$ le topos co-évanescent du morphisme $h_\et\colon \oX^\circ_\et\rightarrow X_\et$ (\cite{agt} VI.3.12) et par
\begin{equation}\label{hmdf7g}
\rho\colon X_\et\gtimes_{X_\et}\oX^\circ_\et\rightarrow \tE
\end{equation}
le morphisme canonique (\cite{agt} VI.10.15). 

On note $\tE_s$ le sous-topos fermé de $\tE$ complémentaire de l'ouvert $\sigma^*(X_\eta)$ \eqref{ahttf12}, 
\begin{equation}\label{hmdf7d}
\delta\colon \tE_s\rightarrow \tE
\end{equation} 
le plongement canonique et 
\begin{equation}\label{hmdf7e}
\sigma_s\colon \tE_s\rightarrow X_{s,\et}
\end{equation} 
le morphisme de topos induit par $\sigma$ \eqref{ahttf12b}. 

On considère aussi les objets et notations analogues pour $f'$, que l'on munit d'un $^\prime$. 

\subsection{}\label{hmdf8}
Pour tout $(V\rightarrow U)\in \ob(E)$, on note $\oU^V$ la fermeture intégrale de $\oU$ dans $V$.
On désigne par $\ocB$ le préfaisceau sur $E$ défini pour tout $(V\rightarrow U)\in \ob(E)$, par 
\begin{equation}\label{hmdf8a}
\ocB((V\rightarrow U))=\Gamma(\oU^V,\co_{\oU^V}).
\end{equation}
C'est un anneau de $\tE$ (\cite{agt} III.8.16). 
D'après (\cite{agt} III.8.17), on a un homomorphisme canonique 
\begin{equation}\label{hmdf8b}
\sigma^*(\hbar_*(\co_\oX))\rightarrow \ocB.
\end{equation}
Sauf mention explicite du contraire, on considère $\sigma$ \eqref{hmdf7b}
comme un morphisme de topos annelés
\begin{equation}\label{hmdf8c}
\sigma\colon (\tE,\ocB)\rightarrow (X_\et,\hbar_*(\co_\oX)).
\end{equation}

Pour tout entier $n\geq 1$, on pose 
\begin{equation}\label{hmdf8e}
\ocB_n=\ocB/p^n\ocB.
\end{equation}
C'est un anneau de $\tE_s$ (\cite{agt} III.9.7). On désigne par
\begin{equation}\label{hmdf8f}
\sigma_n\colon (\tE_s,\ocB_n)\rightarrow (X_{s,\et},\co_{\oX_n})
\end{equation}
le morphisme de topos annelés induit par $\sigma$ \eqref{hmdf8c} (cf. \cite{agt} (III.9.9.4)). 

On considère aussi les objets et notations analogues pour $f'$, que l'on munit d'un $^\prime$.

\subsection{}\label{hmdf11}
Le foncteur 
\begin{equation}\label{hmdf11a}
\Theta^+\colon E\rightarrow E', \ \ \ (V\rightarrow U)\mapsto (V\times_{X^\circ}X'^\rhd\rightarrow U\times_XX')
\end{equation}
est continu et exact à gauche (\cite{agt} VI.10.12). Il définit donc un morphisme de topos 
\begin{equation}\label{hmdf11b}
\Theta\colon \tE'\rightarrow \tE.
\end{equation}
Il résulte aussitôt des définitions que les carrés du diagramme
\begin{equation}\label{hmdf11c}
\xymatrix{
{X'_\et}\ar[d]_{g}&{\tE'}\ar[l]_-(0.5){\sigma'}\ar[d]^{\Theta}\ar[r]^-(0.5){\beta'}&
{\oX'^\rhd_\fet}\ar[d]^{\upgamma}\\
{X_\et}&{\tE}\ar[l]_{\sigma}\ar[r]^{\beta}&{\oX^\circ_\fet}}
\end{equation}
où $\upgamma\colon \oX'^\rhd\rightarrow \oX^\circ$ est le morphisme induit par $g$, sont commutatifs à isomorphismes canoniques près.  

On a un isomorphisme canonique $\Theta^*(\sigma^*(X_\eta))\simeq \sigma'^*(X'_\eta)$ \eqref{hmdf11c}.
En vertu de (\cite{sga4} IV 9.4.3), il existe donc un morphisme de topos
\begin{equation}\label{hmdf11d}
\uptheta\colon \tE'_s\rightarrow \tE_s
\end{equation}
unique à isomorphisme canonique près tel que le diagramme 
\begin{equation}\label{hmdf11e}
\xymatrix{
{\tE'_s}\ar[r]^{\uptheta}\ar[d]_{\delta'}&{\tE_s}\ar[d]^{\delta}\\
{\tE'}\ar[r]^\Theta&{\tE}}
\end{equation}
soit commutatif à isomorphisme près, et même $2$-cartésien. Il résulte de \eqref{hmdf11c} et (\cite{sga4} IV 9.4.3) 
que le diagramme de morphismes de topos 
\begin{equation}\label{hmdf11f}
\xymatrix{
{\tE'_s}\ar[r]^{\uptheta}\ar[d]_{\sigma'_s}&{\tE_s}\ar[d]^{\sigma_s}\\
{X'_{s,\et}}\ar[r]^{g_s}&{X_{s,\et}}}
\end{equation}
est commutatif à isomorphisme canonique près.

\subsection{}\label{hmdf12}
Pour tout $(V\rightarrow U)\in \ob(E)$, posons $(V'\rightarrow U')=\Theta^+(V\rightarrow U)$ de sorte qu'on a un diagramme commutatif 
\begin{equation}\label{hmdf12a}
\xymatrix{\oX'^\rhd\ar[d]&V'\ar@{}[ld]|{\Box}\ar[l]\ar[r]\ar[d]&\oU'\ar[d]\ar[r]\ar@{}[rd]|{\Box}&\oX'\ar[d]\\
\oX^\circ&V\ar[l]\ar[r]&\oU\ar[r]&\oX}
\end{equation}
On en déduit un morphisme 
\begin{equation}\label{hmdf12b}
\oU'^{V'}\rightarrow \oU^V,
\end{equation}
et par suite un homomorphisme d'anneaux de $\tE$
\begin{equation}\label{hmdf12c}
\ocB\rightarrow \Theta_*(\ocB').
\end{equation}
Nous considérons dans la suite $\Theta$ \eqref{hmdf11b} comme un morphisme de topos annelés (respectivement par $\ocB'$ et $\ocB$). 
Nous utilisons pour les modules la notation $\Theta^{-1}$ pour désigner l'image
inverse au sens des faisceaux abéliens et nous réservons la notation 
$\Theta^*$ pour l'image inverse au sens des modules.

Pour tout entier $n\geq 1$, l'homomorphisme canonique $\Theta^{-1}(\ocB)\rightarrow \ocB'$
induit un homomorphisme $\uptheta^*(\ocB_n)\rightarrow \ocB'_n$. 
Le morphisme $\uptheta$ est donc sous-jacent à un morphisme de topos annelés, que l'on note 
\begin{equation}\label{hmdf12d}
\uptheta_n\colon (\tE'_s,\ocB'_n)\rightarrow (\tE_s,\ocB_n).
\end{equation}
Compte tenu de \eqref{hmdf11f}, on vérifie aussitôt que le diagramme de morphismes de topos 
\begin{equation}\label{hmdf12e}
\xymatrix{
{(\tE'_s,\ocB'_n)}\ar[r]^{\uptheta_n}\ar[d]_{\sigma'_n}&{(\tE_s,\ocB_n)}\ar[d]^{\sigma_n}\\
{(X'_{s,\et},\co_{\oX'_n})}\ar[r]^{\ogg_n}&{(X_{s,\et},\co_{\oX_n})}}
\end{equation}
est commutatif à isomorphisme canonique près.

\subsection{}\label{hmdf20}
Soit $U$ un objet de $\Et_{/X}$. 
D'après (\cite{agt} VI.10.14), le topos  $\tE_{/\sigma^*(U)}$, localisé de $\tE$ en $\sigma^*(U)$, 
est canoniquement équivalent au topos de Faltings associé au morphisme $\oU^\circ\rightarrow U$ \eqref{ahttf2e}. 
On désigne par
\begin{equation}\label{hmdf20a}
\jmath_U\colon \tE_{/\sigma^*(U)}\rightarrow \tE
\end{equation}
le morphisme de localisation de $\tE$ en $\sigma^*(U)$, qui s'identifie au morphisme de fonctorialité induit par le morphisme canonique 
$U\rightarrow X$, et par 
\begin{equation}\label{hmdf20b}
\beta_U\colon \tE_{/\sigma^*(U)} \rightarrow \oU^\circ_\fet
\end{equation}
le morphisme canonique \eqref{ahttf2d}. 

On considère les notations analogues pour le topos $\tE'$ et les objets de $\Et_{/X'}$, que l'on munit d'un exposant~$^\prime$. 

Soit $\upmu \colon U'\rightarrow U$ un morphisme au-dessus de $g\colon X'\rightarrow X$ tel que $U'\rightarrow X'$ et 
$U\rightarrow X$ soient étales. On note abusivement $\oupmu\colon \oU'^\rhd\rightarrow \oU^\circ$ le morphisme induit par $\upmu$, 
et on désigne par 
\begin{equation}\label{hmdf20c}
\Theta_\upmu\colon \tE'_{/\sigma'^*(U')}\rightarrow \tE_{/\sigma^*(U)}
\end{equation}
le morphisme de fonctorialité induit par $\upmu$ \eqref{hmdf11}. Les diagrammes 
\begin{equation}\label{hmdf20d}
\xymatrix{
{\tE'_{/\sigma'^*(U')}}\ar[r]^-(0.4){\Theta_\upmu}\ar[d]_{\jmath'_{U'}}&{\tE_{/\sigma^*(U)}}\ar[d]^{\jmath_U}\\
{\tE'}\ar[r]^-(0.5)\Theta&\tE}
\end{equation}
\begin{equation}\label{hmdf20e}
\xymatrix{
{\tE'_{/\sigma'^*(U')}}\ar[r]^-(0.4){\Theta_\upmu}\ar[d]_{\beta'_{U'}}&{\tE_{/\sigma^*(U)}}\ar[d]^{\beta_U}\\
{\oU'^\rhd_\fet}\ar[r]^-(0.5)\oupmu&{\oU^\circ_\fet}}
\end{equation}
sont commutatifs à isomorphismes canoniques près (\cite{agt} (VI.10.12.6)).

\subsection{}\label{ahttfg12}
Soient $\ox'$ un point géométrique de $X'$, $\ox=g(\ox')$, $\uX$ (resp. $\uX'$) le localisé strict de $X$ en $\ox$ (resp. $X'$ en $\ox'$),
$\ug\colon \uX'\rightarrow \uX$ et $\uupgamma\colon \uoX'^\rhd\rightarrow \uoX^\circ$ les morphismes induits par $g$. 
On désigne par $\tuE$ (resp. $\tuE'$) le topos de Faltings associé au morphisme canonique $\uoX^\circ \rightarrow X$
(resp. $\uoX'^\rhd \rightarrow X'$), par 
\begin{equation}\label{ahttfg12a}
\uTheta\colon \tuE'\rightarrow \tuE
\end{equation}
le morphisme de fonctorialité induit par le morphisme $\ug$, et par 
\begin{eqnarray}
\ubeta\colon \tuE\rightarrow  \uoX^\circ_\fet,\label{ahttfg12b}\\
\ubeta'\colon \tuE'\rightarrow \uoX'^\rhd_\fet,\label{ahttfg12c}
\end{eqnarray}
les morphismes canoniques \eqref{ahttf2d}.
D'après (\cite{agt} (VI.10.12.6)), le diagramme 
\begin{equation}\label{ahttfg12d}
\xymatrix{
{\tuE'}\ar[r]^{\uTheta}\ar[d]_{\ubeta'}&{\tuE}\ar[d]^{\ubeta}\\
{\uoX'^\rhd_\fet}\ar[r]^{\uupgamma}&{\uoX^\circ_\fet}}
\end{equation}
est commutatif à isomorphisme canonique près. On désigne par 
\begin{eqnarray}
\theta\colon \uoX^\circ_\fet\rightarrow \tuE,\label{ahttfg12e}\\
\theta'\colon \uoX'^\rhd_\fet\rightarrow \tuE',\label{ahttfg12f}
\end{eqnarray}
les sections canoniques de $\ubeta$ et $\ubeta'$ (\cite{agt} VI.10.23). Le diagramme 
\begin{equation}\label{ahttfg12g}
\xymatrix{
{\uoX'^\rhd_\fet}\ar[r]^-(0.5)\uupgamma\ar[d]_{\theta'}&{\uoX^\circ_\fet}\ar[d]^{\theta}\\
{\tuE'}\ar[r]^-(0.5){\uTheta}&{\tuE}}
\end{equation}
est commutatif à isomorphisme canonique près. En effet, pour tout $\uX$-schéma étale, séparé et de présentation finie  $U$, 
si l'on note $U^\rf$ sa partie $\uX$-finie ({\em i.e.}, la somme disjointe des localisés stricts de $U$ en les points de $U_\ox$), alors 
$U^\rf\times_\uX\uX'$ est la partie $\uX'$-finie de $U\times_\uX\uX'$ (cf. \cite{agt} VI.10.22). 

Par ailleurs, le diagramme 
\begin{equation}\label{ahttfg12h}
\xymatrix{
{\tuE'}\ar[r]^-(0.5){\uTheta}\ar[d]_{\Phi'}&{\tuE}\ar[d]^\Phi\\
{\tE'}\ar[r]^-(0.5){\Theta}&{\tE}}
\end{equation}
où les flèches verticales sont les morphismes de fonctorialité induits par les morphismes canoniques $\uX\rightarrow X$ et $\uX'\rightarrow X'$, 
est commutatif à isomorphisme canonique près. On désigne par 
\begin{eqnarray}
\varphi_\ox\colon \tE\rightarrow \uoX^\circ_\fet,\label{ahttfg12i}\\
\varphi'_{\ox'}\colon \tE'\rightarrow \uoX'^\rhd_\fet,\label{ahttfg12ii}
\end{eqnarray}
les foncteurs composés $\theta^*\circ \Phi^*$ et $\theta'^*\circ \Phi'^*$ (cf. \ref{TFA14}).  On a donc un isomorphisme canonique de foncteurs
\begin{equation}\label{ahttfg12j}
\uupgamma^* \circ \varphi_\ox\stackrel{\sim}{\rightarrow} \varphi'_{\ox'} \circ \Theta^*.
\end{equation}

\subsection{}\label{hmdf17}
On désigne par $G$ le {\em $\mU$-site de Faltings relatif} associé au couple de morphismes 
$(h\colon \oX^\circ \rightarrow X,g\colon X'\rightarrow X)$ (\cite{ag} 3.4.1). Les objets de la catégorie sous-jacente à $G$ sont
les triplets  $(U,U'\rightarrow U, V\rightarrow U)$ formés d'un $X$-schéma $U$ et de deux morphismes $U'\rightarrow U$ et $V\rightarrow U$ 
au dessus de $g$ et $h$ respectivement, c'est-à-dire des diagrammes commutatifs de morphismes de schémas 
\begin{equation}
\xymatrix{
U'\ar[r]\ar[d]&U\ar[d]&V\ar[l]\ar[d]\\
X'\ar[r]^g&X&\oX^\circ\ar[l]_h}
\end{equation}
tels que les morphismes $U\rightarrow X$ et $U'\rightarrow X'$ soient étales et que le morphisme $V\rightarrow \oU^\circ$ soit fini étale; 
un tel objet sera noté $(U'\rightarrow U\leftarrow V)$.  Soient $(U'\rightarrow U\leftarrow V)$ et $(U'_1\rightarrow U_1\leftarrow V_1)$ deux objets de $G$. 
Un morphisme de $(U'_1\rightarrow U_1\leftarrow V_1)$ dans $(U'\rightarrow U\leftarrow V)$ est la donnée de trois morphismes 
$U'_1\rightarrow U'$, $U_1\rightarrow U$ et $V_1\rightarrow V$ au-dessus de $X'$, $X$ et $\oX^\circ$ respectivement, qui rendent le diagramme 
\begin{equation}
\xymatrix{
U'_1\ar[r]\ar[d]&U_1\ar[d]&V_1\ar[l]\ar[d]\\
U'\ar[r]&U&V\ar[l]}
\end{equation}
commutatif. 

On munit $G$ de la topologie {\em co-évanescente} (\cite{ag} 3.4.1), c'est-à-dire la topologie engendrée par les recouvrements 
\[
\{(U'_i\rightarrow U_i\leftarrow V_i)\rightarrow (U'\rightarrow U\leftarrow V)\}_{i\in I}
\] 
des trois types suivants~:
\begin{itemize}
\item[(a)] $V_i=V$, $U_i=U$ pour tout $i\in I$, et $(U'_i\rightarrow U')_{i\in I}$ est une famille couvrante.
\item[(b)] $U'_i=U'$, $U_i=U$ pour tout $i\in I$, et $(V_i\rightarrow V)_{i\in I}$ est une famille couvrante.
\item[(c)] $I=\{1\}$, $U'_1=U'$ et le morphisme $V_1\rightarrow V\times_{U}U_1$ est un isomorphisme
(il n'y a aucune condition sur le morphisme $U_1\rightarrow U$). 
\end{itemize}

On note $\tG$ le {\em $\mU$-topos de Faltings relatif} associé à $(h,g)$,  
c'est-à-dire le topos des faisceaux de $\mU$-ensembles sur $G$. 

On désigne par 
\begin{eqnarray}
\pi\colon \tG \rightarrow X'_\et,\label{hmdf17a}\\
\lambda\colon \tG \rightarrow \oX^\circ_\fet,\label{hmdf17b}
\end{eqnarray}
les morphismes canoniques (\cite{ag} 3.4.4). 
D'après (\cite{ag} 3.4.18), le diagramme commutatif canonique
\begin{equation}\label{hmdf17c}
\xymatrix{
X'\ar[d]_g&{\oX'^\rhd}\ar[d]^{\upgamma}\ar[l]_{h'}\\
X&{\oX^\circ}\ar[l]_h}
\end{equation}
induit des morphismes de topos 
\begin{equation}\label{hmdf17d}
\xymatrix{
\tE'\ar[r]^{\tau}&\tG\ar[r]^{\lgg}&\tE}
\end{equation}
dont le composé est le morphisme $\Theta\colon \tE'\rightarrow \tE$ \eqref{hmdf11b}.
Les triangles et carrés du diagramme de morphismes de topos 
\begin{equation}\label{hmdf17e}
\xymatrix{
&\tE'\ar[d]^{\tau}\ar[r]^{\beta'}\ar[ld]_{\sigma'} & \oX'^\rhd_\fet\ar[d]^{\upgamma_\fet}\\
X'_\et\ar[d]_{g_\et}&\tG\ar[d]_-(0.5){\lgg}\ar[r]^-(0.4){\lambda}\ar[l]_-(0.4){\pi}&\oX^\circ_\fet\\
X_\et&\tE\ar[ur]_{\beta}\ar[l]_{\sigma}&}
\end{equation}
sont commutatifs à isomorphismes canoniques près (\cite{ag} (3.4.18.4)). 

On désigne par $X'_\et\gtimes_{X_\et}\oX^\circ_\et$ le produit orienté des morphismes $g_\et\colon X'_\et\rightarrow X_\et$
et $f_\et\colon \oX^\circ_\et\rightarrow X_\et$ (\cite{agt} VI.3.10) et par 
\begin{equation}\label{hmdf17f}
\varrho\colon X'_\et\gtimes_{X_\et}\oX^\circ_\et\rightarrow \tG
\end{equation}
le morphisme canonique de topos (\cite{ag} (3.4.9.2)). 
On vérifie aussitôt que les carrés du diagramme de morphismes de topos 
\begin{equation}\label{hmdf17g}
\xymatrix{
{X'_\et\gtimes_{X'_\et}\oX'^\rhd_\et}\ar[d]\ar[r]^-(0.5){\rho'}&\tE'\ar[d]^\tau\\
{X'_\et\gtimes_{X_\et}\oX^\circ_\et}\ar[d]\ar[r]^-(0.5){\varrho}&\tG\ar[d]^\lgg\\
{X_\et\gtimes_{X_\et}\oX^\circ_\et}\ar[r]^-(0.5){\rho}&\tE}
\end{equation}
où $\rho$ et $\rho'$ sont les morphismes canoniques \eqref{hmdf7g}, 
sont commutatifs à isomorphismes canoniques près (\cite{agt} VI.4.10 et \cite{ag} 3.4.17).

\subsection{}\label{hmdf26}
Considérons un diagramme commutatif de morphismes de schémas
\begin{equation}\label{hmdf26a}
\xymatrix{
U'\ar[r]\ar[d]&U\ar[d]\\
X'\ar[r]&X}
\end{equation}
tel que les flèches verticales soient des morphismes étales. 
Pour tout préfaisceau $F$ sur $G$, on définit le préfaisceau $F_{U'\rightarrow U}$ sur $\Et_{\rf/\oU^\circ}$ en posant pour tout $V\in \ob(\Et_{\rf/\oU^\circ})$,
\begin{equation}\label{hmdf26b}
F_{U'\rightarrow U}(V)=F(U'\rightarrow U\leftarrow V).
\end{equation}
Si $F$ est un faisceau de $\tG$, alors $F_{U'\rightarrow U}$ est un faisceau de $\oU^\circ_\fet$.

\subsection{}\label{hmdf27}
A tout point géométrique $\ox'$ de $X'$, on associe une catégorie $\fC_{\ox'}$ de la façon suivante. 
Les objets de $\fC_{\ox'}$ sont les diagrammes commutatifs de morphismes de schémas 
\begin{equation}\label{hmdf27a}
\xymatrix{
&U'\ar[r]\ar[d]&U\ar[d]\\
\ox'\ar[r]\ar[ru]&X'\ar[r]&X}
\end{equation}
tels que les morphismes $U'\rightarrow X'$ et $U\rightarrow X$ soient étales. 
Un tel objet sera noté $(\ox'\rightarrow U'\rightarrow U)$. Soient $(\ox'\rightarrow U'\rightarrow U)$, $(\ox'\rightarrow U'_1\rightarrow U_1)$ deux objets de
$\fC_{\ox'}$. Un morphisme de $(\ox'\rightarrow U'_1\rightarrow U_1)$ vers $(\ox'\rightarrow U'\rightarrow U)$ est la donnée d'un $X'$-morphisme
$U'_1\rightarrow U'$ et d'un $X$-morphisme $U_1\rightarrow U$ tels que le diagramme 
\begin{equation}\label{hmdf27b}
\xymatrix{
&U'_1\ar[r]\ar[d]&U_1\ar[d]\\
\ox'\ar[r]\ar[ru]&U'\ar[r]&U}
\end{equation}
soit commutatif. On observera que les produits fibrés sont représentables dans $\fC_{\ox'}$ (cf. la preuve de \cite{agt} VI.10.3).
Les limites projectives finies sont donc représentables dans $\fC_{\ox'}$ (cf. \cite{sga4} I 2.3).
Par suite, la catégorie $\fC_{\ox'}$ est cofiltrante (\cite{sga4} I 2.7.1).

\subsection{}\label{hmdf28}
Soient $(\oy\rightsquigarrow \ox')$ un point de $X'_\et\gtimes_{X_\et}\oX^\circ_\et$ (\cite{ag} 3.4.23), 
$\uX'$ le localisé strict de $X'$ en $\ox'$, $\uX$ le localisé strict de $X$ en $g(\ox')$, 
$\fC_{\ox'}$ la catégorie associée à $\ox'$ dans \ref{hmdf27}. 
On désigne par $u\colon \oy\rightarrow \uX$ le $X$-morphisme qui définit le point $(\oy\rightsquigarrow \ox')$, 
et par $v\colon \oy\rightarrow \uoX^\circ$ le $\oX^\circ$-morphisme induit \eqref{hmdf3b}. 
Pour chaque objet $(\ox'\rightarrow U'\rightarrow U)$ de $\fC_{\ox'}$, on a un $X'$-morphisme canonique 
$\uX'\rightarrow U'$ et un $X$-morphisme canonique $\uX\rightarrow U$ qui s'insèrent dans un diagramme commutatif  
\begin{equation}\label{hmdf28a}
\xymatrix{
&\uX'\ar[r]^g\ar[d]&\uX\ar[d]\\
\ox'\ar[r]\ar[ru]&U'\ar[r]&U}
\end{equation}
où on a encore noté $g$ le morphisme induit par $g$. On en déduit un morphisme $\uoX\rightarrow \oU$. 
Le morphisme $v\colon \oy\rightarrow \uoX^\circ$ induit alors un point géométrique de $\oU^\circ$ que l'on note encore $\oy$. 
Le diagramme 
\begin{equation}\label{hmdf28b}
\xymatrix{
\uX\ar[d]&\oy\ar[l]_-(0.5)u\ar[d]\\
U&\oU^\circ\ar[l]}
\end{equation}
est commutatif. 

On désigne par $\varrho(\oy\rightsquigarrow \ox')$ l'image de $(\oy\rightsquigarrow \ox')$ par  le morphisme $\varrho$ \eqref{hmdf17f}, 
qui est donc un point de $\tG$. D'après (\cite{ag} (6.5.9.5)), pour tout préfaisceau $F$ sur $G$, on a un isomorphisme canonique
\begin{equation}\label{hmdf28e}
F^a_{\varrho(\oy\rightsquigarrow \ox')}\stackrel{\sim}{\rightarrow} \underset{\underset{(\ox'\rightarrow U'\rightarrow U)\in \fC^\circ_{\ox'}}{\longrightarrow}}{\lim}\ 
(F^a_{U'\rightarrow U})_{\rho_{\oU^\circ}(\oy)},
\end{equation}
où $F_{U'\rightarrow U}$ est le préfaisceau sur $\Et_{\rf/\oU^\circ}$ défini dans \eqref{hmdf26b}, l'exposant $^a$ désigne les faisceaux associés et 
$\rho_{\oU^\circ}\colon \oU^\circ_\et \rightarrow \oU^\circ_\fet$ est le morphisme canonique \eqref{notconv10a}.

\subsection{}\label{hmdf29}
Soient $\ox'$ un point géométrique de $X'$, $\uX'$ le localisé strict de $X'$ en $\ox'$, $\uX$ le localisé strict de $X$ en $g(\ox')$. 
On désigne par $\uG$ (resp. $\tuG$) le site (resp. topos) de Faltings relatif associé au couple de morphismes 
$(\uh\colon \uoX^\circ\rightarrow \uX, \ug\colon \uX'\rightarrow \uX)$ induits par $h$ et $g$ \eqref{hmdf3},  par 
\begin{equation}\label{hmdf29a}
\Phi\colon \tuG\rightarrow \tG
\end{equation}
le morphisme de fonctorialité (\cite{ag} (3.4.10.3)) et par 
\begin{equation}\label{hmdf29b}
\vartheta\colon \uoX^\circ_\fet\rightarrow \tuG
\end{equation}
le morphisme défini dans (\cite{ag} (3.4.26.9)). On pose  
\begin{equation}\label{hmdf29c}
\phi_{\ox'}=\vartheta^*\circ \Phi^*\colon \tG\rightarrow \uoX^\circ_\fet.
\end{equation}

D'après (\cite{ag} 3.4.34), pour tout groupe abélien $F$ de $\tG$
et tout entier $q\geq 0$, on a un isomorphisme canonique et fonctoriel
\begin{equation}\label{hmdf29d}
\rR^q\pi_*(F)_{\ox'}\stackrel{\sim}{\rightarrow}\rH^q(\uoX^\circ_\fet,\phi_{\ox'}(F)). 
\end{equation}

On désigne par 
\begin{equation}\label{hmdf29e}
\varphi'_{\ox'}\colon \tE'\rightarrow \uoX'^\rhd_\fet
\end{equation}
le foncteur canonique \eqref{ahttfg12ii} et par 
\begin{equation}\label{hmdf29f}
\uupgamma\colon \uoX'^\rhd\rightarrow \uoX^\circ
\end{equation} 
le morphisme induit par $\ug$.

\begin{prop}[\cite{ag} 6.5.16]\label{hmdf30}
Conservons les hypothèses et notations de \ref{hmdf29}. 
\begin{itemize}
\item[{\rm (i)}] Pour tout groupe abélien $F$ de $\tE'$ et tout entier $q\geq 0$, on a un isomorphisme canonique fonctoriel
\begin{equation}\label{hmdf30a}
\phi_{\ox'}(\rR^q\tau_*(F))\stackrel{\sim}{\rightarrow}\rR^q\uupgamma_{\fet*}(\varphi'_{\ox'}(F)).
\end{equation}
\item[{\rm (ii)}] Pour toute suite exacte de faisceaux abéliens 
$0\rightarrow F'\rightarrow F\rightarrow F''\rightarrow 0$ de $\tE'$ et tout entier $q\geq 0$, le diagramme 
\begin{equation}\label{hmdf30b}
\xymatrix{
{\phi_{\ox'}(\rR^q\tau_*(F''))} \ar[r]\ar[d]&{\phi_{\ox'}(\rR^{q+1}\tau_*(F'))}\ar[d]\\
{\rR^q\uupgamma_{\fet*}(\varphi_{\ox'}(F''))}\ar[r]&{\rR^{q+1}\uupgamma_{\fet*}(\varphi_{\ox'}(F'))}}
\end{equation}
où les flèches verticales sont les isomorphismes canoniques \eqref{hmdf30a} et les flèches horizontales  
sont les bords des suites exactes longues de cohomologie, est commutatif. 
\end{itemize}
\end{prop}

\subsection{}\label{hmdf18}
On désigne par $\tG_s$ le sous-topos fermé de $\tG$ complémentaire de l'ouvert $\pi^*(X'_\eta)$ (\cite{sga4} IV 9.3.5), et par 
\begin{equation}\label{hmdf18a}
\kappa\colon \tG_s\rightarrow \tG
\end{equation} 
le plongement canonique (\cite{sga4} IV 9.3.5) (cf. \cite{ag} 6.5.2). En vertu de (\cite{sga4} IV 9.4.3), il existe un morphisme de topos
\begin{equation}\label{hmdf18b}
\pi_s\colon \tG_s\rightarrow X'_{s,\et}
\end{equation} 
unique à isomorphisme canonique près tel que le diagramme
\begin{equation}\label{hmdf18c}
\xymatrix{
{\tG_s}\ar[r]^-(0.5){\pi_s}\ar[d]_{\kappa}&{X'_{s,\et}}\ar[d]^{a'}\\
{\tG}\ar[r]^-(0.5){\pi}&{X'_\et}}
\end{equation}
où $a'$ est l'injection canonique, soit commutatif à isomorphisme près, et même $2$-cartésien.

On a un isomorphisme canonique $\tau^*(\pi^*(X'_\eta))\simeq \sigma'^*(X'_\eta)$ \eqref{hmdf17e}.
En vertu de (\cite{sga4} IV 9.4.3), il existe donc un morphisme de topos
\begin{equation}\label{hmdf18d}
\tau_s\colon \tE'_s\rightarrow \tG_s
\end{equation}
unique à isomorphisme canonique près tel que le diagramme 
\begin{equation}\label{hmdf18e}
\xymatrix{
{\tE'_s}\ar[r]^{\tau_s}\ar[d]_{\delta'}&{\tG_s}\ar[d]^{\kappa}\\
{\tE'}\ar[r]^\tau&{\tG}}
\end{equation}
soit commutatif à isomorphisme près (cf. \ref{hmdf7}). 

Les foncteurs $\delta'_*$ et $\kappa_*$ étant exacts, 
pour tout groupe abélien $F$ de $\tE'_s$ et tout entier $i\geq 0$, on a un isomorphisme canonique 
\begin{equation}\label{hmdf18k}
\kappa_*(\rR^i\tau_{s*}(F))\stackrel{\sim}{\rightarrow}\rR^i\tau_*(\delta'_*F). 
\end{equation}

Il résulte de \eqref{hmdf17e} et (\cite{sga4} IV 9.4.3) 
que le diagramme de morphismes de topos 
\begin{equation}\label{hmdf18f}
\xymatrix{
{\tE'_s}\ar[r]^{\tau_s}\ar[d]_{\sigma'_s}&{\tG_s}\ar[dl]^{\pi_s}\\
{X'_{s,\et}}&}
\end{equation}
est commutatif à isomorphisme canonique près. 

On a un isomorphisme canonique $\lgg^*(\sigma^*(X_\eta))\simeq \pi^*(X'_\eta)$ \eqref{hmdf17e}.
En vertu de (\cite{sga4} IV 9.4.3), il existe donc un morphisme de topos
\begin{equation}\label{hmdf18g}
\lgg_s\colon \tG_s\rightarrow \tE_s
\end{equation}
unique à isomorphisme canonique près tel que le diagramme 
\begin{equation}\label{hmdf18h}
\xymatrix{
{\tG_s}\ar[r]^{\lgg_s}\ar[d]_{\kappa}&{\tE_s}\ar[d]^{\delta}\\
{\tG}\ar[r]^\lgg&{\tE}}
\end{equation}
soit commutatif à isomorphisme près (cf. \ref{hmdf7}). 

Les foncteurs $\kappa_*$ et $\delta_*$ étant exacts, 
pour tout groupe abélien $F$ de $\tG_s$ et tout entier $i\geq 0$, on a un isomorphisme canonique 
\begin{equation}\label{hmdf18hh}
\delta_*(\rR^i\lgg_{s*}(F))\stackrel{\sim}{\rightarrow}\rR^i\lgg_*(\kappa_*F). 
\end{equation}

Il résulte de \eqref{hmdf17e} et (\cite{sga4} IV 9.4.3) 
que le diagramme de morphismes de topos 
\begin{equation}\label{hmdf18i}
\xymatrix{
{\tG_s}\ar[r]^{\lgg_s}\ar[d]_{\pi_s}&{\tE_s}\ar[d]^{\sigma_s}\\
{X'_{s,\et}}\ar[r]^{g_{s,\et}}&{X_{s,\et}}}
\end{equation}
est commutatif à isomorphisme canonique près.

Il résulte encore de (\cite{sga4} IV 9.4.3) que le composé $\lgg_s\circ \tau_s$ est le morphisme \eqref{hmdf11d}
\begin{equation}\label{hmdf18j}
\uptheta\colon \tE'_s\rightarrow \tE_s.
\end{equation}

\subsection{}\label{hmdf19}
Pour tout objet $(U'\rightarrow U\leftarrow V)$ de $G$, on note $\oU'^V$ la fermeture intégrale de $\oU'$ dans $U'\times_UV$. 
On désigne par $\ocB^!$ le préfaisceau sur $G$ défini pour tout $(U'\rightarrow U\leftarrow V)\in \ob(G)$ par
\begin{equation}\label{hmdf19a}
\ocB^!(U'\rightarrow U\leftarrow V)=\Gamma(\oU'^V,\co_{\oU'^V}).
\end{equation}
Comme $\oX'$ est normal et localement irréductible (\cite{agt} III.4.2(iii)),
$\ocB^!$ est un faisceau pour la topologie co-évanescente de $G$ d'après (\cite{ag} 3.6.4). 

D'après (\cite{ag} 6.5.18), on a des homomorphismes canoniques
\begin{eqnarray}
\ocB^!&\rightarrow&\tau_*(\ocB'),\label{hmdf19b}\\
\ocB&\rightarrow& \lgg_*(\ocB^!),\label{hmdf19c}
\end{eqnarray}
dont le premier est un isomorphisme en vertu de (\cite{ag} 6.5.19).
On a aussi un homomorphisme canonique 
\begin{equation}\label{hmdf19d}
\hbar'_*(\co_{\oX'})\rightarrow \pi_*(\ocB^!).
\end{equation}

Pour tout entier $n\geq 1$, on pose 
\begin{equation}\label{hmdf19f}
\ocB^!_n=\ocB^!/p^n\ocB^!.
\end{equation}
C'est un anneau de $\tG_s$ (\cite{ag} 6.5.27). 
L'homomorphisme canonique $\pi^*(\hbar_*(\co_{\oX'}))\rightarrow \ocB^!$ \eqref{hmdf19d}
induit un homomorphisme $\pi_s^*(\co_{\oX'_n})\rightarrow \ocB^!_n$ de $\tG_s$ (\cite{ag} 6.5.28).
Le morphisme $\pi_s$ \eqref{hmdf18b} est donc sous-jacent à un morphisme de topos annelés, que l'on note 
\begin{equation}\label{hmdf19h}
\pi_n\colon (\tG_s,\ocB^!_n)\rightarrow (\oX'_{s,\et},\co_{\oX'_n}).
\end{equation}

L'homomorphisme canonique $\tau^*(\ocB^!)\rightarrow \ocB'$ \eqref{hmdf19b}
induit un homomorphisme $\tau_s^*(\ocB^!_n)\rightarrow \ocB'_n$ de $\tE'_s$. 
Le morphisme $\tau_s$ \eqref{hmdf18d} est donc sous-jacent à un morphisme de topos annelés, que l'on note 
\begin{equation}\label{hmdf19i}
\tau_n\colon (\tE'_s,\ocB'_n)\rightarrow (\tG_s,\ocB^!_n).
\end{equation}

L'homomorphisme canonique $\lgg^{-1}(\ocB)\rightarrow \ocB^!$ \eqref{hmdf19c}
induit un homomorphisme $\lgg_s^*(\ocB_n)\rightarrow \ocB^!_n$ de $\tG_s$. 
Le morphisme $\lgg_s$ \eqref{hmdf18g} est donc sous-jacent à un morphisme de topos annelés, que l'on note 
\begin{equation}\label{hmdf19j}
\lgg_n\colon (\tG_s,\ocB^!_n)\rightarrow (\tE_s,\ocB_n).
\end{equation}

On vérifie aussitôt que le composé $\lgg_n\circ \tau_n$ est le morphisme \eqref{hmdf12d}
\begin{equation}\label{hmdf19k}
\uptheta_n\colon (\tE'_s,\ocB'_n)\rightarrow (\tE_s,\ocB_n).
\end{equation}

D'après (\cite{ag} 6.5.28.10), le triangle et le carré du diagramme de morphismes de topos annelés
\begin{equation}\label{hmdf19l}
\xymatrix{
{(\tE'_s,\ocB'_n)}\ar[r]^-(0.5){\tau_n}\ar[rd]_{\sigma'_n}&{(\tG_s,\ocB_n^!)}\ar[d]^-(0.5){\pi_n}\ar[r]^-(0.5){\lgg_n}&{(\tE_s,\ocB_n)}\ar[d]^-(0.5){\sigma_n}\\
&{(X'_{s,\et},\co_{\oX'_n})}\ar[r]^{\ogg_n}&{(X_{s,\et},\co_{\oX_n})}}
\end{equation}
sont commutatifs à isomorphismes canoniques près  (\cite{egr1} 1.2.3).

\subsection{}\label{hmdf14}
Reprenons les notations de \ref{notconv13}. On désigne par $\bvocB$ l'anneau $(\ocB_{n+1})_{n\in \mN}$ de $\tE_s^{\mN^\circ}$ \eqref{hmdf8e},
par $\bvocB'$ l'anneau $(\ocB'_{n+1})_{n\in \mN}$ de $\tE'^{\mN^\circ}_s$ et 
par $\bvocB^!$ l'anneau $(\ocB^!_{n+1})_{n\in \mN}$ de $\tG_s^{\mN^\circ}$ \eqref{hmdf19f}. 
On désigne  par $\co_{\bvoX}$ l'anneau $(\co_{\oX_{n+1}})_{n\in \mN}$ de $X_{s,\et}^{\mN^\circ}$ ou $X_{s,\zar}^{\mN^\circ}$,
selon le contexte \eqref{notconv12}, et
par $\txi^{-1}\tOmega^1_{\bvoX/\bvoS}$ le $\co_{\bvoX}$-module $(\txi^{-1}\tOmega^1_{\oX_{n+1}/\oS_{n+1}})_{n\in \mN}$ \eqref{hmdf6d}.
On désigne par $\co_{\bvoX'}$ l'anneau $(\co_{\oX'_{n+1}})_{n\in \mN}$ de $X'^{\mN^\circ}_{s,\et}$ ou $X'^{\mN^\circ}_{s,\zar}$, selon le contexte,
et par $\txi^{-1}\tOmega^1_{\bvoX'/\bvoS}$ (resp. $\txi^{-1}\tOmega^1_{\bvoX'/\bvoX}$) 
le $\co_{\bvoX'}$-module $(\txi^{-1}\tOmega^1_{\oX_{n+1}/\oS_{n+1}})_{n\in \mN}$ (resp. $(\txi^{-1}\tOmega^1_{\oX_{n+1}/\oX_{n+1}})_{n\in \mN}$) \eqref{hmdf6e}. 
On prendra garde de ne pas confondre $\co_{\bvoX}$ et $\co_{\coX}$ (resp. $\co_{\bvoX'}$ et $\co_{\coX'}$) \eqref{hmdf1b}.
Considérons le diagramme de morphismes de topos annelés 
\begin{equation}\label{hmdf14a}
\xymatrix{
{(\tE'^{\mN^\circ}_s,\bvocB')}\ar[rd]^-(0.5){\bvsigma'}\ar[d]_{\bvtau}\ar@/_3pc/[dd]_{\bvuptheta}&\\
{(\tG^{\mN^\circ}_s,\bvocB^!)}\ar[r]^-(0.5){\bvpi}\ar[d]_{\bvlgg}&{(X'^{\mN^\circ}_{s,\et},\co_{\bvoX'})}\ar[d]^{\bvogg}\\
{(\tE^{\mN^\circ}_s,\bvocB)}\ar[r]^-(0.5){\bvsigma}&{(X^{\mN^\circ}_{s,\et},\co_{\bvoX})}}
\end{equation}
induits par $(\uptheta_{n+1})_{n\in \mN}$ \eqref{hmdf12d}, $(\tau_{n+1})_{n\in \mN}$ \eqref{hmdf19i}, 
$(\lgg_{n+1})_{n\in \mN}$ \eqref{hmdf19j}, $(\sigma_{n+1})_{n\in \mN}$ \eqref{hmdf8f}, $(\sigma'_{n+1})_{n\in \mN}$, 
$(\pi_{n+1})_{n\in \mN}$ \eqref{hmdf19h} et $\ogg$ (cf. \cite{agt} III.7.5). 
Les deux triangles et le carré sont commutatifs à isomorphismes canoniques près \eqref{hmdf19l}.

\subsection{}\label{hmdf13}
On pose $\cS=\Spf(\co_C)$ et on désigne par $\fX$ (resp. $\fX'$) le schéma formel complété $p$-adique de $\oX$ (resp. $\oX'$), 
et par $\fgg\colon \fX'\rightarrow \fX$ le morphisme induit par $g$ \eqref{hmdf3a}. Pour tout entier $n\geq 1$, on note
\begin{equation}\label{hmdf13c}
u_n\colon (X_{s,\et},\co_{\oX_n})\rightarrow (X_{s,\zar},\co_{\oX_n})
\end{equation}
le  morphisme canonique \eqref{notconv12j}. 
On désigne par 
\begin{equation}\label{hmdf13d}
\bvu\colon (X_{s,\et}^{\mN^\circ},\co_{\bvoX})\rightarrow (X_{s,\zar}^{\mN^\circ},\co_{\bvoX})
\end{equation}
le morphisme de topos annelés défini par les $(u_{n+1})_{n\in \mN}$ et par
\begin{equation}\label{hmdf13e}
\uplambda\colon (X_{s,\zar}^{\mN^\circ},\co_{\bvoX})\rightarrow (X_{s,\zar}, \co_\fX)
\end{equation}
le morphisme de topos annelés pour lequel  le foncteur $\uplambda_*$ est le foncteur limite projective \eqref{notconv13a}. 
On note
\begin{equation}\label{hmdf13f}
\hupsigma\colon (\tE_s^{\mN^\circ},\bvocB)\rightarrow (X_{s,\zar},\co_{\fX})
\end{equation}
le morphisme composé $\uplambda\circ \bvu\circ \bvsigma$. 
On considère aussi les notations analogues pour $f'$, que l'on munit d'un exposant $^\prime$. 
On note
\begin{equation}\label{hmdf13g}
\huppi\colon (\tG_s^{\mN^\circ},\bvocB^!)\rightarrow (X'_{s,\zar},\co_{\fX'})
\end{equation}
le morphisme composé $\uplambda'\circ \bvu'\circ \bvpi$.
D'après \eqref{hmdf14a}, le diagramme de morphismes de topos annelés
\begin{equation}\label{hmdf13i}
\xymatrix{
{(\tE'^{\mN^\circ}_s,\bvocB')}\ar[r]^-(0.5){\bvtau}\ar[rd]_{\hupsigma'}\ar@/^2pc/[rr]^{\bvuptheta}&{(\tG^{\mN^\circ}_s,\bvocB^!)}
\ar[d]^-(0.5){\huppi}\ar[r]^-(0.5){\bvlgg}&{(\tE^{\mN^\circ}_s,\bvocB)}\ar[d]^-(0.5){\hupsigma}\\
&{(X'_{s,\zar},\co_{\fX'})}\ar[r]^{\fgg}&{(X_{s,\zar},\co_{\fX})}}
\end{equation}
est commutatif à isomorphisme canonique près  (\cite{egr1} 1.2.3).

Pour tout $\co_\fX$-module $\cF$ de $X_{s,\zar}$, on a un isomorphisme canonique 
\begin{equation}\label{hmdf13j}
\hupsigma^*(\cF)\stackrel{\sim}{\rightarrow}(\sigma_n^*(u_n^*(\cF/p^{n+1}\cF)))_{n\in \mN}. 
\end{equation}
De même, pour tout $\co_{\fX'}$-module $\cF'$ de $X'_{s,\zar}$, on a un isomorphisme canonique 
\begin{equation}\label{hmdf13k}
\huppi^*(\cF')\stackrel{\sim}{\rightarrow}(\pi_n^*(u'^*_n(\cF'/p^{n+1}\cF')))_{n\in \mN}. 
\end{equation}
En particulier, $\hupsigma^*(\cF)$ et $\huppi^*(\cF')$ sont adiques (\cite{agt} III.7.18). 

On utilise les notations introduites dans \ref{indsh21} et \ref{indsh22} pour les morphismes de topos annelés.

\subsection{}\label{hmdf40}
On reprend les notations de \ref{ahttf43}, \ref{ahttf40} et \ref{ahttf50} pour les morphismes $f$ et $f'$. 
On a un foncteur pleinement fidèle canonique \eqref{ahttf50c}
\begin{equation}\label{hmdf40a}
\bMod^{\coh}(\co_{\fX}[\frac 1 p])\rightarrow \bIndMod(\co_{\fX}).
\end{equation}
On identifie $\bMod^{\coh}(\co_{\fX}[\frac 1 p])$ à une sous-catégorie pleine de $\bIndMod(\co_{\fX})$ par ce foncteur. 
On considère donc tout $\co_{\fX}[\frac 1 p]$-module cohérent aussi comme un ind-$\co_{\fX}$-module. 
Le foncteur $\rI\hupsigma^*$ \eqref{ahttf43c} induit donc un foncteur que l'on note encore
\begin{equation}\label{hmdf40ii}
\rI\hupsigma^*\colon \bMod^\coh(\co_{\fX}[\frac 1 p])\rightarrow  \bIndMod(\bvocB).
\end{equation}
Par ailleurs, le foncteur 
\begin{equation}\label{hmdf40i}
\hupsigma^*_\mQ\colon \bMod_\mQ(\co_{\fX})\rightarrow  \bMod_\mQ(\bvocB)
\end{equation}
induit, compte tenu de \eqref{ahttf50b}, un foncteur qu'on note encore 
\begin{equation}\label{hmdf40j}
\hupsigma^*_\mQ\colon \bMod^\coh(\co_{\fX}[\frac 1 p])\rightarrow  \bMod_\mQ(\bvocB).
\end{equation} 
On adopte les mêmes conventions pour les $\co_{\fX'}[\frac 1 p]$-modules cohérents.

On désigne par $\bMod(\bvocB^!)$ la catégorie des $\bvocB^!$-modules de $\tG^{\mN^\circ}_s$, par 
$\bMod_\mQ(\bvocB^!)$ la catégorie des $\bvocB^!$-modules de $\tG^{\mN^\circ}_s$ à isogénie près \eqref{indsh11}
et par $\bIndMod(\bvocB^!)$ la catégorie des ind-$\bvocB^!$-modules de $\tG^{\mN^\circ}_s$ \eqref{indsh15}. 
On note
\begin{equation}\label{hmdf40b}
\iota_{\bvocB^!}\colon \bMod(\bvocB^!)\rightarrow \bIndMod(\bvocB^!)
\end{equation}
le foncteur canonique, qui est exact et pleinement fidèle \eqref{indsh15}. 
On identifiera $\bMod(\bvocB^!)$ à une sous-catégorie pleine de $\bIndMod(\bvocB^!)$ 
par ce foncteur qu'on omettra des notations. 

D'après \ref{indsh21}, le morphisme de topos annelés $\huppi$ \eqref{hmdf13g} induit deux foncteurs additifs adjoints
\begin{eqnarray}
\rI \huppi^*\colon \bIndMod(\co_{\fX'}) \rightarrow \bIndMod(\bvocB^!),\label{hmdf40c}\\
\rI \huppi_*\colon \bIndMod(\bvocB^!) \rightarrow \bIndMod(\co_{\fX'}).\label{hmdf40d}
\end{eqnarray}
Le foncteur $\rI \huppi^*$ (resp. $\rI \huppi_*$) est exact à droite (resp. gauche). 
Compte tenu de \eqref{ahttf50c}, le foncteur $\rI \huppi^*$ induit un foncteur qu'on note encore 
\begin{equation}\label{hmdf40dd}
\rI \huppi^*\colon \bMod^\coh(\co_{\fX'}[\frac 1 p])\rightarrow  \bIndMod(\bvocB^!).
\end{equation} 
Le foncteur $\rI \huppi_*$ admet un foncteur dérivé à droite
\begin{equation}\label{hmdf40e}
\rR\rI \huppi_*\colon \bD^+(\bIndMod(\bvocB^!))\rightarrow \bD^+(\bIndMod(\co_{\fX'})).
\end{equation}
On désigne par $\vuppi_*$ le foncteur composé \eqref{ahttf43g}
\begin{equation}\label{hmdf40f}
\vuppi_*=\kappa_{\co_{\fX'}}\circ \rI \huppi_*\colon \bIndMod(\bvocB^!) \rightarrow \bMod(\co_{\fX'}).
\end{equation}
Celui-ci est exact à gauche. Il admet un foncteur dérivé à droite
\begin{equation}\label{hmdf40g}
\rR\vuppi_*\colon \bD^+(\bIndMod(\bvocB^!))\rightarrow \bD^+(\bMod(\co_{\fX'})),
\end{equation}
canoniquement isomorphe à $\kappa_{\co_{\fX'}}\circ \rR\rI \huppi_*$ (\cite{ks2} 13.3.13). 

D'après \ref{indsh11}, on a un foncteur canonique \eqref{indsh11c}
\begin{equation}\label{hmdf40l}
\upalpha_{\bvocB^!}\colon \bMod_\mQ(\bvocB^!)\rightarrow \bIndMod(\bvocB^!),
\end{equation}
qui est pleinement fidèle \eqref{indsh5e} et exact \eqref{indsh16}.
On identifiera $\bMod_\mQ(\bvocB^!)$ à une sous-catégorie pleine de $\bIndMod(\bvocB^!)$ 
par ce foncteur qu'on omettra des notations. On considérera donc tout
$\bvocB^!_\mQ$-module comme un ind-$\bvocB^!$-module.

D'après \ref{indsh20}, le morphisme de topos annelés $\huppi$ \eqref{hmdf13g} induit deux foncteurs additifs adjoints 
\begin{eqnarray}
\huppi^*_\mQ\colon \bMod_\mQ(\co_{\fX'}) \rightarrow \bMod_\mQ(\bvocB^!),\label{hmdf40h}\\
\huppi_{\mQ*}\colon \bMod_\mQ(\bvocB^!) \rightarrow \bMod_\mQ(\co_{\fX'}).\label{hmdf40hh}
\end{eqnarray}
Le foncteur $\huppi^*_\mQ$ (resp. $\huppi_{\mQ*}$) est exact à droite (resp. gauche). 
Compte tenu de \eqref{ahttf50b}, le foncteur $\huppi^*_\mQ$ induit un foncteur qu'on note encore 
\begin{equation}\label{hmdf40k}
\huppi^*_\mQ\colon \bMod^\coh(\co_{\fX'}[\frac 1 p])\rightarrow  \bMod_\mQ(\bvocB^!).
\end{equation} 
On note encore  
\begin{equation}\label{hmdf40kk}
\huppi_{\mQ*}\colon \bMod_\mQ(\bvocB^!) \rightarrow \bMod(\co_{\fX'}[\frac 1 p])
\end{equation}
le composé du foncteur $\huppi_{\mQ*}$ \eqref{hmdf40hh} et du foncteur exact canonique \eqref{ahttf40k}
\begin{equation}\label{hmdf40q}
\bMod_\mQ(\co_{\fX'})\rightarrow \bMod(\co_{\fX'}[\frac 1 p]). 
\end{equation} 
Ces abus de notation n'induisent aucune confusion. 

D'après \ref{indsh14}, le foncteur $\huppi_{\mQ*}$ \eqref{hmdf40hh} admet un foncteur dérivé à droite
\begin{equation}\label{hmdf40m}
\rR \huppi_{\mQ*}\colon \bD^+(\bMod_\mQ(\bvocB^!))\rightarrow \bD^+(\bMod_\mQ(\co_{\fX'})).
\end{equation}
En vertu de \eqref{indsh14g}, le diagramme
\begin{equation}\label{hmdf40n}
\xymatrix{
{\bD^+(\bMod_\mQ(\bvocB^!))}\ar[r]^-(0.5){\rR\huppi_{\mQ*}}\ar[d]_{\upalpha_{\bvocB^!}}&{\bD^+(\bMod_\mQ(\co_{\fX'}))}\ar[d]^{\upalpha_{\co_{\fX'}}}\\
{\bD^+(\bIndMod(\bvocB^!))}\ar[r]^-(0.5){\rR\rI \huppi_*}&{\bD^+(\bIndMod(\co_{\fX'}))}}
\end{equation}
est commutatif à un isomorphisme canonique près. 

Le foncteur $\huppi_{\mQ*}$ \eqref{hmdf40kk} admet un foncteur dérivé à droite
\begin{equation}\label{hmdf40o}
\rR \huppi_{\mQ*}\colon \bD^+(\bMod_\mQ(\bvocB^!))\rightarrow \bD^+(\bMod(\co_{\fX'}[\frac 1 p]))
\end{equation}
qui n'est autre que le composé du foncteur $\rR \huppi_{\mQ*}$ \eqref{hmdf40m} et du foncteur \eqref{hmdf40q}. 
Il résulte de \eqref{ahttf40i} que le diagramme
\begin{equation}\label{hmdf40p}
\xymatrix{
{\bD^+(\bMod_\mQ(\bvocB^!))}\ar[r]^-(0.5){\rR\huppi_{\mQ*}}\ar[d]_{\upalpha_{\bvocB^!}}&{\bD^+(\bMod(\co_{\fX'}[\frac 1 p]))}\ar[d]\\
{\bD^+(\bIndMod(\bvocB^!))}\ar[r]^-(0.5){\rR\vuppi_*}&{\bD^+(\bMod(\co_{\fX'}))}}
\end{equation}
où $\rR\vuppi_*$ est le foncteur \eqref{hmdf40g} et la flèche non libellée est le foncteur canonique, est commutatif à un isomorphisme canonique près.

\subsection{}\label{hmdf410}
On reprend les notations de \ref{nrmh} pour les $\cS$-schémas formels $\fX$ et $\fX'$, en particulier, 
$\txi^{-1}\tOmega^1_{\fX/\cS}$ désigne le complété $p$-adique du $\co_\coX$-module \eqref{hmdf1b} (\cite{egr1} 2.5.1)
\begin{equation}\label{hmdf410a}
\txi^{-1}\tOmega^1_{\coX/\coS}=\txi^{-1}\tOmega^1_{X/S}\otimes_{\co_X}\co_{\coX},
\end{equation} 
et $\txi^{-1}\tOmega^1_{\fX'/\cS}$ le complété $p$-adique du $\co_{\coX'}$-module
\begin{equation}\label{hmdf410b}
\txi^{-1}\tOmega^1_{\coX'/\coS}=\txi^{-1}\tOmega^1_{X'/S}\otimes_{\co_{X'}}\co_{\coX'}.
\end{equation} 
On note, de plus, $\txi^{-1}\tOmega^1_{\fX'/\fX}$ le complété $p$-adique du $\co_{\coX'}$-module
\begin{equation}\label{hmdf410c}
\txi^{-1}\tOmega^1_{\coX'/\coX}=\txi^{-1}\tOmega^1_{X'/X}\otimes_{\co_{X'}}\co_{\coX'}.
\end{equation} 
La suite exacte localement scindée \eqref{hmdf6c} induit une suite exacte 
\begin{equation}\label{hmdf410d}
0\rightarrow \fgg^*(\tOmega^1_{\fX/\cS})\rightarrow \tOmega^1_{\fX'/\cS}\rightarrow \tOmega^1_{\fX'/\fX}\rightarrow 0.
\end{equation}

La filtration de Koszul $\rW^\bullet \tOmega^\bullet_{\fX'/\cS}$ de la $\co_{\fX'}$-algèbre extérieure $\tOmega^\bullet_{\fX'/\cS}$ associée à la suite exacte \eqref{hmdf410d}, 
définie dans \eqref{MH90b}, induit pour tout entier $j\geq 0$ une suite exacte \eqref{MH90e}
\begin{equation}\label{hmdf410e}
0\rightarrow \fgg^*(\tOmega^1_{\fX/\cS})\otimes_{\co_{\fX'}} \tOmega^{j-1}_{\fX'/\fX}\rightarrow \rW^0(\tOmega^\bullet_{\fX'/\cS})/\rW^2(\tOmega^\bullet_{\fX'/\cS})
\rightarrow \tOmega^j_{\fX'/\fX}\rightarrow 0.
\end{equation}
Compte tenu de la formule de projection (\cite{sp} \href{https://stacks.math.columbia.edu/tag/0B54}{0B54}), on en déduit un morphisme de $\bD^+(\bMod(\co_\fX))$ 
\begin{equation}\label{hmdf410f}
\rR \fgg_*(\tOmega^j_{\fX'/\fX})\rightarrow \tOmega^1_{\fX/\cS}\otimes_{\co_{\fX}} \rR\fgg_*(\tOmega^{j-1}_{\fX'/\fX})[+1],
\end{equation}
qu'on appelle {\em application de Kodaira-Spencer de $\fgg$} (cf. \cite{katz2} 1.2).

\subsection{}\label{hmdf41}
On note $\bMH(\co_{\fX'},\txi^{-1}\tOmega^1_{\fX'/\fX})$ la catégorie des $\co_{\fX'}$-modules de Higgs à coefficients dans 
$\txi^{-1}\tOmega^1_{\fX'/\fX}$  \eqref{MH1}. On dit qu'un tel module de Higgs est {\em cohérent} si le $\co_{\fX'}$-module sous-jacent est cohérent.
On note $\bMH^\coh(\co_{\fX'}, \txi^{-1}\tOmega^1_{\fX'/\fX})$ la sous-catégorie pleine de $\bMH(\co_{\fX'},\txi^{-1}\tOmega^1_{\fX'/\fX})$
formée des modules de Higgs cohérents.
On sous-entend par {\em $\co_{\fX'}[\frac 1 p]$-module de Higgs à coefficients dans $\txi^{-1}\tOmega^1_{\fX'/\fX}$}, 
un $\co_{\fX'}[\frac 1 p]$-module de Higgs à coefficients dans $\txi^{-1}\tOmega^1_{\fX'/\fX}[\frac 1 p]$.
On dit qu'un tel module de Higgs est {\em cohérent} si le $\co_{\fX'}[\frac 1 p]$-module sous-jacent est cohérent. 
On désigne par $\bMH(\co_{\fX'}[\frac 1 p], \txi^{-1}\tOmega^1_{\fX'/\fX})$ la catégorie des 
$\co_{\fX'}[\frac 1 p]$-modules de Higgs à coefficients dans $\txi^{-1}\tOmega^1_{\fX'/\fX}$
et par $\bMH^\coh(\co_{\fX'}[\frac 1 p], \txi^{-1}\tOmega^1_{\fX'/\fX})$ la sous-catégorie pleine formée des modules de Higgs cohérents. 
On omettra le champ de Higgs de la notation d'un module de Higgs lorsqu'on n'en a pas explicitement besoin.

On désigne par $\bIH(\co_{\fX'},\txi^{-1}\tOmega^1_{\fX'/\fX})$ la catégorie des $\co_{\fX'}$-isogénies  
de Higgs à coefficients dans $\txi^{-1}\tOmega^1_{\fX'/\fX}$ \eqref{indsh23}, par 
$\bIH^\coh(\co_{\fX'},\txi^{-1}\tOmega^1_{\fX'/\fX})$ la sous-catégorie pleine 
formée des quadruplets $(\cM,\cN,u,\theta)$ tels que $\cM$ et $\cN$ soient des $\co_{\fX'}$-modules cohérents
et par $\bIH_\mQ(\co_{\fX'},\txi^{-1}\tOmega^1_{\fX'/\fX})$ (resp. $\bIH^\coh_\mQ(\co_{\fX'},\txi^{-1}\tOmega^1_{\fX'/\fX})$) 
la catégorie des objets de $\bIH(\co_{\fX'},\txi^{-1}\tOmega^1_{\fX'/\fX})$  (resp. $\bIH^\coh(\co_{\fX'},\txi^{-1}\tOmega^1_{\fX'/\fX})$) 
à isogénie près (\cite{agt} III.6.1.1).
D'après (\cite{agt} III.6.20), on a une équivalence canonique de catégories 
\begin{equation}\label{hmdf41a}
\bIH^\coh_\mQ(\co_{\fX'},\txi^{-1}\tOmega^1_{\fX'/\fX})\stackrel{\sim}{\rightarrow} \bMH^\coh(\co_{\fX'}[\frac 1 p], \txi^{-1}\tOmega^1_{\fX'/\fX}).
\end{equation}

On note $\bIndMH(\co_{\fX'},\txi^{-1}\tOmega^1_{\fX'/\cS})$ (resp. $\bIndMH(\co_{\fX'},\txi^{-1}\tOmega^1_{\fX'/\fX})$) 
la catégorie des ind-$\co_{\fX'}$-modules de Higgs à coefficients dans $\txi^{-1}\tOmega^1_{\fX'/\cS}$ (resp. $\txi^{-1}\tOmega^1_{\fX'/\fX}$) \eqref{indsh30}.
On rappelle qu'on a un foncteur pleinement fidèle \eqref{indmdlb2f}
\begin{equation}\label{hmdf41c} 
\bMH^\coh(\co_{\fX'}[\frac 1 p], \txi^{-1}\tOmega^1_{\fX'/\cS})\rightarrow \bIndMH(\co_{\fX'},\txi^{-1}\tOmega^1_{\fX'/\cS}).
\end{equation}
De même, on a un foncteur pleinement fidèle
\begin{equation}\label{hmdf41b} 
\bMH^\coh(\co_{\fX'}[\frac 1 p], \txi^{-1}\tOmega^1_{\fX'/\fX})\rightarrow \bIndMH(\co_{\fX'},\txi^{-1}\tOmega^1_{\fX'/\fX}).
\end{equation}
On identifiera $\bMH^\coh(\co_{\fX'}[\frac 1 p], \txi^{-1}\tOmega^1_{\fX'/\cS})$ 
(resp. $\bMH^\coh(\co_{\fX'}[\frac 1 p], \txi^{-1}\tOmega^1_{\fX'/\fX})$)
à une sous-catégorie pleine de $\bIndMH(\co_{\fX'},\txi^{-1}\tOmega^1_{\fX'/\cS})$
(resp. $\bIndMH(\co_{\fX'},\txi^{-1}\tOmega^1_{\fX'/\fX})$) par ces foncteurs qu'on omettra des notations.

\begin{defi}\label{hmdf42}
On appelle {\em $\co_{\fX'}[\frac 1 p]$-fibré de Higgs à coefficients dans $\txi^{-1}\tOmega^1_{\fX'/\fX}$} 
tout $\co_{\fX'}[\frac 1 p]$-module de Higgs à coefficients dans $\txi^{-1}\tOmega^1_{\fX'/\fX}$ 
dont le $\co_{\fX'}[\frac 1 p]$-module sous-jacent est localement projectif de type fini \eqref{notconv14}. 
\end{defi}

\begin{prop}\label{hmdf43}
Si le morphisme $g\colon X'\rightarrow X$ est propre, pour tous entiers $i,j\geq 0$, on a les propriétés suivantes. 
\begin{itemize}
\item[{\rm (i)}] Le $\co_{X_\eta}$-module $\rR^ig_*(\tOmega^j_{X'/X})|X_\eta$ est localement libre de type fini.
\item[{\rm (ii)}] On a un isomorphisme canonique 
\begin{equation}\label{hmdf43a}
\rR^i\fgg_*(\tOmega^j_{\fX'/\fX})\stackrel{\sim}{\rightarrow} a^{-1}(\rR^ig_*(\tOmega^j_{X'/X}))\otimes_{a^{-1}(\co_X)}\co_\fX,
\end{equation} 
où $a\colon X_s\rightarrow X$ est l'injection canonique.
\item[{\rm (iii)}] Le $\co_\fX[\frac 1 p]$-module $\rR^i\fgg_*(\tOmega^j_{\fX'/\fX})\otimes_{\co_\fX}\co_\fX[\frac 1 p]$ est localement projectif de type fini \eqref{notconv14}.
\end{itemize}
\end{prop}

(i) On notera d'abord que le morphisme $g$ étant saturé, est exact (\cite{ogus} III 2.5.2), et que les fibres du monoïde $(\cM_X/\co^\times_X)|X_\eta$
sont libres puisque $\cM_X|X_\eta$ est défini par un diviseur à croisements normaux (\cite{agt} III.4.7). 
La proposition résulte alors de (\cite{ikn} 7.2) (cf. \cite{deligne1} 5.5 pour le cas lisse sans structures logarithmiques).

(ii) Cela résulte de (\cite{egr1} 2.5.5(ii) et 2.12.2).

(iii) La question étant locale, on peut supposer $X$ affine. Posons $\cF=\rR^ig_*(\tOmega^j_{X'/X})$.  
D'après (i), il existe un $\co_X$-module cohérent $\cG$, un entier $n\geq 1$ et un morphisme $\co_X$-linéaire $u\colon \cF\oplus\cG\rightarrow \co_X^n$
induisant un isomorphisme sur $X_\eta$. Il existe un entier $m\geq 0$ tel que $p^m$ annule le noyau et le conoyau de $u$. 
D'après (\cite{ag} 2.6.3), il existe alors un morphisme $\co_X$-linéaire $v\colon \co_X^n\rightarrow \cF\oplus\cG$
tel que $u\circ v= p^{2m}\id_{\co_X^n}$ et $v\circ u=p^{2m}\id_{\cF\oplus\cG}$, autrement dit, $u$ est une isogénie \eqref{caip1}. La proposition s'ensuit compte tenu de (ii).

\section{Changement de base}\label{chb}

\subsection{}\label{chb1}
Soient $\cM$ un $\co_{\fX'}$-module, $\cN$ un $\bvocB'$-module, $q$ un entier $\geq 0$. 
Le morphisme d'adjonction $\huppi^*(\cM)\rightarrow \bvtau_*(\hupsigma'^*(\cM))$ \eqref{hmdf13i} et le cup-produit induisent un morphisme bifonctoriel
\begin{equation}\label{chb1a}
\huppi^*(\cM)\otimes_{\bvocB^!}\rR^q\bvtau_*(\cN)\rightarrow \rR^q\bvtau_*(\hupsigma'^*(\cM)\otimes_{\bvocB'}\cN). 
\end{equation}
Celui-ci induit pour tous objets $\cF$ de $\bMod_\mQ(\co_{\fX'})$ et $\cG$ de $\bMod_\mQ(\bvocB')$, un morphisme bifonctoriel de $\bMod_\mQ(\bvocB^!)$
\begin{equation}\label{chb1b}
\huppi^*_\mQ(\cF)\otimes_{\bvocB^!_\mQ}\rR^q\bvtau_{\mQ*}(\cG)\rightarrow \rR^q\bvtau_{\mQ*}(\hupsigma'^*_\mQ(\cF)\otimes_{\bvocB'_\mQ}\cG). 
\end{equation}
Compte tenu de \eqref{indsh21g}, le morphisme \eqref{chb1a} induit aussi pour tous objets $\cF$ de $\bIndMod(\co_{\fX'})$ et $\cG$ de $\bIndMod(\bvocB')$, 
un morphisme bifonctoriel de $\bIndMod(\bvocB^!)$
\begin{equation}\label{chb1c}
\rI\huppi^*(\cF)\otimes_{\bvocB^!}\rR^q\rI\bvtau_*(\cH)\rightarrow \rR^q\rI\bvtau_*(\rI\hupsigma'^*(\cF)\otimes_{\bvocB'}\cH). 
\end{equation}

Compte tenu de \eqref{aspglob1d} (resp. \eqref{hmdf40a}), le morphisme \eqref{chb1b} (resp. \eqref{chb1c}) existe lorsque $\cF$ est un $\co_{\fX'}[\frac 1 p]$-module cohérent.

\begin{lem}\label{chb2}
Soit $q$ un entier $\geq 0$. Alors,
\begin{itemize}
\item[{\rm (i)}] Pour tout $\co_{\fX'}$-module localement libre de type fini $\cM$ et tout $\bvocB'$-module $\cN$, le morphisme canonique \eqref{chb1a}
\begin{equation}
\huppi^*(\cM)\otimes_{\bvocB^!}\rR^q\bvtau_*(\cN)\rightarrow\rR^q\bvtau_*(\hupsigma'^*(\cM)\otimes_{\bvocB'}\cN)
\end{equation}
est un isomorphisme.
\item[{\rm (ii)}] Pour tout $\co_{\fX'}[\frac 1 p]$-module localement projectif de type fini $\cF$ \eqref{notconv14}
et tout $\bvocB'_\mQ$-module $\cG$, le morphisme canonique \eqref{chb1b} 
\begin{equation}
\huppi^*_\mQ(\cF)\otimes_{\bvocB^!_\mQ}\rR^q\bvtau_{\mQ*}(\cG)\rightarrow\rR^q\bvtau_{\mQ*}(\hupsigma'^*_\mQ(\cF)\otimes_{\bvocB'_\mQ}\cG)
\end{equation}
est un isomorphisme.
\item[{\rm (iii)}] Pour tout $\co_{\fX'}$-module localement libre de type fini $\cF$ et tout ind-$\bvocB'$-module $\cG$, le morphisme canonique \eqref{chb1c} 
\begin{equation}
\huppi^*(\cF)\otimes_{\bvocB^!}\rR^q\rI\bvtau_*(\cG)\rightarrow\rR^q\rI\bvtau_*(\hupsigma'^*(\cF)\otimes_{\bvocB'}\cG)
\end{equation}
est un isomorphisme.
\item[{\rm (iv)}] Pour tout $\co_{\fX'}[\frac 1 p]$-module localement projectif de type fini $\cF$ \eqref{notconv14}
et tout ind-$\bvocB'$-module $\cG$, le morphisme canonique \eqref{chb1c} 
\begin{equation}
\rI\huppi^*(\cF)\otimes_{\bvocB^!}\rR^q\rI\bvtau_*(\cG)\rightarrow\rR^q\rI\bvtau_*(\rI\hupsigma'^*(\cF)\otimes_{\bvocB'}\cG)
\end{equation}
est un isomorphisme.
\end{itemize}
\end{lem}

(i) C'est immédiat.  

(ii) Il existe un recouvrement fini $(U'_i)_{i\in I}$ par des ouverts de Zariski de $X'$ tel que pour tout $i\in I$, la restriction 
de $\cF$ à $(U'_i)_s$ soit un facteur direct d'un $(\co_{\fX'}|(U'_i)_s)[\frac 1 p]$-module libre de type fini. 
Pour tout $i\in I$, notons $\bvpi^*((U'_i)_s)$ l'objet constant de $\tG^{\mN^\circ}_s$ de valeur $\pi^*_s((U'_i)_s)$ \eqref{hmdf18b}. 
Les objets $(\bvpi^*((U'_i)_s))_{i\in I}$ forment alors un recouvrement de l'objet final de $\tG^{\mN^\circ}_s$ (cf. \cite{agt} III.7.4 et III.7.5). 
On laissera au lecteur le soin d'écrire l'énoncé analogue pour $\tE'^{\mN^\circ}_s$. 
Compte tenu de (\cite{agt} III.6.7(ii)), on peut alors se borner au cas où $\cF$ est un facteur direct
d'un $\co_{\fX'}[\frac 1 p]$-module libre de type fini, et même au cas où $\cF$ est un 
$\co_{\fX'}[\frac 1 p]$-module libre de type fini, auquel cas l'assertion est évidente. 

(iii) Cela résulte de \ref{indsh47}. On peut aussi le déduire de (i), compte tenu de \eqref{indsh21g} et du fait que le produit tensoriel et 
le foncteur $\rR^q\rI\bvtau_*$ commutent aux petites limites inductives filtrantes d'après \ref{indsh18} et \ref{indsh4}(ii).

(iv) Compte tenu de \ref{indsh39}(iii), on peut se réduire au cas où $\cF$ est un $\co_{\fX'}[\frac 1 p]$-module libre de type fini, cf. la preuve de (ii). 
L'assertion résulte alors du fait que le produit tensoriel et le foncteur $\rR^q\rI\bvtau_*$ commutent aux petites limites inductives filtrantes 
d'après \ref{indsh18}, \eqref{indsh21g} et \ref{indsh4}(ii).

\begin{lem}\label{chb3}
Il existe un entier $N\geq 0$ tel que pour tout entier $n\geq 0$, les propriétés suivantes soient remplies.
\begin{itemize}
\item[{\rm (i)}] Pour toute suite exacte de $\co_{\oX_n}$-module de $X_{s,\et}$, 
$0\rightarrow \cF'\rightarrow \cF\rightarrow \cF''\rightarrow 0$, la suite  
\begin{equation}\label{chb3a}
0\rightarrow \sigma^*_n(\cF')\rightarrow \sigma^*_n(\cF)\rightarrow \sigma^*_n(\cF'')\rightarrow 0
\end{equation}
est $p^N$-exacte \eqref{notconv17}; 
\item[{\rm (ii)}] Pour toute suite exacte de $\co_{\oX'_n}$-module de $X'_{s,\et}$, 
$0\rightarrow \cF'\rightarrow \cF\rightarrow \cF''\rightarrow 0$, la suite  
\begin{equation}\label{chb3b}
0\rightarrow \pi^*_n(\cF')\rightarrow \pi^*_n(\cF)\rightarrow \pi^*_n(\cF'')\rightarrow 0
\end{equation}
est $p^N$-exacte.
\end{itemize}
\end{lem}

Comme $X$ est quasi-compact, on peut supposer que le morphisme $f$ \eqref{hmdf3} admet une carte adéquate (\cite{agt} III.4.7). 
Montrons que si $N$ désigne l'entier fourni par la proposition (\cite{ag} 5.3.9), l'entier $N+1$ convient. 

(i) D'après (\cite{ag} 5.3.9), pour tout point $(\oy\rightsquigarrow \ox)$ de $X_\et\gtimes_{X_\et}\oX^\circ_\et$ \eqref{TFA9} 
tel que $\ox$ soit au-dessus de $s$, 
et toute suite exacte de $\co_{\oX,\ox}$-modules $0\rightarrow M'\rightarrow M\rightarrow M''\rightarrow 0$, la suite 
\begin{equation}
0\rightarrow M'\otimes_{\co_{\oX,\ox}}\ocB_{\rho(\oy\rightsquigarrow \ox)}\rightarrow M\otimes_{\co_{\oX,\ox}}\ocB_{\rho(\oy\rightsquigarrow \ox)}
\rightarrow M''\otimes_{\co_{\oX,\ox}}\ocB_{\rho(\oy\rightsquigarrow \ox)}\rightarrow 0
\end{equation}
est $p^N$-exacte. On notera que $k$ étant algébriquement clos \eqref{hmdf1}, $\ox$ est naturellement un point géométrique de $\oX$.
La propriété (i) s'ensuit compte tenu de (\cite{agt} (VI.10.18.1) et III.9.5).

(ii) Soit $(\oy\rightsquigarrow \ox')$ un point de $X'_\et\gtimes_{X_\et}\oX^\circ_\et$ (\cite{ag} 3.4.23) tel que $\ox'$ soit au-dessus de $s$. 
Posons $\ox=g(\ox')$ et notons $(\oy\rightsquigarrow \ox)$ l'image de $(\oy\rightsquigarrow \ox')$ 
par le morphisme canonique
\begin{equation}
X'_\et\gtimes_{X_\et}\oX^\circ_\et \rightarrow X_\et\gtimes_{X_\et}\oX^\circ_\et.
\end{equation}
D'après (\cite{ag} 6.5.29), pour tout entier $n\geq 0$, l'homomorphisme canonique 
\begin{equation}
(\ocB_{\rho(\oy\rightsquigarrow \ox)}/p^n\ocB_{\rho(\oy\rightsquigarrow \ox)})\otimes_{\co_{\oX,\ox}} \co_{\oX',\ox'} 
\rightarrow \ocB^!_{\varrho(\oy\rightsquigarrow \ox')}/p^n \ocB^!_{\varrho(\oy\rightsquigarrow \ox')}
\end{equation}
est un $\alpha$-isomorphisme \eqref{hmdf17g}. Il résulte alors de (\cite{ag} 5.3.9) que pour toute
suite exacte de $(\co_{\oX',\ox'}/p^n\co_{\oX',\ox'})$-modules $0\rightarrow M'\rightarrow M\rightarrow M''\rightarrow 0$,  la suite 
\begin{equation}
0\rightarrow M'\otimes_{\co_{\oX',\ox'}}\ocB^!_{\varrho(\oy\rightsquigarrow \ox')}\rightarrow M\otimes_{\co_{\oX',\ox'}}\ocB^!_{\varrho(\oy\rightsquigarrow \ox')}
\rightarrow M''\otimes_{\co_{\oX',\ox'}}\ocB^!_{\varrho(\oy\rightsquigarrow \ox')}\rightarrow 0
\end{equation}
est $p^{N+1}$-exacte. La propriété (ii) s'ensuit compte tenu de (\cite{ag} (3.4.23.1) et 6.5.3).

\begin{lem}\label{chb4}
Il existe un entier $N\geq 0$ vérifiant les propriétés suivantes:
\begin{itemize}
\item[{\rm (i)}] pour toute suite exacte de $\co_\fX$-modules $0\rightarrow \cF'\rightarrow \cF\rightarrow \cF''\rightarrow 0$ où $\cF''$ est $\cS$-plat, la suite de $\bvocB$-modules
\begin{equation}\label{chb4a}
0\rightarrow \hupsigma^*(\cF')\rightarrow \hupsigma^*(\cF)\rightarrow \hupsigma^*(\cF'')\rightarrow 0
\end{equation}
est $p^N$-exacte \eqref{notconv17}; 
\item[{\rm (ii)}] pour toute suite exacte de $\co_{\fX'}$-modules $0\rightarrow \cF'\rightarrow \cF\rightarrow \cF''\rightarrow 0$ où $\cF''$ est $\cS$-plat, la suite 
\begin{equation}\label{chb4b}
0\rightarrow \huppi^*(\cF')\rightarrow \huppi^*(\cF)\rightarrow \huppi^*(\cF'')\rightarrow 0
\end{equation}
est $p^N$-exacte. 
\end{itemize}
\end{lem}

Montrons que l'entier $N$ fourni par la proposition \ref{chb3} convient. 
Soit $0\rightarrow \cF'\rightarrow \cF\rightarrow \cF''\rightarrow 0$ une suite exacte de $\co_\fX$-modules telle que $\cF''$ soit $\cS$-plat. Pour tout 
entier $n\geq 0$, la suite de $\co_{\oX_n}$-modules de $X_{s,\zar}$ 
\begin{equation}
0\rightarrow \cF'/p^n\cF'\rightarrow \cF/p^n\cF\rightarrow \cF''/p^n\cF''\rightarrow 0
\end{equation}
est exacte. Comme $u_n$ est plat \eqref{hmdf13c}, la suite de $\ocB_n$-modules
\begin{equation}
0\rightarrow \sigma_n^*(u_n^*(\cF'/p^n\cF'))\rightarrow \sigma_n^*(u_n^*(\cF/p^n\cF))\rightarrow \sigma_n^*(u_n^*(\cF''/p^n\cF''))\rightarrow 0
\end{equation}
est $p^N$-exacte en vertu de \ref{chb3}. La propriété (i) s'ensuit compte tenu de \eqref{hmdf13j} et (\cite{agt} III.7.3(i)).  
La propriété (ii) se démontre de même.

\begin{prop}\label{chb5}
Les foncteurs \eqref{hmdf40j} et \eqref{hmdf40k}
\begin{eqnarray}
\hupsigma^*_\mQ\colon \bMod^\coh(\co_\fX[\frac 1 p])&\rightarrow &\bMod_\mQ(\bvocB),\label{chb5a}\\
\hupsigma'^*_\mQ\colon \bMod^\coh(\co_{\fX'}[\frac 1 p])&\rightarrow &\bMod_\mQ(\bvocB'),\label{chb5b}\\
\huppi^*_\mQ\colon \bMod^\coh(\co_{\fX'}[\frac 1 p])&\rightarrow &\bMod_\mQ(\bvocB^!),\label{chb5c}
\end{eqnarray}
sont exacts.
\end{prop}

En effet, d'après (\cite{agt} III.6.16 et III.6.1.4), toute suite exacte de $\co_\fX[\frac 1 p]$-modules cohérents 
$0\rightarrow \cG'\rightarrow \cG\rightarrow \cG''\rightarrow 0$ s'obtient à partir d'une suite exacte de $\co_\fX$-modules cohérents
$0\rightarrow \cF'\rightarrow \cF\rightarrow \cF''\rightarrow 0$ en inversant $p$. Comme le noyau $\cF''_\tor$ du morphisme canonique 
$\cF''\rightarrow \cF''[\frac 1 p]$ est un $\co_\fX$-module cohérent (\cite{egr1} 2.10.14), 
remplaçant $\cF''$ par $\cF''/\cF''_\tor$, on se réduit au cas où $\cF''$ est $\cS$-plat. Il résulte alors de \ref{chb4} que la suite 
\begin{equation}\label{chb5d}
0\rightarrow \hupsigma^*_\mQ(\cG')\rightarrow \hupsigma^*_\mQ(\cG)\rightarrow \hupsigma^*_\mQ(\cG'')\rightarrow 0
\end{equation}
est exacte. On démontre de même que les foncteurs \eqref{chb5b} et \eqref{chb5c} est exact.

\begin{cor}\label{chb50}
Les foncteurs \eqref{hmdf40ii} et \eqref{hmdf40dd}
\begin{eqnarray}
\rI\hupsigma^*\colon \bMod^\coh(\co_\fX[\frac 1 p])&\rightarrow &\bIndMod(\bvocB),\label{chb50a}\\
\rI\hupsigma'^*\colon \bMod^\coh(\co_{\fX'}[\frac 1 p])&\rightarrow &\bIndMod(\bvocB'),\label{chb50b}\\
\rI\huppi^*\colon \bMod^\coh(\co_{\fX'}[\frac 1 p])&\rightarrow &\bIndMod(\bvocB^!),\label{chb50c}
\end{eqnarray}
sont exacts.
\end{cor}

Cela résulte de \ref{indsh16}(ii) et \ref{chb5}.

\subsection{}\label{chb6}
Pour tout entier $n\geq 0$, le diagramme de morphismes de topos annelés \eqref{hmdf19l} 
\begin{equation}\label{chb6a}
\xymatrix{
{(\tG_s,\ocB_n^!)}\ar[r]^-(0.5){\pi_n}\ar[d]_-(0.5){\lgg_n}&{(X'_{s,\et},\co_{\oX'_n})}\ar[d]^{\ogg_n}\\
{(\tE_s,\ocB_n)}\ar[r]^-(0.5){\sigma_n}&{(X_{s,\et},\co_{\oX_n})}}
\end{equation}
est commutatif à isomorphisme canonique près.  
Pour tout $\co_{\oX'_n}$-module $\cF'$ de $X'_{s,\et}$ et tout entier $q\geq 0$, on a un morphisme canonique  de changement de base (\cite{egr1} (1.2.3.3))
\begin{equation}\label{chb6b}
\sigma_n^*(\rR^q\ogg_{n*}(\cF'))\rightarrow \rR^q\lgg_{n*}(\pi_n^*(\cF')),
\end{equation} 
où $\sigma_n^*$ et $\pi_n^*$ désignent les images inverses au sens des topos annelés. 

\begin{teo}[\cite{ag} 6.5.31]\label{chb7}
Supposons le morphisme $g\colon X'\rightarrow X$ propre. 
Il existe alors un entier $N\geq 0$ tel que pour tous entiers $n\geq 1$ et $q\geq 0$ et 
tout $\co_{\oX'_n}$-module quasi-cohérent $\cF'$ de $X'_{s,\zar}$, 
que l'on considère aussi comme un $\co_{\oX'_n}$-module de $X'_{s,\et}$ \eqref{notconv12}, 
le noyau et le conoyau du morphisme de changement de base \eqref{chb6b}
\begin{equation}\label{chb7a}
\sigma_n^*(\rR^q\ogg_{n*}(\cF'))\rightarrow \rR^q\lgg_{n*}(\pi_n^*(\cF'))
\end{equation}
soient annulés par $p^N$.
\end{teo}

\subsection{}\label{chb8}
Le diagramme de morphismes de topos annelés \eqref{hmdf14a}
\begin{equation}\label{chb8a}
\xymatrix{
{(\tG^{\mN^\circ}_s,\bvocB^!)}\ar[r]^-(0.5){\bvpi}\ar[d]_-(0.5){\bvlgg}&{(X'^{\mN^\circ}_{s,\et},\co_{\bvoX'})}\ar[d]^{\bvogg}\\
{(\tE^{\mN^\circ}_s,\bvocB)}\ar[r]^-(0.5){\bvsigma}&{(X^{\mN^\circ}_{s,\et},\co_{\bvoX})}}
\end{equation}
est commutatif à isomorphisme canonique près.  
Pour tout $\co_{\bvoX'}$-module $\cF'$ de $X'^{\mN^\circ}_{s,\et}$ et tout entier $q\geq 0$, 
on a un morphisme canonique de changement de base (\cite{egr1} (1.2.3.3))
\begin{equation}\label{chb8b}
\bvsigma^*(\rR^q\bvogg_*(\cF'))\rightarrow \rR^q\bvlgg_*(\bvpi^*(\cF')),
\end{equation}
où $\bvsigma^*$ et $\bvpi^*$ désignent les images inverses au sens des topos annelés.

\begin{lem}\label{chb9}
Pour tout $\co_{\bvoX'}$-module $\cF'=(\cF'_n)_{n\geq 0}$ de $X'^{\mN^\circ}_{s,\et}$ 
et tout entier $q\geq 0$, le morphisme de changement de base \eqref{chb8b}
\begin{equation}\label{chb9a}
\bvsigma^*(\rR^q\bvogg_*(\cF'))\rightarrow \rR^q\bvlgg_*(\bvpi^*(\cF'))
\end{equation}
est induit par les morphismes de changement de base \eqref{chb6b}, pour tous les entiers $n\geq 0$,
\begin{equation}\label{chb9b}
\sigma_n^*(\rR^q\ogg_{n*}(\cF'_n))\rightarrow \rR^q\lgg_{n*}(\pi_n^*(\cF'_n)). 
\end{equation}
\end{lem}
En effet, d'après (\cite{agt} III.7.1), pour tout entier $n\geq 0$, il existe deux morphismes de topos 
$a_n\colon \tE_s\rightarrow \tE^{\mN^\circ}_s$ et $b_n\colon \tG_s\rightarrow \tG^{\mN^\circ}_s$ tels que pour tous objets 
$M=(M_n)_{n\geq 0}$ de $\tE^{\mN^\circ}_s$ et $N=(N_n)_{n\geq 0}$ de $\tG^{\mN^\circ}_s$, on ait $a_n^*(M)=M_n$ et $b_n^*(G)=G_n$. 
D'après (\cite{agt} (III.7.5.4)), le diagramme de morphisme de topos
\begin{equation}
\xymatrix{
{\tG_s}\ar[r]^{b_n}\ar[d]_{\lgg_s}&{\tG^{\mN^\circ}_s}\ar[d]^{\bvlgg}\\
{\tE_s}\ar[r]^{a_n}&{\tE^{\mN^\circ}_s}}
\end{equation}
est commutatif à isomorphisme canonique prêt. Par ailleurs, pour tout $\bvocB^!$-module $N$ de $\tG^{\mN^\circ}_s$, le morphisme de changement de base 
\begin{equation}
a_n^*(\rR^q\bvlgg_*(N))\rightarrow\rR^q\lgg_{n*}(b_n^*(N))
\end{equation}
est un isomorphisme (\cite{agt} (III.7.5.5)). La proposition résulte alors de (\cite{egr1} 1.2.4(ii)).

\begin{prop}\label{chb10}
Supposons le morphisme $g\colon X'\rightarrow X$ propre. 
Il existe alors un entier $N\geq 0$ tel que pour tout $\co_{\bvoX'}$-module $\cF'=(\cF'_n)_{n\geq 0}$ de $X'^{\mN^\circ}_{s,\et}$,
où les $\co_{\oX'_n}$-modules $\cF'_n$ sont induits par des $\co_{\oX'_n}$-modules quasi-cohérents de $X'_{s,\zar}$ \eqref{notconv12a},  
et tout entier $q\geq 0$, le noyau et le conoyau du morphisme de changement de base \eqref{chb8b}
\begin{equation}\label{chb10a}
\bvsigma^*(\rR^q\bvogg_*(\cF'))\rightarrow \rR^q\bvlgg_*(\bvpi^*(\cF'))
\end{equation}
sont annulés par $p^N$.
\end{prop}
Cela résulte de \ref{chb7}, \ref{chb9} et (\cite{agt} III.7.3(i)).

\subsection{}\label{chb11}
On vérifie aussitôt que les deux carrés du diagramme de morphismes de topos annelés
\begin{equation}\label{chb11a}
\xymatrix{
{(X'^{\mN^\circ}_{s,\et},\co_{\bvoX'})}\ar[r]^-(0.5){\bvu'}\ar[d]_{\bvogg}
&{(X'^{\mN^\circ}_{s,\zar},\co_{\bvoX'})}\ar[r]^-(0.5){\uplambda'}\ar[d]_{\bvogg}&{(X'_{s,\zar},\co_{\fX'})}\ar[d]^{\fgg}\\
{(X^{\mN^\circ}_{s,\et},\co_{\bvoX})}\ar[r]^-(0.5)\bvu&
{(X^{\mN^\circ}_{s,\zar},\co_{\bvoX})}\ar[r]^-(0.5)\uplambda&{(X_{s,\zar},\co_{\fX})}}
\end{equation}
où $\uplambda$ et $\uplambda'$ (resp. $\bvu$ et $\bvu'$) sont les morphismes de topos annelés canoniques \eqref{notconv13a} 
(resp. \eqref{hmdf13d}), sont commutatifs à isomorphismes canoniques près. 

\begin{lem}[\cite{ag} 6.5.37]\label{chb12}
Supposons le morphisme $g\colon X'\rightarrow X$ séparé et quasi-compact. 
Soient $\cF'=(\cF'_n)_{n\in \mN}$ un $\co_{\bvoX'}$-module de $X'^{\mN^\circ}_{s,\zar}$ tel que pour tout entier $n\geq 0$, 
le $\co_{\oX'_n}$-module $\cF'_n$ soit quasi-cohérent, 
$q$ un entier $\geq 0$. Alors, le morphisme de changement de base relativement au carré de gauche du diagramme \eqref{chb11a}
\begin{equation}\label{chb12a}
\bvu^*(\rR^q\bvogg_{\zar*}(\cF'))\rightarrow \rR^q\bvogg_{\et*}(\bvu'^*(\cF'))
\end{equation}
est un isomorphisme.
\end{lem}

\begin{prop}[\cite{ag} 6.5.38]\label{chb13}
Supposons le morphisme $g\colon X'\rightarrow X$ propre. Soit $\cF'$ un $\co_{\fX'}$-module cohérent de $X'_{s,\zar}$, $q$ un entier $\geq 0$. 
Alors, il existe un entier $N\geq 0$ tel que le noyau et le conoyau du morphisme de changement de base relativement au carré 
de droite du diagramme \eqref{chb11a}
\begin{equation}\label{chb13a}
\uplambda^*(\rR^q\fgg_*(\cF'))\rightarrow \rR^q\bvogg_{\zar*}(\uplambda'^*(\cF'))
\end{equation}
soient annulés par $p^N$.
\end{prop}

\subsection{}\label{chb14}
Le diagramme de morphismes de topos annelés \eqref{hmdf13i} 
\begin{equation}\label{chb14a}
\xymatrix{
{(\tG^{\mN^\circ}_s,\bvocB^!)}\ar[r]^-(0.5){\huppi}\ar[d]_-(0.5){\bvlgg}&{(X'_{s,\zar},\co_{\fX'})}\ar[d]^{\fgg}\\
{(\tE^{\mN^\circ}_s,\bvocB)}\ar[r]^-(0.5){\hupsigma}&{(X_{s,\zar},\co_{\fX})}}
\end{equation}
est commutatif à isomorphisme canonique près.  
Pour tout $\co_{\fX'}$-module $\cF'$ de $X'_{s,\zar}$ et tout entier $q\geq 0$, on a un morphisme canonique de changement de base (\cite{egr1} (1.2.3.3))
\begin{equation}\label{chb14b}
\hupsigma^*(\rR^q\fgg_*(\cF'))\rightarrow \rR^q\bvlgg_*(\huppi^*(\cF')),
\end{equation}
où $\hupsigma^*$ et $\huppi^*$ désignent les images inverses au sens des topos annelés. 
Compte tenu de \eqref{indsh14d}, celui-ci induit pour tout objet 
$F'$ de $\bMod_\mQ(\co_{\fX'})$ et tout entier $q\geq 0$, un morphisme canonique, que l'on appelle aussi morphisme de changement de base,
\begin{equation}\label{chb14c}
\hupsigma^*_\mQ(\rR^q\fgg_{\mQ*}(F'))\rightarrow \rR^q\bvlgg_{\mQ*}(\huppi^*_\mQ(F')). 
\end{equation}

Pour tout $\co_{\fX'}[\frac 1 p]$-module cohérent $\cF'$ de $X'_{s,\zar}$ et tout entier $q\geq 0$
tels que le $\co_\fX[\frac 1 p]$-module $\rR^q\fgg_*(\cF')$ soit cohérent, 
le morphisme \eqref{chb14c} induit un morphisme 
canonique, que l'on appelle aussi morphisme de changement de base,
\begin{equation}\label{chb14d}
\hupsigma^*_\mQ(\rR^q\fgg_*(\cF'))\rightarrow \rR^q\bvlgg_{\mQ*}(\huppi^*_\mQ(\cF')),
\end{equation}
où $\hupsigma^*_\mQ$ et $\huppi^*_\mQ$ désignent les foncteurs \eqref{hmdf40j} et \eqref{hmdf40k}.

\begin{prop}\label{chb15}
Supposons le morphisme $g\colon X'\rightarrow X$ propre. 
Il existe alors un entier $N\geq 0$ tel que pour $\co_{\fX'}$-module cohérent  $\cF'$ et tout entier $q\geq 0$, 
le noyau et le conoyau du morphisme de changement de base \eqref{chb14b}
\begin{equation}\label{chb15a}
\hupsigma^*(\rR^q\fgg_*(\cF'))\rightarrow \rR^q\bvlgg_*(\huppi^*(\cF'))
\end{equation}
soient annulés par $p^N$.
\end{prop}

En effet, on a $\huppi=\bvpi\circ\bvu'\circ \uplambda'$ et $\hupsigma=\bvsigma\circ\bvu\circ \uplambda$. 
La proposition résulte alors de \ref{chb10}, \ref{chb12} et \ref{chb13}, compte tenu de \eqref{notconv12e}, (\cite{ag} 2.6.3) 
et (\cite{egr1} 1.2.4(ii)). 
On notera que le faisceau d'anneaux $\co_{\fX'}$ de $X'_{s,\zar}$  est cohérent (\cite{egr1} 2.8.1).

\begin{cor}\label{chb16}
Supposons le morphisme $g\colon X'\rightarrow X$ propre. 
Soient $\cF'$ un $\co_{\fX'}[\frac 1 p]$-module cohérent, $q$ un entier $\geq 0$. 
Alors, le $\co_\fX[\frac 1 p]$-module $\rR^q\fgg_*(\cF')$ est cohérent et le morphisme de changement de base \eqref{chb14d}
\begin{equation}\label{chb16a}
\hupsigma^*_\mQ(\rR^q\fgg_*(\cF'))\rightarrow \rR^q\bvlgg_{\mQ*}(\huppi^*_\mQ(\cF'))
\end{equation}
est un isomorphisme.
\end{cor}

En effet, la première assertion résulte de (\cite{egr1} 2.10.24 et 2.11.5). 
La seconde assertion résulte alors de \ref{chb15}. On notera que le faisceau d'anneaux $\co_{\fX'}[\frac 1 p]$ de $X'_{s,\zar}$ est cohérent.

\subsection{}\label{chb17}
D'après (\cite{sp} \href{https://stacks.math.columbia.edu/tag/013K}{013K}), 
pour tout complexe borné inférieurement de $\bvocB^!$-modules $\cK^\bullet$ de $\tG_s^{\mN^\circ}$,
il existe un complexe borné inférieurement de $\bvocB^!$-modules injectifs $\cL^\bullet$
et un quasi-isomorphisme $u\colon \cK^\bullet\rightarrow \cL^\bullet$. 
Ce dernier induit un morphisme de $\bD^+(\bMod(\co_{\fX'})))$ 
\begin{equation}\label{chb17a}
\huppi_*(\cK^\bullet)\rightarrow \rR \huppi_*(\cK^\bullet).
\end{equation}
D'après (\cite{sp} \href{https://stacks.math.columbia.edu/tag/05TG}{05TG}), ce morphisme ne depend que de $\cK^\bullet$, mais pas de $u$,  
et il en dépend fonctoriellement.  

De même, comme $\bMod_\mQ(\bvocB^!)$  a assez d'injectifs compte tenu de \ref{indsh16}(iii), 
pour tout complexe borné inférieurement $\cK^\bullet$  de $\bMod_\mQ(\bvocB^!)$, 
il existe un complexe borné inférieurement d'objets injectifs $\cL^\bullet$ de $\bMod_\mQ(\bvocB^!)$
et un quasi-isomorphisme $u\colon \cK^\bullet\rightarrow \cL^\bullet$. 
Ce dernier induit un morphisme de $\bD^+(\bMod_\mQ(\co_{\fX'})))$, 
\begin{equation}\label{chb17b}
\huppi_{\mQ*}(\cK^\bullet)\rightarrow \rR \huppi_{\mQ*}(\cK^\bullet).
\end{equation}
D'après (\cite{sp} \href{https://stacks.math.columbia.edu/tag/05TG}{05TG}), ce morphisme ne depend que de $\cK^\bullet$, mais pas de $u$,  
et il en dépend fonctoriellement.

\subsection{}\label{chb18}
Pour tout complexe de $\co_{\fX'}$-modules $\cF'^\bullet$ tel que $\cF'^i=0$ pour $i<0$, et tout entier $q\geq 0$, 
on a un morphisme canonique de changement de base relativement au diagramme \eqref{chb14a},
\begin{equation}\label{chb18a}
\hupsigma^*(\rR^q\fgg_*(\cF'^\bullet))\rightarrow \rR^q\bvlgg_*(\huppi^*(\cF'^\bullet)),
\end{equation}
où $\huppi^*(\cF'^\bullet)$ désigne l'image inverse de $\cF'^\bullet$ définie terme à terme (non dérivée), et 
$\rR^q\fgg_*(-)$ et $\rR^q\bvlgg_*(-)$ désignent les modules d'hypercohomologie. 
En effet, cela revient à se donner un morphisme 
\begin{equation}\label{chb18b}
\rR^q\fgg_*(\cF'^\bullet)\rightarrow \hupsigma_*(\rR^q\bvlgg_*(\huppi^*(\cF'^\bullet))),
\end{equation}
et on prend le morphisme composé 
\begin{eqnarray}
\lefteqn{\rR^q\fgg_*(\cF'^\bullet)\rightarrow \rR^q\fgg_*(\huppi_*(\huppi^*(\cF'^\bullet)))\rightarrow} \label{chb18c}\\
&&\rR^q(\fgg\circ \huppi)_*(\huppi^*(\cF'^\bullet))\stackrel{\sim}{\rightarrow} \rR^q(\hupsigma\circ \bvlgg)_*(\huppi^*(\cF'^\bullet))
\rightarrow  \hupsigma_*(\rR^q\bvlgg_*(\huppi^*(\cF'^\bullet))),\nonumber
\end{eqnarray}
où la première flèche est induite par le morphisme d'adjonction $\cF'^\bullet\rightarrow \huppi_*(\huppi^*(\cF'^\bullet))$, 
la seconde par \eqref{chb17a}, la troisième par l'isomorphisme sous-jacent à \eqref{chb14a}, 
et la quatrième est l'edge-homomorphisme de la 
seconde suite spectrale d'hypercohomologie du foncteur $\hupsigma_*$ par rapport au complexe $\rR \bvlgg_*(\huppi^*(\cF'^\bullet))$ 
(\cite{ega3} 0.11.4.3) et (\cite{hodge2} 1.4.5 et 1.4.6).

Le morphisme de changement de base \eqref{chb18a} est fonctoriel en le complexe $\cF'^\bullet$. 
On notera toutefois que le foncteur $\huppi^*$ ne transforme pas a priori quasi-isomorphisme en quasi-isomorphisme. 
Le morphisme \eqref{chb18a} ne peut donc pas être étendu à $\cF'^\bullet$ objet de $\bD^+(\bMod(\co_{\fX'}))$.

\begin{lem}\label{chb19}
Pour tout $\co_{\fX'}$-module $\cF'$ de $X'_{s,\zar}$ et tout entier $q\geq 0$, le morphisme de changement de base 
\begin{equation}\label{chb19a}
\hupsigma^*(\rR^q\fgg_*(\cF')))\rightarrow \rR^q\bvlgg_*(\huppi^*(\cF'))
\end{equation}
défini dans \eqref{chb18a} coïncide avec le morphisme de changement de base défini dans \eqref{chb14b}.
\end{lem}

En effet, la seconde suite spectrale d'hypercohomologie du foncteur $\hupsigma_*$ 
par rapport au complexe $\rR \bvlgg_*(\huppi^*(\cF'))$ coïncide avec 
la suite spectrale de Cartan-Leray du foncteur composé $\hupsigma_*\circ \bvlgg_*$ (\cite{sga4} V 5.4). 
Le quatrième morphisme de \eqref{chb18c} coïncide donc avec celui induit par cette suite spectrale de Cartan-Leray.  
De même, le deuxième morphisme \eqref{chb18c} coïncide avec l'edge-homomorphisme de la suite spectrale de Cartan-Leray 
du foncteur composé $\fgg_*\circ \huppi_*$ en vertu de (\cite{ega3} (0.11.3.4.2)), d'où la proposition.

\begin{lem}\label{chb24}
Soit $h\colon T'\rightarrow T$ un morphisme de topos, $K^\bullet$ (resp. $L^\bullet$) un complexe de groupes abéliens de $T$ 
(resp. $T'$) tel que $K^i=0$ (resp. $L^i=0$) pour $i<0$, 
$u\colon h^{-1}(K^\bullet)\rightarrow L^\bullet$ un morphisme de complexes de groupes abéliens, $q$ un entier $\geq 0$. On note 
$v\colon K^\bullet\rightarrow h_*(L^\bullet)$ le morphisme adjoint de $u$, défini terme à terme, $u^q\colon h^{-1}(\rH^q(K^\bullet))\rightarrow \rH^q(L^\bullet)$ 
le morphisme induit par $u$ et $v^q\colon \rH^q(K^\bullet)\rightarrow h_*(\rH^q(L^\bullet))$ son adjoint. Alors, $v^q$ est le composé 
\begin{equation}\label{chb24a}
\xymatrix{
{\rH^q(K^\bullet)}\ar[r]^-(0.5){\rH^q(v)}&{\rH^q(h_*(L^\bullet))}\ar[r]&{\rR^q h_*(L^\bullet)}\ar[r]&{h_*(\rH^q(L^\bullet))}},
\end{equation}
où la seconde flèche est induite par le morphisme canonique $h_*(L^\bullet)\rightarrow \rR h_*(L^\bullet)$ et la troisième flèche est 
l'edge-homomorphisme de la seconde suite spectrale d'hypercohomologie du foncteur $h_*$ par rapport au complexe $L^\bullet$. 
\end{lem}

En effet, la deuxième et la troisième flèches de \eqref{chb24a} étant fonctorielles en $L^\bullet$, on peut se borner au cas où $L=h^{-1}(K^\bullet)$ et $u=\id$. 
Comme $v$ est le morphisme d'adjonction $K^\bullet\rightarrow h_*(h^{-1}(K^\bullet))$ et que $u^q$ est l'identité de $h^{-1}(\rH^q(K^\bullet))$,  
il s'agit alors de montrer que le morphisme composé
\begin{equation}\label{chb24b}
\xymatrix{
{\rH^q(K^\bullet)}\ar[r]^-(0.5){\rH^q(v)}&{\rH^q(h_*(h^{-1}(K^\bullet)))}\ar[r]&{\rR^q h_*(h^{-1}(K^\bullet))}\ar[r]&{h_*(h^{-1}(\rH^q(K^\bullet)))}},
\end{equation}
est le morphisme d'adjonction de $\rH^q(K^\bullet)$. 
Par fonctorialité, considérant le morphisme canonique $\tau_{\leq q}(K^\bullet) \rightarrow K^\bullet$, 
on se ramène au cas où $K^\bullet$ est concentré en degrés $[0,q]$, puis au cas où $K^\bullet=M[-q]$ pour un groupe abélien $M$ de $T$.  
Le premier morphisme de \eqref{chb24b}  s'identifie alors au morphisme d'adjonction $M\rightarrow h_*(h^{-1}(M))$, et le 
deuxième et le troisième morphismes s'identifient au morphisme identique de $h_*(h^{-1}(M))$, d'où la proposition.

\subsection{}\label{chb23}
Considérons le diagramme commutatif de morphismes de topos annelés \eqref{hmdf13i} 
\begin{equation}\label{chb23a}
\xymatrix{
{(\tG^{\mN^\circ}_s,\huppi^{-1}(\co_{\fX'}))}\ar[r]^-(0.5){\huppi}\ar[d]_-(0.5){\bvlgg}&{(X'_{s,\zar},\co_{\fX'})}\ar[d]^{\fgg}\\
{(\tE^{\mN^\circ}_s,\hupsigma^{-1}(\co_\fX))}\ar[r]^-(0.5){\hupsigma}&{(X_{s,\zar},\co_{\fX})}}
\end{equation}
Nous utilisons les notations $\huppi^{-1}$ et $\hupsigma^{-1}$ pour désigner les images inverses au sens des faisceaux abéliens
et nous réservons les notations $\huppi^*$ et $\hupsigma^*$ pour les images inverses au sens des modules par les morphismes de topos annelés 
représentés par les flèches horizontales du diagramme \eqref{chb14a}. 
Calquant la construction de \ref{chb18}, pour tout complexe de $\co_{\fX'}$-modules $\cF'^\bullet$ tel que $\cF'^i=0$ pour $i<0$, et tout entier $q\geq 0$, 
on a un morphisme canonique de changement de base relativement au diagramme \eqref{chb23a}
\begin{equation}\label{chb23b}
u^q\colon \hupsigma^{-1}(\rR^q\fgg_*(\cF'^\bullet))\rightarrow \rR^q\bvlgg_*(\huppi^{-1}(\cF'^\bullet)).
\end{equation}
Nous en donnons une construction alternative.  
Soient $K^{\bullet\bullet}$ (resp. $L^{\bullet\bullet}$) une résolution de Cartan-Eilenberg injective de $\cF'^\bullet$ 
(resp. $\huppi^{-1}(\cF'^\bullet)$) telle que $K^{ij}=0$ (resp. $L^{ij}=0$) pour $i<0$ (\cite{ega3} 0.11.4.2). 
Le morphisme de changement de base relativement au diagramme de \eqref{chb23a} 
appliqué terme à terme induit un morphisme de bicomplexes de $\hupsigma^{-1}(\co_\fX)$-modules
\begin{equation}\label{chb23c}
\hupsigma^{-1}(\fgg_*(K^{\bullet\bullet}))\rightarrow \bvlgg_*(\huppi^{-1}(K^{\bullet\bullet})).
\end{equation}
Il est clair que $\huppi^{-1}(K^{\bullet\bullet})$ est une résolution de Cartan-Eilenberg de $\huppi^{-1}(\cF'^\bullet)$. D'après (\cite{ega3} 0.11.4.2), 
il existe donc un morphisme de bicomplexes de $\huppi^{-1}(\co_{\fX'})$-modules $\huppi^{-1}(K^{\bullet\bullet})\rightarrow L^{\bullet\bullet}$.
Prenant l'image par le foncteur $\bvlgg_*$ et composant avec \eqref{chb23c}, on obtient un morphisme de bicomplexes de $\hupsigma^{-1}(\co_\fX)$-modules
\begin{equation}\label{chb23d}
\hupsigma^{-1}(\fgg_*(K^{\bullet\bullet}))\rightarrow \bvlgg_*(L^{\bullet\bullet}).
\end{equation}
Celui-ci induit entre les $\hupsigma^{-1}(\co_\fX)$-modules de cohomologie des complexes simples associés, un morphisme 
\begin{equation}\label{chb23e}
v^q\colon \hupsigma^{-1}(\rR^q\fgg_*(\cF'^{\bullet}))\rightarrow \rR^q\bvlgg_*(\huppi^{-1}(\cF'^\bullet)),
\end{equation}
qui n'est autre que le morphisme $u^q$ \eqref{chb23b}. En effet, l'adjoint du morphisme \eqref{chb23c} est le composé 
\begin{equation}\label{chb23f}
\fgg_*(K^{\bullet\bullet})\rightarrow \fgg_*(\huppi_*(\huppi^{-1}(K^{\bullet\bullet}))) \rightarrow \fgg_*(\huppi_*(L^{\bullet\bullet}))
\stackrel{\sim}{\rightarrow} \hupsigma_*(\bvlgg_*(L^{\bullet\bullet})),
\end{equation}
où la première flèche est induite par le morphisme d'adjonction $\id \rightarrow \huppi_*\huppi^{-1}$, la seconde flèche par le morphisme
$\huppi^{-1}(K^{\bullet\bullet})\rightarrow L^{\bullet\bullet}$, et la troisième flèche par l'isomorphisme sous-jacent au diagramme \eqref{chb23a}. 
D'après \ref{chb24}, l'adjoint du morphisme $v^q$ est le composé 
\begin{equation}
\rR^q\fgg_*(\cF'^\bullet)\rightarrow \rR^q(\hupsigma\circ \bvlgg)_*(\huppi^{-1}(\cF'^\bullet))\rightarrow \hupsigma_*(\rR^q\bvlgg_*(\huppi^{-1}(\cF'^\bullet))),
\end{equation}
où la première flèche est induite par le morphisme composé \eqref{chb23f} et la seconde flèche est l'edge-homomorphisme de la 
seconde suite spectrale d'hypercohomologie du foncteur $\hupsigma_*$ par rapport au complexe $\rR \bvlgg_*(\huppi^*(\cF'^\bullet))$. 
Par ailleurs, le diagramme  
\begin{equation}
\xymatrix{
{\rR^q\fgg_*(\cF'^\bullet)}\ar[r]^-(0.5){a}\ar[d]&{\rR^q\fgg_*(\huppi_*(\huppi^{-1}(\cF'^\bullet)))}\ar[r]^-(0.5){b}\ar[d]&
{\rR^q\fgg_*(\rR \huppi_*(\huppi^{-1}(\cF'^\bullet)))}\ar[d]\\
{\rR^q\fgg_*(\Tot(K^{\bullet\bullet}))}\ar[r]^-(0.5){a'}\ar[d]&{\rR^q\fgg_*(\huppi_*(\huppi^{-1}(\Tot(K^{\bullet\bullet}))))}\ar[r]^-(0.5){b'}&
{\rR^q\fgg_*(\huppi_*(\Tot(L^{\bullet\bullet})))}\ar[d]\\
{\rH^q(\fgg_*(\Tot(K^{\bullet\bullet})))}\ar[rr]^-(0.5){a''}&&{\rH^q(\Tot(\fgg_*(\huppi_*(L^{\bullet\bullet})))}}
\end{equation}
où $\Tot(-)$ désigne le complexe total associé, $a$ et $a'$ sont induits par le morphisme d'adjonction $\id\rightarrow \huppi_*\huppi^{-1}$,
$b$ est induit par \eqref{chb17a}, $b'$ est induit par le morphisme $\huppi^{-1}(K^{\bullet\bullet})\rightarrow L^{\bullet\bullet}$, 
$a''$ est induit par le composé des deux premières flèches de \eqref{chb23f} et les flèches verticales sont les morphismes canoniques, est commutatif.  
L'égalité $u^q=v^q$ s'ensuit compte tenu de la définition \eqref{chb18c} de l'adjoint du morphisme $u^q$ \eqref{chb23b}. 

\subsection{}\label{chb20}
Pour tout complexe $F'^\bullet$ de $\bMod_\mQ(\co_{\fX'})$ tel que $F'^i=0$ pour $i<0$, et tout entier $q$, 
on a un morphisme canonique de changement de base relativement au diagramme \eqref{chb14a},
\begin{equation}\label{chb20a}
\hupsigma^*_\mQ(\rR^q\fgg_{\mQ*}(F'^\bullet))\rightarrow \rR^q\bvlgg_{\mQ*}(\huppi^*_\mQ(F'^\bullet)),
\end{equation}
où $\huppi^*_\mQ(F'^\bullet)$ désigne l'image inverse de $F'^\bullet$ définie terme à terme \eqref{hmdf40h}, et 
$\rR^q\fgg_{\mQ*}(-)$ et $\rR^q\bvlgg_{\mQ*}(-)$ désignent les modules d'hypercohomologie. 
En effet, cela revient à se donner un morphisme 
\begin{equation}\label{chb20b}
\rR^q\fgg_{\mQ*}(F'^\bullet)\rightarrow \hupsigma_{\mQ*}(\rR^q\bvlgg_{\mQ*}(\huppi^*_\mQ(F'^\bullet))),
\end{equation}
et on prend le morphisme composé 
\begin{eqnarray}
\lefteqn{\rR^q\fgg_{\mQ*}(F'^\bullet)\rightarrow \rR^q\fgg_{\mQ*}(\huppi_{\mQ*}(\huppi^*_\mQ(F'^\bullet)))\rightarrow} \label{chb20c}\\
&&\rR^q(\fgg\circ \huppi)_{\mQ*}(\huppi^*_\mQ(F'^\bullet))\stackrel{\sim}{\rightarrow} \rR^q(\hupsigma\circ \bvlgg)_{\mQ*}(\huppi^*_\mQ(F'^\bullet))
\rightarrow  \hupsigma_{\mQ*}(\rR^q\bvlgg_{\mQ*}(\huppi^*_\mQ(F'^\bullet))),\nonumber
\end{eqnarray}
où la première flèche est induite par le morphisme d'adjonction $F'^\bullet\rightarrow \huppi_{\mQ*}(\huppi^*_\mQ(F'^\bullet))$, 
la seconde par \eqref{chb17b}, la troisième par l'isomorphisme sous-jacent à \eqref{chb14a}, 
et la quatrième est l'edge-homomorphisme de la 
seconde suite spectrale d'hypercohomologie du foncteur $\hupsigma_{\mQ*}$ par rapport au complexe 
$\rR \bvlgg_{\mQ*}(\huppi^*_\mQ(F'^\bullet))$ 
(\cite{ega3} 0.11.4.3)  et (\cite{hodge2} 1.4.5 et 1.4.6). 
On notera que la catégorie $\bMod_\mQ(\bvocB)$ a assez d'injectifs d'après \ref{indsh16}(iii). 
En particulier, tout complexe de $\bMod_\mQ(\bvocB)$ admet une résolution de Cartan-Eilenberg injective. 

Le morphisme de changement de base \eqref{chb20a} est fonctoriel en le complexe $F'^\bullet$. 
On notera toutefois que le foncteur $\huppi^*_\mQ$ ne transforme pas a priori quasi-isomorphisme en quasi-isomorphisme. 
Le morphisme \eqref{chb18a} ne peut donc pas être étendu à $F'^\bullet$ objet de $\bD^+(\bMod_\mQ(\co_{\fX'}))$.

\begin{lem}\label{chb21}\ 
\begin{itemize}
\item[{\rm (i)}] Pour tout complexe de $\co_{\fX'}$-modules $\cF'^\bullet$ tel que $\cF'^i=0$ pour $i<0$, et tout entier $q$, 
le morphisme de changement de base \eqref{chb20a}
\begin{equation}\label{chb21a}
\hupsigma^*_\mQ(\rR^q\fgg_{\mQ*}(\cF'^\bullet_\mQ))\rightarrow \rR^q\bvlgg_{\mQ*}(\huppi^*_\mQ(\cF'^\bullet_\mQ))
\end{equation} 
est l'image canonique du morphisme de changement de base \eqref{chb18a}
\begin{equation}\label{chb21b}
\hupsigma^*(\rR^q\fgg_*(\cF'^\bullet))\rightarrow \rR^q\bvlgg_*(\huppi^*(\cF'^\bullet)).
\end{equation}
\item[{\rm (ii)}] Pour tout objet $F'$ de $\bMod_\mQ(\co_{\fX'})$ et tout entier $q\geq 0$, le morphisme de changement de base 
\begin{equation}\label{chb21c}
\hupsigma^*_\mQ(\rR^q\fgg_{\mQ*}(F')))\rightarrow \rR^q\bvlgg_{\mQ*}(\huppi^*_\mQ(F'))
\end{equation}
défini dans \eqref{chb20a} coïncide avec le morphisme de changement de base défini dans \eqref{chb14c}.
\end{itemize}
\end{lem}

(i) En effet, chacun des morphismes apparaissant dans \eqref{chb20c} 
est l'image canonique du morphisme correspondant apparaissant dans \eqref{chb18c}.  
Ceci est évident pour la première et la troisième flèches compte tenu de \eqref{indsh14d}, et résulte de \ref{indsh16}(iii) pour la seconde flèche.  
Le cas de la quatrième flèche résulte du fait que l'image canonique d'une résolution de Cartan-Eilenberg injective d'un complexe de $\bvocB$-modules 
$\cK^\bullet$ est une résolution de Cartan-Eilenberg injective du complexe $\cK^\bullet_\mQ$.

(ii) Cela résulte de (i) et \ref{chb19}.

\subsection{}\label{chb25}
Calquant la construction de \ref{chb18}, pour tout complexe $F'^\bullet$ de $\bMod_\mQ(\co_{\fX'})$ tel que $F'^i=0$ pour $i<0$, et tout entier $q$,
on a un morphisme canonique de changement de base relativement au diagramme \eqref{chb23a}
\begin{equation}\label{chb25a}
\hupsigma^{-1}_\mQ(\rR^q\fgg_{\mQ*}(F'^\bullet))\rightarrow \rR^q\bvlgg_{\mQ*}(\huppi^{-1}_\mQ(F'^\bullet)).
\end{equation}
Comme dans \ref{chb23}, on peut en donner une autre construction utilisant des résolutions de Cartan-Eilenberg injective de $F'^\bullet$ et $\huppi^{-1}_\mQ(F'^\bullet)$,
qui existent puisque les catégories $\bMod_\mQ(\co_{\fX'})$ et $\bMod_\mQ(\huppi^{-1}(\co_{\fX'}))$ ont assez d'injectifs d'après \ref{indsh16}(iii).

\begin{teo}\label{chb22}
Supposons $g\colon X'\rightarrow X$ propre. 
Soient $\cF'^\bullet$ un complexe de $\co_{\fX'}[\frac 1 p]$-modules cohérents tel que $\cF'^i=0$ pour $i<0$, $q$ un entier $\geq 0$. 
Alors, le $\co_\fX[\frac 1 p]$-module $\rR^q\fgg_*(\cF'^\bullet)$ est cohérent et le morphisme de changement de base \eqref{chb20a}
\begin{equation}\label{chb22a}
\hupsigma^*_\mQ(\rR^q\fgg_*(\cF'^\bullet))\rightarrow \rR^q\bvlgg_{\mQ*}(\huppi^*_\mQ(\cF'^\bullet))
\end{equation}
où $\hupsigma^*_\mQ$ et $\huppi^*_\mQ$ désignent les foncteurs exacts \eqref{hmdf40j} et \eqref{hmdf40k} (cf. \ref{chb5}), 
est un isomorphisme.
\end{teo}

En effet, la première assertion résulte de (\cite{egr1} 2.10.24 et 2.11.5) compte tenu de 
la seconde suite spectrale d'hypercohomologie du foncteur $\fgg_*$ par rapport au complexe $\cF'^\bullet$, 
\begin{equation}\label{chb22b}
\rE_2^{i,j}= \rR^i\fgg_*(\cH^j(\cF'^\bullet))\Rightarrow \rR^{i+j}\fgg_*(\cF'^\bullet).
\end{equation}
Compte tenu de \ref{chb5}, cette dernière induit une suite spectrale 
\begin{equation}\label{chb22c}
\rE_2^{i,j}= \hupsigma^*_\mQ(\rR^i\fgg_*(\cH^j(\cF'^\bullet)))\Rightarrow \hupsigma^*_\mQ(\rR^{i+j}\fgg_*(\cF'^\bullet)).
\end{equation}
Par ailleurs, compte tenu encore de \ref{chb5}, 
la seconde suite spectrale d'hypercohomologie du foncteur $\bvlgg_{\mQ*}$ par rapport au complexe $\huppi^*_\mQ(\cF'^\bullet)$ s'écrit
\begin{equation}\label{chb22d}
\rE_2^{i,j}= \rR^i\bvlgg_{\mQ*}(\huppi^*_\mQ(\cH^j(\cF'^\bullet)))\Rightarrow \rR^{i+j}\bvlgg_{\mQ*}(\huppi^*_\mQ(\cF'^\bullet)).
\end{equation} 
On rappelle que la catégorie $\bMod_\mQ(\bvocB^!)$ a assez d'injectifs d'après \ref{indsh16}(iii). 
Il suffit de montrer que les morphismes de changement de base \eqref{chb14d}
\begin{equation}\label{chb22e}
\hupsigma^*_\mQ(\rR^i\fgg_*(\cH^j(\cF'^\bullet)))\rightarrow \rR^i\bvlgg_{\mQ*}(\huppi^*_\mQ(\cH^j(\cF'^\bullet))),
\end{equation}
et les morphismes de changement de base \eqref{chb22a}
\begin{equation}\label{chb22f}
\hupsigma^*_\mQ(\rR^{i+j}\fgg_*(\cF'^\bullet))\rightarrow \rR^{i+j}\bvlgg_{\mQ*}(\huppi^*_\mQ(\cF'^\bullet))
\end{equation}
définissent un morphisme de la suite spectrale \eqref{chb22c} vers la suite spectrale \eqref{chb22d}. 
On effet, comme les morphismes \eqref{chb22e} sont des isomorphismes en vertu de \ref{chb16}, 
on en déduit que le morphisme \eqref{chb22a} est un isomorphisme pour tout $q\geq 0$. 

Pour tout objet $F'$ de $\bMod_\mQ(\co_{\fX'})$, on a un morphisme canonique de changement de base 
relativement au diagramme \eqref{chb23a},
\begin{equation}\label{chb22h}
\hupsigma^{-1}_\mQ(\rR^q\fgg_{\mQ*}(F'))\rightarrow \rR^q\bvlgg_{\mQ*}(\huppi^{-1}_\mQ(F')).
\end{equation}
Par ailleurs, pour tout complexe $F'^\bullet$ de $\bMod_\mQ(\co_{\fX'})$ tel que $F'^i=0$ pour $i<0$, et tout entier $q\geq 0$, 
on a un morphisme canonique de changement de base \eqref{chb25a} relativement au diagramme \eqref{chb23a},
\begin{equation}\label{chb22i}
\hupsigma^{-1}_\mQ(\rR^q\fgg_{\mQ*}(F'^\bullet))\rightarrow \rR^q\bvlgg_{\mQ*}(\huppi^{-1}_\mQ(F'^\bullet)).
\end{equation}
Par ailleurs, on dispose de deux suites spectrales 
\begin{equation}\label{chb22j}
\rE_2^{i,j}= \hupsigma^{-1}_\mQ(\rR^i\fgg_{\mQ*}(\cH^j(F'^\bullet)))\Rightarrow \hupsigma^{-1}_\mQ(\rR^{i+j}\fgg_{\mQ*}(F'^\bullet)),
\end{equation}
\begin{equation}\label{chb22k}
\rE_2^{i,j}= \rR^i\bvlgg_{\mQ*}(\huppi^{-1}_\mQ(\cH^j(F'^\bullet)))\Rightarrow \rR^{i+j}\bvlgg_{\mQ*}(\huppi^{-1}_\mQ(F'^\bullet)).
\end{equation} 
Montrons d'abord que les morphismes de changement de base \eqref{chb22h}
\begin{equation}\label{chb22l}
\hupsigma^{-1}_\mQ(\rR^i\fgg_{\mQ*}(\cH^j(F'^\bullet)))\rightarrow \rR^i\bvlgg_{\mQ*}(\huppi^{-1}_\mQ(\cH^j(F'^\bullet)))
\end{equation}
et les morphismes \eqref{chb22i} définissent un morphisme de la suite spectrale \eqref{chb22j} vers la suite spectrale \eqref{chb22k}. 

Soient $K^{\bullet\bullet}$ (resp. $L^{\bullet\bullet}$) une résolution de Cartan-Eilenberg injective de $F'^\bullet$ 
(resp.  $\huppi^{-1}_\mQ(F'^\bullet)$) dans la catégorie $\bMod_\mQ(\co_{\fX'})$ (resp. $\bMod_\mQ(\huppi_\mQ^{-1}(\co_{\fX'}))$)
telle que $K^{ij}=0$ (resp. $L^{ij}=0$) pour $i<0$ (\cite{ega3} 0.11.4.2). 
La suite spectrale \eqref{chb22j} (resp. \eqref{chb22k}) est par définition la seconde suite spectrale du bicomplexe 
$\hupsigma_\mQ^{-1}(\fgg_{\mQ*}(K'^{\bullet\bullet}))$ (resp. $\bvlgg_{\mQ*}(L^{\bullet\bullet})$) 
(\cite{ega3} 0.11.3.2). Le morphisme de changement de base 
\eqref{chb22h} définit un morphisme de bicomplexes 
\begin{equation}\label{chb22m}
\hupsigma_\mQ^{-1}(\fgg_{\mQ*}(K^{\bullet\bullet})) \rightarrow \bvlgg_{\mQ*}(\huppi^{-1}_\mQ(K^{\bullet\bullet})).
\end{equation}
Par ailleurs, $\huppi^{-1}_\mQ(K^{\bullet\bullet})$ étant une résolution de Cartan-Eilenberg de $\huppi^{-1}_\mQ(F'^\bullet)$,
il existe un morphisme de bicomplexes $\huppi^{-1}_\mQ(K^{\bullet\bullet})\rightarrow L^{\bullet\bullet}$ (\cite{ega3} 0.11.4.2),
compatible avec les morphismes de $\huppi^{-1}_\mQ(F'^\bullet)$ dans $\huppi^{-1}_\mQ(K^{\bullet0})$ et 
$L^{\bullet0}$. On en déduit un morphisme de bicomplexes 
\begin{equation}\label{chb22n}
\hupsigma_\mQ^{-1}(\fgg_{\mQ*}(K^{\bullet\bullet})) \rightarrow \bvlgg_{\mQ*}(L^{\bullet\bullet}),
\end{equation}
et par suite un morphisme entre les secondes suites spectrales associées \eqref{chb22j} et \eqref{chb22k}.
D'après \ref{chb25}, ce dernier est défini sur les termes initiaux par \eqref{chb22l} et sur les aboutissements par \eqref{chb22i}.  

Prenons maintenant pour $F'^\bullet$ le complexe $\cF'^\bullet$ \eqref{ahttf50b}.
Compte tenu de \ref{chb5}, le morphisme canonique de complexes $\huppi_\mQ^{-1}(\cF'^\bullet)\rightarrow \huppi_\mQ^*(\cF'^\bullet)$ 
induit un morphisme de la suite spectrale \eqref{chb22k} vers la suite spectrale \eqref{chb22d}. 
Composant avec le morphisme défini plus haut, on en déduit 
un morphisme de la suite spectrale \eqref{chb22j} vers la suite spectrale \eqref{chb22d}. 
Par $\bvocB_\mQ$-linéarisation \eqref{chb5}, on en déduit le morphisme recherché de 
la suite spectrale \eqref{chb22c} vers la suite spectrale \eqref{chb22d}, d'où la proposition.

\section{Fonctorialité des algèbres de Higgs-Tate dans les topos de Faltings}\label{ahttfg}

Dans cette section $r$ désigne un nombre rationnel $\geq 0$ et $n$ un entier $\geq 0$. 

\subsection{}\label{ahttfg1}
On désigne par $\bP$ la sous-catégorie pleine de $\Et_{/X}$ formée des schémas {\em affines} $U$ tels que l'une des deux conditions suivantes 
soit remplie:
\begin{itemize}
\item[(i)] le schéma $U_s$ est vide; ou 
\item[(ii)] le morphisme $(U,\cM_X|U)\rightarrow (S,\cM_S)$ induit par $f$ \eqref{hmdf3} admet une carte adéquate \eqref{cad1} (\cite{agt} III.4.4). 
\end{itemize}
On désigne par $\bQ$ la sous-catégorie pleine de $\bP$ formée des schémas affines $U$ tels que l'une des conditions suivantes soit remplie:
\begin{itemize}
\item[(iii)] le schéma $U_s$ est vide; ou 
\item[(iv)] il existe une carte fine et saturée $M\rightarrow \Gamma(U,\cM_X)$ pour $(U,\cM_X|U)$
induisant un isomorphisme 
\begin{equation}\label{ahttfg1a}
M\stackrel{\sim}{\rightarrow} \Gamma(U,\cM_X)/\Gamma(U,\co^\times_X).
\end{equation}
Cette carte est a priori indépendante de la carte adéquate requise dans (ii). 
\end{itemize}

On désigne par 
\begin{eqnarray}
\pi_\bP\colon E_\bP&\rightarrow& \bP,\label{ahttfg1b}\\
\pi_\bQ\colon E_\bQ&\rightarrow& \bQ,\label{ahttfg1c}
\end{eqnarray} 
les sites fibrés déduits de $\pi$ \eqref{hmdf7a} par changement de base par les foncteurs d'injection canoniques de $\bP$ et $\bQ$ dans $\Et_{/X}$. 
On notera que $\bP$ et $\bQ$ sont des sous-catégories topologiquement génératrices de $\Et_{/X}$. 
Par suite, $E_\bP$ et $E_\bQ$ sont des sous-catégories topologiquement génératrices de $E$. 
Par ailleurs, $\bP$ étant stable par produits fibrés, $\pi_\bP$ est un site fibré co-évanescent \eqref{cftf1}; 
mais ce n'est en général pas le cas de $E_\bQ$ (cf. \ref{ahttf4}, \ref{ahttf5} et \ref{ahttf45}). 

On définit de même les sous-catégories pleines $\bP'$ et $\bQ'$ de $\Et_{/X'}$ relativement au morphisme $f'$ \eqref{hmdf3}.

\subsection{}\label{ahttfg100}
Comme $X$ est noethérien et donc quasi-séparé, tout objet de $\bP$ est cohérent sur $X$. 
Considérons le site fibré co-évanescent \eqref{hmdf7f}
\begin{equation}\label{ahttfg100a}
\pi'_\coh\colon E'_\coh\rightarrow \Et_{\coh/X'}.
\end{equation}
On désigne par 
\begin{equation}\label{ahttfg100b}
\varphi\colon \bP\rightarrow \Et_{\coh/X'}
\end{equation}
le foncteur image inverse par $g\colon X'\rightarrow X$, et par 
\begin{equation}\label{ahttfg100c}
\Phi\colon E_\bP\rightarrow E'_{\coh}
\end{equation}
le foncteur induit par $\Theta^+$ \eqref{hmdf11a}. Le diagramme de foncteurs
\begin{equation}\label{ahttfg100d}
\xymatrix{
{E_\bP}\ar[r]^-(0.5){\pi_\bP}\ar[d]_{\Phi}&{\bP}\ar[d]^{\varphi}\\
{E'_\coh}\ar[r]^-(0.5){\pi'_\coh}&{\Et_{\coh/X'}}}
\end{equation}
est strictement commutatif, {\em i.e.} on a $\varphi \circ \pi_\bP=\pi'_\coh\circ \Phi$. 
On désigne par 
\begin{equation}\label{ahttfg100e}
\pi'_\bP\colon E'_\bP\rightarrow \bP
\end{equation} 
la catégorie fibrée déduite de $\pi'_\coh$ par changement de base par $\varphi$ (\cite{sga1} VI § 3), et par 
\begin{equation}\label{ahttfg100f}
\upphi\colon E_\bP\rightarrow E'_\bP
\end{equation}
le $\bP$-foncteur induit par $\Phi$ \eqref{ahttfg100c}. Ce foncteur $\upphi$ est clairement cartésien. 
Le cadre envisagé ci-dessus est donc un cas particulier de celui considéré dans \ref{cftf5}.
Par ailleurs, les hypothèses de \ref{cftf8} sont satisfaites. 
Le foncteur $\Phi$ est donc continu pour les topologies co-évanescentes sur $E_\bP$ et $E'_\coh$ \eqref{cftf5}. 
Il induit un morphisme de topos qui s'identifie canoniquement à $\Theta\colon \tE'\rightarrow \tE$ \eqref{hmdf11b}.

\subsection{}\label{ahttfg103}
On note $\mI$ la catégorie des morphismes $U'\rightarrow U$ au-dessus du morphisme $g\colon X'\rightarrow X$
tels que le morphisme $U'\rightarrow X'$ (resp. $U\rightarrow X$) soit étale de présentation finie (resp. un objet de $\bP$).
Celle-ci s'identifie à la catégorie portant le même nom définie dans \ref{cftf12} relativement au foncteur $\varphi$ \eqref{ahttfg100b}. 
Considérons les foncteurs 
\begin{eqnarray}
\tts\colon \mI\rightarrow \Et_{\coh/X'},&& (U'\rightarrow U)\mapsto U',\label{ahttfg103a}\\
\ttb\colon \mI\rightarrow \bP,&& (U'\rightarrow U)\mapsto U.\label{ahttfg103b}
\end{eqnarray}
Pour tout $U'\in \ob(\Et_{\coh/X'})$, on désigne par $I^{U'}_\varphi$ la catégorie fibre de $\tts$ au-dessus de $U'$, 
autrement dit la catégorie des objets $(\upmu\colon U'\rightarrow U)$ de $\mI$; un tel objet de $I^{U'}_\varphi$ sera aussi noté $(U,\upmu)$.  

On désigne par $\mJ$ la sous-catégorie pleine de $\mI$ formée des morphismes $U'\rightarrow U$ tels que $U'$ soit un objet de $\bQ'$. 
Pour tout objet $(U'\rightarrow U)$ de $\mI$, on note $\mJ_{/(U'\rightarrow U)}$ la catégorie des morphismes de $\mI$ 
d'un objet de $\mJ$ dans $(U'\rightarrow U)$.  

On désigne par $J'$ l'image essentielle de la catégorie $\mJ$ par le foncteur $\tts$. 
Pour tout objet $U'$ de $\Et_{\coh/X'}$, on note $J'_{/U'}$ la catégorie des $X'$-morphismes d'un objet de $J'$ dans $U'$.

\begin{lem}\label{ahttfg6}
\
\begin{itemize}
\item[{\rm (i)}] Le foncteur $\ttb$ est essentiellement surjectif, et la sous-catégorie $J'$ de $\Et_{\coh/X'}$ est topologiquement génératrice. 
\item[{\rm (ii)}] Supposons le schéma $X$ séparé. 
Alors, pour tout objet $(U'\rightarrow U)$ de $\mI$, tout objet $V'$ de $J'$ et tout $X'$-morphisme $u'\colon V'\rightarrow U'$, il existe 
un objet $(V'\rightarrow W)$ de $\mJ$ et un morphisme $(v',v)\colon (V'\rightarrow W)\rightarrow (U'\rightarrow U)$ de $\mI$.
\begin{equation}
\xymatrix{
V'\ar[r]^-(0.5){v'}\ar[d]&U'\ar[d]\\
W\ar[r]^-(0.5)v&U}
\end{equation}
\end{itemize}
\end{lem}

(i) Cela résulte de (\cite{agt} II.5.17). 

(ii) Soient $(U'\rightarrow U)$ un objet de $\mI$, $V'$ un objet de $J'$, $u'\colon V'\rightarrow U'$ un $X'$-morphisme. 
Il existe alors un objet $V'\rightarrow V$ de $\mJ$. Considérons le diagramme commutatif 
\begin{equation}
\xymatrix{
V'\ar[r]^{u'}\ar[d]\ar@/_2pc/[dd]&U'\ar[r]\ar[d]&X'\ar[dd]\\
V\times_XU\ar[r]\ar[d]&U\ar[rd]&\\
V\ar[rr]&&X}
\end{equation}
Comme $X$ est séparé et que $U$ et $V$ sont affines, $V\times_XU$ est affine. Par ailleurs, une carte adéquate pour $(V,\cM_X|V)$ induit 
une carte adéquate pour $(V\times_XU,\cM_X|V\times_XU)$. Par suite, $V\times_XU$ est un objet de $\bP$,   
et le morphisme $V'\rightarrow V\times_XU$ est un objet de $\mJ$; d'où la proposition.

\begin{prop}\label{ahttfg101}
Soit $F=\{U\in \bP^\circ \mapsto F_U\}$ un v-préfaisceau sur $E_\bP$ \eqref{cftf3}. Alors, 
\begin{itemize}
\item[{\rm (i)}] 
Pour tout $U'\in \ob(\Et_{\coh/X'})$, les faisceaux $\oupmu^*(F_U)$, pour $(U,\upmu\colon U'\rightarrow U)\in \ob((I^{U'}_\varphi)^\circ)$, 
où $\oupmu\colon \oU'^\rhd\rightarrow \oU^\circ$ désigne abusivement le morphisme induit par $\upmu$, forment naturellement un système inductif 
de $\oU'^\rhd_\fet$.  Posons 
\begin{equation}\label{ahttfg101a}
F'_{U'}=\underset{\underset{(U,\upmu)\in (I^{U'}_\varphi)^\circ}{\longrightarrow}}{\lim}\ \oupmu^*(F_U).
\end{equation}
\item[{\rm (ii)}] La collection $F'=\{U'\in \Et^\circ_{\coh/X'}\mapsto F'_{U'}\}$ forme naturellement un v-préfaisceau sur $E'$. 
Pour tout morphisme $\uplambda\colon U'_1\rightarrow U'_2$ de $\Et_{\coh/X'}$ et tout objet $\upmu_2\colon U'_2\rightarrow U$ de $\mI$ \eqref{ahttfg103}, 
posant $\upmu_1=\upmu_2\circ \uplambda \colon U'_1\rightarrow U$, le diagramme 
\begin{equation}
\xymatrix{
{\ouplambda^{\rhd*}(\oupmu_2^*(F_U))}\ar[r]\ar[d]&{\oupmu_1^*(F_U)}\ar[d]\\
{\ouplambda^{\rhd*}(F'_{U'_2})}\ar[r]&{F'_{U'_1}}}
\end{equation}
où les flèches verticales sont les morphismes canoniques \eqref{ahttfg101a}, la flèche horizontale supérieure est l'isomorphisme canonique 
et la flèche horizontale inférieure est le morphisme adjoint du morphisme $F'_{U'_2}\rightarrow \ouplambda^\rhd_*(F'_{U'_1})$ 
définissant la structure de préfaisceau sur $F'$, est commutatif.
\item[{\rm (iii)}] On a un isomorphisme canonique fonctoriel 
\begin{equation}\label{ahttfg101b}
\Theta^*(F^a)\stackrel{\sim}{\rightarrow} F'^a,
\end{equation} 
où l'exposant $^a$ désigne les faisceaux associés. 
\item[{\rm (iv)}] Avec les notations de \ref{hmdf20}, pour tout objet $\upmu\colon U'\rightarrow U$ de $\mI$ \eqref{ahttfg103}, le diagramme 
\begin{equation}\label{ahttfg101c}
\xymatrix{
{\oupmu^*(F_U)}\ar[r]^-(0.5){\oupmu^*(a)}\ar[d]_b&{\oupmu^*(\beta_{U*}(\jmath_U^*(F^a)))}\ar[r]^-(0.5)c&
{\beta'_{U'*}(\Theta_\upmu^*(\jmath_U^*(F^a)))}\ar[d]^d\\
{F'_{U'}}\ar[r]^-(0.5){a'}&{\beta'_{U'*}(\jmath'^*_{U'}(F'^a))}\ar[r]^-(0.5)e&{\beta'_{U'*}(\jmath'^*_{U'}(\Theta^*(F^a)))}}
\end{equation}
où $a\colon F_U\rightarrow \beta_{U*}(\jmath_U^*(F^a))$, $a'$ et $b$ sont les morphismes canoniques,  
$c$ est le morphisme de changement de base relativement à \eqref{hmdf20e}, $d$ est l'isomorphisme
sous-jacent à \eqref{hmdf20d} et $e$ est l'isomorphisme induit par \eqref{ahttfg101b}, est commutatif. 
\end{itemize}
\end{prop}

Les propositions (i), (ii) et (iii) sont des cas particuliers de \ref{cftf9}. Montrons la proposition (iv). 
Par localisation \eqref{hmdf20d}, on peut se borner au cas où $\upmu=g$.
Il résulte alors de \ref{cftf15}(ii), \eqref{cftf8c} et \eqref{cftf8f} que le diagramme 
\begin{equation}
\xymatrix{
{F_X}\ar[r]^-(0.5){a}\ar[d]_{b'}&{\beta_*(F^a)}\ar[r]^-(0.5){\beta_*(\ad)}&{\beta_*(\Theta_*(\Theta^*(F^a)))}\ar[d]^{c'}\\
{\upgamma_*(F'_{X'})}\ar[r]^-(0.5){\upgamma_*(a')}&{\upgamma_*(\beta'_*(F'^a))}\ar[r]^-(0.5){e'}&{\upgamma_*(\beta'_*(\Theta^*(F^a)))}}
\end{equation}
où $\upgamma\colon \oX'^\rhd\rightarrow \oX^\circ$ est le morphisme induit par $g$, 
$a$ et $a'\colon F'_{X'}\rightarrow \beta'_*(F'^a)$ sont les morphismes canoniques, $b'$ est l'adjoint du morphisme canonique 
$\upgamma^*(F_X)\rightarrow F'_{X'}$ \eqref{ahttfg101a}, $\ad\colon F^a \rightarrow \Theta_*(\Theta^*(F^a))$
est le morphisme d'adjonction, $e'$ est induit par l'isomorphisme \eqref{ahttfg101b} et $c'$ est l'isomorphisme
sous-jacent à \eqref{hmdf20e}, est commutatif. La proposition s'ensuit par adjonction.

\subsection{}\label{ahttfg19}
Il résulte aussitôt de \ref{ahttfg101} que 
l'homomorphisme canonique $\Theta^{-1}(\ocB)\rightarrow \ocB'$ \eqref{hmdf12c} induit pour tout objet $(\upmu\colon U'\rightarrow U)$ 
de $\mI$ \eqref{ahttfg103}, un homomorphisme de $\oU'^\rhd_\fet$, 
\begin{equation}\label{ahttfg19a}
\oupmu^*(\ocB_U)\rightarrow \ocB'_{U'},
\end{equation}
où $\oupmu\colon \oU'^\rhd\rightarrow \oU^\circ$ désigne abusivement le morphisme induit par $\upmu$. 
Avec la terminologie de \ref{cftf122}, ces homomorphismes forment un 
$\mI$-système de $\Phi$-morphismes compatibles de $\ocB$ dans $\ocB'$ \eqref{ahttfg100c}. 

\subsection{}\label{ahttfg2}
Soient $Y$ un objet de $\bP$ tel que $Y_s$ soit non-vide,  $((P,\gamma),(\mN,\iota),\vartheta)$
une carte adéquate pour le morphisme $f|Y\colon (Y,\cM_X|Y)\rightarrow (S,\cM_S)$ induit par $f$, 
$\oy$ un point géométrique de $\oY^\circ$.  
Le schéma $\oY$ étant localement irréductible d'après (\cite{ag} 4.2.7 et \cite{agt} III.3.3),  
il est la somme des schémas induits sur ses composantes irréductibles. On note $\oY^\star$
la composante irréductible de $\oY$ contenant $\oy$. 
De même, $\oY^\circ$ est la somme des schémas induits sur ses composantes irréductibles
et $\oY^{\star \circ}=\oY^\star\times_{X}X^\circ$ est la composante irréductible de $\oY^\circ$ contenant $\oy$.
On reprend les notations introduites dans
\ref{ahttf6}, en particulier,  la représentation discrète $\oR^\oy_Y$ de $\pi_1(\oY^{\star\circ},\oy)$ définie dans \eqref{ahttf44b},
la suite exacte canonique \eqref{ahttf6f}
\begin{equation}\label{ahttfg2a}
0\rightarrow \hoR^\oy_Y\rightarrow \cF^\oy_Y\rightarrow \txi^{-1}\tOmega^1_{X/S}(Y) \otimes_{\co_X(Y)}\hoR^\oy_Y\rightarrow 0,
\end{equation} 
et la $\hoR^\oy_Y$-algèbre \eqref{ahttf6g} 
\begin{equation}\label{ahttfg2b}
\cC^\oy_Y=\underset{\underset{m\geq 0}{\longrightarrow}}\lim\ \rS^m_{\hoR^\oy_Y}(\cF^\oy_Y).
\end{equation}
On notera que ces représentations dépendent de la carte adéquate $((P,\gamma),(\mN,\iota),\vartheta)$. 
Toutefois, elles n'en dépendent pas si $Y$ est un objet de $\bQ$ d'après \ref{pmh7}. 

On désigne par $\cF_Y^{\oy,(r)}$ 
l'extension de $\hoR_Y^\oy$-modules déduite de $\cF_Y^\oy$ 
\eqref{ahttfg2a} par image inverse par le morphisme de multiplication par $p^r$ sur 
$\txi^{-1}\tOmega^1_{X/S}(Y)\otimes_{\co_X(Y)} \hoR^\oy_Y$,
de sorte qu'on a une suite exacte de $\hoR^\oy_Y$-modules 
\begin{equation}\label{ahttfg2c}
0\rightarrow \hoR^\oy_Y\rightarrow \cF_Y^{\oy,(r)}\rightarrow \txi^{-1}\tOmega^1_{X/S}(Y)\otimes_{\co_X(Y)} \hoR^\oy_Y
\rightarrow 0.
\end{equation}
On désigne par $\cC_Y^{\oy,(r)}$ la $\hoR^\oy_Y$-algèbre \eqref{taht10c}
\begin{equation}\label{ahttfg2d}
\cC_Y^{\oy,(r)}= \underset{\underset{m\geq 0}{\longrightarrow}}\lim\ \rS^m_{\hoR^\oy_Y}(\cF_Y^{\oy,(r)}).
\end{equation}

On considère les notations analogues pour $f'$, que l'on munit d'un exposant $^\prime$.

\subsection{}\label{ahttfg3}
Pour tout objet $U$ de $\Et_{/X}$, on pose $\ocB_{U}=\ocB\circ \iota_{U}$  
\eqref{TFA2d} et $\ocB_{U,n}=\ocB_U/p^n\ocB_U$ \eqref{ahttf3e}.
Suivant \ref{ahttf35}, à tout objet $Y$ de $\bP$ tel que $Y_s$ soit non vide
et à toute carte adéquate $((P,\gamma),(\mN,\iota),\vartheta)$ pour le morphisme $f|Y\colon (Y,\cM_X|Y)\rightarrow (S,\cM_S)$ induit par $f$,
on associe une suite exacte canonique de $\ocB_{Y,n}$-modules de $\oY^\circ_\fet$ \eqref{ahttf35a}
\begin{equation}\label{ahttfg3a}
0\rightarrow \ocB_{Y,n}\rightarrow \cF^{(r)}_{Y,n}\rightarrow 
\txi^{-1}\tOmega^1_{X/S}(Y)\otimes_{\co_X(Y)}\ocB_{Y,n} \rightarrow 0,
\end{equation}
et une $\ocB_{Y,n}$-algèbre de $\oY^\circ_\fet$ \eqref{ahttf35b}
\begin{equation}\label{ahttfg3b}
\cC^{(r)}_{Y,n}=\underset{\underset{m\geq 0}{\longrightarrow}}\lim\ \rS^m_{\ocB_{Y,n}}(\cF^{(r)}_{Y,n}).
\end{equation}
Ces objets sont définis par réduction modulo $p^n$ de ceux définis dans \ref{ahttfg2} (cf. \ref{ahttf35} pour plus de détails). 
On pose $\cF_{Y,n}=\cF^{(0)}_{Y,n}$ et $\cC_{Y,n}=\cC^{(0)}_{Y,n}$ (cf. \ref{ahttf7}). On notera que si $Y$ est un objet de $\bQ$, 
$\cF^{(r)}_{Y,n}$ et $\cC^{(r)}_{Y,n}$ ne dépendent pas de la carte adéquate. 

Pour tout objet $Y$ de $\bP$ tel que $Y_s$ soit vide, on pose $\cC^{(r)}_{Y,n}=\cF^{(r)}_{Y,n}=0$ \eqref{ahttf9}. 

Suivant  \ref{ahttf37} et compte tenu de \ref{cftf10}, on considère les faisceaux associés dans $\tE$
\begin{eqnarray}
\cF^{(r)}_n&=&\{Y\in \bQ^\circ\mapsto \cF^{(r)}_{Y,n}\}^a,\label{ahttfg3c}\\
\cC^{(r)}_n&=&\{Y\in \bQ^\circ \mapsto \cC^{(r)}_{Y,n}\}^a. \label{ahttfg3d}
\end{eqnarray}
D'après \ref{ahttf38}, on a une suite exacte localement scindée canonique de $\ocB_n$-modules 
\begin{equation}\label{ahttfg3e}
0\rightarrow \ocB_n\rightarrow \cF^{(r)}_n\rightarrow 
\sigma_n^*(\txi^{-1}\tOmega^1_{\oX_n/\oS_n})\rightarrow 0,
\end{equation}
et un isomorphisme canonique de $\ocB_n$-algèbres  
\begin{equation}\label{ahttfg3f}
\cC^{(r)}_n \stackrel{\sim}{\rightarrow}\underset{\underset{m\geq 0}{\longrightarrow}}\lim\ \rS^m_{\ocB_n}(\cF^{(r)}_n),
\end{equation}
où les morphismes de transition du système inductif sont induits par \eqref{ahttfg3e}. On pose $\cF_n=\cF^{(0)}_n$ et $\cC_n=\cC^{(0)}_n$ 
(cf. \ref{ahttf37}).

Pour tous nombres rationnels $r\geq r'\geq 0$, on a un morphisme $\ocB_n$-linéaire canonique \eqref{ahttf37d}
\begin{equation}\label{ahttfg3g}
\tta_n^{r,r'}\colon \cF^{(r)}_n\rightarrow \cF_n^{(r')},
\end{equation}
et un homomorphisme canonique de $\ocB_n$-algèbres \eqref{ahttf37e}
\begin{equation}\label{ahttfg3h}
\alpha_n^{r,r'}\colon \cC_n^{(r)}\rightarrow \cC_n^{(r')}.
\end{equation}
Pour tous nombres rationnels $r\geq r'\geq r''\geq 0$, on a
\begin{equation}\label{ahttfg3i}
\tta_n^{r,r''}=\tta_n^{r',r''} \circ \tta_n^{r,r'} \ \ \ {\rm et}\ \ \ \alpha_n^{r,r''}=\alpha_n^{r',r''} \circ \alpha_n^{r,r'}.
\end{equation}
On renvoie à \ref{ahttf38} pour les propriétés de ces morphismes. 

On considère les notations analogues pour $f'$, que l'on munit d'un $^\prime$.

\subsection{}\label{ahttfg4}
Soient $(\upmu\colon Y'\rightarrow Y)$ un objet de $\mJ$ \eqref{ahttfg103} tel que $Y'_s$ soit non vide, 
$((P,\gamma),(\mN,\iota),\vartheta)$ une carte adéquate pour le morphisme $f|Y\colon (Y,\cM_X|Y)\rightarrow (S,\cM_S)$ 
induit par $f$, $\oy'$ un point géométrique de $\oY'^\rhd$ \eqref{hmdf3c}.
On note abusivement $\oupmu\colon \oY'^\rhd\rightarrow \oY^\circ$ le morphisme induit par $\upmu$ et on pose $\oy=\oupmu(\oy')$. 
Reprenons les notations de \ref{ahttfg2} pour $(Y,\oy)$ et $(Y',\oy')$.
D'après \ref{mtht101}, on a un morphisme $\hoR^{\oy}_{Y}$-linéaire et $\pi_1(\oY'^{\star\rhd},\oy')$-équivariant \eqref{mtht101c}
\begin{equation}\label{ahttfg4a}
\cF_Y^\oy\rightarrow \cF'^{\oy'}_{Y'}
\end{equation}
qui s'insère dans un diagramme commutatif 
\begin{equation}\label{ahttfg4b}
\xymatrix{
0\ar[r]&{\hoR^{\oy}_{Y}}\ar[d]\ar[r]&{\cF_{Y}^{\oy}}\ar[d]\ar[r]&{\txi^{-1}\tOmega^1_{X/S}(Y)\otimes_{\co_{X}(Y)}\hoR^{\oy}_{Y}}\ar[r]\ar[d]&0\\
0\ar[r]&{\hoRp_{Y'}^{\oy'}}\ar[r]&{\cF'^{\oy'}_{Y'}}\ar[r]&{\txi^{-1}\tOmega^1_{X'/S}(Y')\otimes_{\co_{X'}(Y')}\hoRp^{\oy'}_{Y'}}\ar[r]&0}
\end{equation}
On en déduit un morphisme $\pi_1(\oY'^{\star\rhd},\oy')$-équivariant de $\hoR^{\oy}_{Y}$-algèbres \eqref{mtht101d}
\begin{equation}\label{ahttfg4c}
\cC_Y^\oy\rightarrow \cC'^{\oy'}_{Y'}
\end{equation}
qui prolonge \eqref{ahttfg4a}.

Le morphisme \eqref{ahttfg4a} induit un morphisme $\hoRp^{\oy'}_{Y'}$-linéaire $\pi_1(\oY'^{\star\rhd},\oy')$-équivariant \eqref{ahttfg2c}
\begin{equation}\label{ahttfg4m}
\cF^{\oy,(r)}\otimes_{\hoR^{\oy}_{Y}}\hoRp^{\oy'}_{Y'}\rightarrow \cF'^{\oy',(r)}
\end{equation}
et par suite un homomorphisme de $\hoR^{\oy}_{Y}$-algèbres $\pi_1(\oY'^{\star\rhd},\oy')$-équivariant \eqref{ahttfg2d}
\begin{equation}\label{ahttfg4l}
\cC^{\oy,(r)}_Y\rightarrow \cC'^{\oy',(r)}_{Y'}.
\end{equation}

On désigne par $\Pi(\oY^{\star \circ})$ et $\Pi(\oY'^{\star \rhd})$ les groupoïdes fondamentaux de $\oY^{\star \circ}$ et 
$\oY'^{\star \rhd}$ et par
\begin{equation}\label{ahttfg4d}
\oupmu_*\colon \Pi(\oY'^{\star \rhd})\rightarrow \Pi(\oY^{\star \circ})
\end{equation}
le foncteur induit par le foncteur image inverse $\Et_{\rf/\oY^{\star \circ}}\rightarrow \Et_{\rf/\oY'^{\star \rhd}}$ par $\oupmu$.
On note 
\begin{equation}
F_{Y,n}\colon \Pi(\oY^{\star \circ})\rightarrow \Ens\ \ \ {\rm et}\ \ \  F'_{Y',n}\colon \Pi(\oY'^{\star \rhd})\rightarrow \Ens
\end{equation} 
les foncteurs associés par (\cite{agt} VI.9.11) aux objets $\cF_{Y,n}|\oY^{\star \circ}$ de $\oY^{\star \circ}_\fet$ et 
$\cF'_{Y',n}|\oY'^{\star \rhd}$ de $\oY'^{\star \rhd}_\fet$, respectivement \eqref{ahttfg3}. 
Le morphisme \eqref{ahttfg4a} induit clairement un morphisme de foncteurs
\begin{equation}\label{ahttfg4e}
F_{Y,n}\circ \oupmu_* \rightarrow F'_{Y',n}. 
\end{equation} 
On en déduit par (\cite{agt} VI.9.11) un morphisme $\oupmu^*(\ocB_{Y,n})$-linéaire de $\oY'^\rhd_\fet$,
\begin{equation}\label{ahttfg4f}
v_{\upmu,n}\colon \oupmu^*(\cF_{Y,n}) \rightarrow \cF'_{Y',n}.
\end{equation} 
Il résulte de \eqref{ahttfg4b} que le diagramme 
\begin{equation}\label{ahttfg4h}
\xymatrix{
0\ar[r]&{\oupmu^*(\ocB_{Y,n})}\ar[r]\ar[d]&{\oupmu^*(\cF_{Y,n})}\ar[r]\ar[d]_{v_{\upmu,n}}&
{\txi^{-1}\tOmega^1_{X/S}(Y)\otimes_{\co_X(Y)}\oupmu^*(\ocB_{Y,n})}\ar[r]\ar[d]&0\\
0\ar[r]&{\ocB'_{Y',n}}\ar[r]&{\cF'_{Y',n}}\ar[r]&
{\txi^{-1}\tOmega^1_{X'/S}(Y')\otimes_{\co_{X'}(Y')}\ocB'_{Y',n}}\ar[r]&0}
\end{equation}
est commutatif. 

Le morphisme \eqref{ahttfg4f} induit un morphisme $\oupmu^*(\ocB_{Y,n})$-linéaire de $\oY'^\rhd_\fet$,
\begin{equation}\label{ahttfg4i}
v^{(r)}_{\upmu,n}\colon \oupmu^*(\cF^{(r)}_{Y,n})\rightarrow \cF'^{(r)}_{Y',n}
\end{equation}
qui s'insère dans un diagramme commutatif 
\begin{equation}\label{ahttfg4j}
\xymatrix{
0\ar[r]&{\oupmu^*(\ocB_{Y,n})}\ar[d]\ar[r]&{\oupmu^*(\cF^{(r)}_{Y,n})}\ar[d]_{v^{(r)}_{\upmu,n}}\ar[r]&
{\txi^{-1}\tOmega^1_{X/S}(Y)\otimes_{\co_X(Y)}\oupmu^*(\ocB_{Y,n})}\ar[r]\ar[d]&0\\
0\ar[r]&{\ocB'_{Y',n}}\ar[r]&{\cF'^{(r)}_{Y',n}}\ar[r]&{\txi^{-1}\tOmega^1_{X'/S}(Y')\otimes_{\co_{X'}(Y')}\ocB'_{Y',n}}\ar[r]&0}
\end{equation}
On en déduit un morphisme de $\oupmu^*(\ocB_{Y,n})$-algèbres de $\oY'^\rhd_\fet$
\begin{equation}\label{ahttfg4k}
w^{(r)}_{\upmu,n}\colon \oupmu^*(\cC^{(r)}_{Y,n})\rightarrow \cC'^{(r)}_{Y',n}.
\end{equation}

On notera que les morphismes \eqref{ahttfg4i} et \eqref{ahttfg4k} dépendent du choix de la carte adéquate $((P,\gamma),(\mN,\iota),\vartheta)$. 
Toutefois, ils n'en dépendent pas si $Y$ est un objet de $\bQ$ d'après \ref{pmh7}.

\begin{rema}\label{ahttfg104}
Supposons que le morphisme $g$ \eqref{hmdf3a} soit étale et strict.
Soit $(\upmu\colon Y'\rightarrow Y)$ un objet de $\mI$ \eqref{ahttfg103} tel que $Y'$ et $Y$ soient des objets de $\bQ$. 
Il résulte de \ref{ahttf36} et \ref{ahttfg4} que le morphisme $\oupmu^*(\ocB_{Y,n})$-linéaire de $\oY'^\circ_\fet$ \eqref{ahttfg4i}
\begin{equation}\label{ahttfg104a}
v^{(r)}_{\upmu,n}\colon \oupmu^*(\cF^{(r)}_{Y,n})\rightarrow \cF'^{(r)}_{Y',n}
\end{equation}
n'est autre que le morphisme \eqref{ahttf36a}, c'est-à-dire l'adjoint du morphisme de transition du préfaisceau  
$\{U\in \bQ^\circ\mapsto \cF^{(r)}_{U,n}\}$ sur $E_\bQ$  \eqref{ahttf37}. 
\end{rema}

\subsection{}\label{ahttfg5}
Supposons que le schéma $X$ soit séparé et que le morphisme $f\colon (X,\cM_X)\rightarrow (S,\cM_S)$ 
admette une carte adéquate $((P,\gamma),(\mN,\iota),\vartheta)$ que l'on fixe. 
Pour tout objet $Y$ de $\bP$, équipant le morphisme $f|Y\colon (Y,\cM_X|Y)\rightarrow (S,\cM_S)$ 
induit par $f$ de la carte adéquate induite par $((P,\gamma),(\mN,\iota),\vartheta)$, 
on définit un $\ocB_{Y,n}$-module $\cF^{(r)}_{Y,n}$ et une $\ocB_{Y,n}$-algèbre $\cC^{(r)}_{Y,n}$ \eqref{ahttfg3}. 
D'après \ref{ahttfg4}, les correspondances 
\begin{equation}\label{ahttfg5k}
\{Y\in \bP^\circ\mapsto \cF^{(r)}_{Y,n} \} \ \ \ {\rm et}\ \ \ \{Y\in \bP^\circ\mapsto \cC^{(r)}_{Y,n}\}
\end{equation} 
définissent des préfaisceaux  sur $E_\bP$ \eqref{ahttfg1b} de modules et d'algèbres, respectivement, 
relativement à l'anneau $\{Y\in \bP^\circ \mapsto \ocB_{Y,n}\}$.
En vertu de \ref{cftf11}(i), on a des isomorphismes canoniques 
\begin{eqnarray}
\cF^{(r)}_n&\stackrel{\sim}{\rightarrow}&\{Y\in \bP^\circ\mapsto \cF^{(r)}_{Y,n}\}^a,\label{ahttfg5i}\\
\cC^{(r)}_n&\stackrel{\sim}{\rightarrow}&\{Y\in \bP^\circ\mapsto \cC^{(r)}_{Y,n}\}^a.\label{ahttfg5j}
\end{eqnarray}

Pour tout objet $(\upmu\colon Y'\rightarrow Y)$ de $\mJ$ \eqref{ahttfg103}, notant abusivement $\oupmu\colon \oY'^\rhd\rightarrow \oY^\circ$
le morphisme induit par $\upmu$, on a un morphisme canonique $\oupmu^*(\ocB_{Y,n})$-linéaire de $\oY'^\rhd_\fet$ \eqref{ahttfg4i}
\begin{equation}\label{ahttfg5a}
v^{(r)}_{\upmu,n}\colon \oupmu^*(\cF^{(r)}_{Y,n})\rightarrow \cF'^{(r)}_{Y',n},
\end{equation}
et un morphisme canonique de $\oupmu^*(\ocB_{Y,n})$-algèbres de $\oY'^\rhd_\fet$ \eqref{ahttfg4k}
\begin{equation}\label{ahttfg5b}
w^{(r)}_{\upmu,n}\colon \oupmu^*(\cC^{(r)}_{Y,n})\rightarrow \cC'^{(r)}_{Y',n}.
\end{equation}

Soit $(m',m)\colon (\upmu_1\colon Y'_1\rightarrow Y_1)\rightarrow (\upmu_2\colon Y'_2\rightarrow Y_2)$ un morphisme de $\mJ$.
\begin{equation}\label{ahttfg5c}
\xymatrix{
Y'_1\ar[r]^{\upmu_1}\ar[d]_{m'}&Y_1\ar[d]^m\\
Y'_2\ar[r]^{\upmu_2}&Y_2}
\end{equation}
Avec les notations de \ref{mtht3}, on vérifie aussitôt que le diagramme 
\begin{equation}\label{ahttfg5d}
\xymatrix{
{(\tmY'_1,\cM'_{\tmY'_1})}\ar[rrr]\ar[ddd]&&&{(\tmY_1,\cM_{\tmY_1})}\ar[ddd]\\
&{(\hmY'_1,\cM_{\hmY'_1})}\ar[r]\ar[d]\ar[lu]_{i'_{Y'_1}}&{(\hmY_1,\cM_{\hmY_1})}\ar[d]\ar[ru]^{i_{Y_1}}&\\
&{(\hmY'_2,\cM_{\hmY'_2})}\ar[r]\ar[ld]_{i'_{Y'_2}}&{(\hmY_2,\cM_{\hmY_2})}\ar[rd]^{i_{Y_2}}&\\
{(\tmY'_2,\cM'_{\tmY'_2})}\ar[rrr]&&&{(\tmY_2,\cM_{\tmY_2})}}
\end{equation}
est commutatif (cf. \eqref{mtht101b}). On en déduit péniblement que le diagramme 
\begin{equation}\label{ahttfg5e}
\xymatrix{
{\om'^{\rhd*}(\oupmu_2^*(\cF^{(r)}_{Y_2,n}))}\ar[rr]^-(0.5){\om'^{\rhd*}(v^{(r)}_{\upmu_2,n})}\ar[d]_c&&{\om'^{\rhd*}(\cF'^{(r)}_{Y'_2,n})}\ar[d]^{a'}\\
{\oupmu_1^*(\om^{\circ*}(\cF^{(r)}_{Y_2,n}))}\ar[r]^-(0.5){\oupmu_1^*(a)}&{\oupmu_1^*(\cF^{(r)}_{Y_1,n})}\ar[r]^-(0.5){v^{(r)}_{\upmu_1,n}}&{\cF'^{(r)}_{Y'_1,n}}}
\end{equation}
où $a\colon \om^{\circ*}(\cF^{(r)}_{Y_2,n})\rightarrow \cF^{(r)}_{Y_1,n}$ (resp. $a'$) est le morphisme de transition du v-préfaisceau 
$\{Y\in \bP^\circ\mapsto \cF^{(r)}_{Y,n}\}$ (resp. $\{Y'\in \bQ'^\circ\mapsto \cF'^{(r)}_{Y',n}\}$) \eqref{cftf3}  et $c$ est l'isomorphisme induit par 
le diagramme commutatif \eqref{ahttfg5c}, est commutatif. On en déduit que le diagramme 
\begin{equation}\label{ahttfg5f}
\xymatrix{
{\om'^{\rhd*}(\oupmu_2^*(\cC^{(r)}_{Y_2,n}))}\ar[rr]^-(0.5){\om'^{\rhd*}(w^{(r)}_{\upmu_2,n})}\ar[d]_c&&{\om'^{\rhd*}(\cC'^{(r)}_{Y'_2,n})}\ar[d]^{b'}\\
{\oupmu_1^*(\om^{\circ*}(\cC^{(r)}_{Y_2,n}))}\ar[r]^-(0.5){\oupmu_1^*(b)}&{\oupmu_1^*(\cC^{(r)}_{Y_1,n})}\ar[r]^-(0.5){w^{(r)}_{\upmu_1,n}}
&{\cC'^{(r)}_{Y'_1,n}}}
\end{equation}
où $b\colon \om^{\circ*}(\cC^{(r)}_{Y_2,n})\rightarrow \cC^{(r)}_{Y_1,n}$ (resp. $b'$) est le morphisme de transition du v-préfaisceau 
$\{Y\in \bP^\circ\mapsto \cC^{(r)}_{Y,n}\}$ (resp. $\{Y'\in \bQ'^\circ\mapsto \cC'^{(r)}_{Y',n}\}$) et $c$ est l'isomorphisme induit par 
le diagramme commutatif \eqref{ahttfg5c}, est commutatif.

Avec la terminologie de \ref{cftf122}, les morphismes $v^{(r)}_{\upmu,n}$ \eqref{ahttfg5a} (resp. $w^{(r)}_{\upmu,n}$ \eqref{ahttfg5b})
forment donc un $\mJ$-système de $\Phi$-morphismes compatibles  de $\{Y\in \bP^\circ\mapsto  \cF^{(r)}_{Y,n}\}$ 
dans $\{Y'\in \bQ'^\circ\mapsto  \cF'^{(r)}_{Y',n}\}$ (resp. $\{Y\in \bP^\circ\mapsto  \cC^{(r)}_{Y,n}\}$ 
dans $\{Y'\in \bQ'^\circ\mapsto  \cC'^{(r)}_{Y',n}\}$) \eqref{ahttfg100c}. 
Par suite, en vertu de \ref{cftf14} et \ref{ahttfg6}, ils  définissent des morphismes de $\tE'$
\begin{eqnarray}
v^{(r)}_n\colon\Theta^{-1}(\cF^{(r)}_n)&\rightarrow& \cF'^{(r)}_n,\label{ahttfg5g}\\
w^{(r)}_n\colon\Theta^{-1}(\cC^{(r)}_n)&\rightarrow& \cC'^{(r)}_n,\label{ahttfg5h}
\end{eqnarray}
où $\Theta$ est le morphisme de topos \eqref{hmdf11b}.

\begin{rema}\label{ahttfg105}
Conservons les hypothèses et notations de \ref{ahttfg5}.
Il résulte de \ref{ahttfg104} que si le morphisme $g$ \eqref{hmdf3a} est étale et strict, le morphisme $v^{(r)}_n$ \eqref{ahttfg5g} 
n'est autre que l'isomorphisme \eqref{ahttf41j}. Par suite, $w^{(r)}_n$ \eqref{ahttfg5h} est aussi un isomorphisme.
\end{rema}

\begin{lem}\label{ahttfg7}
Reprenons les hypothèses et notations de \ref{ahttfg5}. Alors, 
\begin{itemize}
\item[{\rm (i)}] Le morphisme $v^{(r)}_n$ \eqref{ahttfg5g} est $\Theta^{-1}(\ocB_n)$-linéaire. Il  induit donc un morphisme $\ocB'_n$-linéaire de $\tE'_s$
\begin{equation}\label{ahttfg7a}
\upnu^{(r)}_n\colon \uptheta^*_n(\cF^{(r)}_n)\rightarrow \cF'^{(r)}_n,
\end{equation}
où $\uptheta_n$ est le morphisme de topos annelés \eqref{hmdf12d}. 
\item[{\rm (ii)}]  On a un diagramme commutatif 
\begin{equation}\label{ahttfg7b}
\xymatrix{
0\ar[r]&{\ocB'_n}\ar@{=}[d]\ar[r]&{\uptheta^*_n(\cF^{(r)}_n)}\ar[d]_{\upnu^{(r)}_n}\ar[r]&
{\uptheta^*_n(\sigma^*_n(\txi^{-1}\tOmega^1_{\oX_n/\oS_n}))}\ar[r]\ar[d]&0\\
0\ar[r]&{\ocB'_n}\ar[r]&{\cF'^{(r)}_n}\ar[r]&{\sigma'^*_n(\txi^{-1}\tOmega^1_{\oX'_n/\oS_n})}\ar[r]&0}
\end{equation}
où les lignes horizontales sont induites par les suites exactes \eqref{ahttfg3e} et la flèche verticale de droite est induite par le morphisme canonique
\begin{equation}\label{ahttfg7c}
\ogg_n^*(\tOmega^1_{\oX_n/\oS_n})\rightarrow  \tOmega^1_{\oX'_n/\oS_n}
\end{equation}
et le diagramme commutatif \eqref{hmdf12e}. 
\item[{\rm (iii)}] Le morphisme $w^{(r)}_n$ \eqref{ahttfg5h} est un homomorphisme de $\Theta^{-1}(\ocB_n)$-algèbres. 
Il induit donc un morphisme de $\ocB'_n$-algèbres de $\tE'_s$
\begin{equation}\label{ahttfg7d}
\upomega^{(r)}_n\colon \uptheta^*_n(\cC^{(r)}_n)\rightarrow \cC'^{(r)}_n.
\end{equation}
\item[{\rm (iv)}] Le diagramme 
\begin{equation}\label{ahttfg7e}
\xymatrix{
{\uptheta^*_n(\cF^{(r)}_n)}\ar[r]^-(0.5){\upnu^{(r)}_n}\ar[d]&{\cF'^{(r)}_n}\ar[d]\\
{\uptheta^*_n(\cC^{(r)}_n)}\ar[r]^-(0.5){\upomega^{(r)}_n}&{\cC'^{(r)}_n}}
\end{equation}
où les flèches verticales sont les morphismes canoniques, est commutatif. 
\end{itemize}
\end{lem}

(i) Posons $F=\{U\in \bP^\circ\mapsto \cF^{(r)}_{U,n}\}$ qui est un v-préfaisceau sur $E_\bP$. 
D'après \ref{cftf121} et sa preuve, il existe un v-préfaisceau $F'=\{U'\in \Et^\circ_{\coh/X'}\mapsto F'_{U'}\}$ sur $E'_\coh$
et un $\mI$-système de $\Phi$-morphismes compatibles de $F$ dans $F'$,  
\begin{equation}
v'_\upmu\colon \oupmu^*(\cF^{(r)}_{U,n})\rightarrow F'_{U'}, \ \ \ \forall \  (\upmu\colon U'\rightarrow U)\in \ob(\mI),
\end{equation}
où $\oupmu \colon \oU'^\rhd\rightarrow \oU^\circ$ désigne abusivement le morphisme induit par $\upmu$, 
tels que pour tous objets $U'$ de $\bQ'$ et $(\upmu\colon U'\rightarrow U)$ de $\mJ$, 
on ait $F'_{U'}=\cF'^{(r)}_{U',n}$ et $v'_\upmu=v^{(r)}_{\upmu,n}$. 
Par ailleurs, pour tous objets $U'$ de $\Et_{\coh/X'}$ et $(\upmu\colon U'\rightarrow U)$ de $\mI$, on peut supposer que 
que $F'_{U'}$ est un $\ocB'_{U',n}$-module et que le morphisme $v'_\upmu$ est $\oupmu^*(\ocB_{U,n})$-linéaire. 
La proposition s'ensuit compte tenu de \ref{ahttfg101} (cf. \ref{cftf14}). 

(ii) On a un isomorphisme canonique \eqref{ahttf12c}
\begin{equation}
\sigma^*_n(\txi^{-1}\tOmega^1_{\oX_n/\oS_n})\stackrel{\sim}{\rightarrow}
\sigma^{-1}(\txi^{-1}\tOmega^1_{X/S})\otimes_{\sigma^{-1}(\co_X)}\ocB_n.
\end{equation}
Donc en vertu de (\cite{agt} VI.5.34(ii), VI.8.9 et VI.5.17), $\sigma^*_n(\xi^{-1}\tOmega^1_{\oX_n/\oS_n})$ est 
le faisceau de $\tE$ associé au préfaisceau sur $E$ défini par la correspondance
\begin{equation}
\{U\in \Et^\circ_{/X} \mapsto \xi^{-1}\tOmega^1_{X/S}(U)\otimes_{\co_X(U)}\ocB_{U,n}\}.
\end{equation}

D'après \ref{ahttfg19}, pour tout objet $(\upmu\colon U'\rightarrow U)$ de $\mI$ \eqref{ahttfg103}, 
notant abusivement $\oupmu\colon \oU'^\rhd\rightarrow \oU^\circ$ le morphisme induit par $\upmu$,  
on a un morphisme canonique de $\oU'^\rhd_\fet$, 
\begin{equation}
\txi^{-1}\tOmega^1_{X/S}(U)\otimes_{\co_X(U)}\oupmu^*(\ocB_{U,n})\rightarrow \txi^{-1}\tOmega^1_{X'/S}(U')\otimes_{\co_{X'}(U')}\ocB'_{U',n}. 
\end{equation}
Ces morphismes forment un 
$\mI$-système de $\Phi$-morphismes compatibles de 
\[
\{U\in \bP^\circ \mapsto \xi^{-1}\tOmega^1_{X/S}(U)\otimes_{\co_X(U)}\ocB_{U,n}\}\ 
{\rm dans}\  \{U'\in \Et^\circ_{\coh/X'} \mapsto \xi^{-1}\tOmega^1_{X'/S}(U')\otimes_{\co_{X'}(U')}\ocB_{U',n}\}.
\] 
Par suite, en vertu de \ref{cftf11}(i) et  \ref{cftf14}, ils  définissent un morphisme de $\tE'$
\begin{equation}
\Theta^{-1}(\sigma^*_n(\txi^{-1}\tOmega^1_{\oX_n/\oS_n}))\rightarrow 
\sigma'^*_n(\txi^{-1}\tOmega^1_{\oX'_n/\oS_n}),
\end{equation}
qui n'est autre que le morphisme induit par le morphisme canonique \eqref{ahttfg7c}
et le diagramme commutatif \eqref{hmdf12e}. La proposition s'ensuit compte tenu de \eqref{ahttfg4j}. 

(iii) Posons $G=\{U\in \bP^\circ\mapsto \cC^{(r)}_{U,n}\}$ qui est un v-préfaisceau sur $E_\bP$.
D'après \ref{cftf121} et sa preuve, il existe un v-préfaisceau $G'=\{U'\in \Et^\circ_{\coh/X'}\mapsto G'_{U'}\}$ sur $E'_\coh$ 
et un $\mI$-système de $\Phi$-morphismes compatibles de $G$ dans $G'$,  
\begin{equation}
v'_\upmu\colon \oupmu^*(\cC^{(r)}_{U,n})\rightarrow G'_{U'}, \ \ \ \forall \  (\upmu\colon U'\rightarrow U)\in \ob(\mI),
\end{equation}
où $\oupmu \colon \oU'^\rhd\rightarrow \oU^\circ$ désigne abusivement le morphisme induit par $\upmu$, 
tels que pour tous objets $U'$ de $\bQ'$ et $(\upmu\colon U'\rightarrow U)$ de $\mJ$, on ait $G'_{U'}=\cC'^{(r)}_{U',n}$ et 
$v'_\upmu=w^{(r)}_{\upmu,n}$. 
Par ailleurs, pour tous objets $U'$ de $\Et_{\coh/X'}$ et $(\upmu\colon U'\rightarrow U)$ de $\mI$, on peut supposer que 
$G'_{U'}$ est une $\ocB'_{U',n}$-algèbre et que le morphisme $v'_\upmu$ est un homomorphisme 
de $\oupmu^*(\ocB_{U,n})$-algèbres. La proposition s'ensuit compte tenu de \ref{ahttfg101} (cf. \ref{cftf14}). 

(iv) Cela résulte de \ref{ahttfg101} (cf. \ref{cftf14}) compte tenu de la définition de \eqref{ahttfg4k}.

\begin{lem}\label{ahttfg10}
Supposons que le schéma $X$ soit séparé et que le morphisme $f\colon (X,\cM_X)\rightarrow (S,\cM_S)$ 
admette une carte adéquate $((P,\gamma),(\mN,\iota),\vartheta)$ que l'on fixe.   
Soit $(\upmu\colon U'\rightarrow U)$ un objet de $\mI$ \eqref{ahttfg103}.
Notons  abusivement $\oupmu\colon \oU'^\rhd\rightarrow \oU^\circ$ le morphisme induit par $\upmu$ et reprenons les notations de \ref{hmdf20}.
Alors, il existe un morphisme canonique 
\begin{equation}\label{ahttfg10a}
\oupmu^*(\cF^{(r)}_{U,n}) \rightarrow \beta'_{U'*}(\jmath'^*_{U'}(\Theta^{-1}(\cF^{(r)}_n))).
\end{equation}
Celui-ci est indépendant de la carte adéquate $((P,\gamma),(\mN,\iota),\vartheta)$ si $U$ est un objet de $\bQ$. 

Si $(\upmu\colon U'\rightarrow U)$ est un objet de $\mJ$, ce morphisme s'insère dans un diagramme commutatif 
\begin{equation}\label{ahttfg10b}
\xymatrix{
{\oupmu^*(\cF^{(r)}_{U,n})}\ar[rr]^-(0.5){v^{(r)}_{\upmu,n}}\ar[d]&&{\cF'^{(r)}_{U',n}}\ar[d]\\
{\beta'_{U'*}(\jmath'^*_{U'}(\Theta^{-1}(\cF^{(r)}_n)))}\ar[rr]^-(0.5){\beta'_{U'*}(\jmath'^*_{U'}(v^{(r)}_n))}&&{\beta'_{U'*}(\jmath'^*_{U'}(\cF'^{(r)}_n))}}
\end{equation}
où la flèche verticale de droite est le morphisme \eqref{ahttfg9a} et 
$v^{(r)}_{\upmu,n}$ (resp. $v^{(r)}_n$) est le morphisme \eqref{ahttfg5a} (resp. \eqref{ahttfg5g}).
\end{lem}

Posons $F=\{Y\in \bP^\circ\mapsto \cF^{(r)}_{Y,n}\}$ qui est un v-préfaisceau sur $E_\bP$, et pour tout $Y'\in \ob(\Et_{\coh/X'})$,
\begin{equation}\label{ahttfg10c}
F'_{Y'}=\underset{\underset{(Y,\uplambda)\in (I^{Y'}_\varphi)^\circ}{\longrightarrow}}{\lim}\ \ouplambda^*(\cF^{(r)}_{Y,n}),
\end{equation}
où pour tout objet $(Y,\uplambda)$ de $I^{Y'}_\varphi$  \eqref{ahttfg103}, on a abusivement noté $\ouplambda\colon \oY'^\rhd\rightarrow \oY^\circ$ 
le morphisme induit par $\uplambda\colon Y'\rightarrow Y$. 
D'après \ref{ahttfg101}, la collection $F'=\{Y'\in \Et^\circ_{\coh/X'}\mapsto F'_{Y'}\}$ forme naturellement un v-préfaisceau sur $E'$, 
et on a un isomorphisme canonique fonctoriel 
\begin{equation}\label{ahttfg10d}
\Theta^{-1}(\cF^{(r)}_n)\stackrel{\sim}{\rightarrow} F'^a.
\end{equation} 
On prend alors pour morphisme \eqref{ahttfg10a} le morphisme composé 
\begin{equation}
\oupmu^*(\cF^{(r)}_{U,n}) \rightarrow F'_{U'} \rightarrow \beta'_{U'*}(\jmath'^*_{U'}(F'^a)\stackrel{\sim}{\rightarrow} 
\beta'_{U'*}(\jmath'^*_{U'}(\Theta^{-1}(\cF^{(r)}_n))),
\end{equation}
où la première et la deuxième flèches sont les morphismes canoniques et la dernière flèche est induite par l'isomorphisme \eqref{ahttfg10d}.
Il résulte de \ref{ahttfg101}(iv) et \ref{ahttfg9} que ce morphisme est indépendant de la carte adéquate 
$((P,\gamma),(\mN,\iota),\vartheta)$ si $U$ est un objet de $\bQ$. 

La seconde assertion résulte de la définition du morphisme $v^{(r)}_n$ \eqref{ahttfg5g} (cf. \ref{cftf14}).

\subsection{}\label{ahttfg110}
Supposons que le schéma $X$ soit séparé et que le morphisme $f\colon (X,\cM_X)\rightarrow (S,\cM_S)$ 
admette une carte adéquate $((P,\gamma),(\mN,\iota),\vartheta)$ que l'on fixe.  
Soient $\ox'$  un point géométrique de $X'$,  $\uX'$ le localisé strict de $X'$ en $\ox'$. 
On désigne par $\tuE'$ le topos de Faltings associé au morphisme canonique $\uoX'^\rhd\rightarrow \uX'$, par 
\begin{equation}\label{ahttfg110a}
\Phi'\colon \tuE'\rightarrow \tE'
\end{equation}
le morphisme de fonctorialité induit par le morphisme canonique $\uX'\rightarrow X'$, par 
\begin{equation}\label{ahttfg110b}
\ubeta'\colon \tuE'\rightarrow \uoX'^\rhd_\fet
\end{equation}
le morphisme canonique \eqref{ahttf2d} et par 
\begin{equation}\label{ahttfg110c}
\theta'\colon \uoX'^\rhd_\fet\rightarrow \tuE'
\end{equation}
la section canonique de $\ubeta'$ (\cite{agt} VI.10.23). On note 
\begin{equation}\label{ahttfg110d}
\varphi'_{\ox'}\colon \tE'\rightarrow \uoX'^\rhd_\fet
\end{equation}
le foncteur composé $\theta'^*\circ \Phi'^*$.  

On désigne par $\fV'_{\ox'}$ (resp. $\bP'_{\ox'}$, resp. $\bQ'_{\ox'}$) la catégorie des voisinages de $\ox'$ 
dans le site $\Et_{/X'}$ (resp. $\bP'$, resp. $\bQ'$) \eqref{ahttf47}, 
et par $\mI_{\ox'}$ la catégorie des couples de morphismes $(\iota'\colon \ox'\rightarrow U', \upmu\colon U'\rightarrow U)$, 
où $\iota'$ est un $X'$-morphisme et $\upmu$ est un objet de $\mI$ \eqref{ahttfg103}.

Soit $(\iota'\colon \ox'\rightarrow U', \upmu\colon U'\rightarrow U)$ un objet de $\mI_{\ox'}$. 
On désigne encore par $\iota'\colon \uX'\rightarrow U'$ le $X'$-morphisme déduit de $\iota'$ (\cite{sga4} VIII 7.3)
et on note abusivement $\oiota'\colon \uoX'^\rhd\rightarrow \oU'^\rhd$ le morphisme induit.
Reprenons les notations de \ref{hmdf20}. 
Le morphisme $\iota'\colon \uX'\rightarrow U'$ induit par fonctorialité un morphisme 
\begin{equation}\label{ahttfg110e}
\Phi'_{\iota'}\colon \tuE'\rightarrow \tE'_{/\sigma'^*(U')}
\end{equation}  
qui s'insère dans un diagramme commutatif à isomorphisme canonique près
\begin{equation}\label{ahttfg110f}
\xymatrix{
{\tuE'}\ar[r]^-(0.5){\Phi'_{\iota'}}\ar[d]_{\ubeta'}&{\tE'_{/\sigma'^*(U')}}\ar[d]^{\beta'_{U'}}\\
{\uoX'^\rhd_\fet}\ar[r]^-(0.5){\oiota'}&{\oU'^\rhd_\fet}}
\end{equation}
Appliquant le foncteur composé $\theta'^*\circ \Phi'^*_{\iota'}$ au morphisme $\beta'^*_{U'}(\oupmu^*(\cF^{(r)}_{U,n}))
\rightarrow \jmath'^*_{U'}(\Theta^{-1}(\cF_n^{(r)}))$ adjoint de \eqref{ahttfg10a} et tenant compte du diagramme commutatif \eqref{ahttfg110f}, 
on obtient un morphisme 
\begin{equation}\label{ahttfg110g}
a'\colon \oiota'^*(\oupmu^*(\cF^{(r)}_{U,n}))\rightarrow \varphi'_{\ox'}(\Theta^{-1}(\cF^{(r)}_n)).
\end{equation}
Celui-ci est indépendant de la carte adéquate $((P,\gamma),(\mN,\iota),\vartheta)$ si $U$ est un objet de $\bQ$ \eqref{ahttfg10}.

Posons $\ox=g(\ox')$ et notons $\uX$ le localisé strict de $X$ en $\ox$. 
Le $X'$-morphisme $\iota'$ induit un $X$-morphisme 
$\iota=\iota'\circ \upmu\colon \ox\rightarrow U$. Reprenons les notations de \ref{ahttf47} et \ref{ahttfg12}. Comme $(U,\iota)$ est un objet de $\bP_\ox$,
on peut considérer le morphisme \eqref{ahttf47g}
\begin{equation}\label{ahttfg110h}
a\colon \oiota^*(\cF^{(r)}_{U,n})\rightarrow \varphi_\ox(\cF^{(r)}_n).
\end{equation}
Le diagramme \eqref{ahttfg101c} induit un diagramme commutatif 
\begin{equation}\label{ahttfg110i}
\xymatrix{
{\uupgamma^*(\oiota^*(\cF^{(r)}_{U,n}))}\ar[r]^-(0.5){\uupgamma^*(a)}\ar[d]&{\uupgamma^*(\varphi_{\ox}(\cF^{(r)}_n))}\ar[d]\\
{\oiota'^*(\oupmu^*(\cF^{(r)}_{U,n}))}\ar[r]^-(0.5){a'}&{\varphi'_{\ox'}(\Theta^{-1}(\cF^{(r)}_n))}}
\end{equation}
où la flèche verticale de droite (resp. de gauche) est l'isomorphisme \eqref{ahttfg12j} 
(resp. induit par la relation $\oiota\circ \uupgamma=\oupmu\circ \oiota'$).

Si $(\upmu\colon U'\rightarrow U)$ est un objet de $\mJ$ \eqref{ahttfg103}, le digramme \eqref{ahttfg10b} induit un diagramme commutatif
\begin{equation}\label{ahttfg110j}
\xymatrix{
{\oiota'^*(\oupmu^*(\cF^{(r)}_{U,n}))}\ar[rr]^{\oiota'^*(v^{(r)}_{\upmu,n})}\ar[d]_{a'}&&{\oiota'^*(\cF'^{(r)}_{U',n})}\ar[d]^{b'}\\
{\varphi'_{\ox'}(\Theta^{-1}(\cF^{(r)}_n))}\ar[rr]^-(0.5){\varphi'_{\ox'}(v^{(r)}_n)}&&{\varphi'_{\ox'}(\cF'^{(r)}_n)}}
\end{equation}
où $b'$ est le morphisme \eqref{ahttf47g} pour $\cF'^{(r)}_n$, et
$v^{(r)}_n$ (resp.  $v^{(r)}_{\upmu,n}$) est le morphisme \eqref{ahttfg5g} (resp. \eqref{ahttfg5a}).

\begin{prop}\label{ahttfg15}
Supposons que le schéma $X$ soit séparé et que le morphisme $f\colon (X,\cM_X)\rightarrow (S,\cM_S)$ 
admette une carte adéquate $((P,\gamma),(\mN,\iota),\vartheta)$. Alors, le morphisme  
\begin{equation}\label{ahttfg15a}
v^{(r)}_n\colon\Theta^{-1}(\cF^{(r)}_n)\rightarrow \cF'^{(r)}_n
\end{equation}
défini dans \eqref{ahttfg5g} ne dépend pas du choix de la carte adéquate $((P,\gamma),(\mN,\iota),\vartheta)$. 
\end{prop}

Soient $\ox'$ un point géométrique de $X'$, $\ox=g(\ox')$. Reprenons les notations de \ref{ahttf47} et  \ref{ahttfg110}.
On désigne par $\mK_{\ox'}$ la sous-catégorie de $\mI_{\ox'}$ 
formée des objets $(\iota'\colon \ox'\rightarrow U', \upmu\colon U'\rightarrow U)$ tels que $U$ (resp. $U'$) soit un objet de $\bQ$ (resp. $\bQ'$). 
On notera que les limites projectives finies dans $\mI_{\ox'}$ sont représentables (\cite{sga4} I 2.3). 
Par suite, la catégorie $\mI_{\ox'}$ est cofiltrante (\cite{sga4} I 2.7.1).
Le foncteur d'injection canonique $\mK_{\ox'}\rightarrow \mI_{\ox'}$ est 
initial et la catégorie $\mK_{\ox'}$ est cofiltrante en vertu de (\cite{sga4} I 8.1.3(c)) et (\cite{agt} II.5.17). 

Les foncteurs $\tts$ \eqref{ahttfg103a} et $\ttb$ \eqref{ahttfg103b} induisent des foncteurs
\begin{eqnarray}
\tts_{\ox'}\colon \mK_{\ox'}\rightarrow \bQ'_{\ox'},&& (\iota'\colon \ox'\rightarrow U', \upmu\colon U'\rightarrow U)\mapsto (U',\iota'),\label{ahttfg15b}\\
\ttb_{\ox'}\colon \mK_{\ox'}\rightarrow \bQ_\ox,&& (\iota'\colon \ox'\rightarrow U', \upmu\colon U'\rightarrow U)\mapsto (U,\iota=\iota'\circ \upmu\colon \ox\rightarrow U).\label{ahttfg15c}
\end{eqnarray}
Il résulte de (\cite{sga4} I 8.1.3(b)), (\cite{agt} II.5.17) et du fait que les catégories $\bQ'_{\ox'}$ et $\bQ_{\ox}$ sont cofiltrantes 
que les foncteurs $\tts^\circ_{\ox'}$ et $\ttb^\circ_{\ox'}$ sont cofinaux. 
Par suite, compte tenu de \ref{ahttf48}(i), \eqref{ahttfg110i} et \eqref{ahttfg110j}, le morphisme $\varphi'_{\ox'}(v^{(r)}_n)$
s'identifie au morphisme obtenu par passage à la limite inductive sur la catégorie $\mK_{\ox'}^\circ$ du morphisme composé 
\begin{equation}
\xymatrix{
{\uupgamma^*(\oiota^*(\cF^{(r)}_{U,n}))}\ar[r]&{\oiota'^*(\oupmu^*(\cF^{(r)}_{U,n}))}\ar[rr]^-(0.5){\oiota'^*(v^{(r)}_{\upmu,n})}&&{\oiota'^*(\cF'^{(r)}_{U',n})}},
\end{equation}
où la première flèche est l'isomorphisme induit par la relation $\oiota\circ \uupgamma=\oupmu\circ \oiota'$. Comme $U$ (resp. $U'$) est un objet de $\bQ$
(resp. $\bQ'$), le morphisme $v^{(r)}_{\upmu,n}$ ne dépend pas de la carte adéquate $((P,\gamma),(\mN,\iota),\vartheta)$. 
La proposition s'ensuit puisque la famille des foncteurs $\varphi'_{\ox'}$ lorsque $\ox'$ décrit l'ensemble des points géométriques de $X'$ 
est conservative (\cite{agt} VI.10.32).

\subsection{}\label{ahttfg17}
Considérons un diagramme commutatif de morphismes de schémas
\begin{equation}\label{ahttfg17a}
\xymatrix{
U'\ar[r]^{\upmu}\ar[d]_{u'}&U\ar[d]^u\\
X'\ar[r]^g&X}
\end{equation}
tel que $u$ et $u'$ soient étales de présentation finie. On munit $U$ (resp. $U'$) de la structure logarithmique $\cM_{U}$ (resp. $\cM_{U'}$) 
image inverse de $\cM_X$ (resp. $\cM_{X'}$). On lui associe le digramme commutatif de morphismes de schémas 
\begin{equation}\label{ahttfg17b}
\xymatrix{
{(\coU',\cM_{\coU'})}\ar[rrr]^-(0.5){\cou'}\ar[rd]\ar[ddd]_{\coupmu}&&&{(\coX',\cM_{\coX'})}\ar[ddd]^{\cog}\ar[ld]\\
&{(\tU',\cM_{\tU'})}\ar[r]^-(0.5){\tu'}\ar[d]_{\tupmu}&{(\tX',\cM_{\tX'})}\ar[d]^{\tg}&\\
&{(\tU,\cM_{\tU})}\ar[r]^-(0.5){\tu}&{(\tX,\cM_{\tX})}&\\
{(\coU,\cM_{\coU})}\ar[rrr]^-(0.5){\cou}\ar[ru]&&&{(\coX,\cM_{\coX})}\ar[lu]}
\end{equation}
où $\cM_{\coU}$ (resp. $\cM_{\coU'}$) désigne l'image inverse de la structure logarithmique $\cM_{X}$ (resp. $\cM_{X'}$) sur $\coU$ (resp. $\coU'$) \eqref{hmdf1b},
$\tu$ (resp. $\tu'$) est l'unique relèvement étale strict de $\cou$ (resp. $\cou'$) et $\tupmu$ est l'unique relèvement de $\coupmu$ qui rend le carré intérieur 
commutatif. Reprenons les notations de \ref{hmdf20}. 
D'après \ref{ahttf41}, plus précisément \eqref{ahttf41j}, $\jmath_U^*(\cF_n^{(r)})$ s'identifie canoniquement au faisceau analogue pour le morphisme 
adéquat $f\circ u\colon (U,\cM_U)\rightarrow (S,\cM_S)$ muni de la $(\tS,\cM_{\tS})$-déformation $(\tU,\cM_{\tU})$ de $(\coU,\cM_{\coU})$. 
De même, $\jmath'^*_{U'}(\cF'^{(r)}_n)$ s'identifie canoniquement au faisceau analogue pour le morphisme 
adéquat $f'\circ u'\colon (U',\cM_{U'})\rightarrow (S,\cM_S)$ muni de la $(\tS,\cM_{\tS})$-déformation $(\tU',\cM_{\tU'})$ de $(\coU',\cM_{\coU'})$. 

Supposons que le schéma $X$ soit séparé et que le morphisme $f$ admette une carte adéquate. On dispose alors du morphisme \eqref{ahttfg5g}
\begin{equation}\label{ahttfg17c}
v^{(r)}_n\colon\Theta^{-1}(\cF^{(r)}_n)\rightarrow \cF'^{(r)}_n,
\end{equation}
dont nous avons démontré dans \ref{ahttfg15} qu'il ne dépend pas de la carte adéquate pour $f$. 
Supposons de plus que le schéma $U$ soit séparé et que le morphisme $f\circ u$ admette une carte adéquate. 
Le morphisme composé 
\begin{equation}\label{ahttfg17d}
\xymatrix{
{\Theta_{\upmu}^{-1}(\jmath_U^*(\cF_n^{(r)}))}\ar[r]^-(0.5)\sim&{\jmath'^*_{U'}(\Theta^{-1}(\cF^{(r)}_n))}\ar[rr]^-(0.5){\jmath'^*_{U'}(v^{(r)}_n)}&&
{\jmath'^*_{U'}(\cF'^{(r)}_n})},
\end{equation}
où la première flèche est l'isomorphisme sous-jacent au diagramme \eqref{hmdf20d}, s'identifie alors au morphisme analogue à $v^{(r)}_n$ 
défini relativement aux morphismes $\upmu$  et $\tupmu$. En effet, choisissant une carte adéquate pour $f$ et prenant pour $f\circ u$
la carte induite, l'assertion résulte facilement des définitions.

\begin{prop}\label{ahttfg16}
Il existe un morphisme $\Theta^{-1}(\ocB_n)$-linéaire canonique de $\tE'_s$,
\begin{equation}\label{ahttfg16a}
v^{(r)}_n\colon \Theta^{-1}(\cF^{(r)}_n)\rightarrow \cF'^{(r)}_n,
\end{equation}
vérifiant la propriété suivante: pour tout diagramme commutatif de morphismes de schémas
\begin{equation}\label{ahttfg16b}
\xymatrix{
U'\ar[r]^{\upmu}\ar[d]_{u'}&U\ar[d]^u\\
X'\ar[r]^g&X}
\end{equation}
tel que les morphismes $u$ et $u'$ soient étales de présentation finie,
que le schéma $U$ soit séparé et que le morphisme $f|U\colon (U,\cM_X|U)\rightarrow (S,\cM_S)$ induit par $f$ admette une carte adéquate, 
considérant le diagramme \eqref{ahttfg17b} associé à \eqref{ahttfg16b} et reprenant les notations de \ref{hmdf20}, le morphisme composé 
\begin{equation}\label{ahttfg16c}
\xymatrix{
{\Theta_{\upmu}^{-1}(\jmath_U^*(\cF_n^{(r)}))}\ar[r]^-(0.5)\sim&{\jmath'^*_{U'}(\Theta^{-1}(\cF^{(r)}_n))}\ar[rr]^-(0.5){\jmath'^*_{U'}(v^{(r)}_n)}&&
{\jmath'^*_{U'}(\cF'^{(r)}_n})},
\end{equation}
où la première flèche est l'isomorphisme sous-jacent au diagramme \eqref{hmdf20d}, s'identifie grâce à \eqref{ahttf41j},
au morphisme défini dans \eqref{ahttfg5g} relativement aux morphismes $\upmu$  et $\tupmu$ (et à n'importe quelle carte adéquate pour $f|U$).
\end{prop} 

En effet, le morphisme $v^{(r)}_n$ \eqref{ahttfg16a} est défini par descente à partir du morphisme défini dans \eqref{ahttfg5g}, 
compte tenu de \ref{ahttfg15}  et \ref{ahttfg17}.

\subsection{}\label{ahttfg18}
Le morphisme $v^{(r)}_n$ \eqref{ahttfg16a} induit un morphisme $\ocB'_n$-linéaire de $\tE'_s$
\begin{equation}\label{ahttfg18a}
\upnu^{(r)}_n\colon \uptheta^*_n(\cF^{(r)}_n)\rightarrow \cF'^{(r)}_n,
\end{equation}
où $\uptheta_n$ est le morphisme de topos annelés \eqref{hmdf12d}. 
Il résulte aussitôt de \ref{ahttfg7}(ii) qu'on a un diagramme commutatif 
\begin{equation}\label{ahttfg18b}
\xymatrix{
0\ar[r]&{\ocB'_n}\ar@{=}[d]\ar[r]&{\uptheta^*_n(\cF^{(r)}_n)}\ar[d]_{\upnu^{(r)}_n}\ar[r]&
{\uptheta^*_n(\sigma^*_n(\txi^{-1}\tOmega^1_{\oX_n/\oS_n}))}\ar[r]\ar[d]&0\\
0\ar[r]&{\ocB'_n}\ar[r]&{\cF'^{(r)}_n}\ar[r]&{\sigma'^*_n(\txi^{-1}\tOmega^1_{\oX'_n/\oS_n})}\ar[r]&0}
\end{equation}
où les lignes horizontales sont induites par les suites exactes \eqref{ahttfg3e} et la flèche verticale de droite est induite par le morphisme canonique
\begin{equation}\label{ahttfg18c}
\ogg_n^*(\tOmega^1_{\oX_n/\oS_n})\rightarrow  \tOmega^1_{\oX'_n/\oS_n}
\end{equation}
et le diagramme commutatif \eqref{hmdf12e}. 
Compte tenu de \eqref{ahttfg3f}, le morphisme $\upnu^{(r)}_n$ induit un morphisme de $\ocB'_n$-algèbres de $\tE'_s$
\begin{equation}\label{ahttfg18d}
\upomega^{(r)}_n\colon \uptheta^*_n(\cC^{(r)}_n)\rightarrow \cC'^{(r)}_n.
\end{equation}
Dans la suite, on considère $\cC'^{(r)}_n$ comme une $\uptheta^*_n(\cC^{(r)}_n)$-algèbre via $\upomega^{(r)}_n$.

Pour tout nombre rationnel $t\geq r$, on considère la $\uptheta^*_n(\cC^{(r)}_n)$-algèbre
\begin{equation}\label{ahttfg18f}
\cC'^{(t,r)}_n=\cC'^{(t)}_n\otimes_{\uptheta^*_n(\cC^{(t)}_n)}\uptheta^*_n(\cC^{(r)}_n)
\end{equation}
déduite de $\cC'^{(t)}_n$ par changement de base par l'homomorphisme $\uptheta^*_n(\alpha^{t,r}_n)$, où $\alpha^{t,r}_n\colon \cC^{(t)}_n\rightarrow \cC^{(r)}_n$
est l'homomorphisme canonique \eqref{ahttf37e}. 
Reprenant les notations de \eqref{hmdf19l}, l'homomorphisme canonique $\uptheta^*_n(\cC^{(r)}_n)\rightarrow \cC'^{(t,r)}_n$ 
induit un morphisme de $\ocB^!_n$-algèbres de $\tG_s$
\begin{equation}\label{ahttfg18h}
\lgg^*_n(\cC^{(r)}_n)\rightarrow \tau_{n*}(\cC'^{(t,r)}_n).
\end{equation}

Pour tous nombres rationnels  $t,t',r,r'$  tels que $t'\geq t\geq r$ et $t'\geq r'\geq r$, les diagrammes 
\begin{equation}\label{ahttfg18e}
\xymatrix{
{\cC^{(r')}_n}\ar[d]_{\alpha^{r',r}_n}&&
{\cC^{(t')}_n}\ar[ll]_{\alpha^{t',r'}_n}\ar[d]^{\alpha^{t',t}_n}\\
{\cC^{(r)}_n}&&{\cC^{(t)}_n}\ar[ll]^{\alpha^{t,r}_n}}
\ \ \
\xymatrix{
{\uptheta^*_n(\cC^{(t')}_n)}\ar[r]^-(0.5){\upomega^{(t')}_n}\ar[d]_{\uptheta_n^*(\alpha^{t',t}_n)}&{\cC'^{(t')}_n}\ar[d]^{\alpha'^{t',t}_n}\\
{\uptheta^*_n(\cC^{(t)}_n)}\ar[r]_-(0.5){\upomega^{(t)}_n}&{\cC'^{(t)}_n}}
\end{equation}
sont commutatifs. On en déduit un homomorphisme canonique de $\uptheta^*_n(\cC^{(r')}_n)$-algèbres
\begin{equation}\label{ahttfg18g}
\alpha'^{t',t,r',r}_n\colon \cC'^{(t',r')}_n\rightarrow \cC'^{(t,r)}_n.
\end{equation}

Il résulte de \eqref{ahttf37f} que pour tous nombres rationnels  $t,t',t'',r,r',r''$  tels que $t''\geq t'\geq t\geq r\geq 0$, $t''\geq t'\geq r'\geq r$ et $t''\geq r''\geq r'\geq r$, on a 
\begin{equation}\label{ahttfg18i}
\alpha'^{t',t,r',r}_n\circ \alpha'^{t'',t',r'',r'}_n =\alpha'^{t'',t,r'',r}_n.
\end{equation}

\subsection{}\label{ahttfg20}
Reprenons les notations de \ref{hmdf14}.
On désigne par $\bvcF^{(r)}=(\cF^{(r)}_{m+1})_{m\in \mN}$ la $\bvocB$-extension de Higgs-Tate d'épaisseur $r$ 
et par  $\bvcC^{(r)}=(\cC^{(r)}_{m+1})_{m\in \mN}$ la $\bvocB$-algèbre de Higgs-Tate d'épaisseur $r$, associées à $(f,\tX,\cM_\tX)$ (cf. \ref{ahttf14}). 
On considère aussi les objets et notations analogues pour $f'$, que l'on munit d'un $^\prime$. 

D'après (\cite{agt} III.7.5 et III.7.12), les homomorphismes $(\upomega^{(r)}_{m+1})_{m\in \mN}$ \eqref{ahttfg18d} définissent un homomorphisme 
\begin{equation}\label{ahttfg20a}
\bvupomega^{(r)}\colon \bvuptheta^*(\bvcC^{(r)})\rightarrow \bvcC'^{(r)}. 
\end{equation}
Dans la suite, on considère $\bvcC'^{(r)}$ comme une $\bvuptheta^*(\bvcC^{(r)})$-algèbre via $\bvupomega^{(r)}$.

Pour tout nombre rationnel $t\geq r$, on considère la $\bvuptheta^*(\bvcC^{(r)})$-algèbre
\begin{equation}\label{ahttfg20b}
\bvcC'^{(t,r)}=\bvcC'^{(t)}\otimes_{\bvuptheta^*(\bvcC^{(t)})}\bvuptheta^*(\bvcC^{(r)})
\end{equation}
déduite de $\bvcC'^{(t)}$ par changement de base par l'homomorphisme $\bvuptheta^*(\bvalpha^{t,r})$, où $\bvalpha^{t,r}\colon \bvcC^{(t)}\rightarrow \bvcC^{(r)}$
est l'homomorphisme canonique \eqref{ahttf14f}. On a un isomorphisme canonique \eqref{ahttfg18f}
\begin{equation}\label{ahttfg20c}
\bvcC'^{(t,r)}\stackrel{\sim}{\rightarrow}(\cC'^{(t,r)}_{m+1})_{m\in \mN}
\end{equation}     
L'homomorphisme canonique $\bvuptheta^*(\bvcC^{(r)})\rightarrow \bvcC'^{(t,r)}$ induit un homomorphisme
\begin{equation}\label{ahttfg20d}
\bvlgg^*(\bvcC^{(r)})\rightarrow \bvtau_*(\bvcC'^{(t,r)}),
\end{equation}
qui s'identifie à celui induit par les homomorphismes \eqref{ahttfg18h}. 

Pour tous nombres rationnels $t,t',r,r'$ tels que $t'\geq t\geq r\geq 0$ et $t'\geq r'\geq r$, 
on a un homomorphisme canonique de $\bvuptheta^*(\bvcC^{(r')})$-algèbres
\begin{equation}\label{ahttfg20e}
\bvalpha'^{t',t,r',r}\colon \bvcC'^{(t',r')}\rightarrow \bvcC'^{(t,r)},
\end{equation}
qui s'identifie à l'homomorphisme $(\alpha'^{t',t,r',r}_{m+1})_{m\in \mN}$ \eqref{ahttfg18g}.

Il résulte de \eqref{ahttfg18i} que pour tous nombres rationnels  $t,t',t'',r,r',r''$  tels que $t''\geq t'\geq t\geq r\geq 0$, $t''\geq t'\geq r'\geq r$ et $t''\geq r''\geq r'\geq r$, on a 
\begin{equation}\label{ahttfg20j}
\bvalpha'^{t',t,r',r}\circ \bvalpha'^{t'',t',r'',r'}=\bvalpha'^{t'',t,r'',r}.
\end{equation}

On désigne par $\bMod(\bvlgg^*(\bvcC^{(r)}))$ la catégorie des $\bvlgg^*(\bvcC^{(r)})$-modules de $\tG^{\mN^\circ}_s$, par 
$\bMod_\mQ(\bvlgg^*(\bvcC^{(r)}))$ la catégorie des $\bvlgg^*(\bvcC^{(r)})$-modules à isogénie près \eqref{indsh11}
et par $\bIndMod(\bvlgg^*(\bvcC^{(r)}))$ la catégorie des ind-$\bvlgg^*(\bvcC^{(r)})$-modules \eqref{indsh15}. 
On considère les notations analogues pour la catégorie des $\bvuptheta^*(\bvcC^{(r)})$-modules de $\tE'^{\mN^\circ}_s$. 
Le morphisme de topos annelé $\bvtau$ \eqref{hmdf13i} induit un morphisme de topos annelés que nous notons encore abusivement
\begin{equation}\label{ahttfg20f}
\bvtau\colon (\tE'^{\mN^\circ}_s,\bvuptheta^*(\bvcC^{(r)}))\rightarrow (\tG^{\mN^\circ}_s,\bvlgg^*(\bvcC^{(r)}))
\end{equation}
Cette notation n'induit aucun risque d'ambiguïté puisque le foncteur image inverse pour les $\bvlgg^*(\bvcC^{(r)})$-modules coïncide avec 
la restriction du foncteur image inverse pour les $\bvocB^!$-modules \eqref{hmdf14a}. 
On utilise pour $\bvtau$ \eqref{ahttfg20f} les notations introduites dans \ref{indsh21} et \ref{indsh22}. 
On a en particulier deux foncteurs additifs adjoints
\begin{eqnarray}
\rI \bvtau^*\colon \bIndMod(\bvlgg^*(\bvcC^{(r)})) \rightarrow \bIndMod(\bvuptheta^*(\bvcC^{(r)})),\label{ahttfg20g}\\
\rI \bvtau_*\colon \bIndMod(\bvuptheta^*(\bvcC^{(r)})) \rightarrow \bIndMod(\bvlgg^*(\bvcC^{(r)})).\label{ahttfg20h}
\end{eqnarray}
Le foncteur $\rI \bvtau^*$ (resp. $\rI \bvtau_*$) est exact à droite (resp. gauche). 
Le foncteur $\rI \tau_*$ admet un foncteur dérivé à droite
\begin{equation}\label{ahttfg20i}
\rR\rI \bvtau_*\colon \bD^+(\bIndMod(\bvuptheta^*(\bvcC^{(r)})))\rightarrow \bD^+(\bIndMod(\bvlgg^*(\bvcC^{(r)}))).
\end{equation}

\section{Calculs cohomologiques}\label{ccoh}

\subsection{}\label{ccoh0}
Pour les calculs cohomologiques locaux, il est commode d'introduire les notations suivantes. 
Soient $(\oy'\rightsquigarrow \ox')$ un point de $X'_\et\gtimes_{X'_\et}\oX'^\rhd_\et$ (\cite{ag} 3.4.23)
tel que $\ox'$ soit au-dessus de $s$, $\uX'$ le localisé strict de $X'$ en $\ox'$.
D'après (\cite{agt} III.3.7), $\uoX'$ est normal et strictement local (et en particulier intègre); 
on peut donc l'identifier au localisé strict de $\oX'$ en $\oa'(\ox')$ \eqref{hmdf4}. 
Le $X'$-morphisme $\oy'\rightarrow \uX'$
définissant $(\oy'\rightsquigarrow \ox')$ se relève en un $\oX'^\rhd$-morphisme $v'\colon \oy'\rightarrow \uoX'^\rhd$ et 
induit donc un point géométrique de $\uoX'^\rhd$ que l'on note aussi (abusivement) $\oy'$.
On pose $\uDelta'=\pi_1(\uoX'^\rhd,\oy')$. 

Posons $\ox=g(\ox')$ et $\oy=\upgamma(\oy')$ \eqref{hmdf11c} qui sont donc des points géométriques de $X$ et 
$\oX^\circ$ respectivement, et notons $(\oy\rightsquigarrow \ox')$ (resp. $(\oy\rightsquigarrow \ox)$) l'image de $(\oy'\rightsquigarrow \ox')$ 
par le premier (resp. le composé) des morphismes canoniques 
\begin{equation}\label{ccoh0a}
X'_\et\gtimes_{X'_\et}\oX'^\rhd_\et\rightarrow X'_\et\gtimes_{X_\et}\oX^\circ_\et
\rightarrow X_\et\gtimes_{X_\et}\oX^\circ_\et.
\end{equation}
On désigne par $\uX$ le localisé strict de $X$ en $\ox$. 
D'après (\cite{agt} III.3.7), $\uoX$ est normal et strictement local (et en particulier intègre); 
on peut donc l'identifier au localisé strict de $\oX$ en $\oa(\ox)$ \eqref{hmdf4}. 
Le $X$-morphisme $\oy\rightarrow \uX$
définissant $(\oy\rightsquigarrow \ox)$ se relève en un $\oX^\circ$-morphisme $v\colon \oy\rightarrow \uoX^\circ$ et 
induit donc un point géométrique de $\uoX^\circ$ que l'on note aussi (abusivement) $\oy$. 
On pose $\uDelta=\pi_1(\uoX^\circ,\oy)$. On note
\begin{equation}\label{ccoh0b}
\uupgamma\colon \uoX'^\rhd\rightarrow \uoX^\circ
\end{equation} 
le morphisme induit par $g$. Comme $\uupgamma(\oy')=\oy$, celui-ci induit un homomorphisme
\begin{equation}\label{ccoh0c}
\uDelta'\rightarrow \uDelta.
\end{equation}
On désigne par $\uPi$ son noyau. 
On désigne par $\bB_{\uDelta'}$ (resp. $\bB_{\uDelta}$) le topos classifiant du groupe profini $\uDelta'$ (resp. $\uDelta$) et par 
\begin{eqnarray}
\psi_{\uX',\oy'}\colon \uoX'^\rhd_\fet&\stackrel{\sim}{\rightarrow}&\bB_{\uDelta'},\label{ccoh0f1}\\
\psi_{\uX,\oy}\colon \uoX^\circ_\fet&\stackrel{\sim}{\rightarrow}&\bB_{\uDelta},\label{ccoh0f2}
\end{eqnarray}
les foncteurs fibres (\cite{agt}  (VI.9.8.4)). 

On rappelle que la donnée d'un voisinage du point de $X_\et$ associé à $\ox$ 
dans le site $\Et_{/X}$ (resp. $\bP$, resp. $\bQ$ \eqref{ahttfg1})
est équivalente à la donnée d'un $X$-schéma étale $\ox$-pointé (resp. de $\bP$, resp. de $\bQ$) (\cite{sga4} IV 6.8.2). 
Ces objets forment naturellement une catégorie cofiltrante, que l'on note $\fV_\ox$ (resp. $\bP_\ox$,
resp. $\bQ_\ox$). Les catégories $\bP_\ox$ et $\bQ_\ox$ sont $\mU$-petites, et 
les foncteurs d'injection canoniques $\bQ\rightarrow \bP\rightarrow \Et_{/X}$ induisent des foncteurs pleinement 
fidèles et cofinaux $\bQ_\ox\rightarrow \bP_\ox\rightarrow \fV_\ox$.
On définit de même les catégories $\fV'_{\ox'}$, $\bP'_{\ox'}$ et $\bQ'_{\ox'}$. 

Considérons $\fC_{\ox'}$ la catégorie associée à $\ox'$ dans \ref{hmdf27}.
Soit $(\ox'\rightarrow U'\rightarrow U)$ un objet de $\fC_{\ox'}$; on omettra $\ox'$ pour alléger les notations. 
On a un $X'$-morphisme canonique 
$\uX'\rightarrow U'$ et un $X$-morphisme canonique $\uX\rightarrow U$ qui s'insèrent dans un diagramme commutatif  
\begin{equation}\label{ccoh0k}
\xymatrix{
&\uX'\ar[r]^{\ug}\ar[d]&\uX\ar[d]\\
\ox'\ar[r]\ar[ru]&U'\ar[r]&U}
\end{equation}
où $\ug$ est le morphisme induit par $g$. On en déduit des morphismes $\uoX'\rightarrow \oU'$ et $\uoX\rightarrow \oU$. 
Le morphisme $v'\colon \oy'\rightarrow \uoX'^\rhd$ induit un point géométrique de $\oU'^\rhd$ que l'on note aussi $\oy'$. 
De même, le morphisme $v\colon \oy\rightarrow \uoX^\circ$ 
induit un point géométrique de $\oU^\circ$ que l'on note aussi $\oy$. On observera que le diagramme 
\begin{equation}\label{ccoh0l}
\xymatrix{
\oy'\ar[r]\ar[d]&\oy\ar[d]\\
\oU'^\rhd\ar[r]&\oU^\circ}
\end{equation}
est commutatif. 

Les schémas $\oU$ et $\oU'$ étant localement irréductibles d'après (\cite{agt} III.3.3 et III.4.2(iii)),  
ils sont les sommes des schémas induits sur leurs composantes irréductibles. 
On note $\oU^\star$ (resp. $\oU'^\star$)
la composante irréductible de $\oU$ (resp. $\oU'$) contenant $\oy$ (resp. $\oy'$). 
De même, $\oU^\circ$ (resp. $\oU'^\rhd$) est la somme des schémas induits sur ses composantes irréductibles
et $\oU^{\star \circ}=\oU^\star\times_{X}X^\circ$ (resp. $\oU'^{\star \rhd}=\oU'^\star\times_{X'}X'^\rhd$) 
est la composante irréductible de $\oU^\circ$ (resp. $\oU'^\rhd$) contenant $\oy$ (resp. $\oy'$).
On pose $\Delta_U=\pi_1(\oU^{\star\circ},\oy)$ et $\Delta'_{U'}=\pi_1(\oU'^{\star\rhd},\oy')$. 
D'après \eqref{ccoh0l}, on a un homomorphisme canonique $\Delta'_{U'}\rightarrow \Delta_U$. On note $\Pi_{U'\rightarrow U}$ son noyau. 
On désigne par $\bB_{\Delta'_{U'}}$ (resp. $\bB_{\Delta_U}$) le topos classifiant du groupe profini $\Delta'_{U'}$ (resp. $\Delta_U$) et par 
\begin{eqnarray}
\psi_{U',\oy'}\colon \oU'^\rhd_\fet&\stackrel{\sim}{\rightarrow}&\bB_{\Delta'_{U'}},\label{ccoh03f}\\
\psi_{U,\oy}\colon \oU^\circ_\fet&\stackrel{\sim}{\rightarrow}&\bB_{\Delta_U},\label{ccoh0f4}
\end{eqnarray}
les foncteurs fibres. 

Les foncteurs 
\begin{eqnarray}
\fC_{\ox'}&\rightarrow&\fV'_{\ox'},\ \ (\ox'\rightarrow U'\rightarrow U) \mapsto (\ox'\rightarrow U'),\label{ccoh0g1}\\
\fC_{\ox'}&\rightarrow&\fV_{\ox},\ \ \ (\ox'\rightarrow U'\rightarrow U) \mapsto (\ox\rightarrow U),\label{ccoh0g2}
\end{eqnarray}
où le second est défini par composition et le fait que $\ox=g(\ox')$, sont initiaux en vertu de (\cite{sga4} I 8.1.3(b)).
Par suite, les morphismes canoniques
\begin{eqnarray}
\uoX'&\rightarrow& \underset{\underset{(U'\rightarrow U)\in \fC_{\ox'}}{\longleftarrow}}{\lim}\ \oU'^\star,\\
\uoX&\rightarrow& \underset{\underset{(U'\rightarrow U)\in \fC_{\ox'}}{\longleftarrow}}{\lim}\ \oU^\star,
\end{eqnarray}
sont des isomorphismes. D'après (\cite{agt} VI.11.8), les morphismes canoniques 
\begin{eqnarray}
\uDelta'&\rightarrow& \underset{\underset{(U'\rightarrow U)\in \fC_{\ox'}}{\longleftarrow}}{\lim}\ \Delta'_{U'},\\
\uDelta&\rightarrow& \underset{\underset{(U'\rightarrow U)\in \fC_{\ox'}}{\longleftarrow}}{\lim}\ \Delta_U.
\end{eqnarray}
sont donc des isomorphismes. On en déduit que le morphisme canonique 
\begin{equation}\label{ccoh0m}
\uPi\rightarrow \underset{\underset{(U'\rightarrow U)\in \fC_{\ox'}}{\longleftarrow}}{\lim}\ \Pi_{U'\rightarrow U}
\end{equation}
est un isomorphisme.

On désigne par 
\begin{eqnarray}
\varphi'_{\ox'}\colon \tE'&\rightarrow& \uoX'^\rhd_\fet,\label{ccoh0d}\\
\phi_{\ox'}\colon \tG&\rightarrow& \uoX^\circ_\fet,\label{ccoh0e}\\
\varphi_{\ox}\colon \tE&\rightarrow& \uoX^\circ_\fet,\label{ccoh0de} 
\end{eqnarray} 
les foncteurs canoniques définis dans \eqref{ahttfg12ii} et \eqref{hmdf29c}, 
Il résulte de \ref{hmdf30}(i), (\cite{ag} 6.6.7 et (6.5.11.1)),
que pour tout groupe abélien $F$ de $\tE'$ et tout entier $i\geq 0$, on a un isomorphisme canonique
\begin{equation}\label{ccoh0h}
(\rR^i\tau_*(F))_{\varrho(\oy\rightsquigarrow \ox')}\stackrel{\sim}{\rightarrow} \rH^i(\uPi,\psi'_{\uX',\oy'}(\varphi'_{\ox'}(F))),
\end{equation}
où $\varrho$ est le morphisme \eqref{hmdf17f}.

\subsection{}\label{ccoh00}
Conservons les hypothèses et notations de \ref{ccoh0}. 
Pour tout objet $(\ox'\rightarrow U'\rightarrow U)$ de $\fC_{\ox'}$, on pose 
\begin{eqnarray}
\oR'^{\oy'}_{U'}&=&\psi'_{U',\oy'}(\ocB'_{U'}|\oU'^{\star \rhd}),\label{ccoh00a}\\
\oR^{\oy}_{U}&=&\psi_{U,\oy}(\ocB_{U}|\oU^{\star \circ}),\label{ccoh00b}\\
\oR^{!\oy}_{U'\rightarrow U}&=&\psi_{U,\oy}(\ocB^!_{U'\rightarrow U}|\oU^{\star \circ}),\label{ccoh00c}
\end{eqnarray}
où $\ocB_U$ et $\ocB'_{U'}$ (resp. $\ocB^!_{U'\rightarrow U}$) sont les faisceaux définis dans \eqref{TFA2d} (resp. \eqref{hmdf26b}).  

Explicitement, soit $(V_i)_{i\in I}$ le revêtement universel normalisé de $\oU^{\star \circ}$ en $\oy$ (\cite{agt}  VI.9.8).
Pour chaque $i\in I$, $(V_i\rightarrow U)$ (resp. $(U'\rightarrow U\leftarrow V_i)$) est naturellement un objet de $E$ (resp. $G$). On a  alors
\begin{eqnarray}
\oR^{\oy}_{U}&=&\underset{\underset{i\in I}{\longrightarrow}}{\lim}\ \ocB(V_i\rightarrow U),\label{ccoh00d}\\
\oR^{!\oy}_{U'\rightarrow U}&=&\underset{\underset{i\in I}{\longrightarrow}}{\lim}\ \ocB^!(U'\rightarrow U\leftarrow V_i).\label{ccoh00e}
\end{eqnarray} 

L'homomorphisme canonique $\lgg^*(\ocB)\rightarrow \ocB^!$ \eqref{hmdf19c} induit pour tout $i\in I$, un morphisme (fonctoriel en $i$)
\begin{equation}
\ocB(V_i\rightarrow U)\rightarrow \ocB^!(U'\rightarrow U\leftarrow V_i).
\end{equation}
On en déduit par passage à la limite inductive un homomorphisme 
\begin{equation}\label{ccoh00f}
\oR^{\oy}_{U}\rightarrow \oR^{!\oy}_{U'\rightarrow U}.
\end{equation}

Soit $(W_j)_{j\in J}$ le revêtement universel normalisé de $\oU'^{\star \rhd}$ en $\oy'$. 
Pour chaque $j\in J$, $(W_j\rightarrow U')$ est naturellement un objet de $E'$. On a alors 
\begin{equation}\label{ccoh00g}
\oR'^{\oy'}_{U'}=\underset{\underset{j\in J}{\longrightarrow}}{\lim}\  \ocB'(W_j\rightarrow U').
\end{equation}

Pour tout $i\in I$, on a un $\oU^\circ$-morphisme canonique $\oy\rightarrow V_i$.
On en déduit un $\oU'^\rhd$-morphisme $\oy'\rightarrow V_i\times_{\oU^\circ}\oU'^\rhd$.
Le schéma $V_i\times_{\oU^\circ}\oU'^\rhd$ étant localement irréductible,  
il est la somme des schémas induits sur ses composantes irréductibles. 
On note $V'_i$ la composante irréductible de $V_i\times_{\oU^\circ}\oU'^\rhd$ contenant l'image de $\oy'$. 
Les schémas $(V'_i)_{i\in I}$ forment naturellement un système projectif de revêtements étales finis connexes $\oy'$-pointés de $\oU'^{\star \rhd}$.
Le morphisme canonique 
\begin{equation}
\oU'^\rhd\times_{\oU^\circ}V_i\rightarrow U'\times_{(U\times_XX')}(V_i\times_{X^\circ}X'^\rhd)
\end{equation}
est un isomorphisme. On a donc un isomorphisme canonique de $\tE'$ 
\begin{equation}
\tau^*((U'\rightarrow U\leftarrow V_i)^a)\stackrel{\sim}{\rightarrow} (V_i\times_{\oU^\circ}\oU'^\rhd\rightarrow U')^a.
\end{equation}
Comme l'homomorphisme canonique $\ocB^!\rightarrow\tau_*(\ocB')$ \eqref{hmdf19b} est un isomorphisme, 
on en déduit un homomorphisme canonique (fonctoriel en $i$)
\begin{equation}
\ocB^!(U'\rightarrow U\leftarrow V_i)=\ocB'(V_i\times_{\oU^\circ}\oU'^\rhd\rightarrow U')\rightarrow \ocB'(V'_i\rightarrow U').
\end{equation}
Posant
\begin{equation}\label{ccoh00h}
\oR^{\intern\oy'}_{U'\rightarrow U}=\underset{\underset{i\in I}{\longrightarrow}}{\lim}\  \ocB'(V'_i\rightarrow U'),
\end{equation}
on a donc un homomorphisme canonique
\begin{equation}
\oR^{!\oy}_{U'\rightarrow U}\rightarrow \oR^{\intern\oy'}_{U'\rightarrow U}.
\end{equation}

Pour tous $i\in I$ et $j\in J$, il existe au plus un morphisme de $\oU'^{\star \rhd}$-schémas pointés $W_j\rightarrow V'_i$. 
De plus, pour tout $i\in I$, il existe $j\in J$ et un morphisme de $\oU'^{\star \rhd}$-schémas pointés $W_j\rightarrow V'_i$. 
On a donc un homomorphisme
\begin{equation}\label{ccoh00i}
\oR^{\intern\oy'}_{U'\rightarrow U}=\underset{\underset{i\in I}{\longrightarrow}}{\lim}\  \ocB'(V'_i\rightarrow U')\rightarrow 
\underset{\underset{j\in J}{\longrightarrow}}{\lim}\  \ocB'(W_j\rightarrow U')=\oR'^{\oy'}_{U'}. 
\end{equation}

On dispose donc de trois homomorphismes canoniques $\Delta'_{U'}$-équivariants
\begin{equation}\label{ccoh00j}
\oR^\oy_U\rightarrow \oR^{!\oy}_{U'\rightarrow U} \rightarrow \oR^{\intern\oy'}_{U'\rightarrow U}\rightarrow \oR'^{\oy'}_{U'}.
\end{equation}
Ces anneaux et ces homomorphismes étant fonctoriels en $(\ox'\rightarrow U'\rightarrow U)\in \ob(\fC_{\ox'})$, posons
\begin{eqnarray}
\oR'^{\oy'}_{\uX'}&=&\underset{\underset{(\ox'\rightarrow U'\rightarrow U)\in \ob(\fC^\circ_{\ox'})}{\longrightarrow}}{\lim}\ \oR'^{\oy'}_{U'},\label{ccoh00k1}\\
\oR^{\intern\oy}_{\uX'\rightarrow \uX}&=&\underset{\underset{(\ox'\rightarrow U'\rightarrow U)\in \ob(\fC^\circ_{\ox'})}{\longrightarrow}}{\lim}
\oR^{\intern\oy'}_{U'\rightarrow U},\label{ccoh00k2}\\
\oR^{!\oy}_{\uX'\rightarrow \uX}&=&\underset{\underset{(\ox'\rightarrow U'\rightarrow U)\in \ob(\fC^\circ_{\ox'})}{\longrightarrow}}{\lim}
\oR^{!\oy}_{U'\rightarrow U},\label{ccoh00k3}\\
\oR^{\oy}_{\uX}&=&\underset{\underset{(\ox'\rightarrow U'\rightarrow U)\in \ob(\fC^\circ_{\ox'})}{\longrightarrow}}{\lim}\ \oR^\oy_U,\label{ccoh00k4}
\end{eqnarray}
dont les deux premiers sont des anneaux de $\bB_{\uDelta'}$ et les deux autres sont des anneaux de $\bB_{\uDelta}$. 

D'après \eqref{TFA12g}, comme les foncteurs \eqref{ccoh0g1} et \eqref{ccoh0g2} sont initiaux, on a des isomorphismes canoniques
\begin{eqnarray}
\psi_{\uX',\oy'}(\varphi'_{\ox'}(\ocB'))&\stackrel{\sim}{\rightarrow} &\oR'^{\oy'}_{\uX'},\label{ccoh00l3}\\
\psi_{\uX,\oy}(\varphi_{\ox}(\ocB))&\stackrel{\sim}{\rightarrow} &\oR^{\oy}_{\uX}.\label{ccoh00l4}
\end{eqnarray}
Ceux-ci induisent des isomorphismes canoniques \eqref{TFA11c}
\begin{eqnarray}
\ocB'_{\rho'(\oy'\rightsquigarrow \ox')}&\stackrel{\sim}{\rightarrow} &\oR'^{\oy'}_{\uX'},\label{ccoh00l1}\\
\ocB_{\rho(\oy\rightsquigarrow \ox)}&\stackrel{\sim}{\rightarrow} &\oR^{\oy}_{\uX}.\label{ccoh00l2}
\end{eqnarray}
Par ailleurs, d'après \eqref{hmdf28e}, on a un isomorphisme canonique
\begin{equation}\label{ccoh00m}
\ocB^!_{\varrho(\oy\rightsquigarrow \ox')}\stackrel{\sim}{\rightarrow} \oR^{!\oy}_{\uX'\rightarrow \uX}. 
\end{equation}

Les homomorphismes \eqref{ccoh00j} induisent par passage à la limite inductive des homomorphismes canoniques $\uDelta'$-équivariants
\begin{equation}\label{ccoh00n}
\oR^\oy_\uX\rightarrow \oR^{!\oy}_{\uX'\rightarrow \uX} \rightarrow \oR^{\intern\oy'}_{\uX'\rightarrow \uX}\rightarrow \oR'^{\oy'}_{\uX'}.
\end{equation}
En vertu de (\cite{ag} 6.5.26), l'homomorphisme central est un isomorphisme. 
Compte tenu de ce qui précède, le premier homomorphisme et le composé des deux autres s'identifient aux homomorphismes 
\begin{equation}\label{ccoh00o}
\ocB_{\rho(\oy\rightsquigarrow \ox)}\rightarrow \ocB^!_{\varrho(\oy\rightsquigarrow \ox')}\rightarrow \ocB'_{\rho'(\oy'\rightsquigarrow \ox')}
\end{equation}
induits par les adjoints des homomorphismes \eqref{hmdf19c} et \eqref{hmdf19b}, respectivement.

\begin{rema}\label{ccoh01}
Sous les hypothèses de \ref{ccoh00},  pour tout objet $(\ox'\rightarrow U'\rightarrow U)$ de $\fC_{\ox'}$ tel que la restriction 
$(U',\cM_{X'}|U')\rightarrow (U,\cM_X|U)$ de $g$ vérifie les hypothèses de \ref{hypmdef2}, les anneaux 
$\oR^\oy_U\rightarrow \oR^{\intern\oy'}_{U'\rightarrow U}\rightarrow \oR'^{\oy'}_{U'}$ coïncident aves les anneaux 
$\oR\rightarrow \oR^{\intern}\rightarrow \oR'$ définis dans \ref{eccr34}. 
\end{rema}

\subsection{}\label{ccoh02}
Conservons les hypothèses et notations de \ref{ccoh0} et \ref{ccoh00}. 
On note $W$ l'anneau des vecteurs de Witt à coefficients dans $k$ relatif à $p$,
$K_0$ le corps des fractions de $W$ et $\fd$ la différente de l'extension $K/K_0$.
On fixe une uniformisante $\varpi$ de $\co_K$. 
On pose $\delta=0$ dans le cas absolu \eqref{definf10} et $\delta=v(\varpi\fd)$ dans le cas relatif. 
On prendra garde au changement de notation par rapport à \eqref{mph1g} pour éviter la confusion avec le morphisme $\rho$ \eqref{hmdf7g}.  
D'après \eqref{definf16a} et \eqref{definf17c}, on a un isomorphisme $\co_C$-linéaire canonique
\begin{equation}\label{ccoh02d}
\co_C(1)\stackrel{\sim}{\rightarrow} p^{\delta+\frac{1}{p-1}}\txi\co_C.
\end{equation}

On désigne par  $\fC'_{\ox'}$ la sous-catégorie pleine de la catégorie $\fC_{\ox'}$ \eqref{hmdf27}
formée des objets $(\ox'\rightarrow U'\rightarrow U)$ tels que les conditions suivantes soient remplies:
\begin{itemize}
\item[(i)] Les schémas $U$ et $U'$ sont affines et connexes, 
et la restriction $(U',\cM_{X'}|U')\rightarrow (U,\cM_X|U)$ du morphisme $g$ \eqref{hmdf5a} admet une carte relativement adéquate (\cite{ag} 5.1.11).
Ces conditions correspondent à celles (pour $g$) fixées dans \ref{hypmdef2}.
\item[(ii)] Il existe une carte fine et saturée $M\rightarrow \Gamma(U',\cM_{X'})$ pour $(U',\cM_{X'}|U')$
induisant un isomorphisme 
\begin{equation}\label{ccoh02a}
M\stackrel{\sim}{\rightarrow} \Gamma(U',\cM_{X'})/\Gamma(U',\co^\times_{X'}).
\end{equation}
\end{itemize}
La catégorie $\fC'_{\ox'}$ est cofiltrante, et le foncteur d'injection canonique
$\fC'_{\ox'}\rightarrow \fC_{\ox'}$ est initial d'après (\cite{agt} II.5.17) et (\cite{sga4} I 8.1.3(c)).
De même, le foncteur $\fC'_{\ox'}\rightarrow \bQ'_{\ox'}$ défini par $(\ox'\rightarrow U'\rightarrow U)\mapsto (\ox'\rightarrow U')$ est initial \eqref{ccoh0}. 

Pour tout objet $(\ox'\rightarrow U'\rightarrow U)$ de $\fC'_{\ox'}$ et tout entier $n\geq 0$, 
on a une suite exacte canonique de $\oR^{\oy'}_{U'}$-représentations de $\Delta'_{U'}$ 
\begin{equation}\label{ccoh02b}
0\rightarrow \oR'^{\oy'}_{U'}/p^n \oR'^{\oy'}_{U'}\rightarrow \cF'^{\oy'}_{U'}/p^n \cF'^{\oy'}_{U'} 
\rightarrow \txi^{-1}\tOmega^1_{\oX'_n/\oS_n}(\oU'^\star) \otimes_{\co_{\oX'_n}(\oU'^\star)}\oR'^{\oy'}_{U'}\rightarrow 0,
\end{equation} 
déduite de la suite exacte \eqref{ahttf6f} (pour $f'$ et $U'$).  On désigne par 
\begin{equation}\label{ccoh02c}
\partial_{U'\rightarrow U}\colon \txi^{-1}\tOmega^1_{\oX'_n/\oS_n}(\oU'^\star)\otimes_{\co_{\oX'}(\oU'^\star)}\oR^{\intern\oy'}_{U'\rightarrow U}\rightarrow 
\rH^1(\Pi_{U'\rightarrow U},\oR'^{\oy'}_{U'}/p^n\oR'^{\oy'}_{U'})
\end{equation} 
le bord de la suite exacte longue de cohomologie déduite de la suite exacte courte \eqref{ccoh02b}, 
où $\Pi_{U'\rightarrow U}$ est le noyau de l'homomorphisme canonique $\Delta'_{U'}\rightarrow \Delta_U$. 
On a alors un diagramme commutatif
\begin{equation}\label{ccoh02e}
\xymatrix{
{\txi^{-1}\tOmega^1_{\oX'_n/\oS_n}(\oU'^\star)\otimes_{\co_{\oX'}(\oU'^\star)}\oR^{\intern\oy'}_{U'\rightarrow U}}
\ar[dd]_u\ar[r]^-(0.5){\partial_{U'\rightarrow U}}&
{\rH^1(\Pi_{U'\rightarrow U},\oR'^{\oy'}_{U'}/p^n\oR'^{\oy'}_{U'})}\\
&{\rH^1(\Pi_{U'\rightarrow U},p^{\delta+\frac{1}{p-1}}\oR'^{\oy'}_{U'}/p^{\delta+\frac{1}{p-1}+n}\oR'^{\oy'}_{U'})}\ar[u]_c\\
{\txi^{-1}\tOmega^1_{\oX'_n/\oX_n}(\oU'^\star)\otimes_{\co_{\oX'}(\oU'^\star)}\oR^{\intern\oy'}_{U'\rightarrow U}}\ar[r]^-(0.5){a}&
{\rH^1(\Pi_{U'\rightarrow U},\txi^{-1}(\oR'^{\oy'}_{U'}/p^{n}\oR'^{\oy'}_{U'})(1))}\ar[u]_{-b}}
\end{equation}
où $u$ est le morphisme canonique, $a$ est la composante en degré un du morphisme \eqref{eccr51a}, 
$b$ est induit par l'isomorphisme \eqref{ccoh02d} et $c$ est induit par l'injection canonique $p^{\delta+\frac{1}{p-1}}\co_C \rightarrow \co_C$. 
En effet, on peut supposer $U=X$ et $U'=X'$. Le morphisme $\partial_{U'\rightarrow U}$ ne dépend pas de la déformation $(\tX',\cM_{\tX'})$ choisie \eqref{taht8}. 
On peut donc se borner au cas où $(\tX',\cM_{\tX'})$
est la déformation définie par la carte adéquate de $f'$ fixée dans (i) (cf. \ref{pmh8}). Posant $\oR'=\oR'^{\oy'}_{U'}$, 
la même carte détermine une section du torseur de Higgs-Tate sur $\Spec(\hoRp)$ \eqref{mtht6}, 
et par suite un scindage de l'extension \eqref{mtht6a} (cf. \cite{agt} II.10.11). 
L'assertion résulte alors de \ref{pmh90}; cf. (\cite{ag} 6.6.12) qui traite le cas absolu \eqref{definf10}.

\begin{lem}\label{ccoh03}
Conservons les hypothèses et notations de \ref{ccoh02}. Soient, de plus, $(\ox'\rightarrow U'\rightarrow U)$ un objet de $\fC'_{\ox'}$ \eqref{ccoh02},
$n$ un entier $\geq 0$. Alors, 
\begin{itemize}
\item[{\rm (i)}] Le morphisme $\partial_{U'\rightarrow U}$ se factorise à travers $u$ \eqref{ccoh02e} et induit un morphisme
\begin{equation}\label{ccoh03a}
\partial'_{U'\rightarrow U}\colon \txi^{-1}\tOmega^1_{\oX'_n/\oX_n}(\oU'^\star)\otimes_{\co_{\oX'}(\oU'^\star)}\oR^{\intern\oy'}_{U'\rightarrow U}\rightarrow 
\rH^1(\Pi_{U'\rightarrow U},\oR'^{\oy'}_{U'}/p^n\oR'^{\oy'}_{U'}). 
\end{equation} 
\item[{\rm (ii)}] Il existe un et un unique homomorphisme de $\oR^{\intern\oy'}_{U'\rightarrow U}$-algèbres graduées 
\begin{equation}\label{ccoh03b}
\wedge \left(\txi^{-1}\tOmega^1_{\oX'_n/\oX_n}(\oU'^\star)\otimes_{\co_{\oX'}(\oU'^\star)}\oR^{\intern\oy'}_{U'\rightarrow U}\right)
\rightarrow \oplus_{i\geq 0}\rH^i(\Pi_{U'\rightarrow U},\oR'^{\oy'}_{U'}/p^n\oR'^{\oy'}_{U'})
\end{equation}
dont la composante en degré un est le morphisme $\partial'_{U'\rightarrow U}$ \eqref{ccoh03a}. De plus, 
son noyau est annulé par $p^{2\ell(\delta+\frac{1}{p-1})}\fm_\oK$ et son conoyau est annulé par $p^{2\ell\delta+\frac{2\ell+1}{p-1}}\fm_\oK$, où $\ell=\dim(X'/X)$.
\item[{\rm (iii)}] Pour tout $i\geq \ell+1$, $\rH^i(\Pi_{U'\rightarrow U},\oR'^{\oy'}_{U'}/p^n\oR'^{\oy'}_{U'})$ est $\alpha$-nul.
\end{itemize}
\end{lem}

(i) Cela résulte aussitôt du diagramme \eqref{ccoh02e}.  

(ii) En vertu de \ref{eccr51}(i), il existe un et un unique homomorphisme de $\oR^{\intern\oy'}_{U'\rightarrow U}$-algèbres graduées
\begin{equation}\label{ccoh03c}
\wedge (\xi^{-1}\tOmega^1_{\oX'_n/\oX_n}(\oU'^\star)\otimes_{\co_{\oX'}(\oU'^\star)}\oR^{\intern\oy'}_{U'\rightarrow U})
\rightarrow \oplus_{i\geq 0}\rH^i(\Pi_{U'\rightarrow U},\xi^{-i}(\oR'^{\oy'}_{U'}/p^n\oR'^{\oy'}_{U'})(i))
\end{equation}
dont la composante en degré un est le morphisme $a$ du diagramme \eqref{ccoh02e}. Son noyau est $\alpha$-nul et son conoyau est annulé par 
$p^{\frac{1}{p-1}}\fm_\oK$. On en déduit qu'il existe un et un unique homomorphisme de 
$\oR^{\intern\oy'}_{U'\rightarrow U}$-algèbres graduées
\begin{equation}\label{ccoh03d}
\wedge \left(\txi^{-1}\tOmega^1_{\oX'_n/\oX_n}(\oU'^\star)\otimes_{\co_{\oX'}(\oU'^\star)}\oR^{\intern\oy'}_{U'\rightarrow U}\right)
\rightarrow \oplus_{i\geq 0}\rH^i(\Pi_{U'\rightarrow U},\oR'^{\oy'}_{U'}/p^n\oR'^{\oy'}_{U'})
\end{equation}
dont la composante en degré un est le morphisme $\partial'_{U'\rightarrow U}$. 
Une chasse au diagramme \eqref{ccoh02e} 
montre que le noyau de \eqref{ccoh03d} est annulé par $p^{2\ell(\delta+\frac{1}{p-1})}\fm_\oK$.
Compte tenu de (iii) ci-dessous, le conoyau de \eqref{ccoh03d} est annulé par $p^{2\ell\delta+\frac{2\ell+1}{p-1}}\fm_\oK$.

(iii) C'est l'énoncé \ref{eccr51}(ii) mentionné pour rappel.

\subsection{}\label{ccoh1}
Soit $n$ un entier $\geq 0$. 
Reprenons les notations de \ref{hmdf6} et \eqref{hmdf19l} et considérons la suite exacte de $\ocB'_n$-modules 
\begin{equation}\label{ccoh1a}
0\rightarrow \ocB'_n\rightarrow \cF'_n\rightarrow \sigma'^*_n(\txi^{-1}\tOmega^1_{\oX'_n/\oS_n})\rightarrow 0,
\end{equation}
analogue de la suite \eqref{ahttfg3e} pour $f'$ avec $r=0$. Considérons le morphisme $\ocB^!_n$-linéaire de $\tG_s$ 
\begin{equation}\label{ccoh1b}
\pi_n^*(\txi^{-1}\tOmega^1_{\oX'_n/\oS_n})\rightarrow \rR^1\tau_{n*}(\ocB'_n),
\end{equation}
composé du morphisme d'adjonction 
\begin{equation}\label{ccoh1c}
\pi_n^*(\txi^{-1}\tOmega^1_{\oX'_n/\oS_n})\rightarrow \tau_{n*}(\tau^*_n(\pi_n^*(\txi^{-1}\tOmega^1_{\oX'_n/\oS_n})))
\end{equation}
et du cobord de la suite exacte \eqref{ccoh1a} en tenant compte de l'isomorphisme $\sigma'_n\stackrel{\sim}{\rightarrow} \pi_n \tau_n$ \eqref{hmdf19l}.

\begin{prop}\label{ccoh2}
Conservons les hypothèses et notations de \ref{ccoh1}.  Alors, 
\begin{itemize}
\item[{\rm (i)}] Le morphisme \eqref{ccoh1b} se factorise à travers le morphisme canonique surjectif 
\begin{equation}\label{ccoh2a}
\pi_n^*(\txi^{-1}\tOmega^1_{\oX'_n/\oS_n})\rightarrow \pi_n^*(\txi^{-1}\tOmega^1_{\oX'_n/\oX_n}),
\end{equation}
et il induit un morphisme $\ocB^!_n$-linéaire de $\tG_s$ 
\begin{equation}\label{ccoh2b}
\pi_n^*(\txi^{-1}\tOmega^1_{\oX'_n/\oX_n})\rightarrow \rR^1\tau_{n*}(\ocB'_n).
\end{equation}
\item[{\rm (ii)}] Il existe un et un unique homomorphisme de $\ocB^!_n$-algèbres graduées de $\tG_s$
\begin{equation}\label{ccoh2c}
\wedge (\pi_n^*(\txi^{-1}\tOmega^1_{\oX'_n/\oX_n}))\rightarrow \oplus_{i\geq 0}\rR^i\tau_{n*}(\ocB'_n)
\end{equation}
dont la composante en degré un est le morphisme \eqref{ccoh2b}. De plus, 
son noyau est annulé par $p^{2\ell(\delta+\frac{1}{p-1})}\fm_\oK$ et son conoyau est annulé par $p^{2\ell\delta+\frac{2\ell+1}{p-1}}\fm_\oK$, 
où $\delta$ est défini dans \ref{ccoh02} et $\ell=\dim(X'/X)$.  
\item[{\rm (iii)}] Pour tout entier $i\geq \ell+1$, $\rR^i\tau_{n*}(\ocB'_n)$ est $\alpha$-nul. 
\end{itemize}
\end{prop}

Soient $(\oy\rightsquigarrow \ox')$ un point de $X'_\et\gtimes_{X_\et}\oX^\circ_\et$ tel que $\ox'$ soit au-dessus de $s$ (\cite{ag} 3.4.23), 
$\uX'$ le localisé strict de $X'$ en $\ox'$, $\ox=g(\ox')$, $\uX$ le localisé strict de $X$ en $\ox$.
Le morphisme $\uupgamma\colon \uoX'^\rhd\rightarrow \uoX^\circ$ induit par $g$ est fidèlement plat en vertu de (\cite{ag} 2.4.1). 
Il existe donc un point $(\oy'\rightsquigarrow \ox')$ de $X'_\et\gtimes_{X'_\et}\oX'^\rhd_\et$ 
qui relève le point $(\oy\rightsquigarrow \ox')$ de $X'_\et\gtimes_{X_\et}\oX^\circ_\et$ \eqref{hmdf17g}.
Reprenons alors les notations introduites dans \ref{ccoh0}, \ref{ccoh00} et \ref{ccoh02}. 

D'après \eqref{TFA14g} et \ref{ahttf48}, 
l'image de la suite exacte \eqref{ccoh1a} par le foncteur composé $\psi_{\uX',\oy'}\circ \varphi'_{\ox'}$ \eqref{ccoh0}
s'identifie à la limite inductive de la suite \eqref{ccoh02b}, lorsque $(\ox'\rightarrow U'\rightarrow U)$ décrit la catégorie $\fC'^\circ_{\ox'}$ \eqref{ccoh02}.

En vertu de \eqref{hmdf18k} et \eqref{ccoh0h}, pour tout entier $i\geq 0$, 
on a un isomorphisme canonique
\begin{equation}\label{ccoh2d}
\rR^i\tau_{n*}(\ocB'_n)_{\varrho(\oy\rightsquigarrow \ox')}\stackrel{\sim}{\rightarrow}\rH^i(\uPi,\psi_{\uX',\oy'}(\varphi'_{\ox'}(\ocB'_n))).
\end{equation}
Par ailleurs, d'après \eqref{TFA12g}, on a un isomorphisme canonique de $\bB_{\uDelta'}$
\begin{equation}\label{ccoh2e}
\psi_{\uX',\oy'}(\varphi'_{\ox'}(\ocB'))\stackrel{\sim}{\rightarrow} \oR'^{\oy'}_{\uX'}.
\end{equation}
D'après \eqref{ccoh0m} et (\cite{serre1} I prop.~8), le morphisme canonique
\begin{equation}
\underset{\underset{(U'\rightarrow U)\in \fC'^\circ_{\ox'}}{\longrightarrow}}{\lim}\
\rH^i(\Pi_{U'\rightarrow U},\oR'^{\oy'}_{U'}/p^n\oR'^{\oy'}_{U'}) \rightarrow \rH^i(\uPi,\oR'^{\oy'}_{\uX'}/p^n\oR'^{\oy'}_{\uX'})
\end{equation} 
est un isomorphisme.

En vertu de \eqref{ccoh00m} et (\cite{ag} (3.4.23.1) et 6.5.26), on a un isomorphisme canonique 
\begin{equation}\label{ccoh2f}
\pi^*_n(\txi^{-1}\tOmega^1_{\oX'_n/\oS_n})_{\varrho(\oy\rightsquigarrow \ox')}\stackrel{\sim}{\rightarrow}
\txi^{-1}\tOmega^1_{\oX'_n/\oS_n,\ox'}\otimes_{\co_{\oX',\ox'}}\oR^{\intern\oy'}_{\uX'\rightarrow \uX},
\end{equation}
et de même si l'on remplace $\tOmega^1_{\oX'_n/\oS_n}$ par $\tOmega^1_{\oX'_n/\oX_n}$.

D'après \ref{hmdf30}(ii), la fibre du morphisme \eqref{ccoh1b} en $\varrho(\oy\rightsquigarrow \ox')$ 
s'identifie donc à la limite inductive, lorsque $(\ox'\rightarrow U'\rightarrow U)$ décrit la catégorie $\fC'^\circ_{\ox'}$, du morphisme \eqref{ccoh02c}
\begin{equation}\label{ccoh2i}
\partial_{U'\rightarrow U}\colon \txi^{-1}\tOmega^1_{\oX'_n/\oS_n}(\oU'^\star)\otimes_{\co_{\oX'}(\oU'^\star)}\oR^{\intern\oy'}_{U'\rightarrow U} \rightarrow 
\rH^1(\Pi_{U'\rightarrow U},\oR'^{\oy'}_{U'}/p^n\oR'^{\oy'}_{U'}).
\end{equation} 
La proposition résulte alors de \ref{ccoh03} compte tenu de (\cite{ag} 6.5.3).

\subsection{}\label{ccoh16}
Reprenons les hypothèses et notations de \ref{ccoh0} et \ref{ccoh00}. 
Supposons, de plus, que le schéma $X$ est séparé et que le morphisme $f\colon (X,\cM_X)\rightarrow (S,\cM_S)$ 
admet une carte adéquate $((P,\gamma),(\mN,\iota),\vartheta)$, et reprenons les notations de \ref{ahttfg5}. 

On désigne par $\mJ_{\ox'}$ la sous-catégorie pleine de la catégorie $\fC_{\ox'}$ définie dans \ref{hmdf27}
formée des objets $(\iota' \colon \ox'\rightarrow U',\upmu\colon U'\rightarrow U)$ tels que $\upmu$ soit un objet de la catégorie $\mJ$ définie dans \ref{ahttfg103}, 
autrement dit, tels que $U'$ soit un objet de $\bQ'$ et que $U$ soit un objet de $\bP$ \eqref{ahttfg1}. 
Le foncteur d'injection canonique $\mJ_{\ox'}\rightarrow \fC_{\ox'}$ est 
initial et la catégorie $\mJ_{\ox'}$ est cofiltrante en vertu de (\cite{sga4} I 8.1.3(c)) et (\cite{agt} II.5.17). 
Par ailleurs, les foncteurs 
\begin{eqnarray}
\fC_{\ox'}\rightarrow \bQ'_{\ox'},&& (\ox'\rightarrow  U'\rightarrow U)\mapsto (\ox'\rightarrow U'),\label{ccoh16a}\\
\fC_{\ox'}\rightarrow \bP_\ox,&& (\ox'\rightarrow  U'\rightarrow U)\mapsto (\ox\rightarrow U),\label{ccoh16b}
\end{eqnarray}
où les catégories cofiltrantes $\bP_\ox$ et $\bQ'_{\ox'}$ sont définies dans \ref{ccoh0} 
et le second foncteur est défini par composition et le fait que $\ox=g(\ox')$, sont initiaux en vertu de (\cite{sga4} I 8.1.3(b)) et (\cite{agt} II.5.17).

Pour tous nombre rationnel $r\geq 0$ et objet $(\iota' \colon \ox'\rightarrow U'\rightarrow U)$ de $\mJ_{\ox'}$, 
on munit le morphisme $f|U\colon (U,\cM_X|U)\rightarrow (S,\cM_S)$ induit par $f$ de la carte adéquate induite par 
$((P,\gamma),(\mN,\iota),\vartheta)$, et on reprend les notations introduites dans \ref{ahttfg4}.  
On a alors un morphisme canonique $\pi_1(\oU'^{\star\rhd},\oy')$-équivariant de $\hoR^{\oy}_{U}$-algèbres \eqref{ahttfg4l}
\begin{equation}\label{ccoh16c}
\cC^{\oy,(t)}_{U}\rightarrow \cC'^{\oy',(t)}_{U'}.
\end{equation}
On considère la $\cC^{\oy,(r)}_U$-algèbre \eqref{mtht11e}
\begin{equation}\label{ccoh16d}
\cC'^{\oy',(t,r)}_{U'\rightarrow U}=\cC'^{\oy',(t)}_{U'}\otimes_{\cC^{\oy,(t)}_U}\cC^{\oy,(r)}_U
\end{equation}
déduite de $\cC'^{\oy',(t)}_{U'}$ par changement de base par l'homomorphisme canonique $\cC^{\oy,(t)}_U\rightarrow \cC^{\oy,(r)}_U$ \eqref{mtht10}. 

En vertu de \ref{ahttf48}, \eqref{ahttfg12j} et \eqref{ccoh00l3}, l'image de l'homomorphisme 
$\upomega^{(t)}_n\colon \uptheta^*_n(\cC^{(t)}_n)\rightarrow \cC'^{(t)}_n$ \eqref{ahttfg18d}  
par le foncteur composé $\psi_{\uX',\oy'}\circ \varphi'_{\ox'}$ \eqref{ccoh0}
s'identifie à la limite inductive des homomorphismes
\begin{equation}\label{ccoh16e}
(\cC^{\oy,(r)}_{U}/p^n\cC^{\oy,(t)}_{U})\otimes_{\oR^{\oy}_{U}} \oR'^{\oy'}_{U'} \rightarrow \cC'^{\oy',(r)}_{U'}/p^n\cC'^{\oy',(t)}_{U'}
\end{equation}
induits par \eqref{ccoh16c}, lorsque $(\ox'\rightarrow U'\rightarrow U)$ décrit la catégorie $\mJ^\circ_{\ox'}$ (cf. la preuve de \ref{ahttfg15}).
On en déduit un isomorphisme canonique 
\begin{equation}\label{ccoh16f}
\psi_{\uX',\oy'}(\varphi'_{\ox'}(\cC'^{(t,r)}_n)) \stackrel{\sim}{\rightarrow} \underset{\underset{(\ox'\rightarrow U'\rightarrow U)\in \mJ^\circ_{\ox'}}{\longrightarrow}}{\lim}\
\cC'^{\oy',(t,r)}_{U'\rightarrow U}/p^n\cC'^{\oy',(t,r)}_{U'\rightarrow U}. 
\end{equation}

\begin{prop}\label{ccoh3}
Soient $r,t,t'$ trois nombres rationnels tels que $t'>t>r\geq 0$. Alors,
\begin{itemize}
\item[{\rm (i)}] Pour tout entier $n\geq 1$, l'homomorphisme canonique \eqref{ahttfg18h}
\begin{equation}\label{ccoh3a}
\lgg_n^*(\cC^{(r)}_n)\rightarrow \tau_{n*}(\cC'^{(t,r)}_n)
\end{equation}
est $\alpha$-injectif. Notons $\cH^{(t,r)}_n$ son conoyau. 
\item[{\rm (ii)}] Il existe un nombre rationnel $a>0$ tel que pour tout entier $n\geq 1$, 
le morphisme 
\begin{equation}\label{ccoh3b}
\cH^{(t',r)}_n\rightarrow \cH^{(t,r)}_n
\end{equation} 
induit par l'homomorphisme canonique 
$\cC'^{(t',r)}_n\rightarrow \cC'^{(t,r)}_n$ \eqref{ahttfg18g} soit annulé par $p^a$. 
\item[{\rm (iii)}] Il existe un nombre rationnel $b>0$ tel que pour tous entiers $n,q\geq 1$, 
le morphisme canonique
\begin{equation}\label{ccoh3c}
\rR^q\tau_{n*}(\cC'^{(t',r)}_n)\rightarrow \rR^q\tau_{n*}(\cC'^{(t,r)}_n)
\end{equation}
soit annulé par $p^b$.
\end{itemize}
\end{prop}

La question étant locale sur $X$ (\cite{agt} VI.10.14 et \cite{ag} 3.4.12), 
on peut supposer que le schéma $X$ est séparé et que le morphisme $f\colon (X,\cM_X)\rightarrow (S,\cM_S)$ 
admet une carte adéquate $((P,\gamma),(\mN,\iota),\vartheta)$. 
 
Soient $(\oy\rightsquigarrow \ox')$ un point de $X'_\et\gtimes_{X_\et}\oX^\circ_\et$ tel que $\ox'$ soit au-dessus de $s$, 
$\uX'$ le localisé strict de $X'$ en $\ox'$, $\ox=g(\ox')$, $\uX$ le localisé strict de $X$ en $\ox$.
Le morphisme $\uupgamma\colon \uoX'^\rhd\rightarrow \uoX^\circ$ induit par $g$ étant fidèlement plat en vertu de (\cite{ag} 2.4.1), 
il existe un point $(\oy'\rightsquigarrow \ox')$ de $X'_\et\gtimes_{X'_\et}\oX'^\rhd_\et$ 
qui relève le point $(\oy\rightsquigarrow \ox')$ de $X'_\et\gtimes_{X_\et}\oX^\circ_\et$ \eqref{ccoh0a}.
Reprenons alors les notations introduites dans \ref{ccoh0}, \ref{ccoh00} et \ref{ccoh16}. 

En vertu de \eqref{hmdf18k} et \eqref{ccoh0h}, pour tout entier $i\geq 0$, 
on a un isomorphisme canonique
\begin{equation}
\rR^i\tau_{n*}(\cC'^{(t,r)}_n)_{\varrho(\oy\rightsquigarrow \ox')}\stackrel{\sim}{\rightarrow}\rH^i(\uPi,\psi_{\uX',\oy'}(\varphi'_{\ox'}(\cC'^{(t,r)}_n))).
\end{equation}
D'après \eqref{ccoh0m}, \eqref{ccoh16f} et (\cite{serre1} I prop.~8), le morphisme canonique
\begin{equation}
\underset{\underset{(\ox'\rightarrow U'\rightarrow U)\in \mJ^\circ_{\ox'}}{\longrightarrow}}{\lim}\
\rH^i(\Pi_{U'\rightarrow U}, \cC'^{\oy',(t,r)}_{U'\rightarrow U}/p^n\cC'^{\oy',(t,r)}_{U'\rightarrow U}) \rightarrow \rH^i(\uPi,\psi_{\uX',\oy'}(\varphi'_{\ox'}(\cC'^{(t,r)}_n)))
\end{equation} 
est un isomorphisme. Compte tenu de \ref{ahttf48}, \eqref{ccoh00m} et \eqref{ccoh16e}, 
on en déduit que la fibre du morphisme \eqref{ccoh3a} en $\varrho(\oy\rightsquigarrow \ox')$ 
s'identifie donc à la limite inductive du morphisme \eqref{mtht20a}
\begin{equation}
(\cC^{\oy,(r)}_U/p^n\cC^{\oy,(r)}_U)\otimes_{\hoR^{\oy}_{U}}\oR^{\intern \oy'}_{U'\rightarrow U}\rightarrow 
(\cC'^{\oy',(t,r)}_{U'\rightarrow U}/p^n\cC'^{\oy',(t,r)}_{U'\rightarrow U})^{\Pi_{U'\rightarrow U}}
\end{equation}
lorsque $(\ox'\rightarrow U'\rightarrow U)$ décrit la catégorie $\mJ^\circ_{\ox'}$.
La proposition résulte alors de \ref{mtht20} compte tenu de (\cite{ag} 6.5.3).

\begin{cor}\label{ccoh13}
Soient $r,t,t'$ trois nombres rationnels tels que $t'>t>r\geq 0$. Alors,
\begin{itemize}
\item[{\rm (i)}] L'homomorphisme canonique de $\tG^{\mN^\circ}_s$ \eqref{ahttfg20d}
\begin{equation}\label{ccoh13a}
\bvlgg^*(\bvcC^{(r)})\rightarrow \bvtau_*(\bvcC'^{(t,r)})
\end{equation}
est $\alpha$-injectif. Notons $\bvcH^{(t,r)}$ son conoyau. 
\item[{\rm (ii)}] Il existe un nombre rationnel $a >0$ tel que le morphisme 
\begin{equation}\label{ccoh13b}
\bvcH^{(t',r')}\rightarrow \bvcH^{(t,r)}
\end{equation} 
induit par l'homomorphisme canonique $\bvalpha'^{t',t,r,r}$ \eqref{ahttfg20e} soit annulé par $p^a$. 
\item[{\rm (iii)}] Il existe un nombre rationnel $b>0$ tel que pour tout entier $q\geq 1$, 
le morphisme canonique de $\tG^{\mN^\circ}_s$
\begin{equation}\label{ccoh13c}
\rR^q\bvtau_*(\bvcC'^{(t',r)})\rightarrow \rR^q\bvtau_*(\bvcC^{(t,r)})
\end{equation}
soit annulé par $p^b$.
\end{itemize}
\end{cor}

Cela résulte de \ref{ccoh3} et (\cite{agt} III.7.3(i) et (III.7.5.5)). 

\begin{cor}\label{ccoh14}
Soient $r,t,t'$ trois nombres rationnels tels que $t'>t>r\geq 0$. Alors,
\begin{itemize}
\item[{\rm (i)}] Le morphisme canonique de $\bvlgg^*(\bvcC^{(r)})_\mQ$-modules \eqref{ahttfg20d}
\begin{equation}\label{ccoh14a}
\ttu^{t,r}\colon \bvlgg^*(\bvcC^{(r)})_\mQ\rightarrow \bvtau_*(\bvcC'^{(t,r)})_\mQ
\end{equation}
admet un inverse à gauche canonique
\begin{equation}\label{ccoh14b}
\ttv^{t,r}\colon \bvtau_*(\bvcC'^{(t,r)})_\mQ\rightarrow \bvlgg^*(\bvcC^{(r)})_\mQ.
\end{equation}
\item[{\rm (ii)}] Le composé 
\begin{equation}\label{ccoh14c}
\bvtau_*(\bvcC'^{(t',r)})_\mQ\stackrel{\ttv^{t',r}}{\longrightarrow} \bvlgg^*(\bvcC^{(r)})_\mQ
\stackrel{\ttu^{t,r}}{\longrightarrow} \bvtau_*(\bvcC'^{(t,r)})_\mQ
\end{equation}
est l'homomorphisme canonique.
\item[{\rm (iii)}] Pour tout entier $q\geq 1$, le morphisme canonique 
\begin{equation}\label{ccoh14d}
\rR^q\bvtau_*(\bvcC'^{(t',r)})_\mQ\rightarrow \rR^q\bvtau_*(\bvcC'^{(t,r)})_\mQ
\end{equation}
est nul. 
\end{itemize}
\end{cor}

En effet, d'après \ref{ccoh13}(i)-(ii), $\ttu^{t,r}$ est injectif et il existe un et un unique morphisme $\bvlgg^*(\bvcC^{(r)})_\mQ$-linéaire
\begin{equation}
\ttv^{t',t,r}\colon \bvtau_*(\bvcC'^{(t',r)})_\mQ\rightarrow  \bvlgg^*(\bvcC^{(r)})_\mQ
\end{equation}
tel que $\ttu^{t,r}\circ \ttv^{t',t,r}$ soit l'homomorphisme canonique $\bvtau_*(\bvcC'^{(t',r)})_\mQ\rightarrow \bvtau_*(\bvcC'^{(t,r)})_\mQ$. 
Comme on a $\ttu^{t,r}\circ \ttv^{t',t,r}\circ \ttu^{t',r}=\ttu^{t,r}$, on en déduit que $\ttv^{t',t,r}$ est un inverse à gauche de $\ttu^{t',r}$.  
On vérifie aussitôt qu'il ne dépend pas de $t$; d'où les propositions (i) et (ii).  
La proposition (iii) résulte aussitôt de \ref{ccoh13}(iii).

\begin{cor}\label{ccoh15}
Pour tout nombre rationnel $r\geq 0$, l'homomorphisme canonique 
\begin{equation}\label{ccoh15a}
\bvlgg^*(\bvcC^{(r)})_\mQ\rightarrow \underset{\underset{t\in \mQ_{>r}}{\longrightarrow}}{\lim}\  \bvtau_*(\bvcC'^{(t,r)})_\mQ
\end{equation}
est un isomorphisme, et pour tout entier $q\geq 1$,
\begin{equation}\label{ccoh15b}
\underset{\underset{t\in \mQ_{>r}}{\longrightarrow}}{\lim}\ \rR^q\bvtau_*(\bvcC'^{(t,r)})_\mQ=0,
\end{equation}
où les limites sont prises dans $\bMod_\mQ(\bvlgg^*(\bvcC^{(r)}))$. Celles-ci sont en particulier représentables. 
\end{cor}

L'énoncé vaut encore si l'on considère les différents $\bvlgg^*(\bvcC^{(r)})_\mQ$-modules comme des  
ind-$\bvlgg^*(\bvcC^{(r)})$-modules via le foncteur $\alpha_{\bvlgg^*(\bvcC^{(r)})}$ \eqref{indsh20c} et si l'on remplace 
$\underset{\underset{t\in \mQ_{>r}}{\longrightarrow}}{\lim}$ par $\underset{\underset{t\in \mQ_{>r}}{\longrightarrow}}{\mlq\mlq\lim \mrq\mrq}$.

\subsection{}\label{ccoh5}
Soient $t,r$ deux nombres rationnels tels que  $t\geq r\geq 0$, $n$ un entier $\geq 1$. On note
\begin{equation}\label{ccoh5a}
d_{\cC'^{(t)}_n}\colon \cC'^{(t)}_n\rightarrow \sigma'^*_n(\txi^{-1}\tOmega^1_{\oX'_n/\oS_n})\otimes_{\ocB'_n}\cC'^{(t)}_n
\end{equation}
la $\ocB'_n$-dérivation universelle de $\cC'^{(t)}_n$ (cf. \ref{ahttf23}). 
C'est un $\ocB'_n$-champ de Higgs à coefficients dans $\sigma'^*_n(\txi^{-1}\tOmega^1_{\oX'_n/\oS_n})$ d'après \ref{ahttf15}(ii). 
On désigne par
\begin{equation}\label{ccoh5b}
\ud_{\cC'^{(t)}_n}\colon \cC'^{(t)}_n\rightarrow \sigma'^*_n(\txi^{-1}\tOmega^1_{\oX'_n/\oX_n})\otimes_{\ocB'_n}\cC'^{(t)}_n
\end{equation}
la $\ocB'_n$-dérivation induite par $d_{\cC'^{(t)}_n}$. C'est en fait la $\uptheta_n^*(\cC_n^{(t)})$-dérivation universelle de $\cC'^{(t)}_n$. 
Cela résulte de la suite exacte 
\begin{equation}
\sigma'^*_n(\ogg^*_n(\txi^{-1}\tOmega^1_{\oX_n/\oS_n}))\otimes_{\ocB'_n}\cC'^{(t)}_n \rightarrow
\sigma'^*_n(\txi^{-1}\tOmega^1_{\oX'_n/\oS_n})\otimes_{\ocB'_n}\cC'^{(t)}_n\rightarrow \Omega^1_{\cC'^{(t)}_n/\uptheta^*_n(\cC^{(t)}_n)}
\rightarrow 0,
\end{equation}
compte tenu de \ref{ahttf23}, \eqref{hmdf19l} et (cf. \cite{illusie1} II 1.1.2). 
C'est donc un $\uptheta_n^*(\cC_n^{(t)})$-champ de Higgs à coefficients dans $\sigma'^*_n(\txi^{-1}\tOmega^1_{\oX'_n/\oX_n})\otimes_{\ocB'_n}\uptheta_n^*(\cC_n^{(t)})$ 
d'après \ref{MH8}(i).

Le morphisme  
\begin{equation}\label{ccoh5c}
\ud_{\cC'^{(t,r)}_n}\colon \cC'^{(t,r)}_n\rightarrow \sigma'^*_n(\txi^{-1}\tOmega^1_{\oX'_n/\oX_n})\otimes_{\ocB'_n}\cC'^{(t,r)}_n
\end{equation}
déduit de $\ud_{\cC'^{(t)}_n}$ par extension des scalaires, est la $\uptheta_n^*(\cC_n^{(r)})$-dérivation universelle de $\cC'^{(t,r)}_n$ \eqref{ahttfg18f}. 
C'est un $\uptheta_n^*(\cC_n^{(r)})$-champ de Higgs à coefficients dans $\sigma'^*_n(\txi^{-1}\tOmega^1_{\oX'_n/\oX_n})\otimes_{\ocB'_n}\uptheta_n^*(\cC_n^{(r)})$.
On désigne par $\umK^\bullet(\cC'^{(t,r)}_n)$ le complexe de Dolbeault
du $\uptheta_n^*(\cC_n^{(r)})$-module de Higgs $(\cC'^{(t,r)}_n,p^t\ud_{\cC'^{(t,r)}_n})$ \eqref{MH2a}
et par $\utmK^\bullet(\cC'^{(t,r)}_n)$ le complexe de Dolbeault augmenté
\begin{equation}\label{ccoh5d}
\uptheta_n^*(\cC^{(r)}_n)\rightarrow \umK^0(\cC'^{(t,r)}_n)\rightarrow \umK^1(\cC'^{(t,r)}_n)\rightarrow 
\umK^2(\cC'^{(t,r)}_n)\rightarrow \dots,
\end{equation}
où $\uptheta_n^*(\cC^{(r)}_n)$ est placé en degré $-1$ et la première différentielle est l'homomorphisme canonique. 

Pour tous nombres rationnels $t',r'$ tels que $t'\geq t$ et $t'\geq r'\geq r$, on a \eqref{ahttf23c}
\begin{equation}\label{ccoh5e}
p^{t'}(\id \otimes \alpha'^{t',t,r',r}_n) \circ \ud_{\cC'^{(t',r')}_n}=p^t\ud_{\cC'^{(t,r)}_n}\circ \alpha'^{t',t,r',r}_n,
\end{equation}
où $\alpha'^{t',t,r',r}_n$ est l'homomorphisme \eqref{ahttfg18g}.
Par suite, $\alpha'^{t',t,r,r}_n$ induit un morphisme 
\begin{equation}\label{ccoh5f}
\upiota^{t',t,r}_n\colon \utmK^\bullet(\cC'^{(t',r)}_n)\rightarrow \utmK^\bullet(\cC'^{(t,r)}_n).
\end{equation}

\begin{prop}\label{ccoh6}
Pour tous nombres rationnels $r,t,t'$ tels que $t'>t>r\geq 0$, 
il existe un nombre rationnel $a\geq 0$ tel que pour tous entiers $n$ et $q$ avec $n\geq 0$, le morphisme 
\begin{equation}\label{ccoh6a}
\rH^q(\upiota^{t',t,r}_n)\colon \rH^q(\utmK^\bullet(\cC'^{(t',r)}_n))\rightarrow 
\rH^q(\utmK^\bullet(\cC'^{(t,r)}_n))
\end{equation}
soit annulé par $p^a$. 
\end{prop}

La question étant locale sur $X$ (\cite{agt} VI.10.14  et \cite{ag} 3.4.10), 
on peut supposer que le schéma $X$ est séparé et que le morphisme $f\colon (X,\cM_X)\rightarrow (S,\cM_S)$ 
admet une carte adéquate $((P,\gamma),(\mN,\iota),\vartheta)$. 
 
Soient $(\oy'\rightsquigarrow \ox')$ un point de $X'_\et\gtimes_{X'_\et}\oX'^\rhd_\et$  tel que $\ox'$ soit au-dessus de $s$, $n$ un entier $\geq 0$.
Reprenons les notations introduites dans \ref{ccoh0}, \ref{ccoh00} et \ref{ccoh16}. 
Il résulte de (\cite{agt} III.10.30) que 
la fibre de la dérivation $\ud_{\cC'^{(t,r)}_n}$ \eqref{ccoh5c} en $\rho'(\oy' \rightsquigarrow \ox')$ s'identifie à la limite inductive 
des $(\cC^{\oy,(r)}_U/p^n\cC^{\oy,(r)}_U)$-dérivations universelles 
\begin{equation}
\cC'^{\oy',(t,r)}_{U'\rightarrow U}/p^n\cC'^{\oy',(t,r)}_{U'\rightarrow U}\rightarrow \txi^{-1}\tOmega^1_{X'/X}(U')\otimes_{\co_{X'}(U')} 
(\cC'^{\oy',(t,r)}_{U'\rightarrow U}/p^n\cC'^{\oy',(t,r)}_{U'\rightarrow U})
\end{equation}
de $\cC'^{\oy',(t,r)}_{U'\rightarrow U}/p^n\cC'^{\oy',(t,r)}_{U'\rightarrow U}$,
lorsque $(\ox'\rightarrow U'\rightarrow U)$ décrit la catégorie $\mJ^\circ_{\ox'}$ (cf. \ref{cdlbr1}).
La proposition résulte alors de \ref{cdlbr2}(i) compte tenu de (\cite{agt} III.9.5).

\subsection{}\label{ccoh11}
Soient $r,t$ deux nombres rationnels tels que $t\geq r\geq 0$. On note
\begin{equation}\label{ccoh11a}
d_{\bvcC'^{(t)}}\colon \bvcC'^{(t)}\rightarrow \hupsigma'^*(\txi^{-1}\tOmega^1_{\fX'/\cS})\otimes_{\bvocB'}\bvcC'^{(t)}
\end{equation}
la $\bvocB'$-dérivation universelle de $\bvcC'^{(t)}$ \eqref{ahttf30a}, où $\hupsigma'$ est le morphisme de topos annelés défini dans \eqref{hmdf13i}.
C'est un $\bvocB'$-champ de Higgs à coefficients dans $\bvsigma^*(\txi^{-1}\tOmega^1_{\fX'/\cS})$ d'après \ref{ahttf15}(iii). 
On désigne par
\begin{equation}\label{ccoh11b}
\ud_{\bvcC'^{(t)}}\colon \bvcC'^{(t)}\rightarrow \hupsigma'^*(\txi^{-1}\tOmega^1_{\fX'/\fX})\otimes_{\bvocB'}\bvcC'^{(t)}
\end{equation}
la $\bvocB'$-dérivation induite par $d_{\bvcC'^{(t)}}$. 
C'est en fait la $\bvuptheta^*(\bvcC^{(t)})$-dérivation universelle de $\bvcC'^{(t)}$. Cela résulte de la suite exacte 
\begin{equation}\label{ccoh11g}
\hupsigma'^*(\fgg^*(\txi^{-1}\tOmega^1_{\fX/\cS}))\otimes_{\bvocB'}\bvcC'^{(t)} \rightarrow
\hupsigma'^*(\txi^{-1}\tOmega^1_{\fX'/\cS})\otimes_{\bvocB'}\bvcC'^{(t)}\rightarrow \Omega^1_{\bvcC'^{(t)}/\bvuptheta^*(\bvcC^{(t)})}
\rightarrow 0,
\end{equation}
compte tenu de \eqref{ahttf30a} et \eqref{hmdf13i} (cf. \cite{illusie1} II 1.1.2). 
On observera que $\ud_{\bvcC'^{(t)}}$ s'identifie au morphisme $(\ud_{\cC'^{(t)}_{n+1}})_{n\in \mN}$ \eqref{ccoh5b},
et que c'est un $\bvuptheta^*(\bvcC^{(t)})$-champ de Higgs à coefficients dans $\hupsigma'^*(\txi^{-1}\tOmega^1_{\fX'/\fX})\otimes_{\bvocB'}\bvuptheta^*(\bvcC^{(t)})$.

Le morphisme   
\begin{equation}\label{ccoh11c}
\ud_{\bvcC'^{(t,r)}}\colon \bvcC'^{(t,r)}\rightarrow \hupsigma'^*(\txi^{-1}\tOmega^1_{\fX'/\fX})\otimes_{\bvocB'}\bvcC'^{(t,r)},
\end{equation}
déduit par extension des scalaires de $\ud_{\bvcC'^{(t)}}$, est la $\bvuptheta^*(\bvcC^{(r)})$-dérivation universelle de $\bvcC'^{(t,r)}$.  
C'est un $\bvuptheta^*(\bvcC^{(r)})$-champ de Higgs à coefficients dans $\hupsigma'^*(\txi^{-1}\tOmega^1_{\fX'/\fX})\otimes_{\bvocB'}\bvuptheta^*(\bvcC^{(r)})$. 
On désigne par $\umK^\bullet(\bvcC'^{(t,r)})$ le complexe de Dolbeault
du $\bvuptheta^*(\bvcC^{(r)})$-module de Higgs $(\bvcC'^{(t,r)},p^t\ud_{\bvcC'^{(t,r)}})$,
et par $\utmK^\bullet(\bvcC'^{(t,r)})$ le complexe de Dolbeault augmenté
\begin{equation}\label{ccoh11d}
\bvuptheta^*(\bvcC^{(r)})\rightarrow \umK^0(\bvcC'^{(t,r)})\rightarrow \umK^1(\bvcC'^{(t,r)})\rightarrow \dots 
\rightarrow \umK^n(\bvcC'^{(t,r)})\rightarrow \dots,
\end{equation}
où $\bvuptheta^*(\bvcC^{(r)})$ est placé en degré $-1$ et la première différentielle est l'homomorphisme canonique. 

Pour tous nombres rationnels $t',r'$ tels que $t'\geq t$ et $t'\geq r'\geq r$, on a \eqref{ccoh5e}
\begin{equation}\label{ccoh11e}
p^{t'}(\id \otimes \bvalpha'^{t',t,r',r}) \circ \ud_{\bvcC'^{(t',r')}}=p^t\ud_{\bvcC'^{(t,r)}}\circ \bvalpha'^{t',t,,r',r},
\end{equation}
où $\bvalpha'^{t',t,r',r}$ est l'homomorphisme \eqref{ahttfg20e}. Par suite, $\bvalpha'^{t',t,r,r}$ induit un morphisme de complexes 
\begin{equation}\label{ccoh11f}
\bvupiota^{t',t,r}\colon \utmK^\bullet(\bvcC'^{(t',r)})\rightarrow \utmK^\bullet(\bvcC'^{(t,r)}).
\end{equation}

On note $\umK^\bullet_{\mQ}(\bvcC'^{(t,r)})$ et $\utmK^\bullet_{\mQ}(\bvcC'^{(t,r)})$
les images des complexes $\umK^\bullet(\bvcC'^{(t,r)})$ et $\utmK^\bullet(\bvcC'^{(t,r)})$
dans $\bMod_\mQ(\bvuptheta^*(\bvcC^{(r)}))$ \eqref{ahttfg20}.
On considèrera ces complexes aussi comme des complexes de 
la catégorie $\bIndMod(\bvuptheta^*(\bvcC^{(r)}))$ via le foncteur $\alpha_{\bvuptheta^*(\bvcC^{(r)})}$ \eqref{indsh20c}.

\begin{prop}\label{ccoh7}
Pour tous nombres rationnels  $r,t,t'$ tels que $t'>t>r\geq 0$  et tout entier $q$, le morphisme canonique de 
$\bvuptheta^*(\bvcC^{(r)})_\mQ$-modules \eqref{ccoh11f} 
\begin{equation}\label{ccoh7a}
\rH^q(\bvupiota^{t',t,r}_\mQ)\colon \rH^q(\utmK^\bullet_{\mQ}(\bvcC'^{(t',r)}))\rightarrow 
\rH^q(\utmK^\bullet_{\mQ}(\bvcC'^{(t,r)}))
\end{equation}
est nul.
\end{prop}

Cela résulte de \ref{ccoh6} et (\cite{agt} III.7.3(i)).

\begin{cor}\label{ccoh8}
Soient $r,t,t'$ trois nombres rationnels tels que $t'>t>r\geq 0$. Alors,
\begin{itemize}
\item[{\rm (i)}] Le morphisme canonique de $\bvuptheta^*(\bvcC^{(r)})_\mQ$-modules
\begin{equation}\label{ccoh8a}
u^{t,r}\colon \bvuptheta^*(\bvcC^{(r)})_\mQ\rightarrow \rH^0(\umK^\bullet_{\mQ}(\bvcC'^{(t,r)}))
\end{equation}
admet un inverse à gauche canonique 
\begin{equation}\label{ccoh8b}
v^{t,r}\colon  \rH^0(\umK^\bullet_{\mQ}(\bvcC'^{(t,r)})) \rightarrow \bvuptheta^*(\bvcC^{(r)})_\mQ.
\end{equation}
\item[{\rm (ii)}] Le composé 
\begin{equation}
\xymatrix{
{\rH^0(\umK^\bullet_{\mQ}(\bvcC'^{(t',r)}))}\ar[r]^-(0.5){v^{t',r}}&{\bvuptheta^*(\bvcC^{(r)})_\mQ}\ar[r]^-(0.5){u^{t,r}}&
{\rH^0(\umK^\bullet_{\mQ}(\bvcC'^{(t,r)}))}}
\end{equation}
est le morphisme canonique.
\item[{\rm (iii)}] Pour tout entier $q\geq 1$, le morphisme canonique 
\begin{equation}\label{ccoh8c}
\rH^q(\umK^\bullet_{\mQ}(\bvcC'^{(t',r)}))\rightarrow 
\rH^q(\umK^\bullet_{\mQ}(\bvcC'^{(t,r)}))
\end{equation}
est nul.
\end{itemize}
\end{cor}

En effet, considérons le diagramme commutatif canonique (sans la flèche pointillée)
\begin{equation}
\xymatrix{
{\bvuptheta^*(\bvcC^{(r)})_\mQ}\ar@{=}[d]\ar[r]^-(0.5){u^{t',r}}&{\rH^0(\umK^\bullet_{\mQ}(\bvcC'^{(t',r)}))}\ar@{.>}[ld]_-(0.5){v^{t',t,r}}
\ar[d]^-(0.5){\varpi^{t',t,r}}
\ar@{->>}[r]&{\rH^0(\utmK^\bullet_{\mQ}(\bvcC'^{(t',r)}))}\ar[d]^{\rH^0(\bvupiota^{t',t,r}_\mQ)}\\
{\bvuptheta^*(\bvcC^{(r)})_\mQ}\ar[r]_-(0.5){u^{t,r}}& {\rH^0(\umK^\bullet_{\mQ}(\bvcC'^{(t,r)}))}\ar@{->>}[r]&
{\rH^0(\utmK^\bullet_{\mQ}(\bvcC'^{(t,r)}))}}
\end{equation}
Il résulte de \ref{ccoh7} que $u^{t',r}$ et par suite $u^{t,r}$ sont injectifs, et qu'il existe un et un unique morphisme
$v^{t',t,r}$ comme ci-dessus tel que $\varpi^{t',t,r}=u^{t,r}\circ v^{t',t,r}$. Comme on a 
$u^{t,r}\circ v^{t',t,r}\circ u^{t',r}=u^{t,r}$, on en déduit que $v^{t',t,r}$ est un inverse à gauche de $u^{t',r}$.  
On vérifie aussitôt qu'il ne dépend pas de $t$; d'où les propositions (i) et (ii).  
La proposition (iii) résulte aussitôt de \ref{ccoh7}.

\begin{cor}\label{ccoh17}
Pour tout nombre rationnel $r\geq 0$, le morphisme canonique de complexes de ind-$\bvuptheta^*(\bvcC^{(r)})$-modules 
\begin{equation}\label{ccoh17a}
\bvuptheta^*(\bvcC^{(r)})_\mQ[0]\rightarrow \underset{\underset{r\in \mQ_{>0}}{\longrightarrow}}{\mlq\mlq\lim \mrq\mrq}\ \umK^\bullet_{\mQ}(\bvcC'^{(t,r)})
\end{equation}
est un quasi-isomorphisme.
\end{cor}

On rappelle d'abord que $\bIndMod(\bvuptheta^*(\bvcC^{(r)}))$ admet des petites limites inductives et que les petites limites inductives filtrantes sont exactes \eqref{indsh6e}. 
La proposition est donc équivalente au fait que le morphisme canonique 
\begin{equation}\label{ccoh17b}
\bvuptheta^*(\bvcC^{(r)})_\mQ\rightarrow \underset{\underset{t\in \mQ_{>r}}{\longrightarrow}}{\mlq\mlq\lim \mrq\mrq}\ 
\rH^0(\umK^\bullet_{\mQ}(\bvcC'^{(t,r)}))
\end{equation}
est un isomorphisme, et que pour tout entier $q\geq 1$, 
\begin{equation}\label{ccoh17c}
\underset{\underset{t\in \mQ_{>r}}{\longrightarrow}}{\mlq\mlq\lim \mrq\mrq}\ \rH^q(\umK^\bullet_{\mQ}(\bvcC'^{(t,r)}))=0.
\end{equation}
Ces énoncés résultent aussitôt de \ref{ccoh8}.

\begin{rema}\label{ccoh9}
Soit $r$ un nombre rationnel $\geq 0$. 
Bien que les limites inductives filtrantes ne soient pas a priori représentables dans la catégorie $\bMod_\mQ(\bvuptheta^*(\bvcC^{(r)}))$, 
il résulte de \ref{ccoh8} que dans cette catégorie, le morphisme canonique 
\begin{equation}\label{ccoh9a}
\bvuptheta^*(\bvcC^{(r)})_\mQ\rightarrow \underset{\underset{t\in \mQ_{>r}}{\longrightarrow}}{\lim}\ 
\rH^0(\umK^\bullet_{\mQ}(\bvcC'^{(t,r)}))
\end{equation}
est un isomorphisme, et pour tout entier $q\geq 1$, 
\begin{equation}\label{ccoh9b}
\underset{\underset{t\in \mQ_{>r}}{\longrightarrow}}{\lim}\ \rH^q(\umK^\bullet_{\mQ}(\bvcC'^{(t,r)}))=0.
\end{equation}
\end{rema}

\subsection{}\label{ccoh50}
Soient $t,r$ deux nombres rationnels tels que  $t\geq r\geq 0$. On note 
\begin{equation}\label{ccoh50f}
v\colon \bvuptheta^*(\hupsigma^*(\txi^{-1}\tOmega^1_{\fX/\cS})) \rightarrow \hupsigma'^*(\txi^{-1}\tOmega^1_{\fX'/\cS})
\end{equation}
le morphisme induit par le morphisme canonique $\fgg^*(\tOmega^1_{\fX/\cS})\rightarrow \tOmega^1_{\fX'/\cS}$ et le diagramme commutatif \eqref{hmdf13i}. 
D'après \ref{ahttfg7}(ii), le diagramme 
\begin{equation}\label{ccoh50a}
\xymatrix{
{\bvuptheta^*(\bvcC^{(t)})}\ar[rr]^-(0.5){\bvuptheta^*(d_{\bvcC^{(t)}})}\ar[d]_{\bvupomega^{(t)}}&&
{\bvuptheta^*(\hupsigma^*(\txi^{-1}\tOmega^1_{\fX/\cS}))\otimes_{\ocB'_n}\bvuptheta^*(\bvcC^{(t)})}\ar[d]^{v\otimes\bvupomega^{(t)}}\\
{\bvcC'^{(t)}}\ar[rr]^-(0.5){d_{\bvcC'^{(t)}}}&&{\hupsigma'^*(\txi^{-1}\tOmega^1_{\fX'/\cS})\otimes_{\bvocB'}\bvcC'^{(t)}}}
\end{equation}
où $\bvupomega^{(t)}$ est l'homomorphisme \eqref{ahttfg20a}, est commutatif. Compte tenu de \eqref{ahttf30c}, 
il existe donc une $\bvocB'$-dérivation de l'algèbre $\bvcC'^{(t,r)}=\bvcC'^{(t)}\otimes_{\bvuptheta^*(\bvcC^{(t)})}\bvuptheta^*(\bvcC^{(r)})$,
\begin{equation}\label{ccoh50b}
\delta_{\bvcC'^{(t,r)}}\colon \bvcC'^{(t,r)}\rightarrow \hupsigma'^*(\txi^{-1}\tOmega^1_{\fX'/\cS})\otimes_{\bvocB'}\bvcC'^{(t,r)},
\end{equation}
définie par 
\begin{equation}\label{ccoh50c}
\delta_{\bvcC'^{(t,r)}}=p^td_{\bvcC'^{(t)}}\otimes \id_{\bvuptheta^*(\bvcC^{(r)})} 
+(v\otimes_{\bvocB'}\id_{\bvcC'^{(t,r)}})\circ(\id_{\bvcC'^{(t)}} \otimes \bvuptheta^*(p^r d_{\bvcC^{(r)}})). 
\end{equation}
Comme $d_{\bvcC'^{(t)}}$ et $d_{\bvcC^{(r)}}$ sont des champs de Higgs d'après \ref{ahttf15}(ii), 
on vérifie aussitôt que $\delta_{\bvcC'^{(t,r)}}$ est un $\bvocB'$-champ de Higgs à coefficients dans $\hupsigma'^*(\txi^{-1}\tOmega^1_{\fX'/\cS})$.
Le morphisme 
\begin{equation}\label{ccoh50d}
\udelta_{\bvcC'^{(t,r)}}\colon \bvcC'^{(t,r)}\rightarrow \hupsigma'^*(\txi^{-1}\tOmega^1_{\fX'/\fX})\otimes_{\bvocB'}\bvcC'^{(t,r)}
\end{equation}
induit par $\delta_{\bvcC'^{(t,r)}}$ n'est autre que $p^t\ud_{\bvcC'^{(t,r)}}$ \eqref{ccoh11c}. C'est donc une $\bvuptheta^*(\bvcC^{(r)})$-dérivation et un 
$\bvuptheta^*(\bvcC^{(r)})$-champ de Higgs à coefficients dans $\hupsigma'^*(\txi^{-1}\tOmega^1_{\fX'/\fX})\otimes_{\bvocB'}\bvuptheta^*(\bvcC^{(r)})$

Pour tous nombres rationnels $t,t',r,r'$ tels que $t'\geq t\geq r\geq 0$ et $t'\geq r'\geq r$, on a \eqref{ahttf23c}
\begin{equation}\label{ccoh50e}
(\id \otimes \bvalpha'^{t',t,r',r}) \circ \delta_{\bvcC'^{(t',r')}}=\delta_{\bvcC'^{(t,r)}}\circ \bvalpha'^{t',t,r',r},
\end{equation}
où $\bvalpha'^{t',t,r',r}$ est l'homomorphisme \eqref{ccoh11e}.

\subsection{}\label{ccoh51}
Soient $t,r$ deux nombres rationnels tels que  $t\geq r\geq 0$. 
On désigne par $\mK^\bullet(\bvcC'^{(t,r)})$ le complexe de Dolbeault du $\bvocB'$-module de Higgs $(\bvcC'^{(t,r)},\delta_{\bvcC'^{(t,r)}})$ \eqref{ccoh50b} et par
$\umK^\bullet(\bvcC'^{(t,r)})$ le complexe de Dolbeault du $\bvuptheta^*(\bvcC^{(r)})$-module de Higgs $(\bvcC'^{(t,r)},p^t\ud_{\bvcC'^{(t,r)}})$ \eqref{ccoh11c}. 
On munit $\mK^\bullet(\bvcC'^{(t,r)})$ de la filtration de Koszul \eqref{MH90g} associée à l'image par le foncteur $\hupsigma'^*$ de la suite exacte \eqref{hmdf410d}, 
et on note 
\begin{equation}\label{ccoh51a}
\partial^{(t,r)}\colon \umK^\bullet(\bvcC'^{(t,r)})\rightarrow 
\hupsigma'^*(\fgg^*(\txi^{-1}\tOmega^1_{\fX/\cS})) \otimes_{\bvocB'} \umK^\bullet(\bvcC'^{(t,r)})
\end{equation}
le bord associé dans $\bD^+(\bMod(\bvocB'))$ \eqref{MH90j}. Le diagramme de morphismes de $\bD^+(\bMod(\bvocB'))$
\begin{equation}\label{ccoh51b}
\xymatrix{
{\bvuptheta^*(\bvcC^{(r)})[0]}\ar[rr]^-(0.5){p^r\bvuptheta^*(d_{\bvcC^{(r)}})}\ar[d]&&{\bvuptheta^*(\hupsigma^*(\txi^{-1}\tOmega^1_{\fX/\cS})\otimes_{\bvocB}\bvcC^{(r)})[0]}\ar[d]\\
{\umK^\bullet(\bvcC'^{(t,r)})}\ar[rr]^-(0.5){\partial^{(t,r)}}&&{\hupsigma'^*(\fgg^*(\txi^{-1}\tOmega^1_{\fX/\cS})) \otimes_{\bvocB'} \umK^\bullet(\bvcC'^{(t,r)})}}
\end{equation}
où les flèches verticales sont les morphismes canoniques, cf. \eqref{ccoh11d} et \eqref{hmdf13i}, est commutatif. 
Cela résulte de \ref{MH100} en prenant pour suite exacte \eqref{MH90a} l'image par le foncteur $\hupsigma'^*$ de la suite exacte \eqref{hmdf410d}, 
pour $(M,\theta)$ le $\bvocB'$-module de Higgs trivial $(\bvocB',0)$ et pour $(N,\kappa)$ le
$\bvocB'$-module de Higgs $(\bvuptheta^*(\bvcC^{(r)}),\bvuptheta^*(d_{\bvcC^{(r)}}))$.
Avec les notations de \ref{MH98}, on a $\utheta'=0$, de sorte que les différentielles de $\umK'^\bullet$ sont nulles. 
L'assertion recherchée résulte alors de \ref{MH100} par fonctorialité du morphisme bord \eqref{MH90j}.

\section{Cohomologie relative des ind-modules de Dolbeault}\label{crmd}

\subsection{}\label{crindmd1}
On désigne par $\bIndMod^\Dolb(\bvocB')$ la catégorie des ind-$\bvocB'$-modules de Dolbeault \eqref{indmdlb5}
et par $\bMH^\sol(\co_{\fX'}[\frac 1 p], \txi^{-1}\tOmega^1_{\fX'/\cS})$ 
la catégorie des $\co_{\fX'}[\frac 1 p]$-fibrés de Higgs solubles à coefficients dans $\txi^{-1}\tOmega^1_{\fX'/\cS}$,
relativement à la déformation $(\tX',\cM_{\tX'})$ fixée dans \ref{hmdf5}. On note  
\begin{equation}\label{crindmd1a}
\cH'\colon \bIndMod(\bvocB')\rightarrow \bMH(\co_{\fX'}, \txi^{-1}\tOmega^1_{\fX'/\cS})
\end{equation}
le foncteur défini dans \eqref{indmdlb7c}, relativement à la déformation $(\tX',\cM_{\tX'})$.
D'après \ref{indmdlb20}, celui-ci induit une équivalence de catégories que l'on note encore
\begin{equation}\label{crindmd1b}
\cH'\colon \bIndMod^\Dolb(\bvocB')\stackrel{\sim}{\rightarrow} \bMH^\sol(\co_{\fX'}[\frac 1 p], \txi^{-1}\tOmega^1_{\fX'/\cS}).
\end{equation}

On désigne par 
\begin{equation}\label{crindmd1d}
\ucH'\colon \bIndMod(\bvocB')\rightarrow \bMH(\co_{\fX'}, \txi^{-1}\tOmega^1_{\fX'/\fX}). 
\end{equation}
le composé du foncteur $\cH'$ et du foncteur canonique \eqref{hmdf410d}
\begin{equation}\label{crindmd1e}
\bMH(\co_{\fX'}, \txi^{-1}\tOmega^1_{\fX'/\cS})\rightarrow \bMH(\co_{\fX'}, \txi^{-1}\tOmega^1_{\fX'/\fX}). 
\end{equation}

\subsection{}\label{crindmd2}
Soient $t,r$ deux nombres rationnels tels que $t\geq r \geq 0$. 
On associe à $(f',\tX',\cM_{\tX'})$ des objets analogues à ceux associés à $(f,\tX,\cM_{\tX})$ dans \ref{indmdlb1}--\ref{indmdlb3} 
qu'on note par les mêmes symboles affectés d'un exposant $^\prime$. On dispose en particulier de la catégorie $\bIndMC(\bvcC'^{(t)}/\bvocB')$ 
et des foncteurs $\rI\fS'^{(t)}$ \eqref{indmdlb1b}, $\rI\hupsigma'^{(t)*}$ \eqref{indmdlb1f} et $\rI\varepsilon'^{t,t'}$ \eqref{indmdlb3a} 
pour tout nombre rationnel $t'$ tel que $t\geq t'\geq 0$. 

On introduit des variantes relatives. On désigne par $\bIndMC(\bvcC'^{(t,r)}/\bvuptheta^*(\bvcC^{(r)}))$ la catégorie 
des ind-$\bvcC'^{(t,r)}$-modules à $p^t$-connexion intégrable relativement à l'extension
$\bvcC'^{(t,r)}/\bvuptheta^*(\bvcC^{(r)})$  (cf. \ref{indsh28} et \ref{indsh35}). 
Chaque objet de cette catégorie est un ind-$\bvuptheta^*(\bvcC^{(r)})$-module de Higgs
à coefficients dans $\hupsigma'^*(\txi^{-1}\tOmega^1_{\fX'/\fX})\otimes_{\bvocB'}\bvuptheta^*(\bvcC^{(r)})$ d'après \ref{indsh38}(i). 
On peut donc lui associer un complexe de Dolbeault dans $\bIndMod(\bvuptheta^*(\bvcC^{(r)}))$.  

On a le foncteur 
\begin{equation}\label{crindmd2k}
\rI\ufS'^{(t,r)}\colon
\begin{array}[t]{clcr}
\bIndMod(\bvocB')&\rightarrow &\bIndMC(\bvcC'^{(t,r)}/\bvuptheta^*(\bvcC^{(r)}))\\
\cM&\mapsto& (\bvcC'^{(t,r)}\otimes_{\bvocB'}\cM,p^t\ud_{\bvcC'^{(t,r)}}\otimes_{\bvocB'}\id_{\cM}), 
\end{array}
\end{equation}
où $\ud_{\bvcC'^{(t,r)}}$ est la $\bvuptheta^*(\bvcC^{(r)})$-dérivation universelle de $\bvcC'^{(t,r)}$ \eqref{ccoh11c}. 
Compte tenu de \ref{indsh38}(ii) et avec les notations de \ref{hmdf41},
le foncteur $\rI\hupsigma'^*$ \eqref{hmdf13i} induit un foncteur 
\begin{eqnarray}\label{crindmd2n}
\rI\uhupsigma'^{(t,r)*}\colon
\bIndMH(\co_{\fX'},\txi^{-1}\tOmega^1_{\fX'/\fX})\rightarrow \bIndMC(\bvcC'^{(t,r)}/\bvuptheta^*(\bvcC^{(r)})),\\
(\cN,\theta)\mapsto(\bvcC'^{(t,r)}\otimes_{\bvocB'}\rI\hupsigma'^*(\cN),p^t\ud_{\bvcC'^{(t,r)}} \otimes_{\bvocB'}\id+\id \otimes_{\bvocB'} \rI\hupsigma'^*(\theta)).\nonumber
\end{eqnarray}

Soit $(\cF,\nabla)$ un ind-$\bvcC'^{(t)}$-module à $p^t$-connexion intégrable relativement à l'extension $\bvcC'^{(t)}/\bvocB'$. 
On note
\begin{equation}\label{crindmd2d}
\unabla\colon \cF\rightarrow \hupsigma'^*(\txi^{-1}\tOmega^1_{\fX'/\fX}) \otimes_{\bvocB'} \cF
\end{equation}
le morphisme de ind-$\bvocB'$-modules déduit de $\nabla$. Comme la $\bvuptheta^*(\bvcC^{(t)})$-dérivation universelle  $\ud_{\bvcC'^{(t)}}$ de  
$\bvcC'^{(t)}$ \eqref{ccoh11b} est induite par $d_{\bvcC'^{(t)}}$ \eqref{ccoh11a}, on peut considérer $\unabla$ 
canoniquement comme un morphisme de ind-$\bvuptheta^*(\bvcC^{(t)})$-modules 
d'après \ref{indsh33} et \ref{indsh51}, et c'est alors une $p^t$-connexion intégrable sur $\cF$ relativement à l'extension $\bvcC'^{(t)}/\bvuptheta^*(\bvcC^{(t)})$. 
Le morphisme 
\begin{equation}\label{crindmd2e}
\unabla^{(t,r)}\colon \cF\otimes_{\bvuptheta^*(\bvcC^{(t)})}\bvuptheta^*(\bvcC^{(r)})
\rightarrow \hupsigma'^*(\txi^{-1}\tOmega^1_{\fX'/\fX}) \otimes_{\bvocB'} \cF\otimes_{\bvuptheta^*(\bvcC^{(t)})}\bvuptheta^*(\bvcC^{(r)})
\end{equation}
déduit de $\unabla$ par extension des scalaires, est  une $p^t$-connexion intégrable sur $\cF\otimes_{\bvcC'^{(t)}}\bvcC'^{(t,r)}$
 relativement à l'extension $\bvcC'^{(t,r)}/\bvuptheta^*(\bvcC^{(r)})$ \eqref{indsh44}. On obtient ainsi un foncteur 
\begin{equation}\label{crindmd2f}
\rI\ulambda^{t,r}\colon
\begin{array}[t]{clcr} 
\bIndMC(\bvcC'^{(t)}/\bvocB')&\rightarrow &\bIndMC(\bvcC'^{(t,r)}/\bvuptheta^*(\bvcC^{(r)}))\\
(\cF,\nabla)&\mapsto&(\cF\otimes_{\bvcC'^{(t)}}\bvcC'^{(t,r)},\unabla^{(t,r)}).
\end{array}
\end{equation}

On vérifie aussitôt qu'on a un isomorphisme canonique de foncteurs
\begin{equation}\label{crindmd2l}
\rI\ufS'^{(t,r)} \stackrel{\sim}{\rightarrow} \rI\ulambda^{t,r}\circ \rI\fS'^{(t)}.
\end{equation}
On vérifie aussi que le diagramme 
\begin{equation}\label{crindmd2o}
\xymatrix{
{\bIndMH(\co_{\fX'},\txi^{-1}\tOmega^1_{\fX'/\cS})}\ar[r]^-(0.5){\rI\hupsigma'^{(t)*}}\ar[d]&{\bIndMC(\bvcC'^{(t)}/\bvocB')}\ar[d]^{\rI\ulambda^{t,r}}\\
{\bIndMH(\co_{\fX'},\txi^{-1}\tOmega^1_{\fX'/\fX})}\ar[r]^-(0.5){\rI\uhupsigma'^{(t,r)*}}&{\bIndMC(\bvcC'^{(t,r)}/\bvuptheta^*(\bvcC^{(r)}))}}
\end{equation}
où la flèche non libellée est le foncteur canonique, est commutatif à isomorphisme canonique près.

Soient $t',r'$ deux nombres rationnels tels que $t\geq t'\geq r'\geq 0$ et $r\geq r'$. 
D'après \ref{indsh43} et \eqref{ccoh11e}, pour tout ind-$\bvcC'^{(t,r)}$-module à $p^{t}$-connexion intégrable $(\cF,\nabla)$ relativement à l'extension 
$\bvcC'^{(t,r)}/\bvuptheta^*(\bvcC^{(r)})$, l'ind-$\bvcC'^{(t',r')}$-module $\cF\otimes_{\bvcC'^{(t,r)}}\bvcC'^{(t',r')}$ défini par extension des scalaires \eqref{ahttfg20e}
est canoniquement muni d'une $p^{t'}$-connexion intégrable 
$\nabla'$ relativement à l'extension $\bvcC'^{(t',r')}/\bvuptheta^*(\bvcC^{(r')})$. 
On définit ainsi un foncteur 
\begin{equation}\label{crindmd2g}
\rI\uvarepsilon'^{t,t',r,r'}\colon 
\begin{array}[t]{clcr}
\bIndMC(\bvcC'^{(t,r)}/\bvuptheta^*(\bvcC^{(r)}))&\rightarrow& \bIndMC(\bvcC'^{(t',r')}/\bvuptheta^*(\bvcC^{(r')}))\\
(\cF,\nabla)&\mapsto&(\cF\otimes_{\bvcC'^{(t,r)}}\bvcC'^{(t',r')},\nabla').
\end{array}
\end{equation}
On a des isomorphismes canoniques de foncteurs  
\begin{eqnarray}
\rI\uvarepsilon'^{t,t',r,r'}\circ \rI\ufS'^{(t,r)}&\stackrel{\sim}{\rightarrow}& \rI\ufS'^{(t',r')}, \label{crindmd2p}\\
\rI\uvarepsilon'^{t,t',r,r'}\circ \rI\uhupsigma'^{(t,r)*}&\stackrel{\sim}{\rightarrow}& \rI\uhupsigma'^{(t',r')*}.\label{crindmd2q}
\end{eqnarray}

Pour tous nombres rationnels $t',r',t'', r''$ tels $t\geq t'\geq t''\geq r''\geq 0$, $t'\geq r'\geq r''$ et $r\geq r'$, on a un isomorphisme canonique de foncteurs
\begin{equation}\label{crindmd2i}
\rI\uvarepsilon'^{t',t'',r',r''}\circ \rI\uvarepsilon'^{t,t',r,r'}\stackrel{\sim}{\rightarrow}\rI\uvarepsilon'^{t,t'',r,r''}.
\end{equation}

\subsection{}\label{crindmd3}
Soient $t,r$ deux nombres rationnels tels que $t\geq r \geq 0$. 
On associe à $(f',\tX',\cM_{\tX'})$ des objets analogues à ceux associés à $(f,\tX,\cM_{\tX})$ dans \ref{aspglob6}--\ref{dolbff1} 
qu'on note par les mêmes symboles affectés d'un exposant $^\prime$. On dispose en particulier de la catégorie 
$\bIMC_\mQ(\bvcC'^{(t)}/\bvocB')$ et du foncteur $\hupsigma'^{(t)*}$ \eqref{aspglob6g}.

On introduit des variantes relatives. 
On notera simplement $\bIMC(\bvcC'^{(t,r)}/\bvuptheta^*(\bvcC^{(r)}))$ la catégorie des $p^t$-isoconnexions intégrables relativement à l'extension 
$\bvcC'^{(t,r)}/\bvuptheta^*(\bvcC^{(r)})$ (cf. \ref{indsh24}); 
on omet donc l'exposant $p^t$ de la notation  introduite dans \ref{indsh24} considérant qu'il est redondant avec l'exposant de $\bvcC'^{(t,r)}$. 
C'est une catégorie additive.
On désigne par $\bIMC_\mQ(\bvcC'^{(t,r)}/\bvuptheta^*(\bvcC^{(r)}))$ la catégorie des objets de $\bIMC(\bvcC'^{(t,r)}/\bvuptheta^*(\bvcC^{(r)}))$ 
à isogénie près. D'après \ref{indsh26}(i), tout objet de $\bIMC(\bvcC'^{(t,r)}/\bvuptheta^*(\bvcC^{(r)}))$ est une $\bvuptheta^*(\bvcC^{(r)})$-isogénie
de Higgs à coefficients dans $\hupsigma'^*(\txi^{-1}\tOmega^1_{\fX'/\fX})$ \eqref{ccoh11c}. 
En particulier, on peut associer fonctoriellement à tout objet de $\bIMC_\mQ(\bvcC'^{(t,r)}/\bvuptheta^*(\bvcC^{(r)}))$ un complexe de Dolbeault dans 
$\bMod_{\mQ}(\bvuptheta^*(\bvcC^{(r)}))$ (cf. \ref{indsh23}).  

Si $(\cN,\cN',v,\theta)$ est une $\co_{\fX'}$-isogénie de Higgs
à coefficients dans $\txi^{-1}\tOmega^1_{\fX'/\fX}$, 
\begin{equation}\label{crindmd3d}
(\bvcC'^{(t,r)}\otimes_{\bvocB'}\hupsigma'^*(\cN),\bvcC'^{(t,r)}\otimes_{\bvocB'}\hupsigma'^*(\cN'),\id \otimes_{\bvocB'}\hupsigma'^*(v),
p^t\ud_{\bvcC'^{(t,r)}} \otimes\hupsigma'^*(v)+\id \otimes \hupsigma'^*(\theta))
\end{equation}
est un objet de $\bIMC(\bvcC'^{(t,r)}/\bvuptheta^*(\bvcC^{(r)}))$. Avec les notations de \ref{hmdf41}, on obtient ainsi un foncteur 
\begin{equation}\label{crindmd3e}
\uhupsigma'^{(t,r)*}\colon \bIH(\co_{\fX'},\txi^{-1}\tOmega^1_{\fX'/\fX})\rightarrow \bIMC(\bvcC'^{(t,r)}/\bvuptheta^*(\bvcC^{(r)})).
\end{equation}
Compte tenu de \eqref{hmdf41a}, celui-ci induit un foncteur que l'on note encore
\begin{equation}\label{crindmd3f}
\uhupsigma'^{(t,r)*}_\mQ\colon \bMH^\coh(\co_{\fX'}[\frac 1 p], \txi^{-1}\tOmega^1_{\fX'/\fX})\rightarrow \bIMC_\mQ(\bvcC'^{(t,r)}/\bvuptheta^*(\bvcC^{(r)})).
\end{equation}

\subsection{}\label{crindmd15}
Soient $t,r$ deux nombres rationnels tels que $t\geq r \geq 0$, $(\cF,\nabla)$ un ind-$\bvcC'^{(t)}$-module à 
$p^t$-connexion intégrable relativement à l'extension $\bvcC'^{(t)}/\bvocB'$. 
On rappelle que la $\bvocB'$-dérivation $\delta_{\bvcC'^{(t,r)}}$ \eqref{ccoh50b} de l'algèbre $\bvcC'^{(t,r)}$
est un $\bvocB'$-champ de Higgs à coefficients dans $\hupsigma'^*(\txi^{-1}\tOmega^1_{\fX'/\cS})$.
D'après \ref{indsh380}, $\nabla$ et $\delta_{\bvcC'^{(t,r)}}$ définissent sur le $\bvcC'^{(t,r)}$-module $\cF\otimes_{\bvcC'^{(t)}}\bvcC'^{(t,r)}$ 
un $\bvocB'$-champ de Higgs à coefficients dans $\hupsigma'^*(\txi^{-1}\tOmega^1_{\fX'/\cS})$ que l'on note $\nabla^{(t,r)}$. 

Pour $\cZ=\fX$ ou $\cS$, notons $\bIndMH(\bvocB',\hupsigma'^*(\txi^{-1}\tOmega^1_{\fX'/\cZ}))$  
la catégorie des ind-$\bvocB'$-modules de Higgs à coefficients dans $\hupsigma'^*(\txi^{-1}\tOmega^1_{\fX'/\cZ})$. 
On obtient alors un foncteur 
\begin{equation}\label{crindmd15c}
\rI\lambda^{t,r}\colon
\begin{array}[t]{clcr} 
\bIndMC(\bvcC'^{(t)}/\bvocB')&\rightarrow &\bIndMH(\bvocB',\hupsigma'^*(\txi^{-1}\tOmega^1_{\fX'/\cS}))\\
(\cF,\nabla)&\mapsto&(\cF\otimes_{\bvcC'^{(t)}}\bvcC'^{(t,r)},\nabla^{(t,r)}).
\end{array}
\end{equation}
On vérifie aussitôt que le diagramme 
\begin{equation}\label{crindmd15d}
\xymatrix{
{\bIndMC(\bvcC'^{(t)}/\bvocB')}\ar[r]^-(0.5){\rI\lambda^{t,r}}\ar[d]_{\rI\ulambda^{t,r}}&
{\bIndMH(\bvocB',\hupsigma'^*(\txi^{-1}\tOmega^1_{\fX'/\cS}))}\ar[d]\\
{\bIndMC(\bvcC'^{(t,r)}/\bvuptheta^*(\bvcC^{(r)}))}\ar[r]&
{\bIndMH(\bvocB',\hupsigma^*(\txi^{-1}\tOmega^1_{\fX'/\fX}))}}
\end{equation}
où $\rI\ulambda^{t,r}$ est le foncteur \eqref{crindmd2f}, et les foncteurs non libellés sont les foncteurs canoniques, est commutatif.

On désigne par $\rI\hupsigma'^{(t,r)*}$ le foncteur composé $\rI\lambda^{t,r}\circ \rI\hupsigma'^{(t)*}$ \eqref{indmdlb1f}. 
Compte tenu de \ref{indsh380} et \ref{ccoh50}, celui-ci est défini par 
\begin{eqnarray}\label{crindmd15e}
\rI\hupsigma'^{(t,r)*}\colon
\bIndMH(\co_{\fX'},\txi^{-1}\tOmega^1_{\fX'/\cS})\rightarrow \bIndMH(\bvocB',\hupsigma'^*(\txi^{-1}\tOmega^1_{\fX'/\cS})),\\
(\cN,\theta)\mapsto(\bvcC'^{(t,r)}\otimes_{\bvocB'}\rI\hupsigma'^*(\cN),\delta_{\bvcC'^{(t,r)}} \otimes_{\bvocB'}\id+\id \otimes_{\bvocB'} \rI\hupsigma'^*(\theta)).\nonumber
\end{eqnarray}
Il résulte de \eqref{crindmd2o} et \eqref{crindmd15d} que le diagramme 
\begin{equation}\label{crindmd15f}
\xymatrix{
{\bIndMH(\co_{\fX'},\txi^{-1}\tOmega^1_{\fX'/\cS})}\ar[r]^-(0.5){\rI\hupsigma'^{(t,r)*}}\ar[d]&
{\bIndMH(\bvocB',\hupsigma^*(\txi^{-1}\tOmega^1_{\fX'/\cS}))}\ar[dd]\\
{\bIndMH(\co_{\fX'},\txi^{-1}\tOmega^1_{\fX'/\fX})}\ar[d]_-(0.5){\rI\uhupsigma'^{(t,r)*}}&\\
{\bIndMC(\bvcC'^{(t,r)}/\bvuptheta^*(\bvcC^{(r)}))}\ar[r]&{\bIndMH(\bvocB',\hupsigma^*(\txi^{-1}\tOmega^1_{\fX'/\fX}))}}
\end{equation}
où $\rI\uhupsigma'^{(t,r)*}$ est le foncteur \eqref{crindmd2n} et les foncteurs non libellés sont les foncteurs canoniques, est commutatif. 

On désigne par $\rI\fS'^{(t,r)}$ le foncteur composé $\rI\lambda^{t,r}\circ \rI\fS'^{(t)}$ \eqref{indmdlb1b}. 
Compte tenu de \ref{indsh380} et \ref{ccoh50}, celui-ci est défini par 
\begin{eqnarray}\label{crindmd15g}
\rI\fS'^{(t,r)}\colon
\bIndMod(\bvocB')\rightarrow \bIndMH(\bvocB',\hupsigma^*(\txi^{-1}\tOmega^1_{\fX'/\cS})),\\
\cM\mapsto(\bvcC'^{(t,r)}\otimes_{\bvocB'}\cM,\delta_{\bvcC'^{(t,r)}} \otimes_{\bvocB'}\id).\nonumber
\end{eqnarray}

\begin{lem}\label{crindmd4}
Soient  $\cM$ un ind-$\bvocB'$-module rationnel \eqref{ahttf49} et plat \eqref{indsh46}, $r$ un nombre rationnel $\geq 0$, 
$q$ un entier $\geq 0$. On a alors un isomorphisme canonique fonctoriel de ind-$\bvlgg^*(\bvcC^{(r)})$-modules \eqref{hmdf14a}
\begin{equation}\label{crindmd4a}
\rR^q\rI\bvtau_*(\cM\otimes_{\bvocB'}\bvuptheta^*(\bvcC^{(r)}))\stackrel{\sim}{\rightarrow} 
\underset{\underset{t\in \mQ_{>r}}{\longrightarrow}}{\mlq\mlq\lim \mrq\mrq}\
\rR^q\rI\bvtau_*(\cM\otimes_{\bvocB'}\umK^\bullet(\bvcC'^{(t,r)})),
\end{equation}
où $\rI\bvtau_*$ est le foncteur \eqref{ahttfg20h}
et $\umK^\bullet(\bvcC'^{(t,r)})$ est le complexe de Dolbeault du $\bvuptheta^*(\bvcC^{(r)})$-module de Higgs 
$(\bvcC'^{(t,r)},p^t\ud_{\bvcC'^{(t,r)}})$ \eqref{ccoh11}. 
\end{lem}

En effet, en vertu de \ref{ccoh17}, le morphisme canonique de complexes de ind-$\bvuptheta^*(\bvcC^{(r)})$-modules 
\begin{equation}\label{crindmd4b}
\bvuptheta^*(\bvcC^{(r)})_\mQ[0]\rightarrow \underset{\underset{t\in \mQ_{>r}}{\longrightarrow}}{\mlq\mlq\lim \mrq\mrq}\ \umK^\bullet_{\mQ}(\bvcC'^{(t,r)})
\end{equation}
est un quasi-isomorphisme. Par ailleurs, $\cM$ étant rationnel, pour tout $\bvocB'$-module $\cF$, le morphisme canonique 
$\cM\otimes_{\bvocB'}\cF\rightarrow \cM\otimes_{\bvocB'}\cF_\mQ$ \eqref{indsh5g} est un isomorphisme. Comme $\cM$ est plat, 
on en déduit que le morphisme canonique de complexes de ind-$\bvuptheta^*(\bvcC^{(r)})$-modules 
\begin{equation}\label{crindmd4c}
\cM\otimes_{\bvocB'}\bvuptheta^*(\bvcC^{(r)})[0]\rightarrow 
\underset{\underset{t\in \mQ_{>r}}{\longrightarrow}}{\mlq\mlq\lim \mrq\mrq}\ \cM\otimes_{\bvocB'}\umK^\bullet(\bvcC'^{(t,r)})
\end{equation}
est un quasi-isomorphisme. 
Comme $\rR^q\rI\bvtau_*$ commute aux petites limites inductives filtrantes \eqref{indsh9e}, on en déduit l'isomorphisme \eqref{crindmd4a}.

\subsection{}\label{crindmd16} 
On désigne par $\bIH(\bvocB^!,\huppi^*(\txi^{-1}\tOmega^1_{\fX'/\fX}))$ la catégorie des $\bvocB^!$-isogénies de Higgs à coefficients dans 
$\huppi^*(\txi^{-1}\tOmega^1_{\fX'/\fX})$ \eqref{indsh23} et par $\bIH_\mQ(\bvocB^!,\huppi^*(\txi^{-1}\tOmega^1_{\fX'/\fX}))$
la catégorie des objets de $\bIH(\bvocB^!,\huppi^*(\txi^{-1}\tOmega^1_{\fX'/\fX}))$ à isogénie près (\cite{agt} III.6.1.1). 
\`A tout objet $\cF$ de $\bIH_\mQ(\bvocB^!,\huppi^*(\txi^{-1}\tOmega^1_{\fX'/\fX}))$, on associe un complexe de Dolbeault 
$\umK(\cF)$ dans $\bMod_\mQ(\bvocB^!)$ \eqref{indsh23e}. 

Avec les notations de \ref{hmdf41}, le morphisme $\huppi$ \eqref{hmdf13g} induit un foncteur 
\begin{equation}\label{crindmd16a}
\huppi^*\colon  \bIH(\co_{\fX'},\txi^{-1}\tOmega^1_{\fX'/\fX})\rightarrow \bIH(\bvocB^!,\huppi^*(\txi^{-1}\tOmega^1_{\fX'/\fX})),
\end{equation}
et par suite un foncteur 
\begin{equation}\label{crindmd16b}
\huppi^*_\mQ\colon  \bIH_\mQ(\co_{\fX'},\txi^{-1}\tOmega^1_{\fX'/\fX})\rightarrow \bIH_\mQ(\bvocB^!,\huppi^*(\txi^{-1}\tOmega^1_{\fX'/\fX})).
\end{equation}
Compte tenu de \eqref{hmdf41a}, on en déduit un foncteur qu'on note encore 
\begin{equation}\label{crindmd16c}
\huppi^*_\mQ\colon  \bMH^\coh(\co_{\fX'}[\frac 1 p], \txi^{-1}\tOmega^1_{\fX'/\fX})\rightarrow \bIH_\mQ(\bvocB^!,\huppi^*(\txi^{-1}\tOmega^1_{\fX'/\fX})).
\end{equation}
Pour tout objet $\cN$ de $\bMH^\coh(\co_{\fX'}[\frac 1 p], \txi^{-1}\tOmega^1_{\fX'/\fX})$, 
on voit aussitôt que le complexe de Dolbeault $\umK^\bullet(\huppi^*_\mQ(\cN))$ de $\huppi^*_\mQ(\cN)$ se déduit de celui de $\cN$ 
en appliquant le foncteur $\huppi^*_\mQ$ \eqref{hmdf40k} terme à terme.

\begin{lem}\label{crindmd5}
Soient $\cN$ un $\co_{\fX'}[\frac 1 p]$-fibré de Higgs à coefficients dans $\txi^{-1}\tOmega^1_{\fX'/\fX}$ \eqref{hmdf42}, 
$r$ un nombre rationnel $\geq 0$. 
On désigne par $\umK^\bullet(\cN)$ le complexe de Dolbeault de $\cN$ \eqref{MH2a}, que l'on considère comme un complexe de ind-$\co_{\fX'}$-modules \eqref{hmdf40a},
par $\rI\huppi^*(\umK^\bullet(\cN))$ son image par le foncteur $\rI\huppi^*$ \eqref{chb50c}, et 
pour tout nombre rationnel $t\geq r$, par $\rI\uhupsigma'^{(t,r)*}(\cN)$ l'objet de $\bIndMC(\bvcC'^{(t,r)}/\bvuptheta^*(\bvcC^{(r)}))$ associé à $\cN$ \eqref{crindmd2n}
et par $\umK^\bullet(\rI\uhupsigma'^{(t,r)*}(\cN))$ son complexe de Dolbeault \eqref{crindmd2}. 
Alors, pour tout nombre rationnel $t>r$, il existe un morphisme canonique de 
$\bD^+(\bIndMod(\bvlgg^*(\bvcC^{(r)})))$ 
\begin{equation}\label{crindmd5b}
\rI\huppi^*(\umK^\bullet(\cN))\otimes_{\bvocB^!}\bvlgg^*(\bvcC^{(r)})\rightarrow 
\rR\rI \bvtau_*(\umK^\bullet(\rI\uhupsigma'^{(t,r)*}(\cN))),
\end{equation}
où le produit tensoriel à gauche est défini terme à terme (non dérivé) et 
$\rI \bvtau_*$ est le foncteur \eqref{ahttfg20h}. C'est en fait un morphisme de systèmes inductifs indexés par $t\in \mQ_{>r}$, 
où les morphismes de transition du but sont 
induits par les homomorphismes $\bvalpha'^{t',t,r,r}$ \eqref{ahttfg20e} pour $t'\geq t >r$. 
De plus, le morphisme induit de $\bD^+(\bIndMod(\bvlgg^*(\bvcC^{(r)})))$ 
\begin{equation}\label{crindmd5ab}
\rI\huppi^*(\umK^\bullet(\cN))\otimes_{\bvocB^!}\bvlgg^*(\bvcC^{(r)})\rightarrow 
\rR\rI \bvtau_*(\underset{\underset{t\in \mQ_{>r}}{\longrightarrow}}{\mlq\mlq\lim \mrq\mrq}\ \umK^\bullet(\rI\uhupsigma'^{(t,r)*}(\cN)))
\end{equation}
est un isomorphisme. En particulier, pour tout entier $q\geq 0$, le morphisme induit 
\begin{equation}\label{crindmd5a}
\rH^q(\rI\huppi^*(\umK^\bullet(\cN))\otimes_{\bvocB^!}\bvlgg^*(\bvcC^{(r)}))
\rightarrow \underset{\underset{t\in \mQ_{>r}}{\longrightarrow}}{\mlq\mlq\lim \mrq\mrq}\
\rR^q\rI \bvtau_*(\umK^\bullet(\rI\uhupsigma'^{(t,r)*}(\cN)))
\end{equation}
est un isomorphisme. 
\end{lem}

Soit $t$ un nombre rationnel $>r$. 
On désigne par $\uhupsigma'^{(t,r)*}_\mQ(\cN)$ l'objet de $\bIMC_\mQ(\bvcC'^{(t,r)}/\bvuptheta^*(\bvcC^{(r)}))$ associé à $\cN$ \eqref{crindmd3f},
et par $\umK^\bullet(\uhupsigma'^{(t,r)*}_\mQ(\cN))$ son complexe de Dolbeault \eqref{crindmd3}. 
Comme $\ud_{\bvcC'^{(t,r)}}$ est une $\bvuptheta^*(\bvcC^{(r)})$-dérivation \eqref{ccoh11c}, 
on a un morphisme canonique de complexes de  $\bvlgg^*(\bvcC^{(r)})_\mQ$-modules 
\begin{equation}\label{crindmd5i}
\huppi^*_\mQ(\umK^\bullet(\cN))\otimes_{\bvocB^!}\bvlgg^*(\bvcC^{(r)})\rightarrow 
\bvtau_{\mQ*}(\umK^\bullet(\uhupsigma'^{(t,r)*}_\mQ(\cN))).
\end{equation}
D'après (\cite{sp} \href{https://stacks.math.columbia.edu/tag/013K}{013K}), 
comme la catégorie abélienne $\bMod_{\mQ}(\bvuptheta^*(\bvcC^{(r)}))$ a assez d'injectifs \eqref{indsh14},  
il existe un complexe borné inférieurement de $\bvuptheta^*(\bvcC^{(r)})_\mQ$-modules injectifs $\cL^\bullet$
et un quasi-isomorphisme $u\colon \umK^\bullet(\uhupsigma'^{(t,r)*}_\mQ(\cN))\rightarrow \cL^\bullet$. 
Ce dernier induit un morphisme de $\bD^+(\bMod_\mQ(\bvlgg^*(\bvcC^{(r)})))$ 
\begin{equation}\label{crindmd5j}
\bvtau_{\mQ*}(\umK^\bullet(\uhupsigma'^{(t,r)*}_\mQ(\cN)))\rightarrow 
\rR \bvtau_{\mQ*}(\umK^\bullet(\uhupsigma'^{(t,r)*}_\mQ(\cN))),
\end{equation}
qui ne depend que de $\umK^\bullet(\uhupsigma'^{(t,r)*}_\mQ(\cN))$, mais pas de $u$
(\cite{sp} \href{https://stacks.math.columbia.edu/tag/05TG}{05TG}).

Compte tenu de \eqref{indsh13a}, \eqref{indsh14g} et \eqref{indsh20d}, 
on prend pour morphisme \eqref{crindmd5b} l'image par le foncteur $\upalpha_{\bvlgg^*(\bvcC^{(r)})}$ \eqref{indsh20c} 
du morphisme composé de \eqref{crindmd5i} et \eqref{crindmd5j}. 
C'est clairement un morphisme de systèmes inductifs indexés par $t\in \mQ_{>r}$, 
où les morphismes de transition du but sont 
induits par les homomorphismes $\bvalpha'^{t',t,r,r}$ \eqref{ahttfg20e} pour $t'\geq t >r$.

Compte tenu de \ref{chb2}(i), pour tout entier $j\geq 0$, $\ud_{\bvcC'^{(t,r)}}$ \eqref{ccoh11c} induit un morphisme $\bvlgg^*(\bvcC^{(r)})$-linéaire
\begin{equation}\label{crindmd5c}
\theta^{j,(t,r)}\colon \rR^j\bvtau_*(\bvcC'^{(t,r)})\rightarrow \huppi^*(\txi^{-1}\tOmega^1_{\fX'/\fX})\otimes_{\bvocB^!}\rR^j\bvtau_*(\bvcC'^{(t,r)}),
\end{equation}
qui est clairement un $\bvlgg^*(\bvcC^{(r)})$-champ de Higgs sur $\rR^j\bvtau_*(\bvcC'^{(t,r)})$ à coefficients dans 
$\huppi^*(\txi'^{-1}\tOmega^1_{\fX'/\fX})\otimes_{\bvocB^!}\bvlgg^*(\bvcC^{(r)})$. 
Soit $\huppi^*_\mQ(\cN)$ l'objet de $\bIH_\mQ(\bvocB^!,\huppi^*(\txi^{-1}\tOmega^1_{\fX'/\fX}))$ associé à $\cN$ \eqref{crindmd16c}. 
On peut alors considérer $\huppi^*_\mQ(\cN)\otimes_{\bvocB^!_\mQ}\rR^j\bvtau_*(\bvcC'^{(t,r)})_\mQ$, le produit tensoriel \eqref{indsh23h}   
de $\huppi^*_\mQ(\cN)$ et de l'isogénie de Higgs 
\begin{equation}
(\rR^j\bvtau_*(\bvcC'^{(t,r)}),\rR^j\bvtau_*(\bvcC'^{(t,r)}),\id,p^t\theta^{j,(t,r)}),
\end{equation} 
qui est naturellement une $\bvlgg^*(\bvcC^{(r)})$-isogénie de Higgs à coefficients dans 
$\huppi^*(\txi^{-1}\tOmega^1_{\fX'/\fX})\otimes_{\bvocB^!}\bvlgg^*(\bvcC^{(r)})$, définie à isogénie près (cf. \ref{indsh23}).
On note $\umK^\bullet(\huppi^*_\mQ(\cN)\otimes_{\bvocB^!_\mQ}\rR^j\bvtau_*(\bvcC'^{(t,r)})_\mQ)$ son complexe de Dolbeault,
qui est un complexe de $\bMod_{\mQ}(\bvlgg^*(\bvcC^{(r)}))$ \eqref{indsh23e}. 

D'après \ref{chb2}(ii), pour tout entier $i\geq 0$, on a un isomorphisme canonique de $\bvlgg^*(\bvcC^{(r)})_\mQ$-modules
\begin{equation}\label{crindmd5d}
\rR^j\bvtau_{\mQ*}(\umK^i(\uhupsigma'^{(t,r)*}_\mQ(\cN)))\stackrel{\sim}{\rightarrow} 
\umK^i(\huppi^*_\mQ(\cN)\otimes_{\bvocB^!_\mQ}\rR^j\bvtau_{\mQ*}(\bvcC'^{(t,r)}_\mQ)),
\end{equation} 
compatible avec les morphismes induits par les différentielles des deux complexes de Dolbeault.

La catégorie abélienne $\bMod_{\mQ}(\bvuptheta^*(\bvcC^{(r)}))$ ayant assez d'injectifs, 
on a une suite spectrale canonique, fonctorielle en $t$,
\begin{equation}\label{crindmd5e}
{^t\rE}_1^{i,j}=\rR^j\bvtau_{\mQ*}(\umK^i(\uhupsigma'^{(t,r)*}_\mQ(\cN)))\Rightarrow 
\rR^{i+j}\bvtau_{\mQ*}(\umK^\bullet(\uhupsigma'^{(t,r)*}_\mQ(\cN))).
\end{equation}
Celle-ci induit, pour tout entier $q\geq 0$, un morphisme canonique
\begin{equation}\label{crindmd5k}
\rH^q({^t\rE}_1^{\bullet,0})\rightarrow 
\rR^q\bvtau_{\mQ*}(\umK^\bullet(\uhupsigma'^{(t,r)*}_\mQ(\cN))),
\end{equation}
qui n'est autre que le morphisme induit par \eqref{crindmd5j} d'après (\cite{ega3} 0.11.3.4).

Compte tenu de \eqref{indsh13a}, \eqref{indsh14g} et \eqref{indsh20d}, 
appliquant le foncteur exact $\upalpha_{\bvlgg^*(\bvcC^{(r)})}$ à la suite spectrale \eqref{crindmd5e}, 
on obtient une suite spectrale de ind-$\bvlgg^*(\bvcC^{(r)})$-modules
\begin{equation}\label{crindmd5f}
{^t\cE}_1^{i,j}=\rR^j\rI\bvtau_*(\umK^i(\rI\uhupsigma'^{(t,r)*}(\cN)))\Rightarrow 
\rR^{i+j}\rI\bvtau_*(\umK^\bullet(\rI\uhupsigma'^{(t,r)*}(\cN))).
\end{equation}

D'après \ref{ccoh14} et \eqref{crindmd5d}, pour tout $i\geq 0$, on a un isomorphisme canonique 
\begin{equation}\label{crindmd5g}
\underset{\underset{t\in \mQ_{>r}}{\longrightarrow}}{\mlq\mlq\lim \mrq\mrq}\ {^t\cE}_1^{i,0}\stackrel{\sim}{\rightarrow}
\rI\huppi^*(\umK^i(\cN))\otimes_{\bvocB^!}\bvlgg^*(\bvcC^{(r)}),
\end{equation}
et pour tout $j\geq 1$, on a 
\begin{equation}\label{crindmd5h}
\underset{\underset{t\in \mQ_{>r}}{\longrightarrow}}{\mlq\mlq\lim \mrq\mrq}\ {^t\cE}_1^{i,j}=0.
\end{equation}
De plus, les isomorphismes \eqref{crindmd5g}  (pour $i\in \mN$) forment un isomorphisme de complexes.
Comme les petites limites inductives filtrantes existent et sont exactes dans $\bIndMod(\bvlgg^*(\bvcC^{(r)}))$ \eqref{indsh6e}, 
on en déduit que les morphismes \eqref{crindmd5a} sont des isomorphismes.
Par suite, le morphisme \eqref{crindmd5ab} est un isomorphisme puisque 
les foncteurs $\rR^q\rI \bvtau_*$ commutent aux petites limites inductives filtrantes \eqref{indsh9e}.

\begin{prop}\label{crindmd6}
Soient $t, r$ deux nombres rationnels tels que $t>r\geq 0$, $\cM$ un ind-$\bvocB'$-module rationnel \eqref{ahttf49} et plat \eqref{indsh46}, 
$\cN$ un $\co_{\fX'}[\frac 1 p]$-fibré de Higgs à coefficients dans $\txi^{-1}\tOmega^1_{\fX'/\fX}$, 
\begin{equation}
\ualpha\colon \rI\ufS'^{(t,r)}(\cM) \stackrel{\sim}{\rightarrow}\rI\uhupsigma'^{(t,r)*}(\cN)
\end{equation}
un isomorphisme de $\bIndMC(\bvcC'^{(t,r)}/\bvuptheta^*(\bvcC^{(r)}))$, où $\rI\ufS'^{(t,r)}$ est le foncteur \eqref{crindmd2k}
et $\rI\uhupsigma'^{(t,r)*}$ est le foncteur \eqref{crindmd2n}. 
On désigne par $\umK^\bullet(\cN)$ le complexe de Dolbeault de $\cN$ \eqref{MH2a}
et par $\rI\huppi^*(\umK^\bullet(\cN))$ son image par le foncteur $\rI\huppi^*$ \eqref{chb50c}. 
On a alors un isomorphisme canonique de $\bD^+(\bIndMod(\bvlgg^*(\bvcC^{(r)})))$
\begin{equation}\label{crindmd6a}
\rR\rI\bvtau_*(\cM\otimes_{\bvocB'}\bvuptheta^*(\bvcC^{(r)}))\rightarrow
\rI\huppi^*(\umK^\bullet(\cN))\otimes_{\bvocB^!}\bvlgg^*(\bvcC^{(r)}),
\end{equation}
où $\rI \bvtau_*$ est le foncteur \eqref{ahttfg20h} et le produit tensoriel à droite est défini terme à terme (non dérivé). 
\end{prop}

En effet, pour tout nombre rationnel $t'$ tel que $t\geq t' \geq r$, on désigne par
\begin{equation}\label{crindmd6c}
\ualpha^{(t',r)}\colon \rI\ufS'^{(t',r)}(\cM)\stackrel{\sim}{\rightarrow} \rI\uhupsigma'^{(t',r)*}(\cN)
\end{equation}
l'isomorphisme de $\bIndMC(\bvcC'^{(t',r)}/\bvuptheta^*(\bvcC^{(r)}))$ induit par $\rI\uvarepsilon'^{t,t',r}(\ualpha)$ \eqref{crindmd2g} et 
les isomorphismes \eqref{crindmd2p} et \eqref{crindmd2q}. 
Celui-ci induit un isomorphisme de complexes de ind-$\bvuptheta^*(\bvcC^{(r)})$-modules 
\begin{equation}\label{crindmd6d}
\cM\otimes_{\bvocB'}\umK^\bullet(\bvcC'^{(t',r)})\stackrel{\sim}{\rightarrow} \umK^\bullet(\rI\uhupsigma'^{(t',r)*}(\cN)),
\end{equation}
où $\umK^\bullet(\bvcC'^{(t',r)})$ est le complexe de Dolbeault du $\bvuptheta^*(\bvcC^{(r)})$-module de Higgs 
$(\bvcC'^{(t',r)},p^{t'}\ud_{\bvcC'^{(t',r)}})$ \eqref{ccoh11}.
Ces isomorphismes forment un isomorphisme de systèmes inductifs (pour $r<t'< t$) et ils induisent donc un isomorphisme 
de complexes de ind-$\bvuptheta^*(\bvcC^{(r)})$-modules 
\begin{equation}\label{crindmd6e}
\underset{\underset{t'\in \mQ_{>r}}{\longrightarrow}}{\mlq\mlq\lim \mrq\mrq}\ \cM\otimes_{\bvocB'}\umK^\bullet(\bvcC'^{(t',r)})\stackrel{\sim}{\rightarrow}
\underset{\underset{t'\in \mQ_{>r}}{\longrightarrow}}{\mlq\mlq\lim \mrq\mrq}\ \umK^\bullet(\rI\uhupsigma'^{(t',r)*}(\cN)).
\end{equation}
La proposition s'ensuit compte tenu de \eqref{crindmd4c} et \ref{crindmd5}.

\begin{prop}\label{crindmd8}
Soient $\cM$ un ind-$\bvocB'$-module de Dolbeault \eqref{crindmd1}, 
$\ucH'(\cM)$ le $\co_{\fX'}[\frac 1 p]$-fibré de Higgs à coefficients dans $\txi^{-1}\tOmega^1_{\fX'/\fX}$ associé \eqref{crindmd1d}, 
$\umK^\bullet(\ucH'(\cM))$ son complexe de Dolbeault, $\rI\huppi^*(\umK^\bullet(\ucH'(\cM)))$ l'image de ce dernier par le foncteur $\rI\huppi^*$ \eqref{chb50c}. 
Alors, il existe un nombre rationnel $r>0$ et un isomorphisme de $\bD^+(\bIndMod(\bvlgg^*(\bvcC^{(r)})))$
\begin{equation}\label{crindmd8a}
\rR\rI\bvtau_*(\cM\otimes_{\bvocB'}\bvuptheta^*(\bvcC^{(r)}))\rightarrow
\rI\huppi^*(\umK^\bullet(\ucH'(\cM)))\otimes_{\bvocB^!}\bvlgg^*(\bvcC^{(r)}),
\end{equation}
où $\rI \bvtau_*$ est le foncteur \eqref{ahttfg20h} et le produit tensoriel à droite est défini terme à terme (non dérivé).
\end{prop}

En effet, il existe un nombre rationnel $t>0$ et un isomorphisme de $\bIndMC(\bvcC'^{(t)}/\bvocB')$, 
\begin{equation}
\alpha\colon \rI\fS'^{(t)}(\cM)\stackrel{\sim}{\rightarrow} \rI\hupsigma'^{(t)*}(\cH'(\cM)),
\end{equation}
où $\cH'(\cM)$ est le $\co_{\fX'}[\frac 1 p]$-fibré de Higgs à coefficients dans $\txi^{-1}\tOmega^1_{\fX'/\cS}$ associé à $\cM$ \eqref{crindmd1a}. 
Pour tout nombre rationnel $0< r< t$,  l'isomorphisme $\rI\ulambda^{t,r}(\alpha)$ \eqref{crindmd2f} induit,
compte tenu de l'isomorphisme \eqref{crindmd2l} et du diagramme commutatif \eqref{crindmd2o}, 
un isomorphisme de $\bIndMC(\bvcC'^{(t,r)}/\bvuptheta^*(\bvcC^{(r)}))$,
\begin{equation}
\ualpha\colon \rI\ufS'^{(t,r)}(\cM) \stackrel{\sim}{\rightarrow} \rI\uhupsigma'^{(t,r)*}(\ucH'(\cM)). 
\end{equation}
Comme le ind-$\bvocB'$-module $\cM$ est rationnel et plat d'après \ref{indmdlb6},
la proposition résulte alors de \ref{crindmd6}. 

\begin{rema}\label{crindmd9}
Sous les hypothèses de \ref{crindmd8}, il existe un nombre rationnel $r>0$ tel que 
l'isomorphisme \eqref{crindmd8a} soit canonique et qu'il dépende fonctoriellement de $\cM$.
En effet, en vertu de \ref{indmdlb13}, il existe un nombre rationnel $t>0$ et un isomorphisme de $\bIndMC(\bvcC'^{(t)}/\bvocB')$, 
\begin{equation}
\alpha\colon \rI\fS'^{(t)}(\cM)\stackrel{\sim}{\rightarrow} \rI\hupsigma'^{(t)*}(\cH'(\cM)),
\end{equation}
vérifiant les propriétés \ref{indmdlb13}(i)-(ii), 
où $\cH'(\cM)$ est le $\co_{\fX'}[\frac 1 p]$-fibré de Higgs à coefficients dans $\txi^{-1}\tOmega^1_{\fX'/\cS}$ associé à $\cM$ \eqref{crindmd1a}. 
Pour tout nombre rationnel $t'$ tel que $0\leq t'< t$, on désigne par
\begin{equation}
\alpha^{(t')}\colon \rI\fS'^{(t')}(\cM) \stackrel{\sim}{\rightarrow} \rI\hupsigma'^{(t')*}(\cH'(\cM))
\end{equation}
l'isomorphisme de $\bIndMC(\bvcC'^{(t')}/\bvocB')$ induit par $\rI\varepsilon'^{t,t'}(\alpha)$ et 
l'analogue des isomorphismes \eqref{indmdlb3b} et \eqref{indmdlb3c}. D'après la preuve de \ref{indmdlb20}, 
$\alpha^{(t')}$ ne dépend que de $\cM$ (mais pas de $\alpha$) et il en dépend fonctoriellement. 
Pour tous nombres rationnels $t',r$ tel que $0\leq r<t'< t$, on désigne par
\begin{equation}
\ualpha^{(t',r)}\colon \rI\ufS'^{(t',r)}(\cM) \stackrel{\sim}{\rightarrow} \rI\uhupsigma'^{(t',r)*}(\ucH'(\cM))
\end{equation}
l'isomorphisme de $\bIndMC(\bvcC'^{(t',r)}/\bvuptheta^*(\bvcC^{(r)}))$ induit par $\rI\ulambda^{t',r}(\alpha^{(t')})$,
l'isomorphisme \eqref{crindmd2l} et le diagramme commutatif \eqref{crindmd2o}. Fixant $r$ tel que $0< r< t$, 
les isomorphismes $\ualpha^{(t',r)}$ forment un isomorphisme de systèmes inductifs, pour $r<t'< t$. 
L'assertion s'ensuit compte tenu de \ref{crindmd4} et \ref{crindmd5} (cf. la preuve de \ref{crindmd6}). 
\end{rema}

\begin{lem}\label{crindmd10}
Soient $r$ un nombre rationnel $\geq 0$, 
$n$, $q$ deux entiers $\geq 0$, $\cM$ un $\ocB'_n$-module de $\tE'_s$, $\cN$ un $\ocB^!_n$-module de $\tG_s$. 
Alors, les morphismes canoniques 
\begin{eqnarray}
\rR^q\tau_{n*}(\cM)\otimes_{\ocB^!_n}\lgg_n^*(\cC^{(r)}_n)&\rightarrow& \rR^q\tau_{n*}(\cM\otimes_{\ocB'_n}\uptheta_n^*(\cC^{(r)}_n)),\label{crindmd10a}\\
\rR^q\lgg_{n*}(\cN)\otimes_{\ocB_n}\cC^{(r)}_n&\rightarrow& \rR^q\lgg_{n*}(\cN\otimes_{\ocB^!_n}\lgg_n^*(\cC^{(r)}_n)),\label{crindmd10b}\\
\rR^q\uptheta_{n*}(\cM)\otimes_{\ocB_n}\cC^{(r)}_n&\rightarrow& \rR^q\uptheta_{n*}(\cM\otimes_{\ocB'_n}\uptheta_n^*(\cC^{(r)}_n)),\label{crindmd10c}
\end{eqnarray}
où les morphismes de topos annelés $\tau_n$, $\lgg_n$ et $\uptheta_n$ sont définis dans \eqref{hmdf19l} et \eqref{hmdf19k}, sont des isomorphismes. 
\end{lem}

En effet, comme le $\ocB_n$-module $\cF^{(r)}_n$ est localement libre de type fini \eqref{ahttfg3e}, pour tout entier $m\geq 0$, les morphismes canoniques 
\begin{eqnarray}
\rR^q\tau_{n*}(\cM)\otimes_{\ocB^!_n}\lgg_n^*(\rS^m_{\ocB_n}(\cF^{(r)}_n))&\rightarrow& \rR^q\tau_{n*}(\cM\otimes_{\ocB'_n}\uptheta_n^*(\rS^m_{\ocB_n}(\cF^{(r)}_n))),\\
\rR^q\lgg_{n*}(\cN)\otimes_{\ocB_n}\rS^m_{\ocB_n}(\cF^{(r)}_n)&\rightarrow& \rR^q\lgg_{n*}(\cN\otimes_{\ocB^!_n}\lgg_n^*(\rS^m_{\ocB_n}(\cF^{(r)}_n))),\\
\rR^q\uptheta_{n*}(\cM)\otimes_{\ocB_n}\rS^m_{\ocB_n}(\cF^{(r)}_n)&\rightarrow& \rR^q\uptheta_{n*}(\cM\otimes_{\ocB'_n}\uptheta_n^*(\rS^m_{\ocB_n}(\cF^{(r)}_n))),
\end{eqnarray}
sont des isomorphismes \eqref{notconv9}. Par ailleurs, les topos $\tE$ et $\tG$ sont cohérents d'après (\cite{agt} VI.10.5(ii)) et (\cite{ag} 3.4.22(ii)). 
Le morphisme de schémas $g\colon X'\rightarrow X$ étant cohérent, les morphismes de topos $\tau$ et $\lgg$ \eqref{hmdf17d} sont cohérents 
en vertu de (\cite{agt} VI.10.4 et VI.10.5(i)), (\cite{ag} 3.4.21(i) et 3.4.22(ii)) et (\cite{sga4} VI 3.3). 
Par suite, les foncteurs $\rR^q\tau_*$, $\rR^q\lgg_*$ et $\rR^q\Theta_*$ commutent aux limites inductives filtrantes de faisceaux abéliens 
en vertu de (\cite{sga4} VI 5.1). La proposition s'ensuit compte tenu de \eqref{ahttfg3f}, \eqref{hmdf18k} et \eqref{hmdf18hh}.

\begin{lem}\label{crindmd12}
Soient $r$ un nombre rationnel $\geq 0$, 
$q$ un entier $\geq 0$, $\cM$ un ind-$\bvocB'$-module de $\tE'^{\mN^\circ}$, $\cN$ un ind-$\bvocB^!$-module de $\tG^{\mN^\circ}$. 
On a alors des isomorphismes canoniques 
\begin{eqnarray}
\rR^q\rI\bvtau_*(\cM)\otimes_{\bvocB^!}\bvlgg^*(\bvcC^{(r)})&\stackrel{\sim}{\rightarrow}& \rR^q\rI\bvtau_*(\cM\otimes_{\bvocB'}\bvuptheta^*(\bvcC^{(r)})),\label{crindmd12a}\\
\rR^q\rI\bvlgg_*(\cN)\otimes_{\bvocB}\bvcC^{(r)}&\stackrel{\sim}{\rightarrow}& \rR^q\rI\bvlgg_*(\cN\otimes_{\bvocB^!}\bvlgg^*(\bvcC^{(r)})),\label{crindmd12b}\\
\rR^q\rI\bvuptheta_*(\cM)\otimes_{\bvocB}\bvcC^{(r)}&\stackrel{\sim}{\rightarrow}& \rR^q\rI\bvuptheta_*(\cM\otimes_{\bvocB'}\bvuptheta^*(\bvcC^{(r)})),\label{crindmd12c}
\end{eqnarray}
où les morphismes de topos annelés $\bvtau$, $\bvlgg$ et $\bvuptheta$ sont définis dans \eqref{hmdf14a}. 
\end{lem}

En effet, le cas où $\cM$ est un $\bvocB'$-module de $\tE'^{\mN^\circ}$ et $\cN$ est un $\bvocB^!$-module de $\tG^{\mN^\circ}$, 
résulte de \eqref{indsh21d}, \ref{crindmd10} et (\cite{agt} III.7.5). 
Il implique le cas général compte tenu de \ref{indsh18} et \eqref{indsh21g}.

\begin{prop}\label{crindmd13}
Supposons $g\colon X'\rightarrow X$ propre. 
Soient $\cM$ un ind-$\bvocB'$-module de Dolbeault \eqref{crindmd1}, $q$ un entier $q\geq 0$. 
On désigne par $\ucH'(\cM)$ le $\co_{\fX'}[\frac 1 p]$-fibré de Higgs à coefficients dans $\txi^{-1}\tOmega^1_{\fX'/\fX}$ associé à $\cM$ \eqref{crindmd1d}, 
et par $\umK^\bullet(\ucH'(\cM))$ son complexe de Dolbeault. 
Alors, le $\co_\fX[\frac 1 p]$-module $\rR^q\fgg_*(\umK^\bullet(\ucH'(\cM)))$ est cohérent, 
et il existe un nombre rationnel $r>0$, indépendant de $q$, et un isomorphisme de ind-$\bvcC^{(r)}$-modules
\begin{equation}\label{crindmd13a}
\rR^q\rI\bvuptheta_*(\cM)\otimes_{\bvocB}\bvcC^{(r)}\stackrel{\sim}{\rightarrow}
\rI\hupsigma^*(\rR^q\fgg_*(\umK^\bullet(\ucH'(\cM))))\otimes_{\bvocB}\bvcC^{(r)},
\end{equation}
où $\rI \hupsigma^*$ est le foncteur \eqref{hmdf40ii}.
\end{prop}

Cela résulte de \ref{crindmd8}, \ref{crindmd12}, \ref{chb22}, \eqref{indsh14h} et \eqref{indsh22f}.

\subsection{}\label{chb26}
Soit $\cK^\bullet$ un complexe borné inférieurement de $\bMod_\mQ(\bvocB^!)$. 
On note $\bvtau^*_\mQ(\cK^\bullet)$ le complexe image inverse défini terme à terme \eqref{hmdf13i}. 
On a alors un morphisme canonique de complexes de  $\bvocB^!_\mQ$-modules 
\begin{equation}\label{chb26a}
\cK^\bullet\rightarrow \bvtau_{\mQ*}(\bvtau^*_\mQ(\cK^\bullet)).
\end{equation}
D'après (\cite{sp} \href{https://stacks.math.columbia.edu/tag/013K}{013K}), 
comme la catégorie abélienne $\bMod_{\mQ}(\bvocB')$ a assez d'injectifs \eqref{indsh14},  
il existe un complexe borné inférieurement de $\bvocB'_\mQ$-modules injectifs $\cL^\bullet$
et un quasi-isomorphisme $u\colon \bvtau^*_\mQ(\cK^\bullet)\rightarrow \cL^\bullet$. 
Ce dernier induit un morphisme de $\bD^+(\bMod_\mQ(\bvocB^!))$ 
\begin{equation}\label{chb26b}
\bvtau_{\mQ*}(\bvtau^*_\mQ(\cK^\bullet))\rightarrow \rR \bvtau_{\mQ*}(\bvtau^*_\mQ(\cK^\bullet)).
\end{equation}
D'après (\cite{sp} \href{https://stacks.math.columbia.edu/tag/05TG}{05TG}), ce morphisme ne depend que de $\cK^\bullet$, mais pas de $u$,  
et il en dépend fonctoriellement.  
On obtient par composition un morphisme canonique de $\bD^+(\bMod_\mQ(\bvocB^!))$, fonctoriel en $\cK^\bullet$, 
\begin{equation}\label{chb26c}
\cK^\bullet\rightarrow \rR \bvtau_{\mQ*}(\bvtau^*_\mQ(\cK^\bullet)).
\end{equation}
D'après \eqref{indsh13a} et \eqref{indsh14g}, celui-ci induit un morphisme canonique de $\bD^+(\bIndMod(\bvocB^!))$, fonctoriel en $\cK^\bullet$, 
\begin{equation}\label{chb26d}
\cK^\bullet\rightarrow \rR \rI \bvtau_*(\rI\bvtau^*(\cK^\bullet)).
\end{equation}

On notera que les foncteurs $\bvtau^*_\mQ$ et $\rI\bvtau^*$
ne transforment pas a priori quasi-isomorphisme en quasi-isomorphisme. 
Les morphismes \eqref{chb26c} et \eqref{chb26d} ne peuvent donc pas être étendus à $\cK^\bullet$ 
objet de $\bD^+(\bMod_\mQ(\bvocB^!))$. On a toutefois l'énoncé \ref{crindmd140}.

On appliquera les constructions ci-dessus en particulier au cas où $\cK^\bullet=\huppi^*_\mQ(\mK^\bullet)\otimes_{\bvocB^!}\cE$, 
$\mK^\bullet$ étant un complexe borné inférieurement de $\bMod^\coh(\co_{\fX'}[\frac 1 p])$ et $\cE$  un $\bvocB^!$-module plat (cf. \ref{hmdf40}).

\subsection{}\label{crindmd14}
Soit $(\cN,\theta)$ un $\co_{\fX'}[\frac 1 p]$-fibré de Higgs à coefficients dans $\txi^{-1}\tOmega^1_{\fX'/\cS}$. 
Pour alléger les notations, on omettra $\theta$ lorsque cela n'induit aucun risque d'ambiguïté. 
On désigne par
\begin{equation}\label{crindmd14a}
\utheta\colon \cN\rightarrow \txi^{-1}\tOmega^1_{\fX'/\fX}\otimes_{\co_{\fX'}} \cN
\end{equation}
le $\co_{\fX'}[\frac 1 p]$-champ de Higgs induit par $\theta$ et 
par $\mK^\bullet(\cN)$ (resp.  $\umK^\bullet(\cN)$) le complexe de Dolbeault de $(\cN,\theta)$ (resp. $(\cN,\utheta)$) \eqref{MH2c}.
On munit le complexe $\mK^\bullet(\cN)$ de la filtration de Koszul associée à la suite exacte \eqref{hmdf410d}, cf. \eqref{MH90g}, et on note 
\begin{equation}\label{crindmd14b}
\partial\colon \umK^\bullet(\cN)\rightarrow \fgg^*(\txi^{-1}\tOmega^1_{\fX/\cS}) \otimes_{\co_{\fX'}} \umK^\bullet(\cN)
\end{equation}
le bord associé dans $\bD^+(\bMod(\co_{\fX'}[\frac 1 p]))$ \eqref{MH90j}. 

Compte tenu de \ref{chb5}, la filtration de Koszul de $\mK^\bullet(\cN)$ induit une filtration du complexe $\hupsigma'^*_\mQ(\mK^\bullet(\cN))$ \eqref{chb5b} 
et un bord associé 
\begin{equation}\label{crindmd14c}
\hupsigma'^\star_\mQ(\partial)\colon \hupsigma'^*_\mQ(\umK^\bullet(\cN))\rightarrow \hupsigma'^*_\mQ(\fgg^*(\txi^{-1}\tOmega^1_{\fX/\cS}) \otimes_{\co_{\fX'}} \umK^\bullet(\cN))
\end{equation}
dans $\bD^+(\bMod_\mQ(\bvocB'))$, défini de façon analogue à $\partial$. 
Nous utilisons la notation $\hupsigma'^\star_\mQ(\partial)$ pour la distinguer du foncteur $\hupsigma'^*_\mQ$ que nous n'étendons pas aux catégories dérivées envisagées.

On note
\begin{equation}\label{crindmd14d}
\rI\hupsigma'^\star(\partial)\colon \rI\hupsigma'^*(\umK^\bullet(\cN))\rightarrow \rI \hupsigma'^*(\fgg^*(\txi^{-1}\tOmega^1_{\fX/\cS}) \otimes_{\co_{\fX'}} \umK^\bullet(\cN))
\end{equation}
l'image canonique de $\hupsigma'^\star_\mQ(\partial)$ dans $\bD^+(\bIndMod(\bvocB'))$ \eqref{indsh20c}.
Ce morphisme n'est autre que le bord de la filtration de Koszul du complexe de Dolbeault $\rI\hupsigma'^*(\mK^\bullet(\cN))$ 
du ind-$\bvocB'$-module de Higgs $(\rI\hupsigma'^*(\cN),\rI\hupsigma'^*(\theta))$ \eqref{chb50a}, 
associée à l'image par le foncteur $\hupsigma'^*$ de la suite exacte \eqref{hmdf410d}, cf. \eqref{indsh60g}. 
Nous utilisons la notation $\rI\hupsigma'^\star(\partial)$ pour la distinguer du foncteur $\rI\hupsigma'^*$ 
que nous n'étendons pas aux catégories dérivées envisagées.

De même, la filtration de Koszul de $\mK^\bullet(\cN)$ induit une filtration du complexe $\huppi^*_\mQ(\mK^\bullet(\cN))$ \eqref{chb5c} 
et un bord associé 
\begin{equation}\label{crindmd14e}
\huppi^\star_\mQ(\partial)\colon \huppi^*_\mQ(\umK^\bullet(\cN))\rightarrow \huppi^*_\mQ(\fgg^*(\txi^{-1}\tOmega^1_{\fX/\cS}) \otimes_{\co_{\fX'}} \umK^\bullet(\cN))
\end{equation}
dans $\bD^+(\bMod_\mQ(\bvocB^!))$, défini de façon analogue à $\partial$. On note 
\begin{equation}\label{crindmd14f}
\rI\huppi^\star(\partial)\colon \rI\huppi^*(\umK^\bullet(\cN))\rightarrow \rI\huppi^*(\fgg^*(\txi^{-1}\tOmega^1_{\fX/\cS}) \otimes_{\co_{\fX'}} \umK^\bullet(\cN))
\end{equation}
son image canonique dans $\bD^+(\bIndMod(\bvocB^!))$.

Pour tout entier $q\geq 0$, on désigne par 
\begin{equation}\label{crindmd14g}
\kappa^q\colon \rR^q\fgg_*(\umK^\bullet(\cN))\rightarrow \txi^{-1}\tOmega^1_{\fX/\cS} \otimes_{\co_{\fX}} \rR^q\fgg_*(\umK^\bullet(\cN)) 
\end{equation}
le $\co_\fX$-champ de Katz-Oda de la cohomologie de Dolbeault de $\cN$ \eqref{MH96}, qui n'est autre que $\rR^q\fgg_*(\partial)$.

\begin{prop}\label{crindmd141}
Conservons les hypothèses et notations de \ref{crindmd14} supposons de plus que le morphisme $g\colon X'\rightarrow X$ soit propre 
et que le $\co_{\fX'}[\frac 1 p]$-module de Higgs $(\cN,\theta)$ soit nilpotent \eqref{definf22}.
Alors, pour tout entier $q\geq 0$, le $\co_\fX[\frac 1 p]$-module de Higgs $(\rR^q\fgg_*(\umK^\bullet(\cN)),\kappa^q)$, où $\kappa^q$ est le champ de Katz-Oda \eqref{crindmd14g}, 
est nilpotent.   
\end{prop}

La preuve est identique à celle de \ref{MH104} en utilisant \ref{definf231} au lieu de \ref{MH141}. On notera que le foncteur $\rR^q\fgg_*$ 
transforme les $\co_{\fX'}[\frac 1 p]$-modules cohérents en des $\co_\fX[\frac 1 p]$-modules cohérents d'après (\cite{egr1} 2.10.24 et 2.11.5).

\begin{lem}\label{crindmd140}
Conservons les hypothèses et notations de \ref{crindmd14}. Soit, de plus, $\cE$ un $\bvocB^!$-module plat.
Alors, le diagramme de morphismes de $\bD^+(\bIndMod(\bvocB^!))$
\begin{equation}
\text{\tiny {\xymatrix{
{\rI\huppi^*(\umK^\bullet(\cN))\otimes_{\bvocB^!}\cE}\ar[rr]^-(0.5){\rI\huppi^\star(\partial)\otimes\id}\ar[d]&&
{\rI\huppi^*(\fgg^*(\txi^{-1}\tOmega^1_{\fX/\cS}) \otimes_{\co_{\fX'}} \umK^\bullet(\cN))\otimes_{\bvocB^!}\cE}\ar[d]\\
{\rR\rI \bvtau_*(\rI\hupsigma'^*(\umK^\bullet(\cN))\otimes_{\bvocB'}\bvtau^*(\cE))}\ar[rr]^-(0.5){\rR\rI \bvtau_*(\rI\hupsigma'^\star(\partial)\otimes\id)}&&
{\rR\rI \bvtau_*(\rI \hupsigma'^*(\fgg^*(\txi^{-1}\tOmega^1_{\fX/\cS}) \otimes_{\co_{\fX'}} \umK^\bullet(\cN))\otimes_{\bvocB'}\bvtau^*(\cE))}}}}
\end{equation}
où les flèches verticales sont les morphismes \eqref{chb26d}, est commutatif.  
\end{lem}

Pour alléger les notations, on omettra $\cN$ des notations des complexes $\mK^\bullet(\cN)$ et $\umK^\bullet(\cN)$, 
et on notera leurs différentielles $\theta^\bullet$ et $\utheta^\bullet$. Posons $\Omega= \fgg^*(\txi^{-1}\tOmega^1_{\fX/\cS})$ et
\begin{equation}
\mG^\bullet=\mK^\bullet/\rW^2\mK^\bullet,
\end{equation}
où $\rW^\bullet \mK^\bullet$ désigne la filtration de Koszul de $\mK^\bullet$ associée à la suite exacte \eqref{hmdf410d}, cf. \eqref{MH90g}. 
On note encore $\theta^\bullet$ les différentielles de $\mG^\bullet$.
On a une suite exacte canonique de complexes de $\co_{\fX'}[\frac 1 p]$-modules \eqref{MH90i}
\begin{equation}
0\longrightarrow \Omega\otimes_{\co_{\fX'}}\umK^\bullet[-1]\stackrel{u^\bullet}{\longrightarrow} \mG^\bullet 
\stackrel{v^\bullet}{\longrightarrow} \umK^\bullet \longrightarrow 0. 
\end{equation}
On désigne par $\rC^\bullet$ le cône de $u^\bullet$. Pour tout entier $i$, on a donc $\rC^i=(\Omega\otimes_{\co_{\fX'}}\umK^i)\oplus \mG^i$ 
et la différentielle $c^i\colon \rC^i\rightarrow \rC^{i+1}$ est définie par la matrice 
\begin{equation}
\begin{pmatrix}
\id\otimes\utheta^i &0 \\
 u^{i+1} & \theta^i
\end{pmatrix}.
\end{equation}
Notons $\pi^\bullet_1$ et $\pi^\bullet_2$ les projections canoniques de $\rC^\bullet$ sur $\Omega\otimes_{\co_{\fX'}}\umK^\bullet$ et $\mG^\bullet$, 
respectivement. Le composé $v^\bullet \circ \pi^\bullet_2\colon \rC^\bullet\rightarrow \umK^\bullet$ est alors un quasi-isomorphisme,
et $-\partial$ \eqref{crindmd14b} est le composé  dans $\bD^+(\bMod(\co_{\fX'}[\frac 1 p]))$ de l'inverse de $v^\bullet \circ \pi^\bullet_2$ et de $\pi^\bullet_1$
(\cite{sp} \href{https://stacks.math.columbia.edu/tag/09KF}{09KF}). 

D'après \ref{chb50}, la suite de complexes de ind-$\bvocB'$-modules 
\begin{equation}
0\rightarrow \rI\hupsigma'^*(\Omega\otimes_{\co_{\fX'}}\umK^\bullet)[-1]\rightarrow \rI\hupsigma'^*(\mG^\bullet) 
\rightarrow \rI\hupsigma'^*(\umK^\bullet) \rightarrow 0
\end{equation}
est exacte. Le cône de $\rI\hupsigma'^*(u^\bullet)$ s'identifie à $\rI\hupsigma'^*(\rC^\bullet)$, et le morphisme $\rI\hupsigma'^\star(\partial)$ de 
$\bD^+(\bIndMod(\bvocB'))$ est alors défini de façon analogue à $\partial$. On décrit de même le morphisme $\rI\huppi^\star(\partial)$ de 
$\bD^+(\bIndMod(\bvocB^!))$. 

Comme le morphisme \eqref{chb26d} est fonctoriel, les diagrammes de $\bD^+(\bIndMod(\bvocB^!))$
\begin{equation}
\text{\tiny {\xymatrix{
{\rI\huppi^*(\umK^\bullet)\otimes_{\bvocB^!}\cE}\ar[d]&
{\rI\huppi^*(\rC^\bullet)\otimes_{\bvocB^!}\cE}\ar[r]^-(0.5){a_1}\ar[l]_-(0.5){a_2}
\ar[d] &{\rI\huppi^*(\Omega\otimes_{\co_{\fX'}}\umK^\bullet)\otimes_{\bvocB^!}\cE}\ar[d]\\
{\rR\rI\bvtau_*(\rI\hupsigma'^*(\umK^\bullet)\otimes_{\bvocB'}\bvtau^*(\cE))}&{\rR\rI\bvtau_*(\rI\hupsigma'^*(\rC^\bullet)\otimes_{\bvocB'}\bvtau^*(\cE))}
\ar[r]^-(0.5){b_1}\ar[l]_-(0.5){b_2}&
{\rR\rI\bvtau_*(\rI\hupsigma'^*(\Omega\otimes_{\co_{\fX'}}\umK^\bullet)\otimes_{\bvocB'}\bvtau^*(\cE))}}}}
\end{equation}
où les flèches verticales sont les morphismes \eqref{chb26d}, $a_1=\rI\huppi^*(\pi_1^\bullet)\otimes \id$, $a_2=\rI\huppi^*(v^\bullet \circ \pi^\bullet_2)\otimes\id$
$b_1=\rR\rI\bvtau_*(\rI\hupsigma'^*(\pi_1^\bullet)\otimes \id)$ et $b_2=\rR\rI\bvtau_*(\rI\hupsigma'^*(v^\bullet \circ \pi^\bullet_2)\otimes \id)$,
sont commutatifs. Comme $\cE$ est $\bvocB^!$-plat, $a_2$ est un quasi-isomorphisme \eqref{chb50},  
et $-\rI\huppi^\star(\partial)\otimes\id$ est le composé de l'inverse de $a_2$ et de $a_1$.  
De même, $\rI\hupsigma'^*(v^\bullet \circ \pi^\bullet_2)\otimes \id$ est un quasi-isomorphisme, et il en est donc de même de $b_2$, et 
$-\rR\rI \bvtau_*(\rI\hupsigma'^\star(\partial)\otimes\id)$ est le composé de l'inverse de $b_2$ et de $b_1$, d'où la proposition.

\begin{lem}\label{crindmd22}
Sous les hypothèses de \ref{crindmd14} et avec les mêmes notations, pour tout entier $q\geq 0$, le diagramme 
\begin{equation}\label{crindmd22a}
\xymatrix{
{\hupsigma^*_\mQ(\rR^q\fgg_*(\umK(\cN)))}\ar[r]\ar[d]_{\hupsigma^*_\mQ(\kappa^q)}&
{\rR^q\bvlgg_{\mQ*}(\huppi^*_\mQ(\umK(\cN)))}\ar[dd]^{\rR^q\bvlgg_{\mQ*}(\huppi^\star_\mQ(\partial))}\\
{\hupsigma^*_\mQ(\txi^{-1}\tOmega^1_{\fX/\cS} \otimes_{\co_{\fX}}\rR^q\fgg_*(\umK(\cN)))}\ar[d]&\\
{\hupsigma^*_\mQ(\rR^q\fgg_*(\fgg^*(\txi^{-1}\tOmega^1_{\fX/\cS}) \otimes_{\co_{\fX'}}\umK(\cN)))}\ar[r]&
{\rR^q\bvlgg_{\mQ*}(\huppi^*_\mQ(\fgg^*(\txi^{-1}\tOmega^1_{\fX/\cS}) \otimes_{\co_{\fX'}}\umK(\cN)))}}
\end{equation}
où les flèches horizontales sont induites par le morphisme de changement de base \eqref{chb20a} et la flèche verticale non libellée est l'isomorphisme canonique, 
est commutatif. 
\end{lem}

Il revient au même de dire que le diagramme 
\begin{equation}\label{crindmd22b}
\xymatrix{
{\hupsigma^*_\mQ(\rR^q\fgg_*(\umK(\cN)))}\ar[r]\ar[d]_{\hupsigma^*_\mQ(\rR^q\fgg_*(\partial))}&
{\rR^q\bvlgg_{\mQ*}(\huppi^*_\mQ(\umK(\cN)))}\ar[d]^{\rR^q\bvlgg_{\mQ*}(\huppi^\star_\mQ(\partial))}\\
{\hupsigma^*_\mQ(\rR^q\fgg_*(\fgg^*(\txi^{-1}\tOmega^1_{\fX/\cS}) \otimes_{\co_{\fX'}}\umK(\cN)))}\ar[r]&
{\rR^q\bvlgg_{\mQ*}(\huppi^*_\mQ(\fgg^*(\txi^{-1}\tOmega^1_{\fX/\cS}) \otimes_{\co_{\fX'}}\umK(\cN)))}}
\end{equation}
est commutatif. Reprenons les notations de la preuve de \ref{crindmd140}. D'après \ref{chb5}, la suite de complexes de ind-$\bvocB^!$-modules 
\begin{equation}
0\rightarrow \huppi^*_\mQ(\Omega\otimes_{\co_{\fX'}}\umK^\bullet)[-1]\rightarrow \huppi^*_\mQ(\mG^\bullet) 
\rightarrow \huppi^*_\mQ(\umK^\bullet) \rightarrow 0
\end{equation}
est exacte. Le cône de $\huppi^*_\mQ(u^\bullet)$ s'identifie à $\huppi^*_\mQ(\rC^\bullet)$, et le morphisme $\huppi_\mQ^\star(\partial)$ de 
$\bD^+(\bMod_\mQ(\bvocB^!))$ est alors défini de façon analogue à $\partial$. 

Comme le morphisme de changement de base \eqref{chb20a} est fonctoriel, les diagrammes de $\bD^+(\bIndMod(\bvocB))$
\begin{equation}
\text{\tiny {\xymatrix{
{\hupsigma^*_\mQ(\rR^q\fgg_*(\umK^\bullet))}\ar[d]&{\hupsigma^*_\mQ(\rR^q\fgg_*(\rC^\bullet))}\ar[r]^-(0.5){a_1}\ar[l]_-(0.5){a_2}
\ar[d] &{\hupsigma^*_\mQ(\rR^q\fgg_*(\Omega\otimes_{\co_{\fX'}}\umK^\bullet))}\ar[d]\\
{\rR^q\bvlgg_{\mQ*}(\huppi^*_\mQ(\umK^\bullet))}&{\rR^q\bvlgg_{\mQ*}(\huppi^*_\mQ(\rC^\bullet))}
\ar[r]^-(0.5){b_1}\ar[l]_-(0.5){b_2}&
{\rR^q\bvlgg_{\mQ*}(\huppi^*_\mQ(\Omega\otimes_{\co_{\fX'}}\umK^\bullet))}}}}
\end{equation}
où les flèches verticales sont les morphismes \eqref{chb20a}, $a_1=\hupsigma^*_\mQ(\rR^q\fgg_*(\pi_1^\bullet))$, 
$a_2=\hupsigma^*_\mQ(\rR^q\fgg_*(v^\bullet \circ \pi^\bullet_2))$
$b_1=\rR^q\bvlgg_{\mQ*}(\huppi^*_\mQ(\pi_1^\bullet))$ et $b_2=\rR^q\bvlgg_{\mQ*}(\huppi^*_\mQ(v^\bullet \circ \pi^\bullet_2))$,
sont commutatifs. On notera que $a_2$ est un isomorphisme, et que 
$-\hupsigma^*_\mQ(\rR^q\fgg_*(\partial))$ est le composé de l'inverse de $a_2$ et de $a_1$.  
Par ailleurs, $\huppi^*_\mQ(v^\bullet \circ \pi^\bullet_2)$ est un quasi-isomorphisme \eqref{chb5}. Par suite, $b_2$ est un isomorphisme, et 
$-\rR^q\bvlgg_{\mQ*}(\huppi^\star_\mQ(\partial))$ est le composé de l'inverse de $b_2$ et de $b_1$, d'où la proposition.

\subsection{}\label{crindmd145}
Soient $(\cN,\theta)$ un $\co_{\fX'}[\frac 1 p]$-fibré de Higgs à coefficients dans $\txi^{-1}\tOmega^1_{\fX'/\cS}$,
$t,r$ deux nombres rationnels tels que $t\geq r\geq 0$. On désigne par $\rI\hupsigma'^{(t,r)*}(\cN)$ l'image de $(\cN,\theta)$ par le foncteur 
\eqref{crindmd15e}, par $\utheta$ le $\co_{\fX'}[\frac 1 p]$-champ de de Higgs sur $\cN$ à coefficients dans $\txi^{-1}\tOmega^1_{\fX'/\fX}$ induit par $\theta$,
par $\rI\uhupsigma'^{(t,r)*}(\cN)$ l'image de $(\cN,\utheta)$ par le foncteur \eqref{crindmd2n}, et par $\cF^{(r)}$ le ind-$\bvuptheta^*(\bvcC^{(r)})$-module
\begin{equation}\label{crindmd145b}
\cF^{(r)}=\underset{\underset{t\in \mQ_{>r}}{\longrightarrow}}{\mlq\mlq\lim \mrq\mrq}\ \rI\hupsigma'^*(\cN)\otimes_{\bvocB'}\bvcC'^{(t,r)}. 
\end{equation}
Les ind-$\bvocB'$-modules de Higgs $\rI\hupsigma'^{(t,r)*}(\cN)$ à coefficients dans 
$\hupsigma'^*(\xi^{-1}\tOmega^1_{\fX'/\cS})$, pour $t\in \mQ_{>r}$, forment un système inductif \eqref{ccoh50e}. 
On désigne par 
\begin{equation}\label{crindmd145c}
\vartheta^{(r)}\colon \cF^{(r)} \rightarrow \hupsigma'^*(\xi^{-1}\tOmega^1_{\fX'/\cS})\otimes_{\bvocB'}\cF^{(r)}
\end{equation}
le $\bvocB'$-champ de Higgs induit par ceux des $\rI\hupsigma'^{(t,r)*}(\cN)$ pour $t\in \mQ_{>r}$, et par $\mK^\bullet(\cF^{(r)})$ le complexe de Dolbeault de 
$(\cF^{(r)}, \vartheta^{(r)})$. 
De même, les ind-$\bvuptheta^*(\bvcC^{(r)})$-modules de Higgs $\rI\uhupsigma'^{(t,r)*}(\cN)$ à coefficients dans 
$\hupsigma'^*(\xi^{-1}\tOmega^1_{\fX'/\fX})\otimes_{\bvocB'} \bvuptheta^*(\bvcC^{(r)})$, pour $t\in \mQ_{>r}$, forment un système inductif. 
On désigne par 
\begin{equation}\label{crindmd145d}
\uvartheta^{(r)}\colon \cF^{(r)} \rightarrow \hupsigma'^*(\xi^{-1}\tOmega^1_{\fX'/\fX})\otimes_{\bvocB'}\cF^{(r)}
\end{equation}
le $\bvuptheta^*(\bvcC^{(r)})$-champ de Higgs induit par ceux des $\rI\uhupsigma'^{(t,r)*}(\cN)$, pour $t\in \mQ_{>r}$, 
et par $\umK^\bullet(\cF^{(r)})$ le complexe de Dolbeault de $(\cF^{(r)}, \uvartheta^{(r)})$. 
On notera que le morphisme de ind-$\bvocB'$-modules sous-jacent à $\uvartheta^{(r)}$ est induit par $\vartheta^{(r)}$.
On note 
\begin{equation}\label{crindmd145e}
\partial^{(r)}\colon \umK^\bullet(\cF^{(r)})\rightarrow 
\hupsigma'^*(\fgg^*(\txi^{-1}\tOmega^1_{\fX/\cS})) \otimes_{\bvocB'} \umK^\bullet(\cF^{(r)})
\end{equation}
le morphisme de $\bD^+(\bIndMod(\bvocB'))$, bord de la filtration de Koszul de $\mK^\bullet(\cF^{(r)})$ 
associée à l'image par le foncteur $\hupsigma'^*$ de la suite exacte \eqref{hmdf410d}, cf. \eqref{indsh60g}.

D'après \ref{crindmd5}, il existe un isomorphisme canonique de $\bD^+(\bIndMod(\bvlgg^*(\bvcC^{(r)})))$ 
\begin{equation}\label{crindmd145f}
\rI\huppi^*(\umK^\bullet(\cN))\otimes_{\bvocB^!}\bvlgg^*(\bvcC^{(r)})\rightarrow 
\rR\rI \bvtau_*(\umK^\bullet(\cF^{(r)})).
\end{equation}

Compte tenu de \eqref{hmdf13i}, nous identifions $\bvlgg^*(d_{\bvcC^{(r)}})$ \eqref{ahttf30a} à une $\bvocB^!$-dérivation 
\begin{equation}\label{crindmd145g}
\bvlgg^*(d_{\bvcC^{(r)}})\colon \bvlgg^*(\bvcC^{(r)})\rightarrow \huppi^*(\fgg^*(\txi^{-1}\tOmega^1_{\fX/\cS}))\otimes_{\bvocB^!}\bvlgg^*(\bvcC^{(r)}).
\end{equation}

\begin{prop}\label{crindmd17}
Sous les hypothèses de \ref{crindmd145} et avec les mêmes notations, pour tout nombre rationnel $r\geq 0$, le diagramme de $\bD^+(\bIndMod(\bvocB'))$
\begin{equation}
\xymatrix{
{\rI\hupsigma'^*(\umK^\bullet(\cN))\otimes_{\bvocB^!}\bvuptheta^*(\bvcC^{(r)})}\ar[r]^-(0.4)a\ar[d]&
{\rI\hupsigma'^*(\fgg^*(\txi^{-1}\tOmega^1_{\fX/\cS})\otimes_{\co_{\fX'}}\umK^\bullet(\cN))\otimes_{\bvocB'}\bvuptheta^*(\bvcC^{(r)})}\ar[d]\\
{\umK^\bullet(\cF^{(r)})}\ar[r]^-(0.5){\partial^{(r)}}&
{\hupsigma'^*(\fgg^*(\txi^{-1}\tOmega^1_{\fX/\cS})})\otimes_{\bvocB'}\umK^\bullet(\cF^{(r)})}
\end{equation}
où les flèches verticales sont les morphismes canoniques et $a=\rI\hupsigma'^\star(\partial)\otimes\id+p^r\id\otimes \bvuptheta^*(d_{\bvcC^{(r)}})$ \eqref{crindmd14d}, est commutatif.
\end{prop} 

Cela résulte de \ref{indsh62} par fonctorialité du morphisme bord \eqref{indsh60g}. 

\begin{cor}\label{crindmd18}
Sous les hypothèses de \ref{crindmd145} et avec les mêmes notations, pour tout nombre rationnel $r\geq 0$, le diagramme de $\bD^+(\bIndMod(\bvocB^!))$
\begin{equation}
\xymatrix{
{\rI\huppi^*(\umK^\bullet(\cN))\otimes_{\bvocB^!}\bvlgg^*(\bvcC^{(r)})}\ar[r]^-(0.4)a\ar[d]&
{\rI\huppi^*(\fgg^*(\txi^{-1}\tOmega^1_{\fX/\cS})\otimes_{\co_{\fX'}}\umK^\bullet(\cN))\otimes_{\bvocB^!}\bvlgg^*(\bvcC^{(r)})}\ar[d]\\
{\rR\rI \bvtau_*(\umK^\bullet(\cF^{(r)}))}\ar[r]^-(0.5){b}&
{\huppi^*(\fgg^*(\txi^{-1}\tOmega^1_{\fX/\cS}))\otimes_{\bvocB^!}\rR\rI \bvtau_*(\umK^\bullet(\cF^{(r)}))}}
\end{equation}
où les flèches verticales sont induites par l'isomorphisme \eqref{crindmd145f},  
$a=\rI\huppi^\star(\partial)\otimes\id+p^r\id\otimes \bvlgg^*(d_{\bvcC^{(r)}})$ \eqref{crindmd14d} et $b$ est induit par $\rR\rI \bvtau_*(\partial^{(r)})$ et \ref{indsh47}, est commutatif.
\end{cor} 

Cela résulte de \ref{crindmd17}, \ref{crindmd140} et de la fonctorialité du morphisme \eqref{chb26d}.

\subsection{}\label{crindmd19}
Soit $\cM$ un ind-$\bvocB'$-module de Dolbeault \eqref{crindmd1}. On désigne par $\cH'(\cM)$ (resp. $\ucH'(\cM)$) le $\co_{\fX'}[\frac 1 p]$-module de Higgs à coefficients 
dans $\txi^{-1}\tOmega^1_{\fX'/\cS}$ (resp. $\txi^{-1}\tOmega^1_{\fX'/\fX}$) associé à $\cM$ \eqref{crindmd1a} (resp. \eqref{crindmd1d}), 
par $\mK^\bullet(\cH'(\cM))$ (resp. $\umK^\bullet(\ucH'(\cM))$) son complexe de Dolbeault, 
et par $\rI\huppi^*(\umK^\bullet(\ucH'(\cM)))$ l'image de ce dernier par le foncteur $\rI\huppi^*$ \eqref{chb50c}. 
Selon \ref{crindmd14}, la filtration de Koszul de $\mK^\bullet(\cH'(\cM))$ associée à la suite exacte \eqref{hmdf410d} détermine 
un morphisme ``bord'' de $\bD^+(\bMod(\co_{\fX'}[\frac 1 p]))$ \eqref{crindmd14b}
\begin{equation}\label{crindmd19a}
\partial_\cM\colon \umK^\bullet(\ucH'(\cM))\rightarrow \fgg^*(\txi^{-1}\tOmega^1_{\fX/\cS}) \otimes_{\co_{\fX'}} \umK^\bullet(\ucH'(\cM)),
\end{equation}
et un morphisme ``bord'' de $\bD^+(\bIndMod(\bvocB^!))$  \eqref{crindmd14f}
\begin{equation}\label{crindmd19b}
\rI\huppi^\star(\partial_\cM)\colon \rI\huppi^*(\umK^\bullet(\ucH'(\cM)))\rightarrow \rI \huppi^*(\fgg^*(\txi^{-1}\tOmega^1_{\fX/\cS}) \otimes_{\co_{\fX'}} \umK^\bullet(\ucH'(\cM))). 
\end{equation}

Pour tout entier $q\geq 0$, on désigne par 
\begin{equation}\label{crindmd19c}
\kappa^q_\cM\colon \rR^q\fgg_*(\umK^\bullet(\ucH'(\cM)))\rightarrow \txi^{-1}\tOmega^1_{\fX/\cS} \otimes_{\co_{\fX}} \rR^q\fgg_*(\umK^\bullet(\ucH'(\cM))) 
\end{equation}
le $\co_\fX$-champ de Katz-Oda \eqref{MH96}, qui n'est autre que le morphisme $\rR^q\fgg_*(\partial_\cM)$ .

\begin{prop}\label{crindmd20}
Sous les hypothèses de \ref{crindmd19}, il existe un nombre rationnel $r>0$ et un isomorphisme de $\bD^+(\bIndMod(\bvlgg^*(\bvcC^{(r)})))$
\begin{equation}\label{crindmd20a}
\rR\rI\bvtau_*(\cM\otimes_{\bvocB'}\bvuptheta^*(\bvcC^{(r)}))\stackrel{\sim}{\rightarrow}
\rI\huppi^*(\umK^\bullet(\ucH'(\cM)))\otimes_{\bvocB^!}\bvlgg^*(\bvcC^{(r)}),
\end{equation}
où $\rI \bvtau_*$ est le foncteur \eqref{ahttfg20h} et le produit tensoriel à droite est défini terme à terme (non dérivé), tels que le diagramme de 
$\bD^+(\bIndMod(\bvocB^!))$
\begin{equation}\label{crindmd20b}
\text{ {\tiny \xymatrix{
{\rR\rI\bvtau_*(\cM\otimes_{\bvocB'}\bvuptheta^*(\bvcC^{(r)}))}\ar[r]\ar[d]_a&{\rI\huppi^*(\umK^\bullet(\ucH'(\cM)))\otimes_{\bvocB^!}\bvlgg^*(\bvcC^{(r)})}\ar[d]^b\\
{\rR\rI\bvtau_*(\cM\otimes_{\bvocB'}\bvuptheta^*(\bvcC^{(r)}\otimes_{\bvocB}\hupsigma^*(\txi^{-1}\tOmega^1_{\fX/\cS})))}\ar[r]&
{\rI\huppi^*(\fgg^*(\txi^{-1}\tOmega^1_{\fX/\cS})\otimes_{\co_{\fX'}}\umK^\bullet(\ucH'(\cM)))\otimes_{\bvocB^!}\bvlgg^*(\bvcC^{(r)})}}}}
\end{equation}
où les flèches horizontales sont induites par l'isomorphisme \eqref{crindmd20a} et \ref{indsh47}, $a=p^r\rR\rI\bvtau_*(\id\otimes \bvuptheta^*(d_{\bvcC^{(r)}}))$ 
\eqref{ahttf30a} et $b=\rI\huppi^\star(\partial_\cM)\otimes\id+p^r\id\otimes \bvlgg^*(d_{\bvcC^{(r)}})$ \eqref{crindmd19b}, soit commutatif. 
\end{prop}

Pour alléger les notations posons $(\cN,\theta)=\cH'(\cM)$ et $(\cN,\utheta)=\ucH'(\cM)$; 
on omettra les champs de Higgs lorsque cela n'induit aucun risque d'ambiguïté. 
D'après \ref{indmdlb12}, il existe un nombre rationnel $t>0$ et un isomorphisme de $\bIndMC(\bvcC'^{(t)}/\bvocB')$ 
\begin{equation}\label{crindmd20c}
\alpha\colon \rI\fS'^{(t)}(\cM)\stackrel{\sim}{\rightarrow} \rI\hupsigma'^{(t)*}(\cN).
\end{equation} 
Pour tout nombre rationnel $t'$ tel que $0\leq t'\leq  t$, on désigne par
\begin{equation}\label{crindmd20d}
\alpha^{(t')}\colon \rI\fS'^{(t')}(\cM) \stackrel{\sim}{\rightarrow} \rI\hupsigma'^{(t')*}(\cN)
\end{equation}
l'isomorphisme de $\bIndMC(\bvcC'^{(t')}/\bvocB')$ induit par $\rI\varepsilon'^{t,t'}(\alpha)$ \eqref{crindmd2} et 
l'analogue des isomorphismes \eqref{indmdlb3b} et \eqref{indmdlb3c}.
Pour tous nombres rationnels $t',r$ tels que $0\leq r\leq t'\leq  t$, on désigne par
\begin{equation}\label{crindmd20e}
\alpha^{(t',r)}\colon \rI\fS'^{(t',r)}(\cM) \stackrel{\sim}{\rightarrow} \rI\hupsigma'^{(t',r)*}(\cN)
\end{equation}
l'isomorphisme de $\bIndMH(\bvocB',\hupsigma'^*(\txi^{-1}\tOmega^1_{\fX'/\cS}))$ induit par $\rI\lambda^{t',r}(\alpha^{(t')})$ \eqref{crindmd15c},
où $\rI\fS'^{(t',r)}$ est le foncteur \eqref{crindmd15g} et $\rI\hupsigma'^{(t',r)*}$ est le foncteur \eqref{crindmd15e}. 
Soit $r$ un nombre rationnel tel que $0< r< t$. Comme les $\alpha^{(t')}$ \eqref{crindmd20d}, pour $r<t'\leq t$, 
forment un isomorphisme de systèmes inductifs de ind-$\bvocB'$-modules de Higgs à coefficients dans $\hupsigma'^*(\txi^{-1}\tOmega^1_{\fX'/\cS})$,
il en est de même des $\alpha^{(t',r)}$ \eqref{crindmd20e} compte tenu de \eqref{ccoh50e}.
Reprenant les notations de \ref{crindmd145}, on en déduit par passage à la limite inductive dans 
$\bIndMH(\bvocB',\hupsigma'^*(\txi^{-1}\tOmega^1_{\fX'/\cS}))$ un morphisme 
\begin{equation}\label{crindmd20f}
\cM\otimes_{\bvocB'}\bvcC'^{(t',r)}\rightarrow \cF^{(r)},
\end{equation}
où $\cM\otimes_{\bvocB'}\bvcC'^{(t',r)}$ est muni du $\bvocB'$-champ de Higgs $\id\otimes \delta_{\bvcC'^{(t',r)}}$ \eqref{ccoh50b},
et par suite des morphismes de complexes d'ind-$\bvocB'$-modules
\begin{eqnarray}
\cM\otimes_{\bvocB'}\mK^\bullet(\bvcC'^{(t',r)}) &\rightarrow& \mK^\bullet(\cF^{(r)}),\label{crindmd20g1}\\
\cM\otimes_{\bvocB'}\umK^\bullet(\bvcC'^{(t',r)}) &\rightarrow& \umK^\bullet(\cF^{(r)}),\label{crindmd20g2}
\end{eqnarray}
où les complexes de Dolbeault $\mK^\bullet(\bvcC'^{(t',r)})$ et $\umK^\bullet(\bvcC'^{(t',r)})$ sont définis dans \ref{ccoh51}. 
D'après \ref{indmdlb6}, le ind-$\bvocB'$-module $\cM$ est rationnel et plat. 
Prenant les morphismes bord associés aux filtrations de Koszul \eqref{indsh60g} des complexes de Dolbeault apparaissant dans \eqref{crindmd20g1} 
relativement à l'image par le foncteur $\hupsigma'^*$ de la suite exacte \eqref{hmdf410d}, 
on obtient par fonctorialité un diagramme commutatif de $\bD^+(\bIndMod(\bvocB'))$
\begin{equation}\label{crindmd20h}
\xymatrix{
{\cM\otimes_{\bvocB'}\umK^\bullet(\bvcC'^{(t',r)})}\ar[r]\ar[d]_{\id\otimes\partial^{(t',r)}}&{\umK^\bullet(\cF^{(r)})}\ar[d]^{\partial^{(r)}}\\
{\hupsigma'^*(\fgg^*(\txi^{-1}\tOmega^1_{\fX/\cS})) \otimes_{\bvocB'} \cM\otimes_{\bvocB'}\umK^\bullet(\bvcC'^{(t',r)})}\ar[r]&{\hupsigma'^*(\fgg^*(\txi^{-1}\tOmega^1_{\fX/\cS})) \otimes_{\bvocB'} \umK^\bullet(\cF^{(r)})}}
\end{equation}
où les flèches horizontales sont induites par \eqref{crindmd20g2}, et 
$\partial^{(t',r)}$ (resp. $\partial^{(r)}$) est le morphisme bord \eqref{ccoh51a} (resp. \eqref{crindmd145e}). 

D'après \eqref{crindmd4c} et \eqref{crindmd145b}, 
les morphismes \eqref{crindmd20g2} induisent un quasi-isomorphisme 
\begin{equation}\label{crindmd20i}
\cM\otimes_{\bvocB'}\bvuptheta^*(\bvcC^{(r)})[0]\rightarrow \umK^\bullet(\cF^{(r)}). 
\end{equation}
Compte tenu de \eqref{ccoh51b} et \eqref{crindmd20h}, le diagramme 
\begin{equation}\label{crindmd20j}
\xymatrix{
{\cM\otimes_{\bvocB'}\bvuptheta^*(\bvcC^{(r)})[0]}\ar[r]\ar[d]_{p^r\id\otimes \bvuptheta^*(d_{\bvcC^{(r)}})}&{\umK^\bullet(\cF^{(r)})}\ar[d]^{\partial^{(r)}}\\
{\bvuptheta^*(\hupsigma^*(\txi^{-1}\tOmega^1_{\fX/\cS})) \otimes_{\bvocB'} \cM\otimes_{\bvocB'}\bvuptheta^*(\bvcC^{(r)})[0]}\ar[r]&{\hupsigma'^*(\fgg^*(\txi^{-1}\tOmega^1_{\fX/\cS})) \otimes_{\bvocB'} \umK^\bullet(\cF^{(r)})}}
\end{equation}
où les flèches horizontales sont induites par \eqref{crindmd20i} et \eqref{hmdf13i}, est commutatif. 
La proposition s'ensuit en vertu de \ref{crindmd5} et \ref{crindmd18}.

\begin{teo}\label{crindmd21}
Supposons $g\colon X'\rightarrow X$ propre. 
Soient $\cM$ un ind-$\bvocB'$-module de Dolbeault \eqref{crindmd1}, $q$ un entier $\geq 0$.
On désigne par $\ucH'(\cM)$ le $\co_{\fX'}[\frac 1 p]$-module de Higgs à coefficients 
dans $\txi^{-1}\tOmega^1_{\fX'/\fX}$ associé à $\cM$ \eqref{crindmd1d} et par $\umK^\bullet(\ucH'(\cM))$ son complexe de Dolbeault.
On munit le $\co_\fX[\frac 1 p]$-module $\rR^q\fgg_*(\umK^\bullet(\ucH'(\cM)))$ du champ de Katz-Oda $\kappa^q_\cM$ \eqref{crindmd19c}.
Alors, le $\co_\fX[\frac 1 p]$-module $\rR^q\fgg_*(\umK^\bullet(\ucH'(\cM)))$ est cohérent, 
et il existe un nombre rationnel $r>0$, indépendant de $q$, et un isomorphisme de $\bIndMC(\bvcC^{(r)}/\bvocB)$ 
\begin{equation}\label{crindmd21a}
\rI\fS^{(r)}(\rR^q\rI\bvuptheta_*(\cM))\stackrel{\sim}{\rightarrow}
\rI\hupsigma^{(r)*}(\rR^q\fgg_*(\umK^\bullet(\ucH'(\cM))),\kappa^q_\cM),
\end{equation}
où les foncteurs $\rI\fS^{(r)}$ et $\rI\hupsigma^{(r)*}$ sont définis dans \eqref{indmdlb1b} et \eqref{indmdlb1f} compte tenu de \eqref{indmdlb2f}.
\end{teo}

En effet, en vertu de \ref{crindmd13}, le $\co_\fX[\frac 1 p]$-module $\rR^q\fgg_*(\umK^\bullet(\ucH'(\cM)))$ est cohérent, 
et il existe un nombre rationnel $r>0$ et un isomorphisme de ind-$\bvcC^{(r)}$-modules 
\begin{equation}\label{crindmd21b}
\rR^q\rI\bvuptheta_*(\cM)\otimes_{\bvocB}\bvcC^{(r)}\stackrel{\sim}{\rightarrow}
\rI\hupsigma^*(\rR^q\fgg_*(\umK^\bullet(\ucH'(\cM))))\otimes_{\bvocB}\bvcC^{(r)}.
\end{equation}
De plus, il résulte de \ref{crindmd20}, \ref{crindmd12}, \ref{chb22}, \ref{crindmd22}, \eqref{indsh14h} et \eqref{indsh22f} qu'on 
peut trouver un tel isomorphisme compatible avec les champs de Higgs, la source et le but étant munis des champs de Higgs totaux, d'où la proposition. 

\begin{cor}\label{crindmd26}
Soient $\cM$ un ind-$\bvocB'$-module de Dolbeault \eqref{crindmd1}, $q$ un entier $\geq 0$. 
On désigne par $\ucH'(\cM)$ le $\co_{\fX'}[\frac 1 p]$-module de Higgs à coefficients 
dans $\txi^{-1}\tOmega^1_{\fX'/\fX}$ associé à $\cM$ \eqref{crindmd1d} et par $\umK^\bullet(\ucH'(\cM))$ son complexe de Dolbeault.
On munit le $\co_\fX[\frac 1 p]$-module $\rR^q\fgg_*(\umK^\bullet(\ucH'(\cM)))$ du champ de Katz-Oda $\kappa^q_\cM$ \eqref{crindmd19c}.
Supposons le morphisme $g\colon X'\rightarrow X$ propre et le $\co_\fX[\frac 1 p]$-module $\rR^q\fgg_*(\umK^\bullet(\ucH'(\cM)))$ localement projectif de type fini. 
Alors, le ind-$\bvocB$-module $\rR^q\rI\bvuptheta_*(\cM)$ est de Dolbeault, le $\co_\fX[\frac 1 p]$-module de Higgs $(\rR^q\fgg_*(\umK^\bullet(\ucH'(\cM))),\kappa^q_\cM)$ 
est soluble \eqref{indmdlb5}, et on a un isomorphisme de $\co_\fX[\frac 1 p]$-fibrés de Higgs   
\begin{equation}\label{crindmd26a}
\cH(\rR^q\rI\bvuptheta_*(\cM))\stackrel{\sim}{\rightarrow} (\rR^q\fgg_*(\umK^\bullet(\ucH'(\cM))), \kappa_\cM^q),
\end{equation}
où $\cH$ est le foncteur \eqref{indmdlb7c}. 
\end{cor}

Cela résulte de \ref{crindmd21} et \ref{indmdlb10}.

\subsection{}\label{crindmd23}
On reprend les notations de \ref{sld2} pour $f$ et on considère les notations analogues pour $f'$, que l'on munit d'un exposant $^\prime$ \eqref{hmdf3}. 
On a alors le diagramme commutatif de morphismes de topos 
\begin{equation}\label{crindmd23a}
\xymatrix{
{(\oX'^\rhd_\et)^{\mN^\circ}}\ar[d]_{\bvupgamma}\ar[r]^{\bvpsi'}&{\tE'^{\mN^\circ}}\ar[d]^{\bvTheta}\\
{(\oX^\circ_\et)^{\mN^\circ}}\ar[r]^\bvpsi&{\tE^{\mN^\circ}}}
\end{equation}
induits par $\psi$, $\psi'$ \eqref{hmdf7h}, $\upgamma$ et $\Theta$ \eqref{hmdf11c}. 
On les considérera naturellement comme des morphismes de topos annelés par les anneaux $\bvmZ_p$ \eqref{sld1}. 

On désigne par 
\begin{equation}\label{crindmd23b}
\cH'_\mQ\colon \bMod_\mQ(\bvocB')\rightarrow \bMH(\co_{\fX'}[\frac 1 p], \txi^{-1}\tOmega^1_{\fX'/\cS})
\end{equation}
le foncteur défini dans \eqref{aspglob11b}, et par 
\begin{equation}\label{crindmd23c}
\ucH'_\mQ\colon \bMod_\mQ(\bvocB')\rightarrow \bMH(\co_{\fX'}[\frac 1 p], \txi^{-1}\tOmega^1_{\fX'/\fX})
\end{equation}
le composé du foncteur $\cH'_\mQ$ et du foncteur canonique \eqref{crindmd1e}.  
Ces foncteurs sont compatibles avec les foncteurs \eqref{crindmd1a} et \eqref{crindmd1d}; cf. \eqref{aspglob11d}.

\begin{prop}\label{crindmd130}
Supposons $g\colon X'\rightarrow X$ propre. 
Soient $\cM$ un $\bvocB'_\mQ$-module de Dolbeault \eqref{aspglob1}, $q$ un entier $\geq 0$. 
On désigne par $\ucH'_\mQ(\cM)$ le $\co_{\fX'}[\frac 1 p]$-fibré de Higgs à coefficients dans 
$\txi^{-1}\tOmega^1_{\fX'/\fX}$ associé à $\cM$ \eqref{crindmd23c} et par $\umK^\bullet$ son complexe de Dolbeault. 
Alors, le $\co_\fX[\frac 1 p]$-module $\rR^q\fgg_*(\umK^\bullet)$ est cohérent, 
et il existe un nombre rationnel $r>0$, indépendant de $q$, et un isomorphisme de $\bvcC^{(r)}_\mQ$-modules
\begin{equation}\label{crindmd130a}
\rR^q\bvuptheta_{\mQ*}(\cM)\otimes_{\bvocB_\mQ}\bvcC^{(r)}_\mQ\stackrel{\sim}{\rightarrow}
\hupsigma^*_\mQ(\rR^q\fgg_*(\umK^\bullet))\otimes_{\bvocB_\mQ}\bvcC^{(r)}_\mQ,
\end{equation}
où $\hupsigma^*_\mQ$ est le foncteur \eqref{hmdf40j}.
\end{prop}

Cela résulte de \ref{crindmd13}, compte tenu de \eqref{indsh13a}, \eqref{indsh14h}
et du fait que le foncteur $\upalpha_{\bvcC^{(r)}}$ est pleinement fidèle \eqref{indsh11c}.

\begin{cor}\label{crindmd24}
Soit $M=(M_n)_{n\geq 0}$ un $\bvmZ_p$-système local de $(\oX'^\rhd_\et)^{\mN^\circ}$ \eqref{sld3}. Posons $\cM=\bvpsi'_*(M)\otimes_{\bvmZ_p}\bvocB'$, 
et notons $\ucH'_\mQ(\cM_\mQ)$ le $\co_{\fX'}[\frac 1 p]$-fibré de Higgs à coefficients dans $\txi^{-1}\tOmega^1_{\fX'/\fX}$ 
associé à $\cM_\mQ$ \eqref{crindmd23c} et $\umK^\bullet$ son complexe de Dolbeault. 
Supposons les conditions suivantes satisfaites:
\begin{itemize}
\item[{\rm (i)}] le morphisme $g\colon X'\rightarrow X$ est propre;
\item[{\rm (ii)}] le $\bvmZ_{p,\mQ}$-système local $M_\mQ$ de Dolbeault \eqref{sld5}, i.e., le $\bvocB'_\mQ$-module $\cM_\mQ$ de Dolbeault.
\end{itemize}
Alors, les $\co_\fX[\frac 1 p]$-modules $\rR^i\fgg_*(\rH^j(\umK^\bullet))$ sont cohérents pour tous $i,j\geq 0$,
et il existe un nombre rationnel $r>0$ et une suite spectrale de $\bvcC^{(r)}_\mQ$-modules 
\begin{equation}\label{crindmd24a}
\rE_2^{i,j}=\hupsigma^*_\mQ(\rR^i\fgg_*(\rH^j(\umK^\bullet)))\otimes_{\bvocB_\mQ}\bvcC^{(r)}_\mQ
\Rightarrow (\bvpsi_*(\rR^{i+j}\bvupgamma_*(M))\otimes_{\bvmZ_p}\bvcC^{(r)})_\mQ,
\end{equation}
où $\hupsigma^*_\mQ$ est le foncteur \eqref{hmdf40j}.
\end{cor}

En effet, la première assertion résulte de (\cite{egr1} 2.10.24 et 2.11.5).
D'après \ref{chb5}, pour tout nombre rationnel $r\geq 0$, comme $\bvcC^{(r)}$ est $\bvocB$-plat, la seconde suite spectrale d'hypercohomologie 
\begin{equation}
\rE_2^{i,j}=\rR^i\fgg_*(\rH^j(\umK^\bullet))\Rightarrow \rR^{i+j}\fgg_*(\umK^\bullet)
\end{equation}
induit une suite spectrale de $\bvcC^{(r)}_\mQ$-modules 
\begin{equation}
\rE_2^{i,j}=\hupsigma^*_\mQ(\rR^i\fgg_*(\rH^j(\umK^\bullet)))\otimes_{\bvocB_\mQ}\bvcC^{(r)}_\mQ\Rightarrow 
\hupsigma^*_\mQ(\rR^{i+j}\fgg_*(\umK^\bullet))\otimes_{\bvocB_\mQ}\bvcC^{(r)}_\mQ.
\end{equation}
La proposition s'ensuit en vertu de \ref{crindmd130}, (\cite{ag} 5.7.4 et 5.7.6) et (\cite{agt} III.7.5).

\begin{cor}\label{crindmd240}
Supposons le morphisme $g\colon X'\rightarrow X$ propre. 
Alors, il existe un nombre rationnel $r>0$ et pour tout entier $n\geq 0$, un isomorphisme canonique de $\bvcC^{(r)}_\mQ$-modules 
\begin{equation}\label{crindmd240a}
\bvpsi_*(\rR^n\bvupgamma_*(\bvmZ_p))\otimes_{\bvmZ_p}\bvcC^{(r)}_\mQ
\stackrel{\sim}{\rightarrow}\oplus_{0\leq i\leq n}\hupsigma^*(\rR^i\fgg_*(\tOmega^{n-i}_{\fX'/\fX}))\otimes_{\bvocB}\bvcC^{(r)}_\mQ(i-n),
\end{equation}
où $\hupsigma$ est le morphisme de topos annelés \eqref{ahttf13e}.
\end{cor}

Cela résulte de \ref{crindmd130} appliqué à $\cM=\bvocB'_\mQ$, (\cite{ag} 5.7.5 et 5.7.6) et (\cite{agt} III.7.5), puisque  le $\bvocB'_\mQ$-module $\bvocB'_\mQ$ est de Dolbeault et que 
$\ucH'_\mQ(\bvocB'_\mQ)$ est le fibré trivial $\co_{\fX'}[\frac 1 p]$ muni du champ de Higgs nul \eqref{aspglob12}. 

\begin{rema}
Le corollaire \ref{crindmd240} s'applique en particulier en prenant pour $\tg$ dans le cas relatif \eqref{definf10} la déformation triviale \eqref{hmdf50}. 
\end{rema}

\begin{teo}[\cite{ag} 6.7.5]\label{crindmd25}
Si le morphisme $g\colon X'\rightarrow X$ est propre, il existe une suite spectrale canonique de $\bvocB_\mQ$-modules
\begin{equation}\label{crindmd25a}
\rE_2^{i,j}=\hupsigma^*(\rR^i\fgg_*(\tOmega^j_{\fX'/\fX}))\otimes_{\bvocB}\bvocB_\mQ(-j)\Rightarrow \bvpsi_*(\rR^{i+j}\bvupgamma_*(\bvmZ_p))\otimes_{\mZ_p}\bvocB_\mQ,
\end{equation}
où $\hupsigma$ est le morphisme de topos annelés \eqref{ahttf13e}.
\end{teo}

Dans l'énoncé (\cite{ag} 6.7.5), on requiert en fait que $g$ soit projectif. Toutefois, comme signalé sous (\cite{ag} 6.7.6), le résultat vaut sous l'hypothèse 
plus générale que $g$ soit propre. En effet, l'hypothèse de projectivité sur $g$ est utilisée dans la preuve de (\cite{ag} 5.7.4)
qui s'étend également aux morphismes propres (voir \cite{ag} 5.7.6). 

Cette suite spectrale ne requiert  la considération d'aucune déformation \eqref{hmdf5}. 
Elle est $G_K$-équivariante (\cite{ag} 6.7.10) et elle dégénère en $\rE_2$ (\cite{ag} 6.7.13). 
Mais la filtration aboutissement n'est pas scindée en général (cf. \cite{ag} 1.3.2 et 1.3.3). 
Toutefois, on peut vérifier qu'elle se scinde après changement de base de $\bvocB$
à $\bvcC^{(r)}$ pour un nombre rationnel $r>0$, et qu'elle correspond à la décomposition \eqref{crindmd240a}.

\begin{prop}\label{crindmd27}
Supposons $g\colon X'\rightarrow X$ propre. 
Soient $\cM$ un $\bvocB'_\mQ$-module de Dolbeault \eqref{aspglob1}, $q$ un entier $\geq 0$. 
On désigne par $\ucH'_\mQ(\cM)$ le $\co_{\fX'}[\frac 1 p]$-fibré de Higgs à coefficients dans $\txi^{-1}\tOmega^1_{\fX'/\fX}$ associé à $\cM$ \eqref{crindmd23c}, 
par $\umK^\bullet$ son complexe de Dolbeault et par 
\begin{equation}\label{crindmd27a}
\kappa^q_\cM\colon \rR^q\fgg_*(\umK^\bullet)\rightarrow \txi^{-1}\tOmega^1_{\fX/\cS} \otimes_{\co_{\fX}} \rR^q\fgg_*(\umK^\bullet)
\end{equation}
le champ de Katz-Oda \eqref{MH96}. 
Alors, le $\co_\fX[\frac 1 p]$-module $\rR^q\fgg_*(\umK^\bullet)$ est cohérent, 
et il existe un nombre rationnel $r>0$, indépendant de $q$, et un isomorphisme de $\bIMC_\mQ(\bvcC^{(r)}/\bvocB)$ 
\begin{equation}\label{crindmd27b}
\fS^{(r)}(\rR^q\bvuptheta_{\mQ*}(\cM))\stackrel{\sim}{\rightarrow}
\hupsigma^{(r)*}(\rR^q\fgg_*(\umK^\bullet),\kappa^q_\cM),
\end{equation}
où les foncteurs $\fS^{(r)}$ et $\hupsigma^{(r)*}$ sont définis dans \eqref{aspglob6b} et \eqref{aspglob6g}.
\end{prop}

Cela résulte de \ref{crindmd21} au moyen du foncteur pleinement fidèle $\upalpha_{\bvcC^{(r)}/\bvocB}$ \eqref{dolbff1c}, 
compte tenu de \eqref{dolbff1d}, \eqref{dolbff1e} et \eqref{aspglob11d}.

\begin{cor}\label{crindmd28}
Soient $\cM$ un $\bvocB'_\mQ$-module de Dolbeault \eqref{aspglob1}, $q$ un entier $\geq 0$. 
On désigne par $\ucH'_\mQ(\cM)$ le $\co_{\fX'}[\frac 1 p]$-fibré de Higgs à coefficients dans $\txi^{-1}\tOmega^1_{\fX'/\fX}$ associé à $\cM$ \eqref{crindmd23c}, 
par $\umK^\bullet$ son complexe de Dolbeault et par
\begin{equation}\label{crindmd28a}
\kappa^q_\cM\colon \rR^q\fgg_*(\umK^\bullet)\rightarrow \txi^{-1}\tOmega^1_{\fX/\cS} \otimes_{\co_{\fX}} \rR^q\fgg_*(\umK^\bullet)
\end{equation}
le champ de Katz-Oda \eqref{MH96}. 
Supposons le morphisme $g\colon X'\rightarrow X$ propre et le $\co_\fX[\frac 1 p]$-module $\rR^q\fgg_*(\umK^\bullet)$ localement projectif de type fini. 
Alors, le $\bvocB_\mQ$-module $\rR^q\bvuptheta_{\mQ*}(\cM)$ est de Dolbeault, le $\co_\fX[\frac 1 p]$-module de Higgs $(\rR^q\fgg_*(\umK^\bullet),\kappa_\cM^q)$ 
est rationnellement soluble \eqref{aspglob1}, et on a un isomorphisme de $\co_\fX[\frac 1 p]$-fibrés de Higgs   
\begin{equation}\label{crindmd28b}
\cH_\mQ(\rR^q\bvuptheta_{\mQ*}(\cM))\stackrel{\sim}{\rightarrow} (\rR^q\fgg_*(\umK^\bullet), \kappa_\cM^q),
\end{equation}
où $\cH_\mQ$ est le foncteur \eqref{aspglob11b}. 
\end{cor}

Cela résulte de \ref{crindmd27} et \ref{aspglob14}.

\begin{cor}\label{crindmd29}
Supposons le morphisme $g\colon X'\rightarrow X$ propre. Pour tout entier $q\geq 0$, posons $\cM^q=\bvpsi_*(\rR^q\bvupgamma_*(\bvmZ_p))\otimes_{\bvmZ_p}\bvocB$.  
Alors, le $\bvocB_\mQ$-module $\cM^q_\mQ$ est de Hodge-Tate \eqref{sld5}, et on a un isomorphisme de $\co_\fX[\frac 1 p]$-fibrés de Higgs
\begin{equation}\label{crindmd29a}
\cH_\mQ(\cM^q_\mQ)\stackrel{\sim}{\rightarrow} \oplus_{0\leq i\leq q}\rR^i\fgg_*(\txi^{i-q}\tOmega^{q-i}_{\fX'/\fX}) \otimes_{\co_\fX}\co_\fX[\frac 1 p],
\end{equation}
où $\cH_\mQ$ est le foncteur \eqref{aspglob11b}, 
le champ de Higgs sur le terme de droite étant induit par les applications de Kodaira-Spencer de $\fgg$ \eqref{hmdf410f}. 
\end{cor}

En effet, il résulte de (\cite{ag} 5.7.5 et 5.7.6) et (\cite{agt} III.7.5) qu'on a un isomorphisme canonique 
\begin{equation}
\cM^q_\mQ\stackrel{\sim}{\rightarrow} \rR^q\bvuptheta_{\mQ*}(\bvocB'_\mQ).
\end{equation}

D'après \ref{aspglob12}, le $\bvocB'_\mQ$-module $\bvocB'_\mQ$ est de Dolbeault et 
le module de Higgs $\ucH'_\mQ(\bvocB'_\mQ)$ \eqref{crindmd23c} est le fibré trivial $\co_{\fX'}[\frac 1 p]$ muni du champ de Higgs nul. 
Son complexe de Dolbeault $\umK^\bullet$ est donc donné par 
\begin{equation}
\umK^\bullet=\oplus_{i\geq 0} \txi^{-i}\tOmega^i_{\fX'/\fX}\otimes_{\co_{\fX'}}\co_{\fX'}[\frac 1 p][-i].
\end{equation}
Par suite, on a   
\begin{equation}\label{crindmd29c}
\rR^q\fgg_*(\umK^\bullet)=\oplus_{0\leq i\leq q}\rR^i\fgg_*(\txi^{i-q}\tOmega^{q-i}_{\fX'/\fX}) \otimes_{\co_\fX}\co_\fX[\frac 1 p].
\end{equation}
Le champ de Katz-Oda \eqref{MH96}
\begin{equation}\label{crindmd29d}
\kappa^q\colon \rR^q\fgg_*(\umK^\bullet)\rightarrow \txi^{-1}\tOmega^1_{\fX/\cS} \otimes_{\co_{\fX}} \rR^q\fgg_*(\umK^\bullet)
\end{equation}
est induit par les applications de Kodaira-Spencer de $\fgg$ \eqref{hmdf410f}. Il est donc nilpotent \eqref{definf22}.

En vertu de \ref{hmdf43} et \ref{crindmd28}, le $\bvocB_\mQ$-module $\rR^q\bvuptheta_{\mQ*}(\bvocB'_\mQ)$ est de Hodge-Tate, 
et on a un isomorphisme de $\co_\fX[\frac 1 p]$-fibrés de Higgs   
\begin{equation}\label{crindmd29e}
\cH_\mQ(\rR^q\bvuptheta_{\mQ*}(\bvocB'_\mQ))\stackrel{\sim}{\rightarrow} (\rR^q\fgg_*(\umK^\bullet), \kappa^q), 
\end{equation}
d'où la proposition. 

\begin{rema}
Le corollaire \ref{crindmd29} s'applique en particulier en prenant pour $\tg$ dans le cas relatif \eqref{definf10} la déformation triviale \eqref{hmdf50}. 
\end{rema}

\chapter*{\texorpdfstring{Errata and  Addenda to ``The $p$-adic Simpson correspondence''}
{Errata and  Addenda to ``The p-adic Simpson correspondence''}} 
\begin{center}
by A. Abbes, M. Gros and T. Tsuji, \href{http://press.princeton.edu/titles/10779.html}{Ann. of Math. Stud. {\bf 193}}, Princeton Univ. Press (2016)
\end{center}

\vspace{0.2cm}

\begin{center}

{\large \bf A) Misprints}

\end{center}
 
\vspace{0.2cm} 

{\bf (I.1.2).} Line 1,  replace {\em``valuation ring''} by {\em  ``valuation field''}.

\vspace{0.2cm}

{\bf (I.2.1).} Line 1,  replace  {\em``valuation ring''} by  {\em``valuation field''}.

\vspace{0.2cm}

{\bf (I.4.6.1).} One line under (I.4.6.1), replace $\Gamma = \pi_1(X, \oy)$ by $\Gamma = \pi_1(X_{\eta}, \oy)$.

\vspace{0.2cm}

{\bf (I.4.7.4).} Two lines above (I.4.7.4), replace $\alpha^{r,r'} : \cC^{(r')} \rightarrow \cC^{(r)}$ by $\alpha^{r,r'} : \cC^{(r)} \rightarrow \cC^{(r')}$.

\vspace{0.2cm}

{\bf (I.4.7.4).} One line above (I.4.7.4), replace $h_{\alpha}^{r,r'} : \hcC^{(r')} \rightarrow \hcC^{(r)}$ by $h_{\alpha}^{r,r'} : \hcC^{(r)} \rightarrow \hcC^{(r')}$.

\vspace{0.2cm}

{\bf (I.4.7.7).} Replace $p^{r'} ({\rm{id}} \times \alpha^{r,r'}) \circ d_{\cC^{(r')}} =  p^r d_{\cC^{(r)}} \circ  \alpha^{r,r'}$ by 
$p^{r} ({\rm{id}} \times \alpha^{r,r'}) \circ d_{\cC^{(r)}} =  p^{r'} d_{\cC^{(r')}} \circ  \alpha^{r,r'}$.

\vspace{0.2cm}

{\bf (I.4.13.1).} Replace the target $\oplus_{i \in \mZ}{\rm{D}}^{i}(V)\otimes_{\hR} \hRun(-1)$ by $\oplus_{i \in \mZ}{\rm{D}}^{i}(V)\otimes_{\hR} \hRun(-i)$. 

\vspace{0.2cm}

{\bf (I.5.12.7).} One line under (I.5.12.7), replace  {\em``multiplication by  $p^{r'-r}$''} by  {\em``multiplication by  $p^{r-r'}$''}.

\vspace{0.2cm}

{\bf (II.2.1).} Line 1-2,  replace  {\em``valuation ring''} by  {\em``valuation field''}.

\vspace{0.2cm}

{\bf (II.3.9).} Line 5, replace {\em``has a left inverse''}  by   {\em``has a right inverse''}.

\vspace{0.2cm}

{\bf (II.3.12).} Line 1-2, replace  {\em``$M$ a topological  $A$-$G$-module, and $N$ a topological  $A$-$H$-module''}
by  {\em``$M$ a linearly topologized  $A$-$G$-module, and $N$ a linearly topologized  $A$-$H$-module''}.

\vspace{0.2cm}

{\bf (II.3.34).} Line 4-5, replace  {\em``with values in $\id_r+a^m\Mat_r(A/a^qA)$ associated with $N$  and $N'$, respectively''} 
by  {\em``with values in $\id_r+a^m\Mat_r(A/a^qA)$ associated with $\rho$ and $\rho'$, respectively''}.

\vspace{0.2cm}

{\bf (II.5.11).} Line -2, replace  {\em``... under $f$, or, equivalently ...''} by  {\em``... under $f$. If $f$ is strict , then the the 
canonical homomorphism $f^{-1}(\cM^{\sharp}_{Y}) \rightarrow \cM^{\sharp}_{X}$ is an isomorphism and the converse 
holds  if, moreover,  $\cM_{X}$ is $u$-integral in the sense of [58], I. Def. 1.3.1, 3.  (cf. [58], III. Cor. 1.2.11)''}.

\vspace{0.2cm}

{\bf (II.6.20.1).} Replace $E_1^{i,j}$ by $E_2^{i,j}$.

\vspace{0.2cm}

{\bf (II.8.1.10).} Replace the ${\rm{H}}_1$'s by the ${\rm{H}}_0$'s without changing the conclusion.

\vspace{0.2cm}

{\bf (II.8.1.11).} Replace the ${\rm{H}}_0$'s by the ${\rm{H}}_1$'s without changing the conclusion.

\vspace{0.2cm}

{\bf (II.8.21).} Replace in the last line  {\em``whose kernel is annihilated by $p^{\frac{1}{p-1}}$.''} by  {\em``whose 
kernel is annihilated by $\fm_\oK p^{\frac{1}{p-1}}$.''}

\vspace{0.2cm}

{\bf (II.9.3).} Line 5-6, replace the sentence  {\em``For every integer $n\geq 1$, the canonical projection  $\cR_A\rightarrow A/pA$ onto the $(n+1)$th component of the 
inverse system $(A/pA)_{\mN}$ (that is, the component of index $n$)''} by the sentence  {\em``For every integer $n\geq 1$, 
the canonical projection  $\cR_A\rightarrow A/pA$ onto the $n$th component of the inverse system $(A/pA)_{\mN}$ 
(that is, the component of index $n-1$)''}.

\vspace{0.2cm}

{\bf (II.11.1).} Line 4-5,  we could have warn the reader that the notation $\alpha^{r,r'}$ and $h_{\alpha}^{r,r'}$ is not
compatible with that used in (I.4.7). Same remark for (II.12.1.6). 

\vspace{0.2cm}

{\bf (II.11.13).} Line 3, replace {\em``for every $1\leq i \leq n$''} by  {\em``for every $1\leq i \leq d$''}.

\vspace{0.2cm}

{\bf (II.11.14).} Line -3 of the proof, replace   {\em``by virtue of (II.11.14.6) and II.11.12(ii)''} by  {\em``by virtue of (II.11.14.6) and II.11.12(iii)''}.

\vspace{0.2cm}

{\bf (II.12.1.6).} Line -2,  we could have warn the reader that the notation $\alpha^{r,r'}$ and $h_{\alpha}^{r,r'}$ are not compatible with those used in (I.4.7). 

\vspace{0.2cm}

{\bf (II.13.15.7).} Replace $k\varphi(g) = \exp(\sum_{i=1}^d\xi^{-1}\theta_i\otimes \chi_{i}(g))$ by $\varphi(g) = \exp(\sum_{i=1}^d\xi^{-1}\theta_i\otimes \chi_{i}(g))$. 

\vspace{0.2cm}

{\bf (II.14.2).} Line 1, replace {\em``Let $\alpha$ be rational number $>\frac{1}{p-1}$''} by  {\em``Let $\alpha$ be a rational number $>\frac{1}{p-1}$''}.

\vspace{0.2cm}

{\bf (II.14.4.4).} Four lines  before (II.14.4.4),  replace   {\em``By II.14.3 ... and a  $\Delta_{p^\infty}$-equivariant $\oR$-linear isomorphism''} by  
{\em``By II.14.3 .. and a  $\Delta$-equivariant $\oR$-linear isomorphism''}. 

\vspace{0.2cm}

{\bf (II.14.4.4).}  One line under (II.14.4.4), replace the sentences:
 
{\em``By virtue of  II.14.1, for all integers $n\geq m>\alpha$, there exists a unique $\Delta_{p^\infty}$-equivariant $R_1$-linear isomorphism
\begin{equation}
N_n/p^{m-\alpha}N_n\stackrel{\sim}{\rightarrow}N_m/p^{m-\alpha}N_m
\end{equation}
that is compatible with the isomorphisms (II.14.4.4).
Consequently, the $R_1$-modules $(N_n)_{n > \alpha}$
\[
(N_n/p^{n-\alpha-\beta-\frac{1}{p-1}}N_n)_{n>\alpha+\beta+\frac{1}{p-1}}
\] 
form an inverse system...''}

by the sentences: 

{\em``By virtue of  II.14.1, for all integers $n\geq m>\alpha+\beta+\frac{1}{p-1}$, there exists a unique 
$\Delta_{p^\infty}$-equivariant $R_1$-linear isomorphism
\begin{equation}\label{desc4c}
N_n/p^{m-\alpha-\beta-\frac{1}{p-1}}N_n\stackrel{\sim}{\rightarrow}N_m/p^{m-\alpha-\beta-\frac{1}{p-1}}N_m
\end{equation}
that is compatible with the isomorphisms (II.14.4.4). Consequently, the $R_1$-modules
\[
(N_n/p^{n-\alpha-\beta-\frac{1}{p-1}}N_n)_{n>\alpha+\beta+\frac{1}{p-1}}
\] 
form an inverse system...''}

\vspace{0.2cm}

{\bf (III.1).} Line 11, replace  {\em``valuation ring''} by  {\em``valuation field''}.

\vspace{0.2cm}

{\bf (III.2.1).} Line 1-2,  replace  {\em``valuation ring''} by  {\em``valuation field''}.

\vspace{0.2cm}
 
{\bf (III.10.16.11).} Replace the $\times$ by $\otimes$ on each of the 2 horizontal lines.

\vspace{0.2cm}

{\bf (III.10.19.2).} Replace the $\times$ by $\otimes$ on each of the 2 horizontal lines.

\vspace{0.2cm}

{\bf (III.10.30 (iii)).} Line 2, replace {\em``inverse limit''} by {\em``direct limit''}.

\vspace{0.2cm}

{\bf (III.12.7).} Line 8, replace {\em``every object of $\Xi^r$ is a Higgs $\bvocB$-isogeny''} by  {\em``every object of $\Xi^r_{\mQ}$ is a Higgs $\bvocB$-isogeny''}.

\vspace{0.2cm}

{\bf (III.14.7.5).} In the descending left arrow, replace $\Phi^*$ by $\bvPhi^*$.

\vspace{0.2cm}

{\bf (III.14.7.11).} Replace the middle term $ \top'^r_+\circ \bvPhi^*\circ \fS$ by  $ \top'^r_+\circ \bvPhi^*\circ \fS^r$.

\vspace{0.2cm}

{\bf (IV.5.1).} Line 21 of page 384, replace $[\;]\colon \oCA\to W(R_{\oCA})$ by $[\;]\colon R_{\oCA}\to W(R_{\oCA})$. 

\vspace{0.2cm}

{\bf (IV.6.1).} The last line of page 411, replace $\Gamma(\os,\CO_{\os})$ by $\kappa(\os)=\Gamma(\os,\CO_{\os})$.

\vspace{0.2cm}

{\bf (IV.6.2).} Line 20 of page 414, replace $\Gal(\kappa(\os)/\kappa(s^g))$ by $\Gal(\kappa(s)^{\ur}/\kappa(s^g))$. 

\vspace{0.2cm}

{\bf (IV.6.2).} Line 20 of page 414, replace $\Aut_{\CC_{\gpt}}((U,\os))^{\circ}$ by $\Aut_{(\CU_{\oK,\triv})_{\gpt}}(\os)^{\circ}$.

\vspace{0.2cm}

{\bf (V.11.7).} Replace $E_1^{a,b}$ by $E_2^{a,b}$.

\vspace{0.2cm}

{\bf (V.12.1).} Line 3 of page 482, replace {\em``the multiplication by $a$''} by {\em``the multiplication by $\tr_G(b)$''}. 

\vspace{0.2cm}

{\bf (VI.1.13).}  Line 19, replace {\em``valuation ring''} by  {\em``valuation field''}.

\vspace{0.2cm}

{\bf (VI.3.4).} Line 2, Line 3 and 2 times Line -3,  replace $e_Z$ by $e_S$. 

\vspace{0.2cm}

{\bf (VI.3.5).} Line 5, replace 2 times $e_Z$ by $e_S$.  

\vspace{0.2cm}

{\bf (VI.5.9).} Line 4 : replace $(V_{n,m} \rightarrow V_n)_{n \in M_n}$ by $(V_{n,m} \rightarrow V_n)_{m \in M_n}$

\vspace{0.2cm}

{\bf (VI.6.5.1).} Two lines under (VI.6.5.1), replace $, \rightarrow$ by $\rightarrow$.

\vspace{0.2cm}

\begin{center}

{\large \bf B) Errata}

\end{center}

\vspace{0.2cm}

{\bf (II.9.5).}  Replace the proof of the proposition by the following:

Indeed, we clearly have $\theta(\xi)=0$. By ([73]  A.2.3), since $\rW(A^\flat)$ and $\hA$ are $\mZ_p$-flat 
and complete and separated for the $p$-adic topologies, it is enough to prove that the sequence 
\begin{equation*}\label{eip4c}
\tag{II.9.5.3}
\xymatrix{
0\ar[r]&{A^\flat}\ar[r]^-(0.5){\cdot \upp}&{A^\flat}\ar[r]^-(0.5){v_1}&{A/pA}\ar[r]&0},
\end{equation*}
where $v_1$ is the homomorphism induced by the projection on the first component of the projective system $(A/pA)_\mN$ (II.9.3.1),
is exact. By (iii), $v_1$ is surjective.

Let $y=(y_n)_{n\in \mN}\in A^\flat$ such that $\upp y=0$. 
For every $n\geq 0$, let $\ty_n$ be a lift of $y_n$ in $A$. We have $p_n\ty_n\in pA$. Consequently, $\ty_n\in p_n^{p^n-1}A$ 
because $p$ is not a zero divisor in $A$. It follows that 
\begin{equation*}\label{eip4d}
\tag{II.9.5.4}
y_n=y_{n+1}^p=(\ty_{n+1}^p \mod pA)=0
\end{equation*} 
because $p^{n+2}-p\geq p^{n+1}$. Then, $\upp$ is not a zero divisor in $A^\flat$.

It is clear that $v_1(\upp y)=0$ for any $y\in A^\flat$.  
Conversely, let $x=(x_n)_{n\in \mN}\in A^\flat$ such that $x_0=0$. 
For every $n\geq 0$, let $\tx_n$ be a lift of $x_n$ in $A$.
By (i) and (ii), there exists $\ty_n\in A$ such that $\tx_n=p_n \ty_n$. 
From the relation $\tx_{n+1}^p\equiv \tx_n\mod pA$ we deduce that $\ty_{n+1}^p\equiv \ty_n \mod p_n^{p^n-1}A$. 
For every $n\geq 1$, we have 
\begin{equation*}\label{eip4e}
\tag{II.9.5.5}
\ty_{n+1}^{p^2}\equiv \ty_n^p \mod pA,
\end{equation*} 
because $p^{n+1}-p\geq p^{n}$. Hence, $y=(\ty_{n+1}^p \mod pA)_{n\geq 0}\in A^\flat$.  
Since $x_n=x_{n+1}^p=p_n \ty_{n+1}^p \mod p$, we have $x=\upp y$. 
Therefore, the sequence \eqref{eip4c} is exact in the center.

\vspace{0.2cm}

{\bf (II.11.12).} Proof of (ii): The sufficiency of the required exactness of (II.11.12.5) is missing there. One can rather argue as follows:

(ii) Since for every integer $n\geq 0$, we clearly have 
\begin{equation*}
\tag{II.11.12.4}
p^n \cdot {_{i}\cS^{(r)}_{p^\infty}}={_{i}\cS^{(r)}_{p^\infty}}\cap (p^n\cdot {_{(i-1)}\cS^{(r)}_{p^\infty}}),
\end{equation*} 
the canonical morphism $u\colon  {_{i}\hcS^{(r)}_{p^\infty}}  \rightarrow  {_{(i-1)}\hcS^{(r)}_{p^\infty}} $ is injective. 

For every $\nu\in {_{(i-1)}\Xi_{p^\infty}}$, we denote by $(R^{(\nu)}_{p^\infty})^{\wedge}$ and $({_{(i-1)}\cS^{(r)}_{p^\infty}}(\nu))^{\wedge}$
the $p$-adic Hausdorff completions of  $R^{(\nu)}_{p^\infty}$ and ${_{(i-1)}\cS^{(r)}_{p^\infty}}(\nu)$, respectively. 
Note that $({_{(i-1)}\cS^{(r)}_{p^\infty}}(\nu))^\wedge$ is identified with a sub-$R_1$-module of $({_{(i-1)}\cS^{(r)}_{p^\infty}})^{\wedge}$ stable by $\gamma_i$. 
By (II.11.11.3), every element $x$ of $({_{(i-1)}\cS^{(r)}_{p^\infty}})^{\wedge}$ can be written as the sum of a series 
\begin{equation}
\tag{II.11.12.5}
\sum_{\nu\in {_{(i-1)}\Xi_{p^\infty}}}x_\nu,
\end{equation}
where $x_\nu\in {(_{(i-1)}\cS^{(r)}_{p^\infty}}(\nu))^\wedge$ and, for every integer $n\geq 0$, except for finitely many $\nu\in {_{(i-1)}\Xi_{p^\infty}}$, 
$x_\nu\in p^n ({_{(i-1)}\cS^{(r)}_{p^\infty}}(\nu))^\wedge$.
Such an element $x$ is zero if and only if $x_\nu$ is zero for every $\nu\in {_{(i-1)}\Xi_{p^\infty}}$.

Let $\nu\in {_{(i-1)}\Xi_{p^\infty}}$. For every
\begin{equation*}
\tag{II.11.12.6}
z=\sum_{\un\in J_{i-1}}p^{r|\un|}a_\un \uy^\un\in {_{(i-1)}\cS^{(r)}_{p^\infty}}(\nu),
\end{equation*}
we have (II.11.6.1)
\begin{equation*}\label{mtht31gg}
\tag{II.11.12.7}
(\gamma_i-\id)(z)=\sum_{\un\in J_{i-1}}p^{r|\un|}b_\un \uy^\un\in {_{(i-1)}\cS^{(r)}_{p^\infty}}(\nu),
\end{equation*}
where, for every $\un=(n_1,\dots,n_{d'})\in J_{i-1}$, 
\begin{equation*}
b_\un=(\nu(\gamma_i)-1)a_\un+\sum_{\um=(m_1,\dots,m_{d'})
\in J_{i-1}(\un)} p^{r(m_i-n_i)}\binom{m_i}{n_i} \nu(\gamma_i)a_{\um}w^{m_i-n_i},
\end{equation*}
$J_{i-1}(\un)$ denotes the subset $J_{i-1}$ made up of elements $\um=(m_1,\dots,m_{d'})$ such that
$m_j=n_j$ for $j\not=i$ and $m_i>n_i$, 
and $w=\xi^{-1}\log([\zeta])$ is an element of valuation $\frac{1}{p-1}$ of $\co_C$ (II.9.18).

Let $z$ be an element of $({_{(i-1)}\cS^{(r)}_{p^\infty}}(\nu))^\wedge$. Then, $z$ can be written as the sum of a series 
\begin{equation*}
z=\sum_{\un\in J_{i-1}}p^{r|\un|}a_\un \uy^\un,
\end{equation*}
where $a_\un \in (R^{(\nu)}_{p^\infty})^\wedge$ and $a_\un$ tends to $0$ when $|\un|$ tends to infinity. 
As $R^{(\nu)}_{p^\infty}$ is $\co_C$-flat, the same therefore holds for $(R^{(\nu)}_{p^\infty})^\wedge$ (cf. the proof of II.6.14). 
Consequently,  $z$ is zero if and only if $a_\un$ is zero for every $\un\in J_{i-1}$.
We immediately see that $(\gamma_i-\id)(z)$ is also given by the formula \eqref{mtht31gg}.

Suppose that $\gamma_i(z)=z$ and $\nu(\gamma_i)\not=1$.  
As $(R^{(\nu)}_{p^\infty})^\wedge$ is $\co_C$-flat and  $v(\nu(\gamma_i)-1)\leq \frac{1}{p-1}$, 
for every  $\un=(n_1,\dots,n_{d'})\in J_{i-1}$, we have 
\[
a_\un=-(\nu(\gamma_i)-1)^{-1}\sum_{\um=(m_1,\dots,m_{d'})
\in J_{i-1}(\un)} p^{r(m_i-n_i)}\binom{m_i}{n_i} \nu(\gamma_i)a_{\um}w^{m_i-n_i}.
\]
We deduce from this that for every $\alpha\in \mN$ and every $\un\in J_{i-1}$,
we have $a_\un\in p^{r\alpha} (R'^{(\nu)}_{p^\infty})^\wedge$ (this is proved by induction on $\alpha$);
therefore $z=0$ as $(R'^{(\nu)}_{p^\infty})^\wedge$ is separated for the $p$-adic topology. 
Consequently, $\gamma_i-\id$ is injective on  $({_{(i-1)}\cS^{(r)}_{p^\infty}}(\nu))^\wedge$. 

Suppose that $\gamma_i(z)=z$ and $\nu(\gamma_i)=1$, so that $\nu\in {_{(i)}\Xi_{p^\infty}}$. 
Then, for every $\un=(n_1,\dots,n_{d'})\in J_{i-1}$, if we set 
$\un'=(n'_1,\dots,n'_{d'})\in J_{i-1}(\un)$ with $n'_i=n_i+1$, we have
\[
(n_i+1)!a_{\un'}=-\sum_{\um=(m_1,\dots,m_{d'})
\in J_{i-1}(\un')} p^{r(m_i-n_i-1)}m_i! a_{\um}\frac{w^{m_i-n_i-1}}{(m_i-n_i)!}.
\]
We have $w^{m-1}/m!\in \co_C$ for every integer $m\geq 1$. We deduce from this that 
for every $\alpha\in \mN$ and every $\un=(n_1,\dots,n_{d'})\in J_{i-1}$ such that $n_i\geq 1$, 
we have $n_i!a_{\un}\in p^{r\alpha} R'^{(\nu)}_{p^\infty}$ (this is proved by induction on $\alpha$); therefore $a_\un=0$. 
Consequently $z \in ({_{i}\cS^{(r)}_{p^\infty}}(\nu))^\wedge$. We deduce from this that the sequence
\begin{equation*}
\tag{II.11.12.8}
\xymatrix{
0\ar[r]&{({_{i}\cS^{(r)}_{p^\infty}}(\nu))^\wedge}\ar[r]^-(0.45){u_\nu}&{({_{(i-1)}\cS^{(r)}_{p^\infty}}(\nu)^\wedge}\ar[r]^{\gamma_i-\id}&{({_{(i-1)}\cS^{(r)}_{p^\infty}}(\nu))^\wedge}},
\end{equation*}
where $u_\nu$ is the canonical morphism, is exact. The proposition follows. 

\vspace{0.2cm}

{\bf (II.13.9).} Modify the beginning of (II.13.9) as follows: 

Let $M$ be an $\hRun$-module of finite type which is $\co_C$-flat.
By ([1]  1.10.2), $M$ is complete and separeted for the $p$-adic topology. We deduce that the canonical morphism
\begin{equation*}\label{smrep8a}
\tag{II.13.9.a}
\End_{\hRun}(M)\rightarrow \underset{\underset{n\geq 0}{\longleftarrow}}{\lim}\ \Hom_{\hRun}(M,M/p^nM)
\end{equation*}
is an isomorphism. On the other hand, for any integer $n\geq 0$, the exact sequence 
\begin{equation*}\label{smrep8b}
\tag{II.13.9.b}
0\rightarrow M\stackrel{p^n}{\rightarrow} M\rightarrow M/p^nM\rightarrow 0
\end{equation*}
shows that the canonical map 
\begin{equation*}\label{smrep8c}
\tag{II.13.9.c}
\End_{\hRun}(M)/p^n\End_{\hRun}(M)\rightarrow \Hom_{\hRun}(M,M/p^nM)
\end{equation*}
is injective. We deduce from \eqref{smrep8a} and \eqref{smrep8c} that $\End_{\hRun}(M)$ is complete and separated 
for the $p$-adic topology. Moreover, it is $\co_C$-flat and for every rational number $\alpha\geq 0$, the canonical 
homomorphism $p^{\alpha}\End_{\hRun}(M)\rightarrow \Hom_{\hRun}(M,p^{\alpha}M)$ is an isomorphism. Let $u$...

\vspace{0.2cm}

{\bf (III.10.3).} The category $\bP$ should moreover be required to contain the affine schemes $U$ such that the special fiber $U_s$
is empty in order to have covering families.

\vspace{0.2cm}

{\bf (III.10.5).} The category $\bQ$ should moreover be required to contain  the affine schemes $U$ such that the special fiber $U_s$
is empty in order to have covering families.

\vspace{0.2cm}

{\bf (V.12.1).} Lemma V.12.1(2) is wrong. There is an obvious counterexample: $A\neq 0$, $B=A[X]$, 
$G=\{1\}$, $M=B$, $a=1$, $r=1$, $b_1=c_1=1$. The claim $\varphi\circ \psi=a\cdot 1_{B\otimes_AM^G}$ in the proof of Lemma V.12.1(2)
is not true in general unless the homomorphism $A\to B^G$ is surjective. 

Replace Lemma V.12.1(2) and its proof by the following:

{\it Assume that $a\in B$ and $a'\in B^G$ satisfy the following conditions:\par
{\rm (i)} There exist $b_i,c_i\in B$ $(1\leq i\leq r)$ such that
$\sum_{i=1}^rb_ic_i=a$ and $\sum_{i=1}^r b_ig(c_i)=0$ for all $g\in G\backslash\{1\}$.
\par
{\rm (ii)} $a' B^G$ is contained in the image of the homomorphism $A\to B^G$.\par
Then, the kernel and the cokernel of $\psi\colon B\otimes_AM^G\to M;b\otimes x\mapsto bx$
are killed by $aa'$.}

{\sc Proof.} We define a map $\varphi\colon M\to B\otimes_AM^G$ by
$\varphi(m)=\sum_{i=1}^r b_i\otimes a'\trace_G(c_im)$, which is 
$A$-linear. We assert $\varphi\circ \psi=aa'\cdot 1_{B\otimes_AM^G}$
and $\psi\circ\varphi=aa'\cdot 1_M$. The image of $b\otimes m \in B\otimes_AM^G$
$(b\in B,m\in M^G)$ under $\varphi\circ\psi$ is 
$\sum_{i=1}^rb_i\otimes a'\trace_G(c_ibm)=\sum_{i=1}^rb_i\otimes(a'\trace_G(c_ib))m$
by the definition of $\varphi$ and $\psi$, and $m\in M^G$. By the condition (ii), we
obtain
\begin{multline*}
\varphi\circ \psi(b\otimes m)=\sum_{i=1}^r (b_ia'\trace_G(c_ib)\otimes m)
=\biggl(a'\sum_{g\in G}\biggl(\sum_{i=1}^rb_ig(c_i)\biggr)g(b)\biggr)\otimes m=a'ab\otimes m.
\end{multline*}
The second equality is shown as
\begin{equation*}\psi\circ \varphi(m)=\sum_{i=1}^r b_ia'\trace_G(c_im)=
a'\sum_{g\in G}\biggl(\sum_{i=1}^rb_ig(c_i)\biggr)g(m)=a'am.
\qedhere
\end{equation*}
\hfill{$\Box$}

\vspace{0.2cm}

{\bf (V.12.5) and (V.12.6).} Replace the proof of Corollary V.12.6 by the direct proof below. This replacement allows us to 
apply the corrected Lemma V.12.1(2) above in the proof of Proposition V.12.5.

{\sc Proof.} By Lemma V.12.4, $R\to S$ is almost faithfully flat and the homomorphism
$S\otimes_RS\to \prod_{g\in G}S;x\otimes y\mapsto (xg(y))_{g\in G}$ is an
almost isomorphism. By Proposition V.9.1, we have an almost exact sequence 
$R\to S\xrightarrow{d}S\otimes_RS$,  where $d$ is defined by 
$d(s)=1\otimes s-s\otimes 1$. The composition of
$d$ with the almost isomorphism $S\otimes_RS\xrightarrow{\approx}\prod_{g\in G}S$
sends $s$ to $(g(s)-s)_{g\in G}$. Hence $R\to S^G$ is an almost isomorphism.
\hfill{$\Box$}

\vspace{0.2cm}

{\bf (VI.10.40).} The introduction of $X_{(\ox)}$ in the proposition is useless.  A better way to formulate the conclusion is to say that the map 
\begin{equation*}
\cH^i(F) \rightarrow \rR^i\sigma_*(F^a)
\end{equation*} 
is an isomorphism. It follows immediately from the corresponding isomorphism on the stalks given in the proposition and its proof.

\end{document}